\documentclass[reqno,10pt,letterpaper]{amsart}

\usepackage{lipsum}
\usepackage{amsmath}
\usepackage{amssymb}
\usepackage{amsthm}
\usepackage{mathrsfs}
\usepackage{accents}
\usepackage{calc}
\usepackage{arydshln}
\usepackage{upgreek}
\usepackage{slashed}
\usepackage{xifthen}
\usepackage{graphicx}
\usepackage{subcaption}
\usepackage{longtable}
\usepackage[inline]{enumitem}

\usepackage{xr}
\newcommand{\citestab}[1]{\cite[#1]{HintzKerrStab}}

\usepackage{tikz}

\usepackage{xcolor}
\definecolor{winered}{rgb}{0.6,0,0}
\definecolor{lessblue}{rgb}{0,0,0.7}

\usepackage[pdftex,colorlinks=true,linkcolor=winered,citecolor=lessblue,urlcolor=lessblue,breaklinks=true,bookmarksopen=true]{hyperref}

\makeatletter
\newcommand{\myitem}[2]{\item[\rm(#2)]\def\@currentlabel{#2}\label{#1}}
\makeatother

\usepackage{titletoc}

\makeatletter

\def\@tocline#1#2#3#4#5#6#7{
\begingroup
  \par
    \parindent\z@ \leftskip#3 \relax \advance\leftskip\@tempdima\relax
                  \rightskip\@pnumwidth plus 4em \parfillskip-\@pnumwidth
    \ifcase #1 
       \vskip 0.6em \hskip 0em 
       \or
       \or \hskip 0em 
       \or \hskip 1em 
    \fi%
    #6
    \nobreak\relax{\leavevmode\leaders\hbox{\,.}\hfill}
    \hbox to\@pnumwidth {\@tocpagenum{#7}}
  \par
\endgroup
}

 \def\l@section{\@tocline{0}{0pt}{0pc}{}{}}

\renewcommand{\tocsection}[3]{%
  \indentlabel{\@ifnotempty{#2}{ 
    \ignorespaces\bfseries{#2. #3}}}
  \indentlabel{\@ifempty{#2}{\ignorespaces\bfseries{#3}}{}} 
    \vspace{1.5pt}}

\renewcommand{\tocsubsection}[3]{%
  \indentlabel{\@ifnotempty{#2}{
    \ignorespaces#2. #3}}
  \indentlabel{\@ifempty{#2}{\ignorespaces #3}{}}
    \vspace{1.5pt}}

\renewcommand{\tocsubsubsection}[3]{%
  \indentlabel{\@ifnotempty{#2}{
    \ignorespaces#2. #3}}
  \indentlabel{\@ifempty{#2}{\ignorespaces #3}{}}
    \vspace{1.5pt}}

\makeatother

\makeatletter
\def\@nomenstarted{0}
\newlength{\@nomenoldtabcolsep}

\newcommand{\nomenstart}
  {%
    \def\@nomenstarted{1}%
    \setlength{\@nomenoldtabcolsep}{\tabcolsep}%
    \setlength{\tabcolsep}{3.5pt}%
    \begin{longtable}{p{0.11\textwidth} p{0.86\textwidth}}
  }

\newcommand{\nomenitem}[2]{%
    \ifcase\@nomenstarted%
      \or 
      \or \\ 
    \fi%
    #1\,{\leavevmode\leaders\hbox{\,.}\hfill} & #2%
    \def\@nomenstarted{2}%
  }%
\newcommand{\nomenend}
  {\\%
      \end{longtable}%
      \setlength{\tabcolsep}{\@nomenoldtabcolsep}%
      \def\@nomenstarted{0}%
  }
\makeatother

\makeatletter
\newcommand{\BIG}{\bBigg@{3.5}}
\newcommand{\vast}{\bBigg@{4}}
\newcommand{\Vast}{\bBigg@{5}}
\newcommand{\VAST}[1]{\bBigg@{#1}}
\makeatother

\allowdisplaybreaks

\numberwithin{equation}{section}
\numberwithin{figure}{section}
\newtheorem{thm}{Theorem}[section]

\newtheorem{prop}[thm]{Proposition}
\newtheorem{lemma}[thm]{Lemma}
\newtheorem{cor}[thm]{Corollary}

\newtheorem*{thm*}{Theorem}
\newtheorem*{prop*}{Proposition}
\newtheorem*{cor*}{Corollary}
\newtheorem*{conj*}{Conjecture}

\theoremstyle{definition}
\newtheorem{definition}[thm]{Definition}
\newtheorem{notation}[thm]{Notation}

\theoremstyle{remark}
\newtheorem{rmk}[thm]{Remark}
\newtheorem{example}[thm]{Example}

\makeatletter
\newcommand{\fakephantomsection}{%
  \Hy@MakeCurrentHref{\@currenvir.\the\Hy@linkcounter}
  \Hy@raisedlink{\hyper@anchorstart{\@currentHref}\hyper@anchorend}%
  \Hy@GlobalStepCount\Hy@linkcounter%
}
\makeatother

\newcommand{\mc}{\mathcal}
\newcommand{\cA}{\mc A}

\newcommand{\cC}{\mc C}
\newcommand{\cD}{\mc D}
\newcommand{\cE}{\mc E}
\newcommand{\cF}{\mc F}

\newcommand{\cH}{\mc H}

\newcommand{\cK}{\mc K}
\newcommand{\cL}{\mc L}
\newcommand{\cM}{\mc M}
\newcommand{\cN}{\mc N}
\newcommand{\cO}{\mc O}
\newcommand{\cP}{\mc P}

\newcommand{\cR}{\mc R}

\newcommand{\cT}{\mc T}
\newcommand{\cU}{\mc U}
\newcommand{\cV}{\mc V}
\newcommand{\cW}{\mc W}
\newcommand{\cX}{\mc X}
\newcommand{\cY}{\mc Y}

\newcommand{\ms}{\mathscr}
\newcommand{\sA}{\ms A}

\newcommand{\sC}{\ms C}
\newcommand{\sD}{\ms D}

\newcommand{\sG}{\ms G}
\newcommand{\sH}{\ms H}

\newcommand{\scri}{\ms I}
\newcommand{\sscri}{{\!\scri}}

\newcommand{\sP}{\ms P}

\newcommand{\sS}{\ms S}
\newcommand{\sV}{\ms V}
\newcommand{\sW}{\ms W}

\newcommand{\HH}{\mathbb{H}}

\newcommand{\TT}{\mathbb{T}}

\newcommand{\C}{\mathbb{C}}
\newcommand{\N}{\mathbb{N}}
\newcommand{\R}{\mathbb{R}}
\newcommand{\Z}{\mathbb{Z}}

\newcommand{\Sph}{\mathbb{S}}

\newcommand{\sfb}{\mathsf{b}}

\newcommand{\sfl}{\mathsf{l}}
\newcommand{\sfm}{\mathsf{m}}
\newcommand{\sfp}{\mathsf{p}}
\newcommand{\sfr}{\mathsf{r}}
\newcommand{\sfs}{\mathsf{s}}

\newcommand{\sfH}{\mathsf{H}}

\newcommand{\sfM}{\mathsf{M}}

\newcommand{\sfX}{\mathsf{X}}

\newcommand{\bfB}{\mathbf{B}}

\newcommand{\fB}{\mathfrak{B}}

\newcommand{\fk}{\mathfrak{k}}
\newcommand{\fm}{\mathfrak{m}}
\newcommand{\fM}{\mathfrak{M}}

\newcommand{\fq}{\mathfrak{q}}

\newcommand{\ft}{\mathfrak{t}}

\newcommand{\fX}{\mathfrak{X}}

\newcommand{\slg}{\slashed{g}{}}

\newcommand{\slDelta}{\slashed{\Delta}{}}

\newcommand{\slpi}{\slashed{\pi}{}}

\newcommand{\ran}{\operatorname{ran}}
\newcommand{\ann}{\operatorname{ann}}

\newcommand{\End}{\operatorname{End}}

\newcommand{\Hom}{\operatorname{Hom}}
\renewcommand{\Re}{\operatorname{Re}}
\renewcommand{\Im}{\operatorname{Im}}

\newcommand{\ind}{{\operatorname{ind}}}
\newcommand{\mathspan}{\operatorname{span}}
\newcommand{\supp}{\operatorname{supp}}
\newcommand{\sgn}{\operatorname{sgn}}

\newcommand{\tr}{\operatorname{tr}}

\newcommand{\dv}{\operatorname{div}}

\newcommand{\diag}{\operatorname{diag}}

\let\ind\relax
\DeclareMathOperator{\ind}{ind}

\newcommand{\eps}{\epsilon}
\newcommand{\ftrans}{\;\!\wh{\ }\;\!}

\newcommand{\hra}{\hookrightarrow}
\newcommand{\la}{\langle}

\newcommand{\ol}{\overline}
\newcommand{\pa}{\partial}
\newcommand{\dd}{{\mathrm d}}
\newcommand{\ra}{\rangle}
\newcommand{\spec}{\operatorname{spec}}

\newcommand{\weakto}{\rightharpoonup}
\newcommand{\wh}{\widehat}
\newcommand{\wt}{\widetilde}
\newcommand{\xra}{\xrightarrow}

\newcommand{\ubar}[1]{\underaccent{\bar}#1}
\newcommand{\pfstep}[1]{$\bullet$\ \underline{\textit{#1}}}
\newcommand{\pfsubstep}[2]{{\bf#1}\ \textit{#2}}

\newcommand{\bop}{{\mathrm{b}}}

\newcommand{\scop}{{\mathrm{sc}}}
\newcommand{\chop}{{\mathrm{c}\semi}}

\newcommand{\cl}{{\mathrm{cl}}}

\newcommand{\ebop}{{\mathrm{eb}}}
\newcommand{\eop}{{\mathrm{e}}}
\newcommand{\tbop}{{3\mathrm{b}}}
\newcommand{\etbop}{{\mathrm{e}3\mathrm{b}}}
\newcommand{\ebeop}{{\mathrm{e}\mathrm{b}\mathrm{e}}}

\newcommand{\cuop}{{\mathrm{cu}}}
\newcommand{\scbtop}{{\mathrm{sc}\text{-}\mathrm{b}}}
\newcommand{\schop}{{\mathrm{sc},\semi}}

\newcommand{\semi}{\hbar}

\newcommand{\lb}{{\mathrm{lb}}}
\newcommand{\rb}{{\mathrm{rb}}}

\newcommand{\cface}{{\mathrm{cf}}}

\newcommand{\scface}{{\mathrm{scf}}}
\newcommand{\sctface}{{\mathrm{sctf}}}
\newcommand{\sface}{{\mathrm{sf}}}

\newcommand{\tface}{{\mathrm{tf}}}

\newcommand{\ztface}{{\mathrm{ztf}}}
\newcommand{\zface}{{\mathrm{zf}}}

\newcommand{\zeroface}{{\mathrm{0f}}}

\newcommand{\res}{{\mathrm{res}}}

\newcommand{\cp}{{\mathrm{c}}}

\newcommand{\Diff}{\mathrm{Diff}}

\DeclareMathOperator{\Op}{Op}

\newcommand{\Vb}{\cV_\bop}
\newcommand{\Ve}{\cV_\eop}
\newcommand{\Vscbt}{\cV_\scbtop}

\newcommand{\Diffb}{\Diff_\bop}
\newcommand{\Diffe}{\Diff_\eop}
\newcommand{\Diffscbt}{\Diff_\scbtop}

\newcommand{\Psib}{\Psi_\bop}

\newcommand{\Psiscbt}{\Psi_\scbtop}

\newcommand{\Psisc}{\Psi_\scop}

\newcommand{\Vtb}{\cV_\tbop}
\newcommand{\Vetb}{\cV_\etbop}

\newcommand{\Difftb}{\Diff_\tbop}

\newcommand{\Veb}{\cV_\ebop}

\newcommand{\Diffeb}{\Diff_\ebop}

\newcommand{\Tscbt}{{}^\scbtop T}

\newcommand{\Vsc}{\cV_\scop}
\newcommand{\Diffsc}{\Diff_\scop}

\newcommand{\Psih}{\Psi_\semi}

\newcommand{\WF}{\mathrm{WF}}
\newcommand{\Ell}{\mathrm{Ell}}
\newcommand{\Char}{\mathrm{Char}}

\newcommand{\Ellsc}{\mathrm{Ell_\scop}}
\newcommand{\Omegab}{{}^{\bop}\Omega}

\newcommand{\Tb}{{}^{\bop}T}

\newcommand{\Tsc}{{}^{\scop}T}

\newcommand{\Te}{{}^{\eop}T}
\newcommand{\Teb}{{}^{\ebop}T}
\newcommand{\Ttb}{{}^{\tbop}T}
\newcommand{\Tetb}{{}^{\etbop}T}

\newcommand{\Sb}{{}^{\bop}S}
\newcommand{\Se}{{}^{\eop}S}

\newcommand{\Ssc}{{}^{\scop}S}

\newcommand{\Seb}{{}^{\ebop}S}
\newcommand{\Stb}{{}^{\tbop}S}

\newcommand{\sigmasc}{\upsigma_\scop}

\newcommand{\sigmatb}{\upsigma_\tbop}

\newcommand{\sub}{{\mathrm{sub}}}

\newcommand{\loc}{{\mathrm{loc}}}
\newcommand{\CI}{\cC^\infty}

\newcommand{\CIc}{\cC^\infty_\cp}

\newcommand{\CIb}{\cC^\infty_\bop}

\newcommand{\Hb}{H_{\bop}}

\newcommand{\He}{H_{\eop}}

\newcommand{\Heb}{H_{\ebop}}

\newcommand{\Htb}{H_\tbop}

\newcommand{\Hsc}{H_{\scop}}

\newcommand{\bhm}{\fm}

\newcommand{\openbigpmatrix}[1]
  {%
    \def\@bigpmatrixsize{#1}%
    \addtolength{\arraycolsep}{-#1}%
    \begin{pmatrix}%
  }
\newcommand{\closebigpmatrix}
  {%
    \end{pmatrix}%
    \addtolength{\arraycolsep}{\@bigpmatrixsize}%
  }

\newlength{\enummargin}

\newcommand{\usref}[1]{{\upshape\ref{#1}}}

\DeclareGraphicsExtensions{.mps}

\makeatletter
\newcommand*{\fwbw}[1]{\expandafter\@fwbw\csname c@#1\endcsname}
\newcommand*{\@fwbw}[1]{\ifcase #1 \or {\rm fw}\or {\rm bw}\fi}
\AddEnumerateCounter{\fwbw}{\@fwbw}
\makeatother

\begin{document}

\pagenumbering{roman}

\title[(Non-)Linear waves on asymptotically flat spacetimes. II]{(Non-)Linear waves on asymptotically flat spacetimes. II: trapping, bound states, nonlinear applications}

\date{\today}


\author{Peter Hintz}
\address{Department of Mathematics, Pennsylvania State University, 54 McAllister St, State College,\newline PA 16801, United States}
\email{phintz@psu.edu}

\begin{abstract}
  We study wave-type equations on dynamical spacetimes that settle down at suitable rates to a subextremal Kerr black hole spacetime. We prove strong estimates for solutions of (tensorial) linear wave-type equations when the time-translation-invariant model satisfies a spectral assumption of mode stability type. We allow for this model to admit zero energy bound states; besides the scalar wave operator (which has no bound states), examples include the wave operator on 1-forms and the linearization of the Einstein field equations in generalized harmonic gauge. We demonstrate the utility of our estimates by proving the global existence of solutions to some quasilinear wave equations, including in the presence of zero energy bound states. The results proved here are, moreover, crucial ingredients in the proof of the nonlinear stability of subextremal Kerr black holes in \cite{HintzKerrStab}.

  Our key novel linear estimate controls linear waves in weighted $L^2$-based spacetime Sobolev spaces that encode b-regularity, by which we mean regularity with respect to spacetime scaling, spatial scaling (in a hyperboloidal foliation of spacetime), and angular derivatives; this estimate is moreover tame in the b-regularity order, as needed for its applicability in a Nash--Moser iteration scheme. Its proof combines four main ingredients: microlocal propagation estimates in the edge-b-setting near null infinity (as introduced by the author with Vasy) and in the author's 3b-setting in the forward cone; estimates for the stationary model operator; energy estimates on edge-b-spaces on finite time intervals; and commutations with b-vector fields. For the nonlinear applications, we moreover develop a dictionary between decay rates in different asymptotic spacetime regimes on the one hand and weighted low-energy resolvent estimates on the other hand.

  This paper builds on Part I only a broad conceptual level, and is largely self-contained: we only use as black boxes the construction of relevant algebras of pseudodifferential operators, in the setting of scaled bounded geometry structures introduced recently by the author, estimates at the trapped set of dynamical spacetimes, and mode stability results for wave operators on subextremal Kerr spacetimes. We include a detailed description of the null-geodesic dynamics of subextremal Kerr spacetimes following Wunsch--Zworski, Vasy, and Dyatlov. Moreover, we give a streamlined exposition, with complete proofs, of the spectral theory of stationary wave operators on asymptotically flat spacetimes, including at large real and complex frequencies.
\end{abstract}

\maketitle
\thispagestyle{empty}

\clearpage

{\pagestyle{empty}
\renewcommand{\contentsname}{Table of Contents}
\setlength{\parskip}{0.00pt}
\tableofcontents
\setlength{\parskip}{0.03in}
\clearpage
}

\pagenumbering{arabic}

\section{Introduction}
\label{SI}

\begin{figure}[!ht]
\centering
\includegraphics{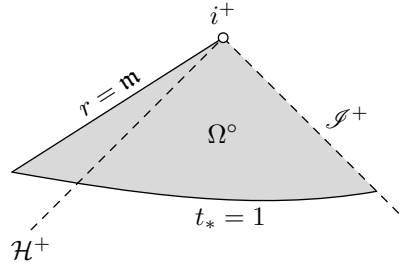}
\caption{Penrose diagrammatic depiction of the domain $\Omega^\circ\subset M^\circ$. The event horizon of the black hole is labeled $\cH^+$, and null infinity is $\scri^+$.}
\label{FigIOmega}
\end{figure}

A subextremal Kerr black hole \cite{KerrKerr,BoyerLindquistKerr} is described by the parameters $\bhm>0$ (mass) and $a\in(-\bhm,\bhm)$ (specific angular momentum). Its (Ricci-flat, Lorentzian) metric is given by
\[
  g_{\bhm,a} = -\frac{\mu}{\varrho^2}(\dd\ft-a\sin^2\theta\,\dd\phi)^2 + \varrho^2\Bigl(\frac{\dd r^2}{\mu}+\dd\theta^2\Bigr) + \frac{\sin^2\theta}{\varrho^2}\bigl((r^2+a^2)\dd\phi-a\,\dd\ft\bigr)^2
\]
on $\R_\ft\times(r_+,\infty)_r\times\Sph^2_{\theta,\phi}$, where $\mu(r)=r^2-2\bhm r+a^2$ and $\varrho^2(r,\theta)=r^2+a^2\cos^2\theta$, and $r_+=\bhm+\sqrt{\bhm^2-a^2}$ (radius of the event horizon) is the larger root of $\mu$. To track waves falling into the black hole or escaping to infinity, it is convenient to replace $\ft,\phi$ by coordinates $t_*,\phi_*$ in which $g_{\bhm,a}$ extends analytically to a Ricci-flat metric on
\[
  M^\circ=\R_{t_*}\times[\bhm,\infty)_r\times\Sph^2_{\theta,\phi_*}.
\]
Concretely, one takes $t_*\approx t-r$, $\phi_*=\phi$ for $r\gg\bhm$, and $t_*\approx t+c\log(r-r_+)$, $\phi_*\approx\phi+c'\log(r-r_+)$ for $r\approx r_+$ with suitable $c,c'>0$; we moreover arrange for the level sets of $t_*$ to be spacelike hypersurfaces of $M^\circ$. (See~\S\ref{SsTsK} for details.) Initial value problems for the wave equation $\Box_{g_{\bhm,a}}u=0$ are then well-posed on the domain
\[
  \Omega^\circ := \{ t_*\geq 1\} \subset M^\circ,
\]
with initial data posed at $t_*=1$. (The boundary hypersurface of $\Omega^\circ$ at $r=\bhm$ is a \emph{final} spacelike hypersurface, and hence no initial data need to be specified there.) See Figure~\ref{FigIOmega}.

\subsection{Description of the main results}
\label{SsIM}

The main result of the paper is a black-box estimate for linear wave equations on dynamical black hole spacetimes, which requires, as an input, suitable spectral properties of the stationary model equation. We first present concrete nonlinear applications in~\S\S\ref{SssIMN}--\ref{SssIMN2}---in which these properties are satisfied---before describing a rough version of the general linear result in~\S\ref{SssIML}.

\subsubsection{Nonlinear applications, I: zero energy bound states}
\label{SssIMN}

In the context of the black hole stability problem, a key challenge (besides issues related to gauges) is to identify the parameters of the ``final'' black hole. Other geometric/analytical difficulties stem from the fact that the null-geodesic flow for a subextremal Kerr metric exhibits (normally hyperbolic) trapping, and that there exist a horizon and an ergoregion. Moreover, the weak inverse polynomial decay of waves is in tension with the hope to regard the nonlinearity (to the maximal extent possible) as perturbative. We demonstrate in a toy problem how to surmount all of these difficulties (except for gauge issues): we consider the 1-form wave operator
\[
  \Box_{g_{\bhm,a}} := \tr\nabla^2 = (\dd+\delta_{g_{\bhm,a}})^2,
\]
which by work of Andersson--H\"afner--Whiting \cite{AnderssonHaefnerWhitingMode} (see \cite[\S{7}]{HaefnerHintzVasyKerr} for the slowly rotating case) is known to satisfy mode stability at nonzero frequencies and to admit a 1-dimensional space of (spatially decaying) stationary solutions, spanned by a divergence-free 1-form $u_0$ given explicitly by~\eqref{EqA2u0}; it satisfies
\begin{equation}
\label{EqIu0}
  u_0=\frac{1}{r}\,\dd\ft+\cO(r^{-2}).
\end{equation}

\begin{thm}[Nonlinear waves with zero energy bound states: simplified version]
\label{ThmI1}
  Let $\eps_\cK\in(0,\frac12)$ and $f_0\in\CIc(\Omega^\circ;T^*_{\Omega^\circ}M^\circ)$, with $f_0$ vanishing for $t_*\leq 1$. Then there exists $\eps_0>0$ such that for all $\eps\in[0,\eps_0]$ and $f:=\eps f_0$, the equation
  \begin{equation}
  \label{EqI1}
    \Box_{g_{\bhm,a}} u = f + |\delta_{g_{\bhm,a}}u|^2 u
  \end{equation}
  has a unique global forward solution $u$ (i.e., $u$ vanishes for $t_*\leq 1$) on $\Omega^\circ$; and there exists $c\in\C$ such that $u$ (that is, each of its components in the frame $\pa_{t_*},\pa_{x^1},\pa_{x^2},\pa_{x^3}$) satisfies the decay estimate
  \begin{equation}
  \label{EqI1Baby}
    |u(t_*,x) - c u_0(x)| \lesssim t_*^{-\eps_\cK}r^{-1}.
  \end{equation}
\end{thm}

A more quantitative version of this result, which also gives more precise asymptotics on $u$, is stated as Theorem~\ref{ThmI1Q} below.

\begin{rmk}[Preliminary comments]
\label{RmkI1Prelim}
  \begin{enumerate}
  \item \textit{Choice of nonlinearity.} The nonlinearity is the simplest one that is 1-form-valued, has sufficient decay at null infinity (i.e., $|t_*|=\cO(1)$, $r\to\infty$) to render it harmless vis-\`a-vis the (weak) null condition, and, most importantly, eliminates the leading-order term $c u_0$ of $u$ (and thus decays even if $u$ itself does not).
  \item \textit{Quasilinear perturbations.} One can replace the Kerr metric $g_{\bhm,a}$ in~\eqref{EqI1} by a metric depending on $u$ in a suitable manner; a simple (and artificial) example is $g_{\bhm,a}+\chi(x)(\delta_g u)u\otimes_s\dd r$ where $\chi\in\CIc(\R^3)$, the relevant property of which is that the metric perturbation has $t_*^{-1-\eps_\cK}$-decay as $t_*\to\infty$. We describe the class of allowed metric perturbations in~\eqref{EqIPmlPert} below.
  \item \textit{Initial value problems.} The same conclusion holds for the initial value problem for $\Box_{g_{\bhm,a}}u=|\delta_{g_{\bhm,a}}u|^2 u$ for initial data posed at a Cauchy hypersurface transversal to $\cH^+$ and tending to $i^0$, provided they are sufficiently decaying (pointwise $r^{-1-\eps_\cK}$- and $r^{-2-\eps_\cK}$-decay for $u$ and its transversal derivative suffice), regular, and small. Solving from such a Cauchy hypersurface up to some finite level set of $t_*$ is an \emph{exterior} problem (i.e., localized to regions $\ft\geq 0$, $t_*\ll -1$, $r\gg 1$ far from the black hole) that can be solved using standard energy methods. (Concrete implementations using the geometric setup of the present paper can be found in \cite[\S{4.1}]{HintzVasyMink4}, \cite[\S{3.2}]{HintzMink4Gauge}, \cite{KadarKehrbergerPhgScatter}.) The continuation of a solution known to exist until some $t_*$-level set can be accomplished by solving a forward problem in $t_*\geq 1$ (see \cite[\S{14}]{HaefnerHintzVasyKerr} for an instance of the relevant argument involving a $t_*$-cutoff). In the remainder of the paper, we shall thus focus exclusively on forward problems in $t_*\geq 1$.
  \end{enumerate}
\end{rmk}

To the author's knowledge, Theorem~\ref{ThmI1} is the first global existence result for a nonlinear wave equation on a black hole spacetime in the presence of zero energy bound states. (Previous proofs of black hole stability results for the Schwarzschild family \cite{DafermosHolzegelRodnianskiSchwarzschildStability} and the slowly rotating Kerr family \cite{KlainermanSzeftelGCM1,KlainermanSzeftelGCM2,KlainermanSzeftelKerr,ShenGCMKerr,GiorgiKlainermanSzeftelStability} use special properties of the Einstein field equations, e.g., the ability to improve decay via careful choices of gauge or approximate reductions to scalar equations (Teukolsky or Regge--Wheeler) without zero energy bound states, and express the equations as a coupled system of hyperbolic, elliptic, and transport equations. While thus of a somewhat different flavor than~\eqref{EqI1}, these works do resolve the aforementioned difficulties for a very particular quasilinear equation on slowly rotating Kerr backgrounds.)

For the linear equation $\Box_{g_{\bhm,a}}u=f$ on 1-forms, the decay estimate~\eqref{EqI1Baby} remains valid. While this has not been explicitly stated in the literature before, this follows already from a simplified version of the arguments in the proof of the linear stability of Kerr by the author with H\"afner and Vasy \cite{HaefnerHintzVasyKerr,HaefnerHintzVasyKerrErratum,HaefnerHintzVasyKerrLarge}. (For the closely related Maxwell equations ($\dd F=G_1$, $\delta_{g_{\bhm,a}}F=G_2$), Andersson--Blue \cite{AnderssonBlueMaxwellKerr} proved decay towards the stationary Coulomb solution, and Metcalfe--Tataru--Tohaneanu \cite{MetcalfeTataruTohaneanuMaxwellSchwarzschild} proved $(t_*+r)^{-1}t_*^{-3}$-decay of the Maxwell field $F$ under conditions on initial data that guarantee that the final charge $c$ is zero. This followed their earlier work \cite{MetcalfeTataruTohaneanuPriceNonstationary} on Price's law \cite{PriceLawI,PriceLawII,PriceBurkoLaw,DonningerSchlagSofferPrice,DonningerSchlagSofferSchwarzschild,TataruDecayAsympFlat,HintzPrice,AngelopoulosAretakisGajicKerr} on nonstationary spacetimes.) The bulk of the present paper is concerned with the proof of linear decay estimates and partial asymptotics for (dynamical perturbations of) $\Box_{g_{\bhm,a}}$ (or, more generally, for structurally similar wave-type operators on Kerr\footnote{A subset of our methods can also be used to study linear and nonlinear wave equations on simpler spacetimes, such as non-trapping asymptotically Minkowski spacetimes, including in the presence of bound states. In a similar vein, our methods apply directly also to subextremal Kerr--Newman spacetimes \cite{NewmanCouchChinnaparedExtonPrakashTorrenceKN}, except the normal hyperbolicity of the trapped set is yet to be proved in the full subextremal range of masses, specific angular momenta, and charges; moreover, a proof of mode stability for tensorial equations such as Maxwell or linearized Einstein remains elusive.}) that can be used in a nonlinear iteration scheme.

Before stating a more precise version of Theorem~\ref{ThmI1}, let us consider a simple ODE analogue of~\eqref{EqI1} that features a zero energy bound state, namely
\begin{equation}
\label{EqIODE}
  P_0 u=f+|\pa_{t_*}u|^2 u,\quad u\colon[1,\infty)_{t_*}\to\C^2,\ P_0:=\begin{pmatrix} \pa_{t_*} & 0 \\ 0 & \pa_{t_*}+1 \end{pmatrix},
\end{equation}
where $f$ is supported in $t_*\geq 1$. The bound state is the $t_*$-independent ``function'' (i.e., constant vector) $u_0=(1,0)$. Consider first the linear equation $P_0 u=f=\cO(t_*^{-\alpha})$. If $\alpha<1$, then $u=\cO(t_*^{-\alpha+1})$ loses one power of $t_*$; for $\alpha>1$, we have $u=c u_0+\cO(t_*^{-\alpha+1})$, so the remainder term sees the same $t_*$-loss. Conversely, in order to close a nonlinear iteration scheme for the solution of~\eqref{EqIODE}, one would like the forward operator $\Phi\colon u\mapsto P_0 u-f-|\pa_{t_*}u|^2 u$ to map $u=c u_0+u^\flat$, $u^\flat=\cO(t_*^{-\alpha+1})$, back to a function of size $\cO(t_*^{-\alpha})$. But the linear part $P_0$ does not map $\cO(t_*^{-\alpha+1})\to\cO(t_*^{-\alpha})$. It is thus important to observe that the part $u^\flat=(u^\flat_1,u^\flat_2)$ of a linear solution of $P_0 u=f$ has more structure: namely, $u^\flat_2=\cO(t_*^{-\alpha})$ has the same decay as $f=(f_1,f_2)$, and only the bound state component $u^\flat_1$ has weaker decay; but it is given by $u^\flat_1(t_*)=\int_1^{t_*} f_1(s)\,\dd s$, and thus gains a power of $t_*$-decay upon differentiation in $t_*$. That is, we have the finer description
\begin{equation}
\label{EqIODEDecomp}
  u(t_*) = (c+a(t_*))u_0 + \tilde u(t_*),\quad a=\cO_1(t_*^{-\alpha+1})\ \text{($\C$-valued)},\ \ \tilde u\in\cO(t_*^{-\alpha})\ \text{($\C^2$-valued)},
\end{equation}
where we write $\cO_1(t_*^{-\alpha+1})$ for a function whose $t_*$-derivative is of size $\cO(t_*^{-\alpha})$. For such $u$ then, we do have $\Phi(u)\in\cO(t_*^{-\alpha})$ and can close a nonlinear iteration scheme for~\eqref{EqIODE} upon using proper function spaces; see~\S\ref{SssIPODE}.

\begin{figure}[!ht]
\centering
\includegraphics{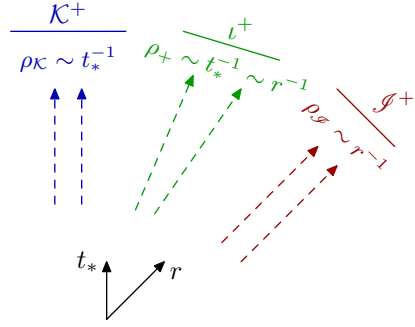}
\caption{Illustration of the three asymptotic regimes $\scri^+$, $\iota^+$, resp.\ $\cK^+$ of $\Omega^\circ$. We also draw ideal boundary hypersurfaces ``at infinity'' (later described via a compactification of $\Omega^\circ$ to a manifold with corners) where $\rho_\sscri$, $\rho_+$, resp.\ $\rho_\cK$ vanish; these contain the ``endpoints'' of curves with constant $t_*$ and $r\to\infty$ ($\scri^+$), constant $\frac{t}{r}\in(0,1)$ (i.e., $\frac{t_*}{r}\in(0,\infty)$) and $t_*\to\infty$ ($\iota^+)$, resp.\ constant $x$ and $t_*\to\infty$ ($\cK^+$).}
\label{FigIRegimes}
\end{figure}

In the context of~\eqref{EqI1}, we use $\alpha=1+\eps_\cK$ as the $t_*$-decay rate in compact spatial regions. But now there is more than one asymptotic regime. We thus briefly discuss what one should expect regarding the requirements for decay rates of source terms $f$ and linear waves $u$ in the various regimes in order to be able to close a nonlinear iteration.
\begin{enumerate}
\item[$(\cK^+)$] \textit{Spatially compact sets.} When $t_*\to\infty$ in spatially compact sets, one expects a decomposition of the form~\eqref{EqIODEDecomp}.
\item[$(\scri^+)$] \textit{Null infinity.} On Minkowski space, a source term $f$ for the scalar wave equation that is of size $\cO(r^{-\beta_\sscri-1})$ and compactly supported in $t_*$ (thus is nontrivial only near null infinity) typically produces a wave with $\cO(t_*^{-\beta_\sscri})$-decay for the solution of $P$ in the forward cone (as follows from simple estimates for the explicit forward solution); expecting a similar relationship to hold between $f$ and the ``regular part'' $\tilde u$ in a decomposition of the type~\eqref{EqIODEDecomp}, one should thus require $\beta_\sscri\geq\alpha$ (and so $\beta_\sscri>1$) for the decay rate of $u$ at $\scri^+$. There is a subtlety, though, since when $\beta_\sscri>1$, then $u$ itself is not of size $\cO(r^{-\beta_\sscri})$ at null infinity, but rather is given by a radiation field $\sim r^{-1}F(t_*,\theta,\phi)$ to leading order as $r\to\infty$ \cite{FriedlanderRadiation}.
\item[$(\iota^+)$] \textit{Punctured timelike infinity.} As already noted in \cite{HintzNonstat}, it turns out to be very convenient for the analysis (and for having a simple description of the asymptotics of metric perturbations in the nonlinear stability problem) to separate null infinity (where the radiation field dominates) and spatially compact regions (which are the origin of the bound state) by a third regime in which $\frac{r}{t}$ lies in a compact subset of $(0,1)$ and $t_*\gg 1$. For source terms $f$ localized in such a region and with $\cO(t_*^{-\beta_+-2})$-decay, the scalar wave sourced by $f$ on Minkowski space typically has $\cO(t_*^{-\beta_+})$-decay.
\end{enumerate}
To capture decay in these three asymptotic regions of $\Omega^\circ$, we use weights which are powers of the functions\footnote{For waves on Minkowski spacetime, the analogues of these weights are $\rho_\sscri=\frac{t-r}{t}$, $\rho_+=\frac{t}{(t-r)(r+1)}$, $\rho_\cK=\frac{r+1}{t}$.}
\begin{equation}
\label{EqIbdfs}
  \rho_\sscri := \frac{t_*}{t_*+r},\quad \rho_+=\frac{t_*+r}{t_* r},\quad \rho_\cK=\frac{r}{t_*+r}.
\end{equation}
Thus, in the first regime, $\rho_\cK$ tends to $0$ as $t_*\to\infty$ while $\rho_\sscri$ and $\rho_+$ remain bounded away from $0$; similarly for the other two regimes. See Figure~\ref{FigIRegimes}.

Furthermore, extrapolating from the regularity properties of $a$ in~\eqref{EqIODEDecomp}, we shall work with function spaces that capture arbitrary amounts of regularity with respect to $t_*\pa_{t_*}$; and for this to be consistent with wave evolution, we also need spatial $r\pa_x$-regularity. (For example, waves on Minkowski space only remain bounded in the forward cone upon applying arbitrary numbers of $t\pa_t$-derivatives if the initial data remain bounded upon applying arbitrary numbers of $\la x\ra\pa_x$-derivatives.) We shall refer to $t_*\pa_{t_*}$ and $r\pa_x$ as \emph{b-derivatives}, and regularity with respect to these as \emph{b-regularity}. (One can equivalently use the vector fields $t_*\pa_{t_*}+r\pa_r$ (spacetime scaling), $r\pa_r$ (weighted outgoing null-derivative), and spherical vector fields.) We thus introduce the weighted norm
\[
  \|u\|_{\Hb^{k,(\beta_\sscri,\beta_+,\beta_\cK)}(\Omega,\mu_\bop)}^2 := \sum_{j+|\alpha|\leq k} \int_{\Omega^\circ} | \rho_\sscri^{-\beta_\sscri}\rho_+^{-\beta_+}\rho_\cK^{-\beta_\cK} (t_*\pa_{t_*})^j (r\pa_x)^\alpha u(t_*,x)|^2\,\dd\mu_\bop,\quad \mu_\bop=\frac{\dd t_*\,\dd x}{t_* r^3}.
\]
We write $\Hb^{k,(\beta_\sscri,\beta_+,\beta_\cK)}(\Omega,\mu_\bop)$ for the space of functions for which this norm is finite; for $k=\infty$, this space is the intersection over all $k$. The insertion of the weight $\frac{1}{t_* r^3}$ in the volume density has the effect of making Sobolev embedding particularly simple, namely
\[
  | u | \leq C\rho_\sscri^{\beta_\sscri}\rho_+^{\beta_+}\rho_\cK^{\beta_\cK} \|u\|_{\Hb^{k,(\beta_\sscri,\beta_+,\beta_\cK)}(\Omega,\mu_\bop)}
\]
for $k\geq 3$, similarly for higher b-derivatives.

\begin{thm}[Nonlinear waves with zero energy bound states: precise version]
\label{ThmI1Q}
  Fix weights
  \[
    1 < \beta_\cK < \beta_+ < \beta_\sscri < \tfrac32.
  \]
  Then there exist $d\in\N$ and $\eps>0$ such that for all 1-forms $f$ on $M^\circ$ with support in $\{t_*\geq 1\}$ that satisfy
  \[
    f \in \Hb^{\infty,(\beta_\sscri+1,\beta_++2,\beta_\cK)}(\Omega,\mu_\bop),\quad
    \|f\|_{\Hb^{d,(\beta_\sscri+1,\beta_++2,\beta_\cK)}(\Omega,\mu_\bop)}<\eps,
  \]
  the equation~\eqref{EqI1} has a (unique) global forward solution $u$ on $\Omega^\circ$. This solution is of the form
  \begin{equation}
  \label{EqI1Qu}
    u(t_*,x) = (c + a(t_*))u_0(x) + \chi_\sscri r^{-1} u_\sscri\Bigl(t_*,\frac{x}{|x|}\Bigr) + \tilde u(t_*,x),
  \end{equation}
  where $a,u_\sscri$, and $\tilde u$ vanish for $t_*\leq 1$, and:
  \begin{enumerate}
  \item $c\in\C$, $a\in\Hb^{\infty,\beta_\cK-1}([1,\infty],|\frac{\dd t_*}{t_*}|)$ (meaning that $t_*^{\beta_\cK-1}a$ and all of its derivatives along $t_*\pa_{t_*}$ lie in $L^2([1,\infty),|\frac{\dd t_*}{t_*}|)$);
  \item each component of $u_\sscri$ in the frame $\dd t_*,\dd x$ lies in $\Hb^{\infty,\beta_+-1}([1,\infty]\times\Sph^2;|\frac{\dd t_*}{t_*}\,\dd g_{\Sph^2}|)$ (regularity with respect to $t_*\pa_{t_*}$ and rotation vector fields), and $\chi_\sscri=\chi_\sscri(\rho_\sscri)$ is supported in a small neighborhood of $\rho_\sscri=0$ (thus, ``near null infinity'');
  \item each component of $\tilde u$ in the frame $\dd t_*,\dd x$ lies in $\Hb^{\infty,(\beta_\sscri,\beta_+,\beta_\cK)}(\Omega,\mu_\bop)$.
  \end{enumerate}
\end{thm}

\begin{figure}[!ht]
\centering
\includegraphics{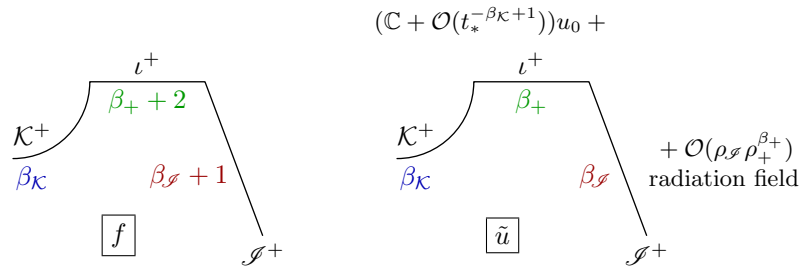}
\caption{Illustration of the decay rates of $f$ and the asymptotic behavior of (the pieces of) $u$ in~\eqref{EqI1Qu}.}
\label{FigI1Q}
\end{figure}

This is restated, using slightly more precise notation, in Theorem~\ref{ThmA2}. See Figure~\ref{FigI1Q} for an illustration. Note that the b-regularity of $a$ is related to the ``$\cO_1$'' behavior noted in~\eqref{EqIODEDecomp}. We discuss further context for the decision to work with b-derivatives in~\S\ref{SssIMb} below.

\begin{rmk}[Decay orders]
\label{RmkI1QDecay}
  The decay rate $\beta_\sscri$ of $\tilde u$ at null infinity is stronger than the $r^{-1}$-decay of the radiation field. The requirement $\beta_+<\beta_\sscri$ arises already for solutions of the wave equation of Minkowski space (the decay rate of sourced waves in the forward cone typically not being any better than at null infinity), similarly for $\beta_\cK<\beta_+$.
\end{rmk}

\begin{rmk}[Sharpness]
\label{RmkI1QSharp}
  The decay rates of $a,u_\sscri$, and $\tilde u$ in Theorem~\ref{ThmI1Q} are sharp, given the assumptions on $f$: indeed, for $u$ of the form~\eqref{EqI1Qu}, we have $\Phi_0(u):=\Box_{g_{\bhm,a}}u-|\delta_{g_{\bhm,a}}u|^2 u\in\Hb^{\infty,(\beta_\sscri+1,\beta_++2,\beta_\cK)}$. Of course, the precise relationship between the decay rates of $f$ and (pieces of) $u$ is not very important for the final result, but it is crucial for the solution operators for linearizations of $\Phi_0(u)$ in a nonlinear iteration scheme (as any mismatch would not allow one to close the iteration).
\end{rmk}

\subsubsection{Nonlinear applications, II: without zero energy bound states}
\label{SssIMN2}

Proving global existence for nonlinear wave equations in the \emph{absence} of zero bound states for the underlying time-translation-invariant linear wave operator is considerably less delicate. We mention here three simple examples, each of which is a nonlinear perturbation of the linear scalar wave equation on subextremal Kerr:

\begin{thm}[Nonlinear waves without zero energy bound states]
\label{ThmINo}
  For sufficiently small $f\colon\Omega^\circ\to\C$ vanishing for $t_*\leq 1$, there exist global forward solutions of the following nonlinear PDE.
  \begin{enumerate}
  \item\label{ItINoPower} {\rm (Power nonlinearity.)} $\Box_{g_{\bhm,a}}u=f+u^p$, $\N\ni p\geq 3$. (See Theorem~\usref{ThmA1Power} for $p\geq 4$ and Theorem~\usref{ThmDp} for $p=3$.)
  \item{\rm (Null-form nonlinearity.)} $\Box_{g_{\bhm,a}}u=f+g_{\bhm,a}(\nabla u,\nabla u)$. (See Theorem~\usref{ThmA1Null}.)
  \item\label{ItINoQ}{\rm (A toy quasilinear equation.)} $\Box_{g(u,\dd u)}u=f$, where $g(t_*,x,u,\dd u)-g_{\bhm,a}(t_*,x)$ has compact $x$-support, vanishes for $u=0$, and has coefficients that are smooth functions of $x,u(t_*,x),\dd u(t_*,x)$. (See Theorem~\usref{ThmA1Q}, and~\S\usref{SsA1Q} for the full setup.)
  \end{enumerate}
  More precisely, in each of these cases except for~\eqref{ItINoPower} with $p=3$, fix $\frac12=\beta_\cK<\beta_+<\beta_\sscri<1$, then there exist $d\in\N$ and $\eps>0$ such that for
  \[
    f\in\Hb^{\infty,(\beta_\sscri+1,\beta_++2,\beta_\cK)}(\Omega,\mu_\bop),\quad \|f\|_{\Hb^{d,(\beta_\sscri+1,\beta_++2,\beta_\cK)}}<\eps,
  \]
  the forward solutions of these equations satisfy
  \[
    u\in\Hb^{\infty,(\beta_\sscri,\beta_+,\beta_\cK)}(\Omega,\mu_\bop).
  \]
  For the case~\eqref{ItINoPower} with $p=3$, we require $1<\beta_\cK<\beta_+<\beta_\sscri<\frac32$ instead, and the solution $u$ is the sum of a radiation field term and an element of $\Hb^{\infty,(\beta_\sscri,\beta_+,\beta_\cK)}(\Omega,\mu_\bop)$ (see~\eqref{EqDpu}).
\end{thm}

These toy models have been (partially) discussed in previous works:
\begin{enumerate}[leftmargin=2em]
\item In the slowly rotating case $|\frac{a}{\bhm}|\ll 1$, global well-posedness was shown for $\Box_{g_{\bhm,a}}u=u^p$ for all (real) exponents $p>1+\sqrt{2}$ by Lindblad--Metcalfe--Sogge--Tohaneanu \cite{LindbladMetcalfeSoggeTohaneanuWangStrauss}. (The same exponent appears in the \emph{Strauss conjecture} on Minkowski spacetime; see, e.g., \cite{JohnBlowupPower,StraussConjecture,SchaefferBoxuup}.) This followed earlier work in the Schwarzschild case by Dafermos--Rodnianski \cite{DafermosRodnianskiNonlinear} ($p>4$, radial data), Blue--Sterbenz \cite{BlueSterbenzSemilinear} ($p>3$), and Marzuola--Metcalfe--Tataru--Tohaneanu \cite{MarzuolaMetcalfeTataruTohaneanuStrichartz} ($p=5$). While our result is the first one explicitly treating this equation in the full subextremal range, the methods of \cite{DafermosHolzegelRodnianskiTaylorQuasilinear,DafermosHolzegelRodnianskiTaylorQuasilinear2} should readily apply and provide an alternative approach.
\item Small data global well-posedness for $\Box_{g_{\bhm,a}}u=|\nabla u|^2$ was shown in \cite{DafermosHolzegelRodnianskiTaylorQuasilinear,DafermosHolzegelRodnianskiTaylorQuasilinear2}, with earlier work by Luk \cite{LukKerrNonlinear} settling the slowly rotating case.
\item Quasilinear equations $\Box_{g(u,\dd u)}u=0$ near a subextremal Kerr metric were studied (in the special case that $g=g(x,u)$ depends only on spatial variables and $u$) by Lindblad--Tohaneanu \cite{LindbladTohaneanuSchwarzschildQuasi,LindbladTohaneanuKerrQuasi,DafermosHolzegelRodnianskiTaylorQuasilinear} on Schwarzschild and slowly rotating Kerr spacetimes; the paper \cite{DafermosHolzegelRodnianskiTaylorQuasilinear2} covers the full subextremal range. (Unlike the references, we prove b-regularity here, but the flipside is that we use more derivatives.)
\end{enumerate}

A further example of a nonlinear toy model discussed in the literature (but not here) is the wave map type equation on slowly rotating Kerr studied by Ionescu--Klainerman \cite{IonescuKlainermanWavemap} and Stogin \cite{StoginKerrWaveMap}.

The reason why the cubic nonlinear wave equation ($p=3$) is more delicate than the case $p\geq 4$ is due to decay considerations, the most delicate place being the intermediate regime $\iota^+$: if $u$ has $t_*^{-\beta_+}$ decay in the regime $\iota^+$ where $\beta_+<1$ is close to $1$, then $u^p$ has $t_*^{-p\beta_+}$ decay, which is stronger than $t_*^{-\beta_+-2}$ only when $p\geq 4$, and thus the $\iota^+$-decay of $u^p$ is compatible with that of $f$. For $p=3$ on the other hand, we must have $\beta_+>1$ orders of $t_*$-decay at $\iota^+$ for $u$ to close a nonlinear iteration, which then also requires $u$ to have stronger-than-$r^{-1}$-decay at $\scri^+$ (modulo the presence of a radiation field); this is why for $p=3$ we need stronger decay of $f$ and capture a radiation field plus a stronger decaying remainder for $u$. (See~\S\S\ref{SsA1Power} and \ref{SsDp} for more detailed heuristics.) For a systematic discussion of the relationship between decay rates and acceptable nonlinearities, we refer the reader to \cite{KadarAsymptotics}.

\begin{rmk}[More general settings]
\label{RmkA1Gen}
  One can also mix all three settings without any difficulty, and consider much more general classes of nonlinear equations in the spirit of \cite{DafermosHolzegelRodnianskiTaylorQuasilinear,DafermosHolzegelRodnianskiTaylorQuasilinear2}. Since our goal is merely to illustrate the utility of Theorem~\ref{ThmF} (the main application being \cite{HintzKerrStab}), we do not aim for any generality here. In particular, the class of metric perturbations considered in Theorem~\ref{ThmINo}\eqref{ItINoQ} is highly artificial; we stress, though, that our methods do cover those metric perturbations which arise in nonlinear stability problems for the Einstein equations, e.g., near null infinity as analyzed in \cite{HintzVasyMink4,HintzMink4Gauge,HintzVasyScrieb}.
\end{rmk}

\subsubsection{Rough version of the linear main result}
\label{SssIML}

The starting points of the proofs of Theorems~\ref{ThmI1Q} and \ref{ThmINo} are precise estimates for forward solutions of linear wave-type equations; the proof of such estimates lies at the heart of the paper. In the context of these nonlinear results, the relevant linear wave operators are the linearizations of the nonlinear maps $\Phi(u)=\Box_{g_{\bhm,a}}u-f-|\delta_{g_{\bhm,a}}u|^2 u$, $\Box_{g_{\bhm,a}}u-f-u^p$, etc.\ around $u$ lying in the respective function spaces in which we seek the nonlinear solutions.

We first state a very rough version of the main result, and restrict attention only to the settings relevant for the aforementioned nonlinear applications:

\begin{thm}[Estimates for linear waves, rough version]
\label{ThmIML}
  Consider one of the following settings:
  \begin{enumerate}
  \item[\rm (i.a)] $P_0=\Box_{g_{\bhm,a}}$ is the scalar wave operator; let $\frac12<\beta_+<\beta_\sscri<1$ and set $\aleph=0$, $\delta=0$.
  \item[\rm (i.b)] $P_0=\Box_{g_{\bhm,a}}$ is the 1-form wave operator; let $\frac12<\beta_+<\beta_\sscri<1$, set $\aleph=1$, and fix any $\delta>\beta_+-\frac12$.
  \end{enumerate}
  Furthermore:
  \begin{enumerate}
  \myitem{ItIML2}{ii} Let $\ell_+>\delta$ and $\ell_\cK>\aleph$. Consider a metric $g=g_{\bhm,a}+h$ on $\Omega^\circ$ where the coefficients of $h$ (and all of its b-derivatives) in the frame $\dd t_*,\dd x$ are of size $\cO(\rho_+^{\ell_+}\rho_\cK^{\ell_\cK})$ away from $\scri^+$ (say, for $\frac{r}{t}<\frac{9}{10}$) and where $h$ has appropriately decaying coefficients near $\scri^+$; we assume that $t_*=1$ is a Cauchy hypersurface for $(\Omega^\circ,g)$.
  \myitem{ItIML3}{iii} Let $P$ be a wave-type operator on $(\Omega^\circ,g)$ (i.e., $P$ has the same principal part as $\Box_g$) such that the coefficients of $P-P_0$, expressed in terms of $r\pa_{t_*}$, $r\pa_x$, are of size $\cO(\rho_+^{\ell_++2}\rho_\cK^{\ell_\cK})$ away from $\scri^+$ and have an appropriate structure near $\scri^+$.
  \end{enumerate}
  Then there exists $d\in\N$ such that the following holds for all $k\in\N_0$: the forward solution $u$ of
  \begin{subequations}
  \begin{equation}
  \label{EqIMLPu}
    P u = f \in \Hb^{k+d,(\beta_\sscri+1,\beta_++2,\frac12)}(\Omega,\mu_\bop)
  \end{equation}
  (with $f$ vanishing for $t_*\leq 1$) satisfies
  \begin{equation}
  \label{EqIMLu}
    u \in \Hb^{k,(\beta_\sscri,\beta_+-\delta,-\aleph+\frac12)}(\Omega,\mu_\bop).
  \end{equation}
  \end{subequations}
  Moreover, $u$ satisfies \emph{tame estimates}, namely,
  \begin{equation}
  \label{EqIMLbTame}
    \|u\|_{\Hb^{k,(\beta_\sscri,\beta_+-\delta,-\aleph+\frac12)}} \leq C_k \Bigl( \|f\|_{\Hb^{k+d,(\beta_\sscri+1,\beta_++2,\frac12)}} + \|P-P_0\|_{k+d} \|f\|_{\Hb^{d,(\beta_\sscri+1,\beta_++2,\frac12)}} \Bigr),
  \end{equation}
  where $\|P-P_0\|_{k+d}$ is a weighted $L^\infty$-norm of the coefficients of $P-P_0$ and their b-derivatives of order $\leq k+d$, and $C_k$ is independent of $P$ provided $P$ is close to $P_0$ in some fixed low-regularity norm.
\end{thm}

See Theorem~\ref{ThmF} for the full statement. We discuss various aspects of this result:

\begin{enumerate}[leftmargin=2em]
\item{\rm (Tensorial equations.)} We work \emph{directly} at the level of tensorial wave-type operators. \emph{At no point do we rely on reductions to scalar problems.} (The only caveat is that, in the literature such as \cite{AnderssonHaefnerWhitingMode}, the verification of mode stability for tensorial stationary models $P_0$---which we use as a black box in the concrete settings under consideration---does use such a reduction.)
\item{\rm (The stationary model $P_0$, I.)} Theorem~\ref{ThmIML} applies more generally, with $\aleph=0$ and $\delta=0$, when the stationary model $P_0\equiv\Box_{g_{\bhm,a}}\bmod\Diff^1$ satisfies mode stability in the closed upper half plane, including at zero energy. The allowed range of values of $\beta_+$ depends on indicial roots of the zero energy operator $\wh{P_0}(0)$. See Theorem~\ref{ThmA1Gen} for the precise statement.
\item{\rm (The stationary model $P_0$, II.)} The general assumption is that $P_0$ be \emph{$\aleph$-admissible with $\sface$-loss $\delta$} (Definition~\ref{DefSSAlephAdm}); in practice, we only encounter integer $\aleph\in\N_0$. We allow for $P_0$ to act on sections of arbitrary vector bundles. Another example is the linearization of the Einstein field equations in a suitable generalized harmonic gauge, which has $\aleph=2$ (and some $\delta\in(0,1)$), as shown in \citestab{Theorem~\ref*{ThmAdm}} following closely related arguments in~\cite[\S\S{3.4--3.5}]{HintzGlueLocIII}. Roughly speaking, the quantity $\aleph\geq 0$ is either $0$, or equal to one plus the maximal order of $t_*$-growth of generalized zero energy states, by which we mean elements in $\ker P_0$ that are polynomials in $t_*$ and have fast enough spatial decay (typically $o(1)$ as $r\to\infty$). No such states exist for the scalar wave operator, while the 1-form wave operator does admit them, namely multiples of $u_0$ (and nothing else); and for the linearized gauge-fixed Einstein equation, deformation tensors of Lorentz boosts have $t_*^1$-growth (see \cite[Theorem~10.4 and Proposition~9.4]{HaefnerHintzVasyKerr} for small $|\frac{a}{\bhm}|$ and \cite[Proposition~6.8 and Lemma~6.9]{HaefnerHintzVasyKerrLarge} for the full subextremal range). Equivalently, $\aleph$ is the order of the pole of the resolvent $\wh{P_0}(\sigma)^{-1}$ at $\sigma=0$ (or, more precisely, due to the failure of meromorphicity, the largest power of $|\sigma|^{-1}$ appearing in its low-energy expansion when acting on $\CIc$ inputs). --- The condition of $\aleph$-admissibility requires mode stability at nonzero frequencies in the closed upper half plane; but due to the wide range of possibilities for the low or zero energy behavior in the presence of zero energy bound states, we do not give a spectral characterization of $\aleph$-admissibility in this paper when $\aleph>0$.
\item{\rm (Structure near $\scri^+$.)} We will state precise assumptions on $h$ and $P$ near $\scri^+$ in and after~\eqref{EqIPmlPert}. They involve anisotropic decay conditions for the metric perturbation $h$ and the coefficients of $P-P_0$. For now, we only mention that these assumptions
  \begin{enumerate}
  \item guarantee that $\scri^+$ is the correct location of null infinity not only for $g_{\bhm,a}$ but also for $g=g_{\bhm,a}+h$ (i.e., $h$ contributes only $o(r^{-1})$ to the metric component $g^{-1}(\dd t_*,\dd t_*)$);
  \item are consistent with the metric perturbations arising in the stability of Minkowski space in generalized harmonic gauge \cite{LindbladRodnianskiGlobalExistence,LindbladRodnianskiGlobalStability,LindbladAsymptotics,HintzVasyMink4,HintzMink4Gauge};
  \item are satisfied for the linearization of the nonlinear wave equations described in~\S\S\ref{SssIMN}--\ref{SssIMN2}.
  \end{enumerate}
\item\label{ItIMLDecay}{\rm (Decay orders of sources and solutions.)} The pointwise bounds implied by~\eqref{EqIMLu} via Sobolev embedding are
  \begin{equation}
  \label{EqIMLuPointwise}
    |u| \lesssim \rho_\sscri^{\beta_\sscri}\rho_+^{\beta_+-\delta}\rho_\cK^{-\aleph+\frac12}\quad\text{on}\ \Omega^\circ,
  \end{equation}
  likewise for up to $k-3$ b-derivatives of $u$. This thus differs from the decay rates of $f$ by one order at $\scri^+$ as usual (as, e.g., in Figure~\ref{FigI1Q}), by $2+\delta$ orders at $\iota^+$ (which is $\delta$ more than usual), and by $\aleph$ orders at $\cK^+$. The latter loss is to be expected already by comparison with the ODE example~\eqref{EqIODE}---which (by analogy) has $\aleph=1$. The $\delta$-loss at $\iota^+$ is more subtle; it arises for the 1-form wave operator $P=P_0=\Box_{g_{\bhm,a}}$ from the lack of decay of a term $a(t_*)u_0$, $a\in\Hb^{k,\frac12}([1,\infty],|\frac{\dd t_*}{t_*}|)$, in the late-time description of $u$ (this term essentially arising from a $t_*$-integral of $f$).\footnote{This refers to~\eqref{EqA2Admau0}. There we use the metric density instead of $\mu_\bop$, and do not resolve null infinity, hence the shift of orders by $(-\frac32,-\frac12)$ and the absence of a $\scri^+$-decay rate.} --- Only in the case $\aleph=0$ does~\eqref{EqIMLuPointwise} entail decay. (For $\beta_\sscri$ and $\beta_+$ close to $1$, for example, it gives $r^{-\frac12+\eps}t^{-\frac12}$-decay on $\Omega^\circ$.) 
\item{\rm (Decay orders of perturbations.)} In light of the losses $\delta$ and $\aleph$ just discussed, the requirement that the perturbation $P-P_0$ decay with rates $\ell_+>\delta$ and $\ell_\cK>\aleph$ allows one to use $P_0$ as an effective late-time model for $P$ (as far as controlling decay is concerned). Let us illustrate the relationship of $\aleph$ and $\ell_\cK$ with simple ODE examples:
  \begin{enumerate}
  \item $P_0=\pa_{t_*}$ (with $\aleph=1$). Then forward solutions of $(P_0+\tilde p(t_*))u=f\in\CIc(\R_{t_*})$ are given at late times by $c\exp(-\int_1^{t_*} \tilde p(s)\,\dd s)$ and thus (in the absence of sign conditions on $\tilde p$) typically grow super-polynomially unless $|\tilde p|=\cO(t_*^{-1})$; and the condition $|\tilde p|=\cO(t_*^{-\ell_\cK})$ with $\ell_\cK>1$ guarantees that the asymptotics of $u$ are qualitatively the same as those of solutions of $P_0$. (By contrast, generic perturbations $\tilde p$ of size $t_*^{-1}$ lead to modified decay rates.) --- Another example is given by~\eqref{EqIODE}: the linearization of $u\mapsto P_0 u-f-|\pa_{t_*}u|^2 u$ around $u$ as in~\eqref{EqIODEDecomp} is given (to leading order in the sense of decay, and for $\R^2$-valued $u$) by $P_u\dot u=P_0\dot u-2 u(\pa_{t_*}u)\cdot\pa_{t_*}\dot u$; the coefficient of $\pa_{t_*}\dot u$ on the right has $\cO(t_*^{-\alpha})$-decay, and we recall that we are taking $\alpha(\,=:\ell_\cK)>1$.
  \item $P_0=\pa_{t_*}^2$ (with $\aleph=2$). Now solutions of $(P_0+\tilde p(t_*))u=f\in\CIc$ have the same asymptotic behavior as those of $P_0$ only when $|\tilde p|=\cO(t_*^{-\ell_\cK})$, $\ell_\cK>2$.
  \end{enumerate}
\end{enumerate}

In the case $\aleph=\delta=0$, one can immediately deduce from the decay estimates of Theorem~\ref{ThmIML} (cf.\ point~\eqref{ItIMLDecay} above) global results for nonlinear wave equations. This involves two ingredients:
\begin{enumerate}
\item \textit{Nash--Moser iteration.} We refer to~\eqref{EqIMLbTame} as \emph{b-tame estimates}, as they are tame in the b-regularity order. (The proof involves also other, weaker, notions of regularity.) This means that a high-regularity norm of $u$ is bounded by a high-regularity norm of $f$ plus the product of a high-regularity norm of the coefficients of the operator $P$ and a low-regularity norm of $f$: there are no products of two high-regularity norms. See \cite{HamiltonNashMoser} for background and examples. We use here the simple version of Nash--Moser developed by Saint Raymond \cite{SaintRaymondNashMoser}. In the simplest case, then, if $\Phi(u)$ is a nonlinear wave operator (such as $\Phi(u)=\Box_{g_{\bhm,a}}u-f-u^p$) such that, for all $u\in\Hb^{\infty,(\beta_\sscri,\beta_+,\frac12)}$ with small low-regularity norm, the linearization $\Psi(u):=D_u\Phi$ satisfies the assumptions of Theorem~\ref{ThmIML}, then the equation $\Phi(u)=0$ admits a solution in $\Hb^{\infty,(\beta_\sscri,\beta_+,\frac12)}$, provided $\Phi(0)=-f$ has sufficiently small low-regularity norm.
\item \textit{Structure of the nonlinearity.} The forward map $\Phi(u)$ must map the space~\eqref{EqIMLu} back into the space of source terms in~\eqref{EqIMLPu}. This is only nontrivial for the nonlinear term:
  \begin{enumerate}
  \item The nonlinearity $u^p$, $p\geq 4$, is acceptable, as are simple quasilinear perturbations. This yields parts of Theorem~\ref{ThmINo}\eqref{ItINoPower} and \eqref{ItINoQ}.
  \item The null-form nonlinearity is acceptable. We only discuss this near $\scri^+$, which is the key region where, on (perturbations of) $(3+1)$-dimensional Minkowski space, the structure of nonlinearities matters, as was first realized by Klainerman \cite{KlainermanNullCondition} and Christodoulou \cite{ChristodoulouGlobalSolutionsSmallData} (see also \cite{LindbladRodnianskiWeakNull,LindbladRodnianskiGlobalStability,KeirWeak}). The key point, on Minkowski space (for notational clarity) and pretending that the nonlinearity is given by $(\pa_t+\pa_r)u\cdot(\pa_t-\pa_r)u$, is that the b-regularity of $u$ entails regularity (i.e., without decay losses) with respect to the \emph{weighted} outgoing null vector field $r(\pa_t+\pa_r)$ (i.e., $r\pa_r$ in the coordinates $t_*,r$), so application of $\pa_t+\pa_r$ \emph{gains one power of $r^{-1}\sim\rho_\sscri\rho_+$}. (By contrast, the incoming null vector field $\pa_t-\pa_r$ is, near $\scri^+$, the product of $t_*^{-1}$ and the b-vector field $t_*(\pa_t-\pa_r)$, and hence its application only gains one power of $t_*^{-1}\sim\rho_+$.)
  \end{enumerate}
  The details are given in~\S\ref{SA1}.
\end{enumerate}

When $\aleph\geq 1$ on the other hand, then~\eqref{EqIMLuPointwise} gives at best a $t_*^{\frac12}$-bound for $u$, which is of course not immediately useful for any (non)linear purposes. Furthermore, unless $\aleph=\delta=0$, the description~\eqref{EqIMLu} is not sharp as far as the relationship of the decay rates of $u$ and $f$ is concerned (meaning that $P$ does not map $u$ of class~\eqref{EqIMLu} back into the space~\eqref{EqIMLPu}). Thus,
\begin{equation}
\label{EqIMLKey}
  \text{\parbox{0.8\textwidth}{\it the key point of Theorem~\usref{ThmIML} is that while it only gives weak polynomial bounds on $u$, it entails \emph{arbitrarily high (b-)regularity} of $u$, with b-tame estimates.}}
\end{equation}
This is then the starting point for the extraction of stronger asymptotics and decay of $u$ (assuming stronger decay of $f$) in the contexts of Theorems~\ref{ThmI1Q} and \ref{ThmINo}\eqref{ItINoPower} with $p=3$, such as radiation fields via integration along approximate characteristics near $\scri^+$ (\S\ref{SsDscri}), or partial asymptotic expansions as $t_*\to\infty$ via rewriting $P u=f$ as
\[
  P_0 u = f - (P-P_0)u
\]
and inverting $P_0$ using spectral theory and (low-energy) resolvent analysis to improve control of $u$ (\S\S\ref{SsDRes}--\ref{SsDp}); we discuss this in~\S\ref{SssIPD}. We remark already here that these methods use up b-derivatives in return for a gain in decay. In the context of Theorem~\ref{ThmI1Q}, where we have already commented on the necessity of capturing at least one order of b-regularity of $a(t_*)$ in~\eqref{EqI1Qu} (i.e., $t_*\pa_{t_*}$-regularity), this demonstrates the considerable utility of having a large surplus of b-derivatives at one's disposal.

\bigskip

Keeping in mind~\eqref{EqIMLKey}, we are, in this paper, not concerned with establishing optimal decay for linear or nonlinear solutions (when the forcing terms or initial data have strong decay, say). In particular, we only push the resolvent analysis of the scalar or 1-form wave operators on subextremal Kerr far enough to prove Theorems~\ref{ThmI1Q} and \ref{ThmINo}, which yields the asymptotic description~\eqref{EqI1Qu} for linear 1-form waves and almost $t_*^{-2}$-decay for scalar waves (see Remark~\ref{RmkDpLinear}). Decay estimates for scalar wave-type equations on Schwarzschild have a long history, starting with works of Wald and Kay \cite{WaldSchwarzschild,KayWaldSchwarzschild} who proved boundedness on Schwarzschild; subsequent and related works include \cite{BlueSofferSchwarzschildDecay} and the aforementioned \cite{DonningerSchlagSofferPrice,DonningerSchlagSofferSchwarzschild} and \cite{AngelopoulosAretakisGajicRNPrice} on Price's law (see also \cite{DafermosRodnianskiPrice} for a nonlinear version in spherical symmetry). On Kerr spacetimes, the first proof of decay in the full subextremal range was given by Dafermos--Rodnianski--Shlapentokh-Rothman \cite{DafermosRodnianskiShlapentokhRothmanDecay} following earlier work by many authors including the pioneering works by Andersson--Blue \cite{AnderssonBlueHiddenKerr} and Tataru--Tohaneanu \cite{TataruTohaneanuKerrLocalEnergy} in the slowly rotating regime; see \cite{MarzuolaMetcalfeTataruTohaneanuStrichartz,TohaneanuKerrStrichartz} for Strichartz estimates. Sharp decay was proved in \cite{TataruDecayAsympFlat,HintzPrice,AngelopoulosAretakisGajicKerr}, and on dynamical perturbations in \cite{MetcalfeTataruTohaneanuPriceNonstationary,LukOhTwoTails}. Besides the aforementioned papers \cite{DafermosHolzegelRodnianskiTaylorQuasilinear,DafermosHolzegelRodnianskiTaylorQuasilinear2}, we also mention recent work by Ma--Szeftel \cite{MaSzeftelEnergyKerr,MaSzeftelTeukolsky} on integrated local energy decay estimates for the scalar wave equation and the Teukolsky equation on dynamical Kerr spacetimes. Radiation fields and scattering isomorphisms were studied, e.g., in \cite{BaskinWangRad,DafermosRodnianskiShlapentokhRothmanScattering}. More relevant to the (linear) stability problem of Kerr black holes is the Teukolsky equation \cite{TeukolskySeparation}, for which sharp decay estimates were proved (for all half-integer spins) in the full subextremal range by Millet \cite{MilletTeukolskyDecay}, with prior results including \cite{MaZhangSchwarzschildDiracSharp,MaZhangTeukolsky}; integrated local energy decay estimates are due to Shlapentokh-Rothman and Teixeira da Costa \cite{ShlapentokhRothmanTeixeiradCTeukolskyI,ShlapentokhRothmanTeixeiradCTeukolskyII}, following the earlier work \cite{DafermosHolzegelRodnianskiTeukolsky}. A goal of some of these works (but not of the present paper) is to obtain an energy boundedness statement without a loss of derivatives; that this is a fairly robust feature in the scalar setting was shown by Shlapentokh-Rothman--Tohaneanu \cite{ShlapentokhRothmanTohaneanuBdd}.

We do not study (perturbations of) extremal Kerr metrics here; even for the scalar wave operator on extremal Kerr, little is known beyond boundedness in axisymmetry \cite{AretakisExtremalKerr}, conditional decay results \cite{GajicExtremalRN}, and mode stability \cite{TeixeiradCModes}. Our methods do, in principle (i.e., up to checking spectral assumptions on the linear stationary wave operators one is interested in, and pending a full verification of the dynamical hypotheses, especially regarding the dynamics near the trapped set), apply directly to subextremal Kerr--Newman black holes, but not to any extremal Kerr--Newman spacetime such as extremal Reissner--Nordstr\"om \cite{AngelopoulosAretakisGajicExtremalRN}. The reason is that the subextremality of the event horizon is crucially used in our regularity theory, both for perturbed wave equations on spacetime and in the analysis of the spectral family of the stationary model.

The main tensorial equations studied on Kerr are the Maxwell equations (and the closely related 1-form wave equation) and linearized Einstein equations; literature on the former includes \cite{AnderssonBlueMaxwellKerr,MetcalfeTataruTohaneanuMaxwellSchwarzschild}, while the latter are an active field of study; the linear stability of subextremal Kerr black holes was shown \cite{AnderssonBackdahlBlueMaKerr} (unconditional for small $|\frac{a}{\bhm}|$, conditional in the full subextremal range) and \cite{HaefnerHintzVasyKerrLarge} (unconditional in the full subextremal range).

\medskip

In previous work \cite{HintzNonstat}, the author studied a large class of wave-type equations on asymptotically flat and dynamical but asymptotically stationary spacetimes and proved decay and \mbox{(b-)}re\-gu\-lar\-i\-ty estimates. Let us thus compare \cite{HintzNonstat} and the present work:
\begin{enumerate}[leftmargin=2em]
\item Unlike \cite{HintzNonstat}, we allow for the presence of trapped null-geodesics here. We exploit the fact that the trapped set of subextremal Kerr is $\sfr$-normally hyperbolic for every $\sfr$ \cite{HirschPughShubInvariantManifolds}; we give a self-contained proof of this fact here, following Wunsch--Zworski \cite{WunschZworskiNormHypResolvent,WunschZworskiNormHypResolventCorrection} and Dyatlov \cite{DyatlovWaveAsymptotics}. Estimates at the trapped set for dynamical perturbations of Kerr were proved by the author in \cite{HintzPolyTrap}.
\item Unlike \cite{HintzNonstat}, we allow for an event horizon (required to be non-degenerate, i.e., subextremal). This is easily handled using a variant of the radial estimates on spacetime proved in \cite{HintzVasySemilinear} following \cite{VasyMicroKerrdS}.
\item Unlike \cite{HintzNonstat}, we allow for the stationary model operator $P_0$ to have zero energy bound states. The precise spectral properties of $P_0$ are black-boxed as a spacetime estimate for forward solutions of $P_0$ (Definition~\ref{DefSSAlephAdm}), which in applications is proved using (in particular low-energy) resolvent estimates (see~\S\S\ref{SsA1Adm} and \ref{SsA2Adm} for the scalar and 1-form wave operators, respectively).
\item In \cite{HintzNonstat}, sharp decay is proved for stationary wave operators under suitable non-degeneracy conditions---which are violated for the scalar operator on Minkowski space, say, but which is satisfied for the wave operator coupled to a generic asymptotically inverse square potential \cite[\S{6.2}]{HintzNonstat}.
\item Both in \cite{HintzNonstat} and here, we prove the b-regularity of solutions of linear wave equations. The two key advances in the present paper in this context are the following.
  \begin{enumerate}
  \item We prove \emph{tame estimates}, which enable us to solve nonlinear wave equations here and in \cite{HintzKerrStab}. Our proofs of tame estimates use the structure of commutators of $P$ with (suitable) b-vector fields; we adapt ideas introduced in a different analytic setting in \cite{HintzGlueLocII}. (See the end of~\S\ref{SssIPb}.) In particular, we do not need to develop a theory of rough-coefficient pseudodifferential operators (unlike in \cite{BealsReedMicroNonsmooth,HintzQuasilinearDS,HintzVasyQuasilinearKdS}), at the minor expense of not getting explicit bounds on the number of b-derivatives required in Theorem~\ref{ThmIML} and in the nonlinear applications.
  \item We prove higher b-regularity in the presence of bound states: this cannot be proved in a transparent manner ``after the fact'' via commutation arguments with $t_*\pa_{t_*}$; we will discuss this issue (and our resolution of it) in~\S\ref{SssIPb}.
  \end{enumerate}
\item Neither here nor in \cite{HintzNonstat} do we handle nonzero modes (purely oscillatory or exponentially growing). These are allowed for in the linear work \cite{MetcalfeSterbenzTataruLED}. Exponentially growing modes are considered unsurmountable obstacles to robust global nonlinear existence results, unless one can eliminate them by modifying the equation (as in \cite{HintzVasyKdSStability,HintzPetersenVasyKdS}) or by working with positive codimension data (as, e.g., in \cite{KriegerNakanishiSchlagCenterStableNLW}, or to some extent in \cite{DafermosHolzegelRodnianskiSchwarzschildStability}).
\item The heart of the analysis both in \cite{HintzNonstat} and the present paper relies on the frameworks for the (microlocal) analysis of waves near null infinity $\scri^+$ developed by the author with Vasy in \cite{HintzVasyScrieb} and in the forward cone by the author in \cite{Hintz3b}.
\end{enumerate}

\subsubsection{On the notion of regularity}
\label{SssIMb}

We give further context for our decision to work with function spaces capturing \emph{b-regularity} in Theorem~\ref{ThmIML} and its nonlinear applications.

First of all, we stress that we work with spaces encoding regularity with respect to an underlying function space \emph{on spacetime} which here we take to be
\[
  L^2(\Omega^\circ)
\]
(or a weighted version thereof). Working thus with \emph{spacetime-$L^2$-spaces} is essential for our analysis in two respects: it is only relative to $L^2$ that the microlocal tools we will employ (in particular, the propagation of regularity) are applicable; and we prove estimates for the stationary model $P_0$ (which are ``glued into'' appropriate estimates for the dynamical operator $P$, for the purpose of controlling decay properties of solutions of $P$) using the Fourier transform in $t_*$ (and partly also using the Mellin transform with respect to spacetime scalings), which is cleanest on $L^2(\R_{t_*},|\dd t_*|)$ by Plancherel's theorem.\footnote{Note that passage from the density $\mu_\bop$ in~\eqref{EqIMLPu} to the metric density $|\dd g_{\bhm,a}|$ causes a shift in the weights, giving $f\in\Hb^{k+d,(\beta_\sscri-\frac12,\beta_+,0)}(\Omega,|\dd g_{\bhm,a}|)$, and the $\cK^+$-order $0$ is what allows for a clean description of such source terms on the Fourier transform side.} (By contrast, we do not rely on Duhamel's principle and dispersive $L^1\to L^\infty$ type estimates.)

Relative to $L^2$ then, one would certainly like to prove regularity with respect to the standard vector fields
\begin{equation}
\label{EqIMbTrans}
  \pa_{t_*},\ \pa_{x^1},\ \pa_{x^2},\ \pa_{x^3}.
\end{equation}
Stronger notions of regularity are important for nonlinear applications. An early instance of this was the observation by Klainerman \cite{KlainermanUniformDecay} that regularity with respect to generators of the Poincar\'e group (i.e., translations~\eqref{EqIMbTrans}, spatial rotations $x^i\pa_{x^j}-x^j\pa_{x^i}$, and Lorentz boosts $x^i\pa_t+t\pa_{x^i}$) and the scaling vector field $t\pa_t+x\cdot\pa_x=t_*\pa_{t_*}+x\cdot\pa_x$ leads to improved pointwise decay estimates. (A more basic example relevant to the present paper is that if $a\in L^2(\R_{t_*},|\dd t_*|)$, $\supp a\subset\{t_*\geq 1\}$, and $t_*\pa_{t_*}a\in L^2$, then $|a|\lesssim t_*^{-\frac12}$; cf.\ also the $\rho_\cK^{\frac12}$-decay in~\eqref{EqIMLuPointwise} for $\aleph=0$.) The fact that these vector fields include (as linear combinations with bounded and smooth coefficients), for instance, $r$-weighted outgoing null derivatives (but not incoming null derivatives) is relevant for structural requirements (e.g., null conditions) on nonlinearities. Regularity with respect to such ``Klainerman vector fields'' is typically the strongest sensible notion on asymptotically Minkowski spaces (absent the possibility of using further conformal Killing vector fields), and indeed holds quite generally \cite{BaskinVasyWunschRadMink,BaskinVasyWunschRadMink2} and plays an important role in nonlinear stability proofs \cite{ChristodoulouKlainermanStability,BieriZipserStability,LindbladRodnianskiGlobalStability,HintzVasyMink4}.

Now, solutions of wave equations cannot be regular with respect to boosts in the presence of nontrivial (asymptotically) stationary perturbations such as potentials $0\neq V=V(x)$ or stationary metrics (other than Minkowski) such as Kerr, since each application of $t\pa_{x^i}({}+x^i\pa_t)$ will produce $t$-growth. One is thus led to drop boosts and retain only the smaller collection
\begin{equation}
\label{EqIMbReg}
\begin{alignedat}{2}
  &\text{rotations:}&&\ x^i\pa_{x^j}-x^j\pa_{x^i}\ \text{(schematically written as $\pa_\omega$)}, \\
  &\text{spacetime scaling:}&&\ t\pa_t+x\pa_x, \\
  &\text{outgoing null regularity:}\ &&\ r(\pa_t+\pa_r)
\end{alignedat}
\end{equation}
of vector fields in the region $r\gtrsim 1$ (in addition to~\eqref{EqIMbTrans}). In the coordinates $(t_*,r)$, $t_*=t-r$, the weighted outgoing null vector field reads $r\pa_r$, and one can instead of~\eqref{EqIMbReg} equivalently use
\begin{equation}
\label{EqIMbReg2}
  t_*\pa_{t_*},\ r\pa_x.
\end{equation}
These are precisely the vector fields we are using in this paper (see the discussion preceding Theorem~\ref{ThmI1Q}), and also in \cite{HintzNonstat}. We also recall that linear waves on subextremal Kerr spacetimes are regular with respect to these, as shown for scalar waves and the Teukolsky equation in \cite{HintzPrice,MilletTeukolskyDecay} and for linearized metric perturbations in \cite{HaefnerHintzVasyKerr,HaefnerHintzVasyKerrLarge}. The double dyadic decomposition in \cite[\S{3.1}]{MetcalfeTataruTohaneanuPriceNonstationary} effectively also amounts to an encoding of b-regularity in the forward cone. Most other linear works only keep track of weaker notions of regularity, e.g., \cite{MetcalfeTataruTohaneanuPriceNonstationary,MorganWunschPrice} drop the outgoing null vector field and \cite{DafermosRodnianskiShlapentokhRothmanDecay} only uses~\eqref{EqIMbTrans} and rotations.

While the vector fields~\eqref{EqIMbTrans}--\eqref{EqIMbReg} yield the strongest sensible notion of regularity for solutions of wave equations on (perturbations) of Kerr spacetimes, it is possible to use weaker notions of regularity to close nonlinear iteration schemes or bootstrap arguments.\footnote{It is a general rule that one should aim to close a nonlinear iteration using the weakest possible function spaces, the expectation being that weaker notions of regularity are easier to recover. We demonstrate in this paper, however, that full b-regularity is highly convenient (especially when dealing with zero energy states, cf.\ the term $a(t_*)$ in~\eqref{EqI1Qu}) and robust, and provable in a conceptually transparent manner.} For example, the works \cite{DafermosHolzegelRodnianskiTaylorQuasilinear,DafermosHolzegelRodnianskiTaylorQuasilinear2} use~\eqref{EqIMbTrans} and rotation vector fields, and encode weighted control needed to capture null conditions via weighted energy estimates (see~\cite[Assumption~4.7.1]{DafermosHolzegelRodnianskiTaylorQuasilinear}). The stability proof of slowly rotating Kerr spacetimes in \cite{KlainermanSzeftelKerr} on the other hand effectively uses the vector fields~\eqref{EqIMbReg} \emph{except for} spacetime scaling; in a splitting of spacetime as (a subset of) $\R_{t_*}\times\R^3$, this amounts to standard $t_*$-regularity and b-regularity (i.e., $r\pa_r,\pa_\omega$) on $\R^3$. This notion of regularity is quite natural also from a spectral point of view and is used, e.g., in \cite[\S{13}]{HaefnerHintzVasyKerr}. It is, however, somewhat delicate: for example, if one adds to the Kerr metric a size $\cO(t_*^{-\alpha})$-perturbation in the transition region $\frac{r}{t}\in(\frac13,\frac23)$, say, then one should only expect regularity with respect to the same vector fields when $\alpha>1$ since, in the case $\alpha<1$, the outgoing light cones get shifted by growing $r$-dependent amounts in the transition region as $r\to\infty$. (The case $\alpha=1$ is borderline.) In such a scenario, then, one would need to adapt the choice of these vector fields to the dynamical geometry (or, equivalently, normalize the dynamical geometry via pullback along a suitable diffeomorphism that approximately straightens out the outgoing light cones). Full b-regularity, by contrast, is more robust, since it does not single out outgoing light cones in this transition region. (Note that the decay rate $\ell_+$ of perturbations in Theorem~\ref{ThmIML} only needs to exceed $\delta$, with $\delta$ being small in our applications.)

\medskip

We end this discussion by relating some of these notions of regularity to b-regularity on manifolds with corners in the sense of \cite{MelroseMendozaB,MelroseAPS}. (This was already discussed in \cite{HintzNonstat} and features in some form also in \cite{WangRadiation,HintzVasySemilinear,BaskinVasyWunschRadMink,SussmanKG}.) Recall that for a manifold $M$ with corners, locally modelled on $[0,\infty)_x^k\times\R^{n-k}_y$ for varying $k$, the space $\Vb(M)\subset\cV(M)=\CI(M;T M)$ of \emph{b-vector fields} is the Lie subalgebra of smooth vector fields that are tangent to each boundary hypersurface of $M$; in a local chart, these are $\CI$-spans of the vector fields $x^i\pa_{x^i}$ and $\pa_{y^j}$. On the \emph{radial compactification} $\ol{\R^4}$ of $\R^4$, given by attaching a sphere $\Sph^3$ ``at infinity'' whose neighborhood in $\ol{\R^4}$ is given by $[0,1)_\varrho\times\Sph^3$ where $\varrho$ is an inverse radius (see~\eqref{EqCRadCpt}), the space $\Vb(\ol{\R^4})$ is spanned by the coordinate vector fields on $\R^4$, scaling $-\varrho\pa_\varrho$, and rotations on $\Sph^3$. Let us use the coordinates $t\in\R$ and $x\in\R^3$ on $\R^4=\R\times\R^3$; then we may pass to the \emph{blow-up} of $\ol{\R^4}$ at the ``light cones at infinity'' given by the closures of $Y_\pm:=\{t\mp r=0\}$, $r=|x|$, in $\ol{\R^4}$ with the boundary $\pa\ol{\R^4}$ at infinity, which is the set $\{\frac{t}{r}\mp 1=0,\ \frac{1}{r}=0\}$; this amounts to replacing $\ol{\R^4}$ by the smooth manifold $[\ol{\R^4};Y_-,Y_+]$ with corners in which polar coordinates around $Y_\pm$ are smooth down to the origin: $Y_\pm\cong\Sph^2$ is thus replaced by the \emph{front face} $\ol\R\times\Sph^2$, with $t\mp r$ being a projective coordinate in the first factor (and $r^{-1}\sqrt{(t\mp r)^2+1}$ and $\arctan(t\mp r)$ are local radial and angular coordinates relative to $Y_\pm$). In particular, while the boundary at infinity of $\{t-r=c\}$ in $\ol{\R^4}$ is equal to $Y_+$ independently of $c\in\R$, this is no longer the case on $[\ol{\R^4};Y_-,Y_+]$. One can then check that
\[
  \Vb([\ol{\R^4};Y_-,Y_+])
\]
is spanned (over the space of smooth functions on this manifold) by the Klainerman vector fields, i.e., spacetime translations, spacetime scaling, rotations, and Lorentz boosts; this essentially follows from \cite[Proposition~2.1]{SussmanKG}. It is thus on the space $[\ol{\R^4};Y_-,Y_+]$ that solutions of wave equations on Minkowski space (and similar spacetimes, e.g., the Lorentzian scattering spaces of \cite{BaskinVasyWunschRadMink}) are b-regular. The front face arising from blowing up $Y_\pm$ is future/past null infinity $\scri^\pm$.

\begin{figure}[!ht]
\centering
\includegraphics{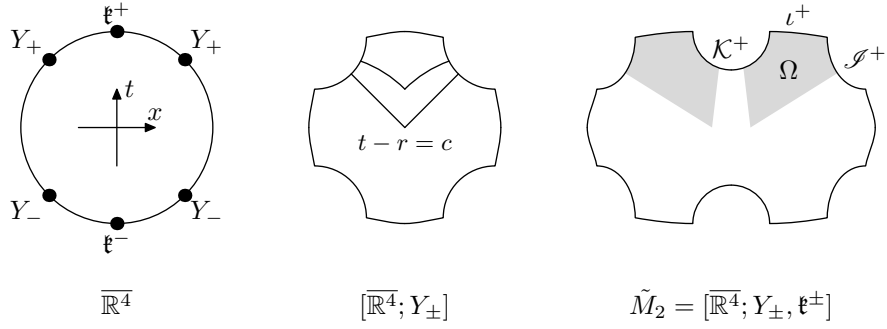}
\caption{Compactifications and blow-ups, and the compactified domain $\Omega$ from~\eqref{EqIOmega}. In the middle, we show two outgoing light cones $t-r=c$. On the domain $\Omega$ on the right, the functions $\rho_\sscri$, $\rho_+$, and $\rho_\cK$ from~\eqref{EqIbdfs} vanish precisely at $\scri^+$, $\iota^+$, and $\cK^+$, respectively (and simply so).}
\label{FigICpt}
\end{figure}

The coefficients of the Kerr metric are not smooth on $[\ol{\R^4};Y_-,Y_+]$ since near the ``north/south pole,'' i.e., the point $\fk^\pm$ given by $\frac{x}{t}=0$, $\frac{1}{t}=0$ at $t=\pm\infty$, they are not smooth functions of $(\frac{x}{t},\frac{1}{t})$. To resolve them, we thus blows up also $\fk^\pm$; this replaces them by a copy of $\ol{\R^3}_x$. On the resulting manifold
\begin{equation}
\label{EqIM1tilde}
  \tilde M_2 := [\ol{\R^4};Y_\pm,\fk^\pm],
\end{equation}
the space of b-vector fields is then spanned by~\eqref{EqIMbReg}.\footnote{The general result, using the terminology of \cite{MelroseDiffOnMwc}, is that when $Y\subset M$ is a p-submanifold, then $\Vb([M;Y])$ is spanned by the lifts of all elements of $\Vb(M)$ that are tangent to $Y$. In the present context, tangency to $\{\fk^\pm\}$ eliminates Lorentz boosts but retains spatial rotations and spacetime scalings.} In the main part of the paper, we will be more economical, as we are only interested in wave equations in the region $t_*\geq 1$ (see Definitions~\ref{DefCMSpacetime} and \ref{DefCMDomain}), so, upon compactifying, in
\begin{equation}
\label{EqIOmega}
  \Omega = \ol{\{t_*\geq 1,\ r\geq\bhm\}} \subset \tilde M_2.
\end{equation}
See Figure~\ref{FigICpt}. 

\subsection{Elements of the proofs}
\label{SsIP}

We now describe the main ideas behind the proofs of Theorem~\ref{ThmIML} and its nonlinear applications.
\begin{enumerate}
\item For a substantial part of the analysis, the crucial notion of regularity---namely, the one for which we can establish appropriate microlocal propagation results for---is \emph{edge-3b-regularity} (or e3b-regularity for short), as introduced in \cite{HintzNonstat}; this is weaker than b-regularity. We motivate and describe this in~\S\ref{SssIPmlG}, and in the process also explain the precise form of the metric perturbations $h$ and operators $P$ allowed in Theorem~\ref{ThmIML}.
\item In~\S\ref{SssIP0}, we explain how to use estimates for the stationary model $P_0$ of $P$ to control the decay of solutions of $P u=f$; we also discuss how to prove these estimates in concrete settings.
\item We explain our method for proving higher b-regularity and b-tame estimates in~\S\ref{SssIPb}.
\item Getting stronger decay and asymptotics requires some post-processing of Theorem~\ref{ThmIML}, which we explain in~\S\ref{SssIPD}.
\item We illustrate our approach to the analysis of wave equations in a simple ODE toy model in~\S\ref{SssIPODE}. (The reader may wish to read this part first.)
\end{enumerate}

\subsubsection{Geometry, dynamics, and microlocal analysis in the e3b-setting}
\label{SssIPmlG}

A large part of our analysis is based on microlocal results on the propagation of regularity in phase space using real principal type propagation results \cite{DuistermaatHormanderFIO2} and uniform (infinite propagation time) versions thereof to deal, in particular, with horizons and trapping. Recall that the Hamiltonian flow associated with (the principal symbol of) a wave operator in the characteristic set---the \emph{null-bicharacteristic flow}---is the lift of the null-geodesic flow to phase space \cite[Proposition~9.9]{HintzMicro}. The first task for the purpose of proving high regularity estimates that are \emph{uniform} on the domain $\Omega^\circ=\{t_*\geq 1,\ r\geq\bhm\}$ is thus, simply put, to find coordinates on $T^*\R^4$ in the three asymptotic regimes of interest ($\scri^+$, $\iota^+$, and $\cK^+$ in Figure~\ref{FigICpt}) in which the null-bicharacteristic flow has a simple uniform description. This is most cleanly accomplished by devising a vector bundle over the compactified space $\Omega$ which over $\Omega^\circ$ is identified with $T^*\R^4$.

Near $\scri^+$ and $\iota^+$ (where thus $r=|x|$ is large), let us discuss the Minkowski metric
\[
  \ubar g = -\dd t^2 + \dd x^2 = -\dd t^2 + \dd r^2 + r^2\,\slg = -\dd t_*^2 - 2\,\dd t_*\,\dd r + r^2\slg,\quad t=t_*+r,
\]
instead of the Kerr metric $g_{\bhm,a}=\ubar g+\cO(r^{-1})$; here $\slg=\dd\theta^2+\sin^2\theta\,\dd\phi^2$ is the standard metric on $\Sph^2$. It is more convenient to work with the dual metric
\[
  \ubar g^{-1} = -\pa_t^2 + \pa_x^2 = -\pa_t^2 + \pa_r^2 + r^{-2}\slg^{-1} = -2\pa_{t_*}\otimes_s\pa_r + \pa_r^2 + r^{-2}\slg^{-1}.
\]

\pfstep{A first attempt: scattering phase space.} The fact that $\ubar g^{-1}$ is a Lorentzian signature quadratic form in the vector fields $\pa_t,\pa_{x^1},\pa_{x^2},\pa_{x^3}$ seems to suggest using the corresponding phase space, i.e., the one for which these vector fields give (via the dual pairing with covectors) linear momentum coordinates, for the analysis of associated wave-type operators; the corresponding notion of regularity is then given by regularity with respect to these vector fields (or equivalently \eqref{EqIMbTrans}) relative to $L^2$. Working on the base manifold $\ol{\R^4}$ this produces the \emph{scattering cotangent bundle}
\[
  \Tsc^*\ol{\R^4} = \ol{\R^4} \times \R^4,
\]
where a point $(z,\zeta)$ with $z=(t,x)\in\R^4\subset\ol{\R^4}$ is identified with the covector $\zeta_\mu\,\dd z^\mu\in T^*_z\R^4$.

Of course, one can replace the base $\ol{\R^4}$ by other compactifications of $\R^4$, such as $\tilde M_2$ defined in~\eqref{EqIM1tilde}; we recall that we must work on\footnote{The resulting vector bundle is then typically called the \emph{3-body scattering cotangent bundle}, following Vasy \cite{VasyThreeBody}.} $\tilde M_2$ in order to be able to work with radiative spacetimes and asymptotically stationary spacetimes (hence performing the blow-ups producing $\scri^+$ and $\cK^+$). Null-bicharacteristics in $\Tsc^*$ are then (reparameterizations of) (limits of families of) curves $s\mapsto (t+\tau s,x+\xi s,\tau,\xi)$; over the boundary of $\tilde M_2$, this includes curves over $\scri^+$ limiting to an endpoint over the future boundary $\scri^+\cap\iota^+$, curves limiting to a point over $\iota^+\cap\cK^+$, resp.\ $\iota^+\cap\scri^+$ in the past, resp.\ future (which capture waves scattered by a stationary metric or potential at $\cK^+$), or curves starting from the interior and limiting to a point over $\scri^+$ (which capture waves escaping to null infinity). In other words, the null-bicharacteristic flow features several invariant manifolds of source, sink, or saddle type over the boundary hypersurfaces or corners of $\tilde M_2$: this is a common feature for hyperbolic PDE on singular (e.g., compactified) spaces \cite{HintzVasySemilinear,BaskinVasyWunschRadMink,HintzVasyScrieb}.

A critical drawback of $\Tsc^*$ (over $\tilde M_2$, say), however, is that the null-bicharacteristic flow has several degeneracies. \emph{First}, in the zero section\footnote{Recall that in the scattering cotangent bundle over a compact manifold with boundary such as $\ol{\R^4}$, one can (micro)localize to finite frequencies $\zeta$ over base infinity, and similarly then over (the interior of) $\iota^+$ in $\tilde M_2$.} (where $\zeta=0$) over $\iota^+\subset\tilde M_2$ (or over $\pa\ol{\R^4}\subset\ol{\R^4}$), the flow is stationary: the Hamiltonian vector field $-2\tau\pa_t+2\xi\cdot\pa_x$ vanishes there, and thus solutions of wave equations cannot be controlled there by propagation from an initial Cauchy hypersurface $t_*=1$, say. (See \cite[Remark~1.12]{HintzNonstat} for a related discussion.) This is analogous to how the Laplace operator $\Delta$ on $\R_x^n$ degenerates as a scattering differential operator over the zero section of $\Tsc^*\ol{\R^n}$ over $\pa\ol{\R^n}$; and indeed $\pa_x$-regularity is not the correct notion for the study of $\Delta$, rather, $\la x\ra\pa_x$-regularity is, i.e., b-regularity on $\ol{\R^n}$. \emph{Second}, at the aforementioned sink for the future-directed flow over $\scri^+$, the flow turns out to have a quadratic degeneracy.

On the positive side, for the null-bicharacteristic flow for the Kerr metric, $\Tsc^*$ \emph{is} a convenient phase space over spatially compact regions, i.e., near $(\cK^+)^\circ\subset\tilde M_2$.\footnote{The locus of microlocal analysis there is fiber infinity of $\Tsc^*$ only: at finite frequencies, one needs to take non-microlocal, global, phenomena (related to spectral information of the stationary wave operator) into account, as done in \cite{VasyThreeBody} in scattering theory and \cite{BaskinDollGellRedmanKG} for the Klein--Gordon equation on asymptotically stationary spacetimes.} Indeed, there is a critical set given by the conormal bundle of the event horizon (over $t_*=\infty$) at which the flow has a nondegenerate saddle point structure (which is closely related to the subextremality of the horizon and the red-shift effect \cite{DafermosRodnianskiRedShift,GannotHorizons}). Moreover, there is an invariant set for the flow over $(\cK^+)^\circ$ which is the \emph{trapped set} (in phase space) $\Gamma$, a smooth submanifold of phase space (in the full subextremal range) at which the null-bicharacteristic flow is $\sfr$-normally hyperbolic for every $\sfr$. A complete description of the flow in this region is given in~\S\S\ref{SssTs3bH}--\ref{SssTs3bSum} by combining \cite[\S{6}]{VasyMicroKerrdS} and \cite{WunschZworskiNormHypResolvent,DyatlovWaveAsymptotics}.

\bigskip

\pfstep{Desingularization: edge-b-phase space near $\scri^+$, 3b-phase space near $\cK^+$.} The crucial observation that enables one to circumvent the degeneracies of the flow in $\Tsc^*$ is that the Minkowski metric (and the Kerr metric) is a (weighted) Lorentzian signature quadratic form also in other frames of vector fields, including ones that are uniformly inequivalent to the frame of coordinate vector fields.\footnote{All phase spaces are canonically identified with $T^*\R^4$ over $\R^4$, but the resulting bundle map from one phase space on $\tilde M_2$ to another may not be smooth up to the boundary of $\tilde M_2$, so the flow over the boundary of $\tilde M_2$ can have different properties, e.g., orders of degeneracy at invariant sets.} Wave-type operators on these spacetimes can, in principle, be analyzed microlocally in each of these phase spaces. The challenge is to identify a phase space on which the flow is non-degenerate in an appropriate sense near each of the various boundary hypersurfaces of the compactified spacetime manifold. This was accomplished on asymptotically flat and dynamical spacetimes by the author in \cite{HintzNonstat} by combining insights from \cite{HintzVasyScrieb,Hintz3b}; we recall this now. (The only features not present in \cite{HintzNonstat} are horizons and trapping.)

A suitable phase space near $\cK^+$ was identified in \cite{Hintz3b} to be the \emph{3b-cotangent bundle} (\S\ref{SssCT3b}): its (uniform down to the boundary at infinity of $\tilde M_2$ near $\cK^+$) momentum coordinates are defined with respect to the vector fields
\begin{equation}
\label{EqIPmlG3b}
  r\pa_{t_*},\ \ r\pa_r,\ \ \pa_\omega.
\end{equation}
Note that in regions $0<\delta<\frac{r}{t}<\delta^{-1}$, these span the space of b-vector fields on $\ol{\R^4}$. The idea of using b-notions to study massless wave equations on asymptotically flat spacetimes goes back to \cite{VasyMicroKerrdS,BaskinVasyWunschRadMink}, and is the Lorentzian analogue of using b-notions to study the ``massless'' Laplace equation on Euclidean space. The Minkowski dual metric is schematically given by
\begin{equation}
\label{EqIPmlG3bMet}
  \ubar g^{-1} = \rho_+^2 \bigl( -2 r\pa_r\otimes_s r\pa_{t_*} + (r\pa_r)^2 + \pa_\omega^2 \bigr),
\end{equation}
where we write $\rho_+=r^{-1}$ for a local defining function\footnote{In the region $\{t_*\geq 1,\ \frac{r}{t}<\frac23\}$, this function is a bounded smooth positive multiple of $\rho_+$ in~\eqref{EqIbdfs}.} of $\iota^+$ near $\cK^+$. This is thus a weighted Lorentzian signature quadratic form in the vector fields~\eqref{EqIPmlG3b}. Microlocal propagation results in 3b-phase space control waves at infinite 3b-frequencies, i.e., in the sense of 3b-regularity, which is regularity with respect to the vector fields~\eqref{EqIPmlG3b}. Note that in spatially compact regions and for $t_*\geq 1$, this is indistinguishable from the scattering setting, while near $\iota^+$ (i.e., for large $r$), this concerns a qualitatively different notion of regularity.

In order to study massless waves in a neighborhood of $\scri^+$, the author and Vasy \cite{HintzVasyScrieb} introduced the \emph{edge-b-cotangent bundle} (``eb'' for short) (\S\ref{SssCTeb}), which combines Mazzeo's edge cotangent bundle \cite{MazzeoEdge} with Melrose's b-cotangent bundle \cite{MelroseTransformation,MelroseAPS}. (In the massive case, Sussman used the closely related double edge-scattering-cotangent bundle \cite{SussmanKG} for microlocalization.) In terms of the local defining functions\footnote{These are bounded smooth positive multiples of the functions~\eqref{EqIbdfs} in the region $\{t_*\geq 1,\ \frac{r}{t}>\frac13\}$.} $\rho_\sscri=\frac{t_*}{r}$ and $\rho_+=\frac{1}{t_*}$ of $\scri^+$ and $\iota^+$, respectively, we set $x_\sscri:=\sqrt{\rho_\sscri}$ and use the vector fields
\begin{equation}
\label{EqIPmlGeb}
  x_\sscri\pa_{x_\sscri}=-2 r\pa_r,\quad \rho_+\pa_{\rho_+}=-(t_*\pa_{t_*}+r\pa_r),\quad x_\sscri\pa_\omega=\sqrt{\frac{t_*}{r}}\pa_\omega
\end{equation}
as a frame; one computes that (schematically, and to leading order in $x_\sscri$)
\begin{equation}
\label{EqIPmlGebMet}
  \ubar g^{-1} = \frac12 \rho_\sscri\rho_+^2 \bigl( x_\sscri\pa_{x_\sscri} \otimes (x_\sscri\pa_{x_\sscri}-2\rho_+\pa_{\rho_+}) + (x_\sscri\pa_\omega)^2 \bigr).
\end{equation}
(See~\eqref{EqCTebMink} for the full expression.) Up to the overall weight $\rho_\sscri\rho_+^2$, this is thus a Lorentzian signature quadratic form in the vector fields~\eqref{EqIPmlGeb}

\bigskip

\begin{figure}[!ht]
\centering
\includegraphics{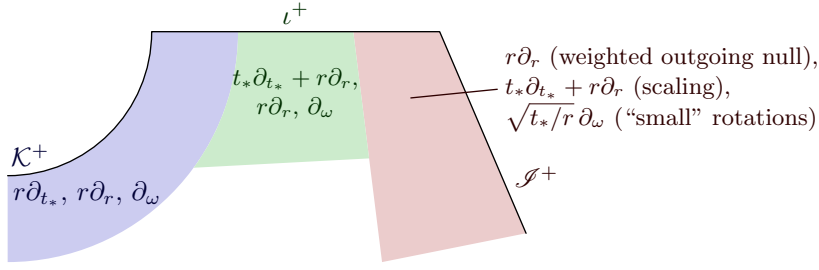}
\caption{The basic vector fields underlying the notion of e3b-regularity, taken from~\eqref{EqIPmlG3b}, \eqref{EqIPmlGeb}, and b-vector fields near $(\iota^+)^\circ$. By contrast, b-regularity amounts to regularity with respect to $t_*\pa_{t_*}$, $r\pa_r$, $\pa_\omega$.}
\label{FigIPmle3bReg}
\end{figure}

\begin{figure}[!ht]
\centering
\includegraphics{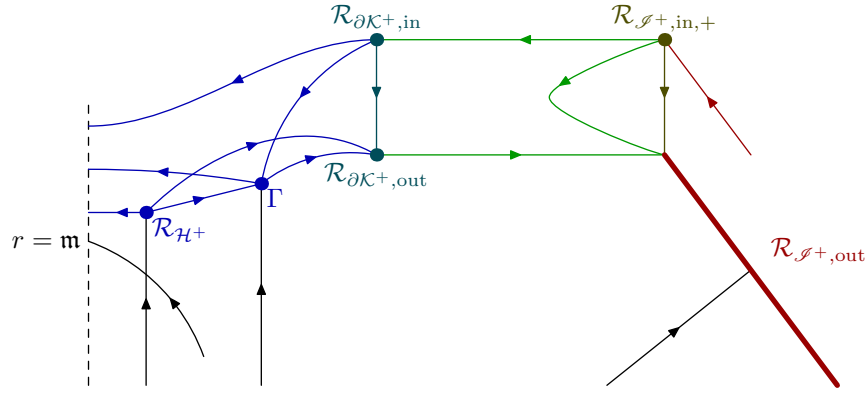}
\caption{Illustration of the future-directed null-bicharacteristic flow for a subextremal Kerr metric in e3b-phase space over the (compactified) domain $\Omega$. The radial sets $\cR_{\scri^+,{\rm out}}$ and $\cR_{\scri^+,{\rm in},+}$ are defined in Definition~\ref{DefTsebRad}, $\cR_{\cK^+,{\rm in/out}}$ in Definition~\ref{DefTs3bRad}, $\cR_{\cH^+}$ in Definition~\ref{DefTs3bH}, and the trapped set $\Gamma$ in~\eqref{EqTs3bGamma} (based on Definition~\ref{DefTs3bOTrap0}).}
\label{FigIPmle3bFlow}
\end{figure}

We combine both phase spaces to the e3b-phase space ${}^\etbop T^*\to\Omega$, which is thus the eb-phase space near $\scri^+$ and the 3b-phase space near $\cK^+$, with both matching up where they overlap over $(\iota^+)^\circ$. The key point is then that the null-geodesic flow lifted to ${}^\etbop T^*$ features only \emph{non-degenerate} saddle points (and trapping, as already discussed); this is verified by explicit computations in~\S\ref{STs}. We summarize the basic e3b-vector fields in Figure~\ref{FigIPmle3bReg}, and the structure of the null-bicharacteristic flow in ${}^\etbop T^*$ over $\Omega$ in Figure~\ref{FigIPmle3bFlow}. (A more schematic version which also shows the initial and final boundary hypersurfaces of $\Omega$ is shown in Figure~\ref{FigTse3bDyn}.)

\begin{rmk}[No direct microlocal control of b-regularity]
\label{RmkIPmlG}
  The Minkowski and Kerr dual metrics are \emph{not} (weighted) non-degenerate Lorentzian signature quadratic forms in the b-vector fields $t_*\pa_{t_*}$, $\la x\ra\pa_x$ (or $r\pa_x$ in $r\gtrsim 1$). Therefore, the corresponding wave operators are highly degenerate when studied in the corresponding b-cotangent bundle, and hence standard microlocal regularity results cannot yield the b-regularity of waves as asserted in Theorem~\ref{ThmIML}. The proof of b-regularity requires, instead, additional commutation-type arguments (\S\ref{SssIPb}).
\end{rmk}

\begin{rmk}[Complexity]
\label{RmkIPmlComplex}
  Over every codimension $2$ corner in a compactification of spacetime, given by the intersection of two adjacent boundary hypersurfaces $H_1$ and $H_2$, one expects at least $2$ radial sets, corresponding to (lifts of) limits of families of null-geodesics that pass from $H_1$ to $H_2$ or vice versa. Thus, even on spacetimes that appear geometrically simple, such as mild stationary perturbations of the Minkowski metric that also feature mild radiation near null infinity, the setup of our phase space analysis necessitates propagation through two times $2$ (for the two corners $\iota^+\cap\scri^+$ and $\iota^+\cap\cK^+$) critical sets (saddles or sources/sinks); see~\S\ref{SssR3R}. The Kerr metric introduces one additional radial set (over the event horizon) and normally hyperbolic trapping (\S\ref{SssR3Tr}). We must, therefore, describe microlocal propagation results near \emph{five} qualitatively different radial sets as well as near the trapped set. Furthermore, in the analysis of various model operators (such as $t_*$-translation-invariant wave operators $P_0$ in~\S\ref{SSp}, or dilation-invariant models at $\iota^+$ in~\S\ref{SsNip}), images of these radial and trapped sets on spacetime appear yet again, necessitating proofs of microlocal propagation estimates also at these ``projected'' versions. The latter are conceptually very similar; but on occasion, additional technicalities arise, e.g., when proving estimates for the spectral family at complex frequencies $\sigma$ with $\Im\sigma>0$ (e.g., Propositions~\ref{PropSpHiPr2}--\ref{PropSpHiPr3}) which correspond to spacetime estimates on exponentially growing spaces---which the spacetime estimates on polynomially weighted e3b-Sobolev spaces do not capture. It is this geometric complexity of the phase space dynamics that is partly responsible for the formidable length of the present paper. --- We note that, by contrast, in the analysis of (non-)linear waves in a neighborhood of the domain of outer communications of asymptotically Kerr--\emph{de~Sitter} spacetimes, there is only \emph{one} asymptotic regime (encoded by the boundary hypersurface $\cong I\times\Sph^2$ of the compactification at $e^{-t_*}=0$, where $I\subset\R$ is a compact radial interval), and only \emph{one} class of radial sets (subextremal horizons). Moreover, the b-perspective suffices for all geometric and analytical purposes there, as the relevant wave operators \emph{are} non-degenerate b-differential operators.
\end{rmk}

\pfstep{e3b-metrics and operators.} For the purpose of analysis in e3b-phase space, what matters most is the principal symbol of wave-type operators as weighted e3b-differential operators, and thus the (dual) metric as a weighted e3b-metric (i.e., expressed in terms of basic e3b-vector fields as in~\eqref{EqIPmlG3bMet} and \eqref{EqIPmlGebMet}). Thus, the precise requirement for the metrics $g=g_{\bhm,a}+h$ in Theorem~\ref{ThmIML} is that
\begin{equation}
\label{EqIPmlPert}
  \text{\parbox{0.8\textwidth}{the coefficients of $h$ when expressed in terms of the dual bases $\tfrac{\dd t_*}{r}$, $\tfrac{\dd r}{r}$, $\dd\omega$, resp.\ $\tfrac{\dd x_\sscri}{x_\sscri}$, $\tfrac{\dd\rho_+}{\rho_+}$, $\tfrac{\dd\omega}{x_\sscri}$ of~\eqref{EqIPmlG3b} (near $\cK^+$), resp.\ \eqref{EqIPmlGeb} (near $\scri^+$), are of size $\cO(\rho_\sscri\rho_+^2\cdot\rho_\sscri^{\ell_\sscri}\rho_+^{\ell_+}\rho_\cK^{\ell_\cK})$ where $\ell_\sscri>0$, $\ell_+>\delta$, and $\ell_\cK>\aleph$, together with all of their b-derivatives.}}
\end{equation}

See also Lemma~\ref{LemmaCTebsc}, and \cite[Lemma~3.23]{HintzMink4Gauge} for the eb-structure of metric perturbations in the exterior stability problem for asymptotically flat spacetimes. Away from $\scri^+$, this is equivalent to $\cO(\rho_+^{\ell_+}\rho_\cK^{\ell_\cK})$-decay of the coefficients of $h$ in the standard basis $\dd t_*,\dd x$ as stated in condition~\eqref{ItIML2} of Theorem~\ref{ThmIML}.

The wave-type operators $P$ relative to such metrics $g$ are correspondingly required to be e3b-differential operators (i.e., built from compositions of e3b-vector fields) with weight $\rho_\sscri\rho_+^2$. The principal symbol being determined by the (dual) metric, this is a non-vacuous requirement only for lower-order terms. As an illustration, consider the scalar wave operator $\Box_{\ubar g}$ on Minkowski space: near $\cK^+$, this is given by $\Box_{\ubar g}=r^{-2}((r\pa_t)^2-(r\pa_r)^2-r\pa_r-\pa_\omega^2)=r^{-2}(-2 r\pa_{t_*}\,r\pa_r-(r\pa_r)^2-r\pa_r-\pa_\omega^2)$ (expressed using~\eqref{EqIPmlG3b}), while near $\scri^+$, we have
\begin{equation}
\label{EqIPmlBoxeb}
  \Box_{\ubar g}\equiv\tfrac12\rho_\sscri\rho_+^2\bigl( -(x_\sscri\pa_{x_\sscri}-2)(x_\sscri\pa_{x_\sscri}-2\rho_+\pa_{\rho_+}) - (x_\sscri\pa_\omega)^2 \bigr)
\end{equation}
to leading order at $x_\sscri=0$ (expressed using~\eqref{EqIPmlGeb}). (The operator $\Box_{g_{\bhm,a}}$ is the same to leading order at $\rho_+=0$ and $x_\sscri=0$.) The term $x_\sscri\pa_{x_\sscri}-2$ is responsible for the $x_\sscri^2$ ($=\cO(r^{-1})$ when $|t_*|=\cO(1)$) radiation field of solutions of the wave equation. We permit the operator $P$ to include, instead, a factor $x_\sscri\pa_{x_\sscri}-2(1+p_1)$ for some $p_1=p_1(t_*,\omega)$ (a bundle map for tensorial operators), leading to potentially different decay rates. The full conditions are stated in Definition~\ref{DefSDWAdm}.

\bigskip
\pfstep{e3b-microlocal analysis.} Given a wave-type operator $P$ of the type described above, we can now analyze it using basic microlocal regularity results (microlocal elliptic regularity, real principal type propagation) as well as propagation estimates at radial sets and at the trapped set (for which we can directly quote \cite{HintzPolyTrap}), \emph{provided} we have an algebra of pseudodifferential operators for the purpose of localizing in (fiber-wise conic) subsets of phase space ${}^\etbop T^*\to\Omega$. This algebra can be defined most simply as the algebra of \emph{uniform} (or \emph{bounded geometry}) \emph{pseudodifferential operators} on $\Omega^\circ$ (ignoring the finite boundary hypersurfaces $t_*^{-1}(1)$ and $r^{-1}(\bhm)$) described by Shubin \cite{ShubinBounded}: a smooth Riemannian e3b-metric for which the basic e3b-vector fields~\eqref{EqIPmlG3b} and \eqref{EqIPmlGeb} have unit length has bounded geometry, and ps.d.o.s can then be defined as standard quantizations on $\R^4$ in unit size balls (``distinguished charts'') with respect to such a metric. A more elementary framework was introduced by the author in \cite{HintzScaledBddGeo}, in which the distinguished charts are the primary object, with uniformly bounded smooth vector fields in the charts being precisely the e3b-vector fields (with e3b-regular coefficients); we recall this in~\S\ref{SsMS}, and describe e3b-ps.d.o.s in~\S\ref{SssMUe3b}. With e3b-ps.d.o.s at hand, standard positive commutator arguments as in \cite[Chapters~8 and 10]{HintzMicro} then yield the desired microlocal real principal type and radial propagation results in e3b-phase space (see~\S\S\ref{SsMBasic} and \ref{SsR3}).

Only near the finite boundary hypersurfaces of $\Omega$ do microlocal arguments lose their potency, as they do not respect boundary conditions and support properties. Instead, near $t_*^{-1}(1)$ and $r^{-1}(\bhm)$, we use standard energy estimates to control waves for $t_*$ near $1$ (via $t_*$-propagation) and $r$ near $\bhm$ (via $r$-propagation), respectively. The former energy estimate is somewhat delicate since $t_*$-level sets extend all the way to null infinity; but energy estimates in the required weighted e3b-Sobolev spaces can be proved using the methods introduced in \cite{HintzVasyMink4,HintzVasyScrieb}. The details are given in~\S\S\ref{SsET}--\ref{SsER}.

At this stage, one can already prove an a priori estimate for forward solutions of $P u=f$ (in the setting of Theorem~\ref{ThmIML}) that controls $u$ in the sense of e3b-regularity, i.e., at high e3b-frequencies; this reads
\begin{equation}
\label{EqIPmlEst}
  \|u\|_{H_\etbop^{s,(\beta_\sscri,\beta_+,\beta_\cK)}(\Omega,\mu_\bop)} := \|\rho_\sscri^{-\beta_\sscri}\rho_+^{-\beta_+}\rho_\cK^{-\beta_\cK}u\|_{H_\etbop^s} \lesssim \|P u\|_{H_\etbop^{s,(\beta_\sscri+1,\beta_++2,\beta_\cK)}} + \|u\|_{H_\etbop^{s-10,(\beta_\sscri,\beta_+,\beta_\cK)}}.
\end{equation}
Here, $H_\etbop^s$ denotes the spacetime $L^2$-space with $s$ additional orders of regularity with respect to e3b-vector fields. Furthermore,
\begin{enumerate}
\item the $\scri^+$-decay rate $\beta_\sscri$ must be below the threshold given by the decay rate towards $\scri^+$ of the radiation field term for solutions of $P$;
\item the $\iota^+$-decay rate $\beta_+$ must be less than $\beta_\sscri$, as already discussed in Remark~\ref{RmkI1QDecay};
\item the regularity order $s$ and the relative decay rate $\beta_\cK-\beta_+$ must satisfy appropriate upper and lower threshold conditions; these arise in the two radial point estimates over $\cK^+\cap\iota^+$ and are closely related to similar conditions in stationary scattering theory for the absence of incoming and the allowance of outgoing radiation (cf.\ \cite{MelroseEuclideanSpectralTheory} and \cite[Proposition~5.28]{VasyMinicourse}).
\end{enumerate}
See\footnote{In the bulk of the paper, we use powers of $x_\sscri=\rho_\sscri^{\frac12}$ as $\scri^+$-weights, and we work with $L^2$ relative to the metric density rather than $\mu_\bop$, which leads to constant shifts in the orders.} Proposition~\ref{PropR3} for the precise estimate. (A technical caveat, again familiar from stationary scattering theory, is that the two threshold conditions on $s$ cannot be simultaneously satisfied by a single real number, and hence we work with variable regularity orders.)

We stress that
\begin{equation}
\label{EqIPmlInsensitive}
  \text{\parbox{0.8\textwidth}{\centering \textit{the estimate~\eqref{EqIPmlEst} is agnostic to spectral properties of the stationary model $P_0$,}}}
\end{equation}
such as mode stability properties and the zero energy behavior. But of course the failure of, say, mode stability would imply that the final norm on the right would, typically, be infinite for forward solutions of $P$ or $P_0$. The next step is thus to control the \emph{decay} of $u$ using spectral information.

\subsubsection{The stationary model and control of decay}
\label{SssIP0}

Consider the estimate~\eqref{EqIPmlEst} and recall from assumption~\eqref{ItIML3} in Theorem~\ref{ThmIML} that the wave-type operator $P$ (on a dynamical Kerr spacetime) is equal to the stationary model $P_0$ (on an exact subextremal Kerr spacetime) to leading order as $t_*\to\infty$. Unconditional estimates for forward solutions of $P_0$ should therefore be useful in controlling the decay of forward solutions of $P$. We proceed to explain this in more detail.

\pfstep{Improving decay, I: ignoring null infinity. Connection to spectral theory.} Conceptually speaking, $t_*$-translation-invariance should mainly be used only near $\cK^+$. (Orbits of the $t_*$-translation action on $\R^4=\R_{t_*}\times\R^3_x$ are parameterized by points in $(\cK^+)^\circ\cong\R^3_x$.) Localized to a neighborhood of $\cK^+$, the norms in~\eqref{EqIPmlEst} are equal to norms of weighted 3b-Sobolev spaces
\[
  \Htb^{s,(\beta_+,\beta_\cK)}(\Omega,\mu_\bop)=\rho_\sface^{\beta_+}\rho_\cK^{\beta_\cK}\Htb^s(\Omega,\mu_\bop),\quad \rho_\sface:=r^{-1},\ \rho_\cK:=\frac{r}{t_*+r},
\]
in which regularity is measured with respect to the 3b-vector fields $r\pa_{t_*}$, $r\pa_r$, $\pa_\omega$ from~\eqref{EqIPmlG3b}. (More precisely, one should consider $\Omega$ here as the closure of $\{t_*\geq 1\}$ inside of $[\ol{\R^4};\fk^\pm]$, which is like $\tilde M_2$ in Figure~\ref{FigICpt} except $\scri^+$ (and $\scri^-$) are not blown up; we write ``$\sface$'' for the boundary $r^{-1}=0$.) The key assumption on $P_0$ is \emph{$\aleph$-admissibility with $\sface$-weight $\beta_\sface$ and $\sface$-loss $\delta$}, which roughly speaking means:
\begin{equation}
\label{EqIP0aleph}
  \text{\parbox{0.8\textwidth}{\it For $f$ with support in $t_*\geq 1$, the forward solution $u$ of $P_0 u=f$ satisfies the quantitative estimate
    \[
      \|u\|_{\Htb^{s,(\beta_\sface-\delta,-\aleph+\frac12)}(\Omega,\mu_\bop)} \leq C\|f\|_{\Htb^{s,(\beta_\sface+2,\frac12)}(\Omega,\mu_\bop)}.
    \]}}
\end{equation}
The full requirements are stated in Definition~\ref{DefSSAlephAdm} (using the metric volume density and thus involving a shift in weights) and in particular require also a higher b-regularity estimate, discussed in~\S\ref{SssIPb} below. That we can allow for a loss $\aleph>0$ here is a crucial novelty (with considerable repercussions also for the proof of b-regularity) of the present paper over \cite{HintzNonstat}.

Before discussing~\eqref{EqIP0aleph} (and how to prove it in concrete cases), let us indicate how it can be used to improve the estimate~\eqref{EqIPmlEst}, still ignoring the $\scri^+$-orders and the edge-nature of the notion of e3b-regularity there: namely, one takes $\beta_+=\beta_\sface-\delta$ and $\beta_\cK=-\aleph+\frac12$ and estimates the final term in~\eqref{EqIPmlEst} by
\[
  \|P_0 u\|_{\Htb^{s-10,(\beta_++2+\delta,\frac12)}(\Omega,\mu_\bop)} \leq \|P u\|_{\Htb^{s-10,(\beta_++2+\delta,\frac12)}} + \|(P-P_0)u\|_{\Htb^{s-10,(\beta_++2+\delta,\frac12)}}.
\]
But since the coefficients of $P-P_0$ have $\ell_+>2+\delta$ and $\ell_\cK>\aleph$ orders of decay (as 3b-differential operators) at $\iota^+$ and $\cK^+$, respectively (see the assumptions stated in Theorem~\ref{ThmIML}\eqref{ItIML2}, \eqref{ItIML3}), the error term $(P-P_0)u$ is bounded by $\|u\|_{\Htb^{s-8,(\beta_+-\eta,-\aleph+\frac12-\eta)}}$, which on $\Omega$ is thus (again ignoring $\scri^+$) weaker than the norm with weights $(\beta_+,-\aleph+\frac12)$ by a weight $t_*^{-\eta}$; so one gets
\[
  \|u\|_{\Htb^{s,(\beta_+,-\aleph+\frac12)}} \lesssim \|P u\|_{\Htb^{s,(\beta_++2+\delta,\frac12)}} + \|t_*^{-\eta}u\|_{\Htb^{s-8,(\beta_+,-\aleph+\frac12)}}.
\]
Since one can control $u$ in terms of $P u$ for any fixed finite $t_*$-interval by the (energy estimate-based) finite-time theory, one can localize to very late times. But then the error term is \emph{small} and can thus be absorbed into the left-hand side, yielding the (at this point still a priori) estimate
\[
  \|u\|_{\Htb^{s,(\beta_+,-\aleph+\frac12)}} \lesssim \|P u\|_{\Htb^{s,(\beta_++2+\delta,\frac12)}}.
\]

\begin{figure}[!ht]
\centering
\includegraphics{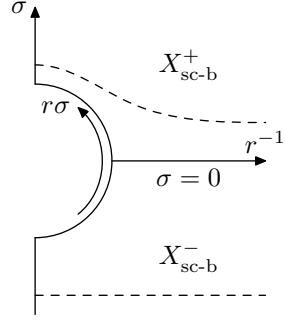}
\caption{The manifold for low-energy resolvent analysis, for which we must work with function spaces that capture different rates of spatial decay at nonzero energies $\sigma$ (the vertical boundary) and at $\sigma=0$ (the front face, shown curved here). Two level sets of $\sigma$ are shown as dashed lines.}
\label{FigIP0scbt}
\end{figure}

\bigskip
\pfstep{Verifying $\aleph$-admissibility.} In the applications here and in \cite{HintzKerrStab}, the estimate~\eqref{EqIP0aleph} for forward solutions of the stationary problem $P_0 u=f$ is proved as a consequence of estimates for the resolvent, i.e., the inverse of the spectral family $\wh{P_0}(\sigma)$ of $P_0$. The starting point is the Fourier transform in $t_*$; passing from the density $\mu_\bop$ to the metric density $\mu=|\dd t_*\,\dd x|$ on Minkowski space (or Kerr) leads to a shift of orders,
\[
  \|f\|_{\Htb^{s,(\beta_\sface+2,\frac12)}(\Omega,\mu_\bop)} = \|f\|_{\Htb^{s,(\alpha_\sface+2,0)}(\Omega,\mu)} = \|r^{\alpha_\sface+2}f\|_{\Htb^s(\Omega,\mu)},\quad \alpha_\sface:=\beta_\sface-\tfrac32.
\]
For $s=0$, this norm can be expressed in terms of the Fourier transform $\cF$, with the convention $\cF f(\sigma,x):=\hat f(\sigma,x)=\int_\R e^{i\sigma t_*}f(t_*,x)\,\dd t_*$, using Plancherel's theorem as
\[
  \|\hat f\|_{L^2\bigl(\R_\sigma; r^{-(\alpha_\sface+2)}L^2(\R^3,|\dd x|)\bigr)}.
\]
(This is where usage of the particular weight $\frac12$ at $\cK^+$ of $f$ in~\eqref{EqIP0aleph} is key.) Since $\cF$ commutes with $r\pa_r$ and $\pa_\omega$, and $\cF\circ r\pa_{t_*}=-i\sigma r\cF$, the space $\Htb^{s,(\alpha_\sface+2,0)}(\Omega,\mu)$ admits a similar characterization, where the $r^{-(\alpha_\sface+2)}L^2(\R^3)$-norm is replaced by the $\sigma$-dependent norm capturing $s$ degrees of regularity with respect to $\sigma r$ (as a multiplication operator) and $r\pa_r$, $\pa_\omega$, or equivalently $r\pa_x$ (for $r\gtrsim 1$). For $\sigma=0$, this is thus a weighted b-Sobolev norm, for $\sigma\neq 0$ a standard norm on a weighted Sobolev space on $\R^3$ (in the geometric singular analysis called \emph{(weighted) scattering Sobolev space} \cite{MelroseEuclideanSpectralTheory}), for unbounded $\sigma$ a semiclassical scattering Sobolev norm (with semiclassical parameter $|\sigma|^{-1}$) \cite{VasyZworskiScl}, and for $|\sigma|\leq 1$ a weighted \emph{scattering-b-transition Sobolev norm} \cite{HintzKdSMS}; this is recalled in Lemma~\ref{LemmaMUetbFT} (following \cite[\S{4.4}]{Hintz3b}). Importantly, on these function spaces, one can prove precise estimates for the spectral family on asymptotically flat spaces, as we recall in detail in~\S\ref{SSp} (see Theorems~\ref{ThmSpB}, \ref{ThmSpHi}, \ref{ThmSp0}, and \ref{ThmSpLo}). As two instances of this, we note here only that the variable nature of the order $s$ descends to the variable decay rates (for incoming vs.\ outgoing spherical waves) in stationary scattering theory at frequencies $\sigma\in\R\setminus\{0\}$; and the range of weights $\alpha_\sface$ for $\aleph$-admissibility is determined by the indicial roots of the zero energy operator $\wh{P_0}(0)$. Furthermore, we point out for later purposes that the manifold with corners on which the low-energy resolvent analysis takes place is a blow-up of $(-1,1)_\sigma\times\ol{\R^3}$ on which the function $\sigma r$ mentioned above is resolved to a projective coordinate on the front face arising from blowing up $\{\sigma=0\}\times\pa\ol{\R^3}$; see Figure~\ref{FigIP0scbt}, and \S\S\ref{SssMUK}, \ref{SsSptf}, and \ref{SsSpLo} for details.

In~\S\ref{SSp}, we prove the Fredholm index $0$ property of $\wh{P_0}(\sigma)$ acting on the appropriate function spaces, as well as precise regularity and asymptotics for elements of $\ker\wh{P_0}(\sigma)$, which are often used in the assumptions of mode stability results such as \cite{WhitingKerrModeStability,ShlapentokhRothmanModeStability,AnderssonMaPaganiniWhitingModeStab,TeixeiradCModes} for the scalar wave equation or the Teukolsky equation on Kerr spacetimes. If $\wh{P_0}(\sigma)$ is invertible for all $\sigma\in\R$, then the function $u$ defined by $\hat u(\sigma)=\wh{P_0}(\sigma)^{-1}\hat f(\sigma)$ satisfies~\eqref{EqIP0aleph}; and if $\wh{P_0}(\sigma)^{-1}$ is invertible in the full half-space $\Im\sigma\geq 0$, then a Paley--Wiener argument implies that this $u$ is, indeed, the \emph{forward} solution of $P_0 u=f$, showing that
\[
  \text{\parbox{0.8\textwidth}{\centering\it $P_0$ is $0$-admissible if it satisfies mode stability in $\Im\sigma\geq 0$.}}
\]
See Theorem~\ref{ThmA1Gen} for the precise result, and Theorem~\ref{ThmA1Adm} for the case of the scalar wave operator on subextremal Kerr spacetimes.

When mode stability fails at $\sigma=0$, as is the case for the 1-form wave operator by \cite[Theorem~5.1]{AnderssonHaefnerWhitingMode}, $\wh{P_0}(\sigma)^{-1}\hat f(\sigma)$ no longer inherits the square-integrability in $\sigma$ (with values in the appropriate spatial Sobolev space) from $\hat f(\sigma)$. We instead  show for the 1-form wave operator on subextremal Kerr spacetimes that one can solve
\begin{equation}
\label{EqIuregsing}
  \wh{P_0}(\sigma) \bigl(\hat u_{\rm reg}(\sigma) + \hat a(\sigma)\sigma^{-1}u_0\bigr) = \hat f(\sigma)
\end{equation}
where $|\hat a(\sigma)|$ and the norm of $\hat u_{\rm reg}(\sigma)$ obey \emph{uniform} bounds in terms of the norm of $\hat f(\sigma)$; here $u_0\in\ker\wh{P_0}(0)$ is the zero energy state (see~\eqref{EqIu0}). We achieve this using a Grushin problem \cite{GrushinICM,SjostrandZworskiLinAlg}, concretely, by using uniform low-energy estimates for the augmented operator
\begin{equation}
\label{EqIP0Grushin}
  \wt{P_0}(\sigma) := \begin{pmatrix} \wh{P_0}(\sigma) & \wh{P_0}(\sigma)(\sigma^{-1}u_0) \\ \lambda(\cdot) & 0 \end{pmatrix}
\end{equation}
where $\lambda$ is a linear functional that does not annihilate $u_0$. Note that $\wh{P_0}(\sigma)(\sigma^{-1}u_0)=\pa_\sigma\wh{P_0}(0)u_0+\cO(|\sigma|)$ is bounded at $\sigma=0$. The point is that $\wt{P_0}(\sigma)$ \emph{is} invertible at (and thus uniformly near) $\sigma=0$; see Lemma~\ref{LemmaA2Adm0Op}. The lack of decay of the forward solution $u$ of $P_0 u=f$ at $\cK^+$, relative to the decay of $f$, arises solely from the singularity of the term $\hat a(\sigma)\sigma^{-1}u_0$ at $\sigma=0$; on the inverse Fourier transform side, this is essentially the $t_*$-integral of $(\cF^{-1}a)u_0$, with $\cF^{-1}a\in L^2(\R_{t_*})$ having the same $t_*$-integrability as $f$. The fact that this $t_*$-integration increases the $t_*$-weight by \textit{1} results in the \textit{1}-admissibility of the 1-form wave operator. See Theorem~\ref{ThmA2Adm} for details.

\bigskip
\pfstep{Improving decay, II: dealing with null infinity.} General edge-3b-Sobolev spaces do not have a Plancherel-type description on the spectral side. They do, however, interact well with the \emph{Mellin transform} with respect to spacetime scalings $(t_*,x)\mapsto(\lambda t_*,\lambda x)$; for $\frac{t_*}{r}\in(0,\infty)$, the interior of $\iota^+$ is naturally identified with a transversal for this dilation action in $t_*>0$. Indeed, one can check that near $\iota^+\cap\cK^+$, the space of e3b-vector fields is spanned by $R\tau\pa_\tau$, $R\pa_R$, $\pa_\omega$ where $R=\frac{r}{t_*}$ and $\tau=\frac{1}{t_*}$, and near $\iota^+\cap\scri^+$ by $\tau\pa_\tau$, $x_\sscri\pa_{x_\sscri}$, and $x_\sscri\pa_\omega$ where $x_\sscri=\sqrt{t_*/r}$; and the Mellin transform in $\tau$ intertwines $\tau\pa_\tau$ with multiplication by a spectral parameter. The (parameter-dependent) function spaces on the spectral side are then b- and semiclassical cone \cite{HintzConicPowers,HintzConicProp} Sobolev spaces near $\iota^+\cap\cK^+$, and 0- \cite{MazzeoMelroseHyp} and semiclassical 0-Sobolev spaces near $\iota^+\cap\scri^+$. On the level of operators, it is the dilation-homogeneous model for $P_0$ at $\iota^+$ that descends to the Mellin-transform side. In the case $P_0=\Box_{g_{\bhm,a}}$, this model is the corresponding (scalar or 1-form) wave operator on Minkowski space; the general formulas are given in Definition~\ref{DefNipOp}.

In~\cite{HintzNonstat}, estimates for this dilation-invariant model---a sharpening of which is proved here in Proposition~\ref{PropNip}---were plugged into the error term of~\eqref{EqIPmlEst} in order to control $u$ to leading order at $\iota^+$ (in the correct e3b-Sobolev space). In the present paper, however, the loss of $\cK^+$-decay due to the possibility of zero energy bound states makes these estimates (which require the $\iota^+\cap\cK^+$-weight to lie in a certain range dictated by the indicial roots of $\wh{P_0}(0)$) difficult to apply directly. We do use them, however, to find an exact local solution of $P_0 u=f$ near $\scri^+$ on e3b-spaces, and mix-and-match such an e3b-solvability theory with the 3b-estimates from the $\aleph$-admissibility assumption to upgrade the estimate~\eqref{EqIP0aleph} to an estimate
\begin{equation}
\label{EqIP0alephe3b}
  \|u\|_{H_\etbop^{s,(\beta_\sscri,\beta_\sface-\delta,-\aleph+\frac12)}(\Omega,\mu_\bop)} \leq C\|P_0 u\|_{H_\etbop^{s,(\beta_\sscri+1,\beta_\sface+2,\frac12)}(\Omega,\mu_\bop)}.
\end{equation}
on e3b-Sobolev spaces; this is done in~\S\ref{SN} (see Theorem~\ref{ThmNFw}).

\bigskip
\pfstep{The final a priori estimate, and how to make it unconditional.} Plugging~\eqref{EqIP0alephe3b} (with regularity order $s-10$) into~\eqref{EqIPmlEst} and using the strong decay of $P-P_0$ as before now yields the a priori estimate
\begin{equation}
\label{EqIP0Est}
  \|u\|_{H_\etbop^{s,(\beta_\sscri,\beta_+,-\aleph+\frac12)}} \lesssim \|f\|_{H_\etbop^{s,(\beta_\sscri+1,\beta_++2+\delta,\frac12)}},\quad f=P u.
\end{equation}
To turn this into an \emph{unconditional} estimate for the forward solution of $P u=f$, we combine two observations.
\begin{enumerate}
\item Forward solutions for $P_0$ \emph{do} satisfy the estimate~\eqref{EqIP0Est}, as just discussed.
\item The estimate~\eqref{EqIP0Est} holds \emph{uniformly} for perturbations of $P$.
\end{enumerate}
To wit, we use~\eqref{EqIP0Est} for operators $P^\eps$ which transition from $P$ to $P_0$ at times $t_*\sim\eps^{-1}$; for such operators,~\eqref{EqIP0Est} is an unconditional estimate. Taking a (weak subsequential) limit as $\eps\to 0$ produces the forward solution of $P u=f$, and the estimate~\eqref{EqIP0Est} remains valid in the limit.\footnote{The solvability argument in~\cite{HintzNonstat}, by contrast, was based on a duality argument, which, however, appears to be difficult to implement in the presence of zero energy bound states.}

\subsubsection{Higher b-regularity and tame estimates}
\label{SssIPb}

We have thus far only discussed the solvability of $P u=f$ in the context of Theorem~\ref{ThmIML} on rather subtle (variable order, hence microlocally defined) e3b-Sobolev spaces. In order to prove Theorem~\ref{ThmIML} as stated on the technically much simpler b-Sobolev spaces, one first upgrades~\eqref{EqIP0Est} to an estimate
\begin{equation}
\label{EqIPb}
  \|u\|_{H_{\etbop;\bop}^{(s;k),(\beta_\sscri,\beta_+,-\aleph+\frac12)}} \lesssim \|P u\|_{H_{\etbop;\bop}^{(s;k),(\beta_\sscri+1,\beta_++2+\delta,\frac12)}}
\end{equation}
on function spaces which, in addition to $s$ degrees of e3b-regularity, encode $k\in\N_0$ degrees of b-regularity. Upon crudely bounding e3b-derivatives using b-derivatives, one then obtains the pure b-statement~\eqref{EqIMLPu}--\eqref{EqIMLu}.

\pfstep{Issues with direct commutation arguments.} A typical approach towards proving~\eqref{EqIPb} would be to commute b-derivatives through the equation $P u=f$, as done in related contexts in, e.g., \cite[\S{5}]{HintzVasyScrieb} and \cite[\S{5.5}]{HintzNonstat}. This, however, is subtle in the presence of zero energy bound states, as we proceed to discuss in the context of proving $t_*\pa_{t_*}$-regularity. Consider the commuted equation for $t_*\pa_{t_*}u$, which reads
\[
  P(t_*\pa_{t_*}u) = t_*\pa_{t_*}f + [P,t_*\pa_{t_*}]u.
\]
Writing $P=P_{(2)}+P_{(1)}\pa_{t_*}+P_{(0)}\pa_{t_*}^2$ where $P_{(j)}$ is a $j$-th order differential operator in the spatial variables, we have
\begin{equation}
\label{EqIPbComm}
  [P,t_*\pa_{t_*}]=P_{(1)}\pa_{t_*}+2 P_{(0)}\pa_{t_*}^2.
\end{equation}
Suppose forward solutions of $P$ lose $\aleph$ powers of $t_*$-decay relative to $f$. Using only the estimate~\eqref{EqIP0Est}, then, $[P,t_*\pa_{t_*}]u$ has the same $\aleph$-loss, and the inversion of $P$ produces yet another $\aleph$-loss, and so on. When $\aleph=0$, there are no losses, and this procedure succeeds. For $\aleph\geq 1$, however, this would result in an unacceptable extra power of $t_*$ with each $t_*\pa_{t_*}$-derivative. Partial relief may be obtained in concrete settings (e.g., the 1-form wave operator) when the $t_*$-decay loss of $u$ arises solely from a $t_*$-integration (as discussed after~\eqref{EqIP0Grushin}) and is thus undone by the presence of $\pa_{t_*}$ in~\eqref{EqIPbComm}. In the presence of a \emph{double} pole of the resolvent at $\sigma=0$ (as is present for the linearized gauge-fixed Einstein equation \cite{HaefnerHintzVasyKerr,HintzKerrStab}), however, this observation appears to recover at best one, but not two, orders of $t_*$-decay.

\begin{rmk}[More careful accounting: an ODE example]
\label{RmkIPbComm}
  Consider an ODE example involving a loss of $2$ powers of $t_*$-decay, namely
  \[
    P_0 u = f,\quad P_0 := \begin{pmatrix} \pa_{t_*} & 1 \\ 0 & \pa_{t_*} \end{pmatrix},\quad t_*\geq 1.
  \]
  Then the commuted equation reads
  \begin{equation}
  \label{EqIPbCommODE}
    P_0(t_*\pa_{t_*} u) = t_*\pa_{t_*}f + \begin{pmatrix} \pa_{t_*} & 0 \\ 0 & \pa_{t_*} \end{pmatrix} u.
  \end{equation}
  For $f\in L^2$ with $t_*\pa_{t_*}f\in L^2$, we only have $u\in t_*^2 L^2$ and $\pa_{t_*}u\in t_* L^2$, which is insufficient for proving $t_*\pa_{t_*}u$-regularity (relative to $t_*^2 L^2$) of $u$ using~\eqref{EqIPbCommODE}. Only when one considers the components of $u=(u_1,u_2)$ separately can one argue successfully: $u_2$ loses one power of $t_*$-decay relative to $f_2$, and $u_1$ loses one power relative to $f_1$ and $u_2$; thus the second term on the right in~\eqref{EqIPbCommODE} amounts to an extra source term $\tilde f$ with $\tilde f_1\in t_* L^2$ and $\tilde f_2\in L^2$, which in this more precise accounting \emph{is} sufficient to conclude the desired $t_*\pa_{t_*}$-regularity of $u$. --- In the wave equation setting, this suggests that one must keep very precise track of which ``components'' of the solution of $P_0 u=f$ grow relative to which ``components'' of $f$. In applications, we can do this efficiently using Grushin problems, as discussed below.
\end{rmk}

A frequently useful observation (e.g., in arguments involving \emph{module regularity} as introduced in \cite{HassellMelroseVasySymbolicOrderZero}) is that one can rewrite~\eqref{EqIPbComm} as an effective contribution to the linear operator $P$, to wit,
\begin{equation}
\label{EqIPbEffective}
  \bigl(P - P_{(1)}t_*^{-1} - 2 P_{(0)}\pa_{t_*} t_*^{-1}\bigr) (t_*\pa_{t_*} u) = t_*\pa_{t_*}f.
\end{equation}
The operator on the left is equal to $P$ up to \emph{decaying} (and subprincipal) $\cO(t_*^{-1})$ terms. But already in the ODE example~\eqref{EqIODE}, $\cO(t_*^{-1})$-perturbations of the model operator $P_0$ are borderline unacceptable, and they are catastrophic when $P_0$ has a second order zero energy bound state (e.g., for $P_0=\pa_{t_*}^2$, with solutions of $P_0+\cO(t_*^{-1})$ typically growing super-polynomially). Thus, the reformulation~\eqref{EqIPbEffective} is not directly useful.

\pfstep{b-regularity for the model problem.} Our general requirement for the stationary model problem is the strengthening
\begin{equation}
\label{EqIPbP0k}
  P_0 u=f,\ k\in\N_0 \implies \|u\|_{H_{\tbop;\bop}^{(s;k),(\beta_\sface-\delta,-\aleph+\frac12)}(\Omega,\mu_\bop)} \leq C_k\|f\|_{H_{\tbop;\bop}^{(s;k),(\beta_\sface+2,\frac12)}(\Omega,\mu_\bop)}
\end{equation}
of~\eqref{EqIP0aleph}, where the norm of $u$ is the sum of weighted $\Htb^s$-norms of $u$ and its up to $k$-fold b-derivatives (see~\eqref{EqIMbReg2}), similarly for $f$. (See again Definition~\ref{DefSSAlephAdm} for details.)

In practice, this is proved via resolvent estimates for $\wh{P_0}(\sigma)^{-1}$; spatial b-derivatives ($r\pa_x$) commute with the Fourier transform $\cF$, while
\begin{equation}
\label{EqIPbFtdt}
  \cF\circ t_*\pa_{t_*}=-\pa_\sigma\sigma\circ\cF,
\end{equation}
i.e., we need to prove the $\sigma\pa_\sigma$-regularity of the resolvent (which in particular entails b-, or conormal, regularity at $\sigma=0$). Spatial b-regularity for the output of $\wh{P_0}(\sigma)^{-1}$ is to be expected since $\wh{P_0}(\sigma)$ is the spectral family of $P_0$ with respect to the foliation of spacetime by level sets of the function $t_*$, not $t$. As an illustration, on Minkowski space, the latter spectral family is $\Delta-\sigma^2$, and the former is $e^{-i\sigma r}(\Delta-\sigma^2)e^{i\sigma r}$, so in the output of the inverse of this operator, the oscillatory part $e^{i\sigma r}$ is factored out, and what remains is the b-regular $\cO(r^{-1})$-coefficient. The advantage of this ``conjugated'' perspective on asymptotically flat resolvents was highlighted by Vasy \cite{VasyLAPLag,VasyLowEnergyLag}. (A closely related point of view is that this spatial b-regularity is the same as \emph{module regularity} at the outgoing radial set \cite{GellRedmanHassellShapiroZhangHelmholtz}.)

Turning to $\sigma\pa_\sigma$-regularity, we discuss the high- and low-energy regimes separately.
\begin{enumerate}[leftmargin=2em]
\item\textit{High energies.} The main tool is the identity
  \begin{equation}
  \label{EqIPbResId}
    \sigma\pa_\sigma\wh{P_0}(\sigma)^{-1} = -\wh{P_0}(\sigma)^{-1}\circ\sigma\pa_\sigma\wh{P_0}(\sigma)\circ\wh{P_0}(\sigma)^{-1}.
  \end{equation}
  Operator norm bounds on $\wh{P_0}(\sigma)$ in the high-energy regime $|\Re\sigma|\to\infty$ involve powers of $|\sigma|$; and the presence of trapping leads to losses of powers $|\sigma|^{\delta_\Gamma}$, $\delta_\Gamma>0$, relative to non-trapping estimates \cite{BurqNontrapping,BonyBurqRamondTrapping}; see Theorem~\ref{ThmSpHi} for precise statements. Iterating the formula~\eqref{EqIPbResId}, this means that $k$ degrees of b-regularity lead to a loss of $\sim k\delta_\Gamma$ powers of $|\sigma|$, which translates into standard regularity losses upon taking inverse Fourier transforms. For linear theory, this is acceptable; but for the precise (tame) b-theory needed for nonlinear applications, it is not, \emph{unless} one can take $\delta_\Gamma$ to be arbitrarily small (e.g., $\delta_\Gamma=\frac{1}{k}$). We prove such almost lossless high-energy estimates as consequences of almost-sharp estimates at normally hyperbolic trapping; see Theorem~\ref{ThmSpHiTr}, which is based on \cite{DyatlovSpectralGaps} and interpolation arguments.
\item\textit{Low energies.} When $\wh{P_0}(\sigma)^{-1}$ has a pole at $\sigma=0$, the identity~\eqref{EqIPbResId}, applied naively, would suggest that $\sigma\pa_\sigma\wh{P_0}(\sigma)^{-1}$ has a pole of higher order (leading to the apparent loss of decay discussed on the physical space side after~\eqref{EqIPbComm}). The resolution, in the concrete case of the 1-form wave operator, is simply to \emph{apply~\eqref{EqIPbResId} to the augmented operator} $\wt{P_0}(\sigma)$ from~\eqref{EqIP0Grushin}, which \emph{does} have a uniformly bounded resolvent at low energies. This Grushin problem reformulation \emph{automatically} takes care of the otherwise subtle bookkeeping hinted at in Remark~\ref{RmkIPbComm}.
\end{enumerate}

\bigskip
\pfstep{Microlocal e3b-estimates relative to b-spaces.} The proof of~\eqref{EqIPb} for $k\geq 1$ proceeds in conceptually the same fashion as in the case $k=0$: we combine~\eqref{EqIPbP0k} (or rather an upgrade of it to an estimate on $(\etbop;\bop)$-spaces of order $(s;k)$) with e3b-microlocal estimates that now take place relative to a (weighted) space $\Hb^k$ encoding $k$ orders of b-regularity; we proceed to discuss these e3b-microlocal estimates. We emphasize from the outset that these e3b-microlocal arguments do not require any spectral information about the stationary model $P_0$ (cf.\ \eqref{EqIPmlInsensitive}), and thus are insensitive to decay losses in the presence of bound states. Now, while such ``mixed'' estimates are typically proved using an algebra of mixed pseudodifferential-differential (here: $\etbop$-$\bop$) operators (see, e.g., \cite{VasyPropagationCorners,MelroseVasyWunschDiffraction} and \cite[\S\S{2.4} and 5.5]{HintzNonstat}), we use here the scaled bounded geometry approach introduced in \cite{HintzScaledBddGeo} which, in the present context, provides us with e3b-ps.d.o.s whose coefficients are b-regular; these operators act boundedly, and e3b-microlocally, on mixed $H_{\etbop;\bop}^{(s;k)}$-spaces. (This is analogous to the approach used in \cite[\S{2.5.4}]{HintzGlueLocII}.)

The basic elliptic and real principal type propagation estimates are proved in~\S\ref{SsMBasic}. The key structural input is the following. Suppose $P u=f$, and one has e3b-regularity control on $u$ in terms of $f$ (as in~\eqref{EqIPmlEst}). Let $X$ be a b-vector field, then $P(X u)=X f+[P,X]u$; so in order for the e3b-regularity theory for $P$ to be applicable, we need $[P,X]$ to be an e3b-operator again. This is not automatic (since b-vector fields are stronger than e3b-vector fields); rather, it requires $X$ to have good commutation properties with e3b-vector fields. (We formalize this in Lemma~\ref{LemmaCTe3bComm}; see also Lemmas~\ref{LemmaCT3bDil} and~\ref{LemmaCTebComm}.) Note that this is a \emph{structural} requirement (on the level of different algebras vector fields and their commutators), not a geometric one (e.g., Killing vector fields play no role here). Near $\cK^+$, we mainly use $X=t_*\pa_{t_*}$ (\S\S\ref{SssRbR}--\ref{SssRbTr}), and near $\scri^+$ we use $t_*\pa_{t_*}$, $r\pa_r$, and $\pa_\omega$ (\S\ref{SssRbComm}). 

It is, in fact, important that we argue more carefully: the operator $[P,X]$ is typically a second order (e3b-)differential operator, while the regularity theory for $P$ gains at most $1$ derivative (since $P$ is a wave operator), and less in the presence of trapping. In nonlinear applications, we cannot afford the resulting (e3b-)derivative losses (which would increase with the desired number of b-derivatives). We must thus, instead, consider a full set of commuted equations $P(X_i u)=X_i f+[P,X_i]u$, $i=1,\ldots,N$, at once, rewrite $[P,X_i]=\sum_j Q_{i j}X_j$ where $Q_{i j}$ is a first order operator, and regard $(P-\sum_j Q_{j i})(X_i u)=X_i f$, $i=1,\ldots,N$, as a single tensorial wave equation; note that the effective operator $(P\delta_{j i}-Q_{j i})_{i,j=1,\ldots,N}$ has the same scalar principal symbol as $P$ itself. This tensorial wave equation, then, does not anymore feature a source term with regularity loss.\footnote{It is not always necessary to use a full set of b-vector fields $X_i$: for example, near $\cK^+$ but away from the ergoregion, the single vector field $t_*\pa_{t_*}$ suffices since its rescaling $r\pa_{t_*}=\frac{r}{t_*}t_*\pa_{t_*}$ is microlocally e3b-elliptic on the characteristic set of $\Box_{g_{\bhm,a}}$ and thus automatically controls a full 3b-derivative, so in particular giving $r\pa_r$- and $\pa_\omega$-control, and thus full b-regularity.} In such a procedure, subtleties arise: one must ensure that the additional subprincipal terms $Q_{j i}$ do not shift threshold conditions at radial sets and do not affect sign conditions on subprincipal symbols at the trapped set. The details are in~\S\ref{SsRb}. We only point out that at the trapped set, the effective operator in the the commuted equation~\eqref{EqIPbEffective} differs from $P$ by a \emph{decaying} perturbation, which therefore does not affect the subprincipal symbol condition (cf.\ Definition~\ref{DefSSTrapAdm}) for the stationary model.

\bigskip
\pfstep{Tame estimates.} The above arguments imply the b-regularity of linear waves stated in~\eqref{EqIMLPu}--\eqref{EqIMLu}. The proof of the tame estimate~\eqref{EqIMLbTame}, in which the constant in the estimate is quantitatively controlled by the coefficients of $P-P_0$, is conceptually very simple, even if the detailed arguments require some care. We first make two observations:
\begin{enumerate}[leftmargin=2em]
\item It suffices to prove a b-tame version of the e3b-microlocal regularity estimate
  \begin{equation}
  \label{EqIPbTame}
    \|u\|_{H_{\etbop;\bop}^{(s;k),\cdots}} \lesssim \|P u\|_{H_{\etbop;\bop}^{(s;k),\cdots}} + \|u\|_{H_{\etbop;\bop}^{(s-10;k),\cdots}}
  \end{equation}
  only, which is the analogue of~\eqref{EqIPmlEst} with $k$ degrees of additional b-regularity. Indeed, the remainder term is then estimated using estimates for the stationary model $P_0$---which has smooth coefficients and in applications varies at worst in some finite-dimensional space of operators (e.g., the linearized gauge-fixed Einstein equation on subextremal Kerr spacetimes).
\item The proof of a e3b-microlocal estimate for a PDE $P u=f$ at e3b-regularity order $s$ requires only some finite (but sufficiently large, depending on $s$) \emph{e3b-regularity} of the coefficients of $P$. (This follows immediately via the bounded geometry perspective from the corresponding statement on $\R^n$; see Lemma~\ref{LemmaMSOpNorm}.)
\end{enumerate}

The idea is thus to prove a b-tame version of~\eqref{EqIPbTame} by applying the usual e3b-microlocal estimates to $k$-times-commuted equations, \emph{always at fixed e3b-regularity order}; in the commuted equations---schematically, $P D^k u=D^k f+[P,D^k]u\sim D^k f+\sum_j [D,[\cdots,[D,P]\cdots]] D^{k-j}u$ with a $j$-fold commutator---we put low-order ($j\leq 2$ or $3$) commutators of $P$ with b-vector fields ($D$ in this schematic notation) into an appropriate effective operator, while the higher-order commutator terms are put on the right-hand side as additional source terms. One must then only prove b-tame estimates for this additional source term (which involves less than $k$ b-derivatives, and indeed less than $k$ derivatives of any sort, acting on $u$, and can thus be estimated by induction on $k$); but these follow almost directly from basic b-tame estimates for products of Sobolev functions (Lemma~\ref{LemmaMTameMult} and its e3b-microlocal version, Proposition~\ref{PropMTameMicr}). See Proposition~\ref{PropMTamePr} for an illustrative example, and Lemma~\ref{LemmaRbRHComm} for an exemplary commuted equation which is utilized in the proof of Proposition~\ref{PropRbRH}. The b-tame estimate at the trapped set is proved in~\S\ref{SssRbTr}.

\subsubsection{Improved decay}
\label{SssIPD}

For some nonlinear applications such as Theorems~\ref{ThmI1} and \ref{ThmINo}\eqref{ItINoPower} with $p=3$, one needs to prove stronger decay for linear waves than what Theorem~\ref{ThmIML} provides.

Asymptotics and decay at null infinity $\scri^+$ are straightforward to obtain via integration along approximate characteristics. For example, the tensor wave operator $P$ on Kerr is given to leading order at $\scri^+$ by the right-hand side of~\eqref{EqIPmlBoxeb}; given the already known b-regularity of a solution of $P u=f$, the term $(x_\sscri\pa_\omega)^2 u$ is of lower order in the sense of decay and can thus be moved to the right-hand side. Upon integration of $x_\sscri\pa_{x_\sscri}-2\rho_+\pa_{\rho_+}$ from $t_*=1$ (where $u$ and $f$ vanish), one therefore obtains an equation for $(x_\sscri\pa_{x_\sscri}-2)u$, the integration of which yields the $x_\sscri^2$-radiation field of $u$ (plus terms with more decay as $x_\sscri\to 0$); see~\S\ref{SsDscri} for details.

\begin{figure}[!ht]
\centering
\includegraphics{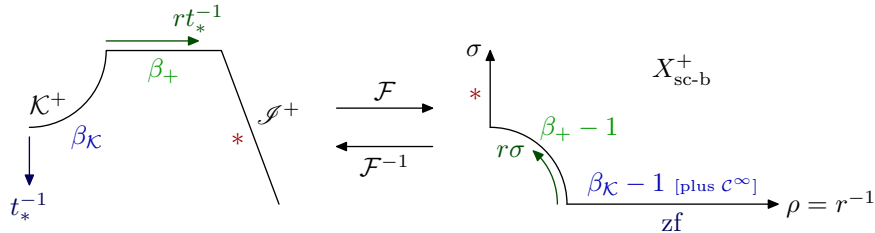}
\caption{Rough illustration of the relationship between decay rates of functions on spacetime (relative to $L^2$ with density $|\frac{\dd t_*}{\la t_*\ra}\,\frac{\dd x}{r^3}|$) and the decay rates (relative to $L^2$ with density $|\frac{\dd\sigma}{\sigma}\,\frac{\dd x}{r^3}|$) of their Fourier transforms at low energies. }
\label{FigIPD}
\end{figure}

We establish improved decay at $\iota^+\cup\cK^+=t_*^{-1}(\infty)$ using the stationary model $P_0$ via the Fourier transform in $t_*$ and improved control on the low-regularity resolvent. (The high-frequency parts of waves have arbitrary $t_*$-decay if one gives up a corresponding amount of b-regularity; see Corollary~\ref{CorDResHiLoc}.) The main technical tool is a precise relationship between the decay rates of functions on spacetime at $\iota^+$ and $\cK^+$ and the decay rates of their Fourier transforms on the manifolds with corners $X_\scbtop^\pm$ from Figure~\ref{FigIP0scbt}; this relationship, which also (roughly) preserves b-regularity (cf.\ the discussion around~\eqref{EqIPbFtdt}), is illustrated in Figure~\ref{FigIPD} and made precise in Propositions~\ref{PropDFTInv} and \ref{PropDFTInv}. It essentially follows from a quantitative variant of the well-known fact that the Fourier transform of $(t_*^{-1})^\alpha$ (appropriately interpreted as a homogeneous or $\pm i 0$ type distribution) is $\sigma^{\alpha-1}$ (see~\S\ref{SsDFT} for precise statements).

Concretely, one rewrites $P u=f$ as $P_0 u=\tilde f:=f-(P-P_0)u$, so $u=P_0^{-1}\tilde f$. This is useful when $P-P_0$ has more $t_*$-decay (or, more precisely, more $\cK^+$- and $\iota^+$-decay) than what the inversion of $P_0$ is expected to lose. On the Fourier transform side, $\cF\tilde f$ then has better decay orders (and better regularity at $\zface$) at the boundary hypersurfaces of $X_\scbtop^\pm$. Uniform bounds on $\wh{P_0}(\sigma)^{-1}$ at low energies imply similarly improved decay orders for $\cF u$ (see Proposition~\ref{PropDResLo}), which translates into stronger decay for $u$ itself. This is implemented for scalar waves in~\S\ref{SsDp}.

For the 1-form wave operator, and thus in the presence of zero energy bound states, one instead uses uniform bounds for the resolvent of the \emph{augmented} operator~\eqref{EqIP0Grushin} and controls the inverse Fourier transform of the singular part $\hat a(\sigma)\sigma^{-1}u_0$ of~\eqref{EqIuregsing} via $t_*$-integration of the inverse Fourier transform of $\hat a$ (which inherits regularity and decay at $\sigma=0$ from $\cF\tilde f$). This is implemented in~\S\S\ref{SssA2B}--\ref{SssA2lot}.

\subsubsection{Analysis of an ODE toy model}
\label{SssIPODE}

Many\footnote{but of course not all: there are no propagation estimates, radial sets, or trapping, and there is only one asymptotic regime---$t_*\to\infty$---instead of three.} aspects of our approach can be illustrated in the ODE toy model
\[
  P_0 u=f+(\pa_{t_*}u\cdot\pa_{t_*}u) u,\quad P_0=\begin{pmatrix}\pa_{t_*} & 0 \\ 0 & \pa_{t_*}+1\end{pmatrix},
\]
taken from~\eqref{EqIODE}, but slightly modified for notational simplicity. We aim to prove the global existence of solutions $u$ in $t_*\geq 1$ when $f$ is sufficiently small and decaying, with the solution being of the form $u(t_*)=(c+a(t_*))u_0+\tilde u(t_*)$ as in~\eqref{EqIODEDecomp}, where $u_0=(1,0)$ is the zero energy state of $P_0$; more precisely, working with weighted b-Sobolev spaces that measure $t_*\pa_{t_*}$-regularity relative to $L^2([1,\infty),|\frac{\dd t_*}{t_*}|)$,
\begin{equation}
\label{EqIPODEStruct}
  c\in\C,\quad
  a \in \Hb^{\infty,\alpha-1}([1,\infty];\C)=t_*^{-(\alpha-1)}\Hb^\infty,\quad
  \tilde u\in\Hb^{\infty,\alpha}([1,\infty];\C^2)=t_*^{-\alpha}\Hb^\infty,
\end{equation}
and $\alpha\in(1,2)$. In a Nash--Moser iteration scheme, we are thus led to studying the linearization of $u\mapsto P_0 u-f-(\pa_{t_*}u)^2 u$ around such $u$, which is given by
\[
  P\dot u = (P_0 + \tilde P)\dot u,\quad \tilde P := -2 u\pa_{t_*}u\cdot\pa_{t_*}(\cdot) - (\pa_{t_*}u\cdot\pa_{t_*}u)(\cdot).
\]
As a cusp operator, i.e., expressed in terms of $\pa_{t_*}$, the coefficients of $\tilde P$ are of class $\Hb^{\infty,\alpha}$ and thus have pointwise $t_*^{-\alpha}$-decay (together with all b-, i.e., $t_*\pa_{t_*}$-, derivatives).

\pfstep{Analysis in the cusp setting.} We use the term ``cusp'' to refer to regularity with respect to $\pa_{t_*}$; this is the ODE analogue of e3b-regularity discussed in~\S\ref{SssIPmlG}; cf.\ the first vector field in~\eqref{EqIPmlG3b}, restricted to, say, $r=1$. (The terminology comes from \cite{MazzeoMelroseHyp} and in the present setting refers to regularity under $\tau^2\pa_\tau$ where $\tau=t_*^{-1}$.) When $u$ is small in $\cC^1$, the operator $P$ is uniformly elliptic, and therefore we immediately obtain, for any fixed $s,\beta\in\R$, the a priori estimate
\begin{equation}
\label{EqIPODEcu}
  \|\dot u\|_{H_\cuop^{s,\beta}} := \| t_*^\beta\dot u \|_{H^s(\R,|\frac{\dd t_*}{t_*}|)} \lesssim \|P\dot u\|_{H_\cuop^{s,\beta}} + \|\dot u\|_{H_\cuop^{-N,\beta}},\quad \supp\dot u\subset\{t_*\geq 1\}.
\end{equation}
(The precise cusp-regularity orders are not important, so we do not keep track of elliptic gains.) For example, one can localize $P\dot u=f$ to $t_*$-intervals of length $1$ and use elliptic estimates in each such interval. This is the analogue of~\eqref{EqIPmlEst}.

\pfstep{Estimates for the stationary model.} Following~\S\ref{SssIP0}, we improve the error term~\eqref{EqIPODEcu} using an unconditional estimate for $P_0$; concretely, we claim that, for any fixed $s\in\R$,
\begin{subequations}
\begin{equation}
\label{EqIPODEP0b}
  \|\dot u\|_{H_\cuop^{s,-\frac12}} \lesssim \|P_0\dot u\|_{H_\cuop^{s,\frac12}},\quad \supp\dot u\subset\{t_*\geq 1\},
\end{equation}
which is the analogue of~\eqref{EqIP0aleph} with $\aleph=1$. Passing to the density $|\dd t_*|$, this reads
\begin{equation}
\label{EqIPODEP0}
  \|\dot u\|_{H_\cuop^{s,-1}(\R,|\dd t_*|)} \lesssim \|f\|_{H_\cuop^s(\R,|\dd t_*|)},\quad f:=P_0 u,\ \ \supp f,\ \supp u\subset\{t_*\geq 1\}.
\end{equation}
\end{subequations}
We prove this using the Fourier transform; we only discuss the case\footnote{Other values of $s$ only cause overall factors of $\la\sigma\ra^{-s}$ on the Fourier transform side, which do not affect our arguments.} $s=0$. We set
\begin{equation}
\label{EqIPODEP0FT}
  \hat u(\sigma) := \wh{P_0}(\sigma)^{-1}\hat f(\sigma),\quad \wh{P_0}(\sigma)=\begin{pmatrix} -i\sigma & 0 \\ 0 & -i\sigma+1 \end{pmatrix},
\end{equation}
initially only for $\sigma\in\C$ with $\Im\sigma\geq 0$ and $\sigma\neq 0$. Since $\hat f\in L^2(\R)$, we have $\hat u|_{\R\setminus[-1,1]}\in L^2$. We analyze the singularity of $\wh{P_0}(\sigma)^{-1}$ at $\sigma=0$ using the augmentation
\[
  \wt{P_0}(\sigma) = \begin{pmatrix} \wh{P_0}(\sigma) & \wh{P_0}(\sigma)(\sigma^{-1}u_0) \\ e_1^* & 0 \end{pmatrix},
\]
where $e_1^*$ is a linear functional mapping $u_0\mapsto 1$. Note that $\wt{P_0}(\sigma)$ \emph{is} invertible at $\sigma=0$, and thus nearby. Thus, the quantities
\begin{equation}
\label{EqIPODEP0Sing}
  (\hat u_{\rm reg}(\sigma),\ \hat a(\sigma)) := \wt{P_0}(\sigma)^{-1} \bigl( \hat f(\sigma), 0 \bigr)
\end{equation}
are controlled in $L^2([-1,1])$ by $\|f\|_{L^2}$. But then $\hat u(\sigma)=\hat u_{\rm reg}(\sigma)+\sigma^{-1}\hat a(\sigma)u_0$. Evaluating the inverse Fourier transform of $\hat u$ on the contour $\Im\sigma=1$, say, and shifting to the real axis, one can then show that $u$ is the sum of a regular piece in $L^2([1,\infty);\C^2)$ and the $t_*$-integral of an $L^2$-function (with compact Fourier support), which thus lies in $t_* L^2$. This gives~\eqref{EqIPODEP0}.

\pfstep{Unconditional estimate for $P$.} We plug~\eqref{EqIPODEP0b} with $s=-N$ into the final term of~\eqref{EqIPODEcu} for $\beta=-\frac12$ and replace $P_0$ by $P$, leading to an error
\[
  \|(P-P_0)\dot u\|_{H_\cuop^{-N,\frac12}} \lesssim \|\dot u\|_{H_\cuop^{-N+1,\frac12-\alpha}}.
\]
Since $\frac12-\alpha<-\frac12$, this error is \emph{small} (compared to the $H_\cuop^{s,-\frac12}$-norm of $\dot u$ when taking $-N+1\leq s$) when $\dot u$ is supported in $t_*\geq T$ for sufficiently large $T$. Controlling $\dot u$ up to time $T$ in terms of $P\dot u$ is trivial by standard finite-time ODE theory. We thus obtain
\begin{equation}
\label{EqIPODEEst}
  \|\dot u\|_{H_\cuop^{s,-\frac12}} \lesssim \|P\dot u\|_{H_\cuop^{s,\frac12}},\quad \supp\dot u\subset\{t_*\geq 1\}.
\end{equation}
This is the analogue of~\eqref{EqIP0Est}. The approximation argument sketched after~\eqref{EqIP0Est} can be used also in the present context and shows: given $f\in H_\cuop^{s,\frac12}$ with support in $t_*\geq 1$, the forward solution $\dot u$ of $P\dot u=f$ satisfies $\dot u\in H_\cuop^{s,-\frac12}$ and the estimate~\eqref{EqIPODEEst}. (This is considerably simpler than the duality arguments in the proof of \cite[Theorem~A.1]{HintzNonstat}.)

\pfstep{Higher b-regularity.} As discussed in~\S\ref{SssIPb}, we need two ingredients. The \emph{first ingredient} is a version of~\eqref{EqIPODEP0b} that includes $k\in\N_0$ additional degrees of b-regularity (i.e., $t_*\pa_{t_*}$) on $\dot u$ and $P_0\dot u$, to wit,
\[
  \|\dot u\|_{H_{\cuop;\bop}^{(s;k),-\frac12}} \lesssim \|P_0\dot u\|_{H_{\cuop;\bop}^{(s;k),\frac12}},\qquad
  \|\dot u\|_{H_{\cuop;\bop}^{(s;k),\beta}} := \sum_{j=0}^k \|(t_*\pa_{t_*})^j t_*^\beta\dot u\|_{H^s}.
\]
On the spectral side, this amounts to controlling $k$-many $\sigma\pa_\sigma$-derivatives of~\eqref{EqIPODEP0FT}. At high frequencies, this follows from the fact that $\|\wh{P_0}(\sigma)^{-1}\|\lesssim|\sigma|^{-1}$ and the identity~\eqref{EqIPbResId} (iterated $k$ times). At low frequencies, one uses the b-regularity (i.e., $\sigma\pa_\sigma$) of $\wt{P_0}(\sigma)^{-1}$ at $\sigma=0$ to deduce the b-regularity of the pieces $\hat u_{\rm reg}$ and $\hat a$ in~\eqref{EqIPODEP0Sing}, which translate into the desired $t_*\pa_{t_*}$-regularity upon inverse Fourier transforming.

The \emph{second ingredient} is the analogue of~\eqref{EqIPODEcu} with extra b-regularity,
\begin{equation}
\label{EqIPODEcub}
  \|\dot u\|_{H_{\cuop;\bop}^{(s;k),\beta}} \lesssim \| P\dot u \|_{H_{\cuop;\bop}^{(s;k),\beta}} + \|\dot u\|_{H_{\cuop;\bop}^{(-N;k),\beta}}.
\end{equation}
Now, this cannot be proved by localizing in unit-length intervals in $t_*$ using, say, a translated version of a cutoff function $\chi\in\CIc((-1,1))$: the issue is that localization using $\chi(\cdot-T)$ destroys uniform b-regularity (since $t_*\pa_{t_*}\chi(t_*-T)\sim T\chi'$). Instead, one can localize in dyadic intervals in $t_*$ using, say, localizers $\chi(\frac{t_*}{T})$ where $\chi\in\CIc((1,10))$ equals $1$ on $[2,9]$ and $T$ is a power of $2$. On $t_*$-intervals $[T,10 T]$, one then uses uniform elliptic estimates relative to a space $H_{-;T}^{(s;k)}$ with norm $\|v\|_{H_{-;T}^{(s;k)}}=\sum_{j=0}^k\|(T\pa_{t_*})^j v\|_{H^s}$. (This is the scaled bounded geometry perspective on mixed spaces, cf.\ \cite[(3.16)]{HintzScaledBddGeo}.) One can replace $T\pa_{t_*}$ by $t_*\pa_{t_*}$ here. The point is that $t_*\pa_{t_*}$ commutes well with $P$ in that $[t_*\pa_{t_*},P]$ has uniformly bounded coefficients as a cusp-operator (i.e., when written in terms of $\pa_{t_*}$); this is an instance of the statement in the line after~\eqref{EqCT3bDil} (upon dropping all spatial variables).

Upon combining these two ingredients and repeating the approximation argument, one obtains the unconditional estimate
\[
  \|\dot u\|_{H_{\cuop;\bop}^{(s;k),-\frac12}} \leq C_k\|P\dot u\|_{H_{\cuop;\bop}^{(s;k),\frac12}}.
\]
To prove a version of this estimate that is tame in the b-regularity order, one repeats the same argument but starting with a b-tame version of~\eqref{EqIPODEcub}. The latter can be proved as follows: using Lemma~\ref{LemmaMTameComm}, write
\[
  P((t_*\pa_{t_*})^k\dot u) = (t_*\pa_{t_*})^k P\dot u - \sum_{j=1}^k ({\rm ad}_{t_*\pa_{t_*}}^j P) (t_*\pa_{t_*})^{k-j}\dot u,\quad {\rm ad}_{t_*\pa_{t_*}}:=[t_*\pa_{t_*},\cdot].
\]
The top order term of the $j=1$ summand is of the schematic form $(t_*\pa_{t_*}p)\pa_{t_*}\circ (t_*\pa_{t_*})^{k-1}\dot u$ where $p$ is a coefficient of $P$, so $t_*\pa_{t_*}p=\cO(t_*^{-\alpha})$, and is thus of the form $\cO(t_*^{-\alpha}\cdot t_*^{-1})(t_*\pa_{t_*})^k\dot u$; move this to the left-hand side. Apply the \emph{basic} estimate~\eqref{EqIPODEcu} (needed only for $\beta=-\frac12$) to the resulting equation; it then remains to bound, in a b-tame manner, the terms
\[
  \|({\rm ad}_{t_*\pa_{t_*}}^j P) (t_*\pa_{t_*})^{k-j}\dot u\|_{H_\cuop^{-N,-\frac12}},\quad j\geq 2,
\]
which involve at most $1+(k-2)\leq k-1$ derivatives on $\dot u$. For $N=0$, this can be estimated using standard Moser-type estimates (in the coordinate $\log t_*$), see, e.g., \cite[Chapter~13, Proposition~3.6]{TaylorPDE3}; these give the desired tame estimates, with the low-high terms involving at most $k-1$ b-derivatives of $\dot u$, which can thus be estimated inductively. See (the proof of) Proposition~\ref{PropMTameMicr} on how to deal with general $N$.

\pfstep{Improving decay.} Finally, to close a nonlinear iteration scheme, we must recover the structure~\eqref{EqIODEDecomp}, \eqref{EqIPODEStruct} of the solution of $P\dot u=f$ when $f\in\Hb^{\infty,\alpha}$. At present, we only have $\dot u\in\Hb^{\infty,-\frac12}$. But then the source term of
\begin{equation}
\label{EqIPODERewrite}
  P_0\dot u = \tilde f:=f - (P-P_0)\dot u \in \Hb^{\infty,\alpha} + t_*^{-\alpha}\Hb^{\infty,-\frac12} = \Hb^{\infty,\frac12+\eps},\quad \eps:=\alpha-1,
\end{equation}
has better-than-$t_*^{-\frac12}$-decay. On the Fourier transform side, this amounts to better-than-$L^2$-regularity at $\sigma=0$. This is preserved by $\wt{P_0}(\sigma)^{-1}$, and upon taking the inverse Fourier transform, one obtains
\[
  \dot u = \dot u_{\rm reg} + \dot a(t_*)u_0,\quad \dot u_{\rm reg}\in\Hb^{\infty,\frac12+\eps}([1,\infty];\C^2),\ \dot a\in\Hb^{\infty,-\frac12+\eps}([1,\infty];\C),
\]
with $\dot a$ arising from the $t_*$-integration of an element of $\Hb^{\infty,\frac12+\eps}$. This improved decay can be plugged into the right-hand side of~\eqref{EqIPODERewrite}; this provides an iterative procedure in which we improve the $t_*$-decay of $\dot u_{\rm reg}$ and $u_0$ at each step by the fixed amount $\eps$. (This is the analogue of the arguments in~\S\ref{SssA2B}.)

A qualitative change occurs when the decay of $\dot u$ has been improved to the extent that now $\tilde f\in\Hb^{\infty,1+\eta}$ for some $\eta>0$: then $\dot a$, arising from integration of an element of $\Hb^{\infty,1+\eta}$, is equal to a constant leading order term plus a remainder in $\Hb^{\infty,\eta}$. Iterating this further yields the desired description of $\dot u=(\dot c+\dot a(t_*))u_0+\tilde{\dot u}(t_*)$ in terms of pieces $\dot c,\dot a,\tilde{\dot u}$ that are of class~\eqref{EqIPODEStruct}. (This is the analogue of the arguments in~\S\S\ref{SssA2Lot}--\ref{SssA2lot}.)

\subsection{Outline}
\label{SsIPlan}

The plan of the paper is as follows.

\begin{itemize}[leftmargin=2em]
\item \S\ref{SC}: \textit{compactifications and phase spaces}. We define the compactified spacetime manifold and introduce the notions of b- and e3b-vector fields and -phase spaces motivated in~\S\ref{SssIPmlG}, as well as the basic commutation properties of b-vector fields needed for our proof of higher b-regularity.
\item \S\ref{SM}: \textit{microlocal tools}. Our analysis of wave-type operators on dynamical spacetimes in e3b-Sobolev spaces being almost entirely microlocal in nature. We define the relevant classes of pseudodifferential operators after a review of the scaled bounded geometry perspective of \cite{HintzScaledBddGeo}. We also describe Plancherel-type characterizations of e3b-spaces on the Fourier transform and Mellin transform sides; and we prove basic e3b-microlocal estimates, including b-tame versions of those. This machinery will be at the heart of our proof of tame e3b-microlocal regularity estimates of the form~\eqref{EqIPbTame} later on.
\item \S\ref{STs}: \textit{null-geodesic dynamics on subextremal Kerr spacetimes.} We describe in full detail the dynamics of the null-geodesic flow on Kerr lifted to e3b-phase space. In spatially compact regions, this collects known results from the literature. We give a self-contained account in the language of the present paper, including at $\iota^+$ and $\scri^+$.
\item \S\ref{SS}: \textit{dynamical metrics and operators.} We describe which perturbations of the stationary Kerr metric and stationary wave-type operators on Kerr backgrounds we can handle; this includes a description of the spectral assumptions on the stationary model (including $\aleph$-admissibility).
\item \S\ref{SDy}: \textit{dynamics on dynamical spacetimes.} We give a rough (but sufficient for our purposes) description of the null-bicharacteristic dynamics on dynamical spacetimes. Over the boundary hypersurfaces of the compactified spacetimes, this is the same as for the Kerr metric. Only the description of the perturbed (un)stable manifold of the trapped set requires some discussion, which we base on \cite{HintzPolyTrap}.
\item \S\ref{SR}: \textit{e3b-microlocal control on spacetime.} We begin the microlocal analysis of dynamical wave-type operators in earnest; in addition to the basic (elliptic, real principal type) estimates from~\S\ref{SM}, this includes proofs of estimates at each of the radial sets (see Figure~\ref{FigIPmle3bFlow}) and at the trapped set. We also prove (b-tame) estimates on spaces encoding additional b-regularity.
\item \S\ref{SE}: \textit{energy estimates.} We prove energy estimates on the same function spaces on which the microlocal analysis takes place, as well as b-tame versions thereof. The basic energy estimates are proved using standard energy methods, while higher regularity is proved using microlocal propagation results.
\item \S\ref{SSp}: \textit{estimates for spectral families.} This section is not needed for the proof of our main linear result (Theorem~\ref{ThmIML}), but it collects all results necessary to verify, in practice, the $\aleph$-admissibility requirement of concrete wave-type operators arising in applications. Many (though not all) of the results are contained in the literature, and this section may serve as a reference point for future investigations of spectral and scattering theory on asymptotically flat spacetimes.
\item \S\ref{SN}: \textit{sharper estimates for the stationary model.} In this quite technical section, we upgrade the (time-translation-invariant) estimate for $P_0$ required by $\aleph$-admissibility (and taking place on 3b-Sobolev spaces) to an estimate on e3b-Sobolev spaces, cf.\ the discussion leading up to~\eqref{EqIP0alephe3b}.
\item \S\ref{SF}: \textit{main linear result.} By combining the e3b-regularity estimate with the estimate for the stationary model, we prove a precise version of Theorem~\ref{ThmIML} on b-tame estimates for forward solutions of linear wave equations on asymptotically subextremal Kerr spacetimes; see Theorem~\ref{ThmF}.
\item \S\ref{SA1}: \textit{first nonlinear applications.} We state the Nash--Moser theorem and supply the necessary smoothing operators in order to deduce, from Theorem~\ref{ThmF}, the small data global solvability of simple nonlinear wave equations (Theorem~\ref{ThmINo}).
\item \S\ref{SD}: \textit{improving decay.} We develop the tools, largely based on a dictionary between spacetime decay and low-energy bounds via the Fourier transform (cf.\ Figure~\ref{FigIPD}), needed to prove sharper decay results for solutions of wave-type equations.
\item \S\ref{SA2}: \textit{nonlinear waves in the presence of zero energy bound states.} We show how to combine Theorem~\ref{ThmF} with the methods of~\S\ref{SD} to obtain a precise late-time description of linear waves in the presence of zero energy bound states and use this to prove Theorem~\ref{ThmI1} on the small data global existence of a suitable nonlinear wave equation with zero energy bound states.
\end{itemize}

\subsection*{Acknowledgments.}

I gratefully acknowledge support from the U.S.\ National Science Foundation under the grant DMS-2554160.

\section{Spacetime compactifications}
\label{SC}

The material in this section is essentially taken from \cite[\S{3.2}]{HintzNonstat}, except here we start from the perspective taken in \cite[\S{2}]{HintzKerrCD} (where, on the other hand, the light cone at future null infinity is not resolved). In short, we combine the compactification near null infinity introduced in Hintz--Vasy \cite{HintzVasyScrieb} with the compactification in the forward timelike cone introduced by the author in \cite{Hintz3b}. The virtue of these compactifications is that the metrics and the wave-type operators under study in this paper have smooth (or at least conormal) coefficients on these compactifications (and this is not true for any `smaller' compactification); thus, these compactifications provide convenient settings for the uniform analysis of regularity and asymptotics. However, in this section we shall use metrics and operators only as \emph{motivations} for the introduction of compactifications and tangent bundles over them; we explain in~\S\ref{STs} how the Kerr metric relates to the notions developed here.

We recall that the radial compactification of $\R^n$, $n\in\N$, is defined by
\begin{equation}
\label{EqCRadCpt}
  \ol{\R^n} := \Bigl[ \R_Z^n \sqcup \Bigl( [0,\infty)_\varrho \times \Sph^{n-1}_\varpi \Bigr) \Bigr] \Big/ \sim,\quad 0\neq Z=R\varpi \sim (\varrho,\varpi)=(R^{-1},\varpi),
\end{equation}
where $R=|Z|$ is the Euclidean norm and $\varpi=\frac{Z}{|Z|}\in\Sph^{n-1}$. Thus, $\ol{\R^n}$ is a smooth manifold with boundary $\pa\ol{\R^n}=\varrho^{-1}(0)\cong\Sph^{n-1}$. If in the region $\{\pm Z^1>0\}$ we introduce the projective coordinates
\[
  \rho_1 := \frac{1}{\pm Z^1},\quad
  \hat Z^i := \frac{Z^i}{Z^1},\ i=2,\ldots,n,
\]
then a chart of $\ol{\R^n}$ covering the half space $\ol{\{\pm Z^1>0\}}\setminus\ol{\{Z^1=0\}}$ is given by $[0,\infty)_{\rho_1}\times\R^{n-1}_{(\hat Z^2,\ldots,\hat Z^n)}$. Linear isomorphisms of $\R^n$ extend uniquely to diffeomorphisms of $\ol{\R^n}$ (see \cite[Lemma~6.50]{HintzMicro}); thus $\ol{\R^n}$ only depends on the linear structure of $\R^n$.

We work on $\R^4=\R_t\times\R^3_x$ and write its standard coordinates as
\begin{equation}
\label{EqCCoord}
  z=(t,x),\quad x=(x^1,x^2,x^3),\quad x=r\omega,\ r=|x|=\Biggl(\;\sum_{i=1}^3(x^i)^2\Biggr)^{\frac12},\ \omega=\frac{x}{|x|}\in\Sph^2.
\end{equation}
We moreover define the inverse radial coordinate
\begin{equation}
\label{EqCCoordRho}
  \rho := \frac{1}{r}.
\end{equation}
We shall mainly be interested in (compactifications of) a neighborhood of the domain of outer communications of the subextremal Kerr spacetime. For now, it suffices to fix a mass parameter
\begin{equation}
\label{EqCMass}
  \bhm \in (0,\infty).
\end{equation}
The event horizon $\cH^+$ of subextremal Kerr with such mass, in Boyer--Lindquist coordinates (recalled in~\S\ref{SsTsK} below), lies at $r=r_+$ where $r_+\in(\bhm,2\bhm]$ (depending on the specific angular momentum).

\begin{definition}[Compactification of the spatial manifold]
\label{DefCSpatial}
  Recalling~\eqref{EqCMass}, we define the compactified spatial manifold by
  \begin{equation}
  \label{EqCSpatial}
    X := \ol{\{r\geq\bhm\}} \subset \ol{\R^3_x} = [0,\bhm^{-1}]_\rho\times\Sph^2_\omega \subset \tilde X := \ol{\R^3}.
  \end{equation}
  We denote its boundary hypersurfaces (by a mild abuse of notation) by\footnote{Analytically, $\pa X$ will be a boundary at infinity where $X$ should be considered singular, i.e., the notion of regularity we will use there will be \emph{b-} or \emph{scattering} regularity, defined below. By contrast, $\pa_\sharp X$ is an artificial boundary which will lie in the interior of the Kerr black hole and is introduced to localize in a neighborhood of its domain of outer communications.}
  \[
    \pa_\sharp X := r^{-1}(\bhm),\quad \pa X := \rho^{-1}(0).
  \]
\end{definition}

\begin{figure}[!ht]
\centering
\includegraphics{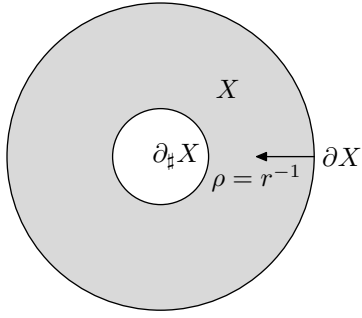}
\caption{The compactified spatial manifold from Definition~\ref{DefCSpatial}, here drawn in $2$ (instead of $3$) dimensions.}
\label{FigCSpatial}
\end{figure}

See Figure~\ref{FigCSpatial}. Thus, $\rho$ is \emph{boundary defining function} for $\pa X$. Recall here that a boundary defining function for an embedded boundary hypersurface $H$ of a manifold with corners $M$ is a smooth function $\rho_H\in\CI(M;[0,\infty))$ such that $H=\rho_H^{-1}(0)$ and $\dd\rho_H\neq 0$ on $H$. We shall also often work with local boundary defining functions $\rho_H\in\CI(U;[0,\infty))$ where $U\subset M$ is open, which means that $H\cap U=\rho_H^{-1}(0)$ and $\dd\rho_H\neq 0$ on $U\cap H$.\footnote{Equivalently, for every compact $K\subset U$, there exists a boundary defining function of $H$ that equals $\rho_H$ on $K$.}

\subsection{Compactifications of the spacetime manifold; bundles}
\label{SsCM}

Smooth functions on the radial compactification $\ol{\R_{(t,x)}^4}$ are, in $\ol{\{t\geq 0\}}\setminus\ol{\{t=0\}}$, smooth functions of $\frac{1}{t}\in[0,\infty)$ and $\frac{x}{t}\in\R^3$. This does not include any non-constant smooth function of $x$ alone. The first step towards defining a compactification suitable for stationary metrics such as Kerr is thus to resolve the `north pole'
\begin{equation}
\label{EqCMkp}
  \fk^+ := \Bigl\{\frac{1}{t}=0,\ \frac{x}{t}=0\Bigr\} \subset \pa\ol{\R^4}
\end{equation}
via a real blow-up; we recall that this produces a new manifold with corners in which polar coordinates around $\fk^+$ are smooth down to the polar coordinate origin. A second resolution is required when discussing radiative spacetimes, where the metric coefficients may depend in a nontrivial manner on $t-r$ near $r^{-1}=0$. We denote the set of endpoints of future null-geodesics on $\R^4$ equipped with the Minkowski metric $\ubar g=-\dd t^2+\dd r^2+r^2\slg$ by
\begin{equation}
\label{EqCMscri}
  Y^+ := \Bigl\{\frac{1}{t}=0,\ \frac{r-t}{t}=0\Bigr\} \subset \pa\ol{\R^4};
\end{equation}
this is the ``light cone at future infinity'' and diffeomorphic to $\Sph^2$. For the following definition, recall that the \emph{lift} of a subset $A\subset M$ to the blow-up $[M;X]$ of $M$ at a submanifold $X$ is the closure of $A\setminus X$ in $[M;X]$ when $A\not\subset X$, and otherwise the preimage of $A$ under the blow-down (polar coordinate) map $[M;X]\to M$. (For an extensive discussion of blow-ups and related notions, see Melrose~\cite{MelroseDiffOnMwc}.)

\begin{definition}[Spacetime compactifications]
\label{DefCMSpacetime}
  We use the notation~\eqref{EqCMkp}, \eqref{EqCMscri}.
  \begin{enumerate}
  \item\label{ItCMSpacetimeM0}{\rm (Non-radiative compactification.)} Let $\fk^-$ be the antipodal point of $\fk^+$ (i.e., the ``south pole'' of $\ol{\R^4}$). We define
    \[
      \tilde M_0 := [\ol{\R^4}; \fk^+,\fk^-],\quad \upbeta_0\colon \tilde M_0\to\ol{\R^4}.
    \]
    Its boundary hypersurfaces are denoted $\sface$ (the lift of $\pa\ol{\R^4}$) and $\cK^\pm$ (the lift of $\fk^\pm$).
  \item\label{ItCMSpacetimeRad}{\rm (Radiative compactification.)} We define
    \[
      \tilde M_1 := [\tilde M_0; Y^+] = [\ol{\R^4}; \fk^+,\fk^-, Y^+],\quad \upbeta_1\colon \tilde M_1\to\ol{\R^4}.
    \]
    Its boundary hypersurfaces are denoted as follows: $\scri^+_1$ is the new front face (i.e., the lift of $Y^+$); $\cK^\pm$ is the lift of $\fk^\pm$ as before; $\iota^+_1$ is the closure of $\{\frac{r}{t}\in(0,1),\ \frac{1}{t}=0\}$; $I^0_1$ is the closure of $\{\frac{t}{r}<1,\ \frac{1}{r}=0\}$. Finally, we define the manifold with corners
    \[
      \tilde M := [\tilde M_1; \scri^+_1, \tfrac12],\quad \upbeta\colon\tilde M\to\ol{\R^4},
    \]
    as the square root blow-up of $\tilde M_1$ at $\scri^+_1$.\footnote{This means that a square root of a defining function of $\scri^+_1$ is a defining function of its lift to $\tilde M$.} Its boundary hypersurfaces are:
    \begin{enumerate}
    \item \emph{null infinity} $\scri^+$, the lift of $\scri^+_1$;
    \item the \emph{Kerr face} $\cK^+$ (and also the past Kerr face $\cK^-$);
    \item \emph{punctured future timelike infinity} $\iota^+$, the lift of $\iota^+_1$;
    \item \emph{spacelike infinity} $I^0$, the lift of $I^0_1$.
    \end{enumerate}
    We denote (local) boundary defining functions of $\scri^+$, $\cK^+\cup\cK^-$, $\iota^+$, and $I^0$ by $x_\sscri$, $\rho_\cK$, $\rho_+$, and $\rho_0\in\CI(\tilde M)$, respectively. (Local) defining functions of $\scri_1^+$ are denoted by $\rho_\sscri$.\footnote{One can thus take $\rho_\sscri=x_\sscri^2$.}
  \item{\rm (Excision of the black hole interior.)} Write $\cl_B(A)$ for the closure of a set $A$ inside of a topological space $B$. Recalling~\eqref{EqCMass}, we then write
    \[
      M_0 := \cl_{\tilde M_0}\{r\geq\bhm\},\quad
      M_1 := \cl_{\tilde M_1}\{r\geq\bhm\},\quad
      M := \cl_M\{r\geq\bhm\}.
    \]
  \end{enumerate}
\end{definition}

The terminology ``spacelike infinity'' is not accurate near $t=-r$, but we shall only be interested in the region $t\geq 0$ for our analysis of forward problems for wave equations (and mostly even just in $t-r\geq 1$, which is disjoint from $I^0$). On a related note, we blow up $\fk^-$ only for notational convenience later on when defining 3b-structures (see the discussion after Definition~\ref{DefCT3b}). We introduce:

\begin{definition}[Domains for waves]
\label{DefCMDomain}
  We define, as submanifolds with corners of $\tilde M$,
  \[
    \tilde\Omega_\sharp := \cl_{\tilde M}\{t\geq 0\},\quad \Omega_\sharp := \tilde\Omega \cap M,\quad
    \Omega := \cl_M\{t-r\geq 1\}.
  \]
\end{definition}

\begin{figure}[!ht]
\centering
\includegraphics{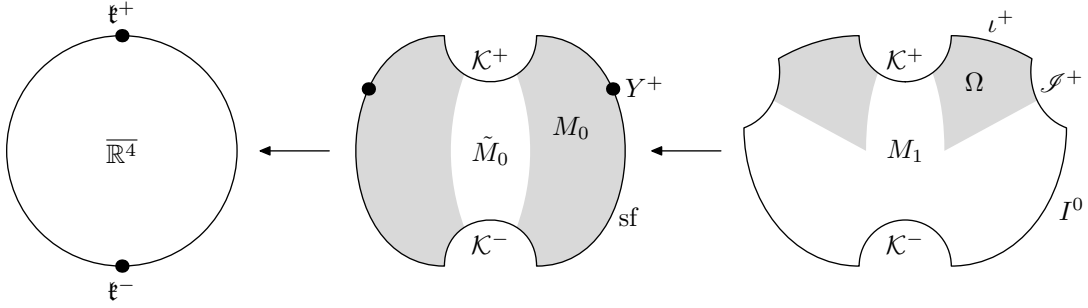}
\caption{The spacetime compactifications from Definition~\ref{DefCMSpacetime} (here shown with spatial dimension $1$ instead of $3$). \textit{On the left:} $\ol{\R^4}$. \textit{In the middle:} $\tilde M_0$ and its submanifold with corners $M_0$ (shaded, connected in $3+1$ dimensions). \textit{On the right:} $\tilde M$ and the domain $\Omega\subset M\subset\tilde M$ from Definition~\ref{DefCMDomain} (shaded, connected in $3+1$ dimensions). --- The arrows are the blow-down maps.}
\label{FigCMSpacetime}
\end{figure}

Coordinate charts covering $\Omega_\sharp\subset M$ are:
\begin{enumerate}[label=(\Alph*)]
\item in $M^\circ$: $t,x$; or $t,r,\omega$; or
  \begin{equation}
  \label{EqCMCoordInt}
    t_*:=t-r\in\R,\quad r\in[\bhm,\infty),\quad \omega\in\Sph^2;
  \end{equation}
\item near the interiors of boundary hypersurfaces:
  \begin{enumerate}
  \item near $(\cK^+)^\circ$: $\frac{1}{t}\geq 0$, $x\in\R^3$ (with $|x|\geq\bhm$); or $\frac{1}{t_*}\geq 0$, $x\in\R^3$;
  \item near $(\iota^+)^\circ$:
    \begin{equation}
    \label{EqCMCoordip}
      v:=\frac{t-r}{r}=\rho t_*\in(0,\infty),\quad \rho_+=\frac{1}{t_*},\quad \omega\in\Sph^2;
    \end{equation}
    or $v,\frac{1}{r},\omega$;
  \item near $(\scri^+)^\circ$: $t-r\in\R$, $r^{-\frac12}\in[0,\infty)$, $\omega\in\Sph^2$;
  \item near $(I^0)^\circ\cap\{t\geq 0\}$: $\frac{r-t}{r}\in(0,1]$, $r^{-1}\in[0,\infty)$, $\omega\in\Sph^2$;
  \end{enumerate}
\item near the corners:
  \begin{enumerate}
  \item near $I^0\cap\scri^+$ (specifically, for $t-r<2$):
    \begin{equation}
    \label{EqCMCoordI0scri}
      \rho_0=\frac{1}{r-t+2},\quad x_\sscri=\sqrt{\frac{r-t+2}{r}},\quad\omega\in\Sph^2;
    \end{equation}
  \item near $\scri^+\cap\iota^+$ (specifically, for $t-r>0$):
    \begin{equation}
    \label{EqCMCoordscriip}
      \rho_+=\frac{1}{t-r},\quad x_\sscri=\sqrt{\frac{t-r}{r}}=\sqrt v,\quad\omega\in\Sph^2;
    \end{equation}
  \item near $\iota^+\cap\cK^+$:
    \begin{equation}
    \label{EqCMCoordipK}
      \rho_+=\frac{1}{r},\quad \rho_\cK=\frac{r}{t},\quad\omega\in\Sph^2.
    \end{equation}
  \end{enumerate}
\end{enumerate}
See Figure~\ref{FigCMSpacetime}. The restriction of the blow-down map $\upbeta$ to $\scri^+$ realizes $\scri^+$ as the total space of a fibration of $\scri^+$ by $\ol\R$,
\begin{equation}
\label{EqCMscriFibr}
  \ol\R - \scri^+ \xra{\upbeta} Y^+=\Sph^2;
\end{equation}
the interiors of the fibers are given by $\{t-r\in\R,\ r^{-\frac12}=0,\ \omega=\omega_0\}$ for fixed $\omega_0\in\Sph^2$.

\begin{rmk}[Explicit defining functions on $\Omega$]
\label{RmkCMExpl}
  Concrete choices of boundary defining functions on $\Omega$ are
  \begin{equation}
  \label{EqCMExpl}
    x_\sscri = \sqrt{\frac{t_*}{t_*+r}},\quad
    \rho_+ = \frac{t_*+r}{t_* r}, \quad
    \rho_\cK = \frac{r}{t_*+r}.
  \end{equation}
  One arrives at these expressions by noting that $t_*$ blows up linearly at $\iota^+$ and $\cK^+$, so $t_*\sim\rho_+^{-1}\rho_\cK^{-1}$ (in that their quotient is a positive smooth function on $\Omega$), similarly $r\sim x_\sscri^{-2}\rho_+^{-1}$ and $t_*+r\sim x_\sscri^{-2}\rho_+^{-1}\rho_\cK^{-1}$, and then solving linear systems of equations for the exponents of $t_*,r,t_*+r$ to isolate $x_\sscri,\rho_+,\rho_\cK$.
\end{rmk}

\subsection{Tangent and cotangent bundles; vector fields}
\label{SsCT}

The manifold $M$ carries a large amount of structure which gives rise to several classes of vector fields and associated tangent, and, dually, cotangent bundles, which will play important roles in our analysis. (This section summarizes part of \cite[\S{2}]{HintzNonstat}.)

\subsubsection{Scattering and b-structures}
\label{SssCTsc}

The Minkowski metric $\ubar g=-\dd t^2+\dd x^2$, having constant coefficients, is, \emph{uniformly} in $z=(t,x)$, a Lorentzian signature quadratic form on the tangent bundle in the frame $\pa_{z^\mu}$, $\mu=0,1,2,3$. Following Melrose \cite{MelroseEuclideanSpectralTheory}, this motivates:

\begin{definition}[Scattering bundles]
\label{DefCTsc}
  The \emph{scattering tangent bundle} $\Tsc\ol{\R^n}\to\ol{\R^n}$ is the trivial rank $n$ bundle $\Tsc\ol{\R^n}=\ol{\R^n}\times\R^n$ which over $\R^n$ is identified with $T\R^n$ via $\Tsc_{\R^n}\ol{\R^n}\ni(z,v)\mapsto v^\mu\pa_{z^\mu}\in T\R^n$. The dual bundle $\Tsc^*\ol{\R^n}\to\ol{\R^n}$ is called the \emph{scattering cotangent bundle}. We denote by $\Vsc(\ol{\R^n}):=\CI(\ol{\R^n};\Tsc\ol{\R^n})$ the space of \emph{scattering vector fields}.
\end{definition}

A scattering vector field is thus a linear combination of $\pa_{z^\mu}$ with coefficients which are of class $\CI(\ol{\R^n})$ (in particular, uniformly bounded). In the case $n=4$, we have
\begin{equation}
\label{EqCTscMink}
  -\dd t^2+\dd x^2 \in \CI(\ol{\R^4};S^2\,\Tsc^*\ol{\R^4}),\quad
  -\pa_t^2+\pa_x^2 \in \CI(\ol{\R^4};S^2\,\Tsc\ol{\R^4}).
\end{equation}
Stationary metrics with non-constant coefficients, such as Kerr, are not smooth sections of $S^2\,\Tsc^*\ol{\R^4}$. To deal with them, we need to allow for more general coefficients than $\CI(\ol{\R^4})$; we do this by working with the pullback of $\Tsc^*\ol{\R^4}$ to resolutions of $\ol{\R^4}$:

\begin{definition}[Pullback bundle]
\label{DefCTscPullback}
  In the notation of Definition~\usref{DefCMSpacetime}, we write
  \[
    \cT \to \tilde M_0,\quad
    \cT \to \tilde M_1,\quad\text{and}\quad
    \cT \to \tilde M
  \]
  for the pullbacks of $\Tsc\ol{\R^4}$ along the blow-down maps $\upbeta_0$, $\upbeta_1$, and $\upbeta$, respectively; and $\cT^*$ denotes the corresponding dual bundles (pullbacks of $\Tsc^*\ol{\R^4})$.
\end{definition}

We shall not distinguish these bundles notationally. The only difference between the spaces of smooth sections of $\cT^*\to\tilde M_0$ and $\cT^*\to\tilde M_1$, say, is that the coefficients of $\pa_{z^\mu}$ in the latter case are smooth functions on $\tilde M_1$, which includes all smooth functions on $\tilde M_0$ but not vice versa. We will realize the Kerr metric as a smooth section of $S^2\cT^*$ over $M_1\subset\tilde M_1$ (and thus also on $M\subset\tilde M$); see~\S\ref{SsTsK}. The bundle $\cT\to\tilde M_0$ had previously been introduced in the context of quantum three-body scattering in Vasy \cite{VasyThreeBody}.

Since $\R^n$ is dense in $\ol{\R^n}$, we have an inclusion map $\Vsc(\ol{\R^n})\hra\cV(\R^n)=\CI(\R^n;T\R^n)$. Local coordinates on $\ol{\R_z^n}$ in the region $\ol{\{z^1\geq 0\}}\setminus\ol{\{z^1=0\}}$ are $\rho_1=\frac{1}{z^1}$, $\hat z^i=\frac{z^i}{z^1}$, $i=2,\ldots,n$, and one computes
\[
  \pa_{z^1} = \rho_1\biggl(-\rho_1\pa_{\rho_1}-\sum_{i=2}^n \hat z^i\pa_{\hat z^i}\biggr),\quad
  \pa_{z^i} = \rho_1\pa_{\hat z^i},\ i=2,\ldots,n.
\]
Therefore, $\rho_1^{-1}\pa_{z^\mu}$ is a smooth vector field in this region which is, moreover, tangent to the boundary $\rho_1=0$. Globally, with $\rho\in\CI(\ol{\R^n})$ denoting a boundary defining function such as $\la z\ra^{-1}$, we thus have
\[
  \Vsc(\ol{\R^n})=\rho\Vb(\ol{\R^n})
\]
where, following Melrose and Mendoza \cite{MelroseMendozaB,MelroseAPS}, we introduce:

\begin{definition}[b-bundles: $\ol{\R^n}$]
\label{DefCTb}
  The \emph{b-tangent bundle} $\Tb\ol{\R^n}\to\ol{\R^n}$ is the trivial rank $n$ bundle $\Tb\ol{\R^n}=\ol{\R^n}\times\R^n$ which over $\R^n$ is identified with $T\R^n$ via $\Tb_{\R^n}\ol{\R^n}\ni(z,v_\bop)\mapsto v^\mu\rho^{-1}\pa_{z^\mu}\in T\R^n$; here $\rho\in\CI(\ol{\R^n})$ is a fixed boundary defining function. The dual bundle $\Tb^*\ol{\R^n}\to\ol{\R^n}$ is called the \emph{b-cotangent bundle}. We denote by $\Vb(\ol{\R^n}):=\CI(\ol{\R^n};\Tb\ol{\R^n})$ the space of \emph{b-vector fields}.
\end{definition}

In local coordinates as above, a smooth section of $\Tb\ol{\R^n}$ over $\ol{\{z^1\geq 0\}}\setminus\ol{\{z^1=0\}}$ is thus spanned by $\rho_1\pa_{\rho_1}$, $\pa_{\hat z^i}$ ($i=2,\ldots,n$). Globally, $\Vb(\ol{\R^n})$ is spanned by vector fields $\pa_{z^\mu}$ and $z^\mu\pa_{z^\nu}$; or, equivalently, by the set of vector fields tangent to $\pa\ol{\R^n}$.

\begin{definition}[b-vector fields in general]
\label{DefCTbGen}
  On a manifold with corners $\tilde M$, we denote by $\Vb(\tilde M)$ the space of all smooth vector fields on $\tilde M$ which are tangent to $\pa\tilde M$, i.e., to each boundary hypersurface of $\tilde M$.
\end{definition}

The strongest notion of regularity for waves $u$ on $\R^4$ which we consider in the present paper is b-regularity on $\tilde M$, i.e., the membership of $u$ in a fixed weighted $L^2$-space upon application of (some number of) b-vector fields on $\tilde M$. Note that the space of b-vector fields is unchanged when passing from $\tilde M_1$ to its square root blow-up $\tilde M$, except one allows for a smooth dependence of the coefficients on $x_\sscri$ (instead of $\rho_\sscri:=x_\sscri^2$, the defining function of $\scri^+_1\subset\tilde M_1$). We thus record:

\begin{example}[b-vector fields on $\tilde M$ from Definition~\ref{DefCMSpacetime}]
\label{ExCTbtildeM}
  We write down spanning sets of b-vector fields (over $\CI(\tilde M)$ or $\CI(\tilde M_1)$ on $\tilde M$ or $\tilde M_1$) in some of the coordinate charts~\eqref{EqCMCoordInt}--\eqref{EqCMCoordipK}. In the coordinates $\rho_\cK=\frac{1}{t}$ and $x\in\R^3$ near $(\cK^+)^\circ$, one can take
  \begin{equation}
  \label{EqCTbtildeM1}
    -\rho_\cK\pa_{\rho_\cK} = t\pa_t,\quad \pa_x\ \text{(or $\pa_r,\pa_\omega$ for $r\geq\bhm$)}.
  \end{equation}
  In the coordinates $\rho_+=\frac{1}{r}$, $\rho_\cK=\frac{r}{t}$, $\omega\in\Sph^2$ from~\eqref{EqCMCoordipK}, so $r=\frac{1}{\rho_+}$, $t=\frac{1}{\rho_+\rho_\cK}$, near $\cK^+\cap\iota^+$, one can take
  \begin{equation}
  \label{EqCTbtildeM2}
    -\rho_+\pa_{\rho_+} = t\pa_t+r\pa_r\ \text{(scaling)},\quad -\rho_\cK\pa_{\rho_\cK}=t\pa_t,\quad\pa_\omega\ \text{(rotations, say)}.
  \end{equation}
  In the coordinates $\rho_+=\frac{1}{t-r}$, $\rho_\sscri=\frac{t-r}{r}$, $\omega\in\Sph^2$ near $\scri_1^+\cap\iota_1^+$ (cf.\ \eqref{EqCMCoordscriip}), so $r=\frac{1}{\rho_+\rho_\sscri}$ and $t=\frac{1}{\rho_+}+\frac{1}{\rho_+\rho_\sscri}$, one can take
  \[
    -\rho_+\pa_{\rho_+} = t\pa_t+r\pa_r\ \text{(scaling)},\quad \rho_\sscri\pa_{\rho_\sscri}-\rho_+\pa_{\rho_+}=(t-r)\pa_t,\quad\pa_\omega.
  \]
  One can replace the scaling vector field here by $-\rho_\sscri\pa_{\rho_\sscri}=r(\pa_t+\pa_r)$; thus b-regularity of a wave captures the usual ``good'' regularity along (but not across) future light cones. The same vector fields also span $\Vb$ near $\scri^+\cap I^0$. A global frame on the domain $\Omega=\cl_{\tilde M}\{r\geq\bhm,\ t_*\geq 1\}$ is, in the coordinates $t_*,x$ (and setting $r=|x|$), given by
  \begin{equation}
  \label{EqCTbtildeMGlobal}
    r\pa_x,\quad t_*\pa_{t_*}.
  \end{equation}
  (This is true also for $\cl_{\tilde M_1}\{r\geq\bhm,\ t_*\geq 1\}\subset\tilde M_1$.) This follows from the earlier expressions by change-of-variables computations. (More conceptually, this follows from Lemma~\ref{LemmaCPXDiffeo} below and the fact that $|x|\pa_x$ and $t_*\pa_{t_*}$ span $\Vb(\ol{\R^3_x})$ and $\Vb(\ol{\R_{t_*}})$ for $|x|>\bhm$ and $t_*>1$, respectively, thus together they span $\Vb(\ol{\R_{t_*}}\times\ol{\R^3_x})$ in this region; and this spanning property persists under blow-ups of boundary faces, i.e., intersections of boundary hypersurfaces, cf.\ \cite[Proposition~5.11.1]{MelroseDiffOnMwc}.)
\end{example}

We make two further remarks. First, $\Vb(\tilde M)$ and $\Vsc(\ol{\R^n})$ are Lie algebras (and Lie subalgebras of $\cV(\R^4)$ and $\cV(\R^n)$, respectively). We can thus consider spaces of differential operators
\begin{equation}
\label{EqCTbDiff}
  \Diffb^m(\tilde M),\quad \Diffsc^m(\ol{\R^n}),
\end{equation}
consisting of up to $m$-fold compositions of elements of $\Vb(\tilde M)$ and $\Vsc(\ol{\R^n})$, respectively, with $0$-fold compositions being multiplication operators with smooth functions on $\tilde M$, resp.\ $\ol{\R^n}$; and the order $m$ here matches the usual order of differential operators on $\tilde M^\circ$ and $\R^n$.

Second, if $\tilde M$ is any manifold with corners and $H\subset\tilde M$ is a compact embedded boundary hypersurface, then there exists a collar neighborhood $[0,\eps)_{\rho_H}\times H$ of $H$ in $\tilde M$, which gives rise to the \emph{b-normal vector field} $\rho_H\pa_{\rho_H}$. A change of variables computation shows that this vector field is, in fact, \emph{independent} of the choice of collar neighborhood modulo elements of $\rho_H\Vb(\tilde M)$ (defined in a neighborhood of $H$).

\subsubsection{3b-structures}
\label{SssCT3b}

Consider the manifold $\tilde M_0=[\ol{\R^4};\fk^+]$ from Definition~\ref{DefCMSpacetime}. The part of its boundary hypersurface $\sface$ where $\frac{r}{t}\in(0,1)$ is the set of endpoints of the timelike curves $r=w t$, $w\in(0,1)$, in Minkowski space. In this regime, every asymptotically flat metric is close to the Minkowski metric, and thus one expects the homogeneity (of degree $-2$) of the latter with respect to spacetime dilations to play an important role. The Minkowski metric can be written as
\[
  \ubar g = -\dd t^2 + \dd r^2 + r^2\slg = r^2 \Bigl[-\Bigl(\frac{\dd t}{r}\Bigr)^2 + \Bigl(\frac{\dd r}{r}\Bigr)^2 + \slg\Bigr],
\]
with $\frac{\dd t}{r}$ and $\frac{\dd r}{r}$ (and spherical 1-forms) being dilation-invariant; the dual vector fields are $r\pa_t$, $r\pa_r$ (and spherical vector fields). Following \cite{Hintz3b}, this suggests introducing:

\begin{definition}[3b-bundles]
\label{DefCT3b}
  The \emph{3b-tangent bundle} $\Ttb\tilde M_0\to\tilde M_0$ is the trivial bundle $\Ttb\tilde M_0=\tilde M_0\times\R^4$ which over $\R^4$ is identified with $T\R^4$ via $\Ttb\tilde M_0\ni(z,v)\mapsto v^\mu\la x\ra\pa_{z^\mu}$, where $z=(t,x)$. (Over $M_0$, one can equivalently, i.e., up to bundle isomorphism, use the identification $(z,v)\mapsto v^\mu\,r\pa_{z^\mu}$.) We write $\Ttb^*\tilde M_0$ for the dual (\emph{3b-cotangent}) bundle, and $\Vtb(\tilde M_0)=\CI(\tilde M_0;\Ttb\tilde M_0)$ for the space of \emph{3b-vector fields}. We define the space $\Difftb^m(\tilde M_0)$ of $m$-th order 3b-differential operators analogously to~\eqref{EqCTbDiff}.
\end{definition}

Note that near $(\cK^+)^\circ$, one can use $\dd t,\dd x$ as a local frame of $\Ttb^*\tilde M_0$; in this fashion, $\Ttb\tilde M_0$ captures the transition from the (asymptotically) translation-invariant regime $(\cK^+)^\circ$ to the (asymptotically) dilation-homogeneous regime $\sface^\circ$. Returning to the Minkowski metric, we thus have
\begin{equation}
\label{EqCT3bMink}
  \ubar g \in r^2\CI(\tilde M_0;S^2\,\Ttb^*\tilde M_0),
\end{equation}
and indeed this is a non-degenerate section in that
\[
  \ubar g^{-1}=r^{-2}(-(r\pa_t)^2+(r\pa_r)^2+\slg^{-1})\in r^{-2}\CI(\tilde M_0;S^2\,\Ttb\tilde M_0).
\]

We have
\[
  \Vtb(\tilde M_0)\subset\Vb(\tilde M_0);
\]
near the corner $\cK^+\cap\sface$, where we can use the coordinates $\rho_\cK=\frac{r}{t}$ and $\rho_+=\frac{1}{r}$, this follows from $r\pa_t=-\rho_\cK^2\pa_{\rho_\cK}$ and $r\pa_r=\rho_\cK\pa_{\rho_\cK}-\rho_+\pa_{\rho_+}$.\footnote{It is here that the blow-up of $\fk^-$ in Definition~\ref{DefCMSpacetime}\eqref{ItCMSpacetimeM0} is convenient: it ensures that $t\mapsto-t$ induces a diffeomorphism of $\tilde M_0$, and thus this computation also applies near $\cK^-\cap\sface$. If we had not blown up $\fk^-$, then typical elements of $\Vtb(\tilde M_0)$ would not be smooth at $\fk^-$. For example, $\la x\ra\pa_t=\frac{\la x\ra}{t}t\pa_t$, with $t\pa_t$ being smooth at $\fk^-$ but $\frac{\la x\ra}{t}=\sqrt{\frac{1}{t^2}+|\frac{x}{t}|^2}$ not.} Since, moreover, $[r\pa_r,r\pa_t]=r\pa_t$ is again a 3b-vector field, it follows easily that $\Vtb(\tilde M_0)$ is a Lie algebra. In the reverse direction, we have
\[
  \rho_\cK\Vb(\tilde M_0)\subset\Vtb(\tilde M_0);
\]
near the corner, this follows for $\rho_\cK=\frac{r}{t}$ from $r\pa_r,\pa_\omega\in\Vtb(\tilde M_0)$ (no need for the extra factor $\rho_\cK$) and $\frac{r}{t}t\pa_t=r\pa_t\in\Vtb(\tilde M_0)$. More generally, for $V\in\Vb(\tilde M_0)$, we have $V\in\Vtb(\tilde M_0)$ if and only if the function $V t$ on $\{t>0\}$, which for general $V\in\Vb(\tilde M_0)$ lies in $t\CI(\tilde M_0\setminus\ol{\{t\leq 0\}})$, lies in the smaller space $t\rho_\cK\CI$. Indeed, this condition restricts the coefficient of $t\pa_t$ to lie in $\rho_\cK\CI$; but $\rho_\cK t\pa_t$ is a smooth multiple of $\la x\ra\pa_t$. (See also \cite[Lemma~3.4(4)]{Hintz3b}.)

These computations also show that the identity map on $T\R^4$ induces a bundle isomorphism\footnote{Both bundles are trivial, but the point is that this isomorphism is \emph{canonical} in that it is the unique smooth extension of the identity map on $T\R^4$.}
\[
  \Ttb_{\tilde M_0\setminus\cK^+}\tilde M_0 = \Tb_{\tilde M_0\setminus\cK^+}\tilde M_0.
\]
Thus, 3b-vector fields are the same as b-vector fields away from $\cK^+$. Moreover, since $\tilde M_0\setminus Y^+=\tilde M_1\setminus\scri^+_1=\tilde M\setminus\scri^+$ (in the sense that the identity map on $\R^4$ induces diffeomorphisms of these manifolds), this implies
\begin{equation}
\label{EqCT3bIsob}
  \Ttb_{\tilde M\setminus(\cK^+\cup\scri^+)}\tilde M = \Tb_{\tilde M\setminus(\cK^+\cup\scri^+)}\tilde M.
\end{equation}

The following result will be important for the proof of the microlocal e3b-regularity, relative to finite order b-regularity, of waves near $\cK^+$ (especially near the trapped set of Kerr):

\begin{lemma}[Commutators with time dilations]
\label{LemmaCT3bDil}
  Let $V\in\Vb(\tilde M_0)$ be such that $V-t\pa_t\in\rho_\cK\Vb(\tilde M_0)$. Let $W\in\Vtb(\tilde M_0)$. Then there exist $w^\flat\in\CI(\tilde M_0)$ and $W^\sharp\in\Vtb(\tilde M_0)$ such that
  \begin{equation}
  \label{EqCT3bDil}
    [V,W] = \rho_\cK w^\flat V + \rho_\cK W^\sharp.
  \end{equation}
  In particular, $[V,W]\in\Vtb(\tilde M_0)$. More generally, if $W\in\rho_\sface^{-\alpha_\sface}\rho_\cK^{-\alpha_\cK}\Difftb^m(\tilde M_0)$, then $[V,W]\in\rho_\sface^{-\alpha_\sface}\rho_\cK^{-\alpha_\cK}\Difftb^m(\tilde M_0)$; and in the special case $W\in\rho_\sface^{-\alpha_\sface}\Difftb^m(\tilde M_0)$, we can write
  \begin{equation}
  \label{EqCT3bDilOp}
    [V,W] = \rho_\cK W^\flat V + \rho_\cK W^\sharp,\quad W^\flat\in\rho_\sface^{-\alpha_\sface}\Difftb^{m-1}(\tilde M_0),\ W^\sharp\in\Difftb^m(\tilde M_0).
  \end{equation}
  The same results hold true for operators acting on sections of a smooth vector bundle $\cE\to\tilde M_0$, except now we require $V\in\Diffb^1(\tilde M_0;\cE)$ and $V-t\pa_t\in\rho_\cK\Diffb^1(\tilde M_0;\cE)$ where $t\pa_t$ is defined as $\nabla_{t\pa_t}$ for any smooth affine connection $\nabla$ on $\cE$, and $W^\flat,W^\sharp$ act on sections of $\cE$ as well.
\end{lemma}

The condition on $W\in\Vtb(\tilde M_0)$ is equivalent to
\begin{equation}
\label{EqCT3bDilEquiv}
  V\equiv -\rho_\cK\pa_{\rho_\cK}\bmod\rho_\cK\Vb(\tilde M_0)
\end{equation}
in any collar neighborhood $[0,\eps)_{\rho_\cK}\times\cK^+$ of $\cK^+$ in $\tilde M_0$; this follows from the fact (see~\eqref{EqCTbtildeM1}--\eqref{EqCTbtildeM2}) that $-t\pa_t$ is a b-normal vector field at $\cK^+$. One possible choice of $W$ on $M\subset\tilde M$ is thus $t_*\pa_{t_*}$ in the coordinates $(t_*,r)=(t-r,r)$ since this equals $t\pa_t-r\pa_t=(1-\rho_\cK)t\pa_t$, $\rho_\cK=\frac{r}{t}$.

\begin{proof}[Proof of Lemma~\usref{LemmaCT3bDil}]
  This is the content of \cite[Lemmas~2.2 and 2.4]{HintzNonstat}; we give a slightly different proof for the convenience of the reader. Consider first the case $W\in\Vtb(\tilde M_0)$. We only give the computations near the corner $\cK^+\cap\sface$, where we take $\rho_\cK=\frac{r}{t}$. If the conclusion holds for $W$, and if $f\in\CI(\tilde M_0)$, consider
  \[
    [V,f W] = f[V,W] + (V f)W.
  \]
  Expressing $[V,W]$ as in~\eqref{EqCT3bDil} shows that $f[V,W]$ is also of the form~\eqref{EqCT3bDil}; moreover, $V f\in\rho_\cK\CI(\tilde M_0)$ (since $V\in\rho_\cK\cV(\tilde M_0)$ by~\eqref{EqCT3bDilEquiv}), so $(V f)W\in\rho_\cK\Vtb(\tilde M_0)$ is of the same form as the second term in~\eqref{EqCT3bDil}. It thus suffices to check~\eqref{EqCT3bDil} for $W=r\pa_t,r\pa_r,\pa_\omega$.

  We next compute
  \[
    [t\pa_t, r\pa_t] = -r\pa_t = -\rho_\cK t\pa_t,\quad
    [t\pa_t, r\pa_r] = 0,\quad
    [t\pa_t, \pa_\omega] = 0,
  \]
  which proves~\eqref{EqCT3bDil} for $V=t\pa_t$ (with $W^\sharp=0$ even). Finally, for $V=t\pa_t+\rho_\cK V^\flat$, we can write $V^\flat=a t\pa_t+V'$ where $a\in\CI(\tilde M_0)$ and $V'$, lying in the $\CI(\tilde M_0)$-span of $r\pa_r,\pa_\omega$, is an element of $\Vtb(\tilde M_0)$. Consider then
  \begin{align*}
    [V,W] &= [t\pa_t,W] + [\rho_\cK a t\pa_t,W] + [\rho_\cK V',W] \\
      &= (1+\rho_\cK a)[t\pa_t,W] - W(\rho_\cK a)t\pa_t + [\rho_\cK V',W].
  \end{align*}
  The third term lies in $\rho_\cK\Vtb(\tilde M_0)$. The second term is of the form $\rho_\cK f t\pa_t$ for $f\in\CI(\tilde M_0)$ and thus equal to $\rho_\cK f V$ modulo $\rho_\cK^2\Vb(\tilde M_0)\subset\rho_\cK\Vtb(\tilde M_0)$. The first term, using~\eqref{EqCT3bDilEquiv} for $t\pa_t$ in place of $V$, can be written as $\rho_\cK w^\flat t\pa_t+\rho_\cK W^\sharp$, and replacing $t\pa_t$ by $V$ in the first term here generates another error of class $\rho_\cK^2\Vb(\tilde M_0)\subset\rho_\cK\Vtb(\tilde M_0)$.

  Turning to the case of differential operators, write $W=w W_0$ where $w=\rho_\sface^{-\alpha_\sface}\rho_\cK^{-\alpha_\cK}$, $W_0\in\Difftb^m(\tilde M_0)$. Then the membership $[V,W_0]\in\Difftb^m(\tilde M_0)$ follows from what we have already shown by writing $W_0$ as the sum of compositions of 3b-vector fields. Since $[V,w]\in w\CI(\tilde M_0)$ (this being true more generally for any b-vector field $V$ on $\tilde M_0$), we get $[V,W]=w[V,W_0]+[V,w]W_0\in w\Diff_\tbop^m(\tilde M_0)$. In the case $\alpha_\cK=0$, let us take $\rho_\sface=\la x\ra^{-1}$ and $V=t\pa_t$. (Additional terms of class $\rho_\cK\Vb(\tilde M_0)$ change the commutator $[V,W]$ by terms in $[\rho_\cK\Vb,w\Diff_\tbop^m]\subset\rho_\cK w\Diff_\tbop^m$, which can be put into $W^\sharp$ in~\eqref{EqCT3bDilOp}.) Since $V$ then commutes with $\rho_\sface$, it suffices to consider the case $\alpha_\sface=0$, i.e., $W\in\Diff_\tbop^m(\tilde M_0)$. The desired conclusion~\eqref{EqCT3bDilOp} now follows by expanding $W$ into products of 3b-vector fields and commuting $V$ through to the rightmost place; for zeroth order terms of $W$, i.e., multiplication operators by $f\in\CI(\tilde M_0)$, we note that $[V,f]=V(f)\in\rho_\cK\CI(\tilde M_0)$ by~\eqref{EqCT3bDilEquiv}, which thus contributes to $W^\sharp$ in~\eqref{EqCT3bDilOp}.

  In the case of bundles, finally, we again only work near $\cK^+\cap\sface$ and note that, in terms of $\rho_\cK=\frac{r}{t}$ and $\rho_\sface=\frac{1}{r}$, we have $t\pa_t=-\rho_\cK\pa_{\rho_\cK}$ and thus $\nabla_{t\pa_t}=-\rho_\cK\nabla_{\pa_{\rho_\cK}}$. In a local trivialization of $\cE$, the operator $\nabla_{t\pa_t}$ is thus given by $t\pa_t$ acting component-wise, plus a zeroth order term of class $\rho_\cK\CI(\tilde M_0;\End(\cE))$. Commutators with the first summand were already discussed, and the commutator of the second summand with an element of $\Difftb^m(\tilde M_0;\cE)$ lie in $\rho_\cK\Difftb^{m-1}(\tilde M_0;\cE)$ and thus can be absorbed into the term $\rho_\cK W^\sharp$ in~\eqref{EqCT3bDilOp}, similarly in the case of weighted operators.
\end{proof}

\begin{cor}[Enough 3b-commutator b-vector fields]
\label{CorCT3bComm}
  The set of \emph{3b-commutator b-vector fields}, i.e., $W\in\Vb(\tilde M_0)$ such that $[V,W]\in\Vtb(\tilde M_0)$ for all $W\in\Vtb(\tilde M_0)$, spans $\Vb(\tilde M_0)$ over $\CI(\tilde M_0)$.
\end{cor}
\begin{proof}
  This is nontrivial only near $\cK^+$. But by Lemma~\ref{LemmaCT3bDil}, $t\pa_t$ is such a vector field; and the 3b-vector fields $r\pa_r$ and $\pa_\omega$ (spherical vector fields), or more globally $\la x\ra\pa_{z^\mu}$, are trivially 3b-commutator b-vector fields. The span of $t\pa_t$ and $\la x\ra\pa_{z^\mu}$ over $\CI(U)$ is equal to $\Vb(U)$ when $U\subset\tilde M_0$ is a small neighborhood of $\cK^+$, as follows from \eqref{EqCTbtildeM1}--\eqref{EqCTbtildeM2}.
\end{proof}

\subsubsection{edge-b-structures}
\label{SssCTeb}

We now consider the structure of the manifold $\tilde M$ in Definition~\ref{DefCMSpacetime} near $\scri^+$. Here, the homogeneity of the Minkowski metric, or more non-trivially of an asymptotically spacetime featuring nontrivial radiation at $\scri^+$, with respect to scalings when $\frac{r}{t}\in(0,1)$ and $r\to\infty$ gives way to a different structure. Using $t_*=t-r$, we first write the Minkowski metric as
\begin{equation}
\label{EqCTebMinkTstar}
  \ubar g = -\dd t^2 + \dd r^2 + r^2\slg = -\dd t_*^2 - 2\,\dd t_*\,\dd r + r^2\slg,\quad
  \ubar g^{-1} = -2\pa_{t_*}\otimes_s\pa_r + \pa_r^2 + r^{-2}\slg^{-1}.
\end{equation}
Passing then to the local boundary defining functions $\rho_+=\frac{1}{t_*}$ and $x_\sscri=\sqrt{\frac{t_*}{r}}$ of $\iota^+$ and $\scri^+\subset\tilde M$ near $\scri^+\cap\iota^+$ from~\eqref{EqCMCoordscriip}, so $t_*=\frac{1}{\rho_+}$ and $r=\frac{1}{\rho_+ x_\sscri^2}$, we compute
\begin{equation}
\label{EqCTebConv}
\begin{alignedat}{2}
  \dd t_*&=-\rho_+^{-1}\frac{\dd\rho_+}{\rho_+},&\quad
  \dd r &= -\rho_+^{-1}x_\sscri^{-2}\Bigl(\frac{\dd\rho_+}{\rho_+} + 2\frac{\dd x_\sscri}{x_\sscri}\Bigr), \\
  \pa_{t_*} &= \rho_+\Bigl(-\rho_+\pa_{\rho_+} + \frac12 x_\sscri\pa_{x_\sscri}\Bigr), &\quad
  \pa_r &= -\frac12\rho_+ x_\sscri^2\cdot x_\sscri\pa_{x_\sscri}.
\end{alignedat}
\end{equation}
Therefore, we have
\begin{equation}
\label{EqCTebMink}
\begin{split}
  \ubar g &= 2 x_\sscri^{-2}\rho_+^{-2}\Bigl( -\frac{\dd\rho_+}{\rho_+} \otimes_s \Bigl[ \Bigl(1+\frac12 x_\sscri^2\Bigr)\frac{\dd\rho_+}{\rho_+} + 2\frac{\dd x_\sscri}{x_\sscri}\Bigr] + \frac{\slg}{x_\sscri^2}\Bigr) \\
  \ubar g^{-1} &= \frac12 x_\sscri^2\rho_+^2\Bigl( x_\sscri\pa_{x_\sscri}\otimes_s \Bigl(-2\rho_+\pa_{\rho_+} + \Bigl(1+\frac12 x_\sscri^2\Bigr)x_\sscri\pa_{x_\sscri}\Bigr) + x_\sscri^2\slg^{-1} \Bigr),
\end{split}
\end{equation}
with similar expressions near the other corner $\scri^+\cap I^0$ of $\scri^+$ (cf.\ \cite[(3.7b), (3.8b)]{HintzVasyScrieb}). Thus, $x_\sscri^{-2}\rho_+^{-2}\ubar g^{-1}$ is a Lorentzian signature quadratic form in $\rho_+\pa_{\rho_+}$, $x_\sscri\pa_{x_\sscri}$, and $x_\sscri\pa_\omega$ (rescaled spherical vector fields); their $\CI(\tilde M)$-span is $\Veb(\tilde M)$, where we introduce:

\begin{definition}[edge-b-bundles]
\label{DefCTeb}
  We denote by $\Veb(\tilde M)$ the space of all b-vector fields on $\tilde M$ that, moreover, are tangent to the fibers of the fibration~\eqref{EqCMscriFibr} of $\scri^+$; these are called \emph{edge-b-vector fields}. By $\Teb\tilde M\to\tilde M$, we denote the vector bundle, identified with $T\R^4$ over $\R^4=\tilde M^\circ$, such that $\Veb(\tilde M)=\CI(\tilde M;\Teb\tilde M)$;\footnote{Here we identify a section of $\Teb\tilde M$ with a section of $T\R^4$ over $\R^4$, i.e., a vector field.} it is called the \emph{edge-b-tangent bundle}. The dual \emph{edge-b-cotangent bundle} is denoted $\Teb^*\tilde M\to\tilde M$. We define the space $\Diffeb^m(\tilde M)$ of $m$-th order edge-b-differential operators in the usual fashion, cf.\ \eqref{EqCTbDiff}.
\end{definition}

Edge-bundles were introduced by Mazzeo \cite{MazzeoEdge}, and the edge-b-combination first appeared in \cite{MelroseVasyWunschDiffraction}; its utility near null infinity was first realized in \cite{HintzVasyScrieb}. In the coordinates $\rho_+,x_\sscri$ as above, and in local coordinates $\omega\in\R^2$ on $\Sph^2$, a local frame of $\Teb\tilde M$, resp.\ $\Teb^*\tilde M$, is given by
\begin{equation}
\label{EqCTebFrame}
  \rho_+\pa_{\rho_+},\ x_\sscri\pa_{x_\sscri},\ x_\sscri\pa_\omega,\quad\text{resp.}\quad
  \frac{\dd\rho_+}{\rho_+},\ \frac{\dd x_\sscri}{x_\sscri},\ \frac{\dd\omega}{x_\sscri}.
\end{equation}
We thus have
\begin{equation}
\label{EqCTebMink2}
  \ubar g \in \rho_0^{-2}x_\sscri^{-2}\rho_+^{-2}\CI(\tilde M\setminus\cK^+;S^2\,\Teb^*\tilde M),\quad
  \ubar g^{-1} \in \rho_0^2 x_\sscri^2\rho_+^2\CI(\tilde M\setminus\cK^+;S^2\,\Teb\tilde M).
\end{equation}
Moreover, the identity map on $T\R^4$ induces a bundle isomorphism
\begin{equation}
\label{EqCTebIsob}
  \Teb_{\tilde M\setminus\scri^+}\tilde M = \Tb_{\tilde M\setminus\scri^+}\tilde M.
\end{equation}

The edge-b-setting is useful for analytic purposes (i.e., regarding wave-type operators as weighted edge-b-operators), whereas it is always the scattering setting (cf.\ \eqref{EqCTscMink} and Definition~\ref{DefCTscPullback}) that is most useful for geometric purposes. We thus record (cf.\ \cite[Lemmas~3.24 and 3.25]{HintzVasyScrieb} and \cite[Lemma~3.13]{HintzMink4Gauge}):

\begin{lemma}[Symmetric 2-tensors: edge-b vs.\ scattering]
\label{LemmaCTebsc}
  Schematically write $\dd\omega$ for a 1-form on $\Sph^2$. Then a local spanning set of $\rho_0^{-2}x_\sscri^{-2}\rho_+^{-2}\CI(\tilde M;S^2\,\Teb^*\tilde M)$ near $\scri^+$ is given by
  \begin{equation}
  \label{EqCTebscStar}
  \begin{alignedat}{3}
    &x_\sscri^{-2}\,\dd t_*^2, &\quad & \dd t_*\otimes_s\dd r, &\quad &x_\sscri^{-1}\,\dd t_*\otimes_s r\,\dd\omega, \\
    &x_\sscri^2\,\dd r^2, &\quad & x_\sscri\,\dd r\otimes_s r\,\dd\omega, &\quad &r\,\dd\omega\otimes_s\,r\,\dd\omega.
  \end{alignedat}
  \end{equation}
  Setting $x^0=t+r$ and $x^1=t-r$, another local spanning set is
  \begin{equation}
  \label{EqCTebsc01}
  \begin{alignedat}{3}
    &x_\sscri^2(\dd x^0)^2, &\quad & \dd x^0\otimes_s\dd x^1, &\quad &x_\sscri\,\dd x_0\otimes_s r\,\dd\omega, \\
    &x_\sscri^{-2}(\dd x^1)^2, &\quad & x_\sscri^{-1}\,\dd x^1\otimes_s r\,\dd\omega, &\quad &r\,\dd\omega\otimes_s\,r\,\dd\omega.
  \end{alignedat}
  \end{equation}
\end{lemma}
\begin{proof}
Let us write $U\subset\tilde M$ for a neighborhood of $\scri^+$ which is disjoint from $\cK^+$ and $\{r=0\}$. Write $\CI=\CI(U;\Teb^*\tilde M)$. Then by~\eqref{EqCTebConv}, we have
  \[
    x_\sscri^{-1}\,\dd t_*,\ x_\sscri\,\dd r,\ r\,\dd\omega \in \rho_+^{-1}x_\sscri^{-1}\CI(U;\Teb^*\tilde M).
  \]
  Similarly, using $\dd x^0=\dd t_*+2\,\dd r$, we have
  \begin{align*}
    x_\sscri\,\dd x^0 &= -\rho_+^{-1}x_\sscri^{-1}\Bigl( (1+2 x_\sscri^2)\frac{\dd\rho_+}{\rho_+} + 2\frac{\dd x_\sscri}{x_\sscri}\Bigr), \\
    x_\sscri^{-1}\,\dd x^1 &= -\rho_+^{-1}x_\sscri^{-1}\frac{\dd\rho_+}{\rho_+}, \\
    r\,\dd\omega &= \rho_+^{-1}x_\sscri^{-1}\cdot x_\sscri^{-1}\,\dd\omega.
  \end{align*}
  Their $\CI(U)$-span is thus equal to $x_\sscri^{-1}\rho_+^{-1}\CI(U;\Teb^*\tilde M)$. Taking symmetric tensor products gives~\eqref{EqCTebscStar}--\eqref{EqCTebsc01}.
\end{proof}

For the proof of b-regularity of waves near $\scri^+$, we have the following analogue of Lemma~\ref{LemmaCT3bDil} (taken from \cite[\S{5.1}]{HintzVasyScrieb}, with a concrete version being indicated in \cite[Lemma~3.16]{HintzMink4Gauge}):

\begin{lemma}[Commutators with suitable b-vector fields]
\label{LemmaCTebComm}
  Let $W\in\Veb(\tilde M)$. Then, in the coordinates $\rho_+,x_\sscri,\omega$ in a neighborhood $U$ of $\scri^+\cap\iota^+$ from~\eqref{EqCMCoordscriip}, we have
  \begin{equation}
  \label{EqCTebComm}
    [ W,\rho_+\pa_{\rho_+} ],\ [ W,x_\sscri\pa_{x_\sscri} ],\ [ W,\pa_\omega ] \in \Veb(U).
  \end{equation}
  More generally, if $V\in\Vb(\tilde M)$ is such that for each $\omega\in\Sph^2$, the pushforward $\upbeta_*(V|_{(t_*,0,\omega)})\in T_\omega\Sph^2$ of $V$ from the point on $\scri^+$ with $t-r=t_*$, $r^{-\frac12}=0$, $\omega$ is independent of $t_*$, then $[V,W]\in\Veb(\tilde M)$ for all $W\in\Veb(\tilde M)$. For any such $V$, and for any $W\in\Diffeb^m(\tilde M)$, we have $[V,W]\in\Diffeb^m(\tilde M)$, similarly when $W$ is a weighted eb-differential operator.
\end{lemma}
\begin{proof}
  Since $\Veb$ is a Lie algebra and $\rho_+\pa_{\rho_+},x_\sscri\pa_{x_\sscri}\in\Veb(U)$, only the membership $[\pa_\omega,W]\in\Veb(U)$ requires an argument. But this is true for $W=\rho_+\pa_{\rho_+}$, $x_\sscri\pa_{x_\sscri}$, and $W=x_\sscri W'$ where $W'$ is a vector field on $\Sph^2$; and since for $f\in\CI$ we have $[f W,\pa_\omega]=f[W,\pa_\omega]-(\pa_\omega f)W$ with $\pa_\omega f\in\CI$, it is true for general $W\in\Veb(\tilde M)$.

  To prove the second part, note that every $V\in\Vb(\tilde M)$ as in the statement can locally be written as $a\rho_+\pa_{\rho_+}+b x_\sscri\pa_{x_\sscri}+c\pa_\omega$ for smooth $a,b,c$ where, moreover, $c|_{\scri^+}$ is constant on the fibers of $\scri^+$. Correspondingly expanding the commutator of $V$ with $W\in\Veb(\tilde M)$ into three summands, the first two summands clearly lie in $\Veb(\tilde M)$, and for the third summand we use that $[W,c\pa_\omega]=(W c)\pa_\omega+c[W,\pa_\omega]$; the second summand here lies in $\Veb(\tilde M)$ by what we have already shown; and so does the first summand since $W c|_{\scri^+}=0$ for $W\in\Veb(\tilde M)$, and hence $W c\in x_\sscri\CI(\tilde M)$.

  The final statement follows by expanding $W$ into products of eb-vector fields.
\end{proof}

\begin{cor}[Enough eb-commutator b-vector fields]
\label{CorCTebComm}
  The set $\cV_{\bop,[\ebop]}(\tilde M)$ of b-vector fields $W$ on $\tilde M$ such that $[V,W]\in\Veb(\tilde M)$ for all $V\in\Veb(\tilde M)$ spans $\Vb(\tilde M)$ over $\CI(\tilde M)$.
\end{cor}
\begin{proof}
  This is only nontrivial near $\scri^+$, where however the vector fields in~\eqref{EqCTebComm} do span $\Vb(\tilde M)$.
\end{proof}

\subsubsection{Combination: e3b-structures}
\label{SssCTe3b}

Since the 3b-(co)tangent bundle on $\tilde M_0$ and the edge-b-(co)tan\-gent bundle on $\tilde M$ agree on $\tilde M\setminus(\cK^+\cup\scri^+)=\tilde M_0\setminus(\fk^+\cup Y^+)$ by~\eqref{EqCT3bIsob} and \eqref{EqCTebIsob}, we can glue them together to an edge-b-3b-(co)tangent bundle. We shall focus our attention mainly on the region $\Omega\subset M$, which is disjoint from $I^0$, and thus the `b'-behavior occurs only at $\iota^+$. This suggests the terminology in the following definition:

\begin{definition}[e3b-bundles]
\label{DefCTe3b}
  The \emph{e3b-tangent bundle} $\Tetb\tilde M\to\tilde M$ is equal to $\Ttb\tilde M_0$ over $\tilde M\setminus\scri^+=\tilde M_0\setminus Y^+$ and equal to $\Teb\tilde M$ over $\tilde M\setminus\cK^+$, with the two bundles identified over $\tilde M\setminus(\cK^+\cup\scri^+)$ via~\eqref{EqCT3bIsob} and \eqref{EqCTebIsob} (i.e., via the extension of the identity map on $T\R^4$). The \emph{e3b-cotangent bundle} is the dual bundle $\Tetb^*\tilde M\to\tilde M$. We write $\Vetb(\tilde M)=\CI(\tilde M;\Tetb\tilde M)$, and $\Diff_\etbop^m(\tilde M)$ for the space of $m$-th order e3b-differential operators.
\end{definition}

Equivalently, $\Tetb\tilde M$ is the vector bundle, identified over $\R^4$ with $T\R^4$, such that the space of its smooth sections is equal to the space $\Vetb(\tilde M)$ of all b-vector fields $V$ on $\tilde M$ which are tangent to the fibers of $\scri^+$ and such that $V t\in t\rho_\cK\CI$ near $\cK^+$. Combining~\eqref{EqCT3bMink} and \eqref{EqCTebMink2}, the Minkowski metric $\ubar g$ thus satisfies, in the domain $\Omega$ from Definition~\ref{DefCMDomain}:
\[
  \ubar g \in \rho_0^{-2}x_\sscri^{-2}\rho_+^{-2}\CI(\Omega;S^2\,\Tetb^*\tilde M),\quad
  \ubar g^{-1} \in \rho_0^2 x_\sscri^2\rho_+^2\CI(\Omega;S^2\,\Tetb\tilde M).
\]
Similarly to \cite{HintzNonstat}, we will analyze wave-type operators on asymptotically Kerr spacetimes microlocally in $\Tetb^*\tilde M$ (see~\S\ref{SR}).

Finally, the combination of Corollaries~\ref{CorCT3bComm} and \ref{CorCTebComm} gives:

\begin{lemma}[Enough e3b-commutator b-vector fields]
\label{LemmaCTe3bComm}
  $\Vb(\tilde M)$ is spanned over $\CI(\tilde M)$ by the set $\cV_{\bop,[\etbop]}(\tilde M)$ of b-vector fields $W$ such that $[V,W]\in\Vetb(\tilde M)$ for all $W\in\Vetb(\tilde M)$. More precisely, there exists a \emph{good spanning set} of $\Vb(\tilde M)$, by which we mean a spanning set
  \[
    \sV = \{V_1,\ldots,V_N\}\subset\cV_{\bop,[\etbop]}(\tilde M)
  \]
  of $\Vb(\tilde M)$ such that for each $j=1,\ldots,N$ there exists a (bounded) weight $w_j\in\CI(\tilde M)$ such that $V_{0,j}:=w_j V_j\in\Vetb(\tilde M)$, and such that $\{V_{0,1},\ldots,V_{0,N}\}$ spans $\Vetb(\tilde M)$ over $\CI(\tilde M)$.
\end{lemma}

The point of a good spanning set is that while every b-vector field $V_j$ is a weighted e3b-vector field, the weight $w_j$ with the property that $V_{0,j}=w_j V_j$ is an unweighted e3b-vector field may need to vanish to fairly high orders at the boundary hypersurfaces of $\tilde M$, thus preventing the vector fields $V_{0,j}$ from spanning $\Vetb(\tilde M)$ due to too many orders of vanishing at certain boundary hypersurfaces. Locally near $\cK^+$, $\{t\pa_t,\,r\pa_r,\,\pa_\omega\}$ is a good spanning set (with $w_1=\frac{r}{t}$, so $V_{0,1}=r\pa_t$), whereas $\{t\pa_t,\,t\pa_t+r\pa_r,\,\pa_\omega\}$ is not (since we need to take, at least, $w_1=\frac{r}{t}$ and $w_2=\frac{r}{t}$, so $V_{0,1}=r\pa_t$ and $V_{0,2}=r\pa_t+r\rho_\cK\pa_r$, $\rho_\cK=\frac{r}{t}$, and thus $r\pa_r$ does not lie in their $\CI$-span).

\begin{proof}[Proof of Lemma~\usref{LemmaCTe3bComm}]
  We only comment on the existence of a \emph{good spanning set}: we may simply take $\sV$ to consist of a finite number of e3b-vector fields together with $\chi_\cK t\pa_t=\rho_\cK^{-1}\cdot\chi_\cK r\pa_t$ where $\chi_\cK\in\CI(\tilde M)$ localizes to a neighborhood of $\cK$ and $\rho_\cK=\frac{r}{t}$, and (schematically) $\chi_\sscri\pa_\omega=x_\sscri^{-1}\cdot\chi_\sscri x_\sscri\pa_\omega$ where $\chi_\sscri\in\CI(\tilde M)$ localizes to a neighborhood of $\scri^+$ and $x_\sscri$ is a defining function of $\scri^+$.
\end{proof}

\subsubsection{Hamiltonian vector fields}
\label{SssCHam}

On a smooth manifold $\cM$ without boundary, consider a smooth function $p\in\CI(T^*\cM)$. Its Hamiltonian vector field $H_p\in\cV(T^*\cM)$ is defined in local coordinates $x\in\R^n$ on $\cM$ and the associated canonical momentum variables $\xi\in\R^n$ as
\[
  H_p = (\pa_\xi p)\cdot\pa_x - (\pa_x p)\cdot\pa_\xi.
\]
Consider now a smooth function $p\in\CI(\Ttb^*\tilde M_0)$ in the notation of Definition~\ref{DefCMSpacetime}; in local coordinates $\rho_\cK=\frac{r}{t_*}$, $\rho_+=\frac{1}{r}$, $\omega\in\Sph^2$ near $\cK^+\cap\iota^+$ and writing 3b-covectors as
\begin{equation}
\label{EqCHam3bCoord}
  \sigma_\tbop\frac{\dd t_*}{r}+\xi_\tbop\frac{\dd r}{r}+\eta_\tbop,\quad \eta_\tbop\in T^*\Sph^2,
\end{equation}
we then have $p=p(\rho_\cK,\rho_+,\omega;\sigma_\tbop,\xi_\tbop,\eta_\tbop)$. A change of variables computation (i.e., passing to coordinates defined via $\sigma_\tbop\frac{\dd t_*}{r}+\xi_\tbop\frac{\dd r}{r}+\eta_\tbop=-\sigma\,\dd t_*+\xi\,\dd r+\eta$) shows that
\begin{equation}
\label{EqCHam3b}
\begin{split}
  H_p &= (\pa_{\sigma_\tbop}p) r\pa_{t_*} + (\pa_{\xi_\tbop}p)(r\pa_r+\sigma_\tbop\pa_{\sigma_\tbop}) + (\pa_{\eta_\tbop}p)\pa_\omega \\
      &\quad - (r\pa_{t_*}p)\pa_{\sigma_\tbop} - \bigl((r\pa_r+\sigma_\tbop\pa_{\sigma_\tbop})p\bigr)\pa_{\xi_\tbop} - (\pa_\omega p)\pa_{\eta_\tbop} \\
      &= -(\pa_{\sigma_\tbop}p)\rho_\cK^2\pa_{\rho_\cK} + (\pa_{\xi_\tbop}p)( \rho_\cK\pa_{\rho_\cK}-\rho_+\pa_{\rho_+} + \sigma_\tbop\pa_{\sigma_\tbop} ) + (\pa_{\eta_\tbop}p)\pa_\omega \\
      &\quad + (\rho_\cK^2\pa_{\rho_\cK}p)\pa_{\sigma_\tbop} - \bigl((\rho_\cK\pa_{\rho_\cK}-\rho_+\pa_{\rho_+}+\sigma_\tbop\pa_{\sigma_\tbop})p\bigr)\pa_{\xi_\tbop} - (\pa_\omega p)\pa_{\eta_\tbop}.
\end{split}
\end{equation}
This is therefore a 3b-vector field on $\Ttb^*\tilde M_0$, i.e., the pushforward of $H_p|_\varpi$, for any $\varpi\in\Ttb^*\tilde M_0$, along the base projection lies in $\Ttb\tilde M_0$. We write this as
\[
  H_p \in \Vtb(\Ttb^*\tilde M_0);
\]
there is a bundle $\Ttb(\Ttb^*\tilde M_0)\to\Ttb^*\tilde M_0$ whose smooth sections are precisely $\Vtb(\Ttb^*\tilde M_0)$. (In local coordinates as above, a local frame for this bundle is $\pa_{\sigma_\tbop},\pa_{\xi_\tbop},\pa_{\eta_\tbop}$, together with the basic 3b-vector fields on $\tilde M_0$.) Moreover, the formula~\eqref{EqCHam3b} shows that the map $p\mapsto H_p$ is a 3b-differential operator:
\begin{equation}
\label{EqCHam3bMap}
  (p\mapsto H_p) \in \Difftb^1(\Ttb^*\tilde M_0;\ubar\C,\Ttb(\Ttb^*\tilde M_0)),
\end{equation}
where $\ubar\C$ is the trivial bundle. The structural reason is that $H_p$ is defined as $\omega(H_p,\cdot)=\dd p$ where $\omega=-\dd\alpha$ is the symplectic form, given as the differential of the canonical 1-form $\alpha_\varpi(V)=(\pi_*\varpi)(V)$ on $T^*\cM$ where $\pi\colon T^*\cM\to\cM$ is the projection. In the 3b-setting then, we have, directly from its definition, that $\alpha\in\CI(\Ttb^*\tilde M_0;\Ttb^*(\Ttb^*\tilde M_0))$. Since $(\dd\alpha)(V,W)=V(\alpha(W))-W(\alpha(V))-\alpha([V,W])$ and since the space of 3b-vector fields on $\Ttb^*\tilde M_0$ is a Lie algebra, we obtain $\omega\in\CI(\Ttb^*\tilde M_0;\Lambda_\tbop^2(\Ttb^*\tilde M_0))$. Combining this with the observation that $\dd\in\Difftb^1(\tilde M_0;\ubar\C,\Ttb^*\tilde M_0)$ yields the above results without computation.

Analogously, then, consider $p\in\CI(\Tetb^*\tilde M)$ in local coordinates $x_\sscri,\rho_+,\omega$ near $\scri^+\cap\iota^+$, and writing eb-covectors as
\begin{equation}
\label{EqCHamebCoord}
  \sigma_\ebop\frac{\dd\rho_+}{\rho_+}+\xi_\ebop\frac{\dd x_\sscri}{x_\sscri}+\eta_\ebop\frac{\dd\omega}{x_\sscri}.
\end{equation}
One then computes
\begin{equation}
\label{EqCHameb}
\begin{split}
  H_p &=(\pa_{\sigma_\ebop}p)\rho_+\pa_{\rho_+} + (\pa_{\xi_\ebop}p)(x_\sscri\pa_{x_\sscri}+\eta_\ebop\pa_{\eta_\ebop}) + (\pa_{\eta_\ebop}p) x_\sscri\pa_\omega \\
    &\quad - (\rho_+\pa_{\rho_+}p)\pa_{\sigma_\ebop} - \bigl((x_\sscri\pa_{x_\sscri}+\eta_\ebop\pa_{\eta_\ebop})p\bigr)\pa_{\xi_\ebop} - (x_\sscri\pa_\omega p)\pa_{\eta_\ebop}.
\end{split}
\end{equation}
Globally, either from these explicit expressions or for the aforementioned structural reasons, we have
\begin{equation}
\label{EqCHame3bMem}
  H_p\in\Vetb(\Tetb^*\tilde M),\quad (p\mapsto H_p) \in \Diff_\etbop^1(\Tetb^*\tilde M;\ubar\C,\Tetb^*(\Tetb^*\tilde M)).
\end{equation}

\subsection{Relation to product compactifications}
\label{SsCP}

In order to prepare for relating dynamical metrics and variable coefficient operators to their model operators---a stationary (Kerr) model at $\cK^+$, and a dilation-invariant (Minkowskian) model at $\iota^+$---we first relate the spacetime compactifications of Definition~\ref{DefCMSpacetime} to other, product-type, compactifications.

\subsubsection{Slicing into copies of \texorpdfstring{$X$}{the spatial manifold}}
\label{SssCPX}

A natural starting point for the analysis of time-translation invariant operators is to slice the spacetime into level sets of $t_*$ (the level sets of which are transversal to $\scri^+$) i.e., to consider $\ol{\R_{t_*}}\times X$ where $X$ is as in~\eqref{EqCSpatial}. Local coordinates near the corner $\{t_*=\pm\infty\}\times\pa X$ are
\[
  \frac{1}{\pm t_*}=\frac{1}{\pm(t-r)},\quad \rho=\frac{1}{r},\quad\omega\in\Sph^2.
\]
The blow-up
\begin{equation}
\label{EqCPXM1p}
  M_1' := \bigl[\ol{\R_{t_*}}\times X; \pa\ol\R\times\pa X\bigr]
\end{equation}
is covered by the following sets of coordinates:
\begin{enumerate}
\item in the interior: $t_*,x$; or $t_*,r,\omega$ as in~\eqref{EqCMCoordInt};
\item near $\pa\ol\R\times X^\circ$: $\rho_\cK:=\frac{1}{\pm t_*}$, $x\in\R^3$ (with $|x|\geq\bhm$), which are also local coordinates near $(\cK^+)^\circ$ in $M$ and $M_1$;
\item near $\R\times\pa X$: $t_*$, $\rho=\frac{1}{r}$, $\omega\in\Sph^2$, which are local coordinates near $(\scri_1^+)^\circ\subset M_1$ (with $\sqrt\rho$ being a local coordinate near $(\scri^+)^\circ\subset M$);
\item near the front face corresponding to the blow-up of $\{+\infty\}\times\pa X$:
  \begin{enumerate}
  \item adjacent to the lift of $\{+\infty\}\times X$: $\rho=\frac{1}{r}$, $\frac{1/t_*}{\rho}=\frac{r}{t-r}$, $\omega\in\Sph^2$; since $\frac{r}{t-r}=\frac{r}{t}(1-\frac{r}{t})^{-1}=\rho_\cK(1-\rho_\cK)$, these match the coordinates~\eqref{EqCMCoordipK} near $\iota^+\cap\cK^+$ up to a smooth coordinate change;
  \item adjacent to the lift of $\ol\R\times\pa X$: $t_*$, $\rho t_*=v$, $\omega\in\Sph^2$, which are also local coordinates near $\scri_1^+\cap\iota_1^+$ (cf.\ \eqref{EqCMCoordscriip});
  \end{enumerate}
\item near the front face corresponding to the blow-up of $\{-\infty\}\times\pa X$, and adjacent to the lift of $\ol\R\times\pa X$: $\rho_0=\frac{1}{-(t_*-2)}=\frac{1}{r-t+2}$, $\rho_\sscri=\frac{\rho}{-1/(t_*-2)}=\frac{r-t+2}{r}$, $\omega\in\Sph^2$, which are also local coordinates near $I_1^0\cap\scri_1^+$ (cf.\ \eqref{EqCMCoordI0scri}).
\end{enumerate}
In the final coordinate system, the hypersurface $t=0$ is given by $\rho_\sscri(1-2\rho_0)=1$, which intersects $\rho_0=0$ at $\rho_\sscri=1$. (This is why we do not consider coordinates near the lift of $\{-\infty\}\times X$.) In summary, these computations prove the following variant of \cite[Lemma~3.18]{HintzNonstat}:

\begin{lemma}[Slicing by $X$]
\label{LemmaCPXDiffeo}
  The map $(t,r,\omega)\mapsto(t_*,r,\omega)=(t-r,r,\omega)$ (i.e., the identity map on $\R^4$) induces a diffeomorphism
  \[
    M_1 \cap \cl_{M_1}\{t\geq 0\} \to M'_1 \cap \cl_{M_1'}\{t_*\geq-r\}
  \]
  which identifies $\cK^+$, $\iota_1^+$, $\scri_1^+$, and $I^0_1\cap\{\frac{t}{r}\geq 0\}$ with the lift of $\{+\infty\}\times X$, the lift of $\{+\infty\}\times\pa X$, the lift of $\ol\R\times\pa X$, and the intersection of the lift of $\{-\infty\}\times\pa X$ with $\{\rho_\sscri\leq 1\}$ in $M'_1$, respectively.
\end{lemma}

\begin{figure}[!ht]
\centering
\includegraphics{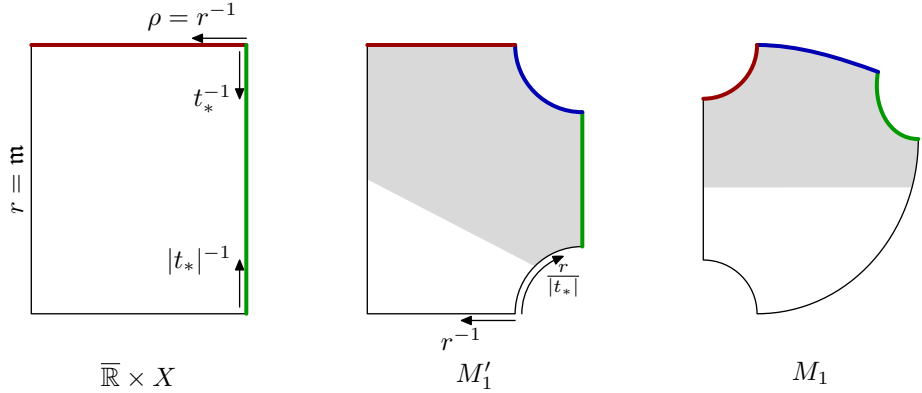}
\caption{\textit{On the left:} $\ol{\R_{t_*}}\times X$ (without the $\Sph^2$ factor) and some local coordinates. \textit{In the middle:} the resolution $M'_1$, defined in~\eqref{EqCPXM1p}. The shaded region is where $t_*\geq-r$. \textit{On the right:} the space $M_1$, shaded where $t\geq 0$. The two shaded regions are naturally diffeomorphic (via the extension of the identity map on $\R^4$). Similarly, the interiors of boundary hypersurfaces with matching colors are identified.}
\label{FigCPXDiffeo}
\end{figure}

See Figure~\ref{FigCPXDiffeo}.

Given a 3b-differential operator $A\in\Diff_\tbop^m(M_0)$, there exists a unique $t_*$-translation-invariant operator $A_0$ on $\R_{t_*}\times X$ such that $\chi(A-A_0)\in\rho_\cK\Diff_\tbop^m(M_0)$ where $\chi\in\CI(M_0)$ is supported in a collar neighborhood of $\cK^+$; this is the \emph{$\cT$-normal operator} of $A$ in the terminology of \cite[Definition~3.19]{Hintz3b}. Indeed, in local coordinates $\rho_\cK=t_*^{-1}$, $x\in\R^3$ near $(\cK^+)^\circ$, one can write
\begin{subequations}
\begin{equation}
\label{EqCPXA}
  A = \sum_{j+|\alpha|\leq m} a_{j\alpha}(\rho_\cK,x) \pa_{t_*}^j \pa_x^\alpha,
\end{equation}
and then $A_0$ is given by freezing coefficients at $\cK^+$, i.e.,
\begin{equation}
\label{EqCPXA0}
  A_0 = \sum_{j+|\alpha|\leq m} a_{j\alpha}(0,x) \pa_{t_*}^j \pa_x^\alpha.
\end{equation}
(The case of operators $A\in\rho_\sface^{-\alpha_\sface}\Diff_\tbop^m(M_0)$ with weights at $\sface$ is completely analogous.) We work here on $M_0$ only since the boundary hypersurface $\scri^+_1\subset M_1$ is disjoint from $\cK^+$ and thus plays no role for this discussion. The structural reasons for the existence of $A_0$ are given in \cite[\S{3.2}]{Hintz3b}, and we give a geometric reason (based on an inspection of unit cells for 3b-vector fields) in~\S\ref{SssMUK}. For later reference, it is useful to give global expressions near $\cK^+$: in terms of local coordinates $t_*,r,\omega$ and local defining functions $\rho_\cK$ and $\rho_+=r^{-1}$ of $\cK^+,\iota^+$, write
\begin{equation}
\label{EqCPXAgl}
  A = \sum_{j+k+|\alpha|\leq m} a_{j k\alpha}(\rho_+,\rho_\cK,\omega) (r\pa_{t_*})^j (r\pa_r)^k \pa_\omega^\alpha;
\end{equation}
then
\begin{equation}
\label{EqCPXA0gl}
  A_0 = \sum_{j+k+|\alpha|\leq m} a_{j k\alpha}(r^{-1},0,\omega) (r\pa_{t_*})^j (r\pa_r)^k \pa_\omega^\alpha.
\end{equation}
\end{subequations}

\subsubsection{Slicing into copies of \texorpdfstring{$\iota^+$}{punctured timelike infinity}}
\label{SssCPip}

The boundary hypersurface $\iota_1^+$ of $\tilde M_1$ is one component of the blow-up of $\sface\subset\tilde M_0$ at its hypersurface $Y^+\subset\sface$. We thus have
\[
  \iota_1^+ = [0,1]_{\rho_\cK} \times \Sph^2,\quad \rho_\cK:=\frac{r}{t}.
\]
By this we mean that the local coordinate map $\tilde M\to[0,2)_{\rho_\cK}\times[0,1)_{\rho_+}\times\Sph^2$ restricts to $\iota_1^+$ to a diffeomorphism of $\iota_1^+$ with $[0,1]\times\{0\}\times\Sph^2$. Equivalently, $\iota_1^+=[0,\infty]_v\times\Sph^2$ where $v=\frac{t-r}{r}$ as in~\eqref{EqCMCoordip}. On the level of $\iota_1^+$, passage to $\tilde M$ amounts to performing a square root blow-up of $\rho_\cK^{-1}(1)=v^{-1}(0)\subset\iota_1^+$; thus,
\begin{equation}
\label{EqCPip}
  \iota^+ = \bigl( [0,2)_{\sqrt{v}} \cup (1,\infty]_v \bigr) \times \Sph^2.
\end{equation}
In a neighborhood of $(\iota^+)^\circ$, we can use the coordinates $v,\frac{1}{t_*},\omega$ as in~\eqref{EqCMCoordip}. Note that $\frac{1}{t_*}$ is not a local defining function of $\iota^+\subset M$, however, but rather a total defining function of $\iota^+\cup\cK^+$. Nonetheless, spacetime dilations act in these coordinates as dilations of $\frac{1}{t_*}$ only. This motivates:

\begin{lemma}[Slicing by $\iota^+$]
\label{LemmaCPip}
  With $\iota^+$ as in~\eqref{EqCPip}, define
  \[
    M_+' := \bigl[ [0,\infty)_\tau \times \iota^+; \{0\}\times v^{-1}(\infty) \bigr].
  \]
  Then the map $(t,r,\omega)\mapsto(\tau,v,\omega)$ where $\tau=\frac{1}{t_*}$, $v=\frac{t_*}{r}$ (with $t_*=t-r$), i.e., the identity map on a subset of $\R^4$, induces a diffeomorphism
  \begin{equation}
  \label{EqCPipDiffeo}
    M \cap \cl_M\{t_*\geq 1,\ r\geq 1\} \to M_+' \cap \cl_{M_+'}\{\tau\leq 1,\ \tau v\leq 1\}
  \end{equation}
  which identifies $\iota^+$, $\cK^+\cap\{r\geq 1\}$, and $\scri^+\cap\{t_*\geq 1\}$ with the lift of $\{0\}\times\iota^+$, the intersection of the front face of $M'_+$ with $\{\tau v\leq 1\}$, and $[0,1]\times v^{-1}(0)$, respectively.
\end{lemma}

\begin{figure}[!ht]
\centering
\includegraphics{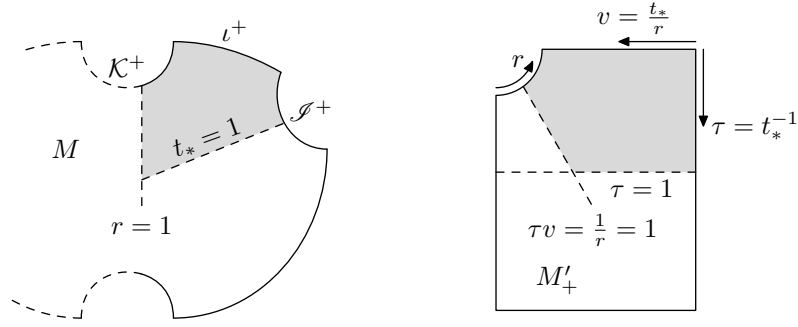}
\caption{\textit{On the left:} $M$. \textit{On the right:} $M'_+$. The regions of interest in~\eqref{EqCPipDiffeo} are shaded in gray.}
\label{FigCPipDiffeo}
\end{figure}

See Figure~\ref{FigCPipDiffeo}.

\begin{proof}[Proof of Lemma~\usref{LemmaCPip}]
  It suffices to show that local coordinates covering the lift of $\{0\}\times\iota^+$ in $M'_+$ are also local coordinates on $M$ near the lift of $\iota^+$. In $v<2$, we can use, on $M'_+$, the coordinates $\sqrt{v}$ and $\tau=\frac{1}{t_*}$, which, on $M$, are local boundary defining functions of $\scri^+$ and $\iota^+$ (see~\eqref{EqCMCoordscriip}). In $v^{-1}<2$, we first pass to the local coordinates $\frac{1}{v}$, $\tau$ on $[0,\infty)\times\iota^+$ and then, upon blowing up $\frac{1}{v}=\tau=0$, to the local coordinates $\frac{1}{v}=\frac{r}{t_*}=\frac{r}{t}(1-\frac{r}{t})$ and $\frac{\tau}{1/v}=\frac{1}{r}$, which are, indeed, local coordinates near $\iota^+\cap\cK^+$ (cf.\ \eqref{EqCMCoordipK}).
\end{proof}

Given an e3b-differential operator $A\in\Diff_\etbop^m(\tilde M)$, there exists a unique dilation-invariant operator $A_0$ on $(0,\infty)_\tau\times\iota^+$ such that $\chi(A-A_0)\in\rho_+\Diff_\etbop^m(\tilde M)$ where $\chi\in\CI(\tilde M)$ is supported in a small collar neighborhood of $\iota^+$ (where Lemma~\ref{LemmaCPip} applies). We check this first near $\iota^+\cap\scri^+$ in local coordinates $x_\sscri=\sqrt{v}$, $\rho_+=\frac{1}{t_*}$, and $\omega\in\Sph^2$: there, we can write
\begin{subequations}
\begin{equation}
\label{EqCPipA}
  A = \sum_{j+k+|\alpha|\leq m} a_{j k\alpha}(x_\sscri,\rho_+,\omega) (x_\sscri\pa_{x_\sscri})^j (\rho_+\pa_{\rho_+})^k (x_\sscri\pa_\omega)^\alpha,
\end{equation}
and then $A_0$ is locally given by
\begin{equation}
\label{EqCPipA0scri}
  A_0 = \sum_{j+k+|\alpha|\leq m} a_{j k\alpha}(x_\sscri,0,\omega) (x_\sscri\pa_{x_\sscri})^j (\rho_+\pa_{\rho_+})^k (x_\sscri\pa_\omega)^\alpha,
\end{equation}
whereas near $\iota^+\cap\cK^+$, we write $\rho_+=\frac{1}{r}$, $\rho_\cK=\frac{r}{t_*}$, and work with the coordinates $R=\frac{r}{t_*}=\rho_\cK=v^{-1}$ and $\tau=\frac{1}{t_*}$ (so $r\pa_{t_*}=-R(R\pa_R+\tau\pa_\tau)$, $r\pa_r=R\pa_R$) and
\begin{equation}
\label{EqCPipA2}
  A = \sum_{j+k+|\alpha|\leq m} a_{j k\alpha}(\rho_+,\rho_\cK,\omega) (R \tau\pa_\tau)^j (R\pa_R)^k \pa_\omega^\alpha,
\end{equation}
we have
\begin{equation}
\label{EqCPipA0scri2}
  A_0 = \sum_{j+k+|\alpha|\leq m} a_{j k\alpha}(0,\rho_\cK,\omega) (R \tau\pa_\tau)^j (R\pa_R)^k \pa_\omega^\alpha.
\end{equation}
\end{subequations}
Carefully note that the coordinates used in~\eqref{EqCPipA} and \eqref{EqCPipA2} are product coordinates on $[0,\infty)_\tau\times\iota^+$. We also note that these computations show that $R\tau\pa_\tau$, $R\pa_R$, $\pa_\omega$ are a local frame for the space of 3b-vector fields near $\iota^+\cap\cK^+$.

\section{Microlocal toolkit}
\label{SM}

In this section, we introduce/recall all of the ps.d.o.\ algebras utilized in our analysis of waves on asymptotically Kerr spacetimes.\footnote{We do assume that the reader is familiar with standard ps.d.o.\ calculus, including in the semiclassical setting, as presented in many standard references such as \cite{ShubinSpectralTheory,GrigisSjostrandBook,HormanderAnalysisPDE3,ZworskiSemiclassical,VasyMinicourse,DyatlovZworskiBook}; we follow the conventions of \cite{HintzMicro}.} These algebras had already been used by the author in \cite{HintzNonstat} for the linear analysis of wave-type operators on asymptotically flat spacetimes, with the important technical caveat that in \cite{HintzNonstat} we only needed algebras of ps.d.o.s with smooth (or conormal) coefficients; these were recalled from a geometric microlocal perspective in \cite[\S{2}]{HintzNonstat}. In the present paper, by contrast, we need to rely on the scaled bounded geometry perspective introduced by the author in \cite{HintzScaledBddGeo} and recalled in~\S\ref{SsMS}, which handles operators whose coefficients have less (uniform) regularity. The reason is that for nonlinear applications, it is crucial to prove estimates for waves in which the constants depend explicitly on (high regularity) norms of the coefficients of the underlying operator, as measured in function spaces in which the estimate itself takes place. This essentially automatic when, say, one utilizes ps.d.o.s with e3b-regular coefficients for the purpose of proving e3b-microlocal estimates; see Proposition~\ref{PropMEll} (with $k$ there taken to be $0$) for a simple example. By contrast, the proof of the boundedness of a 3b-ps.d.o.\ on 3b-Sobolev spaces as proved in \cite{Hintz3b} uses the smoothness (or, at the very least, the conormality) of the coefficients of the ps.d.o.\ (due to the use of pullback/push-forward theorems of Melrose \cite{MelrosePushfwd}).

The e3b-algebra on $\tilde M$, which will be our main workhorse for the uniform analysis of the regularity of waves on spacetime, is defined from the scaled b.g.\ perspective in~\S\ref{SssMUe3b}, and the remaining sections of~\S\ref{SsMU} introduce related algebras arising upon Fourier transforming in $t_*$ or Mellin-transforming along spacetime scalings. In~\S\ref{SsMBasic}, we prove basic microlocal estimates (elliptic regularity, real principal type propagation) on e3b-spaces as well as on space encoding additional integer amounts $k$ of b-regularity. In~\S\ref{SsMTame}, we prove versions which are \emph{tame} in the b-regularity order $k$ when $k$ gets large. (The full regularity analysis on dynamical Kerr spacetimes requires, in addition, radial point and trapping estimates; these are discussed in~\S\ref{SR}.)

\begin{rmk}[Bundles]
\label{RmkMBundles}
  For brevity, we state definitions and results here for complex-valued symbols and distributions and operators acting on such. We leave the notational adaptations when working with sections of vector bundles to the reader; it suffices to require these vector bundles to be uniform vector bundles in the bounded geometry setting (cf.\ \cite[Remarks~3.13 and 3.47]{HintzScaledBddGeo}) or smooth vector bundles on the underlying manifold with corners in the e3b-setting and related settings.
\end{rmk}

\subsection{Scaled bounded geometry structures}
\label{SsMS}

We recall the basic ideas and constructions of \cite{HintzScaledBddGeo}. First, essentially following Shubin \cite{ShubinBounded}, a \emph{bounded geometry} (or \emph{b.g.}) \emph{structure} on an $n$-dimensional manifold without boundary $\cM$ is a (countable) collection $\fB=\{(U_\alpha,\phi_\alpha)\colon\alpha\in\sA\}$ of coordinate charts $\phi_\alpha\colon U_\alpha(\subset\cM)\to (-2,2)^n$ with finite multiplicity (i.e., no point of $\cM$ is contained in more than a fixed number $J<\infty$ of charts $U_\alpha$) such that $\cM=\bigcup_\alpha \phi_\alpha^{-1}((-1,1)^n)$ and such that the transition functions $\tau_{\beta\alpha}:=\phi_\beta\circ\phi_\alpha^{-1}$ are uniformly bounded in $\CI$. We call the sets $U_\alpha$ \emph{unit cells}. A \emph{scaling} is a collection $\rho_\alpha=(\rho_{\alpha,i})_{i=1,\ldots,n}\in(0,1]^n$, $\alpha\in\sA$, such that for all $\gamma\in\N_0^n$ there exists a constant $C_\gamma<\infty$ such that
\begin{equation}
\label{EqMSCuboids}
  |\pa^\gamma(\rho_{\alpha,i}\pa_i\tau_{\beta\alpha}^j(x))|\leq C_\gamma\rho_{\beta,j}\quad\forall\,x\in\phi_\alpha(U_\alpha\cap U_\beta),\ i,j=1,\ldots,n.
\end{equation}
In particular, $\tau_{\beta\alpha}$ maps an $\rho_\alpha$-cuboid (with side length $\rho_{\alpha,i}$ in the $i$-th coordinate direction in the chart $\phi_\alpha$) into a set that contains, resp.\ is contained in a $C^{-1}\rho_\beta$, resp.\ $C\rho_\beta$-cuboid, thus guaranteeing the consistency of the size and shape of these cuboids across charts. Constant scalings, i.e., $\rho_{\alpha,i}=h\in(0,1]$ is independent of $i,\alpha$, always satisfy~\eqref{EqMSCuboids}.

Associated with the \emph{scaled b.g.\ structure} $\fB_\times=\{(U_\alpha,\phi_\alpha,\rho_\alpha)\colon\alpha\in\sA\}$ are the following spaces and notions.
\begin{enumerate}
\item The \emph{coefficient Lie algebra} $\cW$ is the set of all smooth vector fields $V\in\cV(\cM)$ whose coefficients, in the charts $\phi_\alpha$, are uniformly bounded in $\CI$, i.e., writing $(\phi_\alpha)_*V=\sum_{i=1}^n v_\alpha^i\pa_{x^i}$, the functions $v_\alpha^i$ are uniformly bounded in $\CI$. We also write $\cW=\CI_{{\rm uni},\fB}(\cM;T\cM)$.
\item We call a function on $\cM$ \emph{$\cW$-regular} if it remains uniformly bounded under application of any finite number of elements of $\cW$. The space of $\cW$-regular functions is denoted $\CI_{{\rm uni},\fB}(\cM)$.
\item The \emph{operator Lie algebra} $\cV$ consists of all $V\in\cV(\cM)$ such that, upon writing $(\phi_\alpha)_*V=\sum_{i=1}^n V_\alpha^i\,\rho_{\alpha,i}\pa_{x^i}$, the functions $V_\alpha^i$ are uniformly bounded in $\CI$.
\item The \emph{secondary b.g.\ structure} $\rho\fB$ is obtained from $\fB_\times$ by subdividing each unit cell $U_\alpha$ into $\rho_\alpha$-cuboids. The corresponding Lie algebra of uniform vector fields is called the \emph{large operator Lie algebra} $\cV'$. (See \cite[\S{3.2}]{HintzScaledBddGeo} for details.)
\item A \emph{weight} is a function $w>0$ such that $\frac{W w}{w}\in\CI_{{\rm uni},\fB}(\cM)$ for all $W\in\cW$. A weight family is a collection $\{w_\alpha\colon\alpha\in\sA\}$ of nonnegative numbers such that $\frac{w_\alpha}{w_\beta}$ varies over a compact subset of $(0,\infty)$ when $\alpha$ and $\beta$ range over all pairs with $U_\alpha\cap U_\beta\neq\emptyset$. An example is $w_\alpha=\inf_{U_\alpha}w$ or $\sup_{U_\alpha}w$, and one can conversely construct a weight from a weight family using a uniform partition of unity. See \cite[\S{3.3}]{HintzScaledBddGeo}.
\item The \emph{scaling weight} $\bar\rho\in\CI(\cM)$ is a weight that is equivalent to the scaling weight family $\{\max_{i=1,\ldots,n}\rho_{\alpha,i}\colon\alpha\in\sA\}$.
\end{enumerate}

In all settings needed in the present paper, the space $\cW$ will be spanned over $\CI_{{\rm uni},\fB}(\cM)$ by the space $\cW_{[\cV]}$ of \emph{commutator $\cW$-vector fields}, which is defined to consist of all $W\in\cW$ such that $[W,V]\in\cV$ for all $V\in\cV$; see \cite[\S{4.4}]{HintzScaledBddGeo}. (Cf.\ Corollaries~\ref{CorCT3bComm}, \ref{CorCTebComm}, and Lemma~\ref{LemmaCTe3bComm}.)

Sobolev spaces are defined as follows. Fix $\chi\in\CIc((-\frac54,\frac54)^n)$ with $\chi|_{[-1,1]^n}=1$; define the uniform partition of unity $\chi_\alpha:=\frac{1}{S}\phi_\alpha^*\chi$ where $S=\sum_\beta\phi_\beta^*\chi$. Define the scaling map $S_\alpha(x):=(\frac{x^1}{\rho_{\alpha,1}},\ldots,\frac{x^n}{\rho_{\alpha,n}})$. Then
\[
  \|u\|_{H_\cV^s(\cM)}^2 := \sum_{\alpha\in\sA} \| (S_\alpha)_*(\phi_\alpha)_*(\chi_\alpha u) \|_{H^s(\R^n)}^2;
\]
see \cite[Definition~3.14]{HintzScaledBddGeo}. For $s\in\N_0$, this is equivalent to taking the sum of squared $L^2$-norms of up to $s$-fold derivatives of $u$ along the elements of a finite spanning set of $\cV$ over $\CI_{{\rm uni},\fB}(\cM)$. The norm on $w H_\cV^s(\cM)$ is defined in the same way with an extra factor $w_\alpha^{-2}$ on the $\alpha$-th summand. For $k\in\N_0$, we moreover define
\[
  \|u\|_{H_{\cV;\cW}^{(s;k)}(\cM)}^2 := \sum_{j=0}^k \sum_{W_1,\ldots,W_j\in\sW} \|W_1\cdots W_j u\|_{H_\cV^s(\cM)}^2
\]
where $\sW\subset\cW$ is any finite spanning set over $\CI_{{\rm uni},\fB}(\cM)$; see the discussion following \cite[Definition~3.17]{HintzScaledBddGeo} for details and alternative definitions. We then denote by
\begin{equation}
\label{EqMSSobolev}
  H_\cV^s(\cM),\quad H_{\cV;\cW}^{(s;k)}(\cM)
\end{equation}
the (Hilbert) spaces of distributions with finite $\|\cdot\|_{H_\cV^s(\cM)}$- and $\|\cdot\|_{H_{\cV;\cW}^{(s;k)}}$-norms, respectively.

Classes of symbols associated with $\fB_\times$ are defined as follows.

\begin{enumerate}
\item For $m\in\R$, the symbol class $S^m({}^\cV T^*\cM)$ consists of all smooth functions $a$ on $T^*\cM$ such that the symbols $a_\alpha\in S^m((-2,2)^n\times\R^n)$ defined by
  \begin{equation}
  \label{EqMSSymbol}
    a_\alpha(x,\xi) = a\biggl(\phi_\alpha^{-1}(x), \sum_{i=1}^n \xi_i\frac{\dd x^i}{\rho_{\alpha,i}}\biggr),
  \end{equation}
  are uniformly bounded. Here, we use the standard seminorms
  \begin{equation}
  \label{EqMSSymbolSeminorms}
    \|a_\alpha\|_{S^m;q} := \max_{|\beta|,|\gamma|\leq q}\sup_{x\in(-2,2)^n}\sup_{\xi\in\R^n} \la\xi\ra^{-m+|\gamma|}|\pa_x^\beta\pa_\xi^\gamma a_\alpha(x,\xi)|,\quad q\in\N_0.
  \end{equation}
  For $\delta\in(0,\frac12)$, the symbol class $S^m_{1-\delta,\delta}({}^\cV T^*\cM)$ is defined similarly, except the seminorms on $S^m_{1-\delta,\delta}$ are defined using $\la\xi\ra^{-m+|\gamma|-\delta(|\beta|+|\gamma|)}$ instead of $\la\xi\ra^{-m+|\gamma|}$; see \cite[\S{4.2}]{HintzScaledBddGeo}.
\item Fix $\lambda\in S^{-1}({}^\cV T^*\cM)$ with $\lambda^{-1}\in S^1({}^\cV T^*\cM)$. For $\sfm\in S^0({}^\cV T^*\cM)$, the class $S^\sfm({}^\cV T^*\cM)$ of variable-order symbols consists of all elements $a\in\bigcap_{\delta>0}S^{\sup\sfm}_{1-\delta,\delta}$ which are of the form $a=\lambda^{-\sfm}a_0$ with $a_0\in\bigcap_{\delta>0}S^0_{1-\delta,\delta}$.
\item More generally, given a finite collection $0<\rho_1,\ldots,\rho_N\in S^0({}^\cV T^*\cM)$ such that for each $j=1,\ldots,N$ there exists $N_j$ such that $\frac{\bar\rho^{N_j}\lambda^{N_j}}{\rho_j}\in S^0({}^\cV T^*\cM)$ (i.e., $\rho_j$ dominates some power of $\bar\rho\lambda$), we write
  \[
    S^{\sfm,(\sfl_1,\ldots,\sfl_N)}({}^\cV T^*\cM) := \Biggl(\;\prod_{j=1}^N\rho_j^{-\sfl_j}\Biggr)\lambda^{-\sfm} \bigcap_{\delta>0} S^0_{1-\delta,\delta}({}^\cV T^*\cM)
  \]
  for the class of symbols with variable differential order and variable weights.
\end{enumerate}

Given a symbol $a\in S^m({}^\cV T^*\cM)$, we define the associated $\cV$-ps.d.o.\ $\Op_\cV(a)$ as follows. Fix a cutoff function $\psi\in\CIc((-\frac14,\frac14)^n)$ which is equal to $1$ near $0$. Then
\begin{equation}
\label{EqMSQuantFull}
  \Op_\cV(a) \in \Psi_\cV^m(\cM)
\end{equation}
is defined as the sum of quantizations $\Op_\alpha(a_\alpha)$ of the symbol $\chi_\alpha a$, pulled back to $\cM$ via $\phi_\alpha$; here
\begin{equation}
\label{EqMSQuant}
  (\Op_\alpha(a_\alpha)u)(x) = (2\pi)^{-n}\int_{\R^n}\int_{\R^n} \exp\biggl(i\sum_{j=1}^n (x^j-x'{}^j)\frac{\xi_j}{\rho_{\alpha,j}}\biggr)\psi(x-x') a(x,\xi)u(x')\,\dd\xi\,\frac{\dd x'{}^1\cdots\dd x'{}^n}{\rho_{\alpha,1}\cdots\rho_{\alpha,n}}.
\end{equation}
When $a$ is a polynomial in the fibers, then $\Op_\cV(a)$ is a differential operator. Roughly speaking, it is a ``$\cV$-differential operator with $\cW$-regular coefficients'': in local coordinates, it is built from the vector fields $\rho_{\alpha,i}\pa_{x^i}$ but with coefficients that are regular under application of (the stronger vector fields) $\pa_{x^j}$.

The space $\Psi_\cV^m(\cM)$ is defined as the space of all operators of the form $\Op_\cV(a)+R$ where $R$ is a \emph{residual} operator; this means that $R$ and its $L^2$-adjoint with respect to any uniformly positive $\cV$-density are bounded as maps from $w\bar\rho^{-N}H_\cV^{-N}$ to $w\bar\rho^N H_\cV^N$ for all $N\in\R$, and the same remains true for any finite-order commutator of $R$ with elements of $\cW_{[\cV]}$ (see \cite[Definition~3.34, equation~(4.28)]{HintzScaledBddGeo}). The class of all residual operators is denoted $\bar\rho^\infty\Psi_\cV^{-\infty}(\cM)$.\footnote{In \cite{HintzScaledBddGeo}, this is denoted $\bar\rho^\infty\Psi_{\cV,[\,]}^{-\infty}(\cM)$.} The principal symbol map is
\[
  \upsigma_\cV^m \colon \Psi_\cV^m(\cM) \to S^m/\bar\rho S^{m-1}({}^\cV T^*\cM)
\]
and has all the usual properties under addition, composition, and commutators; see \cite[Theorem~3.52]{HintzScaledBddGeo}. That is, for $A_j\in\Psi_\cV^{m_j}(\cM)$, $j=1,2$, with principal symbols $a_1,a_2$,
\[
  \upsigma_\cV^{m_1+m_2}(A_1\circ A_2) = a_1 a_2,\quad
  \upsigma_\cV^{m_1+m_2-1}(i[A_1,A_2]) = H_{a_1}a_2,
\]
where $H_p$ for $p\in\CI(T^*\cM)$ is the Hamiltonian vector field. (Carefully note that in the coordinates $x,\xi$ as used in~\eqref{EqMSQuant}, this is given by $(\pa_{\xi_j}p)\rho_{\alpha,j}\pa_{x^j}-(\rho_{\alpha,j}\pa_{x^j}p)\pa_{\xi_j}$.) In view of our usage of \cite[\S{4.4}]{HintzScaledBddGeo},
\begin{equation}
\label{EqMSComm}
  A\in\Psi^s_\cV(\cM),\ V\in\cW_{[\cV]} \implies [V,A]\in\Psi_\cV^s(\cM),\ \WF'_\cV([V,A])\subset\WF'_\cV(A).
\end{equation}

The definitions for weighted ps.d.o.s $w\Psi_\cV^m(\cM)$ are completely analogous.

\begin{rmk}[Generalizing coefficients]
\label{RmkMSCoeff}
  Note that a $\cW$-regular function is automatically $\cV$-regular (i.e., remains uniformly bounded upon application of any finite number of elements of $\cV$). We can relax the regularity of the coefficients of $\cV$-ps.d.o.s from $\cW$-regularity to $\cV$-regularity by passing to the secondary b.g.\ structure: the large operator Lie algebra $\cV'$ is spanned by $\cV$ over the space $\CI_{{\rm uni},\rho\fB}(\cM)$ of $\cV$-regular functions. Then $\Psi_\cV^m(\cM)\subset\Psi_{\cV'}^m(\cM)$. See also \cite[Definition~1.4(4), Remark~3.49]{HintzScaledBddGeo}. Furthermore, $H_\cV^s(\cM)=H_{\cV'}^s(\cM)$, as argued after \cite[Lemma~3.16]{HintzScaledBddGeo}.
\end{rmk}

Likewise, spaces $\Psi_\cV^\sfm(\cM)$, $\sfm\in S^0({}^\cV T^*\cM)$, of variable order $\cV$-ps.d.o.s are defined in an analogous fashion; the principal symbol is now valued in $S^\sfm/\bigcap_{\delta>0}S^{\sfm-1+2\delta}({}^\cV T^*\cM)$. The cases of operators with (phase space) weights or variable decay orders are similar as well; see \cite[\S\S{4.1} and 4.2]{HintzScaledBddGeo}. One can then use elliptic elements of these spaces of ps.d.o.s to define corresponding scales of Sobolev spaces such as $H_\cV^\sfm(\cM)$ and
\begin{equation}
\label{EqMSMixed}
  H_{\cV;\cW}^{(\sfm;k)}(\cM).
\end{equation}
We recall the that the $L^2$-duals of $H_\cV^\sfm(\cM)$ and $H_\cV^w(\cM)$ are $H_\cV^{-\sfm}(\cM)$ and $H_\cV^{w^{-1}}(\cM)$, respectively; see \cite[Proposition~4.2]{HintzScaledBddGeo} for a proof of the latter duality.

\begin{rmk}[Equivalent b.g.\ structures]
\label{RmkMSEquiv}
  Many choices of scaled b.g.\ structures on $\cM$ give rise to the same spaces of ps.d.o.s (and operator/coefficient Lie algebras and Sobolev spaces): for two scaled b.g.\ structures, this is the case when their union is also a scaled b.g.\ structure, or equivalently if the transition functions across the two structures are uniformly bounded in both the bounded geometry and the scaling sense. See \cite[Definition~3.1, Remark~2.2]{HintzScaledBddGeo} for details.
\end{rmk}

One can use the quantization map~\eqref{EqMSQuantFull} also for symbols $a$ with less regularity, for the following purpose. Consider a scaled b.g.\ structure $\{(U_\alpha,\phi_\alpha,\rho_\alpha)\}$ on $\cM$. For $s,m\in\R$ lying in fixed compact subset $I\subset\R$, recall from \cite[Remark~3.33]{HintzScaledBddGeo} that the boundedness of $\Op_\cV(a)\colon H_\cV^s(\cM)\to H_\cV^{s-m}(\cM)$ only requires the $\cV$-regularity of the symbol $a$ in the base variables. Quantitatively:

\begin{definition}[Finite (mixed) regularity symbols]
\label{DefMSSymbolMixed}
  Given a function $a\colon\cM\to\C$, write
  \[
    a_\alpha(x,\xi)=\chi_\alpha(\phi_\alpha^{-1}(x)) a\Bigl(\phi_\alpha^{-1},\sum_{i=1}^n \xi_i\frac{\dd x^i}{\rho_{\alpha,i}}\Bigr)
  \]
  analogously to~\eqref{EqMSSymbol}, and define the following refinement of~\eqref{EqMSSymbolSeminorms}: for $d_0,k\in\N_0$,
  \begin{equation}
  \label{EqMSSymbolNorms2}
    \|a_\alpha\|_{\cC_{-;\rho_\alpha}^{(d_0;k)}S^m} := \max_{|\beta'|,|\gamma'|\leq d_0} \max_{|\beta|,|\gamma|\leq k} \sup_{x,\xi} \la\xi\ra^{-m+|\gamma|+|\gamma'|} | (\rho_\alpha\pa_x)^{\beta'}\pa_x^\beta\pa_\xi^{\gamma+\gamma'}a_\alpha(x,\xi)|
  \end{equation}
  where $(\rho_\alpha\pa_x)^{\beta'}=\prod_{j=1}^n (\rho_{\alpha,j}\pa_{x^j})^{\beta'_j}$. We then set $\|a\|_{\cC_{\cV;\cW}^{(d_0;k)}S^m}:=\sup_\alpha\|a_\alpha\|_{\cC_{-;\rho_\alpha}^{(d_0;k)}S^m}$, and define
  \[
    \cC_{\cV;\cW}^{(d_0;k)}S^m({}^\cV T^*\cM)
  \]
  for the space of all functions with finite $\|\cdot\|_{\cC_{\cV;\cW}^{(d_0;k)}S^m}$-norm. For $k=0$, we write $\cC_\cV^{d_0}$ instead of $\cC_{\cV;\cW}^{(d_0;0)}$.
\end{definition}

\begin{lemma}[Operator norm bounds]
\label{LemmaMSOpNorm}
  Given $I\Subset\R$, there exist $d_0\in\N$ and $C$ such that for all $s,m\in I$, and $a\in\cC_\cV^{d_0}S^m({}^\cV T^*\cM)$, we have
  \[
    \|\Op_\cV(a)\|_{\cL(H_\cV^s(\cM),H_\cV^{s-m}(\cM))} \leq C\|a\|_{\cC_\cV^{d_0}S^m}.
  \]
  More generally, \emph{for the same value $d_0$}, we have, for all $k\in\N_0$ and $a\in\cC_{\cV;\cW}^{(d_0;k)}S^m({}^\cV T^*\cM)$,
  \[
    \|\Op_\cV(a)\|_{\cL(H_{\cV;\cW}^{(s;k)}(\cM),H_{\cV;\cW}^{(s-m;k)}(\cM))} \leq C\|a\|_{\cC_{\cV;\cW}^{(d_0;k)}S^m}.
  \]
\end{lemma}
\begin{proof}
  We have already argued for the existence of $d_0$ in the case that $k=0$. The case of higher $k$ follows by inspection of the proof of \cite[Proposition~3.32]{HintzScaledBddGeo}; in local coordinates, higher $\xi$-regularity of $a_\alpha$ is not even needed in order to for its quantization to preserve $x$-regularity, cf.\ \cite[equation~(3.29)]{HintzScaledBddGeo}. (We introduce it in the norms~\eqref{EqMSSymbolNorms2} nonetheless, as this enforces compatibility of these norms across coordinate charts.)
\end{proof}

Analogous results hold for weighted operators and spaces, including the case of variable orders when the range of the order functions lies in fixed compact subsets of $\R$. For the standard ps.d.o.\ algebra on $\R^n$, related results, including precise bounds on $d_0$, have been obtained by various authors including \cite{KumanogoNagasePsdo,BealsReedMicroNonsmooth,MarschallPseudo,WittCalculusForNonsmooth,HintzQuasilinearDS}.

\subsubsection{b- and scattering operators}
\label{SssMSbsc}

Consider $\cM=\R^n$ and cover it using the charts
\begin{equation}
\label{EqMSbscCharts}
  U_0 = (-2,2)^n,\quad
  U_{\pm i j k} := \Bigl\{ \pm x^i \in (2^{j-2},2^{j+2}),\ \Bigl(\frac{x^1}{x^i},\ldots,\frac{x^{i-1}}{x^i},\frac{x^{i+1}}{x^i},\ldots,\frac{x^n}{x^i}\Bigr)-k \in (-2,2)^{n-1} \Bigr\}
\end{equation}
where $i\in\{1,\ldots,n\}$, $j\in\N_0$, and $k\in\Z^{n-1}$, $\max|k_i|\leq 10$, with $\phi_0$ the identity and $\phi_{\pm 1,j,k}(x^1,x')=(\log_2(\pm x^1)-j,\frac{x'}{x^1}-k)$ where we write $x=(x^1,x')$; similarly for $\phi_{\pm i j k}$. This defines a b.g.\ structure $\fB$ on $\cM$. The space $\CI_{{\rm uni},\fB}(\cM)$ consists of all functions $u$ which remain uniformly bounded upon application of any finite number of vector fields $\pa_{x^i}$, $x^i\pa_{x^j}$; recalling the discussion after Definition~\eqref{DefCTb}, this means
\[
  \CI_{{\rm uni},\fB}(\cM) = \cC_\bop^\infty(\ol{\R^n}) := \{ u\in\CI(\R^n)\colon V_1\cdots V_k u\in L^\infty(\R^n)\ \forall\,k\in\N_0,\ V_i\in\Vb(\ol{\R^n}) \}.
\]
(The space on the right is typically denoted $\cA^0(\ol{\R^n})$, the space of \emph{bounded conormal functions}.) It then follows that
\begin{equation}
\label{EqMSbscT}
  \CI_{{\rm uni},\fB}(\cM;T\cM) = \CI_\bop(\ol{\R^n})\otimes_{\CI(\ol{\R^n})}\Vb(\ol{\R^n}) = \CI_\bop\Vb(\ol{\R^n}),
\end{equation}
i.e.\ the space of uniform vector fields consists of all linear combinations of b-vector fields with infinitely b-regular (i.e., conormal) coefficients.

If we choose the scalings $\rho_{\alpha,l}=1$ for all $\alpha\in\{0,(\pm i j k)\}$ and $l=1,\ldots,n$, then the coefficient and operator Lie algebras $\cW$ and $\cV$ are both equal to~\eqref{EqMSbscT}, and $\Psi_\cV(\R^n)$ is the algebra of b-pseudodifferential operators with b-regular coefficients, which we shall denote by
\begin{equation}
\label{EqMSbsctildePsi}
  \tilde\Psi_\bop(\ol{\R^n}) = \CI_\bop\Psib(\ol{\R^n});
\end{equation}
weighted versions are
\[
  \tilde\Psi_\bop^{m,l}(\ol{\R^n}) = \rho^{-l}\tilde\Psi_\bop^m(\ol{\R^n}),\quad \rho=\la x\ra^{-1}.
\]
Here and in later examples, the tilde refers to the fact that the operator and coefficient Lie algebras coincide, which happens when all scalings are equal to $1$, i.e., we are in a bounded geometry setting.\footnote{One could thus instead write $\Psi_{{\rm uni},\fB}(\R^n)$. Since we shall need to work with many different b.g.\ structures, we find the notation~\eqref{EqMSbsctildePsi} more convenient, in which the subscript `b' indicates the underlying b.g.\ structure.} An equivalent b.g.\ structure can be defined by working in the projective charts such as $\rho=\frac{1}{x^1}$, $X=(\frac{x^2}{x^1},\ldots,\frac{x^n}{x^1})$, and using as unit cells near $\pa\ol{\R^n}$ sets of the type $(2^{-j-2},2^{-j+2})_\rho\times(k+(-2,2)^{n-1})_X$ where $j\in\N_0$ and $k\in\Z^{n-1}$, $\max|k_i|\leq 10$: indeed, the uniform vector fields are then linear combinations of $\rho\pa_\rho$ and $\pa_X$ with b-regular coefficients, giving again~\eqref{EqMSbscT}.\footnote{Conversely, the space of uniform vector fields, or more generally of the coefficient Lie algebra, determines the underlying b.g.\ structure essentially uniquely; see~\cite[\S\S{2.1} and 3.1]{HintzScaledBddGeo}.} We remark that the double dyadic decomposition in \cite[\S{3.1}]{MetcalfeTataruTohaneanuPriceNonstationary} is precisely a decomposition of the forward cone into b-unit cells.

If, on the other hand, we take $\rho_{0,l}=1$ and $\rho_{(\pm i j k),l}=2^{-j}$ (which is, up to a uniform constant, the size of $\la x\ra$ on $U_{\pm i j k}$), then $\cW$ is still equal to~\eqref{EqMSbscT}, while elements of $\cV$ are $\la x\ra^{-1}$-rescalings of elements of~\eqref{EqMSbscT}; thus, $\cV$ is equal to $\CI_\bop(\ol{\R^n})\Vsc(\ol{\R^n})$, and $\Psi_\cV(\R^n)$ in this case is the space of scattering ps.d.o.s with b-regular coefficients; we denote this by
\[
  \CI_\bop\Psisc(\ol{\R^n});\quad \CI_\bop\Psisc^{m,r}(\ol{\R^n}) = \rho^{-r}\CI_\bop\Psisc^m(\ol{\R^n}).
\]
The principal symbol $\sigmasc^{m,r}$ takes values in $S^{m,r}/S^{m-1,r-1}$. Since a scaling weight in this case is $\bar\rho=\la x\ra^{-1}$, we can work with spaces of scattering ps.d.o.s with variable decay and regularity orders (unlike in the case of $\Psib(\ol{\R^n})$, where only the regularity order may be variable), with the usual modifications in the quotient space which the principal symbol map takes values in.

\subsubsection{Compactifications; locus of microlocal analysis; Sobolev spaces}
\label{SssMSCp}

In the b-setting~\eqref{EqMSbsctildePsi}, the symbol space underlying $\tilde\Psi_\bop^m(\ol{\R^n})$ is $\CI_\bop S^m(\Tb^*\ol{\R^n})$, i.e., the linear combination of b-regular functions on $\ol{\R^n}$ (lifted to phase space) and symbols of order $m$ on the vector bundle $\Tb^*\ol{\R^n}$. (In particular, quantizations of elements of $S^m(\Tb^*\ol{\R^n})$ lie in $\tilde\Psi_\bop^m(\ol{\R^n})$.) If we radially compactify each fiber of $\Tb^*\ol{\R^n}$, then these symbols are smooth functions on the interior of $\ol{\Tb^*}(\ol{\R^n})$. The benefit of this compactification is that the symbol $\lambda=\la\xi\ra^{-1}$ (where $\xi$ is the b-momentum defined by writing covectors as $\xi\cdot\frac{\dd x}{\la x\ra}$), which is elliptic of order $-1$, positive, and invertible,
vanishes at fiber infinity $\Sb^*\ol{\R^n}$ (the \emph{spherical b-conormal bundle}). Since the principal symbol of an element of $\tilde\Psi_\bop^m(\ol{\R^n})$ is valued in $\CI_\bop S^m/\lambda\CI_\bop S^m$ (and one locally has full symbol expansions modulo $\bigcap_N\lambda^N\CI_\bop S^m$), this shows that $\Sb^*\ol{\R^n}$ is a good locus of microlocal objects. Thus, given $a\in\CI_\bop S^m$ and its quantization $A\in\tilde\Psi_\bop^m(\ol{\R^n})$, we define
\[
  \Ell_\bop(a)=\Ell_\bop(A),\ \WF_\bop'(A) \subset \Sb^*\ol{\R^n},\quad \Char_\bop(a)=\Sb^*\ol{\R^n}\setminus\Ell_\bop(a),
\]
in the usual fashion; e.g., $\alpha\in\Sb^*\ol{\R^n}$ does not lie in $\WF_\bop'(A)$ if there exists $\chi\in\CI(\ol{\Tb^*}(\ol{\R^n}))$, $\chi(\alpha)\neq 0$, such that $\chi a\in\CI_\bop S^{-\infty}$. (In the terminology of \cite[Definition~3.62]{HintzScaledBddGeo}, $\ol{\Tb^*}(\ol{\R^n})$ is an \emph{admissible compactification} of ${}^\cV T^*\R^n$, $\cV=\CI_\bop\Vb(\ol{\R^n})$.)

The Sobolev spaces corresponding to this setup are denoted
\[
  \Hb^{s,\alpha}(\ol{\R^n})=\la x\ra^{-\alpha}\Hb^s(\ol{\R^n});
\]
we use the Euclidean density to define the underlying $L^2$-space. The regularity order $s$ may be variable, i.e., lie in $\CI({}^\bop S^*\ol{\R^n})$.

Similarly, the symbol space underlying $\CI_\bop\Psi_\scop^m$ is $\CI_\bop S^m(\Tsc^*\ol{\R^n})$. The functions $\rho=\la x\ra^{-1}$ and $\lambda=\la\xi\ra^{-1}$ (where now $\xi$ is the scattering momentum, i.e., defined by writing covectors as $\xi\cdot\dd x$) are smooth functions on $\ol{\Tsc^*}(\ol{\R^n})$ that vanish simply at fiber infinity $\Ssc^*\ol{\R^n}$ and base infinity $\ol{\Tsc^*_{\pa\ol{\R^n}}}(\ol{\R^n})$, respectively. The principal symbol of $A\in\CI_\bop\Psi_\scop^{m,r}(\ol{\R^n})$ is valued in $\CI_\bop S^m/\rho\lambda\CI_\bop S^m$; correspondingly, its elliptic set is a subset
\[
  \Ellsc(A) \subset \Ssc^*\ol{\R^n} \cup \ol{\Tsc^*_{\pa\ol{\R^n}}}(\ol{\R^n}),
\]
similarly for the characteristic set and operator wave front set. The corresponding Sobolev spaces are denoted $\Hsc^{s,r}(\ol{\R^3})=\la x\ra^{-r}\Hsc^s(\ol{\R^3})$, which for $s,r\in\R$ are the standard weighted Sobolev spaces on $\R^3$ (with norm $\|u\|_{\Hsc^{s,r}(\ol{\R^3})}=\|\la x\ra^r\la D\ra^s u\|_{L^2}$). All orders attached to weights which vanish only at the microlocalization locus can be variable; thus, we can, more generally, work with variable differential and variable decay order scattering Sobolev spaces
\[
  \Hsc^{\sfs,\sfr}(\ol{\R^n}),\quad \sfs\in\CI({}^\scop S^*\ol{\R^n}),\ \sfr\in\CI(\ol{\Tsc^*_{\pa\ol{\R^n}}}\ol{\R^n}).
\]
Similar considerations will apply to all other ps.d.o.\ algebras introduced in~\S\ref{SsMU} below, and we shall thus be brief in the sequel.

\subsubsection{Fourier transforms}
\label{SssMSF}

We recall from \cite[Definition~1.10]{HintzScaledBddGeo} and the discussion following it that \emph{parameterized scaled b.g.\ structures} are families of scaled b.g.\ structures, parameterized by some parameter set $P$, where the bounds on transition functions are also uniform over $P$. In the present paper, they arise via Fourier (and Mellin) transforms as discussed in a general setup in \cite[\S{4.3}]{HintzScaledBddGeo}. Since we will need a slight generalization of the results proved in the reference, we recall the setup here. Given a scaled b.g.\ structure $\fB_{X,\times}=\{(U_\alpha,\phi_\alpha,\rho_\alpha)\colon\alpha\in\sA\}$ on a manifold $\cX$ and weight families $\{\rho_{\alpha,0}\colon\alpha\in\sA\}$ and $\{\tau_\alpha\colon\alpha\in\sA\}$, with $0<\rho_{\alpha,0}\leq 1$, we define a scaled b.g.\ structure $\fB_\times$ on $\cM:=\R_t\times\cX$ using the sets and scalings
\begin{equation}
\label{EqMSFCells}
  U_{(j,\alpha)}:=\Bigl((j-2)\frac{\tau_\alpha}{\rho_{\alpha,0}},(j+2)\frac{\tau_\alpha}{\rho_{\alpha,0}}\Bigr)_t \times U_\alpha,\quad
  \rho_{(j,\alpha),0}=\rho_{\alpha,0},\ \rho_{(j,\alpha),i}=\rho_{\alpha,i}.
\end{equation}
Write $\cW_\cX,\cV_\cX$, resp.\ $\cW,\cV$ for the coefficient and operator Lie algebras of $\fB_{X,\times}$, resp.\ $\fB_\times$. If $\tau$ and $\rho_0$ are weights on $\cX$ equivalent to the weight families $\{\tau_\alpha\}$ and $\{\rho_{0,\alpha}\}$, respectively, then $\cW$ is spanned by $\cW_\cX$ (lifted to $t$-invariant vector fields on $\cM$ tangent to each $t$-level set) and $\rho_0^{-1}\tau\pa_t$ over the space of $\cW$-regular functions; and $\cV$ is spanned by $\cV_\cX$ and $\tau\pa_t$.

Suppose $u\in L^2(\R\times\cX)$ is supported in $\R\times U_\alpha$. In local coordinates $t\in\R$ and $x=\phi_\alpha$ on $\R\times U_\alpha$, the membership $u\in H_\cV^1(\R\times\cX)$ is equivalent to the $L^2$-membership of $u,\tau_\alpha\pa_t u,\rho_{\alpha,i}\pa_{x^i}u$ and thus, by the Plancherel theorem, to $(1+\tau_\alpha|\sigma|)\hat u,\rho_{\alpha,i}\pa_{x^i}\hat u\in L^2$ where\footnote{We use the standard sign convention in spectral theory.}
\[
  \hat u(\sigma,p) = (\cF u)(\sigma,p) := \int_\R e^{-i\sigma t}u(t,p)\,\dd t,\quad \sigma\in\C,\ p\in\cX,
\]
or equivalently,
\[
  \hat u,\ \frac{\rho_{\alpha,i}}{\la\tau_\alpha\sigma\ra}\pa_{x^i}\hat u \in L^2\bigl(\R_\sigma; \la\tau_\alpha\sigma\ra^{-1}L^2(U_\alpha)\bigr).
\]
Following \cite[\S{4.3}]{HintzScaledBddGeo}, we thus introduce the parameterized scaled b.g.\ structure
\begin{equation}
\label{EqMSFParamScbg}
  \hat\fB_{\sigma,\times} = \{ (U_\alpha,\phi_\alpha,\rho_{\sigma,\alpha})\colon\alpha\in\sA\},\ 
  \rho_{\sigma,\alpha,i} := \frac{\rho_{\alpha,i}}{\la\tau_\alpha\sigma\ra},\quad \sigma\in\R,
\end{equation}
on $\cX$, with parameter space $\R$, and obtain an isomorphism
\begin{equation}
\label{EqMSFIsoV}
  \cF \colon H_\cV^s(\R\times\cX) \to L^2\bigl(\R_\sigma; \la\tau\sigma\ra^{-s}H_{\hat\cV_\sigma}^s(\cX)\bigr).
\end{equation}
Here we write $\hat\cV=(\hat\cV_\sigma)_{\sigma\in\R}$ for the operator Lie algebra of $(\cX,(\hat\fB_{\sigma,\times})_{\sigma\in\R})$, the elements of which on $U_\alpha$ are uniformly smooth linear combinations of $\frac{\rho_{\alpha,i}}{\la\tau_\alpha\sigma\ra}\pa_{x^i}$; see \cite[Proposition~4.19]{HintzScaledBddGeo}. For later use, we note that the coefficient Lie algebra is equal to $\cW$, independently of $\sigma$; and as a scaling weight, we may take $\bar{\hat\rho}:=\bar\rho_\cX\la\tau\sigma\ra^{-1}$ where $\bar\rho_\cX$ is any fixed scaling weight on $\cX$.
For the purposes of the present paper, we need to extend~\eqref{EqMSFIsoV} to an isomorphism for mixed spaces $H^{(s;k)}_{\cV;\cW}$. In the above local coordinate setting, $\cW$-regularity of $u$ amounts to regularity with respect to $\frac{\tau_\alpha}{\rho_{0,\alpha}}\pa_t$ and $\pa_{x^i}$, and hence on the Fourier transform side to regularity with respect to $\frac{\tau_\alpha\sigma}{\rho_{0,\alpha}}$ and $\pa_{x^i}$. Note that $\frac{\tau\sigma}{\rho_0}$ is a weight on $(\cX,(\hat\fB_{\sigma,\times})_{\sigma\in\R})$. If we write
\[
  H_{\hat\cV_\sigma;\hat\cW'_\sigma}^{(s;k)}(\cX)
\]
for the space of all $u\in H_{\hat\cV_\sigma}^s(\cX)$ such that $\la\frac{\tau\sigma}{\rho_0}\ra^{-1}W u\in H_{\hat\cV_\sigma}^s(\cX)$ for all $W\in\cW_X$ when $k=1$, similarly for $k\geq 2$, then we get the isomorphism
\begin{equation}
\label{EqMSFIsoVW}
  \cF \colon H_{\cV;\cW}^{(s;k)}(\R\times\cX) \to L^2\Bigl(\R_\sigma; \Big\la\frac{\tau\sigma}{\rho_0}\Big\ra^{-k}\la\tau\sigma\ra^{-s}H_{\hat\cV_\sigma;\hat\cW'_\sigma}^{(s;k)}(\cX)\Bigr).
\end{equation}
Variants of this (with $k=0$) with $t$-independent variable regularity and decay orders are described in~\cite[(4.25a)--(4.25b)]{HintzScaledBddGeo} and will be freely used here. Given a regularity order $\sfs$ on $\R\times\cX$, the regularity order one must use at frequency $\sigma$ is given by the restriction $\sfs_\sigma:=\sfs|_{-\sigma\,\dd t+T^*\cX}$; we call this the \emph{induced order}, similarly for decay orders (when applicable).

Finally, we recall that, given a time-translation invariant (smooth coefficient) differential operator $A\in\Diff_\cV^m(\R\times\cX)$, its spectral family is defined as
\begin{equation}
\label{EqMSFSpecFam}
  \hat A = (\hat A(\sigma))_{\sigma\in\R} \in \la\tau\sigma\ra^m \Diff_{\hat\cV}^m(\cX),\quad (\hat A(\sigma)u)(p)=\bigl(e^{i\sigma t}A(e^{-i\sigma t}u)\bigr)(0,p),\ u\in\CIc(\cX).
\end{equation}
Its $\hat\cV$-principal symbol is $\hat a=(\hat a(\sigma))_{\sigma\in\R}$ where
\begin{equation}
\label{EqMSFPhaseSpace}
  \hat a(\sigma;\varpi) = a(-\sigma\,\dd t+\varpi),\quad\varpi\in T^*\cX;
\end{equation}
here we identify $T^*\cX\cong\ann(\pa_t)\subset T^*_{t^{-1}(0)}\cM$. See \cite[Proposition~4.16]{HintzScaledBddGeo}.

\subsection{Ps.d.o.\ algebras, Sobolev spaces, transforms, normal operators}
\label{SsMU}

We proceed to introduce all ps.d.o.\ algebras utilized in the present paper (in addition to those already introduced in~\S\ref{SssMSbsc}); they all stem from the e3b-algebra introduced in~\S\ref{SssMUe3b} via localization or Fourier/Mellin transforms.

Several of these ps.d.o.\ algebras (e.g., the b- or the scattering-b-transition algebras) are, by themselves, not precise enough to allow for the proof of Fredholm or invertibility properties of differential operators of the same class. For example, what one needs in addition in the b-case on a manifold with boundary $\tilde X$ is information about the \emph{indicial roots}, or boundary spectrum, of a \emph{normal} (model) \emph{operator} at $\pa\tilde X$ \cite[\S{4.15}]{MelroseAPS}. We thus also recall the definitions of the relevant normal operators for smooth coefficient differential operators. (These normal operators themselves are then again differential elements of ps.d.o.\ algebras that are defined using a scaled b.g.\ structure.)

\subsubsection{e3b-setting}
\label{SssMUe3b}

We shall only present unit cells (and scalings) covering the set $\Omega=\ol{\{t_*\geq 1\}}\subset M\subset\tilde M$, $t_*=t-r$, from Definition~\ref{DefCMDomain}, intersected with $M^\circ$; we leave it to the interested reader to write down global covers of $\tilde M^\circ$. The coefficient and operator Lie algebras $\cW$ and $\cV$ will be $\CI_\bop\Vb(\tilde M)$ and $\CI_\bop\Vetb(\tilde M)$, respectively.

We shall use four families of unit cells (and scalings) which are compatible on their pairwise overlaps. We write the scaling of the $t$-coordinate in a unit cell $U_\alpha$ as $\rho_{\alpha,t}$, similarly for the $x^q$-scalings, $q=1,2,3$. The \emph{first family} covers a neighborhood of $\cK^+$. Writing $x^{\hat 1}=(x^2,x^3)$, $x^{\hat 2}=(x^1,x^3)$, and $x^{\hat 3}=(x^1,x^2)$, we define it as
\begin{equation}
\label{EqMUe3bFam1}
\begin{split}
  &U_{1,i} := (2^{i-2}-1,2^{i+2}-1)_t\times(-2,2)^3_x, \\
  &\quad\rho_{(1,i),t}=2^{-i},\ \rho_{(1,i),x^q}=1, \\
  &U_{1,i,\pm j k l} := \Bigl\{ 2^{i-2}-1<t<2^{i+2}-1,\ \pm x^j\in (2^{k-2},2^{k+2}),\ \frac{x^{\hat j}}{x^j} \in (-2,2)^2+l \Bigr\}, \\
  &\quad\rho_{(1,i,\pm j k l),t}=2^{k-i},\ \rho_{(1,i,\pm j k l),x^q}=1,
\end{split}
\end{equation}
where for $U_{1,i}$ we take $i\in\N_0$, while for $U_{1,i,\pm j k l}$ we take $j=1,2,3$, $\N_0\in i\geq k+10$, and $l\in\Z^2$, $|l|\leq 2$.\footnote{The shift of $t$ by $-1$ in the definition of $U_{1,i}$ ensures that $[0,\infty)\times(-2,2)^3_x\subset\bigcup_{i\in\N_0}U_{1,i}$. The lower bound $i\geq k+10$ ensures that the charts $U_{1,i,\pm j k l}$ cover a region with closure contained in $\ol{\{\frac{r}{t}\leq\frac12\}}\subset M$.} These charts are thus products of b-type charts on $\ol\R$ with b-type charts on $\ol{\R^3}$ as in~\eqref{EqMSbscCharts}. In the charts $U_{1,i}$, elements of $\cW$, resp.\ $\cV$ are linear combinations of $\la t\ra\pa_t$, $\pa_x$, resp.\ $\pa_{t_*},\pa_x$, with $\cW$-regular (i.e., conormal) coefficients; and in the charts $U_{1,i,\pm j k l}$, linear combinations of $t_*\pa_{t_*}$, $r\pa_x$, resp.\ $r\pa_{t_*},r\pa_x$, $r=|x|$. (Cf.\ \cite[\S{1.4.4}, Example~(11)]{HintzScaledBddGeo}.)

For the second and third families, we shall work with the coordinates
\[
  \tau = \frac{1}{t_*},\quad \rho_\cK:=\frac{1}{v}=\frac{r}{t_*},\quad x_\sscri:=\sqrt{v}=\sqrt{\frac{t_*}{r}};
\]
here $\rho_\cK$ and $x_\sscri$ are local defining functions of $\cK^+$ and $\scri^+$ in the regions $\rho_\cK<4$ and $x_\sscri<4$, respectively. The \emph{second family} then covers a neighborhood of $\iota^+\cap\{\rho_\cK\leq 3\}$. We define it in polar coordinates. Write $\hat e_1,\hat e_2,\hat e_3\in\Sph^2$ for the standard unit vectors on $\R^3$ and define
\begin{equation}
\label{EqMUe3bS2k}
  \Sph^2_k := \{\omega\in\Sph^2\colon|\omega-\hat e_k|<\pi\}.
\end{equation}
We then set
\begin{equation}
\label{EqMUe3bFam2}
\begin{split}
  &U_{2,i j k} := (2^{-i-2},2^{-i+2})_\tau \times (2^{-j-2},2^{-j+2})_{\rho_\cK} \times \Sph^2_k, \\
  &\quad \rho_{(2,i j k),\tau} = 2^{-j},\ \rho_{(2,i j k),\rho_\cK}=1,\ \rho_{(2,i j k),\omega}=1,
\end{split}
\end{equation}
where $i,j\in\N_0$ with $j\geq i$,\footnote{This means that we only consider charts on which $\rho_\cK\geq 2^{-4}\tau$, i.e., $r=\frac{\rho_\cK}{\tau}\geq 2^{-4}$. Regions of smaller $r$ are covered by~\eqref{EqMUe3bFam1}.} and $k=1,2,3$. In these charts, elements of $\cV$ are linear combinations, with conormal (on $M$) coefficients, of $\rho_\cK\tau\pa_\tau=-\rho_\cK r\pa_r-r\pa_{t_*}$, $\rho_\cK\pa_{\rho_\cK}=r\pa_r$, $\pa_\omega$, or equivalently of the basic 3b-vector fields $r\pa_{t_*}$, $r\pa_r$, $\pa_\omega$; and elements of $\cW$ are linear combinations of the b-vector fields $\tau\pa_\tau$, $\rho_\cK\pa_{\rho_\cK}$, $\pa_\omega$.

The \emph{third family} covers a neighborhood of $\iota^+\cap\{x_\sscri\leq 3\}$; it is
\begin{equation}
\label{EqMUe3bFam3}
\begin{split}
  &U_{3,i j k} := (2^{-i-2},2^{-i+2})_\tau \times (2^{-j-2},2^{-j+2})_{x_\sscri} \times \Sph^2_k, \\
  &\quad \rho_{(3,i j k),\tau} = 1,\ \rho_{(3,i j k),x_\sscri}=1,\ \rho_{(3,i j k),\omega}=2^{-j},
\end{split}
\end{equation}
where $i,j\in\N_0$ and $k=1,2,3$. Elements of $\cV$ are linear combinations, with conormal coefficients, of the basic edge-b-vector fields $\tau\pa_\tau=-t_*\pa_{t_*}-r\pa_r$, $x_\sscri\pa_{x_\sscri}=-2 r\pa_r$, $x_\sscri\pa_\omega$ (cf.\ \eqref{EqCTebFrame} where we can take $\rho_+=\tau$ in present notation), with b-regular coefficients; and $\cW$ is the space of b-vector fields with b-regular coefficients.

The \emph{fourth and final family} consists of a single unit cell which covers the subset of $\Omega\cap M^\circ$ that is not covered by these three families; this subset is compact in $\R^4$, and we can simply take this unit cell to be a sufficiently large cube in $\R^4$.

Complementing these families by analogous families covering the region $t_*\leq 1$ completes the description of a scaled b.g.\ structure on $M$ with $\cV=\CI_\bop\Vetb(\tilde M)$. The constructions recalled in~\S\ref{SsMS} then give rise to the following classes of operators and function spaces:\footnote{We only consider weights at $\scri^+$, $\iota^+$, and $\cK^+$ (with defining functions $x_\sscri,\rho_+$, and $\rho_\cK$), and do not explicitly record the possibility of weights also at $\rho_0$.}
\begin{enumerate}
\item Weighted e3b-ps.d.o.s
  \begin{equation}
  \label{EqMUe3bPsdo}
    \CI_\bop\Psi_\etbop^{\sfs,(2\alpha_\sscri,\alpha_+,\alpha_\cK)}(\tilde M)=x_\sscri^{-2\alpha_\sscri}\rho_+^{-\alpha_+}\rho_\cK^{-\alpha_\cK}\Psi_\etbop^\sfs(\tilde M)
  \end{equation}
  where $\alpha_\sscri,\alpha_+,\alpha_\cK\in\R$, and $\sfs\in\CI({}^\etbop S^*\tilde M)$ is a variable order function. For $\sfs=m\in\N_0$ and without weights for brevity, this contains the space $\CI_\bop\Diff_\etbop^m(\tilde M)$ of e3b-differential operators with b-regular coefficients.
\item Weighted e3b-Sobolev spaces
  \begin{subequations}
  \begin{equation}
  \label{EqMUe3bSpace0}
    H_\etbop^{\sfs,(2\alpha_\sscri,\alpha_+,\alpha_\cK)}(\tilde M) = x_\sscri^{2\alpha_\sscri}\rho_+^{\alpha_+}\rho_\cK^{\alpha_\cK}H_\etbop^\sfs(\tilde M),\quad H_\etbop^\sfs(\tilde M):=H_\etbop^\sfs(\tilde M,|\dd\ubar g|).
  \end{equation}
  (That is, we use the Minkowskian volume density $|\dd\ubar g|=|\dd t\,\dd x|$, or smooth positive multiples thereof such as the Kerr volume density $|\dd g_{\bhm,a}|$, unless otherwise specified.) We also work with mixed $(\etbop;\bop)$-Sobolev spaces which encode additional $k\in\N_0$ degrees of b-regularity; these are special cases of~\eqref{EqMSSobolev} and denoted
  \begin{equation}
  \label{EqMUe3bSpace}
    H_{\etbop;\bop}^{(\sfs;k),(2\alpha_\sscri,\alpha_+,\alpha_\cK)}(\tilde M).
  \end{equation}
  \end{subequations}
\end{enumerate}
As the locus of microlocal analysis, we take the boundary at fiber infinity
\[
  {}^\etbop S^*\tilde M \subset \ol{{}^\etbop T^*}\tilde M.
\]
We shall only need to consider variable order functions $\sfs\in\CI({}^\etbop S^*\tilde M)$; all choices of extensions of such a function to an element of $\CI_\bop S^0({}^\etbop T^*\tilde M)$ lead to the same ps.d.o.\ spaces~\eqref{EqMUe3bPsdo}.

Commutators of e3b-ps.d.o.s and b-vector fields feature prominently in our analysis; in view of Lemma~\ref{LemmaCTe3bComm} and our usage of \cite[\S{4.4}]{HintzScaledBddGeo}, we have, as a special case of~\eqref{EqMSComm},
\begin{equation}
\label{EqMUe3bComm}
  A\in\CI_\bop\Psi_\etbop^s(\tilde M),\ V\in\cV_{\bop,[\etbop]}(\tilde M) \implies [V,A]\in\CI_\bop\Psi_\etbop^s(\tilde M),\ \WF'_\etbop([V,A])\subset\WF'_\etbop(A),
\end{equation}
in the notation of Lemma~\ref{LemmaCTe3bComm}. It is often convenient to be able to expand an e3b-ps.d.o.\ in terms of b-vector fields; for this purpose, we note:

\begin{lemma}[Expanding e3b-ps.d.o.s in terms of b-vector fields]
\label{LemmaMUe3bExpb}
  Let $\sV=\{V_1,\ldots,V_N\}\subset\Vb(\tilde M)$ be such that there exist smooth weights $w_j\in\CI(\tilde M)$ with the property that $V_{0,j}:=w_j V_j\in\Vetb(\tilde M)$ and $\{V_{0,1},\ldots,V_{0,N}\}$ spans $\Vetb(\tilde M)$ over $\CI(\tilde M)$. Put $V_0:=I$. Then for every $A\in\CI_\bop\Psi_\etbop^s(\tilde M)$, there exist $A_j\in\CI_\bop\Psi_\etbop^{s-1}(\tilde M)$, $j=0,\ldots,N$, such that
  \[
    A = \sum_{j=0}^N A_j V_j.
  \]
\end{lemma}

We shall apply this to a \emph{good spanning set} $\sV$ of $\cV_{\bop,[\etbop]}(\tilde M)$ in the terminology of Lemma~\ref{LemmaCTe3bComm}.

\begin{proof}[Proof of Lemma~\usref{LemmaMUe3bExpb}]
  It suffices to show that one can write $A=\sum_{j=0}^N A_{0,j}V_{0,j}$ for suitable $A_{0,j}\in\CI_\bop\Psi_\etbop^s(\tilde M)$ since multiplication by the bounded weights $w_j$ maps $\CI_\bop\Psi_\etbop^s(\tilde M)$ into itself. But since $B:=\sum_{j=1}^N V_{0,j}^*V_{0,j}\in\Diff_\etbop^2(\tilde M)$ is elliptic, we can write $I=Q B+R$ for some $Q\in\CI_\bop\Psi_\etbop^{-2}$, $R\in\CI_\bop\Psi_\etbop^{-\infty}$, and thus
  \[
    A = A Q B + A R = \sum_{j=1}^N A_{0,j} V_{0,j} + A R V_0
  \]
  where $A_{0,j}=A Q V_{0,j}^*\in\CI_\bop\Psi_\etbop^{s-1}$.
\end{proof}

If we use the secondary b.g.\ structure, which is given by subdividing the unit cells~\eqref{EqMUe3bFam1} into $\rho$-cuboids, then we obtain larger classes of weighted e3b-ps.d.o.s with e3b-regular coefficients, denoted
\begin{equation}
\label{EqMUe3bTilde}
  \tilde\Psi_\etbop^{\sfs,(2\alpha_\sscri,\alpha_+,\alpha_\cK)}(\tilde M) = \CI_\etbop\Psi_\etbop^{\sfs,(2\alpha_\sscri,\alpha_+,\alpha_\cK)}(\tilde M).
\end{equation}
(Here $\CI_\etbop$ indicates the regularity of the coefficients of the underlying symbols: they remain bounded under application of any finite number of e3b-differential operators.) Operators of this class still define bounded maps between weighted e3b-Sobolev spaces, and the considerations around Lemma~\ref{LemmaMSOpNorm} apply to them.

\subsubsection{Fourier transform near \texorpdfstring{$\cK^+$}{the Kerr face}: sc-b-transition and semiclassical sc-settings}
\label{SssMUK}

In the sequel, we work in a neighborhood of $\cK^+$, which (over $\tilde M^\circ$) is covered by the unit cells~\eqref{EqMUe3bFam1}. As unit cells for the secondary b.g.\ structure, we can take a subset of the unit cells
\begin{align*}
  U_i &:= (i-2,i+2)_t \times (-2,2)^3_x, \\
  U_{i,\pm j k l} &:= \Bigl\{ (i-2)2^k<t<(i+2)2^k,\ \pm x^j\in(2^{k-2},2^{k+2}),\ \frac{x^{\hat j}}{x^j}\in(-2,2)^2+l \Bigr\},
\end{align*}
where $i\in\Z$, $j=1,2,3$, $k\in\N_0$, $l\in\Z^2$, $|l|\leq 2$. These unit cells on $\R^4=\R_t\times\R_x^3$ are, in turn, of the product form~\eqref{EqMSFCells} relative to the b-structure on $\R^3$ (cf.\ \eqref{EqMSbscCharts}) if we set $\tau_{\pm j k l}=2^k\sim\la x\ra$ and $\rho_{0,(\pm j k l)}=1$. (This is a geometric observation related to the existence of time-translation-invariant models for 3b-operators discussed in~\eqref{EqCPXA}--\eqref{EqCPXA0}.) The Lie algebra of uniform vector fields is equal to $\cC_\tbop^\infty\Vtb(\tilde M_0)$. The associated spaces of pseudodifferential operators as well as Sobolev spaces are denoted by
\begin{equation}
\label{EqMUKPsdo}
\begin{gathered}
  \tilde\Psi_\tbop^{s,(\alpha_\sface,\alpha_\cK)}(\tilde M_0) = \rho_\sface^{\alpha_\sface}\rho_\cK^{\alpha_\cK}\tilde\Psi_\tbop^s(\tilde M_0),\quad \tilde\Psi_\tbop^s=\cC_\tbop^\infty\tilde\Psi_\tbop^s, \\
  H_\tbop^{s,(\alpha_\sface,\alpha_\cK)}(\tilde M_0) = \rho_\sface^{\alpha_\sface}\rho_\cK^{\alpha_\cK}H_\tbop^s(\tilde M_0),
\end{gathered}
\end{equation}
where $\rho_\sface,\rho_\cK$ are defining functions of $\sface$ and $\cK\subset\tilde M_0$, respectively. The underlying $L^2$-space is defined with respect to the Minkowskian volume density, unless otherwise specified.

On the Fourier transform (in $t$) side, this gives rise to a parameterized scaled b.g.\ structure $\{(U_\alpha,\phi_\alpha,\rho_{\sigma,\alpha})\}$ on $\R^3$, with parameter space $\R$, underlying b.g.\ structure $\{(U_\alpha,\phi_\alpha)\}$ the b-structure on $\R^3$, and scalings $\rho_{\sigma,\alpha,i}=\frac{1}{1+|\sigma|\la x\ra}$. The coefficient Lie algebra is $\CI_\bop\Vb(\ol{\R^3})$ for all $\sigma$. The operator Lie algebras $\hat\cV_\sigma$, for various parameters and parameter ranges, are as follows.

{\bf (1)} For $\sigma=0$, we get the b-setting: $\hat\cV_0=\CI_\bop\Vb(\ol{\R^3})$. On the level of a time-translation-invariant 3b-operator $A_0$ as in~\eqref{EqCPXA0gl},\footnote{We only focus on the structure of $\wh{A_0}(\sigma)$ near $r=\infty$ and thus do not pay attention to the fact that the coordinates used in the description of $A_0$ in~\eqref{EqCPXA0gl} are not defined near $r=0$.} this corresponds to considering the spectral family,
  \begin{equation}
  \label{EqMUSpecFamEx}
    \wh{A_0}(\sigma) = \sum_{j+k+|\alpha|\leq m} a_{j k\alpha}(r^{-1},0,\omega)(-i\sigma r)^j(r\pa_r)^k\pa_\omega^\alpha,
  \end{equation}
  and observing that this is a b-differential operator when $\sigma=0$.

{\bf (2)} For $\sigma\neq 0$, we get the scattering setting on $\ol{\R^3}$ with b-regular coefficients: $\hat\cV_\sigma=\CI_\bop\Vsc(\ol{\R^3})$. In particular, the principal symbol map on the associated class of ps.d.o.s captures operators to leading differential and decay (in $\la x\ra^{-1}$) orders. This corresponds to $\wh{A_0}(\sigma)$, upon division by $r^{-m}$, being a scattering differential operator on $\ol{\R^3}$ (i.e., it is built from $\pa_r$, $r^{-1}\pa_\omega$).

\begin{figure}[!ht]
  \centering
  \includegraphics{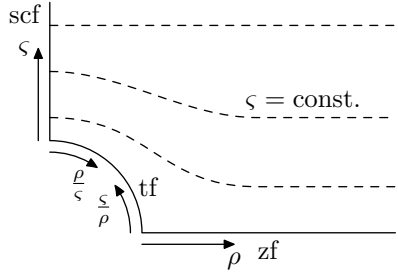}
  \caption{The scattering-b-transition single space $\tilde X_\scbtop$, its boundary hypersurfaces, local coordinates, and some positive level sets of $\varsigma$ (dashed). We only show the region $\rho=r^{-1}\ll 1$ here.}
  \label{FigMUscbt}
  \end{figure}

{\bf (3)} The transition from zero to nonzero $\sigma$ is given by the \emph{scattering-b-transition} setting (with b-regular coefficients), which was first introduced by Guillarmou--Hassell \cite{GuillarmouHassellResI} (under a different name) and Melrose--S\'a Barreto \cite{MelroseSaBarretoLow}, and later extended by the author for low-energy resolvent analysis \cite{HintzKdSMS,Hintz3b,HintzNonstat}. The compactified smooth coefficient version is as follows: writing $\tilde X=\ol{\R^3}$ for brevity, define
  \begin{equation}
  \label{EqMUscbtSingle}
    \tilde X_\scbtop := [ [0,1]_\varsigma \times \tilde X; \{0\}\times\pa\tilde X ],
  \end{equation}
  and write $\scface$, $\tface$, and $\zface$ for the lifts of $[0,1]\times\pa\tilde X$ (the \emph{scattering face}, $\{0\}\times\pa\tilde X$ (the \emph{transition face}, which is thus the front face), and $\{0\}\times\tilde X$ (the \emph{zero face}), respectively. See Figure~\ref{FigMUscbt}.

  Then
  \[
    \Vscbt(\tilde X) := \{ V\in\rho_\scface\Vb(\tilde X_\scbtop)\colon V\ \text{is tangent to every $\varsigma$-level set} \}.
  \]
  Here $\rho_\scface\in\CI(\tilde X_\scbtop)$ is a defining function of $\scface$; a possible choice is $(1+\varsigma\la x\ra)^{-1}$. A spanning set of $\Vscbt(\tilde X)$ is thus given by $\rho_\scface\la x\ra\pa_x$. We then have
  \[
    (\hat\cV_\sigma)_{\sigma\in\pm[0,1]} = \CI_\bop\Vscbt(\tilde X),
  \]
  where we identify $\varsigma=\pm\sigma$. (In this parameterized setting, we write $\CI_\bop$ for uniformly bounded, in $\sigma\in\pm[0,1]$, families of elements of $\CI_\bop(\tilde X)$. There are no regularity requirements in $\sigma$.) The corresponding class of \emph{scattering-b-transition ps.d.o.s} (with uniformly conormal coefficients) is denoted
  \begin{equation}
  \label{EqMUscbt}
    \CI_\bop\Psiscbt^{s,(r,l,q)}(\tilde X) = \rho_\scface^{-r}\rho_\tface^{-l}\rho_\zface^{-q}\CI_\bop\Psiscbt^s(\tilde X);
  \end{equation}
  here $s$ is the regularity order, $r$ the scattering decay order, and $l$ the decay order at the transition face. The orders $s,r$ may be variable. The principal symbol captures an operator to leading order in the regularity and scattering decay senses. One can define a vector bundle $\Tscbt\tilde X\to\tilde X_\scbtop$, the space of smooth sections of which is precisely $\Vscbt(\tilde X)$: a local frame of this bundle is given by $\rho_\scface\la x\ra\pa_x$. The microlocalization locus is thus the union
  \begin{equation}
  \label{EqMUscbtTstar}
    {}^\scbtop S^*\tilde X \cup \ol{\Tscbt^*_\scface}\tilde X
  \end{equation}
  of fiber infinity of $\ol{\Tscbt^*}\tilde X$ and the compactified sc-b-transition phase space of $\scface$. The corresponding Sobolev spaces (with norms depending on the parameter $\varsigma\in[0,1]$) are denoted
  \begin{equation}
  \label{EqMUscbtSob}
    H_{\scbtop,\sigma}^{s,(r,l,q)}(\tilde X)=\rho_\scface^r\rho_\tface^l\rho_\zface^q H_{\scbtop,\sigma}^s(\tilde X).
  \end{equation}
  The orders $s$, $r$ corresponding to the boundary hypersurfaces~\eqref{EqMUscbtTstar} can be variable.

{\bf (4)} The regime of unbounded $\sigma$ is described by the \emph{semiclassical scattering} setting (with uniformly b-regular coefficients). This was first used by Vasy--Zworski \cite{VasyZworskiScl} to prove high-energy resolvent estimates. For the compactified smooth coefficient version, write $\tilde X=\ol{\R^3}$ and set
  \[
    \tilde X_\schop := [0,1]_h \times \tilde X,
  \]
  with boundary hypersurfaces $\sface=\{0\}\times\tilde X$ (the \emph{semiclassical face}) and $\scface=[0,1]\times\tilde X$ (the \emph{scattering face}). Defining functions are $\rho_\sface=h$ and $\rho_\scface=\rho$ (a boundary defining function of $\pa\tilde X$ on $\tilde X$). Let then
  \[
    \cV_\schop(\tilde X) := \{ V\in h\rho_\sface\Vb(\tilde X_\schop) \colon V\ \text{is tangent to every $h$-level set} \}.
  \]
  We then have
  \[
    (\hat\cV_\sigma)_{\sigma\in\pm[1,\infty)} = \CI_\bop\cV_\schop(\ol{\R^3}),
  \]
  where we identify $h=(\pm\sigma)^{-1}$. The corresponding class of \emph{semiclassical scattering ps.d.o.s} (with conormal coefficients) is denoted
  \begin{equation}
  \label{EqMUsch}
    \CI_\bop\Psi_\schop^{s,r,b}(\ol{\R^3}) = \rho^{-r}h^{-b}\CI_\bop\Psi_\schop^s(\ol{\R^3});
  \end{equation}
  here $s$ is the regularity order, $r$ the scattering decay order, and $b$ the semiclassical order. All orders may be variable. The principal symbol captures an operator to leading order in all three senses. The microlocalization locus is the union
  \begin{equation}
  \label{EqMUschTstar0}
    {}^\schop S^*\tilde X \cup \ol{{}^\schop T^*_\scface}\tilde X \cup \ol{{}^\schop T^*_\sface}\tilde X
  \end{equation}
  where ${}^\schop T\tilde X$ is the vector bundle over $\tilde X_{\scop,\semi}$ for which the space of smooth sections is exactly $\cV_\schop(\tilde X)$; and for the corresponding scale of Sobolev spaces (with $h$-dependent norm)
  \begin{equation}
  \label{EqMUschSob}
    H_{\scop,h}^{s,r,b} = \rho^r h^b H_{\scop,h}^s,
  \end{equation}
  all three orders may be variable. (On $\ol{\R^n}$, we have $\|u\|_{H_{\scop,h}^{s,r,b}}=h^{-b}\|\la h D\ra^s(\la x\ra^{-r}u)\|_{L^2}$.) In this paper (as is usual in semiclassical analysis, cf.\ \cite[\S{6.6}]{HintzMicro}) we only work with the third summand of~\eqref{EqMUschTstar0}, which is identified with
  \begin{equation}
  \label{EqMUschTstar}
    \ol{\Tsc^*}\tilde X.
  \end{equation}
  (This identification rescales a semiclassical scattering covector $h\xi$ to $\xi$.) On the level of the spectral family~\eqref{EqMUSpecFamEx}, this amounts to replacing $\sigma$ by $\pm h^{-1}$ and noting that $h^m r^{-m}\wh{A_0}(\pm h^{-1})$ is a semiclassical scattering differential operator (i.e., it is built from $h\pa_r$, $h r^{-1}\pa_\omega$).

\bigskip

It follows from~\eqref{EqMSFSpecFam} that if $A_0\in\rho_\sface^{-l}\Diff_\tbop^m(\tilde M_0)$, $\rho_\sface=\la x\ra^{-1}$, is time-translation-invariant, then for any $\hat\sigma\in\C$, $|\hat\sigma|=1$, and for any $\nu\in\R$, we have
\begin{equation}
\label{EqMUSpecFam}
  (\wh{A_0}(\hat\sigma|\sigma|))_{|\sigma|\in[0,1]} \in \rho_\scface^{-l-m}\rho_\tface^{-l}\Diff_\scbtop^m(\ol{\R^3}),\quad
  (\wh{A_0}(i\nu\pm h^{-1}))_{h\in(0,1]} \in \rho^{-l-m}h^{-m}\Diff_\schop^m(\ol{\R^3});
\end{equation}
here
\begin{equation}
\label{EqMUDiffsch}
  \Diff_\schop^m(\ol{\R^3})
\end{equation}
is the space of semiclassical scattering differential operators, i.e., sums of up to $m$-fold compositions of elements of $\cV_\schop(\ol{\R^3})$. Relatedly, as a special case of~\eqref{EqMSFIsoV}, we can characterize the norm on the weighted 3b-Sobolev space
\[
  \Htb^{\sfs,(\alpha_\sface,0)}(\tilde M_0) = \rho_\sface^{\alpha_\sface}\Htb^\sfs(\tilde M_0,|\dd t\,\dd x|),\quad \rho_\sface=\la x\ra^{-1},
\]
on the Fourier transform side, provided $\sfs\in\CI({}^\tbop S^*\tilde M_0)$ is translation-invariant. Namely,
\begin{align*}
  &\|u\|_{H_\tbop^{\sfs,(\alpha_\sface,0)}(\tilde M_0,|\dd t\,\dd x|)}^2 \\
  &\qquad \sim \sum_\pm \int_{\pm[0,1]} \|\hat u(\sigma)\|_{H_{\scbtop,\sigma}^{\sfs_\sigma,(\sfs_\sigma+\alpha_\sface,\alpha_\sface,0)}(\ol{\R^3},|\dd x|)}^2\,\dd\sigma + \sum_\pm \int_{\pm[1,\infty)} \|\hat u(\sigma)\|_{H_{\scop,|\sigma|^{-1}}^{\sfs_\sigma,\sfs_\sigma+\alpha_\sface,\sfs_\sigma}(\ol{\R^3},|\dd x|)}^2\,\dd\sigma;
\end{align*}
here ``$\sim$'' means that the left-hand side is bounded by a constant (independent of $u$) times the right-hand side and vice versa. (This was also proved in \cite[Proposition~4.24]{Hintz3b}.) The norms on $\hat u(\sigma)$ on the right are weighted $H_{\hat\cV_\sigma}$-norms, though expressed using the notation of Sobolev spaces matching~\eqref{EqMUscbt} and \eqref{EqMUsch}; the orders denoted $\sfs_\sigma$ on the right are those induced by $\sfs$ via the phase space relationship~\eqref{EqMSFPhaseSpace}. The weight at $\scface$ incorporates both the weight $\alpha_\sface$ at $\sface$ and the power of $\la\tau\varsigma\ra\sim 1+|\sigma|\la x\ra$ in~\eqref{EqMSFIsoV}. See also~\S\ref{SsSpTs}.

Now, the 3b-unit cells only coincide with the unit cells of the secondary b.g.\ structure for the e3b-structure defined in~\S\ref{SssMUe3b} in a neighborhood of $\cK^+$ (in fact, outside of any fixed neighborhood of $\scri^+$). Thus:

\begin{lemma}[Fourier transform on 3b-Sobolev spaces]
\label{LemmaMUetbFT}
  For translation-invariant $\sfs\in\CI({}^\tbop S^*\tilde M_0)$, and for all $\alpha_\sface\in\R$, we have
  \begin{align*}
    \|u\|_{H_\tbop^{\sfs,(\alpha_\sface,0)}(\tilde M_0,|\dd t\,\dd x|)}^2 & \sim \sum_\pm \int_{\pm[0,1]} \|\hat u(\sigma)\|_{H_{\scbtop,\sigma}^{\sfs_\sigma,(\sfs_\sigma+\alpha_\sface,\alpha_\sface,0)}(\tilde X,|\dd x|)}^2\,\dd\sigma \\
    &\qquad\quad + \sum_\pm \int_{\pm[1,\infty)} \|\hat u(\sigma)\|_{H_{\scop,|\sigma|^{-1}}^{\sfs_\sigma,\sfs_\sigma+\alpha_\sface,\sfs_\sigma}(\tilde X,|\dd x|)}^2\,\dd\sigma.
  \end{align*}
\end{lemma}

The order $\alpha_\sscri$ plays no role since $\scri^+\cap\supp\chi=\emptyset$. We shall apply this to the Fourier transform with respect to (the slicing of spacetime by level sets of) $t_*$ (for functions supported away from $r=0$), in which case the phase space relationship for variable orders in~\eqref{EqMSFPhaseSpace} with $\dd t_*$ in place of $\dd t$.

We also need spaces with additional degrees of ``b-regularity;'' we use quotes since, as the sole exception to the general notation~\eqref{EqMSSobolev} in this paper, we do \emph{not} mean regularity with respect to b-vector fields on $\tilde M_0$ here, but rather b-vector fields on its resolution $\tilde M_1$. The space of these b-vector fields is spanned by $\la t_*\ra\pa_{t_*}$ and $\la x\ra\pa_x$ as a consequence of Lemma~\ref{LemmaCPXDiffeo} and~\eqref{EqCPXM1p}.

\begin{rmk}[$t$ vs. $t_*$]
\label{RmkMUttstar}
  Since in the present section we work with (the Fourier transform in) $t$, we really use $\la t\ra\pa_t$ here. But when, in the bulk of the paper, we work with (the Fourier transform in) $t_*$, all results stated here hold \emph{mutatis mutandis}, so in particular with $\la t_*\ra\pa_{t_*}$ instead of $\la t\ra\pa_t$ (or, more conveniently, with $t_*\pa_{t_*}$ and $\pa_{t_*}$ instead of $t\pa_t$ and $\pa_t$).
\end{rmk}

Thus, we define
\begin{equation}
\label{EqMUtbb}
  H_{\tbop;\bop}^{(\sfs;k)}(\tilde M_0) := \bigl\{ u\in H_\tbop^\sfs(\tilde M_0) \colon \pa_t^{j_1}(t\pa_t)^{j_2}(\la x\ra\pa_x)^\alpha u\in H_\tbop^\sfs(\tilde M_0)\ \forall\,j_1+j_2+|\alpha|\leq k \bigr\}.
\end{equation}
While b-regularity in the spatial variables commutes with the Fourier transform, b-regularity in $t$ amounts to b-regularity in $\sigma$ in view of $\cF\circ t\pa_t=-\pa_\sigma\sigma\circ\cF$. (This is thus evidently different from~\eqref{EqMSFIsoVW}.) Let us use the notation
\begin{equation}
\label{EqMUscbtbSob}
  H_{(\scbtop,\sigma);\bop}^{(\sfs;k),(\sfr,\sfl,q)}
\end{equation}
for weighted sc-b-transition Sobolev spaces with $k\in\N_0$ additional degrees of b-regularity, i.e., regularity with respect to $\la x\ra\pa_x$. (Without weights, this \emph{is} thus the space $H_{\hat\cV_\sigma;\cW}^{(\sfs;k)}(\R^3)$ in the notation~\eqref{EqMSMixed}).) We furthermore write
\begin{equation}
\label{EqMUHscbh}
  H_{(\scop,|\sigma|^{-1});\bop}^{(\sfs;k),\sfr,\sfb}
\end{equation}
for weighted semiclassical scattering Sobolev spaces with $k$ additional degrees of \emph{non-semiclassical} b-regularity (corresponding to the fact that the coefficient Lie algebra is $\CI_\bop\Vb$).\footnote{Using semiclassical b-regularity might seem more natural, but makes the powers of $|\sigma|^{-1}$ in Lemma~\ref{LemmaMUetbFTb} below less transparent. Since b- and semiclassical b-vector fields differ only by factor of $h$, the translation is straightforward in principle.} We will later on also such spaces for bounded $\sigma$, in which case there is no semiclassical order; we write this space as
\begin{equation}
\label{EqMUHscb}
  H_{\scop;\bop}^{(\sfs;k),\sfr}.
\end{equation}
(This is thus the space $w H_{\cV;\cW}^{(\sfs;k)}$ where $w=\la x\ra^{-\sfr}$, and $\cV$, resp.\ $\cW$ is the Lie algebra of scattering vector fields, resp.\ b-vector fields on $\ol{\R^3}$ with conormal coefficients.) Lemma~\ref{LemmaMUetbFT} now generalizes to:

\begin{lemma}[Fourier transform on 3b-Sobolev spaces with b-regularity]
\label{LemmaMUetbFTb}
  Let $\sfs$, and the induced orders $\sfr,\sfl,\sfb$, be as in Lemma~\usref{LemmaMUetbFT}. Let $\alpha_\sface\in\R$. Then we have
  \begin{align*}
    \|u\|_{H_{\tbop;\bop}^{(\sfs;k),(\alpha_\sface,0)}(\tilde M_0)}^2 &\sim \sum_\pm\sum_{j=0}^k \int_{\pm[0,1]} \|(\sigma\pa_\sigma)^j\hat u(\sigma)\|_{H_{(\scbtop,\sigma);\bop}^{(\sfs_\sigma;k-j),(\sfs_\sigma+\alpha_\sface,\alpha_\sface,0)}(\tilde X,|\dd x|)}^2\,\dd\sigma \\
      &\quad\quad + \sum_\pm\sum_{j=j_1+j_2\leq k} \int_{\pm[1,\infty)} \|\sigma^{j_1}(\sigma\pa_\sigma)^{j_2}\hat u(\sigma)\|_{H_{(\scop,|\sigma|^{-1});\bop}^{(\sfs_\sigma;k-j),\sfs_\sigma+\alpha_\sface,\sfs_\sigma}(\tilde X,|\dd x|)}^2\,\dd\sigma.
  \end{align*}
\end{lemma}

Lemmas~\ref{LemmaMUetbFT} and \ref{LemmaMUetbFTb} are immediately useful in the study of PDEs $P_0 u=f$ near $\cK^+$ when $P_0$ is a stationary, i.e., $t$-translation-invariant, 3b-operator on $\tilde M_0$, as then $\cF(P_0 u)=\wh{P_0}(\sigma)\hat u(\sigma)$ where $\wh{P_0}(\sigma)$ is the spectral family~\eqref{EqMSFSpecFam} (which is thus a sc-b-transition differential operator for $\pm\sigma\in[0,1]$ and a semiclassical scattering differential operator for $\pm\sigma\in[1,\infty)$). The wave-type operators $P$ of interest in this paper are perturbations of stationary operators $P_0$ by operators whose coefficients (as e3b-operators) decay at $\cK^+$, and we will then apply Lemma~\ref{LemmaMUetbFT} to the stationary model problem $P_0 u=f$. For the construction of the stationary $\cK^+$-model in general, see \cite[\S{4.1}]{Hintz3b}.

In the semiclassical scattering algebra, principal symbol based arguments are strong enough to prove Fredholm and invertibility properties for ps.d.o.s. (This algebra is \emph{fully symbolic} and semiclassical, cf.\ \cite[Definition~3.27 and Footnote~14]{HintzScaledBddGeo}.) This is not true in the sc-b-transition algebra uniformly as $|\sigma|\to 0$; instead, one needs to invert \emph{normal operators} at $\zface$ and $\tface$. To define these, consider a weighted sc-b-transition differential operator
\[
  A \in \Diff_\scbtop^{m,(r,l,q)}(\ol{\R^3}) = \rho_\scface^{-r}\rho_\tface^{-l}\rho_\zface^{-q}\Diff_\scbtop^m(\ol{\R^3}).
\]
Working locally near $\pa\ol{\R^3}$ with coordinates $\rho\geq 0$, $\omega\in\R^2$, this is thus the product of $\rho_\scface^{-r}\rho_\tface^{-l}\rho_\zface^{-q}$ and up to $m$-fold compositions of the basic sc-b-transition vector fields $\rho_\scface\rho\pa_\rho$, $\rho_\scface\pa_\omega$; we can locally take $\rho_\scface=\frac{\rho}{\rho+|\sigma|}$. For $q=0$, we can thus restrict the coefficients to $\zface$ (i.e., formally letting $\sigma\to 0$ for all fixed $\rho>0$) and obtain an $m$-fold composition of $\rho\pa_\rho,\pa_\omega$, i.e., a b-differential operator. Since at $\zface$ we have $\rho_\scface\sim 1$ and $\rho_\tface\sim\rho$, this defines the \emph{$\zface$-normal operator}
\begin{subequations}
\begin{equation}
\label{EqMUNzf}
  N_\zface(A) \in \Diffb^{m,l}(\ol{\R^3}) = \rho^{-l}\Diffb^m(\ol{\R^3}),\quad A\in\Diffscbt^{m,(r,l,0)}(\ol{\R^3}).
\end{equation}
For $l=0$ on the other hand, we can restrict the coefficients of $A$ to $\tface$ (i.e.,\ formally letting $\sigma\to 0$ for all fixed ratios $\frac{\rho}{|\sigma|}$); now
\begin{equation}
\label{EqMUtf}
  \tface = [0,\infty]_{\hat\rho}\times\Sph^2,\quad \hat\rho=\frac{\rho}{|\sigma|},
\end{equation}
and thus $\rho_\scface=\frac{\hat\rho}{1+\hat\rho}$ is a defining function of $\sctface:=\tface\cap\scface$. The basic sc-b-transition vector fields are thus $\frac{\hat\rho}{1+\hat\rho}\hat\rho\pa_{\hat\rho}$ and $\frac{\hat\rho}{1+\hat\rho}\pa_\omega$, i.e., scattering vector fields near $\hat\rho=0$ and b-vector fields near $\hat\rho^{-1}=0$. We thus obtain the \emph{$\tface$-normal operator}
\begin{equation}
\label{EqMUNtf}
  N_\tface(A) \in \Diff_{\scop,\bop}^{m,(r,q)}(\tface),\quad A\in\Diffscbt^{m,(r,0,q)}(\ol{\R^3}).
\end{equation}
\end{subequations}
The Sobolev spaces appropriate for the analysis of such scattering-b-operators are
\begin{equation}
\label{EqMUNtfSob}
  H_{\scop,\bop}^{s,(r,q)}(\tface) = \rho_\sctface^r\rho_\ztface^q H_{\scop,\bop}^s(\tface).
\end{equation}
These are a special case of the scaled bounded geometry Sobolev space $w H_\cV(\tface)$ for the weight $w=\rho_\sctface^r\rho_\ztface^q$ and the operator Lie algebra $\cV=\CI_\bop\cV_{\scop,\bop}(\tface)=\rho_\sctface\CI_\bop\Vb(\tface)$ with b-regular coefficients; we may take the differential and scattering decay orders to be variable.

\begin{rmk}[Model problems of model problems]
\label{RmkMUNtfModel}
  Since the b-algebra is not fully symbolic, the analysis of both $N_\zface(A)$ and $N_\tface(A)$ requires the study of a further model problem at $\ztface:=\zface\cap\tface$; we discuss this briefly in~\S\ref{SssMUb}.
\end{rmk}

Applying these considerations to the spectral family of $A_0$ in~\eqref{EqMUSpecFamEx}, one finds (using $\rho=r^{-1}$ and $r\pa_r=-\rho\pa_\rho$)
\begin{align}
  N_\zface\bigl((\wh{A_0}(\sigma))_{\sigma\in\pm[0,1]}\bigr) &= \wh{A_0}(0), \nonumber\\
\label{EqMUNtfExpl}
  N_\tface\bigl((\wh{A_0}(\sigma))_{\sigma\in\pm[0,1]}\bigr) &= \sum_{j+k+|\alpha|\leq m} a_{j k\alpha}(0,0,\omega) (-i\hat\rho)^{-j}(-\hat\rho\pa_{\hat\rho})^k \pa_\omega^\alpha.
\end{align}

Corresponding to the above identifications of operator Lie algebras, sc-b-transition Sobolev spaces are related to b-Sobolev spaces near $\zface$ and to scattering-b-Sobolev spaces near $\tface$. Concretely, if $\chi_\zface$ and $\chi_\tface\in\CI(\tilde X_\scbtop)$ are equal to $1$ near $\zface$ and $\tface$, respectively, and are supported in a small neighborhood thereof, then we have uniform norm equivalences
\begin{equation}
\label{EqMUscbtNormEquiv}
  \|\chi_\zface u\|_{H_{\scbtop,\sigma}^{s,(r,l,q)}(\tilde X)} \sim |\sigma|^{-q}\|\chi_\zface u\|_{H_\bop^{s,l-q}(\tilde X)},\quad
  \|\chi_\tface u\|_{H_{\scbtop,\sigma}^{s,(r,l,q)}} \sim |\sigma|^{\frac32-l}\|\chi_\tface u\|_{H_{\scop,\bop}^{s,(r,q-l)}(\tface)}.
\end{equation}
Here, we use the Euclidean volume densities on $\tilde X$ and $\tface$ to define the underlying $L^2$-spaces. The factor of $|\sigma|^{\frac32}$ in the second norm equivalence arises from the fact that these densities are related, using the projective coordinate $\hat r=r|\sigma|$ on $\tface$, via $r^2\,|\dd r\,\dd\slg|=(|\sigma|^{\frac32})^2\,\hat r^2\,|\dd\hat r\,\dd\slg|$. (See \cite[Appendix~A.4]{HintzKdSMS}, and also \cite[Proposition~2.21]{Hintz3b}.) The first, resp.\ second norm equivalence also holds for variable orders $s,r$, in which case they are then $\sigma$-dependent on the right-hand sides.

\subsubsection{Mellin transform near a boundary in the b-setting}
\label{SssMUb}

On $\cM=\R_x^3$, using the notation~\eqref{EqMUe3bS2k}, and setting $\rho=r^{-1}$, $r=|x|$, we can use the unit cells
\[
  U_{i k} := (2^{-i-2},2^{-i+2})_\rho \times \Sph^2_k,\quad i\in\N_0,\ k=1,2,3,
\]
instead of $U_{\pm i j k}$ in~\eqref{EqMSbscCharts}. Extending the range of $i$ to all of $\Z$ gives a bounded geometry structure on $(0,\infty)\times\Sph^2$, with space of uniform vector fields equal to $\CI_\bop\Vb([0,\infty]\times\Sph^2)$. Passing to $t:=-\!\log_2\rho$, these unit cells are the special case $\tau_\alpha=\rho_{0,\alpha}=1$ of the setup~\eqref{EqMSFCells} where $\{U_\alpha\}=\{\Sph^2_k\colon k=1,2,3\}$. The Fourier transform in $t$ is, up to a constant scaling of the Mellin-dual variable, the Mellin transform in $\rho$, defined here by
\[
  (\cM u)(\sigma,\omega) = \int_0^\infty \rho^{i\sigma}u(\rho,\omega)\,\frac{\dd\rho}{\rho},\quad\sigma\in\C,\ \omega\in\Sph^2.
\]
The construction~\eqref{EqMSFParamScbg} yields a parametrized scaled b.g.\ structure on $\Sph^2$, with parameter $\sigma\in\R$; the operator Lie algebra $(\hat\cV_\sigma)_{\sigma\in\R}$ is, for bounded $\sigma$, simply given by $\cV(\Sph^2)$, while the transition from finite to infinite $\sigma$ is described by $\CI\cV_\semi(\Sph^2)$; here
\[
  \cV_\semi(\Sph^2) := \{ h V \colon V\in\cV([0,1)_h\times\Sph^2) \}
\]
is the space of semiclassical vector fields on $\Sph^2$, and we write $\CI\cV_\semi(\Sph^2)$ for the space of linear combinations of these with uniformly (in $h$) smooth (on $\Sph^2$) functions (with no regularity requirements in the parameter $h$). The corresponding ps.d.o.\ algebras are the standard algebra $\Psi^s(\Sph^2)$ and the semiclassical algebra
\[
  \Psih^{s,b}(\Sph^2) = h^{-b}\Psih^s(\Sph^2);
\]
both orders $s$ and $b$ can be variable. The corresponding semiclassical Sobolev space, at parameter value $h$, is denoted $H_h^{s,b}(\Sph^2)=h^b H_h^s(\Sph^2)$.

Beyond the isomorphism~\eqref{EqMSFIsoV}, we can also consider polynomial weights in $r$ (i.e., exponential weights in $\rho$), which amounts to shifting the Mellin-dual variable by a purely imaginary amount. We thus obtain the following well-known statement (see, e.g., \cite[\S{3.1}]{VasyMicroKerrdS}):

\begin{lemma}[Mellin transform near the boundary of $\ol{\R^3}$]
\label{LemmaMUb}
  Let $\chi\in\CI(\ol{\R^3})$ with $0\notin\supp\chi$. Let $\alpha\in\R$. Then we have an equivalence of norms\footnote{We could write this entirely in terms of semiclassical Sobolev spaces using the semiclassical parameter $\la\sigma\ra^{-1}$. We prefer not to do this for purely aesthetic reasons; after all, small $\sigma$ should not be considered to lie in the semiclassical regime.}
  \begin{align*}
    \|\chi u\|_{\Hb^{s,\alpha}(\ol{\R^3},|\dd x|)}^2 &\sim \int_{-1}^1 \|\cM(\chi u)(\sigma-i(\alpha+\tfrac32))\|_{H^s(\Sph^2)}^2 \\
      &\qquad + \sum_\pm \int_{\pm[1,\infty)} \|\cM(\chi u)(\sigma-i(\alpha+\tfrac32))\|_{H_{|\sigma|^{-1}}^{s,s}(\Sph^2)}^2.
  \end{align*}
\end{lemma}
\begin{proof}
  We only comment on the weights: if we pass from $|\dd x|$ to a b-density such as $\frac{|\dd x|}{\la x\ra^3}$, the weight $\alpha$ is replaced by $\alpha+\frac32$. Near $\rho=0$, any b-density on $\ol{\R^3}$ is a smooth positive multiple of $|\frac{\dd\rho}{\rho}\,\dd\slg|$, and thus of $|\dd t\,\dd\slg|$ if we set $t:=-\log_2\rho$ as before.
\end{proof}

Smooth coefficient b-differential operators $A\in\Diffb^m(\ol{\R^3})$ have a dilation-invariant normal operator at $\pa\ol{\R^3}$: in a collar neighborhood $[0,\eps)_\rho\times\Sph^2$ of $\pa\ol{\R^3}$, write
\[
  A = \sum_{j=0}^m A_j(\rho) (\rho\pa_\rho)^j,\quad A_j\in\CI([0,\eps);\Diff^{m-j}(\Sph^2));
\]
then the normal operator is
\[
  N(A) := \sum_{j=0}^m A_j(0)(\rho\pa_\rho)^j.
\]
Formally conjugating it by the Mellin transform amounts to replacing $\rho\pa_\rho$ by $-i\sigma$ and gives rise to the \emph{indicial family}
\[
  N(A,\sigma) = \sum_{j=0}^m A_j(0)(-i\sigma)^j.
\]
As a special case of~\eqref{EqMSFSpecFam}, we remark that, for any fixed value $\nu\in\R$, the family $(N(A,-i\nu+\sigma))_{\sigma\in\R}$ defines an element $h^{-m}\Diff_\semi^m(\Sph^2)$ with semiclassical parameter $h=\la\sigma\ra^{-1}$.

\begin{definition}[Indicial roots and indicial gaps]
\label{DefMUbInd}
  An \emph{indicial root} of $A$ is a complex number $\lambda$ for which there exists $u\in\CI(\Sph^2)$ such that $N(A)(\rho^\lambda u)=0$. An \emph{indicial gap} of $A$ is a interval $(\beta_-,\beta_+)$ such that $\Re\lambda\notin(\beta_-,\beta_+)$ for all indicial roots $\lambda$ of $A$. For $A\in\rho^{-l}\Diffb^m(\ol{\R^3})$, we define indicial roots and indicial gaps for $A$ as those of $\rho^l A\in\Diffb^m(\ol{\R^3})$.
\end{definition}

For $A$ which are elliptic near $\pa\ol{\R^3}$, the operator $N(A)$ is elliptic, and thus $(N(A,-i\nu+\sigma))_{\sigma\in\R}$ is an elliptic semiclassical operator. The standard semiclassical elliptic parametrix construction thus implies that $N(A,-i\nu+\sigma)$ is invertible as a map $H^s(\Sph^2)\to H^{s-m}(\Sph^2)$ when $|\Re\sigma|\geq C(\nu)$, and hence, by the analytic Fredholm theorem, the inverse $N(A,\sigma)^{-1}$ is meromorphic in $\sigma\in\C$. The indicial roots are then given by $\lambda=-i\sigma$ where $\sigma$ ranges over the poles of $N(A,\sigma)^{-1}$. Therefore, if $\nu$ lies in an indicial gap of $A$, then $N(A,-i\nu+\lambda)$ is invertible for all $\lambda\in\R$, and the inverse is uniformly bounded as a map between the semiclassical Sobolev spaces $h^s H_h^s(\Sph^2)\to h^{s+m}H_h^{s+m}(\Sph^2)$, $h=\la\lambda\ra^{-1}$. If $\alpha+\frac32$ lies in an indicial gap of $A$, it thus follows from Lemma~\ref{LemmaMUb} that for $\chi,\tilde\chi\in\CI(\ol{\R^3})$ with support in a neighborhood of $\pa\ol{\R^3}$, with $\tilde\chi=1$ on $\supp\chi$, we have
\begin{equation}
\label{EqMUbEstNearInfty}
  \|\chi u\|_{\Hb^{s,\alpha}(|\dd x|)} \leq C\Bigl(\|\tilde\chi A u\|_{\Hb^{s-m,\alpha}(|\dd x|)} + \|\tilde\chi u\|_{\Hb^{s,\alpha-1}(|\dd x|)}\Bigr).
\end{equation}
Indeed, this holds without the error term when we replace $\tilde\chi A u$ by $N(A)(\chi u)$. Commuting $N(A)$ through the cutoff and using also that $\chi(N(A)-A)\in\rho\Diffb^m$ gives~\eqref{EqMUbEstNearInfty}. These arguments extend also to the case of variable differential orders $s$ by working with the Mellin transform on the b-Sobolev space whose differential order is the dilation-invariant extension $s_{\rm I}$ of $s|_{\Sb^*\ol{\R^3}}$; since $s$ and $s_{\rm I}$ differ by an arbitrarily small amount in a sufficiently small neighborhood of $\pa\ol{\R^3}$, one obtains~\eqref{EqMUbEstNearInfty} except with the norm $\|\tilde\chi u\|_{\Hb^{s+\eps,\alpha-1}}$ on the right, for any fixed $\eps>0$.

\begin{rmk}[Indicial gaps]
\label{RmkMUbIndGap}
  Since $N(A,\sigma)^{-1}$ only has finitely many poles in any strip of bounded $\Im\sigma$, it follows that the set of $\Re\lambda$ where $\lambda\in\C$ is an indicial root of $A$ is a discrete subset of $\R$, and thus its complement is open. Thus, every $\beta\in\R$ which is not equal to $\Re\lambda$ for all indicial roots $\lambda$ lies in an indicial gap.
\end{rmk}

\subsubsection{Mellin transform at \texorpdfstring{$\iota^+$}{punctured timelike infinity}: \texorpdfstring{$0,\bop$}{0,b} and \texorpdfstring{$0\semi,\chop$-}{semiclassical cone-semiclassical 0-}settings}
\label{SssMUip}

The families~\eqref{EqMUe3bFam2}--\eqref{EqMUe3bFam3} of unit cells covering a neighborhood of $\iota^+$ are an instance of the general setup of~\S\ref{SssMSF} in the following way. Consider the manifold with boundary $\iota^+$, which is given by~\eqref{EqCPip}; we cover its interior by two families of unit cells, namely by the spatial parts
\begin{equation}
\label{EqMUipCells}
\begin{alignedat}{2}
  U_{2,j k} &= (2^{-j-2},2^{-j+2})_{\rho_\cK} \times \Sph^2_k,&\quad &\rho_{(2,j k),\rho_\cK}=1,\ \rho_{(2,j k),\omega}=1, \\
  U_{3,j k} &= (2^{-j-2},2^{-j+2})_{x_\sscri} \times \Sph^2_k,&\quad &\rho_{(3,j k),x_\sscri}=1,\ \rho_{(3,j k),\omega}=2^{-j},
\end{alignedat}
\end{equation}
$i,j\in\N_0$, $k=1,2,3$, of the unit cells~\eqref{EqMUe3bFam2}--\eqref{EqMUe3bFam3}. (The coefficient Lie algebra is thus $\CI_\bop\Vb(\iota^+)$, and the operator Lie algebra is $\CI_\bop\cV_{\bop,\eop}(\iota^+)$ where $\cV_{\bop,\eop}(\iota^+)$ consists of all b-vector fields on $\iota^+$ which are tangent to the fibers of $\iota^+\cap\scri^+$, i.e., subsets of constant angular coordinate, cf.\ \eqref{EqCMscriFibr}; hence the subscript `e' for `edge' \cite{MazzeoEdge}.) In the notation~\eqref{EqMSFCells}, let us set
\begin{equation}
\label{EqMUipScaling}
  \rho_{0,(2,j k)} = \tau_{(2,j k)} = 2^{-j},\quad
  \rho_{0,(2,j k)} = \tau_{(2,j k)} = 1.
\end{equation}
Then the unit cells and scalings~\eqref{EqMSFCells} match~\eqref{EqMUe3bFam2}--\eqref{EqMUe3bFam3} upon identifying $t$ in~\eqref{EqMSFCells} with $-\!\log_2\tau$ in~\eqref{EqMUe3bFam2}--\eqref{EqMUe3bFam3}. (Only the range of $i,j$ is unrestricted in the present discussion.) The Fourier transform in $t$ is thus the same as the Mellin transform in $\tau$, up to a constant scaling of the Mellin-dual variable; we define here
\[
  (\cM u)(\sigma,p) = \int_0^\infty \tau^{i\sigma} u(\tau,p)\,\frac{\dd\tau}{\tau},\quad \sigma\in\C,\ p\in(\iota^+)^\circ.
\]
The operator Lie algebra corresponding to this scaled bounded geometry structure on $(0,\infty)\tau\times(\iota^+)^\circ$ is $\CI_\bop\cV_\ebeop([0,\infty]_\tau\times\iota^+)$ where $\cV_\ebeop([0,\infty]_\tau\times\iota^+)$ is the space of b-vector fields which are tangent to the fibers of $\scri^+$ (i.e., of the projection $[0,\infty]\times(\iota^+\cap\scri^+)\to(\iota^+\cap\scri^+)$) and also to the fibers of the projection $[0,\infty]_\tau\times(\iota^+\cap\cK^+)\to[0,\infty]_\tau$; this is discussed in more detail in~\S\ref{SssNipFw}.

On the Mellin transform side, the construction~\eqref{EqMSFParamScbg} produces a parameterized scaled b.g.\ structure $\{(U_\alpha,\phi_\alpha,\rho_{\sigma,\alpha})\}$ on $(\iota^+)^\circ=(0,\infty)_v\times\Sph^2$, with parameter space $\R$, underlying b.g.\ structure the unit cells~\eqref{EqMUipCells}, and scalings
\[
  \rho_{\sigma,(2,j k),i} = \frac{1}{1+2^{-j}|\sigma|}\sim\frac{1}{1+\rho_\cK|\sigma|},\quad
  \rho_{\sigma,(3,j k),x_\sscri} = \frac{1}{\la\sigma\ra},\ \rho_{\sigma,(3,j k),\omega} = \frac{2^{-j}}{\la\sigma\ra}\sim\frac{x_\sscri}{\la\sigma\ra}.
\]
The coefficient Lie algebra is $\CI_\bop\Vb(\iota^+)$ for all $\sigma$. The operator Lie algebras $\hat\cV_\sigma$ are as follows.

{\bf (1)} For every fixed $\sigma\in\R$, elements of $\hat\cV_\sigma$ are precisely the b-vector fields (with conormal coefficients) near $\iota^+\cap\cK^+$ (i.e., $\rho_\cK=0$) and 0-vector fields near $\iota^+\cap\scri^+$ (i.e., $x_\sscri=0$). Regarding the latter, we recall from Mazzeo--Melrose \cite{MazzeoMelroseHyp} that the space of \emph{0-vector fields} $\cV_0(S)$ on a manifold $S$ with boundary is defined as the space
  \[
    \cV_0(S) := \{ V\in\cV(S) \colon V|_{\pa S}=0 \in \CI(\pa S;T_{\pa S}S) \}
  \]
  of vector fields vanishing at $\pa S$. Applying this with $S=\iota^+\setminus\cK^+$, a local spanning set of $\cV_0(S)$ is thus given by $x_\sscri\pa_{x_\sscri}$, $x_\sscri\pa_\omega$. Writing
  \[
    \cV_{0,\bop}(\iota^+) := \{ V\in\cV(\iota^+) \colon V|_{\iota^+\cap\scri^+}=0,\ V\ \text{is tangent to}\ \iota^+\cap\cK^+\}
  \]
  for the space of vector fields which are of 0-type at $\iota^+\cap\scri^+$ and of b-type at $\iota^+\cap\cK^+$, we thus have
  \[
    \hat\cV_\sigma = \CI_\bop\cV_{0,\bop}(\iota^+).
  \]
  The corresponding classes of $0,\bop$-ps.d.o.s are denoted
  \begin{equation}
  \label{EqMUip0b}
    \CI_\bop\Psi_{0,\bop}^{s,(2\alpha_\sscri,\gamma)}(\iota^+) = x_\sscri^{-2\alpha_\sscri}\rho_\cK^{-\gamma}\CI_\bop\Psi_{0,\bop}^s(\iota^+),
  \end{equation}
  where now $\rho_\cK,x_\sscri\in\CI(\iota^+)$ denote global defining functions of $\cK^+\cap\iota^+$ and $\cK^+\cap\scri^+$. The regularity order $s$ may be variable, while the $\scri^+$- and $\cK^+$-decay orders $2\alpha_\sscri$ and $\gamma$ are real numbers. The principal symbol captures an operator to leading order in the regularity sense. Defining a vector bundle ${}^{0,\bop}T^*\iota^+\to\iota^+$ whose space of sections is precisely $\cV_{0,\bop}(\iota^+)$, we take as the microlocalization locus the boundary
  \[
    {}^{0,\bop}S^*\iota^+
  \]
  of $\ol{{}^{0,\bop}T^*}\iota^+$ at fiber infinity. The associated Sobolev spaces are denoted
  \begin{equation}
  \label{EqMUip0bSob}
    H_{0,\bop}^{s,(2\alpha_\sscri,\gamma)}(\iota^+) = x_\sscri^{2\alpha_\sscri}\rho_\cK^\gamma H_{0,\bop}^s(\iota^+).
  \end{equation}

  For an e3b-differential operator $A_0$ as in~\eqref{EqCPipA0scri}, resp.\ \eqref{EqCPipA0scri2}, which is dilation-invariant in $\rho_+$, its ``spectral family'' (in the literature sometimes called \emph{Mellin-transformed normal operator family}) is
  \begin{equation}
  \label{EqMUipSpecFam}
  \begin{split}
    \wh{A_0}(\sigma) = &\sum_{j+k+|\alpha|\leq m} a_{j k\alpha}(x_\sscri,0,\omega) (-i\sigma)^k (x_\sscri\pa_{x_\sscri})^j (x_\sscri\pa_\omega)^\alpha, \\
    \text{resp.\ } &\sum_{j+k+|\alpha|\leq m} a_{j k\alpha}(0,\rho_\cK,\omega) (-i\sigma R)^j(R\pa_R)^k \pa_\omega^\alpha
  \end{split}
  \end{equation}
  near $\iota^+\cap\scri^+$, resp.\ $\iota^+\cap\cK^+$. For bounded $\sigma\in\C$, this is thus built from the 0-vector fields $x_\sscri\pa_{x_\sscri}$ and $x_\sscri\pa_\omega$, resp.\ the b-vector fields $\rho_\cK\pa_{\rho_\cK}$ and $\pa_\omega$.

{\bf (2)} The transition from finite to infinite $\sigma$ is given near $\cK^+$ by the \emph{semiclassical cone} setting, first introduced by the author \cite{HintzConicPowers,HintzConicProp}, and near $\scri^+$ by the \emph{semiclassical 0-}setting. The smooth setting is as follows. Define the space
  \begin{equation}
  \label{EqMUipSingle}
    \iota^+_{0\semi,\chop} := [ [0,1)_h \times \iota^+; \{0\}\times(\iota^+\cap\cK^+) ].
  \end{equation}
  Label its boundary hypersurfaces as follows:
  \begin{itemize}
  \item $\cface$, the lift of $[0,1)\times(\iota^+\cap\cK^+)$, is the \emph{cone face};
  \item $\tface$ (the front face) is the \emph{transition face};
  \item $\sface$ (the lift of $\{0\}\times\iota^+$) is the \emph{semiclassical face};
  \item $\zeroface$ (the lift of $[0,1)\times(\iota^+\cap\scri^+)$) is the \emph{0-face}.
  \end{itemize}
  Define then the space of \emph{$0\semi,\chop$-vector fields} on $\iota^+$ by
  \begin{align*}
    \cV_{0\semi,\chop}(\iota^+) := \{ V\in\cV(\iota^+_{0\semi,\chop}) &\colon V\ \text{is tangent to all boundary hypersurfaces}, \\
      &\qquad \text{vanishes at}\ \sface\cup\zeroface,\ \text{and is tangent to each $h$-level set} \}.
  \end{align*}
  A spanning set is thus given by $\frac{h}{h+\rho_\cK}V_0$ where $V_0\in\cV_{0,\bop}(\iota^+)$. Identifying $h=|\sigma|^{-1}$ (thus $\frac{1}{1+\rho_\cK|\sigma|}=\frac{h}{h+\rho_\cK}$), we then have
  \[
    (\hat\cV_\sigma)_{\sigma\in\pm[1,\infty)} = \CI_\bop\cV_{0\semi,\chop}(\iota^+).
  \]
  Corresponding classes of pseudodifferential operators are denoted
  \begin{equation}
  \label{EqMUip0hch}
    \CI_\bop\Psi_{0\semi,\chop}^{s,(2\alpha_\sscri,\alpha_\cp,l,b)}(\iota^+) = \rho_\zeroface^{-2\alpha_\sscri}\rho_\cface^{-\alpha_\cp}\rho_\tface^{-l}\rho_\sface^{-b}\CI_\bop\Psi_{0\semi,\chop}^s(\iota^+).
  \end{equation}
  If $\rho_\cK,x_\sscri\in\CI(\iota^+)$ are defining functions of $\cK^+\cap\iota^+$ and $\scri^+\cap\iota^+$, we can take $\rho_\zeroface=x_\sscri$, $\rho_\cface=\frac{\rho_\cK}{h+\rho_\cK}$, $\rho_\tface=h+\rho_\cK$, $\rho_\sface=\frac{h}{h+\rho_\cK}$. The orders $s,b$ can be variable. The principal symbol captures operators to leading order in the regularity order $s$ and the semiclassical order $b$. With ${}^{0\semi,\chop}T^*\iota^+\to\iota^+$ denoting the vector bundle for which the space of smooth sections is $\cV_{0\semi,\chop}(\iota^+)$, we thus use as the microlocalization locus
  \[
    {}^{0\semi,\chop}S^*\iota^+ \cup \ol{{}^{0\semi,\chop}T^*_\sface}\iota^+.
  \]

  For the spectral family~\eqref{EqMUipSpecFam}, this corresponds to
  \begin{subequations}
  \begin{equation}
  \label{EqMUipSpecFam0h}
    \bigl(h^m\wh{A_0}(\pm h^{-1})\bigr)_{h\in(0,1)} \in \Diff_{0\semi}^m
  \end{equation}
  being a semiclassical 0-differential operator near $\iota^+\cap\scri^+$ (i.e., it is built from $h x_\sscri\pa_{x_\sscri}$ and $h x_\sscri\pa_\omega$), resp.
  \begin{equation}
  \label{EqMUipSpecFamch}
    \Bigl(\Bigl((\frac{h}{h+R}\Bigr)^m\wh{A_0}(\pm h^{-1})\Bigr)_{h\in(0,1)} \in \Diff_\chop^m
  \end{equation}
  \end{subequations}
  being a semiclassical cone operator near $\iota^+\cap\cK^+$ (i.e., it is built from $\frac{h}{h+R}R\pa_R$ and $\frac{h}{h+R}\pa_\omega$).

  The corresponding Sobolev spaces are denoted
  \begin{equation}
  \label{EqMUip0hchSob}
    H_{0\semi,\chop}^{s,(2\alpha_\sscri,\alpha_\cface,l,b)} = \rho_\zeroface^{2\alpha_\sscri}\rho_\cface^{\alpha_\cface}\rho_\tface^l\rho_\sface^b H_{0\semi,\chop}^s(\iota^+).
  \end{equation}

\bigskip
The Plancherel type result~\eqref{EqMSFIsoVW} can be applied to $\R\times\iota^+$ equipped with the above scaled b.g.\ structure and gives an isomorphism between Sobolev spaces (without weights in $-\log_2\tau$) associated with the scaled b.g.\ structure~\eqref{EqMSFCells}, \eqref{EqMUipCells}, \eqref{EqMUipScaling} with $L^2(\R_\sigma)$ with values in 0,b- and $0\semi,\chop$-Sobolev spaces. We can allow for polynomial weights $\tau^{\alpha_+}$ by shifting the Mellin dual variable $\sigma$ by $-i\alpha_+$. Recalling now Lemma~\ref{LemmaCPip} on the relationship between $[0,\infty)_\tau\times\iota^+$ and $M$, we obtain:

\begin{lemma}[Mellin transform at $\iota^+$ on e3b-Sobolev spaces]
\label{LemmaMUetbM}
  Let $\chi\in\CI(\tilde M)$ be a function with support in the neighborhood $\cl_{\tilde M}\{r\geq 1,\ t_*\geq 1\}$ of $\iota^+$. Let $\sfs\in\CI({}^\etbop S^*\tilde M)$ be invariant under spacetime dilations $(t_*,v,\omega)\mapsto(\nu t_*,v,\omega)$, $\nu\geq 1$, on $\supp\chi$, and denote the induced orders at the (Mellin-dual) parameter value $\sigma\in\R$ by $\sfs_\sigma$. Denote by $\mu_\bop$ a positive b-density on $\iota^+$. For any $\alpha_\sscri,\alpha_+,\alpha_\cK\in\R$ and $k\in\N_0$, we have
  \begin{align*}
    &\|\chi u\|_{H_{\etbop;\bop}^{(\sfs;k),(2\alpha_\sscri,\alpha_+,\alpha_\cK)}(\tilde M,|\dd t\,\dd x|)}^2 \\
    &\quad \sim \int_{-1}^1 \| \cM(\chi u)(\sigma-i(\alpha_++2)) \|_{H_{0,\bop;\bop}^{(\sfs_\sigma;k-j),(2(\alpha_\sscri+\frac32),\alpha_\cK-\alpha_+-\frac32)}(\iota^+,\mu_\bop)}^2\,\dd\sigma \\
    &\quad\quad + \sum_\pm \int_{\pm[1,\infty)} \sum_{j=0}^k |\sigma|^j\| \cM(\chi u)(\sigma-i(\alpha_++2)) \|_{H_{(0\semi,\chop,|\sigma|^{-1});\bop}^{(\sfs_\sigma;k-j),(2(\alpha_\sscri+\frac32),\alpha_\cK-\alpha_+-\frac32,\alpha_\cK-\alpha_+-\frac32,\sfs_\sigma)}(\iota^+,\mu_\bop)}^2\,\dd\sigma.
  \end{align*}
  For $k=0$, the function spaces on the right are the Sobolev spaces $H_{\hat\cV_\sigma}$, and for $k\geq 1$ the mixed Sobolev spaces $H_{\hat\cV_\sigma;\cW}$ where $\cW=\CI_\bop\Vb(\iota^+)$, i.e., they encode $k$ degrees of (non-semiclassical) b-regularity. The orders of the spaces match the orders of the spaces of ps.d.o.s in~\eqref{EqMUip0b} and \eqref{EqMUip0hch}.
\end{lemma}

This combines \cite[Proposition~4.29]{Hintz3b} and \cite[Proposition~2.10]{HintzNonstat}. For a more explicit description of the phase space relationships and induced orders, see~\S\ref{SssNipPhase}.

\begin{proof}[Proof of Lemma~\usref{LemmaMUetbM}]
  It only remains for us to explain the weights. The isomorphism~\eqref{EqMSFIsoVW} relates $L^2$-based spaces with densities $|\frac{\dd\tau}{\tau}|\mu_\bop$ and $\mu_\bop$. Now, $\mu_\bop$ is a smooth bounded multiple of $|\frac{\dd v}{v}\,\dd\slg|$, and $|\frac{\dd\tau}{\tau}\frac{\dd v}{v}\,\dd\slg|$ is a b-density on the space $[0,\infty)\times\iota^+$ and thus on the blow-up of this space at $\{0\}\times(\iota^+\cap\cK^+)$, so it restricts to $\supp\chi$ as a b-density on $\tilde M$ by Lemma~\ref{LemmaCPip}, and hence is a smooth positive multiple of $|\frac{\dd t_*}{t_*}\,\frac{\dd x}{\la x\ra^3}|$ there. Now, $\frac{1}{t_*\la x\ra^3}$ is a smooth multiple of $x_\sscri^6\rho_+^4\rho_\cK$. For functions supported on $\supp\chi$, the norm on $H_\etbop^{0,(2\alpha_\sscri,\alpha_+,\alpha_\cK)}(\tilde M,|\dd t\,\dd x|)$ is thus equivalent to that on
  \begin{equation}
  \label{EqMUetbMShift}
    H_\etbop^{0,(2(\alpha_\sscri+\frac32),\alpha_++2,\alpha_\cK+\frac12)}\Bigl(\tilde M,\Bigl|\frac{\dd\tau}{\tau}\Bigr|\mu_\bop\Bigr)
  \end{equation}
  Factoring out $\tau^{\alpha_++2}\sim\rho_+^{\alpha_++2}\rho_\cK^{\alpha_++2}$ gives the decay orders $2(\alpha_\sscri+\frac32)$, $0$, $\alpha_\cK-\alpha_+-\frac32$.
\end{proof}

Smooth coefficient (weighted) differential operators $A\in\Diff_{0\semi,\chop}^{m,(2\alpha_\sscri,\alpha_\cp,0,b)}(\iota^+)$, with weight $0$ at $\tface\subset\iota^+_{0\semi,\chop}$, can be restricted to $\tface$. Working in local coordinates $\rho_\cK\geq 0$, $\omega\in\Sph^2$ near $\iota^+\cap\cK^+$, a local spanning set of $\cV_{0\semi,\chop}(\iota^+)$ is given by the vector fields $\frac{h}{h+\rho_\cK}\rho_\cK\pa_{\rho_\cK}$, $\frac{h}{h+\rho_\cK}\pa_\omega$. Passing to the projective coordinate $\hat\rho_\cK:=\frac{\rho_\cK}{h}$ on $\tface$, these become $\frac{1}{1+\hat\rho_\cK}\hat\rho_\cK\pa_{\hat\rho_\cK}$, $\frac{1}{1+\hat\rho_\cK}\pa_\omega$, which are thus b,sc-vector fields on
\begin{equation}
\label{EqMUNtfSpace}
  \tface = [0,\infty]_{\hat\rho_\cK} \times \Sph^2,
\end{equation}
with b-behavior at $\hat\rho_\cK=0$ (i.e., the intersection of $\tface$ with $\cface$) and scattering behavior at $\hat\rho_\cK=\infty$ (i.e., the intersection of $\tface$ with $\sface$). Restriction of coefficients thus defines the \emph{$\tface$-normal operator}
\begin{equation}
\label{EqMUNtfch}
  N_\tface(A) \in \Diff_{\bop,\scop}^{m,(\alpha_\cp,b)}(\tface),\quad A\in\Diff_{0\semi,\chop}^{m,(2\alpha_\sscri,\alpha_\cp,0,b)}(\iota^+).
\end{equation}
In view of this relationship of operator Lie algebras, semiclassical cone Sobolev spaces are related to b-scattering Sobolev spaces on $\tface$. Concretely, if $\chi_\tface\in\CI(\iota^+_{0\semi,\chop})$ (or $\CI_\bop(\iota^+_{0\semi\chop})$, or indeed any function on $(\iota^+_{0\semi,\chop})^\circ$) is supported in a small neighborhood of $\tface$ (e.g., $\chi_\tface=\chi(\rho_\cK+h)$ where $\chi\in\CIc([0,1))$ equals $1$ near $0$), then we have a uniform equivalences
\begin{equation}
\label{EqMUipNormEquiv}
  \|\chi_\tface u\|_{H_{0\semi,\chop,h}^{s,(2\alpha_\sscri,\alpha_\cface,l,b)}(\iota^+)} \sim h^{-l}\|\chi_\tface u\|_{H_{\bop,\scop}^{s,(\alpha_\cface,b-l)}(\tface)}.
\end{equation}
This also holds for variable orders $s,b$, which are then $h$-dependent on the right-hand side.

\subsubsection{Spaces of supported and extendible distributions}
\label{SssMUSupp}

On Kerr, we will work with the spatial manifold $X=\ol{\{r\geq\bhm\}}\subsetneq\tilde X=\ol{\R^3}$. The corresponding Sobolev spaces consist of distributions that are \emph{extendible} across the ``artificial'' boundary at $r=\bhm$ (or of \emph{supported} character there, i.e., vanishing for $r<\bhm$, when studying dual problems). We develop here the relevant notions and (uniform) duality results in the framework of \cite{HintzScaledBddGeo}; see \cite[\S{9.3}]{HintzMicro} for background.

Thus, consider a parameterized scaled b.g.\ structure on the manifold $\cM$ without boundary; denote the parameter (which may be a singleton set) by $\sigma$, and the operator Lie algebra by $\cV=(\cV_\sigma)_\sigma$. We shall only consider unweighted Sobolev spaces $H_{\cV_\sigma}^s(\cM)$ here and leave the notational modifications for weighted spaces to the reader.

\begin{lemma}[Uniform duality]
\label{LemmaMUSuppDual}
  Let $s\in\R$, or let $s$ be a variable order function. Fix a uniformly positive $\cV$-density on $\cM$ to define the $L^2$-spaces $L^2_{\cV_\sigma}(\cM)$.
  \begin{enumerate}
  \item The $L^2_{\cV_\sigma}$-pairing $\la u,v\ra$, $u,v\in\CIc(\cM)$, extends by continuity to $u\in H^s_{\cV_\sigma}(\cM)$, $v\in H^{-s}_{\cV_\sigma}(\cM)$. There exists a constant $C$ such that, for all $u,v,\sigma$,
    \begin{equation}
    \label{EqMUSuppDualCS}
      |\la u,v\ra| \leq C\|u\|_{H^s_{\cV_\sigma}(\cM)}\|v\|_{H^{-s}_{\cV_\sigma}(\cM)}.
    \end{equation}
  \item There exist $c,C>0$ such that, for all $u,\sigma$,
    \begin{equation}
    \label{EqMUSuppDualDual}
      c\|u\|_{H_{\cV_\sigma}^s(\cM)} \leq \sup_{\|v\|_{H_{\cV_\sigma}^{-s}(\cM)}=1} |\la u,v\ra| \leq C\|u\|_{H_{\cV_\sigma}^s(\cM)}.
    \end{equation}
  \end{enumerate}
\end{lemma}
\begin{proof}
  For $s=0$, these are the Cauchy--Schwarz inequality and Riesz' theorem, respectively. For general $s$, fix an elliptic ps.d.o.\ $A=(A_\sigma)_\sigma\in\Psi_\cV^s$ and a parametrix $B=(B_\sigma)_\sigma\in\Psi_\cV^{-s}$ such that $I=B A+R$, $R=(R_\sigma)_\sigma\in\Psi_\cV^{-\infty}$. For $u,v\in\CIc$, we can then write
  \begin{align*}
    \la u,v\ra &= \la A_\sigma u,B_\sigma^*v\ra + \la R_\sigma u,v\ra \\
      &= \la A_\sigma u,B_\sigma^*v\ra + \la R_\sigma u,A_\sigma^*B_\sigma^*v\ra + \la R_\sigma u,R_\sigma^*v\ra \\
      &= \la A_\sigma(I+R_\sigma)u,B_\sigma^*v\ra + \la R_\sigma u,R_\sigma^*v\ra.
  \end{align*}
  In view of the \emph{uniform} boundedness of $A_\sigma\colon H_{\cV_\sigma}^s\to L^2_{\cV_\sigma}$, similarly for $B_\sigma$ and $R_\sigma$, we obtain~\eqref{EqMUSuppDualCS}.

  The upper bound in~\eqref{EqMUSuppDualDual} follows from~\eqref{EqMUSuppDualCS}. To prove the lower bound, we need to show that given $\lambda\in(H_{\cV_\sigma}^{-s})^*$, there exists a (necessarily unique) $u\in H_{\cV_\sigma}^s$ such that $\lambda(v)=\la v,u\ra$, and that the (antilinear) map $\lambda\mapsto u$ is continuous. Consider then the map
  \[
    j \colon L^2_{\cV_\sigma} \oplus L^2_{\cV_\sigma} \ni (v_0,v_1) \mapsto A v_0+R v_1 \in H_{\cV_\sigma}^{-s}.
  \]
  It is surjective, with $v\in H^{-s}_{\cV_\sigma}$ being expressible as $j(B v,v)$. The composition $\lambda\circ j\in (L^2_{\cV_\sigma}\oplus L^2_{\cV_\sigma})^*$ is then given, for suitable $u_0,u_1\in L^2_{\cV_\sigma}$, by $\lambda(j(v_0,v_1))=\la v_0,u_0\ra+\la v_1,u_1\ra$, so
  \[
    \lambda(v) = \lambda(j(B v,v)) = \la B v,u_0\ra + \la R v,u_1\ra = \la v,B^*u_0+R^*u_1\ra = \la v,u\ra,
  \]
  where $u=B^*u_0+R^*u_1\in H_{\cV_\sigma}^s$. Moreover, $\|u_0\|_{L^2_{\cV_\sigma}}+\|u_1\|_{L^2_{\cV_\sigma}}\leq C\|\lambda\|_{(H_{\cV_\sigma}^{-s})^*}$, and thus also $\|u\|_{H_{\cV_\sigma}^s}\leq C'\|\lambda\|_{(H_{\cV_\sigma}^{-s})^*}$.
\end{proof}

Let now $\Omega\subset\cM$ be any open subset.

\begin{definition}[Supported and extendible spaces]
\label{DefMUSupp}
  We define $\dot H_{\cV_\sigma}^s(\bar\Omega):=\{u\in H_{\cV_\sigma}^s(\cM)\colon\supp u\subset\bar\Omega\}$ with the induced norm, and $\bar H_{\cV_\sigma}^s(\Omega):=\{u|_\Omega\colon u\in H_{\cV_\sigma}^s(\cM)\}$ with the quotient norm induced by the isomorphism
  \[
    H_{\cV_\sigma}^s(\cM) / \dot H_{\cV_\sigma}^s(\cM\setminus\Omega) \cong \bar H_{\cV_\sigma}^s(\Omega) 
  \]
  given by restriction to $\Omega$.
\end{definition}

Thus $\|u\|_{\bar H_{\cV_\sigma}^s(\Omega)}=\|\tilde u\|_{H_{\cV_\sigma}^s(\cM)}$ where $\tilde u$ is the unique extension of $u$ with minimal norm. We use the analogous dot/bar notation also for other function spaces. Then:

\begin{lemma}[Uniform duality: extendible/supported case]
\label{LemmaMUExt}
  We use the notation of Lemma~\usref{LemmaMUSuppDual}.
  \begin{enumerate}
  \item\label{ItMUExt1} The $L^2_{\cV_\sigma}$-pairing $\la u,v\ra$, $u\in\bar\cC^\infty(\Omega)$, $v\in\dot\cC^\infty(\bar\Omega)$, extends by continuity to $u\in\bar H_{\cV_\sigma}^s(\Omega)$, $v\in\dot H_{\cV_\sigma}^{-s}(\bar\Omega)$. There exists a constant $C$ such that, for all $u,v,\sigma$,
    \[
      |\la u,v\ra| \leq C\|u\|_{\bar H_{\cV_\sigma}^s(\Omega)}\|v\|_{\dot H_{\cV_\sigma}^{-s}(\bar\Omega)}.
    \]
  \item\label{ItMUExt2} There exist $c,C>0$ such that, for all $u,\sigma$,
    \[
      c\|u\|_{\bar H_{\cV_\sigma}^s(\Omega)} \leq \sup_{\|v\|_{\dot H_{\cV_\sigma}^{-s}(\bar\Omega)}=1} |\la u,v\ra| \leq C\|u\|_{\bar H_{\cV_\sigma}^s(\Omega)};
    \]
    similarly with the roles of $\bar H_{\cV_\sigma}$ and $\dot H_{\cV_\sigma}$ reversed.
  \end{enumerate}
\end{lemma}
\begin{proof}
  Given $u\in\bar H_{\cV_\sigma}^s(\Omega)$, let $\tilde u\in H_{\cV_\sigma}^s(\cM)$ be its minimal norm extension. Then~\eqref{EqMUSuppDualCS} gives
  \[
    |\la u,v\ra| = |\la\tilde u,v\ra| \leq C\|\tilde u\|_{H_{\cV_\sigma}^s(\cM)}\|v\|_{H_{\cV_\sigma}^{-s}(\cM)} = C\|u\|_{\bar H_{\cV_\sigma}^s(\Omega)}\|v\|_{\dot H_{\cV_\sigma}^{-s}(\bar\Omega)}.
  \]
  For part~\eqref{ItMUExt2}, consider $\lambda\in\dot H_{\cV_\sigma}^{-s}(\bar\Omega)^*$. Extend this to $\tilde\lambda\in H_{\cV_\sigma}^{-s}(\cM)^*$ with the same norm (by setting $\tilde\lambda=\lambda$ on $\dot H_{\cV_\sigma}^{-s}(\bar\Omega)$ and $\tilde\lambda=0$ on $\dot H_{\cV_\sigma}^{-s}(\bar\Omega)^\perp\subset H_{\cV_\sigma}^{-s}(\cM)$). By (the proof of) Lemma~\ref{LemmaMUSuppDual}, there exists $\tilde u\in H_{\cV_\sigma}(\cM)$, with norm uniformly controlled by that of $\tilde\lambda$ and thus $\lambda$, such that $\tilde\lambda(v)=\la v,\tilde u\ra$. For $v\in\dot H_{\cV_\sigma}^{-s}(\bar\Omega)$, we thus have
  \[
    \lambda(v) = \la v,\tilde u\ra = \la v,u\ra,\quad u:=\tilde u|_\Omega \in \bar H_{\cV_\sigma}^s(\Omega),
  \]
  and $\|u\|_{\bar H_{\cV_\sigma}^s(\Omega)}\leq\|\tilde u\|_{H_{\cV_\sigma}^s(\cM)}\leq C\|\lambda\|_{(\dot H_{\cV_\sigma}^{-s}(\bar\Omega))^*}$.
\end{proof}

We can now extend the Plancherel isomorphism~\eqref{EqMSFIsoV} to spaces of supported and extendible distributions:

\begin{prop}[Plancherel theorem]
\label{PropMUSuppFT}
  We use the notation of~\eqref{EqMSFCells}--\eqref{EqMSFIsoV}. Let $\Omega\subset\cX$ be open. Then the Fourier transform in the $\R$-factor of $\cM=\R\times\cX$ induces isomorphisms
  \begin{align}
  \label{EqMUSuppFTDot}
    \cF\colon \dot H_\cV^s(\R\times\bar\Omega) &\xra{\cong} L^2\bigl(\R_\sigma;\la\tau\sigma\ra^{-s}\dot H_{\hat\cV_\sigma}^s(\bar\Omega)\bigr), \\
  \label{EqMUSuppFTSupp}
    \cF\colon \bar H_\cV^s(\R\times\Omega) &\xra{\cong} L^2\bigl(\R_\sigma;\la\tau\sigma\ra^{-s}\bar H_{\hat\cV_\sigma}^s(\Omega) \bigr).
  \end{align}
\end{prop}
\begin{proof}
  The isomorphism~\eqref{EqMUSuppFTDot} follows directly from~\eqref{EqMSFIsoV}. To prove~\eqref{EqMUSuppFTSupp}, we write
  \begin{align*}
    \|u\|_{\bar H_\cV^s(\R\times\Omega)}^2 &\leq C\sup_{\|v\|_{\dot H_\cV^{-s}(\R\times\bar\Omega)}=1} |\la u,v\ra|^2 = C\sup_{\|v\|_{\dot H_\cV^{-s}(\R\times\bar\Omega)}=1} \biggl|\int_\R \la\hat u(\sigma),\hat v(\sigma)\ra_{L^2(\cX)}\,\dd\sigma\biggr|^2 \\
      &\leq C C'\sup_{\|\hat v\|_{L^2(\R_\sigma;\la\tau\sigma\ra^s\dot H_{\cV_\sigma}^{-s}(\bar\Omega))}=1}\,\int_\R |\la\hat u(\sigma),\hat v(\sigma)\ra_{L^2(\cX)}|\,\dd\sigma\biggr|^2 \\
      &\leq C C'\sup_{\|\hat v\|_{L^2(\R_\sigma;\la\tau\sigma\ra^s\dot H_{\cV_\sigma}^{-s}(\bar\Omega))}=1} \biggl|\int_\R \|\hat u(\sigma)\|_{\la\tau\sigma\ra^{-s}\bar H_{\cV_\sigma}^s(\Omega)} \|\hat v(\sigma)\|_{\la\tau\sigma\ra^s\dot H_{\cV_\sigma}^{-s}(\bar\Omega)}\,\dd\sigma\biggr|^2 \\
      &\leq C C'C'' \|\hat u\|_{L^2(\R_\sigma;\la\tau\sigma\ra^{-s}\bar H_{\cV_\sigma}^s(\Omega))}^2.
  \end{align*}
  Here, we use Lemma~\ref{LemmaMUExt}\eqref{ItMUExt2} (without parameter dependence) for the first, the isomorphism~\eqref{EqMUSuppFTDot} for the second, Lemma~\ref{LemmaMUExt}\eqref{ItMUExt1} for the third, and Cauchy--Schwarz for the final inequality.
\end{proof}

We thus obtain an analogue of Lemma~\ref{LemmaMUetbFT} for distributions of supported/extendible character, and hence also of Lemma~\ref{LemmaMUetbFTb}.

Finally, for the domain $\Omega=\cl_M\{t-r\geq 1\}\subset M$ of Definition~\ref{DefCMDomain}, we shall write
\begin{subequations}
\begin{equation}
\label{EqMUSuppExt}
  H_\etbop^s(\Omega)^{\bullet,-}
\end{equation}
for the space of restrictions to $\{r>\bhm\}$ of elements of $H_\etbop^s(\tilde M)$ which are supported in $t-r\geq 1$. (This will be the spacetime domain on which we solve forcing problems for wave equations.) The analogous spaces related to $H^{(s;k)}_{\etbop;\bop}(\tilde M)$ are denoted
\begin{equation}
\label{EqMUSuppExt2}
  H_{\etbop;\bop}^{(s;k)}(\Omega)^{\bullet,-},
\end{equation}
\end{subequations}
similarly for weighted spaces.

\subsubsection{\texorpdfstring{$L^\infty$-}{L-infinity-}based spaces and Sobolev embedding}
\label{SssMUC}

While $L^2$-based spaces are convenient for (non-elliptic) microlocal analysis, the \emph{coefficients} of (pseudo)differential operators (especially for nonlinear applications) are best captured in $L^\infty$-based spaces. The $\cV$-quantizations of finite regularity symbols as defined in Definition~\ref{DefMSSymbolMixed} provide a very general class of operators for which operator norm bounds in terms of symbol seminorms are available (Lemma~\ref{LemmaMSOpNorm}). In our applications below, we can, for the most part, work with more restrictive classes of operators which are built out of standard scaled b.g.\ ps.d.o.s and post-multiplication with limited regularity functions. We first define:

\begin{definition}[e3b-continuous functions]
\label{DefMUCe3b}
  Fix boundary defining functions $x_\sscri$, $\rho_+$, and $\rho_\cK$ of $\scri^+$, $\iota^+$, and $\cK^+\subset\tilde M$, respectively. For $d_0\in\N_0$ and $\alpha_\sscri,\alpha_+,\alpha_\cK\in\R$, we then write
  \[
    \cC_\etbop^{d_0,(2\alpha_\sscri,\alpha_+,\alpha_\cK)}(\tilde M)
  \]
  for the space of all $\cC^{d_0}$-functions $u$ on $\tilde M^\circ=\R^4$ such that
  \begin{equation}
  \label{EqMUCe3bC}
    x_\sscri^{-2\alpha_\sscri}\rho_+^{-\alpha_+}\rho_\cK^{-\alpha_\cK}A u \in L^\infty\quad\forall\,A\in\Diff_\etbop^{d_0}(\tilde M).
  \end{equation}
  More generally, for $k\in\N_0$, we write
  \[
    \cC_{\etbop;\bop}^{(d_0;k),(2\alpha_\sscri,\alpha_+,\alpha_\cK)}(\tilde M)
  \]
  when~\eqref{EqMUCe3bC} holds for all $A$ which are linear combinations of products $A_1 A_2$ where $A_1\in\Diffb^k$, $A_2\in\Diff_\etbop^{d_0}$. For $d_0=0$, this space is denoted $\cC_\bop^{k,(2\alpha_\sscri,\alpha_+,\alpha_\cK)}(\tilde M)$ simply.
\end{definition}

Taking the sum of the supremum norms of~\eqref{EqMUCe3bC} where $A$ ranges over a fixed finite spanning set of $\Diff_\etbop^{d_0}(\tilde M)$ over $\CI(\tilde M)$ gives $\cC_\etbop^{d_0,(2\alpha_\sscri,\alpha_+,\alpha_\cK)}(\tilde M)$ the structure of a Banach space; similarly for $\cC_{\etbop;\bop}^{(d_0;k),(2\alpha_\sscri,\alpha_+,\alpha_\cK)}$.

\begin{definition}[e3b-operators with coefficients of limited regularity]
\label{DefMUCe3bOp}
  For $m,s,d_0,k\in\N_0$, write
  \begin{equation}
  \label{EqMUCe3bOp}
    \cC_{\etbop;\bop}^{(d_0;k)}\Diff_\etbop^m(\tilde M),\ \ \text{resp.}\ \ \cC_{\etbop;\bop}^{(d_0;k)}\Psi_\etbop^s(\tilde M)
  \end{equation}
  for the space of finite linear combinations of operators of the form $a P$ where $a\in\cC_{\etbop;\bop}^{(d_0;k)}(\tilde M)$ and $P\in\Diff_\etbop^m(\tilde M)$, resp.\ $P\in\CI_\bop\Psi_\etbop^s(\tilde M)$; similarly for spaces of weighted operators.
\end{definition}

Fixing a finite spanning set $\{P_1,\ldots,P_N\}$ of $\Diff_\etbop^m(\tilde M)$, we can define a norm on the first space in~\eqref{EqMUCe3bOp} by
\begin{equation}
\label{EqMUCe3bOpNorm}
  \|A\|_{\cC_{\etbop;\bop}^{(d_0;k)}\Diff_\etbop^m} := \inf\,\Biggl\{ \max_{1\leq j\leq N} \|a_j\|_{\cC_{\etbop;\bop}^{(d_0;k)}} \colon A = \sum_{j=1}^N a_j P_j \Biggr\}.
\end{equation}
In the pseudodifferential case, there does not exist a finite spanning set anymore. In our applications, we will, however, only encounter the following situation: for some (typically explicitly given) finite set $\sP=\{P_1,\ldots,P_N\}\subset\CI_\bop\Psi_\etbop^s(\tilde M)$, we consider only operators in the $\cC_{\etbop;\bop}^{(d_0;k)}$-span $\cC_{\etbop;\bop}^{(d_0;k)}\sP$ of $P_1,\ldots,P_N$; we can then define
\begin{equation}
\label{EqMUCe3bOpNormPsdo}
  \|P\|_{\cC_{\etbop;\bop}^{(d_0;k)}\Psi_\etbop^s} := \inf\,\Biggl\{ \max_{1\leq j\leq N} \|a_j\|_{\cC_{\etbop;\bop}^{(d_0;k)}} \colon A = \sum_{j=1}^N a_j P_j \Biggr\},\quad P\in\cC_{\etbop;\bop}^{(d_0;k)}\sP.
\end{equation}
For the sake of readability, we do not make the set $\sP$ explicit in the notation for this norm.

We shall prove microlocal estimates for solutions of wave-type equations on mixed $(\etbop;\bop)$-spaces; in nonlinear applications, these solutions also play the simultaneous role of coefficients of wave operators. We thus need Sobolev embedding results. We shall only recall them for b-Sobolev spaces. First, we record:

\begin{lemma}[Mixed and pure spaces]
\label{LemmaMUCe3bTob}
  Let $\sfs$ be a variable e3b-regularity order function and $k\in\N_0$. Let $s_-=\max(0,\lceil-\inf\sfs\rceil)$, and suppose that $k\geq s_-$. Then
  \begin{equation}
  \label{EqMUCe3bTob}
    H_{\etbop;\bop}^{(\sfs;k)}(\tilde M) \subset H_\bop^{k-s_-}(\tilde M).
  \end{equation}
\end{lemma}
\begin{proof}
  If $\sfs\geq 0$, then $s_-=0$ and the inclusion~\eqref{EqMUCe3bTob} is clear. Otherwise, pick finitely many e3b-differential operators $A_1,\ldots,A_N\in\Diff_\etbop^{s_-}(\tilde M)$ such that the intersection of their characteristic sets is empty. We can then write the identity map as
  \[
    I = \sum_{j=1}^N A_j B_j + R,\quad B_j\in\CI_\bop\Psi_\etbop^{-s_-},\ R\in\CI_\bop\Psi_\etbop^{-\infty}.
  \]
  Therefore, any $u\in H_{\etbop;\bop}^{(\sfs;k)}$ can be written as
  \[
    u = A_j u_j + u',\quad u_j := B_j u\in H_{\etbop;\bop}^{(\sfs+s_-;k)}\subset\Hb^k,\ u' := R u\in H_{\etbop;\bop}^{(\infty;k)}\subset\Hb^k.
  \]
  This decomposition gives~\eqref{EqMUCe3bTob}.
\end{proof}

The second step is the following b-Sobolev embedding result:

\begin{lemma}[b-Sobolev embedding]
\label{LemmaMUCe3bSob}
  Let $d\in\N_0$, and let $\alpha_\sscri,\alpha_+,\alpha_\cK\in\R$. Then there exists a constant $C$ such that, for all $u$ with support in $\Omega$ as in Definition~\usref{DefCMDomain},
  \[
    \|u\|_{\cC_\bop^{d,(2(\alpha_\sscri+\frac32),\alpha_++2,\alpha_\cK+\frac12)}} \leq C\|u\|_{\Hb^{d+3,(2\alpha_\sscri,\alpha_+,\alpha_\cK)}(\tilde M;|\dd t\,\dd x|)}.
  \]
  If, on the other hand, $\mu_\bop$ denotes an unweighted b-density on $\tilde M$, then
  \[
    \|u\|_{\cC_\bop^{d,(2\gamma_\sscri,\gamma_+,\gamma_\cK)}} \leq C\|u\|_{\Hb^{d+3,(2\gamma_\sscri,\gamma_+,\gamma_\cK)}(\tilde M,\mu_\bop)}
  \]
  for all $\gamma_\sscri,\gamma_+,\gamma_\cK\in\R$. Conversely, for all $\delta>0$,
  \begin{equation}
  \label{EqMUCe3bSobConv}
    \|u\|_{\Hb^{d,(2\gamma_\sscri,\gamma_+,\gamma_\cK)}} \leq C_\delta\|u\|_{\cC_\bop^{d,(2\gamma_\sscri+\delta,\gamma_++\delta,\gamma_\cK+\delta)}}.
  \end{equation}
\end{lemma}
\begin{proof}
  The shift of weights upon passing from the Minkowskian volume density to a b-density on $\tilde M$ was already discussed around~\eqref{EqMUetbMShift}. It then remains to work in a local coordinate patch $[0,1)_x^k\times\R^{4-k}_y$ near a codimension $k$ corner of $\tilde M$ and the b-density $|\frac{\dd x^1\cdots\dd x^k}{x^1\cdots x^k}\dd y|$ and show that functions in $\Hb^{d+3}$ with compact support lie in $\cC_\bop^d$. Upon passing to logarithmic coordinates $z^j:=-\!\log x^j$, the first claim is then the standard Sobolev embedding $H^{d+3}(\R^4)\hra\cC^d(\R^4)$. In logarithmic coordinates, the final claim amounts to the statement that a function on $\R^4$ that is supported in $[1,\infty)_x^k\times\R^{4-k}$ and has uniformly bounded derivative of order up to $d$ has square integrable derivatives of orders up to $d$ upon division by $e^{\delta z_1}\cdots e^{\delta z_k}$ for any $\delta>0$.
\end{proof}

We also note, as a variant of Lemma~\ref{LemmaMUCe3bSob} (with an analogous proof), that
\begin{equation}
\label{EqMUCbSob}
  \Hb^{d+3,\alpha}(\ol{\R^3},|\dd x|) \hra \cC_\bop^{d,\alpha+\frac32}(\ol{\R^3}).
\end{equation}

\subsection{Basic e3b-microlocal estimates}
\label{SsMBasic}

We only prove microlocal elliptic and real principal type propagation estimates on e3b-spaces here, using the ps.d.o.\ classes introduced in~\eqref{EqMUe3bPsdo} and \eqref{EqMUe3bTilde}; analogous estimates hold also for all the other ps.d.o.\ algebras described in~\S\ref{SsMU}. We shall moreover state estimates only on unweighted spaces and leave the purely notational modifications for weighted estimates to the reader. The proofs are, in essence, adaptations of those in \cite[\S{2.5.4}]{HintzGlueLocII} except for minor differences in the precise assumptions on the operators. For the convenience of the reader, we give slightly expanded proofs here.

\begin{prop}[e3b-microlocal elliptic regularity]
\label{PropMEll}
  Let $m,N\in\R$, $\sfs\in\CI({}^\etbop S^*\tilde M)$. Then there exists $d_0=d_0(m,\sfs,N)\in\N$ such that the following holds for all $k\in\N_0$ and $L\in\cC_{\etbop;\bop}^{(d_0;k)}\Psi_\etbop^m(\tilde M)$. Let $\chi\in\CI(\tilde M)$, and let $B,G\in\CI_\bop\Psi_\etbop^0(\tilde M)$ be such that $B=\chi B\chi$, $G=\chi G\chi$, and $\WF_\etbop'(B)\subset\Ell_\etbop(G)\cap\Ell_\etbop(L)$. Then there exists a constant $C$ such that
  \begin{equation}
  \label{EqMEll}
    \|B u\|_{H_{\etbop;\bop}^{(\sfs;k)}} \leq C\Bigl( \|G L u\|_{H_{\etbop;\bop}^{(\sfs-m;k)}} + \|\chi u\|_{H_{\etbop;\bop}^{(-N;k)}}\Bigr).
  \end{equation}
  This estimate holds locally uniformly for perturbations of $L$ in the class $\cC_{\etbop;\bop}^{(d_0;k)}\sP$ where $\sP\subset\CI_\bop\Psi_\etbop^m(\tilde M)$ is a fixed finite set.
\end{prop}
\begin{proof}
  For $k=0$ and $d_0=\infty$, the estimate~\eqref{EqMEll} follows from a standard elliptic parametrix construction in the algebra $\tilde\Psi_\etbop(\tilde M)$. Since the boundedness on e3b-ps.d.o.s on fixed scales of e3b-Sobolev spaces only requires finite e3b-regularity of the coefficients (see Lemma~\ref{LemmaMSOpNorm}), a sufficiently high but finite e3b-regularity $d_0$ of the coefficients suffices.

  In order to prove~\eqref{EqMEll} for $k\geq 1$, we argue inductively. Let $V\in\cV_{\bop,[\etbop]}(\tilde M)$ be an e3b-commutator b-vector field (see Lemma~\ref{LemmaCTe3bComm}). Apply~\eqref{EqMEll} with $k-1$, $V u$ in place of $k$, $u$. Consider
  \[
    L(V u) = V(L u) - [V,L]u
  \]
  and expand the second term using $L=\sum_{j=1}^N a_j P_j$, $a_j\in\cC_{\etbop;\bop}^{(d_0;k)}$, $P_j\in\sP$. The first term in $[V,a_j P_j]=(V a_j)P_j+a_j[V,P_j]$ is of class $\cC_{\etbop;\bop}^{(d_0;k-1)}\sP$, and the second term is of class $\cC_{\etbop;\bop}^{(d_0;k)}\Psi_\etbop^m$. Therefore, using also that $[G,V]\in\CI_\bop\Psi_\etbop^0$, we get
  \begin{align*}
    \|G L V u\|_{H_{\etbop;\bop}^{(\sfs-m;k-1)}} &\leq \|V G L u\|_{H_{\etbop;\bop}^{(\sfs-m;k-1)}} + \|[G,V]L u\|_{H_{\etbop;\bop}^{(\sfs-m;k-1)}} \\
      &\qquad + \sum_{j=1}^N \|G(V a_j)P_j u\|_{H_{\etbop;\bop}^{(\sfs-m;k-1)}} + \sum_{j=1}^N \|G a_j[V,P_j]u\|_{H_{\etbop;\bop}^{(\sfs-m;k-1)}} \\
      &\leq \|G L u\|_{H_{\etbop;\bop}^{(\sfs-m;k)}} + C\Bigl(\|G^\sharp u\|_{H_{\etbop;\bop}^{(\sfs;k-1)}} + \|\chi u\|_{H_{\etbop;\bop}^{(-N;k-1)}}\Bigr)
  \end{align*}
  for any fixed $G^\sharp\in\CI_\bop\Psi_\etbop^0$ with $\Ell_\etbop(G^\sharp)\supset\WF'_\etbop(G)$. We use here microlocal elliptic regularity on $\WF'_\etbop(G)$ for the already settled b-regularity order $k-1$. We use, moreover, that, by a simple application of the Leibniz rule, multiplication by $V a_j\in\cC_{\etbop;\bop}^{(d_0;k-1)}$ preserves $H_{\etbop;\bop}^{(\sfs-m;k-1)}$ when $d_0$ is sufficiently large (independently of $k$); and we recall that elements of $\CI_\bop\Psi_\etbop$ map boundedly between mixed $(\etbop;\bop)$-Sobolev spaces.

  We can similarly bound $\|B V u\|_{H_{\etbop;\bop}^{(\sfs;k-1)}}\geq\|V B u\|_{H_{\etbop;\bop}^{(\sfs;k-1)}}-C\|\tilde B u\|_{H_{\etbop;\bop}^{(\sfs;k-1)}}$ where the elliptic set of $\tilde B\in\CI_\bop\Psi_\etbop^0$ contains the operator wave front set of $B$. Arguing similarly for the term $\chi V u$ in~\eqref{EqMEll}, we thus obtain~\eqref{EqMEll} as stated (with slightly enlarged operator wave front sets and supports of cutoffs) since the set of e3b-commutator b-vector fields spans $\Vb(\tilde M)$ over $\CI(\tilde M)$.
\end{proof}

While we do not need concrete possible values of $d_0$ in this paper, we remark that reasonable bounds can, in principle, be obtained; see \cite{BealsReedMicroNonsmooth,HintzQuasilinearDS} for simple results of that nature for ps.d.o.s with Sobolev-regular coefficients.

Next, we consider real principal type propagation estimates. We fix a positive homogeneous elliptic symbol
\begin{equation}
\label{EqMPrShift}
  \rho_\infty\in S_{\rm hom}^{-1}(\Tetb^*\tilde M\setminus o)
\end{equation}
to shift orders.

\begin{prop}[e3b-microlocal real principal propagation]
\label{PropMPr}
  Let $m\in\N_0$, $N\in\R$, and $\sfs\in\CI({}^\etbop S^*\tilde M)$. Then there exists $d_0=d_0(m,\sfs,N)\in\N$ such that the following holds for all $\delta>0$, $k\in\N_0$ and
  \begin{equation}
  \label{EqMPrOp}
    L = L_0 + L_1 + L_2,\quad L_0\in\Diff_\etbop^m(\tilde M),\ L_1\in\cC_{\etbop;\bop}^{(d_0;k),(2\delta,\delta,\delta)}\Psi_\etbop^m(\tilde M),\ L_2\in\cC_{\etbop;\bop}^{(d_0;k)}\Psi_\etbop^{m-1}(\tilde M),
  \end{equation}
  where $L_1$ has a homogeneous principal symbol. Assume that the principal symbol $\ell$ of $L$ is real-valued. Write $\sfH_\ell:=\rho_\infty^{m-1}H_\ell$ for its rescaled Hamiltonian vector field.\footnote{Since $H_\ell$ is homogeneous of degree $m-1$ with respect to fiber dilations, $\sfH_\ell$ is homogeneous of degree $0$ and thus induces a vector field on ${}^\etbop S^*\tilde M$.} Let $\chi\in\CI(\tilde M)$, and let $B,E,G\in\CI_\bop\Psi_\etbop^0(\tilde M)$ be such that $B=\chi B\chi$, $E=\chi E\chi$, $G=\chi G\chi$. Suppose that all backwards integral curves of $\sfH_\ell$ starting at $\WF'_\etbop(B)$ reach $\Ell_\etbop(E)$ in uniformly finite time while remaining in $\Ell_\etbop(G)$, and that $\sfH_\ell\sfs\leq 0$ in a neighborhood of $\WF_\etbop'(B)\cap\Char_\etbop(L)$ within $\Char_\etbop(L)\subset{}^\etbop S^*\tilde M$. Suppose also that $\WF_\etbop'(B)\subset\Ell_\etbop(G)$. Then there exists a constant $C$ such that
  \begin{equation}
  \label{EqMPrEst}
    \|B u\|_{H_{\etbop;\bop}^{(\sfs;k)}} \leq C\Bigl(\|G L u\|_{H_{\etbop;\bop}^{(\sfs-m+1;k)}} + \|E u\|_{H_{\etbop;\bop}^{(\sfs;k)}} + \|\chi u\|_{H_{\etbop;\bop}^{(-N;k)}}\Bigr).
  \end{equation}
  This estimate holds locally uniformly also for $\Diff_\etbop^m$, $\cC_{\etbop;\bop}^{(d_0;k)}\Psi_\etbop^{(2\delta,\delta,\delta)}$, and $\cC_{\etbop;\bop}^{(d_0;k)}\Psi_\etbop^{m-1}$ perturbations $L'_0,L'_1,L'_2$ of $L_0,L_1,L_2$ for any fixed class of operators $\cC_{\etbop;\bop}^{(d_0;k)}\sP^{m-j}$, $\sP^{m-j}\subset\Psi_\etbop^{m-j}$ finite, $j=0,1$, that have real principal symbols, and with $\sfH_{\ell'}\sfs\leq 0$ on their characteristic sets over a neighborhood of $\WF'_\etbop(B)$ where $\ell'$ is the principal symbol of $L'=L'_0+L'_1+L'_2$.
\end{prop}

While the conditions on $L$ can be generalized significantly,\footnote{For example, the differential nature of $L_0$ is irrelevant; one can allow $L_0$ to lie in $\CI_\bop\Psi_\etbop^m$, $m\in\R$, and require it to have a smooth principal symbol on $\Tetb^*\tilde M$.} they suffice for our applications. The vanishing of the coefficients of $L_1$ at $\pa\tilde M$ ensures that $H_\ell$ extends by continuity to a vector field on $\Tetb^*\tilde M$, with restriction to $\Tetb^*_{\pa\tilde M}\tilde M$ given by $H_{\ell_0}$ where $\ell_0=\upsigma_\etbop^m(L_0)$. The subprincipal term $L_2$ need not vanish at $\pa\tilde M$; in applications, we will encounter $L_2$ which do vanish at $\iota^+\cap\cK^+$ but not at $\scri^+$.

\begin{proof}[Proof of Proposition~\usref{PropMPr}]
  The case $k=0$, $d_0=\infty$ can be established using a standard positive commutator argument utilizing the algebra $\tilde\Psi_\etbop(\tilde M)$, and it is possible to take $d_0$ to be finite (depending on $m,\sfs,N$) since this argument only requires the boundedness of ps.d.o.s on a fixed range of e3b-Sobolev spaces. Variable orders are discussed in Vasy \cite[Theorem~5.3]{VasyMinicourse}. (Note here that the sharp G\aa{}rding inequality (cf.\ \cite[p.~294]{VasyMinicourse} is available in the b.g.\ setting: it follows easily by using the standard sharp G\aa{}rding inequality (see, e.g., \cite[Theorem~18.1.14]{HormanderAnalysisPDE3}) on each unit cell.) A robust way of choosing a commutant is described in \cite[\S{2}]{DeHoopUhlmannVasyDiffraction} following \cite{MelroseSjostrandSingBVPI}: this does not require straightening out the $\sfH_\ell$-flow and therefore works directly also for perturbations of $L$, thus immediately giving the desired local uniformity statement.

  Instead of working semi-globally using $\tilde\Psi_\etbop(\tilde M)$, one can note that the $\sfH_\ell$-flow, up to any fixed affine parameter, traverses only a subset of phase space over, at most, a \emph{fixed finite number} of e3b-unit cells adjacent to the unit cell containing the starting point; this follows from the fact that $\sfH_\ell$ is an e3b-vector field itself (analogously to~\eqref{EqCHame3bMem}). Therefore,~\eqref{EqMPrEst} follows by summing \emph{standard} real principal type propagation estimates on $(-2,2)^n$ (with variable orders, and thus appealing to sharp G\aa{}rding \emph{on Euclidean space}) for uniformly bounded families of ps.d.o.s on $\R^n$. The usual proofs (e.g., by straightening out the $\sfH_\ell$-flow as in \cite[\S{8.3}]{HintzMicro}) give estimates which are locally uniform upon varying the underlying operator within $\cC^{d_0}\sP$ where $\sP\subset\Psi^m(\R^n)$ is a bounded (not necessarily finite) subset.\footnote{This perspective thus yields a stronger local uniformity statement than the one made in the statement of Proposition~\ref{PropMPr}. However, for general elements of $\tilde\Psi_\etbop^m(\tilde M)$ the Hamiltonian flow does not extend continuously to the boundary of $\tilde M$, which does not mesh well with our insistence, in the present paper, of using ${}^\etbop S^*\tilde M$ as the locus of e3b-microlocal analysis.}

  To treat the case $k\geq 1$, we argue inductively. Fix a \emph{good spanning set} $\{V_1,\ldots,V_N\}\subset\cV_{\bop,[\etbop]}(\tilde M)$ of e3b-commutator b-vector fields (see Lemma~\ref{LemmaCTe3bComm}), and set $V_0:=I$. Using Lemma~\ref{LemmaMUe3bExpb} (and noting that the commutator of $V_j$ with a $\cC_{\etbop;\bop}^{(d_0;k)}$-coefficient of $L_1,L_2$ produces a coefficient of class $\cC_{\etbop;\bop}^{(d_0;k-1)}$), we can then write
  \[
    [L,V_j] = \sum_{k=0}^N A_{j,l}V_l,\quad A_{j,l}\in\CI_\bop\Psi_\etbop^{m-1} + \cC_{\etbop;\bop}^{(d_0;k-1),(2\delta,\delta,\delta)}\Psi_\etbop^{m-1} + \cC_{\etbop;\bop}^{(d_0;k-1)}\Psi_\etbop^{m-2}(\tilde M)
  \]
  Therefore, if we write $L u=f$, we have
  \[
    L^{(1)}u^{(1)}=f^{(1)},\quad u^{(1)}:=(V_j u)_{j=0}^N,\ f^{(1)}:=(V_j f)_{j=0}^N,\quad L^{(1)}=(L\delta_{j,l}-A_{j,l})_{j,l=0,\ldots,N}.
  \]
  Note then that the terms $A_{j,l}$ are subprincipal, so the principal symbol of $L^{(1)}$ is scalar and equal to that of $L$ itself. We can thus apply the inductive hypothesis to this equation, thus obtaining~\eqref{EqMPrEst} for $u^{(1)},L^{(1)},k-1$ in place of $u,L,k$. This gives~\eqref{EqMPrEst} for $B$, resp.\ $G,E$ with slightly shrunk, resp.\ enlarged elliptic sets, upon commuting $V_j$ through $B,G L,E$ similarly to the proof of Proposition~\ref{PropMEll}.
\end{proof}

\subsection{Basic b-tame e3b-microlocal estimates}
\label{SsMTame}

We next prove versions of Propositions~\ref{PropMEll} and \ref{PropMPr} which are \emph{tame} in the b-regularity order $k$, in that the constants in the estimates~\eqref{EqMEll} and \eqref{EqMPrEst} depend only on a fixed low b-regularity norm of the coefficients of $L$. Our approach is closely related to \cite[\S{6.1.1}]{HintzGlueLocII}, though due to the number of required minor modifications (since the function spaces and Lie algebras differ), we present the details here. We first record the following simple result on the structure of commutators (which is the same as \cite[Lemma~3.34]{HintzMink4Gauge} up to sign changes):

\begin{lemma}[Commutators]
\label{LemmaMTameComm}
  Let $R$ be a ring. For $A\in R$, write ${\rm ad}_A=[A,\cdot]\colon B\mapsto A B-B A$. Let $N\in\N$ and $P,X_1,\ldots,X_N\in R$. Then
  \begin{align}
  \label{EqMTameCommR}
    [X_1\cdots X_N,P] &= \sum_{q=1}^N \sum_{1\leq i_1<\cdots<i_q\leq N} ({\rm ad}_{X_{i_1}}\cdots{\rm ad}_{X_{i_q}}P) \prod_{\substack{j=1,2,\ldots,N \\ j\neq i_1,\ldots,i_q}} X_j \\
  \label{EqMTameCommL}
      &= \sum_{q=1}^N (-1)^{q-1} \sum_{1\leq i_1<\cdots<i_q\leq N} \Biggl(\;\prod_{\substack{j=1,2,\ldots,N \\ j\neq i_1,\ldots,i_q}} X_j\Biggr) ({\rm ad}_{X_{i_q}}\cdots{\rm ad}_{X_{i_1}}P).
  \end{align}
\end{lemma}
\begin{proof}
  The case $N=1$ is clear. Assuming~\eqref{EqMTameCommR} and letting $X_{N+1}\in R$, we write
  \[
    [X_1\cdots X_{N+1},P] = [X_1\cdots X_N,P]X_{N+1} + ({\rm ad}_{X_{N+1}}P) X_1\cdots X_N + [X_1\cdots X_N,{\rm ad}_{X_{N+1}}P]
  \]
  Applying the inductive hypothesis to the first and third term produces~\eqref{EqMTameCommR} for $N+1$ in place of $N$. To prove~\eqref{EqMTameCommL} for $N+1$, we apply the inductive hypothesis to
  \[
    [X_1\cdots X_{N+1},P] = X_1[X_2\cdots X_{N+1},P] + X_2\cdots X_{N+1}({\rm ad}_{X_1}P) - [X_2\cdots X_{N+1},{\rm ad}_{X_1}P].\qedhere
  \]
\end{proof}

We will apply this when the $X_i$ are b-vector fields. We introduce the following notation:

\begin{notation}[b-vector fields]
\label{NotMTameb}
  Fix any finite subset $\sV_\sharp\subset\Vb(\tilde M)$ that spans $\Vb(\tilde M)$ over $\CI(\tilde M)$. Fix moreover a \emph{good spanning set} $\sV=\{V_1,\ldots,V_N\}\subset\cV_{\bop,[\etbop]}(\tilde M)$ of e3b-commutator b-vector fields (see Lemma~\usref{LemmaCTe3bComm}). We then write $D_\bop^j$ ($D_\bop^{\leq j}$), resp.\ $D_{[\bop]}^j$ ($D_{[\bop]}^{\leq j}$) for the vector of all (up to) $j$-fold compositions of the elements of $\sV_\sharp$, resp.\ $\sV$, or of $\CI(\tilde M)$-linear combinations thereof.
\end{notation}

Lemma~\ref{LemmaMTameComm} thus reads, schematically,
\[
  [D_\bop^N,P] = \sum_{q=1}^N ({\rm ad}_{D_\bop}^q P) D_\bop^{N-q} = \sum_{q=1}^N D_\bop^{N-q} ({\rm ad}_{D_\bop}^q),
\]
and similarly with $D_{[\bop]}$ in place of $D_\bop$.

\subsubsection{Bounds for products}

Recall that we typically use the Minkowskian volume density to define $L^2$; one can equally well use weighted versions of this, with weights being real powers of $\rho_0,x_\sscri,\rho_+,\rho_\cK$. We then record the following basic multiplication estimate:

\begin{lemma}[Basic multiplication bound]
\label{LemmaMTameMult}
  Given $a,b\in\N_0$, there exists a constant $C_{a,b}$ such that
  \[
    \|(D_\bop^a\ell)(D_\bop^b u)\|_{L^2(\tilde M)} \leq C_{a,b} \Bigl( \|\ell\|_{\cC_\bop^0(\tilde M)}\|u\|_{\Hb^{a+b}(\tilde M)} + \|\ell\|_{\cC_\bop^{a+b}(\tilde M)}\|u\|_{\Hb^3(\tilde M)} \Bigr).
  \]
\end{lemma}
\begin{proof}
  The square of the left-hand side is the sum of $L^2(\R^4)$-norms of $\chi(\phi_\alpha)_*((D_\bop^a\ell)(D_\bop^b u))$ where $\chi\in\CIc((-2,2)^4)$ equals $1$ on $[-1,1]^4$ and the $\phi_\alpha$ are the charts corresponding to the unit cells~\eqref{EqMUe3bFam1}--\eqref{EqMUe3bFam3}. Now, $(\phi_\alpha)_*D_\bop=D$ (schematic notation for a vector of coordinate derivatives, with equality meaning equality up to taking uniformly bounded linear combinations over $\CI$). If $\tilde\chi\in\CIc((-2,2)^4)$ equals $1$ on $\supp\chi$, we may replace $(\phi_\alpha)_*\ell$, $(\phi_\alpha)_*u$ here by $f:=\tilde\chi(\phi_\alpha)_*\ell$, $g:=\tilde\chi(\phi_\alpha)_*u$, and thus need to show
  \[
    \|(D^{\leq a}f)(D^{\leq b}g)\|_{L^2(\R^4)}\lesssim\|f\|_{L^\infty}\|g\|_{H^{a+b}}+\|f\|_{\cC^{a+b}}\|g\|_{H^3}
  \]
  for $f,g$ with compact support in $(-2,2)^4$. This follows from classical Moser estimates (see, e.g., \cite[Chapter~13, Proposition~3.6]{TaylorPDE3}) together with $\|f\|_{H^{a+b}}\leq C\|f\|_{\cC^{a+b}}$ and Sobolev embedding $\|g\|_{L^\infty(\R^4)}\leq C\|g\|_{H^3(\R^4)}$.
\end{proof}

Here and below, we only state estimates on unweighted spaces; analogous on weighted spaces require only notational modifications. The next result is the analogue of \cite[Lemma~6.2]{HintzGlueLocII}:

\begin{prop}[b-tame multiplication, e3b-microlocalized]
\label{PropMTameMicr}
  Let $\sfs,\sfs_0\in\CI({}^\etbop S^*\tilde M)$. Then there exist $d_0,\tilde d\in\N$ such that the following holds for all $B,B^\sharp\in\CI_\bop\Psi_\etbop^0(\tilde M)$, $\chi\in\CI(\tilde M)$ with $B=\chi B\chi$, $B^\sharp=\chi B^\sharp\chi$, and $\WF_\etbop'(B)\subset\Ell_\tbop(B^\sharp)$: for all $k\in\N_0$, there exists a constant $C_k$ such that for all $0\leq j\leq k$,
  \begin{equation}
  \label{EqMTameMicr}
  \begin{split}
    \|B(D_\bop^j\ell)(D_\bop^{k-j}u)\|_{H_\etbop^\sfs(\tilde M)} &\leq C_k\biggl( \|\ell\|_{\cC_{\etbop;\bop}^{(d_0;\tilde d)}(\tilde M)}\Bigl( \|B^\sharp u\|_{H_{\etbop;\bop}^{(\sfs;k)}(\tilde M)} + \|\chi u\|_{H_{\etbop;\bop}^{(\sfs_0;k)}}\Bigr) \\
      &\quad \hspace{1.0em} + \|\ell\|_{\cC_{\etbop;\bop}^{(d_0;k)}(\tilde M)} \Bigl( \|B^\sharp u\|_{H_{\etbop;\bop}^{(\sfs;\tilde d)}(\tilde M)} + \|\chi u\|_{H_{\etbop;\bop}^{(\sfs_0;\tilde d)}}\Bigr)\biggr).
  \end{split}
  \end{equation}
\end{prop}

In our applications, the e3b-orders $\sfs,\sfs_0$ will be chosen (and then fixed) to be compatible with all e3b-microlocal regularity results, while $k$ will be large. The point of~\eqref{EqMTameMicr} is that the products of norms on the right involve only high-low and low-high norms, never high-high norms. 

\begin{proof}[Proof of Proposition~\usref{PropMTameMicr}]
  It suffices to prove~\eqref{EqMTameMicr} with $B^\sharp$ replaced by an operator $B^\flat\in\CI_\bop\Psi_\etbop^0$ such that $\WF'_\etbop(B)\subset\Ell_\etbop(B^\flat)$ and $\WF'_\etbop(B^\flat)\subset\Ell_\etbop(B^\sharp)$; we may in particular take $B^\flat$ such that $\WF'_\etbop(I-B^\flat)\cap\WF'_\etbop(B)=\emptyset$. Relabeling operators, we shall henceforth assume that $\WF'_\etbop(I-B^\sharp)\cap\WF'_\etbop(B)=\emptyset$, i.e., $B^\sharp$ is microlocally the identity operator near $\WF'_\etbop(B)$. We may moreover replace $u$ by $\chi u$.

  \pfstep{$k=0$.} Write $B\ell u=B\ell B^\sharp u+B\ell(I-B^\sharp)u$. For $\ell\in\cC_\etbop^\infty$, multiplication by $\ell$, which constitutes an element of $\tilde\Psi_\etbop^0$, is bounded on $H_\etbop^\sfs$, and hence the same is true for $\ell\in\cC_\etbop^{d_0}$ when $d_0$ is large enough (depending on $\sfs$). Using the boundedness of $B$ on e3b-Sobolev spaces, we thus get
  \[
    \|B\ell B^\sharp u\|_{H_\etbop^s} \leq C\|\ell\|_{\cC_\etbop^{d_0}}\|B^\sharp u\|_{H_\etbop^s}.
  \]
  Again for $\ell\in\cC_\etbop^\infty$, we have $B\ell(I-B^\sharp)\in\tilde\Psi_\etbop^{-\infty}$ since its e3b-operator wave front set is empty; and the operator norm of this operator as a map $H_\etbop^{\sfs_0}\to H_\etbop^\sfs$ is bounded by a finite seminorm of $\ell$. Since, moreover, $B\ell(I-B^\sharp)$ depends linearly on $\ell$, this gives
  \[
    \|B\ell(I-B^\sharp)u\|_{H_\etbop^\sfs} \leq C\|\ell\|_{\cC_\etbop^{d_0}}\|u\|_{H_\etbop^{\sfs_0}}.
  \]
  This proves~\eqref{EqMTameMicr} for $k=0$.
  
  \pfstep{$k\leq\tilde d$,} with the value of $\tilde d$ specified later in the proof. We apply the case $k=0$ to $D_\bop^j\ell$, $D_\bop^{k-j}u$ in place of $\ell$, $u$ to get (using $j,k-j\leq\tilde d$)
  \[
    \|B(D_\bop^j\ell)(D_\bop^{k-j}u)\|_{H_\etbop^\sfs} \leq C\biggl( \|\ell\|_{\cC_{\etbop;\bop}^{(d_0;\tilde d)}}\Bigl(\|B^\sharp(D_\bop^{k-j}u)\|_{H_\etbop^\sfs} + \|u\|_{H_{\etbop;\bop}^{(\sfs_0;\tilde d)}}\Bigr) \biggr).
  \]
  In view of Lemma~\ref{LemmaCTe3bComm}, we can replace $D_\bop$ on the right by $D_{[\bop]}$. Recalling~\eqref{EqMUe3bComm}, we commute $D_{[\bop]}^{k-j}$ through $B^\sharp$; writing ${\rm ad}_{D_{[\bop]}}=[D_{[\bop]},\cdot]$, we thus have
  \[
    B^\sharp D_{[\bop]}^{k-j} = \sum_{m=0}^{k-j} D_{[\bop]}^{k-j-m}({\rm ad}_{D_{[\bop]}}^m B^\sharp)
  \]
  and can furthermore write ${\rm ad}_{D_{[\bop]}}^m B^\sharp=Q_m\tilde B^\sharp+R_m$ using an e3b-microlocal elliptic parametrix, where $\tilde B^\sharp\in\CI_\bop\Psi_\etbop^0$ has elliptic set containing $\WF'_\etbop(B^\sharp)$, and $R_m\in\CI_\bop\Psi_\etbop^{-\infty}$. Using the boundedness of elements of $\CI_\bop\Psi_\etbop$ on mixed $(\etbop;\bop)$-Sobolev spaces thus shows that
  \begin{equation}
  \label{EqMTameMicrSmall}
    \|B^\sharp(D_{[\bop]}^{k-j}u)\|_{H_\etbop^\sfs} \leq C\Bigl(\|\tilde B^\sharp u\|_{H_{\etbop;\bop}^{(\sfs;k-j)}} + \|u\|_{H_{\etbop;\bop}^{(\sfs_0;k-j)}}\Bigr).
  \end{equation}
  This proves~\eqref{EqMTameMicr}, except with $\tilde B^\sharp$ in place of $B^\sharp$---which is easily rectified by beginning the argument with an operator $B^\sharp$ that has smaller operator wave front set.

  \pfstep{$k>\tilde d$.} For $j\leq\tilde d$, we can bound the left-hand side of~\eqref{EqMTameMicr} by using the case $k=0$ and $D_\bop^j\ell$, $D_\bop^{k-j}u$ in place of $\ell,u$, which (using the arguments around~\eqref{EqMTameMicrSmall}) gives a bound by a constant times
  \[
    \|\ell\|_{\cC_{\etbop;\bop}^{(d_0;j)}} \Bigl( \|B^\sharp u\|_{H_{\etbop;\bop}^{(\sfs;k-j)}} + \|u\|_{H_{\etbop;\bop}^{(\sfs_0;k-j)}}\Bigr),
  \]
  which in turn is bounded by the first line on the right in~\eqref{EqMTameMicr}. For $j\geq k-\tilde d$, we again use~\eqref{EqMTameMicr} for $k=0$ with $D_\bop^j\ell$, $D_\bop^{k-j}u$ in place of $\ell,u$; this bounds the left-hand side of~\eqref{EqMTameMicr} by the second line on the right.

  It remains to consider the case $\tilde d<j<k-\tilde d$. The idea is that since the `low' norms on the right-hand side in~\eqref{EqMTameMicr} do involve $\tilde d$ orders of b-regularity, we can be rather wasteful and replace variable order e3b- by integer order b-derivatives. Let $\bar s:=\max(0,\lceil\sup\sfs\rceil)$. Replace $D_\bop$ by $D_{[\bop]}$. We then bound
  \begin{align*}
    \|B(D_{[\bop]}^j\ell)(D_{[\bop]}^{k-j}u)\|_{H_\etbop^\sfs} &\leq C\|B(D_{[\bop]}^j\ell)(D_{[\bop]}^{k-j}u)\|_{\Hb^{\bar s}} \leq \sum_{l=0}^{\bar s} \bigl\|D_{[\bop]}^l \bigl(B(D_{[\bop]}^j\ell)(D_{[\bop]}^{k-j}u)\bigr)\bigr\|_{L^2}.
  \end{align*}
  We then expand $D_{[\bop]}^l B=\sum_{l'=0}^l ({\rm ad}_{D_{[\bop]}}^{l-l'}B) D_{[\bop]}^{l'}$, with ${\rm ad}_{D_{[\bop]}}^{l-l'}B\in\CI_\bop\Psi_\etbop^0$. We only consider the top order term with $l'=l$ (as it involves the highest number of b-derivatives on $\ell,u$), and we use the Leibniz rule to distribute the $l$ b-derivatives onto the two factors involving $\ell,u$. Arguing as before~\eqref{EqMTameMicrSmall}, it then suffices to bound
  \[
    \sum_{l+q\leq\bar s} \|B(D_{[\bop]}^{j+l}\ell)(D_{[\bop]}^{k-j+q}u)\|_{L^2}.
  \]
  We only consider the term with the maximal number of derivatives, i.e., $q=\bar s-l$. Writing $u=B^\sharp u+(I-B^\sharp)u$, and using the $L^2$-boundedness of $B$, we thus need to bound
  \begin{equation}
  \label{EqMTameMicr0}
    \|(D_{[\bop]}^{j+l}\ell)(D_{[\bop]}^{k-j+\bar s-l}(B^\sharp u))\|_{L^2} + \|(D_{[\bop]}^{j+l}\ell)(D_{[\bop]}^{k-j+\bar s-l}((I-B^\sharp)u))\|_{L^2}.
  \end{equation}

  In the first term of~\eqref{EqMTameMicr0}, we use Lemma~\ref{LemmaMTameMult} with $a=j+l-d_1$, $b=k-j+\bar s-l-d_1$ (so $a+b=k+\bar s-2 d_1$) and $D_{[\bop]}^{d_1}\ell$, $D_{[\bop]}^{d_1}(B^\sharp u)$ where $d_1\in\N_0$, $2 d_1\leq\tilde d$ (so $2 d_1<k$), will be chosen momentarily: it is thus bounded by
  \begin{equation}
  \label{EqMTameMicr1}
  \begin{split}
    &\| D_{[\bop]}^{j+l-d_1}(D_{[\bop]}^{d_1}\ell)\cdot D_{[\bop]}^{k-j+\bar s-l-d_1}(D_{[\bop]}^{d_1}B^\sharp u) \|_{L^2} \\
    &\qquad \leq C\Bigl( \|\ell\|_{\cC_\bop^{d_1}}\|B^\sharp u\|_{H_\bop^{k+\bar s-d_1}} + \|\ell\|_{\cC_\bop^{k+\bar s-d_1}} \|B^\sharp u\|_{\Hb^{d_1+3}} \Bigr).
  \end{split}
  \end{equation}
  Let now $\ubar s=\max(0,\lceil-\inf\sfs\rceil)$; then Lemma~\ref{LemmaMUCe3bTob} shows that we can further estimate
  \[
    \|B^\sharp u\|_{\Hb^{k+\bar s-d_1}} \leq C\|B^\sharp u\|_{H_{\etbop;\bop}^{(\sfs;k+\bar s+\ubar s-d_1)}},\quad
    \|B^\sharp u\|_{\Hb^{d_1+3}} \leq C\|B^\sharp u\|_{H_{\etbop;\bop}^{(\sfs;d_1+\ubar s+3)}}.
  \]
  If we choose $d_1$ and then $\tilde d$ such that
  \begin{equation}
  \label{EqMTameMicrCond}
    d_1 \geq \bar s+\ubar s\quad\text{and}\quad 2 d_1\leq\tilde d,\ \ d_1+\ubar s+3\leq\tilde d,
  \end{equation}
  we then have $\|B^\sharp u\|_{\Hb^{k+\bar s-d_1}}\leq C\|B^\sharp u\|_{H_{\etbop;\bop}^{(\sfs;k)}}$ and $\|B^\sharp u\|_{\Hb^{d_1+3}}\leq C\|B^\sharp u\|_{H_{\etbop;\bop}^{(\sfs;\tilde d)}}$. Since then also $k+\bar s-d_1\leq k$, we have shown that the right-hand side of~\eqref{EqMTameMicr1} is bounded by the right-hand side of~\eqref{EqMTameMicr}.

  We can argue similarly for the second term in~\eqref{EqMTameMicr0}, except now we relate b- to $(\etbop;\bop)$-Sobolev spaces with e3b-regularity order $\sfs_0$; thus, letting $\ubar s_0=\max(0,\lceil-\inf\sfs_0\rceil)$, we need to, in addition, require $d_1\geq\bar s+\ubar s_0$ and $2 d_1\leq\tilde d$, $d_1+\ubar s_0+3\leq\tilde d$. For any $\tilde d$ satisfying these conditions as well as~\eqref{EqMTameMicrCond}, we have now proved~\eqref{EqMTameMicr}.
\end{proof}

\subsubsection{Microlocal elliptic regularity and real principal type propagation}
\label{SssMTame}

We first prove a tame version of Proposition~\ref{PropMEll} on microlocal elliptic regularity when the b-regularity order exceeds a certain threshold $d$. (Cf.\ \cite[Lemma~6.4]{HintzGlueLocII}.)

\begin{prop}[b-tame e3b-microlocal elliptic regularity]
\label{PropMTameEll}
  Let $m,N\in\R$, $\sfs\in\CI({}^\etbop S^*\tilde M)$. Then there exist $d_0,d\in\N$ such that the following holds for all $k\in\N_0$ and $L\in\cC_{\etbop;\bop}^{(d_0;d)}\sP$ where $\sP\subset\CI_\bop\Psi_\etbop^m(\tilde M)$ is a finite set. Let $\chi\in\CI(\tilde M)$, $B,G\in\CI_\bop\Psi_\etbop^0(\tilde M)$ with $B=\chi B\chi$, $G=\chi G\chi$, $\WF'_\etbop(B)\subset\Ell_\etbop(G)\cap\Ell_\etbop(L)$. Then there exists a constant $C_k$ such that, in the notation of~\eqref{EqMUCe3bOpNormPsdo},
  \begin{equation}
  \label{EqMTameEll}
  \begin{split}
    \|B u\|_{H_{\etbop;\bop}^{(\sfs;k)}} &\leq C_k\biggl( \|G L u\|_{H_{\etbop;\bop}^{(\sfs-m;k)}} + \|\chi u\|_{H_{\etbop;\bop}^{(-N;k)}} \\
      &\quad\hspace{3em} + \|L\|_{\cC_{\etbop;\bop}^{(d_0;k)}\Psi_\etbop^m} \Bigl( \|G L u\|_{H_{\etbop;\bop}^{(\sfs-m;d)}} + \|\chi u\|_{H_{\etbop;\bop}^{(-N;d)}}\Bigr) \biggr),
  \end{split}
  \end{equation}
  provided $L\in\cC_{\etbop;\bop}^{(d_0;k)}\sP$. The constant $C_k$ can be taken to be locally uniform for perturbations of $L$ in $\cC_{\etbop;\bop}^{(d_0;d)}\sP$.
\end{prop}

Crucially, the local uniformity holds for perturbations of $L$ within a class of operators with a \emph{fixed}, ``small'' amount of b-regularity, and then for \emph{all} $k$ at once.

\begin{proof}[Proof of Proposition~\usref{PropMTameEll}]
  We take $d=\tilde d+1$ in the notation of Proposition~\ref{PropMTameMicr}. We can use Proposition~\ref{PropMEll} for $k\leq d$. For arbitrary $k$, we then note that, schematically,
  \[
    L(D_{[\bop]}^k u) = D_{[\bop]}^k L u + \sum_{j=1}^k ({\rm ad}_{D_{[\bop]}}^j L) D_{[\bop]}^{k-j}u.
  \]
  Applying Proposition~\ref{PropMEll}, with $k=0$ there, to this equation gives
  \[
    \|B D_{[\bop]}^k u\|_{H_\etbop^\sfs} \leq C_k \Bigl( \|G D_{[\bop]}^k L u\|_{H_\etbop^{\sfs-m}} + \|\chi D_{[\bop]}^k u\|_{H_\etbop^{-N}} \Bigr) + C_k \sum_{j=1}^k \|G ({\rm ad}_{D_{[\bop]}}^j L)D_{[\bop]}^{k-j}u \|_{H_\etbop^{\sfs-m}}.
  \]
  In the first term on the right, we can commute the e3b-commutator b-vector fields $D_{[\bop]}$ through $G$ and $\chi$ as before; it is thus bounded by
  \[
    \|G^\sharp L u\|_{H_{\etbop;\bop}^{(\sfs-m;k)}} + \|\tilde\chi u\|_{H_{\etbop;\bop}^{(-N;k)}}
  \]
  where $G^\sharp\in\CI_\bop\Psi_\etbop^0$ has elliptic set containing the operator wave front set of $G$, and $\tilde\chi$ has a slightly enlarged supported compared to $\chi$. As usual, one can revert back to using $G,\chi$ by starting the argument with smaller $G,\chi$; we shall not comment on this anymore in the sequel.

  It remains to prove a tame bound for the final term. Expanding $L=\sum \ell_q A_q$, $\ell_q\in\cC_{\etbop;\bop}^{(d_0;k)}$, $A_q\in\sP$, and using that (iterated) commutators of $D_{[\bop]}$ with the $A_q$ lie in $\CI_\bop\Psi_\etbop^m$, we thus need to estimate, for $A\in\CI_\bop\Psi_\etbop^m$, norms of the form
  \[
    \|G (D_{[\bop]}^j\ell_q) A\,D_{[\bop]}^{k-j}u\|_{H_\etbop^{\sfs-m}},\quad j=1,\ldots,k.
  \]
  We commute $A$ through $D_{[\bop]}^{k-j}$. We may then use Proposition~\ref{PropMTameMicr} for $D_{[\bop]}\ell_q$, $A u$, $j-1$, $k-1$ (in present notation) in place of $\ell,u,j,k$ to bound
  \begin{align*}
    \|G (D_{[\bop]}^{j-1}D_{[\bop]}\ell_q) D_{[\bop]}^{k-j}(A u)\|_{H_\etbop^{\sfs-m}} &\leq C_k \biggl( \|\ell_q\|_{\cC_{\etbop;\bop}^{(d_0;\tilde d+1)}}\Bigl( \|G^\sharp u\|_{H_{\etbop;\bop}^{(\sfs;k-1)}} + \|\chi u\|_{H_{\etbop;\bop}^{(-N;k-1)}}\Bigr) \\
      &\qquad \qquad + \|\ell_q\|_{\cC_{\etbop;\bop}^{(d_0;k)}}\Bigl( \|G^\sharp u\|_{H_{\etbop;\bop}^{(\sfs;\tilde d)}} + \|\chi u\|_{H_{\etbop;\bop}^{(\sfs_0;\tilde d)}}\Bigr) \biggr).
  \end{align*}
  (The b-regularity order on $\ell_q$ here is why we take $d=\tilde d+1$.) Choosing, as we may, the operator wave front set of $G^\sharp$ to be contained in the elliptic set of $L$, the norm $\|G^\sharp u\|_{H_{\etbop;\bop}^{(\sfs;k-1)}}$ can be estimated inductively. This gives~\eqref{EqMTameEll}.
\end{proof}

The tame version of Proposition~\ref{PropMPr} reads as follows. (This is similar to \cite[Lemma~6.5]{HintzGlueLocII}.)

\begin{prop}[b-tame e3b-microlocal real principal propagation estimate]
\label{PropMTamePr}
  Let $m\in\N_0$, $N\in\R$, $\sfs\in\CI({}^\etbop S^*\tilde M)$. Then there exist $d_0,d\in\N$ such that the following holds for all $\delta>0$, $k\in\N_0$, finite sets $\sP^{m-j}\subset\CI_\bop\Psi_\etbop^{m-j}$, $j=0,1$, and
  \[
    L=L_0+L_1+L_2,\quad L_0\in\Diff_\etbop^m,\ L_1\in\cC_{\etbop;\bop}^{(d_0;d),(2\delta,\delta,\delta)}\sP^m,\ L_2\in\cC_{\etbop;\bop}^{(d_0;d)}\cP^{m-1},
  \]
  where $L_1$ has a homogeneous principal symbol. Assume that the principal symbol $\ell$ on $L$ is real-valued, and, recalling~\eqref{EqMPrShift}, write $\sfH=\rho_\infty^{m-1}H_\ell$ for its rescaled Hamiltonian vector field. Let $\chi\in\CI(\tilde M)$, and let $B,E,G\in\CI_\bop\Psi_\etbop^0$ be such that $B=\chi B\chi$, $E=\chi E\chi$, $G=\chi G\chi$. Suppose that all backwards integral curves of $\sfH_\ell$ starting at $\WF'_\etbop(B)$ reach $\Ell_\etbop(E)$ in uniformly finite time while remaining in $\Ell_\etbop(G)$, and that $\sfH_\ell\sfs\leq 0$ in a neighborhood of $\WF'_\etbop(B)\cap\Char_\etbop(L)$ within $\Char_\etbop(L)$. Suppose also that $\WF'_\etbop(B)\subset\Ell_\etbop(G)$. Then there exists a constant $C_k$ such that, in the notation of~\eqref{EqMUCe3bOpNormPsdo},
  \begin{equation}
  \label{EqMTamePrEst}
  \begin{split}
    \|B u\|_{H_{\etbop;\bop}^{(\sfs;k)}} &\leq C_k\biggl( \|G L u\|_{H_{\etbop;\bop}^{(\sfs-m;k)}} + \|E u\|_{H_{\etbop;\bop}^{(\sfs;k)}} + \|\chi u\|_{H_{\etbop;\bop}^{(-N;k)}} \\
      &\quad\hspace{3em} + \|L\|_{\cC_{\etbop;\bop}^{(d_0;k)}\Psi_\etbop^m} \Bigl( \|G L u\|_{H_{\etbop;\bop}^{(\sfs-m;d)}} + \|E u\|_{H_{\etbop;\bop}^{(\sfs;d)}} + \|\chi u\|_{H_{\etbop;\bop}^{(-N;d)}}\Bigr) \biggr),
  \end{split}
  \end{equation}
  provided $L_1\in\cC_{\etbop;\bop}^{(d_0;k),(2\delta,\delta,\delta)}\Psi_\etbop^m$, $L_2\in\cC_{\etbop;\bop}^{(d_0;k)}\Psi_\etbop^{m-1}$. The constant $C_k$ can be taken to be locally uniform for perturbations of $L_0$, $L_1$, and $L_2$ in $\Diff_\etbop^m$, $\cC_{\etbop;\bop}^{(d_0;d),(2\delta,\delta,\delta)}\sP^m$, and $\cC_{\etbop;\bop}^{(d_0;d)}\cP^{m-1}$, respectively.
\end{prop}
\begin{proof}
  We take $d=\tilde d+2$ in the notation of Proposition~\ref{PropMTameMicr}. We can use Proposition~\ref{PropMPr} for $k\leq d$. For larger $k$, we argue inductively. Write $\sV=\{V_1,\ldots,V_N\}$ for a good spanning set of $\cV_{\bop,[\etbop]}(\tilde M)$ (see Lemma~\ref{LemmaCTe3bComm}) and set $V_0=I$. For each multiindex $\alpha\in\N_0^{N+1}$ with $|\alpha|=k$, we then use Lemma~\ref{LemmaMTameComm} and write
  \begin{equation}
  \label{EqMTamePrExp}
    L(\sV^\alpha u) = \sV^\alpha L u + \sum_{\substack{\beta+\gamma=\alpha \\ |\beta|=1}} c_{\beta\gamma} ({\rm ad}_\sV^\beta L) \sV^\gamma u + \sum_{\substack{\beta+\gamma=\alpha \\ |\beta|\geq 2}} c_{\beta\gamma} ({\rm ad}_\sV^\beta L) \sV^\gamma u
  \end{equation}
  where the $c_{\beta\gamma}$ are combinatorial constants. Consider the first sum, which involves only first order b-derivatives of the coefficients of $L$: using Lemma~\ref{LemmaMUe3bExpb}, we may write the $m$-th order operator ${\rm ad}_\sV^\beta L$ as
  \[
    {\rm ad}_\sV^\beta L = \sum_{|\nu|=1} L_{\beta,\nu}\sV^\nu,\quad L_{\beta,\nu}\in\CI_\bop\Psi_\etbop^{m-1} + \cC_{\etbop;\bop}^{(d_0;d-1)}\Psi_\etbop^{m-1}+\cC_{\etbop;\bop}^{(d_0;d-1)}\Psi_\etbop^{m-2}.
  \]
  The vector $u^{(k)}:=(\sV^\alpha u)_{|\alpha|=k}$ thus satisfies the equation
  \[
    L^{(k)}u^{(k)} = f^{(k)} + \tilde f^{(k)},\quad f^{(k)}=(\sV^\alpha L u)_{|\alpha|=k},\ \tilde f^{(k)}:= \Biggl(\;\sum_{\substack{\beta+\gamma=\alpha \\ |\beta|\geq 2}} c_{\beta\gamma}({\rm ad}_\sV^\beta L)\sV^\gamma u\Biggr)_{|\alpha|=k},
  \]
  where $L^{(k)}=(L \delta_{\alpha,\alpha'}-\sum_{\alpha=\alpha'-\beta+\nu,\ |\beta|=|\nu|=1}c_{\beta\gamma}L_{\beta,\nu})_{|\alpha|,|\alpha'|=k}$. (This operator depends, of course, on $k$; notice moreover that it is built from commutators of elements of $\Diff_\etbop^m$, $\sP^m$, and $\sP^{m-1}$ with elements of $\sV$, and hence, even for perturbations of $L$ in the sense stated after~\eqref{EqMTamePrEst}, is built from a fixed finite set of e3b-ps.d.o.s.) Note that $L_{\beta,\nu}$, being of order $m-1$, is subprincipal; moreover, its coefficients involve only $1$ b-derivative, regardless of the value of $k$. We may thus apply Proposition~\ref{PropMTameMicr} to this equation to obtain
  \[
    \|B u^{(k)}\|_{H_\etbop^\sfs} \leq C\Bigl( \|G f^{(k)}\|_{H_\etbop^{\sfs-m+1}} + \|G\tilde f^{(k)}\|_{H_\etbop^{\sfs-m+1}} + \|E u^{(k)}\|_{H_\etbop^\sfs} + \|\chi u^{(k)}\|_{H_\etbop^{-N}} \Bigr).
  \]
  The left-hand side, together with the inductive hypothesis, controls $\|B^\flat u\|_{H_{\etbop;\bop}^{(\sfs;k)}}$ for any $B^\flat\in\CI_\bop\Psi_\etbop^0$ with operator wave front set contained in $\Ell_\etbop(B)$, up to the usual microlocal error term $\|\chi u^{(k)}\|_{H_{\etbop;\bop}^{(-N;k)}}$. On the right-hand side, we can similarly control $\|G f^{(k)}\|_{H_\etbop^{\sfs-m+1}}$ by $\|G^\sharp f\|_{H_{\etbop;\bop}^{(\sfs-m+1;k)}}$ (plus the usual microlocal error term) where $G^\sharp\in\CI_\bop\Psi_\etbop^0$ has elliptic set containing the operator wave front set of $G$; likewise for $E u^{(k)}$.

  It remains to bound $G\tilde f^{(k)}$, i.e., schematically, $\| G ({\rm ad}_{D_{[\bop]}}^j L) D_{[\bop]}^{k-j}u \|_{H_\etbop^{\sfs-m+1}}$, $j\geq 2$. Consider any one term in the sum decomposition of ${\rm ad}_{D_{[\bop]}}^j L$, which we can write as $(D_{[\bop]}^{j'}\ell)A$ where $\ell\in\cC_{\etbop;\bop}^{(d_0;k)}$, $A\in\CI_\bop\Psi_\etbop^m$, and $j'\leq j$. Focusing on the top order term, $j'=j$, and commuting $A$ through $D_{[\bop]}^{k-j}$, we then bound, using Proposition~\usref{PropMTameMicr} for $D_{[\bop]}^2\ell$, $A u$, $j-2$, $k-2$ in place of $\ell,u,j,k$:
  \begin{align*}
    \|G (D_{[\bop]}^{j-2} D_{[\bop]}^2 \ell) (D_{[\bop]}^{k-j}A u) \|_{H_\etbop^{\sfs-m+1}} &\leq C_k\biggl( \|\ell\|_{\cC_{\etbop;\bop}^{(d_0;\tilde d+2)}}\Bigl( \|G^\sharp u\|_{H_{\etbop;\bop}^{(\sfs+1;k-2)}} + \|\chi u\|_{H_{\etbop;\bop}^{(-N;k-2)}} \Bigr) \\
      &\qquad\qquad + \|\ell\|_{\cC_{\etbop;\bop}^{(d_0;k-2)}} \Bigl( \|G^\sharp u\|_{H_{\etbop;\bop}^{(\sfs+1;\tilde d)}} + \|\chi u\|_{H_{\etbop;\bop}^{(\sfs_0;\tilde d)}} \Bigr) \biggr).
  \end{align*}
  Observe then that
  \begin{equation}
  \label{EqMTamePrbToe3b}
    \|G^\sharp u\|_{H_{\etbop;\bop}^{(\sfs+1;k-2)}}\leq C\|G^\sharp u\|_{H_{\etbop;\bop}^{(\sfs;k-1)}}
  \end{equation}
  (which is a consequence of $\Vb(\tilde M)\supset\Vetb(\tilde M)$ in the form of Lemma~\ref{LemmaMUe3bExpb}).\footnote{Thus, it is the loss of $1$ e3b-derivative of the real principal type estimate that prompts us to treat terms in~\eqref{EqMTamePrExp} which involve at least $k-1$ many derivatives on $u$ separately from the remaining terms.} We may then, bound this term inductively (provided the original $G$ was chosen so as to have elliptic and operator wave front set contained in a small neighborhood of the smallest possible subset of ${}^\etbop S^*\tilde M$, given the stated requirements on $G$). This yields~\eqref{EqMTamePrEst}. The local uniformity statement is a consequence of the corresponding statement in Proposition~\ref{PropMPr}.
\end{proof}

\section{The Kerr metric and its phase space dynamics}
\label{STs}

In this section, we explain how the null-geodesic dynamics on a subextremal Kerr spacetime can be analyzed in e3b-phase space ${}^\etbop T^*M\to M$ (Definition~\ref{DefCTe3b}) on the compactification $M$ from Definition~\ref{DefCMSpacetime}. We recall the definition and basic properties of the Kerr metric in~\S\ref{SsTsK} before describing the dynamics (at spatial infinity, near the event horizon, trapping) in detail in~\S\S\ref{SsTs3b}--\ref{SsTse3b}; a basic (but insufficient for our purposes) reference for this part is \cite[Chapter~7, Section~61]{ChandrasekharBlackHoles}. While various aspects of the dynamics have previously been described in \cite{BlueSofferPhaseSpace,VasyMicroKerrdS,DyatlovWaveAsymptotics,MilletTeukolskyDecay,HintzGlueLocII}, we study the dynamics here in a novel (phase space) setting; for example, trapping and radial sets lie entirely over boundary hypersurfaces of $M$. For the sake of completeness and future reference, we thus present the relevant computations in detail. The main results of this section are Propositions~\ref{PropTs3bONHyp} and \ref{PropTse3bDyn}.

\subsection{The Kerr metric}
\label{SsTsK}

We fix the parameters of a \emph{subextremal Kerr black hole}, i.e., its mass $\bhm>0$ and specific angular momentum $a\in\R$ subject to the subextremality condition
\[
  |a| < \bhm.
\]
We recall that in Boyer--Lindquist coordinates \cite{KerrKerr,BoyerLindquistKerr}, the metric and dual metric are given by
\begin{equation}
\label{EqTsBL}
\begin{split}
  g_{\bhm,a} &= -\frac{\mu}{\varrho^2}(\dd\ft-a\sin^2\theta\,\dd\phi)^2 + \varrho^2\Bigl(\frac{\dd r^2}{\mu}+\dd\theta^2\Bigr) + \frac{\sin^2\theta}{\varrho^2}\bigl((r^2+a^2)\dd\phi-a\,\dd\ft\bigr)^2, \\
  \varrho^2 g_{\bhm,a}^{-1} &= -\frac{1}{\mu}\bigl((r^2+a^2)\pa_\ft+a\pa_\phi\bigr)^2 + \mu\pa_r^2 + \pa_\theta^2 + \frac{1}{\sin^2\theta}(a\sin^2\theta\,\pa_\ft+\pa_\phi)^2,
\end{split}
\end{equation}
where
\begin{equation}
\label{EqTsBLFunc}
  \mu(r) = r^2-2\bhm r+a^2,\quad
  \varrho^2(r,\theta) = r^2+a^2\cos^2\theta.
\end{equation}
This form of the metric is valid on the manifold
\begin{equation}
\label{EqTsBLMfd}
  \cM := \R_\ft \times (r_+,\infty)_r \times \Sph^2_{\theta,\phi},\quad r_+:=\bhm+\sqrt{\bhm^2-a^2}.
\end{equation}

\begin{lemma}[Timelike character of $\ft$]
\label{LemmaTsBLTimelike}
  For $r>r_+$, the differential $\dd\ft$ is (past) timelike.
\end{lemma}
\begin{proof}
  We have $\mu>0$ for $r>r_+$. We then compute
  \begin{align*}
    -\mu\varrho^2 g_{\bhm,a}^{-1}(\dd\ft,\dd\ft) &= (r^2+a^2)^2 - a^2\mu\sin^2\theta \\
      &= r^4 + a^2 r^2(1+\cos^2\theta) + a^4\cos^2\theta + 2\bhm a^2 r\sin^2\theta > 0.\qedhere
  \end{align*}
\end{proof}

In terms of the standard metric $\slg=\dd\theta^2+\sin^2\theta\,\dd\phi^2$ on $\Sph^2$, we can write
\begin{align*}
  g_{\bhm,a} &= -\Bigl(1-\frac{2\bhm r}{\varrho^2}\Bigr)\dd\ft^2 + \Bigl(1+\frac{2\bhm r-a^2\sin^2\theta}{\mu}\Bigr)\dd r^2 + \varrho^2\slg \\
    &\quad\hspace{6em} + \Bigl(1+\frac{2\bhm r}{\varrho^2}\Bigr)(a\sin^2\theta\,\dd\phi)^2 - \frac{4 a\bhm r\sin^2\theta}{\varrho^2}\dd\ft\,\dd\phi.
\end{align*}
Since $\sin^2\theta\,\dd\phi=r^{-2}(y\,\dd x-x\,\dd y)$ is a smooth 1-form on $\R^3\setminus\{0\}$, this shows that $g_{\bhm,a}$ is smooth on $\cM$, with $\theta,\phi$ being standard polar coordinates on $\Sph^2$. Furthermore, the metric coefficients, as expressed in terms of $\dd\ft$, $\dd r$, and $r\,\dd\theta$, $r\,\dd\phi$, are smooth $\ft$-independent functions of $r^{-1}$ and independent of $a$ modulo $r^{-2}\CI$.

We shall study the Kerr metric, extended across $r=r_+$, on (subsets of) the compactifications $M_0$ and $M$ defined in Definition~\ref{DefCMSpacetime}. To this end, we first need to relate points $(\ft,r,\theta,\phi)$ to points in $\R^4\subset M_0,M$. We begin by introducing a time function with ``hyperboloidal level sets,'' i.e., level sets that are spacelike and transversal to future null infinity. For this purpose, we only record that
\begin{equation}
\label{EqTsKLMetricKerr}
  g_{\bhm,a} \equiv g_{\bhm,0} = -\Bigl(1-\frac{2\bhm}{r}\Bigr)\,\dd\ft^2 + \Bigl(1-\frac{2\bhm}{r}\Bigr)^{-1}\,\dd r^2 + r^2\slg \bmod r^{-2}\CI.
\end{equation}
Consider then for $r>4\bhm$ the function
\[
  T_\gg(r) := \int_{4\bhm}^r \frac{1+\tilde T_\gg(r)}{1-\frac{2\bhm}{r'}}\,\dd r',
\]
where $\tilde T_\gg(r)$ is a smooth function, specified below, which equals $-1$ for $r\in[4\bhm,5\bhm]$. Setting
\[
  t_* := \ft - T_\gg(r),\quad r>4\bhm,
\]
we then compute
\begin{equation}
\label{EqTsKLMetric}
\begin{split}
  g_{\bhm,0} &= -\Bigl(1-\frac{2\bhm}{r}\Bigr)\,\dd t_*^2 -2(1+\tilde T_\gg)\,\dd t_*\,\dd r - \frac{2+\tilde T_\gg}{1-\frac{2\bhm}{r}}\tilde T_\gg\,\dd r^2 + r^2\slg, \\
  g_{\bhm,0}^{-1} &= \frac{2+\tilde T_\gg}{1-\frac{2\bhm}{r}}\tilde T_\gg\,\pa_{t_*}^2 + \Bigl(1-\frac{2\bhm}{r}\Bigr)\pa_r^2 + 2(1+\tilde T_\gg)\pa_{t_*}\otimes_s\pa_r + r^{-2}\slg^{-1}.
\end{split}
\end{equation}

\begin{lemma}[Spacelike level sets: $r\geq 4\bhm$]
\label{LemmaTsKLSpace}
  There exists a smooth function $\tilde T_\gg$ on $[4\bhm,\infty)$, equal to $-1$ on any fixed bounded subinterval of $[4\bhm,\infty)$ containing $[4\bhm,5\bhm]$, smooth in $r^{-1}$, and vanishing quadratically at $r^{-1}=0$, such that $\dd t_*$ is past timelike with respect to $g_{\bhm,a}$ for all $r\geq 4\bhm$.
\end{lemma}
\begin{proof}
  For $\tilde T_\gg=-1$, the conclusion holds by Lemma~\ref{LemmaTsBLTimelike}. Note that the coefficients of $\dd t_*=\dd\ft+(1-\frac{2\bhm}{r})^{-1}(1+\tilde T_\gg)\,\dd r$ are uniformly bounded when $\tilde T_\gg$ is. For $\tilde T_\gg$ taking values in $[-1,0]$, we can thus use~\eqref{EqTsKLMetric} to estimate, for some $C>0$,
  \[
    g_{\bhm,a}^{-1}(\dd t_*,\dd t_*) \leq g_{\bhm,0}^{-1}(\dd t_*,\dd t_*) + C r^{-2} \leq \tilde T_\gg + C r^{-2}.
  \]
  This is strictly negative if $\tilde T_\gg<-C r^{-2}$. The desired conclusion thus holds if we take $\tilde T_\gg=-1$ on an interval containing $[4\bhm,(C+1)^{-\frac12}+1]$, and transition to $\tilde T_\gg=-(C+1)r^{-2}$ in a smooth and monotonically increasing fashion.
\end{proof}

Turning to $r\leq 4\bhm$, we consider the functions
\begin{equation}
\label{EqTsKLCoordStar}
\begin{alignedat}{2}
  t_* &:= \ft - T_+(r), &\quad T_+(r) &:= \int_r^{4\bhm} \Bigl(\frac{r'{}^2+a^2}{\mu(r')}+\tilde T_+(r')\Bigr)\,\dd r', \\
  \phi_* &:= \phi - \Phi_+(r),&\quad \Phi_+(r) &:= \int_r^{4\bhm} \Bigl(\frac{a}{\mu(r')}+\tilde\Phi_+(r')\Bigr)\,\dd r',
\end{alignedat}
\end{equation}
where the smooth functions $\tilde T_+,\tilde\Phi_+$ on $[\bhm,4\bhm]$ will be chosen momentarily. The dual Kerr metric then takes the form
\begin{align*}
  \varrho^2 g_{\bhm,a}^{-1} &= \mu\pa_r^2 + 2\bigl((r^2+a^2+\mu\tilde T_+)\pa_{t_*}+(a+\mu\tilde\Phi_+)\pa_{\phi_*}\bigr)\otimes_s\pa_r \\
    &\qquad + \bigl((2(r^2+a^2)+\mu\tilde T_+)\pa_{t_*}+(2 a+\mu\tilde\Phi_+)\pa_{\phi_*}\bigr)\otimes_s(\tilde T_+\pa_{t_*}+\tilde\Phi_+\pa_{\phi_*}) \\
    &\qquad + \pa_\theta^2 + \frac{1}{\sin^2\theta}(a\sin^2\theta\,\pa_{t_*}+\pa_{\phi_*})^2,
\end{align*}
which now extends smoothly (and non-degenerately) across $r=r_+$ to $r\geq\bhm$.

\begin{lemma}[Spacelike level sets: $r\leq 4\bhm$]
\label{LemmaTsKLSpace2}
  There exist smooth functions $\tilde T_+,\tilde\Phi_+$ on $[\bhm,4\bhm]$, equal to $-\frac{r^2+a^2}{\mu(r)}$ and $-\frac{a}{\mu(r)}$, respectively, on any fixed compact subinterval of $(r_+,4\bhm]$ containing $[3\bhm,4\bhm]$, such that $\dd t_*$ is past timelike with respect to $g_{\bhm,a}$ for all $r\leq 4\bhm$ and satisfies $g_{\bhm,a}^{-1}(\dd t_*,\dd t_*)\in r^{-2}\CI(X)$. Moreover, $\dd r$ is future timelike for $\bhm\leq r<r_+$.
\end{lemma}
\begin{proof}
  The timelike nature of $\dd t_*$ is equivalent to
  \begin{equation}
  \label{EqTsKLSpace2}
    \varrho^2 g_{\bhm,a}^{-1}(\dd t_*,\dd t_*) =\tilde T_+ \bigl(2(r^2+a^2)+\mu\tilde T_+\bigr) + a^2\sin^2\theta < 0.
  \end{equation}
  In the subset of $\{r>r_+\}$ where $\tilde T_+=-\frac{r^2+a^2}{\mu}$, we indeed have $-\frac{(r^2+a^2)^2}{\mu}+a^2\sin^2\theta=\varrho^2 g_{\bhm,a}^{-1}(\dd\ft,\dd\ft)<0$ (cf.\ the proof of Lemma~\ref{LemmaTsBLTimelike}). For $r>r_+$, where $\mu>0$,~\eqref{EqTsKLSpace2} is equivalent to
  \begin{equation}
  \label{EqTsKLSpace2Intv}
    \Bigl|\tilde T_+ + \frac{r^2+a^2}{\mu}\Bigr| < \frac{1}{\mu}\sqrt{(r^2+a^2)^2 - a^2\mu\sin^2\theta}.
  \end{equation}
  For $r\leq r_+$, where $\mu\leq 0$, the negativity~\eqref{EqTsKLSpace2} holds provided we choose $\tilde T_+<0$ with $2(r^2+a^2)\tilde T_++a^2<0$, i.e.,
  \[
    \tilde T_+ < -\frac{a^2}{2(r^2+a^2)}.
  \]
  We may thus take $\tilde T_+(r)=-\frac{a^2}{2(r^2+a^2)}-1$ for $r\leq r_+$. By continuity,~\eqref{EqTsKLSpace2} then holds also for $r\in[r_+-\delta,r_++\delta]$ for some $\delta>0$. Thus, $\tilde T_+(r_+)$ satisfies~\eqref{EqTsKLSpace2Intv} for $r\in(r_+,r_++\delta]$, and one can then smoothly connect this with $-\frac{r^2+a^2}{\mu}$ for $3\bhm\leq r\leq 4\bhm$.

  Finally, the spacelike nature of $\dd r$ for $r<r_+$ is due to $\varrho^2 g_{\bhm,a}^{-1}(\dd r,\dd r)=\mu<0$ there. Note also that $\dd r$ and $\dd t_*$ are in opposite causal cones since
  \[
    \varrho^2 g_{\bhm,a}^{-1}(\dd r,\dd t_*) = 2(r^2+a^2+\mu\tilde T_+) > 0;
  \]
  here we use that $\mu,\tilde T_+\leq 0$ for $r\leq r_+$.
\end{proof}

\begin{figure}[!ht]
\centering
\includegraphics{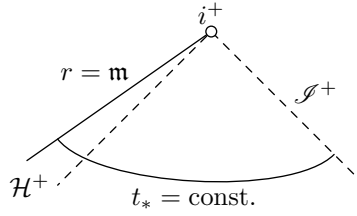}
\caption{Part of the Penrose diagram of the Kerr metric, including level sets of the time function $t_*$, future null infinity $\scri^+$, the future event horizon $\cH^+$, and future timelike infinity $i^+$ (which in the manifold $M$ is resolved to $\iota^+\cup\cK^+$).}
\label{FigTsKLt}
\end{figure}

See Figure~\ref{FigTsKLt}. We now identify $(t_*,r,\theta,\phi)$ with a point in $\R^4\subset M_0$ via $t=t_*+r$. In particular, $t=\ft$ for $r\in[3\bhm,5\bhm]$. We then regard $g_{\bhm,a}$ as a metric on $\R^4\cap\{r\geq\bhm\}$. We summarize and complement our constructions as follows:

\begin{lemma}[Kerr metric on $M$]
\label{LemmaTsKLMetric}
  We write $\ubar g=-\dd t^2+\dd x^2$ for the Minkowski metric, and define $g_{\bhm,a}$ as a metric on $\R_t\times[\bhm,\infty)\times\Sph^2\subset\R^4$ as above. Then $\dd t_*$ is everywhere past timelike, and $\dd r$ is future timelike in $\{\bhm\leq r<r_+\}$. Moreover:
  \begin{enumerate}
  \item\label{ItTsKLMetricSc}{\rm (Scattering metric.)} The Kerr metric $g_{\bhm,a}\in\CI(M;S^2\cT^*)$ is a non-degenerate Lorentzian signature metric, and $g_{\bhm,a}-\ubar g\in r^{-1}\CI(M;S^2\cT^*)$.
  \item\label{ItTsKLMetric3b}{\rm (Weighted 3b-metric.)} $g_{\bhm,a}\in\rho_\sface^{-2}\CI(M_0;S^2\,\Ttb^*M_0)$ is non-degenerate, and $\rho_\sface^2(g_{\bhm,a}-\ubar g)\in\rho_\sface\CI(M_0;S^2\,\Ttb^*M_0)$.
  \item\label{ItTsKLMetrice3b}{\rm (Weighted e3b-metric.)} $g_{\bhm,a}\in \rho_0^{-2}x_\sscri^{-2}\rho_+^{-2}\CI(M;S^2\,\Tetb^*M)$ is non-degenerate, and
    \begin{equation}
    \label{EqTsKLMetrice3b}
      \rho_0^2 x_\sscri^2\rho_+^2(g_{\bhm,a}-\ubar g) \in \rho_0 x_\sscri^2\rho_+\CI(M;S^2\,\Tetb^*M).
    \end{equation}
  \end{enumerate}
\end{lemma}
\begin{proof}
  The Kerr and Minkowski metrics are stationary. Their restrictions to $X=t_*^{-1}(0)$ agree modulo $r^{-1}\CI(X;S^2\cT^*_X)$. This implies part~\eqref{ItTsKLMetricSc} by Lemma~\ref{LemmaCPXDiffeo} and since the blow-down map $M\to M_1$ is smooth.

  Part~\eqref{ItTsKLMetric3b} follows from part~\eqref{ItTsKLMetricSc}, since $\CI(U;S^2\cT^*)=r^2\CI(U;S^2\,\Ttb^*M_0)$ (cf.\ \eqref{EqCT3bMink}), with $r^{-1}$ being a global defining function of $\sface$. Part~\eqref{ItTsKLMetrice3b} follows in a collar neighborhood $U\subset M$ of $\cK^+$ from part~\eqref{ItTsKLMetric3b}, since such $U$ can also be regarded as subsets of $M_0$. In any neighborhood $V$ of $\scri^+$ on the other hand, we combine the memberships~\eqref{EqCTebMink2} for $\ubar g,\ubar g^{-1}$ with two observations. First, Lemma~\ref{LemmaCTebsc} implies $\CI(M;S^2\cT^*)\subset\rho_0^{-2}x_\sscri^{-4}\rho_+^{-2}\CI(M;S^2\,\Tetb^*M)$, and thus
  \begin{equation}
  \label{EqTsKLMetrice3bDiff}
    g_{\bhm,a} - g_{\bhm,0} \in r^{-2}\CI(M;S^2\cT^*) \subset \CI(M;S^2\,\Tetb^*M).
  \end{equation}
  Second, the difference of~\eqref{EqTsKLMetric} and~\eqref{EqCTebMinkTstar} is
  \[
    g_{\bhm,0}-\ubar g = \frac{2\bhm}{r}\,\dd t_*^2 - 2\tilde T_\gg\,\dd t_*\,\dd r- \frac{2+\tilde T_\gg}{1-\frac{2\bhm}{r}}\tilde T_\gg\,\dd r^2.
  \]
  The second and third terms lie in $r^{-2}\CI(V;S^2\cT^*)$ when $V\subset M$ is a neighborhood of $\scri^+$ on which $r\geq 4\bhm$; they thus lie in $\CI(V;S^2\,\Tetb^*M)$ as in~\eqref{EqTsKLMetrice3bDiff}. The first term lies in $\rho_0 x_\sscri^2\rho_+\cdot \rho_0^{-2}\rho_+^{-2}\CI(M;S^2\,\Tetb^*M)=\rho_0^{-1}x_\sscri^2\rho_+^{-1}\CI$ by Lemma~\ref{LemmaCTebsc}. This gives~\eqref{EqTsKLMetrice3b}.
\end{proof}

We proceed to describe the null-bicharacteristic dynamics of $g_{\bhm,a}$ near the Kerr face in $\Ttb^*M$ (\S\ref{SsTs3b}) and then near null infinity in $\Teb^*M$ (\S\ref{SsTseb}). The combination in the e3b-setting is described in~\S\ref{SsTse3b}. We write
\begin{subequations}
\begin{equation}
\label{EqTsChar}
  \Sigma := \{ \zeta\in\Tetb^*M\setminus o \colon g_\etbop^{-1}(\zeta,\zeta)=0 \},\quad g_\etbop^{-1}:=\rho_0^{-2}x_\sscri^{-2}\rho_+^{-2}g_{\bhm,a}^{-1},
\end{equation}
for the collection of dual light cones (without the vertex), and write
\begin{equation}
\label{EqTsChar2}
  \Sigma = \Sigma^+ \sqcup \Sigma^-
\end{equation}
\end{subequations}
as the union of the future (``$+$'') and past (``$-$'') dual light cones; here $\zeta\in\Sigma^+\cap T^*\R^4$ if and only if $g_{\bhm,a}^{-1}(\zeta,\cdot)$ is a future timelike tangent vector. We shall refer to $\Sigma^{(+/-)}$ as the (future/past) characteristic set, as it will be (the conic extension of) the characteristic set of $\Box_{g_{\bhm,a}}$ as a weighted e3b-differential operator. Furthermore, we shall only work in $\Sigma^+$; the dynamics in $\Sigma^-=-\Sigma^+$ are the same up to reversing the affine parameter and an overall minus sign in the momentum variables. We will moreover write
\[
  \pa\Sigma \subset {}^\etbop S^*M
\]
for the boundary of $\Sigma$ at fiber infinity in $\ol{\Tetb^*}M$.

\subsection{Dynamics in 3b-phase space near \texorpdfstring{$\cK^+$}{the Kerr face}}
\label{SsTs3b}

We work in the region $M_+:=\ol{\{r\leq\frac23 t\}}\subset M$ and thus drop the weight and any reference of edge structures at $\scri^+$ from the notation; by Lemma~\ref{LemmaTsKLMetric}, we have $g_{\bhm,a}\in\CI(M_+;S^2\cT^*)$ and $g_{\bhm,a}\in\rho_+^{-2}\CI(M_+;S^2\,\Ttb^*_{M_+}M_0)$. We write
\begin{equation}
\label{EqTs3bGtb}
  G_\tbop(\zeta) := \varrho^2 g_{\bhm,a}^{-1}(\zeta,\zeta).
\end{equation}

\subsubsection{Dynamics near \texorpdfstring{$\cK^+\cap\iota^+$}{the boundary of the Kerr face}}
\label{SssTs3bC}

In this region, only the asymptotically Minkowskian nature of the Kerr metric matters; thus, we can argue as in \cite[\S{5}]{HintzNonstat}. Concretely, we work with the expression
\[
  \ubar g^{-1} = -2\pa_{t_*}\otimes_s\pa_r + \pa_r^2 + r^{-2}\slg,\quad t_*=t-r,
\]
for the Minkowski dual metric from~\eqref{EqCTebMinkTstar}, and recall from Lemma~\ref{LemmaTsKLMetric} that
\[
  g_{\bhm,a}^{-1}-\ubar g^{-1}\in r^{-1}\CI(M_+;S^2\cT)=r\CI(M_+;S^2\,\Ttb M_0).
\]
The rescaled dual metric function on $\Ttb^*M_0$ is
\begin{equation}
\label{EqTs3bCCoord}
  \ubar G_\tbop(\zeta) = r^2\ubar g^{-1}(\zeta,\zeta) = -2\sigma_\tbop\xi_\tbop + \xi_\tbop^2 + |\eta_\tbop|_{\slg^{-1}}^2,\quad \zeta = \sigma_\tbop\frac{\dd t_*}{r} + \xi_\tbop\frac{\dd r}{r} + \eta_\tbop,\ \eta_\tbop\in T^*\Sph^2,
\end{equation}
and $G_\tbop(\zeta)$, defined in~\eqref{EqTs3bGtb}, is equal to this modulo $r^{-1}P^{[2]}(\Ttb M_0)$ where we write $P^{[2]}$ for the space of fiberwise homogeneous quadratic polynomials. Therefore~\eqref{EqCHam3b} and~\eqref{EqCHam3bMap} give, in the coordinates $\rho_\cK=\frac{r}{t_*}$, $\rho_+=\frac{1}{r}$,
\begin{align*}
  H_{G_\tbop} &\equiv H_{\ubar G_\tbop} \bmod r^{-1}\Vtb(\Ttb^*M_0), \\
  H_{\ubar G_\tbop} &= 2(\xi_\tbop-\sigma_\tbop)(\rho_\cK\pa_{\rho_\cK}-\rho_+\pa_{\rho_+}+\sigma_\tbop\pa_{\sigma_\tbop}) + 2\sigma_\tbop\xi_\tbop\pa_{\xi_\tbop} + 2\xi_\tbop\rho_\cK^2\pa_{\rho_\cK} + H_{|\eta_\tbop|_{\slg^{-1}}^2}.
\end{align*}
Since $\ubar G_\tbop=(\xi_\tbop-\sigma_\tbop)^2+|\eta_\tbop|^2-\sigma_\tbop^2$, we have $\sigma_\tbop\neq 0$ on $\Sigma$ over $\iota^+$. Since $\dd t_*$ is past timelike, we moreover have $\sigma_\tbop<0$ on $\Sigma^+\cap{}^\tbop S^*_{\iota^+}M$. We thus pass to projective coordinates near $\pa\Sigma^+$,
\begin{equation}
\label{EqTs3bCCoordProj}
  \rho_\infty := -\frac{1}{\sigma_\tbop},\quad
  \hat\xi_\tbop := \frac{\xi_\tbop}{\sigma_\tbop},\quad
  \hat\eta_\tbop := \frac{\eta_\tbop}{\sigma_\tbop};
\end{equation}
then $\rho_\infty^2\ubar G_\tbop=\hat\xi_\tbop^2-2\hat\xi_\tbop+|\hat\eta_\tbop|^2$ and
\begin{equation}
\label{EqTs3bCHam}
\begin{split}
  \sfH_{G_\tbop} := \rho_\infty H_{G_\tbop} \equiv \rho_\infty H_{\ubar G_\tbop} &= 2(1-\hat\xi_\tbop)(\rho_\cK\pa_{\rho_\cK} - \rho_+\pa_{\rho_+} - \rho_\infty\pa_{\rho_\infty} - \hat\eta_\tbop\pa_{\hat\eta_\tbop}) \\
    &\qquad - 2(2-\hat\xi_\tbop)\hat\xi_\tbop\pa_{\hat\xi_\tbop} - 2\hat\xi_\tbop\rho_\cK^2\pa_{\rho_\cK} - H_{|\hat\eta_\tbop|^2}
\end{split}
\end{equation}
modulo the space of vector 3b-vector fields on $\ol{\Ttb^*}M_0\setminus o$ that are tangent to $\Stb^*M_0$ and vanish simply at $\iota^+$. (We abuse notation here and write $H_{|\hat\eta_\tbop|^2}$ for the Hamiltonian vector field on $T^*\Sph^2$ of the dual metric metric function on $\Sph^2$, evaluated at a point $\hat\eta_\tbop\in T^*\Sph^2$.) At $\iota^+=\rho_+^{-1}(0)$, we have $\sfH_{G_\tbop}=\sfH_{\ubar G_\tbop}$, and this vanishes on $\pa\Sigma^+$ if and only if $\hat\eta_\tbop=0$, $\hat\xi_\tbop\in\{0,2\}$, and $\rho_\cK=0$. We thus introduce:

\begin{definition}[Radial sets over $\pa\cK^+$]
\label{DefTs3bRad}
  In the coordinates~\eqref{EqTs3bCCoord} and~\eqref{EqTs3bCCoordProj} on $\ol{\Ttb^*}M$ near $\Stb^*_{\iota^+}M$, we define the subsets
  \begin{alignat*}{2}
    \cR^+_{\pa\cK^+,{\rm in}} &:= \{ (\rho_\cK,\rho_+,\omega;\sigma_\tbop,\xi_\tbop,\eta_\tbop) \colon \rho_\cK=\rho_+=0,\ \sigma_\tbop<0,\ \xi_\tbop=2\sigma_\tbop,&\ &\eta_\tbop=0 \}, \\
    \cR^+_{\pa\cK^+,{\rm out}} &:= \{ (\rho_\cK,\rho_+,\omega;\sigma_\tbop,\xi_\tbop,\eta_\tbop) \colon \rho_\cK=\rho_+=0,\ \sigma_\tbop<0,\ \xi_\tbop=0,&\ &\eta_\tbop=0 \}
  \end{alignat*}
  of $\Ttb^*_{\pa\cK^+}M\setminus o$, called \emph{outgoing} and \emph{incoming radial sets over $\pa\cK^+$}. We furthermore denote by
  \[
    \pa\cR^+_{\pa\cK^+,{\rm in}}=\{\hat\xi_\tbop=2,\ \hat\eta_\tbop=0\}\cap\Stb^*_{\pa\cK^+}M,\quad
    \pa\cR^+_{\pa\cK^+,{\rm out}}=\{\hat\xi_\tbop=0,\ \hat\eta_\tbop=0\}\cap\Stb^*_{\pa\cK^+}M
  \]
  their boundaries at fiber infinity. We write $(\pa)\cR^-_{\pa\cK^+,{\rm in/out}}$ for the sets obtained from these by fiber-wise multiplication by $-1$, and $(\pa)\cR_{\pa\cK^+,{\rm in/out}}$ for the union of the ``$+$'' and ``$-$'' sets.
\end{definition}

These sets are sources (`in') and sinks (`out') for the $\sfH_{G_\tbop}$-flow over $\cK^+$, and saddle points for over $\cK^+\cup\iota^+$. Since $\sfH_{G_\tbop}$ vanishes at $\pa\cR_{\pa\cK^+,{\rm in}}^+$, its linearization there, denoted $[\sfH_{G_\tbop}]$, is a bundle map on the conormal bundle of $\pa\cR_{\pa\cK^+,{\rm in}}^+$ inside of $\pa\Sigma^+$, mapping $\dd f$ to $\dd(\sfH_{G_\tbop}f)$ where $f$ vanishes on $\pa\cR_{\pa\cK^+,{\rm in}}^+$, similarly at the outgoing radial set (see \cite[Definition~10.8]{HintzMicro}). Since, locally, $\hat\xi_\tbop$ is a function of $\rho_\cK,\rho_+,\hat\eta_\tbop$, we compute
\begin{equation}
\label{EqTs3bLin}
\begin{alignedat}{5}
  &\text{at $\pa\cR_{\pa\cK^+,{\rm in}}^+$},\ &[\sfH_{G_\tbop}]\colon\ \dd\rho_\cK &\mapsto -2\,\dd\rho_\cK, &\ \dd\rho_+ &\mapsto 2\,\dd\rho_+, &\ \dd\hat\eta_\tbop &\mapsto 2\,\dd\hat\eta_\tbop; \\
  &\text{at $\pa\cR_{\pa\cK^+,{\rm out}}^+$},\ &[\sfH_{G_\tbop}]\colon\ \dd\rho_\cK &\mapsto 2\,\dd\rho_\cK, &\ \dd\rho_+ &\mapsto -2\,\dd\rho_+, &\ \dd\hat\eta_\tbop &\mapsto -2\,\dd\hat\eta_\tbop.
\end{alignedat}
\end{equation}

\begin{lemma}[Dynamics near $\pa\cK^+$]
\label{LemmaTs3bDyn}
  Recall $\rho_+=\frac{1}{r}$ and $\rho_\cK=\frac{r}{t_*}$. The $\sfH_{G_\tbop}$-flow in $\pa\Sigma^+\subset\Stb_{\cK^+\cup\iota^+\setminus\scri^+}^*M$ has the following properties.
  \begin{enumerate}
  \item\label{ItTs3bDynKp}{\rm (Convexity in $\cK^+$.)} There exists $\eps>0$ such that if $\sfH_{G_\tbop}\rho_+=0$, $0<\rho_+<\eps$, then $\sfH_{G_\tbop}^2\rho_+<0$.
  \item\label{ItTs3bDynGeod}{\rm (Null-bicharacteristics over $\cK^+$.)} Let $\gamma\colon I\to\pa\Sigma^+$, $I\subseteq\R$, be a maximally extended integral curve of $\sfH_{G_\tbop}$ lying over $\cK^+$ that is not contained in $\pa\cR^+_{\pa\cK^+,{\rm in}}\cup\pa\cR^+_{\pa\cK^+,{\rm out}}$ (otherwise it would be constant). Suppose that $\rho_+(\gamma(0))<\eps$ where $\eps>0$ is as in part~\eqref{ItTs3bDynKp}, and $\frac{\dd}{\dd s}\rho_+(\gamma(s))|_{s=0}\leq 0$, resp.\ $\geq 0$. Then $\gamma(s)$ converges to $\pa\cR^+_{\pa\cK^+,{\rm out}}$ as $s\to+\infty$, resp.\ $\pa\cR^+_{\pa\cK^+,{\rm in}}$ as $s\to-\infty$.
  \item\label{ItTs3bDynip}{\rm (Convexity in $\iota^+$.)} If $\sfH_{G_\tbop}\rho_\cK=0$, $\rho_\cK>0$, then $\sfH_{G_\tbop}^2\rho_\cK>0$.
  \item\label{ItTs3bStable}{\rm ((Un)stable manifolds over $\iota^+$.)} The unstable manifold of $\pa\cR_{\pa\cK^+,{\rm out}}^+$ and stable manifold of $\pa\cR_{\pa\cK^+,{\rm in}}^+$ over $\iota^+$ are
    \begin{equation}
    \label{EqTs3bStable}
    \begin{split}
      \pa\cW_{\rm out}^+ &= \{ \rho_+=0,\ \hat\xi_\tbop=0,\ \hat\eta_\tbop=0 \}, \\
      \pa\cW_{\rm in}^+ &= \{ \rho_+=0,\ \hat\xi_\tbop=2,\ \hat\eta_\tbop=0 \}.
    \end{split}
    \end{equation}
  \end{enumerate}
\end{lemma}

The conic extensions of~\eqref{EqTs3bStable} are
\begin{alignat*}{2}
  \cW_{\rm out}^+ &= \{ (\rho_\cK,\rho_+,\omega;\sigma_\tbop,\xi_\tbop,\eta_\tbop) \colon \rho_+=0,\ \sigma_\tbop<0,\ \xi_\tbop=0,&\ &\eta_\tbop=0 \}, \\
  \cW_{\rm in}^+ &= \{ (\rho_\cK,\rho_+,\omega;\sigma_\tbop,\xi_\tbop,\eta_\tbop) \colon \rho_+=0,\ \sigma_\tbop<0,\ \xi_\tbop=-2\sigma_\tbop,&\ &\eta_\tbop=0 \}.
\end{alignat*}
These are thus the characteristic sets of $r\pa_r$ and $r\pa_{t_*}-2 r\pa_r$, respectively, or, in $(t,r)$ coordinates, of the outgoing, resp.\ incoming null vector fields $r(\pa_t+\pa_r)$, resp.\ $r(\pa_t-\pa_r)$ (for the Minkowski metric).

\begin{proof}[Proof of Lemma~\usref{LemmaTs3bDyn}]
  For part~\eqref{ItTs3bDynKp}, let us write $\sfH_{\tilde G_\tbop}:=\sfH_{G_\tbop}-\sfH_{\ubar G_\tbop}$, then
  \[
    \sfH_{G_\tbop}\rho_+ = -2(1-\hat\xi_\tbop)\rho_+ + \sfH_{\tilde G_\tbop}\rho_+ = 0
  \]
  implies $2(1-\hat\xi_\tbop)=\cO(\rho_+)=\cO(\eps)$, so $\hat\xi_\tbop=1+\cO(\eps)$ and therefore also $|\hat\eta_\tbop|=1+\cO(\eps)$ on $\pa\Sigma^+$. Therefore,
  \begin{equation}
  \label{EqTs3bDynKp}
    \sfH_{G_\tbop}^2\rho_+ = 2\rho_+\sfH_{G_\tbop}\hat\xi_\tbop + \sfH_{G_\tbop}\sfH_{\tilde G_\tbop}\rho_+ = \rho_+\bigl(-2(1+\cO(\eps)) + \rho_+^{-1}\sfH_{G_\tbop}\sfH_{\tilde G_\tbop}\rho_+\bigr).
  \end{equation}
  The second summand in parentheses is bounded in absolute value by $C r^{-1}=\cO(\eps)$. Therefore,~\eqref{EqTs3bDynKp} is negative when $\eps>0$ is small enough.

  For part~\eqref{ItTs3bDynGeod} and assuming $\frac{\dd}{\dd s}\rho_+(\gamma(0))\leq 0$, note that part~\eqref{ItTs3bDynKp} implies that $\rho_+(\gamma(s))$ is strictly increasing for all $s\geq 0$. By~\eqref{EqTs3bLin}, $\gamma(s)$ must remain disjoint from a sufficiently small neighborhood of $\pa\cR_{\pa\cK^+,{\rm in}}^+$. But since over $\pa\cK^+$ we have
  \[
    \sfH_{\ubar G_\tbop}\hat\xi_\tbop = -2(2-\hat\xi_\tbop)\hat\xi_\tbop < 0
  \]
  when $\hat\xi_\tbop$ lies in a compact subset of $(0,2)$, we must have $\hat\xi_\tbop(\gamma(s))\to 0$ as $s\to\infty$. This forces $\gamma(s)$ to converge to $\pa\cR_{\pa\cK^+,{\rm out}}^+$. The arguments in the case $\frac{\dd}{\dd s}\rho_+(\gamma(0))\geq 0$ are completely analogous.

  We next turn to part~\eqref{ItTs3bDynip}. On $\iota^+$, we compute $\sfH_{G_\tbop}\rho_\cK=2(1-\hat\xi_\tbop)\rho_\cK-2\hat\xi_\tbop\rho_\cK^2$. If this vanishes, then $\sfH_{G_\tbop}^2\rho_\cK=-2\rho_\cK(1+\rho_\cK)\sfH_{G_\tbop}\hat\xi_\tbop=4\rho_\cK(1+\rho_\cK)(-\hat\xi_\tbop^2+2\hat\xi_\tbop)=4\rho_\cK(1+\rho_\cK)|\hat\eta_\tbop|^2$ on $\pa\Sigma^+$. This is strictly positive since otherwise, in view of $\rho_\cK>0$, we would have $\hat\eta_\tbop=0$ and thus $\hat\xi_\tbop=0$ or $2$; but then $\sfH_{G_\tbop}\rho_\cK=2\rho_\cK$ or $-2(\rho_\cK+2\rho_\cK^2)$ does not vanish.

  One easily checks that $\sfH_{G_\tbop}$ is tangent to the manifolds~\eqref{EqTs3bStable}; since they have the correct dimension (according to~\eqref{EqTs3bLin}), they are, indeed, the (un)stable manifolds.
\end{proof}

\subsubsection{Dynamics near \texorpdfstring{the event horizon $\cH^+$}{the event horizon}}
\label{SssTs3bH}

Near the event horizon $r=r_+$ of Kerr, it is, computationally, particularly convenient to work with the coordinates
\begin{equation}
\label{EqTs3bHt0phi0}
\begin{alignedat}{2}
  t_0 &:= \ft - T_0(r), &\quad T_0(r)&:=\int_r^{4\bhm}\frac{r'{}^2+a^2}{\mu(r')}\,\dd r', \\
  \phi_0 &:= \phi - \Phi_0(r), &\quad \Phi_0(r)&:=\int_r^{4\bhm}\frac{a}{\mu(r')}\,\dd r',
\end{alignedat}
\end{equation}
which are the same as~\eqref{EqTsKLCoordStar} with $\tilde T_+=\tilde\Phi_+=0$. Thus
\[
  \varrho^2 g_{\bhm,a}^{-1} = \mu\pa_r^2 + 2 \bigl((r^2+a^2)\pa_{t_0}+a\pa_{\phi_0})\otimes_s\pa_r + \pa_\theta^2 + \frac{1}{\sin^2\theta}(a\sin^2\theta\,\pa_{t_0}+\pa_{\phi_0})^2.
\]
Working near $(\cK^+)^\circ$, we introduce fiber-linear coordinates on $\Ttb^*M$ by writing 3b-covectors as
\begin{equation}
\label{EqTs3bHCoord}
  -\sigma_0\,\dd t_0 + \xi_0\,\dd r + \eta_\theta\,\dd\theta + \eta_{\phi_0}\,\dd\phi_0.
\end{equation}
Near the poles of $\Sph^2$, we need pass to smooth local coordinates, to wit,
\begin{equation}
\label{EqTs3bHCoordCart}
  \omega^1=\sin\theta\cos\phi_0,\ \omega^2=\sin\theta\sin\phi_0,\quad \eta_\theta\,\dd\theta + \eta_{\phi_0}\,\dd\phi_0 = \eta_1\,\dd\omega^1 + \eta_2\,\dd\omega^2;
\end{equation}
we record that
\[
  \eta_{\phi_0} = \omega^1\eta_2 - \omega^2\eta_1,\quad
  \slg^{-1} = \pa_{\omega^1}^2+\pa_{\omega^2}^2 - (\omega^1\pa_{\omega^1}+\omega^2\pa_{\omega^2})^2.
\]
We can thus write a rescaling $G_\tbop(\zeta):=\varrho^2 g_{\bhm,a}^{-1}(\zeta,\zeta)$ of the dual metric function as
\begin{equation}
\label{EqTs3bHCarter}
\begin{split}
  G_\tbop &= \mu\xi_0^2 + 2\bigl(-(r^2+a^2)\sigma_0+a\eta_{\phi_0}\bigr)\xi_0 + \sC, \\
  &\quad \sC := \eta_\theta^2 + \frac{1}{\sin^2\theta}(-a\sin^2\theta\,\sigma_0+\eta_{\phi_0})^2 = (\eta_1+a\sigma_0\omega^2)^2 + (\eta_2-a\sigma_0\omega^1)^2 - (\omega^1\eta_1+\omega^2\eta_2)^2.
\end{split}
\end{equation}
(The function $\sC\geq 0$ is \emph{Carter's constant} on the characteristic set.)

As a consequence of the future timelike nature of $\dd r$ for $r<r_+$, we have
\begin{equation}
\label{EqTs3bHMonor}
  H_{G_\tbop}r < 0\quad\text{on}\ \Sigma^+\cap\{r<r_+\}.
\end{equation}
Using the explicit expression
\begin{equation}
\label{EqTs3bHHam}
\begin{split}
  H_{G_\tbop} &= 2(r^2+a^2)\xi_0\pa_{t_0} + 2\bigl(\mu\xi_0 - (r^2+a^2)\sigma_0+a\eta_{\phi_0} \bigr)\pa_r + 2 a\xi_0\pa_{\phi_0} \\
    &\qquad - (\mu'\xi_0 - 4 r\sigma_0)\xi_0\pa_{\xi_0} + H_\sC,
\end{split}
\end{equation}
we can prove:

\begin{lemma}[Monotonicity of $r$]
\label{LemmaTs3bHr}
  Recall the event horizon $\cH^+=r^{-1}(r_+)$. We have $H_{G_\tbop}r<0$ on $\Sigma^+\cap\{r\leq r_+\}\setminus N^*\cH^+$.
\end{lemma}
\begin{proof}
  The inequality $H_{G_\tbop}r\leq 0$ on $\Sigma^+\cap\{r\leq r_+\}$ follows by continuity from~\eqref{EqTs3bHMonor}. At $r=r_+$ then, where $\mu=0$, we have $H_{G_\tbop}r=0$ if and only if $-(r^2+a^2)\sigma_0+a\eta_{\phi_0}=0$, so $G_\tbop=0$ implies $\sC=0$, and hence $\eta_\theta=0$ and $-a\sin^2\theta\,\sigma_0+\eta_{\phi_0}=0$ away from the poles of $\Sph^2$ and $\eta_1=\eta_2=0$ at the poles. Away from the poles, this yields a homogeneous linear system for $\sigma_0,\eta_{\phi_0}$ which admits only the trivial solution $\sigma_0=\eta_{\phi_0}=0$; but these equations together with $r_+=0$ define $N^*\cH^+$. 
\end{proof}

Lemma~\ref{LemmaTs3bHr} implies that future null-geodesics crossing $r=r_+$ in the inward direction pass the artificial interior spacelike boundary hypersurface $r=\bhm$ in finite affine time.

It also follows from~\eqref{EqTs3bHHam} that $N^*\cH^+$ is invariant under the $H_{G_\tbop}$-flow. But in terms of the local defining function $\rho_\cK:=t_0^{-1}$ of $\cK^+$, the quantity $H_{G_\tbop}\rho_\cK=-\rho_\cK^2 H_{G_\tbop}t_0=-2(r^2+a^2)\rho_\cK^2\xi_0$ is strictly negative for $\rho_\cK>0$ and $\xi_0>0$. The appropriate critical set for $H_{G_\tbop}$ thus lies over $\rho_\cK^{-1}(0)$, i.e., $\cK^+$:

\begin{definition}[Generalized radial set over the event horizon at $\cK^+$]
\label{DefTs3bH}
  Recall the coordinates~\eqref{EqTs3bHCoord} and write covectors as $-\sigma_0\,\dd t_0+\xi_0\,\dd r+\eta_0$, $\eta_0\in T^*\Sph^2$. Writing $\rho_\cK=\frac{1}{t_0}$, we then set
  \[
    \cR_{\cH^+}^+ := \{ (\rho_\cK,r,\omega;\sigma_0,\xi_0,\eta_0) \colon \rho_\cK=0,\ r=r_+,\ \sigma_0=0,\ \xi_0>0,\ \eta_0=0 \} \subset \Ttb^*_{(\cK^+)^\circ}M \setminus o.
  \]
  We moreover write $\pa\cR_{\cH^+}^+\subset\Stb^*_{(\cK^+)^\circ}M$ for the boundary of this set at fiber infinity.
\end{definition}

In order to work near $\pa\cR_{\cH^+}^+$, we introduce projective coordinates on the fibers,
\begin{equation}
\label{EqTs3bHCoordProj}
  \rho_\infty := \frac{1}{\xi_0},\ \hat\sigma_0 := \frac{\sigma_0}{\xi_0},\ \hat\eta_\theta := \frac{\eta_\theta}{\xi_0},\ \hat\eta_{\phi_0} := \frac{\eta_{\phi_0}}{\xi_0}.
\end{equation}
We then compute, using $\rho_\cK=\frac{1}{t_0}$,
\begin{equation}
\label{EqTs3bHExpr}
\begin{split}
  \sfH_{G_\tbop} := \rho_\infty H_{G_\tbop} &= -2(r^2+a^2)\rho_\cK^2\pa_{\rho_\cK} + 2(\mu-(r^2+a^2)\hat\sigma_0+a\hat\eta_{\phi_0})\pa_r \\
    &\qquad + 2 a\pa_{\phi_0} + (\mu'-4 r\hat\sigma_0)(\rho_\infty\pa_{\rho_\infty}+\hat\sigma_0\pa_{\hat\sigma_0} + \hat\eta_\theta\pa_{\hat\eta_\theta} + \hat\eta_{\phi_0}\pa_{\hat\eta_{\phi_0}}) + \rho_\infty H_\sC.
\end{split}
\end{equation}
We claim that $\pa\cR_{\cH^+}^+$ is a source for the $\sfH_{G_\tbop}$-flow inside of $\pa\Sigma^+\cap\Stb^*_{(\cK^+)^\circ}M$. To this end, we compute, modulo functions vanishing cubically at $r=r_+$, $\hat\sigma_0=0$, $\hat\eta=0$, and using $\mu'(r_+)=2(r_+-\bhm)$:
\begin{equation}
\label{EqTs3bHQuadDefComp}
  \sfH_{G_\tbop}\hat\sigma_0^2 \equiv 4(r_+-\bhm)\hat\sigma_0^2,\quad 
  \sfH_{G_\tbop}(\rho_\infty^2\sC) \equiv 4(r_+-\bhm)\rho_\infty^2\sC.
\end{equation}
This uses $\pa_r\sC=0$, which implies $\rho_\infty H_\sC(\rho_\infty^2\sC)=\rho_\infty^2\sC H_\sC(\rho_\infty)=0$, and the degree $2$ homogeneity of $\rho_\infty^2\sC=\hat\eta_\theta^2+\frac{1}{\sin^2\theta}(-a\sin^2\theta\,\hat\sigma_0+\hat\eta_{\phi_0})^2$ in $(\hat\sigma_0,\hat\eta_\theta,\hat\eta_{\phi_0})$.

\begin{lemma}[Quadratic defining function of $\pa\cR_{\cH^+}^+$ over $\cK^+$]
\label{LemmaTs3bHQuadDef}
  There exist a neighborhood $\cO$ of $\pa\cR_{\cH^+}^+$ inside of $\pa\Sigma^+\cap{}^\tbop S^*_{(\cK^+)^\circ}M$ such that the function
  \[
    \fq := \hat\sigma_0^2 + \rho_\infty^2\sC
  \]
  is smooth on $\cO$, a quadratic defining function of $\pa\cR_{\cH^+}^+$ inside of $\cO$, and satisfies $\sfH_{G_\tbop}\fq\geq 2(r_+-\bhm)\fq$ on $\cO$.
\end{lemma}
\begin{proof}
  Since $\pa_r G_\tbop=H_{G_\tbop}\xi_0=-\mu'(r_+)\xi_0^2=-2(r_+-\bhm)\xi_0^2$ at $N^*\cH^+$, so $\pa_r(\rho_\infty^2 G_\tbop)\neq 0$ at $\pa\cR_{\cH^+}^+$, one can locally express $r$ as a smooth function of $\hat\sigma_0$ and $(\theta,\phi_0,\hat\eta_\theta,\hat\eta_{\phi_0})\in T^*\Sph^2$; and of course $r=r_+$ at $\pa\cR_{\cH^+}^+$. (This observation is taken from \cite[\S{3}]{HintzHorizons}.) Since $\fq$ has a positive definite Hessian at $\hat\sigma_0=0$ and the zero section of $T^*\Sph^2$ (where $\rho_\infty^2\sC$ vanishes quadratically), this implies that $\fq$ is a quadratic defining function in the stated sense. The remaining claim follows immediately from~\eqref{EqTs3bHQuadDefComp}.
\end{proof}

As for the dynamics in the remaining coordinates $\rho_\cK,\rho_\infty$, we record:
\begin{equation}
\label{EqTs3bHBdfDer}
  \rho_\cK^{-2}\sfH_{G_\tbop}\rho_\cK = -2(r_+^2+a^2) < 0,\ \ 
  \rho_\infty^{-1}\sfH_{G_\tbop}\rho_\infty = 2(r_+-\bhm) > 0\quad\text{at}\ \pa\cR_{\cH^+}^+.
\end{equation}

\subsubsection{Dynamics in the domain of outer communications; trapping}
\label{SssTs3bO}

In the \emph{domain of outer communications}, $\{r_+<r<\infty\}$, we can work in Boyer--Lindquist coordinates used in~\eqref{EqTsBL}. We shall moreover work entirely in the translation-invariant setting; the conclusions for the null-bicharacteristic flow in e3b-phase space will be stated in~\S\ref{SssTs3bSum}.

We write (3b-)covectors as
\begin{equation}
\label{EqTs3bCoord}
  -\sigma\,\dd\ft + \xi\,\dd r + \eta_\theta\,\dd\theta + \eta_\phi\,\dd\phi,
\end{equation}
and introduce smooth local coordinates on $\Sph^2$
\[
  \omega^1=\sin\theta\cos\phi,\ \omega^2=\sin\theta\sin\phi,\quad
  \eta_\theta\,\dd\theta + \eta_\phi\,\dd\phi = \eta_1\,\dd\omega^1 + \eta_2\,\dd\omega^2,
\]
as in~\eqref{EqTs3bHCoordCart}. From~\eqref{EqTsBL}, we can thus read off the rescaled dual metric function $G_\tbop(\zeta):=\varrho^2 g_{\bhm,a}^{-1}(\zeta,\zeta)$ and its Hamiltonian vector field. We introduce
\begin{equation}
\label{EqTs3bOAB}
\begin{alignedat}{2}
  A &:= -(r^2+a^2)\sigma + a\eta_\phi &&=-(r^2+a^2)\sigma + a(\omega^1\eta_2-\omega^2\eta_1), \\
  B &:= -a\sin^2\theta\,\sigma + \eta_\phi &&= \omega^1\eta_2 - \omega^2\eta_1 - a(\omega_1^2+\omega_2^2)\sigma;
\end{alignedat}
\end{equation}
in particular, $\sigma=\varrho^{-2}(-A+a B)$. We moreover define
\begin{equation}
\label{EqTs3bOPsi}
  V := -\mu^{-1}A^2,\quad
  \Psi := -\mu^2\pa_r V = -A(4 r\mu \sigma + A\mu').
\end{equation}
We then have
\begin{equation}
\label{EqTs3bOMetHam}
\begin{split}
  G_\tbop &= \sG + \sC,\quad
  \sG := \mu\xi^2 + V = -\mu^{-1}A^2 + \mu\xi^2,\ \ \sC = \eta_\theta^2 + \frac{1}{\sin^2\theta}B^2, \\
  H_{G_\tbop} &= -\frac{2 A}{\mu}\bigl((r^2+a^2)\pa_\ft + a\pa_\phi\bigr) + 2\mu\xi\pa_r - \Bigl(\mu'\xi^2 - \frac{\Psi}{\mu^2}\Bigr)\pa_\xi + H_\sC.
\end{split}
\end{equation}

\begin{lemma}[Convexity of $r$ near the event horizon]
\label{LemmaTs3bOConvex}
  There exists $\eps>0$ such that the following holds: if $H_{G_\tbop}r=0$ at a point in $\Sigma$ where $r\in(r_+,r_++\eps)$, then $H_{G_\tbop}^2 r<0$.
\end{lemma}
\begin{proof}
  Since $H_{G_\tbop}r=2\mu\xi=0$, we have $\xi=0$. We need to show $H_{G_\tbop}^2 r=2\mu H_{G_\tbop}\xi=\frac{2\Psi}{\mu}<0$, i.e., $\Psi<0$. Note first that $A\neq 0$ since $A=0$ would imply $\sC=0$, and hence $\sigma=\eta_\theta=\eta_\phi=0$, so we would be at the zero section. More quantitatively, $\xi=0$ implies, on $\Sigma$,
  \begin{equation}
  \label{EqTs3bOConvexEst}
    \mu^{-1}A^2 = \sC \geq B^2.
  \end{equation}
  We can then estimate
  \[
    -A^{-2}\Psi = \mu' - \frac{4 r\mu(A-a B)}{\varrho^2 A} \geq \mu' - \frac{4 r\mu}{\varrho^2} - \frac{4 a r\sqrt\mu}{\varrho^2}.
  \]
  Now $\mu'(r_+)=2(r_+-\bhm)>0$, while the remaining two terms are of size $\cO(\sqrt\eps)$ when $r-r_+\in[0,\eps)$. For sufficiently small $\eps>0$, we then get $-A^{-2}\Psi>0$, as desired.
\end{proof}

The following analysis of the dynamics in $r>r_+$ closely follows Dyatlov \cite{DyatlovWaveAsymptotics}; see also Millet \cite{MilletTeukolskyDecay} for a slightly different account. (See also \cite[\S{6.4}]{VasyMicroKerrdS} and \cite{PetersenVasySubextremal} for the Kerr--de~Sitter case.)

\pfstep{Trapping, I: dynamics in $(r,\xi)$.} We note the conservation laws
\[
  H_{G_\tbop}\sigma=0,\quad H_{G_\tbop}\eta_\phi=0,\quad H_{G_\tbop}\sC=0;
\]
the first two correspond to the stationarity and axisymmetry of Kerr, and the third, first observed by Carter~\cite{CarterGlobalKerr}, follows from $H_{G_\tbop}\sC=H_{G_\tbop+\sC}\sC=H_\sG\sC$ and the independence of $\sC$ from $t,\phi,r,\xi$. This independence also yields
\[
  H_\sC r=0,\quad H_\sC\xi=0.
\]
Hamilton's equations for $r,\xi$ thus read $\dot r=\pa_\xi\sG$, $\dot\xi=-\pa_r\sG$, with $\sG$ being a function of $(r,\xi)$ only when $\sigma,\eta_\phi$ are fixed; i.e., they can be studied separately. Let us introduce the notation
\[
  (z,\zeta)\in T^*(\R_\ft\times\Sph^2),\quad z=(t,\theta,\phi),\ \zeta=(\sigma,\eta_\theta,\eta_\phi),
\]
for the remaining variables. We study the null-bicharacteristic starting at any (non-zero) null covector
\[
  \varpi^0 = (r^0,z^0;\xi^0,\zeta^0) \in \Sigma^+ \cap T^*(\R_\ft\times(r_+,\infty)\times\Sph^2)\setminus o,\quad r^0>r_+;
\]
Carter's constant $\sC$ (which only depends on $(z^0,\zeta^0)$) is constant along it, and $V$ in~\eqref{EqTs3bOPsi} becomes a function of $r$ only. The function $\sG$ in~\eqref{EqTs3bOMetHam} is then a function of $r,\xi$ only, with parametric dependence on $(z^0,\zeta^0)$. We denote this by
\begin{equation}
\label{EqTs3bOG0}
  \sG^0_{(z^0,\zeta^0)}(r,\xi) := \mu(r)\xi^2 + \Phi^0_{(z^0,\zeta^0)}(r),\quad
  \Phi^0_{(z^0,\zeta^0)}(r) := (V+\sC)|_{(z,\zeta)=(z^0,\zeta^0)}.
\end{equation}
Let $\pi\colon T^*(\R_\ft\times(r_+,\infty)\times\Sph^2)\to T^*(r_+,\infty)$ denote the projection onto the $(r,\xi)$-components. Denote by $\gamma=e^{s H_{G_\tbop}}\varpi^0$ the integral curve of $H_{G_\tbop}$ with $\gamma(0)=\varpi^0$. Then $\pi\circ\gamma$ is the integral curve of $H_{\sG^0_{(z^0,\zeta^0)}}$ with $(\pi\circ\gamma)(0)=(r^0,\xi^0)$. Since $\varpi^0\in\Sigma$, we have $\sG^0_{(z^0,\zeta^0)}(r^0,\xi^0)=0$ and thus
\[
  \sG^0_{(z^0,\zeta^0)}(r,\xi)=0\quad\text{along}\ \pi\circ\gamma.
\]
In particular, the ``potential function'' $\Phi^0_{(z^0,\zeta^0)}(r)$ is equal to $-\mu(r)\xi^2$ on $\pi\circ\gamma$, and thus $\Phi^0_{(z^0,\zeta^0)}(r)\leq 0$ along $\pi\circ\gamma$. Our first main task will be to elucidate the structure of the zero set of $\sG^0_{(z^0,\zeta^0)}\subset T^*(r_+,\infty)=(r_+,\infty)_r\times\R_\xi$.

\begin{lemma}[Condition for non-trapped $\pi\circ\gamma$]
\label{LemmaTs3bONonTrap}
  Suppose that $\Phi^0_{(z^0,\zeta^0)}(r)=0$ implies $\pa_r\Phi^0_{(z^0,\zeta^0)}(r)\neq 0$. Then $\Sigma_{(z^0,\zeta^0)}:=(\sG^0_{(z^0,\zeta^0)})^{-1}(0)$ is a smooth 1-dimensional submanifold of $T^*(r_+,\infty)$. In both directions along $\pi\circ\gamma$, we have either $r\to r_+$ or $r\to\infty$.
\end{lemma}
\begin{proof}
  The differential $\dd\sG^0_{(z^0,\zeta^0)}=(\mu'\xi^2+\pa_r\Phi^0_{(z^0,\zeta^0)})\,\dd r+2\mu\xi\,\dd\xi$ is nonzero when $\xi\neq 0$. When $\xi=0$ and $\sG^0_{(z^0,\zeta^0)}(r,\xi)=0$, then also $\Phi^0_{(z^0,\zeta^0)}(r)=0$, and therefore $\pa_r\Phi^0_{(z^0,\zeta^0)}(r)\neq 0$, so again $\dd\sG^0_{(z^0,\zeta^0)}\neq 0$. This shows that $\Sigma_{(z^0,\zeta^0)}$ is smooth and 1-dimensional, and $H_{\sG^0_{(z^0,\zeta^0)}}$ does not vanish on it. To finish the proof, we need to show that $\Sigma_{(z^0,\zeta^0)}$ has no compact connected components. (This suffices since the intersection of $\Sigma_{(z^0,\zeta^0)}$ with $T^*_{[a,b]}\R$, $r_+<a<b<\infty$, is compact.) Assuming the contrary, $\Phi^0_{(z^0,\zeta^0)}(r)$ attains a non-positive local minimum on $\Sigma_{(z^0,\zeta^0)}$. This, however, is ruled out by Lemma~\ref{LemmaTs3bOConv} below.
\end{proof}

\begin{lemma}[A convexity result]
\label{LemmaTs3bOConv}
  If $\Phi^0_{(z^0,\zeta^0)}(r)\leq 0$ and $\pa_r\Phi^0_{(z^0,\zeta^0)}(r)=0$, then $\pa_r^2\Phi^0_{(z^0,\zeta^0)}<0$.
\end{lemma}
\begin{proof}
  Since $0\geq\Phi^0_{(z^0,\zeta^0)}(r)=(-\mu^{-1}A^2+\sC)|_{(z,\zeta)=(z^0,\zeta^0)}$, we have $B^2\leq\sC\leq\mu^{-1}A^2$. This implies $A\neq 0$. Indeed, otherwise $A=\sC=B=0$, which implies $\zeta^0=0$, which together with the vanishing of $G_\tbop$ at $\varpi^0=(r^0,z^0;\xi^0,\zeta^0)\in\Sigma$ implies $\xi^0=0$, contradicting the fact that $\varpi^0\neq 0$. Now $\pa_r\Phi^0_{(z^0,\zeta^0)}=\pa_r V=-\mu^{-2}\Psi=0$ gives
  \begin{equation}
  \label{EqTs3bOConvSigma}
    \sigma = -\frac{A\mu'}{4 r\mu}.
  \end{equation}
  Therefore, using also $A'=-2 r\sigma$ and $\mu''=2$, the sign of $\pa_r^2\Phi^0_{(z^0,\zeta^0)}=-\mu^{-2}\pa_r\Psi$ is the opposite of that of
  \begin{equation}
  \label{EqTs3bOConvPsi}
  \begin{split}
    \pa_r\Psi &= -4(r\mu)'\sigma A - 4 r\mu\sigma A' - 2 A A'\mu' - A^2\mu'' \\
      &= \frac{A^2}{2 r\mu}\bigl( (2\mu+r\mu')\mu' - 4 r\mu \bigr) \\
      &= \frac{2 A^2}{r\mu}\bigl((r-\bhm)^3 + \bhm(\bhm^2-a^2)\bigr) > 0.
  \end{split}
  \end{equation}
  This completes the proof.
\end{proof}

\begin{figure}[!ht]
\centering
\includegraphics{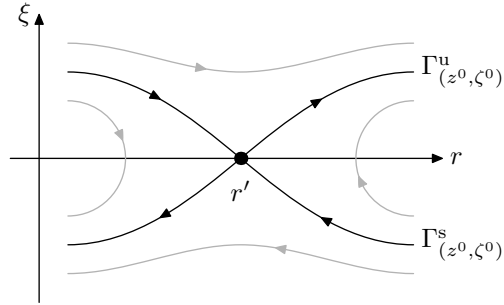}
\caption{Qualitative picture of the null-bicharacteristic dynamics on the unstable/stable trapped sets in $(r,\xi)$-phase space, defined by~\eqref{EqTs3bOGammapm0}. Shown in gray are null-bicharacteristics on $\sG^0_{(z^0,\zeta^0)}$ when $(z^0,\zeta^0)$ satisfies the conditions of Lemma~\usref{LemmaTs3bONonTrap}.}
\label{FigTs3bOTrap}
\end{figure}

Suppose on the other hand that there exists $r'\in(r_+,\infty)$ such that
\begin{equation}
\label{EqTs3bOPhiDouble}
  \Phi^0_{(z^0,\zeta^0)}(r')=\pa_r\Phi^0_{(z^0,\zeta^0)}(r')=0.
\end{equation}
By Lemma~\ref{LemmaTs3bOConv}, we have $\pa_r^2\Phi^0_{(z^0,\zeta^0)}<0$, and therefore $(r-r')\pa_r\Phi^0_{(z^0,\zeta^0)}(r) < 0$ for all $r\in(r_+,\infty)\setminus\{r'\}$ that are close to $r'$. This inequality in fact holds for all $r\in(r_+,\infty)\setminus\{r'\}$, for if it failed, there would be a critical point of $\Phi^0_{(z^0,\zeta^0)}$ closest to $r'$ which would need to be a negative local minimum or saddle point, contradicting Lemma~\ref{LemmaTs3bOConv}. This then implies
\begin{equation}
\label{EqTs3bOPhiSign}
  \Phi^0_{(z^0,\zeta^0)}(r)<0\quad \forall\,r\in(r_+,\infty),\ r\neq r'.
\end{equation}
We conclude that~\eqref{EqTs3bOPhiDouble} implies
\begin{equation}
\label{EqTs3bOGammapm0}
\begin{split}
  &\Sigma_{(z^0,\zeta^0)} = \Gamma_{(z^0,\zeta^0)}^{\rm u} \cup \Gamma_{(z^0,\zeta^0)}^{\rm s}, \\
  &\qquad \Gamma_{(z^0,\zeta^0)}^{\rm u/s} := \Bigl\{ (r,\xi) \colon r>r_+,\ \xi=(+/-)\sgn(r-r')\sqrt{-\Phi^0_{(z^0,\zeta^0)}(r)/\mu(r)}\,\Bigr\},
\end{split}
\end{equation}
where we take the \emph{positive} square root. We have $\Gamma_{(z^0,\zeta^0)}^{\rm u}\cap\Gamma_{(z^0,\zeta^0)}^{\rm s}=\{(r',0)\}$. Furthermore, since $\Phi^0_{(z^0,\zeta^0)}$ has a non-degenerate maximum at $r=r'$, the submanifolds $\Gamma_{(z^0,\zeta^0)}^{\rm u/s}$ are smooth (in fact, analytic). Note also that $H_{\sG^0_{(z^0,\zeta^0)}}|_{(r',0)}=0$, and hence $H_{G_\tbop}$-integral curves starting at a point $\varpi^0$ satisfying~\eqref{EqTs3bOPhiDouble} and $(r^0,\xi^0)=(r',0)$ are \emph{trapped}, i.e., they neither escape to $r=r_+$ nor to $r=\infty$. On $\Gamma_{(z^0,\zeta^0)}^{\rm u}\setminus\{(r',0)\}$, we have
\begin{equation}
\label{EqTs3bOGammapmFlow}
  (r-r')H_{\sG^0_{(z^0,\zeta^0)}}r = 2\mu\xi(r-r') = 2\mu|r-r'|\sqrt{-\Phi^0_{(z^0,\zeta^0)}/\mu} > 0,
\end{equation}
with the sign reversed on $\Gamma^{\rm s}_{(z^0,\zeta^0)}\setminus\{(r',0)\}$. Therefore, future, resp.\ past directed integral curves in $\Gamma^{\rm u}_{(z^0,\zeta^0)}$, resp.\ $\Gamma^{\rm s}_{(z^0,\zeta^0)}$ are either contained in $\{(r',0)\}$ and thus constant, or they escape to $r=r_+$ or $r=\infty$, while in the opposite direction they tend to $r=r'$. This explains the notation `{\it u}nstable' and `{\it s}table.' See Figure~\ref{FigTs3bOTrap}.

Finally, we determine the sign of $\sigma$ on trapped null-bicharacteristics in $\Sigma^+$. If $\gamma(s)=e^{s H_{G_\tbop}}\varpi^0$ satisfies $r(\gamma(s))=r'$ for some $s$, then $\xi(\gamma(s))=0$ since $G_\tbop=0$ along $\gamma$. (And thus $(r,\xi)(\gamma(s'))=(r',0)$ for all $s'$.) Therefore $\sigma^0=\sigma(\gamma(0))=\sigma(\gamma(s))\neq 0$ by~\eqref{EqTs3bOConvSigma} since $\mu'=2(r-\bhm)>0$ for $r\in(r_+,\infty)$. Furthermore, note that
\begin{equation}
\label{EqTs3bOTrapSign}
  \varrho^2 g_{\bhm,a}^{-1}(\dd\ft,\gamma(s)) = -\mu^{-1}(r^2+a^2)A + a B
\end{equation}
and $|a B|\leq |a|\mu^{-\frac12}|A|<\mu^{-1}(r^2+a^2)|A|$ where in the first inequality we use~\eqref{EqTs3bOConvexEst}, and the second inequality follows from $(r^2+a^2)^2-a^2\mu=r^4+a^2 r^2+2\bhm a^2 r>0$; therefore, $g_{\bhm,a}^{-1}(\dd\ft,\gamma(s))$ has the opposite sign of $A$, and thus by~\eqref{EqTs3bOConvSigma} the same sign as $\sigma^0$. Since $\dd\ft$ is past timelike, we conclude that $\gamma(s)\in\Sigma^{\sgn\sigma^0}$, and hence also $\gamma(0)\in\Sigma^{\sgn\sigma^0}$. Since we work in the \emph{future} characteristic set, we thus conclude that for $\varpi^0\in\Sigma^+$, an $r'\in(r_+,\infty)$ with~\eqref{EqTs3bOPhiDouble} can only possibly exist when\footnote{Note that the failure of $\pa_\ft$ to be timelike---or equivalently of $\ann(\pa_\ft)$ to be spacelike---in the \emph{ergoregion} of Kerr with $a\neq 0$ means that, on such spacetimes, there do exist points in $\Sigma^+\cap\{r>r_+\}$ with $\sigma<0$. Thus, the positivity of $\sigma$ on the trapped set is a manifestation of an observation in \cite{DafermosRodnianskiLinear,DafermosRodnianskiShlapentokhRothmanDecay} that ``superradiant frequencies are not trapped''.}
\begin{equation}
\label{EqTs3bOSigmaSign}
  \sigma^0>0.
\end{equation}

\pfstep{Trapping, II: structure of the trapped set.} We remain in the time-translation invariant setting, and introduce:

\begin{definition}[$\ft$-invariant trapped set]
\label{DefTs3bOTrap0}
  Recall $\cM$ from~\eqref{EqTsBLMfd}. The \emph{$\ft$-invariant trapped set} in $\Sigma^+$ is
  \[
    \Gamma_0 := \bigl\{ \varpi\in\Sigma^+\cap T^*\cM \colon r>r_+,\ \xi|_\varpi=H_{G_\tbop}\xi|_\varpi=0 \big\}.
  \]
\end{definition}

(This is motivated by the relationship of~\eqref{EqTs3bOG0}, \eqref{EqTs3bOPhiDouble}, and~\eqref{EqTs3bOPsi}: a point $(t,r,\theta,\phi;\sigma,\xi,\eta_\theta,\eta_\phi)$ lies in $\Gamma_0$ if and only if $\xi=0$ and $\Phi^0_{(z^0,\zeta^0)}$, with $(z^0;\zeta^0)=(t,\theta,\phi;\sigma,\eta_\theta,\eta_\phi)$, has a double zero at $r$.) Equivalently, $\Gamma_0$ can be defined in $\Sigma^+$ by the conditions $H_{G_\tbop}r=0$, $H_{G_\tbop}^2 r=0$, as follows from~\eqref{EqTs3bOMetHam}. By Lemma~\ref{LemmaTs3bDyn}\eqref{ItTs3bDynKp} and Lemma~\ref{LemmaTs3bOConvex}, there exist $r_{\Gamma,-}<r_{\Gamma,+}$ in $(r_+,\infty)$ such that
\begin{equation}
\label{EqTs3bOTrap0Cpt}
  \Gamma_0 \subset r^{-1}((r_{\Gamma,-},r_{\Gamma,+})).
\end{equation}

\begin{lemma}[Regularity of $\Gamma_0$]
\label{LemmaTs3bOTrap0Reg}
  $\Gamma_0$ is a smooth conic codimension $2$ submanifold of $\Sigma^+\cap T^*\cM$ (and thus of codimension $3$ within $T^*\cM$). Furthermore, $\sigma>0$ and $\sC>0$ on $\Gamma_0$.
\end{lemma}
\begin{proof}
  The conic set $\Gamma_0\subset T^*\cM\setminus o$ is defined by the equations $\xi=0$, $H_{G_\tbop}\xi=0$, $G_\tbop=0$. In view of~\eqref{EqTs3bOMetHam}, these are equivalent to $(\xi,\Psi,G_\tbop)=(0,0,0)$. We thus only need to show that at $\Gamma_0$, the differentials of $\xi$, $\Psi$, $G_\tbop$ are linearly independent. But at $\Gamma_0$,
  \begin{equation}
  \label{EqTs3bOTrap0RegDiff}
  \begin{alignedat}{3}
    \pa_\xi\xi&=1, &\qquad \pa_\xi\Psi&=0, &\qquad \pa_\xi G_\tbop&=0, \\
    & &\qquad \pa_r\Psi&\neq 0, &\qquad \pa_r G_\tbop&=0, \\
    & & & &\qquad \pa_\sigma G_\tbop &= -2\varrho^2 g_{\bhm,a}^{-1}(\dd\ft,\cdot)>0;
  \end{alignedat}
  \end{equation}
  here we used~\eqref{EqTs3bOConvSigma}--\eqref{EqTs3bOConvPsi} in the second line, and the discussion of~\eqref{EqTs3bOTrapSign} in the third line.

  The sign of $\sigma$ was already discussed in~\eqref{EqTs3bOSigmaSign}. The vanishing of $\sC$ (thus of $\eta_\theta$ and $B$) and $\xi$ would imply also that of $A$, which cannot hold outside the zero section, and hence $\sC>0$ on $\Gamma_0$.
\end{proof}

We proceed to relate $\Gamma_0$ to the sets~\eqref{EqTs3bOGammapm0}. Define the projection
\begin{equation}
\label{EqTs3bOTrapProj}
  \Pi \colon (r,z;\xi,\zeta) \mapsto (z,\zeta)
\end{equation}
along $T^*(r_+,\infty)$, and set
\begin{equation}
\label{EqTs3bOTrapwtGamma0}
  \tilde\Gamma_0 := \Pi(\Gamma_0) = \bigl\{ (z,\zeta)\in T^*(\R_\ft\times\Sph^2) \colon \exists\,r'>r_+\ \text{such that}\ (r',z;0,\zeta)\in\Gamma_0 \bigr\}.
\end{equation}
Thus, $\tilde\Gamma_0$ consists of all $(z^0,\zeta^0)$ for which there exists $r'\in(r_+,\infty)$ (which is then necessarily unique) satisfying~\eqref{EqTs3bOPhiDouble} (and thus any $(r,\xi)$ such that $(r,z^0;\xi,\zeta^0)\in\Sigma^+$ lies in the unstable or stable trapped set~\eqref{EqTs3bOGammapm0}).

\begin{definition}[$\ft$-invariant stable/unstable trapped sets]
\label{DefTs3bOTrapus0}
  We define the $\ft$-invariant sets
  \begin{align*}
    \Gamma^{\rm u/s}_0 &:= \bigcup_{(z,\zeta)\in\tilde\Gamma_0} \bigl\{ (r,z;\xi,\zeta)\colon (r,\xi)\in\Gamma_{(z,\zeta)}^{\rm u/s} \bigr\} \\
      &= \Bigl\{ (r,z;\xi,\zeta) \colon r>r_+,\ \xi=\pm\sgn(r-r'_{(z,\zeta)})\sqrt{-\Phi_{(z,\zeta)}^0(r)/\mu(r)}\,\Bigr\} \subset \Sigma^+\cap T^*\cM,
  \end{align*}
  where $\Phi^0_{(z,\zeta)}$ is defined by~\eqref{EqTs3bOG0} and $r'_{(z,\zeta)}$ is the unique zero of $\Phi^0_{(z,\zeta)}(r)$ in $r>r_+$; the top sign is for $\Gamma_0^{\rm u}$ and the bottom sign for $\Gamma_0^{\rm s}$.
\end{definition}

Directly from the definitions, we have
\[
  \Gamma_0 = \Gamma_0^{\rm u} \cap \Gamma_0^{\rm s}.
\]
Moreover, since, by construction, the sets $\Gamma_{(z,\zeta)}^{\rm u/s}\subset T^*(r_+,\infty)$ do not depend on the choice of $(z,\zeta)\in\Pi\circ\gamma$ for null-bicharacteristics $\gamma$, the sets $\Gamma_0^{\rm u/s}$ are invariant under the $H_{G_\tbop}$-flow.

\begin{lemma}[Regularity of $\tilde\Gamma_0$ and the (un)stable trapped sets]
\fakephantomsection
\label{LemmaTs3bOTrapwtGamma}
  \begin{enumerate}
  \item The set $\tilde\Gamma_0\subset T^*(\R_\ft\times\Sph^2)$ is a smooth conic codimension $1$ submanifold. The map $\tilde\Gamma_0\ni(z,\zeta)\mapsto r'_{(z,\zeta)}\in(r_+,\infty)$ is smooth and homogeneous of degree $0$ with respect to dilations in $\zeta$.
  \item\label{ItTs3bOTrapwtGamma2} The sets $\Gamma^{\rm u/s}_0$ are smooth conic codimension $1$ submanifolds of $\Sigma^+\cap T^*\cM$ which, as submanifolds of $\Sigma^+\cap T^*\cM$, intersect transversally at $\Gamma_0$.
  \end{enumerate}
\end{lemma}
\begin{proof}
  Let $(z,\zeta)\in\tilde\Gamma_0$ and let $\varpi=(r',z;0,\zeta)\in\Gamma_0$ be its unique preimage in $\Gamma_0$. We first claim that $T_\varpi\Gamma_0 \cap \ker D_\varpi\Pi = \{0\}$. To prove this, recall $T_\varpi\Gamma_0=\ker\dd\xi\cap\ker\dd\Psi\cap\ker\dd G_\tbop$ from the proof of Lemma~\ref{LemmaTs3bOTrap0Reg}; so given $v=a\pa_r+b\pa_\xi\in\ker D_\varpi\Pi$, we get $0=\dd\xi(v)=b$, so $b=0$, and then $\dd\Psi(v)=a\pa_r\Psi=0$ forces $a=0$ by~\eqref{EqTs3bOTrap0RegDiff}, as desired.

  Passing to quotients by fiber-dilations and time-translations, note now that $[\Gamma_0]\subset S^*\cM/\R_\ft$ is compact by~\eqref{EqTs3bOTrap0Cpt}, and we have just shown that $[\Pi]\colon[\Gamma_0]\to S^*(\R_\ft\times\Sph^2)/\R_\ft$ is a smooth immersion which, moreover, is injective due to the uniqueness of $r'$ in~\eqref{EqTs3bOTrapwtGamma0}; therefore, it is an embedding. This shows that $\tilde\Gamma_0$ is a smooth codimension $1$ submanifold. Furthermore, the smoothness of the inverse $(\Pi_{\Gamma_0})^{-1}\colon\tilde\Gamma_0\to\Gamma_0$ gives the smoothness of $r'_{(z,\zeta)}$ in $(z,\zeta)\in\tilde\Gamma_0$. The homogeneity in $\zeta$ is clear.

  Part~\eqref{ItTs3bOTrapwtGamma2} now follows from the definition of $\Gamma_0^{\rm u/s}$, with the transversality being a consequence of the fact that $\Phi^0_{(z,\zeta)}$ has a non-degenerate maximum at $r=r'_{(z,\zeta)}$ by Lemma~\ref{LemmaTs3bOConv}.
\end{proof}

\pfstep{Trapping, III: expansion rates.} Lemma~\ref{LemmaTs3bOTrapwtGamma} enables us to introduce smooth defining functions of $\Gamma_0^{\rm u/s}$:

\begin{definition}[Defining functions]
\label{DefTs3bODefFn}
  We write $\phi_0^{\rm u/s}\in S^1_{\rm hom}(T^*\cM\setminus o)$ for any smooth, stationary, and homogeneous (of degree $1$) extension of the smooth functions
  \begin{equation}
  \label{EqTs3bODefFn}
    \phi_0^{\rm u/s}(r,z;\xi,\zeta) := \xi \mp \sgn(r-r'_{(z,\zeta)}) \sqrt{-\Phi^0_{(z,\zeta)}(r)/\mu(r)},\quad (z,\zeta)\in\tilde\Gamma_0,
  \end{equation}
  defined on $\tilde\Gamma_0\times T^*(r_+,\infty)$; the top, resp.\ bottom sign is for $\phi_0^{\rm u}$, resp.\ $\phi_0^{\rm s}$.
\end{definition}

Note here that $\Phi^0_{(z,\zeta)}$ is homogeneous of degree $2$ since $V,\sC$ are (see~\eqref{EqTs3bOG0} and~\eqref{EqTs3bOPsi}--\eqref{EqTs3bOMetHam}), so~\eqref{EqTs3bODefFn} is homogeneous of degree $1$. Note also that $\Gamma_0^{\rm u/s}$ has codimension $1$ inside of $\tilde\Gamma_0\times T^*(r_+,\infty)$, and $\phi_0^{\rm u/s}$ is a defining function in $\tilde\Gamma_0\times T^*(r_+,\infty)$. Since $\Gamma_0^{\rm u/s}$ is invariant under the $H_{G_\tbop}$-flow, we have $H_{G_\tbop}\phi_0^{\rm u/s}=0$ at $\Gamma_0^{\rm u/s}$ regardless of the choice of extension. Since the codimension of $\Gamma_0^{\rm u/s}\subset\Sigma^+\cap T^*\cM$ is $1$, we can write
\begin{equation}
\label{EqTs3bOnu}
  (\sigma^{-1}H_{G_\tbop}\phi_0^{\rm u/s})|_\Sigma = \mp\nu^{\rm u/s}\phi_0^{\rm u/s}|_\Sigma
\end{equation}
where $\nu^{\rm u/s}$ are smooth functions on $\Sigma\cap T^*\cM$ in a neighborhood of $\Gamma_0^{\rm u/s}$, homogeneous of degree $0$.

\begin{lemma}[Expansion rates]
\label{LemmaTs3bOnu}
  The restrictions $\nu^{\rm u/s}|_{\Gamma_0}\in\CI(\Gamma_0)$ of the expansion rates defined by~\eqref{EqTs3bOnu} are independent of the choice of extensions of $\phi_0^{\rm u/s}$ in Definition~\usref{DefTs3bODefFn}. They are given by the formula
  \begin{equation}
  \label{EqTs3bOnuExpr}
    \nu^{\rm u/s}(r'_{(z,\zeta)},z;0,\zeta) = \nu(z,\zeta) := \sigma^{-1}\sqrt{-2\mu(r'_{(z,\zeta)})\pa_r^2\Phi^0_{(z,\zeta)}(r'_{(z,\zeta)})}.
  \end{equation}
  In particular,
  \begin{equation}
  \label{EqTs3bOnuMin}
    \nu_{\rm min} := \min\Bigl\{\inf_{\Gamma_0}\nu^{\rm u},\,\inf_{\Gamma_0}\nu^{\rm s}\Bigr\}>0.
  \end{equation}
  Moreover, $\{\phi_0^{\rm u},\phi_0^{\rm s}\}|_{\Gamma_0}\in\CI(\Gamma_0)$ is everywhere non-zero.
\end{lemma}
\begin{proof}
  We work inside of $\Sigma^+$. If $\tilde\phi^{\rm u/s}$ is another extension of~\eqref{EqTs3bODefFn}, then for $a^{\rm u/s}:=\tilde\phi^{\rm u/s}-\phi_0^{\rm u/s}$, we have $a^{\rm u/s}=0$ and $H_{G_\tbop}a^{\rm u/s}=0$ on $\Gamma_0^{\rm u}\cup\Gamma_0^{\rm s}$. Since $\Gamma_0^{\rm u/s}$ intersect transversally, we can write $a^{\rm u/s}=\sigma^{-1}\phi_0^{\rm u}\phi_0^{\rm s}a_\sharp^{\rm u/s}$ and $\sigma^{-1}H_{G_\tbop}a^{\rm u/s}=\phi_0^{\rm u}\phi_0^{\rm s}a_\flat^{\rm u/s}$ where $a^{\rm u/s}_{\sharp/\flat}$ are smooth and homogeneous of degree $0$. Therefore, $\tilde\phi^{\rm u/s}=\phi_0^{\rm u/s}(1+\phi_0^{\rm s/u}a_\flat^{\rm u/s})$ and hence
  \begin{align*}
    \sigma^{-1}H_{G_\tbop}\tilde\phi^{\rm u/s} &= \sigma^{-1}H_{G_\tbop}\phi_0^{\rm u/s} + \phi_0^{\rm u}\phi_0^{\rm s}a_\flat^{\rm u/s} \\
      &= \mp\nu^{\rm u/s}(\tilde\phi^{\rm u/s} - \sigma^{-1}\phi_0^{\rm u}\phi_0^{\rm s}a_\sharp^{\rm u/s}) + \phi_0^{\rm u}\phi_0^{\rm s}a_\flat^{\rm u/s} \\
      &= \mp(\nu^{\rm u/s}+\cO(\phi_0^{\rm s/u}))\tilde\phi^{\rm u/s},
  \end{align*}
  which shows that the expansion rates $\tilde\nu^{\rm u/s}$ for $\tilde\phi^{\rm u/s}$ are equal to $\nu^{\rm u/s}$ at $\Gamma_0$.

  When evaluating $\nu^{\rm u/s}$, it thus suffices to compute the limit of $(H_{G_\tbop}\phi_0^{\rm u/s})/(\mp\phi_0^{\rm u/s})$ at $\Gamma_0$ \emph{along $\Gamma_0^{\rm s/u}$}. But at $\Gamma_0^{\rm s/u}$, we have
  \[
    H_{G_\tbop}\phi_0^{\rm u/s} = H_{G_\tbop}(\phi_0^{\rm u/s}+\phi_0^{\rm s/u}) = H_{G_\tbop}(2\xi) = -2\mu'\xi^2 + 2\frac{\Psi}{\mu^2}
  \]
  by~\eqref{EqTs3bOMetHam}. Since $\xi$ vanishes quadratically at $\Gamma_0$, we have $\xi^2/(\mp\phi_0^{\rm u/s})\to 0$ at $\Gamma_0$. To evaluate the second term, we note that on $\Gamma_0^{\rm s/u}$, we have $\mp\phi_0^{\rm u/s}=\mp(\phi_0^{\rm u/s}+\phi_0^{\rm s/u})=\mp 2\xi=2\sgn(r-r'_{(z,\zeta)})\sqrt{-\Phi^0_{(z,\zeta)}(r)/\mu(r)}$ by Definition~\ref{DefTs3bOTrapus0}, and expanding $\Phi^0_{(z,\zeta)}$ in Taylor series around its local non-degenerate maximum at $r=r'_{(z,\zeta)}$ gives
  \begin{equation}
  \label{EqTs3bOnuSqrt}
    \sqrt{-\Phi^0_{(z,\zeta)}(r)/\mu(r)} = \sqrt{(-\pa_r^2\Phi^0_{(z,\zeta)}/(2\mu))|_{r=r'_{(z,\zeta)}}} (r-r'_{(z,\zeta)})\sgn(r-r'_{(z,\zeta)}) + \cO((r-r'_{(z,\zeta)})^2),
  \end{equation}
  so $\mp\phi_0^{\rm u/s}=\sqrt{(-2\pa_r^2\Phi^0_{(z,\zeta)}/\mu)|_{r=r'_{(z,\zeta)}}}(r-r'_{(z,\zeta)})+\cO((r-r'_{(z,\zeta)})^2)$. Since on $(z,\zeta)\times T^*(r_+,\infty)$ we have
  \[
    2\frac{\Psi}{\mu^2} = 2\mu^{-2}\pa_r\Psi(r'_{(z,\zeta)})(r-r'_{(z,\zeta)}) + \cO((r-r'_{(z,\zeta)})^2)
  \]
  in view of $\Psi(r'_{(z,\zeta)})=0$, and since at $r=r'_{(z,\zeta)}$ we have $2\mu^{-2}\pa_r\Psi=2\pa_r(\mu^{-2}\Psi)=-2\pa_r^2 V=-2\pa_r^2\Phi^0_{(z,\zeta)}$, and the formula~\eqref{EqTs3bOnuExpr} follows.

  For the final claim, we again use~\eqref{EqTs3bOnuSqrt} to conclude that, at $(r'_{(z,\zeta)},z;0,\zeta)$,
  \[
    \{\phi_0^{\rm u},\phi_0^{\rm s}\} = H_{\phi_0^{\rm u}}\phi_0^{\rm s} = H_{\phi_0^{\rm u}+\phi_0^{\rm s}}\phi_0^{\rm s} = 2 H_\xi\phi_0^{\rm s}=2\pa_r\phi_0^{\rm s}=2\sqrt{-(\pa_r^2\Phi^0_{(z,\zeta)}/2\mu)|_{r=r'_{(z,\zeta)}}},
  \]
  which is indeed nonzero.
\end{proof}

\pfstep{Trapping, IV: $\sfr$-normal hyperbolicity.} Lemma~\ref{LemmaTs3bOnu} provides all dynamical information necessary for analysis on exact Kerr (in particular, resolvent analysis which near trapping is based on \cite{DyatlovSpectralGaps}). The construction of the unstable trapped set of perturbations of Kerr in \cite{HintzPolyTrap} relies on the considerably stronger $\sfr$-normal hyperbolicity (for every $\sfr$) of the trapped set \cite{HirschPughShubInvariantManifolds}. This was proved for small $|\frac{a}{\bhm}|$ by Wunsch--Zworski \cite{WunschZworskiNormHypResolvent,WunschZworskiNormHypResolventCorrection} and in the full subextremal range by Dyatlov \cite{DyatlovWaveAsymptotics}, whom we follow here.

We consider the rescaling
\[
  \sfH := \frac{1}{H_{G_\tbop}\ft}H_{G_\tbop}
\]
of $H_{G_\tbop}$ near $\Gamma_0\subset\Sigma^+\cap T^*\cM$, which is thus future directed and satisfies $\sfH\ft=1$. Moreover, $\sfH$ is a positive, stationary, homogeneous (of degree $0$) multiple of $\sigma^{-1}H_{G_\tbop}$ since $H_{G_\tbop}\ft=2\varrho^2 g_{\bhm,a}^{-1}(\dd\ft,\cdot)$ is positive and homogeneous (of degree $1$). We may quotient out by $\ft$-translations and fiber-dilations: denote by
\[
  \pa\Sigma^+_{(0)} := \pa\Sigma^+ \cap S^*_{\ft^{-1}(0)}\cM,\quad
  \pa\Gamma_{(0)} := \pa\Gamma_0 \cap S^*_{\ft^{-1}(0)}\cM \subset \pa\Sigma^+_{(0)},
\]
the future characteristic set and the trapped set at fiber infinity $S^*\cM\subset\ol{T^*}\cM$ over $\{\ft=0\}$, respectively. The vector field $\sfH$ on $T^*\cM$ induces a vector field $\sfH_0$ on $\pa\Sigma_{(0)}^+$ (which agrees with $\sfH$ when acting on functions that are independent of $\ft$); and $\sfH_{(0)}$ is tangent to $\pa\Gamma_{(0)}$. The conceptual advantage of $\pa\Gamma_{(0)}$ is that it is a \emph{compact} invariant manifold for $\sfH_{(0)}$, with a compact neighborhood in $\pa\Sigma^+_{(0)}$.

\begin{prop}[$\sfr$-normal hyperbolicity]
\label{PropTs3bONHyp}
  The $\sfH_{(0)}$-flow within $\pa\Sigma^+_{(0)}$ is eventually absolutely $\sfr$-normally hyperbolic at $\pa\Gamma_{(0)}$ for every $\sfr$ in the sense that the time $1$ flow $e^{\sfH_{(0)}}$ satisfies the conditions of \cite[Definition~4]{HirschPughShubInvariantManifolds}. More precisely, there exists a smooth bundle splitting
  \begin{equation}
  \label{EqTs3bONHypSplit}
    T_{\pa\Gamma_{(0)}}\pa\Sigma^+_{(0)} = T\pa\Gamma_{(0)} \oplus \bar N^{\rm u} \oplus \bar N^{\rm s}
  \end{equation}
  with the following properties.
  \begin{enumerate}
  \item\label{ItTs3bONHypPres} For all $\varpi\in\pa\Gamma_{(0)}$ and $s\in\R$, the linearization $D_\varpi e^{s\sfH_{(0)}}$ preserves the splitting~\eqref{EqTs3bONHypSplit}.
  \item\label{ItTs3bONHypNorm} There exist constants $\mu,C>0$ such that
    \begin{equation}
    \label{EqTs3bONHypNorm}
      |D e^{s\sfH_{(0)}}(v)| \leq C e^{-\mu|s|}|v|,\quad v\in\bar N^{\rm u/s},\ \mp s\geq 0.
    \end{equation}
  \item\label{ItTs3bONHypTgt} For all $\eps>0$, there exists a constant $C_\eps>0$ such that
    \begin{equation}
    \label{EqTs3bONHypTgt}
      |D e^{s\sfH_{(0)}}(v)| \leq C_\eps e^{\eps|s|},\quad v\in T\pa\Gamma_{(0)},\ s\in\R.
    \end{equation}
  \end{enumerate}
  Here, $|\cdot|$ denotes any fixed fiber metric on the summands of~\eqref{EqTs3bONHypSplit}.
\end{prop}

We will need the following technical lemma:
\begin{lemma}[Rescalings and linearizations]
\label{LemmaTs3bOResc}
  Let $\cM$ be a manifold with bounded geometry. Let $V_1,V_2\in\cC^1_{\rm uni}(\cM;T\cM)$ be two uniform vector fields, and suppose that $V_2=e^\varphi V_1$ where $\varphi\in\cC^1_{\rm uni}(\cM)$. Let $\cT$ be another manifold with bounded geometry, and let $f\colon T\cM\to T\cT$ be a uniformly bounded bundle map. Denote by $|\cdot|$ uniform norms on the fibers of $T\cM$ and $T\cT$ (i.e., bounded from above and below by the Euclidean norm in the trivializations of $T\cM$ and $T\cT$ induced by the distinguished charts of $\cM$ and $\cT$, respectively). Suppose there exist constants $\lambda\in\R$, $C>0$ such that
  \[
    |f(D e^{s V_1}v)| \leq C e^{\lambda s}|v|,\quad s\geq 0,\ v\in T\cM.
  \]
  Then there exists a constant $C'=C'(C,\lambda,\|D\varphi\|_{\cC^0},f,V_1)$ such that
  \[
    |f(D e^{s V_2}v)| \leq C'(1+s)e^{\lambda s e^{\sup\varphi}}|v|,\quad s\geq 0,\ v\in T\cM.
  \]
\end{lemma}
\begin{proof}
  The flow $\alpha_j(s,p_0):=e^{s V_j}p_0$ is the solution of the ODE $\alpha_j(0,p_0)=p_0$, $\pa_s a_j(s,p_0)=V_j|_{\alpha_j(s,p_0)}$. The flow of $V_2=e^\varphi V_1$ is a reparameterizations of the flow of $V_1$: we have
  \[
    \alpha_2(s,p_0) = \alpha_1(\psi(s,p_0),p_0),
  \]
  provided the function $\psi$ satisfies
  \[
    \psi(0,p_0)=0,\quad \pa_s\psi(s,p_0)=\exp\Bigl( \varphi\bigl(\alpha_1(\psi(s,p_0),p_0)\bigr) \Bigr).
  \]
  In particular, $\psi(s,p_0)\leq s e^{\sup\varphi}$. We then compute the linearization
  \begin{equation}
  \label{EqTs3bOResPf}
  \begin{split}
    D_{p_0}e^{s V_2} = D_2|_{p_0}\alpha_2(s,\cdot) &= D_2|_{p_0}\psi(s,\cdot)(\pa_s\alpha_1)(\psi(s,p_0),p_0) + (D_2|_{p_0}\alpha_1)(\psi(s,p_0),p_0) \\
      &= V_1|_{\alpha_2(s,p_0)}\,D_2|_{p_0}\psi(s,\cdot) + (D_{p_0}e^{s' V_1})|_{s'=\psi(s,p_0)}.
  \end{split}
  \end{equation}
  The linearization $D_2|_{p_0}\psi(s,\cdot)\in\Hom(T_{p_0}\cM,\R)$ satisfies $D_2|_{p_0}\psi(0,\cdot)=0$ and
  \begin{align*}
    \pa_s\bigl(D_2|_{p_0}\psi(s,\cdot)\bigr) &= e^{\varphi(\alpha_2(s,p_0))} D_{\alpha_2(s,p_0)}\varphi \Bigl( D_2|_{p_0}\psi(s,\cdot)\,V_1|_{\alpha_2(s,p_0)} + (D_{p_0}e^{s'V_1})|_{s'=\psi(s,p_0)}\Bigr) \\
      &= (V_2|_{\alpha_2(s,p_0)}\varphi) D_2|_{p_0}\psi(s,\cdot) + e^{\varphi(\alpha_2(s,p_0))} D_{\alpha_2(s,p_0)}\varphi\circ (D_{p_0}e^{s' V_1})|_{s'=\psi(s,p_0)}.
  \end{align*}
  This implies
  \[
    \bigl| \pa_s\bigl(e^{-\varphi(\alpha_2(s,p_0))}D_2|_{p_0}\psi(s,\cdot)\bigr) \bigr| \leq \|D\varphi\|_{\cC^0} \bigl|(D_{p_0}e^{s'V_1})|_{s'=\psi(s,p_0)}\bigr| \leq C' e^{\lambda\psi(s,p_0)} \leq C' e^{\lambda s e^{\sup\varphi}},
  \]
  and therefore $|D_2|_{p_0}\psi(s,\cdot)|\leq s C' e^{\lambda s e^{\sup\varphi}}$. Plugging this into~\eqref{EqTs3bOResPf} gives
  \[
    |f(D_{p_0}e^{s V_2}v)| \leq |f(V_1|_{\alpha_2(s,p_0)})| s C' e^{\lambda s e^{\sup\varphi}}|v| + C e^{\lambda\psi(s,p_0)}|v|,
  \]
  which (for a new constant $C'$) is bounded by $(1+s)C' e^{\lambda s e^{\sup\varphi}}|v|$, as claimed.
\end{proof}

\begin{proof}[Proof of Proposition~\usref{PropTs3bONHyp}]
  It is convenient to project $\Gamma_{(0)}$ not to $S^*\cM$ (which one can regard as one cross section for the fiberwise dilation action) but, instead, to $\{\sigma=\sigma_0\}\subset T^*\cM\setminus o$ for any $\sigma_0>0$ (which is a cross section for the dilation action restricted to $\{\sigma>0\}\supset\Gamma_{(0)}$, cf.\ Lemma~\ref{LemmaTs3bOTrap0Reg}). There is then a one-to-one correspondence, via the differential of the projection $T^*\cM\setminus o\to S^*\cM$, of vector subbundles $\bar N^{\rm u/s}$ of $T_{\pa\Gamma_{(0)}}(S^*\cM)$ in~\eqref{EqTs3bONHypSplit} and vector subbundles $\bar N^{\rm u/s}_1\subset T_{\Gamma_0\cap\ft^{-1}(0)}(T^*\cM)\cap\ker\dd\sigma$ whose fibers are (separately for ${\rm u}$ and ${\rm s}$) identified under fiber-dilations in $T^*\cM$.

  We can then furthermore identify subbundles of $T_{\ft^{-1}(0)}\cM$ with stationary (i.e., invariant under $\ft$-translations) subbundles of $\ker\dd\ft\subset T\cM$, and thus $\bar N^{\rm u/s}_1$ with stationary subbundles $N^{\rm u/s}_0\subset\ker\dd\sigma\cap\ker\dd\ft\subset T_{\Gamma_0}(T^*\cM)$. Note that the $e^{s\sfH}$ acts by translations by $0$ and $s$ in $\sigma$ and $\ft$, respectively, and thus $D e^{s\sfH}$ preserves $\ker\dd\sigma$ and $\ker\dd\ft$. Since $\pa_\ft\in T\Gamma_0$, the original statements are therefore equivalent to analogous statements for
  \begin{equation}
  \label{EqTs3bONHypSplit0}
    T_{\Gamma_0}\Sigma^+ = T\Gamma_0 \oplus N^{\rm u}_0 \oplus N^{\rm s}_0,
  \end{equation}
  where $N^{\rm u/s}_0$ (which we need to find) and the fiber metric used in~\eqref{EqTs3bONHypNorm}--\eqref{EqTs3bONHypTgt} are $\ft$-independent and homogeneous (of degree $0$) with respect to fiber-dilations, and $\dd\sigma=0$, $\dd\ft=0$ on $N^{\rm u/s}_0$.

  Following \cite[\S{2.2}]{DyatlovWaveAsymptotics}, note now that since $\Gamma_0^{\rm u/s}\subset T^*\cM$ has codimension $2$, the symplectic orthocomplement $(T_{\Gamma_0}\Gamma_0^{\rm u/s})^\perp\subset T_{\Gamma_0}(T^*\cM)$ is a rank $2$ subbundle, and it contains $H_{G_\tbop}$ and $H_{\phi_0^{\rm u/s}}$. Since on $\Gamma_0$ we have $H_{G_\tbop}\phi_0^{\rm s/u}=0$ but $H_{\phi_0^{\rm u/s}}\phi_0^{\rm s/u}\neq 0$, we conclude that, in fact, $(T_{\Gamma_0}\Gamma_0^{\rm u/s})^\perp=\mathspan\{H_{G_\tbop},\,H_{\phi_0^{\rm u/s}}\}$. But also $H_{G_\tbop},H_{\phi_0^{\rm u/s}}$ are tangent to $\Gamma_0^{\rm u/s}$ at $\Gamma_0$, so we conclude that $T_{\Gamma_0}\Gamma_0^{\rm u/s}$ is coisotropic. Since, in view of $\sfH\ft=1$, the linear form $\dd\ft$ is injective on $(T_{\Gamma_0}\Gamma_0^{\rm u/s})^\perp$, we can thus define the line subbundle
  \[
    N_0^{\rm u/s} := (T_{\Gamma_0}\Gamma_0^{\rm u/s})^\perp \cap \ker\dd\ft \subset T_{\Gamma_0}\Gamma_0^{\rm u/s}.
  \]
  Explicitly, this is given by the span of $H_{\phi_0^{\rm u/s}}-(H_{\phi_0^{\rm u/s}}\ft)\sfH$, which annihilates $\sigma$ (i.e., is annihilated by $\dd\sigma$). Thus, $T\Gamma_0^{\rm u/s}=T\Gamma_0\oplus N_0^{\rm u/s}$, and hence we obtain a splitting~\eqref{EqTs3bONHypSplit0}.

  \pfstep{Part~\eqref{ItTs3bONHypPres}.} Note that $D e^{s\sfH}$ preserves $T\Gamma_0$ and $T\Gamma_0^{\rm u/s}$ as well as $\ker\dd\ft$. Since on $\Gamma_0$ the vector field $\sfH$ is equal to the \emph{Hamiltonian} vector field $H_{G_\tbop/(H_{G_\tbop}\ft)}$ whose flow is symplectic, also the differential of $e^{s\sfH}$ preserves the splitting~\eqref{EqTs3bONHypSplit0}.

  \pfstep{Part~\eqref{ItTs3bONHypNorm}.} Following \cite[\S{3.2}]{DyatlovWaveAsymptotics}, we define a (stationary) frame $n^{\rm u/s}$ of $N^{\rm u/s}_0$ via the normalization $\dd\phi_0^{\rm s/u}(n^{\rm u/s})=1$. Considering $D e^{s\sfH}(n^{\rm u/s})$ as an $s$-dependent family of sections of $N_0^{\rm u/s}$, we can write its Lie derivative along $\sfH$ as
  \[
    -\cL_\sfH n^{\rm u/s} = \pa_s\bigl(D e^{s\sfH}(n^{\rm u/s})\bigr)\big|_{s=0} =: \pm\tilde\nu^{\rm u/s}n^{\rm u/s}.
  \]
  Using $\sfH\phi_0^{\rm s/u}=\pm\frac{\sigma}{H_{G_\tbop}\ft}\nu^{\rm s/u}\phi_0^{\rm s/u}$ (recalling~\eqref{EqTs3bOnu}), we then compute
  \begin{align*}
    0 &= \sfH\bigl(\dd\phi_0^{\rm s/u}(n^{\rm u/s})\bigr) = \cL_\sfH(\dd\phi_0^{\rm s/u})(n^{\rm u/s}) + \dd\phi_0^{\rm s/u}\bigl(\cL_\sfH n^{\rm u/s}\bigr) = (\sfH\phi_0^{\rm s/u})(n^{\rm u/s}) \mp \tilde\nu^{\rm u/s} \\
      &= \pm\Bigl(\frac{\sigma}{H_{G_\tbop}\ft}\nu^{\rm s/u} - \tilde\nu^{\rm u/s}\Bigr).
  \end{align*}
  Lemma~\ref{LemmaTs3bOnu} thus gives
  \[
    \pa_s\bigl(D e^{s\sfH}(n^{\rm u/s})\bigr) = \pm(\tilde\nu^{\rm u/s}\circ e^{s\sfH})n^{\rm u/s},\quad \tilde\nu^{\rm u/s}=\frac{\sigma}{H_{G_\tbop}\ft}\nu^{\rm s/u} > 0.
  \]
  This implies part~\eqref{ItTs3bONHypNorm}.

  \pfstep{Part~\eqref{ItTs3bONHypTgt}.} This is the content of \cite[Proposition~2.1]{WunschZworskiNormHypResolvent} for $|\frac{a}{\bhm}|\ll 1$ and \cite[Proposition~3.6]{DyatlovWaveAsymptotics} in the full subextremal range. We give the full details here. Let $\eps>0$. By the compactness of $\Gamma_0\cap\{\ft=0,\,\sigma=1\}$, it suffices to show that for each integral curve $\gamma(s)=e^{s\sfH}\gamma(0)$, $\gamma(0)\in\Gamma_0\cap\{\ft=0,\,\sigma=1\}$, of $\sfH$ on $\Gamma_0$, there exists $T>0$ such that
  \begin{equation}
  \label{EqTs3bONHypTgtPf}
    |D_{\gamma(0)} e^{T\sfH}v| \leq e^{\eps T}|v|,\quad v\in T_{\gamma(0)}\Gamma_0.
  \end{equation}
  (Indeed, the same bound then holds, with $2\eps$ in place of $\eps$, for all nearby starting points.) Since $\sigma$ is conserved along the $\sfH$-flow, we have $\dd\sigma(D_{\gamma(0)} e^{T\sfH}v)=\dd\sigma(v)$. Upon subtracting from $v$ a suitable multiple of $\pa_\sigma$, it thus suffices to consider $v\in T_{\gamma(0)}(\Gamma_0\cap\{\sigma=1\})$. For $f=\dd r$, $\dd\xi$, regarded as fiber-linear maps from $T(\Gamma_0\cap\{\sigma=1\})$ to $\R$, we have $f(D_{\gamma(0)} e^{T\sfH}v)=0$ since $r,\xi$ are constant along trapped null-geodesics. Since $\ft$ increases at speed $1$, we moreover have $f(D_{\gamma(0)} e^{T\sfH}v)=f(v)$ for $f=\dd\ft$. It remains to bound the spherical part of $D_{\gamma(0)} e^{T\sfH}v$ when $v\in T(\Gamma_0\cap\{\sigma=1\})$. Lemma~\ref{LemmaTs3bOResc}, applied with $\Gamma_0\cap\{\sigma=1\}=\R_\ft\times(\Gamma_0\cap\{\ft=0,\,\sigma=1\})$ (which is of bounded geometry, with $\Gamma_0\cap\{\ft=0,\,\sigma=1\}$ being compact, and covering $\R_\ft$ by unit length intervals), the vector fields $\sfH$ and $H_{G_\tbop}$, and the projection map $f=\slpi\colon T(\Gamma_0\cap\{\sigma=1\})\subset T^*(\R_\ft\times(r_+,\infty)\times\Sph^2)\to T^*\Sph^2$, shows that it suffices to show that there exists $T>0$ such that
  \begin{equation}
  \label{EqTs3bONHypslpi}
    |\slpi(D_{\gamma(0)} e^{T H_{G_\tbop}}v)| \leq e^{\eps T}|v|,\quad v\in T_{\gamma(0)}(\Gamma_0\cap\{\sigma=1\}).
  \end{equation}

  Due to the time-translation-invariance of $G_\tbop$, it suffices to consider $v$ with vanishing $\pa_\ft$-com\-po\-nent, i.e., $v\in T_{\gamma(0)}(\Gamma_0\cap\{\ft=0,\,\sigma=1\})$. Note now that the projection $\Pi$ in~\eqref{EqTs3bOTrapProj} induces a diffeomorphism of $\Gamma_0\cap\{\ft=0,\,\sigma=1\}$ with $\tilde\Gamma_0\cap\{\ft=0,\,\sigma=1\}$, which in turn can be identified with a codimension $1$ submanifold $\Gamma_{\Sph^2}\subset T^*\Sph^2$. (For example, on Schwarzschild, this is the set $\{|\eta|_{\slg^{-1}}=3\sqrt{3}\bhm\}$, as can be read off from~\cite[(3.35)]{DyatlovWaveAsymptotics}.) The inverse of this identification maps
  \begin{equation}
  \label{EqTs3bONHypId}
    \Gamma_{\Sph^2} \ni (\omega,\eta) \mapsto (\ft,r,\omega;\sigma,\xi,\eta)=(0,r'_{(0,\omega;1,\eta)},\omega;1,0,\eta) \in \Gamma_0.
  \end{equation}
  Therefore, $v$ is of the form
  \[
    v=v_0+\ell(v_0)\pa_r,\quad v_0\in T\Gamma_{\Sph^2},
  \]
  where $\ell\colon T\Gamma_{\Sph^2}\to\R$ is the differential of the map $\Gamma_{\Sph^2}\ni(\omega,\eta)\mapsto r'_{(0,\omega;1,\eta)}$. With $r,\xi$ constant along trapped null-geodesics, we may, for the proof of~\eqref{EqTs3bONHypslpi}, consider only the $T(T^*\Sph^2)$-part of $H_{G_\tbop}$ (i.e., dropping $\pa_\ft,\pa_\sigma,\pa_r,\pa_\xi$), which is
  \[
    -\frac{2 A}{\mu}a\pa_\phi + H_\sC;
  \]
  for the purpose of~\eqref{EqTs3bONHypslpi}, we can regard this as a family
  \[
    H_{G_{r_0}},\quad G_{r_0} = a_{r_0}\eta_\phi + \sC|_{\sigma=1},\ a_{r_0}:=-\frac{2 A}{\mu}\Big|_{\sigma=1,\,r=r_0},
  \]
  of Hamiltonian vector fields on $T^*\Sph^2$; i.e., it suffices to show
  \[
    |D e^{T H_{G_r}}v_0| + |\ell(v_0)\pa_r e^{T H_{G_r}}| \leq e^{\eps T}|v_0|.
  \]
  Since the time $T$ flows of $H_{G_{r_j}}$, $j=1,2$, for $r_1,r_2\in(r_+,\infty)$ are related by a degree $T(a_{r_2}-a_{r_1})$ rotation in $\phi$, we have an even stronger bound, by $C T|v_0|$, for the second summand. For the first summand, it suffices, for the same reason, to consider not $H_{G_{r_0}}$ but
  \[
    H_{\sC_1},\quad \sC_1:=\sC|_{\sigma=1};
  \]
  we need to show the existence of $T>0$ such that
  \begin{equation}
  \label{EqTs3bONHypCarter}
    |D e^{T H_{\sC_1}}v_0| \leq e^{\eps T}|v_0|,\quad v_0\in T\Gamma_{\Sph^2}.
  \end{equation}

  At this point, we finally use Hamiltonian mechanics similarly to \cite{DyatlovWaveAsymptotics}. The quantities $\sC_1$ and $\eta_\phi$ are conserved for the $H_{\sC_1}$-flow, and they are in involution. Let us determine when their differentials are linearly independent. Away from the poles $\theta=0,\pi$, this happens if and only if $\pa_{\eta_\theta}\sC_1=2\eta_\theta=0$ and $\pa_\theta\sC_1=2\cot\theta(a^2\sin^2\theta-\frac{\eta_\phi^2}{\sin^2\theta})=0$, which for $\theta\neq\frac{\pi}{2}$ implies $\frac{\eta_\phi^2}{\sin^2\theta}=a^2\sin^2\theta$. Since we only care about $v_0$ with base points $(\omega;\eta)$ in $\Gamma_{\Sph^2}$, we in addition obtain (cf.\ \eqref{EqTs3bONHypId})
  \begin{equation}
  \label{EqTs3bONHyp0}
    0 = G_\tbop|_{(\ft,r,\omega;\sigma,\xi,\eta)=(0,r'_{(0,\omega;1,\eta)},\omega;1,0,\eta)} = -\mu^{-1}A^2 + \sC_1 = -\mu^{-1}A^2 - 2 a B.
  \end{equation}
  We homogenize this term by inserting $\sigma=1=\varrho^{-2}(-A+a B)$; and then
  \[
    \mu^{-1}A^2 + 2 a B = \frac{A^2}{\mu}\Bigl(1 - \frac{2 a\mu B}{A\varrho^2}\Bigr) + \frac{2 a^2 B^2}{\varrho^2}.
  \]
  Using~\eqref{EqTs3bOConvexEst}, the term in parentheses is bounded from below by $1-\frac{2 a\sqrt\mu}{\varrho^2}$; this is positive, as follows from $\varrho^2\geq r^2$ and $r^4\geq 4 a^2\mu$, the second inequality being equivalent to $(r^2-2 a^2)^2+8 a^2(\bhm r-a^2)>0$, which holds in view of $|a|<\bhm<r$. Therefore~\eqref{EqTs3bONHyp0} would imply $A=B=0$, which (as argued before~\eqref{EqTs3bOConvexEst}) does not happen on $\Gamma_0$. At the poles on the other hand, we work with the expression~\eqref{EqTs3bHCarter} for $\sC$ (with $\sigma_0$ there set to $1$); thus at the pole $(\omega^1,\omega^2)=(0,0)$, we have $\dd\eta_\phi=\eta_2\,\dd\omega^1-\eta_1\,\dd\omega^2$ and $\dd\sC_1=2\eta\,\dd\eta-2 a\,\dd\eta_\phi$, which are linearly independent since $\eta=(\eta_1,\eta_2)\neq 0$ on $\Gamma_{\Sph^2}$. (The latter follows from~\eqref{EqTs3bONHypId} and the fact that $\dd\ft$ is not null.) Thus, $\dd\sC_1$ and $\dd\eta_\phi$ are everywhere linearly dependent on $\Gamma_{\Sph^2}$ except at the conormal bundle of the equator,
  \[
    E := \Gamma_{\Sph^2} \cap \Bigl\{ (\theta,\phi;\eta_\theta,\eta_\phi) \colon \theta=\frac{\pi}{2},\ \eta_\theta=0 \Bigr\}.
  \]

  Consider now an integral curve $\gamma$ of $H_{\sC_1}$ which is not contained in, and hence is disjoint from, $E$. By the Arnold--Liouville theorem, there exists a local symplectomorphism near $\gamma$ from $T^*\Sph^2$ to the cotangent bundle $T^*\TT^2$ of the 2-torus $\TT^2=\R^2_{y^1,y^2}/\Z^2$ which conjugates $H_{\sC_1}$ to $H_f$, $f=f(\eta_1,\eta_2)$, where we now write $\eta_1,\eta_2$ for the fiber-linear coordinates on the fibers of $T^*\TT^2$. Thus, $H_f=(\pa_\eta f)\pa_y$. The linearization $D e^{s H_f}$ evolves a tangent vector $v_y\pa_y+v_\eta\pa_\eta$ at $(y_0,\eta_0)$ according to $\pa_s y(s)=\pa_\eta f(\eta(s))$, $\pa_s\eta(s)=0$ in the base, so $\pa_s y(s)=\pa_\eta f(\eta_0)$, $\eta(s)=\eta_0$, and thus $\pa_s v_\eta(s)=0$, i.e., $v_\eta(s)=v_\eta$, and $\pa_s v_y(s)=\nabla^2 f(\eta_0)v_\eta$ in the fibers. Therefore, $(v_y(s),v_\eta(s))$ grows (at most) linearly in $s$, proving~\eqref{EqTs3bONHypCarter} in this case.

  Suppose on the other hand that $\gamma\subset E$, and consider the evolution $v(s):=D e^{s\sC_1}v_0=v_\theta(s)\pa_\theta+v_\phi(s)\pa_\phi+v_{\eta_\theta}(s)\pa_{\eta_\theta}+v_{\eta_\phi}(s)\pa_{\eta_\phi}$ of a tangent vector $v(0)=v_0\in T_E\Gamma_{\Sph^2}$. Consider first the flow in the base: it is given by
  \begin{equation}
  \label{EqTs3bONHypFlow}
  \begin{alignedat}{2}
    \pa_s\phi(s) &= 2\Bigl(\frac{\eta_\phi(s)}{\sin^2\theta(s)}-a\Bigr), &\quad \pa_s\eta_\phi(s) &= 0, \\
    \pa_s\theta(s) &= 2\eta_\theta(s), &\quad \pa_s\eta_\theta(s) &= -\pa_\theta\sC_1.
  \end{alignedat}
  \end{equation}
  Therefore, $\theta(s)=\frac{\pi}{2}$, $\phi(s)=2 s(\eta_\phi(0)-a)$, $\eta_\theta(s)=0$, $\eta_\phi(s)=\eta_\phi(0)$. Linearizing the system~\eqref{EqTs3bONHypFlow} around this solution gives $\pa_s v_{\eta_\phi}(s)=0$, so $v_{\eta_\phi}(s)=v_{\eta_\phi}(0)$, and $\pa_s v_\phi(s)=2 v_{\eta_\phi}(s)=2 v_{\eta_\phi}(0)$; these two components grow at most linearly in $s$. To compute the evolution equations for $v_\theta,v_{\eta_\theta}$, we use that $\pa_\theta\sC_1$ is independent of $\phi,\eta_\theta$ to find
  \begin{align*}
    \pa_s v_\theta(s) &= 2 v_{\eta_\theta}(s), \\
    \pa_s v_{\eta_\theta}(s) &= -\pa_\theta^2\sC_1(\gamma(s))v_\theta(s) - \pa_{\eta_\phi}\pa_\theta\sC_1(\gamma(s))v_{\eta_\phi}(s) \\
      &= -\pa_\theta^2\sC_1(\gamma(s))v_\theta(s) + 4\eta_\phi(0)v_{\eta_\phi}(0),
  \end{align*}
  and therefore
  \[
    \pa_s^2 v_\theta(s) = -2\pa_\theta^2\sC_1(\gamma(s))v_\theta(s) + 8\eta_\phi(0)v_{\eta_\phi}(0).
  \]
  To finish the proof, it remains to show that along $\gamma(0)$, the (constant) quantity $\pa_\theta^2\sC_1=2(\eta_\phi^2-a^2)$ is positive. To see this, recall from~\eqref{EqTs3bOConvexEst} that $A^2=\mu B^2\neq 0$ at $E$. Therefore, using $\mu<r^2+a^2$, we have $A^2<(r^2+a^2)B^2$, i.e.,
  \[
    (r^2+a^2)^2 - 2 a(r^2+a^2)\eta_\phi + a^2\eta_\phi^2 < (r^2+a^2)(\eta_\phi^2 - 2 a\eta_\phi + a^2),
  \]
  which is equivalent to $\eta_\phi^2>r^2+a^2>\bhm^2+a^2$. Thus $\pa_\theta^2\sC_1>2\bhm^2>0$ indeed. This proves~\eqref{EqTs3bONHypCarter} also for equatorial null-geodesics, and thus finishes the proof of part~\eqref{ItTs3bONHypTgt}.
\end{proof}

\subsubsection{Summary of the dynamics at the Kerr face}
\label{SssTs3bSum}

Recall the radial sets $\cR^+_{\pa\cK^+,{\rm in/out}}$ from Definition~\ref{DefTs3bRad}, $\cR^+_{\cH^+}$ from Definition~\ref{DefTs3bH}, and define\footnote{We use the notation $\bar\Gamma^{\rm u/s}$ for consistency with the notation used in \cite{HintzPolyTrap}.}
\begin{equation}
\label{EqTs3bGamma}
  \Gamma \subset \Ttb^*_{\cK^+}M_0\setminus o,\quad \bar\Gamma^{\rm u/s} \subset \Ttb^*_{\cK^+}M_0\setminus o
\end{equation}
to be the intersections of the closures of $\Gamma_0$ and $\Gamma_0^{\rm u/s}$ (see Definitions~\ref{DefTs3bOTrap0} and \ref{DefTs3bOTrapus0}) inside of $\Ttb^*M\setminus o$ with $\Ttb^*_{\cK^+}M_0\setminus o$.\footnote{Equivalently, $\Gamma=\{\varpi\in\Sigma^+\cap\Ttb^*_{\cK^+}M\setminus o\colon r>r_+,\ \xi|_\varpi=H_{G_\tbop}\xi|_\varpi=0\}$, and $\bar\Gamma^{\rm u/s}$ are defined as in Definition~\ref{DefTs3bOTrapus0} except $z$ there is replaced by points $(\rho_\cK,\omega)$ with $\rho_\cK:=\ft^{-1}$ and $\omega\in\Sph^2$. Another equivalent definition is to write $\Gamma_0=\R_\ft\times\Gamma$ where $\Gamma\subset T^*_{\ft^{-1}(0)}\cM$, and identify $T^*_{\ft^{-1}(0)}\cM$ with $\Ttb^*_{\cK^+\cap\{r_+<r<\infty\}}M$ and thus $\Gamma$ with a subset of the latter.} All of these sets are invariant under the $H_{G_\tbop}$-flow, and they lie over $\cK^+$ (consistent with the idea that, e.g., trapping does not occur on any finite time interval). Writing 3b-covectors near $(\cK^+)^\circ$ as in~\eqref{EqTs3bCoord}, we also recall from~\eqref{EqTs3bOSigmaSign} that
\begin{equation}
\label{EqTs3bSumEll}
  \sigma\ \text{is a symbol (of order $1$ on $\Tb^*M_0$) which is positive (and thus elliptic) at $\Gamma$}.
\end{equation}

We denote by $\pa\cR^+_{\pa\cK^+,{\rm in/out}}$ etc.\ the boundaries of the sets $\cR^+_{\pa\cK^+,{\rm in/out}}$ etc.\ at fiber infinity $\Stb^*_{\cK^+}M$.

\begin{prop}[Dynamics over the Kerr face]
\label{PropTs3bSum}
  Denote by $\rho_\infty\in\CI(\ol{\Tetb^*}M)$ a defining function of ${}^\etbop S^*M$. Let $\gamma\colon I\to\pa\Sigma^+\cap\Stb^*_{\cK^+}M$ denote a maximally extended integral curve of $\sfH:=\rho_\infty H_{G_\tbop}$, with $G_\tbop$ defined in~\eqref{EqTs3bGtb}. Assume that $\gamma$ is not contained in any of the $\sfH$-invariant sets $\pa\cR_{\pa\cK^+,{\rm out}}$, $\pa\cR_{\pa\cK^+,{\rm in}}$, $\pa\cR^+_{\cH^+}$, $\pa\Gamma$. Then:
  \begin{enumerate}
  \item In the forward direction, $\gamma$ either tends to $\pa\Gamma$ or $\pa\cR^+_{\pa\cK^+,{\rm out}}$ as $s\to\infty$, or $\gamma$ reaches $r=\bhm$ in finite time.
  \item In the backward direction, $\gamma$ either tends to $\pa\Gamma$, $\pa\cR^+_{\pa\cK^+,{\rm in}}$, or $\pa\cR^+_{\cH^+}$ as $s\to-\infty$.
  \item $\gamma$ cannot tend to $\pa\Gamma$ in \emph{both} the forward and backward directions.
  \end{enumerate}
\end{prop}
\begin{proof}
  \pfstep{Case~1: $\exists\,s_0\in I$ with $r(\gamma(s_0))\in(r_+,\infty)$.} By Lemma~\ref{LemmaTs3bONonTrap}, $r\circ\gamma$ then leaves every compact subset of $(r_+,\infty)$ in both directions unless $\gamma(s_0)\in\pa\Gamma^{\rm s}\cup\pa\bar\Gamma^{\rm u}$. If $\gamma(s_0)\in\pa\bar\Gamma^{\rm u}$, then $\gamma$ tends to $\pa\Gamma$ in the backward direction and $r\circ\gamma$ leaves every compact subset of $(r_+,\infty)$ in the forward direction by~\eqref{EqTs3bOGammapmFlow}. If, on the other hand, $\gamma(s_0)\in\pa\Gamma^{\rm s}$, then $\gamma$ tends to $\pa\Gamma$ in the forward direction and $r\circ\gamma$ leaves every compact subset of $(r_+,\infty)$ in the backward direction.

  Consider the case that in the forward, resp.\ backward direction, $\limsup r\circ\gamma=\infty$. Lemma~\ref{LemmaTs3bDyn}\eqref{ItTs3bDynGeod} implies that $\gamma(s)$ tends to $\pa\cR^+_{\pa\cK^+,{\rm out}}$ as $s\to\infty$, resp.\ to $\pa\cR^+_{\pa\cK^+,{\rm in}}$ as $s\to-\infty$.

  If in the forward direction starting from $\gamma(s_0)$, the function $r\circ\gamma$ gets arbitrarily close to $r=r_+$, then by Lemma~\ref{LemmaTs3bOConvex}, $r\circ\gamma$ is, eventually, monotonically decreasing. Since $\pa\cR^+_{\cH^+}$ is a source for the $\sfH$-flow by Lemma~\ref{LemmaTs3bHQuadDef}, the only accumulation points of $\gamma$ (while it remains over $r^{-1}((r_+,\infty))$) in the forward direction lie in $\pa\Sigma^+\cap r^{-1}(r_+)\setminus\pa\cR^+_{\cH^+}$. Since $\sfH r<0$ there by Lemma~\ref{LemmaTs3bHr}, $\gamma$ must, in fact, cross $r=r_+$ in finite time. But once $r\circ\gamma$ has attained a value $<r_+$, it is monotonically decreasing by Lemma~\ref{LemmaTs3bHr}, and $\gamma$ reaches $r^{-1}(\bhm)$ in finite time.

  If $r\circ\gamma$ gets arbitrarily close to $r=r_+$ in the \emph{backward} direction, then $\gamma$ must tend to the source $\pa\cR^+_{\cH^+}$ as $s\to-\infty$, for otherwise Lemma~\ref{LemmaTs3bHr} would force $\gamma$ to move into $\{r<r_+\}$ at some argument $s_1<s_0$, which would imply that $r\circ\gamma$ would be less than $r_+$ for $s>s_1$, so in particular at $s=s_0$.

  \pfstep{Case~2: $r(\gamma(s))\leq r_+$ for all $s$.} By Lemma~\ref{LemmaTs3bHr}, $\gamma$ must reach $r^{-1}(\bhm)$ in the forward direction in finite affine time. In the backward direction, $r\circ\gamma$ is increasing while it is less than $r_+$, and thus $\gamma$ must tend to $\pa\cR^+_{\cH^+}$ as $s\to-\infty$, for otherwise $\gamma$ would have an accumulation point in the backward direction lying in $\pa\Sigma^+\cap r^{-1}(r_+)\setminus\pa\cR^+_{\cH^+}$, where Lemma~\ref{LemmaTs3bHr} would apply to show that $\gamma$ would have to cross $r=r_+$ from the domain of outer communications into the black hole interior; so we would be in Case~1 after all.

  \pfstep{Case~3: $r(\gamma(s))=\infty$ for all $s$.} Lemma~\ref{LemmaTs3bDyn} shows that $\gamma$ tends to $\pa\cR^+_{\pa\cK^+,{\rm in}}$ in the backward and $\pa\cR^+_{\pa\cK^+,{\rm out}}$ in the forward direction.
\end{proof}

\subsection{Dynamics in edge-b-phase space near \texorpdfstring{$\scri^+$}{null infinity}}
\label{SsTseb}

We only describe the null-bicharacteristic flow near $\iota^+\cap\scri^+$ following \cite{HintzVasyScrieb}; see \cite[\S{4.1}]{HintzVasyScrieb} for a description of the flow in a full neighborhood of $\scri^+$. Since we work away from $\cK^+$, we shall now drop the `3' from `$\etbop$' from the notation. We use the coordinates
\begin{equation}
\label{EqTsebCoord}
  \rho_+ = \frac{1}{t_*},\quad
  x_\sscri = \sqrt{\frac{t_*}{r}},\quad
  \omega = \frac{x}{|x|}
\end{equation}
in the base and fiber-linear coordinates on $\Teb^*M$ by writing covectors as $\zeta=\sigma_\ebop\frac{\dd\rho_+}{\rho_+}+\xi_\ebop\frac{\dd x_\sscri}{x_\sscri}+\eta_\ebop\frac{\dd\omega}{x_\sscri}$ as in~\eqref{EqCHamebCoord}. Lemma~\ref{LemmaTsKLMetric}\eqref{ItTsKLMetrice3b} implies that the rescaled dual metric function on Kerr,
\[
  G_\ebop(\zeta) := x_\sscri^{-2}\rho_+^{-2}g_{\bhm,a}^{-1}(\zeta,\zeta) \in P^{[2]}(\Teb^*M),
\]
is equal to its Minkowskian analogue $\ubar G_\ebop(\zeta):=x_\sscri^{-2}\rho_+^{-2}\ubar g^{-1}(\zeta,\zeta)$ modulo $x_\sscri^2\rho_+ P^{[2]}(\Teb^*M)$, which in turn, by~\eqref{EqCTebMink}, is equal to
\[
  G^0_\ebop(\zeta) = \frac12\xi_\ebop(\xi_\ebop-2\sigma_\ebop) + |\eta_\ebop|_{\slg^{-1}}^2
\]
modulo $x_\sscri^2 P^{[2]}(\Teb^*M)$. Using~\eqref{EqCHameb} and~\eqref{EqCHame3bMem}, we thus have
\begin{align*}
  H_{G_\ebop}&\equiv H_{G^0_\ebop} \bmod x_\sscri^2\Veb(\Teb^*M), \\
  H_{G^0_\ebop}&=-\xi_\ebop\rho_+\pa_{\rho_+} + (\xi_\ebop-\sigma_\ebop)(x_\sscri\pa_{x_\sscri}+\eta_\ebop\pa_{\eta_\ebop}) - 2|\eta_\ebop|_{\slg^{-1}}^2\pa_{\xi_\ebop} + x_\sscri H_{|\eta_\ebop|_{\slg^{-1}}^2}.
\end{align*}
Since $\dd t_*$ is past causal for the Minkowski metric, and thus $\frac{\dd t_*}{t_*}=-\frac{\rho_+}{\rho_+}$ is past causal for the eb-metric $x_\sscri^2\rho_+^2 g_{\bhm,a}$ at $\iota^+$, we have $\sigma_\ebop<0$ on the past characteristic set. Working in the future characteristic set $\Sigma^+=G_\ebop^{-1}(0)\subset\Teb^*M\setminus o$, we may thus pass to the eb-cosphere bundle via
\begin{equation}
\label{EqTsebCoordFiber}
  \rho_\infty := \frac{1}{\sigma_\ebop},\quad
  \hat\xi_\ebop := \frac{\xi_\ebop}{\sigma_\ebop},\quad
  \hat\eta_\ebop := \frac{\eta_\ebop}{\sigma_\ebop};
\end{equation}
we then compute $\rho_\infty^2 G_\ebop^0=\frac12\hat\xi_\ebop(\hat\xi_\ebop-2)+|\hat\eta_\ebop|^2$ and
\begin{equation}
\label{EqTsebHam}
\begin{split}
  \sfH_{G_\ebop} &:= \rho_\infty H_{G_\ebop} \\
    &\equiv \sfH_{G^0_\ebop} := \rho_\infty H_{G^0_\ebop} = -\hat\xi_\ebop\rho_+\pa_{\rho_+} + (\hat\xi_\ebop-1)(x_\sscri\pa_{x_\sscri}+\hat\eta_\ebop\pa_{\hat\eta_\ebop}) - 2|\hat\eta_\ebop|^2\pa_{\hat\xi_\ebop} + x_\sscri H_{|\hat\eta_\ebop|_{\slg^{-1}}^2}
\end{split}
\end{equation}
modulo the space of eb-vector fields on $\ol{\Teb^*}M\setminus o$ that are tangent to $\Seb^*M$ and vanish at $\scri^+$. (Here we abuse notation and write $H_{|\eta_\ebop|^2_{\slg^{-1}}}$ for the Hamiltonian vector field of $|\eta|_{\slg^{-1}}^2$ on $T^*\Sph^2$ evaluated at $\eta=\hat\eta_\ebop$.) At $\iota^+=\rho_+^{-1}(0)$, we have $\sfH_{G_\etbop}\hat\xi_\ebop=-2|\hat\eta_\ebop|^2+\cO(x_\sscri)$, the vanishing of which forces $|\hat\eta_\ebop|=o(1)$ (as $x_\sscri\to 0$), which on $\Sigma^+$ implies $\hat\xi_\ebop=o(1)$ or $\hat\xi_\ebop=2+o(1)$. In the latter case, we have $\sfH_{G_\ebop}x_\sscri=(1+o(1))x_\sscri$, which is thus nonzero unless $x_\sscri=0$. Thus, over $\iota^+$, the vector field $\sfH_{G_\etbop}$ vanishes only over $\iota^+\cap\scri^+$. Over $\scri^+$, we have $\sfH_{G_\etbop}=\sfH_{G_\ebop^0}$, the vanishing of which forces $\hat\eta_\ebop=0$ and thus $\hat\xi_\ebop\in\{0,2\}$, but when $\hat\xi_\ebop=2$ we have $\sfH_{G_\etbop}\rho_+=-4\rho_+$, which is nonzero unless $\rho_+=0$. We have thus identified the following sets:

\begin{definition}[Radial sets over $\scri^+$]
\label{DefTsebRad}
  In the coordinates~\eqref{EqTsebCoord} and \eqref{EqTsebCoordFiber} on $\ol{\Teb^*}M$ near $\Sigma^+\cap\Seb^*_{\scri^+\cap\iota^+}M$, we define the subsets
  \begin{alignat*}{3}
    \cR^+_{\scri^+,{\rm out}} &:= \{ (\rho_+,x_\sscri,\omega;\sigma_\ebop,\xi_\ebop,\eta_\ebop) \colon &&x_\sscri=0,\ \sigma_\ebop>0,\ \xi_\ebop=0,\ \eta_\ebop=0 \}, \\
    \cR^+_{\scri^+,{\rm in},+} &:= \{ (\rho_+,x_\sscri,\omega;\sigma_\ebop,\xi_\ebop,\eta_\ebop) \colon \rho_+=0,\ &&x_\sscri=0,\ \sigma_\ebop>0,\ \xi_\ebop=2\sigma_\ebop,\ \eta_\ebop=0 \}
  \end{alignat*}
  of $\Teb^*_{\scri^+}M\setminus o$ and $\Teb^*_{\scri^+\cap\iota^+}M\setminus o$, called \emph{outgoing} and \emph{incoming radial sets over $\scri^+$} and $\scri^+\cap\iota^+$, respectively. We furthermore denote by
  \begin{equation}
  \label{EqTsebRadS}
  \begin{split}
    \pa\cR^+_{\scri^+,{\rm out}} &= \{ \hat\xi_\ebop=0,\ \hat\eta_\ebop=0 \} \cap \Seb^*_{\scri^+}M, \\
    \pa\cR^+_{\scri^+,{\rm in},+} &= \{ \hat\xi_\ebop=2,\ \hat\eta_\ebop=0 \} \cap \Seb^*_{\scri^+\cap\iota^+}M
  \end{split}
  \end{equation}
  their boundaries at fiber infinity.
\end{definition}

To compute the linearization of $\sfH_{G_\ebop}$ at the sets~\eqref{EqTsebRadS} inside of $\pa\Sigma^+$, note that $\hat\xi_\ebop$ is, locally, a function of $\rho_+,x_\sscri,\hat\eta_\ebop$, so we only compute
\begin{equation}
\label{EqTsebRadLin}
\begin{alignedat}{5}
  &\text{at}\ \pa\cR^+_{\scri^+,{\rm out}},&&\ [\sfH_{G_\ebop}] \colon &\ \dd x_\sscri &\mapsto -\dd x_\sscri, &\ \dd\hat\eta_\ebop &\mapsto -\dd\hat\eta_\ebop, \\
  &\text{at}\ \pa\cR^+_{\scri^+,{\rm in},+},&&\ [\sfH_{G_\ebop}] \colon \dd\rho_+ \mapsto -2\,\dd\rho_+, &\ \dd x_\sscri &\mapsto \dd x_\sscri, &\ \dd\hat\eta_\ebop &\mapsto \dd\hat\eta_\ebop.
\end{alignedat}
\end{equation}
Thus, $\pa\cR^+_{\scri^+,{\rm out}}$ is a sink, while $\pa\cR^+_{\scri^+,{\rm in},+}$ is a saddle point.

\begin{lemma}[Dynamics near $\scri^+$]
\label{LemmaTsebDyn}
  The $\sfH_{G_\ebop}$-flow in $\pa\Sigma^+$ in $\Seb^*_{\scri^+\cup\iota^+\setminus\cK^+}M$ has the following properties.
  \begin{enumerate}
  \item\label{ItTsebDynChar}{\rm (Null-bicharacteristics over $\scri^+$.)} Let $\gamma\colon I\to\pa\Sigma^+$, $I\subseteq\R$, be a maximally extended integral curve of $\sfH_{G_\ebop}$ lying over $\scri^+\cap\{t_*\geq 1\}$ that is not contained in $\pa\cR^+_{\scri^+,{\rm out}}\cup\pa\cR^+_{\scri^+,{\rm in},+}$ (otherwise it would be constant). Then in the forward direction, $\gamma$ tends to $\pa\cR^+_{\scri^+,{\rm in},+}$ or $\pa\cR^+_{\scri^+,{\rm out}}$, while in the backward direction, $\gamma$ either reaches $t_*^{-1}(1)$ in finite time or tends to $\pa\cR^+_{\scri^+,{\rm in},+}$.
  \item\label{ItTsebDynConv}{\rm (Convexity in $\iota^+$.)} If $H_{G_\ebop}x_\sscri=0$, $x_\sscri>0$, then $\sfH_{G_\ebop}^2 x_\sscri<0$.
  \item\label{ItTsebDynMfd}{\rm ((Un)stable manifolds.)} The stable manifold of $\pa\cR^+_{\scri^+,{\rm out}}$ is
  \begin{equation}
  \label{EqTsebDynStable}
    \pa\cW^+_{\rm out} = \{ \rho_+=0,\ \hat\xi_\ebop=0,\ \hat\eta_\ebop=0 \}.
  \end{equation}
  The stable manifold of $\pa\cR^+_{\scri^+,{\rm in},+}$ over $\scri^+$ is $\pa\cW_{\scri^+}^+:=\{\hat\xi_\ebop=2,\ \hat\eta_\ebop=0\}\cap\Seb^*_{\scri^+}M$.
  \end{enumerate}
\end{lemma}

The manifold~\eqref{EqTsebDynStable} is the same as that defined in~\eqref{EqTs3bStable}, but now expressed in our current coordinate system.

\begin{proof}[Proof of Lemma~\usref{LemmaTsebDyn}]
  For part~\eqref{ItTsebDynChar}, since $\gamma$ is not contained in $\pa\cR^+_{\scri^+,{\rm out}}$, we have $\hat\xi_\ebop>0$; in view of $\sfH_{G_\ebop}\rho_+=-\hat\xi_\ebop\rho_+\leq 0$, $\rho_+$ is thus monotonically decreasing in the forward direction along $\gamma$. If $\gamma$ lies over $(\scri^+)^\circ$, then $\sfH_{G_\ebop}\rho_+=-\hat\xi_\ebop\rho_+<0$ implies that $t_*\circ\gamma$ is decreasing in the backward direction; accumulation points of $\gamma$ in the backward direction cannot lie in $\pa\cR^+_{\scri^+,{\rm out}}$ due to the source nature of this set (cf.\ \eqref{EqTsebRadLin}), and thus $\gamma$ must, in fact, reach $t_*=1$ in finite affine time. In general, we have $\sfH_{G_\ebop}\hat\xi_\ebop=-2|\hat\eta_\ebop|^2=-\hat\xi_\ebop(2-\hat\xi_\ebop)$ on $\pa\Sigma^+$ over $\scri^+$, and thus $\hat\xi_\ebop\to 0$ along $\gamma$ in the forward direction, and thus also $\hat\eta_\ebop\to 0$ on $\pa\Sigma^+$, so $\gamma$ tends to $\pa\cR^+_{\scri^+,{\rm out}}$. Finally, when $\gamma$ lies over $\scri^+\cap\iota^+$, we have $\hat\xi_\ebop\to 2$ in the backward direction, and thus $\hat\eta_\ebop\to 0$, and thus $\gamma$ tends to $\pa\cR^+_{\scri^+,{\rm in},+}$.

  Part~\eqref{ItTsebDynConv} is the same as Lemma~\ref{LemmaTs3bDyn}\eqref{ItTs3bDynip}. In present coordinates, we have $\sfH_{G_\ebop}x_\sscri=x_\sscri(\hat\xi_\ebop-1+o(1))$ (as $x_\sscri\to 0$), the vanishing of which at some point with $x_\sscri>0$ implies $\hat\xi_\ebop=1+o(1)$ and $x_\sscri^{-1}\sfH_{G_\ebop}^2 x_\sscri=\sfH_{G_\ebop}\hat\xi_\ebop+o(1)=-2|\hat\eta_\ebop|^2+o(1)$. But on $\pa\Sigma^+$, $\hat\xi_\ebop=1+o(1)$ implies $|\hat\eta_\ebop|^2=\frac12+o(1)$, which yields the desired conclusion.

  For part~\eqref{ItTsebDynMfd}, it suffices to note that the stated manifolds are invariant under the $\sfH_{G_\ebop}$-flow and have the correct tangent spaces and dimensions at the outgoing and incoming radial sets.
\end{proof}

\subsection{Dynamics in e3b-phase space}
\label{SsTse3b}

We now describe the null-bicharacteristic flow in $\Tetb^*M\setminus o$ over the domain $\Omega=\ol{\{t_*\geq 1\}}$ (see Definition~\ref{DefCMDomain}). We pass to the cosphere bundle by fixing a defining function
\[
  \rho_\infty \in \CI(\ol{\Tetb^*}M)
\]
of ${}^\etbop S^*M$ and considering the rescaled Hamiltonian vector field
\begin{equation}
\label{EqTse3bVF}
  \sfH_G := \rho_\infty x_\sscri^{-2}\rho_+^{-2}H_G
\end{equation}
where $x_\sscri$ and $\rho_+\in\CI(M)$ denote defining functions of $\scri^+$ and $\iota^+\subset M$, respectively. Thus, $\sfH_G$ is an edge-3b-vector field on $\ol{\Tetb^*_\Omega}M$. Note that on the characteristic set $\Sigma$, we have $\sfH_G=\rho_\infty H_{x_\sscri^{-2}\rho_+^{-2}G}$, and $\varrho^2$ is a smooth positive multiple of $\rho_+^{-2}$ near $\cK^+$; thus $\sfH_G$ is (up to smooth positive multiples) equal to the vector field $\sfH_{G_\tbop}$ studied in~\S\ref{SsTs3b}, and equal to $\sfH_{G_\ebop}$ studied in~\S\ref{SsTseb}.

\begin{prop}[Global null-bicharacteristic dynamics on Kerr]
\label{PropTse3bDyn}
  Let $\gamma\colon I\to{}^\etbop S^*_\Omega M$ denote a maximally extended integral curve of $\sfH_G$ inside of the boundary $\pa\Sigma^+$ at fiber infinity of the future characteristic set $(\rho_\infty^2 x_\sscri^{-2}\rho_+^{-2}G)^{-1}(0)\subset\Tetb_\Omega^*M\setminus o$. Let us say that $\gamma$ is \emph{of type $(A,B)$}, where $A$ and $B$ are one of $\ol{\{t_*=1\}}$, $\pa\cR^+_{\cH^+}$ (Definition~\usref{DefTs3bH}), $\pa\cR^+_{\scri^+,{\rm in},+}$, $\pa\cR^+_{\scri^+,{\rm out}}$ (Definition~\usref{DefTsebRad}), $\pa\cR^+_{\pa\cK^+,{\rm in/out}}$ (Definition~\usref{DefTs3bRad}), $\pa\Gamma$ (see~\eqref{EqTs3bGamma}), $\ol{\{r=\bhm\}}$, provided that $\gamma$ tends to $A$ in the backward direction (or reaches $A$ in finite time when $A$ is $\ol{\{t_*=1\}}$) and to $B$ in the forward direction (or reaches $B$ in finite time when $B$ is $\ol{\{r=\bhm\}}$). Then either $\gamma$ is contained in $\pa\cR_{\cH^+}^+$, $\pa\cR^+_{\scri^+,{\rm in},+}$, $\pa\cR^+_{\pa\cK^+,{\rm in}}$, $\pa\Gamma$, $\pa\cR^+_{\pa\cK^+,{\rm out}}$, or $\pa\cR^+_{\scri^+,{\rm out}}$, or $\gamma$ is of one of the following types (and each of these types occurs):
  \begin{equation}
  \label{EqTse3bDynTable}
  \vcenter{\hbox{
  \begin{tabular}{l|l}
    $A$ & $B$ \\
    \hline
    $\ol{\{t_*=1\}}$ & $\ol{\{r=\bhm\}}$, $\pa\cR^+_{\cH^+}$, $\pa\Gamma$, $\pa\cR^+_{\scri^+,{\rm in},+}$, $\pa\cR^+_{\scri^+,{\rm out}}$ \\
    $\pa\cR^+_{\cH^+}$ & $\ol{\{r=\bhm\}}$, $\pa\Gamma$, $\pa\cR^+_{\pa\cK^+,{\rm out}}$ \\
    $\pa\cR^+_{\scri^+,{\rm in},+}$ & $\pa\cR^+_{\pa\cK^+,{\rm in}}$, $\pa\cR^+_{\scri^+,{\rm out}}$ \\
    $\pa\cR^+_{\pa\cK^+,{\rm in}}$ & $\ol{\{r=\bhm\}}$, $\pa\Gamma$, $\pa\cR^+_{\pa\cK^+,{\rm out}}$ \\
    $\pa\Gamma$ & $\ol{\{r=\bhm\}}$, $\pa\cR^+_{\pa\cK^+,{\rm out}}$ \\
    $\pa\cR^+_{\pa\cK^+,{\rm out}}$ & $\pa\cR^+_{\scri^+,{\rm out}}$
  \end{tabular}
  }}
  \end{equation}
\end{prop}

\begin{figure}[!ht]
\centering
\includegraphics{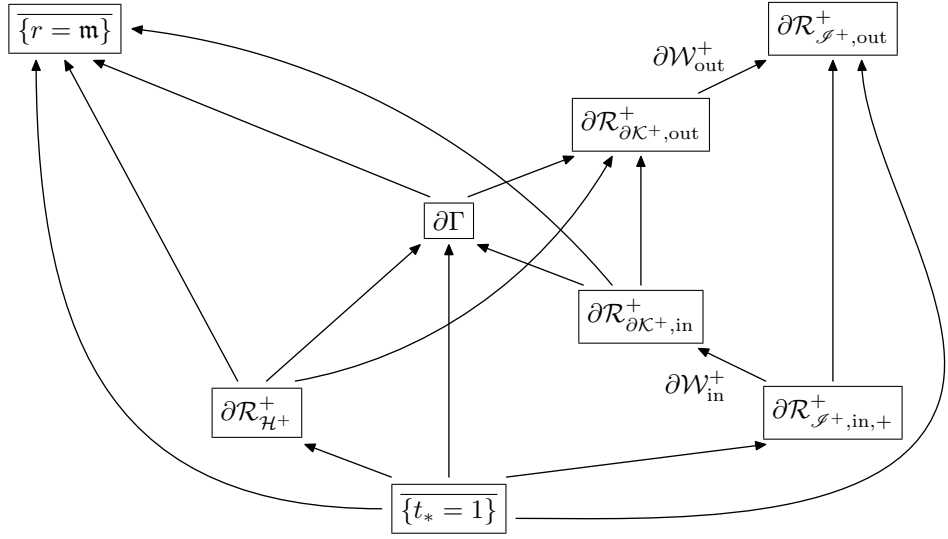}
\caption{Illustration of Proposition~\ref{PropTse3bDyn}. An arrow from a node labeled $A$ to a node labeled $B$ corresponds to a null-bicharacteristic (integral curve of $\sfH_G$ inside of $\pa\Sigma^+$ over $\Omega$) tending to $A$ and $B$ in the backward and forward direction, respectively.}
\label{FigTse3bDyn}
\end{figure}

Note that the value of $A$ in the $j$-th row does not appear as a value of $B$ in the $l$-th row for any $l\geq j$. Thus, one can arrange the radial sets, trapped set, and hypersurfaces in the form of a directed graph, with root node $\ol{\{t_*=1\}}$ (which does not appear in the $B$ column) and final nodes $\ol{\{r=\bhm\}}$ and $\pa\cR^+_{\scri^+,{\rm out}}$ (which do not appear in the $A$ column); see Figure~\ref{FigTse3bDyn}.

\begin{proof}[Proof of Proposition~\usref{PropTse3bDyn}]
  When $\gamma$ lies over $\cK^+$ or $\scri^+$, this follows from Proposition~\ref{PropTs3bSum} and Lemma~\ref{LemmaTsebDyn}.

  Consider the case that $\gamma$ lies over $\iota^+\setminus(\cK^+\cup\scri^+)$. Suppose that $\rho_\cK\circ\gamma$ is increasing at a point $\gamma(s_0)$, then it is monotonically increasing for $s\geq s_0$ since by Lemma~\ref{LemmaTs3bDyn}\eqref{ItTs3bDynip} it cannot have a local maximum; therefore $x_\sscri\circ\gamma$ is decreasing. The only accumulation points of $\gamma$ in the forward direction must then lie over $\iota^+\cap\scri^+$, and due to the source nature of $\pa\cR^+_{\scri^+,{\rm in},+}$ over $\iota^+$, $\gamma$ must tend to $\pa\cR^+_{\scri^+,{\rm out}}$. If, on the other hand, $\rho_\cK\circ\gamma$ is always non-increasing, it must be decreasing, and then (in view of the source nature of $\pa\cR^+_{\pa\cK^+,{\rm out}}$ over $\iota^+$) must tend to $\pa\cR^+_{\pa\cK^+,{\rm in}}$ in the forward direction. (In this case, $\gamma$ lies in $\pa\cW^+_{\rm in}$, see~\eqref{EqTs3bStable}, and tends to $\pa\cR^+_{\scri^+,{\rm in},+}$ in the backward direction.) Similar arguments show that in the backward direction, $\gamma$ must tend either to $\pa\cR^+_{\scri^+,{\rm in},+}$ or $\pa\cR^+_{\pa\cK^+,{\rm out}}$. (In the latter case, $\gamma$ lies in $\pa\cW^+_{\rm out}$ and tends to $\pa\cR^+_{\scri^+,{\rm out}}$ in the forward direction.) (One can, alternatively, compute the $\sfH_G$-flow over $\iota^+$ explicitly in order to reach the same conclusions.)

  It remains to consider $\gamma$ lying over $M^\circ$; such $\gamma$ must reach $t_*=1$ in the backward direction in finite time, and $t_*\circ\gamma$ is monotonically increasing in the forward direction. If $r\circ\gamma$ attains a value in $[\bhm,r_+)$, then $\gamma$ reaches $r=\bhm$ in the forward direction. Otherwise, $r(\gamma(s))\geq r_+$ for all $s\in I$. Every accumulation point of $\gamma$ must then lie over $\scri^+\cup\iota^+\cup\cK^+\cap\ol{\{r\geq r_+\}}$. Since the stable manifolds of $\pa\cR^+_{\pa\cK^+,{\rm out}}$, $\pa\cR^+_{\pa\cK^+,{\rm in}}$, and $\pa\cR^+_{\scri^+,{\rm in},+}$ lie over $\cK^+$, $\iota^+$, and $\scri^+$, respectively, $\gamma$ cannot tend to either of them. The only possibilities are, therefore, that $\gamma$ tends to $\pa\cR^+_{\cH^+}$ (which happens when $\gamma(0)\in S N^*\cH^+=(N^*\cH^+\setminus o)/\R_+$), $\pa\Gamma$ (which happens when $\gamma(0)\in\pa\Gamma_0$), or $\pa\cR^+_{\scri^+,{\rm out}}$ (which happens, e.g., when $\frac{\dd}{\dd s}r(\gamma(0))>0$ and $r(\gamma(0))$ is sufficiently large).
\end{proof}

\begin{rmk}[Smooth dependence on Kerr parameters]
\label{RmkTse3bKerrDep}
  The only sets in Proposition~\ref{PropTse3bDyn} depending on the Kerr parameters $\bhm,a$ are $\pa\cR^+_{\cH^+}$ and $\pa\Gamma$. Both depend smoothly on $(\bhm,a)$ in a neighborhood of fixed subextremal Kerr parameters $(\bhm_0,a_0)$: for $\pa\cR^+_{\cH^+}$, this follows immediately from the formula $\bhm+\sqrt{\bhm^2-a^2}$ for the area radius of the event horizon; for $\pa\Gamma$, it follows via the implicit function theorem from the proof of Lemma~\ref{LemmaTs3bOTrap0Reg}, concretely, the linear independence of the differentials of the defining functions $\xi,\Psi,G_\tbop$ of $\Gamma_0$ in~\eqref{EqTs3bOTrap0RegDiff} and the smooth dependence of these functions on $(\bhm,a)$.
\end{rmk}

\section{Geometric and analytic setup}
\label{SS}

In~\S\ref{SsSS}, we will introduce the class of stationary (i.e., time-translation-invariant) wave-type operators on exact subextremal Kerr spacetimes which we allow for as stationary models for dynamical problems. In~\S\ref{SsSD} then, we first describe the metric perturbations which we consider (\S\ref{SssSDG}) before describing the class of dynamical (but asymptotically stationary) wave-type operators on such perturbed Kerr spacetimes. Throughout, we fix subextremal Kerr parameters
\[
  \bhm > 0,\quad a\in(-\bhm,\bhm),
\]
the (compactified) spacetime manifolds $M_0,M$ from Definition~\ref{DefCMSpacetime} and the domain $\Omega=\ol{\{t_*\geq 1\}}$ (where $t_*=t-r$), with the closure taken in $M_0$ or $M$ depending on the context. We recall the spatial manifold $X=[\bhm,\infty]_r\times\Sph^2$ from Definition~\ref{DefCSpatial} and the area radius $r_+$ of the event horizon from~\eqref{EqTsBLMfd}. We write $\rho:=r^{-1}$ for a boundary defining function of $X$, and we use the notation $x_\sscri$, $\rho_+$, $\rho_\cK$ for (local) boundary defining functions of $\scri^+$, $\iota^+$, $\cK^+$, respectively. We write $\pi_X\colon M\to X$ for the projection $(t,x)\mapsto x$.\footnote{Note here that $(t,x)\mapsto x$ induces a smooth map $\tilde M_0\to\ol{\R^3}$, which then lifts to a smooth map $\tilde M\to\ol{\R^3}$ whose restriction to $M$ is the map $\pi_X$.}

\subsection{The stationary model}
\label{SsSS}

We consider the following class of stationary wave-type operators.

\begin{definition}[Stationary wave-type operators]
\label{DefSSAdm}
  Let $\cE\to M$ be a stationary vector bundle, i.e., the pullback $\cE=\pi_X^*\cE_X$ of a vector bundle $\cE_X\to X$. Then an operator $P_0\in\Diff^2(M^\circ;\cE)$ commuting with translations in $t_*$ is called a \emph{stationary wave-type operator on $(M,g_{\bhm,a})$} if the following conditions hold.
  \begin{enumerate}
  \item{\rm (Principal symbol.)} $\upsigma^2(P_0)(z,\zeta)=G(z,\zeta):=g_{\bhm,a}^{-1}|_z(\zeta,\zeta)$.
  \item{\rm (Structure of lower-order terms.)} Near $\rho=r^{-1}=0$, $P_0$ is of the form
    \begin{subequations}
    \begin{equation}
    \label{EqSSAdmOp}
      P_0 = -2\pa_{t_*}\rho(\rho\pa_\rho-1-S) + \wh{P_0}(0) + Q\pa_{t_*} - g^{0 0}\pa_{t_*}^2,
    \end{equation}
    where $g^{0 0}=g_{\bhm,a}^{-1}(\dd t_*,\dd t_*)\in\rho^2\CI(X)$ and
    \begin{equation}
    \label{EqSSAdmOpPieces}
      \wh{P_0}(0) \in \rho^2\Diffb^2(X;\cE_X), \quad
      S \in \CI(X;\End(\cE_X)), \quad
      Q \in \rho^3\Diffb^1(X;\cE_X).
    \end{equation}
    \end{subequations}
  \end{enumerate}
\end{definition}

The quantity
\begin{equation}
\label{EqSSAdmubarS}
  \ubar S := \inf_{\substack{p\in\pa X \\ \mu\in\spec S(p)}} \Re\mu
\end{equation}
arises frequently in the analysis of $P_0$ and is closely related to the decay rate of solutions of $P_0 u=f$ at $\scri^+$.

\begin{rmk}[Information from the principal symbol ]
\label{RmkSSAdmSymb}
The principal symbols of $\wh{P_0}(0)$ and $Q$ are determined by $g_{\bhm,a}$. In particular, by~\eqref{EqTsKLMetricKerr} and \eqref{EqTsKLMetric}, and the membership $\tilde T_\gg\in r^{-2}\CI$ in Lemma~\ref{LemmaTsKLSpace}, the principal symbol $g_{\bhm,a}^{-1}(\dd t_*,\cdot)+2\la\cdot,\pa_r\ra$ of $Q$ is of class $r^{-2}P^{[1]}(\Tsc^*X)=r^{-3}P^{[1]}(\Tb^*X)$. The only pieces of information in Definition~\ref{DefSSAdm} not determined by $g_{\bhm,a}$ are, therefore, the membership of $S$ and the precise expression for the lower order terms of $\wh{P_0}(0)$ and $Q$. Note also that the terms involving $S$ and $Q$ combine into $2\pa_{t_*}\rho S+Q\pa_{t_*}=(2\rho S+Q)\pa_{t_*}$, so any zeroth order term contributing to $Q$ can equivalently be moved over to $S$.
\end{rmk}

\begin{rmk}[Derivatives of sections of $\cE$]
\label{RmkSSAdmDer}
  We clarify the meaning of $\pa_{t_*}$ and $\rho\pa_\rho$ in~\eqref{EqSSAdmOp} when $\cE$ is not trivial. (Cf.\ the final part of Lemma~\ref{LemmaCT3bDil}.) Fix any smooth affine connection $\nabla^{\cE_X}\in\Diff^1(X;T^*X\otimes\cE_X)$ on $\cE_X$, and denote its pullback along $\pi_X$ by
  \begin{equation}
  \label{EqSSAdmDer}
    \nabla^\cE\in\Diff^1(M;T^*M\otimes\cE).
  \end{equation}
  We then define $\pa_{t_*}$ and $\rho\pa_\rho$ as $\nabla^\cE_{\pa_{t_*}}$ and $\nabla^\cE_{\rho\pa_\rho}$, respectively. They have the following properties.
  \begin{enumerate}
  \item In a trivialization of $\cE$ given by the pullback of a trivialization of $\cE$, the operator $\nabla^\cE_{\pa_{t_*}}$ acts as the component-wise derivative along $\pa_{t_*}$. It is, thus independent of the choice of connection $\nabla^{\cE_X}$.
  \item Note that $\nabla^\cE_{\rho\pa_\rho}=\rho\nabla^\cE_{\pa_\rho}$. Therefore, in any trivialization of $\cE$ and for a section $\vec e=(e_1,\ldots,e_k)$ of $\cE$ (with $k$ denoting the rank of $\cE$), we have $\nabla^\cE_{\rho\pa_\rho}\vec e=\rho\pa_\rho\vec e+\rho A\vec e$ where $\rho\pa_\rho$ acts component-wise and $A$ is a smooth $k\times k$-matrix valued function depending on the choice of connection. In particular, any two choices of connections on $\cE_X$ yield the same operator $\nabla^\cE_{\rho\pa_\rho}$ modulo corrections of class $\rho\CI(X;\End(\cE_X))$ (pulled back along $\pi_X$). Such corrections thus affect sub-leading terms of $S$ at $\pa X$. In particular, $S|_{\pa X}\in\CI(\pa X;\End(\cE_X))$ is well-defined independently of any choices.
  \end{enumerate}
\end{rmk}

\begin{example}[Tensor wave operators on Kerr]
\label{ExSSAdmBox}
  Recalling Lemma~\ref{LemmaTsKLMetric}\eqref{ItTsKLMetric3b}, it follows from the Koszul formula that the Levi-Civita connection of $g_{\bhm,a}$, acting on sections of a 3b-tensor bundle $\Ttb^{p,q}M_0:=(\Ttb M_0)^{\otimes p}(\Ttb^* M_0)^{\otimes q}$ defines an element
  \[
    \nabla \in \Diff_\tbop^1(M_0;\Ttb^{p,q}M_0,\Ttb^*M_0\otimes\Ttb^{p,q}M_0).
  \]
  Replacing $\Ttb^{p,q}M_0$ here by a standard tensor bundle $\cT^{p,q}:=\cT^{\otimes p}\otimes(\cT^*)^{\otimes q}$ (Definition~\eqref{DefCTscPullback}) amounts to conjugation of $\nabla$ by $\rho_\sface^{p-q}$, which preserves the space $\Diff_\tbop^1$. Thus, the tensor wave operator $\Box_{g_{\bhm,a}}^{(p,q)}:=\tr_g\nabla\circ\nabla$ on $(p,q)$-tensors satisfies
  \[
    \Box_{g_{\bhm,a}}^{(p,q)} \in \rho_\sface^2\Diff_\tbop^2(M_0;\cT^{p,q}).
  \]
  Using the frame $r\pa_{t_*}$, $r\pa_r$, $\pa_\omega$ of $\Vtb(M_0)$, this is of the schematic form $r^{-2}\CI(M_0)(r\pa_{t_*},r\pa_r,\pa_\omega)^{\leq 2}$. The terms $r^{-2}\CI(M_0)(r\pa_r,\pa_\omega)^{\leq 2}$ contribute to $\wh{P_0}(0)$ in~\eqref{EqSSAdmOp}, the term $r^{-2}\CI(M_0)(r\pa_{t_*})^2$ is the final term in~\eqref{EqSSAdmOp} (by principal symbol considerations as in Remark~\ref{RmkSSAdmSymb}), and the term $r^{-2}\CI(M_0)r\pa_{t_*}$ is encoded in $S$; finally, the mixed terms $r^{-2}\CI(M_0)r\pa_{t_*}\circ V$, $V=r\pa_r,\pa_\omega$, are encoded by $Q$ (see, again, the principal symbol considerations in Remark~\ref{RmkSSAdmSymb}). Finally, we claim that $S=0$: this follows from the fact that the difference of $\Box_{g_{\bhm,a}}^{(p,q)}$ and the Minkowski analogue $\Box_{\ubar g}^{(p,q)}$ has an extra order of decay at $\sface$ by Lemma~\ref{LemmaTsKLMetric}\eqref{ItTsKLMetric3b}; and $\Box_{\ubar g}^{(p,q)}$ acts component-wise as the scalar wave operator in the bundle trivialization induced by standard coordinates on $\R^4$, which takes the form
  \[
    \Box_{\ubar g} = -2\pa_{t_*}\rho(\rho\pa_\rho-1) + \rho^2\bigl( -(\rho\pa_\rho)^2 + \rho\pa_\rho + \slDelta \bigr)
  \]
  in the coordinates $t_*=t-r$, $\rho=r^{-1}$, and $\omega\in\Sph^2$. Since $S=0$ for this operator, we have $S\in\rho\CI$ for $\Box_{g_{\bhm,a}}^{(p,q)}$. The fact that $\Box_{g_{\bhm,a}}^{(p,q)}-\Box_{\ubar g}^{(p,q)}\in\rho_\sface^3\Difftb^2$ also implies that the spectral families of $\Box_{g_{\bhm,a}}^{(p,q)}$ and $\Box_{\ubar g}^{(p,q)}$ agree to leading order at $\pa X$; this is particularly useful at zero energy for the computation of indicial roots. In summary,
  \begin{equation}
  \label{EqSSAdmBox}
  \begin{split}
    \Box_{g_{\bhm,a}}^{(p,q)}\ &\text{is a stationary wave-type operator, with}\ S\in\rho\CI(X;\End(\cT^{p,q}))\ \text{(and thus $\ubar S=0$)}; \\
    &\text{and}\ \wh{\Box_{g_{\bhm,a}}^{(p,q)}}(0)-\wh{\Box_{\ubar g}^{(p,q)}}(0)\in\rho^3\Diffb^2(X;\cT^{p,q}).
  \end{split}
  \end{equation}
\end{example}

\begin{rmk}[Smoothness]
\label{RmkSSAdmCI}
  The smoothness requirement on the operators~\eqref{EqSSAdmOpPieces} can be relaxed as in \cite[Definition~3.8]{HintzNonstat}: $\wh{P_0}(0)$, $S$, and $Q$ may have $\rho^{2+\delta}\CI_\bop\Diffb^2$, $\rho^\delta\CI_\bop$, and $\rho^{2+\delta}\CI_\bop\Diffb^1$ contributions, respectively, where $\delta\in(0,1]$. This flexibility is not needed in our main application \cite{HintzKerrStab}, and hence we content ourselves with~\eqref{EqSSAdmOpPieces}.
\end{rmk}

\begin{rmk}[Dimension]
\label{RmkSSDim}
  In the present paper, we work in $3+1$ dimensions for concreteness. The constant $1$ in~\eqref{EqSSAdmOp} arises naturally when $P_0=\Box_{g_{\bhm,a}}$, in which case $S|_{\pa X}=0$ then. See \cite[\S{3.1.2}]{HintzNonstat} for matching conventions in general dimension.
\end{rmk}

The following is a variant of \cite[Lemma~3.27]{HintzNonstat}:

\begin{lemma}[$P_0$ as an (e)3b-operator]
\label{LemmaSSOpMem}
  In the notation of Definition~\usref{DefCMSpacetime}, we have
  \begin{equation}
  \label{EqSSOpMem}
    P_0 \in r^{-2}\Diff_\tbop^2(M_0;\cE),\quad
    P_0 \in \rho_0^2 x_\sscri^2\rho_+^2\Diff_\etbop^2(M;\cE).
  \end{equation}
\end{lemma}
\begin{proof}
  The 3b-membership follows from the fact that the space of 3b-vector fields on $M_0$ is spanned over $\CI(M_0)$ by $r\pa_{t_*}=\rho^{-1}\pa_{t_*}$, $r\pa_r=-\rho\pa_\rho$, and $\pa_\omega$ (spherical vector fields); and $\rho\pa_\rho,\pa_\omega$ span $\Vb(X)$ over $\CI(X)$.

  We check the e3b-membership only near $\scri^+\cap\iota^+$. At the principal symbol level, it follows directly from Lemma~\ref{LemmaTsKLMetric}\eqref{ItTsKLMetrice3b}. It remains to note that, in terms of $\rho_+=\frac{1}{t_*}$ and $x_\sscri=\sqrt{\frac{t_*}{r}}$, the vector fields
  \begin{equation}
  \label{EqSSOpMemVF}
  \begin{split}
    \rho\pa_{t_*} &= x_\sscri^2\rho_+^2\Bigl(\frac12 x_\sscri\pa_{x_\sscri} - \rho_+\pa_{\rho_+}\Bigr), \\
    \rho^2\rho\pa_\rho &= \frac12 x_\sscri^4\rho_+^2\cdot x_\sscri\pa_{x_\sscri}, \\
    \rho^2\pa_\omega &= x_\sscri^3\rho_+^2\cdot x_\sscri\pa_\omega
  \end{split}
  \end{equation}
  are of class $x_\sscri^2\rho_+^2\Veb$.
\end{proof}

For later use, we note that~\eqref{EqSSOpMemVF} implies that, near $\scri^+\cap\{\rho_+\leq 1\}$,
\begin{equation}
\label{EqSSOpMemeb}
  2 x_\sscri^{-2}\rho_+^{-2}P_0 \equiv -\bigl(x_\sscri\pa_{x_\sscri}-2(1+S)\bigr)(x_\sscri\pa_{x_\sscri}-2\rho_+\pa_{\rho_+}) + 2 x_\sscri^2\slDelta \bmod x_\sscri\Diffeb^2.
\end{equation}
Here we write $\slDelta$ for any element of $\Diff^2(M;\cE)$ with principal symbol given by the dual metric function of the standard 2-sphere $(\Sph^2,\slg)$ (e.g., $\slDelta$ is a sum of terms $\nabla^\cE_{V_a}\nabla^\cE_{V_b}$ for appropriate $V_a,V_b\in\cV(\Sph^2)$); any two choices lead to the same edge-b-operator $x_\sscri^2\slDelta$ modulo $x_\sscri\Diffeb^1$. To prove~\eqref{EqSSOpMemeb}, note that the principal part matches~\eqref{EqCTebMink} (up to a sign change due to the factors of $i^{-1}$ arising in the computation of principal symbols). Indeed, the terms $\rho^{-1}Q\rho\pa_{t_*}$ (with $\rho^{-1}Q\in\rho^2\Diffb^1\subset x_\sscri^3\rho_+^2\Diffeb^1$) and $g^{0 0}\pa_{t_*}^2=\rho^{-2}g^{0 0}(\rho\pa_{t_*})^2$ as well as lower order contributions $\rho^2\Diffb^1$ and $\rho^3\Diffb^2$ to $\wh{P_0}(0)$ do not contribute, so besides the first term in~\eqref{EqSSAdmOp}, only the principal part $\rho^2(\rho\pa_\rho)^2+\rho^2\slDelta$ of $\wh{P_0}(0)$ enters, with $\rho^2(\rho\pa_\rho)^2$ contributing only to the error term in~\eqref{EqSSOpMemeb}.

The stationary wave-type operators we shall study in this paper need to satisfy additional assumptions:
\begin{enumerate}
\item a condition on the subprincipal symbol of $P_0$ at the trapped set which, roughly speaking, prohibits too much amplitude growth of wave packets localized near trapped null-geodesics (\S\ref{SssSSTrap});
\item spectral assumptions which, roughly speaking, amount to mode stability (absence of purely oscillatory or exponentially growing mode solutions of $P_0$) and appropriate behavior near zero energy (\S\ref{SssSSSpec}).
\end{enumerate}
The full notion of ``admissibility'' which we require for $P_0$ is stated in Definition~\ref{DefSSAlephAdm} below.

\subsubsection{Admissibility, I: trapping}
\label{SssSSTrap}

We phrase the subprincipal symbol condition as in \cite[\S{4.4}]{HintzGlueLocII} using the \emph{subprincipal operator}: in a local trivialization of $\cE_X$, we define
\[
  S_\sub(P_0) := -i H_G + \upsigma_\sub(P_0)
\]
where $H_G$ acts component-wise, and $\upsigma_\sub(P_0)$ is the matrix of subprincipal symbols (see \cite[\S{5.2}]{DuistermaatHormanderFIO2} and \cite[Exercise~5.26]{HintzMicro}) of $P_0$, where we trivialize the 1-density bundle using the metric volume density $|\dd g_{\bhm,a}|$. As shown in \cite{HintzPsdoInner} (which is based on \cite{DenckerPolarization}), this is well-defined as an operator
\[
  S_\sub(P_0) \in \Diff^1(T^*M^\circ\setminus o;\pi^*\cE_X)
\]
where $\pi\colon T^*M^\circ\to M^\circ$ is the base projection. Since $P_0$ is stationary, also $S_\sub(P_0)$ commutes with time-translations.

\begin{definition}[Trapping admissibility]
\label{DefSSTrapAdm}
  We call $P_0$ \emph{trapping admissible} if there exists a stationary positive definite fiber inner product on the pullback bundle $\pi^*\cE_X\to T^*M^\circ\setminus o$, homogeneous of degree $0$ with respect to dilations in the fibers of $T^*M^\circ\setminus o$, such that
  \begin{equation}
  \label{EqSSTrapAdm}
    \rho_\infty\frac{1}{2 i}\bigl(S_\sub(P_0)-S_\sub(P_0)^*\bigr) < \frac12 \nu_{\rm min}\quad\text{at}\ \Gamma_0,\quad \rho_\infty:=\sigma^{-1},
  \end{equation}
  where $\Gamma_0$ and $\nu_{\rm min}$ were defined in Definition~\usref{DefTs3bOTrap0} and~\eqref{EqTs3bOnuMin}, and $\sigma$ is as in~\eqref{EqTs3bCoord}. (The left-hand side is a self-adjoint bundle endomorphism of $\pi^*\cE_X$.) If~\eqref{EqSSTrapAdm} holds with right-hand side $\eps$ for every $\eps>0$, with the fiber inner product allowed to depend on $\eps$, then we say that $P_0$ is \emph{strongly trapping admissible}.
\end{definition}

This is used as follows: fixing any positive definite fiber inner product on $\cE_X\to X$, \cite[Proposition~3.12]{HintzPsdoInner} produces a time-translation-invariant elliptic operator $Q\in\tilde\Psi_\tbop^0(\tilde M_0;\cE)$ such that
\begin{equation}
\label{EqSSTrapConj}
  \rho_\infty\upsigma^1\Bigl(\frac{1}{2 i}\bigl(Q P_0 Q^- - (Q P_0 Q^-)^*\bigr)\Bigr) < \frac12\nu_{\rm min}\quad\text{at}\ \Gamma_0;
\end{equation}
here $Q_-\in\tilde\Psi_\tbop^0(\tilde M_0;\cE)$ denotes an elliptic parametrix of $Q$.

\begin{rmk}[Strong trapping admissibility]
\label{RmkSSTrapAdmStrong}
  The operators $P_0$ arising in all applications here (i.e., the scalar and 1-form wave operators) and in \cite{HintzKerrStab} are strongly trapping admissible (the Einstein case being treated in \citestab{\S\ref*{SsWETr}}); and this fact is used in a crucial manner in the proof of arbitrarily high b-regularity for its solutions in the context of their $\aleph$-admissibility (see Definition~\ref{DefSSAlephAdm} and the discussion after~\eqref{EqA1AdmHiFinal}). An example of an operator that is trapping admissible but not strongly so is the scalar operator $\Box_{g_{\bhm,a}}+c\chi\pa_{t_*}$ where $\chi$ is a bump function near the base projection of the (time-translation-invariant) trapped set and $c$ is a small constant with the ``wrong'' sign ($c<0$). (Readers content with having merely large but finite b-regularity in linear and nonlinear results only need to require~\eqref{EqSSTrapAdm} with right hand side $\eps$ for \emph{some} sufficiently small $\eps>0$.)
\end{rmk}

\subsubsection{Admissibility, II: spectral assumptions}
\label{SssSSSpec}

The spectral family of $P_0$ is defined, analogously to~\eqref{EqMSFSpecFam}, by
\[
  \wh{P_0}(\sigma) = e^{i\sigma t_*}P_0 e^{-i\sigma t_*},\quad \sigma\in\C,
\]
acting on stationary sections of $\cE_X$; thus
\begin{equation}
\label{EqSSSpecFam}
  \wh{P_0}(\sigma) = 2 i\sigma\rho(\rho\pa_\rho-1-S) + \wh{P_0}(0) - i\sigma Q + g^{0 0}\sigma^2.
\end{equation}
It follows from~\eqref{EqMUSpecFam} (with $l=-2$, $m=2$) that, for $\hat\sigma\in e^{i[0,\pi]}$ and $\nu\geq 0$,\footnote{Other values of $\hat\sigma,\nu$ will not be of interest in our analysis.}
\begin{equation}
\label{EqSSSpecFamMem}
  (\wh{P_0}(\hat\sigma|\sigma|))_{|\sigma|\in[0,1]} \in \rho_\tface^2\Diff_\scbtop^2(X;\cE_X), \quad
  (\wh{P_0}(i\nu\pm h^{-1}))_{h\in(0,1]} \in h^{-2}\Diff_\schop^2(X;\cE_X),
\end{equation}
and in particular $\wh{P_0}(0)\in\rho^2\Diffb^2(X;\cE_X)$ and $\wh{P_0}(\sigma)\in\Diffsc^2(X;\cE_X)$, $\sigma\neq 0$. We define the transition face normal operator of the first family in~\eqref{EqSSSpecFamMem} (governing the transition from zero to nonzero frequencies) upon rescaling by $|\sigma|^{-2}$. Concretely, working in the collar neighborhood $[0,\bhm^{-1})_\rho\times\Sph^2$ of $\pa X\cong\Sph^2$, write
\begin{equation}
\label{EqSSSpec0Op}
  \wh{P_0}(0) = \rho^2 P_{(0)}(\rho,\omega,\rho\pa_\rho,\pa_\omega),
\end{equation}
then $N(\rho^{-2}\wh{P_0}(0))=P_{(0)}(0,\omega,\rho\pa_\rho,\pa_\omega)$ and thus
\begin{equation}
\label{EqSStfOp}
\begin{split}
  N_\tface(P_0,\hat\sigma) &:= N_\tface\bigl(|\sigma|^{-2}(\wh{P_0}(\hat\sigma|\sigma|))_{|\sigma|\in[0,1]}\bigr) \\
    &= 2 i\hat\sigma\hat\rho(\hat\rho\pa_{\hat\rho}-1-S) + \hat\rho^2 P_{(0)}(0,\omega,\hat\rho\pa_{\hat\rho},\pa_\omega),\quad \hat\rho:=\frac{\rho}{|\sigma|}.
\end{split}
\end{equation}
(This is an instance of~\eqref{EqMUNtfExpl}.) Thus, if in the notation of~\eqref{EqMUtf}--\eqref{EqMUNtf} we denote by
\[
  \rho_\sctface := \frac{\hat\rho}{1+\hat\rho},\quad
  \rho_\ztface := \frac{1}{1+\hat\rho}
\]
defining functions of $\sctface=\scface\cap\tface$ and $\ztface=\zface\cap\tface\subset\tface$, respectively, then
\begin{equation}
\label{EqSStfOpStruct}
  N_\tface(P_0,\hat\sigma)\in\Diff_{\scop,\bop}^{2,(0,2)}(\tface;\pi_\tface^*\cE_X)=\rho_\ztface^{-2}\Diff_{\scop,\bop}^2(\tface;\pi_\tface^*\cE_X);
\end{equation}
here $\pi_\tface\colon\tface=[0,\infty]_{\hat\rho}\times\Sph^2\to\Sph^2=\pa X$ is the projection. The following is then taken from~\cite[Definition~3.14(3)]{HintzNonstat}:

\begin{definition}[tf-admissibility]
\label{DefSStfAdm}
  We call a stationary wave-type operator $P_0$ \emph{$\tface$-admissible with weight $\beta\in\R$} if the following holds: $\beta$ lies in an indicial gap $(\beta^-,\beta^+)$ for $\wh{P_0}(0)$ (Definition~\usref{DefMUbInd} and Remark~\ref{RmkMUbIndGap}),\footnote{That is, $P_{(0)}(0,\omega,\rho\pa_\rho,\pa_\omega)(\rho^\lambda u(\omega))=0$ (or equivalently $P_{(0)}(0,\omega,\lambda,\pa_\omega)u=0$) implies $\Re\lambda\neq\beta$.} and for all $\alpha\in\R$, the nullspace of $N_\tface(P_0,e^{i\theta})$ on
  \begin{equation}
  \label{EqSStfAdmA}
    \cA^{\alpha,-\beta}(\tface;\pi_\tface^*\cE_X) := \rho_\sctface^\alpha\rho_\ztface^{-\beta}\CI_\bop(\tface;\pi_\tface^*\cE_X)
  \end{equation}
  is trivial for all $\theta\in[0,\pi]$, and the nullspace of $N_\tface(P_0,1)^*$ on
  \begin{equation}
  \label{EqSStfAdmA2}
    e^{-2 i/\rho_\sctface}\cA^{\alpha,-1+\beta}(\tface;\pi_\tface^*\cE_X)
  \end{equation}
  is trivial; here, we use the (Euclidean) density $\hat\rho^{-3}|\frac{\dd\hat\rho}{\hat\rho}\,\dd\slg|$ to define the adjoint.
\end{definition}

This assumption was verified for $\beta\in(0,1)$ for the scalar wave operator in \cite[Theorem~2.22]{HintzPrice}; see also \cite[Proposition~6.1]{HintzNonstat}. It also holds for the wave operator on any tensor bundle. If it holds for $\beta$ in an indicial gap $I\subset\R$, then it holds for all $\beta\in I$; it is thus an open condition in $\beta$. For details and examples, see~\S\ref{SsSptf}.

As in \cite{HintzNonstat}, we shall require mode stability at $\sigma\geq 0$, $\Im\sigma\geq 0$. Unlike in the reference, however, we now wish to allow for the presence of zero energy states; for example, the scalar wave operator has no zero energy states, the Hodge d'Alembertian on 1-forms has a one-dimensional space of bound states (and the resolvent $\wh{P_0}(\sigma)^{-1}$ has, in a certain sense, a simple pole at $\sigma=0$), and, as shown in \cite{HaefnerHintzVasyKerr,HaefnerHintzVasyKerrLarge}, the linearized gauge-fixed Einstein operator has a $7$-dimensional space of bound states (and the resolvent has a double pole at $\sigma=0$). It is difficult to state assumptions that are general enough to encompass all scenarios of interest while keeping the analysis of such general settings sufficiently transparent. We instead choose to directly require the validity of spacetime estimates for forward solutions of $P_0 u=f$. The relevant function spaces are 3b-Sobolev spaces with variable regularity order. The conditions on the regularity order (relative to a choice of weights at $\sface$ and $\cK^+\subset M_0$), ultimately arising from threshold conditions in radial point estimates (see~\S\ref{SsR3}), are as follows:

\begin{definition}[Stationary-$P_0$-admissible 3b-regularity orders]
\label{DefSSOrderAdm}
  We call a time-translation-in\-var\-i\-ant function $\sfs\in\CI(\Stb^*M_0)$ \emph{stationary-$P_0$-admissible with weights $\alpha_\sface,\alpha_\cK$ and margin $\eta\geq 0$} if it satisfies the following conditions.
  \begin{enumerate}
  \item\label{ItSSOrderAdm1}{\rm (Monotonicity.)} $\sfs$ is monotonically non-decreasing along the future null-bicharacteristic flow. (In the notation~\eqref{EqTse3bVF} and~\eqref{EqTsChar2}, this means $\pm\sfH_G\sfs\leq 0$ on $\pa\Sigma^\pm$.)
  \item\label{ItSSOrderAdm2}{\rm (Evenness.)} $\sfs$ is fiberwise even (i.e., pullback under $\zeta\mapsto-\zeta$ preserves $\sfs$).
  \item\label{ItSSOrderAdm3}{\rm (Constancy near the black hole and trapping.)} $\sfs$ is constant for $r\leq r_{\Gamma,+}$ where $r_{\Gamma,+}>r_+$ is such that $r<r_{\Gamma,+}$ over the trapped set.
  \item\label{ItSSOrderAdm4}{\rm (Threshold condition: event horizon.)} With respect to a positive definite fiber inner product on $\cE_X$, and in the coordinates used in~\eqref{EqTs3bHCoord}, define
    \begin{equation}
    \label{EqSSOrderAdm4}
      \vartheta_{\cH^+} := \sup_{N^*\cH^+\setminus o} \spec\Bigl[ \frac{\varrho^2}{2(r_+-\bhm)}\rho_\infty \upsigma^1\Bigl(\frac{P_0-P_0^*}{2 i}\Bigr)\Bigr],\quad \rho_\infty:=\xi_0^{-1}.
    \end{equation}
    (Here ``$\spec$'' denotes the spectrum of a self-adjoint map on a fiber of $\cE_X$.) Then, in a neighborhood of $\pa\cR^+_{\cH^+}$ (Definition~\usref{DefTs3bH}), or equivalently (by stationarity) in a neighborhood of $S N^*\cH^+\setminus o$, $\sfs$ is constant and satisfies
    \[
      \sfs>\frac12 + \vartheta_{\cH^+} + \eta.
    \]
  \item\label{ItSSOrderAdm5}{\rm (Threshold condition: outgoing radial set.)} With respect to a (possibly different) positive definite fiber inner product on $\cE_X$, and in the coordinates used in~\eqref{EqTs3bCCoord}, recall Definition~\ref{DefTs3bRad} and define
    \begin{equation}
    \label{EqSSOrderAdmOut0}
      \vartheta_{\pa\cK^+,{\rm out}} := \inf_{\cR^+_{\pa\cK^+,{\rm out}}} \spec\Bigl[-\rho_\infty\sigmatb^1\Bigl(r^2\frac{P_0-P_0^*}{2 i}\Bigr)\Bigr],\quad \rho_\infty:=-\sigma_\tbop^{-1}.
    \end{equation}
    Then, near $\pa\cR^+_{\pa\cK^+,{\rm out}}$, the function $\sfs$ is constant and satisfies
    \begin{equation}
    \label{EqSSOrderAdmOut}
      \sfs+\alpha_\sface<\frac12(-1+\vartheta_{\pa\cK^+,{\rm out}})+\alpha_\cK - \eta.
    \end{equation}
  \item\label{ItSSOrderAdm6}{\rm (Threshold condition: incoming radial set.)} With respect to a (possibly different) positive definite fiber inner product on $\cE_X$, and with $\rho_\infty=-\sigma_\tbop^{-1}$ as in~\eqref{EqSSOrderAdmOut0}, define
    \[
      \vartheta_{\pa\cK^+,{\rm in}} := \sup_{\cR^+_{\pa\cK^+,{\rm in}}} \spec\Bigl[\rho_\infty\sigmatb^1\Bigl(r^2\frac{P_0-P_0^*}{2 i}\Bigr)\Bigr].
    \]
    Then, near $\pa\cR^+_{\pa\cK^+,{\rm in}}$, the function $\sfs$ is constant and satisfies
    \begin{equation}
    \label{EqSSOrderAdmIn}
      \sfs+\alpha_\sface>\frac12(-1+\vartheta_{\pa\cK^+,{\rm in}})+\alpha_\cK + \eta.
    \end{equation}
  \end{enumerate}
\end{definition}

The threshold condition at the event horizon matches that in \cite[Proposition~4.9]{HintzGlueLocII}, while the threshold conditions at $\pa\cR^+_{\pa\cK^+,{\rm in/out}}$ match those of \cite[Lemma~5.8]{HintzNonstat} (except we phrase the definitions of $\vartheta_{\pa\cK^+,{\rm in/out}}$ more transparently here). The optimal, i.e., smallest possible, value of $\vartheta_{\pa\cK^+,{\rm out}}$ is easily computed (see Lemma~\ref{LemmaSSOrderThr} below). The evenness assumption is made for simplicity (so we do not need to specify conditions at radial sets in the past characteristic set). The margin $\eta\geq 0$ encodes a buffer in the threshold conditions which will be useful for technical purposes (when combining normal operator estimates with global regularity estimates).

\begin{lemma}[Threshold quantity]
\label{LemmaSSOrderThr}
  Recalling $S\in\CI(X;\End(\cE_X))$ from~\eqref{EqSSAdmOp}--\eqref{EqSSAdmOpPieces} and $\ubar S$ from~\eqref{EqSSAdmubarS}.
  For all $\eps>0$, there exists a fiber inner product on $\cE_X$ such that
  \[
    2\ubar S-\eps \leq \vartheta_{\pa\cK^+,{\rm out}} \leq 2\ubar S.
  \]
\end{lemma}
\begin{proof}
  This is related to computations around \cite[equation~(5.21)]{HintzNonstat}. Consider the form~\eqref{EqSSAdmOp} of $P_0$: the terms $Q\pa_{t_*},g^{0 0}\pa_{t_*}^2\in r^{-4}\Diff_\tbop^2$ have more decay than $r^{-2}$ and thus do not contribute to $\vartheta_{\pa\cK^+,{\rm out}}$. Next, $r^2(\wh{P_0}(0)^*-\wh{P_0}(0))$ is of class $\Diffb^1(X)\subset\Diff_\tbop^1(M_0)$, and thus its principal symbol, which is a linear function in the spatial momenta, vanishes at $\cR^+_{\pa\cK^+,{\rm out}}$ where the spatial momenta are zero (see Definition~\ref{DefTs3bRad}). Next, $-2\pa_{t_*}\rho(\rho\pa_\rho-1)\in r^{-2}\Diff_\tbop^2$ is symmetric with respect to the Minkowskian volume density, and hence its imaginary part with respect to the Kerr volume density (also in the presence of bundles, cf.\ Remark~\ref{RmkSSAdmDer}) if of class $r^{-3}\Diff_\tbop^1$. It thus remains to compute, for $\rho_\infty=\frac{1}{-\sigma_\tbop}$,
  \[
    -\rho_\infty\sigmatb^1\Bigl(\frac{1}{2 i}\bigl(2\rho\pa_{t_*} S - (2\rho\pa_{t_*} S)^*\bigr)\Bigr) = S+S^*.
  \]
  (See~\eqref{EqETPPs} for a related computation.) The claim now follows from $\inf_{\pa X}\spec(\frac12(S+S^*))\leq\ubar S$, valid for any fiber inner product on $\cE_X|_{\pa X}$, and the fact that for all $\eps>0$ there exists a fiber inner product such that $\ubar S-\frac{\eps}{2}\leq\inf_{\pa X}\spec(\frac12(S+S^*))$, as shown in \cite[Appendix~B]{HintzNonstat}.
\end{proof}

We recall from~\eqref{EqMUtbb} and Remark~\ref{RmkMUttstar} that the space
\begin{subequations}
\begin{equation}
\label{EqSSHtbb}
  H_{\tbop;\bop}^{(\sfs;k)}(\tilde M_0)
\end{equation}
consists of all $u\in\Htb^\sfs(\tilde M_0)$ (see~\eqref{EqMUKPsdo}) such that $\pa_{t_*}^{j_1}(t_*\pa_{t_*})^{j_2}(\la x\ra\pa_x)^\alpha u\in\Htb^\sfs(\tilde M_0)$ for all $j_1+j_2+|\alpha|\leq k$. The space $H_{\tbop;\bop}^{(\sfs;k)}(\Omega)^{\bullet,-}$ consists of restrictions to $\{r>\bhm\}$ of elements of $H_{\tbop;\bop}^{(\sfs;k)}(\tilde M_0)$ with support in $\{t_*\geq 1\}$. Since $|x|\geq\bhm$ and $t_*\geq 1$ on $\Omega$, it suffices to use $(t_*\pa_{t_*})^j(r\pa_x)^\alpha$, $j+|\alpha|\leq k$, to test for b-regularity here. On $\Omega$, the weights in the space
\begin{equation}
\label{EqSSHtbbOmega}
  H_{\tbop;\bop}^{(\sfs;k),(\alpha_\sface,\alpha_\cK)}(\Omega)^{\bullet,-}=\rho_\sface^{\alpha_\sface}\rho_\cK^{\alpha_\cK}H_{\tbop;\bop}^{(\sfs;k)}(\Omega)^{\bullet,-}
\end{equation}
\end{subequations}
can be taken to be powers of local (on $\Omega$) defining functions, such as $\rho_\sface=r^{-1}$ and $\rho_\cK=\frac{r}{t_*+r}$.

\begin{definition}[$\aleph$-admissibility]
\label{DefSSAlephAdm}
  Let $\aleph\geq 0$ and $\alpha_\sface\in\R$, $\delta\geq 0$. We then call a stationary wave-type operator $P_0$ (Definition~\usref{DefSSAdm}) \emph{$\aleph$-admissible with $\sface$-weight $\alpha_\sface$ and $\sface$-loss $\delta$} such the following holds for all stationary-$P_0$-admissible orders $\sfs$ with weights $\alpha_\sface,-\aleph$ and margin $0$ (Definition~\ref{DefSSOrderAdm}):
  \begin{enumerate}
  \item\label{ItSSAlephAdm1} $P_0$ is trapping admissible (Definition~\ref{DefSSTrapAdm});
  \item\label{ItSSAlephAdm2} $P_0$ is $\tface$-admissible with weight $\alpha_\sface+\frac32$ (Definition~\ref{DefSStfAdm});\footnote{In particular, $\alpha_\sface+\frac32$ lies in the indicial gap $(\beta^-,\beta^+)$ of $\wh{P_0}(0)$.}
  \item\label{ItSSAlephAdm3} for every $k\in\N_0$ and
    \[
      f \in H_{\tbop;\bop}^{(\sfs;k),(\alpha_\sface+2,0)}(\Omega;\cE)^{\bullet,-},
    \]
    the unique forward solution $u$ of $P_0 u=f$ satisfies
    \begin{equation}
    \label{EqSSAlephAdmSol}
    \begin{split}
      &u \in H_{\tbop;\bop}^{(\sfs;k),(\alpha_\sface-\delta,-\aleph)}(\Omega;\cE)^{\bullet,-}, \\
      &\qquad \|u\|_{H_{\tbop;\bop}^{(\sfs;k),(\alpha_\sface-\delta,-\aleph)}(\Omega;\cE)^{\bullet,-}} \leq C_k\|P_0 u\|_{H_{\tbop;\bop}^{(\sfs;k),(\alpha_\sface+2,0)}(\Omega;\cE)^{\bullet,-}}.
    \end{split}
    \end{equation}
    (Here we use the Minkowskian or Kerr density to define the underlying $L^2$-space; and we recall $\Omega$ from Definition~\usref{DefCMDomain}.)
  \end{enumerate}
\end{definition}

\begin{rmk}[Mild relaxation]
\label{RmkSSAlephAdmRelax}
  The following relaxed version of Definition~\ref{DefSSAlephAdm} is easier to verify in some settings (such as \citestab{Theorem~\ref*{ThmAdm}}) and sufficient for all of our purposes: we say that $P_0$ is $\aleph$-admissible with $\sface$-weight $\alpha_\sface$, $\sface$-loss $\delta$, \emph{and margin $d_\aleph\geq 0$} if for all stationary-$P_0$-admissible orders $\sfs$ with weights $\alpha_\sface,-\aleph$ \emph{and margin $d_\aleph$}, conditions~\eqref{ItSSAlephAdm1}--\eqref{ItSSAlephAdm2} hold as stated, while~\eqref{EqSSAlephAdmSol} is replaced by
  \[
    u \in H_{\tbop;\bop}^{(\sfs-d_\aleph;k),(\alpha_\sface-\delta,-\aleph)}(\Omega;\cE)^{\bullet,-},
  \]
  with norm bounded by $\|P_0 u\|_{H_{\tbop;\bop}^{(\sfs;k),(\alpha_\sface+2,0)}}$. --- The additional 3b-regularity loss here does not cause any conceptual difficulties: the only places where it enters is in Theorem~\ref{ThmNFw} below (see Remark~\ref{RmkNFwRelax}), and in the initial choice of 3b-regularity order in the proof of Theorem~\ref{ThmF} (see Remark~\ref{RmkFRelax}).
\end{rmk}

We stress that while 3b-regularity losses will be inconsequential for our analysis, it \emph{is} crucial that the $k$ degrees of b-regularity of $f$ get inherited \emph{without loss} by the forward solution $u$.

In concrete settings of interest to us, the membership~\eqref{EqSSAlephAdmSol} will be a consequence of appropriate resolvent estimates (using also, for the b-regularity part, the conormal regularity of the resolvent at $\sigma=0$).\footnote{The quantitative estimate in~\eqref{EqSSAlephAdmSol} follows from the stated membership for $u$ and the closed graph theorem for $P_0^{-1}$. But in practice, the stated membership is proved as a consequence of the quantitative estimate.} For example, ignoring the higher b-regularity requirement, \cite[Theorem~5.3]{HintzGlueLocII} proves the $0$-admissibility of the scalar wave operator (with $\sface$-weight $\alpha_\sface\in(-\frac32,-\frac12)$ and $\sface$-loss $0$); see Theorem~\ref{ThmA1Adm} for the full result. More generally, we show that mode stability, including at zero energy, implies $0$-admissibility with $\sface$-loss $0$ (Theorem~\ref{ThmA1Gen}).

The quantity $\aleph$ is the loss of decay at $\cK^+$ of the forward solution relative to the forcing term. On the resolvent side, $\aleph$ is, roughly speaking, the order of the pole of $\wh{P_0}(\sigma)^{-1}$ at $\sigma=0$. Allowing for a positive loss $\delta$ is necessary in the presence of zero energy bound states with limited decay; see Theorem~\ref{ThmA2Adm} for a concrete example in the case of the wave operator on 1-forms. As another example, we note that it follows from \cite[Theorem~3.22 and Remark~3.28]{HintzGlueLocIII} that a suitable version of the linearized gauge-fixed Einstein operator is $2$-admissible with $\sface$-weight $\alpha_\sface\in(-\frac32,-\frac12)$ and any $\sface$-loss $\delta>\alpha_\sface+\frac32$.\footnote{Moreover, \cite[Lemma~3.27]{HintzGlueLocIII} \emph{almost} gives the $2$-admissibility of the same operator, with $\sface$-weight $\alpha_\sface\in(-2,-\frac32)$ and $\sface$-loss $0$, except the $\tface$-admissibility fails for such $\alpha_\sface$. The reference in fact shows a stronger statement than~\eqref{EqSSAlephAdmSol} for $\alpha_\sface\in(-\frac32,-\frac12)$, namely one can weaken the $\sface$-decay order of $f$ to $\alpha_\sface+2-\delta$. There are many slight modifications of $\aleph$-admissibility one can consider, and we settle for the one in Definition~\ref{DefSSAlephAdm} since it involves the smallest number of parameters and is easiest to verify in practice.} The Einstein case (albeit for different gauge-fixing and constraint damping) will be discussed in the framework of the present paper in \cite{HintzKerrStab}.

\subsection{Dynamical spacetimes}
\label{SsSD}

We describe the classes of asymptotically Kerr metrics (\S\ref{SssSDG}), and wave-type operators relative to those (\S\ref{SssSDW}), only in the region $\Omega=\ol{\{t_*\geq 1\}}\subset M$. (We recall that the analysis of waves in the ``exterior region,'' where $0\leq t\leq r+1$, was already accomplished in~\cite{HintzVasyScrieb} except for tame estimates---which were later proved in two settings for the Einstein equations in \cite{HintzMink4Gauge,HintzKerrStab}. The analysis there can be done entirely using simple energy estimates and commutator arguments; all novelties in the present paper happen only in $\Omega$.)

\subsubsection{Geometry: asymptotically Kerr metrics}
\label{SssSDG}

In our applications, we exclusively work with small and suitably decaying (in space and time) perturbations of the Kerr metric.\footnote{Large perturbations, as, e.g., allowed for in \cite[Definition~3.22]{HintzNonstat}, can be handled simply by realizing that at sufficiently late times, and relaxing the decay orders slightly, they are small perturbations after all; and on the remaining initial finite time interval, energy estimates as proved in~\S\ref{SsET} control forward solutions of wave-type equations.} We thus introduce:

\begin{definition}[Asymptotically Kerr metrics]
\label{DefSDGMetric}
  Let $\ell_\sscri\in[0,\frac12]$, $\ell_+\in[0,1]$, $\ell_\cK\geq 0$. Let $d_0,k\in\N$. With subextremal Kerr parameters $\bhm,a$ fixed, and recalling Definitions~\usref{DefCTe3b} and~\usref{DefMUCe3b}, we then define
  \begin{align*}
    &\sG_{\etbop;\bop}^{(d_0;k),(2\ell_\sscri,\ell_+,\ell_\cK)} := \bigl\{ g_{\bhm,a}+h \colon h\in\tilde\sG_{\etbop;\bop}^{(d_0;k),(2\ell_\sscri,\ell_+,\ell_\cK)} \bigr\}, \\
    &\qquad \tilde\sG_{\etbop;\bop}^{(d_0;k),(2\ell_\sscri,\ell_+,\ell_\cK)} := x_\sscri^{-2}\rho_+^{-2}\cC_{\etbop;\bop}^{(d_0;k),(2\ell_\sscri,\ell_+,\ell_\cK)}(\Omega;S^2\,\Tetb_\Omega^*M).
  \end{align*}
  We write these spaces as $\sG_{\etbop;\bop}^{(d_0;k),(2\ell_\sscri,\ell_+,\ell_\cK)}(\Omega)$ and $\tilde\sG_{\etbop;\bop}^{(d_0;k),(2\ell_\sscri,\ell_+,\ell_\cK)}(\Omega)$ when we need to make the domain $\Omega$ explicit. We equip $\tilde\sG_{\etbop;\bop}^{(d_0;k),(2\ell_\sscri,\ell_+,\ell_\cK)}$ with the norm of $\cC_{\etbop;\bop}^{(d_0;k),(2\ell_\sscri+2,\ell_++2,\ell_\cK)}$, and the affine space $\sG_{\etbop;\bop}^{(d_0;k),(2\ell_\sscri,\ell_+,\ell_\cK)}$ with the induced metric. For $d_0=0$, resp.\ $k=0$, we omit the subscript ``$\etbop$,'' resp.\ ``$\bop$.''
\end{definition}

Lemma~\ref{LemmaTsKLMetric}\eqref{ItTsKLMetrice3b} implies that if $h\in\tilde\sG_\etbop^{0,(0,0,0)}$ is sufficiently small, then $g_{\bhm,a}+h$ is a non-degenerate time-oriented Lorentzian metric on $\Omega^\circ$, and $x_\sscri^{-2}\rho_+^{-2}(g_{\bhm,a}+h)$ is a non-degenerate Lorentzian signature section of $S^2\,\Tetb^*_\Omega M$. The requirements $2\ell_\sscri,\ell_+\leq 1$ are made for convenience: they imply that smooth (in $x_\sscri,\rho_+$) perturbations of $g_{\bhm,a}$ as $r\to\infty$ are allowed (so one can, e.g., replace $g_{\bhm,a}$ by $\ubar g$ away from $\cK^+$ up to acceptable error terms). For the decay rate $\ell_\cK$ towards stationarity on the other hand, there are no such natural restrictions.

\begin{lemma}[Time function]
\label{LemmaSDGTime}
  Let $\ell_\sscri\in(0,\frac12]$. Then there exist $\eps_\sG>0$ and a function $\ft_*\in\CI(\Omega\setminus(\scri^+\cup\iota^+\cup\cK^+))$ with
  \begin{enumerate}[label={\rm (\alph*)}]
  \item $\ft_*=t_*(1-(\frac{t_*}{r})^{\ell_\sscri})$ for all sufficiently small $\frac{t_*}{r}$,
  \item $\delta t_*\leq\ft_*\leq t_*$ for some $\delta>0$,
  \myitem{ItSDGTime3}{c} $\ft_*=c_0 t_*+c_1$ for $r\leq 2\bhm$, for some $c_0,c_1>0$,
  \end{enumerate}
  such that the following statements hold for all $g=g_{\bhm,a}+h$ where $h\in\tilde\sG_\etbop^{0,(2\ell_\sscri,0,0)}$ has norm $<\eps_\sG$:
  \begin{enumerate}
  \item\label{ItSDGTimeR} $r$ has future timelike differential on $\{\bhm\leq r\leq r_\natural\}$ where $r_\natural:=\frac{\bhm+r_+}{2}$;
  \item\label{ItSDGTimeT} $\ft_*$ has past timelike differential on $\Omega\setminus(\scri^+\cup\iota^+\cup\cK^+)$.
  \end{enumerate}
\end{lemma}

\begin{figure}[!ht]
\centering
\includegraphics{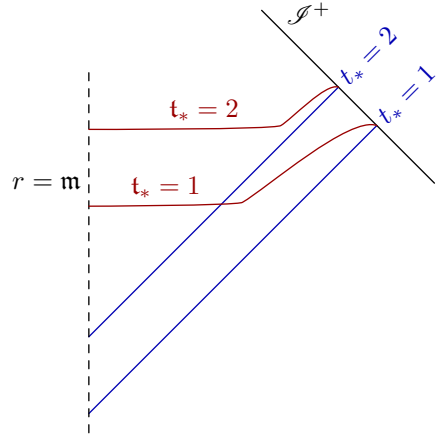}
\caption{Illustration of Lemma~\ref{LemmaSDGTime}, showing level sets of $t_*$ (in blue) and $\ft_*$ (in red).}
\label{FigSDGTime}
\end{figure}

See Figure~\ref{FigSDGTime}.

\begin{rmk}[Other Kerr parameters]
\label{RmkSDGOtherKerr}
  Since $g_{\bhm,a}-\ubar g\in\tilde\sG_\etbop^{0,(1,0,0)}$ by~\eqref{EqTsKLMetrice3b} and thus also $g_{\bhm,a}-g_{\bhm',a'}$, the timelike nature of $\dd\ft_*$ persists also for (perturbations of) nearby Kerr metrics.
\end{rmk}

\begin{proof}[Proof of Lemma~\usref{LemmaSDGTime}]
  For part~\eqref{ItSDGTimeR}, we note that $\sup_{\bhm\leq r\leq r_\natural} g_{\bhm,a}^{-1}(\dd r,\dd r)=:c<0$. For $r\leq r_\natural$, we furthermore have $g^{-1}(\dd r,\dd r)\leq g_{\bhm,a}^{-1}(\dd r,\dd r)+C\eps_\sG$, so $\eps_\sG<|c|/C$ works.

  For part~\eqref{ItSDGTimeT}, we recall Lemma~\ref{LemmaCTebsc}, which implies that the coefficients of $g^{-1}-g_{\bhm,a}^{-1}$ in the basis
  \[
    \rho_\sscri\pa_{t_*}^2,\ \pa_{t_*}\otimes_s\pa_r,\ x_\sscri\pa_{t_*}\otimes r^{-1}\pa_\omega,\ \rho_\sscri^{-1}\pa_r^2,\ x_\sscri^{-1}\pa_r\otimes_s r^{-1}\pa_\omega,\ r^{-1}\pa_\omega\otimes_s r^{-1}\pa_\omega,
  \]
  where $\rho_\sscri=x_\sscri^2$, are bounded in absolute value by $C\eps_\sG\rho_\sscri^{\ell_\sscri}$; by~\eqref{EqTsKLMetrice3b}, this holds also for $g^{-1}-\ubar g^{-1}$, and we recall $\ubar g^{-1}=-2\pa_{t_*}\otimes_s\pa_r+\pa_r^2+r^{-2}\slg$ from~\eqref{EqCTebMinkTstar}. We first work in $r\geq C t_*$ (so $t\leq\frac{C+1}{C}r$) for some large $C$ (to be fixed momentarily), where we can take $\rho_\sscri=\frac{t_*}{r}$; then the inner product of
  \[
    \dd\Bigl[t_*\Bigl(1-\Bigl(\frac{t_*}{r}\Bigr)^{\ell_\sscri}\Bigr)\Bigr] = \bigl(1-(1+\ell_\sscri)\rho_\sscri^{\ell_\sscri}\bigr)\,\dd t_* + \ell_\sscri\rho_\sscri^{\ell_\sscri+1}\,\dd r
  \]
  with itself with respect to $g^{-1}$ is equal to
  \[
    -2\ell_\sscri\bigl(1-(1+\ell_\sscri)\rho_\sscri^{\ell_\sscri}\bigr)\rho_\sscri^{\ell_\sscri+1} + \cO(\rho_\sscri^{2\ell_\sscri+2}) + C\eps_\sG\rho_\sscri^{\ell_\sscri}\bigl(\rho_\sscri\cdot\cO(1) + 1\cdot\cO(\rho_\sscri^{\ell_\sscri+1}) + \rho_\sscri^{-1}\cdot\cO(\rho_\sscri^{2(\ell_\sscri+1)})\bigr)
  \]
  and thus $\leq -\ell_\sscri\rho_\sscri^{\ell_\sscri+1}<0$ when $\eps_\sG,\rho_\sscri\leq C^{-1}$ if $C$ is large enough.

  At $\rho_\sscri=C^{-1}$, we have $r=C t_*$ and thus
  \[
    t_*\Bigl(1-\Bigl(\frac{t_*}{r}\Bigr)^{\ell_\sscri}\Bigr) = C'(t_*+r)\quad\text{at}\ r=C t_*\ \text{for}\ C':=\frac{1-C^{-\ell_\sscri}}{1+C}.
  \]
  Now $t_*+r=t$ has past timelike differential for the Minkowski metric, and thus
  \[
    g^{-1}(\dd t,\dd t) \leq -1+\cO(r^{-1}) + C\eps_\sG \leq -\frac12
  \]
  for $r_0\leq r\leq C t_*$ when $r_0$ is sufficiently large and $\eps_\sG$ is small; here we use that in the region $r\leq C t_*$ (and thus away from $\scri^+$), the metric perturbation $h$ has coefficients of size $\cO(\eps_\sG)$ in the basis given by symmetric tensor products of $\dd t$, $\dd r$, $r\,\dd\omega$.

  Finally, at $r=r_0$ we have $C'(t_*+r)=C'(t_*+r_0)$; and in the region $r\leq r_0$, the differential $\dd(C'(t_*+r_0))=C'\,\dd t_*$ is timelike for $g_{\bhm,a}$ and thus also for $g$ when $\eps_\sG$ is sufficiently small. We may then take $\ft_*$ to be a smoothing of
  \[
    \begin{cases}
      C'(t_*+r_0), & r\leq r_0, \\
      C'(t_*+r), & r_0\leq r\leq C t_*, \\
      t_*\bigl(1-\left(\frac{t_*}{r}\right)^{\ell_\sscri}\bigr), & r\geq C t_*.
    \end{cases}
  \]
  This completes the construction.
\end{proof}

\subsubsection{Wave-type operators}
\label{SssSDW}

The following definition of the class of wave-type operators we shall study is modeled on \cite[Definition~3.28]{HintzNonstat} (and on \cite[Definition~6.2]{HintzVasyScrieb} near $\scri^+$). Besides the fact that in the present paper we choose to work only on $\Omega$ and need to consider operators whose coefficients have limited regularity (matching our analysis below, we assume mixed $(\etbop;\bop)$-regularity), the main new feature is that we need to impose lower bounds on the decay rates towards stationarity (see~\eqref{EqSDWAdmWeights}) to balance the loss of decay of forward solutions of the stationary model in Definition~\ref{DefSSAlephAdm}:

\begin{definition}[Admissible wave-type operators]
\label{DefSDWAdm}
  Let $\aleph,\delta\geq 0$. Let $\ell_\sscri,\ell_+,\ell_\cK>0$ be such that
  \begin{equation}
  \label{EqSDWAdmWeights}
    \ell_\sscri\in(0,\tfrac12],\quad
    \ell_+>\delta,\quad
    \ell_\cK>\aleph.
  \end{equation}
  Let $P_0$ be an $\aleph$-admissible stationary wave-type operator with $\sface$-weight $\alpha_\sface$ and $\sface$-loss $\delta$ (Definition~\usref{DefSSAlephAdm}), acting on section of the bundle $\cE=\pi_X^*\cE_X\to M$ (where $\cE_X\to X$ is a bundle and $\pi_X\colon M\to X$ the spatial projection). We then call an operator $P\in\Diff^2(\Omega^\circ;\cE)$ \emph{admissible (relative to $P_0$) of class $((d_0;k),(2\ell_\sscri,\ell_+,\ell_\cK))$} if the following conditions hold.
  \begin{enumerate}
  \item The principal symbol of $P$ is scalar, and equal to the dual metric function of an asymptotically Kerr metric $g=g_{\bhm,a}+h\in\sG_{\etbop;\bop}^{(d_0;k),(2\ell_\sscri,\ell_+,\ell_\cK)}$ (Definition~\usref{DefSDGMetric}).
  \item\label{ItSDWAdm2} Recalling the coordinates~\eqref{EqCMCoordscriip} near $\scri^+\cap\{\rho_\cK\leq 1\}$, fix a cutoff function $\chi_\sscri=\chi_\sscri(x_\sscri)$. Then there exist
    \begin{equation}
    \label{EqSDWp0p1}
      p_0 \in \rho_+^{\ell_+}\cC_\bop^{k+d_0}(\scri^+;\End(\cE)),\quad
      \tilde p_1 \in \rho_+^{\ell_+}\cC_\bop^{k+d_0}(\scri^+;\End(\cE))
    \end{equation}
    such that (defining derivatives using a connection constructed as in~\eqref{EqSSAdmDer})
    \begin{equation}
    \label{EqSDWAdmOp}
    \begin{split}
      \tilde P&:=P - \Bigl[P_0 + \chi_\sscri \frac12 x_\sscri^2\rho_+^2 \bigl( 2\tilde p_1(x_\sscri\pa_{x_\sscri}-2\rho_+\pa_{\rho_+}) + p_0 \bigr)\Bigr] \\
      &\hspace{8em} \in x_\sscri^2\rho_+^2\cC_{\etbop;\bop}^{(d_0;k),(2\ell_\sscri,\ell_+,\ell_\cK)}\Diff_\etbop^2(\Omega;\cE).
    \end{split}
    \end{equation}
  \item\label{ItSDWAdmSplit} There exists a bundle splitting $\cE_X|_{\pa X}=\bigoplus_{j=1}^J\cE_{X,j}$ in which $p_1:=S|_{\pa X}+\tilde p_1$ (in the notation of~\eqref{EqSSAdmOp}) is lower triangular, with the diagonal entries $p_{1,j j}\in\End(\cE_{X,j})$ having real spectrum; and $p_0$ is strictly lower triangular.
  \end{enumerate}
  We call $P$ \emph{weakly admissible} if the decay orders $\ell_\sscri,\ell_+,\ell_\cK$ only satisfy
  \begin{equation}
  \label{EqSDWAdmWeak}
    \ell_\sscri\in(0,\tfrac12],\quad
    \ell_+>0,\quad
    \ell_\cK>0
  \end{equation}
  instead of the stronger~\eqref{EqSDWAdmWeights}.
\end{definition}

\begin{rmk}[Motivation]
\label{RmkSDWAdmMotiv}
  The requirement~\eqref{EqSDWAdmOp} goes back to \cite[Definition~3.7]{HintzVasyScrieb}. See \cite[\S{3.2}]{HintzVasyScrieb} for examples. In particular, the term $\tilde p_1$ determines decay rates of different components of a wave at $\scri^+$; this arises in the context of the (linearized) gauge-fixed Einstein field equations in \cite[Proposition~3.29]{HintzMink4Gauge}. The strictly lower triangular nature of $p_0$ allows for studying $\scri^+$-decay in a hierarchical manner. See \cite[Remark~6.3]{HintzVasyScrieb} for further discussion of part~\eqref{ItSDWAdmSplit}.
\end{rmk}

The upper bound on $\ell_\sscri$ ensures that different choices of local coordinates preserve the validity of~\eqref{EqSDWAdmOp}. In view of~\eqref{EqSSOpMemeb}, it follows from~\eqref{EqSDWAdmOp} that
\begin{equation}
\label{EqSDWAdmNearScri}
  2 x_\sscri^{-2}\rho_+^{-2}P - \Bigl[-\bigl(x_\sscri\pa_{x_\sscri} - 2(1+p_1)\bigr)(x_\sscri\pa_{x_\sscri}-2\rho_+\pa_{\rho_+}) + 2 x_\sscri^2\slDelta + p_0\Bigr] \in \cC_{\ebop;\bop}^{(d_0;k),(2\ell_\sscri,0)}\Diffeb^2.
\end{equation}
Thus, while on the principal symbol level $P$ is a decaying perturbation of $P_0$ at $\scri^+$, $\iota^+$, and $\cK^+$, we do permit $P$ to have additional subprincipal terms that do not decay at $\scri^+$ (relative to the overall decay rate $x_\sscri^2$ of $P$ as an edge-b-differential operator there). The characteristic exponents of $x_\sscri\pa_{x_\sscri}-2(1+p_1)$ in~\eqref{EqSDWAdmNearScri} are $2(1+\mu)$ where $\mu$ are the eigenvalues of $p_1$. We thus introduce:

\begin{definition}[The quantity $\ubar p_1$]
\label{DefSDWubarp1}
  Setting $p_1:=S|_{\pa X}+\tilde p_1$, we define
  \[
    \ubar p_1 := \inf_{\substack{p\in\scri^+ \\ \mu\in\spec p_1(p)}} \Re\mu.
  \]
\end{definition}

\begin{rmk}[Relaxed orders]
\label{RmkSDWRelax}
  The lower bounds on $\ell_+$ and $\ell_\cK$ in~\eqref{EqSDWAdmWeights} are only important for closing estimates for $P$ by using the mapping properties of the stationary model $P_0$. At prior steps of our analysis (in particular, for e3b-regularity estimates), it suffices to assume only $\ell_+,\ell_\cK>0$ as in~\eqref{EqSDWAdmWeak}. When proving estimates on sets disjoint from $\iota^+,\cK^+$ (resp.\ $\scri^+,\iota^+$), the orders $\ell_+,\ell_\cK$ (resp.\ $\ell_\sscri,\ell_+$) are irrelevant altogether, and we shall thus take them to be $0$ in such cases.
\end{rmk}

\begin{rmk}[Regularity of $p_0,\tilde p_1$]
\label{rmkSDWReg}
  In~\eqref{EqSDWAdmOp}, we implicitly extend $p_0,\tilde p_1$ to $\supp\chi_\sscri$ by requiring them to be independent of $x_\sscri$. (To make sense of this in the presence of $\cE$, one fixes any bundle isomorphism of $\cE$ in the chart with the pullback of $\cE|_{\scri^+}$ along the map $(x_\sscri,\rho_+,\omega)\mapsto(0,\rho_+,\omega)$.) Thus, the terms involving $p_0,\tilde p_1$ in~\eqref{EqSDWAdmOp} are, a fortiori, of class $x_\sscri^2\rho_+^2\cC_{\etbop;\bop}^{(d_0;k),(0,\ell_+,\infty)}$, and different choices of extensions lead to errors in the space~\eqref{EqSDWAdmOp}. By contrast, since the restriction of $\cV_\etbop(M)$ to $\scri^+$ is spanned (over $\CI(\scri^+)$) by the generator $\rho_+\pa_{\rho_+}$ of dilations of the fibers of $\scri^+\setminus I^0$, it would be more natural to assume $p_0,\tilde p_1$ to only have $k$ degrees of b-regularity and $d_0$ additional degrees of regularity in $\rho_+\pa_{\rho_+}$. If, under this weaker assumption, one extends $p_0,\tilde p_1$ to a level set $x_\sscri=s>0$ using the time $s$ heat flow on $\Sph^2$ (with parametric dependence on $\rho_+$), one would obtain extensions with $(\ebop;\bop)$-regularity $(d_0;k)$, but at the cost that the weight $2\ell_\sscri$ in the error space in~\eqref{EqSDWAdmOp} would need to be replaced by mere vanishing at $\scri^+$. All of our microlocal and energy estimates near $\scri^+$ go through with such weaker error terms as well. But since, ultimately, we will be rather wasteful with b-derivatives anyway, we content ourselves with the stronger (but simpler to work with) condition~\eqref{EqSDWp0p1} in this paper.
\end{rmk}

We measure the size of perturbations of a given stationary wave-type operator $P_0$ as follows:

\begin{definition}[Norm]
\label{DefSDWAdmNorm}
  In the notation of Definition~\usref{DefSDWAdm}, we write $\|P-P_0\|_{(d_0;k),(2\ell_\sscri,\ell_+,\ell_\cK)}$ for the sum of the norms of $p_0$, $\tilde p_1$ in the spaces~\eqref{EqSDWp0p1} and of the operator $\tilde P$ in~\eqref{EqSDWAdmOp}; these are thus supremum norms of derivatives of the coefficients of $P-P_0$. If we take the suprema only over a subset $\Omega'\subset\Omega$, we denote the resulting seminorm by $\|P-P_0\|_{(d_0;k),(2\ell_\sscri,\ell_+,\ell_\cK),\Omega'}$.
\end{definition}

Since the principal symbol of $P=P_0+(P-P_0)$ is the dual metric function of $g$, this norm controls also the perturbation of the metric, in that
\[
  \|g-g_{\bhm,a}\|_{\tilde\sG_{\etbop;\bop}^{(d_0;k),(2\ell_\sscri,\ell_+,\ell_\cK)}} \leq C\|P-P_0\|_{(d_0;k),(2\ell_\sscri,\ell_+,\ell_\cK)},
\]
similarly on subsets $\Omega'\subset\Omega$.

\section{Phase space dynamics of asymptotically Kerr metrics}
\label{SDy}

Fix subextremal Kerr parameters $\bhm>0$, $a\in(-\bhm,\bhm)$, and weights
\begin{equation}
\label{EqDyWeights}
  \ell_\sscri \in (0,\tfrac12],\quad
  \ell_+\in(0,1],\quad
  \ell_\cK>0.
\end{equation}
Recall the time function $\ft_*$ from Lemma~\ref{LemmaSDGTime}, and let
\begin{equation}
\label{EqDyDomain}
   \Omega_* := \cl_M\{\ft_*\geq 1\} \subset \Omega.
\end{equation}
In this section, we will describe aspects of the null-geodesic flow of metrics
\begin{equation}
\label{EqDyMetric}
  g = g_{\bhm,a} + h,\quad h\in\tilde\sG_\etbop^{d_0,(2\ell_\sscri,\ell_+,\ell_\cK)}(\Omega_*),\ \|h\|_{\tilde\sG_\etbop^{d_0,(2\ell_\sscri,\ell_+,\ell_\cK)}(\Omega_*)} < \eps_\sG,
\end{equation}
where $2\leq d_0\in\N$ is fixed and $\eps_\sG>0$ is small enough so that the conclusions of Lemma~\ref{LemmaSDGTime} hold. (See Definition~\ref{DefSDGMetric} for the notation.) We denote by $G(\zeta)=g^{-1}(\zeta,\zeta)$ the dual metric function; moreover, fixing defining functions
\begin{subequations}
\begin{equation}
\label{EqDyRhos}
  \rho_\infty \in \CI(\ol{{}^\etbop T^*}M),\quad
  x_\sscri,\,\rho_+\in\CI(M),
\end{equation}
of fiber infinity ${}^\etbop S^*M$ and $\scri^+$, $\iota^+\subset M$, respectively, we consider the rescaled Hamiltonian vector field
\begin{equation}
\label{EqDyHam}
  \sfH_G := \rho_\infty x_\sscri^{-2}\rho_+^{-2}H_G
\end{equation}
\end{subequations}
on the characteristic set $G^{-1}(0)$. We write
\begin{equation}
\label{EqDyChar}
  \Sigma := ({}^\etbop T^*_{\Omega_*}M \setminus o) \cap G^{-1}(0) = \Sigma^+\sqcup\Sigma^-,\quad
  \pa\Sigma^\pm \subset {}^\etbop S^*M,
\end{equation}
where $(z,\zeta)\in\Sigma$ lies in the future characteristic set $\Sigma^+$ if and only if $g^{-1}(\zeta,\dd\ft_*)\leq 0$; and $\Sigma^-=-\Sigma^+$ (fiber-wise scaling by $-1$).

Since the metric perturbation $h$ decays (as a section of $S^2\,{}^\etbop T^*M$) relative to $g_{\bhm,a}$, the rescaled Hamiltonian vector field $\sfH_G$ differs from its Kerr version $\sfH_{G_{\bhm,a}}$ (where $G_{\bhm,a}(\zeta)=g_{\bhm,a}^{-1}(\zeta,\zeta)$) by a continuous b-vector field on ${}^\etbop S^*_{\Omega_*}M$ that vanishes (as a b-vector field) on the e3b-cosphere bundle over the boundary hypersurfaces $\scri^+$, $\iota^+$, $\cK^+$ of $\Omega_*$ at infinity. Thus, the $\sfH_G$-flow and $\sfH_{G_{\bhm,a}}$-flow are identical over $\Omega_*\cap\pa\tilde M$.

\begin{prop}[Global null-bicharacteristic dynamics on an asymptotically Kerr spacetime]
\label{PropDy}
  Let $\gamma\colon I\to{}^\etbop S^*_{\Omega_*}M$ denote a maximally extended integral curve of $\sfH_G$ in $\pa\Sigma^+\cap{}^\etbop S^*_{\Omega^*}M$. Then one of the following (not mutually exclusive) possibilities must occur:
  \begin{itemize}
  \item $\gamma$ lies over $(\Omega_*\cap\scri^+)\cup\iota^+\cup\cK^+$, in which case it is an integral curve of $H_{G_{\bhm,a}}$ (which is then described by Proposition~\usref{PropTse3bDyn});
  \item $\gamma$ crosses $\ol{\{\ft_*=1\}}$ in the backward direction, while in the forward direction either $\gamma$ crosses $\ol{\{r=\bhm\}}$ or $\gamma$ tends to the e3b-cosphere bundle over $(\Omega^*\cap\scri^+)\cup\iota^+\cup\cK^+$.
  \end{itemize}
\end{prop}
\begin{proof}
  When $\gamma$ lies over $M^\circ$, then since $\ft_*$ is a time function for $(\Omega_*,g)$, we have $\sfH_G\ft_*>0$ on $\pa\Sigma^+\cap S^*\tilde M^\circ$; in fact, $\sfH_G\ft_*$ is bounded from below by a positive constant over every fixed compact subset of $\pa\Sigma^*\cap S^*\tilde M^\circ$. If $\ft_*(\gamma(s))$ is unbounded as $s\to\sup I$, then $\gamma$ tends to the set $\ft_*^{-1}(\infty)={}^\etbop S^*_{\iota^+\cup\cK^+}M$. If $\ft_*(\gamma(s))$ remains bounded, the accumulation points of $\gamma(s)$ as $s\to\sup I$ must either lie over $r^{-1}(\bhm)$ or $\scri^+$. In the former case, since $\sfH_G r<0$ on $\pa\Sigma^+\cap S^*_{\Omega_*\cap\{r=\bhm\}}M^\circ$, we conclude that $\gamma$ crosses $\{r=\bhm\}$ in the forward direction.
\end{proof}

It is tempting to formulate the exact analogue of Proposition~\ref{PropTse3bDyn} for the $\sfH_G$-flow, which would thus entail a more precise description of the forward asymptotics of null-bicharacteristics over $M^\circ$ (see the first row of~\eqref{EqTse3bDynTable}). While such an analogue can, indeed, be proved with some effort (requiring, e.g., proofs of convexity properties for $g$ such as Lemma~\ref{LemmaTs3bDyn}\eqref{ItTs3bDynKp}), we do not need such a strong result for our purely \emph{analytic} purposes: roughly speaking, the microlocal regularity results we will employ at the various invariant sets for $\sfH_G=\sfH_{G_{\bhm,a}}$ over $\pa M$ give control of solutions of the PDE $P u=f$ in \emph{open neighborhoods} of these invariant sets and thus, ultimately, in an open neighborhood of $\Omega_*\cap\pa M$; Proposition~\ref{PropDy} then implies that regularity in the complement of a slightly smaller neighborhood can be obtained by propagation starting from the initial hypersurface $\ft_*^{-1}(1)$. (This is an instance of the perturbation stability already pointed out in \cite[\S{2.7}]{VasyMicroKerrdS}.)

\begin{rmk}[Dynamical horizon]
\label{RmkDyHor}
  The event horizon of $(\Omega_*,g)$ is a null hypersurface that is graphical over the event horizon $\cH^+\cap\{\ft_*\geq 1\}$, where $\cH^+=r^{-1}(r_+)$, as shown for metrics of the class studied here in \cite{ChenKlainermanHorizon,HintzHorizons}.
\end{rmk}

\subsection{Trapping}
\label{SsDyTr}

The only place where we need a more precise description of the dynamics of $\sfH_G$ is the trapped set. We will be working near a compact subset of $(\cK^+)^\circ$ and thus away from $\scri^+$ and $\iota^+$, and hence we are not concerned with edge- or b- or 3-body-features: a local spanning set of $\cV_\etbop(M)$ is given by the vector fields $\pa_{t_*}$, $\pa_x$, or, in terms of the local defining function $\rho_\cK=t_*^{-1}$ of $\cK^+$, the vector fields $\rho_\cK^2\pa_{\rho_\cK}$, $\pa_x$. These are the \emph{cusp vector fields} in the terminology of Mazzeo--Melrose \cite{MazzeoMelroseFibred}. Correspondingly (and to match the terminology used in \cite[\S{4}]{HintzPolyTrap}), for submanifolds $D\subset\tilde M\setminus(\scri^+\cup\iota^+)$, we drop weights at $\scri^+,\iota^+$ from the notation and shall write
\begin{subequations}
\begin{equation}
\label{EqDyTrCuObj}
  {}^\cuop T^*D := {}^\etbop T^*_D\tilde M,\quad
  \bar H_\cuop^{s,\alpha_\cK}(D) = \rho_\cK^{\alpha_\cK}\bar H_\cuop^s(D) := \bar H_\etbop^{s,(0,0,\alpha_\cK)}(D)
\end{equation}
for phase space and Sobolev spaces with regularity order $s\in\R$ and decay order $\alpha_\cK\in\R$, and
\begin{equation}
\label{EqDyTrCuObj2}
  \cC_\cuop^{d,\alpha_\cK}(D) = \rho_\cK^{\alpha_\cK}\cC_\cuop^d(D),
\end{equation}
\end{subequations}
$d\in\N_0$, for the space of $\cC^d(D^\circ)$-functions whose up to $d$-fold derivatives along $\pa_{t_*}$ and $\pa_x$ lie in $\rho_\cK^{\alpha_\cK}L^\infty(D)$. (For $D=\R_t\times\cU$ where $\cU\subset\R^3$ is open and precompact, the space $\cC_\cuop^{d,\alpha_\cK}(D)$ is the same as the space $\rho\cC_b^d(D)$, $\rho:=\rho_\cK^{\alpha_\cK}$, in the notation of \cite[\S{2A}]{HintzPolyTrap}, and the norm of elements of $\bar H_\cuop^{s,\alpha_\cK}(D)$ with support (on $\R^4$) contained in $\R_t\times K$, $K\subset\cU$ compact, is equivalent to the norm on the space $H_\cuop^{s,\alpha_\cK}$ in the notation of \cite[\S{3B1}]{HintzPolyTrap}.) We shall also use the notation~\eqref{EqDyTrCuObj2} for subsets $D\subset{}^\cuop S^*D$ where ${}^\cuop S^*D$ is the boundary at fiber infinity of the fiber-wise radial compactification $\ol{{}^\cuop T^*}D$, where the notion of cusp-regularity now also entails regularity in each fiber.

Recall the unstable/stable trapped set $\bar\Gamma^{\rm u/s}\subset\Stb^*_{\cK^+}M\cap\pa\Sigma^+$ for the Kerr metric $g_{\bhm,a}$ from~\eqref{EqTs3bGamma}. Recall also the defining functions $\phi_0^{\rm u/s}\in S^1_{\rm hom}(T^*\cM\setminus o)$ of the stationary extensions $\Gamma^{\rm u/s}_0=\R_\ft\times\bar\Gamma^{\rm u/s}\subset T^*\cM\setminus o$ inside of the future characteristic set
\[
  \Sigma_{\bhm,a}^+ \subset \Sigma_{\bhm,a}:=G_{\bhm,a}^{-1}(0)\cap(T^*\cM\setminus o)
\]
of the exact Kerr metric from Definition~\ref{DefTs3bODefFn}. Recall $r_{\Gamma,\pm}\in(r_+,\infty)$ from~\eqref{EqTs3bOTrap0Cpt}. Set
\begin{equation}
\label{EqDyTrD}
  D := \{ \rho_\cK<1,\ r_{\Gamma,-}<r<r_{\Gamma,+} \} \subset M\setminus(\scri^+\cup\iota^+),\quad \rho_\cK:=t_*^{-1},
\end{equation}
Thus, ${}^\cuop S^*D$ contains an open neighborhood of $\pa\Gamma\subset{}^\cuop S^*_{(\cK^+)^\circ}D$.

\begin{figure}[!ht]
\centering
\includegraphics{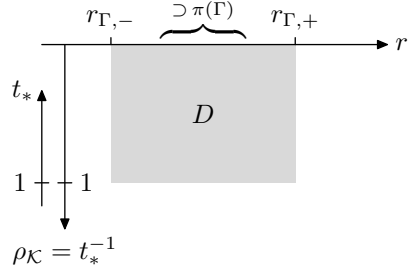}
\caption{The domain $D$ from~\eqref{EqDyTrD} over which Proposition~\ref{PropDyTr} describes the stable and unstable trapped sets of an asymptotically Kerr metric. We write $\pi(\Gamma)$ for the projection of the trapped set to the base.}
\label{FigDyTrD}
\end{figure}

\begin{prop}[Stable and unstable trapped sets of asymptotically Kerr metrics]
\label{PropDyTr}
  For $d_0\geq 2$ and metrics $g=g_{\bhm,a}+h$ where $h\in\tilde\sG_\etbop^{d_0,(0,0,\ell_\cK)}$ (with $\ell_\cK>0$) has sufficiently small norm, there exist functions
  \[
    \tilde\phi^{\rm u/s} \in \cC_\cuop^{d_0-1,\ell_\cK}({}^\cuop S^*D)
  \]
  with norm $\leq 1$ such that, upon setting $\phi^{\rm u/s}:=\phi_0^{\rm u/s}+\tilde\phi^{\rm u/s}$, the vector field $H_G$ is tangent to
  \[
    \Gamma^{\rm u/s}:=\Sigma^+\cap(\phi^{\rm u/s})^{-1}(0).
  \]
  Here $G$ and $\Sigma^+$ are as in~\eqref{EqDyMetric}--\eqref{EqDyChar}. The same statement holds true if we moreover allow for sufficiently small perturbations of $\bhm,a$.
\end{prop}

In view of the decay of $\tilde\phi^{\rm u/s}$, the sets $\Gamma^{\rm u/s}$ are $\cC_\cuop^{d_0-1,\ell_\cK}$-perturbations of $\Gamma_0^{\rm u/s}$ (and thus their intersections with ${}^\cuop T^*_{\cK^+}D$ are equal to $\bar\Gamma^{\rm u/s}\cap{}^\cuop T^*D$).

\begin{proof}[Proof of Proposition~\usref{PropDyTr}]
  After some pre-processing as in \cite[\S{4}]{HintzPolyTrap}, this follows from \cite[Theorem~2.6]{HintzPolyTrap}, with the local uniformity in the subextremal Kerr parameters $\bhm,a$ following from the continuous (in fact, smooth) dependence of the trapped set on Kerr on $(\bhm,a)$ (see Remark~\ref{RmkTse3bKerrDep}). Concretely, denote the characteristic set for the Kerr metric $g_{\bhm,a}$ at $\cK^+\cap D$ by
  \[
    \fX := \pa\Sigma_{\bhm,a}^+ \cap {}^\cuop S_{(\cK^+)^\circ}^*D.
  \]
  Its stationary extension $\fM_0=[0,1)_{\rho_\cK}\times\fX$, is thus a smooth codimension $1$ submanifold of ${}^\cuop S^*D$ with orientable normal bundle. Denote by $\fX\times(-1,1)$ a tubular neighborhood of $\fX$ within ${}^\cuop S^*_{(\cK^+)^\circ}D$ such that the set $\{x\}\times(-1,1)$ for $x\in\fX$ is contained in the fiber of ${}^\cuop S^*D$ containing $x$. Extend this tubular neighborhood to a tubular neighborhood $([0,1)\times\fX)\times(-1,1)$ of $\fM_0$ in ${}^\cuop S^*D$. We then claim that over $D$, the characteristic set $\pa\Sigma^+\cap{}^\cuop S^*D$ (which, in each fiber of ${}^\cuop S^*D$, is a 2-sphere) for the metric $g$ as the graph of a function
  \[
    f \in \cC_\cuop^{d_0,\ell_\cK}\bigl(([0,1)\times\fX) \cap {}^\cuop S^*D; (-1,1)\bigr).
  \]
  Indeed, complementing the defining function $x^0:=\rho_\infty^2 G_{\bhm,a}$ of $\fX$ inside of ${}^\cuop S^*_{(\cK^+)^\circ}D$ with local coordinates $x'$ on $\fX$, and writing $\rho_\infty^2 G=x^0+\rho_\cK^{\ell_\cK}q_0(\rho_\cK,x^0,x')$ where $q\in\cC_\cuop^{d_0}([0,1)_{\rho_\cK}\times(-1,1)_{x^0}\times\R^5_{x'}$, the function $f=f(\rho_\cK,x')$ is the solution $x^0$ of $x^0+\rho_\cK^{\ell_\cK}q_0(\rho_\cK,x^0,x')=0$. Writing $f=\rho_\cK^{\ell_\cK}f_0$, the function $f_0$ is the fixed point $y$ of the map $y\mapsto-q_0(\rho_\cK,\rho_\cK^{\ell_\cK}y,x')$, which indeed exists and is unique in the ball of radius $\|q_0\|_{L^\infty}+1$. Since $q_0$ is $d_0$-times continuously differentiable for $\rho_\cK>0$, so is $f_0$; and the membership $f_0\in\cC_\cuop^{d_0}$ then follows by direct differentiation of the defining equation $f_0(\rho_\cK,x')+q_0(\rho_\cK,\rho_\cK^{\ell_\cK}f_0(\rho_\cK,x'),x')=0$.

  Pulling back the normalized vector field $\sfH_G/(\sfH_G\ft_*)$ along the map $(\rho_\cK,x')\mapsto(\rho_\cK,f(\rho_\cK,x'),x')$ (which maps $\fM_0$ into the characteristic set for $g$), we obtain a vector field $V:=\tilde\sfH_G$ on $\fM_0=[0,1)\times\fX$. Applying this procedure in the case $h=0$ gives the restriction $V_0$ of $\tilde\sfH_{G_{\bhm,a}}=\sfH_{G_{\bhm,a}}/\sfH_{G_{\bhm,a}}\ft_*$ to the characteristic set of $g_{\bhm,a}$. Thus,
  \[
    V = \tilde\sfH_G = V_0 + \tilde V, \quad \tilde V\in\rho_\cK^{\ell_\cK}\cC_\cuop^{d_0-1}(\fM_0;T\fX).
  \]
  (The notation emphasizes that the perturbation $\tilde V$ does not have a $\pa_{\ft_*}$-component: we have $V\ft_*=V_0\ft_*=1$.)

  We can now apply \cite[Theorem~2.6]{HintzPolyTrap} to the vector field $-V$ on $\fM_0$ to obtain a unique stable trapped set $\Gamma^{\rm s}$ and a (non-unique) unstable trapped set $\Gamma^{\rm u}=\bigcup_{s=0}^\infty e^{s V}(\Gamma^{\rm u}_0\cap\{\rho_\cK=1\}$. (The failure of compactness of $\fM_0$ is not an issue since the construction of $\Gamma^{\rm s/u}$ is local near the stationary extension $\Gamma_0$ of the trapped set of Kerr. Note also that stable and unstable manifolds exchange roles when passing from $V$ to $-V$; hence the difference to the statement of \cite[Theorem~2.6]{HintzPolyTrap} in our present discussion. Finally, note that our present definition of $\Gamma^{\rm u}$ matches the initial choice in the proof of \cite[Theorem~2.6]{HintzPolyTrap} which is then propagated to all times by iterating the time $1$ flow $e^V$ in \cite[Part~2 of the proof of Theorem~2.3]{HintzPolyTrap}.) Finally, the proof of \cite[Theorem~2.3]{HintzPolyTrap} (which \cite[Theorem~2.6]{HintzPolyTrap} is reduced to by considering the time $1$ flow $e^V$) constructs the perturbed stable/unstable manifolds \emph{directly} as graphs over the unperturbed (stationary) stable/unstable manifolds, and thus \emph{directly} produces the desired defining functions $\phi^{\rm u/s}$. The claimed uniform boundedness of $\phi^{\rm u/s}$, when $h$ is itself uniformly bounded, follows from the uniformity of all estimates in the proof of \cite[Theorem~2.3]{HintzPolyTrap} on the perturbation $\wt V$ in assumption~(II.2) before \cite[Theorem~2.3]{HintzPolyTrap}.
\end{proof}

\subsection{Monotone order function}
\label{SsDyO}

In~\S\ref{SR}, we will study the regularity of waves on $\Omega_*$ on e3b-Sobolev spaces with variable orders $\sfs$ which must be monotonically decreasing along the future-directed null-bicharacteristic flow. In this section, we recall the construction of appropriate variable order functions. We strengthen the requirements in Definition~\ref{DefSSOrderAdm} to include an appropriate stability under perturbations. Recall $\ell_\sscri,\ell_+,\ell_\cK$ from~\eqref{EqDyWeights} and the notation~\eqref{EqDyRhos}--\eqref{EqDyHam}.

\begin{prop}[Admissible order function]
\label{PropDyO}
  Let $P_0$ be a stationary wave-type operator (Definition~\usref{DefSSAdm}) relative to the Kerr metric $g_{\bhm,a}$. Let $\eta\geq 0$ and $\alpha_+,\alpha_\cK\in\R$. Then there exist an order function $\sfs\in\CI({}^\etbop S^*M)$ and a number $\eps_\sG>0$ such that the following conditions hold.
  \begin{enumerate}
  \item $\sfs|_{{}^\etbop S^*_{\cK^+}M}$, identified with\footnote{Note here that ${}^\tbop S^*M_0$ and ${}^\etbop S^*M$ are equal over $\cK^+$.} a time-translation-invariant element of $\Stb^*M_0$, is stationary-$P_0$-admissible (Definition~\usref{DefSSOrderAdm}) with margin $\eta$ and with $\alpha_\sface=\alpha_+$.
  \item\label{ItDyO2} For all metrics $g=g_{\bhm,a}+h$ with $\|h\|_{\tilde\sG_\etbop^{1,(2\ell_\sscri,\ell_+,\ell_\cK)}(\Omega_*)}<\eps_\sG$, we have $\sfH_G\sfs\leq 0$ on the future characteristic set $\pa\Sigma^+$ (see~\eqref{EqDyChar}).
  \item\label{ItDyO3} $\sfs$ is invariant under dilations $(t_*,r,\omega)\mapsto(\lambda t_*,\lambda r,\omega)$ in a neighborhood of $\iota^+\cap\{\frac{r}{t_*}\in[1,10]\}\subset M$.
  \item The order function $\sfs|_{{}^\etbop S^*_{\cK^+}M}$ is stationary-$P_1$-admissible also for all sufficiently small perturbations $P_1$ of $P_0$ in the class of stationary wave-type operators (including for different but nearby Kerr parameters), and part~\eqref{ItDyO2} is valid also for nearby Kerr parameters.
  \end{enumerate}
\end{prop}

\begin{definition}[Admissible order functions]
\label{DefDyO}
  We call order functions $\sfs$ satisfying the conclusions of Proposition~\usref{PropDyO} \emph{admissible} (for $P_0$) \emph{with weights $\alpha_+,\alpha_\cK$, and margin $\eta$}.
\end{definition}

\begin{proof}[Proof of Proposition~\usref{PropDyO}]
  This follows from \cite[Lemma~5.8 and Remark~5.10]{HintzNonstat}: we may take $\sfs$ to be a large constant over an arbitrarily large compact subset of $(\cK^+)^\circ$, and transition to a sufficiently negative number near $\pa\cR^+_{\pa\cK^+,{\rm out}}\cup\pa\cW^+_{\rm out}\cup\pa\cR^+_{\scri^+,{\rm out}}$ (cf.\ \eqref{EqTs3bStable}). Note that only the dynamics near this latter set, which lies over $\iota^+$, matter in the construction given in \cite[Appendix~C]{HintzNonstat}, and these are the same here as in \cite{HintzNonstat} (i.e., unaffected by trapping and horizons). Part~\eqref{ItDyO3} is arranged in the reference as well: in the notation of \cite[Appendix~C]{HintzNonstat}, the cutoff function $\chi$ used to define $\sfs$ is equal to $\chi_{\iota,\delta}$ for $x_\sscri\geq\eps_\sscri$ where $\eps_\sscri\ll 1$ and $x_\sscri=\sqrt{(t_*+2)/r}$ in the reference; and $\chi_{\iota,\delta}$ is the product of cutoff functions in $\hat\xi_\tbop,\hat\eta_\tbop$ (called $\phi,\psi$ after \cite[(C.9)]{HintzNonstat}), which are thus dilation-invariant, and the function $\tilde\psi$ which, on the support of $\phi$ (where $|\hat\xi_\tbop|\leq\frac12$), is identically $1$ for small $\rho_+$.
\end{proof}

\section{Microlocal regularity control on spacetime}
\label{SR}

Given a weakly admissible wave-type operator $P$ (Definition~\ref{DefSDWAdm}) (relative to $P_0$), we prove precise global results on the edge-3b-regularity of solutions $u$ of $P u=f$ via microlocal elliptic regularity, propagation of singularities, radial point, and trapping estimates on weighted e3b-Sobolev spaces (see~\eqref{EqMUe3bSpace0} and~\eqref{EqMUSuppExt}); see~\S\ref{SsR3}. \emph{These estimates rely only on the principal symbol of $P$ and thus on the dynamical properties of asymptotically Kerr metrics (except for the precise values of certain threshold conditions on weights and regularity orders which are determined by subprincipal terms of $P$).} We then upgrade these to e3b-regularity results on spaces capturing additional $k\in\N_0$ degrees of b-regularity (\S\ref{SsRb}) which are, moreover, tame in the b-regularity order.

The e3b-microlocal results (beyond the elliptic and real principal type results proved in~\S\S\ref{SsMBasic}--\ref{SsMTame}) can almost be quoted from the literature: \cite[\S{4}]{HintzVasyScrieb} near the radial sets $\pa\cR^+_{\scri^+,{\rm in},+}$ and $\pa\cR+_{\scri^+,{\rm out},+}$ over $\scri^+$ (Definition~\ref{DefTsebRad}), \cite[\S{5.1}]{HintzNonstat} near the radial sets $\pa\cR^+_{\pa\cK^+,{\rm in/out}}$ over $\pa\cK^+$ (Definition~\ref{DefTs3bRad}), \cite[\S{4.3}]{HintzGlueLocII} near the radial set $\pa\cR^+_{\cH^+}$ over the event horizon (Definition~\ref{DefTs3bH}), and \cite[\S{3B2}]{HintzPolyTrap} near the trapped set $\pa\Gamma$ (see~\eqref{EqTs3bGamma}). For the sake of completeness, we do present the detailed proofs here (except at the trapped set). One minor novelty here is that we record the stability of these estimates under perturbations of $P$ within the class of admissible wave-type operators with coefficients of limited e3b-regularity.

Similarly, b-regularity results were established in these references at those invariant sets, \emph{except} at the event horizon and the trapped set. But even at the radial sets that are treated in the references, our proof is more robust in that it generalizes easily to estimates that are tame in the b-regularity order (cf.\ Remark~\ref{RmkRbBetter}). Such tame estimates appear here for the first time, and they are based on commutations with e3b-commutator b-vector fields (i.e., Lemma~\ref{LemmaCTe3bComm} and sharper versions in special regimes, such as Lemma~\ref{LemmaCT3bDil}).

\bigskip

Throughout this section, we fix the time function $\ft_*$ from Lemma~\usref{LemmaSDGTime}, and we recall $\Omega_*=\ol{\{\ft_*\geq 1\}}\subset\Omega\subset M$ from~\eqref{EqDyDomain} and Definitions~\usref{DefCMSpacetime} and \usref{DefCMDomain}. Moreover, recall $r_\natural=\frac{\bhm+r_+}{2}\in(\bhm,r_+)$ from Lemma~\ref{LemmaSDGTime} (so $\dd r$ is timelike for $r\leq r_\natural$) and fix
\begin{equation}
\label{EqRRadii}
  \bhm < r_1 < r_1^+ < r_\natural
\end{equation}
as well as a cutoff function
\begin{equation}
\label{EqRCutoff}
\begin{split}
  \chi\in\CI(\Omega_*),\quad \chi&=0\ \text{on}\ \ft_*^{-1}([1,2])\cup r^{-1}([\bhm,r_1]), \\ \chi&=1\ \text{on}\ \{\ft_*\geq 3\} \cap \{r\geq r_1^+\}.
\end{split}
\end{equation}
We study weakly admissible wave-type operators $P$ of class $((d_0;d),(2\ell_\sscri,\ell_+,\ell_\cK))$, acting on sections of a stationary bundle $\cE\to M$, where $d_0,d\in\N_0$ will be sufficiently large (depending on the e3b-regularity order with which one works) but fixed. Recall from Definition~\ref{DefSDWAdm} that this means that the weights are subject only to
\[
  \ell_\sscri\in(0,\tfrac12],\ \ell_+\in(0,1],\ \ell_\cK>0.
\]
(The stronger requirements~\eqref{EqSDWAdmWeights} play no role in the present microlocal regularity considerations.) We require $P$ to be a perturbation of $P_0$ in the norm $\|\cdot\|_{(d_0;0),(2\ell_\sscri,\ell_+,\ell_\cK)}$ (which thus implies the smallness conditions in Lemma~\ref{LemmaSDGTime} and Proposition~\ref{PropDyTr}). We moreover fix weights
\begin{equation}
\label{EqRWeights}
  \alpha_\sscri < -\frac12+\ubar p_1,\quad
  \alpha_+ < \alpha_\sscri - \frac12,\quad \alpha_\cK \in \R,
\end{equation}
where $\ubar p_1$ was defined in Definition~\ref{DefSDWubarp1}, and admissible order functions
\begin{equation}
\label{EqROrderFns}
  \sfs_0,\,\sfs\in\CI({}^\etbop S^*M)
\end{equation}
for $P_0$ with weights $\alpha_+,\alpha_\cK$ (and margin $0$),\footnote{The threshold conditions in Definition~\usref{DefSDWAdm} arise in the radial point estimates at those radial sets that lie over $\cK^+$. The condition on $\alpha_\sscri$ in~\eqref{EqRWeights} arises in the radial point estimate at $\pa\cR^+_{\scri^+,{\rm out}}$. The upper bound on $\alpha_+$ is a requirement in the radial point estimate at $\pa\cR_{\scri^+,{\rm in},+}$.} as produced by Proposition~\ref{PropDyO}, with $\sfs_0<\sfs$.

\subsection{Edge-3b-regularity estimate}
\label{SsR3}

The basic regularity estimate reads as follows.

\begin{prop}[Global e3b-estimate]
\label{PropR3}
  With $\alpha_\sscri,\alpha_+,\alpha_\cK,\sfs_0,\sfs$ as above, there exists a constant $C$ such that
  \begin{equation}
  \label{EqR3Est}
    \|\chi u\|_{H_\etbop^{\sfs,(2\alpha_\sscri,\alpha_+,\alpha_\cK)}} \leq C\Bigl( \|P u\|_{H_\etbop^{\sfs,(2\alpha_\sscri+2,\alpha_++2,\alpha_\cK)}(\Omega_*;\cE)^{\bullet,-}} + \|u\|_{H_\etbop^{\sfs_0,(2\alpha_\sscri,\alpha_+,\alpha_\cK)}(\Omega_*;\cE)^{\bullet,-}}\Bigr)
  \end{equation}
  for $P=P_0$, in the strong sense that if the right-hand side is finite, then so is the left-hand side and the estimate holds. Moreover, there exists $d_0\in\N$ (depending only on $\sfs_0,\sfs$) such that this estimate also holds for all weakly admissible wave-type operators $P$ relative to $P_0$ which are sufficiently small perturbations of $P_0$ as measured in the norm $\|\cdot\|_{(d_0;0),(2\ell_\sscri,\ell_+,\ell_\cK)}$ (Definition~\usref{DefSDWAdmNorm}), with the constant $C$ uniform for such perturbations.
\end{prop}

The gap between having an estimate for the cut-off version $\chi u$ and wanting an estimate for $u$ on all of $\Omega_*$ will be bridged using energy estimates in~\S\ref{SsEReg}.

\begin{rmk}[Uniformity in $P_0$]
\label{RmkR3UnifP0}
  The estimate~\eqref{EqR3Est} holds also when the underlying stationary operator $P_0$ is perturbed within the class of stationary wave-type operators (Definition~\ref{DefSSAdm}), which in particular allows for perturbations of the Kerr parameters $\bhm$ and $a$. This follows from the uniformity of all microlocal estimates entering in its proof below: for elliptic and real principal type estimates, this is shown for more general perturbations in Propositions~\ref{PropMEll}, \ref{PropMPr}, \ref{PropMTameEll}, and \ref{PropMTamePr}, for radial point estimates this follows from the openness of the positivity of symbols in the symbolic (positive commutator) arguments, and likewise for the propagation estimate at the trapped set which, in view of its continuous (in fact, smooth) dependence on the Kerr parameters (cf.\ Remark~\ref{RmkTse3bKerrDep} and the last part of Proposition~\ref{PropDyTr}).
\end{rmk}

\begin{rmk}[Large perturbations of $P_0$]
\label{RmkR3UnifLarge}
  The phrasing of the perturbation-stability of the estimate~\eqref{EqR3Est} is made for convenience and since it suffices for our applications. We do note, however, that if $P$ is an arbitrary weakly admissible wave-type operator relative to $P_0$ for which $\ft_*$ is a global time function and $r$ is a time function near $r=\bhm$, the estimate~\eqref{EqR3Est} still holds by the same proof: the estimates near $\scri^+\cup\iota^+\cup\cK^+$ are unaffected, and propagating control to arbitrarily large times and into the black hole is still accomplished by standard real principal type propagation estimates.
\end{rmk}

The proof strategy is to concatenate radial point and trapping estimates at each of the aforementioned invariant sets, in the order given by Proposition~\ref{PropTse3bDyn}, i.e., going from bottom to top in Figure~\ref{FigTse3bDyn}.

\subsubsection{Radial point estimates}
\label{SssR3R}

We shall only prove estimates in the future characteristic set $\pa\Sigma^+$, as estimates in the past component $\pa\Sigma^-$ are proved in exactly the same manner upon reversing the sign of the Hamiltonian vector field.

We begin with the radial set over the event horizon at $t_*=\infty$. We recall from~\eqref{EqTs3bHCoord} the coordinates $t_0,\phi_0$ and the momentum variables $\sigma_0,\xi_0,\eta_\theta,\eta_{\phi_0}$. Since we are working away from $\scri^+$, we drop all sub- and superscripts referring to the edge-nature of phase space over $\scri^+$; equivalently, we can work with 3b-objects on the manifold $M_0$ from Definition~\ref{DefCMSpacetime} (see~\eqref{EqMUKPsdo} for the notation for 3b-ps.d.o.s.)

\begin{prop}[Radial point estimate near $\pa\cR^+_{\cH^+}$]
\label{PropR3RH}
  There exist operators $B,E,G\in\tilde\Psi_\tbop^0$ such that
  \begin{enumerate}
  \item\label{ItIR3RH1} their operator wave front sets are contained in any fixed neighborhood $\cU$ of $\pa\cR^+_{\cH^+}\subset{}^\tbop S^*M$,
  \item\label{ItIR3RH2} their Schwartz kernels are supported in $K\times K$ for any fixed compact neighborhood $K\subset M$ of $\cH^+\cap\cK^+=\{r=r_+,\ \rho_\cK=0\}$,
  \item\label{ItIR3RH3} $B$ is elliptic at $\pa\cR^+_{\cH^+}$,
  \item\label{ItIR3RH4} the Schwartz kernel of $E$ is supported in $(K\setminus\cK^+)^2$,
  \end{enumerate}
  and such that for all
  \[
    s>s_0>\frac12+\vartheta_{\cH^+}
  \]
  (in the notation of Definition~\usref{DefSSOrderAdm}\eqref{ItSSOrderAdm4}) and $\alpha_\cK\in\R$, the following holds for any cutoff function $\tilde\chi\in\CI(M)$ supported near $K$ and equal to $1$ near $K$. There exists a constant $C$ such that
  \begin{equation}
  \label{EqR3RHEst}
    \|B u\|_{\Htb^{s,\alpha_\cK}} \leq C\Bigl( \|G P u\|_{\Htb^{s-1,\alpha_\cK}} + \|E u\|_{\Htb^{s,\alpha_\cK}} + \|\tilde\chi u\|_{\Htb^{s_0,\alpha_\cK}} \Bigr),
  \end{equation}
  for $P=P_0$, in the strong sense that if the right-hand side is finite, then so is the left-hand side and the estimate holds. Moreover, there exists $d_0\in\N$ (depending only on $s_0,s$) such that this estimate also holds for all weakly admissible wave-type operators $P$ relative to $P_0$ which are sufficiently small perturbations of $P_0$ as measured in the norm $\|\cdot\|_{(d_0;0),(0,0,\ell_\cK)}$ (see Definition~\usref{DefSDWAdmNorm}), with the constant $C$ uniform for such perturbations.
\end{prop}

The estimate~\eqref{EqR3RHEst} allows us to propagate $s$ degrees of 3b-regularity from finite times (the term $E u$) into $\pa\cR^+_{\cH^+}$ (the term $B u$). We remark that a similar result, but for exponentially weighted spaces, is given in \cite[Theorem~6.10]{HintzQuasilinearDS}; unlike in the reference, we do not give an explicit lower bound for the regularity $d_0$ of the coefficients of $P$ but instead only use the abstract observation in Lemma~\ref{LemmaMSOpNorm}.

\begin{proof}[Proof of Proposition~\usref{PropR3RH}]
  For now, we assume that $P$ is of class $((\infty;0),(0,0,\ell_\cK))$. As is usual for positive commutator arguments (cf.\ \cite[Chapter~8]{HintzMicro}), we will first establish the estimate
  \begin{equation}
  \label{EqR3RHEst2}
    \|B u\|_{\Htb^{s,\alpha_\cK}} \leq C\Bigl( \|G P u\|_{\Htb^{s-1,\alpha_\cK}} + \|E u\|_{\Htb^{s,\alpha_\cK}} + \|G u\|_{\Htb^{s-\frac12,\alpha_\cK}} + \|\tilde\chi u\|_{\Htb^{s_0,\alpha_\cK}} \Bigr)
  \end{equation}
  which thus improves on the regularity order of $u$ by (at most) half an order; this estimate also holds in the usual strong sense. Using this estimate for shrunk sets $\cU$ and $K$, one can then estimate the norm of $G u$ inductively as long as $s-\frac12>s_0$, and thus obtain~\eqref{EqR3RHEst} as stated.

  To prove~\eqref{EqR3RHEst2}, we recall from Lemma~\ref{LemmaTs3bHQuadDef} the quadratic defining function $\fq$ of $\pa\cR^+_{\cH^+}$ inside of a small neighborhood $\cO$ of $\pa\cR^+_{\cH^+}$ inside of $\pa\Sigma^+\cap\Stb^*_{(\cK^+)^\circ}M$; extend $\fq$ to a smooth function on a neighborhood of $\cO$ inside of $\Stb^*M$, which we still denote by $\fq$. Let
  \begin{equation}
  \label{EqR3RHPsi}
    \psi\in\CIc((-1,1)),\quad \psi=1\ \text{on}\ [-\tfrac12,\tfrac12],\quad \sqrt{-x\psi(x)\psi'(x)}\in\CI.
  \end{equation}
  With $\rho_\infty=\xi_0^{-1}$ as in~\eqref{EqTs3bHCoordProj} (which is elliptic of order $-1$ near $\pa\cR^+_{\cH^+}$), we then consider the commutant
  \begin{align*}
    &a = \check a^2,\quad \check a:=\rho_\cK^{-\alpha_\cK}\rho_\infty^{-s+\frac12} \psi_\Sigma \psi_\cR \psi_\cK, \\
    &\qquad \psi_\Sigma = \psi(\digamma\rho_\infty^2 G), \ \psi_\cR = \psi(\digamma_\cR\fq), \ \psi_\cK = \psi(\digamma\rho_\cK),
  \end{align*}
  where we shall choose $\digamma,\digamma_\cR>1$ large; note that $\supp a\cap{}^\tbop S^*M\to\pa\cR^+_{\cH^+}$ when $\digamma,\digamma_\cR\to\infty$. To prove the strong version of~\eqref{EqR3RHEst2}, we insert a regularizer,
  \[
    a_\eps = \check a_\eps^2,\ \check a_\eps := \phi_\eps(\rho_\infty)\check a,\quad \phi_\eps(\rho_\infty)=(1+\eps \rho_\infty^{-1})^{-\Theta},\ \Theta := s-s_0>0,
  \]
  so $\{a_\eps\colon\eps\in(0,1]\}$ is uniformly bounded in $S^{2 s-1,2\alpha_\cK}(\Ttb^*M)=\rho_\cK^{-2\alpha_\cK}S^{2 s-1}(\Ttb^*M)$ and lies in $S^{2 s_0-1,2\alpha_\cK}$ for all $\eps>0$. We shall use that
  \[
    \rho_\infty\phi_\eps'(\rho_\infty) = \Theta f_\eps(\rho_\infty)\phi_\eps(\rho_\infty),\quad f_\eps(\rho_\infty):=\frac{\eps\rho_\infty^{-1}}{1+\eps\rho_\infty^{-1}} \in [0,1].
  \]
  We fix the positive definite fiber inner product on $\cE_X$ used in the definition of $\vartheta_{\cH^+}$ in~\eqref{EqSSOrderAdm4}. We use a quantization map $\Op_\tbop$ producing 3b-ps.d.o.s (with 3b-regular coefficients) and introduce
  \[
    A_\eps = \check A_\eps^*\check A_\eps,\quad \check A_\eps=\Op_\tbop(\psi_\infty\check a_\eps);\ \ (A_\eps)_{\eps\in(0,1]} \in L^\infty\bigl((0,1]_\eps;\tilde\Psi_\tbop^{s-\frac12,\alpha_\cK}\bigr);
  \]
  here $\psi_\infty=\psi(\rho_\infty)$ cuts away the singularity at the zero section. Note that $A_\eps\in\tilde\Psi_\tbop^{s_0-\frac12,\alpha_\cK}$ for all $\eps\in(0,1]$. We then compute the $L^2$-pairing
  \begin{equation}
  \label{EqR3RHCommCalc}
    2 \Im\la\check A_\eps P u,\check A_\eps u\ra = 2 \Im \la P u,A_\eps u\ra = \la \sC_\eps u,u\ra,\quad \sC_\eps := i[P,A_\eps] + 2\frac{P-P^*}{2 i}A_\eps.
  \end{equation}
  Write $p_\sub:=\upsigma_\tbop^1(\frac{P-P^*}{2 i})$ (which is a symbol of order $1$ valued in self-adjoint bundle endomorphisms of $\cE$) and $\sfp_\sub:=\rho_\infty p_\sub$. Note that since $P$ equals $P_0$ plus a 3b-differential operator $P-P_0$ whose coefficients decay as $\rho_\cK\to 0$, the formula~\eqref{EqSSOrderAdm4} for $\vartheta_{\cH^+}$ can equivalently be written as
  \begin{equation}
  \label{EqR3RHThreshold}
    \vartheta_{\cH^+} = \sup_{\pa\cR^+_{\cH^+}} \spec\Bigl[\frac{\varrho^2}{2(r_+-\bhm)}\rho_\infty\upsigma_\tbop^1\Bigl(\frac{P-P^*}{2 i}\Bigr)\Bigr] = \sup_{\pa\cR^+_{\cH^+}} \spec\Bigl[\frac{\varrho^2}{2(r_+-\bhm)}\sfp_\sub\Bigr].
  \end{equation}
  Put $\sfH_G:=\rho_\infty H_G$ (which is homogeneous of degree $0$ in the fibers of $\Stb^*M$). Let us write $\psi'_\Sigma=\psi'(\digamma\rho_\infty^2 G)$, $\psi'_\cR=\psi'(\digamma_\cR\fq)$, and $\psi'_\cK=\psi'(\digamma\rho_\cK)$. We then compute the principal symbol of $\sC_\eps$ to be
  \begin{subequations}
  \begin{equation}
  \label{EqR3RHSymb}
    H_G a_\eps + 2 p_\sub a_\eps = -\delta\rho_\infty^{-1}\check a_\eps^2 - b_\eps^2 - b_{\cR,\eps}^2 + e_\eps^2 + j_\eps G,
  \end{equation}
  where
  \begin{equation}
  \label{EqR3RHSymb2}
  \begin{split}
    b_\eps &:= \phi_\eps\rho_\infty^{-s}\rho_\cK^{-\alpha_\cK} \psi_\Sigma\psi_\cR\psi_\cK \Bigl[ (\rho_\infty^{-1}\sfH_G\rho_\infty)(2 s-1- 2\Theta f_\eps) - \delta - 2\sfp_\sub + 2 (\rho_\cK^{-2}\sfH_G\rho_\cK)\alpha_\cK\rho_\cK \Bigr]^{\frac12}, \\
    b_{\cR,\eps} &:= \phi_\eps\rho_\infty^{-s}\rho_\cK^{-\alpha_\cK}\psi_\Sigma\psi_\cK\bigl[ -2 \digamma_\cR(\sfH_G\fq)\psi_\cR\psi'_\cR\bigr]^{\frac12}, \\
    e_\eps &:= \phi_\eps\rho_\infty^{-s}\rho_\cK^{-\alpha_\cK+\frac12}\psi_\Sigma\psi_\cR\bigl[ 2 \digamma (\rho_\cK^{-2}\sfH_G\rho_\cK) \rho_\cK\psi_\cK\psi'_\cK \bigr]^{\frac12}, \\
    j_\eps &:= 2\phi_\eps^2\rho_\infty^{-2 s+2}\rho_\cK^{-2\alpha_\cK}\psi_\cR^2\psi_\cK^2 \digamma\psi_\Sigma\psi'_\Sigma (\rho_\infty^{-2}\sfH_G\rho_\infty^2).
  \end{split}
  \end{equation}
  \end{subequations}
  The first three terms have the following properties:
  \begin{itemize}
  \item By~\eqref{EqTs3bHExpr} or \eqref{EqTs3bHBdfDer} (and recalling that $G_\tbop=\varrho^2 G_{\bhm,a}$, and $G=G_{\bhm,a}$ over $\cK^+$), we have $\rho_\infty^{-1}\sfH_G\rho_\infty=\varrho^{-2}\mu'(r_+)=\frac{2(r_+-\bhm)}{\varrho^2}$ at $\pa\cR^+_{\cH^+}$. Therefore, if $\supp a$ is contained in a sufficiently small neighborhood of $\pa\cR^+_{\cH^+}$, the term of $b_\eps$ in square brackets is strictly positive on $\supp a$ when $\delta>0$ is sufficiently small, since $\frac{2(r_+-\bhm)}{\varrho^2}(s_0-\frac12)>\sfp_\sub$. Note here that $2 s-2\Theta f_\eps\geq 2 s_0$ for all $\eps>0$. (The weight $\alpha_\cK$ does not matter due to the additional factor of $\rho_\cK$, which vanishes at $\rho_\cK=0$.) Thus, $b_\eps$ is well-defined, of class $S^{s,\alpha_\cK}$ uniformly in $\eps$, and of class $S^{s_0,\alpha_\cK}$ for every $\eps>0$.
  \item For $e_\eps$, we note that $\rho_\cK^{-2}\sfH_G\rho_\cK<0$ at $\pa\cR^+_{\cH^+}$ by~\eqref{EqTs3bHExpr} or \eqref{EqTs3bHBdfDer}, and thus also $e_\eps$ is well-defined when $\supp a$ is a sufficiently small neighborhood of $\pa\cR^+_{\cH^+}$.
  \item Finally, regarding $b_{\cR,\eps}$, note that if $\digamma_\cR$ is sufficiently large, we have $\sfH_G\fq\geq 2 c$ for some constant $c>0$ on the support of $\psi'_\cR$ (where $\fq\in[\frac{1}{2\digamma_\cR},\frac{1}{\digamma_\cR}]$) intersected with $\rho_\cK^{-1}(0)$; this follows from Lemma~\ref{LemmaTs3bHQuadDef}. Then also $\sfH_G\fq\geq c$ on $\supp a\cap\supp\psi'_\cR$ if we localize sufficiently closely to the characteristic set and to $\rho_\cK=0$, which we can do by choosing $\digamma$ large (depending on the choice of $\digamma_\cR$). 
  \end{itemize}

  One now passes to quantizations $B_\eps=\Op_\tbop(\psi_\infty b_\eps)$, $B_{\cR,\eps}$, $E_\eps$, $J_\eps$ of the corresponding lower case symbols. Note that one can make their Schwartz kernels have support in the square of any small neighborhood of projection of $\supp a$ to the base. Then
  \[
     \sC_\eps = -\delta(\Lambda^{\frac12}\check A_\eps)^*(\Lambda^{\frac12}\check A_\eps) - B_\eps^*B_\eps - B_{\cR,\eps}^*B_{\cR,\eps} + E_\eps^*E_\eps + J_\eps P + R_\eps,
  \]
  where $\Lambda^{\frac12}\in\tilde\Psi_\tbop^{\frac12}$ has principal symbol $\rho_\infty^{-\frac12}$ near $\supp a$, and where $B_\eps,B_{\cR,\eps},E_\eps$ are uniformly bounded in $\tilde\Psi_\tbop^{s,\alpha_\cK}$ and of order $s_0$ for every $\eps>0$, while $J_\eps$ and $R_\eps$ are uniformly bounded in $\tilde\Psi_\tbop^{2 s-2,2\alpha_\cK}$ and $\tilde\Psi_\tbop^{2 s-1,2\alpha_\cK}$, respectively, with orders reduced to $2 s_0-2$ and $2 s_0-1$ for $\eps>0$. The term $R_\eps$ collects the subprincipal contributions to $\sC_\eps$ not captured by the above principal symbol computation. Plugging this into~\eqref{EqR3RHCommCalc} and using the Peter--Paul inequality gives the $L^2$-estimate
  \begin{equation}
  \label{EqR3RHCommEst}
  \begin{split}
    \delta\|\Lambda^{\frac12}\check A_\eps u\|^2 + \|B_\eps u\|^2 + \|B_{\cR,\eps}u\|^2 &\leq \delta\|\Lambda^{\frac12}\check A_\eps u\|^2 + \frac{C}{\delta}\|\Lambda^{-\frac12}\check A_\eps P u\|^2 + \|E_\eps u\|^2 \\
      &\qquad + |\la J_\eps u,P u\ra| + |\la R_\eps u,u\ra|
  \end{split}
  \end{equation}
  some constant $C$ depending only on the choice $\Lambda^{-\frac12}\in\tilde\Psi_\tbop^{-\frac12}$ of an operator that is elliptic near $\supp a$. The first terms on the left- and right-hand sides cancel. If $G\in\tilde\Psi_\tbop^0$ is elliptic on $\supp a$ and has Schwartz kernel with support close to that of $A_\eps$, we can uniformly (in $\eps$) estimate, using elliptic regularity,
  \begin{align*}
    |\la J_\eps u,P u\ra| &\leq C\Bigl(\|G u\|_{H_\tbop^{s-1,\alpha_\cK}}^2 + \|G P u\|_{H_\tbop^{s-1,\alpha_\cK}}^2 + \|\tilde\chi u\|_{H_\tbop^{s_0,\alpha_\cK}}^2\Bigr), \\
    |\la R_\eps u,u\ra| &\leq C\Bigl(\|G u\|_{H_\tbop^{s-\frac12,\alpha_\cK}}^2 + \|\tilde\chi u\|_{H_\tbop^{s_0,\alpha_\cK}}^2\Bigr).
  \end{align*}
  Plugging these into~\eqref{EqR3RHCommEst}, one gets the uniform bound
  \[
    \|B_\eps u\| \leq C\Bigl( \|G P u\|_{H_\tbop^{s-1,\alpha_\cK}} + \|E u\| + \|G u\|_{H_\tbop^{s-\frac12,\alpha_\cK}} + \|\tilde\chi u\|_{H_\tbop^{s_0,\alpha_\cK}} \Bigr).
  \]
  Letting $\eps\searrow 0$ and using the weak compactness of the unit ball in $L^2$ as well as the distributional convergence $B_\eps u\to B_0 u$, we conclude that $B_0 u\in L^2$ and the same estimate for $\|B_0 u\|$. This establishes~\eqref{EqR3RHEst2} (upon shifting the orders of $B,E$ relative to present notation) and finishes the proof; note here that $B_0$ is elliptic at $\pa\cR^+_{\cH^+}$, so we have indeed proved the strong version of~\eqref{EqR3RHEst2}.

  As mentioned before, the stability of the estimate~\eqref{EqR3RHEst} under perturbations follows from the fact that the mapping properties of 3b-ps.d.o.s used in the above argument between fixed pairs of 3b-Sobolev spaces only require a sufficiently large finite amount of 3b-regularity of the underlying symbols (see Lemma~\ref{LemmaMSOpNorm}).
\end{proof}

\begin{rmk}[Pseudodifferential lower order terms]
\label{RmkR3RHPsdolot}
  Proposition~\ref{PropR3RH} holds also for operators of class
  \[
    P = P_0 + \tilde P,\quad \tilde P \in \tilde\Psi_\tbop^1;
  \]
  more generally still, it suffices to assume $\tilde P\in\cC_\ebop^{d_0}\Psi_\tbop^1$ for sufficiently large $d_0$; and the resulting estimate is locally uniform in $\tilde P$ when the pseudodifferential terms involve only ps.d.o.s in any fixed finite subset of $\tilde\Psi_\tbop^1$ (cf.\ \eqref{EqMUCe3bOpNormPsdo}). The proof goes through without any changes since subprincipal terms of $P$, such as $\tilde P$, only enter the argument through $\sfp_\sub$ and thus through the symbol of $\frac{1}{2 i}(\tilde P-\tilde P^*)$. Thus, the only place where $\tilde P$ enters is in the determination of the threshold condition. In our application of this pseudodifferential generalization of Proposition~\ref{PropR3RH} in~\S\ref{SssRbR} below, the principal symbol $\upsigma_\tbop^1(\tilde P)$ is continuous (in fact, smooth plus decaying conormal) down to $\cK^+$, and thus the contribution of $\upsigma_\tbop^1(\frac{1}{2 i}(\tilde P-\tilde P^*))$ at $N^*\cH^+\setminus o$ to $\vartheta_{\cH^+}$ in~\eqref{EqR3RHThreshold} is well-defined. (If the symbol of $\tilde P$ is only conormal down to $\cK^+$, one needs to define $\vartheta_{\cH^+}$ as the $\limsup$ of the right-hand side in~\eqref{EqR3RHThreshold} as one approaches $\pa\cR^+_{\cH^+}$; and then Proposition~\ref{PropR3RH} is again valid, with the same proof.)
\end{rmk}

The proofs of the other radial point estimates are very similar (including the regularization, and including the possibility of allowing for additional first order pseudodifferential terms as in Remark~\ref{RmkR3RHPsdolot}); we shall thus be brief.

\begin{prop}[Radial point estimate near $\pa\cR^+_{\pa\cK^+,{\rm in}}$]
\label{PropR3RKI}
  There exist operators $B,E,G\in\tilde\Psi_\tbop^0$ such that
  \begin{enumerate}
  \item their operator wave front sets are contained in any fixed neighborhood of $\pa\cR^+_{\pa\cK^+,{\rm in}}$,
  \item their Schwartz kernels are supported in $K\times K$ where $K$ is any fixed compact neighborhood of $\pa\cK^+$,
  \item $B$ is elliptic at $\pa\cR^+_{\pa\cK^+,{\rm in}}$,
  \item the Schwartz kernel of $E$ is supported in $(K\setminus\cK^+)^2$,
  \end{enumerate}
  and such that for all $s,s_0,\alpha_+,\alpha_\cK\in\R$ with $s>s_0$ and
  \begin{equation}
  \label{EqR3RKIThr}
    s_0+\alpha_+>\frac12(-1+\vartheta_{\pa\cK^+,{\rm in}})+\alpha_\cK
  \end{equation}
  (in the notation of Definition~\usref{DefSSOrderAdm}\eqref{ItSSOrderAdm6}), the following holds for any cutoff function $\tilde\chi\in\CI(M)$ supported near $K$ and equal to $1$ near $K$. There exists a constant $C$ such that
  \begin{equation}
  \label{EqR3RKIEst}
    \|B u\|_{\Htb^{s,(\alpha_+,\alpha_\cK)}} \leq C\Bigl( \|G P u\|_{\Htb^{s-1,(\alpha_++2,\alpha_\cK)}} + \|E u\|_{\Htb^{s,(\alpha_+,\alpha_\cK)}} + \|\tilde\chi u\|_{\Htb^{s_0,(\alpha_+,\alpha_\cK)}} \Bigr)
  \end{equation}
  holds in the strong sense for $P=P_0$ and all of its sufficiently small perturbations as measured in the norm $\|\cdot\|_{(d_0;0),(0,\ell_+,\ell_\cK)}$, where $d_0\in\N$ is sufficiently large (depending only on $s,s_0$).
\end{prop}
\begin{proof}
  This is a minor variation of \cite[Proposition~5.5]{HintzNonstat}. We use local defining functions $\rho_\cK=\frac{r}{t_*}$, $\rho_+=\frac{1}{r}$ of $\cK^+$ and $\iota^+$, respectively, and use the 3b-momentum variables $\sigma_\tbop,\xi_\tbop,\eta_\tbop$ from~\eqref{EqTs3bCoord}; recalling that $\sigma_\tbop<0$ at $\pa\cR^+_{\pa\cK^+,{\rm in}}$, we moreover set $\rho_\infty=-\frac{1}{\sigma_\tbop}$ and $(\hat\xi_\tbop,\hat\eta_\tbop)=(\frac{\xi_\tbop}{\sigma_\tbop},\frac{\eta_\tbop}{\sigma_\tbop})$ as in~\eqref{EqTs3bCCoordProj}, so $\pa\cR^+_{\pa\cK^+,{\rm in}}$ is defined inside of $\Stb^*M$ by the equations $\rho_\cK=\rho_+=0$, $\hat\xi_\tbop=2$, $\hat\eta_\tbop=0$. We set
  \[
    \fq := |\hat\eta_\tbop|^2 + (\hat\xi_\tbop-2)^2
  \]
  and compute using~\eqref{EqTs3bCHam} that at $\iota^+$,
  \[
    \sfH_{G_\tbop}\fq = 4(\hat\xi_\tbop-1)|\hat\eta_\tbop|^2 + 4\hat\xi_\tbop(\hat\xi_\tbop-2)^2 \geq 2\fq
  \]
  near $\pa\cR^+_{\pa\cK^+,{\rm in}}$. We will use a positive commutator estimate with the commutant
  \begin{equation}
  \label{EqR3RKIComm}
  \begin{split}
    &a = \check a^2,\quad \check a=\rho_\infty^{-s+\frac12}\rho_+^{-\alpha_+-1}\rho_\cK^{-\alpha_\cK}\psi_\Sigma\psi_\cR\psi_+\psi_\cK, \\
    &\qquad \psi_\Sigma = \psi(\digamma\rho_\infty^2\rho_+^{-2}G),\ \psi_\cR=\psi(\digamma_\cR\fq),\ \psi_+=\psi_+(\digamma\rho_+),\ \psi_\cK=\psi(\digamma\rho_\cK),
  \end{split}
  \end{equation}
  where $\psi$ is as in~\eqref{EqR3RHPsi} and $\digamma,\digamma_\cR>1$ will be chosen large in order to localize $\supp a$ near $\pa\cR^+_{\pa\cK^+,{\rm in}}$. We then compute the principal symbol of $i[P,A]+2\frac{P-P^*}{2 i}A$, where $A=\check A^*\check A$, $\check A=\Op_\tbop(\check a)$, and we use the fiber inner product on $\cE_X$ used in the definition of $\vartheta_{\pa\cK^+,{\rm in}}$: it is given by $H_G a+2 p_\sub a$ where $p_\sub=\sigmatb^1(\frac{P-P^*}{2 i})$, so in terms of $\sfp_\sub:=\rho_+^{-2}\rho_\infty^{-1}p_\sub$ and $\sfH_G:=\rho_+^{-2}\rho_\infty H_G$,\footnote{The multiplication by $\rho_+^{-2}$ takes out the quadratically vanishing weight at $\iota^+$, cf.\ the passage from $G_{(\bhm,a)}$ to $G_\tbop$ in~\eqref{EqTs3bGtb}.} we have
  \[
    H_G a + 2 p_\sub a = -\delta\rho_\infty^{-1}\rho_+^2\check a^2 -b^2 - b_\cR^2 - b_\cK^2 + e_+^2 + j G,
  \]
  where the functions
  \begin{equation}
  \label{EqR3RKISymb}
  \begin{split}
    b &:= \rho_\infty^{-s}\rho_+^{-\alpha_+}\rho_\cK^{-\alpha_\cK}\psi_\Sigma\psi_\cR\psi_+\psi_\cK \\
      &\qquad \times \Bigl[ (\rho_\infty^{-1}\sfH_G\rho_\infty)(2 s-1) + (\rho_+^{-1}\sfH_G\rho_+)(2\alpha_++2) + 2(\rho_\cK^{-1}\sfH_G\rho_\cK)\alpha_\cK - \delta - 2\sfp_\sub \Bigr]^{\frac12}, \\
    b_\cR &:= \rho_\infty^{-s}\rho_+^{-\alpha_+}\rho_\cK^{-\alpha_\cK}\psi_\Sigma\psi_+\psi_\cK\bigl[ -2\digamma_\cR \psi_\cR\psi_\cR' (\sfH_G\fq) \bigr]^{\frac12}, \\
    b_\cK &:=\rho_\infty^{-s}\rho_+^{-\alpha_+}\rho_\cK^{-\alpha_\cK}\psi_\Sigma\psi_\cR\psi_\cK\bigl[ -2(\rho_+^{-1}\sfH_G\rho_+)\digamma\rho_+\psi_+\psi'_+\bigr]^{\frac12}, \\
    e_+ &:=\rho_\infty^{-s}\rho_+^{-\alpha_+}\rho_\cK^{-\alpha_\cK}\psi_\Sigma\psi_\cR\psi_+\bigl[ 2(\rho_\cK^{-1}\sfH_G\rho_\cK)\digamma\rho_\cK\psi_\cK\psi'_\cK\bigr]^{\frac12}, \\
    j &:= 2\rho_\infty^{-2 s+2}\rho_+^{-2\alpha_+-2}\rho_\cK^{-2\alpha_\cK} \psi_\cR^2\psi_+^2\psi_\cK^2 \digamma\psi_\Sigma\psi'_\Sigma(\rho_\infty^{-2}\rho_+^2\sfH_G\rho_\infty^2\rho_+^{-2})
  \end{split}
  \end{equation}
  have the following properties.
  \begin{itemize}
  \item At $\pa\cR^+_{\pa\cK^+,{\rm in}}$ and using~\eqref{EqTs3bCHam} (or~\eqref{EqTs3bLin}), the term in square brackets in the expression for $b$ equals $2(2 s-1)+2(2\alpha_++2)-2\cdot 2\alpha_\cK-\delta-2\sfp_\sub$, and hence is positive for small $\delta>0$ (and also for $s_0$ in place of $s$). This persists in a sufficiently small neighborhood of $\pa\cR^+_{\pa\cK^+,{\rm in}}$, and $b$ is a symbol of class $S^{s,(\alpha_+,\alpha_\cK)}$ when $\digamma,\digamma_\cR$ are sufficiently large.
  \item For all sufficiently large $\digamma_\cR$, we have $\sfH_G\fq\geq 2 c>0$ on $\iota^+\cap\supp\psi'_\cR$, and thus $\sfH_G\fq\geq c>0$ nearby; thus $b_\cR$ is well-defined if $\digamma$ are chosen sufficiently large (depending on $\digamma_\cR$).
  \item The symbols $b_\cK$ and $e_+$ are well-defined for all sufficiently large $\digamma,\digamma_\cR$ since $\rho_\cK^{-1}\sfH_G\rho_\cK=-2<0$ and $\rho_+^{-1}\sfH_G\rho_+=2>0$ at $\pa\cR^+_{\pa\cK^+,{\rm in}}$. Note that $\supp e_+\cap\Stb^*_{\cK^+}M=\emptyset$.
  \end{itemize}

  Upon (regularizing and) passing to quantizations as in the proof of Proposition~\ref{PropR3RH}, we obtain the estimate
  \[
    \|B u\|_{L^2}^2 \leq C\Bigl(\|G P u\|_{\Htb^{s-1,\alpha_++2,\alpha_\cK}}^2 + \|E_+ u\|_{L^2}^2 + \|G u\|_{\Htb^{s-\frac12,\alpha_+,\alpha_\cK}}^2 + \|\tilde\chi u\|_{\Htb^{s_0,\alpha_+,\alpha_\cK}}^2 \Bigr)
  \]
  (which holds in the strong sense), from which~\eqref{EqR3RKIEst} follows by iteration.
\end{proof}

\begin{prop}[Radial point estimate near $\pa\cR^+_{\pa\cK^+,{\rm out}}$]
\label{PropR3RKO}
  There exist $B,E,G\in\tilde\Psi^0_\tbop$ such that
  \begin{enumerate}
  \item their operator wave front sets are contained in any fixed neighborhood of $\pa\cR^+_{\pa\cK^+,{\rm out}}$,
  \item their Schwartz kernels are supported in $K\times K$ where $K$ is any fixed compact neighborhood of $\pa\cK^+$,
  \item $B$ is elliptic at $\pa\cR^+_{\pa\cK^+,{\rm out}}$,
  \item the operator wave front set of $E$ is disjoint from the unstable manifold $\pa\cW_{\rm out}^+$ of $\pa\cR^+_{\pa\cK^+,{\rm out}}$ (see~\eqref{EqTs3bStable}),
  \end{enumerate}
  and such that for all $s,s_0,\alpha_+,\alpha_\cK\in\R$ with $s>s_0$ and
  \begin{equation}
  \label{EqR3RKOThr}
    s+\alpha_+<\frac12(-1+\vartheta_{\pa\cK^+,{\rm out}})+\alpha_\cK
  \end{equation}
  (in the notation of Definition~\usref{DefSSOrderAdm}\eqref{ItSSOrderAdm5}), the estimate~\eqref{EqR3RKIEst} holds (for any $\tilde\chi\in\CI(M)$ supported near $K$ and equal to $1$ near $K$) in the strong sense for $P=P_0$ and all of its sufficiently small perturbations in the norm $\|\cdot\|_{(d_0;0),(0,\ell_+,\ell_\cK)}$ when $d_0\in\N$ is large enough and fixed (depending only on $s,s_0$).
\end{prop}
\begin{proof}
  This is a minor variation of \cite[Proposition~5.6]{HintzNonstat}. In the coordinates used already in the proof of Proposition~\ref{PropR3RKI}, we recall that $\pa\cR^+_{\pa\cK^+,{\rm out}}$ is defined inside of $\Stb^*M$ by the equations $\rho_\cK=\rho_+=0$, $\hat\xi_\tbop=0$, $\hat\eta_\tbop=0$. For $\fq:=|\hat\eta_\tbop|^2+\hat\xi_\tbop^2$, we now compute using~\eqref{EqTs3bCHam} that at $\iota^+$,
  \begin{equation}
  \label{EqR3RKOfq}
    \sfH_{G_\tbop}\fq =  - 4(2-\hat\xi_\tbop)\hat\xi_\tbop^2 - 4(1-\hat\xi_\tbop)|\hat\eta_\tbop|^2 \leq -2\fq
  \end{equation}
  near $\pa\cR^+_{\pa\cK^+,{\rm out}}$. For the commutant $a$ given by~\eqref{EqR3RKIComm} (but with our current definition of $\fq$), we now have, in the notation used in the proof of Proposition~\ref{PropR3RKI},
  \[
    H_G a+2 p_\sub a = -\delta\rho_\infty^{-1}\rho_+^2\check a^2 - b^2 - b_+^2 + e_\cR^2 + e_\cK^2 + j G,
  \]
  where $b,j$ are as in~\eqref{EqR3RKISymb}, while
  \begin{align*}
    e_\cR &:= \rho_\infty^{-s}\rho_+^{-\alpha_+}\rho_\cK^{-\alpha_\cK}\psi_\Sigma\psi_+\psi_\cK\bigl[ 2\digamma_\cR \psi_\cR\psi_\cR' (\sfH_G\fq) \bigr]^{\frac12}, \\
    e_\cK &:=\rho_\infty^{-s}\rho_+^{-\alpha_+}\rho_\cK^{-\alpha_\cK}\psi_\Sigma\psi_\cR\psi_\cK\bigl[ 2(\rho_+^{-1}\sfH_G\rho_+)\digamma\rho_+\psi_+\psi'_+\bigr]^{\frac12}, \\
    b_+ &:=\rho_\infty^{-s}\rho_+^{-\alpha_+}\rho_\cK^{-\alpha_\cK}\psi_\Sigma\psi_\cR\psi_+\bigl[ -2(\rho_\cK^{-1}\sfH_G\rho_\cK)\digamma\rho_\cK\psi_\cK\psi'_\cK\bigr]^{\frac12}
  \end{align*}
  are defined like $b_\cR,b_\cK,e_+$ in~\eqref{EqR3RKISymb} except for sign switches under the square roots, necessitated by the sign switches of~\eqref{EqR3RKOfq}, $\rho_\cK^{-1}\sfH_G\rho_\cK=2>0$, and $\rho_+^{-1}\sfH_G\rho_+=-2<0$ at $\pa\cR^+_{\pa\cK^+,{\rm out}}$ compared to the previous proof. The terms $e_\cR,e_\cK$ necessitate the a priori control term $E u$ in~\eqref{EqR3RKIEst}; note that $\supp e_\cR$ is disjoint from $\{\hat\xi_\tbop=\hat\eta_\tbop=0\}$ and thus from $\pa\cW^+_{\rm out}$, and $\supp e_\cK$ is disjoint from $\iota^+$ altogether. Moreover, at $\pa\cR^+_{\pa\cK^+,{\rm out}}$, the term in square brackets in the expression for $b$ in~\eqref{EqR3RKISymb} equals $-2(2 s-1)-2(2\alpha_++2)+2\cdot 2\alpha_\cK-\delta-2\sfp_\sub$ and is thus positive (for small $\delta>0$) under the condition~\eqref{EqR3RKOThr}. From here on, the proof of the proposition proceeds as before, \textit{mutatis mutandis}.
\end{proof}

We next turn to the radial sets over $\scri^+$ (Definition~\ref{DefTsebRad}). We now drop sub- and superscripts referring to the 3-body nature of e3b-phase space over $\cK^+$. Correspondingly, we write $\tilde\Psi_\ebop^{s,(2\alpha_\sscri,\alpha_+)}=x_\sscri^{2\alpha_\sscri}\rho_+^{\alpha_+}\tilde\Psi_\ebop^s$ for spaces of eb-ps.d.o.s with eb-regular coefficients, which are the same as e3b-ps.d.o.s with e3b-regular coefficients whose Schwartz kernels are supported away from $\cK^+$ in both factors. Here we write $x_\sscri,\rho_+\in\CI(\tilde M)$ for defining functions of $\sscri^+,\iota^+$.

\begin{prop}[Radial point estimate near $\pa\cR^+_{\scri^+,{\rm in},+}$]
\label{PropR3RScriI}
  There exist $B,E,G\in\tilde\Psi^0_\ebop$ such that
  \begin{enumerate}
  \item their operator wave front sets are contained in any fixed neighborhood of $\pa\cR^+_{\scri^+,{\rm in}}$,
  \item their Schwartz kernels are supported in $K\times K$ where $K$ is any fixed compact neighborhood of $\scri^+\cap\iota^+$,
  \item $B$ is elliptic at $\pa\cR^+_{\scri^+,{\rm in}}$,
  \item the Schwartz kernel of $E$ is supported in $(K\setminus\iota^+)^2$,
  \end{enumerate}
  and such that for all $s,s_0,\alpha_\sscri,\alpha_+\in\R$ with $s>s_0$ and
  \begin{equation}
  \label{EqR3RScriIThr}
    \alpha_+ < \alpha_\sscri - \frac12,
  \end{equation}
  the following holds for any cutoff function $\tilde\chi\in\CI(M)$ supported near $K$ and equal to $1$ near $K$. There exists a constant $C$ such that
  \begin{equation}
  \label{EqR3RScriIEst}
    \|B u\|_{\Heb^{s,(2\alpha_\sscri,\alpha_+)}} \leq C\Bigl( \|G P u\|_{\Heb^{s-1,(2\alpha_\sscri+2,\alpha_++2)}} + \|E u\|_{\Heb^{s,(2\alpha_\sscri,\alpha_+)}} + \|\tilde\chi u\|_{\Heb^{s_0,(2\alpha_\sscri,\alpha_+)}} \Bigr)
  \end{equation}
  for $P=P_0$ holds in the strong sense, and also for all sufficiently small perturbations of $P_0$ in the norm $\|\cdot\|_{(d_0;0),(2\ell_\sscri,\ell_+,0)}$ when $d_0\in\N$ is large enough and fixed (depending only on $s,s_0$).
\end{prop}
\begin{proof}
  This is essentially the content of \cite[Lemma~4.10(1)]{HintzVasyScrieb}. We again only present the unregularized symbolic computation for the positive commutator argument. We use the coordinates $\rho_+=\frac{1}{t_*}$, $x_\sscri=\sqrt{\frac{t_*}{r}}$ from~\eqref{EqTsebCoord}, the eb-momenta $\sigma_\ebop,\xi_\ebop,\eta_\ebop$ from~\eqref{EqCHamebCoord}, and their projective versions $(\rho_\infty,\hat\xi_\ebop,\hat\eta_\ebop)=(\frac{1}{\sigma_\ebop},\frac{\xi_\ebop}{\sigma_\ebop},\frac{\eta_\ebop}{\sigma_\ebop})$ from~\eqref{EqTsebCoordFiber}. Thus, $\pa\cR^+_{\scri^+,{\rm in},+}$ is defined by $x_\sscri=\rho_+=0$, $\hat\xi_\ebop=2$, $\hat\eta_\ebop=0$. Recalling the overall weight $x_\sscri^2\rho_+^2$ of the dual metric as an edge-b-metric near $\scri^+$, we then consider the commutant
  \begin{align*}
    &a = \check a^2,\quad
    \check a = \rho_\infty^{-s+\frac12}x_\sscri^{-2\alpha_\sscri-1}\rho_+^{-\alpha_+-1}\psi_\Sigma\psi_\cR\psi_\sscri\psi_+, \\
    &\qquad \psi_\Sigma=\psi(\digamma\rho_\infty^2 x_\sscri^{-2}\rho_+^{-2}G),\ \psi_\cR=\psi(\digamma_\cR\fq),\ \psi_\sscri=\psi(\digamma x_\sscri),\ \psi_+=\psi(\digamma\rho_+),
  \end{align*}
  where $\psi$ is as in~\eqref{EqR3RHPsi}, $\fq=(\hat\xi_\ebop-2)^2+|\hat\eta_\ebop|^2$, and $\digamma,\digamma_\cR>1$ will be chosen large. Using~\eqref{EqTsebHam}, we compute over $\scri^+$, where the characteristic set is defined by $|\hat\eta_\ebop|^2=\frac12\hat\xi_\ebop(2-\hat\xi_\ebop)$,
  \begin{equation}
  \label{EqR3RScriIq}
    \sfH_{G_\ebop}\fq = 2(\hat\xi_\ebop-1)|\hat\eta_\ebop|^2 - 2|\hat\eta_\ebop|^2(\hat\xi_\ebop-2) = 2(\hat\xi_\ebop-1)|\hat\eta_\ebop|^2 + 2\hat\xi_\ebop(\hat\xi_\ebop-2)^2 \geq \fq
  \end{equation}
  on the characteristic set near $\pa\cR^+_{\scri^+,{\rm in},+}$. We then write
  \[
    H_G a + 2 p_\sub a = -\delta\rho_\infty^{-1}x_\sscri^2\rho_+^2\check a^2 - b^2 - b_\cR^2 - b_+^2 + e_\sscri^2 + j G,
  \]
  where we set $p_\sub=\upsigma_\ebop(\frac{P-P^*}{2 i})$ and, setting $\sfH_G:=\rho_\infty x_\sscri^{-2}\rho_+^{-2}H_G$ and $\sfp_\sub:=\rho_\infty x_\sscri^{-2}\rho_+^{-2}p_\sub$, define the symbols
  \begin{align*}
    b &:= \rho_\infty^{-s}x_\sscri^{-2\alpha_\sscri}\rho_+^{-\alpha_+}\psi_\Sigma\psi_\cR\psi_\sscri\psi_+ \\
      &\qquad \times \Bigl[ (x_\sscri^{-1}\sfH_G x_\sscri)(4\alpha_\sscri+2) + (\rho_+^{-1}\sfH_G\rho_+)(2\alpha_++2) + (\rho_\infty^{-1}\sfH_G\rho_\infty)(2 s-1) - \delta - 2\sfp_\sub \Bigr]^{\frac12}, \\
    b_\cR &:= \rho_\infty^{-s}x_\sscri^{-2\alpha_\sscri}\rho_+^{-\alpha_+}\psi_\Sigma\psi_\sscri\psi_+\bigl[ -2\digamma_\cR\psi_\cR\psi'_\cR(\sfH_G\fq) \bigr]^{\frac12}, \\
    b_+ &:= \rho_\infty^{-s}x_\sscri^{-2\alpha_\sscri}\rho_+^{-\alpha_+}\psi_\Sigma\psi_\cR\psi_+\bigl[ -2(x_\sscri^{-1}\sfH_G x_\sscri)\digamma x_\sscri\psi_\sscri\psi_\sscri' \bigr]^{\frac12}, \\
    e_\sscri &:= \rho_\infty^{-s}x_\sscri^{-2\alpha_\sscri}\rho_+^{-\alpha_+}\psi_\Sigma\psi_\cR\psi_\sscri\bigl[ 2(\rho_+^{-1}\sfH_G\rho_+)\digamma\rho_+\psi_+\psi'_+ \bigr]^{\frac12}, \\
    j &:= 2\rho_\infty^{-2 s+2}x_\sscri^{-4\alpha_\sscri-2}\rho_+^{-\alpha_+-2}\psi_\cR\psi_\sscri\psi_+ \digamma\psi_\Sigma\psi'_\Sigma (\rho_\infty^{-2}x_\sscri^2\rho_+^2\sfH_G\rho_\infty^2 x_\sscri^{-2}\rho_+^{-2}).
  \end{align*}
  Note here that at $\pa\cR^+_{\scri^+,{\rm in},+}$ and using~\eqref{EqTsebHam} (or~\eqref{EqTsebRadLin}), we have $x_\sscri^{-1}\sfH_G x_\sscri=1>0$, $\rho_+^{-1}\sfH_G\rho_+=-2<0$, and $\rho_\infty^{-1}\sfH_G\rho_\infty=0$; this implies that $b_+$ and $e_\sscri$ are well-defined when $\digamma,\digamma_\cR$ are large enough, and also $b_\cR$ is well-defined by~\eqref{EqR3RScriIq} provided we first choose $\digamma_\cR$ and then $\digamma$ sufficiently large. Regarding $b$, the first four summands in the expression in square brackets sum to $4\alpha_\sscri-4\alpha_+-2-\delta$ at $\pa\cR^+_{\scri^+,{\rm in},+}$, which is positive for sufficiently small $\delta>0$ since $\alpha_\sscri>\alpha_++\frac12$. It remains to show that $\sfp_1=0$ at $\pa\cR^+_{\scri^+,{\rm in},+}$. For later purposes, we compute, more generally,
  \begin{equation}
  \label{EqR3RScripsub}
    \sfp_\sub=(\Re p_1)(\hat\xi_\ebop-2)\ \text{over}\ \scri^+,
  \end{equation}
  where we recall $p_1:=S|_{\pa X}+\tilde p_1$ (in the notation of~\eqref{EqSSAdmOp} and \eqref{EqSDWAdmOp}). To prove~\eqref{EqR3RScripsub}, observe first that the operator $\frac12 x_\sscri^2\rho_+^2(x_\sscri\pa_{x_\sscri}-2)(x_\sscri\pa_{x_\sscri}-2\rho_+\pa_{\rho_+})$ (see~\eqref{EqSDWAdmNearScri}), is symmetric with respect to the Minkowskian volume density $r^2|\dd t_*\,\dd r\,\dd\slg|=2\rho_+^{-4}x_\sscri^{-6}|\frac{\dd\rho_+}{\rho_+}\frac{\dd x_\sscri}{x_\sscri}\,\dd\slg|$. (Note here that over $\scri^+$ we can replace the Kerr by the Minkowski metric in this computation.) The imaginary part of $x_\sscri^4\slDelta$ is of class $x_\sscri^4\Diff^1(\Sph^2)\subset x_\sscri^3\Diff_\ebop^1$, so $x_\sscri^{-2}$ times it vanishes at $\scri^+$. Next, the term $p_0$ in~\eqref{EqSDWAdmOp} is sub-subprincipal and thus does not contribute to $\sfp_\sub$. It remains to compute the imaginary part of $x_\sscri^2\rho_+^2 p_1(x_\sscri\pa_{x_\sscri}-2\rho_+\pa_{\rho_+})$, which is easily found to have symbol $(\Re p_1)(\xi_\ebop-2\sigma_\ebop)$; this gives~\eqref{EqR3RScripsub}. At $\pa\cR^+_{\scri^+,{\rm in},+}$ then, \eqref{EqR3RScripsub} vanishes.
\end{proof}

Finally, we need a radial point estimate at $\pa\cR^+_{\scri^+,{\rm out}}$ which involves localization in $t_*$. (This is not stated in the required form in \cite[Lemma~4.11]{HintzVasyScrieb}.)

\begin{prop}[Radial point estimate near $\pa\cR^+_{\scri^+,{\rm out}}$]
\label{PropR3RScriO}
  Let $t_{*,0}\leq t_{*,-}<t_{*,+}\leq\infty$. Then there exist $B,E,G\in\tilde\Psi_\ebop^0$ such that
  \begin{enumerate}
  \item their operator wave front sets are contained in any fixed neighborhood of the set $\pa\cR^+_{\scri^+,{\rm out}}\cap t_*^{-1}([t_{*,-},t_{*,+}])$,
  \item their Schwartz kernels are supported in $K\times K$ where $K$ is any fixed compact neighborhood of $\scri^+\cap t_*^{-1}([t_{*,-},t_{*,+}])$,
  \item $B$ is elliptic at $\pa\cR^+_{\scri^+,{\rm out}}\cap t_*^{-1}([t_{*,-},t_{*,+}])$,
  \item the operator wave front set of $E$ is disjoint from $\pa\cR^+_{\scri^+,{\rm out}}$,
  \end{enumerate}
  and such that for all $s,s_0,\alpha_\sscri,\alpha_+\in\R$ with $s>s_0$ and
  \begin{equation}
  \label{EqR3RScriO}
    \alpha_\sscri < -\frac12 + \ubar p_1
  \end{equation}
  (in the notation of Definition~\usref{DefSDWubarp1}), the estimate~\eqref{EqR3RScriIEst} holds (for any $\tilde\chi\in\CI(M)$ supported near $K$ and equal to $1$ near $K$) in the strong sense for $P=P_0$ and all of its sufficiently small perturbations in the norm $\|\cdot\|_{(d_0;0),(2\ell_\sscri,\ell_+,0)}$ when $d_0\in\N$ is large enough and fixed (depending only on $s,s_0$).
\end{prop}
\begin{proof}
  We use the same notation as in the proof of Proposition~\ref{PropR3RScriI}. It is convenient to use the time function $\ft_*$ from Lemma~\ref{LemmaSDGTime}, which in terms of the coordinates $x_\sscri,\rho_+$ used in the preceding proof is equal to $t_*(1-x_\sscri^{2\ell_\sscri})=\rho_+^{-1}(1-x_\sscri^{2\ell_\sscri})$ near $\scri^+$. Set $\fq:=\hat\xi_\ebop^2+|\hat\eta_\ebop|^2$; this vanishes quadratically at $\pa\cR^+_{\scri^+,{\rm out}}$ and, by~\eqref{EqTsebHam}, satisfies
  \[
    \sfH_{G_\ebop}\fq = -2(1-\hat\xi_\ebop)|\hat\eta_\ebop|^2 - 4|\hat\eta_\ebop|^2\hat\xi_\ebop = -2(1-\hat\xi_\ebop)|\hat\eta_\ebop|^2 - 4(2-\hat\xi_\ebop)\hat\xi_\ebop^2 \leq -\fq
  \]
  on the characteristic set near $\pa\cR^+_{\scri^+,{\rm out}}$. Fix a function $\eta\in\CI(\R)$ which equals $0$ on $(-\infty,-\frac12]$, $1$ on $[0,\infty)$, and satisfies $\eta'\geq 0$ and $\sqrt{\eta\eta'}\in\CI$, and consider (with $\psi$ as in~\eqref{EqR3RHPsi}) the commutant
  \begin{align*}
    &a = \check a^2,\quad \check a=\rho_\infty^{-s+\frac12}x_\sscri^{-2\alpha_\sscri-1}\rho_+^{-\alpha_+-1} \psi_\Sigma\psi_\cR\psi_\sscri\psi_-\psi_+, \\
    &\qquad \psi_\Sigma=\psi(\digamma\rho_\infty^2 x_\sscri^{-2}\rho_+^{-2}),\ \psi_\cR=\psi(\digamma_\cR\fq),\ \psi_\sscri=\psi(\digamma x_\sscri), \\
    &\qquad \psi_-=\eta(\digamma(t_{*,-}-\ft_*)),\ \psi_+=\eta(\digamma(\ft_*-t_{*,+})).
  \end{align*}
  (When $t_{*,+}=\infty$, we take $\psi_+=1$.) We then write
  \[
    H_G a + 2 p_\sub a = -\delta\rho_\infty^{-1}x_\sscri^2\rho_+^2\check a^2 - b^2 - b_+^2 + e_\cR^2 + e_\sscri^2 + e_-^2 + j G,
  \]
  where we introduce
  \begin{align*}
    b &:= \rho_\infty^{-s}x_\sscri^{-2\alpha_\sscri}\rho_+^{-\alpha_+}\psi_\Sigma\psi_\cR\psi_\sscri\psi_-\psi_+ \\
      &\qquad \times \Bigl[ 2(x_\sscri^{-1}\sfH_G x_\sscri)(2\alpha_\sscri+1) + 2(\rho_+^{-1}\sfH_G\rho_+)\alpha_+ + (\rho_\infty^{-1}\sfH_G\rho_\infty)(2 s-1) - \delta - 2\sfp_\sub \Bigr]^{\frac12}, \\
    b_+ &:= \rho_\infty^{-s}x_\sscri^{-2\alpha_\sscri}\rho_+^{-\alpha_+}\psi_\Sigma\psi_\cR\psi_\sscri\psi_- \bigl[ -2(\sfH_G\ft_*) \digamma\psi_+\psi'_+ \bigr]^{\frac12}, \\
    e_- &:= \rho_\infty^{-s}x_\sscri^{-2\alpha_\sscri}\rho_+^{-\alpha_+}\psi_\Sigma\psi_\cR\psi_\sscri\psi_+ \bigl[ -2(\sfH_G\ft_*) \digamma\psi_-\psi'_- \bigr]^{\frac12}, \\
    e_\cR &:= \rho_\infty^{-s}x_\sscri^{-2\alpha_\sscri}\rho_+^{-\alpha_+}\psi_\Sigma\psi_\sscri\psi_-\psi_+ \bigl[ 2(\sfH_G\fq) \digamma_\cR\psi_\cR\psi_\cR' \bigr]^{\frac12}, \\
    e_\sscri &:= \rho_\infty^{-s}x_\sscri^{-2\alpha_\sscri}\rho_+^{-\alpha_+}\psi_\Sigma\psi_\cR\psi_-\psi_+\bigl[ 2(x_\sscri^{-1}\sfH_G x_\sscri) \digamma x_\sscri \psi_\sscri\psi'_\sscri \bigr]^{\frac12}.
  \end{align*}
  The symbols $b_+,e_-,e_\cR$, and $e_\sscri$ are well-defined if we first choose $\digamma_\cR$ large enough and then $\digamma$ sufficiently large (depending on $\digamma_\cR$); this uses that $\sfH_G\ft_*\geq 0$, and $x_\sscri^{-1}\sfH_G x_\sscri=-1$ at $\pa\cR^+_{\scri^+,{\rm out}}$ by~\eqref{EqTsebHam}.

  Turning to $b$, we note that $\rho_+^{-1}\sfH_G\rho_+=\rho_\infty^{-1}\sfH_G\rho_\infty=0$ at $\pa\cR^+_{\scri^+,{\rm out}}$. Now,~\eqref{EqR3RScripsub} implies that $\sfp_\sub=-2(\Re p_1)$; let us choose a fiber inner product on $\cE_X$ such that $\Re p_1>\ubar p_1-\delta'$ where we fix any $\delta'>0$. (This holds for a diagonal inner product in the splitting $\cE_X=\bigoplus_{j=1}^J \cE_{X,j}$ in the notation of Definition~\ref{DefSDWAdm}\eqref{ItSDWAdmSplit}, with diagonal entries $1,\eta,\ldots,\eta^{J-1}$ for $0<\eta\ll 1$.) Then $\sfp_\sub<-2(\ubar p_1-\delta')$ and thus, at $\pa\cR^+_{\scri^+,{\rm out}}$, the term in square brackets in the expression for $b$ is bounded from below by $-4\alpha_\sscri-2-\delta+4(\ubar p_1-\delta')$; in view of~\eqref{EqR3RScriO}, this is positive if we choose $\delta,\delta'>0$ sufficiently small.

  The remainder of the proof (including regularization) is, again, the same as in Proposition~\ref{PropR3RH}, \textit{mutatis mutandis}.
\end{proof}

\subsubsection{Estimate at normally hyperbolic trapping}
\label{SssR3Tr}

The microlocal propagation estimate near the trapped set $\pa\Gamma$ from~\eqref{EqTs3bGamma} which we need here was proved in \cite{HintzPolyTrap} using cusp ps.d.o.s and associated weighted Sobolev spaces near $(\cK^+)^\circ$ which, as discussed at the beginning of~\S\ref{SsDyTr}, are the same as e3b-ps.d.o.s and e3b-Sobolev spaces upon localization to the set $D=\{t_*>1,\ r_{\Gamma,-}<r<r_{\Gamma,+}\}$ from~\eqref{EqDyTrD}. Setting $\rho_\cK:=t_*^{-1}$, we thus write $\tilde\Psi_\cuop^{s,\alpha_\cK}=\rho_\cK^{-\alpha_\cK}\tilde\Psi_\cuop^0$ for the space of weighted cusp-ps.d.o.s with cusp-regular coefficients; these are the same as elements of $\tilde\Psi_\etbop^{s,(0,0,\alpha_\cK)}$ whose Schwartz kernels are supported in $K\times K$ where $K\subset\tilde M\setminus(\scri^+\cup\iota^+)$. We similarly write $H_\cuop^{s,\alpha_\cK}=\rho_\cK^{\alpha_\cK}H_\cuop^s$ for weighted cusp Sobolev spaces (relative to the Minkowski or Kerr volume density); and we regard the trapped set $\Gamma$ and the dynamical stable and unstable manifolds $\Gamma^{\rm s/u}$ as subsets of ${}^\cuop S^*D$.

\begin{prop}[Estimate near the trapped set]
\label{PropR3Tr}
  There exist operators $B,E,G\in\tilde\Psi_\cuop^0$ such that
  \begin{enumerate}
  \item\label{ItR3Tr1} their operator wave front sets are contained in any fixed neighborhood $\cU\subset{}^\cuop S^*D$ of $\pa\Gamma$,
  \item\label{ItR3Tr2} their Schwartz kernels are supported in $K\times K$ for any fixed compact neighborhood $K\subset M$ of the projection of $\cU$ to the base,
  \item\label{ItR3Tr3} $B$ is elliptic at $\pa\Gamma$,
  \item\label{ItR3Tr4} the operator wave front set of $E$ is disjoint from $\pa\Gamma^{\rm u}$,
  \end{enumerate}
  and such that for all $s>s_0$ and $\alpha_\cK\in\R$, the following holds for any cutoff function $\tilde\chi\in\CI(M)$ supported near $K$ and equal to $1$ near $K$. There exists a constant $C$ such that
  \begin{equation}
  \label{EqR3TrEst}
    \|B u\|_{H_\cuop^{s,\alpha_\cK}} \leq C\Bigl(\|G P u\|_{H_\cuop^{s,\alpha_\cK}} + \|E u\|_{H_\cuop^{s+1,\alpha_\cK}} + \|\tilde\chi u\|_{H_\cuop^{s_0,\alpha_\cK}} \Bigr)
  \end{equation}
  for $P=P_0$, in the strong sense that if the right-hand side is finite, then so is the left-hand side and the estimate holds; and this estimate holds uniformly also for all weakly admissible wave-type operators $P$ relative to $P_0$ which are sufficiently small perturbations of $P_0$ as measured in the norm $\|\cdot\|_{(d_0;0),(0,0,\ell_\cK)}$ (see Definition~\usref{DefSDWAdmNorm}).
\end{prop}
\begin{proof}
  Consider first $P=P_0$ or perturbations of class $((\infty;0),(0,0,\ell_\cK))$; we wish to apply \cite[Theorem~3.9]{HintzPolyTrap} and need to check its assumptions. The assumptions \cite[(P.1)--(P.5)]{HintzPolyTrap} on the stationary model follows from the results in~\S\ref{SssTs3bO} and the observation that $H_{G_{\bhm,a}}\sigma=0$ in the notation of~\eqref{EqTs3bCoord} (so one can take $\wh\rho=|\sigma|^{-1}$ in assumption~(P.2)), while assumption \cite[(P.6')]{HintzPolyTrap} is the content of Proposition~\ref{PropDyTr}; note here that the cusp conormality required in assumption~(P.6') is the same as infinite order cusp regularity of the coefficients of $P-P_0$. The subprincipal symbol condition \cite[(3-24)]{HintzPolyTrap} is almost the content of Definition~\ref{DefSSTrapAdm}; if we work with $Q P Q^-$ instead of $P$ (in the notation of~\eqref{EqSSTrapConj}), then it \emph{is} satisfied. Here we note that~\eqref{EqSSTrapConj} is equivalent to (using ``$\cuop$'' in place of ``$\tbop$'')
  \begin{equation}
  \label{EqR3TrSubpr}
    \rho_\infty \upsigma_\tbop\Bigl(\frac{1}{2 i}\bigl(Q P Q^- - (Q P Q^-)^*\bigr)\Bigr) < \frac12\nu_{\rm min}\quad\text{at}\ \pa\Gamma.
  \end{equation}

  The estimate~\eqref{EqR3TrEst} is then an immediate consequence of \cite[(3-25)]{HintzPolyTrap} with $Q P Q^-$, $B$, $E$, $\alpha_\cK$, $Q u$, and $s_0$ in place of $P$, $B_0$, $B_1$, $r$, $v$, and $-N$, respectively, and with the error term $\|v\|_{H_\cuop^{-N,r}}$ in~\cite[(3-25)]{HintzPolyTrap} replaced by $\|\tilde\chi u\|_{H_\cuop^{s_0,\alpha_\cK}}$ in present notation. This uses two simply observations: first, $Q P Q^-(Q u)=Q P u+Q P(I-Q^- Q)u$, and for $G Q P(I-Q^- Q)\in\tilde\Psi_\cuop^{-\infty}$ we have $\|G Q P(I-Q^- Q)u\|_{H_\cuop^{s,\alpha_\cK}}\leq\|\tilde\chi u\|_{H_\cuop^{s_0,\alpha_\cK}}$. Second, the localizer $\tilde\chi$ can be inserted easily in the error term by applying~\cite[(3-25)]{HintzPolyTrap} to a localization $\chi^\sharp u$ where $\chi^\sharp=1$ near $K$ and $\supp\chi^\sharp\subset\tilde\chi^{-1}(1)$, using that the operator wave front sets of $B_0,B_1,G$ are disjoint from phase space over $\supp\dd\chi^\sharp$ (which implies that, for example, $[E,\chi^\sharp]$ is a residual operator and thus $\|[E,\chi^\sharp]u\|_{H_\cuop^{s,\alpha_\cK}}\leq C\|\tilde\chi u\|_{H_\cuop^{s_0,\alpha_\cK}}$.

  The proof of \cite[Theorem~3.9]{HintzPolyTrap} goes through without any changes if one uses only cusp-ps.d.o.s with cusp-regular coefficients (i.e., using a quantization map for $\tilde\Psi_\cuop=\cC_\cuop^\infty\Psi_\cuop$), as it only utilizes the cusp principal symbol (but no normal operator arguments). Therefore, for fixed $s,s_0$, it applies also for perturbations of $P_0$ by operators with sufficiently high but finite cusp-regularity $d_0$ (cf.\ Lemma~\ref{LemmaMSOpNorm}).
\end{proof}

\subsubsection{Proof of Proposition~\usref{PropR3}}
\label{SssR3Pf}

The proof of Proposition~\ref{PropR3} combines microlocal elliptic and real principal type propagation results with the radial point and trapping estimates from~\S\S\ref{SssR3R}--\ref{SssR3Tr}. The details are as follows; we drop the bundle $\cE$ from the notation, use the notation and assumptions of Proposition~\ref{PropR3}, and recall Figure~\ref{FigTse3bDyn}.

\medskip

\pfstep{Step~1. Estimate for finite $\ft_*$.} For any $\ft_*^\sharp\in(1,\infty)$, we claim that there exists a constant $C=C(\ft_*^\sharp)$ such that
\begin{equation}
\label{EqR3PfWeak}
\begin{split}
  &\|\chi u\|_{H_\etbop^{\sfs+1,(2\alpha_\sscri,\alpha_+,\alpha_\cK)}(\Omega_*\cap\{\ft_*\leq\ft_*^\sharp\})^{\bullet,-}} \\
  &\qquad \leq C \Bigl( \|P u\|_{H_\etbop^{\sfs,(2\alpha_\sscri+2,\alpha_++2,\alpha_\cK)}(\Omega_*)^{\bullet,-}} + \|u\|_{H_\etbop^{\sfs_0,(2\alpha_\sscri,\alpha_+,\alpha_\cK)}(\Omega_*)^{\bullet,-}}\Bigr)
\end{split}
\end{equation}
(and this holds in the strong form). This is thus an estimate for $\chi u$ (see~\eqref{EqRCutoff}) up to any finite time. (We work with the regularity order $\sfs+1$ here in order to counteract a loss of an extra derivative in the trapping estimate in Step~4 below.)

We first prove a weaker version of~\eqref{EqR3PfWeak}. Let $B\in\tilde\Psi_\etbop^0$ be any operator with operator wave front set disjoint from $\pa\cR^+_{\scri^+,{\rm out}}$, and with Schwartz kernel localized sufficiently close to the diagonal that $\chi B u$ is well-defined for $u$ with supported, resp.\ extendible character at the initial, resp.\ final boundary hypersurface of $\Omega_*$. Then $\|\chi B u\|_{H_\etbop^{\sfs+1,(2\alpha_\sscri,\alpha_+,\alpha_\cK)}(\Omega_*\cap\{\ft_*\leq\ft_*^\sharp+1\})^{\bullet,-}}$ is bounded by the right hand side of~\eqref{EqR3PfWeak} where now $C$ also depends on $B$: this follows from microlocal elliptic regularity away from the characteristic set of $P$ (Proposition~\ref{PropMEll}) as well as real principal type propagation for $P$ (Proposition~\ref{PropMPr}), using the fact that all past-directed null-bicharacteristics starting on $\WF'_\cuop(\chi B)$ cross the initial hypersurface $\ft_*^{-1}(1)$ of $\Omega_*$ in finite time and thus enter the region $\ft_*<1$ where $u$ vanishes (so the a priori control term $E u$ in~\eqref{EqMPrEst} is equal to $0$ where one may take $E$ to be a bump function in $\ft_*$ supported in $\{\ft_*<1\}$).

To deduce~\eqref{EqR3PfWeak} from an estimate for $\chi B u$, we apply Proposition~\ref{PropR3RScriO} with $\ft_{*,-}=0$, $\ft_{*,+}=\ft_*^\sharp$, and with $\sfs+1$ in place of $\sfs$ to propagate weighted eb-regularity (with orders $\sfs+1,2\alpha_\sscri$) into $\pa\cR^+_{\scri^+,{\rm out}}$; note that the a priori control term $E u$ in the estimate~\eqref{EqR3RScriIEst} is controlled by what we currently call $\chi B u$. This proves~\eqref{EqR3PfWeak}.

\medskip

\pfstep{Step~2. Radial sets $\#1$ and $\#2$.} The a priori control terms $E u$ in Lemmas~\ref{PropR3RH} (concerning $\pa\cR^+_{\cH^+}$) and \ref{PropR3RScriI} (concerning $\pa\cR^+_{\scri^+,{\rm in},+}$) can be estimated by the left hand side of~\eqref{EqR3PfWeak}, provided $\ft_*^\sharp$ is sufficiently large; this is due to the operators $E$ in these Lemmas having Schwartz kernel supported away from $\ft_*=\infty$. 

\medskip

\pfstep{Step~3. Radial set $\#3$.} Consider the a priori control term $E u$ in Proposition~\ref{PropR3RKI} (concerning $\pa\cR^+_{\pa\cK^+,{\rm in}}$). Past null-bicharacteristics $\gamma$ starting at a point in $\WF_\etbop'(E)$ have $\frac{\dd}{\dd s}\rho_\cK(\gamma(s))>0$; thus, by Lemma~\ref{LemmaTs3bDyn}\eqref{ItTs3bDynip}, $\gamma$ must either tend to $\scri^+$ (if it lies over $\iota^+$) and hence to $\pa\cR^+_{\scri^+,{\rm in},+}$ (note that $\pa\cR^+_{\scri^+,{\rm out}}$ does not feature in the ``$A$'' column of~\eqref{EqTse3bDynTable}), or it crosses $\ft_*^{-1}(1)$ in finite time. Reversing this, we thus conclude that the norm of $E u$ in~\eqref{EqR3RKIEst} can be controlled using the propagation of regularity (Proposition~\ref{PropMPr}) starting in the region where, by Step~2, we already have control on $u$. Therefore, we obtain $H_\etbop^{\sfs+1,(\alpha_+,\alpha_\cK)}$-control of $u$ microlocally near $\pa\cR^+_{\pa\cK^+,{\rm in}}$.

\medskip

\pfstep{Step~4. Trapping.} In Steps~1--3, we have obtained control on $u$ for $\ft_*\leq\ft_*^\sharp$ (for any finite $\ft_*^\sharp$) and on a neighborhood $\cV\subset{}^\etbop S^*M$ of $\pa\cR^+_{\cH^+}\cup\pa\cR^+_{\pa\cK^+,{\rm in}}$ (and $\pa\cR^+_{\scri^+,{\rm in},+}$, which, however, plays no role in the current step). Recall the unstable trapped set $\bar\Gamma^{\rm u}\subset{}^\etbop S^*_{(\cK^+)^\circ}M$ from~\eqref{EqTs3bGamma}. We claim that there exists a neighborhood $\cU\subset{}^\etbop S^*M$ of $\pa\Gamma\subset{}^\etbop S^*_{(\cK^+)^\circ}M$ with the following property: every past-directed null-bicharacteristic starting at a point $\varpi\in\cU\cap\pa\Sigma^+\setminus\bar\Gamma^{\rm u}$ either enters $\cV$ or $\{\ft_*\leq\ft_*^\sharp\}$ in finite time. If $\ft_*(\varpi)<\infty$, the latter happens. If $\varpi$ lies over $(\cK^+)^\circ$, then in the terminology of Proposition~\ref{PropTse3bDyn}, the future-directed null-bicharacteristic $\gamma$ starting at $\varpi$ is of type $(A,B)$ where $A$ cannot be $\pa\Gamma$ since $\varpi\notin\bar\Gamma^{\rm u}$, and both $A,B$ must be subsets of ${}^\etbop S^*_{\cK^+}M$. Therefore $A$ must be one of $\pa\cR^+_{\cH^+}$, $\pa\cR^+_{\pa\cK^+,{\rm in}}$, and hence $\gamma$ enters $\cV$ in the backward direction, indeed.

We can therefore apply Proposition~\ref{PropR3Tr}: the a priori control term $E u$ in~\eqref{EqR3TrEst} is controlled by propagation of regularity starting from points where previous steps already established microlocal $H_\etbop^{\sfs+1,(2\alpha_\sscri,\alpha_+,\alpha_\cK)}$-bounds on $u$. Thus, we conclude weighted $H_\etbop^\sfs$-bounds for $u$ microlocally near $\pa\Gamma$.

\medskip

\pfstep{Step~5. Radial set $\#4$.} We claim that now the a priori control term $E u$ in Proposition~\ref{PropR3RKO} is controlled by the previous steps and the propagation of regularity. This follows from a simple notational modification of the argument in the previous step: we only consider a null-bicharacteristic lying over $\cK^+\cup\iota^+$ and starting in a punctured neighborhood of $\pa\cR^+_{\pa\cK^+,{\rm out}}$ but not in $\pa\cW^+_{\rm out}$: in the notation of~\eqref{EqTse3bDynTable}. If it lies over $\cK^+\cap\iota^+$, then it is of type $(\pa\cR^+_{\pa\cK^+,{\rm in}},\pa\cR^+_{\pa\cK^+,{\rm out}})$ and are thus controlled from Step~3. If it lies over $(\cK^+)^\circ$, then due to the sink nature of $\pa\cR^+_{\pa\cK^+,{\rm out}}$ over $\cK^+$, it must be of type $(A,B)$ with $B=\pa\cR^+_{\pa\cK^+,{\rm out}}$, and thus the only possibilities for $A$ are $\pa\cR^+_{\cH^+}$, $\pa\cR^+_{\pa\cK^+,{\rm in}}$, and $\pa\Gamma$ by~\eqref{EqTse3bDynTable}, near each of which we have control by Step~2, 3, and 4, respectively. Finally, if it lies over $(\iota^+)^\circ$, then $A$ can neither be $\pa\cR^+_{\pa\cK^+,{\rm out}}$ (since the null-bicharacteristic does not start on the unstable manifold $\pa\cW^+_{\rm out}$ of $\pa\cR^+_{\pa\cK^+,{\rm out}}$ over $\iota^+$) nor $\pa\cR^+_{\pa\cK^+,{\rm in}}$ (since the latter set is a repeller for the backwards flow over $\iota^+$), so we must have $A=\pa\cR^+_{\scri^+,{\rm in},+}$, in a neighborhood of which we do have control on $u$ by Step~2. We therefore obtain microlocal $H_\etbop^{\sfs,(2\alpha_\sscri,\alpha_+,\alpha_\cK)}$-control of $u$ near $\pa\cR^+_{\pa\cK^+,{\rm out}}$.

\medskip

\pfstep{Step~6. Radial set $\#5$.} We finally propagate the obtained control into $\pa\cR^+_{\scri^+,{\rm out}}$: every backward null-bicharacteristic $\gamma$ starting in a sufficiently small punctured neighborhood thereof must of be of type $(A,B)$ where $B=\pa\cR^+_{\scri^+,{\rm out}}$ (due to the sink nature of this set) and thus either $A=\pa\cR^+_{\scri^+,{\rm in},+}$, $A=\pa\cR^+_{\pa\cK^+,{\rm out}}$, or $\gamma$ crosses $\ft_*^{-1}(1)$ in finite time. Thus, the term $E u$ in Proposition~\ref{PropR3RScriO} (with $t_{*,-}=t_{*,0}$ and $t_{*,+}=\infty$) is controlled by the previous steps and the propagation of regularity.

\pfstep{Conclusion.} We can propagate $H_\etbop^{\sfs,(2\alpha_\sscri,\alpha_+,\alpha_\cK)}$-control on $u$ also into the black hole to every compact subset of $r>\bhm$. This (together with microlocal elliptic estimates) gives full control on $\chi u$ in $H_\etbop^{\sfs,(2\alpha_\sscri,\alpha_+,\alpha_\cK)}$ and proves the estimate~\eqref{EqR3Est} (in the strong form, since all individual steps entailed the propagation of regularity, too). The proof of Proposition~\ref{PropR3} is complete.

\subsection{Relative b-regularity and tame estimates}
\label{SsRb}

We now prove estimates for weakly admissible wave-type operators $P$ of class $((d_0;k),(2\ell_\sscri,\ell_+,\ell_\cK))$ (Definition~\ref{DefSDWAdm}), i.e., their coefficients have $k$ additional degrees of b-regularity. As in~\S\ref{SsR3}, we only need to assume $\ell_\sscri\in(0,\frac12]$ and $\ell_+,\ell_\cK>0$ here. We use the notation of Definition~\ref{DefSDWAdmNorm}.

\begin{prop}[Global $(\etbop;\bop)$-estimate]
\label{PropRb}
  We use weights and order functions $\alpha_\sscri,\alpha_+,\alpha_\cK,\sfs_0,\sfs$ satisfying the same assumptions as in Proposition~\usref{PropR3}, and we recall $\chi\in\CI(\Omega_*)$ from~\eqref{EqRCutoff}. Then for all $k\in\N_0$ there exists a constant $C_k$ such that
  \begin{equation}
  \label{EqRb}
    \|\chi u\|_{H_{\etbop;\bop}^{(\sfs;k),(2\alpha_\sscri,\alpha_+,\alpha_\cK)}} \leq C_k\Bigl( \|P u\|_{H_{\etbop;\bop}^{(\sfs;k),(2\alpha_\sscri+2,\alpha_++2,\alpha_\cK)}(\Omega_*;\cE)^{\bullet,-}} + \|u\|_{H_{\etbop;\bop}^{(\sfs_0;k),(2\alpha_\sscri,\alpha_+,\alpha_\cK)}(\Omega_*;\cE)^{\bullet,-}} \Bigr)
  \end{equation}
  holds for $P=P_0$ in the strong sense, and also for all weakly admissible wave-type operators $P$ relative to $P_0$ which are sufficiently small perturbations as measured in the norm $\|\cdot\|_{(d_0;k),(2\ell_\sscri,\ell_+,\ell_\cK)}$, where $d_0\in\N$ is sufficiently large and depends only on $\sfs_0,\sfs$. Furthermore, there exists $d\in\N$ such that for $P=P_0$ and all of its sufficiently small perturbations as measured in $\|\cdot\|_{(d_0;d),(2\ell_\sscri,\ell_+,\ell_\cK)}$, one has, for all $k\in\N_0$, the b-tame estimate
  \begin{equation}
  \label{EqRbTame}
  \begin{split}
    &\|\chi u\|_{H_{\etbop;\bop}^{(\sfs;k),(2\alpha_\sscri,\alpha_+,\alpha_\cK)}} \\
    &\quad \leq C_k\biggl( \|P u\|_{H_{\etbop:\bop}^{(\sfs;k),(2\alpha_\sscri+2,\alpha_++2,\alpha_\cK)}(\Omega_*;\cE)^{\bullet,-}} + \|u\|_{H_{\etbop;\bop}^{(\sfs_0;k),(2\alpha_\sscri,\alpha_+,\alpha_\cK)}(\Omega_*;\cE)^{\bullet,-}} \\
    &\quad \quad \hspace{3em} + \|P-P_0\|_{(d_0;k),(2\ell_\sscri,\ell_+,\ell_\cK),\Omega_*} \\
    &\quad \quad \hspace{4.5em} \times \Bigl( \| P u \|_{H_{\etbop;\bop}^{(\sfs;d),(2\alpha_\sscri+2,\alpha_++2,\alpha_\cK)}(\Omega_*;\cE)^{\bullet,-}} + \|u\|_{H_{\etbop;\bop}^{(\sfs_0;d),(2\alpha_\sscri+2,\alpha_++2,\alpha_\cK)}(\Omega_*;\cE)^{\bullet,-}} \Bigr) \biggr),
  \end{split}
  \end{equation}
  provided $P$ is of class $((d_0;k),(2\ell_\sscri,\ell_+,\ell_\cK))$.
\end{prop}

As usual, one uses~\eqref{EqRb} for $k<d$ and~\eqref{EqRbTame} for $k\geq d$. The proof of Proposition~\ref{PropRb} is completely analogous to that of Proposition~\ref{PropR3}, and thus we shall not spell it out. The only difference is that now we use the (b-tame) e3b-microlocal elliptic and real principal type propagation estimates proved in~\S\S\ref{SsMBasic} and~\ref{SssMTame} as well as (b-tame) versions of the radial point and trapping estimates of~\S\S\ref{SssR3R}--\ref{SssR3Tr} which are proved below. After preliminary considerations regarding commutators with b-vector fields near $\scri^+$ in~\S\ref{SssRbComm}, which will be used also in~\S\ref{SE}, we prove the required radial point and trapping estimates in~\S\S\ref{SssRbR}--\ref{SssRbTr} (see Propositions~\ref{PropRbRH} and \ref{PropRbTr}).

\subsubsection{Commutators with b-vector fields near \texorpdfstring{$\scri^+$}{null infinity}}
\label{SssRbComm}

We now work near $\scri^+\cap\iota^+$, and thus only record the edge-b-nature of operators for notational simplicity. Since near $\scri^+\setminus\iota^+$ we will use (non-microlocal) energy estimates (see~\S\ref{SE}), we need to describe commutators of wave-type operators with b-vector fields using only vector fields (i.e., without using proper \emph{pseudo}differential operators).

In the coordinates~\eqref{EqCMCoordscriip} near $\scri^+\cap\iota^+$, set
\begin{subequations}
\begin{equation}
\label{EqRbCommV}
  V_1 := x_\sscri\pa_{x_\sscri}-2\rho_+\pa_{\rho_+} = 2 t_*\pa_{t_*},\quad
  V_2 := x_\sscri\pa_{x_\sscri} = -2 r\pa_r,
\end{equation}
and pick vector fields $V_3,V_4,V_5\in\cV(\Sph^2)$ spanning $\cV(\Sph^2)$ over $\CI(\Sph^2)$ (e.g., rotation vector fields). (For comparison with the discussion above, note that $\pa_{t_*}=\frac12\rho_+(x_\sscri\pa_{x_\sscri}-2\rho_+\pa_{\rho_+})$.) Observe that, in the notation of Corollary~\ref{CorCTebComm} (but writing $\eop$ instead of $\ebop$ here)
\[
  V_i \in \cV_{\bop,[\eop]},
\]
i.e., commutators of $V_i$ with edge vector fields are, again, edge vector fields. Recalling $\nabla^\cE$ from~\eqref{EqSSAdmDer}, set then
\begin{equation}
\label{EqRbCommX}
  X_i := \nabla^\cE_{V_i}\ (i=1,\ldots,5);\quad
  X_1,X_2\in\Diffeb^1(M;\cE),\ X_a\in\Diffb^1(M;\cE)\ (a=3,4,5).
\end{equation}
\end{subequations}
Let us introduce $X_6:=I$ and write
\begin{equation}
\label{EqRbCommXRules}
  \vec X = (X_1,\ldots,X_5,X_6),\quad \vec X^\alpha = X_1^{\alpha_1}\cdots X_5^{\alpha_5}X_6^{\alpha_6},\ \alpha\in\N_0^6,\quad
  {\rm ad}_{\vec X}^\alpha = {\rm ad}_{X_1}^{\alpha_1}\cdots{\rm ad}_{X_5}^{\alpha_5}{\rm ad}_{X_6}^{\alpha_6}.
\end{equation}
(Commutators with $X_6$ vanish, of course.) We moreover introduce
\[
  \vec X_{\rm shift} := x_\sscri^2\rho_+^2 \vec X x_\sscri^{-2}\rho_+^{-2} = (X_1+2,X_2-2,X_3,X_4,X_5,X_6).
\]
Since $[V_1,V_2]=[V_1,V_a]=[V_2,V_a]=0$ for $a=3,4,5$, we have
\begin{equation}
\label{EqRbCommXComm}
  [X_1,X_2],\ [X_1,X_a],\ [X_2,X_a] \in x_\sscri\CI;
\end{equation}
for the second commutator here we use that $\nabla^\cE_{\rho_+\pa_{\rho_+}}$ is, at $\scri^+$, differentiation along the fibers of $\scri^+$ (and independent of the underlying choice of connection on $\cE_X$). 

\begin{lemma}[Commutators near $\scri^+$]
\label{LemmaRbComm}
  Let $k\in\N$, and write $\cA_k:=\{\alpha\in\N_0^6\colon|\alpha|=k\}$. If $P u=f$, then upon setting $u^{(k)}:=(\vec X^\alpha u)_{\alpha\in\cA_k}$ and $f^{(k)}:=(\vec X_{\rm shift}^\alpha f)_{\alpha\in\cA_k}$, we have
  \begin{equation}
  \label{EqRbCommEq}
  \begin{split}
    P^{(k)}u^{(k)} &= f^{(k)} + \tilde f^{(k)},\quad \tilde f^{(k)} = \Biggl(\; \sum_{\substack{\beta+\gamma=\alpha \\ |\beta|\geq 2}} c_{\beta\gamma} x_\sscri^2\rho_+^2 \bigl({\rm ad}_{\vec X}^\beta(x_\sscri^{-2}\rho_+^{-2}P)\bigr) \vec X^\gamma u + \tilde f^{(k),0}_\alpha \Biggr)_{\alpha\in\cA_k}, \\
    &\qquad \tilde f^{(k),0}_\alpha=\sum_{i,\gamma} x_\sscri^2\rho_+^2\bigl( c_\alpha [X_i,p_0] + c'_\alpha([X_i,p_1]X_1 + p_1[X_i,X_1])\bigr)\vec X^\gamma u,
  \end{split}
  \end{equation}
  where the $c_{\beta\gamma},c_\alpha,c'_\alpha\in\N_0$ are combinatorial constants, $\alpha=\gamma+(\delta_{i j})_{j=1,\ldots,6}$, and we recall $[X_i,X_1]\in x_\sscri\CI$ from~\eqref{EqRbCommXComm}; and the $(\alpha,\alpha')$ matrix element of $P^{(k)}$ is of the form
  \begin{equation}
  \label{EqRbCommLot}
  \begin{split}
    &P^{(k)}_{\alpha\alpha'} = \delta_{\alpha\alpha'}P + P^{\sharp,(k)}_{\alpha\alpha'} + \tilde P^{(k)}_{\alpha\alpha'}, \\
    &\qquad P^{\sharp,(k)}_{\alpha\alpha'} \in x_\sscri\cdot x_\sscri^2\rho_+^2\Diffeb^1(M;\cE),\quad \tilde P^{(k)}_{\alpha\alpha'}\in x_\sscri^2\rho_+^2\cC_{\ebop;\bop}^{(d_0;k-1),(2\ell_\sscri,\ell_+)}\Diffeb^1(M;\cE);
  \end{split}
  \end{equation}
  here $\delta_{\alpha\alpha'}$ is the Kronecker delta, and if $\tilde P$ is defined by~\eqref{EqSDWAdmOp}, then the coefficients of $\tilde P^{(k)}_{\alpha\alpha'}$ arise by differentiation of those of $\tilde P$ along a b-differential operator of order $\leq 1$, while those of $P^{\sharp,(k)}_{\alpha\alpha'}$ only depend on $P_0$ (in the notation of Definition~\usref{DefSDWAdm}).
\end{lemma}

\begin{rmk}[Schematic content, I: $P^{(k)}_{\alpha\alpha'}$]
\label{RmkRbCommSchema1}
  The equation for $u^{(k)}$ is a wave-type equation; the diagonal leading order part of $P^{(k)}$ is given by $P$ on $|\cA_k|$ many copies of $\cE$, while the terms $P^{\sharp,(k)}_{\alpha\alpha'}$ and $\tilde P^{(k)}_{\alpha\alpha'}$ are of lower order not only in the differential sense (thus do not affect the underlying Lorentzian metric) but also \emph{in the decay sense at $\scri^+$} (thus do not affect the leading order description~\eqref{EqSDWAdmOp}). The latter feature is crucial because unmodified leading order terms of $P^{(k)}$ at $\scri^+$ mean that the decay rates for solutions of $P$ and $P^{(k)}$ are the same (i.e., threshold quantities in energy and radial point estimates are equal for $P$ and $P^{(k)}$). The origin of the structure~\eqref{EqRbCommLot} is thus that the operators $X_1,\ldots,X_5,X_6$ commute with the leading order part of $P$ in an appropriate sense.\footnote{For the term $x_\sscri^2\slDelta$, see the discussion of~\eqref{EqRbCommSph}.}
\end{rmk}

\begin{rmk}[Schematic content, II: $\tilde f^{(k)}$]
\label{RmkRbCommSchema2}
  Much as in the proof of Proposition~\ref{PropMTamePr} (see the discussion around~\eqref{EqMTamePrbToe3b}), the term $\tilde f^{(k)}$ will be controlled via induction in $k$. It is schematically (and using Notation~\ref{NotMTameb}) of the form $(D_{[\bop]}^q p)D_\ebop^2 D_{[\bop]}^{k-q}u$ where $p$ denotes a coefficient of $P$ (expanded into a frame of eb-vector fields) and $2\leq q=|\beta|\leq k$ as well as $0\leq k-q=|\gamma|\leq k-2$, except for the term $\tilde f_\alpha^{(k),0}$ which is of the form $(D_{[\bop]}^{\leq 1}p)D_\ebop^{\leq 1}D_{[\bop]}^{k-1}u$; therefore, these terms are of the form
  \begin{equation}
  \label{EqRbCommSchema}
    (D_\bop^{\leq 1}p)D_\bop^{\leq k-1}(D_\ebop^{\leq 1}u), \quad (D_\bop^{\leq j-2}D_\bop^{\leq 2}p)D_\bop^{\leq(k-2)-(j-2)}(D_\ebop^{\leq 2}u),\ 2\leq j\leq k.
  \end{equation}
  (We commuted $D_\ebop$ through $D_{[\bop]}$ and used that $D_\ebop^j D_{[\bop]}^k=D_{[\bop]}^{\leq k}D_\ebop^{\leq j}$, and then switched back to arbitrary b-derivatives for better readability.) Weakening the second half of~\eqref{EqRbCommSchema} to
  \begin{equation}
  \label{EqRbCommSchemaWeak}
    (D_\bop^{\leq j-2}D_\bop^{\leq 2}p)D_\bop^{\leq (k-1)-(j-2)}(D_\ebop^{\leq 1}u),\quad 2\leq j\leq k,
  \end{equation}
  each term thus involves at most $k-1$ b-derivatives of $u$ and one eb-derivative (the latter being the amount of regularity that $P^{-1}f$ gains over $f$, given the wave nature of $P$ as an edge-b-differential operator). We use here that spaces of eb- and b-differential operators are invariant under conjugations by products of powers of boundary defining functions (here $x_\sscri^2\rho_+^2$).
\end{rmk}

\begin{proof}[Proof of Lemma~\usref{LemmaRbComm}]
  Motivated by the weights in~\eqref{EqSDWAdmNearScri}, we write
  \[
    x_\sscri^{-2}\rho_+^{-2}P \vec X^\alpha u=\vec X^\alpha(x_\sscri^{-2}\rho_+^{-2}P u)-[\vec X^\alpha,x_\sscri^{-2}\rho_+^{-2}P]u,
  \]
  and therefore
  \begin{equation}
  \label{EqRbCommEqPf}
    P \vec X^\alpha u = x_\sscri^2\rho_+^2 \vec X^\alpha x_\sscri^{-2}\rho_+^{-2} f - \frac12 x_\sscri^2\rho_+^2[\vec X^\alpha,2 x_\sscri^{-2}\rho_+^{-2}P]u.
  \end{equation}
  Let us abbreviate
  \begin{align*}
    P' := 2 x_\sscri^{-2}\rho_+^{-2}P &= -\bigl(X_2-2(1+p_1)\bigr)X_1 + 2 x_\sscri^2\slDelta + p_0 + \tilde P', \\
      &\quad\hspace{6em} \tilde P' \in \cC_{\ebop;\bop}^{(d_0;k),(2\ell_\sscri,\ell_+)}\Diffeb^2(M;\cE).
  \end{align*}
  Lemma~\ref{LemmaMTameComm} gives
  \begin{equation}
  \label{EqRbCommPp}
    [\vec X^\alpha,P'] = -\sum_{\substack{\beta+\gamma=\alpha \\ |\beta|\geq 1}} 2 c_{\beta\gamma}({\rm ad}_{\vec X}^\beta P')\vec X^\gamma
  \end{equation}
  where the $c_{\beta\gamma}$ are combinatorial constants; the normalization ensures that they are the same as in~\eqref{EqRbCommEq}. The terms with $|\beta|\geq 2$ yield the term $\tilde f^{(k)}_\alpha$, and it remains to consider the terms with $|\beta|=1$. The commutator with $X_i\in\Diffb^1$ with $\tilde P'$ lies in $\cC_{\ebop;\bop}^{(d_0;k-1),(2\ell_\sscri,\ell_+)}\Diffeb^2$, as follows from Lemma~\ref{LemmaCTebComm}; since $V_1,V_2$ and $x_\sscri V_a$, $a=3,4,5$, span $\Veb$, we can express $[X_i,\tilde P']=\sum_{j=1}^6 Q_{i j}X_j$ for some $Q_{i j}\in\cC_{\ebop;\bop}^{(d_0;k-1),(2\ell_\sscri,\ell_+)}\Diffeb^1(M;\cE)$.

  Next, we study the leading order term of $P'$. By~\eqref{EqRbCommXComm}, we see that for all $i=1,\ldots,6$, the commutator $[X_i,X_2 X_1]\in x_\sscri\Diffb^1$ can be written as
  \begin{equation}
  \label{EqRbCommxDb1}
    \sum_{j=1}^6 Q_j X_j,\quad Q_j\in x_\sscri\Diffeb^1,
  \end{equation}
  in fact with $Q_j\in x_\sscri\CI$. Furthermore, writing $\slDelta$ modulo $x_\sscri\Diffeb^2$ as a sum of terms $X_a X_b$ (with $\CI$ coefficients) where $a,b\in\{3,4,5\}$, we note that $[X_1,x_\sscri^2 X_a X_b]\in x_\sscri^2\CI\cdot X_a X_b+x_\sscri^2 x_\sscri\Diffb^1$; the first term is $ x_\sscri(x_\sscri\CI\cdot X_a) X_b$, so this can be written in the form~\eqref{EqRbCommxDb1} again; similarly for $[X_2,x_\sscri^2 X_a X_b]$. For $c\in\{3,4,5\}$, we have
  \begin{equation}
  \label{EqRbCommSph}
    [X_c,x_\sscri^2 X_a X_b]=x_\sscri\bigl(x_\sscri[X_c,X_a]\circ X_b+x_\sscri[X_c,X_b]\circ X_a+x_\sscri[X_a,[X_c,X_b]]),
  \end{equation}
  which is thus a sum of terms of the form $x_\sscri\Diffeb^1\circ X_d$, $d=3,4,5,6$; so also this commutator can be expressed in the form~\eqref{EqRbCommxDb1}.

  We finally study the lower order (in the differential sense) terms $2 p_1 X_1$ and $p_0$ of $P'$: if $\beta=(\delta_{i j})_{j=1,\ldots,6}$, then the $\beta$-th term in~\eqref{EqRbCommPp} contributes $x_\sscri^2\rho_+^2 [X_i,p_0]\vec X^\gamma u=x_\sscri^2\rho_+^2[X_i,p_0]\vec X^\gamma u$ times a constant to~\eqref{EqRbCommEqPf}, with $|\gamma|=k-1$. For the term $2 p_1 X_1$, we note $[X_i,p_1 X_1]=[X_i,p_1]X_1+p_1[X_i,X_1]$; but $X_1\in\Diffeb^1$ and $[X_i,X_1]\in x_\sscri\CI$.
\end{proof}

\subsubsection{Radial point estimate near \texorpdfstring{$\cR^+_{\cH^+}$}{the event horizon}}
\label{SssRbR}

The statements of radial point estimates on $H_{\etbop;\bop}$-spaces are similar to the statement Proposition~\ref{PropMTamePr}, but with $B,E,G$ as in the results of~\S\ref{SssR3R}: in order to ensure the boundedness of the e3b-microlocalizers on $H_{\etbop;\bop}$-spaces, we now need to work with spaces $\CI_\bop\Psi_\etbop$ of e3b-ps.d.o.s that have b-regular coefficients. For estimates away from $\scri^+$ and $\cK^+$, we drop the subscripts ``$\eop$'' and ``$3$,'' respectively. (We could also write ``$\cuop$'' instead of ``$\tbop$'' throughout this section since we work near compact subsets of $(\cK^+)^\circ$, but do not do so for consistency with the notation in~\S\ref{SsTs3b}.)

We proceed to prove an analogue of Proposition~\ref{PropR3RH} (concerning $\pa\cR^+_{\pa\cH^+}$). We use a variant of an idea introduced in \cite[\S{4.3.1}]{HintzGlueLocII}: it suffices to commute $t_*\pa_{t_*}$ and $\pa_r$ through the equation $P u=f$ since $\pa_r$ is microlocally elliptic (as a 3b-vector field) near $\pa\cR^+_{\pa\cH^+}$. Thus, 3b-microlocal control of both $\pa_r u$ and $t_*\pa_{t_*}u$ near $\pa\cR^+_{\pa\cH^+}$ gives microlocal control of $u$ in (a weighted version of) $H_{\etbop;\bop}^{(s;1)}$, similarly for higher b-regularity orders. The details are as follows. In order to accommodate the bundle $\cE$, fix a smooth affine connection $\nabla^\cE$ on $\cE$ (such as~\eqref{EqSSAdmDer}) and define
\[
  X_1:=\nabla^\cE_{\pa_r},\quad X_2:=\nabla^\cE_{t_*\pa_{t_*}},\quad X_3:=I
\]
in the coordinates $t_*=t-r$, $r$, and $\omega\in\Sph^2$, write $\vec X=(X_1,X_2,X_3)$, and define $\vec X^\alpha,{\rm ad}_{\vec X}^\alpha$ analogously to~\eqref{EqRbCommXRules}.

\begin{lemma}[Mixed norms using $X_1$ and $X_2$]
\label{LemmaRbRMix}
  Let $\cU\subset\Stb^*M$ be a sufficiently small neighborhood of $\pa\cR^+_{\cH^+}$. Let $K\subset M$ be a compact set containing the base projection of $\cU$ (thus, $K$ is a neighborhood of $\cK^+\cap\cH^+=\{r=r_+,\,t_*^{-1}=0\}$). Suppose the operator wave front sets of $B,\tilde B\in\CI_\bop\Psi_\tbop^0$ are contained in $\cU$, with $\tilde B$ elliptic on $\WF'_\tbop(B)$, and suppose that the Schwartz kernels of $B,\tilde B$ are supported in $K\times K$. Let $\tilde\chi\in\CI(M)$ be supported near and equal to $1$ near $K$. Then for all orders $s,s_0$ and $k\in\N_0$, there exists a constant $C$ such that
  \begin{align}
  \label{EqRbRMix1}
    \|B u\|_{H_{\tbop;\bop}^{(s;k)}} &\leq C\biggl( \|\tilde B u\|_{H_{\tbop;\bop}^{(s;k-1)}} + \sum_{|\alpha|=k} \|B\vec X^\alpha u\|_{H_\tbop^s} + \|\tilde\chi u\|_{H_{\tbop;\bop}^{(s_0;k)}} \biggr), \\
  \label{EqRbRMix2}
    \sum_{|\alpha|=k}\|B\vec X^\alpha u\|_{H_\tbop^s} &\leq C\Bigl( \|\tilde B u\|_{H_{\tbop;\bop}^{(s;k)}} + \|\tilde\chi u\|_{H_{\tbop;\bop}^{(s_0;k)}} \Bigr).
  \end{align}
\end{lemma}
\begin{proof}
  Since the space of b-vector fields is spanned by the space of 3b-vector fields and $t_*\pa_{t_*}$, we have (using schematic notation analogous to Notation~\ref{NotMTameb})
  \[
    \|B u\|_{H_{\tbop;\bop}^{(s;k)}} \sim \sum_{j=0}^k \| D_\tbop^{\leq k-j} X_2^j B u \|_{H_\tbop^s}
  \]
  (i.e., the left-hand side is bounded by a constant multiple times the right-hand side and vice versa). Let now $Q\in\CI_\bop\Psi_\tbop^{-(k-j)}$ be a microlocal elliptic parametrix for $X_1^{k-j}$ on $\bar\cU$, so $I=Q X_1^{k-j}+R$ where $R\in\CI_\bop\Psi_\tbop^0$, $\WF'_\tbop(R)\cap\bar\cU=\emptyset$; we can make the Schwartz kernels of $Q,R$ have support in any fixed neighborhood (as measured by b-unit cells) of the diagonal. Given $A\in\Diff_\tbop^{k-j}$, we then have
  \begin{equation}
  \label{EqRbRMixPx}
    A X_2^j B u = A Q X_1^{k-j}X_2^j B u + A R X_2^j B u.
  \end{equation}
  Consider the first term. Since $A Q\in\CI_\bop\Psi_\tbop^0$ is bounded on $\Htb^s$, we have
  \[
    \|A Q X_1^{k-j}X_2^j B u\|_{\Htb^s} \leq C\|X_1^{k-j}X_2^j B u\|_{\Htb^s}.
  \]
  We then commute $X_1^{k-j}X_2^j$ through $B$ and use Lemma~\ref{LemmaMTameComm} to describe the commutator; this gives
  \[
    \|X_1^{k-j}X_2^j B u\|_{H_\tbop^s} \leq C\biggl( \| B X_1^{k-j}X_2^j u\|_{H_\tbop^s} + \sum_{\substack{|\beta|+|\gamma|=k \\ |\beta|\leq k-1}} \|\vec X^\beta({\rm ad}_{\vec X}^\gamma B)u\|_{H_\tbop^s} \biggr).
  \]
  But $\|\vec X^\beta({\rm ad}_{\vec X}^\gamma B)u\|_{H_\tbop^s}\leq C\|({\rm ad}_{\vec X}^\gamma B)u\|_{H_{\tbop;\bop}^{(s;k-1)}}$. Since by~\eqref{EqMUe3bComm}, (iterated) commutators of $X_1$ and $X_2$ with $B$ are still of class $\CI_\bop\Psi_\tbop^0$ and have operator wave front sets contained in $\WF'_\tbop(B)$, this can in turn be estimated using a microlocal elliptic parametrix for $\tilde B$ by a constant times $\|\tilde B u\|_{H_{\tbop;\bop}^{(s;k-1)}}+\|\tilde\chi u\|_{H_{\tbop;\bop}^{(s_0;k-1)}}$.

  For the second term in~\eqref{EqRbRMixPx}, we commute $X_2^j$ through $B$ and again use Lemma~\ref{LemmaMTameComm} to describe the commutator. The fact that $\WF'_\tbop(R)\cap\WF'_\tbop(B)=\emptyset$ now implies that we can estimate $\|A R X_2^j B u\|_{H_\tbop^s}\leq C\sum_{j'=0}^j\|\chi^\flat X_2^{j'}u\|_{H_{\tbop;\bop}^{s_0}}$ where $\chi^\flat\in\CI(M)$ with $\chi^\flat=1$ near $K$ and $\supp\chi^\flat\subset\tilde\chi^{-1}(1)$; and then we can commute $X_2^{j'}$ through $\chi^\flat$ to further bound this by $C'\|\tilde\chi u\|_{H_{\tbop;\bop}^{(s_0;k)}}$. This concludes the proof of~\eqref{EqRbRMix1}.

  To prove~\eqref{EqRbRMix2}, we commute $\vec X^\alpha$ through $B$ and thus bound
  \[
    \|B\vec X^\alpha u\|_{H_\tbop^s} \leq \|\vec X^\alpha B u\|_{H_\tbop^s} + C \sum_{\substack{|\beta|+|\gamma|=k \\ |\beta|\leq k-1}} \|\vec X^\beta({\rm ad}_{\vec X}^\gamma B)u \|_{H_\tbop^s}.
  \]
  The first term is bounded by $\|B u\|_{H_{\tbop;\bop}^{(s;k)}}$ and thus, using a microlocal elliptic parametrix for $\tilde B$, by $\|\tilde B u\|_{H_{\tbop;\bop}^{(s;k)}} + \|\tilde\chi u\|_{H_{\tbop;\bop}^{(s_0;k)}}$. The same argument applies to all other summands, using again that $\WF'_\tbop({\rm ad}_{\vec X}^\gamma B)\subset\WF'_\tbop(B)$ is contained in the elliptic set of $\tilde B$.
\end{proof}

Now, $t_*\pa_{t_*}$ is a b-vector field with good commutation properties by Lemma~\ref{LemmaCT3bDil}; and commutators with $\pa_r$ lead to weaker threshold 3b-regularities due to the ``enhanced red-shift effect'' \cite{DafermosRodnianskiKerrBoundedness}. More precisely:

\begin{lemma}[Commutators with $t_*\pa_{t_*}$ and $\pa_r$]
\label{LemmaRbRHComm}
  Set $\cA_k:=\{\alpha\in\N_0^3\colon|\alpha|=k\}$. For all sufficiently small neighborhoods $\cU\subset\Stb^*M$ of $\pa\cR^+_{\cH^+}$, the following holds. If $P u=f$ where $P$ is a weakly admissible wave-type operator of class $((d_0;k),(2\ell_\sscri,\ell_+,\ell_\cK))$, then $u^{(k)}:=(\vec X^\alpha u)_{\alpha\in\cA_k}$ satisfies the equation
  \begin{equation}
  \label{EqRbRHCommEq}
  \begin{split}
    &P^{(k)}u^{(k)} = f^{(k)} + \tilde f^{(k)} + \tilde f^{(k)}_{\rm micro}, \\
    &\qquad\qquad f^{(k)} := (\vec X^\alpha f)_{\alpha\in\cA_k},\ \ 
      \tilde f^{(k)} := \Biggl(\;\sum_{\substack{\beta+\gamma=\alpha \\ |\beta|\geq 2}} c_{\beta\gamma} ({\rm ad}_{\vec X}^\beta P)\vec X^\gamma u \Biggr)_{\alpha\in\cA_k},
  \end{split}
  \end{equation}
  where the $c_{\beta\gamma}\in\N_0$ are combinatorial constants; furthermore,
  \begin{enumerate}
  \item\label{ItRbRHComm1} the components of $\tilde f_{\rm micro}^{(k)}$ are sums of terms of the form
    \begin{equation}
    \label{EqRbRHCommMicro}
      (D_\bop^{\leq 1}p) R \vec X^\beta u,\quad |\beta|\leq k-1,\ R\in\CI_\bop\Psi_\tbop^2,\ \WF'_\tbop(R)\cap\bar\cU=\emptyset,
    \end{equation}
    where $p=p^\sharp+\tilde p$, $p^\sharp\in\CI$, $\tilde p\in\rho_\cK^{\ell_\cK}\cC_{\tbop;\bop}^{(d_0;k)}$, denotes an $\End(\cE)$-valued coefficient of $P$ when expressed near $\cK^+\cap\cH^+$ as a sum of terms of the form $p A$, $A\in\Diff_\tbop^2$;\footnote{One can, for example, express $P$ \emph{uniquely} as a linear combination of $I,\nabla^\cE_{\pa_{t_*}},\nabla^\cE_{\pa_{x^i}}$ ($i=1,2,3$), and their up to two-fold products with coefficients $p$ that are sections of $\End(\cE)$.}
  \item\label{ItRbRHComm2} the $(\alpha,\alpha')$ matrix element of $P^{(k)}$ is of the form
    \begin{equation}
    \label{EqRbRHCommP}
      P^{(k)}_{\alpha \alpha'} = \delta_{\alpha\alpha'} ( P + P^{(k),\sharp}_\alpha ) + \tilde P^{(k)}_{\alpha\alpha'},
    \end{equation}
    where
    \begin{itemize}
    \item $P^{(k),\sharp}_\alpha\in\CI_\bop\Psi_\tbop^1$ has principal symbol $-i p^\sharp_\alpha\xi_0$ at $\cR^+_{\cH^+}$ where $p^\sharp_\alpha\in\CI(\Sph^2)$ is a non-negative function on $\cH^+\cap\cK^+=\Sph^2_\omega$, and $\xi_0$ was defined\footnote{and is positive at $\cR^+_{\cH^+}$ by Definition~\ref{DefTs3bH}} in~\eqref{EqTs3bHCoord};
    \item $\tilde P^{(k)}_{\alpha\alpha'}$ is a finite sum of terms of the form $(D_\bop^{\leq 1}\tilde p)\CI_\bop\Psi_\tbop^1$ where $\tilde p\in\rho_\cK^{\ell_\cK}\cC_{\tbop;\bop}^{(d_0;k)}$ denotes a (decaying) coefficient of $P-P_0$ as before, and terms of class $t_*^{-1}\CI_\bop\Psi_\tbop^1$ depending only on $P_0$.
    \end{itemize}
  \end{enumerate}
\end{lemma}

Since we will only need Lemma~\ref{LemmaRbRHComm} microlocally near $\pa\cR^+_{\cH^+}$, the term $\tilde f^{(k)}_{\rm micro}$ should be regarded as an error term: upon applying a 3b-microlocalizer to a neighborhood of $\pa\cR^+_{\cH^+}$ (more precisely, with operator wave front set contained in $\cU$) to it, it becomes regularizing in the 3b-differential order sense. The fact that the imaginary part of $P_\alpha^{(k),\sharp}$ has a sign at $\cR^+_{\cH^+}$ (and that the remaining term $\tilde P^{(k)}_{\alpha\alpha'}$ has decaying coefficients as $t_*\to\infty$) makes the threshold quantity $\vartheta_{\cH^+}^{(k)}$ for $P^{(k)}$ the same as for $P$ itself.\footnote{With the wrong sign, one would need to increase the minimal 3b-regularity order in Proposition~\ref{PropR3RH} with every additional b-derivative one wishes to control.}

\begin{proof}[Proof of Lemma~\usref{LemmaRbRHComm}]
  For $\alpha=(j,l,q)\in\cA_k$, Lemma~\ref{LemmaMTameComm} gives
  \begin{equation}
  \label{EqRbRHCommPrel}
    P\vec X^\alpha u = \vec X^\alpha f + \Bigl( j[X_1,P] X_1^{j-1}X_2^l u + l[X_2,P] X_1^j X_2^{l-1}u\Bigr) + \sum_{\substack{\beta+\gamma=\alpha \\ |\beta|\geq 2}} c_{\beta\gamma} ({\rm ad}_{\vec X}^\beta P) \vec X^\gamma u
  \end{equation}
  for combinatorial constants $c_{\beta\gamma}\in\N_0$. Here, the first, resp.\ second term in parentheses is absent for $j=0$, resp.\ $l=0$.

  \pfstep{Commutators with $\pa_r$.} Let us study the term $[X_1,P]X_1^{j-1}X_2^l u$ when $j\geq 1$. Write $P$ near $\cK^+\cap\{r=r_+\}$ as a sum of terms of the form $p A$ where $p\in\CI+t_*^{-\ell_\cK}\cC_{\tbop;\bop}^{(d_0;k)}$ (valued in endomorphisms of $\cE$), while $A\in\Diff_\tbop^2$ has smooth coefficients. Then
  \[
    [X_1,p A] = [X_1,p] A + p[X_1,A].
  \]
  When $\bar\cU$ is disjoint from the characteristic set of $\pa_r$ near $\pa\cR^+_{\cH^+}$, we can use a microlocal elliptic parametrix of $X_1\in\Diff_\tbop^1(M;\cE)$ near $\bar\cU$ to write
  \begin{equation}
  \label{EqRbRHCommPx}
    I=Q X_1+R,\quad Q\in\CI_\bop\Psi_\tbop^{-1},\ R\in\CI_\bop\Psi_\tbop^0,\ \WF'_\tbop(R)\cap\bar\cU=\emptyset,
  \end{equation}
  and therefore
  \[
    [X_1,p A] = \underbrace{\bigl( [X_1,p] A Q + p[X_1,A] Q \bigr) X_1}_{=:\,{\rm I}} + \underbrace{\bigl( [X_1,p]R A + p [X_1,A]R \bigr)}_{=:\,{\rm II}}.
  \]
  We absorb the term ${\rm II}$, acting on $X_1^{j-1}X_2^l X_3^q u$, into $\tilde f^{(k)}_{\rm micro}$. The term ${\rm I}$ on the other hand is of the schematic form $(D_\bop^{\leq 1}p)\Psi_\ebop^1\circ X_1$ and thus contributes to $j[X_1,P]X_1^{j-1}X_2^l u$ the term
  \[
    j(D_\bop^{\leq 1}p) \Psi_\ebop^1 \circ \vec X^\alpha u
  \]
  (where we inserted $X_3^q=I$ and recalled $\vec X^\alpha=X_1^j X_2^l X_3^q$), which can thus be moved to the left hand side of~\eqref{EqRbRHCommPrel}. The sum of all resulting contributions $-j[X_1,p A]Q$ put on the left hand side (which act on $\vec X^\alpha u$) is $j[P,X_1]Q=i j[P,i^{-1}X_1]Q$. Splitting $P=P_0+(P-P_0)$, the contributions from $P-P_0$ have decaying (as $t_*\to\infty$) coefficients. In the coordinates~\eqref{EqTs3bHCoord}, using~\eqref{EqTs3bHHam}, and recalling that $G_{\bhm,a}=\varrho^{-2}G_\tbop$ over $\cK^+$, we compute the 3b-principal symbol of the contribution $i j[P_0,i^{-1}X_1]Q$ from $P_0$ to be
  \begin{equation}
  \label{EqRbRHCommPrShift}
    i j(H_{G_{\bhm,a}}\xi_0) \xi_0^{-1} = -i j \varrho^{-2}\mu'(r_+) \xi_0\ \text{at}\ \cR^+_{\cH^+}.
  \end{equation}
  In the notation of the statement, we thus set $p_\alpha^\sharp=j\varrho^{-2}\mu'(r_+)$.

  \pfstep{Commutators with $t_*\pa_{t_*}$.} We now study the term $[X_2,P]X_1^j X_2^{l-1}u$ when $l\geq 1$. Write $P$ as a sum of terms $p A$, $p=p^\sharp+\tilde p$, $p^\sharp\in\CI$, $\tilde p\in t_*^{-\ell_\cK}\cC_{\tbop;\bop}^{(d_0;k)}$. Write
  \begin{equation}
  \label{EqRbRHCommX20}
    [X_2,p A] = [X_2,p]A + p[X_2,A].
  \end{equation}
  In the first summand, we insert~\eqref{EqRbRHCommPx} to get, upon acting on $X_1^j X_2^{l-1}u$,
  \[
    \underbrace{[X_2,p]A Q\bigl(X_1^{j+1}X_2^{l-1}X_3^q u\bigr)}_{=:\,{\rm I}} + \underbrace{[X_2,p]A R X_1^j X_2^{l-1}u}_{=:\,{\rm II}}.
  \]
  We absorb the term ${\rm II}$ into $\tilde f^{(k)}_{\rm micro}$. The term ${\rm I}$ is b-principal in that it features $(j+1)+(l-1)+q=|\alpha|=k$ derivatives along $\vec X$, but the coefficients
  \begin{equation}
  \label{EqRbRHCommX2}
    [X_2,p^\sharp] \in t_*^{-1}\CI,\quad
    [X_2,\tilde p] \in t_*^{-\ell_\cK}\cC_{\tbop;\bop}^{(d_0;k-1)}
  \end{equation}
  decay as $t_*\to\infty$. The second membership in~\eqref{EqRbRHCommX2} is clear, and the first one follows by noting that, in a local trivialization of $\cE$, the operator $X_2$ is the sum of $t_*\pa_{t_*}=-\rho_\cK\pa_{\rho_\cK}$ acting component-wise plus a zeroth order term of class $t_*^{-1}\CI$; the commutator of the latter term with $p^\sharp\in\CI$ is of class $t_*^{-1}\CI$ still, while the commutator of the former term with $p^\sharp$ is $-\rho_\cK[\pa_{\rho_\cK},p^\sharp]$ and thus, due to the smoothness of $p^\sharp$ in $\rho_\cK$, of class $\rho_\cK\CI$ indeed. (In fact, in the coordinates $(t_*,r,\omega)$, $p^\sharp$ is \emph{independent} of $t_*$, so $[t_*\pa_{t_*},p^\sharp]=0$.)

  The second summand on the right-hand side of~\eqref{EqRbRHCommX20} can be expressed using Lemma~\ref{LemmaCT3bDil} as
  \[
    p[X_2,A] = p t_*^{-1}A^\flat X_2 + t_*^{-1}p A^\sharp,\quad A^\flat\in\Diff_\tbop^1,\ A^\sharp\in\Diff_\tbop^2.
  \]
  Acting on $X_1^j X_2^{l-1}u$, the first term is b-principal but with coefficient of class $t_*^{-1}\CI+t_*^{-1-\ell_\cK}\cC_{\tbop;\bop}^{(d_0;k)}$. We expand the second term, acting on $X_1^j X_2^{l-1}u$, using~\eqref{EqRbRHCommPx} into
  \[
    t_*^{-1}p A^\sharp Q\bigl(X_1^{j+1}X_2^{l-1}u\bigr) + t_*^{-1}p A^\sharp R X_1^j X_2^{l-1}u,
  \]
  with the first, resp.\ second summand absorbed into $\tilde P^{(k)}_{\alpha\alpha'}$, resp.\ $\tilde f^{(k)}_{\rm micro}$. This completes the proof.
\end{proof}

\begin{prop}[Radial point estimate near $\pa\cR^+_{\cH^+}$ with b-regularity]
\label{PropRbRH}
  Let $s>s_0>\frac12+\vartheta_{\cH^+}$, $k\in\N_0$. Then there exist operators $B,E,G\in\CI_\bop\Psi_\tbop^0$ satisfying the properties~\eqref{ItIR3RH1}--\eqref{ItIR3RH4} in Proposition~\usref{PropR3RH} such that for all $\alpha_\cK\in\R$, and for all $\tilde\chi\in\CI(M)$ supported near $K$ and equal to $1$ near $K$, we have
  \begin{equation}
  \label{EqRbRH}
    \|B u\|_{H_{\tbop;\bop}^{(s;k),\alpha_\cK}} \leq C_k\Bigl( \|G P u\|_{H_{\tbop;\bop}^{(s-1;k),\alpha_\cK}} + \|E u\|_{H_{\tbop;\bop}^{(s;k),\alpha_\cK}} + \|\tilde\chi u\|_{H_{\tbop;\bop}^{(s_0;k),\alpha_\cK}} \Bigr),
  \end{equation}
  which holds uniformly for $P=P_0$ and all of its sufficiently small perturbations as measured in the norm $\|\cdot\|_{(d_0;k),(0,0,\ell_\cK)}$ when $d_0\in\N$ is sufficiently large (depending only on $s_0,s$). Furthermore, there exists $d\in\N$ such that for $P=P_0$ and all of its sufficiently small perturbations as measured in $\|\cdot\|_{(d_0;d),(0,0,\ell_\cK)}$, one has the b-tame estimate
  \begin{equation}
  \label{EqRbRHTame}
  \begin{split}
    &\|B u\|_{H_{\tbop;\bop}^{(s;k),\alpha_\cK}} \\
    &\qquad \leq C_k\biggl( \|G P u\|_{H_{\tbop;\bop}^{(s-1;k),\alpha_\cK}} + \|E u\|_{H_{\tbop;\bop}^{(s;k),\alpha_\cK}} + \|\tilde\chi u\|_{H_{\tbop;\bop}^{(s_0;k),\alpha_\cK}} \\
    &\qquad \quad \hspace{3em} + \|P-P_0\|_{(d_0;k),(0,0,\ell_\cK),K'} \Bigl( \|G P u\|_{H_{\tbop;\bop}^{(s-1;d),\alpha_\cK}} + \|E u\|_{H_{\tbop;\bop}^{(s;d),\alpha_\cK}} + \|\tilde\chi u\|_{H_{\tbop;\bop}^{(s_0;d),\alpha_\cK}} \Bigr) \biggr)
  \end{split}
  \end{equation}
  where $K'=\supp\tilde\chi$, provided $P$ is of class $((d_0;k),(0,0,\ell_\cK))$.
\end{prop}
\begin{proof}
  \pfstep{Proof of~\eqref{EqRbRH}.} We prove~\eqref{EqRbRH} by induction on $k$, the case $k=0$ being Proposition~\ref{PropR3RH}. Having proved~\eqref{EqRbRH} for $k-1$ in place of $k$, we consider the equation~\eqref{EqRbRHCommEq} for $u^{(k)}$. The operator $P^{(k)}$, acting on the bundle $\cE^{(k)}:=\bigoplus_{\alpha\in\cA_k}\cE$, has the same (scalar) principal symbol $G$ as $P$ itself. If on each summand of $\cE^{(k)}$ we use the same inner product as in Definition~\ref{DefSSOrderAdm}\eqref{ItSSOrderAdm4}, then the threshold quantity~\eqref{EqR3RHThreshold} for $P^{(k)}$ is $\leq\vartheta_{\cH^+}$ (in fact, equal to $\vartheta_{\cH^+}$) since the term $P^{(k),\sharp}_\alpha$ in~\eqref{EqRbRHCommP} gives a non-positive contribution (in fact, no contribution at all in case $\alpha=(j,l,q)$ with $j=0$, cf.\ \eqref{EqRbRHCommPrShift}). We now apply Proposition~\ref{PropR3RH} to $P^{(k)}$ and $u^{(k)}$ in place of $P$ and $u$, with unchanged values of $s,s_0$, and for a sufficiently small neighborhood $\cU$ of $\pa\cR^+_{\cH^+}$ (so that $\pa_r$ is elliptic on $\cU$); the first order pseudodifferential terms of $P^{(k)}$ are straightforwardly handled (see Remark~\ref{RmkR3RHPsdolot}). We claim that this implies~\eqref{EqRbRH} for the value $k$, with slightly adjusted $B,G,E,\tilde\chi$.

  In a nutshell, the term $\tilde f^{(k)}$ in~\eqref{EqRbRHCommEq} involve at most $2$ 3b-derivatives and $k-2$ b-derivatives acting on $u$, so a fortiori $1$ 3b- and $k-1$ b-derivatives, and thus they can be controlled inductively; and the term $\tilde f^{(k)}_{\rm micro}$ is microlocalized away from $\pa\cR^+_{\cH^+}$ and hence only contributes 3b-residual errors (the term $\tilde\chi u$ in~\eqref{EqRbRH}). In more detail, let us consider each of the terms in the estimate~\eqref{EqR3RHEst}; we omit the weight $\alpha_\cK$ for brevity. (The following arguments are similar to those in the proof of Proposition~\ref{PropMTamePr}, but since the details differ, mainly due to the presence of the error term $\tilde f^{(k)}_{\rm micro}$ with special microlocal properties, we shall present the complete proof.)

  First, we use~\eqref{EqRbRMix1} to bound
  \[
    \sum_{|\alpha|=k}\|B\vec X^\alpha u\|_{H_\tbop^s}\gtrsim\|B u\|_{H_{\tbop;\bop}^{(s;k)}}-\|\tilde B u\|_{H_{\tbop;\bop}^{(s;k-1)}}-\|\tilde\chi u\|_{H_{\tbop;\bop}^{(s_0;k)}}.
  \]
  That is, $\sum_{|\alpha|=k}\|B\vec X^\alpha u\|_{H_\tbop^s}$ controls the $H_{\tbop;\bop}^{(s;k)}$-norm of $B u$ up to error terms that are lower order in either the b-regularity or 3b-regularity sense. Next, using~\eqref{EqRbRMix2}, we can estimate $\|E\vec X^\alpha u\|_{H_\tbop^s}$ for $|\alpha|=k$ by $\|\tilde E u\|_{H_{\tbop;\bop}^{(s;k)}}+\|\tilde\chi u\|_{H_\tbop^{(s_0:k)}}$ where $\tilde E$ is elliptic on $\WF'_\tbop(E)$. Moreover, $\|\tilde\chi\vec X^\alpha u\|_{\Htb^{s_0}}\leq C\|\chi^\sharp u\|_{H_{\tbop;\bop}^{(s_0;k)}}$ for any fixed $\chi^\sharp\in\CI(M)$ which equals $1$ on $\supp\tilde\chi$.

  It remains to control, for $|\alpha|=k$, the $H_\tbop^{s-1}$-norms of $G$ applied to the each of the three summands on the right-hand side of~\eqref{EqRbRHCommEq}. If $\chi^\sharp\in\CI(M)$ equals $1$ on $\supp\tilde\chi$, then
  \begin{equation}
  \label{EqRbRH1}
    \|G \vec X^\alpha f\|_{H_\tbop^{s-1}} \leq C\Bigl( \|\tilde G f\|_{H_{\tbop;\bop}^{(s-1;k)}} + \|\tilde\chi f\|_{H_{\tbop;\bop}^{(s_0-2;k)}} \Bigr)
  \end{equation}
  by~\eqref{EqRbRMix2} with $f,s_0-2$ in place of $u,s_0$, and with $\tilde G\in\CI_\bop\Psi_\tbop^0$ elliptic on $\WF'_\tbop(G)$. Write $f=P u$, then $\tilde\chi f=P\tilde\chi\chi^\sharp u+[P,\tilde\chi]\chi^\sharp u$. Using the boundedness of $P\circ\tilde\chi$, $[P,\tilde\chi]\colon H_{\tbop;\bop}^{(s_0;k)}\to H_{\tbop;\bop}^{(s_0-2;k)}$ (which only uses $\cC_{\tbop;\bop}^{(d_0;k)}$-regularity of the coefficients of $P-P_0$, with $d_0$ sufficiently large depending only on $s_0$), we can thus estimate
  \[
    \|\tilde\chi f\|_{H_{\tbop;\bop}^{(s_0-2;k)}} \leq C\|\chi^\sharp u\|_{H_{\tbop;\bop}^{(s_0;k)}}.
  \]

  Turning to $\tilde f^{(k)}$ in~\eqref{EqRbRHCommEq}, we need to bound $\|G({\rm ad}_{\vec X}^\beta P)\vec X^\gamma u\|_{\Htb^{s-1}}$ where $|\gamma|\leq k-2$. But $P^{(\beta)}:={\rm ad}_{\vec X}^\beta P$ is, a fortiori, of class $\Difftb^2+\cC_\tbop^{d_0}\Difftb^2$. Choosing $G^\flat\in\CI_\bop\Psi_\tbop^0$ such that $\WF'_\tbop(G^\flat-I)\cap\WF'_\tbop(G)=\emptyset$, we can estimate
  \begin{align*}
    \|G P^{(\beta)}\vec X^\gamma u\|_{\Htb^{s-1}} &\leq \|G P^{(\beta)}G^\flat\vec X^\gamma u\|_{\Htb^{s-1}} + \|G P^{(\beta)}(I-G^\flat)\vec X^\gamma u\|_{\Htb^{s-1}} \\
      &\leq C\Bigl(\|G^\flat\vec X^\gamma u\|_{\Htb^{s+1}} + \|\tilde\chi u\|_{H_{\tbop;\bop}^{(s_0;k-2)}}\Bigr)
  \end{align*}
  since, for $d_0=\infty$, $G P^{(\beta)}(I-G^\flat)\in\tilde\Psi_\tbop^{-\infty}$, and thus $G P^{(\beta)}(I-G^\flat)\colon\Htb^{s_0}\to\Htb^{s-1}$ is bounded when $d_0$ is large enough. Using~\eqref{EqRbRMix2}, we moreover bound
  \[
    \|G^\flat\vec X^\gamma u\|_{\Htb^{s+1}} \leq C\Bigl( \|\tilde G u\|_{H_{\tbop;\bop}^{(s+1;k-2)}} + \|\tilde\chi u\|_{H_{\tbop;\bop}^{(s_0;k-2)}} \Bigr)
  \]
  if $\WF'_\tbop(G^\flat)\subset\Ell_\tbop(\tilde G)$; and finally we note that $\|\tilde G u\|_{H_{\tbop;\bop}^{(s+1;k-2)}}\leq C\|\tilde G u\|_{H_{\tbop;\bop}^{(s;k-1)}}$.

  Finally, consider $\tilde f^{(k)}_{\rm micro}$ in~\eqref{EqRbRHCommEq}: since $D_\bop^{\leq 1}p$ is, a fortiori, of class $\CI+\cC_\tbop^{d_0}$, and since $\WF'_\tbop(R)\cap\WF'_\tbop(G)=\emptyset$, we have, for $|\beta|\leq k-1$ and using the notation from~\eqref{EqRbRHCommMicro},
  \begin{equation}
  \label{EqRbRH3}
    \|G(D_\bop^{\leq 1}p)R\vec X^\beta u\|_{\Htb^{s-1}} \leq C\|\tilde\chi u\|_{H_{\tbop;\bop}^{(s_0;k-1)}}.
  \end{equation}

  In summary, we have now shown the estimate~\eqref{EqRbRH} for slightly adjusted $B,G,E,\tilde\chi$, except for the presence of an additional term $\|\tilde G u\|_{H_{\tbop;\bop}^{(s;k-1),\alpha_\cK}}$ on the right (where we choose $\tilde G$ to have elliptic set containing $\WF'_\tbop(\tilde B)$ in the above arguments). This term, in turn, can be estimated using the inductive hypothesis (provided the operator wave front sets of all operators in the current step were chosen to be sufficiently small neighborhoods of $\pa\cR^+_{\cH^+}$).

  \pfstep{Proof of~\eqref{EqRbRHTame}.} For $k\leq d$, the estimate~\eqref{EqRbRHTame} is implied by~\eqref{EqRbRH}. For $k>d$, we use induction; so assume the validity of~\eqref{EqRbRHTame} for the value $k-1$. Note that regardless of the value of $k$, the coefficients of the operator $P^{(k)}$ (to which we apply Proposition~\ref{PropR3RH}), as described in Lemma~\ref{LemmaRbRHComm}\eqref{ItRbRHComm2}, are obtained from those of $P$ by applying at most one b-derivative. In order to prove~\eqref{EqRbRHTame}, we only need to revisit the estimate for $\|G f^{(k)}\|_{\Htb^{s-1}}$ (again omitting the weight $\alpha_\cK$ for brevity), as this is the only term involving $P$.

  In~\eqref{EqRbRH1}, we again write $\tilde\chi f=P\tilde\chi\chi^\sharp u+[P,\tilde\chi]\chi^\sharp u$. We only estimate the $H_{\tbop;\bop}^{(s_0-2;k)}$-norm of the first summand, the estimate for the second summand being completely analogous. We use schematic notation as in Notation~\ref{NotMTameb} and consider a term $p D_\tbop^2$, $p=p^\sharp+\tilde p$, $p^\sharp\in\CI$, $\tilde p\in\rho_\cK^{\ell_\cK}\cC_{\tbop;\bop}^{(d_0;k)}$, of $P$. Then $\|p^\sharp D_\tbop^2\tilde\chi\chi^\sharp u\|_{H_{\tbop;\bop}^{(s_0-2;k)}}\leq C\|\chi^\sharp u\|_{H_{\tbop;\bop}^{(s_0;k)}}$ where $C$ only depends on $p^\sharp$ (and thus on $P^0$), while for the contribution of $\tilde p$ we must estimate
  \[
    \bigl\|D_{[\bop]}^{\leq k}\bigl(\tilde p D_\tbop^2(\tilde\chi\chi^\sharp u)\bigr)\bigr\|_{H_\tbop^{s_0-2}}.
  \]
  We use the Leibniz rule to distribute the derivatives $D_{[\bop]}^{\leq k}$ to $\tilde p$ and $D_\tbop^2(\tilde\chi\chi^\sharp u)$, and commute the derivatives falling on the second factor through to $\chi^\sharp$ (and using that commutators of $D_{[\bop]}$ with $D_\tbop$ are again of type $D_\tbop$). We moreover commute $D_\tbop^2$ to the front (which leads to at most $2$ 3b-derivatives falling onto $\tilde p$). We then need to estimate $\|(D_\bop^j D_\tbop^{\leq 2}\tilde p)(D_\bop^{\leq k-j}(\chi^\sharp u))\|_{\Htb^{s_0}}$ for $0\leq j\leq k$, which we do using Proposition~\ref{PropMTameMicr} (with $I,s_0,D_\tbop^{\leq 2}\tilde p,\chi^\sharp u$ in place of $B,s,\ell,u$ there).

  For the term $\tilde f^{(k)}$ in~\eqref{EqRbRHCommEq}, we need to estimate, in schematic notation,
  \begin{equation}
  \label{EqRbRHTame1}
    \| G(D_\bop^j p)D_\tbop^2(D_\bop^{\leq k-j}u) \|_{\Htb^{s-1}},\quad 2\leq j\leq k.
  \end{equation}
  We commute $D_\tbop^2$ through $D_\bop^j p$ (which again leads to at most $2$ 3b-derivatives on $p$) and through $G\in\CI_\bop\Psi_\tbop^0$; for $\tilde G\in\CI_\bop\Psi_\tbop^0$ that is elliptic on $\WF'_\tbop(G)$, it then remains to estimate
  \begin{equation}
  \label{EqRbRHTame2}
    \|\tilde G (D_\bop^j D_\tbop^{\leq 2}p)(D_\bop^{\leq k-j}u) \|_{\Htb^{s+1}}.
  \end{equation}
  Tameness considerations only enter for the contribution of $\tilde p$ to $p$. We use Proposition~\ref{PropMTameMicr} with $D_\bop^2 D_\tbop^{\leq 2}\tilde p$ and $k-2$ in place of $\ell$ and $k$, respectively. Let $G^\sharp\in\CI_\bop\Psi_\tbop^0$ be an operator that is elliptic on $\WF'_\tbop(\tilde G)$, and write $\tilde d_0$ for the value of $d_0$ in Proposition~\ref{PropMTameMicr}, then~\eqref{EqRbRHTame2} is bounded by
  \begin{align*}
    &C_k\biggl( \|\tilde p\|_{\cC_{\tbop;\bop}^{(\tilde d_0+2,\tilde d+2)}} \Bigl( \|G^\sharp u\|_{H_{\tbop;\bop}^{(s+1;k-2)}} + \|\tilde\chi u\|_{H_{\tbop;\bop}^{(s_0;k-2)}} \Bigr) \\
    &\qquad + \|\tilde p\|_{\cC_{\tbop;\bop}^{(\tilde d_0+2;k)}}\Bigl( \|G^\sharp u\|_{H_{\tbop;\bop}^{(s+1;\tilde d)}} + \|\tilde\chi u\|_{H_{\tbop;\bop}^{(s_0;\tilde d)}} \Bigr) \biggr),
  \end{align*}
  which for $d_0\geq\tilde d_0+2$ and $d\geq\tilde d+2$ we can further estimate using
  \begin{equation}
  \label{EqRbRHTame3}
    H_{\tbop;\bop}^{(s+1;k-2)} \subset H_{\tbop;\bop}^{(s;k-1)}
  \end{equation}
  by
  \begin{align*}
    &C_k\biggl( \|\tilde p\|_{\cC_{\tbop;\bop}^{(d_0,d)}} \Bigl( \|G^\sharp u\|_{H_{\tbop;\bop}^{(s;k-1)}} + \|\tilde\chi u\|_{H_{\tbop;\bop}^{(s_0;k-2)}} \Bigr) \\
    &\qquad + \|\tilde p\|_{\cC_{\tbop;\bop}^{(d_0;k)}}\Bigl( \|G^\sharp u\|_{H_{\tbop;\bop}^{(s;d)}} + \|\tilde\chi u\|_{H_{\tbop;\bop}^{(s_0;d)}} \Bigr) \biggr),
  \end{align*}
  If we arrange for $G^\sharp$ to have operator wave front set in a small neighborhood of $\pa\cR^+_{\cH^+}$ still, the norms on $G^\sharp u$ can be estimated by the inductive hypothesis (when $k-1>d$) or by~\eqref{EqRbRH} (when $k-1\leq d$).

  The estimate for the term $\tilde f^{(k)}_{\rm micro}$ given in~\eqref{EqRbRH3} is sufficient: it only uses the sufficiently large degree of 3b-regularity of the coefficients of $D_\bop^{\leq 1}p$, i.e., it does not need high b-regularity of $p$.

  In summary, we have shown~\eqref{EqRbRHTame} for slightly adjusted $B,G,E,\tilde\chi$, except for the presence of an additional term $\|G^\sharp u\|_{H_{\tbop;\bop}^{(s;k-1)}}$ on the right-hand side, which is controlled by the inductive hypothesis.
\end{proof}

\begin{rmk}[Improvement over \cite{HintzVasyScrieb}]
\label{RmkRbBetter}
  In the proof of \cite[Corollary~6.6]{HintzVasyScrieb}, the estimate~\eqref{EqRbRH} for orders $(s;k)$ is shown to follow easily from the case $(s+1;k-1)$. Note, however, that the case $(s+1;k-1)$ would then follow from $(s+2;k-2)$, etc., and ultimately from $(s+k;0)$; the 3b-regularity requirements on the coefficients of $P$ would then correspondingly need to increase with $k$, which would render a proof of a (b-)tame estimate such as~\eqref{EqRbRHTame} along such lines impossible.
\end{rmk}

\subsubsection{Other radial point estimates}
\label{SssRbRO}

We now turn to the proof of generalizations of the other radial point estimates in~\S\ref{SssR3R} to estimates on spaces capturing $k\in\N_0$ degrees of b-regularity (including tame bounds in $k$). The proofs will be completely analogous to the proof of Proposition~\ref{PropRbRH} above once one has an analogue of Lemma~\ref{LemmaRbRHComm}, i.e., a suitable collection of commutator b-vector fields for which one can apply the radial point estimates from~\S\ref{SssR3R} to the commuted equation with the same threshold e3b-regularity.

\medskip
\textbf{Radial sets over $\scri^+$.} Lemma~\ref{LemmaRbComm} gives such a collection for the radial point estimates near $\pa\cR^+_{\scri^+,{\rm in},+}$ and $\pa\cR^+_{\scri^+,{\rm out}}$: the operator $P^{(k)}$ featuring in the commuted equation~\eqref{EqRbCommEq} is equal to $P$ acting component-wise (on each of the $|\cA_k|$ copies of $\cE$) plus terms (given by $P_{\alpha\alpha'}^{\sharp,(k)}$ and $\tilde P^{(k)}_{\alpha\alpha'}$ in~\eqref{EqRbCommLot}) which are of lower order both in the sense of eb-differential order and decay; therefore, the threshold quantity $\ubar p_1$ (see Definition~\ref{DefSDWubarp1} and~\eqref{EqR3RScriO}) is the same for $P^{(k)}$ as for $P$.

Compared to the equation~\eqref{EqRbRHCommEq}, the commuted equation~\eqref{EqRbCommEq} features no microlocalizers, but it does feature the term $\tilde f^{(k),0}_\alpha$ which is a sum of terms of the schematic form $x_\sscri^2\rho_+^2(D_\bop^{\leq 1}p_j)D_\ebop^{\leq 1} D_\bop^{k-1}u$, $j=0,1$, which involves $k-1$ many b-derivatives of $u$. Let $G$ be as in Proposition~\ref{PropR3RScriI} or \ref{PropR3RScriO} except now we require $G\in\CI_\bop\Psi_\ebop^0$, and let $\tilde G\in\CI_\bop\Psi_\ebop^0$ be elliptic on $\WF'_\ebop(G)$, with the Schwartz kernels of $G,\tilde G$ supported in the interior of $(\supp\tilde\chi)^2$. When the eb-regularity of the coefficients $p_j$ is sufficiently large (depending only on $s$), we then have eb-microlocal estimates
\begin{align*}
  &\|G x_\sscri^2\rho_+^2(D_\bop^{\leq 1}p_j)D_\ebop^{\leq 1} D_\bop^{k-1}u\|_{H_\ebop^{s-1,(2\alpha_\sscri+2,\alpha_++2)}} \\
  &\qquad \leq \|G x_\sscri^2\rho_+^2(D_\bop^{\leq 1}p_j)D_\ebop^{\leq 1}\tilde G D_\bop^{k-1}u\|_{H_\ebop^{s-1,(2\alpha_\sscri+2,\alpha_++2)}} \\
  &\qquad\qquad + \|G x_\sscri^2\rho_+^2(D_\bop^{\leq 1}p_j)D_\ebop^{\leq 1}(I-\tilde G)D_\bop^{k-1}u\|_{H_\ebop^{s-1,(2\alpha_\sscri+2,\alpha_++2)}} \\
  &\qquad \leq C\Bigl( \|\tilde G D_\bop^{k-1}u\|_{\Heb^{s-1,(2\alpha_\sscri,\alpha_+)}} + \|\tilde\chi D_\bop^{k-1}u\|_{\Heb^{s_0,(2\alpha_\sscri,\alpha_+)}} \Bigr).
\end{align*}
(Since this only requires $p_j$ to have one b-derivative, there are no tameness considerations here.) This (upon replacing $D_\bop^{k-1}$ by $D_{[\bop]}^{k-1}$---recalling that commutator b-vector fields $D_{[\bop]}$ span the space of all b-vector fields $D_\bop$ over $\CI$---and commuting $D_{[\bop]}^{k-1}$ through $\tilde G$ and $\tilde\chi$) can be estimated inductively since it involves only $(k-1)$-many b-derivatives on $u$. \textit{This implies the analogues of Propositions~\usref{PropR3RScriI}--\usref{PropR3RScriO} of the type of Proposition~\usref{PropRbRH}.}

\medskip
\textbf{Radial sets over $\pa\cK^+$.} It remains to describe the commutation of suitable b-vector fields through $P$ near $\pa\cR_{\pa\cK^+,{\rm in/out}}$. This was addressed for single commutators in \cite[Lemma~5.31]{HintzNonstat}. We give a simpler construction here\footnote{We take full advantage of the microlocal ellipticity of time-dilations at these radial sets (see \eqref{EqRbROEll}) but produce a pseudodifferential commuted equation. The reference on the other hand produces, with some extra effort, a commuted equation which is differential.} and describe also iterated commutators. We introduce the operators
\begin{equation}
\label{EqRbROX12}
  X_1 := \nabla^\cE_{t_*\pa_{t_*}},\quad X_2 := I.
\end{equation}
Recall that, therefore, in a local trivialization of $\cE$, the operator $X_1$ acts component-wise as the time dilation vector field $t_*\pa_{t_*}\equiv-\rho_\cK\pa_{\rho_\cK}\bmod\rho_\cK\Vb$ (cf.\ Lemma~\ref{LemmaCT3bDil} and the discussion after~\eqref{EqCT3bDilEquiv}) plus a zeroth order term of class $\rho_\cK\CI$; here we fix
\[
  \rho_\cK := \frac{r}{t_*},\quad \rho_+:=\frac{1}{r}
\]
as the local defining functions of $\cK^+$ and $\iota^+$ near $\pa\cK^+$. Recalling the linear 3b-momentum $\sigma_\tbop$ (dual to $r\pa_{t_*}$) from~\eqref{EqTs3bCCoord} and the fact that the first-order symbol $\sigma_\tbop$ is elliptic at both radial sets $\cR^+_{\pa\cK^+,{\rm in}}$ and $\cR^+_{\pa\cK^+,{\rm out}}$ (see Definition~\ref{DefTs3bRad}), note then that
\begin{equation}
\label{EqRbROEll}
  \upsigma_\tbop^1(\rho_\cK X_1) = \upsigma_\tbop^1(r\pa_{t_*}) = i\sigma_\tbop
\end{equation}
is elliptic at the radial sets. This gives the following variant of Lemma~\ref{LemmaRbRMix}, where we now write $\vec X=(X_1,X_2)$.

\begin{lemma}[Mixed norms using $X_1$]
\label{LemmaRbROMix}
  Let $\cU\subset\Stb^*M$ be a sufficiently small neighborhood of $\pa\cR^+_{\pa\cK^+,{\rm in/out}}$. Let $K$ be a compact set containing the base projection of $\cU$ (thus, $K$ is a neighborhood of $\pa\cK^+$). Suppose the operator wave front sets of $B,\tilde B\in\CI_\bop\Psi_\tbop^0$ are contained in $\cU$, with $\tilde B$ elliptic on $\WF'_\tbop(B)$, and suppose that the Schwartz kernels of $B,\tilde B$ are supported in $K\times K$. Let $\tilde\chi\in\CI(M)$ be supported near and equal to $1$ near $K$. Then for all orders $s,s_0$ and $k\in\N_0$, there exists a constant $C$ such that the estimates~\eqref{EqRbRMix1} and \eqref{EqRbRMix2} hold, now with $\vec X=(X_1,X_2)$ and~\eqref{EqRbROX12}.
\end{lemma}

Roughly speaking, this means that microlocalized $H_{\tbop;\bop}^{(s;k)}$-norms of distributions $u$ can be computed as the sum of the $H_\tbop^s$-norms of $X_1^j u$ when $j$ ranges from $0$ to $k$.

\begin{proof}[Proof of Lemma~\usref{LemmaRbROMix}]
  The key point is that the space of b-vector fields on $M$ is spanned (locally near $\pa\cK^+$) by $V_1:=t_*\pa_{t_*}$ and the 3b-vector fields $r\pa_r$, $\pa_\omega$; writing $I=Q\rho_\cK^{k-j}X_1^{k-j}+R$ (with $j=0,\ldots,k$) where $Q\in\CI_\bop\Psi_\tbop^{-(k-j)}$, $R\in\CI_\bop\Psi_\tbop^0$, $\WF'_\tbop(R)\cap\bar\cU=\emptyset$, we can then write, for $A\in\Diff_\tbop^{k-j}$,
  \[
    A X_1^j B u = A Q\rho_\cK^{k-j} X_1^k B u + A R X_1^j B u.
  \]
  The $\Htb^s$-norm of the first term is bounded by that of $X_1^k B u$ since $A Q \rho_\cK^{k-j}\subset\rho_\cK^{k-j}\CI_\bop\Psi_\tbop\subset\CI_\bop\Psi_\tbop$ is bounded on every 3b-Sobolev space; and the norm of the final term is bounded by $C\|\tilde\chi u\|_{\Htb^{s_0}}$ since $\WF'_\tbop(R)\cap\WF'_\tbop(B)=\emptyset$. This gives~\eqref{EqRbRMix1}; and~\eqref{EqRbRMix2} follows from this as in the proof of Lemma~\ref{LemmaRbRMix}.
\end{proof}

Iterated commutators with $X_1$ are described in the following variant of Lemma~\ref{LemmaRbRHComm}:

\begin{lemma}[Commutators with $t_*\pa_{t_*}$]
\label{LemmaRbROComm}
  Set $\cA_k:=\{\alpha\in\N_0^2\colon|\alpha|=k\}$. For all sufficiently small neighborhoods $\cU\subset\Stb^*M$ of $\pa\cR^+_{\pa\cK^+,{\rm in/out}}$, the following holds. If $P u=f$ where $P$ is a weakly admissible wave-type operator of class $((d_0;k),(2\ell_\sscri,\ell_+,\ell_\cK))$, then $u^{(k)}:=(\vec X^\alpha u)_{\alpha\in\cA_k}$ satisfies the equation~\eqref{EqRbRHCommEq} where now
  \begin{enumerate}
  \item the components of $\tilde f^{(k)}_{\rm micro}$ are sums of terms of the form
    \[
      (D_\bop^{\leq 1}p)R\vec X^\beta u,\quad |\beta|\leq k-1,\ R\in\CI_\bop\Psi_\tbop^2,\ \WF'_\tbop(R)\cap\bar\cU=\emptyset,
    \]
    where $p=p^\sharp+\tilde p$, $p^\sharp\in\rho_+^2\CI$, $\tilde p\in\rho_+^{2+\ell_+}\rho_\cK^{\ell_\cK}\cC_{\tbop;\bop}^{(d_0;k)}$, denotes a coefficient of $P$;
  \item the $(\alpha,\alpha')$ matrix element of $P^{(k)}$ is of the form
    \[
      P^{(k)}_{\alpha\alpha'} = \delta_{\alpha\alpha'}P + \tilde P^{(k)}_{\alpha\alpha'},
    \]
    where $\tilde P^{(k)}_{\alpha\alpha'}$ is a finite sum of terms of the form $(D_\bop^{\leq 1}\tilde p)\CI_\bop\Psi_\tbop^1$ (with $\tilde p$ a decaying coefficient of $P-P_0$) and terms of class $\rho_\cK\CI_\bop\Psi_\tbop^1$ depending only on $P_0$.
  \end{enumerate}
\end{lemma}

Crucially, the term $\tilde P^{(k)}_{\alpha\alpha'}$ is of lower order both in the 3b-differential order sense and in the sense of decay at $\cK^+$; therefore, when applying Proposition~\ref{PropR3RKI}/\ref{PropR3RKO} to the equation~\eqref{EqRbRHCommEq}, the threshold condition on $s,\alpha_+,\alpha_\cK$ in~\eqref{EqR3RKIThr}/\eqref{EqR3RKOThr} is the same as for $P$ itself.

\begin{proof}[Proof of Lemma~\usref{LemmaRbROComm}]
  For $\alpha=(j,q)\in\cA_k$, we have
  \[
    P\vec X^\alpha u = \vec X^\alpha f + j[X_1,P]X_1^{j-1}u + \sum_{\substack{\beta+\gamma=\alpha \\ |\beta|\geq 2}} c_{\beta\gamma}({\rm ad}_{\vec X}^\beta P)\vec X^\gamma u.
  \]
  To prove the lemma, we only need to describe the term $j[X_1,P]X_1^{j-1}u$ for $j\geq 1$. This is very similar to the second part of the proof of Lemma~\ref{LemmaRbRHComm}. We use a microlocal parametrix for $\rho_\cK X_1$ near $\pa\cR^+_{\pa\cK^+,{\rm in/out}}$, so
  \[
    I = Q\rho_\cK X_1 + R,\quad Q\in\CI_\bop\Psi_\tbop^{-1},\ R\in\CI_\bop\Psi_\tbop^0,\ \WF'_\tbop(R)\cap\bar\cU=\emptyset.
  \]
  For a term $p A$ of $P$, with $p\in\rho_+^2\CI+\rho_+^{\ell_+}\rho_\cK^{\ell_\cK}\cC_{\tbop;\bop}^{(d_0;k)}$ and $A\in\Difftb^2$, we then write
  \begin{equation}
  \label{EqRbROCommPf}
    [X_1,p A] = [X_1,p]A + p[X_1,A]
  \end{equation}
  and further
  \[
    [X_1,p]A(X_1^{j-1}u) = \underbrace{[X_1,p]A Q\rho_\cK(X_1^j u)}_{=:\,{\rm I}} + \underbrace{[X_1,p]R(X_1^{j-1}u)}_{=:\,{\rm II}}.
  \]
  We absorb ${\rm II}$ into $\tilde f^{(k)}_{\rm micro}$, and we move the term ${\rm I}$ to the left-hand side of the equation for $P(\vec X^\alpha u)$ as a (diagonal) term of class $\rho_\cK(D_\bop^{\leq 1}p)\Psi_\tbop^1$. For the second term on the right in~\eqref{EqRbROCommPf}, we use Lemma~\ref{LemmaCT3bDil} to write $[X_1,A]=\rho_\cK A^\flat X_1+\rho_\cK A^\sharp$ with $A^\flat\in\rho_+^2\Difftb^1$ and $A^\sharp\in\rho_+^2\Difftb^2$. Thus,
  \[
    p[X_1,A]X_1^{j-1}u = \rho_\cK p A^\flat X_1^j u + \rho_\cK p A^\sharp Q\rho_\cK X_1^j u + \rho_\cK p A^\sharp R X_1^{j-1} u;
  \]
  we move the first two terms to the left-hand side of the equation for $P(\vec X^\alpha)$ where they thus yield further contributions to $\tilde P^{(k)}_{\alpha\alpha}$. The third term on the right is 3b-microlocalized away from $\bar\cU$ and thus put into $\tilde f^{(k)}_{\rm micro}$.
\end{proof}

As argued before, \textit{this implies the analogues of Propositions~\usref{PropR3RKI}--\usref{PropR3RKO} of the type of Proposition~\usref{PropRbRH}.}

\subsubsection{Estimate at normally hyperbolic trapping}
\label{SssRbTr}

Our final task in this section is to prove a version of Proposition~\ref{PropR3Tr} incorporating additional $k\in\N_0$ degrees of b-regularity. Recall from~\eqref{EqTs3bSumEll} that the $\ft$-momentum $\sigma$ is elliptic on the trapped set. Therefore, we shall use the vector field $\ft\pa_\ft$ (in Boyer--Lindquist coordinates as introduced in~\eqref{EqTsBL}) to test for b-regularity; in the presence of the bundle $\cE$, we use a smooth affine connection $\nabla^\cE$ of $\cE\to M$ and use
\[
  X_1 := \nabla^\cE_{\ft\pa_\ft},\quad X_2 := I
\]
much as in~\eqref{EqRbROX12}; the operator $\rho_\cK X_1$, where we use the local defining function
\[
  \rho_\cK := \frac{1}{\ft}
\]
of $\cK^+$ near the base projection of the trapped set, is then elliptic at $\pa\Gamma$. Working near a compact subset of $(\cK^+)^\circ$, we now use the terminology ``$\cuop$'' instead of ``$\tbop$.'' Lemma~\ref{LemmaRbROMix} is then valid, with purely notational changes (namely, replacing ``$\tbop$'' by ``$\cuop$''), when $\cU\subset{}^\cuop S^*D$, with $D$ as in~\eqref{EqDyTrD}, is a sufficiently small neighborhood of $\pa\Gamma$, and $K$ is a compact set containing the base projection of $\cU$.

Similarly, Lemma~\ref{LemmaRbROComm} is valid as well with the same notational change, and the weights in $\rho_+$ can be dropped now since we are proving estimates near trapping, which is disjoint from $\iota^+$. However, this does not suffice for an inductive proof of tame estimates: since the trapping estimate loses an extra (cusp-)derivative compared to radial point estimates, one needs to control not the $H_\cuop^{s-1}$- but the $H_\cuop^s$-norm of expressions such as $\tilde G(D_\bop^j D_\cuop^{\leq 2}p)(D_\bop^{\leq k-j}u)$ where $p\in\CI+\cC_{\cuop;\bop}^{(d_0;k)}$ is a coefficient of $P$, $\tilde G\in\CI_\bop\Psi_\cuop^0$ is elliptic on $\pa\Gamma$ (cf.\ \eqref{EqRbRHTame2}), and $j\geq 2$. But this requires the $H_{\cuop;\bop}^{(s+2;k-j)}$-norm of $u$ which analogously to~\eqref{EqRbRHTame3} one would like to estimate using the $H_{\cuop;\bop}^{(s+3-j;k-1)}$-norm of $u$. This is acceptable for $j\geq 3$, but not for $j=2$ due to the increase in the cusp regularity order. To remedy this, we must therefore refine Lemma~\ref{LemmaRbROComm} to also collect double commutators of $X_1$ with $P$ into the operator $P^{(k)}$. (Similar considerations already appear in the proof of \cite[Proposition~6.7]{HintzGlueLocII}.) Concretely:

\begin{lemma}[Commutators with $\ft\pa_\ft$ near trapping]
\label{LemmaRbTrComm}
  Set $\cA_k:=\{\alpha\in\N_0^2\colon|\alpha|=k\}$. For all sufficiently small neighborhoods $\cU\subset{}^\cuop S^*D$ of $\pa\Gamma$, the following holds. If $P u=f$ where $P$ is a weakly admissible wave-type operator of class $((d_0;k),(2\ell_\sscri,\ell_+,\ell_\cK))$, then $u^{(k)}:=(\vec X^\alpha u)_{\alpha\in\cA_k}$ satisfies
  \begin{equation}
  \label{EqRbTrCommEq}
  \begin{split}
    &P^{(k)}u^{(k)} = f^{(k)} + \tilde f^{(k)} + \tilde f^{(k)}_{\rm micro}, \\
    &\qquad\qquad f^{(k)}=(\vec X^\alpha f)_{\alpha\in\cA_k}, \ \ \tilde f^{(k)} := \Biggl(\;\sum_{\substack{\beta+\gamma=\alpha \\ |\beta|\geq 3}} c_{\beta\gamma}({\rm ad}_{\vec X}^\beta P)\vec X^\gamma u\Biggr)_{\alpha\in\cA_k},
  \end{split}
  \end{equation}
  where the $c_{\beta\gamma}\in\N_0$ are combinatorial constants; furthermore,
  \begin{enumerate}
  \item the components of $\tilde f^{(k)}_{\rm micro}$ are sums of terms of the form
    \[
      (D_\bop^{\leq 2}p)R\vec X^\beta u,\quad |\beta|\leq k-1,\ R\in\CI_\bop\Psi_\cuop^2,\ \WF'_\cuop(R)\cap\bar\cU=\emptyset,
    \]
    where $p=p^\sharp+\tilde p$, $p^\sharp\in\CI$, $\tilde p\in\rho_\cK^{\ell_\cK}\cC_{\cuop;\bop}^{(d_0;k)}$, denotes a coefficient of $P$ when expressed near $(\cK^+)^\circ$ as a sum of terms of the form $p A$, $A\in\Diff_\cuop^2$;
  \item the $(\alpha,\alpha')$ matrix element of $P^{(k)}$ is of the form
    \begin{equation}
    \label{EqRbTrCommPiece}
      P^{(k)}_{\alpha\alpha'} = \delta_{\alpha\alpha'}P + \tilde P^{(k)}_{\alpha\alpha'},
    \end{equation}
    where $\tilde P^{(k)}_{\alpha\alpha'}$ is a finite sum of terms of the form $(D_\bop^{\leq 2}\tilde p)\CI_\bop\Psi_\cuop^1$, with $\tilde p$ as before, and terms of class $\rho_\cK\CI_\bop\Psi_\cuop^1$ depending only on $P_0$.
  \end{enumerate}
\end{lemma}

The key point, similarly to before, is that $\tilde P^{(k)}_{\alpha\alpha'}$ in~\eqref{EqRbTrCommPiece} is of lower order both in the sense of the (cusp-)differential order and in the sense of decay at $\cK^+$, and thus the subprincipal sign condition for $P$ as stated in~\eqref{EqR3TrSubpr} is valid for $P^{(k)}$ as well.

\begin{proof}[Proof of Lemma~\usref{LemmaRbTrComm}]
  For $\alpha=(j,q)\in\cA_k$, we now use Lemma~\ref{LemmaMTameComm} to expand the commutator $[P,\vec X^\alpha]$ to one more order, to wit,
  \begin{equation}
  \label{EqRbTrCommExp}
    P\vec X^\alpha u = \vec X^\alpha f + c_j[X_1,P]X_1^{j-1}u + c'_j[X_1,[X_1,P]]X_1^{j-2} + \sum_{\substack{\beta+\gamma=\alpha \\ |\beta|\geq 3}} c_{\beta\gamma}({\rm ad}_{\vec X}^\beta P)\vec X^\gamma u;
  \end{equation}
  here $c_j=0$ for $j=0$ and $c'_j=0$ for $j=0,1$. Using Lemma~\ref{LemmaCT3bDil} to write
  \begin{equation}
  \label{EqRbTrCommX1A}
    [X_1,A]=\rho_\cK A^\flat X_1+\rho_\cK A^\sharp,\quad A^\sharp\in\Diff_\cuop^1,\ A^\sharp\in\Diff_\cuop^2.
  \end{equation}
  Consider a term $p A$ of $P$, with $p$ as in the statement of the lemma and $A\in\Diff_\cuop^2$, and write
  \begin{equation}
  \label{EqRbTrCommX1pA}
    [X_1,p A] = [X_1,p]A + p[X_1,A] = [X_1,p]A + \rho_\cK p A^\flat X_1 + \rho_\cK p A^\sharp.
  \end{equation}
  When $j\geq 1$, we use a microlocal elliptic (cusp-)parametrix for $\rho_\cK X_1$,
  \begin{equation}
  \label{EqRbTrCommPx1}
    I=Q\rho_\cK X_1+R,\quad Q\in\CI_\bop\Psi_\cuop^{-1},\ R\in\CI_\bop\Psi_\cuop^0,\ \WF'_\cuop(R)\cap\bar\cU=\emptyset,
  \end{equation}
  to write, as before,
  \[
    [X_1,p A]X_1^{j-1}u = \bigl([X_1,p]A Q\rho_\cK + \rho_\cK p A^\flat\bigr)(X_1^j u) + [X_1,p]A R X_1^{j-1}u + \rho_\cK p A^\sharp(X_1^{j-1}u).
  \]
  Since $[X_1,p]\in\rho_\cK\CI+\rho_\cK^{\ell_\cK}\cC_{\cuop;\bop}^{(d_0;k-1)}$, and also $\rho_\cK p$ lies (a fortiori) in this space, we can move the first term on the right to the left-hand side of the equation for $P\vec X^\alpha$; it thus contributes to the term $\tilde P^{(k)}_{\alpha\alpha}$ in~\eqref{EqRbTrCommPiece}. We absorb the second term into $\tilde f^{(k)}_{\rm micro}$. In the third term, we insert~\eqref{EqRbTrCommPx1} yet again to rewrite it as
  \[
    \rho_\cK p A^\sharp Q\rho_\cK(X_1^j u) + \rho_\cK p A^\sharp R(X_1^{j-1}u);
  \]
  the first, resp.\ second term yields a further contribution to $\tilde P^{(k)}_{\alpha\alpha}$, resp.\ $\tilde f^{(k)}_\alpha$.

  For $j\geq 2$, we next analyze the double commutator in~\eqref{EqRbTrCommExp}, with $p A$ in place of $P$. We use~\eqref{EqRbTrCommX1pA} and expand $[X_1,A]$ using~\eqref{EqRbTrCommX1A} to obtain
  \begin{align*}
    [X_1,[X_1,p A]]X_1^{j-2}u &= [X_1,[X_1,p]]A X_1^{j-2}u + [X_1,p]\rho_\cK A^\flat X_1^{j-1}u + [X_1,p]\rho_\cK A^\sharp X_1^{j-2}u \\
      &\qquad + [X_1,\rho_\cK p A^\flat]X_1^{j-1}u + [X_1,\rho_\cK p A^\sharp]X_1^{j-2}u.
  \end{align*}
  In the first term, we insert
  \[
    I = Q'\rho_\cK^2 X_1^2 + R',\quad Q'\in\CI_\bop\Psi_\cuop^{-2},\ R\in\CI_\bop\Psi_\cuop^0,\ \WF'_\cuop(R)\cap\bar\cU=\emptyset,
  \]
  and obtain $(D_\bop^2 p)A Q'\rho_\cK^2 X_1^j u$ (which, due to the $\rho_\cK$-factors, we can put into $\tilde P^{(k)}_{\alpha\alpha}$) plus the term $(D_\bop^2 p)A R' X_1^{j-2}u$ (which we put into $\tilde f^{(k)}_{\rm micro}$); similarly for the third and fifth terms. In the second term, we insert~\eqref{EqRbTrCommPx1} to get $(D_\bop p)\rho_\cK A^\flat Q\rho_\cK X_1^j u$ (put into $\tilde P^{(k)}_{\alpha\alpha}$) plus $(D_\bop p)\rho_\cK A^\flat R X_1^{j-1}u$ (put into $\tilde f^{(k)}_{\rm micro}$); similarly for the fourth term.
\end{proof}

\begin{prop}[Estimate at the trapped set with b-regularity]
\label{PropRbTr}
  There exist operators $B,E,G\in\CI_\bop\Psi_\cuop^0$ satisfying the properties~\eqref{ItR3Tr1}--\eqref{ItR3Tr4} in Proposition~\usref{PropR3Tr} such that for all $s>s_0$ and $\alpha_\cK\in\R$, and for all $\tilde\chi\in\CI(M)$ supported near $K$ and equal to $1$ near $K$, we have
  \begin{equation}
  \label{EqRbTr}
    \|B u\|_{H_{\cuop;\bop}^{(s;k),\alpha_\cK}} \leq C_k\Bigl( \|G P u\|_{H_{\cuop;\bop}^{(s;k),\alpha_\cK}} + \|E u\|_{H_{\cuop;\bop}^{(s+1;k),\alpha_\cK}} + \|\tilde\chi u\|_{H_{\cuop;\bop}^{(s_0;k),\alpha_\cK}} \Bigr),
  \end{equation}
  which holds uniformly for $P=P_0$ and all of its sufficiently small perturbations as measured in the norm $\|\cdot\|_{(d_0;k),(0,0,\ell_\cK)}$ when $d_0\in\N$ is sufficiently large (depending on $s_0,s$). Furthermore, there exists $d\in\N$ such that for $P=P_0$ and all of its sufficiently small perturbations as measured in $\|\cdot\|_{(d_0;d),(0,0,\ell_\cK)}$, one has the b-tame estimate
  \begin{equation}
  \label{EqRbTrTame}
  \begin{split}
    &\|B u\|_{H_{\cuop;\bop}^{(s;k),\alpha_\cK}} \\
    &\qquad \leq C_k\biggl( \|G P u\|_{H_{\cuop;\bop}^{(s;k),\alpha_\cK}} + \|E u\|_{H_{\cuop;\bop}^{(s+1;k),\alpha_\cK}} + \|\tilde\chi u\|_{H_{\cuop;\bop}^{(s_0;k),\alpha_\cK}} \\
    &\qquad \quad \hspace{3em} + \|P-P_0\|_{(d_0;k),(0,0,\ell_\cK),K'} \Bigl( \|G P u\|_{H_{\cuop;\bop}^{(s;d),\alpha_\cK}} + \|E u\|_{H_{\cuop;\bop}^{(s+1;d),\alpha_\cK}} + \|\tilde\chi u\|_{H_{\cuop;\bop}^{(s_0;d),\alpha_\cK}} \Bigr) \biggr)
  \end{split}
  \end{equation}
  where $K'=\supp\tilde\chi$, provided $P$ is of class $((d_0;k),(0,0,\ell_\cK))$.
\end{prop}
\begin{proof}
  We apply Proposition~\ref{PropR3Tr} to the commuted equation~\eqref{EqRbTrCommEq}. The $H_\cuop^{s,\alpha_\cK}$-norms of $G f^{(k)}$ and $G\tilde f^{(k)}_{\rm micro}$ can be estimated exactly as in the proof of Proposition~\ref{PropRbRH} but with $s$ in place of $s-1$; likewise for the $H_\cuop^{s,\alpha_\cK}$-norm of $B\vec X^\alpha u$ and the $H_\cuop^{s+1,\alpha_\cK}$-norm of $E\vec X^\alpha$ (the latter requiring replacing $s$ by $s+1$ in the argument for $E\vec X^\alpha u$ in the proof of Proposition~\ref{PropRbRH}). The only difference to previous arguments is thus the estimate of $\|G\tilde f^{(k)}\|_{H_\cuop^{s,\alpha_\cK}}$, which we schematically (and dropping the weight $\alpha_\cK$ from the notation) write as
  \[
    \|G (D_\bop^j p)D_\cuop^2 D_\bop^{\leq k-j}u\|_{H_\cuop^s},\quad 3\leq j\leq k.
  \]
  (Unlike in~\eqref{EqRbRHTame1}, $j$ now starts at $3$.) Arguing as there, we commute $D_\cuop^2$ through $D_\bop^j p$ and through $G$ so that it remains to estimate, for $\tilde G\in\CI_\bop\Psi_\cuop^0$ which is elliptic on $\WF'_\cuop(G)$, the norm
  \begin{equation}
  \label{EqRbTrTameNorm}
    \|\tilde G(D_\bop^j D_\cuop^{\leq 2}p)(D_\bop^{\leq k-j}u)\|_{H_\cuop^{s+2}}.
  \end{equation}

  Let $G^\sharp\in\CI_\bop\Psi_\cuop^0$ be elliptic on $\WF'_\cuop(\tilde G)$. For the purpose of proving~\eqref{EqRbTr}, we note that $D_\bop^j p$ is, a fortiori, of class $\CI+\cC_\cuop^{d_0}$, and thus this norm is bounded by $\|G^\sharp u\|_{H_{\cuop;\bop}^{(s+2;k-3)}}\leq C\|G^\sharp u\|_{H_{\cuop;\bop}^{(s;k-1)}}$ plus $\|\tilde\chi u\|_{H_{\cuop;\bop}^{(s_0;k-3)}}$; this can be controlled inductively.

  For the tame estimate, we write $p=p^\sharp+\tilde p$ in~\eqref{EqRbTrTameNorm} as in Lemma~\ref{LemmaRbTrComm}; only the term $\tilde p$ matters for tameness considerations. We use Proposition~\ref{PropMTameMicr} with $D_\bop^3 D_\cuop^{\leq 2}p$ and $k-3$ in place of $\ell$ and $k$, respectively. Write $\tilde d_0$ for the value of $d_0$ in Proposition~\ref{PropMTameMicr}, then~\eqref{EqRbTrTameNorm} is bounded by
  \begin{align*}
    &C_k\biggl( \|\tilde p\|_{\cC_{\cuop;\bop}^{(\tilde d_0+2,\tilde d+3)}} \Bigl( \|G^\sharp u\|_{H_{\cuop;\bop}^{(s+2;k-3)}} + \|\tilde\chi u\|_{H_{\cuop;\bop}^{(s_0;k-3)}} \Bigr) \\
    &\qquad + \|\tilde p\|_{\cC_{\cuop;\bop}^{(\tilde d_0+2;k)}}\Bigl( \|G^\sharp u\|_{H_{\cuop;\bop}^{(s+2;\tilde d)}} + \|\tilde\chi u\|_{H_{\cuop;\bop}^{(s_0;\tilde d)}} \Bigr) \biggr),
  \end{align*}
  which for $d_0\geq\tilde d_0+2$ and $d\geq\tilde d+3$ we can further estimate using $H_{\cuop;\bop}^{(s+2;k-3)} \subset H_{\cuop;\bop}^{(s;k-1)}$ by
  \begin{align*}
    &C_k\biggl( \|\tilde p\|_{\cC_{\cuop;\bop}^{(d_0,d)}} \Bigl( \|G^\sharp u\|_{H_{\cuop;\bop}^{(s;k-1)}} + \|\tilde\chi u\|_{H_{\cuop;\bop}^{(s_0;k-3)}} \Bigr) \\
    &\qquad + \|\tilde p\|_{\cC_{\cuop;\bop}^{(d_0;k)}}\Bigl( \|G^\sharp u\|_{H_{\cuop;\bop}^{(s;d)}} + \|\tilde\chi u\|_{H_{\cuop;\bop}^{(s_0;d)}} \Bigr) \biggr).
  \end{align*}
  The norms of $G^\sharp u$ appearing here can be controlled inductively (for the \emph{same} value of $s$), which thus gives~\eqref{EqRbTrTame}.
\end{proof}

\section{Energy estimates}
\label{SE}

We need energy estimates for two purposes: first, to control the solution of wave-type equations on finite time intervals (\S\ref{SsET}); second, to ``cap off'' microlocal regularity estimates, which, due to the presence of error terms with larger supports (as, e.g., in~\eqref{EqR3Est}) prevent one from closing any estimate. (In the second capacity, energy estimates were introduced for this purpose in~\cite{HintzVasySemilinear}.) While the energy estimate in $r<r_+$ is essentially standard, and an energy estimate for bounded $t_*$-intervals, localized near $\scri^+$, is a special case of the results in~\cite[\S{6}]{HintzVasyScrieb}, a novelty here is that we need to prove estimates which are \emph{global} (spatially) on bounded $t_*$-intervals and \emph{tame} in the b-regularity order; this requires new arguments which we base on a commutation result (Lemma~\ref{LemmaRbComm} in~\S\ref{SssRbComm}). The main results of this section are:
\begin{itemize}
\item Propositions~\ref{PropETSolvb} and \ref{PropETSolvbTame}: energy estimates for bounded $\ft_*$-intervals;
\item Propositions~\ref{PropERSolvb}--\ref{PropERSolvbTame}: energy estimates for bounded $r$-intervals in the black hole interior;
\item Theorem~\ref{ThmEReg}: full control of forward solutions of $P u=f$ at large e3b-frequencies.
\end{itemize}

\textit{Throughout this section, we shall only consider metrics which are small $\tilde\sG_\etbop^{0,(2\ell_\sscri,0,0)}$-perturbations of $g_{\bhm,a}$ in the sense of Lemma~\usref{LemmaSDGTime}; we denote by $\ft_*$ the time function defined there.}

\subsection{Finite-time energy estimates in \texorpdfstring{$\ft_*$}{t-star}}
\label{SsET}

We prove estimates for forward problems $P u=f$ on bounded slabs $\ft_*^{-1}([T_0,T_1])$; here
\[
  \inf_\Omega\ft_*+1\leq T_0<T_1<\infty,
\]
and $P$ is a weakly admissible wave-type operator (Definition~\ref{DefSDWAdm}). Since the only boundary hypersurface of $M$ intersecting
\begin{equation}
\label{EqETDomain}
  \Omega_{T_0,T_1} := \ol{\{T_0\leq\ft_*\leq T_1\}} \subset \Omega
\end{equation}
is $\scri^+$, only the edge-nature of edge-3b-Sobolev spaces and the decay rate at $\scri^+$ are relevant here; we shall thus work with spaces
\[
  H_\eop^{s,2\alpha_\sscri} = x_\sscri^{2\alpha_\sscri}H_\eop^s,
\]
with the underlying $L^2$-space defined using the metric volume density $|\dd\ubar g|$ or, equivalently, the metric density $|\dd g|$. For mixed spaces which capture also additional amounts of b-regularity, we use the notation
\begin{equation}
\label{EqETMixed}
  H_{\eop;\bop}^{(s;k),2\alpha_\sscri}.
\end{equation}
In order for our estimates to complement the microlocal estimates in~\S\ref{SR}, we must allow for $s$ to be variable. Note that $\Omega_{T_0,T_1}$ has three finite (i.e., not contained in $\pa M$) boundary hypersurfaces, all of which are spacelike: $r^{-1}(\bhm)$, $\ft_*^{-1}(T_0)$, and $\ft_*^{-1}(T_1)$. The first and third are final hypersurfaces of $\Omega_{T_0,T_1}$, in that future timelike vectors at them point out of $\Omega_{T_0,T_1}$, whereas the second one is an initial hypersurface. We shall correspondingly work with spaces
\[
  H_\eop^{s,2\alpha_\sscri}(\Omega_{T_0,T_1})^{\bullet,-},\quad
\]
of distributions which are of supported, resp.\ extendible character at the initial, resp.\ final hypersurfaces, and with their dual spaces $H_\eop^{-s,-2\alpha_\sscri}(\Omega_{T_0,T_1})^{-,\bullet}$ (cf.\ \S\ref{SssMUSupp}); we define $H_{\eop;\bop}^{(s;k),2\alpha_\sscri}(\Omega_{T_0,T_1})^{\bullet,-}$ similarly.

We use the (global on $\Omega_{T_0,T_1})$ coordinates
\[
  t_*,\ x_\sscri = r^{-\frac12},\ \omega\in\Sph^2,
\]
so the basic edge vector fields are $\pa_{t_*}$, $x_\sscri\pa_{x_\sscri}$, and $x_\sscri\pa_\omega$;\footnote{Since $|\dd\ubar g|=r^3|\dd t_*\,\frac{\dd r}{r}\,\dd\slg|$ is a smooth positive multiple of $x_\sscri^{-6}|\dd t_*\,\frac{\dd x_\sscri}{x_\sscri}\,\dd\slg|$, we thus have, for $s\in\N_0$, \[\|u\|_{\He^{s,2\alpha_\sscri}}^2=\sum_{j+k+|\beta|\leq s}\int_{\Omega_{T_0,T_1}^\circ} x_\sscri^{-2\alpha_\sscri}|\pa_{t_*}^j(x_\sscri\pa_{x_\sscri})^k(x_\sscri\pa_\omega)^\beta u|^2\,x_\sscri^{-6}\,\dd t_*\,\frac{\dd x_\sscri}{x_\sscri}\,\dd\slg.\]} and (cf.\ \eqref{EqSSAdmOp} and \eqref{EqSDWAdmOp})\footnote{Here we write $x_\sscri\pa_{x_\sscri}$ for $\nabla^\cE_{x_\sscri\pa_{x_\sscri}}$ where $\nabla^\cE$ is the pullback of a smooth connection on $\cE_X$, similarly for $\pa_{t_*}$ and $\slDelta$. The terms in square brackets in~\eqref{EqETedgeExpr} change by $x_\sscri\Diffe^2$ if one passes from one connection to another.}
\begin{equation}
\label{EqETedgeExpr}
\begin{split}
  g^{-1} &\equiv x_\sscri^2\bigl[ x_\sscri\pa_{x_\sscri}\otimes_s\pa_{t_*} + x_\sscri^2\slg^{-1} \bigr] \bmod x_\sscri^2\cC_{\eop;\bop}^{(d_0;k),2\ell_\sscri}(\Omega_{T_0,T_1};S^2\,\Te^*M), \\
  P &\equiv \frac12 x_\sscri^2\Bigl[-2\bigl(x_\sscri\pa_{x_\sscri}-2(1+p_1)\bigr)\pa_{t_*} + 2 x_\sscri^2\slDelta + p_0\Bigr] \bmod x_\sscri^2\cC_{\eop;\bop}^{(d_0;k),2\ell_\sscri}\Diffe^2.
\end{split}
\end{equation}
Here, $\ell_\sscri\in(0,\frac12]$. We note that, as a consequence of Definition~\ref{DefSDWAdm}\eqref{ItSDWAdmSplit}, there exists, for each $\eps>0$, a smooth positive definite fiber inner product on $\cE_X$---which induces a $t_*$-translation-invariant fiber inner product on $\cE$---such that
\begin{equation}
\label{EqETProd}
  \ubar p_1(h_\cE) := \inf_{p\in\scri^+} \spec\Bigl(\frac{p_1(p)+p_1(p)^*}{2}\Bigr) > \ubar p_1-\eps, \quad
  -\eps < \frac{p_1-p_1^*}{2 i} < \eps,\quad \|p_0\|<\eps,
\end{equation}
with adjoints and operator norms defined with respect to $h_\cE$; here we recall $\ubar p_1$ from Definition~\ref{DefSDWubarp1}. (We recall that this holds when $h_\cE$ is diagonal in the splitting $\cE=\bigoplus_{j=1}^J\cE_j$, with diagonal entries $1,\eta,\ldots,\eta^{J-1}$ for $0<\eta\ll 1$.)

\subsubsection{Estimates on edge-Sobolev spaces}

We begin with the following basic energy estimate.

\begin{lemma}[First order energy estimate]
\label{LemmaETEnergy}
  In the notation of Definition~\usref{DefSDWubarp1}, let $\alpha_\sscri<-\frac12+\ubar p_1$. Then there exists a constant $C$ such that for all $u\in\He^{2,2\alpha_\sscri}(\Omega_{T_0,T_1};\cE)^{\bullet,-}$, we have
  \begin{equation}
  \label{EqETEnergy}
    \|u\|_{\He^{1,2\alpha_\sscri}(\Omega_{T_0,T_1};\cE)^{\bullet,-}} \leq C\|P u\|_{\He^{0,2\alpha_\sscri+2}(\Omega_{T_0,T_1};\cE)^{\bullet,-}}.
  \end{equation}
  (The left-hand side controls $x^{-2\alpha_\sscri}u$ and its derivatives along $\pa_{t_*}$, $x_\sscri\pa_{x_\sscri}=-2 r\pa_r$, and $x_\sscri\pa_\omega$ in $L^2$ with the metric volume density.) The constant $C$ can be taken to be uniform for sufficiently small perturbations of $P$ as measured in the norm $\|\cdot\|_{(d_0;0),(2\ell_\sscri,0,0)}$ (Definition~\usref{DefSDWAdmNorm}) for $d_0=2$.
\end{lemma}
\begin{proof}
  It suffices to prove~\eqref{EqETEnergy} for $u$ supported in an arbitrarily small (but fixed) neighborhood of $\scri^+$. Indeed, computations as in the proof of Lemma~\ref{LemmaSDGTime} imply that the level sets of
  \begin{equation}
  \label{EqETtstar}
    \ft^*:=t_*+r^\alpha
  \end{equation}
  are spacelike near $\Omega_{T_0,T_1}\cap\scri^+$ for any fixed $\alpha\in(0,\ell_\sscri)$. Taking $\alpha=\frac12\ell_\sscri$, say, we have $\Omega_{T_0,T_1}\cap(\ft^*)^{-1}(T)\to\Omega_{T_0,T_1}\cap\scri^+$ as $T\to\infty$. In the region
  \begin{equation}
  \label{EqETtstarRegion}
    \Omega_{T_0,T_1}^{\leq T}:=\Omega_{T_0,T_1}\cap\{\ft^*\leq T\},
  \end{equation}
  for any fixed finite $T$ (see Figure~\ref{FigETtstar}), a standard energy estimate gives
  \begin{equation}
  \label{EqETEnergyLoc}
    \|u\|_{H^1(\Omega_{T_0,T_1}^{\leq T};\cE)} \leq C\|P u\|_{L^2(\Omega_{T_0,T_1}^{\leq T};\cE)}
  \end{equation}
  for all $u$ with supported character at $\ft_*=T_0$ and extendible character at $r=\bhm$, $\ft_*=T_1$, and $\ft^*=T$. Assuming we have proved~\eqref{EqETEnergy} for $u$ vanishing for $\ft^*\leq\frac{T}{2}$, one then obtains~\eqref{EqETEnergy} in general as follows: pick a cutoff function $\chi=\chi(\ft^*)$ which equals $1$ for $\ft^*\leq\frac{T}{2}$ and $0$ for $\ft^*\geq T$. Applying~\eqref{EqETEnergy} on $\Omega_{T_0,T_1}\cap\{\ft^*\geq\frac{T}{2}\}$ to $P((1-\chi)u)=P u-[P,\chi]u$, one obtains an estimate for $(1-\chi)u$, with right hand side involving $\|[P,\chi]u\|_{H^1}\leq C\|u\|_{H^1(\Omega_{T_0,T_1}^{\leq T})}$, which is, in turn, controlled by~\eqref{EqETEnergyLoc}.

  \begin{figure}[!ht]
  \centering
  \caption{Illustration of the level sets of $\ft^*$ from~\eqref{EqETtstar} and the domain $\Omega_{T_0,T_1}^{\leq T}$ from~\eqref{EqETtstarRegion}.}
  \label{FigETtstar}
  \end{figure}

  \pfstep{Vector field multiplier.} In the case that $P$ is the scalar wave operator $\Box_g$ and $u$ is supported in a sufficiently small neighborhood of $\scri^+$, we run an energy estimate with the vector field multiplier
  \begin{equation}
  \label{EqETVF}
    X = e^{-2\digamma\ft_*} X_0,\quad X_0:=x_\sscri^{-4\alpha_\sscri-2}(-x_\sscri\pa_{x_\sscri} + \pa_{t_*}) = r^{2\alpha_\sscri+1} ( 2 r\pa_r + \pa_{t_*} ),
  \end{equation}
  which is thus a weighted edge vector field; here, $2\digamma$ will be chosen large.\footnote{Comparing with \cite[(7.5)]{HintzVasyScrieb} (which proves a more delicate estimate near $\scri^+\cap\iota^+$), we do not insert a small constant in front of $-x_\sscri\pa_{x_\sscri}$. We do, however, insert the positivity-inducing prefactor $e^{-2\digamma\ft_*}$.} Since $-x_\sscri\pa_{x_\sscri}$ and $\pa_{t_*}$ are both future null edge tangent vectors by~\eqref{EqETedgeExpr}, $X_0$ is uniformly future timelike \emph{as an edge tangent vector}. Crucially, then, the associated $K$-current (cf.\ \cite[Theorem~9.16]{HintzMicro})
  \begin{align*}
    &{}^{(X)}\!K := {}^{(X)}\pi\cdot T = -\frac12(\cL_X g^{-1})^{\mu\nu}T(\pa_\mu,\pa_\nu) , \\
    &\qquad {}^{(X)}\pi:=\frac12\cL_X g,\quad T(X,Y):=X\otimes_s Y-\frac12 g(X,Y)g^{-1},
  \end{align*}
  is $x_\sscri^{-4\alpha_\sscri}$ times a positive definite quadratic form on $\Te^*M$ over $\Omega_{T_0,T_1}$ near $\scri^+$: the condition $\alpha_\sscri<-\frac12$ provides positivity on $\dd t_*\otimes\dd t_*$ by direct computation (cf.\ the $\pa_{t_*}^2$ term in~\eqref{EqETKCurrent}), while $-\nabla\ft_*$ being future timelike (and approximately null as $r\to\infty$) implies that the second term of
  \[
    {}^{(X)}\!K = e^{-2\digamma\ft_*}\bigl({}^{(X_0)}\!K + 2\digamma T(X_0,-\nabla\ft_*)\bigr)
  \]
  gives arbitrary amounts of positive definiteness also on the second symmetric tensor power of the span of $\frac{\dd x_\sscri}{x_\sscri}$, $x_\sscri^{-1}\,T^*\Sph^2$. Concretely, modulo $x_\sscri^{-4\alpha_\sscri+2\ell_\sscri}\cC_\eop^{d_0-1}(\Omega_{T_0,T_1};S^2\,\Te M)$,
  \begin{equation}
  \label{EqETKCurrent}
  \begin{split}
    {}^{(X_0)}\!K &\equiv x_\sscri^{-4\alpha_\sscri}\bigl( -(2\alpha_\sscri+1)\pa_{t_*}^2 - 2\pa_{t_*}\otimes_s x_\sscri\pa_{x_\sscri} - 2(\alpha_\sscri+1)x^2\slg^{-1} \bigr), \\
    T(X_0,-\nabla\ft_*) &\equiv \frac12 x_\sscri^{-4\alpha_\sscri}\bigl((x_\sscri\pa_{x_\scri})^2+x_\sscri^2\slg^{-1}\bigr),
  \end{split}
  \end{equation}
  The term $p_1$ in~\eqref{EqETedgeExpr} will give an additional contribution to the $\pa_{t_*}^2$ term which thus shifts the threshold condition on $\alpha_\sscri$.

  We will quantify positive definiteness of a (weighted) section $K$ of $S^2\,\Te M$ by comparing it to the Riemannian edge-metric
  \begin{equation}
  \label{EqETgR}
    g_R := \dd t_*^2 + \frac{\dd x_\sscri^2}{x_\sscri^2} + x_\sscri^{-2}\slg;
  \end{equation}
  we then say that $K\geq c$ if $K(\zeta,\zeta)\geq c g_R(\zeta,\zeta)$ for all $\zeta\in\Te^*M$.

  \pfstep{Top order estimate. Case I: $u=0$ outside of $\Omega_{T_0,T_1}$.} The details, for general $P$, are as follows. With $\eps>0$ to be fixed later, we work with the fiber inner product $h_\cE$ on $\cE_X$ such that~\eqref{EqETProd} holds. Denote by
  \begin{equation}
  \label{EqETX01Ops}
    \sfX_0 \in x_\sscri^{-4\alpha_\sscri-2}\Diffe^1(\Omega_{T_0,T_1};\cE)
  \end{equation}
  a differential operator whose principal part is the vector field $X_0$ from~\eqref{EqETVF}. More precisely, we choose $\sfX_0$ such that in a product trivialization of $\cE$ near a fiber of $\scri^+$ over a point $p\in\Sph^2=\pa X$ induced by a trivialization of $\cE_X$ near $p$, it satisfies $\sfX_0(e_1,\ldots,e_k)-(X_0 e_1,\ldots,X_0 e_k)=x_\sscri^{-4\alpha_\sscri-2}x_\sscri A(e_1,\ldots,e_k)$ where $k$ is the rank of $\cE_X$, $(e_1,\ldots,e_k)$ is a section of $\cE_X$, and $A\in\CI(\Omega_{T_0,T_1};\End(\cE))$; this condition is consistent upon passing to another smooth trivialization of $\cE_X$, since the resulting trivializations of $\cE$ are invariant under $t_*$-translations, and since $x_\sscri\pa_{x_\sscri}$ is $x_\sscri$ times a smooth vector field on $X_0$.\footnote{A similar construction is performed in the proof of \cite[Theorem~6.4]{HintzVasyScrieb}, but expressed using the notion of edge normal operators.} Since $X_0^*=-X_0-\dv_g X_0$, we then have (relative to any smooth choice of fiber inner product on $\cE_X$, pulled back to a stationary inner product on $\cE$) that
  \[
    \sfX_0^* = -\sfX_0 - \dv_g X_0 + \tilde f_0,\quad \tilde f_0\in x_\sscri\CI(\Omega_{T_0,T_1};\End(\cE)).
  \]
  Define the edge differential operator
  \[
    \sfX := e^{-2\digamma\ft_*}\sfX_0.
  \]
  For $u\in\dot H_\eop^{2,2\alpha_\sscri}(\Omega_{T_0,T_1};\cE)$ (i.e., $u$ is of supported character at $\ft_*=T_0,T_1$ and $r=\bhm$), we then compute the $L^2(\Omega_{T_0,T_1},|\dd g|;\cE,h_\cE)$-inner product
  \begin{equation}
  \label{EqETComm}
  \begin{split}
    2\Re\la P u,\sfX u\ra &= \la (P^*\sfX+\sfX^*P)u,u\ra \\
      &= \big\la \bigl([P,\sfX]-(\dv_g X)P - (P-P^*)\sfX + (\sfX^*+\sfX+\dv_g X)P \bigr)u, u\big\ra.
  \end{split}
  \end{equation}
  We claim that, to leading order, this is a weighted positive definite quadratic form in edge-derivatives of $u$. This follows from a variant of \cite[(9.60)--(9.63)]{HintzMicro}: denoting covectors by $\zeta$ and the dual metric function of $g$ by $G(\zeta)$, the principal symbol of $[P,\sfX]-(\dv_g X)P$ is equal to
  \[
    H_G (\sfX^\mu\zeta_\mu) - (\dv_g X)G(\zeta) = (X^{\mu;\nu}+X^{\nu;\mu}-X^\lambda{}_{;\lambda}g^{\mu\nu})\zeta_\mu\zeta_\nu = 2\,{}^{(X)}\!K(\zeta,\zeta).
  \]
  For the computation of (the principal symbol) of the first order operator $P-P^*$, we use~\eqref{EqETedgeExpr}, the $t_*$-independence of the fiber inner product on $\cE$, and the fact that the metric volume density is, to leading order at $x_\sscri=0$, equal to $2 x_\sscri^{-6}|\dd t_*\frac{\dd x_\sscri}{x_\sscri}\,\dd\slg|$; this gives
  \begin{equation}
  \label{EqETPPs}
    P-P^* \equiv 2 x_\sscri^2(p_1\pa_{t_*}+\pa_{t_*}p_1^*) + \frac12 x_\sscri^2(p_0-p_0^*) \bmod x_\sscri^{2+2\ell_\sscri}\cC_\eop^{d_0-2}\Diffe^1.
  \end{equation}
  A short calculation gives
  \begin{equation}
  \label{EqETXadj}
    \sfX^*+\sfX+\dv_g X = e^{-2\digamma\ft_*} \tilde f_0.
  \end{equation}
  The principal symbol of $P^*\sfX+\sfX^*P$ is thus equal to ${}^{(X,p_1)}\!K(\zeta,\zeta)$ where\footnote{Note that the principal symbol of $\pa_{t_*}\circ\sfX=i^2(i^{-1}\pa_{t_*})\circ(i^{-1}\sfX)$ is $-\pa_{t_*}(\zeta)X(\zeta)$ where $-1=i^2$.}
  \begin{align}
  \label{EqETK}
    {}^{(X,p_1)}\tilde K &\equiv (2\,{}^{(X)}\!K - 2(p_1+p_1^*)\pa_{t_*}\otimes_s X) + e^{-2\digamma\ft_*}\tilde f_0 g^{-1} \\
  \label{EqETKt1}
      &= 2 e^{-2\digamma\ft_*} \Bigl[ {}^{(X_0)}\!K+(p_1+p_1)^*\pa_{t_*}\otimes_s X_0 + 2\digamma T(X_0,-\nabla\ft_*)\Bigr] \\
  \label{EqETKt2}
      &\quad + e^{-2\digamma\ft_*}\tilde f_0 g^{-1}
  \end{align}
  modulo a section of $S^2\,\Te M\otimes\End(\cE)$ of the form $x_\sscri^{2+2\ell_\sscri}\cC_\eop^{d_0-2}(\Omega_{T_0,T_1};\End(\cE)\otimes\Te M)\otimes_s X$, which thus in particular has an extra factor of $x_\sscri^{2\ell_\sscri}$ relative to the main term~\eqref{EqETKt1}. (This differs from ${}^{(X)}\!K$ only through the presence of the terms involving $p_1$, which arise from $P-P^*$.) Now, recalling~\eqref{EqETKCurrent}, the $\pa_{t_*}\otimes\pa_{t_*}$ component of ${}^{(X_0)}\!K+(p_1+p_1)^*\pa_{t_*}\otimes_s X_0$ is equal to $x_\sscri^{-4\alpha_\sscri}(-2\alpha_\sscri-1+p_1+p_1^*)$ and thus positive if $\alpha_\sscri<-\frac12+\ubar p_1-\eps$ (cf.\ \eqref{EqETProd}). Upon choosing $\digamma$ large enough (and recalling~\eqref{EqETKCurrent}), the term in square brackets in~\eqref{EqETKt1} is thus bounded from below over $\Omega_{T_0,T_1}\cap\scri^+$ by $3 c_0 x_\sscri^{-4\alpha_\sscri}$ (in the sense of comparing with~\eqref{EqETgR}, as explained there), where
  \[
    c_0:=\frac14(-2\alpha_\sscri-1+2\ubar p_1) > 0,
  \]
  for all $\digamma\geq\digamma_0$ if $\digamma_0$ is sufficiently large, and thus by $2 c_0 x_\sscri^{-4\alpha_\sscri}$ for $x_\sscri\leq\delta_0$ if we fix $\delta_0>0$ sufficiently small (depending only on $\digamma_0$). Since the term~\eqref{EqETKt2} vanishes (relative to $x_\sscri^{-4\alpha_\sscri}$) at $\scri^+$ and is thus small near $x_\sscri=0$, it can also be dominated by the first term upon shrinking $\delta_0$ if necessary. In summary (and relabeling constants), there exist $\delta_0,c_0>0$ and $\digamma$ (which we will fix for the remainder of the proof) such that, relative to~\eqref{EqETgR},
  \begin{equation}
  \label{EqETKPos}
    {}^{(X,p_1)}\tilde K \geq c_0 x_\sscri^{-4\sscri}\quad\text{on}\ \Omega_{T_0,T_1}\cap\{x_\sscri\leq\delta_0\}.
  \end{equation}
  Importantly, this holds, with the same value of $c_0$, for all sufficiently small choices of $\eps>0$ and inner products $h_\cE$ on $\cE_X$ satisfying~\eqref{EqETProd}.

  In order to describe the lower order terms arising in~\eqref{EqETComm}, we use a connection $\nabla^\cE$ constructed as in~\eqref{EqSSAdmDer}. A particularly convenient choice of $\sfX$ in~\eqref{EqETX01Ops} is then
  \begin{equation}
  \label{EqETNablaX}
    \sfX = \nabla^\cE_X
  \end{equation}
  in the notation of~\eqref{EqETVF}. (That $\sfX$ satisfies the condition stated after~\eqref{EqETX01Ops} follows from the fact that $\nabla^\cE$ comes from a smooth connection on $\cE_X$, down to $x_\sscri=0$.) We denote the adjoint of $\nabla^\cE$ with respect to $|\dd g|$, $g_R$ (on $\Te^*M$), and $h_\cE$ (on $\cE$) by $(\nabla^\cE)^*\in\Diff_\eop^1+x_\sscri^{2\ell_\sscri}\cC_\eop^{d_0-1}(\Omega_{T_0,T_1};\Hom(\Te^*M\otimes\cE,\cE))$. (The zeroth order term involves one edge-derivative of the metric density $|\dd g|$.) We identify $\Te M$ and $\Te^*M$ using $g_R$, and thus also $\Te M\otimes\Te M$ with $\End(\Te^*M)$. Consider first the case that $\cE$ is trivial, $h_\cE$ is constant in the trivialization of $\cE$, $P=\Box_g$ acts component-wise as the scalar wave operator, and $\nabla^\cE_V$, for any vector field $V$, acts component-wise as $V$. Then $P^*\sfX+\sfX^*P$ and $(\nabla^\cE)^*{}^{(X)}\!K\nabla^\cE$ have the same principal symbol, so differ by a first order operator that is symmetric and has real coefficients, and hence must in fact be a zeroth order operator which, since it annihilates constant sections of $\cE$, vanishes identically; so
  \begin{equation}
  \label{EqETKEquality}
    (P^*\sfX+\sfX^*P)-(\nabla^\cE)^*{}^{(X)}\!K\nabla^\cE = 0
  \end{equation}
  in this case. If $P$ includes, in addition, a first order term $a\nabla^\cE_V$, where $a$ is a bundle endomorphism of $\cE$ and $V\in x_\sscri^2(\CI+x_\sscri^{2\ell_\sscri}\cC_\eop^{d_0})\Ve$ is a real vector field, the term $P^*\sfX+\sfX^*P$ includes, in addition, $(\nabla^\cE)^*(a^*V\otimes X+a X\otimes V)\nabla^\cE$. Replacing ${}^{(X)}\!K$ by ${}^{(X,a,V)}\!K:={}^{(X)}\!K+\frac{a+a^*}{2}(X\otimes V+V\otimes X)$ (which is still symmetric), we thus have
  \begin{align*}
    &(P^*\sfX+\sfX^*P)-(\nabla^\cE)^*\,{}^{(X,a,V)}\!K\,\nabla^\cE \\
    &\quad = (\nabla^\cE)^* \Bigl[i\Im(a)(X\otimes V-V\otimes X)\Bigr] \nabla^\cE \\
    &\quad = (\nabla^\cE_X)^* i \Im(a) \nabla^\cE_V - (\nabla^\cE_V)^* i\Im(a) \nabla^\cE_X.
  \end{align*}
  (One expand this further to see that this is a first order differential operator.) We shall apply this with $a=2 x_\sscri^2 p_1$ and $V=\pa_{t_*}$, in which case $\Im(a)=\cO(\eps)$ by~\eqref{EqETProd}. Adding to $P$ a zeroth order term $b$ changes $P^*\sfX+\sfX^*P$ via addition of $b^*\sfX+\sfX^*b$; this term is small when $b$ is small, and we will apply this with $b=\frac12 x_\sscri^2 p_0=\cO(\eps)$ (cf.\ again \eqref{EqETProd}).

  These considerations apply in any chart. In such a chart, passing to a general connection $\nabla^\cE$ (arising from a smooth connection on $\cE_X$ as before, now in a chart on $X$) and fiber inner product $h_\cE$ results in the addition of zeroth order terms to $\nabla^\cE$, resp.\ $(\nabla^\cE)^*$ which are sections of $\Hom(\cE,\Te^*M\otimes\cE)$ of class $x_\sscri\CI$, resp.\ $x_\sscri^{2\ell_\sscri}\cC_\eop^{d_0-1}$. In the term $P^*\sfX+\sfX^*P$, this produces (through the resulting zeroth order modification of~\eqref{EqETNablaX}) additional second order, thus principal, terms, which are captured by the term $e^{-2\digamma\ft_*}\tilde f_0$ in~\eqref{EqETXadj} and must be added to ${}^{(X,a,V)}\!K$. In the case of the operator $P$ under study, this yields ${}^{(X,p_1)}\!\tilde K$ from~\eqref{EqETK}.

  Altogether, then, we obtain from~\eqref{EqETComm} and these considerations the estimate
  \begin{align*}
    &\|{}^{(X,p_1)}\!\tilde K(\nabla^\cE u,\nabla^\cE u)\|_{L^2(\Omega_{T_0,T_1};\Te^*M\otimes\cE,g_R\otimes h_\cE)}^2 \\
    &\qquad \leq \eps_0\|x_\sscri^2 x_\sscri^{2\alpha_\sscri}\sfX u\|_{L^2(\Omega_{T_0,T_1};\cE,h_\cE)}^2 + \frac{1}{2\eps_0}\|x_\sscri^{-2}x_\sscri^{-2\alpha_\sscri}P u\|_{L^2(\Omega_{T_0,T_1};\cE,h_\cE)}^2 \\
    &\qquad \quad + C\Bigl(\|(\sqrt\eps+\delta_0^{2\ell_\sscri})x_\sscri^{-2\alpha_\sscri}u\|_{L^2(\Omega_{T_0,T_1};\cE,h_\cE)}^2 \\
    &\qquad \quad \hspace{4em} + \|(\sqrt\eps+\delta_0^{2\ell_\sscri})x_\sscri^{-2\alpha_\sscri}\nabla^\cE u\|_{L^2(\Omega_{T_0,T_1};\Te^*M\otimes\cE,g_R\otimes h_\cE)}^2 \Bigr)
  \end{align*}
  for all $\eps_0>0$ and $u$ with support in $\Omega_{T_0,T_1}\cap\{x_\sscri\leq\delta_0\}$; the constant $C$ is independent of $\eps,\eps_0,\delta_0$. In view of~\eqref{EqETKPos}, when $\eps_0,\eps,\delta_0$ are sufficiently small, we obtain
  \begin{equation}
  \label{EqETEst}
    \| \nabla^\cE u \|_{x_\sscri^{2\alpha_\sscri}L^2(\Omega_{T_0,T_1})}^2 \leq C\Bigl( \|P u\|_{x_\sscri^{2\alpha_\sscri+2}L^2(\Omega_{T_0,T_1})}^2 + \|(\sqrt\eps+\delta_0^{2\ell_\sscri})u\|_{L^2(\Omega_{T_0,T_1})}^2 \Bigr).
  \end{equation}

  \pfstep{Top order estimate. Case II: general $u$.} The only difference to the case already discussed is that $u$ is no longer required to vanish at the final boundary hypersurface $\ft_*=T_1$ of $\Omega_{T_0,T_1}$. (We still work with $u$ that vanish when $x_\sscri>\delta_0$ or $\ft_*<T_0$.) We enforce integration only over $\Omega_{T_0,T_1}$ by working with $\chi\sfX$ in place of $\sfX$ where $\chi=H(T_1-\ft_*)$, with $H$ denoting the Heaviside function. Starting with~\eqref{EqETKEquality} (which holds for $\chi\sfX$ in place of $\sfX$), one thus finds that if $\cE$ is trivial and $\nabla^\cE$ acts component-wise, then the only additional term one incurs arises at the principal symbol level, namely from the addition to ${}^{(X,p_1)}\tilde K$ of the term $T(X,\nabla\chi)=\delta(\ft_*-T_1)T(X,-\nabla\ft_*)$. Since $T(X,-\nabla\ft_*)$ is positive definite, this term has the same sign as the main term ${}^{(X,p_1)}\tilde K$ and can thus be dropped; thus, we obtain~\eqref{EqETEst} also in this case. For general $\nabla^\cE$, the additional $\cO(x_\sscri)\cdot X$ zeroth order terms arising from~\eqref{EqETNablaX} create additional error terms of size $C\delta_0\|u\|_{x_\sscri^{2\alpha_\sscri}L^2(\ft_*^{-1}(T_1))}\|\nabla^\cE u\|_{x_\sscri^{2\alpha_\sscri}L^2(\ft_*^{-1}(T_1))}$. (Schematically, these arise from principal terms $D^2$ of $P$ interacting with $x_\sscri\cdot\chi$, e.g., $D^2 x_\sscri\chi=D(x_\sscri\chi)D+D[D,x_\sscri\chi]$, with the second summand being a $\delta$-distribution at $\ft_*=T_1$.) Bounding this by $\frac{C\delta_0}{2}(\|u\|_{x_\sscri^{2\alpha_\sscri}L^2(\ft_*^{-1}(T_1))}^2+\|\nabla^\cE u\|_{x_\sscri^{2\alpha_\sscri}L^2(\ft_*^{-1}(T_1))}^2)$, we can absorb the second term into the left-hand side, and thus arrive at~\eqref{EqETEst} with $\delta_0^{\min(2\ell_\sscri,\frac12)}$ in place of $\delta_0^{2\ell_\sscri}$.

  \pfstep{Control of zeroth order terms.} Integrating in $t_*$ from $T_0$ to $T_1$ gives, for $u$ vanishing in $\ft_*\leq T_0$,
  \[
    \|u\|_{x_\sscri^{2\alpha_\sscri}L^2(\Omega_{T_0,T_1};\cE,h_\cE)}^2 \leq C\|\pa_{t_*}u\|_{x_\sscri^{2\alpha_\sscri}L^2(\Omega_{T_0,T_1};\cE,h_\cE)}^2,
  \]
  where $C$ is independent of the value $\eps$ entering in the choice of $h_\cE$. (We write $\pa_{t_*}$ here to emphasize that this follows from a component-wise estimate for the components of $u$ in trivializations of $\cE_X$ in a system of coordinate charts.) Plugging this into~\eqref{EqETEst} and taking $\eps,\delta_0$ to be sufficiently small yields~\eqref{EqETEnergy}.
\end{proof}

A minor extension of the above proof shows, moreover, that in the notation of~\eqref{EqETtstar}--\eqref{EqETtstarRegion}, one has a \emph{uniform} estimate
\begin{equation}
\label{EqETUnifEst}
  \|u\|_{H_\eop^{1,2\alpha_\sscri}(\Omega_{T_0,T_1}^{\leq T};\cE)^{\bullet,-}} \leq C\|P u\|_{H_\eop^{0,2\alpha_\sscri+2}(\Omega_{T_0,T_1}^{\leq T};\cE)^{\bullet,-}}
\end{equation}
for all sufficiently large $T$, where the supported character (``$\bullet$'') refers to $\ft_*=T_0$ as before, and the extendible character (``$-$'') to the hypersurfaces $\ft_*=T_0$, $r=\bhm$, and $\ft^*=T$. Indeed, upon inserting a cutoff $H(T-\ft^*)$ into the vector field multiplier, one obtains a favorable boundary term at $\ft^*=T$ since $-\nabla\ft^*$ is future timelike there.

\begin{cor}[Solvability on $H_\eop^1$]
\label{CorETSolv1}
  Let $\alpha_\sscri<-\frac12+\ubar p_1$. Then there exists a constant $C$ (uniform for $\|\cdot\|_{(d_0;0),(2\ell_\sscri,0,0)}$-perturbations of $P$ where $d_0=2$) such that for all $f\in H_\eop^{0,2\alpha_\sscri+2}(\Omega_{T_0,T_1};\cE)^{\bullet,-}$, there exists a unique forward solution $u$ of $P u=f$, and it satisfies $u\in H_\eop^{1,2\alpha_\sscri}(\Omega_{T_0,T_1};\cE)^{\bullet,-}$ and the estimate~\eqref{EqETEnergy}.
\end{cor}
\begin{proof}
  When $P$ has smooth coefficients, then on $\Omega_{T_0,T_1}^{\leq T}$, existence and uniqueness of a forward solution $u_T$ of $P u_T=f$ follows from the (standard) finite-time solvability of wave equations (see, e.g., \cite[\S{9.4}]{HintzMicro}); and the estimate~\eqref{EqETUnifEst} holds. It remains to take the limit as $T\to\infty$, using the uniformity of $C$. For general $P$, apply these arguments to smoothed-out (on the level of coefficients, and in edge unit cells on $M$ near $\Omega_{T_0,T_1}$) versions $P_\eps$, $\eps>0$, of $P$, let $\eps\to 0$, and take a weak subsequential limit of the (uniformly bounded by~\eqref{EqETUnifEst}) forward solutions $u_\eps$ to obtain $u$.
\end{proof}

Using the propagation of edge-regularity and duality arguments, we can now prove:

\begin{prop}[Solvability and uniqueness on edge-Sobolev spaces]
\label{PropETSolv}
  Let $\sfs$ be an admissible order function for $P_0$ (Definition~\usref{DefDyO}). Then there exists $d_0\in\N$ such that, for $\inf_\Omega\ft_*+1\leq T_0<T_1<\infty$ and $\alpha_\sscri<-\frac12+\ubar p_1$ in the notation of Definition~\usref{DefSDWubarp1}, and for $P$ and all perturbations thereof which are sufficiently small as measured in the norm $\|\cdot\|_{(d_0;0),(2\ell_\sscri,0,0)}$ (Definition~\usref{DefSDWAdmNorm}), the following holds. There exists a constant $C$ such that for all $f\in H_\eop^{\sfs-1,2\alpha_\sscri+2}(\Omega_{T_0,T_1};\cE)^{\bullet,-}$, the unique forward solution $u$ of $P u=f$ satisfies $u\in H_\eop^{\sfs,2\alpha_\sscri}(\Omega_{T_0,T_1};\cE)^{\bullet,-}$ and the estimate
  \begin{equation}
  \label{EqETSolvEst}
    \|u\|_{H_\eop^{\sfs,2\alpha_\sscri}(\Omega_{T_0,T_1};\cE)^{\bullet,-}} \leq C\|f\|_{H_\eop^{\sfs-1,2\alpha_\sscri+2}(\Omega_{T_0,T_1};\cE)^{\bullet,-}}.
  \end{equation}
\end{prop}

In order to create some physical space to accommodate microlocal error terms, we will use extension/restriction arguments as in \cite[Chapter~9]{HintzMicro}. Since in practice we are always given an operator $P$ on a full domain $\Omega\cap\{\ft_*\geq T_0\}$, extending past $\ft_*=T_1$ is not necessary; however, extending past $r=\bhm$ \emph{is}. In preparation for later sections which are global in $t_*$, we thus introduce:

\begin{definition}[Spaces for extension operators]
\label{DefETExt}
  Recalling $\tilde M$ from Definition~\usref{DefCMSpacetime}, $t_*=t-r$, and the Kerr black hole mass $\bhm>0$, set
  \begin{equation}
  \label{EqETExtDomains}
    \tilde\Omega := \ol{\{t_*\geq 1\}} \subset \tilde M, \quad
    \Omega_{\rm int} := \ol{\{\bhm\leq r\leq 2\bhm\}} \cap \Omega,\quad
    \tilde\Omega_{\rm int} := \ol{\{r\leq 2\bhm\}} \cap \tilde\Omega.
  \end{equation}
  For $\alpha\in\R$ and $d_0,k\in\N_0$, denote by $\cC_{\cuop;\bop}^{(d_0;k),\alpha}(\Omega_{\rm int})=t^{-\alpha}\cC_{\cuop;\bop}^{(d_0;k)}(\Omega_{\rm int})$ the space of all functions $u=u(t,x)\in\cC^{d_0+k}(\Omega_{\rm int}\setminus\cK^+)$ such that
  \begin{equation}
  \label{EqETExtNorm}
    \|u\|_{t^{-\alpha}\cC_{\cuop;\bop}^{(d_0;k)}(\Omega_{\rm int})} := \max_{\substack{j+|\beta|\leq d_0 \\ j'+|\beta'|\leq k}}\sup_{\Omega_{\rm int}^\circ} \Bigl(t^\alpha |\pa_t^j (t\pa_t)^{j'}\pa_x^{\beta+\beta'} u(t,x)|\Bigr) < \infty.
  \end{equation}
  Define the normed space $t^{-\alpha}\cC_{\cuop;\bop}^{(d_0;k)}(\tilde\Omega_{\rm int})$ in the same fashion. We moreover define
  \[
    H_{\cuop;\bop}^{(s;k),\alpha}(\Omega_{\rm int}) = t^{-\alpha}H_{\cuop;\bop}^{(s;k)}(\Omega_{\rm int})
  \]
  for $s\in\N_0$ using the squared norm
  \[
    \|u\|_{t^{-\alpha}H_{\cuop;\bop}^{(d_0;k)}(\Omega_{\rm int})}^2 := \sum_{\substack{j+|\beta|\leq s \\ j'+|\beta'|\leq k}} \int_{\Omega_{\rm int}} t^{2\alpha} |\pa_t^j(t\pa_t)^{j'}\pa_x^{\beta+\beta'}u(t,x)|^2\,\dd t\,\dd x,
  \]
  similarly for $t^{-\alpha}H_{\cuop;\bop}^{(s;k)}(\tilde\Omega_{\rm int})$. For $k=0$, we drop the subscript ``$\bop$'' and simply write $t^{-\alpha}\cC_\cuop^{d_0}$ and $t^{-\alpha}H_\cuop^s$.
\end{definition}

The subscript ``$\cuop$'' stands for ``cusp'' in the sense of \cite{MazzeoMelroseFibred}: the norm~\eqref{EqETExtNorm} tests for $d_0$ degrees of regularity with respect to the so-called cusp vector fields $\rho_\cK^2\pa_{\rho_\cK}$, $\pa_x$ where $\rho_\cK=t^{-1}$. On $\Omega_{\rm int}$, $(\cuop;\bop)$-regularity is the same as $(\etbop;\bop)$-regularity (Definition~\ref{DefMUCe3b}). For $s\in\R$, we can thus equivalently define $H_{\cuop;\bop}^{(s;k)}(\Omega_{\rm int})$ and its weighted versions as the spaces $\bar H_{\etbop;\bop}^{(s;k)}(\Omega_{\rm int})$ of extendible distributions on (the interior of) $\Omega_{\rm int}$.

\begin{lemma}[Extension operator]
\label{LemmaETExt}
  There exists an operator
  \[
    E \colon \CIc(\Omega_{\rm int}\setminus\cK^+) \to \CIc(\tilde\Omega_{\rm int}\setminus\cK^+)
  \]
  with the following properties.
  \begin{enumerate}
  \item{\rm (Extension.)} $(E u)|_{\Omega_{\rm int}^\circ}=u$ for all $u$.
  \item{\rm (Continuity on $L^\infty$-based spaces.)} For all $\alpha\in\R$ and $d_0,k\in\N_0$, $E$ extends to a bounded linear map
    \begin{equation}
    \label{EqETExt}
      E \colon t^{-\alpha}\cC_{\cuop;\bop}^{(d_0;k)}(\Omega_{\rm int}) \to t^{-\alpha}\cC_{\cuop;\bop}^{(d_0;k)}(\tilde\Omega_{\rm int}).
    \end{equation}
  \item{\rm (Continuity on $L^2$-based spaces.)} For all $\alpha\in\R$ and $s,k\in\N_0$, $E$ extends to a bounded linear map
    \begin{equation}
    \label{EqETExt2}
      E \colon t^{-\alpha}H_{\cuop;\bop}^{(s;k)}(\Omega_{\rm int}) \to t^{-\alpha}H_{\cuop;\bop}^{(s;k)}(\tilde\Omega_{\rm int}).
    \end{equation}
  \end{enumerate}
\end{lemma}
\begin{proof}
  This is a variant of Seeley's theorem \cite{SeeleyExtension}. Fix a cutoff function $\chi\in\CIc([0,\bhm))$ which equals $1$ on $[0,\frac{\bhm}{2}]$. We will work in the coordinates $t_*$, $z=|x|-\bhm$, $\omega=\frac{x}{|x|}$, and define $E$ by the formula
  \[
    (E u)(t_*,z,\omega) = \begin{cases} u(t_*,z,\omega), & 0\leq z\leq\bhm, \\ \sum_{l=0}^\infty c_l(\chi u)\bigl(t_*,-\frac{z}{\eps_l},\omega\bigr), & -\bhm<z<0, \end{cases}
  \]
  where we set $\eps_l=3^{-l}$ and define $c_l$ via $\sin(\frac{\pi}{2}w)=\sum_{l=0}^\infty c_l w^l$. These choices ensure that $C_j:=\sum_{l=0}^\infty|c_l|\eps_l^{-j}<\infty$ and $\sum_{l=0}^\infty c_l\eps_l^{-j}=(-1)^j$ for all $j\in\N_0$, and hence $E u$ is smooth across $z=0$. Note that $E u$ vanishes for $z<-\frac{\bhm}{2}$, i.e., for $r<\frac{\bhm}{2}$, and since $t_*,\omega$ are merely parameters, $E$ preserves regularity in $\omega$ and along $\pa_{t_*}$ and $t_*\pa_{t_*}$. To prove~\eqref{EqETExt}, it thus suffices to consider $z$-derivatives of $E u$, and to note that for $z<0$, the absolute value of
  \[
    \pa_z^j(E u)(t_*,z,\omega) = \sum_{l=0}^\infty (-1)^j c_l\eps_l^{-j} (\pa_z^j(\chi u))\Bigl(t_*,-\frac{z}{\eps_l},\omega\Bigr)
  \]
  is bounded by $C_j\sup_{0\leq z\leq\bhm} |\pa_z^j(\chi u)(t_*,z,\omega)|$. The proof of~\eqref{EqETExt2} is similar; the underlying $L^2$-estimate is, for $v:=\chi u$,
  \[
    \biggl(\int_{-\bhm}^0 \Bigl|c_l v\Bigl(-\frac{z}{\eps_l}\Bigr)\Bigr|^2\,\dd z \biggr)^{\frac12} = |c_l|\eps_l^{\frac12} \biggl(\int_0^\bhm |v(z)|^2\,\dd z\biggr)^{\frac12}.\qedhere
  \]
\end{proof}

Given an admissible wave-type operator $P$, with $\cC_{\etbop;\bop}^{(d_0;k)}$-regular coefficients, as in Definition~\ref{DefSDWAdm} acting on complex-valued functions, one can then express $P$ near $r=\bhm$ uniquely as a sum of operators of the form $u A$ where $A\in\{I,\pa_{z^\mu},\pa_{z^\mu}\pa_{z^\nu}\}$, $0\leq\mu\leq\nu\leq 3$, $z=(t_*,x)$ and $u=u_0+\tilde u$ where $u_0\in\CI(\R^3\cap\{r\geq\bhm\})$ is $t_*$-independent, while $\tilde u\in t_*^{-\ell_\cK}\cC_{\cuop;\bop}^{(d_0;k)}(\Omega_{\rm int})$; extending $u_0$ to a smooth function on $\R^3$ and $\tilde u$ to an element of $t_*^{-\ell_\cK}\cC_{\cuop;\bop}^{(d_0;k)}(\tilde\Omega_{\rm int})$, one thus obtains an extension of $P$ to a wave-type operator on $\tilde\Omega\cap\{r\geq r_-\}$ which agrees with $P$ on $\Omega$; here $r_-\in(0,\bhm)$ is close enough to $\bhm$ so that the principal symbol of $P$ still comes from a Lorentzian metric for which $\dd\ft_*$ (with $\ft_*$ given by the formula in Lemma~\ref{LemmaSDGTime}\eqref{ItSDGTime3} in the extension region) is everywhere timelike, and $\dd r$ is timelike in the region $\{r_-\leq r<r_+\}$, where we recall $r_+$ from~\eqref{EqTsBLMfd}. When $P$ acts on sections of a stationary bundle $\cE=\pi_X^*\cE_X\to M$, this construction can be applied in local trivializations of $\cE_X$ and then patched together using a partition of unity.

\begin{proof}[Proof of Proposition~\usref{PropETSolv}]
  This is analogous to the proof of \cite[Theorem~6.4]{HintzVasyScrieb}. We drop the bundle $\cE$ from the notation. Write
  \begin{equation}
  \label{EqETSolvOmegas}
    \Omega^\flat:=\Omega_{T_0,T_1}\subset\Omega,\quad
    \Omega^\sharp:=\ol{\{T_0\leq\ft_*\leq T_1+1,\ r\geq r^\sharp\}}\subset\tilde M,
  \end{equation}
  where $r^\sharp\in(0,\bhm)$ is close to $\bhm$ (as in the preceding paragraph). See Figure~\ref{FigETSolvOmegas}. We use the quotient norm on $H_\eop^{s-1,2\alpha_\sscri+2}(\Omega^\flat)^{\bullet,-}$ induced by the restriction map from $\Omega^\sharp$. Lemma~\ref{LemmaETExt} (and property~\eqref{ItSDGTime3} of $\ft_*$ in Lemma~\ref{LemmaSDGTime}) produces an extension of the operator $P$ to a wave-type operator $\tilde P$ on $\Omega^\sharp$.

  \begin{figure}[!ht]
  \centering
  \includegraphics{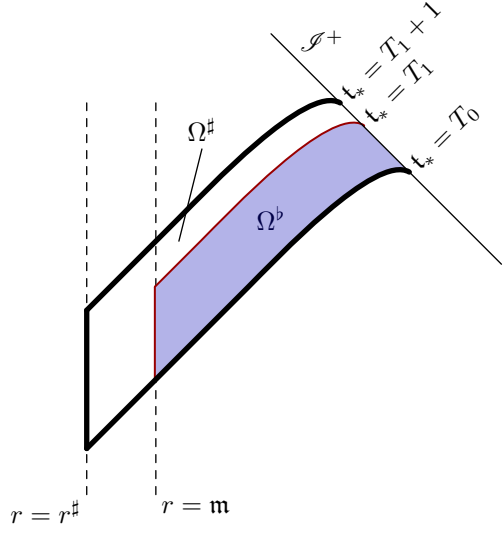}
  \caption{The domain $\Omega^\flat$ (in blue) and its extension $\Omega^\sharp$ (with thick boundary) beyond the final hypersurfaces (in red) from~\eqref{EqETSolvOmegas}.}
  \label{FigETSolvOmegas}
  \end{figure}

  \pfstep{$1\leq s\in\R$.} Extend a given $f\in H_\eop^{s-1,2\alpha_\sscri+2}(\Omega^\flat)^{\bullet,-}$ to $\tilde f\in H_\eop^{s-1,2\alpha_\sscri+2}(\Omega^\sharp)^{\bullet,-}$ with the same norm. Solve the forward problem for $\tilde P\tilde u=\tilde f$ on $\Omega^\sharp$; this gives $\tilde u\in H_\eop^{0,2\alpha_\sscri}(\Omega^\sharp)^{\bullet,-}$ with norm bounded by that of $\tilde f$. Microlocal elliptic regularity, propagation of edge-regularity (starting from $\ft_*<T_0$ where $u$ vanishes), and the localized radial sink estimate, Proposition~\ref{PropR3RScriO}, imply
  \[
    \|\tilde u|_{\Omega^\flat}\|_{H_\eop^{s,2\alpha_\sscri}(\Omega^\flat)} \leq C\Bigl( \|\tilde f\|_{H_\eop^{s-1,2\alpha_\sscri+2}(\Omega^\sharp)} + \|\tilde u\|_{H_\eop^{1,2\alpha_\sscri}(\Omega^\sharp)}\Bigr).
  \]
  Estimating $\tilde u$ on the right by $\|\tilde f\|_{H_\eop^{0,2\alpha_\sscri+2}(\Omega^\sharp)}$, we conclude that $u=\tilde u|_{\Omega^\flat}$ satisfies~\eqref{EqETSolvEst}.

  \pfstep{Existence and uniqueness for $P^*$, $s^*\leq 0$.} Dualization of~\eqref{EqETSolvEst}, now regarded as an a priori estimate for $u\in H_\eop^{s+1,2\alpha_\sscri}(\Omega^\flat)^{\bullet,-}$, say, yields (via the Hahn--Banach theorem) the \emph{solvability} of $P^*u^*=f^*$, given $f^*\in H_\eop^{s^*-1,-2\alpha_\sscri}(\Omega^\flat)^{-,\bullet}$, $s_*:=-s+1\leq 0$, with $u^*\in H_\eop^{s^*,-2\alpha_\sscri-2}(\Omega^\flat)^{-,\bullet}$ having norm bounded by $C$ times the norm of $f^*$. The solvability of $P u=f$ for all $f\in H_\eop^{s,2\alpha_\sscri+2}(\Omega^\flat)^{\bullet,-}$ moreover implies the uniqueness of $u^*$: indeed, if $u^*\in H_\eop^{s^*,-2\alpha_\sscri-2}(\Omega^\flat)^{\bullet,-}$ satisfies $P^*u^*=0$, then for $f\in\CIc((\Omega^\flat)^\circ)$, let $u\in H_\eop^{s+1,2\alpha_\sscri}(\Omega^\flat)^{\bullet,-}$ be the forward solution of $P u=f$; we then compute the $L^2$-pairing
  \[
    \la u^*,f\ra = \la u^*,P u\ra = \la P^*u^*,u\ra = 0.
  \]
  Since $f$ is arbitrary, this gives $u^*=0$.

  \pfstep{Existence and uniqueness for $P^*$, $s^*\in\R$.} An extension/restriction argument, now extending to the closure of $\{T_0-1\leq\ft_*\leq T_1,\ r\geq\bhm\}$ in $\tilde M$, combined with the microlocal edge-regularity results for $P^*$ (in particular, a localized estimate at the radial source $\pa\cR^+_{\scri^+,{\rm out}}$ for $P^*$-propagation analogous to Proposition~\ref{PropR3RKO}, whose proof we leave to the reader) for $u^*$ vanishing when $\ft_*>T_1$ or $r<\bhm$, analogously to the first step of the proof now shows: for all $s^*\in\R$, there exists a constant $C$ (uniform for perturbations of $P$) such that for all $f^*\in H_\eop^{s^*-1,-2\alpha_\sscri}(\Omega^\flat)^{-,\bullet}$ there exists a unique backward solution $u^*\in H_\eop^{s^*,-2\alpha_\sscri-2}(\Omega^\flat)^{-,\bullet}$ of $P^*u^*=f^*$, with norm bounded by $C$ times the norm of $f^*$.

  \pfstep{Proof for general $s\in\R$.} This follows by duality from the previous step for $s^*=-s+1$.

  \pfstep{Proof for general admissible orders $\sfs$.} We again use an extension/restriction argument as in the first step: from what we have already shown, we can solve the extended equation $\tilde P\tilde u=\tilde f\in H_\eop^{\sfs-1,2\alpha_\sscri+2}(\Omega^\sharp)^{\bullet,-}$ on $\Omega^\sharp$ with $\tilde u\in H_\eop^{s_0,2\alpha_\sscri}(\Omega^\sharp)^{\bullet,-}$ where $s_0:=\inf_{\Se^*_{\Omega^\sharp}M}\sfs$. Microlocal regularity results then imply that $u=\tilde u|_{\Omega^\flat}$ lies in $H_\eop^{\sfs,2\alpha_\sscri}(\Omega^\flat)^{\bullet,-}$.
\end{proof}

\subsubsection{Estimates on spaces with additional b-regularity}

We next generalize Proposition~\ref{PropETSolv} to the solvability (with estimates) on spaces encoding additional integer amounts of b-regularity; recall the notation~\eqref{EqETMixed}. To test for b-regularity, we use the operators $X_1,\ldots,X_5$ and $X_6=I$ from~\eqref{EqRbCommV}--\eqref{EqRbCommX}.

\begin{prop}[Solvability and uniqueness on $(\eop;\bop)$-Sobolev spaces]
\label{PropETSolvb}
  Let $\sfs\in\R$ be an admissible order function for $P_0$ (Definition~\usref{DefDyO}). Let $T_0,T_1,\alpha_\sscri$ be as in Proposition~\usref{PropETSolv}. Then there exists $d_0$ such that the following holds for all $k\in\N_0$: there exists a constant $C$ such that the unique forward solution of $P u=f$ on $\Omega_{T_0,T_1}$ satisfies the estimate
  \[
    \|u\|_{H_{\eop;\bop}^{(s;k),2\alpha_\sscri}(\Omega_{T_0,T_1};\cE)^{\bullet,-}} \leq C\|f\|_{H_{\eop;\bop}^{(s-1;k),2\alpha_\sscri+2}(\Omega_{T_0,T_1};\cE)^{\bullet,-}},
  \]
  provided $P$ is of class $((d_0;k),(2\ell_\sscri,0,0))$ (Definition~\usref{DefSDWAdm}). This estimate holds uniformly for sufficiently small perturbations of $P$ as measured in the norm $\|\cdot\|_{(d_0;k),(2\ell_\sscri,0,0)}$.
\end{prop}
\begin{proof}
  For $k=0$, this is the content of Proposition~\ref{PropETSolv}. Let $k\geq 1$, and suppose we have proved the result for $k-1\geq 0$ in place of $k$. Then by Lemma~\ref{LemmaRbComm} and using its notation, $u^{(k)}=(\vec X^\alpha u)_{\alpha\in\cA_k}$ satisfies a wave-type equation $P^{(k)}u^{(k)}=f^{(k)}+\tilde f^{(k)}$ where the quantity $\ubar p_1$ for $P^{(k)}$ equals that for $P$ (in view of the decay of $P^{\sharp,(k)}_{\alpha\alpha'}$ and $\tilde P^{(k)}_{\alpha\alpha'}$ in~\eqref{EqRbCommLot} at $\scri^+$); therefore, setting $\cE^{(k)}:=\cE\oplus\cdots\oplus\cE$ (with $|\cA_k|$ many summands), we have
  \[
    \|u^{(k)}\|_{H_\eop^{s,2\alpha_\sscri}(\Omega_{T_0,T_1};\cE^{(k)})^{\bullet,-}} \leq C_k\Bigl(\|f^{(k)}\|_{H_\eop^{s-1,2\alpha_\sscri+2}(\Omega_{T_0,T_1};\cE^{(k)})^{\bullet,-}} + \|\tilde f^{(k)}\|_{H_\eop^{s-1,2\alpha_\sscri+2}(\Omega_{T_0,T_1};\cE^{(k)})^{\bullet,-}}\Bigr).
  \]
  The inductive hypothesis and the left-hand side together control $\|u\|_{H_{\eop;\bop}^{(s;k),2\alpha_\sscri}(\Omega_{T_0,T_1};\cE)^{\bullet,-}}$, and the norm of $f^{(k)}$ is controlled by $\|f\|_{H_{\eop;\bop}^{(s-1;k),2\alpha_\sscri+2}(\Omega_{T_0,T_1};\cE)^{\bullet,-}}$. For the term $\tilde f^{(k)}$, we use the schematic description~\eqref{EqRbCommSchema}--\eqref{EqRbCommSchemaWeak} of its structure. A fortiori, we thus see that both terms are of the form $(D_\bop^{\leq j+1}p)D_\bop^{\leq k-j}(D_\eop^{\leq 1}u)$ for $j=1,\ldots,k-1$, and we can then estimate
  \begin{equation}
  \label{EqETSolvbEst}
  \begin{split}
    \|(D_\bop^{\leq j+1} p)D_\bop^{\leq k-j}(D_\eop^{\leq 1}u)\|_{H_\eop^{s-1,2\alpha_\sscri+2}(\Omega_{T_0,T_1};\cE)^{\bullet,-}} &\leq C\|D_\bop^{\leq k-j}(D_\eop^{\leq 1}u)\|_{H_\eop^{s-1,2\alpha_\sscri}(\Omega_{T_0,T_1};\cE)^{\bullet,-}} \\
      &\leq C'\|u\|_{H_{\eop;\bop}^{(s;k-j),2\alpha_\sscri}(\Omega_{T_0,T_1};\cE)^{\bullet,-}}.
  \end{split}
  \end{equation}
  In the first line, we use that multiplication by a differentiated coefficient $D_\bop^{j+1}p\in x_\sscri^2\cC_{\eop;\bop}^{(d_0;k-j-1)}\subset x_\sscri^2\cC_\eop^{d_0}$ of $P$ defines a bounded map on $H_\eop^{s-1,2\alpha_\sscri-2}(\Omega_{T_0,T_1};\cE)^{\bullet,-}\to H_\eop^{s-1,2\alpha_\sscri}(\Omega_{T_0,T_1};\cE)^{\bullet,-}$ when $d_0$ is large enough (depending only on $s$); see Lemma~\ref{LemmaETMult} below.
\end{proof}

Mirroring Remark~\ref{RmkRbBetter}, we note that an alternative proof of Proposition~\ref{PropETSolvb}, in which the case of orders $(s;k)$ is rather straightforwardly deduced from the case $(s+1;k-1)$, is given in \cite[Corollary~6.6]{HintzVasyScrieb}. That argument, however, does not generalize well when one needs tame estimates in the b-regularity order $k$. --- To complete the proof, we need to establish the boundedness of multiplication operators on Sobolev spaces with extendible character. The following version of Lemma~\ref{LemmaETExt} will be useful for this purpose:

\begin{lemma}[Another extension operator]
\label{LemmaETExtOther}
  Let $r^\sharp\in(0,\bhm)$. Setting
  \begin{equation}
  \label{EqETExtOther}
  \begin{split}
    \Omega_0 &:= \Omega_{T_0,T_1} = \ol{\{T_0\leq\ft_*\leq T_1,\ r\geq\bhm\}} \subset \tilde M, \\
    \Omega_1 &:= \ol{\{T_0-1\leq\ft_*\leq T_1+1,\ r\geq r^\sharp\}} \subset \tilde M,
  \end{split}
  \end{equation}
  there exists a continuous linear extension operator
  \[
    E \colon \CIc(\Omega_0\setminus\scri^+) \to \CIc(\Omega_1\setminus\scri^+)
  \]
  which, for all $d_0,s,k\in\N_0$ and $\alpha\in\R$, is continuous as a linear map
  \begin{equation}
  \label{EqETExtOther2}
  \begin{alignedat}{2}
    E &\colon& x_\sscri^{2\alpha}\cC_{\eop;\bop}^{(d_0;k)}(\Omega_0) &\to x_\sscri^{2\alpha}\cC_{\eop;\bop}^{(d_0;k)}(\Omega_1), \\
    E &\colon& x_\sscri^{2\alpha}H_{\eop;\bop}^{(s;k)}(\Omega_0)^{\bullet,-} &\to x_\sscri^{2\alpha}H_{\eop;\bop}^{(s;k)}(\Omega_1)^{\bullet,-}.
  \end{alignedat}
  \end{equation}
\end{lemma}
\begin{proof}
  Extend across $\ft_*=T_0,T_1$ with parametric dependence on the remaining coordinates $x_\sscri=r^{-\frac12}$, $\omega\in\Sph^2$, and extend across $r=\bhm$ as in the proof of Lemma~\ref{LemmaETExt}. Note here that one can indeed use $\pa_{\ft_*}$, $x_\sscri\pa_{x_\sscri}$, and $x_\sscri\pa_\omega$, resp.\ $\pa_\omega$ to test for edge-, resp.\ b-regularity, as follows from a change of variables computation using the form of $\ft_*$ near $\scri^+$ given in Lemma~\ref{LemmaSDGTime}\eqref{ItSDGTime3}.
\end{proof}

\begin{lemma}[Multiplication on extendible/supported spaces]
\label{LemmaETMult}
  Write $\Omega_0:=\Omega_{T_0,T_1}$. Let $\sfs\in\CI({}^\etbop S^*\tilde M)$. Then there exists $d_0\in\N$ such that the following statement holds.
  \begin{enumerate}
  \item\label{ItETMult1} Multiplication by $\ell\in\cC_\eop^{d_0}(\Omega_0)$ defines a bounded linear map $H_\eop^\sfs(\Omega_0)^{\bullet,-}\to H_\eop^\sfs(\Omega_0)^{\bullet,-}$.
  \item\label{ItETMult2} More generally, for every $k\in\N_0$, multiplication by $\ell\in\cC_{\eop;\bop}^{(d_0;k)}(\Omega_0)$ defines a bounded linear map $H_{\eop;\bop}^{(\sfs;k)}(\Omega_0)^{\bullet,-}\to H_{\eop;\bop}^{(\sfs;k)}(\Omega_0)^{\bullet,-}$.
  \item\label{ItETMult3} There exists $\tilde d\in\N_0$ such that for all $k\in\N_0$ and using Notation~\usref{NotMTameb}, we have the (b-tame) estimate
    \begin{equation}
    \label{EqETMultTame}
      \|(D_\bop^j\ell)(D_\bop^{k-j}u)\|_{H_\eop^\sfs(\Omega_0)^{\bullet,-}} \leq C_k\Bigl( \|\ell\|_{\cC_{\eop;\bop}^{(d_0;\tilde d)}(\Omega_0)}\|u\|_{H_{\eop;\bop}^{(\sfs;k)}(\Omega_0)^{\bullet,-}} + \|\ell\|_{\cC_{\eop;\bop}^{(d_0;k)}(\Omega_0)}\|u\|_{H_{\eop;\bop}^{(\sfs;\tilde d)}(\Omega_0)^{\bullet,-}}\Bigr).
    \end{equation}
  \end{enumerate}
  Analogous statements hold, mutatis mutandis, for weighted $\ell$ and weighted Sobolev spaces.
\end{lemma}
\begin{proof}
  \pfstep{Part~\eqref{ItETMult1}.} Define $\Omega_1$ as in~\eqref{EqETExtOther}, and set $\ell^\sharp=E\ell$ where $E$ is the extension operator from Lemma~\ref{LemmaETExtOther}. We use the quotient norms on $H_{\eop;\bop}^{(\sfs;k)}(\Omega_0)^{\bullet,-}$ induced by the restriction map from $\dot H_{\eop;\bop}^{(\sfs;k)}(\Omega^\sharp)$. Let $\chi\in\CIc(\tilde M)$ be equal to $1$ on $\Omega_0$ and equal to $0$ near $\tilde M\setminus(\Omega_1^\circ\cup\scri^+)$.

  We need to show that $\chi \ell^\sharp u\in H_\eop^\sfs$ when $u\in\dot H_\eop^\sfs(\Omega^\sharp)$. But when $\ell\in\cC_\eop^\infty$, multiplication by $\chi\ell^\sharp$ is an element of\footnote{Recall our present convention to drop any mention of b- or 3b-notions since we are currently working on domains such as $\Omega_1$ which intersect $\pa\tilde M$ only in $(\scri^+)^\circ$ (where the operator Lie algebra of interest is the edge-algebra, cf.\ \eqref{EqMUe3bFam3}, where we now restrict to for bounded $i$ and thus bounded $t_*=\tau^{-1}$). This is why here we write $\tilde\Psi_\eop^0$ instead of $\tilde\Psi_\etbop^0$, see~\eqref{EqMUe3bTilde}; for present purposes, one may simply define $\tilde\Psi_\eop^0$ as the subspace of $\tilde\Psi_\etbop^0$ consisting of all operators whose Schwartz kernels are supported in $\Omega_1\times\Omega_1$.} $\tilde\Psi_\eop^0$ and thus bounded on $H_\eop^\sfs$. The same then holds true when the edge-regularity of $\ell^\sharp$ is only finite but sufficiently large (depending on $\sfs$); see Lemma~\ref{LemmaMSOpNorm}.

  \pfstep{Part~\eqref{ItETMult2}.} This follows from the first part and the Leibniz rule.

  \pfstep{Part~\eqref{ItETMult3}.} The proof proceeds similarly to that of Proposition~\ref{PropMTameMicr}.\footnote{One may be tempted to directly apply Proposition~\ref{PropMTameMicr} to extensions of $\ell$ and $u$. Unfortunately, we are unable to construct a single extension operator that respects all \emph{variable} differential orders $\sfs$, which would be needed to implement such an argument.} We prove the analogue of Lemma~\ref{LemmaMTameMult} and consider the case $\sfs=0$. Let $u^\sharp=E u$. By multiplying $\ell^\sharp$ and $u^\sharp$ by smooth cutoff functions that equal $1$ on $\Omega_0$, we may assume that they are supported in a neighborhood of $\Omega_0$. We can then use Lemma~\ref{LemmaMTameMult} to obtain
  \begin{align*}
    \|(D_\bop^j\ell)(D_\bop^{k-j}u)\|_{L^2(\Omega_0)^{\bullet,-}} &\leq C\|\chi(D_\bop^j\ell^\sharp)(D_\bop^{k-j}u^\sharp)\|_{L^2(\Omega_1)} \\
      &\leq C_k\Bigl( \|\ell^\sharp\|_{\cC_\bop^0(\Omega_1)}\|u^\sharp\|_{H_\bop^k(\Omega_1)^{\bullet,-}} + \|\ell^\sharp\|_{\cC_\bop^k(\Omega_1)}\|u^\sharp\|_{H_\bop^{(0;3)}(\Omega_1)^{\bullet,-}}\Bigr).
  \end{align*}
  By the continuity properties of $E$ asserted in Lemma~\ref{LemmaETExtOther}, this implies
  \begin{equation}
  \label{EqETMultL2}
    \|(D_\bop^j\ell)(D_\bop^{k-j}u)\|_{L^2(\Omega_0)^{\bullet,-}} \leq C'_k\Bigl(\|\ell\|_{\cC_\bop^0(\Omega_0)}\|u\|_{H_\bop^k(\Omega_0)^{\bullet,-}} + \|\ell\|_{\cC_\bop^k(\Omega_0)}\|u\|_{H_\bop^3(\Omega_0)^{\bullet,-}}\Bigr).
  \end{equation}

  Returning to the proof of~\eqref{EqETMultTame} for general orders $\sfs$, consider first the case $k=0$. Part~\eqref{ItETMult1} gives
  \begin{equation}
  \label{EqETMultTame0}
    \|\ell u\|_{H_\eop^\sfs(\Omega_0)^{\bullet,-}} \leq C\|\ell\|_{\cC_\eop^{d_0}(\Omega_0)}\|u\|_{H_\eop^\sfs(\Omega_0)^{\bullet,-}}.
  \end{equation}
  For $1\leq k\leq\tilde d$, with $\tilde d$ specified later, part~\eqref{ItETMult1} gives, for $0\leq j\leq k$,
  \[
    \|(D_\bop^j\ell)(D_\bop^{k-j}u)\|_{H_\eop^\sfs(\Omega_0)^{\bullet,-}} \leq C\|D_\bop^j\ell\|_{\cC_\eop^{d_0}(\Omega_0)}\|D_\bop^{k-j}u\|_{H_\eop^\sfs(\Omega_0)^{\bullet,-}} \leq C\|\ell\|_{\cC_{\eop;\bop}^{(d_0;\tilde d)}}\|u\|_{H_{\eop;\bop}^{(\sfs;\tilde d)}(\Omega_0)^{\bullet,-}},
  \]
  which is the estimate~\eqref{EqETMultTame}. Consider next $k>\tilde d$. When $j\leq\tilde d$ or $k-j\leq\tilde d$, we use~\eqref{EqETMultTame0} with $0$, $D_\bop^j\ell$, $D_\bop^{k-j}u$ in place of $k$, $\ell$, $u$, to bound
  \[
    \|(D_\bop^j\ell)(D_\bop^{k-j}u)\|_{H_\eop^\sfs(\Omega_0)^{\bullet,-}} \leq C \|\ell\|_{\cC_{\eop;\bop}^{(d_0;j)}(\Omega_0)}\|u\|_{H_{\eop;\bop}^{(\sfs;k-j)}(\Omega_0)^{\bullet,-}},
  \]
  which in turn implies the estimate~\eqref{EqETMultTame}. For $\tilde d<j<k-\tilde d$, we crudely bound the $H_\eop^\sfs$-norm by the $H_\bop^{\bar s}$-norm where $\bar s:=\max(0,\lceil\sup\sfs\rceil)$ and thus need to bound
  \[
    \|(D_\bop^j\ell)(D_\bop^{k-j}u)\|_{H_\eop^\sfs(\Omega_0)^{\bullet,-}} \leq C\sum_{l=0}^{\bar s} \bigl\|D_\bop^l\bigl( (D_\bop^j\ell)(D_\bop^{k-j}u) \bigr)\|_{L^2(\Omega_0)} \leq C'\sum_{l=0}^{\bar s}\| (D_\bop^{j+l}\ell)(D_\bop^{k-j+\bar s-l}u)\|_{L^2}.
  \]
  For $d_1\leq\frac{\tilde d}{2}$ (so $2 d_1<k$) chosen momentarily, we bound this using~\eqref{EqETMultL2} with $j+l-d_1$, $k-j+\bar s-l-d_1$, $D_\bop^{d_1}\ell$, $D_\bop^{d_1}u$ in place of $j$, $k$, $\ell$, $u$, and thus by a constant times
  \[
    \|\ell\|_{\cC_\bop^{d_1}(\Omega_0)}\|u\|_{H_\bop^{k+\bar s-d_1}(\Omega_0)^{\bullet,-}} + \|\ell\|_{\cC_\bop^{k+\bar s-d_1}(\Omega_0)}\|u\|_{H_\bop^{d_1+3}(\Omega_0)^{\bullet,-}}
  \]
  Writing $\ubar s=\max(0,\lceil-\inf\sfs\rceil)$, this is further bounded by
  \[
    \|\ell\|_{\cC_\bop^{d_1}(\Omega_0)}\|u\|_{H_{\eop;\bop}^{(\sfs;k+\bar s+\ubar s-d_1)}(\Omega_0)^{\bullet,-}} + \|\ell\|_{\cC_\bop^{k+\bar s-d_1}(\Omega_0)}\|u\|_{H_{\eop;\bop}^{(\sfs;d_1+\ubar s+3)}(\Omega_0)^{\bullet,-}}.
  \]
  Choosing $d_1$ and then $\tilde d$ such that $d_1\geq\bar s+\ubar s$ (thus also $k+\bar s-d_1\leq k$) and $2 d_1\leq\tilde d$, $d_1+\ubar s+3\leq\tilde d$, this is, finally, bounded by the right-hand side of~\eqref{EqETMultTame}.
\end{proof}

For large b-regularity orders $k$, we then have the following tame version of Proposition~\ref{PropETSolvb}.

\begin{prop}[Solvability and uniqueness with b-tame estimates]
\label{PropETSolvbTame}
  Let $\sfs\in\N$ be an admissible order function for $P_0$ (Definition~\usref{DefDyO}). Let $T_0,T_1,\alpha_\sscri$ be as in Proposition~\usref{PropETSolv}. Then there exist $d_0,d\in\N$ such that the following holds for all $k\in\N_0$: there exists a constant $C_k$ such that the unique forward solution of $P u=f$ on $\Omega_{T_0,T_1}$ satisfies the tame estimate
  \begin{equation}
  \label{EqETSolvbTame}
  \begin{split}
    \|u\|_{H_{\eop;\bop}^{(\sfs;k),2\alpha_\sscri}(\Omega_{T_0,T_1};\cE)^{\bullet,-}} &\leq C_k\Bigl( \|f\|_{H_{\eop;\bop}^{(\sfs-1;k),2\alpha_\sscri+2}(\Omega_{T_0,T_1};\cE)^{\bullet,-}} \\
      &\quad \hspace{3em} + \|P-P_0\|_{(d_0;k),(2\ell_\sscri,0,0)}\|f\|_{H_{\eop;\bop}^{(\sfs-1;d),2\alpha_\sscri+2}(\Omega_{T_0,T_1};\cE)^{\bullet,-}} \Bigr),
  \end{split}
  \end{equation}
  provided $P$ is of class $((d_0;k),(2\ell_\sscri,0,0))$ (Definition~\usref{DefSDWAdm}). This estimate holds uniformly for sufficiently small perturbations of $P$ as measured in the fixed low regularity norm $\|\cdot\|_{(d_0;d),(2\ell_\sscri,0,0))}$. In~\eqref{EqETSolvbTame}, the norm on $P-P_0$ can, in fact, be taken to be the norm\footnote{That is, the sup norm of the coefficients of $P-P_0$, as in Definition~\usref{DefSDWAdmNorm}, and its up to $d_0$-fold edge- and $k$-fold b-derivatives.} over $\Omega_{T_0,T_1}$ only.
\end{prop}

We remark that the variable nature of $\sfs$ will be used in our application since we shall use Proposition~\ref{PropETSolvbTame} on $t_*$-intervals which get arbitrarily large as $k$ increases; and we recall from Definition~\ref{DefSSOrderAdm} that $\sfs$ typically needs to be variable at $\iota^+$ and thus nearby.

\begin{proof}[Proof of Proposition~\usref{PropETSolvbTame}]
  We drop the bundle $\cE$ from the notation. With $d\geq 2$ to be fixed later, the estimate~\eqref{EqETSolvbTame} follows for $k\leq d$ from Proposition~\ref{PropETSolvb}; consider thus the case $k\geq d+1$. We then argue as in the proof of Proposition~\ref{PropETSolvb}, except we replace estimate~\eqref{EqETSolvbEst} by a tame version as follows. We estimate the two types of terms in~\eqref{EqRbCommSchema} separately. For the first term, we estimate
  \[
    \|(D_\bop^{\leq 1}p)D_\bop^{\leq k-1}(D_\ebop^{\leq 1}u)\|_{H_\eop^{\sfs-1,2\alpha_\sscri+2}(\Omega_{T_0,T_1})^{\bullet,-}} \leq \|p\|_{\cC_{\eop;\bop}^{(d_0;1)}}\|u\|_{H_{\eop;\bop}^{(\sfs;k-1),2\alpha_\sscri}(\Omega_{T_0,T_1})^{\bullet,-}}
  \]
  using Lemma~\ref{LemmaETMult}\eqref{ItETMult1}. For the second class of terms in~\eqref{EqRbCommSchema}, we use Lemma~\ref{LemmaETMult}\eqref{ItETMult3} for $\ell=D_\bop^{\leq 2}p$ and $s-1,j-2,k-2$ in place of $s,j,k$, by
  \begin{align*}
    &\|(D_\bop^{\leq j-2}D_\bop^{\leq 2}p)D_\bop^{\leq(k-2)-(j-2)}(D_\eop^{\leq 2}u)\|_{H_\eop^{\sfs-1,2\alpha_\sscri+2}(\Omega_{T_0,T_1})^{\bullet,-}} \\
    &\quad \leq C_k\Bigl( \|D_\bop^{\leq 2}p\|_{\cC_{\eop;\bop}^{(d_0;\tilde d)}} \|u\|_{H_{\eop;\bop}^{(\sfs+1;k-2),2\alpha_\sscri}(\Omega_{T_0,T_1})^{\bullet,-}} + \|D_\bop^{\leq 2}p\|_{\cC_{\eop;\bop}^{(d_0;k-2)}}\|u\|_{H_{\eop;\bop}^{(\sfs+1;\tilde d),2\alpha_\sscri}(\Omega_{T_0,T_1})^{\bullet,-}}\Bigr) \\
    &\quad \leq C'_k\Bigl( \|p\|_{\cC_{\eop;\bop}^{(d_0;\tilde d+2)}} \|u\|_{H_{\eop;\bop}^{(\sfs;k-1),2\alpha_\sscri}(\Omega_{T_0,T_1})^{\bullet,-}} + \|p\|_{\cC_{\eop;\bop}^{(d_0;k)}}\|u\|_{H_{\eop;\bop}^{(\sfs;\tilde d+1),2\alpha_\sscri}(\Omega_{T_0,T_1})^{\bullet,-}}\Bigr).
  \end{align*}
  In the passage to the final line, we use $H_{\eop;\bop}^{(a-1;b+1)}\subset H_{\eop;\bop}^{(a;b)}$ (which follows from $\Ve\subset\Vb$). If we set $d=\tilde d+2$, this gives~\eqref{EqETSolvbTame} upon estimating $\|u\|_{H_{\eop;\bop}^{(\sfs;k-1),2\alpha_\sscri}(\Omega_{T_0,T_1})^{\bullet,-}}$ using induction in $k$, and $\|u\|_{H_{\eop;\bop}^{(\sfs;\tilde d+1),2\alpha_\sscri}(\Omega_{T_0,T_1})^{\bullet,-}}\leq C\|f\|_{H_{\eop;\bop}^{(\sfs-1;\tilde d+1),2\alpha_\sscri+2}(\Omega_{T_0,T_1})^{\bullet,-}}$ from Proposition~\ref{PropETSolvb}. The uniformity of~\eqref{EqETSolvbTame} for perturbations of $P$ uses, moreover, the fact that the coefficients of the operator $P^{(k)}$ in Lemma~\ref{LemmaRbComm} involve only up to first order b-derivatives of those of $P$.
\end{proof}

\subsection{Energy estimates near the interior hypersurface \texorpdfstring{$r=\bhm$}{r=m}}
\label{SsER}

Using the timelike nature of $\dd r$ for $\bhm\leq r\leq r_\natural=\frac{\bhm+r_+}{2}\in(\bhm,r_+)$, we now prove energy estimates in the region $\{r\leq r_\natural\}$ in the black hole interior. These estimates are considerably easier to prove than those in~\S\ref{SsET} since exponential weights in $r$ give arbitrary amounts of positivity of $K$-currents, unlike the weights in $\ft_*$ in the proof of Lemma~\ref{LemmaETEnergy} which only gave partial positivity (cf.\ the second line of~\eqref{EqETKCurrent}).

We use the notation of Definition~\usref{DefETExt}, also for domains such as $\{\ft_*\geq T_0,\ r\in I\}$ for bounded intervals $I\subset[\bhm,\infty)$. Concretely, for
\[
  T_0 \geq 1,\quad 0 < r_1 < r_2 \leq r_\natural,
\]
we write
\begin{equation}
\label{EqERDomain}
  \Omega_{T_0;r_1,r_2} := \ol{\{ \ft_*\geq T_0,\ r_1\leq r\leq r_2 \}} \subset \tilde M;
\end{equation}
when $r_1\geq\bhm$, this is a submanifold with corners of $\Omega$, with initial hypersurfaces $\ft_*=T_0$ and $r=r_2$ and final hypersurface $r=r_1$. For $s,\alpha_\cK\in\R$,\footnote{We do not need variable orders in this part.} we correspondingly write
\[
  H_\cuop^{s,\alpha_\cK}(\Omega_{T_0;r_1,r_2})^{\bullet,-},\quad H_{\cuop;\bop}^{(s;k),\alpha_\cK}(\Omega_{T_0;r_1,r_2})^{\bullet,-}
\]
for spaces of distributions with supported, resp.\ extendible character at the initial, resp.\ final hypersurfaces of $\Omega_{T_0;r_1,r_2}$. (These spaces are \emph{the same} as $H_\etbop^{s,\alpha_\cK}(\Omega_{T_0;r_1,r_2})^{\bullet,-}$ and $H_{\etbop;\bop}^{(s;k),\alpha_\cK}(\Omega_{T_0;r_1,r_2})^{\bullet,-}$.) The analogue of Lemma~\ref{LemmaETEnergy} reads:

\begin{lemma}[First order energy estimate]
\label{LemmaEREnergy}
  Let $\alpha_\cK\in\R$. Then there exists a constant $C$ such that for all $u\in H_\cuop^{2,\alpha_\cK}(\Omega_{T_0;r_1,r_2})^{\bullet,-}$, we have
  \begin{equation}
  \label{EqEREnergy}
    \|u\|_{H_\cuop^{1,\alpha_\cK}(\Omega_{T_0;r_1,r_2})^{\bullet,-}} \leq C\|P u\|_{H_\cuop^{0,\alpha_\cK}(\Omega_{T_0;r_1,r_2})^{\bullet,-}}.
  \end{equation}
  (The left-hand side controls $t_*^{-\alpha_\cK}u$ and its derivatives along $\pa_{t_*}$ and $\pa_x$ in $L^2$ with the metric volume density.) The constant $C$ can be taken to be uniform for sufficiently small perturbations of $P$ as measured in the norm $\|\cdot\|_{(d_0;0),(0,0,\ell_\cK)}$ (Definition~\usref{DefSDWAdmNorm}, with $\ell_\cK>0$, cf.\ also Remark~\usref{RmkSDWRelax}) for $d_0=2$.
\end{lemma}
\begin{proof}
  Write $\rho_\cK:=t_*^{-1}$ for a local defining function of $\cK^+$. If $P=\Box_g$ is the scalar wave operator, we obtain~\eqref{EqEREnergy} immediately as from energy estimate with the vector field multiplier
  \[
    X = e^{2\digamma r}X_0,\quad
    X_0 = \rho_\cK^{-2\alpha_\cK}\nabla r.
  \]
  Note that $X_0$ is future timelike. The $K$-current of $X$ is
  \[
    {}^{(X)}\!K = e^{2\digamma r}\bigl( {}^{(X_0)}\!K + 2\digamma T(\nabla r,\nabla r) \bigr),
  \]
  and hence positive definite---in fact, $\gtrsim\digamma$ compared to the Riemannian metric $g_R=\dd t_*^2+\dd x^2$---when $\digamma$ is sufficiently large. The boundary terms at $r=r_2$ and $\ft_*=T_0$ vanish by the supported character of $u$, while the boundary term at $r=r_1$ (due to the final nature of this spacelike hypersurface) has the same sign as the bulk term from ${}^{(X)}\!K$ and can thus be dropped. As in the proof of Lemma~\ref{LemmaETEnergy}, we can regard the boundary terms as parts of the $K$-current of the sharply cut-off vector field $H(\ft_*-T_0)H(r-r_1)H(r_2-r)X$.

  For general wave-type operators $P$, one works with $\sfX=e^{2\digamma r}\sfX_0$ and $\sfX_0=\nabla^\cE_{X_0}$ where $\nabla^\cE$ is a smooth connection on $\cE$. Starting with~\eqref{EqETKEquality} (when $\cE$ is trivial, $P$ acts component-wise as the scalar wave operator, and $\nabla^\cE$ is the trivial connection) and then generalizing to general connections, an arbitrary but fixed positive definite fiber inner product $h_\cE$ on $\cE$, and adding lower order terms to $P$, one obtains (in keeping with our current convention to write ``$\cuop$'' in place of ``$\etbop$'')
  \begin{align*}
    &\|{}^{(X)}K(\nabla^\cE u,\nabla^\cE u)\|_{L^2(\Omega_{T_0;r_1,r_2};{}^\cuop T^*M\otimes\cE,g_R\otimes h_\cE)}^2 \\
    &\qquad \leq \|\rho_\cK^{\alpha_\cK}e^{-\digamma r}\sfX u\|_{L^2(\Omega_{T_0;r_1,r_2};\cE,h_\cE)}^2 + \frac12\|\rho_\cK^{-\alpha_\cK}e^{\digamma r}P u\|_{L^2(\Omega_{T_0;r_1,r_2};\cE,h_\cE)}^2 \\
    &\qquad \qquad + C\Bigl(\|\rho_\cK^{-\alpha_\cK}e^{\digamma r}u\|_{L^2(\Omega_{T_0;r_1,r_2};\cE)}^2 + \|\rho_\cK^{-\alpha_\cK}e^{\digamma r}\nabla^\cE u\|_{L^2(\Omega_{T_0;r_1,r_2};{}^\cuop T^*M\otimes\cE,g_R\otimes h_\cE)}^2\Bigr).
  \end{align*}
  The left-hand side controls $c\digamma\|\nabla^\cE u\|_{L^2}^2$ for some $\digamma$-independent constant $c>0$. For all sufficiently large $\digamma$, we thus have
  \begin{equation}
  \label{EqEREnergyPf}
    \|e^{\digamma r}\nabla^\cE u\|_{\rho_\cK^{\alpha_\cK}L^2(\Omega_{T_0;r_1,r_2})}^2 \leq C\digamma^{-1}\Bigl( \|e^{\digamma r}P u\|_{\rho_\cK^{\alpha_\cK}L^2(\Omega_{T_0;r_1,r_2})}^2 + \|e^{\digamma r}u\|_{\rho_\cK^{\alpha_\cK}L^2(\Omega_{T_0;r_1,r_2})}^2 \Bigr).
  \end{equation}
  We estimate the second term on the right as follows: using
  \[
    \dv_g(-\rho_\cK^{-2\alpha_\cK}e^{2\digamma r}\nabla r)=e^{2\digamma r}\bigl(-\dv_g(\rho_\cK^{-2\alpha_\cK}\nabla r)-2\digamma \rho_\cK^{-2\alpha_\cK} g^{-1}(\nabla r,\nabla r)\bigr)\geq c\digamma \rho_\cK^{-2\alpha_\cK}e^{2\digamma r},
  \]
  for some $c>0$ and for all sufficiently large $\digamma$, and applying the divergence theorem on $\Omega_{T_0;r_1,r_2}$, one obtains, for $u\in H_\cuop^{1,\alpha_\cK}(\Omega_{T_0;r_1,r_2})^{\bullet,-}$,
  \[
    \digamma\|e^{\digamma r}u\|_{\rho_\cK^{\alpha_\cK}L^2(\Omega_{T_0;r_1,r_2})}^2 \leq C\|e^{\digamma r}\nabla^\cE u\|_{\rho_\cK^{\alpha_\cK}L^2(\Omega_{T_0;r_1,r_2})}^2
  \]
  when $\digamma$ is sufficiently large. Plugging this into~\eqref{EqEREnergyPf} and fixing $\digamma$ to be sufficiently large completes the proof of~\eqref{EqEREnergy}.
\end{proof}

Writing $\Omega_{T_0,T_1;r_1,r_2}:=\Omega_{T_0;r_1,r_2}\cap\{\ft_*\leq T_1\}$, with $\ft_*=T_1$ being a final hypersurface of $\Omega_{T_0,T_1;r_1,r_2}$, one similarly proves the uniform estimate
\[
  \|u\|_{H_\cuop^{1,\alpha_\cK}(\Omega_{T_0,T_1;r_1,r_2};\cE)^{\bullet,-}} \leq C\|P u\|_{H_\cuop^{0,\alpha_\cK}(\Omega_{T_0,T_1;r_1,r_2};\cE)^{\bullet,-}}.
\]
Straightforward modifications of the arguments used in the proofs of Corollary~\ref{CorETSolv1} and Proposition~\ref{PropETSolv} now show:

\begin{prop}[Solvability and uniqueness on cusp-Sobolev spaces]
\label{PropERSolv}
  Let $s\in\R$. Then there exists $d_0\in\N$ such that for $P$ and all perturbations thereof which are sufficiently small as measured in the norm $\|\cdot\|_{(d_0;0),(0,0,\ell_\cK)}$, the following holds. There exists a constant $C$ such that for all $f\in H_\cuop^{s-1,\alpha_\cK}(\Omega_{T_0;r_1,r_2};\cE)^{\bullet,-}$, the unique forward solution $u$ of $P u=f$ satisfies $u\in H_\cuop^{s,\alpha_\cK}(\Omega_{T_0;r_1,r_2};\cE)^{\bullet,-}$ and the estimate
  \[
    \|u\|_{H_\cuop^{s,\alpha_\cK}(\Omega_{T_0;r_1,r_2};\cE)^{\bullet,-}} \leq C\|P u\|_{H_\cuop^{s-1,\alpha_\cK}(\Omega_{T_0;r_1,r_2};\cE)^{\bullet,-}}.
  \]
\end{prop}

The only microlocal regularity results needed here are microlocal elliptic and real principal type propagation (to propagate regularity in the direction of decreasing, resp.\ increasing $r$ for $P$, resp.\ $P^*$). Applying Proposition~\ref{PropERSolv} to $\chi u$ where $\chi=\chi(r)$ equals $1$ on $[r_1,r_1^+]$ and $0$ on $[r_2^-,r_2]$ where $r_1<r_1^+<r_2^-<r_2$ and noting that $\|[P,\chi]u\|_{\bar H_\cuop^{s-1,\alpha_\cK}(\Omega_{T_0;r_1,r_2};\cE)^{\bullet,-}}\leq C\|u\|_{\bar H_\cuop^{s,\alpha_\cK}(\Omega_{T_0;r_1^+,r_2^-};\cE)^{\bullet,-}}$ yields the continuation estimate
\begin{equation}
\label{EqERSolvCont}
  \|u\|_{\bar H_\cuop^{s,\alpha_\cK}(\Omega_{T_0;r_1,r_2};\cE)} \leq C\Bigl( \|P u\|_{\bar H_\cuop^{s-1,\alpha_\cK}(\Omega_{T_0;r_1,r_2};\cE)} + \|u\|_{\bar H_\cuop^{s,\alpha_\cK}(\Omega_{T_0;r_1^+,r_2};\cE)}\Bigr).
\end{equation}

The proof of higher b-regularity, i.e., regularity with respect to $t_*\pa_{t_*}$ and $\pa_x$, relies on Lemma~\ref{LemmaCT3bDil}, which implies in particular that $[t_*\pa_{t_*},P]\in(\CI+\cC_{\cuop;\bop}^{(d_0;k-1),\ell_\cK})\Diff_\cuop^2$ in $r\leq r_+$ when $P$ is of class $((d_0;k),(0,0,\ell_\cK))$. Concretely, fix operators
\[
  X_1,\ldots,X_4 \in \Diff_\cuop^1
\]
with principal parts given by the vector fields $V_1=t_*\pa_{t_*}$, $V_{1+j}=\pa_{x^j}$, e.g., $X_i=\nabla^\cE_{V_i}$. Let $X_5:=I$, and write $\vec X=(X_1,\ldots,X_4,X_5)$. For $k\in\N$, set $\cA_k:=\{\alpha\in\N_0^5\colon|\alpha|=k\}$. If $P u=f$, then, analogously to (but simpler than) Lemma~\ref{LemmaRbComm}, the vector $u^{(k)}:=(\vec X^\alpha u)_{\alpha\in\cA_k}$ satisfies the equation
\begin{equation}
\label{EqERCommuted}
  P^{(k)}u^{(k)} = f^{(k)} + \tilde f^{(k)}
\end{equation}
where $f^{(k)}=(\vec X^\alpha f)_{\alpha\in\cA_k}$ and, schematically,
\[
  P^{(k)}_{\alpha\alpha'} = \delta_{\alpha\alpha'}P + (D_\bop^{\leq 1}p)D_\cuop^{\leq 1}
\]
where $p\in\CI+\cC_{\cuop;\bop}^{(d_0;k),\alpha_\cK}$ denotes the coefficients of $P$ expressed in terms of the basic vector fields $\pa_{t_*},\pa_x$; and $\tilde f^{(k)}$ has the schematic form
\begin{equation}
\label{EqERtildef}
  \tilde f^{(k)} = \sum_{j=2}^k (D_\bop^{\leq j}p)D_\cuop^2 D_\bop^{\leq k-j}u.
\end{equation}
We then have the following generalizations of~\eqref{EqERSolvCont}.

\begin{prop}[Solvability and uniqueness on $(\cuop;\bop)$-Sobolev spaces]
\label{PropERSolvb}
  Let $s,\alpha_\cK\in\R$, and let $T_0$, $r_1<r_1^+<r_2\leq r_\natural$. Then there exists $d_0$ such that the following holds for all $k\in\N_0$: there exists a constant $C$ such that every solution of $P u=f$ on $\Omega_{T_0;r_1,r_2}$ satisfies the estimate
  \[
    \|u\|_{\bar H_{\cuop;\bop}^{(s;k),\alpha_\cK}(\Omega_{T_0;r_1,r_2};\cE)} \leq C\Bigl( \|P u\|_{\bar H_{\cuop;\bop}^{(s-1;k),\alpha_\cK}(\Omega_{T_0;r_1,r_2};\cE)} + \|u\|_{\bar H_{\cuop;\bop}^{(s;k),\alpha_\cK}(\Omega_{T_0;r_1^+,r_2};\cE)}\Bigr),
  \]
  provided $P$ is of class $((d_0;k),(0,0,\ell_\cK))$ (Definition~\usref{DefSDWAdm}). This estimate holds uniformly for sufficiently small perturbations of $P$ as measured in the norm $\|\cdot\|_{(d_0;k),(0,0,\ell_\cK)}$.
\end{prop}
\begin{proof}
  The estimate~\eqref{EqERSolvCont} applies also to the equation~\eqref{EqERCommuted}. It remains to bound the norm of $\tilde f^{(k)}$; but using~\eqref{EqERtildef}, we have, for $j=2,\ldots,k$,
  \[
    \|(D_\bop^{\leq j}p)D_\cuop^2 D_\bop^{\leq k-j}u\|_{\bar H_\cuop^{s-1,\alpha_\cK}(\Omega_{T_0;r_1,r_2};\cE)} \leq \|p\|_{\cC_{\cuop;\bop}^{(d_0;k)}}\|u\|_{\bar H_{\cuop;\bop}^{(s+1;k-2),\alpha_\cK}(\Omega_{T_0;r_1,r_2};\cE)}
  \]
  for sufficiently large $d_0$ (depending only on $s$), as follows from an analogue of Lemma~\ref{LemmaETMult} (which, in turn, is proved using extension operators across $r=r_1$ and $r=r_2$ as constructed\footnote{albeit for particular values of $r_1,r_2$ which is, however, irrelevant for the proof given there} in Lemma~\ref{LemmaETExt}). The norm on $u$ here is bounded by the $\bar H_{\cuop;\bop}^{(s;k-1),\alpha_\cK}$-norm, which can be controlled inductively.
\end{proof}

\begin{prop}[Solvability and uniqueness with b-tame estimates]
\label{PropERSolvbTame}
  Let $s,\alpha_\cK\in\R$ and $T_0,r_1<r_1^+<r_2\leq r_\natural$. Then there exist $d_0,d\in\N$ such that the following holds for all $k\in\N_0$: there exists a constant $C_k$ such that every solution of $P u=f$ on $\Omega_{T_0;r_1,r_2}$ satisfies the estimate
  \begin{equation}
  \label{EqERSolvbTame}
  \begin{split}
    &\|u\|_{\bar H_{\cuop;\bop}^{(s;k),\alpha_\cK}(\Omega_{T_0;r_1,r_2};\cE)} \\
    &\qquad \leq C_k\biggl( \|f\|_{\bar H_{\cuop;\bop}^{(s-1;k),\alpha_\cK}(\Omega_{T_0;r_1,r_2};\cE)} + \|u\|_{\bar H_{\cuop;\bop}^{(s;k),\alpha_\cK}(\Omega_{T_0;r_1^+,r_2};\cE)} \\
    &\qquad\hspace{3.5em} + \|P-P_0\|_{(d_0;k),(0,0,\ell_\cK)} \Bigl( \|f\|_{\bar H_{\cuop;\bop}^{(s-1;d),\alpha_\cK}(\Omega_{T_0;r_1,r_2};\cE)} + \|u\|_{\bar H_{\cuop;\bop}^{(s;d),\alpha_\cK}(\Omega_{T_0;r_1^+,r_2};\cE)}\Bigr) \biggr),
  \end{split}
  \end{equation}
  provided $P$ is of class $((d_0;k),(0,0,\ell_\cK))$. This estimate holds uniformly for sufficiently small perturbations of $P$ as measured in the fixed low regularity norm $\|\cdot\|_{(d_0;d),(0,0,\ell_\cK)}$. In~\eqref{EqERSolvbTame}, the norm on $P-P_0$ can be taken to be the norm over $\Omega_{T_0;r_1,r_2}$ only.
\end{prop}
\begin{proof}
  We omit the bundle $\cE$ from the notation. Set $\Omega_0:=\Omega_{T_0;r_1,r_2}$. We have, analogously to~\eqref{EqETMultTame} (and also with a similar proof), the tame product estimate
  \begin{equation}
  \label{EqERSolvbTameProd}
    \|(D_\bop^j\ell)(D_\bop^{k-j}u)\|_{\bar H_\cuop^s(\Omega_0)} \leq C_k\Bigl( \|\ell\|_{\cC_{\cuop;\bop}^{(d_0;\tilde d)}(\Omega_0)}\|u\|_{\bar H_{\cuop;\bop}^{(s;k)}(\Omega_0)} + \|\ell\|_{\cC_{\cuop;\bop}^{(d_0;k)}(\Omega_0)}\|u\|_{\bar H_{\cuop;\bop}^{(s;\tilde d)}(\Omega_0)}\Bigr)
  \end{equation}
  when $d_0,\tilde d$ are large enough (depending only on $s$). Applying the estimate~\eqref{EqERSolvCont} to the equation~\eqref{EqERCommuted}, we then estimate the $j$-th summand of $\tilde f^{(k)}$ in~\eqref{EqERtildef} using~\eqref{EqERSolvbTameProd} with $s-1$, $D_\bop^{\leq 2}p$, $D_\cuop^{\leq 2}u$, $k-2$, $j-2$ in place of $s$, $\ell$, $u$, $k$, $j$ by
  \begin{align*}
    &\|(D_\bop^{\leq j}p)D_\bop^{\leq k-j}(D_\cuop^{\leq 2}u)\|_{\bar H_\cuop^{s-1,\alpha_\cK}(\Omega_0)} \\
    &\qquad \leq C_k\Bigl( \|D_\bop^{\leq 2}p\|_{\cC_{\cuop;\bop}^{(d_0;\tilde d)}(\Omega_0)} \|D_\cuop^{\leq 2}u\|_{\bar H_{\cuop;\bop}^{(s-1;k-2)}(\Omega_0)} + \|D_\bop^{\leq 2}p\|_{\cC_{\cuop;\bop}^{(d_0;k-2)}(\Omega_0)}\|D_\cuop^{\leq 2}u\|_{\bar H_{\cuop;\bop}^{(s-1;\tilde d)}(\Omega_0)}\Bigr) \\
    &\qquad \leq C'_k\Bigl( \|p\|_{\cC_{\cuop;\bop}^{(d_0;d)}(\Omega_0)}\|u\|_{\bar H_{\cuop;\bop}^{(s;k-1)}(\Omega_0)} + \|p\|_{\cC_{\cuop;\bop}^{(d_0;k)}(\Omega_0)}\|u\|_{\bar H_{\cuop;\bop}^{(s;d)}(\Omega_0)}\Bigr)
  \end{align*}
  for $d:=\tilde d+2$. The $\bar H_{\cuop;\bop}^{(s;d)}$-norm of $u$ is controlled by Proposition~\ref{PropERSolvb}, and the $\bar H_{\cuop;\bop}^{(s;k-1)}$-norm can be controlled via induction on $k$.
\end{proof}

\subsection{Global e3b-regularity estimate}
\label{SsEReg}

Propositions~\ref{PropETSolvb} and \ref{PropETSolvbTame}, applied on a $\ft_*$-interval $[1,4]$, provide (b-tame) estimates for the forward solution $u$ of $P u=f$ on $\Omega_*=\cl_M\{\ft_*\geq 1\}$ near the initial hypersurface $\ft_*^{-1}(1)$ of $\Omega_*$. Proposition~\ref{PropRb} provides (b-tame) estimates for $u$ on $\chi^{-1}(1)\supset\cl_M\{\ft_*\geq 3,\ r\geq r_1^+\}$ (in the notation of~\eqref{EqRRadii} and Proposition~\ref{PropRb}). Finally, Propositions~\ref{PropERSolvb} and \ref{PropERSolvbTame}, applied with $r_1=\bhm$, $r_1^+$ from~\eqref{EqRRadii}, and $r_2=r_\natural$ propagates the (b-tame) estimates for $u$ on $\Omega_*\cap r^{-1}([r_1^+,r_\natural])$ to $\Omega_*\cap r^{-1}([\bhm,r_\natural])$. This proves:

\begin{thm}[Global b-tame e3b-regularity estimate]
\label{ThmEReg}
  Let $\alpha_\sscri,\alpha_+,\alpha_\cK\in\R$, and let $\sfs_0,\sfs\in\CI({}^\etbop S^*M)$ with $\sfs>\sfs_0$. Let $P$ be a weakly admissible wave-type operator (relative to a stationary wave-type operator $P_0$) in the sense of Definition~\usref{DefSDWAdm}. Suppose that $\alpha_\sscri<-\frac12+\ubar p_1$ in the notation of Definition~\usref{DefSDWubarp1}, and that $\sfs,\sfs_0$ are admissible order functions for $P_0$ with weights $\alpha_+,\alpha_\cK$ and margin $0$ (Definition~\usref{DefDyO}). Then for all $k\in\N_0$ there exists a constant $C_k$ such that
  \begin{equation}
  \label{EqEReg}
  \begin{split}
    &\|u\|_{H_{\etbop;\bop}^{(\sfs;k),(2\alpha_\sscri,\alpha_+,\alpha_\cK)}(\Omega_*;\cE)^{\bullet,-}} \\
    &\qquad \leq C_k\Bigl( \|P u\|_{H_{\etbop;\bop}^{(\sfs;k),(2\alpha_\sscri+2,\alpha_++2,\alpha_\cK)}(\Omega_*;\cE)^{\bullet,-}} + \|u\|_{H_{\etbop;\bop}^{(\sfs_0;k),(2\alpha_\sscri,\alpha_+,\alpha_\cK)}(\Omega_*;\cE)^{\bullet,-}} \Bigr)
  \end{split}
  \end{equation}
  holds for $P=P_0$ in the strong sense, and also for all weakly admissible wave-type operators $P$ relative to $P_0$ which are sufficiently small perturbations as measured in the norm $\|\cdot\|_{(d_0;k),(2\ell_\sscri,\ell_+,\ell_\cK)}$ (see Definition~\usref{DefSDWAdmNorm}) where $d_0\in\N$ is sufficiently large and depends only on $\sfs_0,\sfs$. Furthermore, there exists $d\in\N$ such that for $P=P_0$ and all of its sufficiently small perturbations as measured in $\|\cdot\|_{(d_0;d),(2\ell_\sscri,\ell_+,\ell_\cK)}$, one has, for all $k\in\N_0$, the b-tame estimate
  \begin{equation}
  \label{EqERegTame}
  \begin{split}
    &\|u\|_{H_{\etbop;\bop}^{(\sfs;k),(2\alpha_\sscri,\alpha_+,\alpha_\cK)}(\Omega_*;\cE)^{\bullet,-}} \\
    &\quad \leq C_k\biggl( \|P u\|_{H_{\etbop;\bop}^{(\sfs;k),(2\alpha_\sscri,\alpha_+,\alpha_\cK)}(\Omega_*;\cE)^{\bullet,-}} + \|u\|_{H_{\etbop;\bop}^{(\sfs_0;k),(2\alpha_\sscri,\alpha_+,\alpha_\cK)}(\Omega_*;\cE)^{\bullet,-}} \\
    &\quad \quad \hspace{3em} + \|P-P_0\|_{(d_0;k),(2\ell_\sscri,\ell_+,\ell_\cK),\Omega_*} \\
    &\quad \quad \hspace{4.5em} \times \Bigl( \|P u\|_{H_{\etbop;\bop}^{(\sfs;d),(2\alpha_\sscri,\alpha_+,\alpha_\cK)}(\Omega_*;\cE)^{\bullet,-}} + \|u\|_{H_{\etbop;\bop}^{(\sfs_0;d),(2\alpha_\sscri,\alpha_+,\alpha_\cK)}(\Omega_*;\cE)^{\bullet,-}}\Bigr) \biggr),
  \end{split}
  \end{equation}
  provided $P$ is of class $((d_0;k),(2\ell_\sscri,\ell_+,\ell_\cK))$. This holds also when $P_0$ is replaced by sufficiently small perturbation within the class of stationary wave-type operators (including for different but nearby Kerr parameters).
\end{thm}

This theorem provides \emph{full control of $u$ at high e3b-frequencies.} Our main outstanding task in our analysis of admissible wave-type operators is to improve (i.e., weaken) the norms on $u$ on the right-hand sides of~\eqref{EqEReg}--\eqref{EqERegTame} in the sense of decay at $\iota^+\cup\cK^+=t_*^{-1}(\infty)$; this is where the spectral assumptions on the stationary operator $P_0$ enter. We turn to this in~\S\S\ref{SN}--\ref{SF}.

\section{General-purpose estimates for the spectral family}
\label{SSp}

While we have black-boxed the requirements on the mapping properties of $P_0$ in Definition~\ref{DefSSAlephAdm}\eqref{ItSSAlephAdm3} rather than stating assumptions on its spectral family $\wh{P_0}(\sigma)$ (see~\eqref{EqSSSpecFam}), it is useful to record estimates that the spectral family $\wh{P_0}(\sigma)$ satisfies and which follow solely from the dynamical properties of the Kerr metric. (One reason for this is that estimates for the transition face normal operators in~\S\ref{SsSptf} will be needed in our analysis of $P_0$ near $\iota^+$ in~\S\ref{SN} below.) These are:

\begin{enumerate}
\item the Fredholm and index $0$ property of $\wh{P_0}(\sigma)$ at nonzero energies $\sigma$, $\Im\sigma\geq 0$ (\S\ref{SsSpB}), and regularity and decay estimates for elements of $\ker\wh{P_0}(\sigma)$;
\item invertibility and semiclassical estimates for $\wh{P_0}(\sigma)$ at high frequencies (\S\ref{SsSpHi});
\item the Fredholm property of the zero energy operator $\wh{P_0}(0)$ (\S\ref{SsSp0}), and a criterion for it to have index $0$ (Corollary~\ref{CorSpLoInd0});
\item the Fredholm property of the $\tface$- (transition face) normal operators $N_\tface(P_0,e^{i\theta})$ (\S\ref{SsSptf});
\item assuming the invertibility of the $\tface$-normal operators and the zero energy operator: uniform low-energy estimates in the closed upper half plane (\S\ref{SsSpLo}).
\end{enumerate}

We also prove results on the regularity of the resolvent $\wh{P_0}(\sigma)^{-1}$ in $\sigma$ (more precisely, in $\sigma\pa_\sigma$). These results are not enough for the proof of $\aleph$-admissibility for any $\aleph$: one needs, in addition, to require mode stability at $\sigma\neq 0$ and detailed information on the low-energy behavior (including the invertibility of the $\tface$-normal operators). Given this, the estimates in the present section yield the $\aleph$-admissibility of operators $P_0$ of interest in applications rather quickly; see~\S\ref{SA1} for the scalar wave equation, \S\ref{SA2} for the Maxwell equation, and \citestab{\S\ref*{SAdm}} for the linearized gauge-fixed Einstein equation.

\begin{rmk}[Reference point]
\label{RmkSpRef}
  Many of the results stated here exist, at least in some related form, in the literature (see the references in the respective sections). In order to provide a unified point of reference and streamlined proofs for future work on spectral theory on asymptotically flat spacetimes, we present the arguments in some detail (and generality) here, the main focus being on the underlying microlocal propagation estimates. In particular, we give a systematic treatment of (microlocal propagation, radial point, and energy) estimates when $\Im\sigma\geq 0$ (including when $\Im\sigma\gg 1$); and we also give a self-contained proof of an almost sharp semiclassical propagation estimate at normally hyperbolic trapping (Theorem~\ref{ThmSpHiTr}).
\end{rmk}

\begin{notation}[Ps.d.o.\ algebras]
\label{NotSpPsdo}
  The spectral family $\wh{P_0}(\sigma)$ has smooth coefficients, and hence we do not need to prove tame estimates. We shall thus simply write $\Psib$, $\Psisc$, $\Psi_{\scop,\bop}$, for algebras of ps.d.o.s on $\ol{\R^3}$ with b-regular coefficients.
\end{notation}

In order to emphasize the structural properties of $P_0$ relevant for various steps in the analysis of $\wh{P_0}(\sigma)$, we shall repeatedly refer to the following two rather general settings (which stationary wave-type operators in the sense of Definition~\ref{DefSSAdm} fit into):

\begin{definition}[General stationary wave operator]
\label{DefSpGen}
  Let $\cX$ be an open $n$-dimensional manifold ($n\geq 1$). Let $g$ be a smooth stationary time-oriented Lorentzian metric on $\cM=\R_t\times\cX$ such that $\dd t$ is everywhere past timelike. Denote the dual metric function by
  \[
    G(t,x;\sigma,\xi)=g_{(t,x)}^{-1}(-\sigma\,\dd t+\xi\,\dd x,-\sigma\,\dd t+\xi\,\dd x).
  \]
  A \emph{general stationary wave operator} $P_0$ on $\cM$ is then a $t$-translation-invariant operator $P_0\in\Diff^2(\cM)$ with principal symbol $G$.
\end{definition}

\begin{definition}[Stationary asymptotically Minkowski wave operator]
\label{DefSpGenMink}
  Let $\cX$ be an open neighborhood of the boundary $\pa\ol{\R^n}$ of the radial compactification $\ol{\R^n}$ of $\R^n$; write $r=|x|$ and $\rho=r^{-1}$. Then a general stationary wave operator $P_0$ on $\R\times\cX^\circ$, relative to a metric $g$, is a \emph{stationary asymptotically Minkowski wave operator} if $g^{-1}\equiv -2\pa_{t_*}\otimes_s\pa_r + \pa_r^2 + r^{-2}\slg^{-1}$ modulo $\rho\CI(\cX;S^2\,\Tsc\cX)$ (cf.\ \eqref{EqCTebMinkTstar}) and\footnote{These conditions are the same as in~\eqref{EqSSAdmOp}--\eqref{EqSSAdmOpPieces}.}
  \[
    P_0 = -2\pa_{t_*}\rho\Bigl(\rho\pa_\rho-\frac{n-1}{2}-S\Bigr) + \wh{P_0}(0) + Q\pa_{t_*} - g^{0 0}\pa_{t_*}^2,
  \]
  where $g^{0 0}\in\rho^2\CI(\cX)$, $\wh{P_0}(0)\in\rho^2\Diffb^2(\cX)$, $S\in\CI(X)$, $Q\in\rho^3\Diffb^1(\cX)$.
\end{definition}

\textit{We will only consider scalar wave operators; the purely notational modifications required to treat operators acting on vector bundles are left to the reader.}

\subsection{Phase space relationships}
\label{SsSpTs}

The function spaces involved in estimates for the spectral family are those arising from 3b-Sobolev spaces via the Fourier transform and Plancherel's theorem as in Lemma~\ref{LemmaMUetbFT}. Concretely, a stationary-$P_0$-admissible order function
\begin{equation}
\label{EqSpOrderStart}
  \sfs\in\CI({}^\tbop S^*M_0)\ \text{with weights}\ \alpha_+,0
\end{equation}
(Definition~\ref{DefSSOrderAdm}) induces an order function on the spectral side, given at a $\sigma$-level set by $\sfs_\sigma$, via the phase space relationship~\eqref{EqMSFPhaseSpace}. This gives rise to (variable) orders on each of the boundary hypersurfaces of the b- ($\sigma=0$), scattering ($\sigma\neq 0$), sc-b-transition ($\pm\sigma\in[0,1]$), and semiclassical scattering ($\sigma=\pm h^{-1}$) phase spaces at which the corresponding ps.d.o.\ algebras are symbolic; see~\eqref{EqMUscbtTstar} and~\eqref{EqMUschTstar} for the latter two, and also Lemma~\ref{LemmaMUetbFT}. Moreover, also the characteristic set and Hamiltonian vector field of $P_0$ are related to those of $\wh{P_0}(\sigma)$ via the same phase space relationship.

Let us make this explicit so as to illustrate the emergence of the different phase spaces for different $\sigma$. (An abbreviated discussion is given in \cite[\S{2.2}]{HintzNonstat}; see in particular \cite[Figure~2.3]{HintzNonstat} and also \cite[Figure~3.3, Proposition~4.29]{Hintz3b}.) In the coordinates~\eqref{EqTs3bCCoord}, $\sfs$ is a function
\[
  \sfs = \sfs\Bigl(r,\omega; \sigma_\tbop\frac{\dd t_*}{r}+\xi_\tbop\frac{\dd r}{r}+\eta_\tbop\Bigr)
\]
which is homogeneous of degree $0$ in $(\sigma_\tbop,\xi_\tbop,\eta_\tbop)\in(\R\oplus\R\oplus T_\omega^*\Sph^2)\setminus o$. Therefore,
\begin{equation}
\label{EqSpOrderFT}
  \sfs_\sigma\Bigl( r,\omega; \xi\frac{\dd r}{r} + \eta\Bigr) = \sfs\Bigl(r,\omega; -r\sigma\frac{\dd t_*}{r} + \xi\frac{\dd r}{r} + \eta\Bigr),
\end{equation}
where we express spatial covectors in the (b-)form $\xi\frac{\dd r}{r}+\eta$, with $\xi\in\R$, $\eta\in T^*\Sph^2$. We proceed to consider the behavior of $\sfs_\sigma$ in various regimes.

\subsubsection{Nonzero real frequencies}
\label{SssSpTssc}

For $\pm\sigma>0$, we pass to spatial scattering covectors
\[
  \xi_\scop\,\dd r+r\eta_\scop
\]
and use the homogeneity of $\sfs$ to obtain from~\eqref{EqSpOrderFT} the function
\begin{align*}
  &\Tsc^*X \ni \varpi \mapsto \sfs(\Phi_\sigma\varpi)\in\R, \\
  &\qquad \Phi_\sigma \colon (r,\omega;\xi_\scop\,\dd r+r\eta_\scop) \mapsto \Bigl(r,\omega; -\sigma\frac{\dd t_*}{r} + \xi_\scop\frac{\dd r}{r} + \eta_\scop\Bigr);
\end{align*}
here $\Phi_\sigma$ extends to a smooth map $\ol{\Tsc^*}X\to\ol{\Ttb^*_{\cK^+}}M_0$. The function $\sfs\circ\Phi_\sigma$ restricts over $\pa X=r^{-1}(\infty)$, resp.\ at fiber infinity to a scattering decay, resp.\ scattering regularity order
\begin{equation}
\label{EqSpOrdersc}
  \sfr_\sigma \in \CI(\ol{\Tsc^*_{\pa X}}X),\quad\text{resp.}\quad \sfs_\scop \in \CI({}^\scop S^*X);
\end{equation}
note here that $\Phi_\sigma|_{{}^\scop S^*X}$ is independent of $\sigma$.

The principal symbol $G_\sigma$ of
\begin{equation}
\label{EqSpTsscSpecFam}
  \wh{P_0}(\sigma) = 2 i\sigma\rho(\rho\pa_\rho-1-S) + \wh{P_0}(0) - i\sigma Q + g^{0 0}\sigma^2,\quad
  \wh{P_0}(0)=\rho^2 P_{(0)}(\rho,\omega,\rho\pa_\rho,\pa_\omega)
\end{equation}
(cf.\ \eqref{EqSSSpec0Op}), is the preimage under $\Phi_\sigma$ of the 3b-principal symbol of $P_0$. Over $X^\circ$, it consists of all nonzero covectors with vanishing $t_*$-momentum and thus lies over the ergoregion and the black hole horizon and interior; the trapped set is disjoint from the range of $\Phi_\sigma$ (and indeed trapping only plays a role in the high-energy regime). In particular, the order functions $\sfr_\sigma$ and $\sfs_\sigma$ are monotonically decreasing along the $\pm\rho^{-1} H_{G_\sigma}$-flow. Denote the characteristic set by $\Sigma_\sigma\subset\ol{\Tsc^*_{\pa X}}X\cup\,{}^\scop S^*X$. (We will also use the notation $\Sigma_\sigma$ for the fiberwise conic extension of $\Sigma_\sigma$ over $X^\circ$.) Over $\pa X$, we have $G_\sigma=-\sigma^2+(\xi_\scop+\sigma)^2+|\eta_\scop|^2$ (as follows from~\eqref{EqTs3bCCoord}, or directly from~\eqref{EqSSSpecFam}). The radial sets $\cR^\pm_{\pa\cK^+,{\rm in/out}}$ (with $\cR^-_{\pa\cK^+,{\rm out/in}}$ defined as in Definition~\ref{DefTs3bRad} but with $\sigma_\tbop>0$) lie in the image of $\Phi_\sigma$ for $\pm\sigma>0$; their preimages, in the coordinates $\rho=r^{-1}$, $\omega\in\Sph^2$ near $\pa X$, lie over $\pa X$ and are given by the subsets
\begin{equation}
\label{EqSpOrderRinout}
\begin{alignedat}{2}
  \cR_{\sigma,\rm in} &= \{ (\rho,\omega;\xi_\scop,\eta_\scop) \colon \rho=0,\ \xi_\scop=-2\sigma,&\ &\eta_\scop=0 \}, \\
  \cR_{\rm out} &= \{ (\rho,\omega;\xi_\scop,\eta_\scop) \colon \rho=0,\ \xi_\scop=0,&\ &\eta_\scop=0 \}
\end{alignedat}
\end{equation}
of $\Sigma_\sigma$. The null-bicharacteristic flow of $\pm\rho^{-1}H_{G_\sigma}$ in $\Sigma_\sigma$ has a source ($\cR_{\sigma,{\rm in}})$ to sink ($\cR_{\sigma,{\rm out}}$) structure, being the pullback along $\Phi_\sigma$ of the $H_{G_{\bhm,a}}$-flow described in Lemma~\ref{LemmaTs3bDyn}.

The preimages of the radial sets $\cR^\pm_{\cH^+}$ over the event horizon (Definition~\ref{DefTs3bH}, with $\cR^-_{\cH^+}$ being given by the same expression but with $\xi_0<0$) are $\sigma$-independent. We denote their conic extensions in $T^*X^\circ\setminus o$ (by an abuse of notation) by
\[
  \cR^\pm_{\cH^+} = \{ (r,\omega;\xi_0,\eta_0) \colon r=r_+,\ \pm\xi_0>0,\ \eta_0=0 \},
\]
where we write covectors as $\xi_0\,\dd r+\eta_0$, $\eta_0\in T^*_\omega\Sph^2$; and $\cR^\pm_{\cH^+}\subset\Sigma_\sigma^\pm$ (the preimage of the future/past characteristic set $\Sigma^\pm$ under $\Phi_\sigma$). The dynamics of the future null-bicharacteristic flow (i.e., the flow of $\pm H_{G_\sigma}$) can be read off from that of $P_0$, restricted to null-geodesics with $t_*$-momentum $-\sigma$. Besides the source-to-sink flow over $\pa X$ discussed above, there is non-trivial dynamics near the black hole: all null-bicharacteristics that are not contained in $\cR^\pm_{\cH^+}$ escape through $r^{-1}(\bhm)$ in the future direction and tend to $\cR^\pm_{\cH^+}$ in the past direction.

We will then be able to prove Fredholm estimates for $\wh{P_0}(\sigma)$ in~\S\ref{SsSpB} by concatenating radial point and real principal type (and microlocal elliptic) estimates, much as (but simpler than in)~\S\ref{SssR3Pf}.

We will also study complex $\sigma$ with $\Im\sigma>0$. (This is related to studying $P_0$ on spaces with exponentially growing weights in $t_*$, and thus no longer has any direct relationship with the microlocal analysis in~\S\ref{SR}.) The map $\Phi_\sigma$ no longer maps $\Tsc^*_{\pa X}X$ into the (real) phase space $\Ttb^*M_0$. The characteristic set over $X^\circ$ (which thus only lives at fiber infinity) is independent of $\sigma$. Over $\pa X$, lying in the zero set of $G_\sigma=-2\sigma\xi_\scop+\xi_\scop^2+|\eta_\scop|^2$ requires $\xi_\scop=0$ and thus $\eta_\scop=0$; therefore,
\begin{equation}
\label{EqSpOrderCx}
  \Sigma_\sigma \cap \Tsc^*_{\pa X}X = \cR_{\rm out},\quad \Im\sigma>0.
\end{equation}
These considerations remain valid for every stationary asymptotically Minkowski wave operator (Definition~\ref{DefSpGenMink}).

\subsubsection{Large frequencies}
\label{SssSpSemi}

There are two high-energy regimes we must study. The first regime concerns real $\sigma$ with $|\Re\sigma|\to\infty$, or more generally $\sigma$ with imaginary part $\Im\sigma\geq 0$ bounded from above and $|\Re\sigma|\to\infty$; this is relevant for control of the resolvent on the real axis, or more generally in strips above the real axis (the latter being important for contour shifting type arguments). The second regime concerns $\sigma$ with $\Im\sigma\gg 1$ (which is related to the absence of modes with strong exponential growth). Notice, though, that for any fixed $\Im\sigma=C\gg 1$, the rescaled quantity $z=|\sigma|^{-1}\sigma$ converges to $\pm 1$ when $\Re\sigma\to\pm\infty$, though it always does satisfy $\Im z\geq C h$ (up to $\cO(h^2)$ terms) in this limit. Altogether, then, we need to study $\sigma=h^{-1}z$ where $\Im z\geq 0$ and $|z|=1$ (or, slightly more generally, $|z|\in[\frac12,2]$ since for large $\pm\Re\sigma$ it is occasionally convenient to be able to take $h=|\Re\sigma|^{-1}$ and thus $z=\pm 1+i h\Im\sigma$). We shall then study the semiclassical rescaling
\begin{equation}
\label{EqSpPhz}
  P_{h,z} := h^2\wh{P_0}(h^{-1}z),\quad (P_{h,z})_{h\in(0,1)} \in \Diff_{\scop,\semi}^2(X)
\end{equation}
(see~\eqref{EqMUDiffsch}), for such $z$ by separating three cases (with the second regime above split into two cases). Before describing them, recall that the phase space for microlocal analysis in the semiclassical scattering setting is all of $\ol{\Tsc^*}X$, and thus in particular includes finite covectors over $X^\circ$. Put $\Re P_{h,z}=\frac12(P_{h,z}+P_{h,z}^*)$ and $\Im P_{h,z}=\frac{1}{2 i}(P_{h,z}-P_{h,z}^*)$.
\begin{enumerate}
\item When $|\Im z|=\cO(h)$ and $z=\pm 1+\cO(h)$, the principal symbol of $P_{h,z}$ is real and equal to $G_{\pm 1}$ (see below); the imaginary part of $P_{h,z}$ plays only a minor role in shifting threshold quantities in radial point estimates (consistently with the bounded-$\sigma$ regime in Proposition~\ref{PropSpBHor} below).
\item When we require $h\ll\Im z\ll 1$ and $\Re z=\pm 1$, we write
  \begin{equation}
  \label{EqSphDec}
    P_{h,z} = \Re P_{h,z} + i \Im P_{h,z}.
  \end{equation}
  The term $\Im P_{h,z}$ is not uniformly subprincipal in the semiclassical sense (as it has size $\Im z$), but it has a sign compatible with future-directed propagation of regularity. We still regard $\Re P_{h,z}$ as the main term (e.g., determining the semiclassical characteristic set) and $\Im P_{h,z}$ as a complex absorbing potential.\footnote{Recall here from \cite[\S{5.4.5}]{VasyMinicourse} (see also \cite{NonnenmacherZworskiQuantumDecay,WunschZworskiNormHypResolvent}) that for $P-i Q$ with $P,Q\in\Psi_\cl^m$ having real-valued principal symbols, one can propagate regularity of solutions $u$ of $(P-i Q)u=f$ in the forward direction along $H_{\upsigma^m(P)}$ on $\Char(P)$ when $\upsigma^m(Q)\geq 0$, and in the backward direction when $\upsigma^m(Q)\leq 0$. Thus, the sign of the imaginary part forces one direction of propagation. In the present context, for $\Im z\geq 0$, this direction is the \emph{future} direction (i.e., along the Hamiltonian vector field in the future characteristic set, and along the negative of the Hamiltonian vector field in the past characteristic set), both for $\Re z=+1$ and $\Re z=-1$. The reader may be used to being forced to propagate along the Hamiltonian flow everywhere. But when studying, say, the spectral family $\Delta-\sigma^2$ of the Minkowski wave operator, the semiclassical characteristic set for $\sigma=h^{-1}z$, $z=\pm 1+i\eps h$, $\eps\geq 0$, only has one component, namely the future (for $z=1+\cO(h)$) or past (for $z=-1+\cO(h)$). It is only for spectral families of wave operators associated with spacetime metrics featuring an ergoregion (where $\pa_t$ ceases to be timelike) that the semiclassical characteristic set has \emph{two} components.} When $h^{-1}\Im z$ is large enough, then the term $\Im P_{h,z}$ dominates in commutator arguments.
\item When $\Im z\geq c>0$ (in which case $\Re z$ need not have a sign anymore), we can exploit the definite sign of $\Im z$ to prove the ellipticity of $P_{h,z}$ at finite frequencies (except over $\pa X$). The characteristic set at fiber infinity is independent of the spectral parameter. We again use the decomposition~\eqref{EqSphDec}, with $\Im P_{h,z}$ dominating in positive commutator arguments.
\end{enumerate}

For $z=\pm 1+\cO(h)$, we pass to spatial semiclassical scattering covectors
\[
  \xi_\schop\frac{\dd r}{h} + \frac{r\eta_\schop}{h}
\]
and obtain from~\eqref{EqSpOrderFT} for $\sigma=h^{-1}z=\pm h^{-1}+\cO(1)$ the function $\Tsc^*X\ni\varpi\mapsto\sfs(\Phi_{\pm 1}\varpi)=\sfs(r,\omega;\mp\frac{\dd t_*}{r}+\xi_\schop\frac{\dd r}{r}+\eta_\schop)$. This extends to the full compactified semiclassical scattering phase space $\ol{{}^\schop T^*}X$ (see the discussion after~\eqref{EqMUschTstar0}) and induces (by restriction to the boundary hypersurfaces of this phase space listed in~\eqref{EqMUschTstar0}, intersected with $h=0$) a (sign-independent) semiclassical scattering regularity order
\begin{subequations}
\begin{equation}
\label{EqSpOrderHis}
  \sfs_\scop \in \CI({}^\scop S^*X)
\end{equation}
(where we identify fiber infinity of $\ol{{}^\schop T^*}X$ over $h=0$ with fiber infinity ${}^\scop S^*X$ of $\ol{{}^\scop T^*}X$), a (sign-dependent) semiclassical scattering decay order
\begin{equation}
\label{EqSpOrderHir}
  \sfr_\pm \in \CI(\ol{{}^\scop T^*_{\pa X}}X),
\end{equation}
and a (sign-dependent) semiclassical order (powers of $h$)
\begin{equation}
\label{EqSpOrderHib}
  \sfb_\pm \in \CI(\ol{{}^\scop T^*}X).
\end{equation}
\end{subequations}
In fact, we have $\sfr_\pm=\sfr_{\pm 1}$ in the notation of~\eqref{EqSpOrdersc}, and the functions $\sfs_\scop$ here and in~\eqref{EqSpOrdersc} are identical.

The semiclassical characteristic set is given by
\begin{equation}
\label{EqSpOrderSigmah}
  \Sigma_{\pm,\semi} := \ol{\Tsc^*}X \cap \{ \rho_\infty^2 G_{\pm 1}=0 \},
\end{equation}
where $\rho_\infty\in S^{-1}(\Tsc^*X)$ is elliptic and positive. The rescaled Hamiltonian vector field $\rho^{-1}\rho_\infty H_{G_{\pm 1}}$ is a smooth b-vector field on $\ol{\Tsc^*}X$, and we refer to its integral curves inside of $\Sigma_{\pm,\semi}$ as null-bicharacteristics.

The future-directed null-bicharacteristic flow over $\pa X$ has a source-to-sink structure with radial sets $\cR_{\pm 1,{\rm in}}$ and $\cR_{\rm out}$ in the notation of~\eqref{EqSpOrderRinout}. The preimage of $\cR^+_{\cH^+}$ (which lies in the annihilator of $\pa_{t_*}$) under $\Phi_{\pm 1}$ lies at fiber infinity only, i.e., it is given by
\[
  \pa\cR^+_{\cH^+} \subset S^*X^\circ,
\]
similarly for the past component $\pa\cR^-_{\cH^+}$.

Finally, consider the trapped set. The intersection of the time-translation-invariant trapped set $\Gamma_0$ (Definition~\ref{DefTs3bOTrap0}) with the level sets $\mp\dd t_*+T^*X^\circ$ of the $t_*$-momentum coordinate (which is the same as the Boyer--Lindquist $\ft$-momentum) over $t_*^{-1}(0)$, say, is compact by \eqref{EqTs3bOTrap0Cpt} and~\eqref{EqTs3bSumEll}, and therefore the preimage of the trapped set $\Gamma\subset\Ttb^*_{\cK^+}M_0\setminus o$ inside of the future characteristic set $\Sigma^+$ of $P_0$ (see~\eqref{EqTs3bGamma}) under $\Phi_{+1}$ is a compact subset which by an abuse of notation we denote by
\[
  \Gamma \subset T^*X^\circ,
\]
while the preimage under $\Phi_{-1}$ is empty (cf.\ the sign in~\eqref{EqTs3bSumEll}). Analogous considerations apply to the forward/backward trapped sets $\bar\Gamma^{\rm u/s}$ (see also Definition~\ref{DefTs3bOTrapus0}) inside of $\Sigma^+$, whose preimage under $\Phi_{+1}$ are smooth codimension $2$ submanifolds
\[
  \Gamma^{\rm u/s} \subset T^*X^\circ
\]
intersecting transversally at $\Gamma$; restriction of the defining functions from Definition~\ref{DefTs3bODefFn} yields defining functions $\phi^{\rm u/s}$ of $\Gamma^{\rm u/s}$ inside of the characteristic set of $(P_{h,+1})_{h\in(0,1)}$ such that
\[
  H_{G_{+1}}\phi^{\rm u/s} = \mp\nu^{\rm u/s}\phi^{\rm u/s},\quad
  \{\phi^{\rm u},\phi^{\rm s}\}\neq 0
\]
near $\Gamma$. For the (forward/backward) trapped sets inside of the past characteristic set of $P_0$, analogous statements hold with the roles of $\Phi_{+1}$ and $\Phi_{-1}$ reversed, and with $-H_{G_{-1}}$ in place of $H_{G_{+1}}$. (See also \cite[\S{2.3}]{DyatlovWaveAsymptotics} for the relationship between the spacetime and spatial phase spaces near the trapped set.)

When $\Im z$ is not of size $\cO(h)$, the above phase space relationship between $\Tsc^*X$ and $\Ttb^*_{\cK^+}M_0$ breaks down. For $\Im z\gtrsim 1$, the characteristic set over $\pa X$ is given by $\cR_{\rm out}$, exactly as in~\eqref{EqSpOrderCx}. The characteristic set over $X^\circ$ is described in a general setting in Lemma~\ref{LemmaSpChar} below.

\subsubsection{Zero energy}

The relevant phase space at zero energy is the b-cotangent bundle, as discussed in~\S\ref{SssMUK}. Concretely, the order induced by the stationary-$P_0$-admissible order function $\sfs\in\CI({}^\tbop S^*M_0)$ in~\eqref{EqSpOrderFT} is
\begin{equation}
\label{EqSpOrderb}
  \sfs_0=\sfs\circ\Phi_0,\qquad
  \Phi_0 \colon \Bigl(r,\omega; \xi_\bop\frac{\dd r}{r}+\eta \Bigr) \ni \Tb^*X \mapsto \Bigl(r,\omega; \xi_\bop\frac{\dd r}{r}+\eta \Bigr) \in {}^\tbop T^*_{\cK^+}M_0,
\end{equation}
which extends to a smooth map $\ol{\Tb^*}X\to\ol{{}^\tbop T^*_{\cK^+}}M_0$. Since the microlocalization locus of the b-algebra is fiber infinity $\Sb^*X$ (see~\S\ref{SssMSCp}), the characteristic set of $\wh{P_0}(0)$ is the preimage under $\Phi_0$ of the intersection of the characteristic set $\pa\Sigma\subset{}^\tbop S^*M_0$ of $P_0$ with the closure of $\ann\pa_t\subset T^*M_0^\circ$; this is thus non-empty only over the ergoregion, and the only critical set for the null-bicharacteristic flow is $\pa\cR_{\cH^+}$ (Definition~\ref{DefTs3bH}). In particular, $\wh{P_0}(0)\in\rho^2\Diffb^2(X)$ is elliptic (as a weighted b-operator) near $\pa X$.\footnote{As discussed in \cite{MelroseSaBarretoLow,GuillarmouHassellResI,VasyLowEnergy,VasyLowEnergyLag}, the scattering perspective on $\wh{P_0}(0)$, as obtained from~\S\ref{SssSpTssc} by taking $\sigma=0$, leads to a quadratic degeneracy at the scattering zero section over $\pa X$, which is resolved by working in the b-setting instead.}

\subsubsection{Transition from zero to nonzero energies}
\label{SssSpOrderscbt}

In the low-energy $\pm\sigma\in[0,1]$, the order function~\eqref{EqSpOrderFT} (with $\sigma$ being a coordinate on the ``total'' phase space which includes a factor of $\R_\sigma$) gives rise to sc-b-transition regularity and scattering decay orders
\begin{equation}
\label{EqSpOrderscbt}
  \sfs \in \CI({}^\scbtop S^*X),\quad
  \sfr_\pm \in \CI(\ol{{}^\scbtop T^*_\scface}X),
\end{equation}
in the notation of~\eqref{EqMUscbtTstar}. (The order $\sfs$ is independent of the choice of sign; see below.) These orders, in turn, can be restricted to ${}^\scbtop S^*_\tface X={}^{\scop,\bop}S^*\tface$ and $\ol{{}^\scbtop T^*_{\tface\cap\sctface}}X=\ol{{}^{\scop,\bop}T^*_\scface}\tface$; this produces scattering-b-regularity and scattering decay orders which (by a minor abuse of notation) we also denote by
\begin{equation}
\label{EqSpOrdertf}
  \sfs \in \CI({}^{\scop,\bop}S^*\tface),\quad \sfr_\pm \in \CI(\ol{{}^\scbtop T^*_\sctface}\tface).
\end{equation}

Let us make the latter construction explicit. For $\pm\sigma\in[0,1]$, we use $\hat\rho=\frac{1}{\pm r\sigma}=\frac{\rho}{|\sigma|}$ as the projective coordinate on $\tface=[0,\infty]_{\hat\rho}\times\Sph^2$ as in~\eqref{EqMUtf}. A frame of ${}^{\scop,\bop}T^*\tface$ is then given by the covectors
\begin{equation}
\label{EqSpOrderscbtFrame}
  (\hat\rho^{-1}+1)\frac{\dd\hat\rho}{\hat\rho},\ (\hat\rho^{-1}+1)T^*\Sph^2
\end{equation}
(which interpolate between the scattering 1-forms $\frac{\dd\hat\rho}{\hat\rho^2}$, $\hat\rho^{-1}T^*\Sph^2$ near $\sctface=\{\hat\rho=0\}$ and b-1-forms $\frac{\dd\hat\rho}{\hat\rho}$, $T^*\Sph^2$ near $\ztface=\{\hat\rho=\infty\}$). Therefore, exploiting the homogeneity of $\sfs$, we obtain from~\eqref{EqSpOrderFT} the functions
\begin{align*}
  &{}^{\scop,\bop}T^*\tface \ni \varpi \mapsto \sfs(\Phi_\pm\varpi) \in \R, \\
  &\hspace{4em} \Phi_\pm \colon \Bigl(\hat\rho,\omega; \xi_{\scop,\bop}(\hat\rho^{-1}+1)\frac{\dd\hat\rho}{\hat\rho} + (\hat\rho^{-1}+1)\eta_{\scop,\bop}\Bigr) \\
  &\hspace{4em}\hspace{4em} \mapsto \Bigl(\infty,\omega;\mp\frac{\dd t_*}{r}+(1+\hat\rho)\xi_{\scop,\bop}\frac{\dd r}{r} + (1+\hat\rho)\eta_{\scop,\bop}\Bigr) \in \Ttb^*_{\pa\cK^+}M_0.
\end{align*}
The restriction of this function to $\hat\rho=0$, resp.\ to fiber infinity (in $(\xi_{\scop,\bop},\eta_{\scop,\bop})$) yields the orders~\eqref{EqSpOrdertf}. Note that the restriction to fiber infinity does not depend on the choice of sign.

The map $\Phi_\pm$ also relates characteristic sets and Hamiltonian vector fields for $P_0$ over $\pa\cK^+$ with those of $N_\tface(P_0,\pm 1)$. The characteristic set of $N_\tface(P_0,\pm 1)$ lies over $\sctface$ and is given by
\[
  \Sigma_\tface = {}^{\scop,\bop}T^*_\sctface\tface \cap G_\tface^{-1}(0),\quad G_\tface=(\xi_{\scop,\bop}\mp 1)^2 + |\eta_{\scop,\bop}|^2 - 1.
\]
(This also follows from the fact that the principal symbol of~\eqref{EqSStfOp} over $\hat\rho=0$ is given by that of $\hat\rho^2 P_{(0)}(0,\omega,\hat\rho\pa_{\hat\rho},\pa_\omega)-2\hat\sigma\hat\rho^2 D_{\hat\rho}$ with $\hat\sigma=\pm 1$, which in turn is given by the dual metric function of the Minkowski metric restricted to the annihilator of $\pa_{t_*}$; see then~\eqref{EqCTebMinkTstar}.) Fiber infinity of ${}^{\scop,\bop}T^*\tface$, by contrast, is mapped to $\ol{\{\sigma_\tbop=0\}}\subset\ol{\Ttb^*_{\pa\cK^+}}M_0$ by $\Phi_\pm$, which is disjoint from the characteristic set of $P_0$. The radial sets $\cR^\pm_{\pa\cK^+,{\rm out/in}}$  lie in the image of $\Phi_\pm$ but not of $\Phi_\mp$; their preimages are
\begin{equation}
\label{EqSpOrderRad}
\begin{alignedat}{2}
  \cR_{\pm 1,\rm in} &= \{ (\hat\rho,\omega;\xi_{\scop,\bop},\eta_{\scop,\bop}) \colon \hat\rho=0,\ \xi_{\scop,\bop}=\pm 2,&\ &\eta_{\scop,\bop}=0 \}, \\
  \cR_{\rm out} &= \{ (\hat\rho,\omega;\xi_{\scop,\bop},\eta_{\scop,\bop}) \colon \hat\rho=0,\ \xi_{\scop,\bop}=0,&\ &\eta_{\scop,\bop}=0 \}.
\end{alignedat}
\end{equation}

Note that~\eqref{EqSpOrderscbtFrame} is a local frame of ${}^\scbtop T^*X$ near $\tface$ (since $(\hat\rho^{-1}+1)^{-1}=\frac{\rho}{\rho+|\sigma|}$ is a defining function of $\scface\subset X_\scbtop$, cf.\ the discussion around~\eqref{EqMUscbtSingle}). Thus, the same maps $\Phi_\pm$ also serve to map $\ol{{}^\scbtop T^*}X$ near $\tface$ to $\ol{\Ttb^*_{\cK^+}}M_0$. Therefore, in ${}^\scbtop T^*X$, the sets~\eqref{EqSpOrderRad} are radial sets (and lie over $\scface$) for the null-bicharacteristic flow of $(\wh{P_0}(\sigma))_{\sigma\in\pm[0,1]}\in\rho_\tface^2\Diff_\scbtop^2(X)$ (cf.\ \eqref{EqSSSpecFamMem}), which has the usual source-to-sink structure over $\scface$; and this operator is elliptic near $\tface\setminus\scface$. Its characteristic set away from $\scface$ lies only over the ergoregion and the black hole interior, as in the bounded $\sigma$ setting discussed in~\S\ref{SssSpTssc}.

\subsection{Bounded nonzero frequencies}
\label{SsSpB}

For $\sigma\in\R\setminus\{0\}$, we recall the orders $\sfr_\sigma$ and $\sfs_\scop$ from~\eqref{EqSpOrdersc}, which are induced by a stationary-$P_0$-admissible order function $\sfs\in\CI({}^\tbop S^*M_0)$ with weights $\alpha_+,0$; here $\alpha_+\in\R$. By the discussion in~\S\ref{SssSpTssc}, the orders satisfy the threshold conditions
\begin{subequations}
\begin{alignat}{2}
\label{EqSpBThrHor}
  &\text{at}\ \cR^\pm_{\cH^+}:\ & \sfs_\scop &> \frac12+\vartheta_{\cH^+}, \\
\label{EqSpBThrIn}
  &\text{at}\ \cR_{\sigma,{\rm in}}:\ & \sfr_\sigma+\alpha_+ &> \frac12(-1+\vartheta_{\pa\cK^+,{\rm in}}), \\
\label{EqSpBThrOut}
  &\text{at}\ \cR_{\rm out}:\ & \sfr_\sigma+\alpha_+ &< -\frac12+\ubar S,
\end{alignat}
\end{subequations}
where we recall $\vartheta_{\cH^+}$ from~\eqref{EqSSOrderAdm4}.

\begin{thm}[Spectral family at nonzero frequencies]
\fakephantomsection
\label{ThmSpB}
  \begin{enumerate}
  \item\label{ItSpBReal}{\rm (Fredholm estimates, I.)} Let $I\subset\pm(0,\infty)$ be a compact set, and let $\nu\geq 0$, $N\in\R$. Then there exists a constant $C$ such that for all $\sigma=\sigma_0+i\sigma_1$ with $\sigma_0\in I$ and $\sigma_1\in[0,\nu]$, and writing $\sfr=\sfr_{\sigma_0}$, we have the estimates
    \begin{equation}
    \label{EqSpBEst}
    \begin{split}
      \|u\|_{\bar H_\scop^{\sfs_\scop,\sfr+\alpha_+}(X;\cE_X)} &\leq C\Bigl( \|\wh{P_0}(\sigma)u\|_{\bar H_\scop^{\sfs_\scop-1,\sfr+\alpha_++1}(X;\cE_X)} + \|u\|_{\bar H_\scop^{-N,-N}(X;\cE_X)} \Bigr), \\
      \|u^*\|_{\dot H_\scop^{-\sfs_\scop+1,-\sfr-\alpha_+-1}(X;\cE_X)} &\leq C\Bigl( \|\wh{P_0}(\sigma)^*u^*\|_{\dot H_\scop^{-\sfs_\scop,-\sfr-\alpha_+}(X;\cE_X)} + \|u^*\|_{\dot H_\scop^{-N,-N}(X;\cE_X)} \Bigr).
    \end{split}
    \end{equation}
    The extendible, resp.\ supported character of these spaces refers to the boundary hypersurface $\pa_\sharp X=r^{-1}(\bhm)$ of $X$,\footnote{These spaces arise from 3b-Sobolev spaces on $M_0$ with extendible, resp.\ supported character at $r^{-1}(\bhm)$ via the Fourier transform, cf.\ Proposition~\ref{PropMUSuppFT}.} and we use the Euclidean volume density on $X$ for the underlying $L^2$-spaces.
  \item\label{ItSpBImag}{\rm (Fredholm estimates, II.)} Let $I\subset\R$ be a compact set and $0<\nu_\flat<\nu$. Let $\sfr\in\CI(\ol{\Tsc^*_{\pa X}}X)$ satisfy~\eqref{EqSpBThrOut}. Then there exists a constant $C$ such that the estimates~\eqref{EqSpBEst} hold for all $\sigma\in I+i[\nu_\flat,\nu]$.
  \item\label{ItSpBFred} {\rm (Fredholm index.)} In the settings of~\eqref{ItSpBReal} and \eqref{ItSpBImag}, the Fredholm index of
    \begin{equation}
    \label{EqSpBMap}
      \wh{P_0}(\sigma) \colon \cX_\sigma := \bigl\{ u\in\bar H_\scop^{\sfs_\scop,\sfr+\alpha_+} \colon \wh{P_0}(\sigma)u\in\bar H_\scop^{\sfs_\scop-1,\sfr+\alpha_++1} \bigr\} \to \cY_\sigma := \bar H_\scop^{\sfs_\scop-1,\sfr+\alpha_++1}
    \end{equation}
    is independent of $\sigma$ (and, in fact, equal to\footnote{This will be proved in~\S\ref{SsSpHi} below.} $0$).
  \item\label{ItSpBb} {\rm (Higher b-regularity.)} If $u\in\cX_\sigma$ and $\wh{P_0}(\sigma)u\in\bar H_{\scop;\bop}^{(\sfs_\scop-1;k),\sfr+\alpha_++1}$ for some $k\in\N_0$,\footnote{See~\eqref{EqMUHscb} for the definition of this function space.} then also $u\in\bar H_{\scop;\bop}^{(\sfs_\scop;k),\sfr+\alpha_+}$, and we have an estimate
    \begin{equation}
    \label{EqSpBbEst}
      \|u\|_{\bar H_{\scop;\bop}^{(\sfs_\scop;k),\sfr+\alpha_+}} \leq C\Bigl( \|\wh{P_0}(\sigma)u\|_{\bar H_{\scop;\bop}^{(\sfs_\scop-1;k),\sfr+\alpha_++1}} + \|u\|_{\bar H_\scop^{\sfs_\scop,\sfr+\alpha_+}}\Bigr).
    \end{equation}
  \item\label{ItSpBNull}{\rm (Regularity of kernel elements.)} In the settings of~\eqref{ItSpBReal} and \eqref{ItSpBImag}, every $u\in\bar H_\scop^{\sfs_\scop,\sfr+\alpha_+}(X;\cE_X)$ with $\wh{P_0}(\sigma)u=0$ satisfies
    \begin{equation}
    \label{EqSpBNullA}
      u\in\cA^{1+\ubar S-\eps}(X;\cE_X)=\rho^{1+\ubar S-\eps}\CI_\bop(X;\cE_X)\quad\forall\,\eps>0.
    \end{equation}
    (Conversely, the space $\bigcap_{\eps>0}\cA^{1+\ubar S-\eps}(X;\cE_X)$ is contained in the domain of~\eqref{EqSpBMap}.) If for all $p\in\pa X$ and $\lambda,\mu\in\spec S(p)$, we have $\lambda-\mu\notin\Z$,\footnote{When this condition is violated, but $S|_{\pa X}$ has constant spectrum (which is the case in all applications we shall consider), then~\eqref{EqSpBNullspace} needs to be amended by the possibility of logarithmic lower-order terms.} then, in fact,
    \begin{equation}
    \label{EqSpBNullspace}
      u\in\rho^{1+S}\CI(X;\cE_X),
    \end{equation}
  \item\label{ItSpBInvReg}{\rm (Regularity in the spectral parameter.)} If the map~\eqref{EqSpBMap} is invertible for $\sigma=\sigma_0\in\C$, then the estimates~\eqref{EqSpBEst} and~\eqref{EqSpBbEst} hold without the second terms on the right. Moreover, for all $j\in\N_0$, $\wh{P_0}(\sigma)^{-1}$ is $j$ times continuously differentiable in $\sigma$ as a map
    \begin{equation}
    \label{EqSpBInvReg0}
      \wh{P_0}(\sigma)^{-1} \colon \bar H_{\scop;\bop}^{(\sfs_\scop+\eps-1;j),\sfr+\eps+\alpha_++1}\to\bar H_{\scop;\bop}^{(\sfs_\scop-\eps;0);\sfr-\eps+\alpha_+},
    \end{equation}
    where the orders are those for $\sigma=\sigma_0$, for any fixed $\eps>0$ and for all $\sigma$ in a neighborhood (depending only on $\eps$) of $\sigma_0$. In particular, when $\Im\sigma_0>0$, then $\wh{P_0}(\sigma)^{-1}$ is holomorphic near $\sigma=\sigma_0$. The $j$-th derivative is moreover bounded as a map\footnote{The scattering regularity gain of $j$ orders in the output of~\eqref{EqSpBInvReg} will not be used in this paper.}
    \begin{equation}
    \label{EqSpBInvReg}
      \pa_\sigma^j\wh{P_0}(\sigma_0)^{-1} \colon \bar H_{\scop;\bop}^{(\sfs_\scop-1;k+j),\sfr+\alpha_++1} \to \bar H_{\scop;\bop}^{(\sfs_\scop+j;k),\sfr+\alpha_+}\quad\forall\,k\in\N_0.
    \end{equation}
  \end{enumerate}
\end{thm}

\begin{rmk}[Strong estimates]
\label{RmkSpBStrong}
  If one replaces the error term in the first estimate in~\eqref{EqSpBEst} by the norm of $\bar H_\scop^{\sfs_\scop^\flat-1,\sfr^\flat+\alpha_++1}$ where $\sfs_\scop^\flat,\sfr^\flat$ are the orders induced by an admissible order function $\sfs^\flat<\sfs$, then the resulting estimate holds in the strong sense that if the right-hand side is finite, then so is the left-hand side and the estimate is valid. (An analogous statement holds for the adjoint estimate.) This follows from the fact that each microlocal result used in the proof of~\eqref{EqSpBEst} holds in such a strong sense.
\end{rmk}

\begin{rmk}[Related literature]
\fakephantomsection
\label{RmkSpBLit}
  \begin{enumerate}[label=(\roman*)]
  \item Parts~\eqref{ItSpBReal}--\eqref{ItSpBImag} follow by concatenating radial point, real principal type (with mild complex absorption when $\sigma_1>0$), and microlocal elliptic estimates. At $\cR^\pm_{\cH^+}$, the relevant radial point estimates were first proved in Vasy \cite[\S{2.4}]{VasyMicroKerrdS}. At the radial sets over $\pa X$, the relevant estimates for the ``standard'' spectral family $\Delta-\sigma^2$ of asymptotically Euclidean spaces are due to Melrose \cite{MelroseEuclideanSpectralTheory}; see also \cite[\S{5.4.7}]{VasyMinicourse}. We are, in effect, considering a conjugation (roughly $e^{-i\sigma r}(\Delta-\sigma^2)e^{i\sigma r}$) of this ``standard'' spectral family by weights which are exponential in $r$ when $\Im\sigma>0$, and thus the radial point estimates for complex $\sigma$ are more delicate. (For example, $\Delta-\sigma^2$ is scattering elliptic for $\Im\sigma>0$, but the conjugated spectral family always has characteristic set at the scattering zero section over $\pa X$.) Vasy \cite{VasyLAPLag} introduced a second microlocal framework for this purpose, but we work solely with scattering spaces here, following \cite[\S{4.1.2}]{HintzNonstat}.
  \item The essence of part~\eqref{ItSpBFred} is already contained in \cite[\S{2.7}]{VasyMicroKerrdS}. A proof in a special case is discussed in the proof of \cite[Theorem~6.1]{HaefnerHintzVasyKerr}.
  \item The b-regularity in part~\eqref{ItSpBb} appears under the name of \emph{module regularity} as introduced in \cite{HassellMelroseVasySymbolicOrderZero} (and also in an \textit{ad hoc} manner in \cite{MelroseEuclideanSpectralTheory}) and used frequently since, e.g., in \cite[\S{4}]{BaskinVasyWunschRadMink} and \cite{GellRedmanHassellShapiroZhangHelmholtz}; see Remark~\ref{RmkSpBInftyModule}.
  \item Extracting asymptotics~\eqref{EqSpBNullspace} from the b- (i.e., conormal) regularity asserted in part~\eqref{ItSpBb} is essentially a standard matter and already discussed in \cite{MelroseEuclideanSpectralTheory}. For the present ``conjugated'' spectral family (which factors out the oscillatory $e^{i\sigma r}$ term of outgoing mode solutions), this is particularly simple; the details are given below.
  \end{enumerate}
\end{rmk}

\begin{rmk}[No pure b-estimates]
\label{RmkSpNob}
  Unlike in \cite[\S{4}]{HintzNonstat}, we do not prove precise estimates on b-Sobolev spaces here, but rather stay within the framework of (variable order) scattering (and, in later subsections, b-, scattering-b-transition, and semiclassical scattering Sobolev spaces) with integer amounts of additional b-regularity, thus streamlining the analysis and making it directly compatible with the microlocal e3b-theory with variable regularity order (cf.\ Theorem~\ref{ThmEReg}). In particular, we do not rely on \cite{VasyLAPLag} anymore, but rather work entirely using the analytic framework introduced in \cite{HintzKdSMS}. One benefit of direct b-estimates in the spirit of \cite{VasyLAPLag} is that they imply estimates for forward solutions that are less lossy as far as the b-regularity order is concerned. But since in the present paper we are working in a b-tame setting and allow ourselves arbitrary, albeit fixed, amounts of b-regularity losses, this is of no concern to us; we will discuss (rather lossy) pure b-estimates in~\S\ref{SsDRes}.
\end{rmk}

We discuss the radial point estimates underlying parts~\eqref{ItSpBReal}, \eqref{ItSpBImag}, and \eqref{ItSpBb} for the spectral family of rather general wave-type operators. Estimates near horizons are discussed in~\S\ref{SssSpBHor}, and those at spatial infinity in~\S\ref{SssSpBInfty}. The proof of Theorem~\ref{ThmSpB} is then given in~\S\ref{SssSpBPf}.

\subsubsection{Radial point estimates near horizons}
\label{SssSpBHor}

We consider a general stationary wave operator $P_0$ on $\cM=\R_t\times\cX$ as in Definition~\ref{DefSpGen} with the following additional structure:
\begin{enumerate}
\item\label{ItSpBHor1} $\cH\subset\cX$ is a smooth compact submanifold of codimension $1$ with orientable normal bundle;
\item $\cR^+$ is one of the two connected components of $N^*\cH\setminus o$, and $\cR^+\subset\Sigma_0:=\{(x,\xi)\in T^*\cX\setminus o\colon G_0(x,\xi)=0\}$ where $G_0(x,\xi)=G(0,x;0,\xi)$;
\item $\rho_\infty\in S^{-1}_{\rm hom}(T^*\cX\setminus o)$ is positive, homogeneous of degree $-1$, and elliptic at $\cR^+$, and
  \begin{equation}
  \label{EqSpBHorSigns}
    \beta_0:=\rho_\infty^{-1}(\rho_\infty H_G)\rho_\infty > 0,\ \ \beta := \inf_{\cR^+}\bigl( \beta_0^{-1}\rho_\infty H_G t\bigr) > 0;
  \end{equation}
\item\label{ItSpBHor4} there exists a quadratic defining function\footnote{Note that $\dd(G|_{T^*\cX\setminus o})\neq 0$, as follows from $H_G\rho_\infty\neq 0$.} $\fq\in\CI(\pa\Sigma_0)$ (where $\pa\Sigma_0\subset S^*\cX$ is the boundary at fiber infinity) of $\pa\cR^+\subset S^*\cX$ such that in a neighborhood of $\pa\cR^+$, we have $\rho_\infty H_G\fq\geq c\fq$ for some $c>0$.
\end{enumerate}

These conditions are satisfied near the event horizon of Kerr, with $t,\rho_\infty$ replaced by $t_0,\xi_0^{-1}$ (in the notation of~\S\ref{SssTs3bH}), and with $\fq$ given by the restriction $\rho_\infty^2\sC$ to $\{\sigma=0\}$ of the function in Lemma~\ref{LemmaTs3bHQuadDef}. Gannot \cite{GannotHorizons} gives rather general conditions which guarantee the validity of~\eqref{EqSpBHorSigns} (the first condition of which captures the non-degenerate red-shift).

\begin{prop}[Radial point estimate]
\label{PropSpBHor}
  Let $\sigma\in\C$. Set $\beta_1:=\sup_{\cR^+}(\beta_0^{-1}\rho_\infty\upsigma^1(\frac{1}{2 i}(\wh{P_0}(0)-\wh{P_0}(0)^*)))$. Let $\cU\subset S^*\cX$ be a neighborhood of $\pa\cR^+$, and let $K\subset\cX$ be a compact neighborhood of $\cH$. Let $\chi\in\CIc(\cX)$, with $\chi=1$ near $K$.
  \begin{enumerate}
  \item{\rm (Direct estimate.)} Let $s>s_0>\frac12-\beta\Im\sigma+\beta_1$. There exist operators $B,G\in\Psi^0(\cX)$, with Schwartz kernels supported in $K\times K$ and operator wave front sets in $\cU$, and with $B$ elliptic at $\pa\cR^+$, such that
    \[
      \|B u\|_{H^s} \leq C\bigl( \|G\wh{P_0}(\sigma)u\|_{H^{s-1}} + \|\chi u\|_{H^{s_0}}\bigr)
    \]
    holds in the strong sense.
  \item{\rm (Adjoint estimate.)} Let $s>\frac12-\beta\Im\sigma+\beta_1$ and $-N<-s+1$. There exist operators $B,G,E\in\Psi^0(\cX)$ with Schwartz kernels supported in $K\times K$ and operator wave front sets in $\cU$, and with $B$ elliptic at $\pa\cR^+$ and $\WF'(E)\cap\pa\cR^+=\emptyset$, such that
    \[
      \|B u\|_{H^{-s+1}} \leq C\bigl( \|G\wh{P_0}(\sigma)^*u\|_{H^{-s}} + \|E u\|_{H^{-s+1}} + \|\chi u\|_{H^{-N}}\bigr)
    \]
    holds in the strong sense.
  \end{enumerate}
\end{prop}
\begin{proof}
  This is proved using a standard positive commutator argument as in \cite[\S{2.4}]{VasyMicroKerrdS}, \cite[Chapter~10]{HintzMicro}, and~\S\ref{SssR3R}. We only recall the symbolic computations and omit their regularization (needed to prove the strong version of the estimates) and quantization. Note that the principal symbol of $\wh{P_0}(\sigma)$ is $\sigma$-independent and given by $G_0$. Fix
  \begin{equation}
  \label{EqSpBHorPsi}
    \psi\in\CIc((-1,1)),\quad \psi=1\ \text{on}\ [-\tfrac12,\tfrac12],\quad \sqrt{-x\psi(x)\psi'(x)}\in\CI.
  \end{equation}
  For $A=\check A^*\check A$ where $\check A=\Op(\check a)=\check A^*$ with
  \begin{equation}
  \label{EqSbPHorChecka}
    \check a = \rho_\infty^{-s+\frac12}\psi_\Sigma\psi_\cR,\quad \psi_\Sigma:=\psi(\digamma\rho_\infty^2 G_0),\ \psi_\cR:=\psi(\digamma_\cR\fq),
  \end{equation}
  where $\digamma,\digamma_\cR>1$ are to be chosen later, we write $2\Im\la\check A\wh{P_0}(\sigma)u,\check A u\ra_{L^2}=\la\sC u,u\ra$ and compute the principal symbol of
  \begin{equation}
  \label{EqSbPHorsC}
    \sC := i[\wh{P_0}(\sigma),A] + 2\frac{\wh{P_0}(\sigma)-\wh{P_0}(\sigma)^*}{2 i}A
  \end{equation}
  in terms of $\sfH_{G_0}:=\rho_\infty H_{G_0}$ and
  \begin{equation}
  \label{EqSbPHorp1}
    \sfp_1(\sigma):=\rho_\infty\upsigma^1\Bigl(\frac{1}{2 i}(\wh{P_0}(\sigma)-\wh{P_0}(\sigma)^*)\Bigr)
  \end{equation}
  to be
  \begin{align*}
    H_{G_0}(\check a^2&) + 2\rho_\infty^{-1}\sfp_1\check a^2 = -\delta\rho_\infty^{-1}\check a^2 - b^2 - b_\cR^2 + j G_0 \\
    b &:= \rho_\infty^{-s}\psi_\Sigma\psi_\cR\bigl[ (\rho_\infty^{-1}\sfH_{G_0}\rho_\infty)(2 s-1) - \delta - 2\sfp_1(\sigma) \bigr]^{\frac12}, \\
    b_\cR &:= \rho_\infty^{-s}\psi_\Sigma\bigl[ -2(\sfH_{G_0}\fq)\digamma_\cR\psi_\cR\psi'_\cR \bigr]^{\frac12}, \\
    j &:= 2\rho_\infty^{-2 s+2}\psi_\cR \digamma(\rho_\infty^{-2}\sfH_{G_0}\rho_\infty^2)\psi_\Sigma\psi'_\Sigma.
  \end{align*}
  Choosing $\digamma_\cR$ large and then $\digamma$ large ensures that $b_\cR\in S^s$ is well-defined. For small $\delta>0$ and all sufficiently large $\digamma_\cR,\digamma$, we have $b\in S^s$, and $b$ is elliptic at $\cR^+$ since the term in square brackets in the expression for $b$ evaluates at $\cR^+$ and for $\delta=0$ to $2\beta_0$ times
  \begin{equation}
  \label{EqSbPHorThr}
    s-\frac12-\beta_0^{-1}\sfp_1(0) - \beta_0^{-1}(\sfp_1(\sigma)-\sfp_1(0)).
  \end{equation}
  Indeed, for $\sigma=0$ this is $\geq s-\frac12-\beta_1>0$. For general $\sigma$ we note that $\wh{P_0}(\sigma)=\wh{P_0}(0)+\sigma P_1+\sigma^2 P_2$ where $P_1\in\Diff^1$, $P_2\in\Diff^0$ are $\sigma$-independent; differentiating $\wh{P_0}(\sigma)=e^{i\sigma t}P_0 e^{-i\sigma t}$ at $\sigma=0$ gives
  \begin{equation}
  \label{EqSbPpasigma}
    \pa_\sigma\wh{P_0}(0)=-i[P_0,t]
  \end{equation}
  and thus $\upsigma^1(\pa_\sigma\wh{P_0}(0))=-H_{G_0}t$; therefore thus have
  \[
    \sfp_1(\sigma)-\sfp_1(0) = (\rho_\infty H_{G_0}t)\Bigl(\frac{1}{2 i} (-\sigma + \bar\sigma)\Bigr) = -(\Im\sigma)(\rho_\infty H_{G_0}t),
  \]
  and thus $-\beta_0^{-1}(\sfp_1(\sigma)-\sfp_1(0))\geq\beta\Im\sigma$.

  For the adjoint estimate, one uses $-s+1$ in place of $s$; this necessitates to a relative sign switch between $b^2$ and $b_\cR^2$, with $b$ still being the main term (elliptic at $\pa\cR^+$), but with $b_\cR$ now being a priori control term supported on $\supp\check a\cap\supp\dd\fq$, and thus in a punctured neighborhood of $\pa\cR^+$. We leave the detailed symbolic computation to the reader.
\end{proof}

\subsubsection{Radial point estimates near infinity}
\label{SssSpBInfty}

We now consider a stationary asymptotically Minkowski wave operator $P_0$ (Definition~\ref{DefSpGenMink}). Its spectral family is given by
\begin{equation}
\label{EqSpBInftyOp}
  \wh{P_0}(\sigma) = 2 i\sigma\rho\Bigl(\rho\pa_\rho-\frac{n-1}{2}-S\Bigr) + \wh{P_0}(0) - i\sigma Q + \sigma^2 g^{0 0} \in \Diff_\scop^2(\cX),
\end{equation}
where $\wh{P_0}(0)\in\rho^2\Diffb^2(\cX)$ has principal symbol equal to the Euclidean dual metric function plus lower order terms, and we recall $Q\in\rho^3\Diffb^1(\cX)$, $g^{0 0}\in\rho^2\CI(\cX)$. Writing scattering covectors as $\xi_\scop\,\dd r+r\eta_\scop$, $\eta_\scop\in T^*\Sph^2$, the scattering principal symbol of $\wh{P_0}(\sigma)$ is equal to
\[
  G_\sigma \equiv 2\sigma\xi_\scop + \xi_\scop^2 + |\eta_\scop|^2 \bmod \rho P^{[2]}(\Tsc^*\cX).
\]
Therefore, by a change of variables computation,
\[
  \sfH_{G_\sigma} := \rho^{-1}H_{G_\sigma} \equiv -2(\sigma+\xi_\scop)(\rho\pa_\rho+\eta_\scop\pa_{\eta_\scop}) + 2|\eta_\scop|^2\pa_{\xi_\scop} + \sfH_{|\eta_\scop|^2}
\]
modulo vector fields vanishing over $\pa\cX$; here $\sfH_{|\eta_\scop|^2}:=2\eta_\scop\cdot\pa_\omega$ at the center of geodesic normal coordinates $\omega$ around a point in $\Sph^2$. (For $\pm\sigma>0$, the vector field $\pm\sfH_{G_\sigma}$ can be integrated over $\pa\cX$ to yield the source-to-sink dynamics from $\cR_{\sigma,{\rm in}}=\{\xi_\scop=-2\sigma,\ \eta_\scop=0\}\cap\Tsc^*_{\pa\cX}\cX$ to $\cR_{\rm out}=\{\xi_\scop=\eta_\scop=0\}\cap\Tsc^*_{\pa\cX}\cX$ in the notation of~\eqref{EqSpOrderRinout}.) For $\Im\sigma>0$, we recall that $G_\sigma^{-1}(0)\cap\Tsc_{\pa\cX}^*\cX=\cR_{\rm out}$; thus, we only need to prove an estimate near $\cR_{\sigma,{\rm in}}$ when $\sigma$ is almost real.

\begin{prop}[Radial point estimate near $\cR_{\sigma,{\rm in}}$]
\label{PropSpBInftyIn}
  Let $0\neq\sigma_0\in\R$ and $\nu>0$. Let $s,N\in\R$. Define
  \begin{equation}
  \label{EqSpBInftyInThr}
    \vartheta_{\rm in}:=\sup_{\pa\cX} \Bigl( \upsigma_\bop^1\bigl(\Im\wh{P_0}(0)\bigr)\big|_{-2\frac{\dd r}{r}} - 2\Re S \Bigr).
  \end{equation}
  Let $\cU\subset\Tsc^*_{\pa\cX}\cX$ be a compact neighborhood of $\cR_{\sigma_0,{\rm in}}$, and let $K\subset\cX$ be a compact neighborhood of $\pa\cX$. Let $\chi\in\CIc(\cX)$, with $\chi=1$ near $K$.
  \begin{enumerate}
  \item{\rm (Direct estimate.)} Let $r>r_0>\frac12(-1+\vartheta_{\rm in})$. There exist operators $B,G\in\Psi_\scop^{0,0}(\cX)$ with Schwartz kernels supported in $K\times K$ and operator wave front sets in $\cU$, and with $B$ elliptic at $\cR_{\sigma_0,{\rm in}}$, such that
    \begin{equation}
    \label{EqSpBInftyIn}
      \|B u\|_{\Hsc^{s,r}} \leq C\Bigl( \|G\wh{P_0}(\sigma)u\|_{\Hsc^{s-2,r+1}} + \|u\|_{\Hsc^{-N,r_0}} \Bigr),\quad \sigma\in\sigma_0+i[0,\nu],
    \end{equation}
    holds in the strong sense.
  \item{\rm (Adjoint estimate.)} Let $r>\frac12(-1+\vartheta_{\rm in})$. There exist operators $B,E,G\in\Psi_\scop^{0,0}(\cX)$ with Schwartz kernels supported in $K\times K$ and operator wave front sets in $\cU$, and with $B$ elliptic at $\cR_{\sigma_0,{\rm in}}$ and $\WF'_\scop(E)\cap\cR_{\sigma,{\rm in}}=\emptyset$, such that
    \[
      \|B u\|_{\Hsc^{-s+2,-r-1}} \leq C\Bigl( \|G\wh{P_0}(\sigma)^*u\|_{\Hsc^{-s,-r}} + \|E u\|_{\Hsc^{-s+2,-r-1}} + \|u\|_{\Hsc^{-N,-N}} \Bigr),\quad \sigma\in\sigma_0+i[0,\nu],
    \]
    holds in the strong sense.
  \end{enumerate}
\end{prop}
\begin{proof}
  We only prove the direct estimate, which is a minor variation of the standard incoming scattering radial point estimate of Melrose \cite{MelroseEuclideanSpectralTheory}. (The proof of the adjoint estimate is very similar to the proof of the direct estimate at $\cR_{\rm out}$ given below.) Since $\cR_{\sigma,{\rm in}}$ is disjoint from ${}^\scop S^*\cX$, the regularity order is irrelevant. If $\nu'\in(0,\nu)$, then for $\sigma\in\sigma_0+i[\nu',\nu]$, the operator $\wh{P_0}(\sigma)$ is elliptic at $\cR_{\sigma_0,{\rm in}}$, and thus~\eqref{EqSpBInftyIn} follows from an elliptic estimate (and the orders $s-2,r+1$ on the right can be weakened to $s-2,r$).

  We moreover only consider the case $\Re\sigma=\sigma_0>0$ and leave the sign changes required to treat the case $\sigma_0<0$ to the reader. For a suitable commutant $A=\check A^*\check A$, $\check A=\check A^*\in\Psi_\scop^{-\infty,r+\frac12}$, we write $2\Im\la\check A\wh{P_0}(\sigma)u,\check A u\ra_{L^2}=\la\sC u,u\ra$ where
  \begin{align}
    \sC &= i[\Re\wh{P_0}(\sigma),A] - [\Im\wh{P_0}(\sigma),\check A^2] + 2(\Im\wh{P_0}(\sigma))\check A^2 \nonumber\\
  \label{EqSpBInftyC}
      &= i[\Re\wh{P_0}(\sigma),A] + 2\check A(\Im\wh{P_0}(\sigma))\check A + [[\Im\wh{P_0}(\sigma),\check A],\check A].
  \end{align}
  With $\psi$ as in~\eqref{EqSpBHorPsi}, we take $\check A$ to be a quantization of
  \begin{equation}
  \label{EqSpBInftyChecka}
    \check a = \rho^{-r-\frac12}\psi_\pa\psi_\cR, \quad \psi_\cR=\psi(\digamma_\cR\fq),\ \psi_\pa=\psi(\digamma\rho),
  \end{equation}
  where $\fq:=(\xi_\scop+2\Re\sigma)^2+|\eta_\scop|^2$.

  For the first term of~\eqref{EqSpBInftyC}, we note that Hamiltonian vector field of the principal symbol of $\Re\wh{P_0}(\sigma)$ is
  \[
    \sfH_{\Re G_\sigma} = \Re\sfH_{G_\sigma} = -2(\xi_\scop+\Re\sigma)(\rho\pa_\rho+\eta_\scop\pa_{\eta_\scop}) + 2|\eta_\scop|^2\pa_{\xi_\scop} + \sfH_{|\eta_\scop|^2}.
  \]
  Over $\pa\cX$, we then compute
  \[
    \sfH_{\Re G_\sigma}\fq = 2|\eta_\scop|^2(\xi_\scop+2\Re\sigma) - 2(\xi_\scop+\Re\sigma)|\eta_\scop|^2,
  \]
  which on the characteristic set $\{(\xi_\scop+2\Re\sigma)\xi_\scop+|\eta_\scop|^2=0\}$ of $\wh{P_0}(\sigma)$ over $\pa\cX$ equals $-2\xi_\scop(\xi_\scop+2\Re\sigma)^2-2(\xi_\scop+\Re\sigma)|\eta_\scop|^2$; over $\pa\cX$ but near $\cR_{\Re\sigma,{\rm in}}$, where $\xi_\scop=-2\Re\sigma$, we thus see that
  \[
    \sfH_{\Re G_\sigma}\fq\geq (\Re\sigma)\fq.
  \]

  The third term of~\eqref{EqSpBInftyC} is a double commutator, with $\Im\wh{P_0}(\sigma)\in\Diff_\scop^1(\cX)$, and thus of class $\Psi_\scop^{-\infty,2 r-1}$ and hence subprincipal. For the second term of~\eqref{EqSpBInftyC}, we compute
  \begin{equation}
  \label{EqSpBInftyIm}
    \Im\wh{P_0}(\sigma) = \Im\wh{P_0}(0) + \Im\bigl(\sigma\pa_\sigma\wh{P_0}(0)\bigr) + \Im\Bigl(\frac{\sigma^2}{2}\pa_\sigma^2\wh{P_0}(0)\Bigr).
  \end{equation}
  The last summand is of class $\rho^2\CI$ and thus contributes a subprincipal term (of class $\Psi_\scop^{-\infty,2 r-1}$) to $\sC$. The second term is
  \begin{equation}
  \label{EqSpBInftySplit}
    \Re(\sigma)\Im\bigl(\pa_\sigma\wh{P_0}(0)\bigr) + \Im(\sigma)\Re\bigl(\pa_\sigma\wh{P_0}(0)\bigr).
  \end{equation}
  Now, $\pa_\sigma\wh{P_0}(0)=2 i\rho(\rho\pa_\rho-\frac{n-1}{2}-S)-i Q$ is symmetric with respect to the Euclidean volume density $r^{n-1}|\dd r\,\dd\slg|=\rho^{-n+1}|\frac{\dd\rho}{\rho^2}\,\dd\slg|$ when $S=0$ and $Q=0$; and since $Q\in\rho^2\Diffsc^1$, only $S$ contributes to $\Im(\pa_\sigma\wh{P_0}(0))$ to leading order. We thus find
  \begin{subequations}
  \begin{equation}
  \label{EqSpBInftyIm1}
    \upsigma_\scop^1\Bigl(\rho^{-1}\Re(\sigma)\Im\bigl(\pa_\sigma\wh{P_0}(0)\bigr)\Bigr)\Big|_{\cR_{\sigma,{\rm in}}} = -2\Re\sigma\Re S.
  \end{equation}
  (We analyze the second term in~\eqref{EqSpBInftySplit} below.) In view of $\Im\wh{P_0}(0)\in\rho^2\Diffb^1$, the scattering principal symbol of $\rho^{-1}\Im\wh{P_0}(0)$ is fiber-linear over $\pa\cX$, and therefore
  \begin{equation}
  \label{EqSpBInftyIm2}
    \upsigma_\scop^1\Bigl(\rho^{-1}\Im\wh{P_0}(0)\Bigr)\Big|_{\cR_{\sigma,{\rm in}}} = (\Re\sigma)\sfp_1,\quad \sfp_1:=\upsigma_\scop^1\Bigl(\rho^{-1}\Im\wh{P_0}(0)\Bigr)\Big|_{\cR_{1,{\rm in}}}.
  \end{equation}
  \end{subequations}

  Altogether, then, we can write, for sufficiently small $\delta>0$,
  \begin{equation}
  \label{EqSbPInftyComm}
    \sC = \Op_\scop\bigl( -\delta\rho\check a^2 - b^2 - b_\cR^2 \bigr) - 2(\Im\sigma)\check A\Re\bigl(-\pa_\sigma\wh{P_0}(0)\bigr)\check A + R,
  \end{equation}
  where $R\in\Psi_\scop^{-\infty,2 r-1}$ and\footnote{Differentiation of $\psi_\pa$ yields a symbol of class $S^{-\infty,-\infty}$, and hence does not enter here. For the same reason, we could as well have omitted it from the definition of $\check a$, much as we did not include cutoffs away from the zero section in~\eqref{EqSbPHorChecka}, say; but we did include the $\rho$-cutoff here for the sake of clarity.}
  \begin{align*}
    b &:= \rho^{-r}\psi_\pa\psi_\cR \bigl[ (\rho^{-1}\sfH_{\Re G_\sigma}\rho) (2 r+1) - \delta + (\Re\sigma) (4\Re S-2\sfp_1) \bigr]^{\frac12}, \\
    b_\cR &:= \rho^{-r}\psi_\pa\bigl[ -2(\sfH_{\Re G_\sigma}\fq)\digamma_\cR\psi_\cR\psi'_\cR \bigr]^{\frac12}.
  \end{align*}
  At $\cR_{\Re\sigma,{\rm in}}$ where $\xi_\scop=-2\Re\sigma$, we have $\rho^{-1}\sfH_{\Re G_\sigma}\rho=2\Re\sigma$. The threshold condition for $b$ to be well-defined is thus $4 r+2+4\Re S-2\sfp_1>0$, as stated. To complete the proof, it remains to prove that the $\Im\sigma$-term has a sign, to wit,
  \begin{equation}
  \label{EqSbPInftyLower}
    \big\la \check A\Re\bigl(-\pa_\sigma\wh{P_0}(0)\bigr)\check A u,u\big\ra_{L^2} \geq -C\|\chi u\|_{H_\scop^{-N,r_0}}^2.
  \end{equation}
  But $\upsigma_\scop^1(-\Re(\pa_\sigma\wh{P_0}(0)))=-2\xi_\scop=4\Re\sigma$ at $\cR_{\Re\sigma,{\rm in}}$, which is thus positive on $\supp\check a$; thus we can write $\check A\Re(-\pa_\sigma\wh{P_0}(0))\check A=F^*F+R'$ where $F\in\Psi_\scop^{-\infty,r+\frac12}$ and $R'\in\Psi_\scop^{-\infty,-\infty}$, which implies~\eqref{EqSbPInftyLower}.
\end{proof}

Unlike the incoming radial set, the outgoing radial set over $\pa\cX$ persists also for complex $\sigma$. It is then convenient to emphasize the outgoing vector field $-\rho\pa_\rho$ by considering the rescaling
\begin{equation}
\label{EqSbPInftyOp}
  \tilde P_0(\sigma) := \sigma^{-1}\wh{P_0}(\sigma) = 2 i\rho\Bigl(\rho\pa_\rho-\frac{n-1}{2}-S\Bigr) + \frac{\bar\sigma}{|\sigma|^2}\wh{P_0}(0) - i Q + \sigma g^{0 0},
\end{equation}
so its principal symbol $\tilde G_\sigma$ and rescaled Hamiltonian vector field $\sfH_{\tilde G_\sigma}=\rho^{-1}H_{\tilde G_\sigma}$ satisfy
\begin{align*}
  \tilde G_\sigma &\equiv 2\xi_\scop + \frac{\bar\sigma}{|\sigma|^2}(\xi_\scop^2+|\eta_\scop|^2), \\
  \Re\sfH_{\tilde G_\sigma} &\equiv -2\Bigl(1+\frac{\Re\sigma}{|\sigma|^2}\xi_\scop\Bigr)(\rho\pa_\rho+\eta_\scop\pa_{\eta_\scop}) + 2\frac{\Re\sigma}{|\sigma|^2}|\eta_\scop|^2\pa_{\xi_\scop} + \frac{\Re\sigma}{|\sigma|^2}\sfH_{|\eta_\scop|^2}.
\end{align*}
Moreover, using again the symmetry of $2 i\rho(\rho\pa_\rho-\frac{n-1}{2})$, we have
\begin{equation}
\label{EqSpBInftyImTilde}
  \Im\tilde P_0(\sigma) \equiv -2\rho\Re S + (\Re\sigma)|\sigma|^{-2}\Im\wh{P_0}(0) - (\Im\sigma)|\sigma|^{-2}\Re\wh{P_0}(0) \bmod \rho^2\Diffsc^1.
\end{equation}

\begin{prop}[Radial point estimate near $\cR_{\rm out}$]
\label{PropSpBInftyOut}
  Let $\cD$ be a compact subset of $\{\sigma\in\C\colon \sigma\neq 0,\ \Im\sigma\geq 0\}$. Let $s,N\in\R$. Recall the quantity $\ubar S$ from~\eqref{EqSSAdmubarS}. Let $\cU\subset\Tsc^*_{\pa\cX}\cX$ be a compact neighborhood of $\cR_{\rm out}$, and let $K\subset\cX$ be a compact neighborhood of $\pa\cX$. Let $\chi\in\CIc(\cX)$, with $\chi=1$ near $K$.
  \begin{enumerate}
  \item\label{ItSpBInftyOutDir}{\rm (Direct estimate.)} Let $r<-\frac12+\ubar S$. There exist operators $B,E,G\in\Psi_\scop^{0,0}(\cX)$ with Schwartz kernels supported in $K\times K$ and operator wave front sets in $\cU$, and with $B$ elliptic at $\cR_{\rm out}$ and $\WF'_\scop(E)\cap\cR_{\rm out}=\emptyset$, such that
    \begin{equation}
    \label{EqSpBInftyOut}
      \|B u\|_{\Hsc^{s,r}} \leq C\Bigl( \|G\wh{P_0}(\sigma)u\|_{\Hsc^{s-2,r+1}} + \|E u\|_{\Hsc^{s,r}} + \|\chi u\|_{\Hsc^{-N,r_0}} \Bigr),\quad\sigma\in\cD,
    \end{equation}
    holds in the strong sense. More generally, for every $k\in\N_0$, we have the (strong) estimate
    \begin{equation}
    \label{EqSpBInftyOutb}
      \|B u\|_{H_{\scop;\bop}^{(s;k),r}} \leq C\Bigl( \|G\wh{P_0}(\sigma)u\|_{H_{\scop;\bop}^{(s-2;k),r+1}} + \|E u\|_{H_{\scop;\bop}^{(s;k),r}} + \|\chi u\|_{\Hsc^{-N,-N}} \Bigr),\quad\sigma\in\cD.
    \end{equation}
  \item\label{ItSpBInftyOutAdj}{\rm (Adjoint estimate.)} Let $r<r_0<-\frac12+\ubar S$. There exist operators $B,G\in\Psi_\scop^{0,0}(\cX)$ with Schwartz kernels supported in $K\times K$ and operator wave front sets in $\cU$, and with $B$ elliptic at $\cR_{\rm out}$, such that
    \[
      \|B u\|_{\Hsc^{-s+2,-r-1}} \leq C\Bigl( \|G\wh{P_0}(\sigma)^*u\|_{\Hsc^{-s,-r}} + \|\chi u\|_{\Hsc^{-N,-r_0-1}}\Bigr),\quad\sigma\in\cD,
    \]
    holds in the strong sense.
  \end{enumerate}
\end{prop}

\begin{rmk}[b-regularity]
\label{RmkSpBInftyModule}
  The mixed norms in~\eqref{EqSpBInftyOutb} are the sums of $\Hsc$-norms of a distribution and its up to $k$-fold derivatives along $\rho\pa_\rho=\rho^{-1}\rho^2\pa_\rho$ and $\pa_\omega=\rho^{-1}\rho\pa_\omega$, which span, microlocally near $\cR_{\rm out}$ and over $\Psi_\scop^{0,0}$, the space of elements of $\Psi_\scop^{1,1}$ with principal symbols vanishing at $\cR_{\rm out}$. Thus,~\eqref{EqSpBInftyOutb} is an estimate featuring \emph{module regularity}; see \cite[\S{2.4}]{GellRedmanHassellShapiroZhangHelmholtz}.
\end{rmk}

The structural input for the proof of b-regularity is the following variant of the considerations in~\S\ref{SssRbComm} and \cite[\S{2.4}]{GellRedmanHassellShapiroZhangHelmholtz}.

\begin{lemma}[Commutators of $\wh{P_0}(\sigma)$]
\label{LemmaSpBInftyOutComm}
  Let $V_1=\rho\pa_\rho$ and fix vector fields $V_2,\ldots,V_N\in\cV(\Sph^{n-1})$ spanning $\cV(\Sph^{n-1})$ over $\CI(\Sph^{n-1})$. Write $\sV=(V_1,\ldots,V_N)$ and $\sV^\alpha=\prod_{j=1}^N V_j^{\alpha_j}$ for $\alpha\in\N_0^N$. If $\wh{P_0}(\sigma)u=f$, then
  \[
    u^{(k)}:=(\sV^\alpha u)_{\alpha\in\cA_k},\quad \cA_k:=\{\alpha\in\N_0^N\colon|\alpha|=k\},
  \]
  and $f^{(k)}:=(\sV^\alpha f)_{\alpha\in\cA_k}$ satisfy
  \[
    P_\sigma^{(k)}u^{(k)} = f^{(k)} + \tilde f^{(k)},\quad \tilde f^{(k)} = \Biggl(\;\sum_{|\beta|\leq k-1} R_{\alpha\beta}\sV^\beta u \Biggr)_{\alpha\in\cA_k},
  \]
  where $R_{\alpha\beta}\in\rho\Diffsc^2$, while the $(\alpha,\alpha')$-matrix element of $P_\sigma^{(k)}$ is equal to
  \begin{equation}
  \label{EqSpBInftyOutComm}
    (P^{(k)}_\sigma)_{\alpha\alpha'} = \delta_{\alpha\alpha'}\wh{P_0}(\sigma) - 2 i k\sigma \rho \delta_{\alpha,(k,0,\ldots,0)} + \tilde P_{\alpha\alpha'}
  \end{equation}
  where $\tilde P_{\alpha\alpha'}\in\rho\Diffsc^1$ and $\upsigma_\scop^{1,0}(\rho^{-1}\tilde P_{\alpha\alpha'})|_{\cR_{\rm out}}=0$.
\end{lemma}
\begin{proof}
  We need to study the commutator $[\wh{P_0}(\sigma),\sV^\alpha]$ for $\alpha\in\cA_k$. Terms in $[\wh{P_0}(\sigma),\sV^\alpha]$ arising from contributions to $\wh{P_0}(\sigma)$ of class $\rho\CI$, $\rho^2\Diffb^1$, resp.\ $\rho^3\Diffb^2$ are of class $\rho\Diffb^{k-1}$, $\rho^2\Diffb^k$, resp.\ $\rho^3\Diffb^{k+1}$. Since $\sV$ spans $\Vb([0,1)_\rho\times\Sph^{n-1})$ over $\CI([0,1)_\rho\times\Sph^{n-1})$, such terms can be written as sums of terms of the form $\rho\Diffb^0\circ\sV^{\leq k-1}$, $\rho^2\Diffb^1\circ\sV^{\leq k-1}$, resp.\ $\rho^3\Diffb^2\circ\sV^{\leq k-1}$, and thus a fortiori of the form $\sum_{|\beta|\leq k-1}R_{\alpha\beta}\sV^\beta$ where $R_{\alpha\beta}\in\rho\Diffsc^2$; we put such terms into $\tilde f^{(k)}$. This applies to all terms in~\eqref{EqSpBInftyOp} except for $2 i\sigma\rho^2\pa_\rho$ and the term $\rho^2((\rho D_\rho)^2+\slDelta)$ of $\wh{P_0}(0)$.

  We have $[\sV^\alpha,2 i\sigma\rho^2\pa_\rho]=0$ unless $\alpha=(k,0,\ldots,0)$, in which case
  \[
    [(\rho\pa_\rho)^k,2 i\sigma\rho\,\rho\pa_\rho] \equiv 2 i k\sigma\rho(\rho\pa_\rho)^k \bmod \rho\Diffb^{k-1};
  \]
  we put the remainder term into $\tilde f^{(k)}$ and the main term into $P^{(k)}_\sigma$. Finally, we can write
  \[
    [\sV^\alpha,\rho^2((\rho D_\rho)^2+\slDelta)]\in\rho^2\Diffb^{k+1}=\rho^2\Diffb^1\circ\Diffb^k
  \]
  as a sum of terms $\tilde P_{\alpha\beta}\sV^\beta$ where $|\beta|\leq k$ and $\tilde P_{\alpha\beta}\in\rho^2\Diffb^1\subset\rho\Diffsc^1$; thus the principal symbol of $\rho^{-1}\tilde P_{\alpha\beta}\in\Diffsc^1$ vanishes at the scattering zero section $\cR_{\rm out}$ over $\pa\cX$. The terms with $|\beta|=k$, resp.\ $|\beta|\leq k-1$ are put into $P^{(k)}_\alpha$, resp.\ $\tilde f^{(k)}$.
\end{proof}

\begin{proof}[Proof of Proposition~\usref{PropSpBInftyOut}]
  We only prove the direct estimate.

  \pfstep{Step 2.~Basic estimate.} We start with~\eqref{EqSpBInftyOut} (but omit the regularization needed to prove the strong nature of this estimate). We take $A=\check A^2$ where $\check A=\Op_\scop(\check a)=\check A^*\in\Psi_\scop^{-\infty,r+\frac12}$, with
  \begin{equation}
  \label{EqSpBInftyOutChecka}
    \check a = \rho^{-r-\frac12}\psi_\pa\psi_\cR,\quad \psi_\cR=\psi(\digamma_\cR\fq),\ \psi_\pa=\psi(\digamma\rho),
  \end{equation}
  where now $\fq=\xi_\scop^2+|\eta_\scop|^2$, and we will choose $\digamma_\cR$ large and then $\digamma$ large. We study the equation for $\tilde P_0(\sigma)$ and write
  \[
    2\Im\la\check A\tilde P_0(\sigma)u,\check A u\ra=\la\sC u,u\ra
  \]
  where, analogously to~\eqref{EqSpBInftyC},
  \begin{equation}
  \label{EqSpBInftyOutsC}
    \sC = i[\Re\tilde P_0(\sigma),A] + 2\check A(\Im\tilde P_0(\sigma))\check A + [[\Im\tilde P_0(\sigma),\check A],\check A].
  \end{equation}
  The double commutator is again subprincipal, being of class $\Psi_\scop^{-\infty,2 r-1}$.

  \pfsubstep{(1.1)}{Almost real $\sigma$.} Consider $\sigma\in\sigma_0+i[0,\nu]$ where $\sigma_0>0$ (the case $\sigma_0<0$ being completely analogous) and $\nu>0$. Then $\Re\sigma=\sigma_0\neq 0$, and thus over $\pa\cX$, the equation $\Re\tilde G_\sigma=0$ implies $|\eta_\scop|^2=-\xi_\scop(2\frac{|\sigma|^2}{\Re\sigma}+\xi_\scop)$. Therefore,
  \begin{equation}
  \label{EqSpBInftyOutq}
    \sfH_{\Re\tilde G_\sigma}\fq = -4\Bigl(1+\frac{\Re\sigma}{|\sigma|^2}\xi_\scop\Bigr)|\eta_\scop|^2 + 4\frac{\Re\sigma}{|\sigma|^2}|\eta_\scop|^2\xi_\scop \leq -2\fq
  \end{equation}
  over $\pa\cX$ near $\cR_{\rm out}$.

  As a preparation for the treatment of the $\Im\tilde P_0(\sigma)$ term as computed in~\eqref{EqSpBInftyImTilde}, observe that we can write $\rho^{-1}\Re\wh{P_0}(0)\rho^{-1}\equiv(\rho D_\rho)^2+\slDelta\bmod\Diffb^1+\rho\Diffb^2$ as a sum of squares of finitely many b-vector fields to leading order, and hence
  \begin{equation}
  \label{EqSpBInftyReP0}
    \Re\wh{P_0}(0) = \sum_j (\rho V_j)^*(\rho V_j) + R_1 + R_2,\quad V_j\in\Diffb^1(\cX),\ R_1\in\rho^2\Diffb^1(\cX),\ R_2\in\rho^3\Diffb^2(\cX).
  \end{equation}
  (This observation was already made by Vasy in \cite[Lemma~4.19]{VasyLAPLag}.) Thus, $\sfr_1:=\upsigma_\scop^{1,0}(\rho^{-1}R_1)$ and $\sfr_2:=\upsigma_\scop^{2,0}(\rho^{-1}R_2)$ are homogeneous of degree $1$ and $2$ in the fibers of $\Tsc^*_{\pa\cX}\cX$, respectively, and in particular vanish at the zero section $\cR_{\rm out}$. Furthermore, since $\Im\wh{P_0}(0)\in\rho^2\Diffb^1$, also the scattering principal symbol $\sfr_3$ of $\rho^{-1}\Im\wh{P_0}(0)$ vanishes at the zero section.

  Recalling~\eqref{EqSpBInftyImTilde}, we can thus write, for sufficiently small $\delta>0$,
  \begin{equation}
  \label{EqSpBInftyOutSymb}
    \sC = \Op_\scop( -\delta\rho\check a^2 - b^2 + e_\cR^2 ) - 2(\Im\sigma)|\sigma|^{-2}\check A\biggl(\sum_j (\rho V_j)^*(\rho V_j)\biggr)\check A + R
  \end{equation}
  where $R\in\Psi_\scop^{-\infty,2 r-1}$ and
  \begin{align}
  \label{EqSpBInftyOutSymb2}
    b &:= \rho^{-r}\psi_\pa\psi_\cR\bigl[ (\rho^{-1}\sfH_{\Re\tilde G_\sigma}\rho)(2 r+1) + 4\Re S + (\Im\sigma)|\sigma|^{-2}(\sfr_1+\sfr_2) - (\Re\sigma)|\sigma|^{-2}\sfr_3 - \delta\bigr]^{\frac12}, \\
    e_\cR &:= \rho^{-r}\psi_\pa\bigl[ 2(\sfH_{\Re\tilde G_\sigma}\fq) \digamma_\cR\psi_\cR\psi'_\cR \bigr]^{\frac12}. \nonumber
  \end{align}
  Note then that at $\cR_{\rm out}$ we have $\rho^{-1}\sfH_{\Re\tilde G_\sigma}\rho=-2$, so $b$ is well-defined and elliptic since $-4 r-2+4\Re S>0$; and $e_\cR$ is well-defined for large $\digamma_\cR$ (and then large $\digamma$) by~\eqref{EqSpBInftyOutq}.

  Since the term
  \begin{equation}
  \label{EqSpBInftyOutSign}
    -2(\Im\sigma)|\sigma|^{-2}\la \check A(\rho V_j)^*\rho V_j\check A u,u\ra = -2(\Im\sigma)|\sigma|^{-2}\|\rho V_j\check A u\|_{L^2}^2
  \end{equation}
  is $\leq 0$ since $\Im\sigma\geq 0$ and thus has the same term as the main term $-b^2$ in the symbolic computation, we now obtain~\eqref{EqSpBInftyOut} in the usual fashion; the a priori control term $E u$ arises from the term $e_\cR^2$ in~\eqref{EqSpBInftyOutSymb}.

  \pfsubstep{(1.2)}{Non-real $\sigma$.} When $\Im\sigma\geq\nu>0$, the operator $\tilde P_0(\sigma)$ is elliptic at $\Tsc^*_{\pa\cX}\cX\setminus\cR_{\rm out}$, and thus we do not need to record any sign for $\sfH_{\Re\tilde G_\sigma}\fq$ anymore. (In fact, the computation~\eqref{EqSpBInftyOutq} does not apply anymore.) We run the same commutator argument as before, except now we replace the term $\Op_\scop(e_\cR^2)$ in~\eqref{EqSpBInftyOutSymb} by $\Op_\scop(j)\tilde P_0(\sigma)$ where
  \[
    j := \rho^{-2 r}\psi_\pa^2 (\sfH_{\Re\tilde P_0(\sigma)}\fq)\digamma_\cR\psi_\cR\psi'_\cR \tilde G_\sigma^{-1}
  \]
  is well-defined when $\digamma$ (localizing near $\rho=0$) is sufficiently small, as then $\tilde G_\sigma$ is elliptic on $\supp\psi_\pa\psi'_\cR$.

  \pfstep{Step 2.~Estimate with b-regularity.} We prove~\eqref{EqSpBInftyOutb} by induction on $k$, the case $k=0$ being~\eqref{EqSpBInftyOut}. In the notation of Lemma~\ref{LemmaSpBInftyOutComm}, we apply the estimate~\eqref{EqSpBInftyOut} to $\sigma^{-1}P^{(k)}_\sigma u^{(k)}=\sigma^{-1}f^{(k)}$. The term $\tilde P_{\alpha\alpha'}$ in~\eqref{EqSpBInftyOutComm} yields a vanishing contribution to the principal symbol computation for $\sC$ at $\cR_{\rm out}$. Consider the subprincipal contribution
  \[
    -2 k i\rho\delta_{\alpha,(k,0,\ldots,0)}
  \]
  to $\sigma^{-1}P^{(k)}_\sigma$: it enters as a contribution to $2\check A(\Im P^{(k)}_\sigma)\check A$ (cf.\ \eqref{EqSpBInftyOutsC}) in the form
  \[
    -2 k \check A \rho\delta_{\alpha,(k,0,\ldots,0)} \check A,
  \]
  which thus, upon acting on $u$ and pairing with $u$, has the correct sign (namely, that of the main term $-b^2$ in~\eqref{EqSpBInftyOutSymb}) for it to be dropped in the commutator estimate. To obtain~\eqref{EqSpBInftyOutb}, it then remains to observe, firstly, that for $|\beta|\leq k-1$ and $R_{\alpha\beta}\in\rho\Diffsc^2$,
  \[
    \|G R_{\alpha\beta}\sV^\beta u\|_{\Hsc^{s-2,r+1}} \leq C\Bigl(\|\tilde G u\|_{H_{\scop;\bop}^{(s;k-1),r}} + \|\tilde\chi u\|_{H_{\scop;\bop}^{(-N;k-1),r_0}}\Bigr)
  \]
  if $\tilde G\in\Psi_\scop^{0,0}$ is elliptic on $\WF'_\scop(G)$ and has Schwartz kernel contained in $\tilde\chi^{-1}(1)\times\tilde\chi^{-1}(1)$ where $\tilde\chi\in\CIc(\cX)$ equals $1$ near $\supp\chi$. (Here, we use that commutators of b-vector fields with scattering ps.d.o.s with b-regular coefficients are, again, scattering ps.d.o.s of the same order, cf.\ \eqref{EqMSComm}.) The norm of $\tilde G u$ can then be estimated using the inductive hypothesis. Secondly,
  \[
    \|\chi\sV^\beta u\|_{\Hsc^{-N,r_0}} \leq C\|\tilde\chi u\|_{\Hsc^{-N+k-1,r_0+k-1}};
  \]
  but since $N,r_0$ are arbitrary, this is of the required form in~\eqref{EqSpBInftyOutb}.
\end{proof}

\begin{rmk}[b-regularity away from the scattering zero section]
\label{RmkSpBNob}
  For every $\varpi\in\Tsc^*\cX\setminus{}^\scop o_{\pa\cX}$ (with ${}^\scop o_{\pa\cX}\subset\Tsc^*_{\pa\cX}\cX$ denoting the zero section) there exists a b-vector field $V\in\Vb(\cX)$ such that $\rho^{-1}V\in\Vsc(\cX)$ is elliptic at $\varpi$. Therefore, $H_{\scop;\bop}^{(s;k),r}$ and $H_\scop^{s+k,r+k}$ are microlocally the same away from ${}^\scop o_{\pa\cX}$; more precisely, if $B\in\Psi_\scop^{0,0}$ satisfies $\WF'_\scop(B)\cap{}^\scop o_{\pa\cX}=\emptyset$ and has Schwartz kernel supported in $K\times K$, then for $\chi\in\CIc(\cX)$ equal to $1$ near $K$, we have
  \begin{equation}
  \label{EqSpBNobNorm}
  \begin{split}
    \|B u\|_{H_{\scop;\bop}^{(s;k),r}} &\leq C\Bigl(\|B u\|_{H_\scop^{s+k,r+k}} + \|\chi u\|_{H_\scop^{-N,-N}}\Bigr), \\
    \|B u\|_{H_\scop^{s+k,r+k}} &\leq C\Bigl(\|B u\|_{H_{\scop;\bop}^{(s;k),r}} + \|\chi u\|_{H_\scop^{-N,-N}}\Bigr).
  \end{split}
  \end{equation}
  This is why we do not need to prove estimates on $H_{\scop;\bop}$-spaces near any radial set other than that near $\cR_{\rm out}$. See also Lemma~\ref{LemmaDRes} below for a closely related result.
\end{rmk}

\begin{rmk}[Second microlocal spaces]
\label{RmkSpB2nd}
  Instead of the direct commutation argument in Step~2 above, one could alternatively work with the second microlocal function spaces introduced by Vasy \cite{VasyLAPLag}, using which one could prove estimates with $k\geq 0$ degrees of b-regularity directly: this amounts to working on second microlocal Sobolev spaces where the scattering decay order and scattering regularity are increased by $k$, while the b-decay order (which should be thought of as the scattering decay order at the zero section) is unchanged.
\end{rmk}

\subsubsection{Proof of Theorem~\usref{ThmSpB}}
\label{SssSpBPf}

\pfstep{Parts~\eqref{ItSpBReal}--\eqref{ItSpBImag} of Theorem~\usref{ThmSpB}.} We first concatenate the radial point estimates at the sources for the future null-bicharacteristic flow (Propositions~\ref{PropSpBHor} and \ref{PropSpBInftyIn}) with real principal type propagation estimates and the radial point estimate at the sink over $\pa X$ (Proposition~\ref{PropSpBInftyOut}). For $\chi\in\CI(X)$ which equals $1$ for $r\geq r_1^+\in(\bhm,r_+)$, and $0$ for $r\leq r_1\in(\bhm,r_1^+)$, this gives the (strong) estimate of $\|\chi u\|_{H_\scop^{\sfs_\scop,\sfr+\alpha_+}}$ by the right-hand side of~\eqref{EqSpBEst}. Control of $u$ on $r^{-1}([\bhm,r_1^+])$ then follows by propagation from $r^{-1}([r_1^+,\bhm])$ using a standard energy estimate (see \cite[Corollary~9.38]{HintzMicro}).

\pfstep{Part~\eqref{ItSpBFred}.} The Fredholm property of the map~\eqref{EqSpBMap} follows from the compactness of the inclusions $\bar H_\scop^{\sfs_\scop,\sfr+\alpha_+}\hra\bar H_\scop^{-N,-N}$ and $\dot H_\scop^{-\sfs_\scop+1,-\sfr-\alpha_+-1}\hra\dot H_\scop^{-N,-N}$ in~\eqref{EqSpBEst}. Write $\dot\HH^+:=\{\sigma\in\C\colon \sigma\neq 0,\ \Im\sigma\geq 0\}$. Given $\sigma_0\in\dot\HH$, consider the nullspace $\cN$ and a complementary subspace $\cR\subset\CIc(X^\circ)$ in the domain, resp.\ range of~\eqref{EqSpBMap}. Define the map
\begin{equation}
\label{EqSpBPfwtP}
  \wt{P_0}(\sigma) \colon \cX_\sigma \oplus \C^{\dim\cR} \to \cY_\sigma \oplus \C^{\dim\cN},\quad
  \wt{P_0}(\sigma) = \begin{pmatrix} \wh{P_0}(\sigma) & E_{0 1} \\ E_{1 0} & 0 \end{pmatrix},
\end{equation}
where $E_{0 1}\colon\C^{\dim\cR}\to\cR$ is an isomorphism and $E_{1 0}\colon\cX_\sigma\to\C^{\dim\cN}$ restricts to an isomorphism on $\cN$ when $\sigma=\sigma_0$; we take $E_{1 0}$ to be an $\dim\cN$-tuple of linear functionals of the form $\la\cdot,v\ra_{L^2}$ where $v\in\CIc(X^\circ)$. Now $\wt{P_0}(\sigma_0)$ is invertible, and we have uniform Fredholm estimates~\eqref{EqSpBEst} for $\wt{P_0}(\sigma)$ when $\sigma\in\dot\HH$ lies in a neighborhood of $\sigma_0$, where we use norms on $\bar H_\scop^{\sfs_\scop,\sfr+\alpha_+}\oplus\C^{\dim\cR}$ for $u$ and on $\bar H_\scop^{\sfs_\scop-1,\sfr+\alpha_++1}\oplus\C^{\dim\cN}$ for $\wt{P_0}(\sigma)$. As shown in \cite[\S{2.7}]{VasyMicroKerrdS}, this implies the invertibility of $\wt{P_0}(\sigma)$ for $\sigma$ near $\sigma_0$. We recall the argument for $\ker\wt{P_0}(\sigma)=0$, the triviality of the cokernel being similar: if $\wt{P_0}(\sigma_j)(u_j,c_j)=0$ where $\dot\HH\ni\sigma_j\to\sigma$ and $u_j\in\cX_{\sigma_j}$, $c_j\in\C^{\dim\cR}$ with $(u_j,c_j)\neq 0$, then we may normalize $(u_j,c_j)$ to have $\bar H_\scop^{-N,-N}\oplus\C^{\dim\cR}$-norm $1$. The estimate~\eqref{EqSpBEst} then yields a uniform bound for the $\bar H_\scop^{\sfs_\scop,\sfr+\alpha_+}\oplus\C^{\dim\cR}$-norm of $(u_j,c_j)$. (We can work here with $\sigma$-independent variable orders.) Upon passing to a subsequence, $(u_j,c_j)$ thus converges weakly in this space to some limit $(u,c)$, and thus in norm in $\bar H_\scop^{-N,-N}\oplus\C^{\dim\cR}$, with the limit $(u,c)$ having norm $1$. But taking (distributional) limits in $0=\wt{P_0}(\sigma_j)(u_j,c_j)=\wh{P_0}(\sigma_j)u_j+E_{0 1}(c_j)$ yields $\wt{P_0}(\sigma)(u,c)=0$, contradicting the injectivity of $\wt{P_0}$. But the invertibility of~\eqref{EqSpBPfwtP} implies
\[
  \ind\wh{P_0}(\sigma) = \dim\cN-\dim\cR = \ind\wh{P_0}(\sigma_0).
\]

\pfstep{Part~\eqref{ItSpBb}.} In view of the microlocal norm equivalence~\eqref{EqSpBNobNorm}, only the requisite outgoing radial point estimate is different compared to the proof of parts~\eqref{ItSpBReal} and \eqref{ItSpBImag}: we need to use the version of Proposition~\ref{PropSpBInftyOut} with b-regularity, i.e., the (strong) estimate~\eqref{EqSpBInftyOutb}.

\pfstep{Part~\eqref{ItSpBNull}.} Given $u\in\bar H_\scop^{\sfs_\scop,\sfr+\alpha_+}\cap\ker\wh{P_0}(\sigma)$, we first obtain $u\in\bar H_{\scop;\bop}^{(\sfs_\scop;\infty),\sfr+\alpha_++1}$ and thus $u\in\bar H_\bop^{\infty,r}$ for all $r<-\frac12+\ubar S$ (cf.\ \eqref{EqSpBThrOut}). Sobolev embedding (see~\eqref{EqMUCbSob}) gives $u\in\cC_\bop^{\infty,1+\ubar S-\eps}(X)$ for all $\eps>0$. We then rewrite the equation for $u$ as
\begin{equation}
\label{EqSpBPfImpr}
  \Bigl(\rho\pa_\rho-1-S\Bigr)u = -\frac{1}{2 i\sigma\rho}\bigl(\wh{P_0}(0)-i\sigma Q+g^{0 0}\sigma^2\bigr)u.
\end{equation}
The right-hand side lies in $\cC_\bop^{\infty,2+\ubar S-\eps}$. Integrating this ODE from $\rho=\delta>0$ towards $\rho=0$ gives $u=\rho^{1+S}u_0+\tilde u$ where $u_0\in\CI(\pa X)$ and $\tilde u\in\cC_\bop^{2+\ubar S-\eps}$. The right-hand side of~\eqref{EqSpBPfImpr} thus lies in $\rho^{1+S+1}\CI(\pa X)+\cC_\bop^{3+\ubar S-\eps}$. When no two elements of $\spec S$ differ by an integer, integration yields $u\in\sum_{j=0}^1\rho^{1+S+j}\CI(\pa X)+\cC_\bop^{3+\ubar S-\eps}$; continuing in this fashion completes the proof of~\eqref{EqSpBNullspace}.

For the parenthetical statement of part~\eqref{ItSpBNull}, let $u\in\bigcap_{\eps>0}\cC_\bop^{\infty,1+\ubar S-\eps}(X)=\bigcap_{\eps>0}\bar H_\bop^{\infty,1+\ubar S-\eps}(X,\mu_\bop)$ (where we use the unweighted b-density $\mu_\bop:=\frac{|\dd x|}{|x|^3}$ to define $L^2$), which in turn implies $u\in\bar H_\scop^{\infty,\sfr'}(X,|\dd x|)$ for all $\sfr'\in\CI(\ol{\Tsc^*_{\pa X}}X)$ such that $(1+\ubar S)-\frac32>\sup_{\cR_{\rm out}}\sfr'$, essentially due to the observation in Remark~\ref{RmkSpBNob}; see Lemma~\ref{LemmaDRes} below for the detailed statement and proof (which can be read independently). This implies the claim in view of~\eqref{EqSpBThrOut}.

\pfstep{Part~\eqref{ItSpBInvReg}.} Let us first \emph{assume} the differentiability of~\eqref{EqSpBInvReg0}. For later purposes, it is more convenient to study derivatives of the resolvent along $\sigma\pa_\sigma$; note then that~\eqref{EqSpBInftyOp}, now again in the Kerr setting for definiteness, gives
\begin{equation}
\label{EqSpBInvRegDiff}
  \sigma\pa_\sigma\wh{P_0}(\sigma) = 2 i\sigma\rho(\rho\pa_\rho-1-S) - i\sigma Q + 2\sigma^2 g^{0 0},
\end{equation}
which thus lies in $\rho\Diff_\bop^1$. The reason behind the validity of~\eqref{EqSpBInvReg} is as follows: consider
\[
  \sigma\pa_\sigma\wh{P_0}(\sigma)^{-1}=-\wh{P_0}(\sigma)^{-1}\circ\sigma\pa_\sigma\wh{P_0}(\sigma)\circ\wh{P_0}(\sigma)^{-1},
\]
then each application of $\wh{P_0}(\sigma)^{-1}$ \emph{loses} one order of decay (and gaining a scattering derivative), while $\sigma\pa_\sigma\wh{P_0}(\sigma)$ \emph{gains} one order of decay at the expense of one b-derivative. More precisely, for any $\ell\geq 0$ and $q\in\N_0$, we have
\begin{equation}
\label{EqSpBInvRegPf}
  \sigma\pa_\sigma\wh{P_0}(\sigma)\circ\wh{P_0}(\sigma)^{-1} \colon \bar H_{\scop;\bop}^{(\sfs_\scop+\ell-1;q),\sfr+\alpha_++1} \xra{\wh{P_0}(\sigma)^{-1}} H_{\scop;\bop}^{(\sfs_\scop+\ell;q),\sfr+\alpha_+} \xra{\sigma\pa_\sigma\wh{P_0}(\sigma)} H_{\scop;\bop}^{(\sfs_\scop+\ell;q-1),\sfr+\alpha_++1}.
\end{equation}
(The fact that we can take the regularity order to be $\sfs_\scop$ plus a non-negative constant $\ell$ follows from the fact that the only requirements on the regularity order in~\eqref{EqSpBEst} and \eqref{EqSpBbEst} are its monotonicity along the null-bicharacteristic flow and its lower bound at the radial set $\cR_{\cH^+}$ over the event horizon.) Applying~\eqref{EqSpBInvRegPf} with $(\ell,q)=(0,k+j)$, $(1,k+j-1)$, $\ldots$, $(j-1,k+1)$, and applying $\wh{P_0}(\sigma)^{-1}\colon\bar H_{\scop;\bop}^{(\sfs_\scop+j-1;k),\sfr+\alpha_++1}\to\bar H_{\scop;\bop}^{(\sfs_\scop+j;k),\sfr+\alpha_+}$ one last time gives~\eqref{EqSpBInvReg}.

To complete the argument, we need to prove the $j$-fold continuous differentiability of $\wh{P_0}(\sigma)^{-1}$ as a map~\eqref{EqSpBInvReg0} for $\sigma$ near $\sigma_0$. This is closely related to arguments in \cite[\S{2.7}]{VasyMicroKerrdS}. We begin with $j=0$, and first prove the invertibility of $\wh{P_0}(\sigma)$ for $\sigma$ near $\sigma_0$; by part~\eqref{ItSpBFred}, it suffices to prove its injectivity. If this failed, there would exist sequences $\sigma_j\to\sigma_0$ and $u_j\in\cA^{1+\ubar S-\eps}$ such that $\wh{P_0}(\sigma_j)u_j=0$ and $\|u_j\|_{\bar H_\scop^{\sfs_\scop,\sfr_{\sigma_j}+\alpha_+}}=1$, where we make the $\sigma$-dependence of $\sfr$ explicit . By~\eqref{EqSpBEst}, this gives $\|u_j\|_{\bar H_\scop^{-N,-N}}\geq C^{-1}$. A weak subsequential limit $u\in\bar H_\scop^{\sfs_\scop,\sfr_{\sigma_0}-\eps+\alpha_+}$ of $u_j$, thus a norm limit in $\bar H_\scop^{-N,-N}$ (into which $\bar H_\scop^{\sfs_\scop,\sfr_{\sigma_0}-\eps+\alpha_+}$ embeds compactly), has $\bar H_\scop^{-N,-N}$-norm $\geq C^{-1}$ and thus is non-zero; but we have distributional convergence $\wh{P_0}(\sigma_j)u_j\to\wh{P_0}(\sigma_0)u$ along this subsequence, so $\wh{P_0}(\sigma_0)u=0$, contradicting the injectivity of $\wh{P_0}(\sigma_0)$ on~\eqref{EqSpBMap} (with scattering decay order reduced by $\eps$).

Let now $f\in\bar H_\scop^{\sfs_\scop+\eps-1,\sfr+\eps+\alpha_++1}$; we need to show the continuity in $\sigma$ of
\[
  u(\sigma) := \wh{P_0}(\sigma)^{-1}f \in \bar H_\scop^{\sfs_\scop-\eps,\sfr_{\sigma_0}-\eps+\alpha_+}
\]
at $\sigma=\sigma_0$. But $\|u(\sigma)\|_{\bar H_\scop^{\sfs_\scop,\sfr_{\sigma_0}+\alpha_+}}$ is uniformly bounded (since~\eqref{EqSpBEst} without the error term, and with orders increased by $\eps/2$, holds uniformly for $\sigma$ near $\sigma_0$), so for any sequence $\sigma_j\to\sigma_0$ we can pass to a subsequence converging weakly to some $u_0\in\bar H_\scop^{\sfs_\scop,\sfr_{\sigma_0}+\alpha_+}$, and thus in norm in $\bar H_\scop^{\sfs_\scop-\eps,\sfr_{\sigma_0}-\eps+\alpha_+}$; and $\wh{P_0}(\sigma_0)u_0$ equals $f$, being the distributional limit of $\wh{P_0}(\sigma_j)u(\sigma_j)$. Since the solution of $\wh{P_0}(\sigma_0)u=f$ in $\bar H_\scop^{\sfs_\scop-\eps,\sfr_{\sigma_0}-\eps+\alpha_+}$ is unique for all sufficiently small $\eps\geq 0$, we conclude that $u_0$ is independent of the (sub)sequence, and thus $u(\sigma)\to u_0=\wh{P_0}(\sigma_0)^{-1}f=u(\sigma_0)$, as desired.

The differentiability of $\wh{P_0}(\sigma)^{-1}$ in $\sigma$ follows by an inspection of finite difference quotients
\[
  (\sigma-\sigma_0)^{-1}\bigl(\wh{P_0}(\sigma)-\wh{P_0}(\sigma_0)\bigr) = -\wh{P_0}(\sigma)^{-1}\circ\frac{\wh{P_0}(\sigma)-\wh{P_0}(\sigma_0)}{\sigma-\sigma_0}\circ\wh{P_0}(\sigma_0)^{-1}.
\]
The mapping properties of the second factor on the right, which is $\pa_\sigma\wh{P_0}(\sigma_0)+\frac{\sigma-\sigma_0}{2}\pa_\sigma^2\wh{P_0}(\sigma_0)$, with $\pa_\sigma^2\wh{P_0}(\sigma_0)\in\rho^2\CI$, are the same as those for $\pa_\sigma\wh{P_0}(\sigma_0)$ used above; the already established continuity of the resolvent allows one to take the limit $\sigma\to\sigma_0$. Higher regularity then follows by induction on $j$.

The proof of Theorem~\ref{ThmSpB} is complete.

\subsection{High-energy estimates}
\label{SsSpHi}

As in~\eqref{EqSpPhz}, we write $P_{h,z}=h^2\wh{P_0}(h^{-1}z)$. We will use semiclassical scattering ps.d.o.s (see~\eqref{EqMUsch}), and recall the notation $H_{\scop,h}^{s,r,b}=\rho^r h^b H_{\scop,h}^s$ for the associated Sobolev spaces from~\eqref{EqMUschSob}. We recall the orders $\sfs_\scop,\sfr_\pm,\sfb_\pm$ from~\eqref{EqSpOrderHis}--\eqref{EqSpOrderHib}. We only need semiclassical versions of Theorem~\ref{ThmSpB}\eqref{ItSpBReal}, \eqref{ItSpBImag}, \eqref{ItSpBb}, and \eqref{ItSpBInvReg}. We use the notation
\begin{equation}
\label{EqSpHiNormb}
  \|u\|_{\bar H_{(\scop,h);\bop^+}^{(\sfs;k),\sfr,\sfb}} := \sum_{j=0}^k h^{-j}\|u\|_{\bar H_{(\scop,h);\bop}^{(\sfs;k-j),\sfr,\sfb}},
\end{equation}
in which, thus, powers of $h^{-1}$ are considered part of $k$ orders of ``extended'' b-regularity. (On the spacetime side, powers of $h^{-1}=|\sigma|$ amount to regularity in $t_*$; see also Lemma~\ref{LemmaMUetbFTb}.)

\begin{thm}[Spectral family at high frequencies]
\fakephantomsection
\label{ThmSpHi}
  \begin{enumerate}
  \item\label{ItSpHiI}{\rm (High-energy estimates, I).} Let $\nu\geq 0$. Then there exist constants $C$ and $h_0>0$ such that for all $z\in\pm 1+i[0,\nu]$ and all $h\in(0,h_0)$, we have
    \begin{equation}
    \label{EqSpHi}
    \begin{split}
      \|u\|_{\bar H_{\scop,h}^{\sfs_\scop,\sfr+\alpha_+,\sfb}} &\leq C h^{-2}\|P_{h,z}u\|_{\bar H_{\scop,h}^{\sfs_\scop-1,\sfr+\alpha_++1,\sfb}}, \\
      \|u^*\|_{\dot H_{\scop,h}^{-\sfs_\scop+1,-\sfr-\alpha_+-1,-\sfb}} &\leq C h^{-2}\|P_{h,z}^*u^*\|_{\dot H_{\scop,h}^{-\sfs_\scop,-\sfr-\alpha_+,-\sfb}},
    \end{split}
    \end{equation}
    where $\sfr=\sfr_\pm$ and $\sfb=\sfb_\pm$.
  \item\label{ItSpHiII}{\rm (High-energy estimates, II.)} Let $\nu'>0$. Then there exist constants $C$ and $h_0>0$ such that for all $z\in\C$, $|z|\in[\frac12,2]$, with $\Im z\geq\nu'$ and all $h\in(0,h_0)$, the estimates~\eqref{EqSpHi} hold for all order functions $\sfr\in\CI(\ol{\Tsc^*_{\pa X}}X)$ such that $\sfr+\alpha_+<-\frac12+\ubar S$ at $\cR_{\rm out}$, and constant $\sfb$.\footnote{It suffices to assume that $\sfb\in\CI(\ol{\Tsc^*}X)$, and $\sfb$ is monotonically decreasing along the future null-bicharacteristic flow at fiber infinity and constant near $\pa\cR_{\cH^+}$ and near $r=\bhm$. Since we use large $\Im(h^{-1}z)$ estimates only for Paley--Wiener purposes, flexibility in the choices of $\sfr,\sfb$ is not important here, however.}
  \item\label{ItSpHiInv}{\rm (Invertibility.)} There exists $C_0>0$ such that for all $\sigma\in\C$ with $\Im\sigma\geq 0$ and $|\sigma|\geq C_0$, the operator $\wh{P_0}(\sigma)$ is invertible as a map~\eqref{EqSpBMap}.
  \item\label{ItSpHib}{\rm (Smaller $h$-losses; higher b-regularity.)} For $z,h_0$ as in parts~\eqref{ItSpHiI}--\eqref{ItSpHiII}, and for all $k\in\N_0$ and $\eps>0$, there exists a constant $C$ such that
    \begin{equation}
    \label{EqSpHib}
      \|u\|_{\tilde H_{(\scop,h);\bop^+}^{(\sfs_\scop;k),\sfr+\alpha_+,\sfb}} \leq C h^{-1-\delta_\Gamma(\Im(h^{-1}z)-\eps)}\|P_{h,z}u\|_{\tilde H_{(\scop,h);\bop^+}^{(\sfs_\scop-1;k),\sfr+\alpha_++1,\sfb}},
    \end{equation}
    where the norms are those on the space~\eqref{EqMUHscbh}, and the function $\delta_\Gamma$ (which takes values in $[0,1]$) is defined in~\eqref{EqSpHiTrLoss}.\footnote{Thus, they control a weighted semiclassical scattering norm of a distribution together with up to $k$-many of its non-semiclassical b-derivatives.}
  \item\label{ItSpHiInvReg}{\rm (Regularity in the spectral parameter.)} For $\sigma$ as in part~\eqref{ItSpHiInv}, the resolvent $\wh{P_0}(\sigma)^{-1}$ is $j$ times differentiable along $\sigma\pa_\sigma$ for all $j\in\N_0$, with the derivative being uniformly bounded as a map
    \begin{equation}
    \label{EqSpHiInvReg}
      (\sigma\pa_\sigma)^j\wh{P_0}(\sigma)^{-1} \colon \bar H_{(\scop,|\sigma|^{-1});\bop^+}^{(\sfs_\scop-1;k+j),\sfr+\alpha_++1,\sfb} \to |\sigma|^{-1+(j+1)\delta_\Gamma(\Im(\sigma)-\eps)} \bar H_{(\scop,|\sigma|^{-1});\bop^+}^{(\sfs_\scop+j;k),\sfr+\alpha_+,\sfb}
    \end{equation}
    for all $\eps>0$.
  \end{enumerate}
\end{thm}

\begin{rmk}[b-regularity]
\label{RmkSpHib}
  The estimate~\eqref{EqSpHib} is somewhat unusual in two regards. First, it mixes semiclassical scattering and classical b-regularity; but this is what arises upon Fourier transforming $H_{\tbop;\bop}$-spaces (cf.\ Lemma~\ref{LemmaMUetbFTb}). Second (and relatedly), every b-derivative results in an additional maximal loss of $h^{-1}$. On the inverse Fourier transform/spacetime side, this means that control of $k$ spatial b-derivatives of $u$ requires control of (roughly speaking) $k-j$ spatial b-derivatives and $j$ $t_*$-derivatives of $P_0 u$, which is weaker than $k$ orders of spacetime b-regularity (and thus acceptable in the context of estimates of the form~\eqref{EqSSAlephAdmSol}). See the proof of Theorem~\ref{ThmA1Adm} for details in a concrete setting. The estimates~\eqref{EqSpHib} are, of course, equivalent to purely semiclassical estimates (i.e., involving \emph{semiclassical} b-regularity) without such additional losses, to wit,
  \[
    \|u\|_{\bar H_{(\scop;\bop),h}^{(\sfs_\scop;k),\sfr+\alpha_+,\sfb}}\leq C h^{-2}\|P_{h,z}u\|_{\bar H_{(\scop;\bop),h}^{(\sfs_\scop-1;k),\sfr+\alpha_++1,\sfb}}.
  \]
\end{rmk}

\begin{rmk}[Related literature]
\fakephantomsection
\label{RmkSpHiLit}
  \begin{enumerate}[label=(\roman*)]
  \item Semiclassical radial point estimates near the event horizon were first proved in \cite[Proposition~2.10]{VasyMicroKerrdS} for $|\Im z|=\cO(h)$ and in \cite[Proposition~7.1]{VasyMicroKerrdS} for $h\ll\Im z\lesssim 1$. High-energy scattering estimates (i.e., estimates near the incoming and outgoing radial sets over $\pa X$) are due to Vasy--Zworski \cite{VasyZworskiScl} for ``unconjugated'' spectral families (such as $h^2\Delta-z^2$ on Euclidean space) at real energies $z=\pm 1$. The corresponding estimates for the ``conjugated'' spectral family, as studied here, are more delicate when $z$ becomes complex; as far as the author knows, the results we prove here (Propositions~\ref{PropSpHiInftyIn} and \ref{PropSpHiInftyOut}) have not appeared in this form before, though the arguments are strongly inspired by \cite{VasyLAPLag}. The main difference to the case of bounded frequencies, however, is the need to deal with trapping. We recall that the semiclassical estimates proved by Wunsch--Zworski \cite{WunschZworskiNormHypResolvent} and Dyatlov \cite{DyatlovSpectralGaps} are microlocalized near the trapped set by placing complex absorption away from the trapped set. (The reduction of general classes of frequency-dependent operators to the particular form studied in \cite{WunschZworskiNormHypResolvent,DyatlovSpectralGaps} is performed in \cite[\S{4}]{DyatlovResonanceProjectors}.) A more convenient formulation as a microlocal propagation estimate was given in \cite[Theorem~4.7]{HintzVasyQuasilinearKdS}. We give a self-contained proof of such a microlocal estimate in~\S\ref{SssSpHiTr} below, and in fact prove a sharpening as far as the semiclassical loss is concerned (by means of a complex interpolation argument).
  \item A version of the estimates~\eqref{EqSpHi} and \eqref{EqSpHib} on second microlocal Sobolev spaces was proved by Vasy \cite{VasyLAPLag}.
  \item Semiclassical energy estimates for $\Im z=\cO(h)$ are standard (see, e.g., \cite[Proposition~3.7]{VasyMicroKerrdS}, \cite[Theorem~9.44]{HintzMicro}, or \cite[Appendix~E.5.2]{DyatlovZworskiBook}). We present the details for $h\lesssim\Im z\lesssim 1$ in~\S\ref{SssSpHiEn}.
  \end{enumerate}
\end{rmk}

We first prove real principal type propagation estimates in~\S\ref{SssSpHiPr} before turning to radial sets in~\S\ref{SssSpHiHor} (near the event horizon) and \S\ref{SssSpHiInfty} (near spatial infinity); trapping is discussed in~\S\ref{SssSpHiTr}. Energy estimates are discussed in~\S\ref{SssSpHiEn}, and the proof of Theorem~\ref{ThmSpHi} is given in~\S\ref{SssSpHiPf}.

\subsubsection{Real principal type propagation estimates and complex absorption}
\label{SssSpHiPr}

As a warm-up, we discuss the propagation of semiclassical regularity by separately considering each of the three cases described in~\S\ref{SssSpSemi}. We work with a general stationary wave operator $P_0$ on $\cM=\R_t\times\cX$ as in Definition~\ref{DefSpGen}, and recall $P_{h,z}=h^2\wh{P_0}(h^{-1}z)$; here $z\in\C$, $|z|\in[\frac12,2]$.

For $z=\pm 1+\cO(h)$, the characteristic set of $(P_{h,z})_{h\in(0,1)}\in\Diff^2_\semi(\cX)$ is equal to
\[
  \Sigma_{\pm,\semi} := \{ \rho_\infty^2 G_{\pm 1}=0 \} \cap \ol{T^*}\cX,\quad G_z:=G(-z\,\dd t+\cdot) \in P^2(T^*\cX),
\]
where $\rho_\infty\in S^{-1}(T^*\cX)$ is elliptic and positive; and we call integral curves of $\rho_\infty H_{G_{\pm 1}}$ inside of $\Sigma_{\pm,\semi}$ null-bicharacteristics. The future/past component of this is denoted $\Sigma_{\pm,\semi}^\pm$. The zero set of $\rho_\infty^2 G(-z\,\dd t+\cdot)$ over $S^*\cX$ is independent of $z$: it consists of the points in $S^*\cX$ in the direction of $\hat\xi\in T^*\cX\setminus o$ with $G(\hat\xi)=0$ (since $\rho_\infty^2 G(-z\,\dd t+\rho_\infty^{-1}\hat\xi)=G(-z\rho_\infty\,\dd t+\hat\xi)\to G(\hat\xi)$ as $\rho_\infty\to 0$). We denote it by
\begin{equation}
\label{EqSpHiChar}
  \pa\Sigma_\semi \subset S^*\cX,\quad \pa\Sigma_\semi = \pa\Sigma_\semi^+ \sqcup \pa\Sigma_\semi^-.
\end{equation}

\begin{prop}[Propagation of regularity, I: $|\Im z|=\cO(h)$]
\label{PropSpHiPr1}
  Let $z=\pm 1+\cO(h)$, and let $\gamma\colon[0,1]\to\Sigma_{\pm,\semi}\subset\ol{T^*}\cM$ be a null-bicharacteristic of $(P_{h,z})_{h\in(0,1)}$. Fix any neighborhoods $\cU_0$ of $\gamma(0)$ and $\cU$ of $\gamma([0,1])$, and fix a compact neighborhood $K$ of the base projection of $\gamma$. Let $\chi\in\CIc(\cX)$ be equal to $1$ near $K$. Then there exist $h_0>0$ and operators $B,E,G\in\Psi_\semi^0(\cX)$ with Schwartz kernels supported in $K\times K$ and operator wave front sets in $\cU$ such that $\WF_\semi'(E)\subset\cU_0$, $\WF_\semi'(B)\subset\cU$, and $B$ is elliptic on $\gamma([0,1])$, such that, for all $s,N\in\R$,
  \begin{equation}
  \label{EqSpHiPr1}
    \|B u\|_{H_h^s} \leq C\Bigl( h^{-1}\|G P_{h,z}u\|_{H_h^{s-1}} + \|E u\|_{H_h^{s-1}} + h^N\|\chi u\|_{H_h^{-N}} \Bigr),\quad h\in(0,h_0).
  \end{equation}
  The same estimate holds with $P_{h,z}^*$ in place of $P_{h,z}$.
\end{prop}
\begin{proof}
  Since $(P_{h,z})_{h\in(0,1)}$ has real principal symbol $G_{\pm 1}$, this is the standard propagation of regularity, as proved, e.g., in \cite[\S{8.7}]{HintzMicro} and \cite[Appendix~E.4]{DyatlovZworskiBook}.
\end{proof}

As a preparation for the case $\Im z\gg h$, we discuss the real and imaginary parts of the semiclassical principal symbol (Lemma~\ref{LemmaSpChar}) and the imaginary part of $P_{h,z}=\Re P_{h,z}+i\Im P_{h,z}$ (Lemma~\ref{LemmaSpImag}).

\begin{lemma}[Dual metric function]
\label{LemmaSpChar}
  For $z\in\C$ and $\xi\in T^*\cX$, the decomposition of $G(-z\,\dd t+\xi)$ into real and imaginary parts is
  \begin{equation}
  \label{EqSpChar}
    G(-z\,\dd t+\xi) = \Bigl( G\bigl(-(\Re z)\,\dd t+\xi\bigr) - (\Im z)^2 G(\dd t) \Bigr) - 2 i\Im(z) g^{-1}\bigl(-(\Re z)\,\dd t+\xi,\dd t\bigr).
  \end{equation}
  Fix a positive elliptic symbol $\rho_\infty\in S^{-1}(T^*\cX)$, and set $\Sigma_z:=\ol{T^*}\cX\cap(\rho_\infty^2 G(-z\,\dd t+\cdot))^{-1}(0)$. For $z=\pm 1$, write $\Sigma_z^\pm\subset\Sigma_z$ for the future (``$+$''), resp.\ past (``$-$'')  component of the characteristic set.\footnote{For $z=\pm 1$, the set $\Sigma_z^\mp$ is nonempty only over points in $\cX$ where $\pa_t$ is not timelike and thus $T^*\cX$ is not spacelike (so null or Lorentzian).} Then:
  \begin{enumerate}
  \item\label{ItSpChar1} If $\Im z>0$, then $\Sigma\subset S^*\cX$.
  \item\label{ItSpChar2} For $z=\pm 1$, the quantity $\rho_\infty g^{-1}(-(\Re z)\,\dd t+\xi,\dd t)$ is positive, resp.\ negative for $\xi\in\Sigma_z^+$, resp.\ $\xi\in\Sigma_z^-$ (including at fiber infinity).
  \end{enumerate}
\end{lemma}
\begin{proof}
  For part~\eqref{ItSpChar1}, we need to show that $G(-z\,\dd t+\xi)\neq 0$ for all $\xi\in T^*\cX$. Now, $\Re G(-z\,\dd t+\xi)=0$ implies
  \[
    G\bigl(-(\Re z)\,\dd t+\xi\bigr) = (\Im z)^2 G(\dd t) < 0.
  \]
  Therefore, $-(\Re z)\,\dd t+\xi$ is timelike, and hence its inner product with the timelike covector $\dd t$ is nonzero; so $\Im G(-z\,\dd t+\xi)\neq 0$, as desired. Part~\eqref{ItSpChar2} follows for finite $\xi$ from the past timelike nature of $\dd t$. At fiber infinity, the claim is equivalent to $g^{-1}(\xi,\dd t)>0$ when $\xi\in T^*\cX$ is a unit covector (with respect to any norm on the fiber) lying in the future characteristic set of $P_0$ intersected with $T^*\cX\subset T^*_{t^{-1}(0)}\cM$; this, too, follows from the past timelike nature of $\dd t$.
\end{proof}

\begin{lemma}[Imaginary part of $P_{h,z}$]
\label{LemmaSpImag}
  We have
  \begin{equation}
  \label{EqSpImag}
  \begin{split}
    &h^{-1}\Im P_{h,z} = \bigl(h\Im\wh{P_0}(0) + \Re(z)Q_0\bigr) - h^{-1}\Im(z) Q_{1,h}, \\
    &\qquad Q_0:=\Im(\pa_\sigma\wh{P_0}(0))\in\Diff^0(\cX),\quad Q_{1,h}=-\Bigl(h\Re\bigl(\pa_\sigma\wh{P_0}(0)\bigr) + 2\Re(z)G(\dd t)\Bigr),
  \end{split}
  \end{equation}
  here $Q_0\in\Diff^0(\cX)$ is real, $h\Im\wh{P_0}(0)\in\Diff_\semi^1(\cX)$, and $(Q_{1,h})_{h\in(0,1)}\in\Diff_\semi^1(\cX)$ has semiclassical principal symbol given by $T^*\cX\ni\xi\mapsto 2 g^{-1}(-(\Re z)\,\dd t+\xi,\dd t)$.
\end{lemma}

The first term is (symmetric and) of class $\Diff_\semi^1(\cX)$, and is thus harmless in positive commutator arguments (at radial points, it would affect the threshold conditions). The term $Q_{1,h}$ multiplying $h^{-1}\Im(z)$ (which thus dominates when $h^{-1}\Im(z)\gg 1$) is of class $\Diff_\semi^1(\cX)$ and symmetric. Its semiclassical principal symbol can be read off from~\eqref{EqSpChar}; on $\Sigma_{\pm,\semi}^+$, resp.\ $\Sigma_{\pm,\semi}^-$, it is positive, resp.\ negative.

\begin{proof}[Proof of Lemma~\usref{LemmaSpImag}]
  Write $\wh{P_0}(\sigma)=\wh{P_0}(0)+\sigma\pa_\sigma\wh{P_0}(0)+\frac{\sigma^2}{2}\pa_\sigma^2\wh{P_0}(0)$. The operator $\pa_\sigma\wh{P_0}(0)\in\Diff^1(\cX)$ has the real principal symbol $(-H_G t)|_{T^*\cX}$ (as follows from \eqref{EqSbPpasigma}), which is thus equal to the principal symbol of
  \[
    Q_1:=\Re(\pa_\sigma\wh{P_0}(0))\in\Diff^1(\cX).
  \]
  But $-H_G t=\pa_\sigma G$, which in view of
  \[
    G(-\sigma\,\dd t+\xi\,\dd x)=\sigma^2 g^{-1}(\dd t,\dd t)-2\sigma g^{-1}(\dd t,\xi\,\dd x)+g^{-1}(\xi\,\dd x,\xi\,\dd x)
  \]
  evaluates at $T^*\cX=\{\sigma=0\}$ to $-2 g^{-1}(\dd t,\cdot)$. Write
  \[
    \Im\bigl(\sigma \pa_\sigma\wh{P_0}(0)\bigr) = \Im(\sigma)Q_1 + \Re(\sigma)Q_0.
  \]
  We combine this with the term $\pa_\sigma^2\wh{P_0}(0)$, which is a zeroth order operator and equal to $\pa_\sigma^2 G=2 G(\dd t)$. Using that $\Im(\sigma^2)=2\Im(\sigma)\Re(\sigma)$, we then combine the contributions from $\Im(\sigma)Q_1$ and $\frac{\sigma^2}{2}\pa_\sigma^2\wh{P_0}(0)$ to $h^{-1}\Im P_{h,z}$, obtaining
  \[
    h^{-1} h^2\bigl( \Im(h^{-1}z)Q_1 + \Im(h^{-1}z)\Re(h^{-1}z) G(\dd t) \bigr) = h^{-1}\Im(z) \bigl( h Q_1 + 2\Re(z) G(\dd t) \bigr).
  \]
  It then remains to observe that the semiclassical principal symbol of the term in parentheses on the right is $-2 g^{-1}(\dd t,-(\Re z)\,\dd t+\cdot)$.
\end{proof}

\begin{prop}[Propagation of regularity, II: $h\ll\Im z\ll 1$]
\label{PropSpHiPr2}
  Let $\nu\in\R$. Let $\gamma\colon[0,1]\to\Sigma_{\pm,\semi}\subset\ol{T^*}\cM$ be a \emph{future} null-bicharacteristic of $(P_{h,\pm 1})_{h\in(0,1)}$; that is, on $\Sigma_{\pm,\semi}^+$, resp.\ $\Sigma_{\pm,\semi}^-$, $\gamma$ is an integral curve of $\rho_\infty H_{G_{+1}}$, resp.\ $-\rho_\infty H_{G_{-1}}$. Then for $\cU_0,\cU,K,\chi$ as in Proposition~\usref{PropSpHiPr1}, there exist $c,h_0>0$ and $B,E,G\in\Psi_\semi^0(\cX)$ as in Proposition~\usref{PropSpHiPr1} such that the estimate~\eqref{EqSpHiPr1} holds for all $z\in\C$ with $\Re z=\pm 1$ and $h\nu\leq\Im z\leq c$, $h\in(0,h_0)$. The same holds for $P_{h,z}^*$ in place of $P_{h,z}$ if we replace ``future null-bicharacteristic'' by ``past null-bicharacteristic.''
\end{prop}
\begin{proof}
  For any fixed $\nu'>\nu$, this follows from Proposition~\ref{PropSpHiPr1} when $\Im z\in[\nu,\nu']$. We may therefore arbitrarily increase $\nu$ in the course of the proof, as long as $\nu$ remains independent of $h$. Furthermore, we only discuss the case that $\gamma\subset\Sigma_{\pm,\semi}^+$ lies in the future characteristic set.

  The proof is again a positive commutator argument, starting as around~\eqref{EqSbPHorsC}: for $A=\check A^*\check A\in\Psih^{2 s-1}(\cX)$, with $\check A=\Op_\semi(\check a)=\check A^*$ a symmetric semiclassical quantization of a suitable symbol $\check a\in S^{s-\frac12}(T^*\cX)$, we consider the $L^2$-pairing
  \begin{subequations}
  \begin{equation}
  \label{EqSpHiPr2Pair}
    2 h^{-1}\Im\la\check A P_{h,z}u,\check A u\ra = \la\sC u,u\ra,
  \end{equation}
  where, analogously to~\eqref{EqSpBInftyC}, we now write
  \begin{align}
    \sC &= \frac{i}{h}[\Re P_{h,z},A] - \frac{1}{h}[\Im P_{h,z},\check A^2] + 2 h^{-1}(\Im P_{h,z})\check A^2 \nonumber\\
  \label{EqSpHiPr2C}
      &= \frac{i}{h}[\Re P_{h,z},A] + 2 h^{-1}\check A(\Im P_{h,z})\check A + h^{-1}[[\Im P_{h,z},\check A],\check A].
  \end{align}
  \end{subequations}
  The first term lies in $\Psi_\semi^{2 s}$ and has principal symbol given by $H_{\upsigma_\semi^2(\Re P_{h,z})}\check a^2$. Now, the principal symbols $\upsigma_\semi^2(\Re P_{h,z})$ and $\upsigma_\semi^2(\Re P_{h,\pm 1})=G_{\pm 1}$ differ by $(\Im z)^2\CI(\cX)$ (see~\eqref{EqSpChar}), and thus their characteristic sets and null-bicharacteristics are close when $\Im z$ is small. For appropriate (standard) choices of $\check a$ (see \cite[(8.104a) and (8.104b)]{HintzMicro}), we can thus arrange that
  \begin{equation}
  \label{EqSpHiPr2Symb}
    \upsigma_\semi^{2 s}\Bigl(\frac{i}{h}[\Re P_{h,z},A]\Bigr) = -\rho_\infty^{-1}\check a^2 - b^2 + e
  \end{equation}
  where $b\in S^s$ is elliptic on $\gamma([0,1])$ and $e\in S^{2 s}$ has essential support contained in a neighborhood of $\gamma(0)$. The third (double commutator) term in~\eqref{EqSpHiPr2C} lies in $h\Psi_\semi^{2 s-1}$, and hence is subprincipal.

  For the second term in~\eqref{EqSpHiPr2C}, we use Lemma~\ref{LemmaSpImag}, which we rephrase as
  \begin{equation}
  \label{EqSpHiPr2Gard}
  \begin{split}
    &h^{-1}\Im P_{h,z} = -h^{-1}\Im(z)\tilde Q_{1,z,h}, \\
    &\qquad \tilde Q_{1,z,h}=Q_{1,h} - \Bigl(\frac{\Im z}{h}\Bigr)^{-1}\bigl(h\Im\wh{P_0}(0)+\Re(z)Q_0\bigr) \in \Diff_\semi^1(\cX).
  \end{split}
  \end{equation}
  The point is that since $Q_{1,h}$ has a positive principal symbol near $\gamma([0,1])$, the same is true for $\tilde Q_{1,z,h}$ when $\frac{\Im z}{h}\geq\nu$ is sufficiently large, so $\upsigma_\semi^1((\tilde Q_{1,z,h})_{h\in(0,1)})\geq 2 c\rho_\infty^{-1}$, for some $c>0$, in a fixed neighborhood of $\gamma([0,1])$. Thus, when $\WF'_\semi(\check A)$ is contained in this neighborhood, then by the G\aa{}rding inequality, we can estimate
  \begin{equation}
  \label{EqSpHiPr2Gard2}
  \begin{split}
    \la h^{-1}\check A(\Im P_{h,z})\check A u,u\ra &= -\la h^{-1}\Im(z)\tilde Q_{1,z,h}\check A u,\check A u\ra \\
      &\leq -c h^{-1}\Im(z)\|\check A u\|_{L^2}^2 + C_N h^N\|\chi u\|_{H_h^{-N}}^2,
  \end{split}
  \end{equation}
  which (up to the trivial error term involving $\chi u$) thus has the same (negative) sign as the main term $-b^2$ in~\eqref{EqSpHiPr2Symb}. The proof is thus completed in the usual fashion by quantizing~\eqref{EqSpHiPr2Symb}.
\end{proof}

Finally, we study $z$ with $\Im z$ bounded away from $0$. Recall from~\eqref{EqSpHiChar} that the characteristic set $\pa\Sigma_\semi$ of $(P_{h,z})_{h\in(0,1)}$ now lies over $S^*\cX$ and is independent of $z$; likewise, the null-bicharacteristic flow does not depend on $z$.

\begin{prop}[Propagation of regularity, III: $\Im z\gtrsim 1$]
\label{PropSpHiPr3}
  Consider $z\in\C$ with $|z|\in[\frac12,2]$ and $\Im z\geq c>0$. Let $\gamma\colon[0,1]\to\pa\Sigma_\semi^\pm$ be a future, resp.\ past null-bicharacteristic when considering $P_{h,z}$, resp.\ $P_{h,z}^*$. Then the conclusions of Proposition~\usref{PropSpHiPr2} hold.
\end{prop}
\begin{proof}
  We again start with~\eqref{EqSpHiPr2Pair}--\eqref{EqSpHiPr2C}. The second term on the right-hand side of~\eqref{EqSpHiPr2C} has a sign by~\eqref{EqSpHiPr2Gard}--\eqref{EqSpHiPr2Gard2}, and the third term is subprincipal. We thus only need to revisit the first term. The point is that the semiclassical principal symbol at fiber infinity, $(\rho_\infty^2\upsigma_\semi^2(\Re P_{h,z}))|_{S^*\cX}$, is equal to $\rho_\infty^2 G(\cdot)$ and thus independent of $z$. (Only the subprincipal part depends on $\Im z$.) Thus, we can again make a standard choice of $\check a$ (e.g., by straightening out the Hamiltonian flows of $\Re P_{h,z}$ for all $z$ of present interest and taking~\cite[(8.104b)]{HintzMicro}) in order to arrange for~\eqref{EqSpHiPr2Symb}, and then conclude the proof via quantization as usual.
\end{proof}

\begin{rmk}[Propagation along fiber infinity]
\label{RmkSpHiComb}
  For propagation along null-bicharacteristics in $S^*\cX$, which are independent of $z$, the proof of Proposition~\ref{PropSpHiPr2} applies directly also to all values of $z$ considered in Proposition~\ref{PropSpHiPr3}.
\end{rmk}

\subsubsection{Radial point estimates near horizons}
\label{SssSpHiHor}

We return to the general setting of~\S\ref{SssSpBHor} (see Definition~\ref{DefSpGen} and items~\eqref{ItSpBHor1}--\eqref{ItSpBHor4} in~\S\ref{SssSpBHor}) and prove a semiclassical version of Proposition~\ref{PropSpBHor}. Such estimates were first discussed in \cite[\S{7.1}]{VasyMicroKerrdS}.

\begin{prop}[Semiclassical radial point estimate]
\label{PropSpHiHor}
  We use the notation of Proposition~\usref{PropSpBHor}. Let $\nu\in\R$.
  \begin{enumerate}
  \item{\rm (Direct estimate.)} Let $s>s_0>\frac12-\beta\nu+\beta_1$. There exist $h_0>0$ and semiclassical ps.d.o.s $B,G\in\Psi_\semi^0(\cX)$, with Schwartz kernels supported in $K\times K$ and operator wave front sets in $\cU$, and with $B$ semiclassically elliptic at $\pa\cR^+\subset S^*\cX$, such that for all $z\in\C$ with $|z|\in[\frac12,2]$ and $\Im z\geq h\nu$, we have the estimate
    \begin{equation}
    \label{EqSpHiHor}
      \|B u\|_{H_h^s} \leq C\Bigl( h^{-1}\|G P_{h,z}u\|_{H_h^{s-1}} + h^N\|\chi u\|_{H^{s_0}}\Bigr),\quad h\in(0,h_0).
    \end{equation}
  \item{\rm (Adjoint estimate.)} Let $s>\frac12-\beta\nu+\beta_1$ and $-N<-s+1$. Then there exist $h_0>0$ and semiclassical ps.d.o.s $B,G,E\in\Psi_\semi^0(\cX)$, with Schwartz kernels supported in $K\times K$ and operator wave front sets in $\cU$, and with $B$ elliptic at $\pa\cR^+$ and $\WF_\semi'(E)\cap\pa\cR^+=\emptyset$, such that
    \[
      \|B u\|_{H_h^{-s+1}} \leq C\Bigl( h^{-1}\|G P_{h,z}^*u\|_{H_h^{-s}} + \|E u\|_{H_h^{-s+1}} + h^N\|\chi u\|_{H_h^{-N}} \Bigr),\quad h\in(0,h_0).
    \]
  \end{enumerate}
\end{prop}
\begin{proof}
  We only prove the first part. Given any $\nu'>\nu$, we first prove~\eqref{EqSpHiHor} for $\Re z=\pm 1$ and $h\nu\leq\Im z\leq h\nu'$ by a minor modification of the proof of Proposition~\ref{PropSpBHor}: to wit, we now need to localize also near fiber infinity. Thus, recalling the cutoff $\psi$ from~\eqref{EqSpBHorPsi}, we let $a=\check a^2\in S^{2 s-1}$ where
  \[
    \check a = \rho_\infty^{-s+\frac12}\psi_\Sigma\psi_\infty\psi_\cR,\quad \psi_\Sigma:=\psi(\digamma\rho_\infty^2 G_0),\ \psi_\infty:=\psi(\digamma\rho_\infty),\ \psi_\cR:=\psi(\digamma_\cR\fq),
  \]
  where we will choose $\digamma_\cR>1$ large and then $\digamma>1$ large. Consider $\sC$ in~\eqref{EqSpHiPr2Pair}: writing it as
  \[
    \sC = \frac{i}{h}[P_{h,z},A] + 2\frac{P_{h,z}-P_{h,z}^*}{2 i h}A \in \Psi_\semi^{2 s},
  \]
  its semiclassical principal symbol is $H_{G_{\pm 1}}a+2\rho_\infty^{-1}\sfp_1(z)$ in the notation of~\eqref{EqSbPHorp1}. In terms of $\sfH_{G_{\pm 1}}:=\rho_\infty H_{G_{\pm 1}}$, we can write this as
  \begin{align}
  \label{EqSpHiHorExpr}
    \upsigma_\semi^{2 s}(\sC) =& -\delta\rho_\infty^{-1}\check a^2 - b^2 - b_\cR^2 - b_\infty^2 + j G_{\pm 1}, \\
    b &:= \rho_\infty^{-s}\psi_\Sigma\psi_\infty\psi_\cR\bigl[ (2 s-1)\rho_\infty^{-1}\sfH_{G_{\pm 1}}\rho_\infty - \delta - 2\sfp_1(z) \bigr]^{\frac12}, \nonumber\\
    b_\cR &:= \rho_\infty^{-s}\psi_\Sigma\psi_\infty\bigl[ -2(\sfH_{G_{\pm 1}}\fq)\digamma_\cR\psi_\cR\psi'_\cR \bigr]^{\frac12}, \nonumber\\
    b_\infty &:= \rho_\infty^{-s}\psi_\Sigma\psi_\cR\bigl[ -2(\rho_\infty^{-1}\sfH_{G_{\pm 1}}\rho_\infty)\digamma\rho_\infty\psi_\infty\psi'_\infty \bigr]^{\frac12}, \nonumber\\
    j &:= 2\rho_\infty^{-2 s+2}\psi_\cR \digamma(\rho_\infty^{-2}\sfH_{G_{\pm 1}}\rho_\infty^2)\psi_\Sigma\psi'_\Sigma. \nonumber
  \end{align}
  These symbols are well-defined, with $b$ elliptic at $\pa\cR^+$, since $\rho_\infty^{-1}\sfH_{G_{\pm 1}}\rho_\infty>0$ at, and thus near, $\pa\cR^+$. (This quantity is equal to $\rho_\infty^{-1}\sfH_G\rho_\infty$ for either choice of sign, and indeed for $G_z$ in place of $G$ for any $z\in\C$.)

  It remains to prove~\eqref{EqSpHiHor} for $|z|\in[\frac12,2]$ with $h\ll\Im z\leq 2$. (We do not need to separate two cases here corresponding to whether $\Im z\ll 1$ or not; cf.\ Remark~\ref{RmkSpHiComb}.) To this end, we write $\sC$ as in~\eqref{EqSpHiPr2C} and use the positivity of the principal symbol of $\tilde Q_{1,z,h}$ in~\eqref{EqSpHiPr2Gard}; we can borrow some positivity from $\tilde Q_{1,z,h}$ by writing
  \begin{equation}
  \label{EqSpHiHorExpr2}
  \begin{split}
    \sC &\equiv \Bigl[\frac{i}{h}[\Re P_{h,z},A] - 2\Op_\semi(C\rho_\infty^{-1}\check a^2)\Bigr] \\
      &\qquad - 2 h^{-1}\Im(z)\check A\Bigl(\tilde Q_{1,z,h}-\Bigl(\frac{\Im z}{h}\Bigr)^{-1}\Op_\semi(C\rho_\infty^{-1})\Bigr)\check A \bmod h\Psi_\semi^{2 s-1}
  \end{split}
  \end{equation}
  for any fixed constant $C$. (Here, it is convenient to take $\rho_\infty^{-1}$ to be fiber-linear near $\pa\cR^+$, as one may.) In view of the $z$-independence of the principal symbol of $\Re P_{h,z}$ at fiber infinity, the term in square brackets can be expressed as in~\eqref{EqSpHiHorExpr}, with a large positive choice of $C$ shifting the threshold regularity below any desired real number; when $\frac{\Im z}{h}$ is then sufficiently large, the term in the second line of~\eqref{EqSpHiHorExpr2} can then be estimated as in~\eqref{EqSpHiPr2Gard2}.
\end{proof}

\subsubsection{Radial point estimates near infinity}
\label{SssSpHiInfty}

We shall next prove semiclassical analogues of Propositions~\ref{PropSpBInftyIn} and \ref{PropSpBInftyOut}, again in the general setting of a stationary asymptotically Minkowski wave operator $P_0$ (Definition~\ref{DefSpGenMink}). We thus study the semiclassical rescaling
\begin{equation}
\label{EqSpHiInftyOp}
  P_{h,z} = 2 i z \rho h\Bigl(\rho\pa_\rho-\frac{n-1}{2}-S\Bigr) + h^2\wh{P_0}(0) - i h z Q + z^2 g^{0 0}.
\end{equation}
the semiclassical principal symbol of $(P_{h,z})_{h\in(0,1)}\in\Diff_{\scop,\semi}^{2,0}(\cX)$ is
\[
  G_z \equiv 2 z\xi_\schop+\xi_\schop^2+|\eta_\schop|^2 \bmod \rho P^2(\Tsc^*\cX),
\]
where we write semiclassical scattering covectors as $\xi_\schop\frac{\dd r}{h}+\frac{r\eta_\schop}{h}$, $\eta_\schop\in T^*\Sph^2$. For $z=\pm 1+\cO(h)$, this has the two radial sets
\[
  \cR_{\pm 1,{\rm in}} = \{ \xi_\schop=\mp 2,\ \eta_\schop=0 \} \cap \Tsc^*_{\pa\cX}\cX,\quad
  \cR_{\rm out} = \{ \xi_\schop=0,\ \eta_\schop=0 \} \cap \Tsc^*_{\pa\cX}\cX.
\]
When $\Im z>0$, $G_z$ is elliptic over $\pa\cX$ except at the zero section $\cR_{\rm out}$.

We will only explicitly work with semiclassical scattering ps.d.o.s and Sobolev spaces whose semiclassical orders are equal to $0$ (denoted $b$ in~\eqref{EqMUsch}); thus, we write $\Psi_\schop^{s,r}(\cX)=\rho^{-r}\Psi_\schop^s(\cX)$ and $H_{\scop,h}^{s,r}=\rho^r H_\scop^s$.

\begin{prop}[Semiclassical radial point estimate near $\cR_{\pm 1,{\rm in}}$]
\label{PropSpHiInftyIn}
  Let $\nu>0$ and $s,N\in\R$. Define $\vartheta_{\rm in}$ by~\eqref{EqSpBInftyInThr}. Let $\cU\subset\Tsc^*_{\pa\cX}\cX$ be a compact neighborhood of $\cR_{\pm 1,{\rm in}}$, and let $K\subset\cX$ be a compact neighborhood of $\pa\cX$. Let $\chi\in\CIc(\cX)$, with $\chi=1$ near $K$.
  \begin{enumerate}
  \item{\rm (Direct estimate.)} Let $r>r_0>\frac12(-1+\vartheta_{\rm in})$. There exist $h_0>0$ and operators $B,G\in\Psi_{\scop,\semi}^{0,0,0}(\cX)$ with Schwartz kernels supported in $K\times K$ and operator wave front sets in $\cU$, and with $B$ elliptic at $\cR_{\pm,{\rm in}}$, such that
    \begin{equation}
    \label{EqSpHiInftyIn}
      \|B u\|_{H_{\scop,h}^{s,r}} \leq C\Bigl( h^{-1}\|G P_{h,z}u\|_{H_{\scop,h}^{s-2,r+1}} + h^N\|u\|_{H_{\scop,h}^{-N,r_0}}\Bigr),\quad z\in\pm 1+i[0,\nu],\ \ h\in(0,h_0).
    \end{equation}
  \item{\rm (Adjoint estimate.)} Let $r>\frac12(-1+\vartheta_{\rm in})$. There exist $h_0>0$ and operators $B,E,G\in\Psi_\schop^{0,0}(\cX)$ with Schwartz kernels supported in $K\times K$ and operator wave front sets in $\cU$, and with $B$ elliptic at $\cR_{\pm,{\rm in}}$ and $\WF'_\schop(E)\cap\cR_{\pm,{\rm in}}=\emptyset$, such that
    \begin{align*}
      &\|B u\|_{H_{\scop,h}^{-s+2,-r-1}} \leq C\Bigl( h^{-1}\|G P_{h,z}^*u\|_{H_{\scop,h}^{-s,-r}} + \|E u\|_{H_{\scop,h}^{-s+2,-r-1}} + h^N\|u\|_{H_{\scop,h}^{-N,-N}}\Bigr), \\
      &\qquad z\in\pm 1+i[0,\nu],\ \ h\in(0,h_0).
    \end{align*}
  \end{enumerate}
\end{prop}
\begin{proof}
  The proof is very similar to that of Proposition~\ref{PropSpBInftyIn} for $\Re\sigma=\pm 1$. We only consider the direct estimate in the case $\Re\sigma=1$. For a suitable $A=\check A^*\check A$, $\check A=\check A^*\in\Psi_\schop^{-\infty,r+\frac12}$, we now consider~\eqref{EqSpHiPr2Pair}--\eqref{EqSpHiPr2C}. Apart from notational adjustments to reflect the semiclassical nature of our present analysis, there are only two novel aspects. \emph{First}, terms arising from the differentiation of $\psi_\pa$ (in the symbolic analysis of $\sC$) in~\eqref{EqSpBInftyChecka} are no longer subprincipal: they are principal in the semiclassical order sense. Thus, in the first term of~\eqref{EqSbPInftyComm} (with $\Op_{\scop,h}$ in place of $\Op_\scop$), we have an additional term $-b_\pa^2$ where
  \[
    b_\pa := \rho^{-r}\psi_\cR\bigl[ -2(\rho^{-1}\sfH_{G_{+1}}\rho) \digamma\psi_\pa\psi'_\pa \bigr]
  \]
  is well-defined since $\rho^{-1}\sfH_{G_{+1}}\rho=2$ at $\cR_{+1,{\rm in}}$. This term has the same sign as the main term $-b^2$.

  \emph{Second}, there are terms of $P_{h,z}$ that contribute to $h^{-1}\Im P_{h,z}$ at leading \emph{semiclassical} order which only contributed subprincipal (in the sense of scattering decay) terms in the discussion of~\eqref{EqSpBInftyIm}. In the expression~\eqref{EqSpImag}, the term $h\Im\wh{P_0}(0)+\Re(z)Q_0$ lies in $\rho\Diff_\schop^1$, with the principal symbol being $\rho$ times $(\Re z)(\sfp_1-2\Re S)$ (cf.\ \eqref{EqSpBInftyIm1}--\eqref{EqSpBInftyIm2}), and it becomes part of the symbolic computation and leads to the same threshold condition on the decay order as in Proposition~\ref{PropSpBInftyIn}. In the term $-h^{-1}\Im(z)Q_{1,h}$, the semiclassical principal symbol of $(Q_{1,h})_{h\in(0,1)}\in\Diff_\schop^1(\cX)$ at $\cR_{+1,{\rm in}}$ is $2 g^{-1}(-\dd t-2\,\dd r,\dd t)=4>0$; therefore, the G\aa{}rding inequality implies
  \begin{equation}
  \label{EqSpHiInftyInIm}
    \big\la h^{-1}\check A \bigl( -h^{-1}\Im(z)Q_{1,h} \bigr)\check A u,u\big\ra \leq C_N h^N\|\chi u\|_{H_{\scop,h}^{-N,-N}}
  \end{equation}
  for all $N$. (The point is that this is non-positive up to a ``trivial'' error term, similarly to~\eqref{EqSpHiPr2Gard2}.)

  Given these observations, the proof of~\eqref{EqSpHiInftyIn} then proceeds in the standard fashion.
\end{proof}

\begin{prop}[Semiclassical radial point estimate near $\cR_{\rm out}$]
\label{PropSpHiInftyOut}
  Let $s,N\in\R$. Recall the quantity $\ubar S$ from~\eqref{EqSSAdmubarS}. Let $\cU\subset\Tsc^*_{\pa\cX}\cX$ be a compact neighborhood of $\cR_{\rm out}$, and let $K\subset\cX$ be a compact neighborhood of $\pa\cX$. Let $\chi\in\CIc(\cX)$, with $\chi=1$ near $K$.
  \begin{enumerate}
  \item{\rm (Direct estimate.)} Let $r<-\frac12+\ubar S$. There exist operators $B,E,G\in\Psi_\schop^{0,0}(\cX)$ with Schwartz kernels supported in $K\times K$ and operator wave front sets in $\cU$, and with $B$ elliptic at $\cR_{\rm out}$ and $\WF'_\schop(E)\cap\cR_{\rm out}=\emptyset$, such that
    \begin{equation}
    \label{EqSpHiInftyOut}
      \|B u\|_{H_{\scop,h}^{s,r}} \leq C\Bigl( h^{-1}\|G P_{h,z}u\|_{H_{\scop,h}^{s-2,r+1}} + \|E u\|_{H_{\scop,h}^{s,r}} + h^N\|\chi u\|_{H_{\scop,h}^{-N,-N}} \Bigr)
    \end{equation}
    holds for all $z\in\C$ with $|z|\in[\frac12,2]$ and $\Im z\geq 0$. More generally, for every $k\in\N_0$, we have
    \begin{equation}
    \label{EqSpHiInftyOutb}
      \|B u\|_{H_{(\scop,h);\bop^+}^{(s;k),r}} \leq C\Bigl( h^{-1}\|G P_{h,z}u\|_{H_{(\scop,h);\bop^+}^{(s-2;k),r+1}} + \|E u\|_{H_{(\scop,h);\bop^+}^{(s;k),r}} + h^N\|\chi u\|_{H_{\scop,h}^{-N,-N}} \Bigr),
    \end{equation}
    where we use notation as in~\eqref{EqSpHiNormb}.
  \item{\rm (Adjoint estimate.)} Let $r<r_0<-\frac12+\ubar S$. There exist operators $B,G\in\Psi_\schop^{0,0}(\cX)$ with Schwartz kernels supported in $K\times K$ and operator wave front sets in $\cU$, and with $B$ elliptic at $\cR_{\rm out}$, such that
    \[
      \|B u\|_{H_{\scop,h}^{-s+2,-r-1}} \leq C\Bigl( h^{-1}\|G P_{h,z}^*u\|_{H_{\scop,h}^{-s,-r}} + h^N\|\chi u\|_{H_{\scop,h}^{-N,-r_0-1}}\Bigr)
    \]
    holds for all $z\in\C$ with $|z|\in[\frac12,2]$ and $\Im z\geq 0$.
  \end{enumerate}
\end{prop}
\begin{proof}
  \pfstep{Step 1.~Basic estimate.} We only consider the direct estimate. 

  \pfsubstep{(1.1)}{Almost real $z$.} The necessary modifications of the proof of Proposition~\ref{PropSpBInftyOut} are similar to the proof of Proposition~\ref{PropSpHiInftyIn} when $z\in\pm 1+[0,\nu]$ where $\nu>0$. We only indicate the necessary modifications and additional terms one must deal with. We now study
  \begin{equation}
  \label{EqSpHiInftyOutTilde}
    \tilde P_{h,z}:=z^{-1}P_{h,z},
  \end{equation}
  which thus equals $z^{-1}h^2\wh{P_0}(h^{-1}z)=h\tilde P_0(h^{-1}z)$ in the notation of~\eqref{EqSbPInftyOp}. We consider the ``commutator'' $2 h^{-1}\Im\la\check A\tilde P_{h,z}u,\check A u\ra=\la\sC u,u\ra$,
  \[
    \sC = \frac{i}{h}[\Re\tilde P_{h,z},A] + 2\check A h^{-1}(\Im\tilde P_{h,z})\check A + h^{-1}[[\Im\tilde P_{h,z},\check A],\check A],
  \]
  with $\check A=\Op_{\scop,\semi}(\check a)=\check A^*\in\Psi_\schop^{-\infty,r+\frac12}$ where $\check a=\rho^{-r-\frac12}\psi_\pa\psi_\cR$ is given by~\eqref{EqSpBInftyChecka}. Now, differentiation of the cutoff $\psi_\pa$ yields another a priori control term; indeed, this contributes $+e_\pa^2$ to the principal symbol of $\sC$ in the analogue of~\eqref{EqSpBInftyOutSymb} where, in the case of the ``$+$'' sign,
  \[
    e_\pa := \rho^{-r}\psi_\cR\bigl[ 2(\rho^{-1}\sfH_{\Re\tilde G_z}\rho) \digamma\psi_\pa\psi'_\pa \bigr]^{\frac12}.
  \]
  (This is well-defined when $\Im z\leq\nu$ is sufficiently small since $\rho^{-1}\sfH_{G_{+1}}\rho=-2$ at $\cR_{\rm out}$.) The quantizations of the two terms $e_\cR$ and $e_\pa$ (acting on $u$, and then paired with $u$) are controlled by the term $E u$ in~\eqref{EqSpHiInftyOut}.

  We furthermore need to study the semiclassical analogue of~\eqref{EqSpBInftyImTilde}. Using~\eqref{EqSbPInftyOp}, we find
  \begin{align*}
    h^{-1}\Im\tilde P_{h,z} = \Im\tilde P_0(h^{-1}z) &= -2\rho\Re S + (\Re z)|z|^{-2} h\Im\wh{P_0}(0) - (\Im z)|z|^{-2} h^2\Re\wh{P_0}(0) \\
      &\qquad - \Re Q + h^{-1}(\Im z)g^{0 0}.
  \end{align*}
  Completely analogously to~\eqref{EqSpBInftyOutSymb2}, the first term shifts the decay threshold; the second term (of class $\rho\,h\rho\Diff_\bop$, so the principal symbol of $\rho^{-1}$ times it vanishes at $\cR_{\rm out}$) yields a contribution to $b$ in~\eqref{EqSpBInftyOutSymb2} that vanishes at $\cR_{\rm out}$ and thus is arbitrarily small nearby. The term involving $\Re\wh{P_0}(0)$ similarly yields arbitrarily small contributions near $\cR_{\rm out}$ (cf.\ $\sfr_1,\sfr_2$ in~\eqref{EqSpBInftyOutSymb}), plus a term with a good sign arising from the sum in~\eqref{EqSpBInftyReP0}, analogously to~\eqref{EqSpBInftyOutSign}. Next, note that $Q$ has a purely imaginary principal symbol (since the principal symbol of the second-order operator $Q\pa_{t_*}$ in Definition~\ref{DefSpGenMink} is real), and thus $\Re Q\in\rho^3$; therefore, $\Re Q$ yields a further (semiclassically principal but scattering decay subprincipal) contribution to the expression under the square root in~\eqref{EqSpBInftyOutSymb2}. Finally, $g^{0 0}=G(\dd t)\leq 0$ has a sign, and thus
  \[
    \la 2\check A h^{-1}(\Im z)g^{0 0}\check A u,u\ra \leq 0
  \]
  has a good sign as well.

  \pfsubstep{(1.2)}{Non-real $z$.} When $\Im z\geq\nu>0$, we shall argue completely analogously to Step~(1.2) in the proof of Proposition~\ref{PropSpBInftyOut}. Since we now need to control the derivative of $\psi_\pa$, we only need to recall that $P_{h,z}$ is elliptic on $\supp(\psi_\cR\dd\psi_\pa)$ in the notation of~\eqref{EqSpBInftyOutChecka}, as follows from Lemma~\ref{LemmaSpChar}\eqref{ItSpChar1}.

  \pfstep{Step 2.~Estimate with b-regularity.} The estimate~\eqref{EqSpHiInftyOutb} follows by applying~\eqref{EqSpHiInftyOut} to the commuted equation of Lemma~\ref{LemmaSpBInftyOutComm}, multiplied by $h^2$; note that $R_{\alpha\beta,h}:=h^2 R_{\alpha\beta}\in\rho\Diff_\schop^2$, and $h^2(P^{(k)}_\sigma)_{\alpha\alpha'}=\delta_{\alpha\alpha'}P_{h,z}-2 k i h z\rho\delta_{\alpha,(k,0,\ldots,0)}+h^2\tilde P_{\alpha\alpha'}$, with $h^2\tilde P_{\alpha\alpha'}\in h\rho\Diff_\schop^1$ having vanishing principal symbol at $\cR_{\rm out}$; the imaginary part of $2 k i h z\rho\delta_{\alpha,(k,0,\ldots,0)}$ has the correct sign when $\Im z\geq 0$, as in the proof of Proposition~\ref{PropSpBInftyOut}. We estimate the norms of $h^2\tilde f^{(k)}$ via induction on $k$; note that
  \[
    h^{-1}\|G R_{\alpha\beta,h}\sV^\beta u\|_{H_{\scop,h}^{s-2,r+1}}\leq C\Bigl(h^{-1}\|\tilde G u\|_{H_{(\scop,h);\bop}^{(s-2;k-1),r+1}}+h^N\|\tilde\chi u\|_{H_{(\scop,h);\bop}^{(-N;k-1),-N}}\Bigr)
  \]
  if $\tilde G\in\Psi_\schop^{0,0}$ is elliptic on $\WF'_\schop(G)$ and has Schwartz kernel supported on $\tilde\chi^{-1}(1)\times\tilde\chi^{-1}(1)$ where $\tilde\chi\in\CIc(\cX)$. The first term can then be estimated by induction on $k$; notice the factor $h^{-1}$ incurred here. The second term is $\sim h^{N-(k-1)}\|\tilde\chi u\|_{H_{\scop,h}^{-N+k-1,-N+k-1}}$ and thus residual since $N$ is arbitrary.
\end{proof}

\subsubsection{Estimates at the trapped set}
\label{SssSpHiTr}

Let $P_0$ be a general stationary wave operator on $\cM=\R_t\times\cX$ (Definition~\ref{DefSpGen}). We will consider trapped sets for the null-bicharacteristic flow for $P_{h,z}=h^2\wh{P_0}(h^{-1}z)$ over compact subsets of $T^*\cX$. Since, by Lemma~\ref{LemmaSpChar}, for $\Im z>0$ the characteristic set of $(P_{h,z})_{h\in(0,1)}$ lies over $S^*\cX$, we only need to study almost real $z$. For the sake of definiteness, we study $z\in\C$ with $\Re z=1$. Write $G_1(x,\xi):=g_{(0,x)}^{-1}(-\dd t+\xi\cdot\dd x)$ for the semiclassical characteristic function of $P_{h,1}$ and $\Sigma_1:=G_1^{-1}(0)\cap T^*\cX$ for the semiclassical characteristic set at finite momenta; we moreover write $\Sigma_1^+\subset\Sigma_1$ for the future component. We then make the following assumptions:
\begin{enumerate}
\myitem{ItSpHiTrGamma}{$\Gamma$.1} $\Gamma^{\rm u/s}$ are codimension $1$ submanifolds of $\Sigma_1^+$ which intersect transversally in the compact set
  \[
    \Gamma := \Gamma^{\rm u} \cap \Gamma^{\rm s}.
  \]
  Moreover, $H_{G_1}$ is tangent to $\Gamma^{\rm u/s}$ (and thus to $\Gamma$).
\myitem{ItSpHiTrDef}{$\Gamma$.2} There exist smooth defining functions $\phi^{\rm u/s}\in\CI(\Sigma_1^+)$ of $\Gamma^{\rm u/s}$ near $\Gamma$ (that is, in a neighborhood $\cO\subset T^*\cX$ of $\Gamma$, we have $\Gamma^{\rm u/s}\cap\cO=\{\phi^{\rm u/s}=0\}\cap\cU$, and $\dd\phi^{\rm u/s}\neq 0$ on $\Gamma^{\rm u/s}$) such that
  \[
    H_{G_1}\phi^{\rm u} = -\nu^{\rm u}\phi^{\rm u},\quad
    H_{G_1}\phi^{\rm s} = \nu^{\rm s}\phi^{\rm s}
  \]
  on $\Sigma_1^+$ near $\Gamma$ where $\nu^{\rm u},\nu^{\rm s}$ are smooth functions and $\nu_{\rm min}:=\min(\inf_\Gamma\nu^{\rm u},\inf_\Gamma\nu^{\rm s})>0$.
\myitem{ItSpHiTrSym}{$\Gamma$.3} We have $\{\phi^{\rm u},\phi^{\rm s}\}>0$ and $\dd G_1\neq 0$ on $\Gamma$.
\myitem{ItSpHiTrSub}{$\Gamma$.4} Define\footnote{The sign of $\beta_0$ is due to the compactness of $\Gamma$ and the fact that $-\dd t+\xi$ is future null for all $\xi\in\Gamma$, while $\dd t$ is past timelike.}
  \begin{equation}
  \label{EqSpHiTrp1}
    \beta_0 := \min_{\xi\in\Gamma} g^{-1}(-\dd t+\xi,\dd t)>0,\quad
    \sfp_1 := \upsigma_\semi\Bigl(\frac{P_{h,1}-P_{h,1}^*}{2 i h}\Bigr),
  \end{equation}
  and set
  \begin{equation}
  \label{EqSpHiTrpGammas}
    \gamma_- := \frac{\sup_\Gamma\sfp_1-\frac12\nu_{\rm min}}{2\beta_0},\quad
    \gamma_+ := \frac{\sup_\Gamma\sfp_1}{2\beta_0}.
  \end{equation}
  For $\eta\geq\gamma_-$ and $\eps\geq 0$, define the \emph{semiclassical loss function} by
  \begin{equation}
  \label{EqSpHiTrLoss}
    \delta_\Gamma(\eta) := \max\Bigl(\frac{\gamma_+-\eta}{\gamma_+-\gamma_-}, 0\Bigr) \in [0,1].
  \end{equation}
\end{enumerate}

The function~\eqref{EqSpHiTrLoss} interpolates between $1$ (for $\eta=\gamma_-$) and $0$ (for $\eta\geq\gamma_+$); it quantifies the semiclassical (powers of $h^{-1}$) loss that we incur in the high-energy estimates for frequencies $\sigma$ with $\Im\sigma=\eta$. (In our applications, $\gamma_+\geq 0$ will be zero or, for tensor-valued wave equations, less than any specified positive number for a suitable choice of fiber inner product, and thus $\gamma_-=\gamma_+-\frac{\nu_{\rm min}}{4\beta_0}$ will be strictly negative; thus $\delta_\Gamma(0)$ will be $0$ or arbitrarily small in this case, resulting in an arbitrarily small $h^{-\eps}$-loss below.)

Since we work over $T^*\cX$, the regularity orders are irrelevant; we take them to be $0$.

\begin{thm}[Estimate at normally hyperbolic trapping]
\label{ThmSpHiTr}
  Let $\beta>\gamma_-$. Let $\nu>0$ and $N\in\R$. Let $\cU\subset T^*\cX$ be an open neighborhood of $\Gamma$, and let $K\subset\cX$ be a compact neighborhood of the base projection of $\Gamma$. Let $\chi\in\CIc(\cX)$ be equal to $1$ near $K$.
  \begin{enumerate}
  \item{\rm (Direct estimate.)} There exist semiclassical ps.d.o.s $B,E,G\in\Psi_\semi(\cX)$ with Schwartz kernels supported in $K\times K$ and operator wave front sets in $\cU$, with $B$ semiclassically elliptic at $\pa\cR^+$ and $\WF'_\semi(E)$ disjoint from $\Gamma^{\rm u}$, and $C>0$ such that
    \begin{equation}
    \label{EqSpHiTr}
    \begin{split}
      &\|B u\|_{L^2} \leq C_\eps\Bigl( h^{-1-\delta_\Gamma(\Im(h^{-1}z)-\eps)}\|G P_{h,z}u\|_{L^2} + h^{-\delta_\Gamma(\Im(h^{-1}z)-\eps)}\|E u\|_{L^2} + h^N\|\chi u\|_{H_h^{-N}}\Bigr) \\
      &\qquad \forall\,z\in 1+i[h\beta,\nu],\ \ 0<h<1,\ \ \eps>0.
    \end{split}
    \end{equation}
    More precisely, this holds provided $\Gamma\cup\WF_\semi'(B)\subset\Ell_\semi(G)$ and all backward null-bi\-char\-ac\-ter\-is\-tics starting at a point in $\WF_\semi'(B)\setminus\Gamma^{\rm u}$ reach $\Ell_\semi(E)$ in finite time while remaining in $\Ell_\semi(G)$.
  \item{\rm (Adjoint estimate.)} There exist semiclassical ps.d.o.s $B,E,G\in\Psi_\semi(\cX)$ with Schwartz kernels supported in $K\times K$ and operator wave front sets in $\cU$, with $B$ semiclassically elliptic at $\pa\cR^+$ and $\WF'_\semi(E)$ disjoint from $\Gamma^{\rm s}$, and $C>0$ such that
    \begin{align*}
      &\|B u\|_{L^2} \leq C_\eps\Bigl( h^{-1-\delta_\Gamma(\Im(h^{-1}z)-\eps)}\|G P_{h,z}^*u\|_{L^2} + h^{-\delta_\Gamma(\Im(h^{-1}z)-\eps)}\|E u\|_{L^2} + h^N\|\chi u\|_{H_h^{-N}}\Bigr) \\
      &\qquad \forall\,z\in 1+i[-h\beta,\nu],\ \ 0<h<1,\ \ \eps>0.
    \end{align*}
    More precisely, this holds provided $\Gamma\cup\WF_\semi'(B)\subset\Ell_\semi(G)$ and all forward null-bi\-char\-ac\-ter\-is\-tics starting at a point in $\WF_\semi'(B)\setminus\Gamma^{\rm s}$ reach $\Ell_\semi(E)$ in finite time while remaining in $\Ell_\semi(G)$.
  \end{enumerate}
\end{thm}

Our strategy for the proof of Theorem~\ref{ThmSpHiTr} is as follows:
\begin{enumerate}
\item Note that for $\Im(h^{-1}z)>\gamma_+$, the estimates~\eqref{EqSpHiTr} are lossless compared to real principal type or radial point estimates. We first prove this lossless estimate in this case using a simple positive commutator argument, with $\Im P_{h,z}$ providing the main positivity in the symbolic computation.
\item We next prove Theorem~\ref{ThmSpHiTr} with the loss $\delta_\Gamma$ replaced by the constant $1$. In this case, the estimate is essentially deduced in \cite[Theorem~4.7]{HintzVasyQuasilinearKdS} from \cite{DyatlovSpectralGaps} when $|\Im z|=\cO(h)$. We give a self-contained proof here of this weaker result using a microlocalized version of semiclassical defect measures.
\item To obtain the smaller loss $\delta_\Gamma$, we use an interpolation argument. This is most conveniently carried out for an invertible problem with complex absorption, which must then be translated back into a propagation estimate; we do this using the arguments in \cite{HintzVasyQuasilinearKdS}.
\end{enumerate}

\begin{rmk}[$h$ away from $0$]
\label{RmkSpHiNonscl}
  For any $h_0>0$, the estimate~\eqref{EqSpHiTr} holds for $h\in[h_0,1)$ with the terms involving $P_{h,z}$ and $E$ on the right, as $\|B_h u_h\|_{L^2}\leq C\|\chi u\|_{H_h^{-N}}$ follows in this case simply from the fact that $B_h$ is a standard (i.e., not semiclassical) residual ps.d.o. Thus, it suffices to study $h<h_0$, with $h_0$ arbitrarily small.
\end{rmk}

\begin{proof}[Proof of Theorem~\usref{ThmSpHiTr} for $\Im(h^{-1}z)>\gamma_+$]
  Let $\gamma>\gamma_+$. We need to show that for $\Im z\geq h\gamma$, we have
  \begin{equation}
  \label{EqSpHiTrHi}
    \|B u\|_{L^2} \leq C\Bigl( h^{-1}\|G P_{h,z}u\|_{L^2} + \|E u\|_{L^2} + h^N\|\chi u\|_{H_h^{-N}}\Bigr),\quad z\in 1+i[h\gamma,\nu].
  \end{equation}
  For any $\nu'>0$, we have the even stronger elliptic estimate $\|B u\|\leq C(\|G P_{h,z}u\|_{L^2}+h^N\|\chi u\|_{H_h^{-N}})$, so it suffices to study the regime $h\gamma\leq\Im z\leq\nu'$ for any fixed $\nu'>0$.

  We prove~\eqref{EqSpHiTrHi} using a positive commutator argument with $A=\check A^*\check A$, $\check A=\check A^*=\Op_\semi(\check a)$,
  \[
    \check a = \psi_\Sigma\psi^{\rm u}\psi^{\rm s},\quad \psi_\Sigma:=\psi(\digamma_\Sigma G_1),\ \psi^{\rm u/s}:=\psi(\digamma\phi^{\rm u/s}),
  \]
  where $\psi$ is as in~\eqref{EqR3RHPsi} and $\digamma,\digamma_\Sigma>1$; here we write $\phi^{\rm u/s}\in\CI(T^*\cX)$ for arbitrary extensions of the eponymous functions in assumption~\eqref{ItSpHiTrDef} above. When $\digamma$ and $\digamma_\Sigma$ are large enough, this is supported in any fixed neighborhood of $\Gamma$. We again consider the ``commutator''~\eqref{EqSpHiPr2Pair}--\eqref{EqSpHiPr2C} and need to describe the terms $\frac{i}{h}[\Re P_{h,z},A]+2 h^{-1}\check A(\Im P_{h,z})\check A$.

  First, we note that
  \[
    H_{G_1}\phi^{\rm u}=-\nu^{\rm u}\phi^{\rm u}+j^{\rm u}G_1,\quad
    H_{G_1}\phi^{\rm s}=\nu^{\rm s}\phi^{\rm s}+j^{\rm s}G_1
  \]
  for some smooth functions $j^{\rm u/s}$ near $\Gamma$. Therefore, writing $(\psi^{\rm u})'=\psi'(\digamma\phi^{\rm u})$ (similarly for the other cutoff functions), we have
  \[
    H_{G_1}(\psi^{\rm u})^2 = -2\nu^{\rm u}\digamma\phi^{\rm u}(\psi^{\rm u})' \psi^{\rm u} + 2\digamma j^{\rm u}G_1(\psi^{\rm u})'\psi^{\rm u} = -2\digamma\phi^{\rm u}(\psi^{\rm u})' \psi^{\rm u} \Bigl( \nu^{\rm u} - \frac{j^{\rm u}G_1}{\phi^{\rm u}}\Bigr).
  \]
  Given any $\eta>0$, note then that on $\supp(\psi^{\rm u})'\cap\supp\psi_\Sigma$, we have $|G_1|\leq\digamma_\Sigma^{-1}$ and $\phi^{\rm u}\geq\frac{1}{2\digamma}$, so upon choosing $\digamma_\Sigma$ large enough, the term in parenthesis is $\geq(1-\eta)\nu^{\rm u}$ on $\supp\check a$. In order to take into account the possibility of non-real $\sigma$, note that by~\eqref{EqSpChar}, we can write $H_{G_1}-H_{\Re G_z}=(\Im z)\tilde H_z$ where the vector field $\tilde H_z$ has coefficients of size $1$. Since $(\Im z)\tilde H_z(\psi^{\rm u})^2 = 2(\Im z)(\tilde H_z\phi^{\rm u})(\psi^{\rm u})'\psi^{\rm u}$, we have
  \[
    H_{\Re G_z}(\psi^{\rm u})^2 = -2\digamma\phi^{\rm u}(\psi^{\rm u})'\psi^{\rm u}\Bigl(\nu^{\rm u}-\frac{j^{\rm u}G_1}{\phi^{\rm u}} - \frac{\Im z}{2}\frac{\tilde H_z\phi^{\rm u}}{\phi^{\rm u}}\Bigr);
  \]
  the third term in parentheses is bounded by $C|\Im z|\digamma$ in absolute value and thus smaller than $\eta\nu^{\rm u}$ when $|\Im z|$ is small enough depending on $\digamma$.\footnote{This is an instance of the fact that positivity is preserved under small perturbations.} Completely analogous considerations apply to
  \[
    H_{\Re G_z}(\psi^{\rm s})^2 = 2\digamma\phi^{\rm s}(\psi^{\rm s})'\psi^{\rm s}\Bigl(\nu^{\rm s}+\frac{j^{\rm s}G_1}{\phi^{\rm s}}+\frac{\Im z}{2}\frac{\tilde H_z\phi^{\rm s}}{\phi^{\rm s}}\Bigr).
  \]
  Finally, on $\supp\dd\psi_\Sigma\cap\supp\check a$ we have $|G_1|\geq\frac{1}{2\digamma_\Sigma}$, and thus also $G_z$ is elliptic there when $|\Im z|$ is small enough depending on $\digamma_\Sigma$.

  We use Lemma~\ref{LemmaSpImag} to write
  \begin{subequations}
  \begin{equation}
  \label{EqSpHiImag1}
    h^{-1}\Im P_{h,z}= h^{-1}\Im P_{h,1} - h^{-1}\Im(z)Q_{1,h};
  \end{equation}
  recalling~\eqref{EqSpHiTrp1}, we thus have
  \begin{equation}
  \label{EqSpHiImag2}
    \upsigma_\semi(h^{-1}\Im P_{h,1})=\sfp_1,\quad \upsigma_\semi(Q_{1,h})\geq 2\beta_0.
  \end{equation}
  \end{subequations}
  Roughly speaking, then, when $h^{-1}\Im(z)\geq\gamma$, then $h^{-1}\Im P_{h,z}$ is $\leq\sfp_1-2\beta_0\gamma$ on the symbolic level, and thus negative at $\Gamma$ under the condition~\eqref{EqSpHiTrHi}. We now make this precise. Fix $\gamma_0\in(\frac{\inf_\Gamma\sfp_1}{2\beta_0},\gamma)$. If $\digamma$ is large enough, then $\digamma_\Sigma$ is chosen large enough, and then $\nu'>0$ (and thus $|\Im z|$) is small enough, we can write
  \[
    \sC \equiv \Op_\semi(-b^2 - b_{\rm u}^2 + e_{\rm s}^2 + j G_z) - 2 h^{-1}\Im(z)\check A\Bigl(Q_{1,h}-2\beta_0\frac{\gamma_0}{h^{-1}\Im z}\Bigr)\check A \bmod h\Psi_\semi,
  \]
  where the elements
  \begin{align*}
    b &:= \psi_\Sigma\psi^{\rm u}\psi^{\rm s}\bigl[ -2\sfp_1 + 4\beta_0\gamma_0 \bigr]^{\frac12}, \\
    b_{\rm u} &:= \psi_\Sigma\psi^{\rm u}\bigl[ -H_{\Re G_z}(\psi^{\rm s})^2 \bigr]^{\frac12}, \\
    e_{\rm s} &:= \psi_\Sigma\psi^{\rm s}\bigl[ H_{\Re G_z}(\psi^{\rm u})^2 \bigr]^{\frac12}, \\
    j &:= 2(\psi^{\rm u})^2(\psi^{\rm s})^2 \digamma_\Sigma\psi_\Sigma'\psi_\Sigma (H_{\Re G_z}G_1) / G_z
  \end{align*}
  of $\CIc(T^*\cX)$ are well-defined by the above considerations. The term of $\sC$ involving $\Im(z)$ yields a negative contribution (modulo $h^N\|\chi u\|_{H_h^{-N}}$) when acting on $u$ and being paired with $u$ since the principal symbol of $Q_{1,h}-2\beta_0\frac{\gamma_0}{h^{-1}\Im z}$ is $\geq 2\beta_0-2\beta_0\frac{\gamma_0}{\gamma}>0$. The term $-2\check a^2$ yields control at $\Gamma$ given a priori control on $\supp e_{\rm s}$, which is indeed disjoint from $\Gamma^{\rm u}$.
\end{proof}

For parts of the proof of Theorem~\ref{ThmSpHiTr} for $|\Im z|=\cO(h)$, we use a microlocalized version of semiclassical defect measures on $\cU$. (One could alternatively adapt the arguments of \cite[\S{3}]{HintzPolyTrap} to the semiclassical setting, i.e., use a direct commutator argument and the G\aa{}rding inequality.)

\begin{lemma}[Microlocalized semiclassical defect measures]
\label{LemmaSpHiDef}
  Let $\cU\subset T^*\cX$ be open and precompact, and let $\chi\in\CIc(\cX)$ be equal to $1$ near the base projection of $\cU$. Consider a sequence $(h_j)_{j\in\N}$ in $(0,1)$ converging to $0$, and let $u_{h_j}$, $j\in\N$, be distributions on $\cX$ such that
  \begin{subequations}
  \begin{equation}
  \label{EqSpHiDefchi}
    \|\chi u_{h_j}\|_{H_h^{-N_0}}\leq C h^{-N_0}
  \end{equation}
  for some $N_0$; assume moreover that there exists $G=(G_h)_{h\in(0,1)}\in\Psi_\semi(\cX)$, elliptic on $\cU$, such that
  \begin{equation}
  \label{EqSpHiDefG}
    \sup\|G_{h_j}u_{h_j}\|_{L^2}<\infty.
  \end{equation}
  \end{subequations}
  Then there exist a subsequence $(u_{h_{j_k}})_{k\in\N}$ and a Radon measure $\mu$ on $\cU$ such that for all $B=(B_h)_{h\in(0,1)}\in\Psi_\semi(\cX)$ with Schwartz kernel supported in $K\times K$, $\WF'_\semi(B)\subset\cU$ compact, and with semiclassical principal symbol $b=\upsigma_\semi(B)\in\CIc(\cU)$, we have
  \begin{equation}
  \label{EqSpHiDef}
    \lim_{k\to\infty} \la B_{h_{j_k}}u_{h_{j_k}},u_{h_{j_k}}\ra_{L^2} = \int_\cU b\,\dd\mu.
  \end{equation}
\end{lemma}

We call a sequence $(u_{h_j})_{j\in\N}$ of distributions on $\cX$ \emph{$h$-tempered on $\cU,K$} if~\eqref{EqSpHiDefchi}--\eqref{EqSpHiDefG} hold. We moreover call $\mu$ a \emph{(semiclassical) defect measure of $(u_{h_{j_k}})_{k\in\N}$ on $\cU$}. For the proof, we record that
\begin{equation}
\label{EqSpHiBound}
\begin{split}
  &A\in\Psi_\semi,\ \WF'_\semi(A)\subset T^*\cX,\ \text{with Schwartz kernel supported in $K\times K$} \\
  &\qquad \implies \|A u\|_{L^2} \leq 2\|a\|_\infty\|\chi u\|_{L^2} + C h^{\frac12}\|\chi u\|_{H_h^{-N}};
\end{split}
\end{equation}
indeed, either $a=\upsigma_\semi(A)=0$ and thus $A\in h\Psi_\semi$, so $|\la A u,u\ra|\leq C h^{\frac12}\|\chi u\|_{H_h^{-N}}$ (the regularity order being arbitrary), or $\|a\|_\infty=\sup|a|>0$, in which case a square root construction gives $2\|a\|_\infty^2 - A^*A=B^*B+R$ with $B\in\Psi_\semi$ and $R\in h^\infty\Psi_\semi$, so $\|A u\|^2\leq 2\|a\|_\infty^2\|u\|^2+C_N h^N$.

\begin{proof}[Proof of Lemma~\usref{LemmaSpHiDef}]
  This is a microlocalization of the arguments in \cite[Chapter~5.1]{ZworskiSemiclassical}. The details are as follows. Fix any compact subset $\cK\subset\cU$. Given any $B\in\Psi_\semi$ with $\WF'_\semi(B)\subset\cK$, let $Q\in\Psi_\semi$ be a quantization of $b/g$; then $Q G=B+R$ where $R\in h\Psi_\semi$. Adding to $Q$ a quantization of $h\upsigma_\semi(h^{-1}R)/g$, we can improve the remainder to $h^2\Psi_\semi$; continuing in this fashion, we get a new operator $Q\in\Psi_\semi$, with $q=\upsigma_\semi(Q)=b/g$ still and $\WF'_\semi(Q)\subset\cK$, such that $Q G=B+R$ with $R\in h^\infty\Psi_\semi$. Therefore,
  \[
    |\la B u,u\ra| \leq |\la G u,Q^* u\ra| + |\la R u,u\ra| \leq C\|Q^*u\| + C_N h^N.
  \]
  Let $J\in\Psi_\semi$ with $\WF'_\semi(J-I)\cap\cK=\emptyset$ and $\WF'_\semi(J)\subset\cU$, then we can write $Q^*=Q^*J+Q^*(I-J)$; but $Q^*(I-J)\in h^\infty\Psi_\semi$, while we can write $J=\tilde Q G+\tilde R$ where $\tilde Q\in\Psi_\semi$ has a principal symbol $\tilde q$ with $\sup|\tilde q|\leq C_{0,\cK}:=\sup_\cK|g|^{-1}$ and $\tilde R\in h^\infty\Psi_\semi$, so
  \[
    \|Q^*u\|\leq \|Q^*\tilde Q G u\| + C_N h^N.
  \]
  We then use~\eqref{EqSpHiBound} (or \cite[Theorem~5.1]{ZworskiSemiclassical}, or the sharp G\aa{}rding inequality) to bound this by $2\sup|q\tilde q| \|G u\|_{L^2}+C'h^{\frac12}$. Since $|q|\leq C_{0,\cK}|b|$, we have, altogether, shown that there exist constants $C_\cK$, $C_{\cK,B}$ such that
  \begin{equation}
  \label{EqSpHiDefEst}
    |\la B u,u\ra| \leq C_\cK\sup|b| + C_{\cK,B}h^{\frac12},\quad B\in\Psi_\semi,\ \WF'_\semi(B)\subset\cK.
  \end{equation}

  Fixing a countable dense subset $\{b_l\colon l\in\N\}$ of $\CIc(\cU)$, a diagonal argument shows that, upon replacing $(u_{h_j})$ by a subsequence, we have convergence $\lim_{j\to\infty}\la (B_l)_{h_j}u_{h_j},u_{h_j}\ra=:\alpha_l\in\C$ for all $l$; and for $b_l\in\CIc(\cK)$ we have $|\alpha_l|\leq C_\cK|b_l|$ by~\eqref{EqSpHiDefEst}. Thus, the map $b_l\mapsto\alpha_l$ extends to a continuous map $\cC^0_{\rm c}(\cU)\to\C$ which can be represented in the form~\eqref{EqSpHiDef} by the Riesz representation theorem, with $\mu$ a possibly complex measure.

  To conclude, we only need to show that if $B\in\Psi_\semi$, $\WF'_\semi(B)\subset\cK\Subset\cU$, and $b\geq 0$, then $\lim\la B_{h_j}u_{h_j},u_{h_j}\ra\geq 0$. Let $J$ be as above, with principal symbol $j$ equal to $1$ on $\cK$; then given $\eps>0$, we have $b_\eps:=b+\eps|j|^2>0$ on $\cK$, and $b_\eps=\eps|j|^2$ on $\cU\setminus\cK$; so $\sqrt{b_\eps}=q_\eps\in\CIc(\cU)$, and thus for $Q_\eps=\Op_\semi(q_\eps)$, $B+\eps J^*J=Q_\eps^*Q_\eps+h R$ where $R\in\Psi_\semi$, $\WF'_\semi(R)\Subset\cU$. Therefore, the $\liminf$ of
  \[
    \la B u,u\ra = \|Q_\eps u\|^2 + h\la R u,u\ra - \eps\|J u\|^2
  \]
  as $h_j\to 0$ is $\geq -C\eps$. Since $\eps>0$ is arbitrary, we are done.
\end{proof}

\begin{lemma}[Elliptic regularity]
\label{LemmaSpHiDefEll}
  Let $(u_{h_j})_{j\in\N}$ be $h$-tempered on $\cU\subset T^*\cX$, $K\subset\cK$, with $\bar\cU$ and $K$ compact. Let $P\in\Psi_\semi(\cX)$ have Schwartz kernel supported in $K\times K$; denote its principal symbol by $p$. Suppose that $\|G P u\|_{L^2}=o(1)$ as $h\to 0$. If $\mu$ is a defect measure for $(u_{h_j})$ on $\cU$, then $\supp\mu\subset p^{-1}(0)$.
\end{lemma}
\begin{proof}
  Let $b\in\CIc(\cU)$ with $\supp b\cap p^{-1}(0)=\emptyset$. Then $P$ is elliptic on $B=\Op_\semi(b)$, so we can write $B=G P u+R u$ with $R\in h^\infty\Psi_\semi^{-\infty}$, and therefore
  \[
    \int_\cU |b|^2\,\dd\mu = \lim_{h\to 0}\la B u,B u\ra = \la G P u+R u,G P u+R u\ra_{L^2} = 0
  \]
  since $\|G P u\|_{L^2}$, $\|R u\|_{L^2}=o(1)$ as $h\to 0$.
\end{proof}

The following is a microlocal analogue of \cite[Lemma~2.3]{DyatlovSpectralGaps}.

\begin{lemma}[Propagation of regularity]
\label{LemmaSpHiDefPr}
  Let $\mu$ be a microlocal semiclassical defect measure of the $h$-tempered family $(u_{h_j})_{j\in\N}$ on $\cU,K$. Let $\chi\in\CIc(\cX)$ be equal to $1$ near $K$. Let $P\in\Psi_\semi$, with Schwartz kernel supported in $K\times K$, have real principal symbol $p$, and set $\sfp_1:=\upsigma_\semi(\frac{1}{2 i h}(P-P^*))$. Then for all $a\in\CIc(\cU)$ and $Y\in\Psi_\semi$, with Schwartz kernel supported in $K\times K$, which satisfy $\WF'_\semi(Y)\subset\cU$ and $\WF'_\semi(I-Y)\cap\supp a=\emptyset$, we have
  \begin{equation}
  \label{EqSpHiDefPr}
    \biggl| \int_\cU (H_p+2\sfp_1)a\,\dd\mu\biggr| \leq 4\|a\|_\infty \limsup_{j\to\infty} \Bigl( h_j^{-1}\|Y P u\|_{L^2}\|Y u\|_{L^2} \Bigr)
  \end{equation}
  In particular, if $\|Y P u\|_{L^2}=o(h_j)$ as $j\to\infty$, then $\int_\cU (H_p+2\sfp_1)a\,\dd\mu=0$; and in this case, $\supp\mu$ is invariant under the $H_p$-flow in $\cU$.
\end{lemma}
\begin{proof}
  By considering the real and imaginary parts of $a$ separately, we may assume that $a$ is real. Let $A=A^*\in\Psi_\semi$ be a quantization of $a$ with $\WF'_\semi(A)=\supp a$. Then
  \[
    2 h^{-1}\Im\la P u,A u\ra = \la\sC u,u\ra,\quad \sC=\frac{i}{h}[\Re P,A]+h^{-1}\bigl( (\Im P)A+A(\Im P) \bigr).
  \]
  Since $\upsigma_\semi(\sC)=H_p a+2\sfp_1 a$, the right-hand side converges to $\int_\cU (H_p+2\sfp_1)a\,\dd\mu$. We can estimate the left-hand side by
  \[
    2 h^{-1} \bigl(\| Y P u \| + C h^N\|\chi u\|_{H_h^{-N}}\bigr)\bigl(\|A Y u\|_{L^2} + C h^N\|\chi u\|_{H_h^{-N}}\bigr),
  \]
  and in turn $\|A Y u\|_{L^2} \leq 2(\|a\|_\infty\|Y u\|_{L^2} + h^{\frac12}\|\chi u\|_{H_h^{-N}})$ by~\eqref{EqSpHiBound}. Taking the $\limsup$ as $h\to 0$, one obtains the right-hand side of~\eqref{EqSpHiDefPr}.

  If $\varpi\notin\supp\mu$ and indeed $\cV\cap\supp\mu$ for a neighborhood $\cV$ of $\varpi$, then we can take $a$ to be a function localized near a segment $\gamma([0,1])$ of an integral curve of $H_p$ with $\gamma(0)=\varpi$ such that $(H_p+2\sfp_1)a\leq 0$ outside of $\cV$, and this is strictly negative $\gamma([0,1])\setminus\cV$. This gives $\gamma([0,1])\cap\supp\mu=\emptyset$.
\end{proof}

We now turn to the proof of the trapping estimate.

\begin{proof}[Proof of Theorem~\usref{ThmSpHiTr} for $\gamma_-<\Im(h^{-1}z)\lesssim 1$, with loss $1$]
  We demand that \emph{all ps.d.o.s chosen in the course of the proof have Schwartz kernels supported in $K\times K$.} Let $\gamma>\gamma_+$. The semiclassical principal symbol of $(P_{h,z})_{h\in(0,1)}$ for $\Re z=1$ and $|\Im z|=\cO(h)$ is $G_1$. Suppose the estimate~\eqref{EqSpHiTr} does not hold for any constant, then we obtain sequences $u_{h_j}\in\sD'(\cX)$, $h_j\to 0$, and $z_j$ with $\Re z=1$, $-\beta\leq h^{-1}\Im z_j\leq\gamma$, such that
  \begin{equation}
  \label{EqSpHiTrContr}
    \|B u_{h_j}\|=1,\ \ 
    \|G P_{h_j,z_j}u_{h_j}\|=o(h_j^2),\ \ 
    \|E u_{h_j}\|=o(h_j),\ \ 
    h_j^N\|\chi u_{h_j}\|_{H_{h_j}^{-N}} = o(1).
  \end{equation}
  Upon passing to a subsequence, we may assume that $\frac{\Im z_j}{h_j}\to\eta\in[-\beta,\gamma]$ and $\mu$ is a microlocal semiclassical defect measure in $\cU_1:=\Ell_\semi(B)$, $K$. In the remainder of the proof, we drop the subscript ``$j$'' (and often also ``$h$''). Thus $f=P_{h,z}u$ satisfies $\|G f\|=o(h^2)$. Since every point in $\WF'_\semi(B)\setminus\Gamma^{\rm u}$ is connected to a point in $\Ell_\semi(E)$ via a null-bicharacteristic remaining in $\Ell_\semi(G)$, we have $\supp\mu\subset\Gamma^{\rm u}$. On the other hand, null-bicharacteristic starting on $\cU\cap\Gamma^{\rm u}$ reach any fixed neighborhood of $\Gamma$ while remaining in $\Ell_\semi(G)$ (due to $H_p\phi^{\rm s}\neq 0$ on $\Gamma^{\rm u}\setminus\Gamma$); therefore, if $E'\in\Psi_\semi$ is elliptic at $\Gamma$, then
  \[
    \|B u\| \leq C\Bigl( h^{-1}\|G P_{h,z}u\| + \|E u\| + \|E' u\| + h^N\|\chi u\|_{H_h^{-N}}\Bigr)
  \]
  In order to obtain a contradiction to~\eqref{EqSpHiTrContr}, it thus suffices to show that $\|E' u\|=o(1)$ for some such $E'$. Since we may take $E'$ to have $\WF'_\semi(E')\subset\cU_1$ (and thus $\|E' u\|^2\to\int_\cU |e'|^2\,\dd\mu=\int_{\cU\cap\Gamma^{\rm u}} |e'|^2\,\dd\mu$ where $e'=\upsigma_\semi(E')$), we conclude that it suffices to show that, for any arbitrarily small $\delta>0$, we have
  \begin{equation}
  \label{EqSpHiTrGoal}
    \mu(\delta)=0,\quad \mu(x):=\mu\bigl(\Gamma^{\rm u}\cap\{|\phi^{\rm s}|<x\}\bigr).
  \end{equation}

  We replace $P_{h,z}$ by $Y P_{h,z}$ where $Y\in\Psi_\semi$ has $\WF'_\semi(Y)\subset\cU_1$ and $\WF'_\semi(I-Y)$ is disjoint from a fixed neighborhood $\cU_0$ of $\Gamma$, with $\ol{\cU_0}\subset\cU_1$. In all computations of operators below, we write $R_{\rm away}$ for every semiclassical ps.d.o.\ with $\WF'_\semi(R_0)\subset\cU_1$ and $\WF'_\semi(R_0)\cap\ol{\cU_0}=\emptyset$. (Thus, $R_{\rm away}$ may change from line to line.)

  \pfstep{Step~1: auxiliary equation.} Let $\Phi^{\rm u},N^{\rm u}\in\Psi_\semi$ be quantizations of $\phi^{\rm u}$, $\nu^{\rm u}$, multiplied with a microlocalizer that is equal to $I$ modulo $h^\infty\Psi_\semi$ near $\ol{\cU_0}$ and has operator wave front set in supported in $\cU_1$. In the notation of~\eqref{EqSpHiTrp1}, let moreover $J^{\rm u}\in\Psi_\semi$ be a quantization of $j^{\rm u}$ (cut off in the same fashion). Then
  \[
    \frac{i}{h}[P_{h,z},\Phi^{\rm u}] = -N^{\rm u}\Phi^{\rm u} + J^{\rm u}P_{h,z} + h R + R_{\rm away}u,\quad R\in\Psi_\semi,\ \WF'_\semi(R)\subset\cU,
  \]
  and hence $P_{h,z}\Phi^{\rm u}u = \Phi^{\rm u}f + [P_{h,z},\Phi^{\rm u}]u$ implies
  \begin{equation}
  \label{EqSpHiTrAux}
    (P_{h,z}-i h N^{\rm u})(\Phi^{\rm u}u) = (\Phi^{\rm u}-i h J^{\rm u})f - i h^2 R u + R_{\rm away}u.
  \end{equation}
  Recalling~\eqref{EqSpHiImag1}--\eqref{EqSpHiImag2} and $h^{-1}\Im z\geq-\beta$, we have $\upsigma_\semi(h^{-1}\Im(P_{h,z}-i h N^{\rm u}))\leq(\sfp_1+2\beta_0\beta)-\nu^{\rm u}<0$ at $\Gamma$. Therefore, the estimate~\eqref{EqSpHiTrHi} applies, in a small neighborhood of $\Gamma$, to $P_{h,z}-i h N^{\rm u}$ and $\Phi^{\rm u}u$ in place of $P_{h,z}$ and $u$. Notice that the right-hand side of~\eqref{EqSpHiTrAux}, microlocalized to $\cU_0$, is of size $\cO(h^2)$ in $L^2$. Upon shrinking $\cU_0$ if necessary, we thus obtain
  \begin{equation}
  \label{EqSpHiTrAux2}
    \|\Phi^{\rm u}u\|_{L^2} \leq C h.
  \end{equation}
  for some constant $C$.

  \pfstep{Step~2: partial regularity of the defect measure along $\Gamma^{\rm u}$.} We now follow \cite{DyatlovSpectralGaps} with only minor changes. We apply the estimate~\eqref{EqSpHiDefPr} to the equation~\eqref{EqSpHiTrAux2}, thus getting
  \begin{equation}
  \label{EqSpHiTrAux3}
    \biggl|\int_{\cU\cap\Gamma^{\rm u}} H_{\phi^{\rm u}}a\,\dd\mu\biggr|\leq C\|a\|_\infty\quad\forall\,a\in\CIc(\cU_0).
  \end{equation}
  We then take
  \begin{equation}
  \label{EqSpHiTrAuxa}
    a=a_0(\phi^{\rm u})a_1(\phi^{\rm s})\psi(G_1),\quad a_0,a_1,\psi\in\CIc(\R)
  \end{equation}
  where $a_0,a_1,\psi$ are supported in a fixed small neighborhood $(\delta_0,\delta_0)$ of $0$ (so that $\supp a\subset\cU_0$), furthermore $\|a_0\|_\infty$, $\|a_1\|_\infty\leq 1$, further $a_0=1$ and $\psi=1$ near $0$, and for some small $\delta>0$,
  \[
    a'_1\geq\frac{1}{2\delta}\ \text{on}\ [-\delta,\delta], \quad
    a'_1\geq 0\ \text{on}\ [-2\delta,2\delta], \quad
    |a'_1|\leq 2\delta_0^{-1}\ \text{on}\ [-\delta_0,\delta_0]\setminus[-2\delta,2\delta].
  \]
  Since $H_{\varphi^{\rm u}}a=a_0 a_1'(H_{\phi^{\rm u}}\phi^{\rm s})$, condition~\eqref{ItSpHiTrSym} implies that~\eqref{EqSpHiTrAux3} yields $\delta^{-1}\mu(\delta) \leq C\mu(\delta_0)$, i.e.,
  \begin{equation}
  \label{EqSpHiTrAux4}
    \mu(\delta)\leq C\delta\mu(\delta_0).
  \end{equation}

  \pfstep{Step~3: propagation for logarithmic times.} We now use a propagation estimate for the equation $P_{h,z}u=f$ to control $\mu(\delta_0)$ in turn by $\mu(\delta)$. We apply~\eqref{EqSpHiDefPr} to the equation $P_{h,z}u=o(h^2)$. Set $\sfp_z:=\upsigma_\semi(h^{-1}\Im P_{h,z})$ and recall from~\eqref{EqSpHiImag1}--\eqref{EqSpHiImag2} that $\sfp_z\leq\sfp_1-2(h^{-1}\Im z)\beta_0\leq\sfp_1+2\beta\beta_0$. We then have
  \[
    \int_{\cU\cap\Gamma^{\rm u}} (H_{G_1}+2\sfp_z)a\,\dd\mu = 0.
  \]
  We apply this with $a$ of the form~\eqref{EqSpHiTrAuxa} and $a_0,\psi$ as there, but now $a_1$ is supported in a small neighborhood of $[-\delta_0,\delta_0]$ and is of the form $a_1(\phi^{\rm s})=|\phi^{\rm s}|^{-2\lambda}f(|\phi^{\rm s}|)$ where we fix $\lambda\in(\sup_\Gamma\frac{\sfp_z}{\nu^{\rm s}},\frac12)$ (this interval being non-empty in view of assumption~\eqref{ItSpHiTrSub}) and require $f\in\CIc((\frac{\delta}{2},\infty))$, $f\geq 0$, and
  \[
    f|_{[\delta,\delta_0]}=1, \quad
    f'\leq 0\ \text{on}\ [\delta,\infty), \quad
    |x f'(x)| \leq 4\ \text{for}\ x\in[\tfrac{\delta}{2},\delta].
  \]
  Thus $(H_{G_1}+2\sfp_z)a=\nu^{\rm s}|\phi^{\rm s}|^{-2\lambda} ( -2(\lambda-\frac{\sfp_z}{\nu^{\rm s}})f + |\phi^{\rm s}|f'(|\phi^{\rm s}|))$ is non-positive for $|\phi^{\rm s}|\geq\delta$ and less than a negative constant on $\Gamma^{\rm u}\cap\{\delta\leq|\phi^{\rm s}|\leq\delta_0\}$, while for $|\phi^{\rm s}|\in[\frac{\delta}{2},\delta]$ it is bounded in absolute value by $C\delta^{-2\lambda}$. Therefore, $\mu(\delta_0)-\mu(\delta) \leq C\delta^{-2\lambda}\mu(\delta)$, which gives
  \begin{equation}
  \label{EqSpHiTrLog}
    \mu(\delta_0) \leq C\delta^{-2\lambda}\mu(\delta).
  \end{equation}

  \pfstep{Step~4: conclusion.} Plugging~\eqref{EqSpHiTrAux4} into~\eqref{EqSpHiTrLog} gives $\mu(\delta)\leq C'\delta^{1-2\lambda}\mu(\delta)$ and thus $\mu(\delta)=0$ when $\delta>0$ is sufficiently small. This establishes~\eqref{EqSpHiTrGoal} and thus finishes the proof of Theorem~\ref{ThmSpHiTr} with constant loss $\delta_\Gamma=1$.
\end{proof}

To set up the interpolation argument, we now modify $P_{h,z}$ microlocally away from $\Gamma$ by adding complex absorption. Some care is required to ensure the holomorphicity of the modified operator in the parameter $\sigma=h^{-1}z$. For small $\delta_0>0$, we define the compact neighborhoods
\[
  \Gamma(\delta) := \{ |\phi^{\rm u}|\leq\delta,\ |\phi^{\rm s}|\leq\delta,\ |G_1|\leq\delta \} \subset T^*\cX
\]
of $\Gamma$ for $\delta\leq\delta_0$. Let $q_0\in\CI(T^*\cX)$ be equal to $0$ on $\Gamma(\frac{\delta_0}{2})$ and $1$ outside $\Gamma(\frac{3}{4}\delta_0)$. For $\sigma\in\R$ and $\xi\in T_x^*\cX$, set then
\[
  q_1(x;\sigma,\xi) = \begin{cases} q_0\bigl(x;\frac{\xi}{\sigma}\bigr), & \sigma>0,\ \frac{\xi}{\sigma}\in\Gamma(\delta_0), \\ 1 & \text{otherwise}. \end{cases}
\]
Fix $\chi_0\in\CIc((-1,1))$, equal to $1$ on $[-\frac12,\frac12]$, and define
\[
  q_2(x;\sigma,\xi) := \bigl(1-\chi_0(|(\sigma,\xi)|)\bigr) (\sigma^2+|\xi|^2) q_1(x;\sigma,\xi),
\]
where we use a Riemannian metric on $\cX$ to define $|\xi|^2$. Thus, $q_2$ is positively homogeneous of degree $2$ away from the zero section of $\ubar\R\oplus T^*\cX$. Denote its quantization by $Q_2=\Op(q_2)\in\Psi_\cl^2(\R\times\cX)$; its Schwartz kernel takes the form $K_{Q_2}(t,x,t',x')=k_{Q_2}(t-t';x,x')$, in local coordinates $K_{Q_2}(t,x,t',x')=(2\pi)^{-n-1}\chi_0(|x-x'|)\iint e^{i(\xi\cdot(x-x')-\sigma(t-t'))}q_2(x;\sigma,\xi)\,\dd\sigma\,\dd\xi$. We localize near $t-t'$ using a further cutoff and consider the operator $Q_3\in\Psi_\cl^2(\R\times\cX)$ with Schwartz kernel $K_{Q_3}(t,x,t',x')=\chi_0(t-t')k_{Q_2}(t-t';x,x')$. This is thus $t$-translation-invariant, and therefore given by the quantization of a symbol
\[
  q_3 = q_3(x;\sigma,\xi) \in S_\cl^2(\ubar\R\oplus T^*\cX).
\]
Note that $e^{-\eta t}Q_3 e^{\eta t}\in\Psi_\cl^2$, $\eta\in\C$, has the same principal symbol and operator wave front set as $Q_3$; indeed, in local coordinates its full symbol is $q_3(x;\sigma+i\eta,\xi)$. Thus, the spectral family $\wh{Q_3}(\sigma)$ is an entire family, in $\sigma\in\C$, of elements of $\Psi_\cl^2(\cX)$; and for $z=1+i\eta$, $\eta\in\R$, the semiclassical rescaling $h^2\wh{Q_3}(h^{-1}z)$, with full semiclassical symbol $(h,x;\xi_\semi)\mapsto h^2 q_3(x;h^{-1}z,h^{-1}\xi_\semi)=h^2 q_3(x;h^{-1}(1+i\eta),h^{-1}\xi_\semi)$, has $q_0$.

We next modify $P_0-i Q_3$ in order to place it the product of $\R$ with a compact manifold without boundary. Thus, let $\chi\in\CIc(\cX)$ be equal to $1$ near $\pi(\Gamma(\delta_0))$ where $\pi\colon T^*\cX\to\cX$ is the projection, and embed an open neighborhood of $\chi$ diffeomorphically into a closed Riemannian manifold $\cX'$. Set then
\[
  P' := \chi(P_0-i Q_3) - i(1-\chi)(D_t^2+\Delta_{\cX'}) \in \Psi_\cl(\R\times\cX').
\]
The semiclassical rescaling
\[
  P'_{h,z} := h^2\wh{P'}(h^{-1}z)
\]
of its spectral family is well-defined for all $z\in\C$. For $z=1+\cO(h)$, $(P'_{h,z})_{h\in(0,1)}\in\Psi_{\semi,\cl}^2(\cX')$ is semiclassical elliptic outside of $\Gamma(\frac34\delta_0)$, including at fiber infinity.

\begin{lemma}[High-energy estimates for $P'$]
\label{LemmaSpHiTrPp}
  Recall~\eqref{EqSpHiTrLoss}. Let $\gamma_-<\beta<\beta'$. Then there exists $h_0>0$ such that for all $z\in 1+i[h\beta,h\beta']$ and $0<h<h_0$, the operator $P'_{h,z}\colon H^2_h(\cX')\to L^2(\cX')$ is invertible, and we have
  \begin{equation}
  \label{EqSpHiTrPp}
    \|(P'_{h,z})^{-1}u\|_{L^2} \leq C_\eps h^{-1-\delta_\Gamma(\Im(h^{-1}z)-\eps)} \|u\|_{L^2},\quad 0<h<h_0,\ \ \eps>0.
  \end{equation}
\end{lemma}
\begin{proof}
  Fix any $\beta>\gamma_-$ and $\beta'>\gamma_+$. We first claim that, for all $h\in(0,1)$ and $z\in 1+i[h\beta,h\beta']$,
  \begin{align*}
    \|u\|_{L^2} &\leq C h^{-2}\Bigl(\|P'_{h,z}u\|_{H_h^{-2}} + h^N\|u\|_{H_h^{-N}}\Bigr), \\
    \|u\|_{L^2} &\leq C h^{-1}\Bigl(\|P'_{h,1+i h\beta'}u\|_{H_h^{-2}} + h^N\|u\|_{H_h^{-N}}\Bigr).
  \end{align*}
  To prove this, we first use semiclassical microlocal elliptic estimates on the elliptic set of the complex absorption, which is the union of $\supp\chi\cap\Ell(q_0)$ and $\supp(1-\chi)$. In the remaining set, which is a subset of $\Gamma(\delta_0)$, all backwards null-bicharacteristics starting at a point in $\WF_\semi'(P_{h,z})$ either tend to $\Gamma$ or reach the complex absorption in finite time. The version of Theorem~\ref{ThmSpHiTr} we have already proved thus controls $u$ microlocally near $\Gamma$ in the stated fashion.

  For sufficiently small $h>0$ then, we can absorb the error term in these estimates into the left-hand side. Arguing similarly for the adjoint problem, we thus obtain the invertibility of $P'_{h,z}\colon L^2\to H_h^{-2}$ when $h>0$ is sufficiently small. Rewriting this in terms of the (holomorphic) parameter $\sigma=h^{-1}z$ gives
  \begin{alignat*}{2}
    \|\wh{P'}(\sigma)^{-1}f\|_{L^2} &\leq C\|f\|_{L^2},&\quad \Re\sigma&>h_0^{-1},\ \ \Im\sigma\in[\beta,\beta'], \\
    \|\wh{P'}(\sigma)^{-1}f\|_{L^2} &\leq C|\sigma|^{-1}\|f\|_{L^2},&\quad \Re\sigma&>h_0^{-1},\ \ \Im\sigma=\beta'.
  \end{alignat*}
  (We strengthen the norm on $f$ on the right from $H_{|\Re\sigma|^{-1}}^{-2}$ to $L^2$ in order to make the norms here independent of $\sigma$.)

  Let now $0<h<h_0$ and consider $\sigma_0$ the holomorphic function $F(\sigma):=\wh{P'}(\sigma+i\eta)^{-1}f$ for $\sigma\in[\sigma_0,5\sigma_0]+i[-\delta_-,\delta_+]$ where $\delta_-=\beta-\eta$, $\delta_+=\beta'-\eta$. We apply the three line theorem \cite[Lemma~D.1]{DyatlovZworskiBook} with $M_+=M=C\|f\|_{L^2}$ and $M_-=C\sigma_0^{-1}$; note that the width $4\sigma_0$ of the rectangle exceeds $\log(\frac{M}{\min(M_-,M_+)})=C\log\sigma_0$ for all sufficiently large $\sigma_0$.\footnote{The subharmonicity of $\log\|F(\sigma)\|$ for holomorphic Hilbert space-valued $F$, used in the proof of \cite[Lemma~D.1]{DyatlovZworskiBook}, follows from $\|F(z_0)\|=|\la F(z_0),\psi\ra|$ for an appropriate unit vector $\psi$, the subharmonicity of $\log|h(z)|$, $h(z)=\la F(z),\psi\ra$, and $|h(z)|\leq\|A(z)\|$.} to get
  \[
    \|\wh{P'}(\sigma)^{-1}f\|_{L^2} \leq C|\sigma|^{-1+\frac{\Im\sigma-\beta}{\beta'-\beta}}\|f\|_{L^2},\quad \sigma\in[2\sigma_0,4\sigma_0]+i[-\delta_-,\delta_+],
  \]
  with the constant $C$ independent of $\sigma_0$. As $\beta\searrow\gamma_-$ and $\beta'\searrow\gamma_+$, the exponent tends to $-1+\delta_\Gamma(\Im\sigma)$ from above. This implies~\eqref{EqSpHiTrPp}.
\end{proof}

\begin{proof}[Proof of Theorem~\usref{ThmSpHiTr} with precise loss]
  We now turn~\eqref{EqSpHiTrPp} into the propagation estimate~\eqref{EqSpHiTr} by following the proof of \cite[Theorem~4.7]{HintzVasyQuasilinearKdS}, which we recall here (using improved notation) for completeness. \textit{We require all auxiliary ps.d.o.s introduced in the argument below to have Schwartz kernels supported in $\chi^{-1}(1)^2$.} Consider $B,G,B^\sharp\in\Psi_\semi$ such that:
  \begin{itemize}
  \item $G$ is elliptic on $\Gamma(\delta_0)$;
  \item $B$ is elliptic at $\Gamma$, and $\WF'_\semi(B)$ is contained in a small neighborhood of $\Gamma$;
  \item $\WF'_\semi(I-B^\sharp)\cap\WF'_\semi(B)=\emptyset$, and $\WF'_\semi(B^\sharp)\subset\Gamma(\frac{\delta_0}{2})$ (which is thus disjoint from the operator wave front set of the complex absorption).
  \end{itemize}
  Writing $P'_{h,z}=P_{h,z}-i Q_{h,z}$, we then have
  \begin{align}
    B u &= B B^\sharp u + B(I-B^\sharp)u \nonumber\\
      &= B(P'_{h,z})^{-1} (P_{h,z}-i Q_{h,z}) B^\sharp u + B(I-B^\sharp)u \nonumber\\
  \label{EqSpHiTrSplit}
      &= B(P'_{h,z})^{-1} B^\sharp P_{h,z}u + B(P'_{h,z})^{-1}[P_{h,z},B^\sharp]u - i B(P'_{h,z})^{-1}Q_{h,z}B^\sharp u + B(I-B^\sharp)u.
  \end{align}
  Since $Q_{h,z}B^\sharp,B(I-B^\sharp)\in h^\infty\Psi_\semi^{-\infty}$, the last two terms have $L^2$-norms bounded by $C h^N\|\chi u\|_{H_h^{-N}}$. By Lemma~\ref{LemmaSpHiTrPp}, the first term is bounded by
  \[
    C h^{-1-\delta_\Gamma} \| B^\sharp P_{h,z}u \|_{L^2} \leq C\Bigl( h^{-1-\delta_\Gamma}\|G P_{h,z}u\|_{L^2} + C_N h^N\|\chi u\|_{H_h^{-N}} \Bigr),
  \]
  where (for any fixed $\eps>0$) we write $\delta_\Gamma=\delta_\Gamma(\Im(h^{-1}z)-\eps)$; this is of the required form for the right-hand side of~\eqref{EqSpHiTr}. In the second term, note that $[P_{h,z},B^\sharp]u$ is microlocalized away from $\Gamma$. The part localized away from $\Gamma^{\rm u}$ can be controlled by real principal type propagation with a priori control term $E$ as in the statement of Theorem~\ref{ThmSpHiTr}: fix $X^{\rm u}\in\Psi_\semi$ with $\WF'_\semi(X^{\rm u})\subset\Gamma(\delta_0)$ and $\WF'(I-X^{\rm u})\cap\Gamma^{\rm u}\cap\Gamma(\frac{\delta_0}{2})=\emptyset$, and such that all backward null-bicharacteristics from $\WF'(I-X^{\rm u})$ enter $\Ell_\semi(E)$ in finite time; then the standard propagation estimate for $P_{h,z}$ yields
  \[
    h\cdot \|(I-X^{\rm u})h^{-1}[P_{h,z},B^\sharp]u\|_{L^2} \leq h\cdot C\Bigl( h^{-1}\|G P_{h,z}\|_{L^2} + \|E u\|_{L^2} + C_N h^N\|\chi u\|_{H_h^{-N}}\Bigr),
  \]
  and therefore Lemma~\ref{LemmaSpHiTrPp} gives
  \begin{equation}
  \label{EqSpHiTrIXu}
    \| B(P'_{h,z})^{-1}(I-X^{\rm u})[P_{h,z},B^\sharp] \|_{L^2} \leq h\cdot C h^{-1-\delta_\Gamma}\Bigl( h^{-1}\|G P_{h,z}\|_{L^2} + \|E u\|_{L^2} + C_N h^N\|\chi u\|_{H_h^{-N}}\Bigr),
  \end{equation}
  which is of the required form for the right-hand side of~\eqref{EqSpHiTr}.

  \begin{figure}[!ht]
  \centering
  \includegraphics{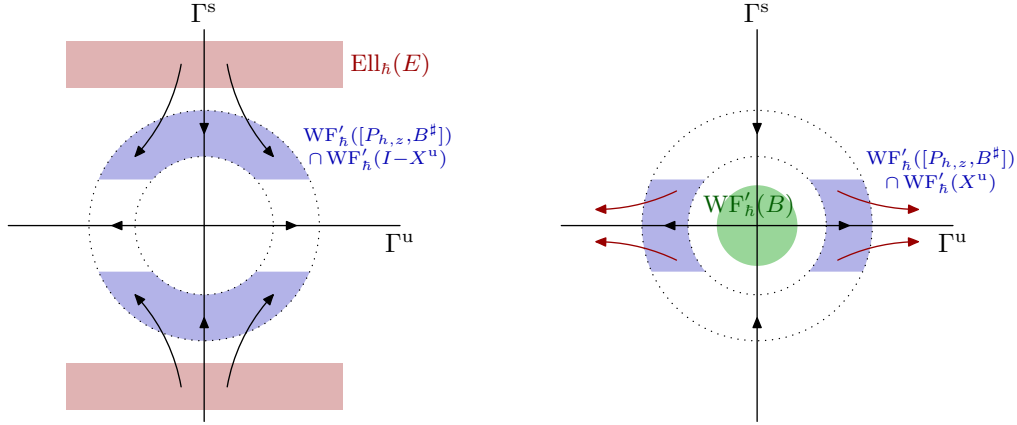}
  \caption{\textit{On the left:} estimating the $L^2$-norm of $B(P'_{h,z})^{-1}X^{\rm u}[P_{h,z},B^\sharp]u$ (see~\eqref{EqSpHiTrIXu}). \textit{On the right:} estimating the $L^2$-norm of $B(P'_{h,z})^{-1}(I-X^{\rm u})[P_{h,z},B^\sharp]u$ (via~\eqref{EqSpHiTrXu}); the red arrows indicate the direction in which semiclassical the wave front set propagates upon applying $(P'_{h,z})^{-1}$ to a function microlocalized in the blue region.}
  \label{FigSpHiTr}
  \end{figure}

  Finally, for the part $B(P'_{h,z})^{-1}X^{\rm u}[P_{h,z},B^\sharp]u$ of the second term of~\eqref{EqSpHiTrSplit} that is localized near $\Gamma^{\rm u}\setminus\Gamma$, we will use that $(P'_{h,z})^{-1}$ propagates semiclassical wave front sets forward. To demonstrate this, we introduce $-i\tilde Q$, where $\tilde Q$ is elliptic near $\Gamma$ and has a non-negative principal symbol, as an additional complex absorption near $\Gamma$, and note that $\|u\|_{L^2}\leq h^{-1}\|(P'_{h,z}-i\tilde Q)u\|_{L^2}$ for small $h>0$ since all backward null-bicharacteristics starting on the characteristic set of $P_{h,z}$ in $\Gamma(\delta_0)$ reach the elliptic set of $Q_{h,z}+\tilde Q$ in finite time. The function
  \[
    v := (P'_{h,z})^{-1}X^{\rm u}[P_{h,z},B^\sharp]u - (P'_{h,z}-i\tilde Q)^{-1}X^{\rm u}[P_{h,z},B^\sharp]u
  \]
  satisfies $P'_{h,z}v=-i\tilde Q(P'_{h,z}-i\tilde Q)^{-1} X^{\rm u}[P_{h,z},B^\sharp]u$; by the direction of propagation for $(P'_{h,z}-i\tilde Q)^{-1}$, this has trivial semiclassical wave front set, so $\|P'_{h,z}v\|\leq C_N h^N\|\chi u\|_{H_h^{-N}}$ and thus also $\|v\|\leq C_N h^N\|\chi u\|_{H_h^{-N}}$ by~\eqref{EqSpHiTrPp}. It remains to observe that
  \begin{equation}
  \label{EqSpHiTrXu}
    \|B (P'_{h,z}-i\tilde Q)^{-1}X^{\rm u}[P_{h,z},B^\sharp]u\|_{L^2} \leq C_N h^N\|\chi u\|_{H_h^{-N}}
  \end{equation}
  if we arrange, as we may, that no backwards null-bicharacteristic starting on $\WF'_\semi(B)$ reaches the operator wave front set of $X^{\rm u}[P_{h,z},B^\sharp]$. (See Figure~\ref{FigSpHiTr}.) This completes the proof.
\end{proof}

\begin{rmk}[Second microlocalization]
\label{RmkSpHiTr2nd}
  One may attempt to prove Theorem~\ref{ThmSpHiTr} more directly by working with an auxiliary equation for ``$\Op_\semi(|\phi^{\rm u}|^\lambda)u$'' for a suitable fractional $\lambda\in(0,1)$; this amounts to second microlocalizing at $\Gamma^{\rm u}$. An analogue of such an argument was successfully carried out by Jia~\cite{JiaTrapping} in a spacetime context, yielding a sharpening of the spacetime trapping estimate \cite[Theorem~3.9]{HintzPolyTrap}.
\end{rmk}

\subsubsection{Energy estimates}
\label{SssSpHiEn}

We finally turn to the semiclassical energy estimates required to ``cap off'' the microlocal estimates in $r>\bhm$ near the interior hypersurface at $r=\bhm$, which (given \cite[Proposition~3.7]{VasyMicroKerrdS}, \cite[Theorem~9.44]{HintzMicro}, \cite[Appendix~E.5.2]{DyatlovZworskiBook}) we only need to prove for $h\ll\Im z\lesssim 1$. We prove a general result here for a general stationary wave operator $P_0$ on $(\cM,g)$ (see Definition~\ref{DefSpGen}), $\cM=\R_t\times\cX$, under the following additional assumption:
\begin{itemize}
\item the function $r\in\CI(\cX)$ is proper, and $r_-<r_+$ are such that $\dd r$ is future timelike on $r^{-1}([r_-,r_+])$.
\end{itemize}
Recall that $\dd t$ is past timelike.

\begin{prop}[High-energy estimate]
\label{PropSpHiEn}
  Let $r_\natural\in(r_-,r_+)$. Let $s,\nu\in\R$. There exists a constant $C$ such that for all $z\in\C$ with $|z|\in[\frac12,2]$ and $\Im z\geq h\nu$, and for all $h\in(0,1)$, we have the direct estimate
  \begin{equation}
  \label{EqSpHiEnDir}
    \|u\|_{\bar H_h^s(r^{-1}([r_-,r_+]))} \leq C\Bigl(h^{-1}\|P_{h,z}u\|_{\bar H_h^{s-1}(r^{-1}([r_-,r_+]))} + \|u\|_{\bar H_h^s(r^{-1}([r_\natural,r_+]))}\Bigr)
  \end{equation}
  and the adjoint estimate
  \begin{equation}
  \label{EqSpHiEnAdj}
    \|u\|_{\bar H_h^s(r^{-1}([r_-,r_+]))} \leq C\Bigl(h^{-1}\|P_{h,z}^*u\|_{\bar H_h^{s-1}(r^{-1}([r_-,r_+]))} + \|u\|_{\bar H_h^s(r^{-1}([r_-,r_\natural]))}\Bigr).
  \end{equation}
\end{prop}

As usual, the direct and the adjoint estimate differ in the direction of propagation.

\begin{proof}[Proof of Proposition~\usref{PropSpHiEn}]
  Write $\Omega=r^{-1}([r_-,r_+])\subset\cX$. For any fixed $\nu'$, the estimates~\eqref{EqSpHiEnDir}--\eqref{EqSpHiEnAdj} follow from \cite[Theorem~9.44]{HintzMicro} when $h\nu\leq\Im z\leq h\nu'$, and thus we only need to consider $h^{-1}\Im z\geq\nu'$ for some arbitrarily large but fixed $\nu'$. We shall, moreover, only prove the a priori estimate
  \begin{equation}
  \label{EqSpHiEnPf}
    \|u\|_{H_h^1(\Omega)^{\bullet,-}} \leq C h^{-1}\|P_{h,z}u\|_{H_h^0(\Omega)^{\bullet,-}},\qquad u\in\cC^2(\Omega)^{\bullet,-},
  \end{equation}
  where ``$\bullet$'' and ``$-$'' indicate the supported character at $r^{-1}(r_+)$ and the extendible character at $r^{-1}(r_-)$, respectively. (Dualization and propagation of regularity arguments as in \cite[Chapter~9.5]{HintzMicro} then yield $\|u\|_{H_h^s(\Omega)^{\bullet,-}}\leq C h^{-1}\|P_{h,z}u\|_{H_h^{s-1}(\Omega)^{\bullet,-}}$; applying this to $\chi u$ where $\chi\in\CI(\Omega)$ equals $1$ on $r^{-1}([r_-,r_\natural])$ and $0$ near $\{r=r_+\}$ yields~\eqref{EqSpHiEnDir}.)

  Writing $T[u](X,Y)=\Re(X u\,\ol{Y u}-\frac12 g(X,Y)g^{-1}(\dd u,\dd\bar u))$ and ${}^{(X)}\!K[u]=\frac12\cL_X g)\cdot T[u]$ (which satisfies ${}^{(w X)}\!K[u]=w({}^{(X)}\!K[u]+T[u](\frac{\nabla w}{w},X))$ by \cite[Proposition~9.13]{HintzMicro}), we apply the energy identity \cite[Theorem~9.16]{HintzMicro} to
  \[
    u' := e^{-i z t/h}u,\quad
    X' := e^{-2(\Im z)t/h}X,\ X=\nabla r
  \]
  on the domain $[0,1]_t\times\Omega\subset\cM$. Note that $X$ is future timelike. With $\nu_t:=-\nabla t/\sqrt{-g(\nabla t,\nabla t)}$ and $\nu_r=\nabla r/\sqrt{-g(\nabla r,\nabla r)}$ denoting the future unit normals at $[0,1]\times r^{-1}(r_-)$ and $\{t\}\times\Omega^\circ\subset\cM$, respectively, this identity reads
  \begin{align*}
    &\int_{\Omega\times[0,1]} {}^{(X')}\!K[u']\,|\dd g| + B(1) + \int_{[0,1]\times r^{-1}(r_-)} T[u'](X',\nu_r)\,\dd\sigma = \int_{\Omega\times[0,1]} \Re(\Box_g u'\,X'u')\,|\dd g| + B(0), \\
    &\qquad B(t) := \int_{\Omega\times\{t\}} T[u'](X',\nu_t)\,\dd\sigma.
  \end{align*}
  Since $g^{-1}(\dd u',\dd\ol{u'})$ and
  \begin{align*}
    \Re(X'u'\,\ol{\nu_t u'}) &= \Re\Bigl[ e^{-2(\Im z)t/h} X( e^{-i z t/h}u ) \ol{\nu_t( e^{-i z t/h}u )}\,\Bigr] \\
      &= \Re\biggl[ \Bigl(X-\frac{i z}{h}(X t)\Bigr)u\cdot\ol{\Bigl(\nu_t-\frac{i z}{h}(\nu_t t)\Bigr)u}\;\biggr]
  \end{align*}
  are independent of $t$, $B(t)$ is constant. Since $T[u](X,\nu_r)$ is a positive definite quadratic form in $\dd u$ at $[0,1]\times r^{-1}(r_0)$, we can drop this term to obtain
  \begin{align}
    \int_{\Omega\times[0,1]} {}^{(X')}\!K[u']\,|\dd g| &\leq \int_{\Omega\times[0,1]} \Re\bigl( e^{-2 i(\Re z)t/h} \wh{\Box_g}(h^{-1}z)u\cdot\hat X(h^{-1}z)u\bigr)\,|\dd g| \nonumber\\
  \label{EqSpHiEnBasic}
      &\leq \int_{\Omega\times[0,1]} \bigl|\wh{\Box_g}(h^{-1}z)u\bigr|\,|\hat X(h^{-1}z)u|\,|\dd g|,
  \end{align}
  where we define $\hat X(\sigma)=e^{i\sigma t}X e^{-i\sigma t}\in\Diff^1(\cX)$ to be the spectral family of $X$; thus $\hat X(h^{-1}z)=\hat X(0)-\frac{i z}{h}(X t)$.

  Note now that
  \[
    {}^{(X')}\!K[u'] = e^{-2(\Im z)t/h}\Bigl( {}^{(X)}\!K[u'] + 2\frac{\Im z}{h} T[u'](X,-\nabla t) \Bigr),
  \]
  with the expression in parentheses being $\geq c\frac{\Im z}{h}$ for some $c>0$ as a quadratic form in $\dd u'$, with respect to any fixed smooth $t$-independent Riemannian metric $R$ on $T^*\cM$, when $\frac{\Im z}{h}$ is large enough. Since $\dd u'=e^{-i z t/h}(-\frac{i z}{h}u\,\dd t+\dd_\cX u)$ where $\dd_\cX$ is the exterior derivative on $\cX$, we get the pointwise bound
  \begin{equation}
  \label{EqSpHiEnKBound}
    {}^{(X')}\!K[u'] \geq c\frac{\Im z}{h}\bigl(|h^{-1}u|^2 + |\dd_\cX u|_R^2\bigr).
  \end{equation}

  As for the right-hand side of~\eqref{EqSpHiEnBasic}, replacing $\Box_g$ by $P_0=\Box_g+V+a$, where $V\in\Diff^1(\cM)$ and $a\in\CI(\cM)$ are $t$-translation-invariant, generates error terms which obey pointwise bounds
  \[
    |\hat V(h^{-1}z)u||\hat X(h^{-1}z)u| \leq C(|h^{-1}u|^2+|\dd_\cX u|_R^2),\quad
    |a u||\hat X(h^{-1}z)u| \leq C h(|h^{-1}u|^2+|\dd_\cX u|_R^2);
  \]
  these can thus be absorbed into~\eqref{EqSpHiEnKBound} for sufficiently large $\frac{\Im z}{h}$.

  In summary, we have shown that if $\nu'$ is sufficiently large, then
  \[
    \|h^{-1}u\|_{L^2}^2 + \|d_\cX u\|^2 \leq C\frac{h}{\Im z} h^{-4} \|P_{h,z}u\|_{L^2}^2,\quad \Im z\geq h\nu'.
  \]
  This yields~\eqref{EqSpHiEnPf} (and indeed a mild strengthening thereof, namely, with an additional factor of $\sqrt{\frac{h}{\Im z}}$ on the right-hand side).
\end{proof}

\subsubsection{Proof of Theorem~\usref{ThmSpHi}}
\label{SssSpHiPf}

The first (direct) estimate in~\eqref{EqSpHi} follows by first using Propositions~\ref{PropSpHiHor} and \ref{PropSpHiInftyIn} (to get microlocal control near $\pa\cR_{\cH^+}$ and the ingoing radial set), then Theorem~\ref{ThmSpHiTr} (to propagate the now obtained control on the future trapped set into the trapped set itself), and finally Proposition~\ref{PropSpHiInftyOut} (to get microlocal control near the outgoing radial set), with real principal type propagation in between. Recall here that Proposition~\ref{PropTse3bDyn} shows, via the phase space relationship discussed in~\S\ref{SssSpSemi}, that these propagation results suffice to give microlocal control in semiclassical scattering Sobolev spaces everywhere on $\ol{\Tsc^*}X\setminus r^{-1}(\bhm)$. To get control near $r^{-1}(\bhm)$, we use Proposition~\ref{PropSpHiEn}. This yields the estimate~\eqref{EqSpHi} except with an additional error term $C_N h^N\|\chi u\|_{\bar H_{\scop,h}^{\sfs_0,\sfr+\alpha_+,\sfb}}$ on the right-hand side, where $\sfs_0<\sfs_\scop$, $\sfr_0<\sfr$. But for fixed $N>0$ and sufficiently small $h$, this error can be absorbed into the left-hand side, yielding~\eqref{EqSpHi} as stated. The proof of the adjoint estimate is similar (with reversed direction and order of propagation). This proves part~\eqref{ItSpHiI}. The proof of part~\eqref{ItSpHiII} is the same.

The invertibility statement in part~\eqref{ItSpHiInv} follows immediately from the estimates~\eqref{EqSpHi}, which imply the triviality of the nullspace of $P_{h,z}$ (and thus $\wh{P_0}(h^{-1}z)$) on $\bar H_\scop^{\sfs_\scop,\sfr+\alpha_+}$ and of $P_{h,z}^*$ (and thus $\wh{P_0}(h^{-1}z)^*$) on $\dot H_\scop^{-\sfs_\scop+1,-\sfr-\alpha_+-1}$ when $|z|\in[\frac12,2]$, $\Im z\geq 0$, and $h>0$ is sufficiently small.

Microlocally away from $\cR_{\rm out}$, control of $k\in\N_0$ degrees of additional b-regularity (part~\eqref{ItSpHib}) is the same as control of $k$ additional scattering regularity, scattering decay, and semiclassical orders each; this is completely analogous to Remark~\ref{RmkSpBNob}. Only at $\cR_{\rm out}$ do we need to use a separate estimate, which is provided by~\eqref{EqSpHiInftyOutb}. This finishes the proof of Theorem~\ref{ThmSpHi}.

For the proof of~\eqref{EqSpHiInvReg}, note that undoing the semiclassical rescaling in~\eqref{EqSpHib} yields the uniform boundedness of
\[
  \wh{P_0}(\sigma)^{-1} \colon \bar H_{(\scop,|\sigma|^{-1});\bop^+}^{(\sfs_\scop-1;k),\sfr+\alpha_++1,\sfb)} \to |\sigma|^{-1+\delta_\Gamma(\Im(\sigma)-\eps)}\bar H_{(\scop,|\sigma|^{-1});\bop^+}^{(\sfs_\scop;k),\sfr+\alpha_+,\sfb)}
\]
for any $\eps>0$. Now, we read off from~\eqref{EqSpBInvRegDiff} that
\[
  \sigma\pa_\sigma\wh{P_0}(\sigma) \in \sigma\rho\Diff_\bop^1 + \sigma^2\rho^2\CI.
\]
Since the norms~\eqref{EqSpHiNormb} treat powers of $\sigma$ on an equal footing with b-derivatives, this is $\sigma\rho$ times an ``extended'' b-derivative. The semiclassical replacement for~\eqref{EqSpBInvRegPf} is thus, for $\ell\geq 0$ and $q\in\N_0$, and setting $\delta:=\delta_\Gamma(\Im(\sigma)-\eps)$,
\begin{alignat*}{2}
  \sigma\pa_\sigma\wh{P_0}(\sigma)\circ\wh{P_0}(\sigma)^{-1} \colon \bar H_{(\scop,|\sigma|^{-1});\bop^+}^{(\sfs_\scop+\ell-1;q),\sfr+\alpha_++1,\sfb} &\;\xra{\wh{P_0}(\sigma)^{-1}}&\ & |\sigma|^{-1+\delta}\bar H_{(\scop,|\sigma|^{-1});\bop^+}^{(\sfs_\scop+\ell;q),\sfr+\alpha_+,\sfb} \\
    &\xra{\sigma\pa_\sigma\wh{P_0}(\sigma)} &\ & |\sigma|^\delta\bar H_{(\scop,|\sigma|^{-1});\bop^+}^{(\sfs_\scop+\ell;q-1),\sfr+\alpha_++1,\sfb}.
\end{alignat*}
We apply this with $(\ell,q)=(0,k+j)$, $(1,k+j-1)$, $\ldots$, $(j-1,k+1)$, and apply $\wh{P_0}(\sigma)$ one final time to get~\eqref{EqSpHiInvReg}.

\subsection{Fredholm property of the zero energy operator}
\label{SsSp0}

We recall the order $\sfs_0\in\CI(\Sb^*X)$ from~\eqref{EqSpOrderb}. The only relevant features of $\sfs_0$ for the analysis in this section are that $\sfs_0$ is monotonically decreasing along the future null-bicharacteristic flow for $\wh{P_0}(0)$, and $\sfs_0>\frac12+\vartheta_{\cH^+}$ at $\pa\cR_{\cH^+}$ (as in~\eqref{EqSpBThrHor}).

\begin{thm}[Spectral family at zero energy]
\label{ThmSp0}
  Let $\alpha\in\R$ be such that $\alpha+\frac32$ lies in an indicial gap $I\subset\R$ of $\wh{P_0}(0)$ (Definition~\usref{DefMUbInd}).
  \begin{enumerate}
  \item\label{ItSp0Fred}{\rm (Fredholm estimates.)} Let $N\in\R$. Then we have the estimates
    \begin{equation}
    \label{EqSp0Est}
    \begin{split}
      \|u\|_{\bar H_\bop^{\sfs_0,\alpha}} &\leq C\Bigl(\|\wh{P_0}(0)u\|_{\bar H_\bop^{\sfs_0-1,\alpha+2}} + \|u\|_{\bar H_\bop^{-N,-N}}\Bigr), \\
      \|u^*\|_{\dot H_\bop^{-\sfs_0+1,-\alpha-2}} &\leq C\Bigl(\|\wh{P_0}(0)^*u^*\|_{\dot H_\bop^{-\sfs_0,-\alpha-2}} + \|u^*\|_{\dot H_\bop^{-N,-N}}\Bigr).
    \end{split}
    \end{equation}
    In particular, the operator
    \begin{equation}
    \label{EqSp0Map}
      \wh{P_0}(0) \colon \{ u\in\bar H_\bop^{\sfs_0,\alpha} \colon \wh{P_0}(0)u\in\bar H_\bop^{\sfs_0-1,\alpha+2} \} \to \bar H_\bop^{\sfs_0-1,\alpha+2}
    \end{equation}
    is Fredholm.\footnote{Its Fredholm index is $0$ when $P_0$ is $\tface$-admissible with weight $\alpha+\frac32$; see Corollary~\ref{CorSpLoInd0}.}
  \item\label{ItSp0Reg}{\rm (Regularity of kernel elements.)} Let $u\in\bar H_\bop^{\sfs_0,\alpha}$ with $\wh{P_0}(0)u=0$. Then we have $u\in\cA^\beta(X)$ for all $\beta<\sup I$.
  \end{enumerate}
\end{thm}

\begin{rmk}[Improvements]
\label{RmkSp0Impr}
  If one replaces the orders of the error term in the first estimate in~\eqref{EqSp0Est} by $\sfs^\flat,\alpha^\flat$, where $\sfs^\flat$ is induced by an admissible order function $<\sfs$ and $\alpha^\flat<\alpha$ are such that $\alpha^\flat+\frac32$ and $\alpha+\frac32$ lie in the same indicial gap, then the estimate holds in the strong form. An analogous comment holds for the second (adjoint) estimate. One can strengthen part~\eqref{ItSp0Reg} to the full polyhomogeneity of $u$ using standard b-analysis techniques (see, e.g., \cite[Proposition~5.61]{MelroseAPS}).
\end{rmk}

\begin{proof}[Proof of Theorem~\usref{ThmSp0}]
  The first step in the proof of the first (direct) estimate~\eqref{EqSp0Est} is to concatenate a radial point estimate at the event horizon (Proposition~\ref{PropSpBHor} with $\sigma=0$), real principal type propagation in the rest of the ergoregion, and a standard energy estimate near $r=\bhm$. Together with elliptic regularity estimates, this gives
  \[
      \|u\|_{\bar H_\bop^{\sfs_0,\alpha}} \leq C\Bigl(\|\wh{P_0}(0)u\|_{\bar H_\bop^{\sfs_0-1,\alpha+2}} + \|u\|_{\bar H_\bop^{-N,\alpha}}\Bigr).
  \]
  We then use~\eqref{EqMUbEstNearInfty} (which uses that $\alpha+\frac32$ lies in an indicial gap) to bound the second, error, term by $\|\wh{P_0}(0)u\|_{\bar H_\bop^{-N-2,\alpha+2}}+\|\chi u\|_{H_\bop^{-N,\alpha-1}}$ where $\chi\in\CI(X)$ is $1$ near, and supported near, $\pa X$. The simple interpolation inequality $\|\chi u\|_{H_\bop^{-N,\alpha-1}}\leq \eps\|\chi u\|_{H_\bop^{-N,\alpha}}+C_\eps\|\chi u\|_{H_\bop^{-N,-N}}$ yields~\eqref{EqSp0Est}.

  To prove part~\eqref{ItSp0Reg}, one first uses that the aforementioned microlocal estimates hold in the strong form, and thus imply that every kernel element lies in $\bar H_\bop^{\infty,\alpha}$. The improvement of the weight follows from a standard contour shifting argument: rewriting then $\wh{P_0}(0)u=0$ as
  \begin{equation}
  \label{EqSp0Impr}
    N(\rho^{-2}\wh{P_0}(0))(\chi u) = -\bigl(\rho^{-2}\wh{P_0}(0)-N(\rho^{-2}\wh{P_0}(0))\bigr)(\chi u) + \bigl[\rho^{-2}\wh{P_0}(0),\chi\bigr]u \in \bar H_\bop^{\infty,\alpha+1},
  \end{equation}
  one passes to the Mellin transform and uses the absence of poles of $N(\rho^{-2}\wh{P_0}(0),\lambda)$ for $\Re\lambda\in[\alpha+\frac32,\min(\alpha+\frac52,\beta)]$ to improve the decay of $\chi u$ to $\min(\alpha+1,\beta-\frac32)$. After finitely many iterations, one obtains $u\in\bar H_\bop^{\infty,\beta-\frac32}\subset\cA^\beta$ by Sobolev embedding (see~\eqref{EqMUCbSob}).
\end{proof}

\subsection{Transition face normal operator}
\label{SsSptf}

For a stationary asymptotically Minkowski wave operator $P_0$ on $\R\times\cX^\circ$, with $\cX\subset\ol{\R^n}$ a neighborhood of $\pa\ol{\R^n}$, we now study the transition face normal operator of its spectral family, as defined in~\eqref{EqSSSpecFam}--\eqref{EqSStfOp}. For $\hat\sigma\in e^{i[0,\pi]}$, this is given by
\begin{equation}
\label{EqSptfOp0}
  N_\tface(P_0,\hat\sigma) = 2 i\hat\sigma\hat\rho\Bigl(\hat\rho\pa_{\hat\rho}-\frac{n-1}{2}-S\Bigr) + \hat\rho^2 P_{(0)}(0,\omega,\hat\rho\pa_{\hat\rho},\pa_\omega) \in \rho_\ztface^{-2}\Diff_{\scop,\bop}^2(\tface),
\end{equation}
where $P_0(\rho,\omega,\rho\pa_\rho,\pa_\omega)=\rho^2 P_{(0)}(\rho,\omega,\rho\pa_\rho,\pa_\omega)$, and $\hat\rho=\frac{\rho}{|\sigma|}$ is a projective coordinate on $\tface=[0,\infty]_{\hat\rho}\times\Sph^2$; moreover, we set $\rho_\ztface:=\frac{1}{1+\hat\rho}$ and $\rho_\sctface:=\frac{\hat\rho}{1+\hat\rho}$. We recall the scattering-b-Sobolev spaces
\[
  H_{\scop,\bop}^{s,(r,q)}(\tface) = \rho_\sctface^r\rho_\ztface^q H_{\scop,\bop}^s(\tface)
\]
from~\eqref{EqMUNtfSob}, with the orders $s\in\CI({{}^{\scop,\bop}T^*}\tface)$ and $r\in\CI(\ol{{}^{\scop,\bop}T^*_\sctface}\tface)$ allowed to be variable. We write $H_{\scop,\bop;\bop}^{(s;k),(r,q)}$ for the space of distributions with $k\in\N_0$ additional degrees of b-regularity.

Recall the orders $\sfs,\sfr_\pm$ from~\eqref{EqSpOrdertf}. In the result below, due to the ellipticity of $N_\tface(P_0,\hat\sigma)$ at fiber infinity, the regularity order can, in fact, be taken to be arbitrary; we denote it by $\sfs$; and for the decay order, only the monotonicity and threshold conditions at the radial sets in~\eqref{EqSpOrderRad}, i.e., the analogues of~\eqref{EqSpBThrIn}--\eqref{EqSpBThrOut} (with $\alpha_+$ as in~\S\ref{SsSpTs}), are needed.

\begin{thm}[Fredholm property of the tf-normal operator]
\label{ThmSptf}
  Let $q\in\R$ be such that $-q+\frac{n}{2}$ lies in an indicial gap $(\beta^-,\beta^+)\subset\R$ of $\wh{P_0}(0)$.
  \begin{enumerate}
  \item\label{ItSptfReal}{\rm (Fredholm estimates, I.)} Let $\theta\in[0,\pi]$ be such that $|e^{i\theta}\pm 1|\leq\frac12$. Let $N\in\R$, and write $\sfr=\sfr_\pm$. Then there exists a constant $C$ such that
    \begin{equation}
    \label{EqSptfEst}
    \begin{split}
      \|u\|_{H_{\scop,\bop}^{\sfs,(\sfr+\alpha_+,q)}} &\leq C\Bigl( \|N_\tface(P_0,e^{i\theta})u\|_{H_{\scop,\bop}^{\sfs-2,(\sfr+\alpha_++1,q-2)}} + \|u\|_{H_{\scop,\bop}^{-N,(-N,-N)}} \Bigr), \\
      \|u\|_{H_{\scop,\bop}^{-\sfs+2,(-\sfr-\alpha_+-1,-q+2)}} &\leq C\Bigl( \|N_\tface(P_0,e^{i\theta})^*u\|_{H_{\scop,\bop}^{-\sfs,(-\sfr-\alpha_+,-q)}} + \|u\|_{H_{\scop,\bop}^{-N,(-N,-N)}} \Bigr).
    \end{split}
    \end{equation}
  \item\label{ItSptfImag}{\rm (Fredholm estimates, II.)} Let $I\subset(0,\pi)$ be a compact set. Let $N\in\R$. Let $\sfr\in\CI(\ol{{}^{\scop,\bop}T^*_\sctface}\tface)$ be such that $\sfr+\alpha_+<-\frac12+\ubar S$ at $\cR_{\rm out}$, where $\ubar S=\inf_{\pa\cX}S$. Then the estimates~\eqref{EqSptfEst} hold uniformly for $\theta\in I$.
  \item\label{ItSptfIndex}{\rm (Fredholm index.)} In the settings of~\eqref{ItSptfReal} and~\eqref{ItSptfImag}, the Fredholm index of
    \begin{equation}
    \label{EqSptfOp}
    \begin{split}
      N_\tface(P_0,e^{i\theta}u) &\colon \cX_{\tface,\theta} := \bigl\{ u\in H_{\scop,\bop}^{\sfs,(\sfr+\alpha_+,q)} \colon N_\tface(P_0,e^{i\theta})u \in H_{\scop,\bop}^{\sfs-2,(\sfr+\alpha_++1,q-2)} \bigr\} \\
        &\quad \to \cY_{\tface,\theta} := H_{\scop,\bop}^{\sfs-2,(\sfr+\alpha_++1,q-2)}.
    \end{split}
    \end{equation}
    is independent of $\theta\in[0,\pi]$ and $q\in(\frac{n}{2}-\beta^+,\frac{n}{2}-\beta^-)$. In fact, its nullspace as well as the nullspace of the adjoint on $H_{\scop,\bop}^{-s+2,(-\sfr-\alpha_+-1,-q+2)}$ are independent of $q$ in this interval.
  \item\label{ItSptfb}{\rm (Higher b-regularity.)} If $u\in\cX_{\tface,\theta}$ and $N_\tface(P_0,e^{i\theta})u\in H_{\scop,\bop;\bop}^{(\sfs-2;k),(\sfr+\alpha_++1,q-2)}$, then $u\in H_{\scop,\bop;\bop}^{(\sfs;k),(\sfr+\alpha_+,q)}$, and we have an estimate
    \begin{equation}
    \label{EqSptfbEst}
      \|u\|_{H_{\scop,\bop;\bop}^{(\sfs;k),(\sfr+\alpha_+,q)}} \leq C\Bigl( \|N_\tface(P_0,e^{i\theta})u\|_{H_{\scop,\bop;\bop}^{(\sfs-2;k),(\sfr+\alpha_++1,q-2)}} + \|u\|_{H_{\scop,\bop}^{\sfs,(\sfr+\alpha_+,q)}}\Bigr).
    \end{equation}
  \item\label{ItSptfKer}{\rm (Regularity of kernel elements.)} In the settings of~\eqref{ItSptfReal} and \eqref{ItSptfImag}, every $u\in H_{\scop,\bop}^{\sfs,(\sfr+\alpha_+,q)}$ with $N_\tface(P_0,e^{i\theta})u=0$ satisfies
    \begin{equation}
    \label{EqSptfKer}
      u \in \cA^{\frac{n-1}{2}+\ubar S-\eps,-\beta}(\tface)\quad\forall\,\eps>0,\ \beta\in(\beta^-,\beta^+).
    \end{equation}
    (Conversely, the space $\bigcap_{\eps>0}\cA^{\frac{n-1}{2}+\ubar S-\eps,-\beta^--\eps}(\tface)$ is included in the domain of~\eqref{EqSptfOp}.) Furthermore, every $u^*\in H_{\scop,\bop}^{-\sfs+2,(-\sfr-\alpha_+-1,-q+2)}$ with $N_\tface(P_0,1)^*u^*=0$ satisfies
    \begin{equation}
    \label{EqSptfKerAdj}
      u^* \in e^{-2 i/\rho_\sctface}\cA^{\alpha,-n+2+\beta}(\tface)
    \end{equation}
    for all $\beta\in(\beta^-,\beta^+)$ and for some $\alpha\in\R$. (Conversely, every element of the space $\bigcap_{\eps>0} e^{-2 i/\rho_\sctface}\cA^{\alpha,-n+2+\beta^+-\eps}$ is included in the domain $H_{\scop,\bop}^{-\sfs+2,(-\sfr-\alpha_+-1,-q+2)}$ of $N_\tface(P_0,1)^*$ when $\sfr|_{\cR_{\rm in,1}}+\alpha_+>\frac{n-2}{2}-\alpha$, i.e., when \eqref{EqSpBThrIn} holds with sufficient margin.)
  \end{enumerate}
\end{thm}
\begin{proof}
  We combine the analysis near infinity in the bounded nonzero frequency regime (\S\ref{SsSpB}) with the b-analysis of the zero energy operator (\S\ref{SsSp0}). Thus, we shall be brief.

  The proofs of parts~\eqref{ItSptfReal} and \eqref{ItSptfImag} are the same. We only consider the direct estimate. The operator $N_\tface(P_0,e^{i\theta})$ is not elliptic (as an operator of class $\rho_\ztface^{-2}\Diff_{\scop,\bop}(\tface)$) only over $\sctface$, while over $\sctface$ the characteristic set and null-bicharacteristic flow are as in the case of frequency $e^{i\theta}$ as discussed in~\S\ref{SsSpB}, radial point estimates (Propositions~\ref{PropSpBInftyIn} and \ref{PropSpBInftyOut}) and real principal type propagation yield the first estimate in~\eqref{EqSptfEst} except the error term is $\|u\|_{H_{\scop,\bop}^{-N,(-N,q)}}$.

  Now, $\hat r:=\hat\rho^{-1}$ is a boundary defining function of $\ztface\subset\tface$, and the b-normal operator of $N_\tface(P_0,e^{i\theta})$ at $\ztface$ is $\hat r^{-2}P_{(0)}(0,\hat r^{-1},-\hat r\pa_{\hat r},\pa_\omega)$. Therefore, $\lambda$ is an indicial root for it if and only if $-\lambda$ is an indicial root of $\wh{P_0}(0)$ (i.e., of $P_{(0)}$). Under the stated condition on $q$, we can thus use the normal operator estimate~\eqref{EqMUbEstNearInfty} to improve the error term to $\|u\|_{H_{\scop,\bop}^{-N,(-N,q-1)}}$, and then~\eqref{EqSptfEst} follows by an interpolation equality (cf.\ the proof of Theorem~\ref{ThmSp0}).

  Part~\eqref{ItSptfIndex} is proved using the same arguments as around~\eqref{EqSpBPfwtP} (with the analogues of the operators $E_{1 0}$ and $E_{0 1}$ now being tuples of elements of $\CIc(\tface^\circ)$ and of $L^2$-pairings with elements of $\CIc(\tface^\circ)$, respectively). The independence of the kernel and cokernel of $q$ follow from part~\eqref{ItSptfKer}, proved below.

  The proof of part~\eqref{ItSptfb} is the same as that of parts~\eqref{ItSptfReal}--\eqref{ItSptfImag}, now using the b-regular version of Proposition~\ref{PropSpBInftyOut}.

  The first half of part~\eqref{ItSptfKer} follows from the infinite order b-regularity of $u$ (which is a consequence of part~\eqref{ItSptfb}) and normal operator arguments: near $\sctface$ as around~\eqref{EqSpBPfImpr}, near $\ztface$ as around~\eqref{EqSp0Impr}. Note that Sobolev embedding reads $\Hb^{\infty,q}([0,1)_{\hat\rho};\hat r^{n-1}\,|\dd\hat r|)\subset\cA^{q-\frac{n}{2}}([0,1)_{\hat\rho})$. The parenthetical remark is proved similarly to that in Theorem~\ref{ThmSpB}\eqref{ItSpBNull}. For the second half, i.e., the membership~\eqref{EqSptfKerAdj}, we first characterize the indicial roots of $N_\tface(P_0,1)^*$. Note that the adjoint here is taken with respect to the Euclidean density $\mu_{\rm e}:=\hat r^{n-1}\,|\dd\hat r\,\dd\slg|=\hat\rho^{-n}\mu_\bop$, $\mu_\bop:=|\frac{\dd\hat\rho}{\hat\rho}\,\dd\slg|$, so making the densities explicit,
  \[
    N_\tface(P_0,1)^{*,\mu_{\rm e}} = \hat\rho^n N_\tface(P_0,1)^{*,\mu_\bop}\hat\rho^{-n},
  \]
  the indicial roots of which are $-1$ times those of
  \[
    \hat\rho^n P_{(0)}^{*,\mu_\bop}\hat\rho^2 \hat\rho^{-n}.
  \]
  If $\lambda$ is an indicial root of this, then $\lambda-n+2$ is an indicial root of $P_{(0)}^{*,\mu_\bop}$, and thus $-\ol{\lambda-n+2}=-\bar\lambda+n-2$ is an indicial root of $P_{(0)}$. That is, the indicial roots of $N_\tface(P_0,1)$ are given by $-n+2+\bar\lambda$ where $\lambda$ runs over the indicial roots of $\wh{P_0}(0)$. In view of the infinite order scattering-b-regularity of $u^*$ (which follows from the ellipticity of $N_\tface(P_0,1)^*$ at fiber infinity), we thus obtain $u^*|_{\{\hat r<\infty\}}\in\cA^{-n+2+\beta}$. It remains to analyze $u^*$ near $\sctface$. In view of the radial point estimate in Proposition~\ref{PropSpBInftyOut}\eqref{ItSpBInftyOutAdj}, $u^*$ has infinite scattering decay microlocally near $\cR_{\rm out}$ and thus, by propagation of regularity (here: scattering decay), microlocally away from $\cR_{1,{\rm in}}$. We need to show that
  \begin{equation}
  \label{EqSptfAdjCon}
    \tilde u:=e^{2 i/\hat\rho}u^*\ \text{is conormal at}\ \hat\rho=0.
  \end{equation}
  We have
  \[
    0 = e^{2 i/\hat\rho}N_\tface(P_0,1)^* e^{-2 i/\hat\rho}\tilde u.
  \]
  (On the scattering phase space level, conjugation by $e^{2 i\hat r}$ amounts to pullback along $\xi\mapsto\xi-2\,\dd\hat r$, which thus shifts the incoming radial set $\cR_{1,{\rm in}}$ to the zero section.) We claim that
  \begin{equation}
  \label{EqSptfAdjOp}
    e^{2 i/\hat\rho}N_\tface(P_0,1)^* e^{-2 i/\hat\rho} = -2 i\hat\rho\Bigl(\hat\rho\pa_{\hat\rho}-\frac{n-1}{2}-\tilde S\Bigr) + \hat\rho^2\tilde P_{(0)}(0,\omega,\hat\rho\pa_{\hat\rho},\pa_\omega)
  \end{equation}
  where $\tilde S\in\CI(\pa X)$ and $\Diff_\bop^2\ni\tilde P_{(0)}\equiv(\hat\rho D_{\hat\rho})^2+\slDelta\bmod\Diff_\bop^1$; this thus has the same structure as $N_\tface(P_0,-1)$. To prove this, note that $(\hat\rho\pa_{\hat\rho})^{*,\mu_{\rm b}}=-\hat\rho\pa_{\hat\rho}$ and thus
  \begin{align*}
    e^{2 i/\hat\rho} \Bigl[ 2 i\hat\rho\Bigl(\hat\rho\pa_{\hat\rho}-\frac{n-1}{2}-S\Bigr)\Bigr]^{*,\mu_{\rm e}} e^{-2 i/\hat\rho} &= -e^{2 i/\hat\rho}2 i\hat\rho\,\hat\rho^{n-1}\Bigl(\hat\rho\pa_{\hat\rho}-\frac{n-1}{2}-S\Bigr)^{*,\mu_\bop}\hat\rho^{-n+1}e^{-2 i/\hat\rho} \\
      &= 2 i\hat\rho\Bigl(\hat\rho\pa_{\hat\rho}-\frac{n-1}{2}+\bar S\Bigr) - 4.
  \end{align*}
  Furthermore, the contribution of $\hat\rho^2\slDelta$ to~\eqref{EqSptfAdjOp} is $\hat\rho^2\slDelta$ again; and a lower order term of the schematic form $\hat\rho^2(\hat\rho\pa_{\hat\rho},\pa_\omega,1)^{\leq 1}$ contributes $\hat\rho^2(\hat\rho\pa_{\hat\rho}+2 i\hat\rho^{-1},\pa_\omega,1)^{\leq 1}$, so a linear combination (with $\CI(\Sph^{n-1})$ coefficients) of $\hat\rho^2\hat\rho\pa_{\hat\rho}$, $\hat\rho^2\pa_\omega$, $\hat\rho^2$ (put into $\hat\rho^2\tilde P_{(0)}$), and $2 i\hat\rho$ (put into $\tilde S$). Finally, the term $(\hat\rho D_{\hat\rho})^2=-(\hat\rho\pa_{\hat\rho})^2$ of $P_{(0)}$ contributes
  \begin{align*}
    -e^{2 i/\hat\rho} \hat\rho^n(\hat\rho\pa_{\hat\rho})^2\hat\rho^{-n}\hat\rho^2 e^{-2 i/\hat\rho} &= -\hat\rho^2(\hat\rho\pa_{\hat\rho}-n+2+2 i\hat\rho^{-1})^2 \\
      &= -\hat\rho^2(\hat\rho\pa_{\hat\rho}-n+2)^2 - 2 i\hat\rho(\hat\rho\pa_{\hat\rho}-n+1) + 4.
  \end{align*}
  Summing yields~\eqref{EqSptfAdjOp}. The radial point estimates with b-regularity in Proposition~\ref{PropSpBInftyOut}\eqref{ItSpBInftyOutDir} imply~\eqref{EqSptfAdjCon}.
\end{proof}

\begin{cor}[Consequences of tf-admissibility]
\label{CorSptfAdm}
  Suppose $P_0$ is $\tface$-admissible with weight $\beta=-q+\frac{n}{2}\in(\beta^-,\beta^+)$ in the sense of Definition~\usref{DefSStfAdm} except with $-1+\beta$ in~\eqref{EqSStfAdmA2} replaced by $-n+2+\beta$. Then it is $\tface$-admissible for \emph{all} weights in $(\beta^-,\beta^+)$. Moreover, the operator~\eqref{EqSptfOp} is invertible for all $\theta\in[0,\pi]$ and $q\in(\frac{n}{2}-\beta^+,\frac{n}{2}-\beta^-)$, and the estimates~\eqref{EqSptfEst} and~\eqref{EqSptfbEst} hold without the compact error terms: the direct estimates read
  \begin{equation}
  \label{EqSptfAdm}
    \|u\|_{H_{\scop,\bop;\bop}^{(\sfs;k),(\sfr+\alpha_+,q)}} \leq C_k\|N_\tface(P_0,e^{i\theta})u\|_{H_{\scop,\bop;\bop}^{(\sfs-2;k),(\sfr+\alpha_++1,q-2)}},\quad k\in\N_0.
  \end{equation}
\end{cor}
\begin{proof}
  First, we show that~\eqref{EqSptfOp} is invertible for $\theta=0$. Every element of its nullspace lies in the space~\eqref{EqSStfAdmA} in view of~\eqref{EqSptfKer} and thus vanishes. Every element $u^*$ in the $L^2$-orthocomplement of its range satisfies $N_\tface(P_0,e^{i\theta})*u^*=0$ and $u^*\in H_{\scop,\bop}^{-\sfs+2,-\sfr-\alpha_+-1,-q+2}$ and thus lies in $e^{-2 i/\rho_\sctface}\cA^{\alpha,-n+2+\beta}$ by~\eqref{EqSptfKerAdj}. Next, the invertibility of~\eqref{EqSptfOp} for all $\theta\in[0,\pi]$ follows from two facts: its Fredholm index is $0$ by Theorem~\ref{ThmSptf}\eqref{ItSptfIndex}, and it is injective in view of~\eqref{EqSptfKer} and the assumption of $\tface$-admissibility.

  A simple functional analytic argument (using the compactness of the inclusion $H_{\scop,\bop}^{\sfs,(\sfr+\alpha_+,q)}\hra H_{\scop,\bop}^{-N,(-N,-N)}$) shows that the direct estimate~\eqref{EqSptfEst} and the triviality of the nullspace of $N_\tface(P_0,e^{i\theta})$ imply~\eqref{EqSptfAdm} for $k=0$. For general $k\in\N$, one combines~\eqref{EqSptfbEst} with~\eqref{EqSptfAdm} for $k=0$.
\end{proof}

An important example of a $\tface$-admissible operator is the wave operator $P_0=\Box_g$ on an asymptotically Minkowski spacetime (such as the Minkowski spacetime itself, with $g=-\dd t_*^2-2\,\dd t_*\,\dd r+r^2\slg$): in this case,
\[
  \hat P(\hat\sigma) := e^{i\hat\sigma\hat r}N_\tface(\Box_g,\hat\sigma)e^{-i\hat\sigma\hat r} = \Delta_{\rm e} - \hat\sigma^2 = -\pa_{\hat r}^2 - \frac{n-1}{\hat r}\pa_{\hat r} + \hat r^{-2}\slDelta - \hat\sigma^2,
\]
where $0\leq\Delta_{\rm e}$ is the Euclidean Laplacian and $0\leq\slDelta$ is the non-negative Laplacian on $\Sph^{n-1}$.

\begin{prop}[tf-admissibility: asymptotically Minkowski case]
\label{PropSptfAdm}
  Suppose the spatial dimension satisfies $n\geq 3$. Then scalar wave operator $\Box_g$, or indeed the tensor wave operator acting on sections of any tensor bundle, on an asymptotically Minkowski spacetime is $\tface$-admissible with all weights $\beta\in(0,n-2)$.
\end{prop}
\begin{proof}
  Since the tensor wave operator on Minkowski space acts component-wise, in the trivialization of the tensor bundle given by coordinate derivatives and differentials, as the scalar wave operator, it suffices to consider the scalar case. We recall the proof in this case from \cite[Proposition~6.1]{HintzNonstat}. Let $u\in\ker\hat P(\hat\sigma)$ where $\hat\sigma=e^{i\theta}$. Upon replacing $u$ by its projection to degree $\ell\in\N_0$ spherical harmonics, we have
  \[
    -u'' - \frac{n-1}{\hat r}u' + \frac{\lambda_\ell}{\hat r^2}u - e^{2 i\theta}u = 0,\quad \lambda_\ell:=\ell(\ell+n-2)\in\spec\slDelta.
  \]
  This is a Bessel ODE, so the solution is a linear combination of Hankel functions,
  \[
    u = \sum_{j=1}^2 c_j \hat r^{-\frac{n-2}{2}} H_{\nu_\ell}^{(j)}(e^{i\theta}\hat r),\quad \nu_\ell:=\frac{n-2}{2}+\ell.
  \]
  Suppose first that $e^{-i\hat\sigma\hat r}u\in\cA^{\alpha,-\beta}(\tface)$ for some $\alpha\in\R$. If $\theta\in(0,\pi)$, then $u$ is exponentially decaying as $\hat r\to\infty$, and hence $c_2=0$ since $H_{\nu_\ell}^{(j)}(e^{i\theta}\hat r)$ is exponentially growing for $j=2$ but decaying for $j=1$ by \cite[\S{7.4.1}, (4.03)--(4.04)]{OlverSpecial}. But since, near $\hat r=0$, the bound $|\hat r^{-\frac{n-2}{2}}H_{\nu_\ell}^{(1)}(e^{i\theta}\hat r)|\gtrsim \hat r^{-\frac{n-2}{2}-\nu_\ell}\gtrsim\hat r^{-(n-2)}$ from \cite[\S{7.4.2}, (4.12) and \S{2.9.3}, (9.09)]{OlverSpecial} is incompatible with $\cO(\hat r^{-\beta})$ asymptotics as $\hat r\to 0$, we must have $c_1=0$, so $u=0$. When $\theta=0$, then again $c_2=0$ since $H_{\nu_\ell}^{(j)}(\hat r)$ has $\hat r^{-\frac{n-1}{2}}e^{\pm i\hat r}$ asymptotics as $\hat r\to\infty$, with the ``$+$'' sign for $j=1$ and the ``$-$'' sign for $j=2$; and then again $c_1=0$ by considering the $\hat r\to 0$ asymptotics. The case $\theta=\pi$ is similar.

  Finally, suppose that $e^{i\hat r}u\in\cA^{\alpha,-n+2+\beta}(\tface)$ for some $\alpha\in\R$ and $u\in\ker\hat P(1)^*=\ker\hat P(1)$. Then $u$ is again of the above form, and now we conclude that $c_1=0$ by considering the large $\hat r$ asymptotics, and then $c_2=0$ by considering the small $\hat r$ asymptotics.
\end{proof}

\subsection{Uniform low-energy estimates}
\label{SsSpLo}

Recall the order~\eqref{EqSpOrderStart} and the induced orders $\sfs$ and $\sfr_\pm$ from~\eqref{EqSpOrderscbt}. We use the function spaces~\eqref{EqMUscbtSob} and~\eqref{EqMUscbtbSob}.

\begin{thm}[Uniform low-energy estimates]
\label{ThmSpLo}
  Let $\alpha\in\R$ be such that the zero energy operator $\wh{P_0}(0)$ is invertible as a map~\eqref{EqSp0Map}. Suppose moreover that $P_0$ is $\tface$-admissible with weight $\alpha+\frac32$. Then there exists $c>0$ such that for all $0\neq\sigma\in\C$ with $\Im\sigma\geq 0$ and $|\sigma|\leq c$, the operator $\wh{P_0}(\sigma)$ is invertible as a map~\eqref{EqSpBMap}. Furthermore, for such $\sigma$, and for all $k\in\N_0$, we have the uniform estimates
  \begin{equation}
  \label{EqSpLoEst}
    \|u\|_{H_{(\scbtop,|\sigma|);\bop}^{(\sfs;k),(\sfr+\alpha_+,\alpha,0)}} \leq C\|\wh{P_0}(\sigma)u\|_{H_{(\scbtop,|\sigma|);\bop}^{(\sfs-1;k),(\sfr+\alpha_++1,\alpha+2,0)}}
  \end{equation}
  where $\sfr=\sfr_+$, resp.\ $\sfr_-$ when $\arg\sigma\in[0,\frac{\pi}{4}]$, resp.\ $[\frac{3\pi}{4},\pi]$, and $\sfr$ is arbitrary except for satisfying the threshold condition~\eqref{EqSpBThrOut} at the outgoing radial set over $\pa X$. Finally, $j$-fold conormal derivatives of the resolvent in $\sigma$ are uniformly bounded as maps
  \begin{equation}
  \label{EqSpLoInvReg}
    (\sigma\pa_\sigma)^j\wh{P_0}(\sigma)^{-1} \colon H_{(\scbtop,|\sigma|);\bop}^{(\sfs-1;k+j),(\sfr+\alpha_++1,\alpha+2,0)} \to H_{(\scbtop,|\sigma|);\bop}^{(\sfs+j;k),(\sfr+\alpha_+,\alpha,0)}.
  \end{equation}
  These conclusions hold for \emph{all} $\alpha$ lying in the interval $(-\frac32+\beta^-,-\frac32+\beta^+)$, uniformly when $\alpha$ lies in a compact subinterval thereof.
\end{thm}
\begin{proof}
  We first discuss the case $k=0$, which is a variation of \cite[Lemma~7.5]{Hintz3b} and \cite[Proposition~4.9]{HintzNonstat}. Consider
  \[
    \sigma=\varsigma e^{i\theta},\quad \theta\in[0,\pi].
  \]

  Let $\chi_\zface=\chi_0(\frac{\varsigma}{\rho})$ where $\chi_0\in\CIc([0,1))$ equals $1$ near $0$; thus $\chi_\zface\in\CI(X_\scbtop)$ is a cutoff to a neighborhood of $\zface\subset X_\scbtop$. Then, using the first norm equivalence in~\eqref{EqMUscbtNormEquiv},
  \begin{align*}
    \|\chi_\zface u\|_{H_{\scbtop,\varsigma}^{\sfs,(\sfr+\alpha_+,\alpha,0)}} &\sim \|\chi_\zface u\|_{H_\bop^{\sfs,\alpha}} \\
      &\leq C\|\wh{P_0}(0)(\chi_\zface u)\|_{H_\bop^{\sfs-1,\alpha+2}} \sim \|\wh{P_0}(0)(\chi_\zface u)\|_{H_{\scbtop,\varsigma}^{\sfs-1,(-N,\alpha_++2,0)}}
  \end{align*}
  for any $N$. We then commute $\chi_\zface$ through $\wh{P_0}(0)$ and replace $\wh{P_0}(0)$ by $\wh{P_0}(\sigma)$ and bound the resulting error terms using
  \[
    [\wh{P_0}(0),\chi_\zface]\in\rho_\scface^\infty\rho_\zface^\infty\rho_\tface^2\Diff_\scbtop^1,\quad
    \wh{P_0}(\sigma)-\wh{P_0}(0) \in \rho_\tface^2\rho_\zface\Diff_\scbtop^1.
  \]
  (The first order differential nature in the second membership follows from the independence of the principal symbol of $\wh{P_0}(\sigma)$ over $X^\circ$ on $\sigma$.) Thus,
  \[
    \|\wh{P_0}(0)(\chi_\zface u)\|_{H_{\scbtop,\varsigma}^{\sfs-1,(-N,\alpha_++2,0)}} \leq C\Bigl(\|\wh{P_0}(\sigma)u\|_{H_{\scbtop,\varsigma}^{\sfs-1,(-N,\alpha_++2,0)}} + \|u\|_{H_{\scbtop,\varsigma}^{\sfs,(-N,\alpha_+,-\eta)}}\Bigr)
  \]
  for any fixed $\eta\in(0,1]$. Since $(1-\chi_\zface)u$ is supported away from $\zface$, we can thus conclude that
  \begin{equation}
  \label{EqSpLo1}
    \|u\|_{H_{(\scbtop,\varsigma);\bop}^{\sfs,(\sfr+\alpha_+,\alpha,0)}} \leq C\Bigl(\|\wh{P_0}(\sigma)u\|_{H_{(\scbtop,\varsigma);\bop}^{\sfs-1,(\sfr+\alpha_++1,\alpha+2,0)}} + \|u\|_{H_{(\scbtop,\varsigma);\bop}^{\sfs,(\sfr+\alpha_+,\alpha,-\eta)}}\Bigr).
  \end{equation}

  Let now $\chi_\tface=\chi_0(\varsigma)\chi_0(\rho)$ be a cutoff to a neighborhood of $\tface\subset X_\scbtop$. Using now the second norm equivalence in~\eqref{EqMUscbtNormEquiv} as well as the $\tface$-normal operator estimate~\eqref{EqSptfAdm}, which is applicable when $\eta>0$ is small enough (cf.\ Remark~\ref{RmkMUbIndGap})
  \begin{align*}
    \|\chi_\tface u\|_{H_{\scbtop,\varsigma}^{\sfs,(\sfr+\alpha_+,\alpha,-\eta)}(X)} &\sim |\sigma|^{\frac32-\alpha} \|\chi_\tface u\|_{H_{\scop,\bop}^{\sfs,(\sfr+\alpha_+,-\alpha-\eta)}(\tface)} \\
      &\leq C|\sigma|^{\frac32-\alpha}\|N_\tface(P_0,e^{i\theta})(\chi_\tface u)\|_{H_{\scop,\bop}^{\sfs-2,(\sfr+\alpha_++1,-\alpha-\eta-2)}(\tface)} \\
      &\lesssim \|\varsigma^2 N_\tface(P_0,e^{i\theta})(\chi_\tface u)\|_{H_{\scbtop,\varsigma}^{\sfs-2,(\sfr+\alpha_++1,\alpha+2,-\eta}(X)}.
  \end{align*}
  We commute $\chi_\tface$ through $\varsigma^2 N_\tface(P_0,e^{i\theta})$ and replace $\varsigma^2 N_\tface(P_0,e^{i\theta})$ by $\wh{P_0}(\sigma)$; note that
  \begin{equation}
  \label{EqSpLoNtfDiff}
    [\varsigma^2 N_\tface(P_0,e^{i\theta}),\chi_\tface] \in \rho_\scface\rho_\tface^\infty\Diff_\scbtop^1,\quad
    \wh{P_0}(\sigma)-\varsigma^2 N_\tface(P_0,e^{i\theta}) \in \rho_\scface\rho_\tface^3\Diff_\scbtop^2,
  \end{equation}
  with the factor of $\rho_\scface$ in the second membership arising by comparison of~\eqref{EqSpTsscSpecFam} and \eqref{EqSptfOp0} (with $\hat\rho=\frac{\rho}{\varsigma}$ and $\hat\sigma=e^{i\theta}$ in present notation). Thus,
  \begin{align*}
    &\|\varsigma^2 N_\tface(P_0,e^{i\theta})(\chi_\tface u)\|_{H_{\scbtop,\varsigma}^{\sfs-2,(\sfr+\alpha_++1,\alpha+2,-\eta)}} \\
    &\qquad\leq C\Bigl(\|\wh{P_0}(\sigma)u\|_{H_{\scbtop,\varsigma}^{\sfs-2,(\sfr+\alpha_++1,\alpha+2,-\eta)}} + \|u\|_{H_{\scbtop,\varsigma}^{\sfs,(\sfr+\alpha_+,\alpha-1,-\eta)}}\Bigr).
  \end{align*}
  Since $(1-\chi_\tface)u$ is supported away from $\tface$, we have now improved~\eqref{EqSpLo1} to
  \[
    \|u\|_{H_{(\scbtop,\varsigma);\bop}^{\sfs,(\sfr+\alpha_+,\alpha,0)}} \leq C\Bigl(\|\wh{P_0}(\sigma)u\|_{H_{(\scbtop,\varsigma);\bop}^{\sfs-1,(\sfr+\alpha_++1,\alpha+2,0)}} + \|u\|_{H_{(\scbtop,\varsigma);\bop}^{\sfs,(\sfr+\alpha_+,\alpha-1,-\eta)}}\Bigr).
  \]
  But\footnote{We are being rather economical here: we only use estimates for the two models $\wh{P_0}(0)$ and $N_\tface(P_0,e^{i\theta})$, but no symbolic estimates. The price to pay is that for fixed $\varsigma$, this estimate does not even imply the semi-Fredholm property of $\wh{P_0}(\sigma)$ since the regularity and scattering decay orders of the final norm on the right are not improved relative to those of the left-hand side; but this is of no concern here, as we can use the smallness of $\sigma$ to conclude.} the norm of the error term is $\leq C\varsigma^\eta$ times the left-hand side and can thus be absorbed when $\varsigma=|\sigma|$ is sufficiently small; this gives~\eqref{EqSpLoEst}.

  For general $k\in\N$, the same arguments apply, now using the radial point and $\tface$-normal operator estimates featuring $k$ degrees of b-regularity (and using the zero energy bound with b-regularity order $\sfs+k$). The upper bound $\delta_k>0$ required for $|\sigma|$ so that one obtains~\eqref{EqSpLoEst} via the above argument may, a priori, depend on $k$, and indeed without loss of generality $\delta_k$ is non-increasing with $k$; but for any fixed value of $k\in\N_0$, one can upgrade~\eqref{EqSpLoEst} for $|\sigma|\in[\delta_k,\delta_0]$ from $k=0$ to the desired value of $k$ using Theorem~\ref{ThmSpB}\eqref{ItSpBb}. This establishes~\eqref{EqSpLoEst} for a $k$-independent neighborhood (of $0$) of frequencies $\sigma$.

  Finally, the proof~\eqref{EqSpLoInvReg} is analogous to that of~\eqref{EqSpBInvReg}: we now use that
  \[
    \sigma\pa_\sigma\wh{P_0}(\sigma) \in \sigma\rho\Diffb^1+\sigma^2\rho^2\CI
  \]
  is uniformly bounded as a map $H_{(\scbtop,|\sigma|);\bop}^{(\sfs+\ell;q),(\sfr+\alpha_+,\alpha,0)}\to H_{(\scbtop,|\sigma|);\bop}^{(\sfs+\ell;q-1),(\sfr+\alpha_++1,\alpha+2,0)}$. Indeed, the gain of one order at $\scface$ arises from the factors of $\rho$, while the gain of two orders at $\tface$ arises from the fact that $\sigma\rho$ vanishes quadratically there.\footnote{There is also a gain of one order at $\zface$, which we do not record here as it does not (typically) persist when applying $\wh{P_0}(\sigma)^{-1}$.}
\end{proof}

\begin{rmk}[Augmentations]
\label{RmkSpLoAug}
  In some applications, e.g., the Maxwell equations and the linearized gauge-fixed Einstein equations discussed in~\S\ref{SA2} and \cite{HintzKerrStab}, respectively, the zero energy operator $\wh{P_0}(0)$ is not invertible. We will then instead study appropriate ``augmentations'' of $\wh{P_0}(\sigma)$ whose zero energy operators \emph{are} invertible. The simplest type of augmentation is
  \begin{equation}
  \label{EqSpLoAug}
    \wt{P_0}(\sigma) := \begin{pmatrix} \wh{P_0}(\sigma) & E_{0 1} \\ E_{1 0} & 0 \end{pmatrix},
  \end{equation}
  where $E_{0 1},E_{1 0}$ are finite rank operators, with $E_{1 0}$, resp.\ $E_{0 1}$ being a vector of elements of $\CIc(X^\circ)$, resp.\ of $L^2$-pairings with elements of $\CIc(X^\circ)$ as in~\eqref{EqSpBPfwtP}. If $\wt{P_0}(0)$ is invertible (on the direct sums of the function spaces in~\eqref{EqSp0Map} with finite-dimensional $\C$-vector spaces), then the above arguments go through with purely notational changes to yield uniform estimates for $\wt{P_0}(\sigma)^{-1}$. Here, the $\tface$-normal operator of $\wt{P_0}(\sigma)$ is taken to be $N_\tface(\wh{P_0}(\sigma),e^{i\theta})$ simply, corresponding to the vanishing of $E_{0 1}$ and $E_{1 0}$ near $\tface$. It suffices, in fact, to require $E_{0 1}$ and $E_{1 0}$ to vanish more than quadratically (and thus faster than $\wh{P_0}(\sigma)$ as a sc-b-transition operator) at $\tface$.
\end{rmk}

\begin{cor}[Index $0$ property of the zero energy operator]
\label{CorSpLoInd0}
  Suppose $P_0$ is $\tface$-admissible with weight $\alpha+\frac32$. Then the Fredholm index of $\wh{P_0}(0)$ as a map~\eqref{EqSp0Map} equals $0$.
\end{cor}
\begin{proof}
  Let $d_0$ and $d_1$ be the dimension of the kernel and cokernel of~\eqref{EqSp0Map}, respectively, and fix $E_{0 1}\colon\C^{d_1}\to\CIc(X^\circ)$, $E_{1 0}\colon\cX_0\to\C^{d_0}$, such that $\ran E_{0 1}$ is complementary to the range of~\eqref{EqSp0Map}, and $E_{1 0}|_{\ker\wh{P_0}(0)}$ is invertible. Then Theorem~\ref{ThmSpLo} applies to the augmented operator~\eqref{EqSpLoAug} and implies the invertibility of $\wt{P_0}(\sigma)$ (as a map~\eqref{EqSpBMap}) for $\sigma=i\eps$ when $\eps>0$ is sufficiently small. But the index of~\eqref{EqSpBMap} is $0$, and thus $d_0-d_1=0$, as claimed.
\end{proof}

\section{Upgrading \texorpdfstring{$\aleph$-admissibility}{aleph-admissibility} to bounds on weighted e3b-Sobolev spaces}
\label{SN}

We now return to the task of improving the estimates in Theorem~\ref{ThmEReg}. The most systematic approach towards upgrading the estimate~\eqref{EqEReg}---i.e., weakening the decay orders $\alpha_+$ and $\alpha_\cK$ of the $H_{\etbop;\bop}^{(\sfs_0;k),(2\alpha_\sscri,\alpha_+,\alpha_\cK)}$-norm of $u$ on the right-hand side---would proceed via the inversion of the two normal operators of $P$ at $\cK^+$ and $\iota^+$.\footnote{One can also weaken the decay order $2\alpha_\sscri$ at $\scri^+$, which is, however, not necessary when one localizes to sufficiently large $t_*$; cf.\ Remark~\ref{RmkFNoScriNorm}.} The normal operator at the former, resp.\ the latter being the time-translation-invariant operator
\begin{equation}
\label{EqNP0}
\begin{split}
  &P_0 = -2\pa_{t_*}\rho(\rho\pa_\rho-1-S) + \rho^2 P_{(0)}(\rho,\omega,\rho\pa_\rho,\pa_\omega) + Q\pa_{t_*} - g^{0 0}\pa_{t_*}^2, \\
  &\hspace{5em} P_{(0)} \equiv \rho^2\bigl((\rho D_\rho)^2 + \slDelta\bigr) \bmod \rho^2\Diffb^1+\rho^3\Diffb^2,
\end{split}
\end{equation}
in the notation of Definition~\ref{DefSSAdm} and~\eqref{EqSSSpec0Op} (cf.\ \S\ref{SssCPX}), resp.\ a spacetime-dilation-invariant operator (cf.\ \S\ref{SssCPip}), this would involve the Fourier transform in $t_*$ and the spectral family of $P_0$, resp.\ the Mellin transform in $\frac{1}{t_*}$ and the Mellin-transformed normal operator family of $P_0$ (which is the same as that of $P$) at $\iota^+$ (cf.\ \eqref{EqCPXA}--\eqref{EqCPXA0gl} and~\S\ref{SssMUK}, resp.\ \eqref{EqCPipA}--\eqref{EqCPipA0scri2} and~\S\ref{SssMUip}). This was the approach pursued in \cite[\S\S{5.3}--{5.4}]{HintzNonstat} for wave-type operators without zero energy bound states (and without trapping), and for elliptic 3b-operators in \cite[\S\S{7.1}--{7.2}]{Hintz3b}.

We will argue differently here and use more strongly that $P_0$ is a good model for $P$ near \emph{all of} $\iota^+\cup\cK^+$. One consequence of Definition~\ref{DefSDWAdm} is that the difference between an admissible wave-type operator $P$ and its stationary model $P_0$ decays at $\iota^+\cup\cK^+$ (as an e3b-differential operator with weight $\rho_+^2$ at $\iota^+$). (Correspondingly, \emph{the current section entirely concerns the stationary model $P_0$}.) While this suggests using the Fourier transform in $t_*$, it is only for 3b-Sobolev spaces (on $\tilde M_0$) that the Fourier transform provides a Plancherel-type isomorphism (see Lemma~\ref{LemmaMUetbFT}), whereas this is not the case for e3b-Sobolev spaces.\footnote{Only function spaces defined via testing with $t_*$-translation-invariant vector fields or ps.d.o.s admit a Plancherel-type characterization on the Fourier transform side, as only such vector fields or ps.d.o.s are diagonal when acting on functions with time dependence $e^{-i\sigma t_*}$, $\sigma\in\C$. Thus, if one uses 3b-vector fields, say $r\pa_{t_*},r\pa_r,\pa_\omega$, near $\cK^+\cap\{r\geq r_0\}$, $r_0>\bhm$, to test for regularity, Plancherel-type passage to the Fourier transform forces one to use the same vector fields everywhere in $\{r\geq r_0\}$ and thus, in particular, near $\iota^+\cup\scri^+$. Thus one cannot transition from 3b- to edge-b-vector fields near $\scri^+$ and retain compatibility with the Fourier transform.} This is why in \cite{Hintz3b,HintzNonstat}, the Mellin transform is used near $\iota^+$ to analyze $P_0$, with a Plancherel-type result being available (see Lemma~\ref{LemmaMUetbM}). This is also the reason why the $\aleph$-admissibility for $P_0$ in Definition~\ref{DefSSAlephAdm} is phrased in terms of 3b-Sobolev spaces (and, in practice, proved via the Fourier transform and resolvent estimates).

A novel issue arises in the present paper, however, due to possibility of zero energy bound states, which in particular leads to the loss of decay at $\cK^+$ upon inverting $P_0$ (cf.\ the parameter $\aleph$ in Definition~\ref{DefSSAlephAdm}): the weights $-\aleph$ and $\alpha_\sface$ at $\cK^+$ and $\iota^+$ for which the stationary estimates~\eqref{EqSSAlephAdmSol} on 3b-Sobolev spaces hold for solutions $u$ of $P_0 u=f$ equation are not (directly) compatible with the weights for which estimates for $P_0$ on e3b-Sobolev spaces near $\iota^+$ hold.\footnote{Roughly speaking, the relative order $\alpha_\sface-\alpha_\cK+\frac32$, with $\alpha_\cK=0$ in Definition~\ref{DefSSAlephAdm}, needs to lie in the indicial gap $(\beta^-,\beta^+)$ in order for 3b-estimates to hold (in our applications this is due to the requirements of the zero energy problem in Theorem~\ref{ThmSp0}). But then, typically, the relative order for the output $u$, which is roughly $\alpha_\sface-(\alpha_\cK-\aleph)$, falls outside of that interval, rendering $\iota^+$-normal operator estimates and thus e3b-estimates impossible; cf.\ the weight $\gamma_\cface$ in Proposition~\ref{PropNip}.} We resolve this issue as follows. Since e3b- and 3b-Sobolev spaces (the latter lifted from $M_0$ to $M$) differ only near $\scri^+$, we shall invert $P_0$ on e3b-spaces, i.e., solve $P_0 u=f$, in several steps so as to isolate the behavior near $\scri^+$:
\begin{enumerate}
\item we first solve $P_0 u=f$ near $\iota^+\cap\scri^+$ on edge-b-spaces by inverting not $P_0$ but a modification $P_+$ of $P_0$ which agrees with $P_0$ at $\iota^+$ and near $\scri^+$ but which cuts away the Kerr behavior at $\cK^+$ (see~\S\ref{SsNMod} for details);
\item we then solve away the remaining part of $f$ (now localized near $\cK^+$ and thus away from $\scri^+$) using Definition~\ref{DefSSAlephAdm}, which thus produces a solution in a 3b-Sobolev space on $M_0$;
\item we show that this 3b-solution in fact lies in an appropriate e3b-Sobolev space on $M$ by relating it near $\scri^+$ to an equation involving the modified operator $P_+$ and exploiting the uniqueness and regularity properties of solutions of $P_+$.
\end{enumerate}

The main result of this section is:

\begin{thm}[Forward solutions for $P_0$: e3b-spaces]
\label{ThmNFw}
  Let $P_0$ be a stationary wave-type operator (Definition~\usref{DefSSAdm}). Let $\alpha_\sscri\in\R$ with
  \begin{equation}
  \label{EqNFwThr}
    \alpha_\sscri<-\frac12+\ubar S,\quad \alpha_+<\alpha_\sscri-\frac12,
  \end{equation}
  where we recall $\ubar S=\inf_{\pa X}\Re\spec S$ from~\eqref{EqSSAdmubarS}. Assume that $P_0$ is $\aleph$-admissible with $\sface$-weight $\alpha_+$ and $\sface$-loss $\delta$ (Definition~\usref{DefSSAlephAdm}),\footnote{or the relaxed version in Remark~\usref{RmkSSAlephAdmRelax}, discussed in Remark~\usref{RmkNFwRelax} below} and let $\sfs\in\CI({}^\etbop S^*M)$ be an order function which is admissible with weights $\alpha_+,-\aleph$ and margin $0$ (Definition~\usref{DefDyO}), and also with weights $\alpha_+,0$ and margin $1$;\footnote{Both conditions hold if $\sfs$ is admissible with weights $\alpha_+,-\aleph$ and margin $\max(1,\aleph)$.} we require moreover that
  \begin{equation}
  \label{EqNFwThr2}
    \alpha_+ < \alpha_\sscri - \frac12.
  \end{equation}
  Then for all $k\in\N_0$, there exists a constant $C_k$ such that the unique forward solution of $P_0 u=f$ on $\Omega=\ol{\{t_*\geq 1\}}\subset M$ satisfies
  \begin{equation}
  \label{EqNFW}
    \|u\|_{H_{\etbop;\bop}^{(\sfs;k),(2\alpha_\sscri,\alpha_+-\delta,-\aleph)}(\Omega;\cE)^{\bullet,-}} \leq C_k\|P_0 u\|_{H_{\etbop;\bop}^{(\sfs;k),(2\alpha_\sscri+2,\alpha_++2,0)}(\Omega;\cE)^{\bullet,-}}.
  \end{equation}
\end{thm}

\begin{rmk}[Threshold condition on $\alpha_\sscri$]
\label{RmkNFwThr}
  The conditions~\eqref{EqNFwThr} and \eqref{EqNFwThr2} are implied by~\eqref{EqRWeights} since, in the notation of Definition~\ref{DefSDWubarp1}, we have $p_1(p)=S(p)$ for $p\in\scri^+\cap\iota^+\cong\pa X$ and thus $\ubar p_1\leq\ubar S$. For $P=P_0$, we have $p_1=S|_{\pa X}$ and thus $\ubar p_1=\ubar S$.
\end{rmk}

Before introducing and analyzing the modified operator $P_+$ in~\S\ref{SsNMod}, we prove preliminary estimates for the dilation-invariant model of $P_0$ at $\iota^+$ (which is related to the transition face normal operators $N_\tface(P_0,e^{i\theta})$ via the Fourier transform) in~\S\ref{SsNip}. The proof of Theorem~\ref{ThmNFw} is then given in~\S\ref{SsNPf}. \emph{We drop the bundles $\cE_X$ and $\cE$ from the notation.}

\subsection{Inversion of the \texorpdfstring{$\iota^+$-normal operator}{normal operator at punctured timelike infinity}}
\label{SsNip}

We shall prove a solvability and uniqueness result for the wave-type operator $N_{\iota^+}(P_0)$, defined below, on the space $[0,\infty)_{t_*^{-1}}\times\iota^+$; see Proposition~\ref{PropNipFwd} below, which is strongly used in~\S\ref{SsNMod}. The parts in this section concerning the Mellin-transformed normal operator family $N_{\iota^+}(P_0,\lambda)$ of $P_0$ at $\iota^+$ (see Definition~\ref{DefNipOp} below) are minor variations on \cite[\S{5.3}]{HintzNonstat} for bounded $\Im\lambda$. (See also \cite[\S{8}]{HintzVasyScrieb} for localized estimates near $\iota^+\cap\scri^+$.) The novelties here are, first, that we prove estimates involving $k\in\N_0$ degrees of b-regularity, which requires us to revisit the arguments in the reference, and, second, that we also prove estimates for large $\Im\lambda$.

Recalling Lemma~\ref{LemmaCPip} (see also \S\ref{SssMUip}), we work in the local coordinates
\begin{equation}
\label{EqNipCoords}
  \tau=\frac{1}{t_*},\quad R=\frac{r}{t_*}=\frac{1}{\rho t_*}
\end{equation}
near $(\iota^+)^\circ$, so $R=v^{-1}$ where $v=\frac{t_*}{r}=\rho t_*$. (Thus $R$ and $v^{\frac12}=R^{-\frac12}$ are local defining functions of $\iota^+\cap\cK^+$ and $\iota^+\cap\scri^+$, respectively, near these boundary hypersurfaces of $\iota^+$.)

\begin{definition}[$\iota^+$-normal operator]
\label{DefNipOp}
  The \emph{$\iota^+$-normal operator} of $P_0$ is the differential operator
  \begin{equation}
  \label{EqNipOp}
  \begin{split}
    N_{\iota^+}(P_0) &= t_*^2 N_{\iota^+}^0(P_0), \\
    N_{\iota^+}^0(P_0) &:= -2 R^{-1}(\tau\pa_\tau+R\pa_R)(R\pa_R+1+S|_{\pa X}) + R^{-2}P_{(0)}(0,\omega,-R\pa_R,\pa_\omega),
  \end{split}
  \end{equation}
  on $[0,\infty)_\tau\times\iota^+$. The \emph{$\iota^+$-normal operator family} of $P_0$ is
  \begin{equation}
  \label{EqNipOpMT}
  \begin{split}
    N_{\iota^+}^0(P_0,\sigma) &:= -2 R^{-1}(i\sigma+R\pa_R)(R\pa_R+1+S|_{\pa X}) + R^{-2}P_{(0)}(0,\omega,-R\pa_R,\pa_\omega) \\
      &= x_\sscri^2\Bigl[-\frac12(x_\sscri\pa_{x_\sscri}-2 i\sigma)(x_\sscri\pa_{x_\sscri}-2(1+S|_{\pa X})) + x_\sscri^2 P_{(0)}(0,\omega,\tfrac12 x_\sscri\pa_{x_\sscri},\pa_\omega) \Bigr],
  \end{split}
  \end{equation}
  where we set $x_\sscri=R^{-\frac12}$ in the second line.
\end{definition}

\begin{rmk}[Not the b-normal operator]
\label{RmkNipOpNotb}
  Since $t_*^{-1}=\tau$ is a \emph{joint} defining function for $\iota^+\cap\cK^+$, this is \emph{not} the same as the normal operator (family) of $P_0$ when regarded as a weighted b-operator on $M$: the latter would involve $\rho_+^{-2}P_0$ and the Mellin transform in $\rho_+$. The ``total'' (i.e., in $\tau$) Mellin transform is, however, what interacts cleanly with 3b-operators and -spaces near $\iota^+\cap\cK^+$. See also \cite[\S{3.3}]{Hintz3b}.
\end{rmk}

Note that $N_{\iota^+}^0(P_0)$ is $\tau$-dilation-invariant. The operator $N_{\iota^+}(P_0)$ is homogeneous of degree $-2$ with respect to dilations in $t_*$; and it is a wave-type operator relative to the metric $g_{\iota^+}$ whose dual is given by
\begin{equation}
\label{EqNipDualMet}
\begin{split}
  g_{\iota^+}^{-1} &:= 2 R^{-1}(\tau\pa_\tau+R\pa_R)\otimes_s R\pa_R+R^{-2}((R\pa_R)^2+\slg^{-1}) \\
    &= x_\sscri^2\Bigl(\frac12(x_\sscri\pa_{x_\sscri} - 2\tau\pa_\tau)\otimes_s x_\sscri\pa_{x_\sscri} + x_\sscri^2\bigl(\tfrac14(x_\sscri\pa_{x_\sscri})^2 + \slg^{-1}\bigr) \Bigr).
\end{split}
\end{equation}
(This is of course the same as the rescaling $\tau^{-2}\ubar g^{-1}=t_*^2(-2\pa_{t_*}\otimes_s\pa_r+\pa_r^2+r^{-2}\slg^{-1})$ of the Minkowski dual metric, expressed in the coordinates $\tau=\frac{1}{t_*}$, $R=\frac{r}{t_*}$, and $\tau=\frac{1}{t_*}$, $x_\sscri=\sqrt{\frac{t_*}{r}}$, respectively.) Moreover, $N_{\iota^+}(P_0)$ equals $P_0\in x_\sscri^2\rho_+^2\Diff_\etbop^2$ (see Lemma~\ref{LemmaSSOpMem}) to leading order at $\iota^+$ in the following sense:

\begin{lemma}[$\iota^+$-normal operator]
\label{LemmaNipOpDiff}
  Let $\chi\in\CI(M)$ be equal to $1$ near $\iota^+$ and supported in a small collar neighborhood of $\iota^+$. Then $\chi N_{\iota^+}(P_0)\in x_\sscri^2\rho_+^2\Diff_\etbop^2(M)$ and
  \begin{equation}
  \label{EqNipOpDiff}
    \chi\bigl(P_0 - N_{\iota^+}(P_0)\bigr) \in x_\sscri^2\rho_+^3\Diff_\etbop^2(M).
  \end{equation}
\end{lemma}

This can be regarded as an instance of~\eqref{EqCPipA}--\eqref{EqCPipA0scri2}.

\begin{proof}[Proof of Lemma~\usref{LemmaNipOpDiff}]
  The membership of $\chi N_{\iota^+}(P_0)$ is checked as in Lemma~\ref{LemmaSSOpMem}. The vector fields $t_*\pa_{t_*}$ and $\rho\pa_\rho$ in the coordinates $(t_*,\rho)$ become $-\tau\pa_\tau-R\pa_R$ and $-R\pa_R$, respectively, in the coordinates~\eqref{EqNipCoords}. The terms in~\eqref{EqNipOp} then arise from the first term in~\eqref{EqNP0} and from $\rho^2 P_{(0)}(0,\omega,\rho\pa_\rho,\pa_\omega)$. We need to show that the remaining terms lie in the space~\eqref{EqNipOpDiff}; near $\scri^+$, we use $x_\sscri=\sqrt{v}$ as a local defining function of $\scri^+$. The difference of $S$ and $S|_{\pa X}$ lies in $\rho\CI$; and $t_*^2\cdot \pa_{t_*}\rho^2=-\tau^{-1}(\tau\pa_\tau+R\pa_R)\tau^2 R^{-2}$ vanishes at $\tau=0$ and to order $v^2=x_\sscri^4$ at $\scri^+$. The difference of $\rho^2 P_{(0)}(\rho,\omega,\rho\pa_\rho,\pa_\omega)$ and its normal operator at $\rho=0$ lies in $\rho^3\Diff_\bop^2(X)$, so $t_*^2$ times this is of the schematic form $R^{-3}\tau(R\pa_R,\pa_\omega)^{\leq 2}$, so a fortiori $x_\sscri^4\tau(x_\sscri\pa_{x_\sscri},x_\sscri\pa_\omega)^{\leq 2}$, which lies in~\eqref{EqNipOpDiff} as well. Similarly, $t_*^2 Q\pa_{t_*}\in t_*^2\rho^3\Diffb^1(X)\pa_{t_*}$ is of the schematic form $R^3\tau^2(\tau\pa_\tau+R\pa_R)^{\leq 1}(R\pa_R,\pa_\omega)^{\leq 1}$ (so vanishes at $\iota^+$ and to fifth order at $\scri^+$), and finally $t_*^2 g^{0 0}\pa_{t_*}^2$ is of the form $(\tau\pa_\tau+R\pa_R)^{\leq 2}R^{-2}\tau^2$ since $g^{0 0}\in\rho^2\CI(X)$; this vanishes at $\iota^+$ and to fourth order at $\scri^+$.
\end{proof}

By the general theory in~\S\ref{SssMUip} or by inspection of~\eqref{EqNipOp}, we have
\[
  N_{\iota^+}^0(P_0,\sigma)\in \Diff_{0,\bop}^{2,(-2,2)}(\iota^+) = x_\sscri^2\rho_\cK^{-2}\Diff_{0,\bop}^2(\iota^+), \quad \sigma\in\C,
\]
where $x_\sscri$ and $\rho_\cK$ are defining functions of $\iota^+\cap\scri^+$ and $\iota^+\cap\cK^+$, respectively. One can, e.g., take
\[
  \rho_\cK = \frac{R}{R+1},\quad
  x_\sscri = \sqrt{\frac{v}{v+1}}.
\]
Indeed, if $\chi\in\CIc([0,3))$ equals $1$ on $[0,2]$, then a spanning set of the space $\cV_{0,\bop}(\iota^+)$ of 0-b-vector fields is $\chi(R)R\pa_R$, $\chi(R)\pa_\omega$, $\chi(R^{-1})R\pa_R$, $\chi(R^{-1})R^{-\frac12}\pa_\omega$. In the high-energy regime, we introduce the semiclassical rescaling
\begin{equation}
\label{EqNipSclResc}
  P_{\iota^+,h,z} := h^2 N_{\iota^+}^0(P_0,h^{-1}z);
\end{equation}
then for bounded $z$,
\[
  (P_{\iota^+,h,z})_{h\in(0,1)} \in \Diff_{0\semi,\chop}^{2,(-2,0,0,2)}(\iota^+) = x_\sscri^2\rho_\cface^{-2}\Diff_{0\semi,\chop}^2(\iota^+),
\]
where, in the notation of~\eqref{EqMUipSingle}, $\rho_\cface\in\CI(\iota^+_{0\semi,\chop})$ is a defining function of $\cface\subset\iota^+_{0\semi,\chop}$ such as $\rho_\cface=\frac{\rho_\cK}{h+\rho_\cK}$ (or, away from $\scri^+$: $\frac{R}{h+R}$). (The orders ``$0$'' refer to the orders at $\sface$ and $\tface$, respectively.) This membership can either be checked directly from the expressions~\eqref{EqNipOpch}--\eqref{EqNipOp0h} below, or from~\eqref{EqMUipSpecFam0h}--\eqref{EqMUipSpecFamch} which gives the membership of $(P_{\iota^+,h,z})_{h\in(0,1)}$ in $x_\sscri^2 h^2\cdot h^{-2}\Diff_{0\semi}^2=x_\sscri^2\Diff_{0\semi}^2$ near $\iota^+\cap\scri^+$ and in $h^2\cdot R^{-2}(\frac{h}{h+R})^{-2}\Diff_\chop^2=(\frac{R}{h+R})^{-2}\Diff_\chop^2$.

For the analysis of this operator family, we use the corresponding classes of ps.d.o.s and Sobolev spaces. Since $N_{\iota^+}^0(P_0,\sigma)$ has smooth coefficients, we follow Notation~\ref{NotSpPsdo} and do not make the coefficient class of ps.d.o.s on $\iota^+$ explicit in the notation; we shall always work with conormal coefficients (i.e., $\CI_\bop(\iota^+)$), and write $\Psi_{0,\bop}$, $\Psi_{0\semi,\chop}$ as in~\S\ref{SssMUip} but omitting ``$\CI_\bop$''; we moreover write
\[
  \Psi_0,\ \text{resp.}\ \Psi_\bop,\quad\text{and}\quad\Psi_{0\semi},\ \text{resp.}\ \Psi_\chop
\]
for operators of class $\Psi_{0,\bop}(\iota^+)$ and $\Psi_{0\semi,\chop}$ which are microlocalized near $\iota^+\cap\scri^+$, resp.\ $\iota^+\cap\cK^+$, unless specified otherwise. We similarly use the notation\footnote{Since the weights in these spaces are related to the weights at the boundary hypersurfaces of $M$ by shifts, as in Lemma~\ref{LemmaMUetbM}, we use the notation $\gamma_\bullet$ instead of $\alpha_\bullet$ in this section.}
\[
  H_{0,\bop}^{s,(2\gamma_\sscri,\gamma)}(\iota^+),\quad
  H_{0\semi,\chop}^{s,(2\gamma_\sscri,\gamma_\cface,l,b)} = x_\sscri^{2\gamma_\sscri}\rho_\cface^{\gamma_\cface}\rho_\tface^l\rho_\sface^b H_{0\semi,\chop}^s
\]
for the associated Sobolev spaces as in~\eqref{EqMUip0bSob} and \eqref{EqMUip0hchSob}, with the density underlying the $L^2$-space being any smooth positive b-density on $\iota^+$ such as $|\frac{\dd R}{R}\,\dd\slg|$ (as also used in Lemma~\ref{LemmaMUetbM}). Spaces with additional $k\in\N_0$ orders of b-regularity are denoted $H_{0,\bop;\bop}^{(s;k)}$ and $H_{0\semi,\chop;\bop}^{(s;k)}$, with the norms denoted $\|\cdot\|_{H_{0,\bop;\bop}^{(s;k)}}$ and $\|\cdot\|_{H_{(0\semi,\chop,h);\bop}^{(s;k)}}$ (the latter norm depending parametrically on $h\in(0,1)$), similarly for weighted versions; and analogously to~\eqref{EqSpHiNormb}, we write
\begin{equation}
\label{EqNipNormbplus}
  \|u\|_{H_{(0\semi,\chop,h);\bop^+}^{(s;k),(2\gamma_\sscri,\gamma_\cface,l,b)}} := \sum_{j=0}^k \|u\|_{H_{(0\semi,\chop,h);\bop}^{(\sfs;k-j),(2\gamma_\sscri,\gamma_\cface,l,b+j)}}.
\end{equation}
Thus, powers of $\rho_\sface^{-1}$, which in $R<2$ are powers of $R h^{-1}$, are regarded as ``extended'' b-derivatives. (In view of the factor of $R$, this sharper than using powers of $h^{-1}$.)

Recalling Definition~\ref{DefDyO}, we use the orders induced by an admissible order function
\begin{equation}
\label{EqNipOrder}
  \sfs\in\CI({}^\etbop S^*M)\ \text{with weights}\ \alpha_+,0\ \text{and margin}\ 0,
\end{equation}
as in Lemma~\ref{LemmaMUetbM}; we make this concrete in~\S\ref{SssNipPhase} below. We denote these orders by
\begin{equation}
\label{EqNipOrders}
\begin{alignedat}{2}
  \sfs_{0,\bop} &\in \CI({}^{0,\bop}S^*\iota^+)&\quad&\text{(differential order)}, \\
  \sfs_{0\semi,\chop} &\in \CI({}^{0\semi,\chop}S^*\iota^+)&\quad&\text{(semiclassical differential order)}, \\
  \sfb_\pm&\in\CI(\ol{{}^{0\semi,\chop}T^*_\sface}\iota^+)&\quad&\text{(semiclassical order when $\pm\sigma\to\infty$)},
\end{alignedat}
\end{equation}
where $\sface$ was defined after~\eqref{EqMUipSingle}.

\begin{prop}[Estimates for the $\iota^+$-normal operator family]
\label{PropNip}
  Recall the upper bound $\beta^+\in\R$ of the indicial gap for $\wh{P_0}(0)$ from Definition~\usref{DefSStfAdm} in which also $\alpha_++\frac32$ lies (cf.\ Definition~\usref{DefSSAlephAdm}\eqref{ItSSAlephAdm2}). Let $\gamma_\sscri<1+\ubar S$ and $\gamma_\cface=-(\alpha_++\frac32)$.
  \begin{enumerate}
  \item\label{ItNipB}{\rm (Bounded frequencies.)} For all $s$ (real or variable, in particular for $s=\sfs_{0,\bop}$) and $k\in\N_0$, the operator
    \begin{equation}
    \label{EqNip}
      N_{\iota^+}^0(P_0,\sigma) \colon H_{0,\bop;\bop}^{(s;k),(2\gamma_\sscri,\gamma_\cface)}(\iota^+) \to H_{0,\bop;\bop}^{(s-2;k),(2\gamma_\sscri+2,\gamma_\cface-2)}(\iota^+)
    \end{equation}
    is Fredholm of index $0$ for all $\sigma\in\C$ satisfying
    \[
      \Im\sigma>-\gamma_\sscri,
    \]
    and invertible if, in addition, $\Im\sigma>-1-\beta^+$.
  \item\label{ItNipHi}{\rm (High-energy estimates.)} Let $\gamma_+<\gamma_\sscri$. Then there exists $h_0>0$ such that for all $z\in\C$ with $|z|\in[\frac12,2]$ and $\Im z\geq -h\gamma_+$, the operator~\eqref{EqNip} is invertible for $\sigma=h^{-1}z$, $0<h<h_0$, with the inverse satisfying the estimate
    \begin{equation}
    \label{EqNipb}
      \|u\|_{H_{(0\semi,\chop,h);\bop^+}^{(\sfs;k),(2\gamma_\sscri,\gamma_\cface,\gamma_\cface,\sfb)}(\iota^+)} \leq C_k\|P_{\iota^+,h,z}u\|_{H_{(0\semi,\chop,h);\bop^+}^{(\sfs-2;k),(2\gamma_\sscri+2,\gamma_\cface-2,\gamma_\cface,\sfb+1)}(\iota^+)}
    \end{equation}
    for all $k\in\N_0$. Here $\sfs$ is arbitrary (e.g., equal to $\sfs_{0\semi,\chop}$), and $\sfb$ is given by $\sfb_\pm$ when $|z\mp 1|\leq\frac12$ and an arbitrary constant otherwise.
  \end{enumerate}
\end{prop}

We first this for large frequencies using semiclassical techniques in~\S\ref{SssNipHi}, and then for bounded frequencies using a parametrix construction in~\S\ref{SssNipB}. Our main interest in Proposition~\ref{PropNip} stems from its relevance for the study of $N_{\iota^+}(P_0)$ as a wave-type operator in its own right; this is discussed in~\S\ref{SssNipFw}.

Since decay as $|\Re\sigma|\to\infty$ corresponds to regularity under $\tau\pa_\tau$, we do not need to study the mapping properties of $\sigma$-derivatives of $N_{\iota^+}^0(P_0,\sigma)^{-1}$ here (by contrast to Theorem~\ref{ThmSpB}\eqref{ItSpBInvReg}, say).

\subsubsection{Phase space relationships}
\label{SssNipPhase}

We make the relationships between the spacetime phase space ${}^\etbop T^*_{\iota^+}M$ and the phase spaces for the analysis of $N_{\iota^+}^0(P_0,\sigma)$ explicit, and thereby also explain the orders induced by $\sfs|_{{}^\etbop S^*_{\iota^+}M}$ and used in Proposition~\ref{PropNip}. In $R<2$, we write (e)3b-covectors as $\sigma_\tbop\frac{\dd t_*}{r}+\xi_\tbop\frac{\dd r}{r}+\eta_\tbop$, $\eta_\tbop\in T^*\Sph^2$, while in $v<2$, we write e(3)b-covectors as $\sigma_\ebop\frac{\dd\rho_+}{\rho_+}+\xi_\ebop\frac{\dd x_\sscri}{x_\sscri}+\frac{\eta_\ebop}{x_\sscri}$, $\eta_\ebop\in T^*\Sph^2$, where $x_\sscri=R^{-\frac12}$ (cf.\ \eqref{EqCHam3bCoord} and \eqref{EqCHamebCoord}). Thus, $\sfs$ is a homogeneous degree $0$ function in $(\sigma_\tbop,\xi_\etbop,\eta_\tbop)$ and $(\sigma_\ebop,\xi_\ebop,\eta_\ebop)$, respectively.

On the spectral (Mellin transform) side, we write b-covectors on $\iota^+$ in $R<2$ as $\xi_\bop\frac{\dd R}{R}+\eta_\bop$, $\eta_\bop\in\Sph^2$, and 0-covectors in $v<2$ as $\xi_0\frac{\dd x_\sscri}{x_\sscri}+\frac{\eta_0}{x_\sscri}$. Note that $\frac{\dd R}{R}=\frac{\dd r}{r}-R\frac{\dd t_*}{r}$. For $\sigma\in\R$, the phase space relationship~\eqref{EqMSFPhaseSpace} (with $t$ there being $\log t_*=-\log\rho_+$, $\rho_+:=t_*^{-1}$, in present notation) thus produces from $\sfs$ the function $\sfs\circ\Phi_\sigma$, where
\begin{alignat*}{3}
  \Phi_\sigma\colon\Bigl(R,\omega;\xi_\bop\frac{\dd R}{R}+\eta_\bop\Bigr)&\ni\Tb^*_{\{R<2\}}\iota^+ &&\mapsto \Bigl(R,\omega;-R(\xi_\bop+\sigma)\frac{\dd t_*}{r}+\xi_\bop\frac{\dd r}{r}+\eta_\bop\Bigr) &&\in {}^\etbop T^*_{\iota^+\cap\{R<2\}}M, \\
    \Bigl(x_\sscri,\omega;\xi_0\frac{\dd x_\sscri}{x_\sscri}+\frac{\eta_0}{x_\sscri}\Bigr) &\ni {}^0 T^*_{\{v<2\}}\iota^+ &&\mapsto \Bigl(x_\sscri,\omega;\sigma\frac{\dd\rho_+}{\rho_+}+\xi_0\frac{\dd x_\sscri}{x_\sscri}+\frac{\eta_0}{x_\sscri}\Bigr) &&\in {}^\etbop T^*_{\iota^+\cap\{v<2\}}M.
\end{alignat*}
The restriction of $\Phi_\sigma$ to fiber infinity is independent of $\sigma$, and thus $\sfs\circ\Phi_\sigma$ restricts to ${}^{0,\bop}S^*\iota^+$ as the $\sigma$-independent function $\sfs_{0,\bop}$ in~\eqref{EqNipOrders}. From the expressions~\eqref{EqTs3bCCoord} and \eqref{EqTsebCoord}, it moreover follows that the image of fiber infinity under $\Phi_\sigma$ is disjoint from the characteristic set of $P_0$, and therefore the operator $N_{\iota^+}^0(P_0,\sigma)\in x_\sscri^2\rho_\cK^{-2}\Diff_{0,\bop}^2(\iota^+)$ is elliptic. This of course also follows directly from~\eqref{EqNipOpMT}, whose principal symbol is $R^{-2}((\xi_\bop+R)^2+|\eta_\bop|^2)=x_\sscri^2(\frac12\xi_0^2+|\eta_0|^2)$.

For the high-energy regime, consider $\sigma=\pm h^{-1}$, $h\in(0,1)$. Working near $\iota^+\cap\cK^+$, write semiclassical cone covectors on $\iota^+$ as
\begin{equation}
\label{EqNipCovecchop}
  \xi_\chop\,\frac{h+R}{h}\frac{\dd R}{R}+\frac{h+R}{h}\eta_\chop,
\end{equation}
then $\sfs\circ\Phi_{\pm h^{-1}}$ is given by
\[
  \Bigl(R,\omega;\xi_\chop\,\frac{h+R}{h}\frac{\dd R}{R}+\frac{h+R}{h}\eta_\chop\Bigr) \mapsto \sfs\Bigl(R,\omega; -\Bigl(R\xi_\chop\pm\frac{R}{h+R}\Bigr)\frac{\dd t_*}{r}+\xi_\chop\frac{\dd r}{r}+\eta_\chop \Bigr),
\]
where we exploited the homogeneity of $\sfs$. The restriction to $\frac{h}{h+R}=0$, given by $\sfs\circ\Phi_{\chop,\pm}$ where
\[
  \Phi_{\chop,\pm}\Bigl(R,\omega;\xi_\chop\frac{h+R}{h}\frac{\dd R}{R}+\frac{h+R}{h}\eta_\chop\Bigr) = \Bigl(R,\omega; (\mp 1+R\xi_\chop)\frac{\dd t_*}{r}+\xi_\chop\frac{\dd r}{r}+\eta_\chop\Bigr),
\]
is the semiclassical order $\sfb_\pm$ in~\eqref{EqNipOrders} on $\{R<2\}$. Near $\iota^+\cap\scri^+$ on the other hand, we write semiclassical 0-covectors as $\xi_{0,\semi}\,\frac{\dd x_\sscri}{h x_\sscri}+\frac{\eta_{0,\semi}}{h x_\sscri}$; then $\sfs\circ\Phi_{\pm h^{-1}}$ is
\[
  \Bigl(x_\sscri,\omega;\xi_{0,\semi}\frac{\dd x_\sscri}{h x_\sscri}+\frac{\eta_{0,\semi}}{h x_\sscri}\Bigr) \mapsto \sfs\Bigl( x_\sscri,\omega; \pm\frac{\dd\rho_+}{\rho_+} + \xi_{0,\semi}\frac{\dd x_\sscri}{x_\sscri}+\frac{\eta_{0,\semi}}{x_\sscri} \Bigr),
\]
and thus, upon identifying semiclassical and standard 0-covectors via multiplication by $h$, the pullback of $\sfs$ along $\Phi_{\pm 1}$. This restricts to $h=0$ to give the semiclassical order $\sfb_\pm$ in~\eqref{EqNipOrders} on $\{v<2\}$. The restrictions to fiber infinity on the other hand are independent of the sign of $\sigma$, and are given by the function $\sfs_{0\semi,\chop}$ in~\eqref{EqNipOrders}; the latter function moreover restricts to ${}^{0\semi,\chop}S^*_\tface\iota^+={}^{\bop,\scop}S^*\tface$ (in the notation introduced around~\eqref{EqMUNtfSpace}) to a b-scattering-regularity order
\begin{equation}
\label{EqNipOrdertf}
  \sfs_{\bop,\scop} \in \CI({}^{\bop,\scop}S^*\tface).
\end{equation}

The radial sets $\pa\cR^\pm_{\pa\cK^+,{\rm in/out}}$ (see Definition~\ref{DefTs3bRad}, with the sets for the ``$-$'' being obtained by fiber-wise multiplication by $-1$) lie only in the image of $\Phi_{\chop,\pm}$ (but not of $\Phi_{\chop,\mp}$), and their preimages are the subsets
\begin{alignat*}{2}
  \cR_{\cface,\pm 1,{\rm in}} &:= \{ (R,\omega; \xi_\chop,\eta_\chop) \colon R=0,\ \xi_\chop=\mp 2,&\ &\eta_\chop=0 \}, \\
  \cR_{\cface,{\rm out}} &:= \{ (R,\omega; \xi_\chop,\eta_\chop) \colon R=0,\ \xi_\chop=0,&\ &\eta_\chop=0 \}
\end{alignat*}
of ${}^{0\semi,\chop} T^*_\sface\iota^+$ (lying over $\sface\cap R^{-1}(0)$). The analogous statement holds for the radial sets $\pa\cR^\pm_{\scri^+,{\rm in,+/out}}$ (see Definition~\ref{DefTsebRad}), whose preimages under $\Phi_{\pm h^{-1}}$ are
\begin{equation}
\label{EqNipOrderRadScri}
\begin{alignedat}{2}
  \cR_{\scri^+,\pm 1,{\rm in}} &:= \{ (x_\sscri,\omega; \xi_{0\semi},\eta_{0\semi}) \colon x_\sscri=0,\ \xi_{0\semi}=\pm 2,&\ &\eta_{0\semi}=0 \}, \\
  \cR_{\scri^+,{\rm out}} &:= \{ (x_\sscri,\omega; \xi_{0\semi},\eta_{0\semi}) \colon x_\sscri=0,\ \xi_{0\semi}=0,&\ &\eta_{0\semi}=0 \}.
\end{alignedat}
\end{equation}
The structure of the future null-bicharacteristic flow can be read off from Proposition~\ref{PropTse3bDyn} via pullback: null-bicharacteristics are of type $(A,B)$ where the allowed pairs are:
\begin{itemize}
\item $A=\cR_{\scri^+,\pm 1,{\rm in}}$ and $B=\cR_{\scri^+,{\rm out}}$ or $\cR_{\cface,\pm 1,{\rm in}}$;
\item $A=\cR_{\cface,\pm 1,{\rm in}}$ and $B=\cR_{\cface,{\rm out}}$;
\item $A=\cR_{\cface,{\rm out}}$ and $B=\cR_{\scri^+,{\rm out}}$.
\end{itemize}
The semiclassical characteristic set is a compact subset of ${}^{0\semi,\chop}T^*_\sface\iota^+$, i.e., disjoint from fiber infinity.

As for the operator $P_{\iota^+,h,z}=h^2 N_{\iota^+}^0(P_0,h^{-1}z)$ itself, we compute it, in the coordinates $R=\frac{r}{t_*}$ (for $R<2$) and $x_\sscri=R^{-\frac12}$ (for $v=R^{-1}<2$), to be
\begin{align}
\label{EqNipOpch}
\begin{split}
  P_{\iota^+,h,z} &= -2 R^{-1}(i z+h R\pa_R)(h R\pa_R + h(1+S|_{\pa X})) \\
    &\quad \hspace{7em} + R^{-2}h^2 P_{(0)}(0,\omega,-R\pa_R,\pa_\omega)
\end{split} \\
\label{EqNipOp0h}
    &= x_\sscri^2 \Bigl[ -\frac12(h x_\sscri\pa_{x_\sscri}-2 i z)(h x_\sscri\pa_{x_\sscri}-2 h(1+S|_{\pa X})) + h^2 x_\sscri^2 P_{(0)}(0,\omega,\tfrac12 x_\sscri\pa_{x_\sscri},\pa_\omega)\Bigr]
\end{align}
Its $(0\semi,\chop)$-principal symbol is therefore
\begin{equation}
\label{EqNipOpSymb}
\begin{split}
  G_{\iota^+,z} &= 2(z+R\xi_\chop)\xi_\chop + \xi_\chop^2+|\eta_\chop|^2 \\
    &= x_\sscri^2\Bigl[ \frac12(\xi_{0\semi}-2 z)\xi_{0\semi} + |\eta_{0\semi}|^2 \Bigr].
\end{split}
\end{equation}
(For $z=\pm 1$, one can then check the above statements about the null-bicharacteristic flow by direct computation.)

\begin{figure}[!ht]
\centering
\includegraphics{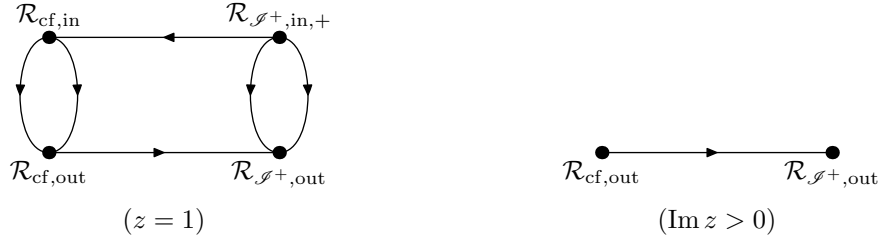}
\caption{\textit{On the left:} characteristic set and null-bicharacteristic dynamics for $(P_{\iota^+,h,z})_{h\in(0,1)}$ in ${}^{0\semi,\chop}T^*_\sface\iota^+$ for $z=1$. \textit{On the right:} characteristic set and null-bicharacteristic dynamics when $\Im z>0$.}
\label{FigNipHiFlow}
\end{figure}

We will also need to study complex $z$; but when $\Im z\neq 0$, then $\Im G_{\iota^+,z}=2(\Im z)\xi_\chop=0$ implies $\xi_\chop=0$, and then $\Re G_{\iota^+,z}=0$ forces $\eta_\chop=0$, similarly in the semiclassical 0-phase space, so the characteristic set is equal to the zero section, with the flow of the Hamiltonian vector field of the real part of
\begin{equation}
\label{EqNipComplex}
  z^{-1}G_{\iota^+,z}
\end{equation}
going from $\cR_{\cface,{\rm out}}$, which is now a source, to the sink $\cR_{\scri^+,{\rm out}}$. See Figure~\ref{FigNipHiFlow}.

\subsubsection{High-energy estimates; proof of Proposition~\usref{PropNip}\eqref{ItNipHi}}
\label{SssNipHi}

In this section, we prove part~\eqref{ItNipHi} of Proposition~\ref{PropNip}. We closely follow \cite[\S{5.3}]{HintzNonstat} but extend it to cover also large $\Im\sigma$ as well as additional degrees b-regularity.

\pfstep{Step~1. Symbolic estimate.} Using radial point, real principal type propagation, and microlocal elliptic estimates, we first prove the estimate
\begin{equation}
\label{EqNipHiSymb}
  \|u\|_{H_{0\semi,\chop,h}^{s,(2\gamma_\sscri,\gamma_\cface,\gamma_\cface,\sfb)}(\iota^+)} \leq C\Bigl( \|P_{\iota^+,h,z}u\|_{H_{0\semi,\chop,h}^{s-2,(2\gamma_\sscri+2,\gamma_\cface-2,\gamma_\cface,\sfb+1)}(\iota^+)} + \|u\|_{H_{0\semi,\chop,h}^{-N,(2\gamma_\sscri,\gamma_\cface,\gamma_\cface,\sfb_\flat)}(\iota^+)} \Bigr)
\end{equation}
for any fixed $N$, where $\sfb_\flat<\sfb$ is induced by an admissible order function $\sfs_\flat<\sfs$ when $z$ is near $\pm 1$, and is a constant when $\Im z>0$. Since the underlying microlocal estimates are very similar to those proved in~\S\ref{SsSpHi}, we shall discuss them rather briefly and leave the statements of standalone propagation estimates to the reader.

\pfsubstep{(1.1)}{Real principal type propagation.} When $|\Im z|=\cO(h)$, Proposition~\ref{PropSpHiPr1} applies. Consider then the case $h\ll\Im z\lesssim 1$: we cannot appeal to Propositions~\ref{PropSpHiPr2}--\ref{PropSpHiPr3} since $\frac{\dd\tau}{\tau}$ is only null, not timelike for the dual metric $g_{\iota^+}^{-1}$ from~\eqref{EqNipDualMet}. Recall that the semiclassical characteristic set of $(P_{\iota^+,h,z})_{h\in(0,1)}$ is the zero section. We instead rescale the operator $P_{\iota^+,h,z}$ similarly to~\eqref{EqSbPInftyOp}, \eqref{EqSpBInftyImTilde}, and~\eqref{EqSpBInftyReP0}--\eqref{EqSpBInftyOutSymb}. We briefly describe this for
\begin{equation}
\label{EqNipHizResc}
\begin{split}
  z^{-1}x_\sscri^{-2}P_{\iota^+,h,z} &= -h x_\sscri D_{x_\sscri} + \frac{\bar z}{|z|^2}\Bigl(\frac12(h x_\sscri D_{x_\sscri})^2+h^2 x_\sscri^2 P_{(0)}(0,\omega,\tfrac12 x_\sscri\pa_{x_\sscri},\pa_\omega)\Bigr) \\
    &\qquad + h z^{-1}(h x_\sscri\pa_{x_\sscri}-2 i z)(1+S).
\end{split}
\end{equation}
(We use~\eqref{EqNipOp0h} to obtain this expression.) In a standard positive commutator argument then, the real part of this operator has symbol $-\xi_{0,\semi}$ plus corrections that vanish quadratically at the zero section $\{\xi_{0,\semi}=\eta_{0,\semi}=0\}$, and hence at the zero section we will propagate along the Hamiltonian vector field $-x_\sscri\pa_{x_\sscri}$ of $-\xi_{0,\semi}$. The second line of~\eqref{EqNipHizResc} is subprincipal. In the term multiplying $\frac{\bar z}{|z|^2}$, note first that $\Im\frac{\bar z}{|z|^2}=-|z|^{-2}\Im z\ll-h$; moreover, the operator in parentheses is equal to $h^2$ times a sum of squares, up to semiclassically subprincipal terms, and thus its overall contribution to the commutator calculation has the ``correct'' sign (namely, the same sign as the main term in the usual positive commutator setup), i.e., is negative as in~\eqref{EqSpBInftyOutSymb} and \eqref{EqSpBInftyOutSign}. (Note that, as usual, the sign of this term determines the direction of propagation.)

\pfsubstep{(1.2)}{Radial point estimates near $\iota^+\cap\scri^+$.} The radial point estimates at $\cR_{\scri^+,\pm 1,{\rm in}}$ and $\cR_{\scri^+,{\rm out}}$ are closely related to Propositions~\ref{PropR3RScriI} and \ref{PropR3RScriO}, with $-\Im\sigma$ and $\gamma_\sscri$ playing the roles of the decay rates $\alpha_++2$ and $\alpha_\sscri+\frac32$ (the shifts arising from passage to a b-density on spacetime, cf.\ the proof of Lemma~\ref{LemmaMUCe3bSob}), respectively, so~\eqref{EqR3RScriIThr} becomes $\Im\sigma>-\gamma_\sscri$, while~\eqref{EqR3RScriO} (with $\ubar p_1=\ubar S$ for $P_0$) becomes $\gamma_\sscri<1+\ubar S$.

\pfsubstep{(1.2.1)}{Incoming radial set.} We first discuss the estimate near $\cR_{\scri^+,+1,{\rm in}}$ for $z\in 1+[-h\gamma_+,\nu]$ in more detail, roughly following the proof of Proposition~\ref{PropSpHiInftyIn}.\footnote{By contrast with Proposition~\ref{PropSpHiInftyIn}, the large $\Im\sigma$ regime in the present discussion is related to estimates for $N_{\iota^+}(P_0)$ on function spaces with large \emph{polynomial} growth $\tau^{i\sigma}=t_*^{-i\sigma}$, which is different than the strong \emph{exponential} growth of $e^{-i\sigma t_*}$.} Write $G_{\iota^+}(\xi)=g_{\iota^+}^{-1}(\xi,\xi)$. The analogue of Lemma~\ref{LemmaSpImag} is that
\[
  h^{-1}\Im P_{\iota^+,h,z} = \bigl( h \Im(N_{\iota^+}^0(P_0,0)) + \Re(z)Q_0 \bigr) - h^{-1}\Im(z) Q_{1,h};
\]
using~$\pa_\sigma N_{\iota^+}^0(P_0,0)=i x_\sscri^2(x_\sscri\pa_{x_\sscri}-2(1+S))$, we have
\begin{equation}
\label{EqNipHiQ0symb}
  Q_0:=\Im\bigl(\pa_\sigma N_{\iota^+}^0(P_0,0)\bigr) = -x_\sscri^2(3+2\Re S),
\end{equation}
where we use the b-density $|\frac{\dd x_\sscri}{x_\sscri}\,\dd\slg|$ to compute adjoints. Concerning $\Im(N_{\iota^+}^0(P_0,0))$, note that, to leading order at $x_\sscri=0$ as a 0-differential operator, we can replace $x_\sscri^4 P_{(0)}$ by $x_\sscri^4 P_{(0)}(0,\omega,0,\pa_\omega)$, whose imaginary part has the schematic form $x_\sscri^4\pa_\omega^{\leq 1}$, and thus the semiclassical 0-principal symbol of $h\Im(x_\sscri^4 P_{(0)})$, so that of $x_\sscri^3\,h x_\sscri\pa_\omega^{\leq 1}$, vanishes cubically, and thus more than quadratically, at $x_\sscri=0$. The principal symbol of $h$ times the imaginary part of the term $-\frac12 x_\sscri^2\,x_\sscri\pa_{x_\sscri}(x_\sscri\pa_{x_\sscri}-2(1+S))$ of $N_{\iota^+}^0(P_0,0)$ is $x_\sscri^2(2+\Re S)\xi_{0\semi}$, the sum of which with~\eqref{EqNipHiQ0symb} is
\[
  x_\sscri^2\bigl( 2\xi_{0\semi}-3 + (\Re S)(\xi_{0\semi}-2) \bigr).
\]
(Cf.\ the arguments after~\eqref{EqR3RScripsub}.) Lastly, the operator
\[
  Q_{1,h}:=-h\Re(\pa_\sigma N_{\iota^+}^0(P_0,0))
\]
has semiclassical principal symbol $\upsigma_{0\semi}(Q_{1,h})=x_\sscri^2\xi_{0\semi}$.\footnote{This equals $2 g_{\iota^+}^{-1}((\Re z)\frac{\dd\tau}{\tau}+\xi,-\frac{\dd\tau}{\tau})=-2 g_{\iota^+}^{-1}(\xi,\frac{\dd\tau}{\tau})$, matching Lemma~\ref{LemmaSpImag}.} This equals $2 x_\sscri^2$ at $\cR_{\scri^+,+1,{\rm in}}$ by~\eqref{EqNipOrderRadScri}, and is a strictly positive multiple of $x_\sscri^2$ nearby.

Consider first the case of $z\in 1+[-h\gamma_+,h\nu']$ where $\nu'\in\R$ is arbitrary but fixed. The Hamiltonian vector field of $x_\sscri^{-2}G_{\iota^+,1}=\frac12(\xi_{0\semi}-2)\xi_{0\semi}+|\eta_{0\semi}|^2$, with points in phase space given by $\xi_{0\semi}\frac{\dd x_\sscri}{x_\sscri}+\eta_{0\semi}\frac{\dd\omega}{x_\sscri}$, is
\[
  H_{x_\sscri^{-2}G_{\iota^+,1}} = (\xi_{0\semi}-1)(x_\sscri\pa_{x_\sscri}+\eta_{0\semi}\pa_{\eta_{0\semi}}) - 2|\eta_{0\semi}|^2\pa_{\xi_{0\semi}} + x_\sscri H_{|\eta_{0\semi}|^2}
\]
where $H_{|\eta_{0\semi}|^2}:=2\eta_{0\semi}\cdot\pa_\omega$ at the origin of geodesic normal coordinates on $\Sph^2$. We then run a positive commutator argument with $A=\check A^2$ where $\check A=\check A^*\in x_\sscri^{-1-2\gamma_\sscri}\Psi_{0\semi}$ is the quantization of
\begin{equation}
\label{EqNipHiChecka}
  \check a = x_\sscri^{-1-2\gamma_\sscri}\psi(\digamma x_\sscri)\psi(\digamma_\cR\fq),
\end{equation}
where $\psi\in\CIc([0,1))$ is equal to $1$ near $0$ and satisfies $\sqrt{-\psi\psi'}\in\CI$, and $\fq:=(\xi_{0\semi}-2)^2+|\eta_{0\semi}|^2$ is a quadratic defining function of $\cR_{\scri^+,+1,{\rm in}}$ over $\{x_\sscri=0\}$; and we shall first choose $\digamma_\cR$ large and then $\digamma$ large. Analogously to~\eqref{EqSpHiPr2C}, we then consider the $L^2(\iota^+,|\frac{\dd x_\sscri}{x_\sscri}\,\dd\slg|)$-pairing $2 h^{-1}\Im\la\check A P_{\iota^+,h,z} u,\check A u\ra=\la\sC u,u\ra$, where
\begin{equation}
\label{EqNipHisC}
  \sC = \frac{i}{h}[\Re P_{\iota^+,h,z},A] + 2 h^{-1}\check A(\Im P_{\iota^+,h,z})\check A + h^{-1}[[\Im P_{\iota^+,h,z},\check A],\check A],
\end{equation}
with the third summand being subprincipal. On the characteristic set of $P_{\iota^+,h,z}$, the principal symbol of the first term is $x_\sscri^2 H_{x_\sscri^{-2}G_{\iota^+,1}}\check a^2$. Differentiation of $\psi(\digamma_\cR\fq)$ yields a non-positive square at $x_\sscri=0$ when $\digamma_\cR$ is large, and thus on $\supp\check a$ when $\digamma$ is large; and differentiation of $\psi(\digamma x_\sscri)$ yields a non-positive square as well (recalling that $\xi_{0\semi}=2$ at the radial set). Differentiation of the weight $x_\sscri^{-2-4\gamma_\sscri}$ of $\check a^2$ yields $-(4\gamma_\sscri+2)x_\sscri^{-4\gamma_\sscri}$ at $\cR_{\scri^+,+1,{\rm in}}$, and the second term of $\sC$ gives an additional contribution at the radial set given by $2 x_\sscri^{-2-4\gamma_\sscri}x_\sscri^2\cdot(2\cdot 2-3-2 h^{-1}\Im(z))$, for a grand total of $-4 x_\sscri^{-4\gamma_\sscri}(\gamma_\sscri+h^{-1}\Im z)$, which is negative since $h^{-1}\Im z\geq-\gamma_+>-\gamma_\sscri$. We thus get a semiclassical 0-estimate for $u$ near $\cR_{\scri^+,+1,{\rm in}}$, without any a priori control term.

When $\Im z\in[h\nu',\nu]$, then if $\nu'$ is sufficiently large, the positivity of the principal symbol of $Q_{1,h}$ at $\cR_{\scri^+,+1,{\rm in}}$ implies that from the term multiplying $-h^{-1}\Im(z)\leq-\nu'$, one can borrow any amount of ``positivity'' (i.e., negativity given the signs in the above commutator calculation) at the radial set, plus an additional term that is negative by the G\aa{}rding inequality analogously to~\eqref{EqSpHiPr2Gard2}. Thus, one gets a radial point estimate at $\cR_{\scri^+,+1,{\rm in}}$ also in this case.

\pfsubstep{(1.2.2)}{Outgoing radial set.} In order to handle complex $z$, we work with $\tilde P_{\iota^+,h,z}:=z^{-1}P_{\iota^+,h,z}$ (see~\eqref{EqNipHizResc}) analogously to~\eqref{EqSpHiInftyOutTilde}. Consider first $z\in 1+[-h\gamma_+,\nu]$ for small $\nu>0$. Then the semiclassical principal symbol of $\tilde P_{\iota^+,h,z}$ is $G_{\iota^+,1}+\cO(\nu)$. Consider the commutant~\eqref{EqNipHiChecka}, now with the quadratic defining function $\fq:=\xi_{0\semi}^2+|\eta_{0\semi}|^2$ of the radial set $\cR_{\scri^+,{\rm out}}$ over $\{x_\sscri=0\}$. Differentiation of either cutoff now yields a non-negative square when $\nu$ is small enough, while the principal symbol of $\frac{i}{h}[\Re\tilde P_{\iota^+,h,z},A]$ at the radial set is $(4\gamma_\sscri+2)x_\sscri^{-4\gamma_\sscri}$. We split the second term in~\eqref{EqNipHisC} into the sum-of-squares contribution of the second term in~\eqref{EqNipHizResc} as in Step~(1.1) above, the contribution of $2\check a^2 h^{-1}\Im(h z^{-1}x_\sscri^2(h x_\sscri\pa_{x_\sscri}-2 i z)(1+S))$, whose principal symbol at $\cR_{\scri^+,{\rm out}}$ is $-4 x_\sscri^{-4\gamma_\sscri}(1+\Re S)$, and the contribution of $2\check a^2 h^{-1}\Im(-x_\sscri^2 h x_\sscri D_{x_\sscri})$, with principal symbol at the radial set given by $-2 x_\sscri^{-4\gamma_\sscri}$. Under the threshold condition $(4\gamma_\sscri+2)-4(1+\Re S)-2<0$, i.e., $\gamma_\sscri<1+\Re S$, we can thus propagate semiclassical 0-control of $u$ from a punctured neighborhood of $\cR_{\scri^+,{\rm out}}$ into $\cR_{\scri^+,{\rm out}}$.

When $\Im z\geq\nu>0$, the derivative of the cutoff in $\fq$ is now supported in the elliptic set of $\tilde P_{\iota^+,h,z}$. The remaining terms are estimated as before. Thus, we again get a radial sink estimate in which now semiclassical 0-regularity is propagated along the zero section into $\cR_{\scri^+,{\rm out}}$.

\pfsubstep{(1.3)}{Radial point estimates near $\iota^+\cap\cK^+$.} Radial point estimates in the semiclassical cone calculus were first discussed in \cite[\S{4}]{HintzConicProp}. The radial point estimates near $\cR_{\cface,\pm 1,{\rm in}}$ and $\cR_{\cface,{\rm out}}$ are related via the Mellin transform in $\tau=t_*^{-1}$ to Propositions~\ref{PropR3RKI} and \ref{PropR3RKO}, with the following identifications: the 3b-regularity is $s=\sfb$; and denoting the decay orders relative to $L^2(\R^4,|\dd g_{\bhm,a}|)$ by $\tilde\alpha_+$ and $\tilde\alpha_\cK$ for clarity, and in view of $R^{\tilde\alpha_\cK+\frac12}(\frac{\tau}{R})^{\tilde\alpha_++2}=R^{\tilde\alpha_\cK-\tilde\alpha_+-\frac32}\tau^{\tilde\alpha_++2}$, we set $\tilde\alpha_++2=-\Im\sigma$ and $\gamma_\cface=\tilde\alpha_\cK-\tilde\alpha_+-\frac32$, so $\tilde\alpha_\cK=\tilde\alpha_+-\alpha_+$ since $\gamma_\cface=-(\alpha_++\frac32)$; the threshold condition~\eqref{EqR3RKIThr}, with $s_0=s$, then reads $\sfb+\tilde\alpha_+>\frac12(-1+\vartheta_{\pa\cK^+,{\rm in}})+\tilde\alpha_\cK$, or equivalently,
\begin{equation}
\label{EqNipHicfInThr}
  \sfb + \alpha_+ > \frac12(-1+\vartheta_{\pa\cK^+,{\rm in}}),
\end{equation}
which indeed holds for the induced semiclassical order $\sfb$ by~\eqref{EqSSOrderAdmIn}; similarly for the threshold condition~\eqref{EqR3RKOThr}.\footnote{In other words, since $-\Im\sigma$ is the $\tau$-decay rate, and $\tau=\rho_+\rho_\cK$ is a \emph{joint} defining function of $\iota^+$ and $\cK^+\subset M$, the relative order at $\cK^+$ and $\iota^+$, i.e., $\tilde\alpha_\cK-\tilde\alpha_+$ in present notation, is independent of $\Im\sigma$.} The details are as follows.

\pfsubstep{(1.3.1)}{Incoming radial set.} We first work near $\cR_{\cface,+1,{\rm in}}$, with $z\in 1+[-h\gamma_+,\nu]$. Consider first the case $z\in 1+[-h\gamma_+,h\nu']$ where $\nu'\in\R$ is arbitrary but fixed. The Hamiltonian vector field of $G_{\iota^+,1}$ in~\eqref{EqNipOpSymb}, with points in phase space given by $\xi_\chop\frac{\dd R}{h}+\eta_\chop\cdot\frac{R\,\dd\omega}{h}$ (which matches~\eqref{EqNipCovecchop} at $\frac{h}{h+R}=0$, i.e., when setting $h=0$ for $R>0$), is given by
\[
  H_{G_{\iota^+,1}} = \hat h \Bigl( 2\bigl((1+2 R)\xi_\chop+1\bigr)(R\pa_R-\eta_\chop\pa_{\eta_\chop}) + 2(|\eta_\chop|^2-R\xi_\chop^2)\pa_{\xi_\chop} + H_{|\eta_\chop|^2}\Bigr),\quad\hat h:=\frac{h}{R},
\]
where we write $H_{|\eta_\chop|^2}:=2\eta_\chop\cdot\pa_\omega$ at the origin of normal coordinates on $\Sph^2$; this follows from~\eqref{EqNipOpSymb} via a change of variables computation. Note that $R$, resp.\ $\hat h$ are local defining functions of $\tface$, resp.\ $\sface\subset\iota^+_{0\semi,\chop}$ near $\sface\cap\tface$. We recall that $\cR_{\cface,+1,{\rm in}}=\{R=0,\ \xi_\chop=-2,\ \eta_\chop=0\}$. We then run a positive commutator argument with $A=\check A^2$, where $\check A=\check A^*\in R^{-\gamma_\cface}\hat h^{-\sfb-\frac12}\Psi_\chop^{-\infty}$ is a quantization of
\begin{equation}
\label{EqNipHicfInChecka}
  \check a := R^{-\gamma_\cface}\hat h^{-\sfb-\frac12}\psi(\digamma R)\psi(\digamma_\cR\fq) = h^{-\sfb-\frac12} R^{\sfb-\gamma_\cface+\frac12}\psi(\digamma R)\psi(\digamma_\cR\fq),
\end{equation}
where $\fq:=(\xi_\chop+2)^2+|\eta_\chop|^2$, and $\digamma,\digamma_\cR$ are large. We now consider $2\Im\la\check A P_{\iota^+,h,z}u,\check A u\ra=\la\sC u,u\ra$, where now
\[
  \sC = i[\Re P_{\iota^+,h,z},A] + 2\check A(\Im P_{\iota^+,h,z})\check A + [[\Im P_{\iota^+,h,z},\check A],\check A].
\]
Since commutators in the semiclassical cone calculus gain a power of $\hat h$, the third term is of class $R^{-\gamma_\cface}\hat h^{-2\sfb+1}\Psi_\chop$ and thus subprincipal in $R^{-\gamma_\cface}\hat h^{-2\sfb}\Psi_\schop$. Upon differentiation of $\check a^2$ along $H_{G_{\iota^+,1}}$, the term $\psi(\digamma R)$ yields a non-negative square, while $\psi(\digamma_\cR\fq)$ yields a non-positive square since $\sface\cap\tface\cap\fq^{-1}(0)=\cR_{\cface,+1,{\rm in}}$ is a source for the $H_{G_{\iota^+,1}}$-flow over $\sface\cap\tface$. Differentiation of the $R$-weight of $\check a^2$ yields $-R^{-2\gamma_\cface}\hat h^{-2\sfb}(4\sfb-4\gamma_\cface+2)$ at the radial set.

We next compute
\begin{align*}
  &\Im P_{\iota^+,h,z} = \bigl( h^2\Im(N_{\iota^+}^0(P_0,0))+\Re(z) Q_{0,h} \bigr) - \Im(z)Q_{1,h}, \\
  &\qquad Q_{0,h} := \Im\bigl(h \pa_\sigma N_{\iota^+}^0(P_0,0)\bigr),\ \ Q_{1,h}=\Re\bigl(-h\pa_\sigma N_{\iota^+}^0(P_0,0)\bigr);
\end{align*}
we are interested in the symbols to order $\cO(\hat h)$ (which thus contribute at the same level as the term $i[\Re P_{\iota^+,h,z},A]$ just discussed), while terms of class $\cO(\hat h^2)$ and $\cO(\hat h R)$ are to be considered subprincipal and can be neglected. Recall that near $\sface$, semiclassical cone vector fields are spanned by $h\pa_R$ and $h R^{-1}\pa_\omega$. Recall now~\eqref{EqNipOpMT}, and note that the term $-2 h^2 R^{-1}R\pa_R(R\pa_R+1+S)$ of $h^2 N_{\iota^+}^0(P_0,0)$ lies in $R\Diff_\chop^2$ and has real principal symbol $R\xi_\chop^2$, so its imaginary part is of class $R\hat h\Diff_\chop^1$ and thus subprincipal. The term $h^2\Im(R^{-2}P_{(0)}(0,\omega,-R\pa_R,\pa_\omega))$ is of class $h^2 R^{-2}\Diffb^1$ and thus of class $\hat h\Diff_\chop^1$, with principal symbol being $\hat h$ times a fiber-linear function; this is equal to the symbol of imaginary part of $h^2 R^{-2}P_{(0)}$ computed with the density $R^3\,|\frac{\dd R}{R}\,\dd\slg|$ plus a correction term given by the symbol of $-h^2\frac{i}{2}[R^{-2}P_{(0)},R^3]R^{-3}$, which is $-\frac32 h^2 \pa_{\xi_\bop}(R^{-2}\xi_\bop^2)$ (where we write $\xi_\bop\frac{\dd R}{R}+\eta_\bop\cdot\dd\omega$ for b-covectors) and thus $-3\hat h\xi_\chop$ (since $\xi_\bop=\frac{R}{h}\xi_\chop$). We add all of this to the symbol of $Q_{0,h}$; since the imaginary part of $\pa_\sigma N_{\iota^+}^0(P_0,0)=-2 i R^{-1}(R\pa_R+1+S)$ is $-R^{-1}(3+2\Re S)$, we have $\upsigma_\chop(Q_{0,h})=-\hat h(3+2\Re S)$. This gives $\hat h$ times $-3-3\xi_\chop-\frac{\xi_\chop}{2}\vartheta_{\pa\cK^+,{\rm in}}$ (the prefactor $-\frac{\xi_\chop}{2}$ being equal to $1$ at $\cR_{\cface,+1,{\rm in}}$). The grand total at the radial set thus far is $R^{-2\gamma_\cface}\hat h^{-2\sfb}(-4\sfb+4\gamma_\cface-2+2\cdot(-3-3\cdot(-2)+\vartheta_{\pa\cK^+,{\rm in}}))$, which (upon plugging in $\gamma_\cface=-(\alpha_++\frac32)$) is negative under the condition~\eqref{EqNipHicfInThr}.

Finally, the $\chop$-principal symbol of $Q_{1,h}$ is $-2\xi_\chop$, so equal to $4$ at $\cR_{\cface,+1,{\rm in}}$, and hence positive nearby. But since $4\Im(z)=4\Im(h^{-1}z)\hat h R$, it has no effect on the threshold condition at the radial set. Altogether, under the condition~\eqref{EqNipHicfInThr}, we can thus propagate semiclassical cone regularity from $\supp\psi'(\digamma R)\psi(\digamma_\cR\fq)$ (which is contained in a neighborhood of the stable manifold of $\cR_{\cface,+1,{\rm in}}$) into the radial set.

When $z\in 1+[h\nu',\nu]$, the arguments are the same; the only difference is that we now estimate the term $-2\Im(z)\check A Q_{1,h}\check A$ using the G\aa{}rding inequality, using the positivity of the principal symbol of $Q_{1,h}$: upon applying this to $u$ and pairing with $u$, one obtains a non-positive contribution up to the usual microlocal error term (of arbitrarily negative $\sface$-order).

\pfsubstep{(1.3.2)}{Outgoing radial set.} Analogously to step (1.2.2), we now work with
\begin{align*}
  \tilde P_{\iota^+,h,z} &= z^{-1}P_{\iota^+,h,z} \\
    &= 2 h D_r + \frac{\bar z}{|z|^2}\Bigl( R^{-2}h^2 P_{(0)}(0,\omega,-R\pa_R,\pa_\omega) + 2 R\cdot R^{-2}h^2(R D_R)^2\Bigr) \\
    &\qquad - 2(1+S)(h z^{-1}h\pa_R + i R^{-1}h).
\end{align*}
The term multiplying $\frac{\bar z}{|z|^2}$ has principal part $h^2$ times $(1+2 R)(D_R)^2+(R^{-1})^2\slDelta$, which is a sum of squares of vector fields. In the $L^2$-pairing in the commutator computation, this thus gives a non-positive contribution for $\Im\bar z=-\Im z<0$. We run a positive commutator argument for $\tilde P_{\iota^+,h,z}$, with commutant~\eqref{EqNipHicfInChecka} where now $\fq:=\xi_\chop^2+|\eta_\chop|^2$ defines $\cR_{\cface,{\rm out}}$ over $\tface\cap\sface=\{R=\frac{h}{R}=0\}$. The second line is subprincipal at $\sface=\{\frac{h}{R}=0\}$. When $z\in 1+[-h\gamma_+,\nu]$ for small $\nu>0$, differentiation of the cutoff in $\fq$ along the Hamiltonian vector field of $\Re\tilde P_{\iota^+,h,z}$ (which is $\nu$-close to $H_{2\xi_\chop}=2\hat h R\pa_R=2 h\pa_R$, plus an additional term that vanishes as a $\chop$-vector field at the zero section) yields a non-negative square (necessitating an a priori control term), while differentiation of the cutoff in $R$ yields a non-positive square. Differentiation of the $R$-weight of $\check a^2$ gives $R^{-2\gamma_\cface}\hat h^{-2\sfb}(4\sfb-4\gamma_\cface+2)$. Twice the imaginary part of $2 h D_r-2(1+S)(h z^{-1}h\pa_R+i R^{-1}h)$ contributes a further term $R^{-2\gamma_\cface}\hat h^{-2\sfb}(-2-(4+4\Re S))$ at the radial set, for a grand total of $R^{-2\gamma_\cface}\hat h^{-2\sfb}$ times $4\sfb+4\alpha_++2+4\Re S$. This is negative when $\sfb+\alpha_+<-\frac12+\Re S$; in view of Lemma~\ref{LemmaSSOrderThr}, this matches the threshold condition~\eqref{EqSSOrderAdmOut} (with $\alpha_\cK=\eta=0$, cf.\ \eqref{EqNipOrder}). Under this condition, we can thus propagate $\chop$-regularity from the stable manifold of $\cR_{\cface,{\rm out}}$ over $\sface\cap\tface$ into the radial set (and then out towards $\cR_{\scri^+,{\rm out}}$).

The argument for $\Im z\geq\nu$ is similar, the only difference being that now the derivative of the cutoff in $\fq$ is supported in the elliptic set of $\tilde P_{\iota^+,h,z}$.

\pfstep{Step~2. Improvement at the transition face.} The error term in~\eqref{EqNipHiSymb} is not yet small as $h\to 0$ since the order $\gamma_\cface$ at the transition face is the same as on the left-hand side. We improve this by inverting the $\tface$-normal operator of $(P_{\iota^+,h,z})_{h\in(0,1)}$. Let $\chi\in\CI([0,1))$ be equal to $1$ on $(-\infty,\frac12]$, and, recalling~\eqref{EqMUipSingle}, consider the cutoff
\[
  \chi_\tface=\chi(h+R)\chi\Bigl(\frac{h+R}{h/(h+R)}\Bigr)
\]
which localizes to a neighborhood of $\tface\subset\iota^+_{0\semi,\chop}$ and the set where $h+R$ (which defines $\tface$ near $\tface\cap\sface$) is smaller than $\frac{h}{h+R}$ (which defines $\sface$ near $\tface\cap\sface$); note that $\chi_\tface\in\CI_\bop(\iota^+_{0\semi,\chop})$ (in fact, this is smooth on the blow-up of $\iota^+_{0\semi,\chop}$ at $\sface\cap\tface$); the idea of using such a conormal cutoff (which gives useful flexibility in~\eqref{EqNipHiNtoP} below) is taken from the proof of \cite[Theorem~4.10]{HintzConicProp}. Using~\eqref{EqMUipNormEquiv}, we then have
\[
  \|\chi_\tface u\|_{H_{0\semi,\chop,h}^{-N,(2\gamma_\sscri,\gamma_\cface,\gamma_\cface,\sfb_\flat)}(\iota^+)} \lesssim h^{-\gamma_\cface}\|\chi_\tface u\|_{H_{\bop,\scop}^{-N,(\gamma_\cface,\sfb_\flat-\gamma_\cface)}(\tface)}.
\]
Using the projective coordinate $\hat R=\frac{R}{h}$ on $\tface=[0,\infty]_{\hat R}\times\Sph^2$, let us compute the $\tface$-normal operator
\[
  N_{\tface,z}\in\Diff_{\bop,\scop}^{2,(-2,0)}(\tface)=\rho_\cface^{-2}\Diff_{\bop,\scop}^2(\tface);
\]
of $(P_{\iota^+,h,z})_{h\in(0,1)}$; here we can use as boundary defining functions of $\tface\cap\cface$ and $\tface\cap\sface$ the functions $\rho_\cface=\frac{\hat R}{\hat R+1}$ and $\rho_\sface=\frac{1}{\hat R+1}$. We find $N_{\tface,z}$ by expressing~\eqref{EqNipOpch} in the coordinates $\hat R,\omega$ and setting $h=0$ (thereby restricting to $\tface$), as discussed in general in~\eqref{EqMUNtfch}; this gives
\begin{align*}
  N_{\tface,z} &= -2 i z\hat R^{-1}(\hat R\pa_{\hat R}+1+S|_{\pa X}) + \hat R^{-2}P_{(0)}(0,\omega,-\hat R\pa_{\hat R},\pa_\omega) \\
    &= 2 i z\hat\rho(\hat\rho\pa_{\hat\rho}-1-S|_{\pa X}) + \hat\rho^2 P_{(0)}(0,\omega,\hat\rho\pa_{\hat\rho},\pa_\omega),\quad \hat\rho:=\hat R^{-1}.
\end{align*}
This is thus the same as the transition face normal operator $N_\tface(P_0,z)$ for the low-energy spectral family introduced in~\eqref{EqSStfOp} and analyzed in~\S\ref{SsSptf}. (The equality of the two transition face normal operators holds for general 3b-operators, as shown in \cite[Propositions~3.28 and 4.17]{Hintz3b}.) In order to make the function spaces compatible, we pass from the presently used b-density $|\frac{\dd\hat R}{\hat R}\,\dd\slg|$ on $\tface$ to the ``Euclidean'' density $\hat R^3\,|\frac{\dd\hat R}{\hat R}\,\dd\slg|=\rho_\cface^3\rho_\sface^{-3}\,|\frac{\dd\hat R}{\hat R}\,\dd\slg|$, so
\[
  \|\chi_\tface u\|_{H_{\bop,\scop}^{-N,(\gamma_\cface,\sfb_\flat-\gamma_\cface)}(\tface,|\frac{\dd\hat R}{\hat R}\,\dd\slg|)} = \|\chi_\tface u\|_{H_{\bop,\scop}^{-N,(\gamma_\cface+\frac32,\sfb_\flat-\gamma_\cface-\frac32)}(\tface,\hat R^3\,|\frac{\dd\hat R}{\hat R}\,\dd\slg|)}.
\]
We estimate this using Corollary~\ref{CorSptfAdm}. (The order of the subscripts ``$\bop$'' and ``$\scop$'' here is reversed relative to the notation in Corollary~\ref{CorSptfAdm}, but the b-behavior and the scattering behavior refers to the boundary at $\hat R=\hat\rho^{-1}=0$ and $\hat\rho=0$, respectively, exactly as in Corollary~\ref{CorSptfAdm}.) Note that $\gamma_\cface+\frac32=-\alpha_+$ and recall that $P_0$ is, by assumption, $\tface$-admissible with weight $\alpha_++\frac32$; thus, Corollary~\ref{CorSptfAdm} indeed applies (with $n=3$, $q=\gamma_\cface+\frac32=-\alpha_+$, $\sfr=\sfb_\flat$) and gives
\begin{align*}
  \|\chi_\tface u\|_{H_{\bop,\scop}^{-N,(\gamma_\cface+\frac32,\sfb_\flat-\gamma_\cface-\frac32)}(\tface,\hat R^3\,|\frac{\dd\hat R}{\hat R}\,\dd\slg|)} &\leq C\|N_{\tface,z}(\chi_\tface u)\|_{H_{\bop,\scop}^{-N-2,(\gamma_\cface-\frac12,\sfb_\flat-\gamma_\cface-\frac12)}(\tface,\hat R^3\,|\frac{\dd\hat R}{\hat R}\,\dd\slg|)} \\
    &= C\|N_{\tface,z}(\chi_\tface u)\|_{H_{\bop,\scop}^{-N-2,(\gamma_\cface-2,\sfb_\flat-\gamma_\cface+1)}(\tface,|\frac{\dd\hat R}{\hat R}\,\dd\slg|)}.
\end{align*}
Therefore, since the order at $\iota^+\cap\scri^+$ is irrelevant when working near $\tface$,
\[
  \|\chi_\tface u\|_{H_{0\semi,\chop,h}^{-N,(2\gamma_\sscri,\gamma_\cface,\gamma_\cface,\sfb_\flat)}(\iota^+)} \leq C\|N_{\tface,z}(\chi_\tface u)\|_{H_{0\semi,\chop,h}^{-N-2,(-N,\gamma_\cface,\gamma_\cface-2,\sfb_\flat+1)}(\iota^+)}.
\]
We now replace $N_{\tface,z}$ by $P_{\iota,h,z}$, the difference of the two operators lying in $\rho_\cface^{-2}\rho_\tface\Diff_\chop^2$ near $\tface$. Note that $\chi_\tface\in\rho_\sface^\eta\rho_\tface^{-\eta}\CI_\bop$ for all $\eta\geq 0$; indeed, if we take $\rho_\sface=\frac{h}{h+R}$ and $\rho_\tface=h+R$ as defining functions, then on the support of $\chi(\frac{h+R}{h/(h+R)})=\chi(\frac{\rho_\tface}{\rho_\sface})$, we have $\frac{\rho_\sface}{\rho_\tface}\gtrsim 1$. Thus,
\begin{equation}
\label{EqNipHiNtoP}
\begin{split}
  &\|N_{\tface,z}(\chi_\tface u)\|_{H_{0\semi,\chop,h}^{-N-2,(-N,\gamma_\cface,\gamma_\cface-2,\sfb_\flat+1)}} \\
  &\qquad \leq C\Bigl( \|\chi_\tface P_{\iota^+,h,z}u\|_{H_{0\semi,\chop,h}^{-N-2,(-N,\gamma_\cface,\gamma_\cface-2,\sfb_\flat+1)}} + \|u\|_{H_{0\semi,\chop,h}^{-N,(-N,\gamma_\cface,\gamma_\cface-1+\eta,\sfb_\flat+1-\eta)}}\Bigr),
\end{split}
\end{equation}
where we also use that $[P_{\iota^+,h,z},\chi_\tface]\in\rho_\cface^\infty\rho_\tface^{-\eta'}\rho_\sface^{1+\eta'}\CI_\bop\Diff_\chop^1$ for all $\eta'\in\R$ near $\tface$, here with $\eta'=-1+\eta$. If we use this with $\eta=1-\delta$, the error term here is
\[
  \|u\|_{H_{0\semi,\chop,h}^{-N,(-N,\gamma_\cface,\gamma_\cface-\delta,\sfb_\flat+\delta)}(\iota^+)} \leq \|u\|_{H_{0\semi,\chop,h}^{-N,(-N,\gamma_\cface,\gamma_\cface-\delta,\sfb-\delta)}}
\]
when $\delta>0$ is sufficiently small. This, in turn, is bounded by $C h^\delta$ times the left-hand side of~\eqref{EqNipHiSymb}, and can therefore be absorbed when $h>0$ is small enough. We have thus shown that there exists $h_0>0$ such that
\begin{equation}
\label{EqNipHiFinal}
  \|u\|_{H_{0\semi,\chop,h}^{s,(2\gamma_\sscri,\gamma_\cface,\gamma_\cface,\sfb)}(\iota^+)} \leq C\|P_{\iota^+,h,z}u\|_{H_{0\semi,\chop,h}^{s-2,(2\gamma_\sscri+2,\gamma_\cface-2,\gamma_\cface,\sfb+1)}(\iota^+)},\quad 0<h<h_0,
\end{equation}
which is the estimate~\eqref{EqNipb} for $k=0$.

An analogous estimate for the adjoint $P_{\iota^+,h,z}^*$ on dual spaces is proved similarly; together with the estimate~\eqref{EqNipHiFinal}, this implies the invertibility of $P_{\iota^+,h,z}$ for small $h>0$. Now, the space $H_{0\semi,\chop,h}^{s,(2\gamma_\sscri,\gamma_\cface,\gamma_\cface,\sfb)}(\iota^+)$ is equal to $H_{0,\bop}^{s,(2\gamma_\sscri,\gamma_\cface)}$ (only its norm depends on $h$ and the remaining orders), so we obtain the invertibility of~\eqref{EqNip} for $\sigma=h^{-1}z$.

\pfstep{Step~3. Higher b-regularity.} To prove the estimate~\eqref{EqNipb} for general $k\in\N_0$, we only need to establish analogues of the microlocal estimates used in Step~1 which feature $k$ additional degrees of b-regularity. (In Step~2 then, one uses~\eqref{EqSptfAdm} with $k$ b-derivatives.) We only discuss the radial point estimates. The key observation for the estimates near $\cR_{\scri^+,\pm 1,{\rm in}}$ and $\cR_{\scri^+,{\rm out}}$ is that the structure of the commuted equation in Lemma~\ref{LemmaRbComm} descends to the $\iota^+$-normal operator family, with the vector fields $V_1,\ldots,V_5$, and the corresponding operators $X_1,\ldots,X_5$ and $X_6=I$, replaced by their respective normal operator families $\wh{X_j}(\sigma)$, obtained near $\iota^+\cap\scri^+$ by replacing $\rho_+\pa_{\rho_+}$ in~\eqref{EqRbCommV} by $-i\sigma$; thus the radial point estimates on semiclassical 0-Sobolev spaces discussed in Step~1 are valid, with the same threshold condition, for the commuted equation satisfied by the vector of $k$-fold b-derivatives of $u$ along $\wh{X_j}(\sigma)$, $j=1,\ldots,6$. Note that the right-hand side of the commuted equation involves $k-1$ b-derivatives, and the high-energy estimate loses a power of $h^{-1}$; thus, one pays one power of $h^{-1}$ for each b-derivative, as captured by the norm~\eqref{EqNipNormbplus}. (This is completely analogous to Step~2 in the proof of Proposition~\ref{PropSpHiInftyOut}.)

Near $\iota^+\cap\cK^+$, we note that away from the zero section of ${}^\chop T^*_\sface(\iota^+\setminus\scri^+)$, $k$ orders of b-regularity are, microlocally, the same as $k$ orders each of semiclassical cone regularity and semiclassical orders each. (Indeed, $R\pa_R=\frac{h+R}{h}\cdot\frac{h}{h+R}R\pa_R$ and $\pa_\omega=\frac{h+R}{h}\cdot\frac{h}{h+R}\pa_\omega$, with the symbols of $\frac{h}{h+R}R\pa_R$ and $\frac{h}{h+R}\pa_\omega$ spanning the space of fiber-linear functions on ${}^\chop T^*$, and with $\frac{h+R}{h}$ being the inverse of a defining function of $\sface$.) Notice that this is the \emph{same} as the observation in Remark~\ref{RmkSpBNob} in view of the relationship between semiclassical cone (near $\iota^+\cap\cK^+$) and b-scattering (on $\tface$) phase spaces recalled before~\eqref{EqNipOrdertf}.

We therefore only need to prove a b-regularity result near $\cR_{\cface,{\rm out}}$. The arguments involving Lemma~\ref{LemmaRbROMix} do not descend to the Mellin transform side. Instead, we use commutators with $V_1=R\pa_R$ and a spanning set $V_a\in\cV(\Sph^2)$, $a=2,3,4$, of $\cV(\Sph^2)$ over $\CI(\Sph^2)$ (or, in the case of operators acting on sections of a bundle, the operators $X_i=\nabla^\cE_{V_i}$, $i=1,\ldots,4$, where $\nabla^\cE$ is as in~\eqref{EqSSAdmDer}). The following result is then closely related to Lemma~\ref{LemmaSpBInftyOutComm}.\footnote{In fact, Lemma~\ref{LemmaSpBInftyOutComm} can be obtained from Lemma~\ref{LemmaNipHiComm} by restricting each computation in the proof to $\tface$.}

\begin{lemma}[Commutators near $\tface$]
\label{LemmaNipHiComm}
  Let $V_5:=I$ and set $\sV:=(V_1,V_2,\ldots,V_5)$ and $\sV_{\rm shift}=R^{-2}\sV R^2=(V_1+2,V_2,\ldots,V_5)$. Set $\cA_k:=\{\alpha\in\N_0^5\colon|\alpha|=k\}$. Write $\rho_\cface=\frac{R}{h+R}$, $\rho_\tface=h+R$ and $\rho_\sface=\frac{h}{h+R}$ for local defining functions of $\cface$, $\tface$, and $\sface\subset\iota^+_{0\semi,\chop}$, respectively. If $P_{\iota^+,h,z}u=f$, then upon setting $u^{(k)}:=(\sV^\alpha u)_{\alpha\in\cA_k}$ and $f^{(k)}:=(\sV_{\rm shift}^\alpha f)$, we have
  \begin{equation}
  \label{EqNipHiComm}
    P^{(k)}_{h,z}u^{(k)} = f^{(k)} + \tilde f^{(k)},\quad
    \tilde f^{(k)} = \Biggl(\;\sum_{|\beta|\leq k-1} R_{\alpha\beta}\sV^\beta u \Biggr)_{\alpha\in\cA_k},
  \end{equation}
  where $R_{\alpha\beta}\in\rho_\cface^{-2}\rho_\sface\Diff_\chop^2$, and the $(\alpha,\alpha')$-matrix element of $P^{(k)}_{h,z}$ is given by
  \begin{equation}
  \label{EqNipHiCommMtx}
    (P^{(k)}_{h,z})_{\alpha\alpha'} = \delta_{\alpha\alpha'}P_{\iota^+,h,z} - 2 i k z h R^{-1}\delta_{\alpha,(k,0,\ldots,0)} + \tilde P_{\alpha\alpha'}
  \end{equation}
  where $\tilde P_{\alpha\alpha'}\in\rho_\cface^{-2}\rho_\sface\Diff_\chop^1$ and $\upsigma_\chop(\rho_\sface^{-1}\tilde P_{\alpha\alpha'})|_{\cR_{\cface,{\rm out}}}=0$.
\end{lemma}
\begin{proof}
  Set $P'_{h,z}:=R^2 P_{\iota^+,h,z}$. For $\alpha\in\cA_k$, we then have
  \[
    P_{\iota^+,h,z}(\sV^\alpha u) = R^{-2}\sV^\alpha(R^2 P_{\iota^+,h,z}u) - R^{-2}[\sV^\alpha,R^2 P_{\iota^+,h,z}]u = \sV_{\rm shift}^\alpha f - R^{-2}[\sV^\alpha,P'_{h,z}]u.
  \]
  Using~\eqref{EqNipOpch}, we have
  \begin{equation}
  \label{EqNipHiCommPp}
    P'_{h,z} = h^2 P_{(0)}(0,\omega,-R\pa_R,\pa_\omega) - 2 R(i z+h R\pa_R)(h R\pa_R+h(1+S|_{\pa X})).
  \end{equation}

  Commutators of $\sV^\alpha$ with terms of $P_{(0)}$ of class $\Diffb^1$ lie in $\Diffb^k$ and can thus be written in the schematic form $\Diffb^1\circ\sV^{\leq k-1}$; since $R\pa_R=\rho_\sface^{-1}\rho_\sface R\pa_R$ and $\pa_\omega=\rho_\sface^{-1}\rho_\sface\pa_\omega$, with $\rho_\sface R\pa_R$ and $\rho_\sface\pa_\omega$ being a frame of $\cV_\chop$, we can thus absorb $R^{-2}h^2\Diffb^1\subset\rho_\cface^{-2}\rho_\sface\Diff_\chop^1$ into the operators $R_{\alpha\beta}$. Commutators of $\sV^\alpha$ with terms of $P_{(0)}$ of class $R\,\Diffb^2$ lie in $R\,\Diffb^{k+1}\subset R\,\Diffb^1\,\sV^{\leq k}$; we then absorb $h^2 R^{-2}\,R\,\Diffb^1\subset\rho_\cface^{-1}\rho_\tface\rho_\sface\Diff_\chop^1$ into $\tilde P_{\alpha\alpha'}$. Consider the remaining piece $(R D_R)^2+\slDelta$ of $P_{(0)}$: this commutes with $R\pa_R$, while commutators with $k$ spherical vector fields are of class $\Diff^{k+1}(\Sph^2)\subset\Diff^1(\Sph^2)\circ\sV^{\leq k}$; and the $\chop$-principal symbol of an element of $R^{-2}h^2\Diff^1(\Sph^2)\subset R^{-2}h^2\Diffb^1\subset\rho_\cface^{-2}\rho_\sface\Diff_\chop^1$ is a linear function of the spherical momentum variables $\eta_{0\semi}$, which vanish at $\cR_{\cface,{\rm out}}$; we put such operators into $\tilde P_{\alpha\alpha'}$.

  Next, commutators of spherical vector fields with the second summand of~\eqref{EqNipHiCommPp} are only sensitive to $S|_{\pa X}$; note then that $R^{-2}\cdot 2 R(i z+h R\pa_R)h\CI(\Sph^2)\circ\sV^{\leq k-1}\subset\rho_\cface^{-1}\rho_\sface\Diff_\chop^0\circ\sV^{\leq k-1}$ is a sum of terms of the form $R_{\alpha\beta}\sV^\beta$. It remains to compute
  \[
    R^{-2}[ (R\pa_R)^k, -2 R(i z+h R\pa_R)h R\pa_R ] \equiv -2 k R^{-1}(i z+h R\pa_R)h(R\pa_R)^k \bmod R^{-2}R(h R\pa_R)^{\leq 2}\Diffb^{k-1}.
  \]
  The term $-2 k R^{-1}h R\pa_R\,h\,(R\pa_R)^k\subset\rho_\cface^{-1}\rho_\tface\rho_\sface\Diff_\chop^1\circ\sV^{\leq k}$ contributes to $\tilde P_{\alpha\alpha'}$ (the factor of $\rho_\tface$ guaranteeing the vanishing condition on the principal symbol at $\cR_{\cface,{\rm out}}$), and the error term can be put into $\tilde P_{\alpha\alpha'}$ as well since it is of the form $R^{-2}R(h R\pa_R)^{\leq 2}\Diffb^{k-1}\subset R^{-1}(h R\pa_R)h\,\Diffb^k$, with $R^{-1}(h R\pa_R)h\in\rho_\cface^{-1}\rho_\tface\rho_\sface\Diff_\chop^1$ again. The only term remaining is $-2 i k z h R^{-1}(R\pa_R)^k$, which gives rise to the second term in~\eqref{EqNipHiCommMtx}.
\end{proof}

We may now apply the already established radial point estimate at $\cR_{\cface,{\rm out}}$ to the commuted equation~\eqref{EqNipHiComm}. Note that each term in $\tilde f^{(k)}$ regains the one power of $\rho_\sface$ which the radial point estimate loses (cf.\ the final order on $P_{\iota^+,h,z}u$ in~\eqref{EqNipHiSymb}). Moreover, the additional term $-2 i k z\hat h\delta_{\alpha,(k,0,\ldots,0)}$ in~\eqref{EqNipHiCommMtx} (where, as before, we write $\hat h=\frac{h}{R}$ for a local defining function of $\sface$) enters the positive commutator argument through its imaginary part, which is $\leq 0$ and thus does not affect the threshold condition. Arguing via induction on $k$ much as in the proof of Proposition~\ref{PropSpHiInftyOut}, we thus obtain a higher b-regularity version of the radial point estimate at $\cR_{\cface,{\rm out}}$.

This finishes the proof of Proposition~\ref{PropNip}\eqref{ItNipHi}.

\begin{rmk}[Second microlocal alternative]
\label{RmkNipHi2nd}
  An alternative to this direct commutator argument is to use second microlocal spaces as introduced in \cite[\S{3.4}]{HintzConicProp} and used in \cite[Theorem~4.13]{HintzConicProp}; this would be analogous to Remark~\ref{RmkSpB2nd}.
\end{rmk}

\subsubsection{Estimates for bounded frequencies; proof of Proposition~\usref{PropNip}\eqref{ItNipB}}
\label{SssNipB}

We now follow the arguments for \cite[Proposition~5.16]{HintzNonstat} and \cite[\S{8.1}]{HintzVasyScrieb} closely. We first show:

\begin{lemma}[Fredholm property]
\label{LemmaNipBFred}
  Let $\gamma_\sscri,\gamma_\cface$ be as in Proposition~\usref{PropNip}. Then for all $\sigma\in\C$ with $\Im\sigma>-\gamma_\sscri$, and for all $s\in\R$, the operator
  \begin{equation}
  \label{EqNipBFred}
    N_{\iota^+}^0(P_0,\sigma) \colon H_{0,\bop}^{s,(2\gamma_\sscri,\gamma_\cface)}(\iota^+) \to H_{0,\bop}^{s-2,(2\gamma_\sscri+2,\gamma_\cface-2)}(\iota^+)
  \end{equation}
  is Fredholm of index $0$. Moreover, if $u\in H_{0,\bop}^{s,(2\gamma_\sscri,\gamma_\cface)}$ satisfies
  \[
    N_{\iota^+}^0(P_0,\sigma)u=f\in H_{0,\bop;\bop}^{(s'-2;k),(2\gamma_\sscri+2,\gamma_\cface-2)}
  \]
  where $s'\geq s$ and $k\in\N_0$, then $u\in H_{0,\bop;\bop}^{(s';k),(2\gamma_\sscri,\gamma_\cface)}$.
\end{lemma}
\begin{proof}
  The ellipticity of $N_{\iota^+}^0(P_0,\sigma)\in\Diff_{0,\bop}^{2,(-2,2)}(\iota^+)$ (discussed in~\S\ref{SssNipPhase}) implies for all $\sigma\in\C$ the estimate
  \begin{equation}
  \label{EqNipBFredPf}
    \|u\|_{H_{0,\bop}^{s,(2\gamma_\sscri,\gamma_\cface)}} \leq C\Bigl(\|N_{\iota^+}^0(P_0,\sigma)u\|_{H_{0,\bop}^{s-2,(2\gamma_\sscri+2,\gamma_\cface-2)}} + \|u\|_{H_{0,\bop}^{-N,(2\gamma_\sscri,\gamma_\cface)}} \Bigr),
  \end{equation}
  where we take $-N<s$. Note next that the normal operator of $R^2 N_{\iota^+}^0(P_0,\sigma)$ at $R=0$ is equal to $P_{(0)}(0,\omega,-R\pa_R,\pa_\omega)$, and thus $\lambda$ is an indicial root for it if and only if $-\lambda$ is an indicial root for $\wh{P_0}(0)$. Since $-\gamma_\cface=\alpha_++\frac32$ lies in the indicial gap $(\beta^-,\beta^+)$ for $\wh{P_0}(0)$, the number $\gamma_\cface$ lies in the indicial gap $(-\beta^+,-\beta^-)$ of $N_{\iota^+}^0(P_0,\sigma)$ at $R=0$. We can thus use a normal operator argument (as in~\eqref{EqMUbEstNearInfty}, but now with b-densities) to estimate, for $\chi\in\CIc([0,3))$ equal to $1$ on $[0,2]$ and for $\tilde\chi\in\CIc([0,4))$ equal to $1$ on $[0,3]$,
  \[
    \|\chi(R)u\|_{H_{0,\bop}^{-N,(2\gamma_\sscri,\gamma_\cface)}} \leq C\Bigl( \| \tilde\chi(R) N_{\iota^+}^0(P_0,\sigma)u \|_{H_{0,\bop}^{-N-2,(2\gamma_\sscri+2,\gamma_\cface-2)}} + \|\tilde\chi(R) u\|_{H_{0,\bop}^{-N,(2\gamma_\sscri,\gamma_\cface-1)}}\Bigr).
  \]
  We plug this into~\eqref{EqNipBFredPf}, and use also that $1-\chi(R)$ has support disjoint from $R^{-1}(0)$, to improve the error term to
  \begin{equation}
  \label{EqNipBFredPf2}
    \|u\|_{H_{0,\bop}^{-N,(2\gamma_\sscri,\gamma_\cface-1)}}.
  \end{equation}

  We next work near $x_\sscri^{-1}(0)$ and use an elliptic parametrix construction in the 0-calculus to weaken the weight at $\scri^+$; we use \cite{Hintz0Px} for this purpose. Using the expression in~\eqref{EqNipOpMT}, the operator $x_\sscri^{-2}N_{\iota^+}^0(P_0,\sigma)$ is given by
  \[
    -\frac12(x_\sscri\pa_{x_\sscri}-2 i\sigma)(x_\sscri\pa_{x_\sscri}-2(1+S|_{\pa X})) + x_\sscri^2\slDelta
  \]
  modulo $x_\sscri\Diff_0^2$. The 0-normal operator is obtained from this by freezing coefficients at the boundary: in the center of geodesic normal coordinates $\omega$ around a point $\omega_0\in\Sph^2$, this gives $-\frac12(x_\sscri\pa_{x_\sscri}-2 i\sigma)(x_\sscri\pa_{x_\sscri}-2(1+S|_{\pa X}))+x_\sscri^2 D_\omega^2$. One then passes to the Fourier transform in $\omega$ and introduces\footnote{We use the notation $t$ for consistency with \cite{Hintz0Px}; this is of course not the same as the time function $t$ on spacetime.} $t:=x_\sscri|\omega|$, $\hat\omega:=\frac{\omega}{|\omega|}$; this yields
  \[
    \hat N = -\frac12(t\pa_t-2 i\sigma)(t\pa_t-2(1+S|_{\omega_0})) + t^2 \in \Diff_{\bop,\scop}^{2,(0,2)}([0,\infty]_t),
  \]
  independently of $\hat\omega$, where the weights $0$ and $2$ are the inverse powers of the defining functions $\frac{t}{1+t}$ and $\frac{1}{1+t}$ of $\{0\}$ and $\{\infty\}$, respectively. The indicial roots at $t=0$ are $2 i\sigma$ and $\spec(2(1+S|_{\omega_0}))$, and hence $2\gamma_\sscri$ lies in the indicial gap $(-2\Im\sigma,2(1+\ubar S))$. Following \cite[Definition~1.2]{Hintz0Px}, with weight $2\gamma_\sscri$, we check that
  \[
    \hat N\colon H_{\bop,\scop}^{s,(2\gamma_\sscri,r)}([0,\infty],|\tfrac{\dd t}{t}|)\to H_{\bop,\scop}^{s-2,(2\gamma_\sscri,r-2)}([0,\infty],|\tfrac{\dd t}{t}|)
  \]
  is invertible; $s,r\in\R$ are arbitrary here by elliptic regularity at fiber infinity and at $t=\infty$ (in the scattering algebra). The Fredholm property follows by elliptic regularity and a normal operator argument at $t=0$. In a suitable basis of the bundle $\cE$, the endomorphism $S|_{\omega_0}$ is lower triangular, and can be continuously deformed to a diagonal matrix without changing the diagonal entries; it thus suffices to show for any one number $\mu\in\spec(1+S|_{\omega_0})$ (so $2\mu>2\gamma_\sscri>-2\Im\sigma=\Re(2 i\sigma)$) that
  \[
    -(t\pa_t-2 i\sigma)(t\pa_t-2\mu) + 2 t^2 \colon H_{\bop,\scop}^{s,(2\gamma_\sscri,r)} \to H_{\bop,\scop}^{s-2,(2\gamma_\sscri,r-2)}
  \]
  is invertible. Multiplying this operator on the left, resp.\ right by $t^{-(i\sigma+\mu)}$, resp.\ $t^{i\sigma+\mu}$, we must prove the invertibility of
  \[
    \hat N' := -(t\pa_t)^2 + (\mu-i\sigma)^2 + 2 t^2 \colon H_{\bop,\scop}^{s,(2\tilde\gamma_\sscri,r')} \to H_{\bop,\scop}^{s-2,(2\tilde\gamma_\sscri,r'-2)}
  \]
  where $2\tilde\gamma_\sscri:=2\gamma_\sscri+\Im\sigma-\mu$ and $r'\in\R$; note that $2\tilde\gamma_\sscri\in(-\mu-\Im\sigma,\mu+\Im\sigma)=(-\Re(\mu-i\sigma),\Re(\mu+i\sigma))$ lies between the two indicial roots (i.e., characteristic exponents) of the normal operator at $t=0$. Since $\mu+\Im\sigma>0$, an element $u$ of the nullspace thus satisfies $|u|,|t\pa_t u|=o(1)$ as $t\to 0$; and $u$ is Schwartz as $t\to\infty$ by ellipticity. If $\mu-i\sigma\in\R$ (which is thus positive), then $0=\la\hat N' u,u\ra_{L^2((0,\infty),|\frac{\dd t}{t}|)}=\|t\pa_t u\|^2+\|(\mu-i\sigma)u\|^2+2\|t u\|^2$ implies $u=0$; otherwise, we conclude $u=0$ from $0=\Im\la\hat N' u,u\ra=\Re(\mu-i\sigma)\Im(\mu-i\sigma)\|u\|^2=0$. The surjectivity of $\hat N'$ is a consequence of the injectivity of its adjoint, which in turn follows from the same arguments with $-\bar\sigma$ in place of $\sigma$.

  It now follows from \cite[Theorem~1.5]{Hintz0Px}, applied on $\tilde\iota:=\iota^+\cap\{x_\sscri<4\}$, that $x_\sscri^{-2}N_{\iota^+}^0(P_0,\sigma)$ admits a left parametrix $Q$ near $x_\sscri=0$, which is the sum
  \[
    Q=Q_0+Q_1
  \]
  of an operator $Q_0\in\Psi_0^{-2}$ and an operator
  \begin{equation}
  \label{EqNipBFredQ1}
    Q_1\in\Psi_0^{-\infty,(\alpha_\lb,-2,\alpha_\rb)}(\tilde\iota;\pi_R^*\Omegab\tilde\iota),
  \end{equation}
  where we use the notation of \cite[\S{2}]{Hintz0Px}, except we work with the pullback of the b-density bundle $\Omegab\tilde\iota$ over $\tilde\iota$ along the right projection $\pi_R\colon\tilde\iota\times\tilde\iota\to\tilde\iota$, and we trivialize the 0-density bundle here using $|\frac{\dd x_\sscri}{x_\sscri}|\,\frac{|\dd\slg|}{x_\sscri^2}$; in~\eqref{EqNipBFredQ1}, $\alpha_\lb=2(1+\ubar S)-\eps>2\gamma_\sscri$ and $\alpha_\rb=2\Im\sigma-\eps>-2\gamma_\sscri$ for any fixed $\eps>0$. (The usage of a right b-density causes a shift of $2=\dim(\tilde\iota)-1$ in the second and third orders of~\eqref{EqNipBFredQ1} relative to the reference.) The operator $Q_0$ is bounded on every 0-Sobolev space, and $Q_1\colon H_{0,\loc}^{-N-2,2\gamma_\sscri}\to H_{0,\loc}^{\infty,2\gamma_\sscri}$; likewise on spaces with additional b-regularity. The local left parametrix property means that $Q x_\sscri^{-2}N_{\iota^+}^0(P_0,\sigma)=I+R$, where $R\in\Psi_0^{-\infty,(\infty,\infty,\alpha_\rb)}$ maps $H_{0,\loc}^{-N-2,\beta}\to H_{0,\loc}^{\infty,\infty}$ for all $\beta>-2\Im\sigma$, which includes values of $\beta$ below $2\gamma_\sscri$; fix $\beta\in(-2\Im\sigma,2\gamma_\sscri)$. Therefore, if $\chi_j\in\CIc([0,j+1))$ equals $1$ on $[0,j]$, then, writing $\chi_j=\chi_j(x_\sscri)$,
  \[
    \chi_1 u = \chi_1\chi_2 u= \chi_1 Q x_\sscri^{-2} N_{\iota^+}^0(P_0,\sigma)(\chi_2 u) + \chi_1 R \chi_2 u;
  \]
  taking $H_0^{-N,2\gamma_\sscri}$-norms yields
  \[
    \|\chi_1 u\|_{H_0^{-N,2\gamma_\sscri}} \leq C\Bigl( \|\chi_2 N_{\iota^+}^0(P_0,\sigma)u\|_{H_0^{-N-2,2\gamma_\sscri+2}} + \|\chi_3 u\|_{H_0^{-N-1,-N'}} + \|\chi_2 u\|_{H_0^{-N,\beta}}\Bigr);
  \]
  here $N'$ is arbitrary, with the second term on the right arising from $[N_{\iota^+}^0(P_0,\sigma),\chi_2]\chi_3$. Plugging this into~\eqref{EqNipBFredPf2}, and using also that $1-\chi_1(x_\sscri)$ has support disjoint from $x_\sscri^{-1}(0)$, we have improved the error term in~\eqref{EqNipBFredPf} to $\|u\|_{H_{0,\bop}^{-N,(\beta,\gamma_\cface-1)}}$. Since now $\beta<2\gamma_\sscri$ and $\gamma_\cface-1<\gamma_\cface$, this is a relatively compact error term.

  Completely analogous arguments for the adjoint $N_{\iota^+}^0(P_0,\sigma)^*$ yield an estimate for the adjoint on dual spaces; this completes the proof of the Fredholm property of~\eqref{EqNipBFred}. The index $0$ property follows from its invertibility for large $\Im\sigma>1$, which was already proved in the previous section.

  Finally, the higher 0-b-regularity statement follows from the ellipticity of $N_{\iota^+}^0(P_0,\sigma)$, and also the additional b-regularity follows away from $\iota^+\cap\scri^+$ from the ellipticity of $N_{\iota^+}^0(P_0,\sigma)$ as a b-differential operator there. Near $\iota^+\cap\scri^+$, we apply the parametrix $Q$ to recover $u$ from $N_{\iota^+}^0(P_0,\sigma)$; we use then that the parametrix $Q$ acts boundedly between 0-Sobolev spaces with additional b-regularity: this is clear for $Q_0$ (cf.\ \eqref{EqMSComm}), and follows for $Q_1$ from the following observation. If $K(x_\sscri,\omega;x_\sscri',\omega')\,|\frac{\dd x_\sscri'}{x_\sscri'}\,\dd\omega'|$ is the Schwartz kernel of $Q_1$ in local coordinates, then the Schwartz kernel of $[x_\sscri\pa_{x_\sscri},Q_1]$ is given by $V K$ where $V=x_\sscri\pa_{x_\sscri}+x'_\sscri\pa_{x'_\sscri}$ is a b-vector field on the 0-double space\footnote{It is, in fact, tangent to the lifted diagonal, and thus preserves also the 0-differential order even when it is not $-\infty$; this is the geometric reason why 0-ps.d.o.s preserve b-regularity.} $[\tilde\iota\times\tilde\iota;\diag_{\pa\tilde\iota}]$, as can be easily checked by a direct computation, and thus preserves the operator class~\eqref{EqNipBFredQ1}; similarly for commutators with $\pa_\omega$, where now $V=\pa_\omega+\pa_{\omega'}$.
\end{proof}

As a corollary of Lemma~\ref{LemmaNipBFred} and the invertibility of~\eqref{EqNipBFred} for large $\Im\sigma$, we conclude that $N_{\iota^+}^0(P_0,\sigma)^{-1}$ is meromorphic in $\Im\sigma>-\gamma_\sscri$. To prove its holomorphicity for $\Im\sigma>-1-\beta^+$, it remains to establish:

\begin{lemma}[Trivial nullspace]
\label{LemmaNipBker}
  Let $\gamma_\sscri,\gamma_\cface$ be as in Proposition~\usref{PropNip}. Then for all $\sigma\in\C$ with $\Im\sigma>-\gamma_\sscri$ and $\Im\sigma>-1-\beta^+$, the kernel of~\eqref{EqNipBFred} is trivial, and hence~\eqref{EqNipBFred} is invertible.
\end{lemma}
\begin{proof}
  Let $u$ be in the kernel of~\eqref{EqNipBFred}; then $u$ is conormal by the final part of Lemma~\ref{LemmaNipBFred}. Working with the smooth structure of $\iota_1^+=[0,\infty]_R\times\Sph^2$ from Definition~\ref{DefCMSpacetime}\eqref{ItCMSpacetimeRad} (i.e., declaring $R^{-1}$, not $R^{-\frac12}$, to be a local defining function of $\{R=\infty\}$), this means
  \[
    u\in\cA^{\gamma_\sscri,\gamma_\cface}(\iota_1^+)=\rho_\sscri^{\gamma_\sscri}\rho_\cK^{\gamma_\cface}\CI_\bop(\iota_1^+),
  \]
  where we can take $\rho_\sscri=\frac{1}{1+R}=\frac{v}{1+v}$ (with $v=\frac{1}{R}=\frac{t_*}{r}$) and $\rho_\cK=\frac{R}{1+R}=\frac{1}{1+v}$ as boundary defining functions.

  We analyze the equation
  \[
    -2 R^{-1}(i\sigma+R\pa_R)(R\pa_R+1+S|_{\pa X}) + R^{-2}P_{(0)}(0,\omega,-R\pa_R,\pa_\omega)u = 0
  \]
  using the ideas from the proof of \cite[Proposition~5.16]{HintzNonstat}, which in turns builds on \cite[Remark~3.27]{Hintz3b}. Thus, we first relate this equation to an equation involving the b-normal operator of $P_0$ at $\sface\subset M_0$, defined as the normal operator of $r^2 P_0=R^2 t_*^2 P_0$ at $\sface$, and indeed its normal operator family, defined as the formal conjugation by the Mellin transform in $r^{-1}$, i.e., considering
  \[
    N^1(P_0,\sigma) := r^{i\sigma} r^2 P_0 r^{-i\sigma} = r^{i\sigma}R^2 t_*^2 P_0 r^{-i\sigma} = R^{i\sigma} R^2 t_*^{i\sigma}t_*^2 P_0 t_*^{-i\sigma} R^{-i\sigma} = R^{i\sigma} R^2 N_{\iota^+}^0(P_0,\sigma) R^{-i\sigma}
  \]
  acting on functions of $v=R^{-1}\in\R$ and $\omega\in\Sph^2$; thus,
  \[
    N^1(P_0,\sigma)(v^{-i\sigma}u) = 0.
  \]
  In the coordinates $v=\frac{t_*}{r}$, $\omega$, we compute
  \[
    N^1(P_0,\sigma) = -2\pa_v(v\pa_v+i\sigma-1-S|_{\pa X}) + P_{(0)}(0,\omega,v\pa_v+i\sigma,\pa_\omega).
  \]
  The $t_*$-translation-invariance of $P_0$ suggests passing to the Fourier transform in $v\in\R$, so we first need to extend $v^{-i\sigma}u$ from $v>0$ to $v\in\R$. Now,
  \begin{equation}
  \label{EqNipBkerCon}
    v^{-i\sigma}u\in\rho_\sscri^{\gamma_\sscri+\Im\sigma}\rho_\cK^{\gamma_\cface-\Im\sigma}\CI_\bop([0,\infty]_v\times\Sph^2),
  \end{equation}
  and $\gamma':=\gamma_\sscri+\Im\sigma>0$. So $v^{-i\sigma}u$ and all of its $v\pa_v$-derivatives are of size $\cO(v^{\gamma'})$ as $v\to 0$. Extending $v^{-i\sigma}u$ by $0$ to $(-\infty,0)_v$ thus yields a function
  \[
    u^\sharp(v,\omega) = \begin{cases} v^{-i\sigma}u, & v>0, \\ 0, & v\leq 0, \end{cases}
  \]
  with the property that $(v\pa_v)^j u^\sharp\in H_\loc^{\frac12+\gamma''}((-1,1)_v\times\Sph^2)$ for all $j\in\N_0$ and $\gamma''<\gamma'$, as follows from elementary bounds on the Fourier transform of $u^\sharp$ (spelled out in~\eqref{EqNipBFT} below); let us fix $\gamma''=\frac12\gamma'>0$. Therefore, $N^1(P_0,\sigma)u^\sharp=0$ everywhere except possibly for a $\delta$-distributional singularity at $v=0$, which however cannot be there since $N^1(P_0,\sigma)u^\sharp\in H_\loc^{-\frac12+\gamma''}$.

  With the convention $(\cF u^\sharp)(\hat r,\omega):=\int_\R e^{i\hat r v}u^\sharp(v,\omega)\,\dd v$, we have $\cF\circ\pa_v=-i\hat r\circ\cF$ and $\cF\circ v=-i\pa_{\hat r}\circ\cF$, so
  \begin{equation}
  \label{EqNipBkerOp}
    \Bigl(2 i\hat r(-\hat r\pa_{\hat r}+i\sigma-2-S|_{\pa X}) + P_{(0)}(0,\omega,-\hat r\pa_{\hat r}-1+i\sigma,\pa_\omega)\Bigr) (\cF u^\sharp) = 0.
  \end{equation}
  But from the conormal bounds~\eqref{EqNipBkerCon}, we deduce that
  \begin{equation}
  \label{EqNipBFT}
  \begin{split}
    w_\pm(\hat r,\omega) &:= (\cF u^\sharp)(\pm\hat r,\omega) \\
      &\in \CIc([0,1)_{\hat r}) + (1+\hat r)^{-1-(\gamma_\sscri+\Im\sigma)+\eps}\Bigl(\frac{\hat r}{1+\hat r}\Bigr)^{-1+\gamma_\cface-\Im\sigma-\eps}\CI_\bop([0,\infty]_{\hat r}\times\Sph^2) \\
      &\subset (1+\hat r)^{-1-\Im\sigma-\gamma_\sscri+\eps} \Bigl(\frac{\hat r}{1+\hat r}\Bigr)^{\min(-1-\Im\sigma+\gamma_\cface,0)-\eps}\CI_\bop([0,\infty]_{\hat r}\times\Sph^2)
  \end{split}
  \end{equation}
  for all $\eps>0$. (See \cite[Lemma~2.25]{Hintz3b} for a proof of this standard result.) We clean up the shifts of $\hat r\pa_{\hat r}$ in~\eqref{EqNipBkerOp} by introducing
  \[
    w_\pm^\flat(\hat r,\omega) := \hat r^{1-i\sigma}w_\pm(\hat r,\omega) \in (1+\hat r)^{-\gamma_\sscri+\eps}\Bigl(\frac{\hat r}{1+\hat r}\Bigr)^{\min(\gamma_\cface,1+\Im\sigma)-\eps}\CI_\bop([0,\infty]_{\hat r}\times\Sph^2),
  \]
  which thus satisfies, upon setting $\hat\rho=\hat r^{-1}$ and multiplying~\eqref{EqNipBkerOp} by $\hat\rho^2$,
  \[
    0 = \Bigl(\pm 2 i\hat\rho(\hat\rho\pa_{\hat\rho}-1-S|_{\pa X}) + \hat\rho^2 P_{(0)}(0,\omega,\hat\rho\pa_{\hat\rho},\pa_\omega)\Bigr) w_\pm^\flat = N_\tface(P_0,\pm 1)w_\pm^\flat;
  \]
  identifying $[0,\infty]_{\hat\rho}\times\Sph^2=\tface$, we have $w_\pm^\flat\in\cA^{\alpha,-\beta}(\tface)$ where $\alpha=\gamma_\sscri-\eps$ and $\beta=-\min(\gamma_\cface,1+\Im\sigma)+\eps\in(\beta^-,\beta^+)$ for small enough $\eps>0$; we use here that $-\gamma_\cface=\alpha_++\frac32$ lies in the indicial gap $(\beta^-,\beta^+)$ of $\wh{P_0}(0)$, and $1+\Im\sigma<\beta^+$. By the $\tface$-admissibility of $P_0$, we thus deduce that $w_\pm^\flat=0$, therefore $w_\pm=0$ and thus $\cF u^\sharp(\hat r,\omega)=0$ for $\hat r\neq 0$. Therefore, $\cF u^\sharp$ is a sum of differentiated $\delta$-distributions at $\hat r=0$, and hence $u^\sharp$ is a polynomial in $v$; but since $\supp u^\sharp\subset\{v\geq 0\}$, this is only possible if $u^\sharp=0$, and thus $u=0$.
\end{proof}

The proof of Proposition~\ref{PropNip} is complete.

\subsubsection{Forward solutions of the \texorpdfstring{$\iota^+$-normal operator}{normal operator at punctured timelike infinity}}
\label{SssNipFw}

We continue writing $\tau=t_*^{-1}$. To capture the structure of $N_{\iota^+}(P_0)$ in Definition~\ref{DefNipOp}, we first introduce:

\begin{definition}[ebe-vector fields]
\label{DefNipebe}
  We denote by $\cV_\ebeop([0,\infty]\times\iota^+)$ the space of all b-vector fields on $[0,\infty]_\tau\times\iota^+$ which, moreover, are tangent to the fibers of the fibrations
  \begin{alignat*}{4}
    &[0,\infty]\times{}&&(\iota^+\cap\scri^+) &&\to \Sph^2,&\quad& (\tau,\omega)\mapsto\omega, \\
    &[0,\infty]\times{}&&(\iota^+\cap\cK^+) &&\to [0,\infty],&\quad& (\tau,\omega)\mapsto\tau.
  \end{alignat*}
  We write ${}^\ebeop T([0,\infty]\times\iota^+)$ and ${}^\ebeop T^*([0,\infty]\times\iota^+)$ for the corresponding tangent and cotangent bundles, respectively.
\end{definition}

By a minor abuse of notation, we shall write $\scri^+$ and $\iota^+$ for the boundary hypersurfaces $[0,\infty]_\tau\times(\iota^+\cap\scri^+)$ and $\{0\}\times\iota^+$ of $[0,\infty]_\tau\times\iota^+$, respectively. We moreover denote by
\[
  \cface := [0,\infty]_\tau \times (\iota^+\cap\cK^+)
\]
the third boundary hypersurface of $[0,\infty]\times\iota^+$. (The fourth boundary hypersurface $\{\infty\}\times\iota^+$ plays no role below.) We denote defining functions of these boundary hypersurfaces by $\rho_+$, $x_\sscri$, and $\rho_\cface$, respectively; possible choices are
\[
  \rho_+=\frac{\tau}{1+\tau},\quad x_\sscri=\sqrt{\frac{v}{1+v}},\quad \rho_\cface=\frac{R}{1+R}.
\]
Other choices are more useful for local computations, such as $\tau$ near $\iota^+$, $\sqrt{v}$ near $\scri^+$, and $R$ near $\cface$. Let $\chi\in\CIc([0,3))$ be equal to $1$ on $[0,2]$. Then $\cV_\ebeop([0,\infty]\times\iota^+)$ is spanned over $\CI$ by the vector fields
\begin{gather}
\label{EqNipebeVF}
  \chi(R) R\pa_R,\ \chi(R)\pa_\omega,\ \chi(R)R\tau\pa_\tau, \\
  \chi(v) v\pa_v=-\chi(R^{-1})R\pa_R,\ \chi(v)\sqrt{v}\pa_\omega=\chi(R^{-1})R^{-\frac12}\pa_\omega,\ \chi(v)\tau\pa_\tau=\chi(R^{-1})\tau\pa_\tau. \nonumber
\end{gather}
By inspection of $\rho_+^{-2}N_{\iota^+}(P_0)$ in~\eqref{EqNipOp}, we thus have
\begin{equation}
\label{EqNipOpebe}
  N_{\iota^+}(P_0) \in x_\sscri^2\rho_+^2\rho_\cface^{-2}\Diff_{\eop\bop\eop}^2([0,\infty]_\tau\times\iota^+).
\end{equation}
Note moreover that since the vector fields~\eqref{EqNipebeVF} are also a local frame for $\Vtb(M_0)$ near $\pa\cK^+$, we have a identification
\begin{equation}
\label{EqNipTebeIdent}
  {}^\ebeop T^*_{\iota^+\cap\cface^+}([0,\infty]\times\iota^+) = {}^\tbop T^*_{\iota^+\cap\cK^+}M_0
\end{equation}

The space $\CI_\bop\cV_\ebeop([0,\infty]\times\iota^+)$ of ebe-vector fields with conormal coefficients was, in fact, already described from a scaled bounded geometry perspective in~\S\ref{SssMUip}: its unit cells are given by~\eqref{EqMSFCells} relative to the spatial (i.e., on $\iota^+$) unit cells~\eqref{EqMUipCells} and the additional weight families~\eqref{EqMUipScaling}. We thus have associated spaces
\[
  \CI_\bop\Psi_\ebeop^{s,(2\gamma_\sscri,\gamma_+,\gamma_\cface)} = x_\sscri^{-2\gamma_\sscri}\rho_+^{-\gamma_+}\rho_\cface^{-\gamma_\cface}\CI_\bop\Psi_\ebeop^s
\]
of ebe-ps.d.o.s with conormal coefficients, and Sobolev spaces
\begin{equation}
\label{EqNipebeSob}
  H_\ebeop^{s,(2\gamma_\sscri,\gamma_+,\gamma_\cface)}([0,\infty]\times\iota^+) = x_\sscri^{2\gamma_\sscri}\rho_+^{\gamma_+}\rho_\cface^{\gamma_\cface}H_\ebeop^s([0,\infty]\times\iota^+);
\end{equation}
we shall always use a positive smooth (unweighted) b-density to define the underlying $L^2$-space. The regularity order $s$ may be variable, while the decay orders $\gamma_\sscri,\gamma_+,\gamma_\cface$ must be constant. We denote spaces encoding additional $k\in\N_0$ degrees of b-regularity (i.e., regularity with respect to $\tau\pa_\tau$, $R\pa_R$, $\pa_\omega$) by
\begin{equation}
\label{EqNipebebSob}
  H_{\ebeop;\bop}^{(s;k),(2\gamma_\sscri,\gamma_+,\gamma_\cface)}.
\end{equation}
Lemma~\ref{LemmaMUetbM}, with purely notational modifications, gives a Plancherel-type description of these spaces, namely
\begin{equation}
\label{EqNipebeEquiv}
\begin{split}
  \|u\|_{H_{\ebeop;\bop}^{(s;k),(2\gamma_\sscri,\gamma_+,\gamma_\cface)}}^2 &\sim \int_{-1}^1 \| (\cM u)(\sigma-i\gamma_+)) \|_{H_{0,\bop;\bop}^{(s;k-j),(2\gamma_\sscri,\gamma_\cface)}(\iota^+)}^2\,\dd\sigma \\
  &\qquad + \sum_\pm\sum_{j=0}^k \int_{\pm[1,\infty)} |\sigma|^j\| (\cM u)(\sigma-i\gamma_+) \|_{H_{(0\semi,\chop,|\sigma|^{-1});\bop}^{(s;k-j),(2\gamma_\sscri,\gamma_\cface,\gamma_\cface,s)}(\iota^+)}^2\,\dd\sigma,
\end{split}
\end{equation}
similarly for variable orders $s\in\CI({}^\ebeop S^*([0,\infty]\times\iota^+))$ that are $\tau$-dilation-invariant, using the induced orders on the right-hand side.

In order to encode causality, we write
\begin{equation}
\label{EqNipebeDot}
  \dot H_\ebeop^s([0,\tau_0]\times\iota^+) := \bigl\{ u\in H_\ebeop^s([0,\infty]\times\iota^+) \colon \supp u\subset\tau^{-1}([0,\tau_0]) \bigr\}
\end{equation}
for $\tau_0>0$, similarly for weighted spaces and in the presence of additional b-regularity. Note that $\dd(-\log\tau)=-\frac{\dd\tau}{\tau}$ is (past) null for the metric~\eqref{EqNipDualMet}, so it is reasonable to expect that forward problems for $N_{\iota^+}(P_0)u=f$ are well-posed on $[0,\tau_0]\times\iota^+$ in function spaces encoding appropriate weights at $\cface$. This expectation is correct, in the following strong sense:

\begin{prop}[Forward solutions of the $\iota^+$-normal operator]
\label{PropNipFwd}
  Let $\sfs\in\CI({}^\ebeop S^*([0,\infty]_\tau\times\iota^+))$ be an order function given by the $\tau$-dilation-invariant extension of the restriction to $\iota^+$ of an admissible order function $\sfs$ with weights $\alpha_+,0$ and margin $0$ (as in~\eqref{EqNipOrder}). Let
  \[
    \gamma_\sscri<1+\ubar S,\ \ \gamma_+<\min(\gamma_\sscri,1+\beta^+),\ \ \gamma_\cface:=-(\alpha_++\tfrac32),
  \]
  where we recall $\beta^+$ from Definition~\usref{DefSStfAdm}. Let $\tau_0>0$ and $k\in\N_0$. Then:
  \begin{enumerate}
  \item\label{ItNipFwd1}{\rm (Existence.)} For all source terms $f\in\dot H_{\ebeop;\bop}^{(\sfs-2;k),(2\gamma_\sscri+2,\gamma_++2,\gamma_\cface-2)}([0,\tau_0]\times\iota^+)$, there exists a unique solution $u\in\dot H_{\ebeop;\bop}^{(\sfs;k),(2\gamma_\sscri,\gamma_+,\gamma_\cface)}([0,\tau_0]\times\iota^+)$ of the equation
    \begin{equation}
    \label{EqNipFwdEq}
      N_{\iota^+}(P_0)u=f.
    \end{equation}
  \item\label{ItNipFwd2}{\rm (Strong uniqueness.)} If $\tilde u$ is a distribution on $(0,\infty)\times(\iota^+)^\circ$ with support in $\tau\leq\tau_0$ that solves $N_{\iota^+}(P_0)\tilde u=0$ and such that, moreover,\footnote{By this notation, we mean that $\chi\tilde u\in H_\ebeop^{\sfs,(2\gamma_\sscri,\gamma_+,\gamma_\cface)}$ for all $\chi\in\CIc((0,\infty)\times\{R<R_0\})$; since such $\chi$ have supports disjoint from $\scri^+$ and $\iota^+$, we drop the orders $2\gamma_\sscri$ and $\gamma_+$ and the subscripts ``$\eop$'' and ``$\bop$'' corresponding to these boundaries.} $\tilde u\in H_{\eop,\loc}^{\sfs,\gamma_\cface}((0,\infty)_\tau\times\{R<R_0\})$ for some $R_0>0$, then $\tilde u=0$.
  \end{enumerate}
\end{prop}
\begin{proof}
  \pfstep{Part~\eqref{ItNipFwd1}.} Equation~\eqref{EqNipFwdEq} is equivalent to $N_{\iota^+}^0(P_0)u=t_*^2 f=\tau^{-2}f$. For the proof of existence, we shall thus define $u$ via its Mellin transform $(\cM u)(\sigma,\cdot)=\int \tau^{-i\sigma}u(\tau,\cdot)\,\frac{\dd\tau}{\tau}$ by
  \[
    (\cM u)(\sigma) := N_{\iota^+}^0(P_0,\sigma)^{-1}\bigl( \cM(\tau^{-2}f) \bigr),\quad \Im\sigma=-\gamma_+;
  \]
  this is well-defined by Proposition~\ref{PropNip}. The estimate~\eqref{EqNipb}, phrased as
  \begin{align*}
    \|(\cM u)(\sigma)\|_{H_{(0\semi,\chop,|\sigma|^{-1});\bop^+}^{(\sfs;k),(2\gamma_\sscri,\gamma_\cface,\gamma_\cface,\sfb)}} &\leq C_k|\sigma|^{-2}\|N_{\iota^+}^0(P_0,\sigma)(\cM u)(\sigma)\|_{H_{(0\semi,\chop,|\sigma|^{-1});\bop^+}^{(\sfs-2;k),(2\gamma_\sscri+2,\gamma_\cface-2,\gamma_\cface,\sfb+1)}} \\
      &= C_k\|\cM(\tau^{-2}f)(\sigma)\|_{H_{(0\semi,\chop,|\sigma|^{-1});\bop^+}^{(\sfs-2;k),(2\gamma_\sscri+2,\gamma_\cface-2,\gamma_\cface-2,\sfb-1)}} \\
      &\leq C'_k\|\cM(\tau^{-2}f)(\sigma)\|_{H_{(0\semi,\chop,|\sigma|^{-1});\bop^+}^{(\sfs-1;k),(2\gamma_\sscri+2,\gamma_\cface-2,\gamma_\cface-2,\sfb-1)}},
  \end{align*}
  and the uniform bounds on the inverse of~\eqref{EqNip} on compact subsets of $\{\sigma\in\C\colon\Im\sigma=-\gamma_+\}$, imply, using the norm equivalence~\eqref{EqNipebeEquiv}, that
  \begin{equation}
  \label{EqNipFwdBd}
    \|u\|_{H_\ebeop^{(\sfs;k),(2\gamma_\sscri,\gamma_+,\gamma_\cface)}([0,\infty]\times\iota^+)} \leq C\|\tau^{-2}f\|_{H_\ebeop^{(\sfs-1;k),(2\gamma_\sscri,\gamma_++2,\gamma_\cface-2)}([0,\infty]\times\iota^+)}.
  \end{equation}

  We next use the Paley--Wiener theorem to prove that $u$ vanishes for $\tau>\tau_0$. In view of the bound~\eqref{EqNipFwdBd}, it suffices to prove this for a dense set of source terms $f$, and thus for $f\in\CIc((0,\tau_0)\times(\iota^+)^\circ)$; if $\supp f\subset[0,\tau_-]\times K$ where $\tau_-<\tau_0$ and $K\subset(\iota^+)^\circ$ is compact, then $(\cM f)(\sigma)\in\CIc((\iota^+)^\circ)$ is supported in $K$ as well. For such $f$, $\cM(\tau^{-2}f)(\sigma,\cdot)$ is an entire $\CIc((\iota^+)^\circ)$-valued function. Let $\phi\in\CIc((\iota^+)^\circ)$; then $\la(\cM u)(\sigma,\cdot),\phi\ra_{L^2(\iota^+)}$ is holomorphic for $\Im\sigma>-\min(\gamma_\sscri,1+\beta^+)$. Now, the estimate~\eqref{EqNipb} then implies
  \[
    |\la (\cM u)(\sigma,\cdot),\phi\ra_{L^2(\iota^+)}| \leq C\la\sigma\ra^N \|\cM(\tau^{-2}f)(\sigma,\cdot)\|_{H^N}\quad \forall\,\sigma\in\C,\ \Im\sigma\geq-\gamma_+,
  \]
  for some $N$, where we crudely estimate $H_{0\semi,\chop}$-norms by standard Sobolev norms times sufficiently high powers of $\la\sigma\ra$. But since $\tau\leq\tau_-<\tau_0$ on $\supp f$, we have
  \[
    \|\cM(\tau^{-2}f)(\sigma,\cdot)\|_{H^N} = \la\sigma\ra^{-2 L} \bigl\|\cM\bigl[(1+(\tau D_\tau)^2)^L\tau^{-2}f\bigr](\sigma,\cdot)\bigr\|_{H^N} \leq C_L\la\sigma\ra^{-2 L} \tau_-^{\Im\sigma}
  \]
  for all $L\in\N_0$; we fix $L$ such that $N-2 L<-1$. Write
  \[
    \la u(\tau,\cdot),\phi\ra_{L^2(\iota^+)} = \frac{1}{2\pi} \int_{\Im\sigma=-\gamma} \tau^{i\sigma} \la(\cM u)(\sigma,\cdot),\phi\ra_{L^2(\iota^+)}\,\dd\sigma,
  \]
  initially for $\gamma=\gamma_+$. If $\tau\geq\tau_0$, we shift the contour by increasing $\gamma$; the integrand is then bounded in absolute value by $\tau^{-\Im\sigma}\la\sigma\ra^{N-2 L}\tau_-^{\Im\sigma}$, the integral of which, for $\Im\sigma=\gamma$ is bounded by $C(\frac{\tau_-}{\tau_0})^\gamma$, which tends to $0$ as $\gamma\to\infty$. Since $\phi$ was arbitrary, this shows that $\supp u\subset\tau^{-1}([0,\tau_0])$, as desired.

  \pfstep{Part~\eqref{ItNipFwd2}.} We first note that the arguments thus far imply the uniqueness of $u$ in the class $H_\ebeop^{\sfs,(2\gamma_\sscri,\gamma_+,\gamma_\cface)}([0,\infty]\times\iota^+)$.

  One way to prove the stronger uniqueness claim is to use \cite[Theorem~3.18]{HintzConicWave}, the main task being the verification of its hypotheses, specifically, the spectral admissibility property of $N_{\iota^+}(P_0)$ in the sense of \cite[Definition~3.12]{HintzConicWave}. We present an alternative and self-contained, albeit rather \emph{ad hoc}, argument. Geometrically speaking, we are presently working in a slicing of the ``spacetime'' $(0,\infty)_\tau\times(\iota^+)^\circ$ by \emph{null} hypersurfaces $\{\tau={\rm const.}\}=\{t_*={\rm const.}\}$. Arguments based on finite speed of propagation and Paley--Wiener arguments, as used in \cite{HintzConicWave}, are thus not immediately available since supports can spread with infinite speed on $t_*$-level sets and immediately (as measured by $t_*$) reach $\scri^+$. The idea is thus to exploit the energy estimates near $\scri^+$ in $t_*$-slabs established in Proposition~\ref{PropETSolv} to extend $\tilde u$ to a full $t_*$-slab, and then use the uniqueness property for $N_{\iota^+}(P_0)$ noted above.

  The details are as follows. Writing $t_*=\tau^{-1}$, we may shift $t_*$ to arrange that $t_*\geq\frac{4}{R_0}$ on the support of the distribution $\tilde u\in H_{\eop,\loc}^{\sfs,\gamma_\cface}((0,\infty)_\tau\times\{R<R_0\})$ satisfying $N_{\iota^+}(P_0)\tilde u=0$. We now pass to the coordinates $t_*$ and $r=R t_*$; thus $\tilde u$ is defined when $r=R t_*<4$. Shifting $t_*$ again, we may now assume that $\tilde u$ vanishes for $t_*\leq 1$; we shall prove its vanishing for $t_*\leq 2$, which by iteration implies $\tilde u=0$ globally. For $j=1,2,3$, let now $\chi_{t_*}\in\CI(\R)$ be equal to $1$ near $(-\infty,2]$ and equal to $0$ near $[3,\infty)$, and let $\chi_r\in\CIc([0,3))$ be equal to $1$ near $[0,2]$. We then have
  \[
    N_{\iota^+}(P_0)(\chi_{t_*}\chi_r \tilde u) = [N_{\iota^+}(P_0),\chi_{t_*}](\chi_r \tilde u) + \chi_{t_*}[N_{\iota^+}(P_0),\chi_r]\tilde u \in H_\ebeop^{\sfs-1,(\infty,\infty,\gamma_\cface-2)}.
  \]
  The second summand is supported in the set $\{1\leq t_*<3,\ 2<r<3\}$. We may then use standard hyperbolic theory to find a solution $u_{\rm ext}\in H^\sfs_\loc(\{t_*<3,\ 1<r\})$ of the wave equation
  \[
    N_{\iota^+}(P_0)u_{\rm ext} = \chi_{t_*}[N_{\iota^+}(P_0),\chi_r]\tilde u,
  \]
  with $u_{\rm ext}$ vanishing for $r<2-(t-3)$ and also for $t_*\leq 1$ by finite speed of propagation, where we recall the Minkowski time coordinate $t=t_*+r$, i.e., its support is contained in $r>1$. Proposition~\ref{PropETSolv} (the proof of which applies without change to the operator $N_{\iota^+}(P_0)$) in fact gives uniform control of $u_{\rm ext}$ as $r\to\infty$; thus
  \[
    \chi_{t_*}u_{\rm ext} \in H_\ebeop^{\sfs,(2\gamma_\sscri,\infty,\infty)}.
  \]
  (Here $\alpha_\sscri$ in~\eqref{EqETSolvEst} and $\gamma_\sscri$ are related via $\gamma_\sscri=\alpha_\sscri+\frac32$ since we presently work with unweighted b-densities.) See Figure~\ref{FigNipFwdUniq}.

  \begin{figure}[!ht]
  \centering
  \includegraphics{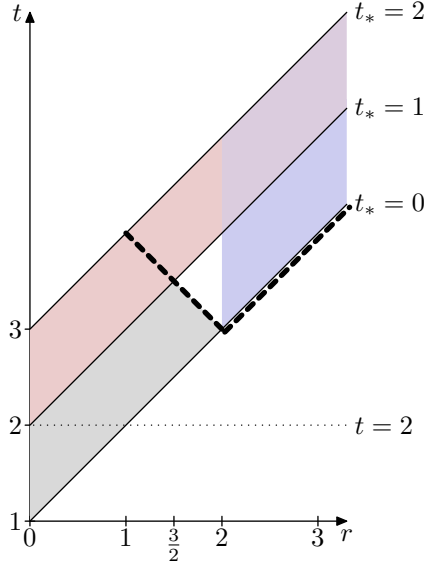}
  \caption{Setup for the uniqueness argument. The cut-off wave $\chi_{t_*}\chi_r\tilde u$ is generated by a source term supported in the union of the red and blue regions (their overlap being purple). The support of $u_{\rm ext}$, which is sourced by a term supported in the blue and purple regions, lies above the dashed lines. We first show that $\tilde u$ vanishes in the gray region; it then vanishes for all $t_*<2$ since it has trivial Cauchy data at $t=2$ (dotted).}
  \label{FigNipFwdUniq}
  \end{figure}

  Let $\chi_{t_*}^\flat\in\CI(\R)$ be equal to $1$ on $(-\infty,2]$ and supported in $\{\chi_{t_*}=1\}$. For the function
  \[
    w := \chi_{t_*}\chi_r \tilde u - \chi_{t_*}^\flat u_{\rm ext} \in \dot H_\ebeop^{\sfs,(2\gamma_\sscri,\infty,\gamma_\cface)}([0,1]_\tau\times\iota^+),
  \]
  we then compute
  \[
    N_{\iota^+}(P_0)w = [N_{\iota^+}(P_0),\chi_{t_*}](\chi_r \tilde u) + \chi_{t_*}(1-\chi_{t_*}^\flat)[N_{\iota^+}(P_0),\chi_r]\tilde u - [N_{\iota^+}(P_0),\chi_{t_*}^\flat](\chi_{t_*}u_{\rm ext}) =: f.
  \]
  We have $t_*\geq 2$ (so $\tau\leq\frac12$) on the support of each summand of $f$, and therefore
  \[
    f\in\dot H_\ebeop^{\sfs-1,(2\gamma_\sscri+2,\infty,\gamma_\cface-2)}([0,\tfrac12]_\tau\times\iota^+).
  \]
  By the first part of the proof, we can therefore solve $N_{\iota^+}(P_0)w'=f$ with $w'\in\dot H_\ebeop^{\sfs,(2\gamma_\sscri,\gamma_+,\gamma_\cface)}([0,\tfrac12]_\tau\times\iota^+)$. But then $N_{\iota^+}(P_0)(w-w')=0$, and as noted above, uniqueness of solutions of this equation in the space $\dot H_\ebeop^{\sfs,(2\gamma_\sscri,\gamma_+,\gamma_\cface)}([0,1]_\tau\times\iota^+)\ni w-w'$ now implies that $w-w'=0$, so $w=w'$ is supported in $\{\tau\leq\frac12\}$. By the support properties of $u_{\rm ext}$, this implies that $\tilde u$ vanishes for $t_*\leq 2$ when $r>\frac32$. Having shown the vanishing of $\tilde u|_{\{t_*<2\}}$ near $r=0$, the vanishing in all of $\{t_*<2\}$ now follows from standard uniqueness for solutions of wave equations, as $\tilde u$ satisfies the homogeneous equation $N_{\iota^+}(P_0)\tilde u=0$ with trivial Cauchy data on $\{t=2,\ r>0\}$.
\end{proof}

\subsection{Modified wave operator}
\label{SsNMod}

The modification of $P_0$ which we will study next interpolates between $P_0$ near $\scri^+$ and $N_{\iota^+}(P_0)$ (Definition~\ref{DefNipOp}) near $R=\frac{r}{t_*}=0$. Concretely:

\begin{definition}[Modified wave operator]
\label{DefNModOp}
  Let $\chi\in\CIc([0,2))$, $0\leq\chi\leq 1$, be equal to $1$ near $[0,1]$; abbreviate $\chi=\chi(R)$. On the manifold $\{t_*>1\}\times\iota^+$, we then define the operator
  \[
    P_+ := \chi N_{\iota^+}(P_0) + (1-\chi) P_0.
  \]
\end{definition}

We use the notation of~\S\ref{SssNipFw}. In view of the edge-b-nature of $P_0$ near $\scri^+\subset M$ (cf.\ Lemma~\ref{LemmaSSOpMem}) and the b-edge-nature of $N_{\iota^+}(P_0)$ near $\cface\subset[0,\infty]\times\iota^+$ (cf.\ \eqref{EqNipOpebe}), we have
\[
  P_+ \in x_\sscri^2\tau^2\rho_\cface^{-2}\Diff_\ebeop^2([0,1)_\tau\times\iota^+).
\]
As the defining function of $\{0\}\times\iota^+$, we can take $\tau=\frac{1}{t_*}$ here. Since $t_*^2 N_{\iota^+}(P_0)$ is the b-normal operator of $t_*^2 P_0$ at $\iota^+$, we have
\begin{equation}
\label{EqNModNormOp}
  P_+ - N_{\iota^+}(P_0) \in x_\sscri^2\tau^3\rho_\cface^{-2}\Diff_\ebeop^2,
\end{equation}
i.e., $N_{\iota^+}(P_0)$ is the normal operator of $P_+$ at $\tau=0$. Unlike, $N_{\iota^+}(P_0)$, however, the operator $P_+$ is \emph{not} homogeneous with respect to $\tau$-dilations (since the Kerr metric itself is not).

Furthermore, $P_+$ is principally scalar, with principal symbol equal to the dual metric function of $g_\chi$ where
\[
  g_\chi^{-1} := \chi\Bigl(\frac{r}{t_*}\Bigr)\ubar g^{-1}+\Bigl(1-\chi\Bigl(\frac{r}{t_*}\Bigr)\Bigr)g_{\bhm,a}^{-1}.
\]
Since $\dd t_*$ is null for $\ubar g$ and (by Lemma~\ref{LemmaTsKLMetric}) timelike for $g_{\bhm,a}$, we have $g_\chi^{-1}(\dd t_*,\dd t_*)\leq 0$, so $\dd t_*$ is past causal for $g_\chi$. We may thus study forward problems for $P_+$ on domains $[0,\tau_0]_\tau \times \iota^+$ where $\tau_0\in(0,\tfrac12]$. We use the edge-b-edge-Sobolev spaces introduced in~\eqref{EqNipebeSob}--\eqref{EqNipebebSob} and \eqref{EqNipebeDot}. (We remark that Baskin--Marzuola \cite{BaskinMarzuolaCone} study time-translation-invariant wave operators on \emph{exact} cones; this thus bears some similarities with $P_+$. However, we do not use b-Sobolev spaces near the timelike curve of cone points here, as edge Sobolev spaces are less precise and thus more flexible.)

\begin{prop}[Forward solutions of $P_+$]
\label{PropNMod}
  Let $\sfs\in\CI({}^\etbop S^*M)$ be an admissible order function with weights $\alpha_+,0$ and margin $1$. Define an order function $\sfs'\in\CI({}^\ebeop S^*([0,1)\times\iota^+))$ by $\sfs'=\chi(R)\sfs_{\rm I}+(1-\chi(R))\sfs$, where $\sfs_{\rm I}\in\CI({}^\ebeop S^*([0,\infty]\times\iota^+))$ is the $\tau$-dilation-invariant extension of $\sfs|_{{}^\etbop S_{\iota^+}M}$. Then there exists $\tau_0\in(0,\frac12]$ such that, for all
  \[
    \gamma_\sscri<1+\ubar S,\quad \gamma_+<\min(\gamma_\sscri,1+\beta^+),\quad \gamma_\cface:=-(\alpha_++\tfrac32),
  \]
  and for all $k\in\N_0$, the following statements hold.
  \begin{enumerate}
  \item{\rm (Existence.)} For all $f\in\dot H_{\ebeop;\bop}^{(\sfs'-1;k),(2\gamma_\sscri+2,\gamma_++2,\gamma_\cface-2)}([0,\tau_0]\times\iota^+)$, there exists a forward solution
    \begin{equation}
    \label{EqNModSol}
      u\in\dot H_{\ebeop;\bop}^{(\sfs';k),(2\gamma_\sscri,\gamma_+,\gamma_\cface)}([0,\tau_0]\times\iota^+)
    \end{equation}
    of the equation $P_+ u=f$.
  \item\label{ItNModUniq}{\rm (Strong uniqueness.)} If $\tilde u$ is a distribution on $(0,\infty)\times(\iota^+)^\circ$ with support in $\tau\leq\tau_0$ that solves $P_+\tilde u=0$ and satisfies $\tilde u\in H_{\eop,\loc}^{\sfs',\gamma_\cface}((0,\infty)_\tau\times\{R<R_0\})$ for some $R_0>0$, then $\tilde u=0$.
  \end{enumerate}
\end{prop}
\begin{proof}
  By Proposition~\ref{PropDyO}\eqref{ItDyO3}, the order functions $\sfs'$ and $\sfs$ agree near $\iota^+\cap\supp\dd\chi$. The monotonicity of $\sfs$ along the null-bicharacteristic flow over $\iota^+$ then implies the same property also for $\sfs'$ on $\{\chi=1\}$ by the dilation-invariance of $P_+$ in this region. On the transition region $\supp\chi$, where $\sfs=\sfs'$ for small $\tau=t_*^{-1}$, the metric $g_\chi$ is a small perturbation of $g_{\bhm,a}$ (indeed, differs from it by a b-metric of class $\tau^3\CI$ since $g_{\bhm,a}$ and $\ubar g$ agree to leading order on $\iota^+$, with $\tau^3$ being one order more than the $\tau^2$ weight of the Kerr metric as a b-metric near $(\iota^+)^\circ$, cf.\ Lemma~\ref{LemmaTsKLMetric}\eqref{ItTsKLMetrice3b}) in the sense of Proposition~\ref{PropDyO}\eqref{ItDyO2}, with $\ell_+=\frac12$, say, provided we localize to a sufficiently small neighborhood of $\iota^+$; thus, the monotonicity of $\sfs=\sfs'$ holds also for $g_\chi$. We can thus pick $\tau_0$ such that all monotonicity and threshold conditions for propagation and radial point estimates hold for $\sfs'$ on $[0,\tau_0]\times\iota^+$.

  \pfstep{Finite-time solvability.} We first show how to solve $P_+ u=f$ in the desired function spaces on finite time intervals. For better readability, we pass to $t_*=\tau^{-1}$ and let $t_{*,0}=\tau_0^{-1}\geq 2$. We first define $u_\cface$ as the forward solution of
  \[
    N_{\iota^+}(P_0)u_\cface = \chi f,
  \]
  which by Proposition~\ref{PropNipFwd} satisfies
  \[
    \chi u_\cface\in\dot H_{\ebeop;\bop}^{(\sfs_{\rm I};k),(\infty,\gamma_+,\gamma_\cface)}([0,\tau_0]\times\iota^+).
  \]
  Let $\chi_\flat=\chi_\flat(R)$ be equal to $1$ near $[0,1]$ and supported in $\chi^{-1}(1)$. Then
  \begin{align*}
    f' &:= f - P_+(\chi_\flat u_\cface) = f-\chi_\flat P_+ u_\cface - [P_+,\chi_\flat]u_\cface \\
      &= (1-\chi_\flat)f - [P_+,\chi_\flat](\chi u_\cface) \in \dot H_{\ebeop;\bop}^{(\sfs'-1;k),(2\gamma_\sscri+2,\gamma_++2,\infty)}([0,\tau_0]\times\iota^+),
  \end{align*}
  where we use that $\chi_\flat P_+=\chi_\flat N_{\iota^+}(P_0)$; and $R\geq 1$ on $\supp f'$, so $r=R t_*\geq t_{*,0}$. Since the metric underlying $P_+$ in $R\leq 1$ is the Minkowski metric, we may then solve the wave equation
  \[
    P_+ u_{\rm mid}=f'
  \]
  for the $t_*$-interval $[t_{*,0},2 t_{*,0}]$ without encountering the singularity at $R=0$. (Indeed, since $\supp f'$ lies in the future of the point $(t,r)=(t_*+r,r)=(2 t_{*,0},t_{*,0})$, the solution $u_{\rm mid}$ vanishes for $t+r\geq 3 t_{*,0}$, and thus for $t_*=t-r\in[t_{*,0},2 t_{*,0}]$ for $r\geq\frac12 t_{*,0}$.) Proposition~\ref{PropETSolvb} moreover gives control of $u_{\rm mid}$ at $\scri^+\cap t_*^{-1}([t_{*,0},2 t_{*,0}])$; so if $\chi_{t_*}=\chi_{t_*}(t_*)\in\CI(\R)$ equals $1$ for $t_*\leq\frac32 t_{*,0}$ and $0$ for $t_*\geq 2 t_{*,0}$, then
  \[
    \chi_{t_*}u_{\rm mid}\in\dot H_{\ebeop;\bop}^{(\sfs';k),(2\gamma_\sscri,\infty,\infty)}([0,\tau_0]\times\iota^+).
  \]
  Since $P_+(\chi_{t_*}u_{\rm mid})=f'$ for $t_*\in[t_{*,0},\frac32 t_{*,0}]$, we conclude that
  \[
    f'' := f' - P_+(\chi_{t_*}u_{\rm mid}) \in \dot H_{\ebeop;\bop}^{(\sfs'-1;k),(2\gamma_\sscri+2,\gamma_++2,\gamma_\cface-2)}([0,\tfrac23\tau_0]\times\iota^+).
  \]
  In other words, $\chi_\flat u_\cface+\chi_{t_*}u_{\rm mid}$ is a local forward solution of $P_+ u=f$ for $t_*\in[t_{*,0},\frac32 t_{*,0}]$. Iterating this argument $N$ times, and using that $(\frac32)^N\to\infty$ as $N\to\infty$, yields a forward solution $u$ on $[\tau_-,\tau_0]\times\iota^+$ for any fixed $\tau_->0$ (though at this point we do not get quantitative control as $\tau_-\to 0$ yet), which we record as
  \begin{equation}
  \label{EqNModExistLoc}
    u \in \dot H_{\ebeop;\bop,\loc}^{(\sfs';k),(2\gamma_\sscri+2,\gamma_+,\gamma_\cface)}((0,\tau_0]\times\iota^+).
  \end{equation}

  \pfstep{Strong uniqueness.} As in the proof of Proposition~\ref{PropNipFwd}, the only non-standard part of the argument is local uniqueness near the timelike curve $R=0$ of conic singularities. The argument is almost the same, and hence we shall be brief: in the notation used in the proof of Proposition~\ref{PropNipFwd}, we have $N_{\iota^+}(P_0)\chi_r=P_+\chi_r$ for $t_*\geq t_{*,0}$, provided $\{t_*\geq t_{*,0},\ r\leq 3\}\subset\{R<1\}\subset\chi^{-1}(1)$, which holds for all $t_{*,0}>3$. Thus
  \[
    N_{\iota^+}(P_0)(\chi_r\tilde u) = P_+(\chi_r\tilde u) = [P_+,\chi_r]\tilde u.
  \]
  We can then solve $N_{\iota^+}(P_0)u_{\rm mid}=-[P_+,\chi_r]\tilde u$ (with the right-hand side supported in $r^{-1}((2,3))$) in the $t_*$-interval $[t_{*,0},t_{*,0}+1]$, with $u_{\rm mid}$ vanishing near $R=0$ and controlled near $\scri^+$ by Proposition~\ref{PropETSolvb}. But then $N_{\iota^+}(P_0)(\chi_r\tilde u+u_{\rm mid})=0$ for $t_*\leq t_{*,0}+1$. If $\chi_{t_*}\in\CI(\R)$ equals $1$ for $t_*\leq t_{*,0}+\frac12$ and vanishes for $t_*\geq t_{*,0}+1$, then $N_{\iota^+}(P_0)(\chi_{t_*}(\chi_r\tilde u+u_{\rm mid}))$ is supported in $\{t_*\geq t_{*,0}+\frac12\}$; therefore, by the uniqueness part of Proposition~\ref{PropNipFwd}, $\chi_{t_*}(\chi_r\tilde u+u_{\rm mid})$ itself must be supported in this region. Since for $t_*\in[t_{*,0},t_{*,0}+1]$, we have $u_{\rm mid}=0$ for $r\leq 2-(t-(t_{*,0}+2))=4+t_{*,0}-t=4+t_{*,0}-(t_*+r)$, i.e., $r\leq 2-\frac12(t_*-t_{*,0})$, by the finite speed of propagation, we conclude that also $\tilde u$ vanishes there. One then concludes $\tilde u=0$ in the whole slab $t_*^{-1}([t_{*,0},t_{*,0}+1])$ by standard uniqueness for $P_+\tilde u=0$, given the vanishing of $\tilde u$ for $t_*\leq t_{*,0}$ and $r\leq 4+t_{*,0}-t$.

  Having thus shown the vanishing of $\tilde u$ for $t_*\leq t_{*,0}+1$, one can iterate this argument to conclude that $\tilde u=0$ everywhere.

  \pfstep{Quantitative bounds near $\{0\}\times\iota^+$: a priori estimate} We next upgrade~\eqref{EqNModExistLoc} to~\eqref{EqNModSol}. The ``standard'' functional analytic approach would be to establish this first for $k=0$ as a consequence of estimates for the backwards problem for the adjoint $P_+^*$ on dual spaces, then using the Hahn--Banach theorem to produce a solution, and then using a commutation argument to prove higher b-regularity. We leave it to the interested reader to carry this out. We shall use a different argument that avoids the need for dualization and that is, instead, based on an approximation argument which was already sketched after~\eqref{EqIP0Est}. (The same type of argument will be used in~\S\ref{SF} below, where a dualization approach is not available.)

  Let $\sfs_\flat$ be an admissible order function with weights $\alpha_+,0$ and margin $0$ satisfying $\sfs_\flat\leq\sfs-1$, and define $\sfs'_\flat$ relative to $\sfs_\flat$ exactly like $\sfs'$ is defined relative to $\sfs$. For $j=0,1$, let $\chi_j\in\CI(\R)$ be equal to $0$ on $(-\infty,t_{*,0}+j]$ and equal to $1$ on $[t_{*,0}+j+1,\infty)$. We shall first establish the symbolic a priori estimate
  \begin{equation}
  \label{EqNModSymb}
    \|\chi_2 u\|_{H_{\ebeop;\bop}^{(\sfs';k),(2\gamma_\sscri,\gamma_+,\gamma_\cface)}} \leq C\Bigl( \|\chi_1 P_+ u\|_{H_{\ebeop;\bop}^{(\sfs'-1;k),(2\gamma_\sscri+2,\gamma_++2,\gamma_\cface-2)}} + \|\chi_1 u\|_{H_{\ebeop;\bop}^{(\sfs'_\flat;k),(2\gamma_\sscri,\gamma_+,\gamma_\cface)}} \Bigr)
  \end{equation}
  via microlocal elliptic, real principal type propagation, and radial point estimates in the ebe-phase space ${}^\ebeop T^*([0,\infty]\times\iota^+)$. In view of the already established finite-time existence result, it suffices to prove this for arbitrarily high but fixed values of $t_{*,0}$.

  Note that $P_+$ inherits the radial sets $\cR_{\scri^+,{\rm in},+}$ and $\cR_{\scri^+,{\rm out}}$ (Definition~\ref{DefTsebRad}) from $P_0$. Furthermore, using the identification~\eqref{EqNipTebeIdent}, the radial sets $\cR_{\pa\cK^+,{\rm in/out}}$ (Definition~\ref{DefTs3bRad}) can be viewed as conic subsets of ${}^\ebeop T^*([0,\infty]\times\iota^+)$ lying over $\{\tau=R=0\}$; and since we extend $N_{\iota^+}(P)$ by dilation-invariance in $\tau$ for $R<1$ to obtain $P_+$, the $\tau$-invariant extensions
  \[
    \cR_{\cface,{\rm in/out}} \subset {}^\ebeop T^*_{R^{-1}(0)}([0,\infty]\times\iota^+)\setminus o
  \]
  of $\cR_{\pa\cK^+,{\rm in/out}}$ are radial sets for $P_+$. (These are the analogues of the incoming and outgoing radial sets $\cR_{\rm in/out}^\pm$ in \cite[\S{2.1}]{HintzConicWave}.) Only the dynamics near $\cR_{\pa\cK^+,{\rm in/out}}$ over $\iota^+$ get inherited by $\cR_{\cface,{\rm in/out}}$, so $\cR_{\cface,{\rm in}}$ is a normal saddle for the flow, with stable manifold near $R=0$ given by radially infalling (towards $R=0$) null-geodesics lifted to phase space, and unstable manifold lying over $R=0$, which is in turn the stable manifold of $\cR_{\cface,{\rm out}}$ whose unstable manifold is given by the lifts of radially (from $R=0$) outgoing null-geodesics.

  The null-bicharacteristic flow for $P_+$ over $\{\infty\}\times\iota^+$ is described by Proposition~\ref{PropTse3bDyn}. Recall moreover that $t_*$, having past causal differential with respect to the metric $g_\chi$, is monotonically non-decreasing along the future null-bicharacteristic flow. It follows that in order to prove~\eqref{EqNModSymb}, it suffices to combine microlocal propagation results in the following order:
  \begin{enumerate}
  \item real principal type propagation along $\scri^+$ to control $u$ on the stable manifold of $\cR_{\scri^+,{\rm in},+}$ (see Lemma~\ref{LemmaTsebDyn}\eqref{ItTsebDynMfd});
  \item a radial (saddle) point estimate to propagate control on $u$ through $\cR_{\scri^+,{\rm in},+}$ (Proposition~\ref{PropR3RScriI});
  \item a radial (saddle) point estimate to propagate control through $\cR_{\cface,{\rm in}}$, with control on the stable manifold arising from propagation from $t_*<t_{*,0}$ or, over $\iota^+$, from the unstable manifold of $\cR_{\scri^+,{\rm in},+}$;
  \item a radial (saddle) point estimate to propagate through $\cR_{\cface,{\rm out}}$;
  \item a radial (sink) estimate to propagate into $\cR_{\scri^+,{\rm out}}$.
  \end{enumerate}
  See Figure~\ref{FigNModFlow} for an illustration. Only the radial point estimates near $\cR_{\cface,{\rm in/out}}$ have not yet been discussed. But in fact they are minimal modifications of Propositions~\ref{PropR3RKI} and \ref{PropR3RKO}: the only difference is that we no longer need to localize near $\iota^+\subset M$, but only near $\tau=0$, so the localizer $\psi_+(\digamma\rho_+)$ in~\eqref{EqR3RKIComm} must be replaced by a cutoff $\psi_+(\tau)$ to a neighborhood of $\tau=0$, with $\psi_+=1$ near $0$ and $\sqrt{-\psi_+\psi'_+}\in\CI$. It is slightly more convenient to use not $\tau$ but rather a strictly timelike function $\tilde\tau$ whose level sets interpolate between $t$-level sets near $R=0$ and $t_*$-level sets for $R>2$. In the commutator argument with $\psi_+(\tilde\tau)$ in place of $\psi_+(\digamma\rho_+)$ in~\eqref{EqR3RKIComm}, then, the Hamiltonian vector field $H_{G_\chi}$ of $P_+$ (with $G_\chi$ denoting the dual metric function of $g_\chi$) falling on $\psi_+^2$ produces $2\psi_+\psi'_+ H_{G_\chi}\tilde\tau$, which on the future characteristic set is a non-negative square and thus gives rise to an a priori control term---which, however, vanishes when controlling $u$ that vanish on $\supp\dd\psi_+$.

  \begin{figure}[!ht]
  \centering
  \includegraphics{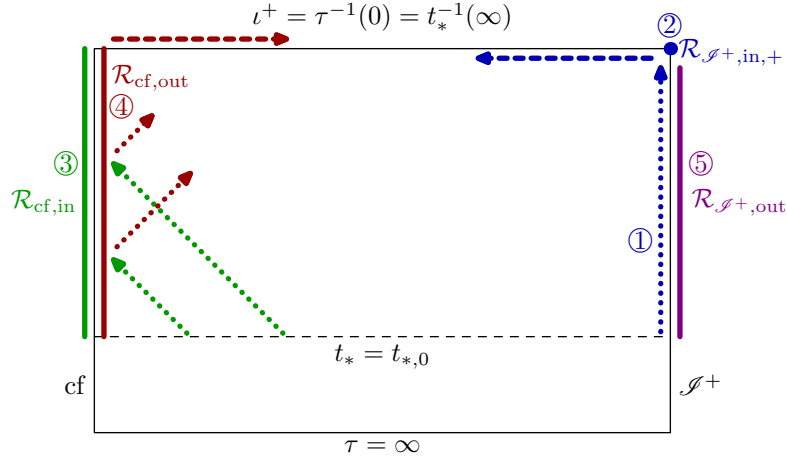}
  \caption{Illustration of the null-bicharacteristic flow of $P_+$ in edge-b-edge-phase space. The numbers indicate the order of propagation in the proof of~\eqref{EqNModSymb}. The dotted, resp.\ dashed arrows indicate the direction of real principal type propagation preceding, resp.\ following the radial point estimate with the matching color.}
  \label{FigNModFlow}
  \end{figure}

  \pfstep{Uniform a priori estimate for approximations of $P_+$.} Having established~\eqref{EqNModSymb}, note first that $u$ can be estimated for $t_{*,0}\leq t_*\leq t_{*,0}+2$ using~\eqref{EqNModExistLoc}; we thus obtain on $[0,\tau_0]\times\iota^+$ the estimate
  \begin{equation}
  \label{EqNModEst}
    \|u\|_{\dot H_{\ebeop;\bop}^{(\sfs';k),(2\gamma_\sscri,\gamma_+,\gamma_\cface)}} \leq C\Bigl( \|P_+^0 u\|_{\dot H_{\ebeop;\bop}^{(\sfs'-1;k),(2\gamma_\sscri+2,\gamma_++2,\gamma_\cface-2)}} + \|u\|_{\dot H_{\ebeop;\bop}^{(\sfs'_\flat;k),(2\gamma_\sscri,\gamma_+,\gamma_\cface)}} \Bigr)
  \end{equation}
  for $P_+^0$. The estimates used thus far (being based on positive commutator arguments, and thus relying only on the positivity of certain functions, which is a stable property under perturbations) are stable under small principally scalar perturbations of $P_+$ in the class $x_\sscri^2\tau^2\rho_\cface^{-2}\CI_\bop\Diff_\ebeop^2$. Thus, the estimate~\eqref{EqNModEst} holds, for a uniform constant $C$, also for all operators
  \[
    P_+^\eps := \chi_0\Bigl(\frac{\tau}{\eps}\Bigr)N_{\iota^+}(P_0) + \Bigl(1-\chi_0\Bigl(\frac{\tau}{\eps}\Bigr)\Bigr)P_+
  \]
  where $\chi_0\in\CIc([0,2))$ equals $1$ on $[0,1]$, and $\eps>0$ is sufficiently small. (Thus, $P_+^\eps$ is dilation-invariant for $\tau\leq\eps$, i.e., for sufficiently late times $t_*$.) Note here that $P_+-P_+^\eps=\chi_0(\frac{\tau}{\eps})(P_+-N_{\iota^+}(P_0))$ is the product of $\chi_0(\frac{\tau}{\eps})$, which is uniformly bounded in $\CI_\bop([0,\infty]_\tau)$ and converges to $0$ in $\tau\CI_\bop([0,\infty])$, and $P_+-N_{\iota^+}(P_0)\in x_\sscri^2\tau^3\rho_\cface^{-2}\Diff_\ebeop^2$; so $P_+^\eps\to P_+$ in $x_\sscri^2\tau^2\rho_\cface^{-2}\CI_\bop\Diff_\ebeop^2$ as $\eps\to 0$, indeed. For later use, we also note that
  \begin{equation}
  \label{EqNModDiff}
    \{ P_+^\eps-N_{\iota^+}(P_0) \colon \eps\in(0,1] \}\ \text{is bounded in}\ x_\sscri^2\tau^3\rho_\cface^{-2}\CI_\bop\Diff_\ebeop^2.
  \end{equation}

  \pfstep{Inversion of the normal operator.} We next estimate the error term in~\eqref{EqNModEst} using the estimate for the normal operator $N_{\iota^+}(P_0)$ provided by Proposition~\ref{PropNipFwd}; using~\eqref{EqNModDiff}, this gives the uniform (in $\eps$) estimate
  \begin{align*}
    \|u\|_{\dot H_{\ebeop;\bop}^{(\sfs'_\flat;k),(2\gamma_\sscri,\gamma_+,\gamma_\cface)}} &\leq C \|N_{\iota^+}(P_0)u\|_{\dot H_{\ebeop;\bop}^{(\sfs'_\flat-1;k),(2\gamma_\sscri,\gamma_+,\gamma_\cface)}} \\
      &\leq C\Bigl( \|P_+^\eps u\|_{\dot H_{\ebeop;\bop}^{(\sfs'_\flat-1;k),(2\gamma_\sscri+2,\gamma_++2,\gamma_\cface-2)}} + \|u\|_{\dot H_{\ebeop;\bop}^{(\sfs'_\flat+1;k),(2\gamma_\sscri,\gamma_+-1,\gamma_\cface)}}\Bigr).
  \end{align*}
  Plugging this into~\eqref{EqNModEst} for $P_+^\eps$, we have thus proved
  \begin{equation}
  \label{EqNModEst2}
    \|u\|_{\dot H_{\ebeop;\bop}^{(\sfs';k),(2\gamma_\sscri,\gamma_+,\gamma_\cface)}} \leq C\Bigl( \|P_+^\eps u\|_{\dot H_{\ebeop;\bop}^{(\sfs'-1;k),(2\gamma_\sscri+2,\gamma_++2,\gamma_\cface-2)}} + \|u\|_{\dot H_{\ebeop;\bop}^{(\sfs';k),(2\gamma_\sscri,\gamma_+-1,\gamma_\cface)}} \Bigr);
  \end{equation}
  here we also use that $\sfs'_\flat+1\leq\sfs'$.

  We can exploit the weakened $\tau$-decay order by introducing another ($\eps$-independent) cutoff in $\tau$. To wit, observe that for fixed orders $\sfs',k,\gamma_\sscri,\gamma_+,\gamma_\cface$, multiplication by $\chi_0(\frac{\tau}{\delta})$ tends to $0$ as a map
  \[
    \dot H_{\ebeop;\bop}^{(\sfs';k),(2\gamma_\sscri,\gamma_+,\gamma_\cface)} \to \dot H_{\ebeop;\bop}^{(\sfs';k),(2\gamma_\sscri,\gamma_+-1,\gamma_\cface)}
  \]
  when $\delta\to 0$; this follows from the uniform boundedness of multiplication by $\frac{\tau}{\delta}\chi_0(\frac{\tau}{\delta})$, which in turn follows from the uniform boundedness of this family of multipliers in $\CIb([0,\infty]_\tau)$. We can, therefore, fix $\delta>0$ small enough so that $C\|\chi_0(\frac{\tau}{\delta})u\|_{\dot H_{\ebeop;\bop}^{(\sfs';k),(2\gamma_\sscri,\gamma_+-1,\gamma_\cface)}}\leq\frac12\|u\|_{\dot H_{\ebeop;\bop}^{(\sfs';k),(2\gamma_\sscri,\gamma_+,\gamma_\cface)}}$ can be absorbed into the left-hand side of~\eqref{EqNModEst2}. Since the norm of $(1-\chi_0(\frac{\tau}{\delta}))u$ is already controlled by the norm of $P_+ u=P_+^\eps u$ (for $\eps<\frac{\delta}{2}$) by the \emph{finite-time} existence and uniqueness theory established at the beginning of the proof, we have thus proved the uniform (for all small $\eps$) a priori estimate
  \begin{equation}
  \label{EqNModEst3}
    \|u\|_{\dot H_{\ebeop;\bop}^{(\sfs';k),(2\gamma_\sscri,\gamma_+,\gamma_\cface)}([0,\tau_0]\times\iota^+)} \leq C \|P_+^\eps u\|_{\dot H_{\ebeop;\bop}^{(\sfs'-1;k),(2\gamma_\sscri+2,\gamma_++2,\gamma_\cface-2)}([0,\tau_0]\times\iota^+)}.
  \end{equation}

  \pfstep{Conclusion of the approximation argument.} Let now $f\in\dot H_{\ebeop;\bop}^{(\sfs'-1;k),(2\gamma_\sscri+2,\gamma_++2,\gamma_\cface-2)}([0,\tau_0]\times\iota^+)$. We claim that the forward solutions $u^\eps$ of the equations
  \[
    P_+^\eps u^\eps = f,
  \]
  for all small $\eps>0$, satisfy
  \begin{equation}
  \label{EqNModueps}
    u^\eps\in\dot H_{\ebeop;\bop}^{(\sfs';k),(2\gamma_\sscri,\gamma_+,\gamma_\cface)}([0,\tau_0]\times\iota^+).
  \end{equation}
  To prove this, let us first solve this equation from $\tau=\tau_0$ up to $\tau=\frac{\eps}{2}$ using finite-time existence, yielding an ``initial'' solution $u^\eps_{\rm ini}$ on $\tau^{-1}([\frac{\eps}{2},\tau_0])$ with regularity $(\sfs';k)$ and weights $2\gamma_\sscri,\gamma_\cface$; and then, writing $\chi^{\eps/2}=\chi_0(\frac{\tau}{\eps/2})$, we have
  \[
    u^\eps = (1-\chi^{\eps/2})u^\eps_{\rm ini} + N_{\iota^+}(P_0)^{-1}\bigl(\chi^{\eps/2} f + [N_{\iota^+}(P_0),\chi^{\eps/2}]u_{\rm ini}^\eps \bigr),
  \]
  where $N_{\iota^+}(P_0)^{-1}$ is the forward solution operator; we use here that $P_+^\eps=N_{\iota^+}(P_0)$ on $\supp\chi^{\eps/2}$. Proposition~\ref{PropNipFwd} thus gives~\eqref{EqNModueps}.

  But then the uniform a priori estimate~\eqref{EqNModEst3}, with constant right-hand side given by $C$ times the norm of $f$, applies to the solutions $u^\eps$ and implies their uniform boundedness. We can thus find a sequence $u_{\eps_j}$, $\eps_j\to 0$, which converges weakly to some $u\in\dot H_{\ebeop;\bop}^{(\sfs';k),(2\gamma_\sscri,\gamma_+,\gamma_\cface)}([0,\tau_0]\times\iota^+)$ (with norm still bounded by that of $f$). Since $P_+^\eps=P_+$ for $\tau\geq 2\eps$, this $u$ satisfies $P_+ u=f$ for $\tau\geq 2\eps$, and thus everywhere since $\eps>0$ is arbitrary. This completes the proof.
\end{proof}

\begin{rmk}[The choice of $P_+$]
\label{RmkNModOther}
  The choice of $P_+$ is motivated by three considerations. First, it must equal $P_0$ in a neighborhood of $\scri^+$ (or, at least, of $\scri^+\cap\iota^+$). Second, we demand that it have the same normal operator at $\iota^+$ as $P_0$, as it is only for this normal operator that we can guarantee the absence of resonances (i.e., poles of $N_{\iota^+}^0(P_0,\sigma)^{-1}$ in Proposition~\ref{PropNip}), as a consequence of the $\tface$-admissibility assumption, and thus the desired polynomial bound $\tau^{\gamma_+}$, $\gamma_+<\gamma_\sscri$, for its forward solutions. (Other choices, such as introducing an artificial timelike boundary at $R=1$ with Dirichlet boundary conditions and working in $\{R>1\}$, or smoothing out the normal operator near $R=0$ so that it has smooth coefficients across $R\omega=0\in\R^3$, might introduce spurious resonances.) Third, we want to avoid introducing further asymptotic regimes (such as $\cK^+$) which would complicate the asymptotic analysis of $P_+$; this is why near $R=0$ we extend $N_{\iota^+}(P_0)$ by homogeneity with respect to $\tau$-dilations. The price to pay here is that $P_+$ has a timelike curve of conic singularities at $R=0$, which, however, is easily analyzed using the Mellin transform and Proposition~\ref{PropNip} as well as ideas originating in \cite{HintzConicWave}.
\end{rmk}

\begin{rmk}[Normal operator estimates]
\label{RmkFNoScriNorm}
  Observe that we are able to prove Theorem~\ref{ThmNFw}, concerning forward solutions for $P_+$ on a manifold with three relevant ``ideal'' boundary hypersurfaces ($\cface$, $\scri^+$, $\tau^{-1}(0)$), by combining symbolic analysis (to control high ebe-frequencies) with only \emph{one} normal operator estimate (namely for $N_{\iota^+}(P_0)$, corresponding to the boundary hypersurface $\tau^{-1}(0)$). This is unexpected from the perspective of general PDEs (and Fredholm theory): for an elliptic PDE on a noncompact or singular space, one typically needs to control normal operators at \emph{every} boundary hypersurface. It is the \emph{hyperbolic} nature of $P_+$, which leads to our ability to solve $P_+ u=f$ rather easily on arbitrarily long (but finite) $t_*$-intervals, that allows to effectively localize near $\tau=0$. In particular, we do not need any normal operator estimates at $\cface$ or at $\scri^+$ (the latter being discussed in~\cite[\S{7}]{HintzVasyScrieb}).
\end{rmk}

\subsection{Proof of Theorem~\usref{ThmNFw}}
\label{SsNPf}

As discussed before the statement of Theorem~\ref{ThmNFw}, we solve
\begin{equation}
\label{EqNPf}
  P_0 u=f\in H_{\etbop;\bop}^{(\sfs;k),(2\alpha_\sscri+2,\alpha_++2,0)}(\Omega)^{\bullet,-}.
\end{equation}
in several steps. First, using the finite-time solvability theory of Proposition~\ref{PropETSolvb}, we can solve this equation up to any finite value of $t_*$, and hence we may assume that
\[
  \supp f\subset \{t_*\geq\tau_0^{-1} \}
\]
for any fixed large constant $\tau_0^{-1}$, $\tau_0\in(0,1]$, on $\supp f$. For $j=1,2,3,4$, fix cutoff functions
\[
  \chi_j \in \CIc([0,j+1)),\ \ \chi_j=1\ \text{near}\ [0,j],
\]
which we regard as functions on $\Omega$ via $(t_*,r,\omega)\mapsto\chi_j(R)$, $R=\frac{r}{t_*}$. We may take $\chi_1=\chi$ in the notation of Definition~\ref{DefNModOp}.

\pfstep{A first solution near $\scri^+$.} Note first that $(1-\chi_2)f$, being supported in $R\geq 2$, defines a function on the space $[0,\infty]_\tau\times\iota^+$, $\tau=t_*^{-1}$, which vanishes near $R=0$; to wit,
\[
  (1-\chi_2)f \in \dot H_{\ebeop;\bop}^{(\sfs;k),(2(\alpha_\sscri+\frac32)+2,(\alpha_++2)+2,\infty)}([0,1]_\tau\times\iota^+),
\]
where the shifts of $\frac32$ and $2$ are due to our convention to use unweighted b-volume densities on $[0,\infty]\times\iota^+$, as opposed to metric volume densities on $\Omega$ (cf.\ \eqref{EqMUetbMShift}). We have $\gamma_\sscri:=\alpha_\sscri+\frac32<1+\ubar S$ by~\eqref{EqNFwThr} and $\gamma_+:=\alpha_++2<\alpha_\sscri+\frac32=\gamma_\sscri$ by~\eqref{EqNFwThr2}; and since $\gamma_+-\frac12=\alpha_++\frac32\in(\beta^-,\beta^+)$ (see Definition~\ref{DefSSAlephAdm}\eqref{ItSSAlephAdm3} and Definition~\ref{DefSStfAdm}), we also have $\gamma_+<\beta^++\frac12<\beta^++1$,\footnote{The extra buffer of $\frac12$ here can be explained as follows. Microlocally, estimates for $u$ in $\Htb^{\sfs,(\alpha_+,0)}$ impose the threshold condition $\sfs+\alpha_+<-\frac12$ (up to a shift by $\frac12\vartheta_{\pa\cK^+,{\rm out}}$) at the radial set $\cR_{\pa\cK^+,{\rm out}}$, and thus one does not get full b-regularity on its flow-out $\cW_{\rm out}$; what one \emph{can} get is ``module regularity'' at $\cW_{\rm out}$, which in particular entails b-regularity in $r,\omega$ (in the coordinates $t_*,r,\omega$). Sobolev embedding in the spatial variables (which lie in a $3$-dimensional space) only thus gives $\frac32$ degrees of $L^\infty$-decay, i.e., $\alpha_++\frac32$ pointwise decay at $\iota^+$. Only if one has additional b-regularity also in $t_*$ does one get $\alpha_++2$ pointwise decay. The reason of the $\frac12$ buffer is thus the discrepancy between the microlocal propagation results that limit the raw 3b-regularity order, and the pointwise decay at $\iota^+$ which is governed by the resonances---that are excluded in $\Im\sigma>-1-\beta^+$ (rather than merely $-\frac12-\beta^+$) by Proposition~\ref{PropNip}.} so the weights $\gamma_\sscri$, $\gamma_+$, and $\gamma_\cface:=-(\alpha_++\frac32)$ satisfy the requirements of Proposition~\ref{PropNMod}; and we take the regularity order in Proposition~\ref{PropNMod} to be $\sfs+1$. (This is acceptable by our requirement that $\sfs$ have margin $\geq 1$.) Let us thus define $u_\sscri$ as the forward solution of
\begin{equation}
\label{EqNPfuscri}
  P_+ u_\sscri = (1-\chi_2)f,\qquad u_\sscri\in\dot H_{\ebeop;\bop}^{(\sfs'+1;k),(2(\alpha_\sscri+\frac32),\alpha_++2,\gamma_\cface)}([0,\tau_0]_\tau\times\iota^+)
\end{equation}
where $\sfs'$ is defined in Proposition~\ref{PropNMod} and equals $\sfs$ for $R\geq 2$. Since $P_+=P_0$ for $R\geq 2$ and thus on $\supp(1-\chi_2)$, we then have
\begin{align*}
  f - P_0\bigl((1-\chi_3)u_\sscri\bigr) &= f-(1-\chi_3)f + [P_0,\chi_3]u_\sscri \\
    &= \chi_3 f+[P_0,\chi_3]u_\sscri =: f' \in H_{\etbop;\bop}^{(\sfs;k),(2\alpha_\sscri,\alpha_+,0)}(\Omega)^{\bullet,-},
\end{align*}
where $f'$ vanishes for $R\geq 4$ and thus near $\scri^+$; hence, we can regard $f'$ as a function on $M_0$, which then satisfies
\begin{equation}
\label{EqNPffp}
  f' \in H_{\tbop;\bop}^{(\sfs;k),(\alpha_++2,0)}(\Omega)^{\bullet,-},\qquad \supp f'\subset\Bigl\{ t_*\geq\tau_0^{-1},\ R=\frac{r}{t_*}\leq 4 \Bigr\}.
\end{equation}

\pfstep{Control of the solution sourced by $f'$.} Define next $u'$ as the forward solution of
\begin{equation}
\label{EqNPfup}
  P_0 u' = f',\quad u'\in H_{\tbop;\bop}^{(\sfs;k),(\alpha_+-\delta,-\aleph)}(\Omega)^{\bullet,-};
\end{equation}
we use here the $\aleph$-admissibility of $P_0$, specifically~\eqref{EqSSAlephAdmSol}. We must show that, when regarded as a function on $M$, the solution $u'$ is e3b-regular near $\scri^+$; this only concerns the behavior of $u'$ for large $R$. Consider thus
\[
  P_0\bigl((1-\chi_4)u'\bigr) = (1-\chi_4)f' - [P_0,\chi_4]u' =: f'';
\]
the first term vanishes by~\eqref{EqNPffp}, and the second term is supported in $\{4\leq R\leq 5\}$, and thus away from $\scri^+$; therefore, upon passing to $[0,\infty]\times\iota^+$,
\[
  f'' \in \dot H_{\ebeop;\bop}^{(\sfs-1;k),(\infty,(\alpha_+-\delta+2)+2,\infty)}([0,\tau_0]\times\iota^+),\quad \supp f''\subset\Bigl\{ t_*\geq\tau_0^{-1},\ 4\leq R\leq 5 \Bigr\}.
\]
But on $\supp(1-\chi_4)$, we have $R\geq 4$ and thus $P_0=P_+$; we may thus compare $(1-\chi_4)u'$ to the forward solution
\[
  u'' \in \dot H_{\ebeop;\bop}^{(\sfs';k),(2(\alpha_\sscri+\frac32),\alpha_+-\delta+2,\gamma_\cface)}([0,\tau_0]\times\iota^+)
\]
of $P_+ u''=f''$ controlled using Proposition~\ref{PropNMod} (now with $\gamma_+$ reduced to $\alpha_+-\delta+2$). The difference $\tilde u:=u''-(1-\chi_4)u'$ satisfies $P_+\tilde u=f''-f''=0$; and since $(1-\chi_4)u'$ vanishes near $R=0$, we can apply Proposition~\ref{PropNMod}\eqref{ItNModUniq} to conclude that $\tilde u=0$, so taking into account the vanishing of $1-\chi_4$ for $R\leq 2$,
\begin{equation}
\label{EqNPfupp}
  (1-\chi_4)u' = u'' \in H_{\etbop;\bop}^{(\sfs;k),(2\alpha_\sscri,\alpha_+-\delta,\infty)}(\Omega)^{\bullet,-}.
\end{equation}

\pfstep{Conclusion.} We have now shown that the forward solution of~\eqref{EqNPf} can be written as
\[
  u = (1-\chi_3)u_\sscri - u' = (1-\chi_3)u_\sscri - \chi_4 u' - u''.
\]
In view of~\eqref{EqNPfuscri}, we have $(1-\chi_3)u_\sscri\in H_{\etbop;\bop}^{(\sfs+1;k),(2\alpha_\sscri,\alpha_+,\infty)}(\Omega)^{\bullet,-}$; and moreover~\eqref{EqNPfup} implies $\chi_4 u'\in H_{\etbop;\bop}^{(\sfs;k),(\infty,\alpha_+-\delta,-\aleph)}(\Omega)^{\bullet,-}$. Together with~\eqref{EqNPfupp}, this implies $u\in H_{\etbop;\bop}^{(\sfs;k),(2\alpha_\sscri,\alpha_+-\delta,-\aleph)}(\Omega)^{\bullet,-}$. Since all constructions involved in this argument come with quantitative norm bounds, we also obtain the estimate~\eqref{EqNFW}. The proof of Theorem~\ref{ThmNFw} is complete.

\begin{rmk}[Modifications for relaxed $\aleph$-admissibility]
\label{RmkNFwRelax}
  Continuing Remark~\ref{RmkSSAlephAdmRelax}, if $P_0$ is $\aleph$-admissible with margin $d_\aleph$, then the estimate~\eqref{EqNFW} is valid when $\sfs$ is in addition admissible with weights $\alpha_+,-\aleph$ \emph{and margin $d_\aleph$} and also with weights $\alpha_+,0$ and margin $1+d_\aleph$, except the norm of $u$ must be replaced by the norm of $H_{\etbop;\bop}^{(\sfs-d_\aleph;k),(2\alpha_\sscri,\alpha_+-\delta,-\aleph)}$. Indeed, the only place in the proof where the $\aleph$-admissibility of $P_0$ is used is in~\eqref{EqNPfup}, where in the current relaxed setting one would only obtain the 3b-regularity $\sfs-d_\aleph$ instead of $\sfs$. This loss then simply remains in place throughout the remainder of the proof.
\end{rmk}

\section{Forward solutions of wave-type operators and tame estimates}
\label{SF}

We are now ready to state and prove the main result of this paper concerning forward solutions of admissible wave-type operators on dynamical perturbations of Kerr spacetimes and b-tame estimates for them. The regularity assumptions and function spaces here \emph{only} concern b-regularity on $M$, or equivalently on $M_1$;\footnote{Recall here that a spanning set of $\Vb(\tilde M_1)$ lifts to a spanning set of $\Vb(\tilde M)$, with an explicit such spanning set being given in~\eqref{EqCTbtildeMGlobal}.} thus, we shall use powers of $\rho_\sscri=x_\sscri^2$ as weights at $\scri^+$, rather than $x_\sscri$. (This leads only to some factors of $2$ relative to earlier estimates.) Recall that we use the volume density of the Kerr (or, equivalently, the Minkowski) metric to define Sobolev spaces.
\begin{equation}
\label{EqFSob}
  \Hb^{k,(\alpha_\sscri,\alpha_+,\alpha_\cK)}(\tilde M) := \rho_\sscri^{\alpha_\sscri}\rho_+^{\alpha_+}\rho_\cK^{\alpha_\cK}\Hb^k(\tilde M),
\end{equation}
where $\rho_\sscri,\rho_+,\rho_\cK$ are defining functions of $\sscri_1^+,\iota^+,\cK^+\subset\tilde M_1$, respectively (Definition~\ref{DefCMSpacetime}); elements of $\Hb^k(\tilde M)$ lie in $L^2$ together with their up to $k$-fold b-derivatives. For domains such as $\Omega=\cl_{\tilde M}\{r\geq\bhm,\ t_*\geq 1\}$ in Definition~\ref{DefCMDomain}, we recall from~\eqref{EqMUSuppExt}--\eqref{EqMUSuppExt2} the notation $\Hb^k(\Omega)^{\bullet,-}$ for the space of restrictions to $\{r>\bhm\}$ of elements of $\Hb^k(\tilde M)$ with support in $t_*\geq 1$; similarly for weighted spaces. Explicit defining functions and frames of b-vector fields are given in Remark~\ref{RmkCMExpl} and Example~\ref{ExCTbtildeM}.

\begin{thm}[Tame estimates for forward solutions]
\label{ThmF}
  On the subextremal Kerr spacetime with parameters $\bhm>0$, $a\in(-\bhm,\bhm)$, let $P_1$ be an $\aleph$-admissible stationary wave-type operator (Definition~\usref{DefSSAlephAdm}) with $\sface$-weight $\alpha_+$ and $\sface$-loss $\delta$,\footnote{or the relaxed version in Remark~\usref{RmkSSAlephAdmRelax}, discussed in Remark~\ref{RmkFRelax} below} acting on sections of the bundle $\cE\to M$. Recall the domain $\Omega_*=\cl_M\{\ft_*\geq 1\}$ from~\eqref{EqDyDomain}. Let $\alpha_\sscri<-\frac12+\ubar S$ where $\ubar S$ was defined in~\eqref{EqSSAdmubarS}, and assume that $\alpha_+<\alpha_\sscri-\frac12$. There exist $\eps_0,d\in\N$ such that the following holds. Let $\ell_\sscri\in(0,\frac12]$, $\ell_+>\delta$, and $\ell_\cK>\aleph$. Let $P$ be an admissible wave-type operator of class $((0;k+d),(2\ell_\sscri,\ell_+,\ell_\cK))$ (Definition~\usref{DefSDWAdm}) relative to $P_0$ where $P_0$ is equal to $P_1$ or a small perturbation thereof in the class of stationary wave-type operators (including for different but nearby Kerr parameters), and such that $\|P-P_0\|_{(0;d),(2\ell_\sscri,\ell_+,\ell_\cK),\Omega_*}<\eps_0$ (Definition~\usref{DefSDWAdmNorm}). Then for all $k\in\N_0$ and all source terms
  \begin{equation}
  \label{EqFf}
    f \in \Hb^{k+d,(\alpha_\sscri+1,\alpha_++2,0)}(\Omega_*;\cE)^{\bullet,-},
  \end{equation}
  the unique forward solution $u$ of $P u=f$ satisfies
  \[
    u \in \Hb^{k,(\alpha_\sscri,\alpha_+-\delta,-\aleph)}(\Omega_*;\cE)^{\bullet,-}
  \]
  and the tame estimate
  \begin{equation}
  \label{EqFTame}
  \begin{split}
    &\|u\|_{\Hb^{k,(\alpha_\sscri,\alpha_+-\delta,-\aleph)}(\Omega_*;\cE)^{\bullet,-}} \\
    &\qquad \leq C_k\Bigl( \|f\|_{\Hb^{k+d,(\alpha_\sscri+1,\alpha_++2,0)}(\Omega_*;\cE)^{\bullet,-}} \\
    &\qquad \hspace{4em} + \|P-P_0\|_{(0;k+d),(2\ell_\sscri,\ell_+,\ell_\cK),\Omega_*} \|f\|_{\Hb^{d,(\alpha_\sscri+1,\alpha_++2,0)}(\Omega_*;\cE)^{\bullet,-}} \Bigr).
  \end{split}
  \end{equation}
  If the principal part of $P$ is equal to that of the Kerr wave operator, this holds for $\Omega=\cl_M\{t_*\geq 1\}$ in place of $\Omega_*$. Finally, the map $(P,f)\mapsto u=P^{-1}f$ is continuous if we use the $((0;k+d),(2\ell_\sscri,\ell_+,\ell_\cK))$-norm on $P-P_0$, the $\Hb^{k+d,(\alpha_\sscri+1,\alpha_++2,0)}$-norm on $f$, and, for any fixed $\eps>0$, the $\Hb^{k-1,(\alpha_\sscri-\eps,\alpha_+-\delta-\eps,-\aleph-\eps)}$-norm on $u$.
\end{thm}

It may be easier to parse~\eqref{EqFTame} when using an unweighted b-density, which is given by $\mu_\bop:=\rho_\sscri^3\rho_+^4\rho_\cK|\dd t\,\dd x|$ (where $\rho_\sscri=x_\sscri^2$) for any choice of boundary defining functions on $\tilde M_1$ such as the weights
\begin{equation}
\label{EqFWeights}
  \rho_\sscri = \frac{t_*}{t_*+r}, \quad
  \rho_+ = \frac{t_*+r}{t_* r},\quad
  \rho_\cK = \frac{r}{t_*+r}
\end{equation}
from~\eqref{EqCMExpl} (in which case $\rho_\sscri^3\rho_+^4\rho_\cK=\frac{1}{t_* r^3}$). The assumptions now read that $P_1$ be $\aleph$-admissible with $\sface$-weight $\beta_+-2$ and $\sface$-loss $\delta$, further $\beta_\sscri<1+\ubar S$ and $\beta_+<\beta_\sscri$; here we use the identifications
\begin{equation}
\label{EqFGammas}
  \beta_\sscri = \alpha_\sscri+\frac32,\quad
  \beta_+ = \alpha_++2.
\end{equation}
The estimate~\eqref{EqFTame} reads (dropping the bundle from the notation)
\begin{subequations}
\begin{equation}
\label{EqFTame20}
\begin{split}
  &\|u\|_{\Hb^{k,(\beta_\sscri,\beta_+-\delta,-\aleph+\frac12)}(\Omega_*,\mu_\bop)^{\bullet,-}} \\
  &\qquad \leq C_k\Bigl(\|f\|_{\Hb^{k+d,(\beta_\sscri+1,\beta_++2,\frac12)}(\Omega_*,\mu_\bop)^{\bullet,-}} \\
  &\qquad \hspace{4em} + \|P-P_0\|_{(0;k+d),(2\ell_\sscri,\ell_+,\ell_\cK),\Omega_*}\|f\|_{\Hb^{d,(\beta_\sscri+1,\beta_++2,\frac12)}(\Omega_*,\mu_\bop)^{\bullet,-}}\Bigr),
\end{split}
\end{equation}
i.e., writing out the norms,
\begin{equation}
\label{EqFTame2}
\begin{split}
  &\iint_{\Omega_*} | \rho_\sscri^{-\beta_\sscri}\rho_+^{-\beta_++\delta}\rho_\cK^{\aleph-\frac12}(t_*\pa_{t_*},r\pa_x)^{\leq k}u |^2\,\frac{\dd t_*\,\dd x}{t_* r^3} \\
  &\qquad \leq C_k \Biggl(\; \iint_{\Omega_*} |\rho_\sscri^{-\beta_\sscri-1}\rho_+^{-\beta_+-2}\rho_\cK^{-\frac12}(t_*\pa_{t_*},r\pa_x)^{\leq k+d}f |^2\,\frac{\dd t_*\,\dd x}{t_* r^3} \\
  &\qquad \hspace{4em} + \|P-P_0\|_{(0;k+d),(2\ell_\sscri,\ell_+,\ell_\cK),\Omega_*} \iint_{\Omega_*} |\rho_\sscri^{-\beta_\sscri-1}\rho_+^{-\beta_+-2}\rho_\cK^{-\frac12}(t_*\pa_{t_*},r\pa_x)^{\leq d}f|^2\,\frac{\dd t_*\,\dd x}{t_* r^3}\Biggr),
\end{split}
\end{equation}
\end{subequations}
under the assumption that $f$ vanish for $\ft_*\leq 1$ (i.e., in the past of the Cauchy hypersurface of $\Omega_*$). Furthermore, Sobolev embedding (Lemma~\ref{LemmaMUCe3bSob}) gives the pointwise bound
\begin{equation}
\label{EqFTame3}
  | (t_*\pa_{t_*},r\pa_x)^{\leq k-3} u | \leq \text{(R.H.S.~of~\eqref{EqFTame2})} \times \rho_\sscri^{\beta_\sscri}\rho_+^{\beta_+-\delta}\rho_\cK^{-\aleph+\frac12}.
\end{equation}

\begin{proof}[Proof of Theorem~\usref{ThmF}]
  \emph{Since the proof uses e3b-estimates, we shall use powers of $x_\sscri$ as weights at $\scri^+$.} We drop the bundle $\cE$ from the notation. We shall construct $u$ by a variant of the approximation argument used in the proof of Proposition~\ref{PropNMod}. Thus, we start with the regularity estimate of Theorem~\ref{ThmEReg}, applied to operators $P^\eps$ which transition from $P$ to the stationary model $P_0$ at times $t_*\gtrsim\eps^{-1}$ (where we will take $\eps\searrow 0$) and estimate the error term in~\eqref{EqEReg} and \eqref{EqERegTame} using the estimate for the stationary model $P_0$ proved in Theorem~\ref{ThmNFw}; we then obtain a uniform a priori estimate for $P^\eps$. Since $P^\eps u^\eps=f$, where
  \begin{equation}
  \label{EqFTamef}
    f \in H_{\etbop;\bop}^{(\sfs;k),(2\alpha_\sscri+2,\alpha_++2,0)}(\Omega_*)^{\bullet,-},
  \end{equation}
  can be solved (the stationary theory being applicable for very late times for $P^\eps$), a weak compactness argument produces the desired solution $u$.

  The details are as follows. Fix
  \[
    \chi_0 \in \CIc([0,2)),\ \chi_0|_{[0,1]}=1,\quad
    \chi^\eps := \chi_0\Bigl(\frac{1}{\eps t_*}\Bigr),
  \]
  which equals $1$ in a $\cO(\eps)$-neighborhood of $\iota^+\cup\cK^+$, and set
  \[
    P^\eps := \chi^\eps P_0 + (1-\chi^\eps)P,
  \]
  which is thus exactly $t_*$-translation-invariant for $t_*\geq\eps^{-1}$. Note that
  \begin{equation}
  \label{EqFPepsP0}
    P^\eps-P_0=(1-\chi^\eps)(P-P_0),
  \end{equation}
  so since $1-\chi^\eps$ is uniformly bounded in $\CI_\bop(\Omega_*)$, $P^\eps$ is a uniformly (in $\eps$) small perturbation of $P_0$, i.e., it satisfies the same assumptions as $P$ itself.

  Let $\sfs\in\CI({}^\etbop S^*M)$ be an admissible order function (Definition~\ref{DefDyO}) for $P_0$ with weights $\alpha_+-\delta,-\aleph$ and margin $2+\delta+\max(1,\aleph)$, and set $\sfs_0:=\sfs-2$. Thus, $\sfs_0$ is admissible with weights $\alpha_+,-\aleph$ and margin $\max(1,\aleph)$. Write $\tilde d_0,\tilde d$ for the quantities denoted by $d_0,d$ in Theorem~\usref{ThmEReg}.

  \pfstep{Step~1. Small b-regularity orders.} We first consider ``small'' b-regularity orders $k$, i.e., $k\leq\tilde d$. 

  \pfsubstep{(1.1)}{Uniform a priori regularity estimate.} We start with the uniform (in $\eps$) a priori estimate~\eqref{EqEReg} with $\alpha_+-\delta$ and $-\aleph$ in place of $\alpha_+$ and $\alpha_\cK$. Then
  \begin{equation}
  \label{EqFRegEst}
  \begin{split}
    &\|u\|_{H_{\etbop;\bop}^{(\sfs;k),(2\alpha_\sscri,\alpha_+-\delta,-\aleph)}(\Omega_*)^{\bullet,-}} \\
    &\qquad \leq C_k\Bigl( \|P^\eps u\|_{H_{\etbop;\bop}^{(\sfs;k),(2\alpha_\sscri+2,\alpha_+-\delta+2,-\aleph)}(\Omega_*)^{\bullet,-}} + \|u\|_{H_{\etbop;\bop}^{(\sfs_0;k),(2\alpha_\sscri,\alpha_+-\delta,-\aleph)}(\Omega_*)^{\bullet,-}} \Bigr).
  \end{split}
  \end{equation}
  This estimate requires $P$, and thus $P^\eps$, to be of class $((\tilde d_0;k),(2\ell_\sscri,\ell_+,\ell_\cK))$, which for $k\leq\tilde d$ follows if $P$ is of class $((0;d),(2\ell_\sscri,\ell_+,\ell_\cK))$ with $d\geq\tilde d_0+\tilde d$.

  \pfsubstep{(1.2)}{Inversion of $P_0$.} Now, $\sfs_0$ satisfies the conditions of Theorem~\ref{ThmNFw}, which yields
  \begin{equation}
  \label{EqFP0Est}
    \|u\|_{H_{\etbop;\bop}^{(\sfs_0;k),(2\alpha_\sscri,\alpha_+-\delta,-\aleph)}(\Omega_*)^{\bullet,-}} \leq C_k\|P_0 u\|_{H_{\etbop;\bop}^{(\sfs_0;k),(2\alpha_\sscri+2,\alpha_+,0)}(\Omega_*)^{\bullet,-}}.
  \end{equation}
  (The domain $\Omega_*$ used here is contained in the domain $\Omega$ used in Theorem~\ref{ThmNFw}.) We wish to replace $P_0$ by $P^\eps$; note then that by~\eqref{EqFPepsP0} and~\eqref{EqSDWp0p1}--\eqref{EqSDWAdmOp},
  \begin{equation}
  \label{EqFPepsP0Diff}
    P^\eps - P_0 \in x_\sscri^2\rho_+^2 \cC_{\etbop;\bop}^{(\tilde d_0;k),(2\ell_\sscri,\ell_+,\ell_\cK)}\Diff_\etbop^2
  \end{equation}
  is uniformly bounded, and therefore
  \begin{equation}
  \label{EqFPepsP0Bd}
    \|(P^\eps-P_0)u\|_{H_{\etbop;\bop}^{(\sfs_0;k),(2\alpha_\sscri+2,\alpha_+,0)}(\Omega_*)^{\bullet,-}} \leq C'_k\|u\|_{H_{\etbop;\bop}^{(\sfs_0+2;k),(2\alpha_\sscri,\alpha_+-\ell_+,-\ell_\cK)}(\Omega_*)^{\bullet,-}}
  \end{equation}
  for some $\eps$-independent constant $C'_k$. Upon localization to a small neighborhood of $t_*^{-1}(\infty)=\iota^+\cup\cK^+$, the right-hand side of~\eqref{EqFPepsP0Bd} is \emph{small} compared to the left-hand side of~\eqref{EqFRegEst}. Concretely, set $\delta':=\min(\ell_+-\delta,\ell_\cK-\aleph)$, then
  \begin{equation}
  \label{EqFchietau}
  \begin{split}
    \|\chi^\eta u\|_{H_{\etbop;\bop}^{(\sfs_0+2;k),(2\alpha_\sscri,\alpha_+-\ell_+,-\ell_\cK)}} &\leq C''_k\|\chi^\eta t_*^{-\delta'} u\|_{H_{\etbop;\bop}^{(\sfs_0+2;k),(2\alpha_\sscri,\alpha_+-\delta,-\aleph)}} \\
      &\leq C''_k C'''_k\eta^{\delta'}\|u\|_{H_{\etbop;\bop}^{(\sfs_0+2;k),(2\alpha_\sscri,\alpha_+-\delta,-\aleph)}},
  \end{split}
  \end{equation}
  where in the passage to the second line we use that $\chi^\eta=t_*^{-\delta'}\eta^{\delta'}\chi_1(\frac{1}{\eta t_*})$ where the family $\chi_1(\frac{1}{\eta t_*})=(\frac{1}{\eta t_*})^{\delta'}\chi_0(\frac{1}{\eta t_*})$ is uniformly bounded in $\CI_\bop([1,\infty]_{t_*})$ for $\eta\in(0,1]$. Fix then $\eta>0$ so that $C_k C'_k C''_k C'''_k\eta^{\delta'}=\frac12$. Plugging these estimates into~\eqref{EqFP0Est}, using $\sfs_0+2=\sfs$, and absorbing the term $C_k C'_k C''_k C'''_k\eta^{\delta'}\|u\|_{H_{\etbop;\bop}^{(\sfs;k),(2\alpha_\sscri,\alpha_+-\delta,-\aleph)}}$ into the left-hand side of~\eqref{EqFRegEst}, we obtain
  \[
    \|u\|_{H_{\etbop;\bop}^{(\sfs_0;k),(2\alpha_\sscri,\alpha_+-\delta,-\aleph)}(\Omega_*)^{\bullet,-}} \leq C_k\Bigl( \|P^\eps u\|_{H_{\etbop;\bop}^{(\sfs_0;k),(2\alpha_\sscri+2,\alpha_+,0)}(\Omega_*)^{\bullet,-}} + \|(1-\chi^\eta)u\|_{H_{\etbop;\bop}^{(\sfs;k),(2\alpha_\sscri,0,0)}}\Bigr),
  \]
  with the $\iota^+$- and $\cK^+$-orders of the error term being arbitrary since $\supp(1-\chi^\eta)$ is disjoint from $\iota^+\cup\cK^+$.

  But since $t_*$ is bounded on $\supp(1-\chi^\eta)$, we can control $u$ on it using the \emph{finite-time} estimate from Proposition~\ref{PropETSolvb}; note that for sufficiently small $\eps>0$, we have $P^\eps=P$ on $\supp(1-\chi^\eta)$. We therefore obtain the uniform (in $\eps>0$) a priori estimate
  \begin{equation}
  \label{EqFUnifApriori}
    \|u\|_{H_{\etbop;\bop}^{(\sfs;k),(2\alpha_\sscri,\alpha_+-\delta,-\aleph)}(\Omega_*)^{\bullet,-}} \leq C_k \|P^\eps u\|_{H_{\etbop;\bop}^{(\sfs;k),(2\alpha_\sscri+2,\alpha_+,0)}(\Omega_*)^{\bullet,-}}.
  \end{equation}

  \pfsubstep{(1.3)}{Conclusion of the approximation argument.} We claim that the forward solution $u^\eps$ of $P^\eps u^\eps=f$, with $f$ as in~\eqref{EqFTamef}, satisfies
  \begin{equation}
  \label{EqFuepsMem}
    u^\eps\in H_{\etbop;\bop}^{(\sfs;k),(2\alpha_\sscri,\alpha_+-\delta,\aleph)}(\Omega_*)^{\bullet,-}.
  \end{equation}
  Letting $u^\eps_{\rm ini}$ be the finite-time solution of $P^\eps u^\eps_{\rm ini}=f$ for $t_*\leq\frac{2}{\eps}$ produced by Proposition~\ref{PropETSolvb}, we have
  \[
    u^\eps = (1-\chi^{\eps/2})u^\eps_{\rm ini} + P_0^{-1}\bigl(\chi^{\eps/2}f + [P_0,\chi^{\eps/2}]u^\eps_{\rm ini}\bigr)
  \]
  since $P^\eps=P_0$ on $\supp\chi^{\eps/2}$. The first term lies in $H_{\etbop;\bop}^{(\sfs+1;k),(2\alpha_\sscri,0,0)}(\Omega_*)^{\bullet,-}$ by Proposition~\ref{PropETSolvb} (the orders at $\iota^+$ and $\cK^+$ being irrelevant since $\supp(1-\chi^{\eps/2})$ is disjoint from $\iota^+\cup\cK^+$). Since the argument of $P_0^{-1}$ in the second term therefore lies in $H_{\etbop;\bop}^{(\sfs;k),(2\alpha_\sscri+2,\alpha_++2,0)}(\Omega_*)^{\bullet,-}$, Theorem~\ref{ThmNFw} shows that the second term is an element of $H_{\etbop;\bop}^{(\sfs;k),(2\alpha_\sscri,\alpha_+-\delta,-\aleph)}$. This establishes~\eqref{EqFuepsMem}. Moreover, the estimate~\eqref{EqFUnifApriori} holds for $u^\eps$ in place of $u$.

  But then the estimate~\eqref{EqFUnifApriori} shows that $u^\eps$ is uniformly bounded by the norm of $f$. We may thus pass to a subsequence $u^{\eps_j}$, $\eps_j\searrow 0$, which converges weakly to some $u\in H_{\etbop;\bop}^{(\sfs;k),(2\alpha_\sscri,\alpha_+-\delta,-\aleph)}(\Omega_*)^{\bullet,-}$. Since $P u=P^\eps u^\eps=f$ for $t_*\leq\frac{1}{2\eps}$ and $\eps$ is arbitrary, we conclude that $P u=f$; the estimate~\eqref{EqFUnifApriori} also applies to the equation $P u=f$, yielding
  \begin{equation}
  \label{EqFUnif}
    \|u\|_{H_{\etbop;\bop}^{(\sfs;k),(2\alpha_\sscri,\alpha_+-\delta,-\aleph)}(\Omega_*)^{\bullet,-}} \leq C_k \|f\|_{H_{\etbop;\bop}^{(\sfs;k),(2\alpha_\sscri+2,\alpha_+,0)}(\Omega_*)^{\bullet,-}},\quad k\leq\tilde d.
  \end{equation}

  \pfstep{Step~2. Large b-regularity orders and tame estimates.} For $k>\tilde d$, the arguments are very similar. We begin with the b-tame estimate~\eqref{EqERegTame} for $P^\eps$, again with $\alpha_+-\delta$ and $-\aleph$ in place of $\alpha_+$ and $\alpha_\cK$, and estimate the final norm on $u$ in terms of $P^\eps u$ using~\eqref{EqFUnifApriori}; this gives the a priori estimate
  \begin{align*}
    &\|u\|_{H_{\etbop;\bop}^{(\sfs;k),(2\alpha_\sscri,\alpha_+-\delta,-\aleph)}(\Omega_*)^{\bullet,-}} \\
    &\qquad \leq C_k\biggl( \|P^\eps u\|_{H_{\etbop;\bop}^{(\sfs;k),(2\alpha_\sscri+2,\alpha_+-\delta+2,-\aleph)}(\Omega_*)^{\bullet,-}} + \|u\|_{H_{\etbop;\bop}^{(\sfs_0;k),(2\alpha_\sscri,\alpha_+-\delta,-\aleph)}(\Omega_*)^{\bullet,-}} \\
    &\qquad \quad \hspace{6em} + \|P^\eps-P_0\|_{(\tilde d_0;k),(2\ell_\sscri,\ell_+,\ell_\cK),\Omega_*} \|P^\eps u\|_{H_{\etbop;\bop}^{(\sfs_0;\tilde d),(2\alpha_\sscri,\alpha_+,0)}(\Omega_*)^{\bullet,-}} \biggr).
  \end{align*}
  (The smallness of $\|P-P_0\|_{(0;d),(2\ell_\sscri,\ell_+,\ell_\cK)}$ implies the smallness of $P-P_0$ required in Theorem~\ref{ThmEReg} provided we take $d\geq\tilde d_0+\tilde d$.) We again estimate the norm of $u$ on the right by~\eqref{EqFP0Est}, and replace $P_0$ on the right-hand side of~\eqref{EqFP0Est} by $P^\eps$; instead of~\eqref{EqFPepsP0Bd} however, we now use a b-tame bound on the action of $P^\eps-P_0$ (see~\eqref{EqFPepsP0Diff}) on $u$ instead of~\eqref{EqFPepsP0Bd}, which reads
  \begin{equation}
  \label{EqFPeP0u}
  \begin{split}
    &\|(P^\eps-P_0)u\|_{H_{\etbop;\bop}^{(\sfs_0;k),(2\alpha_\sscri+2,\alpha_+,0)}(\Omega_*)^{\bullet,-}} \\
    &\qquad \leq C_k\Bigl( \|u\|_{H_{\etbop;\bop}^{(\sfs_0+2;k),(2\alpha_\sscri,\alpha_+-\ell_+,-\ell_\cK)}(\Omega_*)^{\bullet,-}} \\
    &\qquad \quad \hspace{6em} + \|P^\eps-P_0\|_{(\tilde d_0;k),(2\ell_\sscri,\ell_+,\ell_\cK)}\|u\|_{H_{\etbop;\bop}^{(\sfs_0;\tilde d),(2\alpha_\sscri,\alpha_+-\ell_+,-\ell_\cK)}(\Omega_*)^{\bullet,-}} \Bigr);
  \end{split}
  \end{equation}
  we use that the low-regularity norm $\|P^\eps-P_0\|_{(0;d),(2\ell_\sscri,\ell_+,\ell_\cK),\Omega_*}$ is bounded $C\eps_0$, as noted after~\eqref{EqFPepsP0}. (The proof of this estimate, which is completely analogous to similar estimates used in the proofs of Propositions~\ref{PropRb} and~\ref{PropETSolvbTame}; concretely, one uses Lemma~\ref{LemmaETMult} near the initial (and similarly near the final) boundary hypersurface of $\Omega_*$, and Proposition~\ref{PropMTameMicr} away from these.) We again estimate the final norm of $u$ on the right of~\eqref{EqFPeP0u} using~\eqref{EqFUnifApriori}. Upon splitting the first norm on the right of~\eqref{EqFPeP0u} into the norm of $\chi^\eta u$ and $(1-\chi^\eta)u$, the norm of $\chi^\eta u$ is again small when $\eta>0$ is sufficiently small (depending \emph{only} on $k$) in view of~\eqref{EqFchietau} and can thus be absorbed; and for the remaining piece $\|(1-\chi^\eta)u\|_{H_{\etbop;\bop}^{(\sfs_0+2;k),(2\alpha_\sscri,0,0)}}$, we use the b-tame bound for $u$ near $\supp(1-\chi^\eta)$ (which is contained in a \emph{finite} $\ft_*$-slab) given by Proposition~\ref{PropETSolvbTame}. We have now established the uniform (in $\eps$) a priori estimate
  \begin{equation}
  \label{EqFPepsTame}
  \begin{split}
    &\|u\|_{H_{\etbop;\bop}^{(\sfs;k),(2\alpha_\sscri,\alpha_+-\delta,-\aleph)}(\Omega_*)^{\bullet,-}} \\
    &\qquad\leq C_k\biggl( \|P^\eps u\|_{H_{\etbop;\bop}^{(\sfs;k),(2\alpha_\sscri+2,\alpha_++2,0)}(\Omega_*)^{\bullet,-}} \\
    &\qquad \quad \hspace{6em} + C'_k\|P-P_0\|_{(\tilde d_0;k),(2\ell_\sscri,\ell_+,\ell_\cK),\Omega_*} \|P^\eps u\|_{H_{\etbop;\bop}^{(\sfs_0;\tilde d),(2\alpha_\sscri,\alpha_+,0)}(\Omega_*)^{\bullet,-}} \biggr),
  \end{split}
  \end{equation}
  where we bounded the norm of $P^\eps-P_0$ by that of $P-P_0$ using~\eqref{EqFPepsP0} and the uniform boundedness of multiplication by $1-\chi^\eps$ on weighted $\cC_{\etbop;\bop}^{(\tilde d_0;k)}(\Omega_*)$-spaces.

  But recall from~\eqref{EqFuepsMem} that the forward solution $u^\eps$ of $P^\eps u^\eps=f$ (with $f$ as in~\eqref{EqFTamef}) does lie in the space $H_{\etbop;\bop}^{(\sfs;k),(2\alpha_\sscri,\alpha_+-\delta,-\aleph)}(\Omega_*)^{\bullet,-}$, and hence the estimate~\eqref{EqFPepsTame} applies to $u^\eps$ and $f$ in place of $u$ and $P^\eps u$. The right-hand side is now independent of $\eps$, giving a uniform norm bound for $u^\eps$; and a weak sequential compactness argument produces, as before, a weak subsequential limit $u\in H_{\etbop;\bop}^{(\sfs;k),(2\alpha_\sscri,\alpha_+-\delta,-\aleph)}(\Omega_*)^{\bullet,-}$ of $u^\eps$ as $\eps\searrow 0$, which solves $P u=f$ and satisfies the estimate~\eqref{EqFPepsTame} with $P$ in place of $P^\eps$.

  \pfstep{Step~3. Pure b-estimates.} Passage from the above estimates on mixed $(\etbop;\bop)$-spaces to pure b-Sobolev spaces is straightforward: one uses that b-regularity implies e3b-regularity on (subsets of) $\tilde M$ (since $\CI_\bop\cV_\etbop\subset\CI_\bop\cV_\bop$), so
  \begin{alignat*}{2}
    \|u\|_{\Hb^{k,(2\alpha_\sscri,\alpha_+-\delta,-\aleph)}} &\leq C_k\|u\|_{H_{\etbop;\bop}^{(\sfs;k+d'),(2\alpha_\sscri,\alpha_+-\delta,-\aleph)}},&\quad& d':=\max(0,-\lfloor\sfs\rfloor), \\
    \|P u\|_{H_{\etbop;\bop}^{(\sfs;k+d'),(2\alpha_\sscri+2,\alpha_+,0)}} &\leq C_k\|P u\|_{\Hb^{k+d'+d'',(2\alpha_\sscri+2,\alpha_+,0)}},&\quad &d'':=\max(0,\lceil\sfs\rceil).
  \end{alignat*}
  If we require $d\geq d'+d''$ (in addition to the requirements stated earlier), we thus obtain~\eqref{EqFTame} (upon halving the $\scri^+$-weight due to passing to powers of $\rho_\sscri=x_\sscri^2$).

  The penultimate statement follows from the local-in-time solvability of $P u=f$ using standard energy estimates together with Proposition~\ref{PropETSolvb} to bridge the gap between $\Omega$ and $\Omega_*$.

  \pfstep{Step~4. Continuous dependence of $u$ on $(P,f)$.} If $P_{(j)}\to P$ and $f_{(j)}\to f$ in the respective norms, then the uniformity of the estimate~\eqref{EqFTame} shows that the forward solutions $u_{(j)}$ of
  \begin{equation}
  \label{EqFPjuj}
    P_{(j)}u_{(j)}=f_{(j)}
  \end{equation}
  are uniformly bounded in the Hilbert space $\Hb^{k,(\alpha_\sscri,\alpha_+-\delta,-\aleph)}$. Every subsequence of $u_{(j)}$ then has a further subsequence that converges weakly to some limit $u$. Taking the distributional limit of~\eqref{EqFPjuj} as $j\to\infty$ yields $P u=f$. But forward solutions of this equation are unique; therefore, we have convergence $u_{(j)}\weakto u$ of the full sequence. Since the inclusion $\Hb^{k,(\alpha_\sscri,\alpha_+-\delta,-\aleph)}\hra\Hb^{k-1,(\alpha_\sscri-\eps,\alpha_+-\delta-\eps,-\aleph-\eps)}$ is compact, the claim follows.
\end{proof}

\begin{rmk}[Modifications for relaxed $\aleph$-admissibility]
\label{RmkFRelax}
  Continuing Remarks~\ref{RmkSSAlephAdmRelax} and \ref{RmkNFwRelax}, we note that an additional 3b-regularity loss of order $d_\aleph$ incurred from the inversion of $P_0$, or the resulting e3b-regularity loss of order $d_\aleph$ in Theorem~\ref{ThmNFw}, is easily absorbed in the proof above. We now require $\sfs$ to be admissible with the increased margin $2+\delta+\max(1,\aleph)+d_\aleph$ and choose $\sfs_0:=\sfs-2-d_\aleph$. In the estimate~\eqref{EqFP0Est}, the norm on $P_0 u$ now uses the e3b-regularity order $\sfs_0+d_\aleph$ instead of $\sfs_0$, which can ultimately be absorbed into the left-hand side of~\eqref{EqFRegEst} since $\sfs_0+d_\aleph+2=\sfs$.
\end{rmk}

\begin{rmk}[Usage of $(\etbop;\bop)$-Sobolev spaces]
\label{RmkFetbb}
  Note that it is crucial that the symbolic (i.e., e3b-microlocal) estimate of Theorem~\ref{ThmEReg} be lossless as far as the b-regularity order is concerned: both $u$ and $P u$ are estimated in function spaces featuring $k$ orders of b-regularity. If one needed, say, $k+1$ b-orders on $P u$, our argument above would not close anymore, as the right-hand side of~\eqref{EqFPepsP0Bd} would now also require $k+1$ b-orders, and then~\eqref{EqFchietau}, with $k+1$ in place of $k$, could no longer be absorbed into the left-hand side of~\eqref{EqFRegEst}. But in the absence of a viable b-microlocal theory for $P$, we are forced to work with weaker (e3b-)spaces for the symbolic analysis of $P$; and to close estimates for forward solutions of $P u=f$, we can only exploit properties of $P$ (and its stationary model $P_0$) on e3b-spaces (relative to $\Hb^k$ with any fixed value of $k$). This is why the proof of Theorem~\ref{ThmF} uses the full force of the e3b-microlocal (and b-tame) estimates developed in earlier parts of the paper. However, once the b-tame estimate~\eqref{EqFTame}---which is very lossy in the b-regularity sense---is in place, the e3b-microlocal theory is no longer needed at all: it is the b-regularity of $u$, the strongest sensible notion of regularity on a singular space (here: the compactified space $M$, or rather its subset $\Omega_*$), that is most convenient when trying to establish more precise properties of $u$, such as stronger asymptotics. We demonstrate this in (the proofs of) Theorems~\ref{ThmDp} and \ref{ThmA2} below. (See also Remark~\ref{RmkDpLinear}.)
\end{rmk}

\begin{rmk}[Conditions on the weights]
\label{RmkFWeights}
  The conditions on $\alpha_\sscri$ and $\alpha_+$ imply that $\alpha_+<-1+\ubar S$; but also $\alpha_++\frac32$ must lie in the indicial gap $(\beta^-,\beta^+)$ of $\wh{P_0}(0)$ (Definitions \ref{DefSSAlephAdm}\eqref{ItSSAlephAdm2} and~\ref{DefSStfAdm}). Thus, in order for Theorem~\ref{ThmF} to be applicable, we must have $\beta^-<\frac12+\ubar S$ where $(\beta^-,\beta^+)$ is the indicial gap for $\wh{P_0}(0)$. This condition is satisfied in all applications of interest to us here and in \cite{HintzKerrStab}. (If it were violated, it may still be the case that one can prove Theorem~\ref{ThmF} for some $\alpha_+<-1+\ubar S$, provided a strengthening of~\eqref{EqSSAlephAdmSol}, in which one allows $P_0 u$ to have weight $\alpha_\sface+2-\delta$ at $\sface$, holds, where one moreover requires $\alpha_\sface-\delta<-1+\ubar S$. See \cite[Lemma~3.27]{HintzGlueLocIII} for such a strengthening of $2$-admissibility for a linearized gauge-fixed Einstein operator.)
\end{rmk}

\section{Application I: nonlinear scalar waves}
\label{SA1}

We fix a subextremal Kerr black hole metric
\[
  g=g_{\bhm,a},\quad \bhm>0,\ a\in(-\bhm,\bhm).
\]
We study the small data global existence for three ``toy'' scalar semi-/quasi-linear wave equations that have already been partially discussed in the literature (see the discussion in~\S\ref{SssIMN2}): power nonlinearities in~\S\ref{SsA1Power}, a simple quasilinear equation in~\S\ref{SsA1Q}, and a null-form nonlinearity in~\S\ref{SsA1Null}. In all cases, the relevant linearized operators asymptote at large $t_*$ to the scalar wave operator $\Box_g$ on Kerr, which does not have zero energy bound states. Upon establishing the admissibility of $\Box_g$ (\S\ref{SsA1Adm}) and recalling a suitable version of the Nash--Moser iteration scheme in~\S\ref{SsA1NM}, the small data global existence in each of these three case is then an immediate consequence of Theorem~\ref{ThmF} upon applying the Nash--Moser iteration scheme.

\begin{rmk}[Decay]
\label{RmkA1Decay}
  In this section, we are only concerned with the basic problem of global existence. We shall not discuss (almost) sharp decay rates for solutions; we content ourselves with pointwise $t_*^{-\frac12}$-decay in compact spatial regions, which will follow directly from Theorem~\ref{ThmF} (cf.\ \eqref{EqFTame3} with $\aleph=0$). (The decay rates that we establish for the solutions $u$ of our nonlinear equations \emph{are} sharp as far as the decay rates for the source term $f$ are concerned, i.e., $u$ decays $1$, $2$, and $0$ orders faster than $f$ at $\scri^+$, $\iota^+$, and $\cK^+$, respectively.) In~\S\S\ref{SsDp} and \ref{SA2}, however, we will need stronger (but still far from sharp) decay for $u$ (thus strengthening also the requirements on $f$) in order to close the nonlinear iteration; we develop the necessary tools in~\S\ref{SD}.
\end{rmk}

\subsection{Admissibility of the scalar wave operator}
\label{SsA1Adm}

We recall $\Omega=\cl_M\{t_*\geq 1\}$ from Definition~\ref{DefCMDomain}, and recall from Lemma~\ref{LemmaTsKLMetric} that the $t_*$-level sets are spacelike for $g=g_{\bhm,a}$ and transversal to future null infinity. Recalling Definition~\ref{DefSSAlephAdm}, we shall prove:

\begin{thm}[$0$-admissibility of $\Box_g$]
\label{ThmA1Adm}
  $\Box_g$ is $0$-admissible with any $\sface$-weight $\alpha_\sface\in(-\frac32,-\frac12)$ and $\sface$-loss $0$.
\end{thm}
\begin{proof}
  First of all, $\Box_g$ is a stationary wave-type operator by Example~\ref{ExSSAdmBox}. Trapping admissibility (Definition~\ref{DefSSTrapAdm}) holds because $\Box_g$ is symmetric with respect to the metric $L^2$-inner product, so the left-hand side of~\eqref{EqSSTrapAdm} is equal to $0$. The $\tface$-admissibility of $\Box_g$ with weight $\beta\in(\beta^-,\beta^+):=(0,1)$ was checked in Proposition~\ref{PropSptfAdm}.

  It remains to check~\eqref{EqSSAlephAdmSol} for $\aleph=\delta=0$. 

  \pfstep{Basic properties of $\wh{\Box_g}(\sigma)$.} Given $f\in H_{\tbop;\bop}^{(\sfs;k),(\alpha_\sface+2,0)}(\Omega)^{\bullet,-}$ (with $\Omega=\cl_M\{t_*\geq 1\}$), we wish to define
  \begin{equation}
  \label{EqA1Admu}
    u := \cF^{-1}\hat u,\quad \hat u(\sigma) = \wh{\Box_g}(\sigma)^{-1}\hat f(\sigma),
  \end{equation}
  for $\sigma\in\R$, and more generally for $\sigma\in\C$ with $\Im\sigma\geq 0$. This is well-defined for sufficiently large $|\sigma|$ by Theorem~\ref{ThmSpHi}. Furthermore, the mode stability of the Kerr metric on the real axis proved in \cite{ShlapentokhRothmanModeStability,AnderssonMaPaganiniWhitingModeStab} implies that $\wh{\Box_g}(\sigma)^{-1}$ is well-defined as the inverse of the map~\eqref{EqSpBMap} also for real $\sigma$; this uses the regularity~\eqref{EqSpBNullspace} of kernel elements and a separation into fully separated mode solutions (using spheroidal harmonics) that is possible for real $\sigma\neq 0$, and the index $0$ property of~\eqref{EqSpBMap}. For $\sigma=0$, the triviality of $\ker\wh{\Box_g}(0)$ on the space of Theorem~\ref{ThmSp0}\eqref{ItSp0Reg} was checked in \cite[Lemma~3.19]{HintzKdSMS}; this gives the invertibility for small $\sigma$ by Theorem~\ref{ThmSpLo}. For complex $\sigma$, one may use \cite{WhitingKerrModeStability,FinsterSmollerSpheroidalCompleteness} to prove the invertibility of $\wh{\Box_g}(\sigma)$ or, alternatively, the more robust continuity argument from \cite[Theorem~1.7, \S{3.9}]{HintzKdSMS}. The summary of the latter argument is as follows: the high-energy estimates of Theorem~\ref{ThmSpHi} and the real mode stability (including near $0$) imply that for some large $C>1$, $\wh{\Box_{g_{s\bhm,s a}}}(\sigma)$ is invertible when $|\sigma|\leq C^{-1}$ or $|\sigma|\geq C$, for all $s\in[0,1]$. Using Rouch\'e's theorem, one can then show that a putative pole $\sigma(s)$ of $\wh{\Box_{g_{s\bhm,s a}}}(\sigma)^{-1}$ for $s=1$ persists for nearby $s$; and the compactness of $\{\sigma\in\C\colon\Im\sigma\geq 0,\ |\sigma|\in[C^{-1},C]\}$ implies that the set of $s$ such that a pole exists for all parameters in $[s,1]$ is closed; thus it must be all of $[0,1]$. But for $s=0$, i.e., for the Schwarzschild metric, there are no poles; this is a contradiction.

  \pfstep{Low-energy estimates.} Returning to~\eqref{EqA1Admu}, we use Lemma~\ref{LemmaMUetbFTb} (with $\chi$ there equal to $1$ on $\supp f$) to pass to the Fourier transform; note then that for $j=0,\ldots,k$, we have
  \begin{equation}
  \label{EqA1AdmLou}
    (\sigma\pa_\sigma)^j \hat u(\sigma) = \sum_{0\leq j_1+j_2\leq j} \binom{j}{j_1} \bigl((\sigma\pa_\sigma)^{j_1}\wh{\Box_g}(\sigma)^{-1}\bigr) \bigl((\sigma\pa_\sigma)^{j_2}\hat f(\sigma)\bigr).
  \end{equation}
  For low frequencies $|\sigma|\leq 1$, we use~\eqref{EqSpLoInvReg} with $\alpha=\alpha_\sface$, and write $\sfs_\sigma,\sfs_\sigma$ in place of $\sfs,\sfr$, to estimate\footnote{We give up the gain of $j_1$ orders of sc-b-transition-regularity in the output of~\eqref{EqSpLoInvReg}.}
  \begin{align*}
    &\bigl\|\bigl((\sigma\pa_\sigma)^{j_1}\wh{\Box_g}(\sigma)^{-1}\bigr)\bigl((\sigma\pa_\sigma)^{j_2}\hat f(\sigma)\bigr) \bigr\|_{H_{(\scbtop,|\sigma|);\bop}^{(\sfs_\sigma;k-j),(\sfs_\sigma+\alpha_\sface,\alpha_\sface,0)}} \\
    &\qquad \leq C\|(\sigma\pa_\sigma)^{j_2}\hat f(\sigma)\|_{H_{(\scbtop,|\sigma|);\bop}^{(\sfs_\sigma-1;k-j),(\sfs_\sigma+\alpha_\sface+1,\alpha_\sface+2,0)}}.
  \end{align*}
  Upon squaring and integrating over $\sigma\in[-1,1]$, we thus obtain
  \begin{equation}
  \label{EqA1AdmLoEst}
    \int_{-1}^1 \|(\sigma\pa_\sigma)^j\hat u(\sigma)\|_{H_{(\scbtop,|\sigma|);\bop}^{(\sfs_\sigma;k-j),(\sfs_\sigma+\alpha_\sface,\alpha_\sface,0)}}^2\,\dd\sigma \leq C_k\|f\|_{H_{\tbop;\bop}^{(\sfs-1;k),(\alpha_\sface+2,0)}}^2.
  \end{equation}
  (The weight $\alpha_\sface+2$ thus arises from the $\tface$-order in~\eqref{EqSpLoInvReg}, with the two power difference between the orders of $u$ and $f$ being due to the mapping properties of the zero energy operator $\wh{\Box_g}(0)$ in Theorem~\ref{ThmSp0}.)

  \pfstep{High-energy estimates.} For high frequencies $|\sigma|\geq 1$, we note that $\sfp_1=0$ in~\eqref{EqSpHiTrp1} and thus $\gamma_-<\gamma_+=0$ in~\eqref{EqSpHiTrpGammas}. Therefore, the semiclassical loss function~\eqref{EqSpHiTrLoss} satisfies $\delta_\Gamma(\Im\sigma-\eps)\leq\eps$ for $\Im\sigma\geq 0$ and all $\eps>0$. We then use~\eqref{EqSpHiInvReg} to estimate, for $j_1+j_2+j_3\leq j$ and denoting the orders $\sfs_\scop,\sfr,\sfb$ from Theorem~\ref{ThmSpHi} by $\sfs_\sigma,\sfs_\sigma,\sfs_\sigma$,
  \begin{equation}
  \label{EqA1AdmHi}
  \begin{split}
    &\bigl\||\sigma|^{j_3}\bigl((\sigma\pa_\sigma)^{j_1}\wh{\Box_g}(\sigma)^{-1}\bigr)\bigl((\sigma\pa_\sigma)^{j_2}\hat f(\sigma)\bigr)\bigr\|_{H_{(\scop,|\sigma|^{-1});\bop}^{(\sfs_\sigma;k-j);\sfs_\sigma+\alpha_\sface,\sfs_\sigma}} \\
    &\qquad \leq C_\eps|\sigma|^{j_3-1+(j_1+1)\eps} \|(\sigma\pa_\sigma)^{j_2}\hat f(\sigma)\|_{H_{(\scop,|\sigma|^{-1});\bop^+}^{(\sfs_\sigma-1;k-j+j_1),\sfs_\sigma+\alpha_\sface+1,\sfs_\sigma}} \\
    &\qquad \leq C'_\eps \sum_{j_4=0}^{k-j+j_1} |\sigma|^{j_3-1+(j_1+1)\eps} \| |\sigma|^{j_4}(\sigma\pa_\sigma)^{j_2}\hat f(\sigma)\|_{H_{(\scop,|\sigma|^{-1});\bop}^{(\sfs_\sigma-1;k-j+j_1-j_4),\sfs_\sigma+\alpha_\sface+1,\sfs_\sigma}},
  \end{split}
  \end{equation}
  where in the passage to the last line we used the definition~\eqref{EqSpHiNormb}. Let us choose $\eps=\frac{1}{k+1}$, then the total $|\sigma|$-power in the last line is $\leq j_3+j_4$, and therefore the total number of powers of $(|\sigma|,\sigma\pa_\sigma)$ is $\leq j_2+j_3+j_4$; since $(k-j+j_1-j_4)+(j_2+j_3+j_4)\leq (k-(j_1+j_2+j_3)+j_1-j_4)+(j_2+j_3+j_4)=k$, we can thus bound this further by
  \[
    C'_k \sum_{j=j_1+j_2\leq k} \|(\sigma\pa_\sigma)^{j_1}|\sigma|^{j_2}\hat f(\sigma)\|_{H_{(\scop,|\sigma|^{-1});\bop}^{(\sfs_\sigma-1;k-j),\sfs_\sigma+\alpha_\sface+1,\sfs_\sigma}}.
  \]
  Using Lemma~\ref{LemmaMUetbFTb}, we thus obtain
  \begin{equation}
  \label{EqA1AdmHiFinal}
    \sum_\pm \sum_{j=j_3+j_4\leq k}\int_{\pm[1,\infty)} \bigl\| \sigma^{j_3}(\sigma\pa_\sigma)^{j_4}\bigl(\wh{\Box_g}(\sigma)^{-1}\hat f(\sigma)\bigr)\bigr\|_{H_{(\scop,|\sigma|^{-1});\bop}^{(\sfs_\sigma;k-j);\sfs_\sigma+\alpha_\sface,\sfs_\sigma}}^2\,\dd\sigma \leq C''_k\|f\|_{H_{\tbop;\bop}^{(\sfs;k),(\alpha_\sface+1,0)}}^2.
  \end{equation}
  (That we need to use the 3b-differential order $\sfs$ here is due to the fact that the semiclassical order in~\eqref{EqA1AdmHi} is $\sfs_\sigma$, which in turn comes from the $|\sigma|$-power loss in the preceding estimate due to trapping.\footnote{It is here that having only an arbitrarily small loss in high-energy estimates is used in a crucial way. If one could not choose $\eps$ arbitrarily small, one would be faced with $|\sigma|$-power losses which increase with $k$, resulting in being only able to estimate the high-frequency part of $u$ by $\|\pa_{t_*}^{\leq c k}f\|_{H_{\tbop;\bop}^{(\sfs;k),(\alpha_\sface+1,0)}}$ for some $c>0$. The regularity gain in the output of~\eqref{EqSpHiInvReg} does not help since one could only exploit it by working with the 3b-regularity function $\sfs+k$ instead of $\sfs$, say, which however would not satisfy the required threshold conditions at the radial sets anymore unless $\sfs$ itself satisfied them with a $k$-dependent margin, i.e., $\sfs$ would need to be $k$-independent, which, in turn, would break the b-tame analysis developed in this paper.})

  \pfstep{Quantitative bound. Support property.} Using Lemma~\ref{LemmaMUetbFTb} again, we conclude that $u$ defined in~\eqref{EqA1Admu} satisfies $u\in H_{\tbop;\bop}^{(\sfs;k),(\alpha_\sface,0)}(\tilde M_0)$ and
  \begin{equation}
  \label{EqA1AdmBd}
    \|u\|_{H_{\tbop;\bop}^{(\sfs;k),(\alpha_\sface,0)}(\tilde M_0)} \leq C_k \|f\|_{H_{\tbop;\bop}^{(\sfs;k),(\alpha_\sface+2,0)}}.
  \end{equation}
  To complete the proof, we need to show that $t_*\geq 1$ on $\supp u$, which we do using a Paley--Wiener argument similarly to the proof of Proposition~\ref{PropNipFwd}. In view of the quantitative bound~\eqref{EqA1AdmBd}, it suffices to prove this for a dense set of source terms $f$. Let us thus assume that $f\in\CIc((t_+,\infty)_{t_*}\times X^\circ)$ where $t_+>1$ and $X^\circ=\{x\in\R^3\colon|x|\geq\bhm\}$, and let $\phi\in\CIc(X^\circ)$. Then Theorems~\ref{ThmSpLo} and \ref{ThmSpHi} give
  \[
    |\la\hat u(\sigma,\cdot),\phi\ra_{L^2(X^\circ)}| \leq C\la\sigma\ra^N \| \hat f(\sigma,\cdot) \|_{H^N}\quad\forall\sigma\in\C,\ \Im\sigma\geq 0,
  \]
  for some $N$, where we crudely estimate the norms in those theorems by the $H^N$-norm times a large power of $\la\sigma\ra$. Since $t_*\geq t_+>1$ on $\supp f$, we have
  \[
    \|\hat f(\sigma,\cdot)\|_{H^N} = \la\sigma\ra^{-2 L}\|\cF((1+D_t^2)^L f)(\sigma,\cdot)\|_{H^N} \leq C_L\la\sigma\ra^{-2 L}e^{-t_+\Im\sigma}
  \]
  for all $L\in\N_0$; we fix $L$ such that $N-2 L<-1$. Using the continuity of $\la\hat u(\sigma,\cdot),\phi\ra_{L^2(X^\circ)}$ in $\{\sigma\in\C\colon \sigma\neq 0,\ \Im\sigma\geq 0\}$ which follows from Theorem~\ref{ThmSpB}\eqref{ItSpBInvReg} and its uniform boundedness for $|\sigma|\leq 1$ which follows from Theorem~\ref{ThmSpLo}, the contour
  \[
    \la u(t_*,\cdot),\phi\ra_{L^2(X^\circ)} = \frac{1}{2\pi}\int_{\Im\sigma=\gamma} e^{-i\sigma t_*}\la\hat u(\sigma,\cdot),\phi\ra_{L^2(X^\circ)}\,\dd\sigma
  \]
  can be shifted from $\gamma=0$ to $\gamma>0$ without changing its value. The integrand is bounded in absolute value by $e^{t_*\Im\sigma}e^{-t_+\Im\sigma}\la\sigma\ra^{N-2 L}$, and hence the integral over $\Im\sigma=\gamma$ is bounded by $C e^{-(t_+-t_*)\gamma}$. For $t_*\leq 1<t_+$, this tends to $0$ as $\gamma\to\infty$. Since $\phi$ was arbitrary, this implies $u=0$ for $t_*\leq 1$, as desired.
\end{proof}

Theorem~\ref{ThmF} is thus applicable with $\aleph=0$ and $\delta=0$. For $P_0=\Box_g$, we have $S=0$ in~\eqref{EqSSAdmOp} and thus $\ubar S=0$. The weights in Theorem~\ref{ThmF} are thus subject to the conditions $\alpha_\sscri<-\frac12$, $\alpha_+\in(-\frac32,\alpha_\sscri-\frac12)$; we may thus take
\begin{equation}
\label{EqA1Admalphas}
  \alpha_\sscri=-\frac12-\eps_\sscri,\ \ \alpha_+=-1-\eps_+;\quad 0<\eps_+<\eps_\sscri<\frac12.
\end{equation}
Translated into pointwise bounds using~\eqref{EqFGammas} and \eqref{EqFTame3}, we thus obtain pointwise bounds
\begin{subequations}
\begin{equation}
\label{EqA1AdmPwu}
  |D_\bop^{\leq k}u|\lesssim\rho_\sscri^{1-\eps_\sscri}\rho_+^{1-\eps_+}\rho_\cK^{\frac12}
\end{equation}
for forward solutions $u$ of $\Box_g u=f$ as well as for their $k$-fold b-derivatives (i.e., derivatives along $t_*\pa_{t_*}$ and $|x|\pa_x$), given suitably decaying $f$ as in~\eqref{EqFf} or equivalently \eqref{EqFTame2}. For example, by~\eqref{EqMUCe3bSobConv}, it suffices if $f$ satisfies the bounds
\begin{equation}
\label{EqA1AdmPwf}
  |D_\bop^{\leq k+d}f| \lesssim \rho_\sscri^{2-\eps_\sscri+\delta}\rho_+^{3-\eps_++\delta}\rho_\cK^{\frac12+\delta}
\end{equation}
\end{subequations}
for some $\delta>0$. The same results hold true for non-stationary perturbations of $\Box_g$ which are admissible wave-type operators of class $((0;k+d),(2\ell_\sscri,\ell_+,\ell_\cK))$ for any fixed $\ell_\sscri\in(0,\frac12]$, $\ell_+>\delta$, and $\ell_\cK>0$.

\bigskip

Observe that the only property of $\Box_g$ used in the proof of Theorem~\ref{ThmA1Adm} (besides the fact that its principal symbol is the dual metric function of a subextremal Kerr metric) is its mode stability in $\Im\sigma\geq 0$. Thus, the same arguments show:

\begin{thm}[Mode stability implies $0$-admissibility]
\label{ThmA1Gen}
  Let $P_0$ be a stationary wave-type operator on the subextremal Kerr spacetime $(M,g_{\bhm,a})$ (Definition~\usref{DefSSAdm}). Suppose that
  \begin{enumerate}
  \item $P_0$ is $\tface$-admissible with weight $\beta:=\alpha_\sface+\frac32$ (Definition~\usref{DefSStfAdm});
  \item for all $0\neq\sigma\in\C$, $\Im\sigma\geq 0$, the nullspace of $\wh{P_0}(\sigma)$ on $\cA^{1+\ubar S-\eps}(X)$ is trivial for some $\eps>0$;
  \item the nullspace of $\wh{P_0}(0)$ on $\cA^\beta(X)$ is trivial.
  \end{enumerate}
  Then $P_0$ is $0$-admissible with $\sface$-weight $\alpha_\sface$ and $\sface$-loss $0$ (Definition~\usref{DefSSAlephAdm}).
\end{thm}

\subsection{Smoothing operators and the Nash--Moser theorem}
\label{SsA1NM}

In order to solve nonlinear equations (with the tame estimates of Theorem~\ref{ThmF} as the starting point), we will use a version of the Nash--Moser theorem proved by Saint Raymond \cite{SaintRaymondNashMoser}. Since smoothing operators enlarge supports, a bit of care is needed; we follow \cite[\S{5.2}]{HintzVasyQuasilinearKdS} and \cite[\S{6.1.3}]{HintzGlueLocII} here. We recall \cite[Theorem~1]{SaintRaymondNashMoser} here in the form given in \cite[Theorem~6.11]{HintzGlueLocII}:

\begin{thm}[Nash--Moser]
\label{ThmA1NM}
  Let $(B^k,|\cdot|_k)$ and $(\bfB^k,\|\cdot\|_k)$ be Banach spaces for $k\in\N_0$. Assume that $B^s\hra B^t$ is a continuous inclusion when $s\geq t$, and set $B^\infty:=\bigcap_{k\in\N_0}B^k$; and for $\eta\in[0,1]$, let $B_\eta^k\subset B^k$ be a linear subspace, with $B_\eta^k\subseteq B_{\eta'}^k$ when $\eta\geq\eta'$ and $B_0^k=B^k$; similarly for $\bfB^k$. We require the existence of operators $S_\theta\colon B^\infty_{\theta^{-\frac12}}\to B^\infty_0$, $\theta>1$, which satisfy
  \begin{gather}
  \label{EqA1NMSmooth}
    S_\theta \colon B^\infty_\eta \to B^\infty_{\eta-\theta^{-\frac12}}\quad\forall\,\eta\in[\theta^{-\frac12},1], \\
  \label{EqA1NMSmooth2}
    |S_\theta v|_s \leq C_{s,t}\theta^{s-t}|v|_t\quad\forall\,s\geq t,\qquad
    |v-S_\theta v|_s \leq C_{s,t}\theta^{s-t}|v|_t\quad\forall\,s\leq t.
  \end{gather}
  Let $\delta>0$, $d\in\N$, and suppose $\Phi\colon\{u\in B^\infty\colon|u|_{3 d}<\delta\}\to\bfB^\infty$ is a $\cC^2$ map such that $u\in B_\eta^\infty$ implies $\Phi(u)\in\bfB_\eta^\infty$. Suppose moreover that there exist constants $C_1,C_2$, and $C_k$, $k\geq d$, such that for all $u,v,w\in B^\infty$ with $|u|_{3 d}<\delta$,
  \begin{gather}
  \label{EqA1NMPhi}
    \|\Phi(u)\|_k \leq C_k(1+|u|_{k+d})\quad\forall\,k\geq d, \\
    \|\Phi'(u)v\|_{2 d} \leq C_1|v|_{3 d},\quad
    \|\Phi''(u)(v,w)\|_{2 d} \leq C_2|v|_{3 d}|w|_{3 d}. \nonumber
  \end{gather}
  Suppose, finally, that for all $u\in B^\infty$ with $|u|_{3 d}<\delta$ there exists a linear map $\Psi(u)\colon\bfB^\infty\to B^\infty$ mapping $\bfB^\infty_\eta\to B^\infty_\eta$ for all $\eta\in[0,1]$, satisfying $\Phi'(u)\Psi(u)f=f$ for all $f\in\bfB^\infty$, and the tame estimate
  \begin{equation}
  \label{EqA1NMTame}
    |\Psi(u)f|_k \leq C_k\bigl(\|f\|_{k+d} + |u|_{k+d}\|f\|_{2 d}\bigr)\quad\forall\,k\geq d,\ f\in\bfB^\infty.
  \end{equation}
  Then there exists a constant $\eps=\eps(\delta,\eta_0,\{C_s,\ C_{s,t}\colon s,t\leq 16 d^2+43 d+24\})$ such that if $\Phi(0)\in\bfB^\infty_{\eta_0}$ with $\|\Phi(0)\|_{2 d}<\eps$, there exists $u\in B^\infty$ with $|u|_{3 d}<\delta$ and $\Phi(u)=0$.
\end{thm}

\begin{rmk}[Norms]
\label{RmkA1NMNorms}
  The references require the inequality $|\cdot|_t\leq|\cdot|_s$ of norms for $s\geq t$, similarly for the norms $\|\cdot\|_s$. We do not require this explicitly here: this can always be arranged simply by replacing any given norms $|\cdot|_s$, $s\in\N_0$, by the new norms $\max_{t\leq s}|\cdot|_t$. (This replacement only has the effect of modifying the constants $C_s$ and $C_{s,t}$ in the statement of the theorem.)
\end{rmk}

The simplest setting of interest below is
\begin{equation}
\label{EqA1NMSob}
  B^k=\Hb^{k,(\alpha_\sscri,\alpha_+,0)}(\Omega)^{\bullet,-},\quad
  \bfB^k=\Hb^{k,(\alpha_\sscri+1,\alpha_++2,0)}(\Omega)^{\bullet,-},
\end{equation}
with elements of the subspaces $B^k_\eta$ required to be supported in the subset $\{t_*\geq 1+\eta\}$ of $\Omega$. For the null-form nonlinearity, the space $B^\infty$ must also encode the radiation field at $\scri^+$ and thus includes as a second summand a b-Sobolev space on $\scri^+$. For quasilinear equations, the domain $\Omega_*$ is more appropriate, but it is not a smooth submanifold (with corners) of $\tilde M$ anymore since the function $\ft_*$ defining its initial Cauchy hypersurface only satisfies $\ft_*-t_*\in x_\sscri^{2\ell_\sscri}\CI_\bop$ near $(\scri^+)^\circ$ (see Lemma~\ref{LemmaSDGTime}). We proceed to discuss smoothing operators in the required degree of generality. The setup is as follows.

\begin{itemize}
\item $\sfM$ is a compact manifold with corners; we require its boundary hypersurfaces to be embedded submanifolds.
\item The set of boundary hypersurfaces is $\sH=\{H_1,\ldots,H_N\}$. Consider a decomposition $\sH=\sH_\bop\sqcup\bar\sH$ (with $\sH_\bop$ or $\bar\sH$ allowed to be empty).
\item Fix a smooth positive b-density $\mu$ on $\sfM$. Write $\bar H_\bop^k(\sfM)\subset L^2(\sfM^\circ,\mu)$ for the subspace of $L^2$ consisting of all functions such that $A u\in L^2$ where $A$ is an up to $k$-fold composition of smooth vector fields on $\sfM$ that are tangent to every element of $\sH_\bop$, and equip $\bar H_\bop^k(\sfM)$ with the (complete) norm which is the sum of $L^2$-norms of $A u$ for $A$ ranging over a finite spanning set (over $\CI(\sfM)$) of the space of such operators.
\item Write $\rho_H\in\CI(\sfM)$ for a boundary defining function of $H\in\sH$, and write $\CI_\bop(\sfM)$ for the space of bounded conormal functions on $\sfM^\circ$ (i.e., all $u\in\CI(\sfM^\circ)$ which remain uniformly bounded under application of operators $A$ as in the previous part, for all $k$).
\item Let $\Sigma\subset\sfM$ be a smooth hypersurface such that for all $H\in\sH_\bop$, $\Sigma\cap H$ is either empty or contained in $H^\circ$. Let $t\in\CI(\sfM^\circ)$ be smooth in a neighborhood of $\Sigma$, with $\Sigma=t^{-1}(0)$, and such that $\dd t\neq 0$ on $t^{-1}(0)$. Let moreover $\ft\in\CI(\sfM\setminus\bigcup_{H\in\sH_\bop}H)$ be such that $\dd\ft\neq 0$ on $\sfM^\circ\cap\ft^{-1}(0)$ and $(t-\ft)\chi\in\prod_{H\in\sH_\bop}\rho_H^\eps\CI_\bop(\sfM)$ for some $\eps>0$ where $\chi\in\CI(\sfM)$ is equal to $1$ near $t^{-1}(0)$, and suppose that for all $H\in\bar\sH$, the intersection $\ft^{-1}(0)\cap H$ is either empty or contained in $H^\circ$.
\item Set $\Omega_\eta:=\{\ft\geq\eta\}$, and write
  \[
    \Hb^k(\Omega_\eta)^{\bullet,-} := \{ u\in\bar H_\bop^k(\sfM) \colon \supp u\subset\Omega_\eta \}.
  \]
\end{itemize}

\begin{prop}[Smoothing operators]
\label{PropA1NMSmooth}
  Under the above assumptions, there exists $c_0>0$ such that, upon defining
  \[
    B^k_\eta := \Hb^k(\Omega_{c_0\eta})^{\bullet,-},\quad k\in\N_0\cup\{\infty\},
  \]
  there exist smoothing operators $S_\theta\colon B^\infty_{\theta^{-\frac12}}\to B^\infty_0$, $\theta>1$, which satisfy the properties~\eqref{EqA1NMSmooth}--\eqref{EqA1NMSmooth2}. The analogous statement holds more generally for weighted spaces $B^k_\eta=w\Hb^k(\Omega_{c_0\eta})^{\bullet,-}$ with weights $w=\prod_{H\in\sH_\bop}\rho_H^{\alpha_H}$, $\alpha_H\in\R$.
\end{prop}

\begin{example}[Special cases]
\label{ExA1NMSpecial}
  We may take $\sfM=\cl_M\{t_*\geq\frac12\}$, and then $\sH_\bop=\{\cK^+,\iota^+,\scri^+\cap\sfM\}$ and $\sH=\{\sfM\cap r^{-1}(\bhm),\sfM\cap t_*^{-1}(\frac12)\}$. We moreover take $t$ to be the function $t_*-1$, and $\ft$ to be either $t_*-1$ or $\ft_*-1$ in the notation of Lemma~\ref{LemmaSDGTime}. Then Proposition~\ref{PropA1NMSmooth} provides smoothing operators on the spaces used in Theorem~\ref{ThmF}. If instead one takes $\sfM=\scri^+\cap\{t_*\geq\frac12\}$, $\sH_\bop=\{\scri^+\cap\iota^+\}$, $\sH=\{\scri^+\cap t_*^{-1}(\frac12)\}$, and $t=\ft:=t_*-1$, one gets smoothing operators for (weighted b-Sobolev) functions on $\scri^+$ with support in $t_*\geq 1$.
\end{example}

\begin{proof}[Proof of Proposition~\usref{PropA1NMSmooth}]
  \pfstep{Step~1. Smoothing operators on $\R^n$.} Following \cite[Lemma~6.12]{HintzGlueLocII}, we first construct smoothing operators $S_{\theta,0}\colon\sD'(\R^n)\to\CI(\R^n)$ satisfying~\eqref{EqA1NMSmooth2} with $|\cdot|_s:=\|\cdot\|_{H^s(\R^n)}$ such that $\supp(S_{\theta,0}u)$ is contained in a $\theta^{-\frac12}$-neighborhood of $\supp u$ for all $u\in\sD'(\R^n)$. To this end, we fix $\phi\in\CIc(\R^n)$ with $\phi(\xi)=1$ for $|\xi|<1$, further $\chi:=\cF^{-1}\phi\in\sS(\R^n)$, and then
  \[
    (S_{\theta,0}v)(x) := \int_{\R^n} \phi(\theta^{\frac12}y)\theta^n\chi(\theta y)v(x-y)\,\dd y.
  \]
  Without the localizer $\phi(\theta^{\frac12}y)$, this is a Fourier multiplier by $\phi(\frac{\xi}{\theta})$, and the estimates~\eqref{EqA1NMSmooth2} then follow from $\la\xi\ra^s|\phi(\frac{\xi}{\theta})|\leq C_{s,t}\theta^{s-t}\la\xi\ra^t$ for $s\geq t$ and $\la\xi\ra^s|1-\phi(\frac{\xi}{\theta})|\leq C_{s,t}\theta^{s-t}\la\xi\ra^t$ for $s\leq t$. Inserting the localizer yields an error term $\int_{\R^n}(1-\phi(\theta^{\frac12}y))\theta^n\chi(\theta y)v(x-y)\,\dd y$ whose $x$-derivatives can be written as $y$-derivatives falling on $v$ under the integral sign, so upon integrating by parts we only need to note that since $|\theta y|\geq\max(\theta^{\frac12},|y|)$ on $\supp\phi(\theta^{\frac12}y)$,
  \[
    |\theta^N(1-\phi(\theta^{\frac12}y))\theta^n\chi(\theta y)| \leq C_N \theta^N\cdot C_L\theta^n(1+|\theta y|)^{-L-n-1} \leq C_{N,K}\theta^{-K}\la y\ra^{-n-1}
  \]
  for any desired value of $K$, provided $N+n-2 L\leq-K$; this is thus square integrable in $y$, with $L^2$-norm bounded by $\theta^{-K}$ for any $K$.

  \pfstep{Step~2. Smoothing operators on orthants.} Write $\R^n_l:=[0,\infty)_l\times\R^{n-l}$, $l\geq 1$. We next construct smoothing operators $S_{\theta,l}\colon\bar\sD'(\R^n_l)\to\bar\cC^\infty(\R^n_l)$ which enlarge supports by $\leq\theta^{-\frac12}$ and satisfy the properties~\eqref{EqA1NMSmooth2} with $|\cdot|_s:=\|\cdot\|_{\bar H^s(\R^n_l)}$. We use Seeley's theorem \cite{SeeleyExtension} as in the proof of Lemma~\ref{LemmaETExt} to obtain continuous linear extension operators $E_j\colon\bar\cC^\infty(\R^n_j)\to\bar\cC^\infty(\R^n_{j-1})$, $j=l,l-1,\ldots,1$, which extend to continuous maps $\bar H^s(\R^n_j)\to\bar H^s(\R^n_{j-1})$ and such that $\supp E_j\subset\{x_j\geq-1\}$. We then set $S_{\theta,l}=R\circ S_{\theta,0}\circ E_1\circ\cdots\circ E_l$ where $R\colon\R^n\to\R^n_l$ is the restriction map. The estimates for $S_{\theta,0}$ then imply those for $S_{\theta,l}$.

  \pfstep{Step~3. A good cover of $\Omega_0$.} Near a point in $\ft^{-1}(0)\cap H$, $H\in\bar\sH$, the function $\ft$ is smooth, and we can thus use it as one of the local coordinates; we thus get a chart $(-2 c_0,2 c_0)_\ft\times[0,2)_x\times(-2,2)^{n-2}_y$ for some $c_0>0$, with $H$ in this chart being $x^{-1}(0)$. Near a point in $t^{-1}(0)\cap H$, $H\in\sH_\bop$, on the other hand, we get a chart $(-4 c_0,4 c_0)_t\times[0,2)_x\times(-2,2)^{n-2}_y$, which over $x>0$ we subdivide into the countably many charts
  \begin{equation}
  \label{EqA1NMCharts}
    (-4 c_0,4 c_0)_t \times (-2,2)_{x_{(j)}} \times (-2,2)_y^{n-2},\quad x_{(j)} := -\log_2 x - j,\ \ 3\leq j\in\N,
  \end{equation}
  which we call \emph{b-type charts} in this proof. (These are essentially unit cells for b-vector fields on $\sfM$.) In these charts, any fixed $\CI$-seminorm of $\ft-t$ is bounded by a constant times $2^{-\eps j}$. Thus, for all $j\geq j_0$ where $j_0\geq 3$ is sufficiently large, we can use $\ft$ instead of $t$ as a local coordinate, thus getting charts
  \[
    (-2 c_0,2 c_0)_\ft \times (-2,2)_{x_{(j)}} \times (-2,2)_y^{n-2},\quad j\geq j_0.
  \]
  The complement (in $\sfM\setminus\bigcup_{H\in\sH_\bop}H$) of the charts just discussed intersected with $\Sigma$, which is contained in a compact subset of $\Sigma^\circ$, can be covered by finitely many charts of the form $(-2 c_0,2 c_0)_\ft\times(-2,2)^{n-1}_y$.

  Around points in $\ft>0$, we pick charts in $\{\ft>0\}$ which are standard charts on manifolds with corners, i.e., of the form $[0,2)_x^l\times(-2,2)_y^{n-l}$; for those boundary defining functions $x_j$, $j\in\{1,\ldots,l\}$, which correspond to elements of $\sH_\bop$, we pass to countably many b-type charts by subdividing $-\log_2 x_j$ into intervals of length $4$ as in~\eqref{EqA1NMCharts}. By compactness of $\Omega_{-c_0}$ (for sufficiently small $c_0>0$), we have thus covered $\Omega_{-c_0}$ by finitely many families of b-type charts and finitely many standard charts. See Figure~\ref{FigA1NMCharts}.

  \begin{figure}[!ht]
  \centering
  \includegraphics{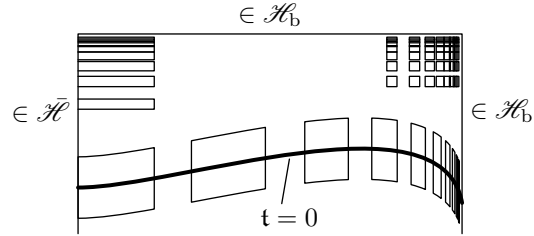}
  \caption{Illustration of the charts constructed in Step~3 of the proof of Proposition~\usref{PropA1NMSmooth}.}
  \label{FigA1NMCharts}
  \end{figure}

  \pfstep{Step~4. Combination.} We can rescale $c_0$ to $1$ by dividing $t$ and $\ft$ by $c_0$. We then fix a partition of unity $\chi_\alpha$ subordinate to the cover by these charts $(\phi_\alpha,U_\alpha)$, with $\phi_\alpha$ being a diffeomorphism from $U_\alpha\subset\sfM$ to $[0,2)^l\times(-2,2)^{n-l}$ with $l$ depending on $\alpha$. Fix another set of cutoff functions $\tilde\chi_\alpha\in\CIc(U_\alpha)$ which equal $1$ on $\supp\chi_\alpha$. Define $S_\theta u=\sum_\alpha \tilde\chi_\alpha\phi_\alpha^* S_\theta((\phi_\alpha)_*(\chi_\alpha u))$, where $S_\theta$ is given by $S_{\theta,l}$ in the chart. That $S_\theta$ satisfies all the requirements follows from the fact that the norm $\|u\|_{\bar H_\bop^k(\Omega_{-1})}$ is equivalent to the square root of $\sum_\alpha \|(\phi_\alpha)_*(\chi_\alpha u)\|_{\bar H^k(\R^n_l)}^2$. For the weighted case, we only need to note that the weight $w$ is approximately constant (i.e., its supremum and infimum are bounded by each other by a universal constant factor) on each chart; setting $w_\alpha:=\sup_{U_\alpha}w$, we then have the norm equivalence $\|u\|_{w\bar H_\bop^k(\Omega_{-1})}\sim(\sum_\alpha w_\alpha^{-1}\|(\phi_\alpha)_*(\chi_\alpha u)\|_{\bar H^k(\R^n_l)}^2)^{\frac12}$.
\end{proof}

The tame estimate~\eqref{EqA1NMPhi} for the ``direct'' operator $\Phi$ is, in applications, typically the consequence of standard Moser-type estimates for products of functions in Sobolev spaces. For the convenience of the reader, we record simple (and far from optimal) such results here, as we will need them in the slightly nonstandard setting of (partially b-)Sobolev spaces on domains with corners. We only consider a single chart, with extensions to more general settings, such as the one introduced above, following from simple partition of unity arguments.

\begin{lemma}[Tame bounds for nonlinear expressions]
\label{LemmaA1Nonlin}
  Let $l,m,n\in\N_0$ and equip the interior of $\cM:=[0,\infty)_\rho^l\times[0,\infty)_x^m\times\R_y^n$ with the volume density $|\frac{\dd\rho_1}{\rho_1}\cdots\frac{\dd\rho_l}{\rho_l}\,\dd x\,\dd y|$. For $k\in\N_0$, write $\bar H_\bop^k(\cM)\subset L^2(\cM)$ for the space of functions $u$ such that $(\rho_1\pa_{\rho_1})^{\alpha_1}\cdots(\rho_l\pa_{\rho_l})^{\alpha_l}\pa_x^\beta\pa_y^\gamma u\in L^2(\cM)$ for all $\alpha=(\alpha_1,\ldots,\alpha_l)\in\N_0^l$, $\beta\in\N_0^m$, $\gamma\in\N_0^n$ with $|\alpha|+|\beta|+|\gamma|\leq k$. Define the norm $\|\cdot\|_k$ as the sum of the $L^2$-norms of $u$ and its derivatives. Let $\N\ni d>\frac{l+m+n}{2}$.
  \begin{enumerate}
  \item\label{ItA1NonlinSob}{\rm (Sobolev embedding.)} There exists a constant $C_{\rm Sob}>0$ such that $\|u\|_{L^\infty}\leq C_{\rm Sob}\|u\|_d$.
  \item\label{ItA1NonlinProd}{\rm (Products.)} For all $k\geq d$, there exists $C_k$ such that for all $u,v\in\bar H_\bop^k(\cM)$, we have
    \begin{equation}
    \label{EqA1NonlinProd}
      u v\in H_\bop^k(\cM),\quad \|u v\|_k \leq C_k\bigl( \|u\|_d\|v\|_k + \|u\|_{k+d}\|v\|_0\bigr).
    \end{equation}
  \item\label{ItA1Nonlin}{\rm Nonlinear expressions.)} Let $\delta>0$, and let $F\colon[-\delta,\delta]\to\R$ be a smooth function with $F(0)=0$. Then for all $k\geq d$ and $u\in\bar H_\bop^k(\cM)$ with $\|u\|_d\leq \delta/C_{\rm Sob}$, we have
    \begin{equation}
    \label{EqA1Nonlin}
      F(u)\in H_\bop^k(\cM),\quad \|F(u)\|_k \leq C_k\|u\|_{k+d}.
    \end{equation}
  \end{enumerate}
\end{lemma}
\begin{proof}
  We follow the ideas of \cite[Lemma~3.33]{HintzMink4Gauge}. Upon passing to logarithmic coordinates $-\log\rho_1$, $\ldots$, $-\log\rho_l$, b-regularity in $\rho_1$ is the same as standard regularity in $-\log\rho_1$, so we may replace $y$ by $(y,-\log\rho_1,\ldots,-\log\rho_l)$; relabeling $n+l$ as $n$, we have thus reduced to the case $l=0$. Next, by applying Seeley extension operators across the hypersurfaces $x^1=0$, $\ldots$, $x^m=0$, we can extend given functions $u,v$ to all of $\R^{m+n}$; once we have proved the lemma for $\cM=\R^{n+m}$, it then follows for $\cM=[0,\infty)_x^m\times\R_y^n$ by restriction. Relabeling $m+n$ as $n$, we have now reduced to the case $l=m=0$.

  We then recall the proofs of the claimed tame estimates on $\R_x^n$; by a density argument, we only need to consider $u,v\in\CIc(\R^n)$. Sobolev embedding is standard. Iterating $\int|\pa u|^2\,\dd x=\int u\ol{\pa^2 u}\,\dd x\leq\|u\|_{L^2}\|\pa^2 u\|_{L^2}$ gives $\|\pa^p u\|_{L^2}\leq C_{p k}\|u\|_{L^2}^{\frac{k-p}{k}}\|\pa^k u\|_{L^2}^{\frac{p}{k}}$ for $0\leq p\leq k$, and thus the $L^2$-norm of $\pa^k(u v)$ can be estimated, after applying the Leibniz rule, using $\|(\pa^p u)(\pa^{k-p}v)\|_{L^2}\leq\|\pa^p u\|_{L^\infty}\|\pa^{k-p}v\|_{L^2}\lesssim\|\pa^p(\pa^{\leq d}u)\|_{L^2}\|\pa^{k-p}v\|_{L^2}\lesssim\|\pa^{\leq d}u\|_{L^2}^{\frac{k-p}{k}}\|\pa^{\leq k+d}u\|_{L^2}^{\frac{p}{k}}\cdot\|v\|_{L^2}^{\frac{p}{k}}\|\pa^k v\|_{L^2}^{\frac{k-p}{k}}$, which gives~\eqref{EqA1NonlinProd}. The proof of~\eqref{EqA1Nonlin} follows similarly by direct differentiation of $F(u)$ and using the $L^\infty$-norm plus Sobolev embedding for all but one factor of (derivatives of) $u$, with this one factor being estimated in $L^2$.
\end{proof}

\subsection{Power nonlinearity}
\label{SsA1Power}

In order to make the weights below more transparent, we use (as in~\eqref{EqFTame20}) an unweighted b-density such as
\begin{equation}
\label{EqA1Powermub}
  \mu_\bop := \frac{1}{t_* r^3}\mu,\quad \mu:=|\dd g|\ \text{(or $\mu=|\dd\ubar g|=|\dd t_*\,\dd x|$)}
\end{equation}
to define $L^2$-spaces on spacetime. The benefit is that $L^2$- and pointwise bounds are compatible in that $\Hb^{\infty,(\beta_\sscri,\beta_+,\beta_\cK)}(\Omega,\mu_\bop)\hra\rho_\sscri^{\beta_\sscri}\rho_+^{\beta_+}\rho_\cK^{\beta_\cK}\cC_\bop^\infty$, while the reverse inclusion holds upon reducing the orders on the $\Hb$-space by an arbitrarily small amount.

We use the notation from~\S\ref{SsA1Adm}. Using only the basic bounds provided by Theorem~\ref{ThmF}, we can now prove:

\begin{thm}[Small data global existence, $p\geq 4$]
\label{ThmA1Power}
  Let $p\geq 4$. Fix $\beta_\sscri=1-\eps_\sscri$ and $\beta_+=1-\eps_+$, where $\eps_\sscri\in(0,\frac12)$ and $0<\eps_+<\min(\frac13,\eps_\sscri)$. Then there exist $d\in\N$ and $\eps>0$ such that for all $f\in\Hb^{\infty,(\beta_\sscri+1,\beta_++2,\frac12)}(\Omega,\mu_\bop)^{\bullet,-}$ with
  \begin{equation}
  \label{EqA1Powerf}
    \|f\|_{\Hb^{d,(\beta_\sscri+1,\beta_++2,\frac12)}}^2 = \iint_\Omega |\rho_\sscri^{-\beta_\sscri-1}\rho_+^{-\beta_+-1}\rho_\cK^{-\frac12}(t_*\pa_{t_*},r\pa_x)^{\leq d}f|^2\,\dd\mu_\bop < \eps^2,
  \end{equation}
  where $\rho_\sscri=\frac{t_*}{t_*+r}$ and $\rho_+=\frac{t_*+r}{t_* r}$ as in~\eqref{EqFWeights}, the equation
  \begin{equation}
  \label{EqA1PowerEq}
    \Box_g u = f + u^p
  \end{equation}
  has a unique forward solution $u\in\Hb^{\infty,(\beta_\sscri,\beta_+,\frac12)}(\Omega)^{\bullet,-}$.
\end{thm}

See~\eqref{EqA1AdmPwu} for the resulting pointwise decay estimate for $u$. Heuristically, the condition $p\geq 4$ arises as follows: if $u$ satisfies the pointwise bounds~\eqref{EqA1AdmPwu}, then $|u^p|\lesssim\rho_\sscri^{p(1-\eps_\sscri)}\rho_+^{p(1-\eps_+)}\rho_\cK^{\frac{p}{2}}$, which is consistent with the bounds~\eqref{EqA1AdmPwf} (ignoring the b-differentiability) provided $p(1-\eps_\sscri)>2-\eps_\sscri$ (so $p\geq 3$), further
\begin{equation}
\label{EqA1PowerHeur}
  p(1-\eps_+)>3-\eps_+, \quad \text{so}\quad p\geq 4,
\end{equation}
and finally $\frac{p}{2}>\frac12$ (so $p\geq 2$). In particular, in order to prove a version of Theorem~\ref{ThmA1Power} for $p=3$, one must improve the $\iota^+$-decay of $u$ to $\rho_+^{1+\delta}$ for some $\delta>0$. This is done in~\S\ref{SsDp}, see Theorem~\ref{ThmDp}.

\begin{proof}[Proof of Theorem~\usref{ThmA1Power}]
  Let $\Phi(u):=\Box_g u-u^p-f$. By shifting $t_*$ by $1$, we may assume that $t_*\geq 2$ on $\supp f$. We claim that $\Phi$ satisfies the hypotheses of Theorem~\ref{ThmA1NM} for the spaces
  \begin{equation}
  \label{EqA1PowerB}
  \begin{alignedat}{2}
    B^k_\eta &:= \{ u\in\Hb^{k,(\beta_\sscri,\beta_+,\frac12)}(\Omega,\mu_\bop)^{\bullet,-} &&\colon t_*\geq 1+\eta\ \text{on}\ \supp u \}, \\
    \bfB^k_\eta &:= \{ f\in\Hb^{k,(\beta_\sscri+1,\beta_++2,\frac12)}(\Omega,\mu_\bop)^{\bullet,-} &&\colon t_*\geq 1+\eta\ \text{on}\ \supp f \},
  \end{alignedat}
  \end{equation}
  where the $\scri^+$-weight is given in terms of powers of $\rho_\sscri=x_\sscri^2$. Following the notation in Theorem~\ref{ThmA1NM}, we denote the norms on these spaces by $|\cdot|_k$ and $\|\cdot\|_k$, respectively.

  Note first that if $u\in B^{k+d}_0$ where $d\geq 3$, then by Sobolev embedding,
  \[
    u \in \rho_\sscri^{\beta_\sscri}\rho_+^{\beta_+}\rho_\cK^{\frac12}\cC_\bop^{k+d-3},
  \]
  which implies
  \begin{equation}
  \label{EqA1Powerpm1}
    u^{p-1}\in\rho_\sscri^{(p-1)\beta_\sscri}\rho_+^{(p-1)\beta_+}\rho_\cK^{\frac{p-1}{2}}\cC_\bop^k,\quad
    u^p \in \Hb^{k,(p\beta_\sscri,p\beta_+,\frac{p}{2})}(\Omega,\mu_\bop)^{\bullet,-}.
  \end{equation}
  Since
  \begin{equation}
  \label{EqA1PowerWeights}
    p\beta_\sscri>\beta_\sscri+1,\quad p\beta_+>\beta_++2
  \end{equation}
  (which are equivalent to $(p-1)(1-\eps_\sscri)\geq 1$ and $(p-1)(1-\eps_+)\geq 2$), we conclude that $u^p\in\bfB_0^k$, which together with $\Box_g u\in\bfB_0^{k+d-2}$ implies $\Phi(u)\in\bfB_0^k$. The tame mapping property~\eqref{EqA1NMPhi} follows from the tame bound $\|u^p\|_k\leq C(|u|_d)|u|_k$, valid for sufficiently large $d$: without weights, this is an instance of Lemma~\ref{LemmaA1Nonlin}, while the weights were just discussed.

  We next verify the properties of the forward solution operator $\Psi(u)$ for the linearized equation $\Phi'(u)v=(\Box_g-p u^{p-1})v$. First of all, $\Psi(u)$ preserves the property of being supported in $\{t_*\geq\eta\}$ for all $\eta\in[0,1]$. Secondly, $P_0:=\Box_g$ is a stationary wave-type operator in the sense of Definition~\ref{DefSSAdm}, as shown in Example~\ref{ExSSAdmBox}. Finally, we note, using~\eqref{EqA1Powerpm1}, that $P:=\Box_g-p u^{p-1}$ is admissible (relative to $P_0$) of class $((0;\infty),(2\ell_\sscri,\ell_+,\ell_\cK))$ provided the conditions $\ell_\sscri+1\leq (p-1)\beta_\sscri$, $2+\ell_+\leq(p-1)\beta_+$, and $\ell_\cK\leq\frac{p-1}{2}$ hold (cf.\ \eqref{EqSDWAdmOp} with $p_0=\tilde p_1=0$); but by~\eqref{EqA1PowerWeights}, we can indeed pick $\ell_\sscri\in(0,\frac12]$, $\ell_+>0$, and $\ell_\cK>0$ that satisfy these conditions. Therefore, Theorem~\ref{ThmF} gives the required tame estimate for $\Psi(u)$.

  Theorem~\ref{ThmA1NM} produces a solution of~\eqref{EqA1PowerEq} with support in $t_*\geq 1$; but by local uniqueness, we in fact have $t_*\geq 2$ on $\supp u$. Shifting $t_*$ back yields the claim.
\end{proof}

\begin{rmk}[Initial value problems]
\label{RmkA1IVP}
  We leave it to the reader to state and prove an analogous result for the initial value problem for $\Box_g u=u^p$ when initial data are given on a Cauchy hypersurface that is, say, equal to a Boyer--Lindquist-$\ft$-level set (see~\eqref{EqTsBL}) for large $r$ and transversal to the future event horizon; the missing ingredient is the solvability of $\Box_g u=u^p$ in an exterior region ($\ft\geq 0$, $t_*\leq 1$, in a region with large $r$), which is essentially standard; cf.\ \cite{HintzMink4Gauge,KadarKehrbergerPhgScatter}, based in parts on \cite{HintzVasyMink4}. The upshot is that one has global existence of solutions for initial data with pointwise $r^{-1+\eps_0}$ and $r^{-2+\eps_0}$ decay for $u$ and its time derivative, respectively, and for a large finite number of their b-derivatives (i.e., derivatives along $r\pa_x$); here $0<\eps_0<\frac12$, and then one needs $\eps_\sscri<\eps_0$ in~\eqref{EqA1Admalphas}.
\end{rmk}

\subsection{A simple quasilinear equation}
\label{SsA1Q}

In order to avoid the need for \emph{ad hoc} conditions on quasilinear perturbations near $\scri^+$ related to the requirements for metrics encoded in Definition~\ref{DefSDGMetric}, we only consider metric perturbations that are supported in spatially compact sets. (We were inspired by \cite{DafermosHolzegelRodnianskiTaylorQuasilinear,DafermosHolzegelRodnianskiTaylorQuasilinear2} to make this simplification. In the Kerr stability problem, the metric itself is the unknown, and it will be shown to satisfy the assumptions at $\scri^+$ throughout the nonlinear iteration scheme.) We introduce the relevant notation: for vector bundles $\cE,\cF\to\cM$, we write
\[
  \CI_\cM(\cE,\cF) \subset \CI(\cE,\cF)
\]
for the space of all smooth maps $F\colon\cE\to\cF$ which respect the fibers, i.e., $F|_{\cE_p}\colon\cE_p\to\cF_p$ for all $p\in\cM$. Fix a cutoff function $\chi\in\CIc(\R_x^3)$. For real-valued smooth functions $u\colon\cD\subset\R^4\to\R$, We then consider metrics of the form
\begin{align*}
  &g(t_*,x,u,\dd u) = g_{\bhm,a}(t_*,x) + \chi(x) h(t_*,x,u,\dd u), \\
  &\qquad h \in \CI_{\ol{\R_{t_*}}\times\R^3}(\ubar\R\oplus\cT^*;S^2\cT^*),\ \ h(t_*,x,0,0)=0.
\end{align*}
where we write $\cT^*\to\ol{\R_{t_*}}\times\R^3$ for a bundle with local frame $\dd t_*,\dd x$ analogously to Definition~\ref{DefCTscPullback}. We write out the $(t_*,x)$-coordinates for clarity here; later on, regarding $(u,\dd u)$ as a section of $\ubar\R\oplus\cT^*$ over $\cD$, we simply write $g(u,\dd u)$. An example of a metric $g$ of this form over the Kerr spacetime manifold is
\[
  g_{\bhm,a} + \chi(x)\biggl(u\,\dd t_*^2 - r^5 u^2\,\dd t_*\,\dd r + \frac{r^7|\dd u|_{g_{\bhm,a}^{-1}}^2}{1+u^2}\,\dd r^2 \biggr).
\]
We work with an unweighted b-density $\mu_\bop$ as in~\eqref{EqA1Powermub}.

\begin{thm}[Small data global existence]
\label{ThmA1Q}
  Fix $\beta_\sscri=1-\eps_\sscri$ and $\beta_+=1-\eps_+$ where $0<\eps_+<\eps_\sscri<\frac12$. Then there exist $d\in\N$ and $\eps>0$ such that for all $f\in\Hb^{\infty,(\beta_\sscri+1,\beta_++2,\frac12)}(\Omega,\mu_\bop)^{\bullet,-}$ with $\|f\|_{\Hb^{d,(\beta_\sscri+1,\beta_++2,\frac12)}}<\eps$, the equation
  \[
    \Box_{g(u,\dd u)}u = f
  \]
  has a unique forward solution $u\in\Hb^{\infty,(\beta_\sscri,\beta_+,\frac12)}(\Omega)^{\bullet,-}$.
\end{thm}

See~\eqref{EqA1Powerf} for the explicit form of the smallness condition on $f$, and~\eqref{EqA1AdmPwu} for the resulting pointwise decay estimates for $u$. Since we only consider spatially compactly supported metric perturbations, the level sets of $t_*$ will be spacelike also with respect to $g(u,\dd u)$ when $u$ is small enough; this is why we can still work with $\Omega$ here.

\begin{proof}[Proof of Theorem~\usref{ThmA1Q}]
  We define $B_\eta^k$ and $\bfB_\eta^k$ by~\eqref{EqA1PowerB} and check the hypotheses of Theorem~\ref{ThmA1NM} for the map $\Phi(u)=\Box_{g(u,\dd u)}u-f$. Note that if the $L^\infty$-norms of $u$ and $\dd u$ are sufficiently small, then $g(u,\dd u)$ is a well-defined Lorentzian metric. For $u\in B^{k+d}_0$ with $|u|_{3 d}$ (where $d$ is some fixed large number) sufficiently small, then, the difference $g_{\bhm,a}-g(u,\dd u)$ is then of class $t_*^{-\frac12}\chi(x)\cC_\bop^{k+d-4}(M)$, where we use the pointwise $t_*^{-\frac12}$ decay from Sobolev embedding (Lemma~\ref{LemmaMUCe3bSob}). It then follows, by writing out $\Box_{g(u,\dd u)}$ in $(t_*,x)$-coordinates and using Lemma~\ref{LemmaA1Nonlin}, that $\Phi(u)\in\bfB^{k+d-4}_0\subset\bfB^k$ provided $d\geq 4$, and the tame estimates~\eqref{EqA1NMPhi} hold. Furthermore, the linearization $\Phi'(u)$ differs from $\Box_{g_{\bhm,a}}$ by an operator which is small in $t_*^{-\frac12}\chi(x)\cC_\bop^{k+d-5}(M)\Diff_\etbop^2(M)$, and thus it is admissible in the sense of Definition~\ref{DefSDWAdm} and of class $((0;k),(2\ell_\sscri,\ell_+,\ell_\cK))$ provided $d\geq 5$, where we can take $\ell_\sscri=\frac12$ and $\ell_+=1$ (these are arbitrary), and $\ell_\cK=\frac12$. Theorem~\ref{ThmF} then shows that the forward solution operator $\Psi(u)$ for $\Phi'(u)$ satisfies the tame estimates~\eqref{EqA1NMTame} (for some sufficiently large $d$). Since $\Psi(u)$ moreover preserves the property of being supported in $t_*\geq\eta$ for $\eta\in[0,1]$, an application of Theorem~\ref{ThmA1NM} (and a local uniqueness argument as in the proof of Theorem~\ref{ThmA1Power}) produces the desired global forward solution $u$.
\end{proof}

\subsection{Null-form nonlinearity}
\label{SsA1Null}

As is well-known since the classical work of Klainerman~\cite{KlainermanNullCondition}, the null condition (or its weak versions as studied, e.g., in \cite{LindbladRodnianskiWeakNull,LindbladRodnianskiGlobalExistence,KeirWeak}) plays an important role in controlling nonlinear interactions near null infinity; it requires that in nonlinear expressions $\pa u\cdot\pa u$, not both derivatives are transversal to the outgoing light cones (i.e., $\pa_t-\pa_r$). Note that near $\scri^+$, the other relevant derivatives,
\[
  \pa_t+\pa_r = r^{-1}\cdot r(\pa_t+\pa_r),\quad
  r^{-1}\pa_\omega = r^{-1}\cdot\pa_\omega,
\]
are $r^{-1}=\rho_\sscri\rho_+$ times b-vector fields on $M$ (cf.\ Example~\ref{ExCTbtildeM}), whereas
\[
  \pa_t-\pa_r = (t-r)^{-1}\cdot(t-r)(\pa_t-\pa_r)
\]
is only $(t-r)^{-1}=\rho_+$ times a b-vector field on $M$, so does not produce an extra power of $\rho_\sscri$ upon application to $u$. If $u$ and its b-derivatives have pointwise $\rho_\sscri^{1-\eps_\sscri}\rho_+^{1-\eps_+}$-decay (cf.\ \eqref{EqA1AdmPwu}), then the decay of $\pa u\cdot\pa u$ where $\pa\in\{\pa_t\pm\pa_r,r^{-1}\pa_\omega\}$, with not both $\pa$ being $\pa_t-\pa_r$, is $\rho_\sscri^{3-2\eps_\sscri}\rho_+^{4-2\eps_+}$, which thus improves over $u$ by (more than) one power of $\rho_\sscri$ and two powers of $\rho_+$. One thus expects the forward solution operator, upon acting on such $\pa u\cdot\pa u$, to output a function with decay rates consistent with those of the original $u$.

Note that according to this heuristic discussion we need not capture the radiation field of $u$ in this analysis; rather, it is again the basic b-tame solvability and regularity theory provided by Theorem~\ref{ThmF} that will (without any post-processing) give the following result; we again work with the unweighted b-density $\mu_\bop$ from~\eqref{EqA1Powermub}.

\begin{thm}[Small data global existence]
\label{ThmA1Null}
  Fix $\beta_\sscri=1-\eps_\sscri$ and $\beta_+=-1-\eps_+$ where $0<\eps_+<\eps_\sscri<\frac12$. Then there exist $d\in\N$ and $\eps>0$ such that for all $f\in\Hb^{\infty,(\beta_\sscri+1,\beta_++2,\frac12)}(\Omega,\mu_\bop)^{\bullet,-}$ with $\|f\|_{\Hb^{d,(\beta_\sscri+1,\beta_++2,\frac12)}}<\eps$, the equation
  \[
    \Box_g u = f + g^{-1}(\dd u,\dd u),\qquad g:=g_{\bhm,a},
  \]
  has a unique forward solution $u\in\Hb^{\infty,(\beta_\sscri,\beta_+,\frac12)}(\Omega,\mu_\bop)^{\bullet,-}$.
\end{thm}
\begin{proof}
  We again use the spaces~\eqref{EqA1PowerB} for the Nash--Moser iteration.

  We first check that $\Phi(u):=\Box_g u-g^{-1}(\dd u,\dd u)-f$ maps $u\in B^\infty_0$ into the space $\bfB^\infty_0$. If we replace $g$ by $\ubar g$, then, according to our discussion above (and recalling that the $\scri^+$-weights in our Sobolev spaces are powers of $x_\sscri=\rho_\sscri^{\frac12}$),
  \begin{equation}
  \label{EqA1NullMink}
  \begin{split}
    \ubar g^{-1}(\dd u,\dd u) &= -|\pa_t u|^2 + |\pa_r u|^2 + r^{-2}|\pa_\omega u|^2 \\
      &= -\bigl((\pa_t-\pa_r)u\bigr)\bigl((\pa_t+\pa_r)u\bigr) + r^{-2}|\pa_\omega u|^2 \\
      & \in \Hb^{\infty,(\beta_\sscri,\beta_++1,\frac12)} \cdot \Hb^{\infty,(\beta_\sscri+1,\beta_++1,\frac12)} + \Hb^{\infty,(\beta_\sscri+1,\beta_++1,\frac12)}\cdot\Hb^{\infty,(\beta_\sscri+1,\beta_++1,\frac12)}.
  \end{split}
  \end{equation}
  The weights for $\CI_\bop$-spaces, to which we pass in each first factor, are equal to those for the $\Hb$-spaces, so we obtain
  \begin{equation}
  \label{EqA1Nullubarg}
    \ubar g^{-1}(\dd u,\dd u) \in \Hb^{\infty,(2\beta_\sscri+1,\ 2\beta_++2,\ 1 )} \subset \bfB^\infty_0;
  \end{equation}
  for the membership in $\bfB^\infty_0$, we use that $2\beta_\sscri+1\geq\beta_\sscri+1$ (i.e., $\beta_\sscri\geq 0$) and $2\beta_++2\geq\beta_++2$ (i.e., $\beta_+\geq 0$). The difference of $g_{\bhm,a}^{-1}$ and $\ubar g^{-1}$ is of class $r^{-1}\CI(M_0;S^2\,\cT)$ by Lemma~\ref{LemmaTsKLMetric}\eqref{ItTsKLMetricSc}. Therefore, $g_{\bhm,a}^{-1}(\dd u,\dd u)-\ubar g^{-1}(\dd u,\dd u)$ is a sum of terms of the schematic form $r^{-1}\pa u\cdot\pa u$ where $\pa\in\{\pa_{t_*},\pa_x\}$, so $\pa\in\rho_+\Vb$. But such terms lie in
  \begin{equation}
  \label{EqA1NullMinkErr}
    r^{-1}\Hb^{\infty,(\beta_\sscri,\beta_++1,\frac12)}\cdot\Hb^{\infty,(\beta_\sscri,\beta_++1,\frac12)} \subset \Hb^{\infty,( 2\beta_\sscri+1,\ 2\beta_++2,\ 1)}
  \end{equation}
  (which is the same space as in~\eqref{EqA1Nullubarg}) and thus again in $\bfB^\infty_0$. The tame estimates for $\Phi(u)$ itself are again a simple consequence of Lemma~\ref{LemmaA1Nonlin}.

  Consider next linearized operator $\Phi'(u)=\Box_g-2 g^{-1}(\dd u,\dd(\cdot))$ for $u\in B^{k+d}_0$; we claim that
  \begin{equation}
  \label{EqA1NullErr}
    \Phi'(u)-\Box_g\in\cC_\bop^{k,(2\ell_\sscri+2,\ell_++2,\ell_\cK)}\Diff_\etbop^1
  \end{equation}
  for $\ell_\sscri=\frac12$, $\ell_+=1-\eps_+$, $\ell_\cK=\frac12$ (cf.\ \eqref{EqSDWAdmOp} with $p_0=\tilde p_1=0$), where we use powers of $x_\sscri$ for the $\scri^+$-weight. By Lemma~\ref{LemmaTsKLMetric}\eqref{ItTsKLMetrice3b}, $g^{-1}(\dd u,\dd(\cdot))$ is of class
  \begin{equation}
  \label{EqA1Null}
    x_\sscri^2\rho_+^2\CI(M;S^2\,\Tetb M) (\dd u,\dd(\cdot)).
  \end{equation}
  Now, $\dd\in\Diff^1_\etbop(M;\ubar\C,\Tetb^*M)\subset\Diff_\bop^1(M;\ubar\C,\Tetb^*M)$, so
  \[
    \dd u\in\Hb^{k+d-1,(2\beta_\sscri,\beta_+,\frac12)}(\Omega,\mu_\bop;\Tetb^*_\Omega M)^{\bullet,-} \subset \cC_\bop^{k+d-4,(2\beta_\sscri,\beta_++2,\frac12)},
  \]
  and therefore~\eqref{EqA1Null} is
  \begin{equation}
  \label{EqA1NullMem}
    g^{-1}(\dd u,\dd(\cdot)) \in x_\sscri^2\rho_+^2\cC_\bop^{k+d-4,(2\beta_\sscri,\beta_+,\frac12)}\Diff_\etbop^1 = \cC_\bop^{k+d-4,(2\beta_\sscri+2,\beta_++2,\frac12)}\Diff_\etbop^1.
  \end{equation}
  This gives~\eqref{EqA1NullErr} provided $\ell_\sscri\leq\beta_\sscri=1-\eps_\sscri$, $\ell_+\leq\beta_+=1-\eps_+$, and $\ell_\cK\leq\frac12$.

  Theorem~\ref{ThmF} now provides the tame estimates for the forward solution operator for $\Psi'(u)$ which, in conjunction with Theorem~\ref{ThmA1NM}, prove the theorem.
\end{proof}

\begin{rmk}[Simplification]
\label{RmkA1NullSimpl}
  Having established~\eqref{EqA1NullMem}, one can plug in $u$ and thus deduce~\eqref{EqA1Nullubarg} directly, without the need for the Minkowskian calculation~\eqref{EqA1NullMink} and error estimate~\eqref{EqA1NullMinkErr}. The underlying reason is that the weight of $g^{-1}$ as a dual edge-b-metric near $\scri^+$ being $x_\sscri^2$ already encodes the relevant null structure (i.e., the absence of terms $(\pa_t-\pa_r)\otimes(\pa_t-\pa_r)$ in the expression for $g^{-1}$).
\end{rmk}

\section{Improving decay}
\label{SD}

The standard procedure to extract stronger asymptotics for a wave-type equation $P u=f$ (where $P$ is as in Theorem~\ref{ThmF}, say) is to replace $P$ by its stationary model $P_0$, so $P_0 u=f-(P-P_0)u$, and inverting the stationary operator $P_0$ using the Fourier transform; the point, roughly speaking, is that $(P-P_0)u$ gains $\ell_\cK>\aleph$ orders of $t_*$-decay (at $\cK^+$) over the a priori decay bounds on $u$, and thus the inversion of $P_0$ should yield control of $u$ modulo errors that gain $\ell_\cK-\aleph>0$ orders of $t_*$-decay over the a priori decay bounds on $u$.

To carry out this procedure, we need to control Fourier transforms of functions in b-Sobolev spaces with $\cK^+$-weights other than $0$ and inverse Fourier transforms of functions which have nontrivial weights at $\sigma=0$ (or, more precisely, at $\zface\in X_\scbtop$ in the notation of~\eqref{EqMUscbtSingle}); we do this in~\S\ref{SsDFT} (see Propositions~\ref{PropDFTFourier} and \ref{PropDFTInv}). In~\S\ref{SsDRes}, we record the resolvent bounds on b-Sobolev spaces which follow from the results in~\S\ref{SSp} (see Propositions~\ref{PropDResLo} and \ref{PropDResHi}). The decay rate of waves at $\scri^+$ imposes an upper bound on the amount of decay one can prove at $\iota^+$ and $\cK^+$, and hence strong decay at $\scri^+$ is a prerequisite for strong decay at $\iota^+\cup\cK^+$; we thus discuss in~\S\ref{SsDscri} how to extract the radiation field and decay at $\scri^+$. In~\S\ref{SsDp}, we then show, as an application, how to use stronger $t_*$-decay estimates to extend the range of exponents $p$ in Theorem~\ref{ThmA1Power} to $p=3$; see Theorem~\ref{ThmDp}.

\begin{rmk}[Nonlinear variants]
\label{RmkDNonlin}
  Such a procedure can often also be carried out at the nonlinear level: once a solution $u$ of, say, the nonlinear equation $\Box_g u=f+u^p$ in~\eqref{EqA1PowerEq} has been proved to exist (Theorem~\ref{ThmA1Power}), but the source term $f$ there in fact has stronger decay than required for~\eqref{EqA1Powerf}, then inversion of the linear operator $\Box_g$ on $f+u^p$ can yield stronger decay for $u$ itself. Such an approach based on spectral theory was implemented by Looi--Xiong \cite{LooiXiongSemilinearAsymp} in some generality (though their results do not cover the case of Kerr spacetimes). Sharpening asymptotics directly at the level of the \emph{linearized} equations (in a Nash--Moser iteration scheme, say) has the significant advantage of leading to stronger global existence results for nonlinear wave equations, as stronger decay means that weaker nonlinearities (e.g., smaller exponents $p$ in power-type nonlinearities) can still be regarded as perturbative.
\end{rmk}

\subsection{Fourier transforms and weighted b-Sobolev spaces}
\label{SsDFT}

We first study the Fourier transform on the line (\S\ref{SssDFTL}) before turning to the spacetime setting (\S\ref{SssDFT}). Our convention for the (inverse) Fourier transform is
\[
  (\cF u)(\sigma) = \hat u(\sigma) = \int_\R e^{i\sigma t}u(t)\,\dd t,\quad
  (\cF^{-1}v)(t) = \check v(t) = \frac{1}{2\pi}\int_\R e^{-i\sigma t}v(\sigma)\,\dd\sigma.
\]

\subsubsection{Fourier transform on the line}
\label{SssDFTL}

For $k\in\N_0$ and $\beta\in\R$, write
\[
  \Hb^{k,\beta}(\ol\R)=\la t\ra^{-\beta}\Hb^k(\ol\R,|\tfrac{\dd t}{\la t\ra}|).
\]
Thus, $\Hb^{0,\beta}(\ol\R)=\la t\ra^{-\beta+\frac12}L^2(\R,|\dd t|)$, and the space $\Hb^{k,\beta}$ tests of $k$-fold regularity relative to $\Hb^{0,\beta}$ with respect to $\la t\ra\pa_t$, or equivalently $(\pa_t,t\pa_t)$. On the Fourier transform side, this means that
\begin{equation}
\label{EqDFTL}
  \cF \colon \Hb^{k,\beta}(\ol\R) \xra{\cong} \bigl\{ v \colon \sigma^j(\sigma\pa_\sigma)^l v\in H^{\beta-\frac12}(\R,|\dd\sigma|)\ \forall\,j,l\in\N_0,\ j+l\leq k \bigr\}.
\end{equation}
The powers of $\sigma$ encode decay as $|\sigma|\to\infty$. The conormal regularity at $\sigma=0$ is more easily captured by relating $H^{\beta-\frac12}(\R_\sigma)$, or rather the spaces $\bar H^{\beta-\frac12}(\pm(0,\infty))$ of restrictions to $\pm(0,\infty)$, with weighted b-Sobolev spaces on the half lines $\pm[0,\infty)$. We first recall some standard results from \cite[Chapter~4]{TaylorPDE1}.

\begin{lemma}[Isomorphism of supported and extendible spaces]
\label{LemmaDFTLIso}
  For $s\in(-\frac12,\frac12)$, the restriction map
  \begin{equation}
  \label{EqDFTLIso}
    \dot H^s([0,\infty)) = \{ u\in H^s(\R) \colon \supp u\subset[0,\infty) \} \xra{(\cdot)|_{(0,\infty)}} \bar H^s((0,\infty)) = \{u|_{(0,\infty)} \colon u\in H^s(\R)\}
  \end{equation}
  is an isomorphism.
\end{lemma}
\begin{proof}
  If $u\in\dot H^s([0,\infty))$ lies in the kernel of~\eqref{EqDFTLIso}, then $\supp u\subseteq\{0\}$, so $u$ is a sum of differentiated $\delta$-distributions at $0$ and therefore does not lie in $H^{-\frac12}(\R)$ unless $u=0$. Thus,~\eqref{EqDFTLIso} is injective for $s\geq-\frac12$. To prove surjectivity for $s\in(-\frac12,\frac12)$, let $u\in\bar H^s((0,\infty))$ and pick $\tilde u\in H^s(\R)$ with $\tilde u|_{(0,\infty)}=u$. Then $u_0(\sigma):=H(\sigma)\tilde u(\sigma)$ (where $H$ is the Heaviside function) lies in $H^s(\R)$ since multiplication by $H(\sigma)$ is continuous on $H^s(\R)$; this follows for $s\in[0,\frac12)$ from \cite[Chapter~4, Proposition~5.3]{TaylorPDE1}, and for $s\in(-\frac12,0]$ by duality.
\end{proof}

\begin{cor}[Relationship between extendible and supported spaces]
\label{CorDFTLRel}
  Let $s\in(-\frac12,\frac12)+j$, $j\in\N_0$. Let $u\in\bar H^s((0,\infty))$. Then $u\in\dot H^s([0,\infty))$ (meaning that $u$ lies in the range of~\eqref{EqDFTLIso}) if and only if
  \begin{equation}
  \label{EqDFTLRel}
    u(0) = u'(0) = \cdots = u^{(j-1)}(0) = 0.
  \end{equation}
\end{cor}
\begin{proof}
  For $j=0$, this is the content of Lemma~\ref{LemmaDFTLIso}. Consider thus $j\geq 1$. By Sobolev embedding, $H^s(\R)\hra\cC^{j-1}(\R)$, so $\bar H^s((0,\infty))\hra\cC^{j-1}([0,\infty))$. Therefore,~\eqref{EqDFTLRel} is a necessary condition. To prove that it is sufficient, we argue by induction. Consider $u_1:=u'\in\bar H^{s-1}((0,\infty))$, which satisfies~\eqref{EqDFTLRel} with $j$ replaced by $j-1$; thus, extension by $0$ gives $u_1\in\dot H^{s-1}([0,\infty))$. Let then
  \begin{equation}
  \label{EqDFTLRelPf}
    \tilde u(\sigma) := \int_0^\sigma u_1(\tau)\,\dd\tau\in\dot H_\loc^s([0,\infty)).
  \end{equation}
  Since $(u-\tilde u)'=u_1-u_1=0$ in the sense of distributions on $(0,\infty)$, we have $u=c+\tilde u$ on $(0,\infty)$ for some $c\in\C$. Evaluating this at $\sigma>0$ and taking the limit as $\sigma\to 0$ gives $0=c+0$, so $c=0$ and hence $u=\tilde u$. In particular, $\tilde u$ lies in $H^s$ near $\infty$, so from~\eqref{EqDFTLRelPf} we obtain $u=\tilde u\in\dot H^s([0,\infty))$.
\end{proof}

Turning to b-Sobolev spaces, we first show:

\begin{lemma}[Standard and b-Sobolev spaces]
\label{LemmaDFTLb}
  For $s>-\frac12$, the inclusion $\CIc((0,\infty))\subset\dot H^s([0,\infty))$ into $\Hb^{s,s-\frac12}([0,\infty),|\frac{\dd\sigma}{\sigma}|)$ extends by continuity to an isomorphism
  \begin{equation}
  \label{EqDFTLbIso}
    \dot H^s([0,\infty)) \cong \Hb^{s,s-\frac12}([0,\infty),|\tfrac{\dd\sigma}{\sigma}|) := \sigma^{s-\frac12}\Hb^s([0,\infty),|\tfrac{\dd\sigma}{\sigma}|).
  \end{equation}
\end{lemma}
\begin{proof}
  This is trivial for $s=0$. Consider the case $s\in\N$. Since $\pa_\sigma=\sigma^{-1}\sigma\pa_\sigma$, we certainly have $\Hb^{s,s-\frac12}([0,\infty),|\frac{\dd\sigma}{\sigma}|)=\Hb^{s,s}([0,\infty),|\dd\sigma|)\hra\dot H^s([0,\infty))$. For the converse, we first claim that
  \begin{equation}
  \label{EqDFTLbHardy}
    k\in\N,\ u\in\dot H^k([0,\infty)) \implies \sigma^{-k}u \in L^2([0,\infty)).
  \end{equation}
  This is a Hardy-type inequality, which for $u\in\CIc((0,\infty))$ follows from
  \begin{align*}
    \sigma^{-k}u(\sigma) &= \frac{1}{\sigma^k}\int_0^\sigma \int_0^{\sigma_1}\cdots\int_0^{\sigma_{k-1}} u^{(k)}(\sigma_k)\,\dd\sigma_k\cdots\dd\sigma_2\,\dd\sigma_1 \\
      &= \int_0^\sigma u^{(k)}(\sigma_k) \biggl(\int_{\sigma_k}^\sigma\cdots\int_{\sigma_2}^\sigma\,\frac{\dd\sigma_1}{\sigma}\cdots\frac{\dd\sigma_{k-1}}{\sigma}\biggr)\,\frac{\dd\sigma_k}{\sigma} \\
      &= \int_0^1 u^{(k)}(\sigma z_k) \biggl( \int_{z_k}^1 \cdots \int_{z_2}^1\,\dd z_1\cdots\dd z_{k-1}\biggr)\,\dd z_k \\
      &= \int_0^1 u^{(k)}(\sigma t)\frac{(1-t)^{k-1}}{(k-1)!}\,\dd t,
  \end{align*}
  and thus
  \begin{align*}
    \|\sigma^{-k}u\|_{L^2([0,\infty))} &\leq \int_0^1 \|u^{(k)}(t\cdot)\|_{L^2([0,\infty))}\frac{(1-t)^{k-1}}{(k-1)!}\,\dd t \\
      &= \|u^{(k)}\|_{L^2([0,\infty))} \int_0^1 t^{-\frac12}\frac{(1-t)^{k-1}}{(k-1)!}\,\dd t \\
      &\leq C_k\|u\|_{\dot H^k([0,\infty))}.
  \end{align*}
  Given $u\in\dot H^s([0,\infty))$, we conclude from~\eqref{EqDFTLbHardy} that
  \[
    \sigma^{-s}u,\ \sigma^{-(s-1)}\pa_\sigma u,\ \cdots,\ \sigma^{-1}\pa_\sigma^{s-1}u,\ \pa_\sigma^s u \in L^2([0,\infty)),
  \]
  i.e., $\sigma^j\pa_\sigma^j u\in\sigma^j L^2([0,\infty))$ for $j=0,\ldots,s$, which is equivalent to $u\in\sigma^s\Hb^s([0,\infty),|\dd\sigma|)$, as claimed. For real $s\geq 0$, the claim now follows by interpolation. For $s\in(-\frac12,0)$, we use Lemma~\ref{LemmaDFTLIso} to deduce
  \[
    \dot H^s([0,\infty)) = \bigl( \bar H^{-s}((0,\infty)) \bigr)^* = \bigl(\dot H^{-s}([0,\infty)) \bigr)^* = \bigl(\Hb^{-s,-s}([0,\infty),|\dd\sigma|)\bigr)^* = \Hb^{s,s}([0,\infty),|\dd\sigma|).
  \]
  This completes the proof.
\end{proof}

For $s\in(-\frac12,\frac12)$, a direct proof that avoids complex interpolation is given in \cite[Lemma~3.12]{LiInternalCorners}.

In order to capture the Taylor coefficients of $u\in\bar H^s((0,\infty))$ at $0$, we introduce the space
\begin{equation}
\label{EqDFTLHbphg}
  \Hb^{s,((0,0),\beta)}([0,\infty)) := \Biggl\{ \chi(\sigma)\sum_{\ell=0}^{\lfloor s-\frac12\rfloor} \sigma^\ell u_\ell + \tilde u \colon u_\ell\in\C,\ \tilde u\in\Hb^{s,\beta}([0,\infty)) \Biggr\},
\end{equation}
where $\chi\in\CIc((-2,2))$ equals $1$ on $[-1,1]$. Note that the sum is empty for $s<\frac12$.

\begin{cor}[Extendible spaces and b-Sobolev spaces]
\label{CorDFTLbExt}
  Let $s\in(-\frac12,\frac12)+j$, $j\in\N_0$. Then we have an isomorphism
  \[
    \bar H^s((0,\infty)) \cong \Hb^{s,((0,0),s-\frac12)}([0,\infty),|\tfrac{\dd\sigma}{\sigma}|).
  \]
\end{cor}
\begin{proof}
  Elements of $\Hb^{s,((0,0),s-\frac12)}([0,\infty))$ are sums of $\CIc$-functions and elements of $\Hb^{s,s-\frac12}([0,\infty))\cong\dot H^s([0,\infty))$, and thus define elements of $\bar H^s((0,\infty))$. Conversely, if $u\in\bar H^s((0,\infty))$, then Corollary~\ref{CorDFTLRel} shows that
  \[
    u(\sigma) - \sum_{\ell=0}^{\lfloor s-\frac12\rfloor} u^{(\ell)}(0)\frac{\sigma^\ell}{\ell!} \in \dot H^s([0,\infty)),
  \]
  and then Lemma~\ref{LemmaDFTLb} finishes the proof.
\end{proof}

Carefully note that the weight in the b-Sobolev spaces in~\eqref{EqDFTLbIso} and \eqref{EqDFTLHbphg} is a weight both at $\sigma=0$ and $\sigma=\infty$. As an immediate consequence of~\eqref{EqDFTL} and Corollary~\ref{CorDFTLbExt}, we deduce:

\begin{lemma}[Fourier transform on the line]
\label{LemmaDFTLine}
  Let $\beta\in(0,\infty)\setminus\N$. Then the Fourier transform defines an isomorphism
  \begin{align*}
    \Hb^{k,\beta}(\ol\R,|\tfrac{\dd t}{\la t\ra}|) &\to \bigl\{ (v_-,v_+) \colon \la\sigma\ra^j(\sigma\pa_\sigma)^l v_\pm \in \Hb^{\beta-\frac12,((0,0),\beta-1)}(\pm[0,\infty),|\tfrac{\dd\sigma}{\sigma}|) \\
      &\quad \hspace{7em} \forall\,j,l\in\N_0,\ j+l\leq k;\ \text{and}\ v_-^{(l)}(0)=v_+^{(l)}(0)\ \forall\,0\leq l<\beta-1 \bigr\}.
  \end{align*}
\end{lemma}

In our applications, the weight $\beta$ will lie in some compact subset of $\R$ (and $k$ will be arbitrarily large), and we shall therefore not keep careful track of the b-regularity order but allow for bounded losses.\footnote{When we also take spatial b-regularity into account in~\eqref{SssDFT}, it is moreover technically convenient to only work with integer b-regularity orders.} (The sharpness of the $\sigma$-weight at $\sigma=0$ is crucial, however.) Upon defining the space
\begin{equation}
\label{EqDFTHbplus}
  H_{\bop^+}^{k,((0,0),\beta)}([0,\infty)) := \{ v \colon \sigma^j(\sigma\pa_\sigma)^l v\in \Hb^{0,((0,0),\beta)}([0,\infty))\ \forall\,j,l\in\N_0,\ j+l\leq k \},
\end{equation}
which puts powers of $\sigma$ on the same footing as b-derivatives, we thus deduce:

\begin{cor}[Fourier transform: lossy version]
\label{CorDFTLossy}
  Let $\beta\in(0,\infty)\setminus\N$, and define $\ell_\pm(\beta):=\lceil(\beta-\frac12)_\pm\rceil$. Let $k\in\N_0$. Then
  \begin{align*}
    \cF &\colon \Hb^{k+\ell_-(\beta),\beta}(\ol\R,|\tfrac{\dd t}{\la t\ra}|) \to \cX^{k,\beta-1} := \bigl\{ (v_-,v_+) \colon v_\pm\in H_{\bop^+}^{k,((0,0),\beta-1)}(\pm[0,\infty),|\tfrac{\dd\sigma}{\sigma}|), \\
      &\quad \hspace{20em} v_-^{(l)}(0)=v_+^{(l)}(0)\ \forall\,0\leq l<\beta-1 \bigr\} \\
    \cF^{-1} &\colon \cX^{k+\ell_+(\beta),\beta-1} \to \Hb^{k,\beta}(\ol\R,|\tfrac{\dd t}{\la t\ra}|).
  \end{align*}
\end{cor}

\subsubsection{Fourier transform on resolved product spaces}
\label{SssDFT}

\begin{figure}[!ht]
\centering
\includegraphics{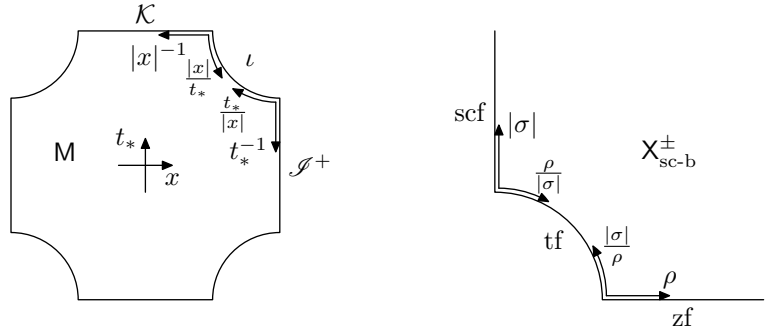}
\caption{\textit{On the left:} the space $\sfM$ (for $n=1$), its boundary hypersurfaces, and some local boundary defining functions. \textit{On the right:} the space $\sfX_\scbtop^+$, its boundary hypersurfaces, and some local boundary defining functions; here $\rho=|x|^{-1}$.}
\label{FigDFTSpaces}
\end{figure}

For $n\in\N$, define
\begin{equation}
\label{EqDFTSpace}
  \sfM := [\,\ol{\R_{t_*}}\times\sfX; \{-\infty,\,+\infty\}\times\pa\sfX\,],\quad \sfX := \ol{\R_x^n}.
\end{equation}
We write $\rho_\sscri$, $\rho_\iota$, $\rho_\cK$ for defining functions of $\scri$ (the lift of $\ol\R\times\pa\sfX$), the front face $\iota$ (i.e., the lift of $\pa\ol\R\times\pa\sfX$), and $\cK$ (the lift of $\pa\ol\R\times\sfX$), respectively. (Note that $\iota$ and $\cK$ have two connected components.) Possible choices are
\begin{equation}
\label{EqDFTbdf}
  \rho_\sscri = \frac{\la t_*\ra}{\la t_*\ra+\la x\ra},\quad
  \rho_\iota = \frac{\la t_*\ra+\la x\ra}{\la t_*\ra\la x\ra}, \quad
  \rho_\cK = \frac{\la x\ra}{\la t_*\ra+\la x\ra}.
\end{equation}
We fix smooth unweighted b-densities on $\ol\R$ and $\sfX$ such as $|\frac{\dd t_*}{\la t_*\ra}|$ and $|\frac{\dd x}{\la x\ra^n}|$ to define $L^2$-spaces. We then denote weighted b-Sobolev spaces by
\begin{equation}
\label{EqDFTHbM}
  \Hb^{k,(\beta_\sscri,\beta_\iota,\beta_\cK)}(\sfM) := \rho_\sscri^{\beta_\sscri}\rho_\iota^{\beta_\iota}\rho_\cK^{\beta_\cK}\Hb^k(\sfM);
\end{equation}
these test for $k$-fold regularity with respect to the vector fields $\la t_*\ra\pa_{t_*}$ (or equivalently $\pa_{t_*}$ and $t_*\pa_{t_*}$) and $\la x\ra\pa_x$. See Figure~\ref{FigDFTSpaces}.

We wish to study the Fourier transform $\hat u(\sigma,x)=\int e^{i\sigma t_*}u(t_*,x)\,\dd t_*$ of elements of such spaces; at low frequencies $\sigma$, this will require working on the resolved space
\begin{equation}
\label{EqDFTXscbtpm}
  \sfX_\scbtop^\pm = \bigl[\pm[0,1]_\sigma \times \sfX; \{0\}\times\pa\sfX\bigr],
\end{equation}
whose boundary hypersurfaces are the \emph{scattering face} $\scface$ (the lift of $\pm[0,1]\times\pa\sfX$), the \emph{transition face} $\tface$ (i.e., the lift of $\{0\}\times\pa\sfX$), and the \emph{zero face} $\zface$ (the lift of $\{0\}\times\sfX$). (The boundary hypersurface at $\sigma=\pm 1$ is arbitrarily chosen and irrelevant.) The terminology ``$\scbtop$'' is taken from~\eqref{EqMUscbtSingle}, but we emphasize that we are presently concerned exclusively with b-regularity. (This manifold was denoted $\sfX_\res^\pm$ in \cite{HintzPrice}.) We work with unweighted b-densities on $[0,1]$ and $\sfX$ (and thus on $\sfX_\scbtop^\pm$) to define $L^2$-spaces, so on $\sfX_\scbtop^\pm$ we fix the (unweighted b-)density $|\frac{\dd\sigma}{|\sigma|}\,\frac{\dd x}{\la x\ra^n}|$. Setting $\rho:=\la x\ra^{-1}$, we may take as boundary defining functions of $\sfX_\scbtop^\pm$ the functions
\begin{equation}
\label{EqDFTbdfs}
  \rho_\scface = \frac{\rho}{\rho+|\sigma|},\quad
  \rho_\tface = |\sigma|+\rho,\quad
  \rho_\zface = \frac{|\sigma|}{\rho+|\sigma|};
\end{equation}
and we consider weighted b-Sobolev spaces
\begin{equation}
\label{EqDFTHbX}
  \Hb^{k,(\gamma_\scface,\gamma_\tface,\gamma_\zface)}(\sfX_\scbtop^\pm) := \rho_\scface^{\gamma_\scface}\rho_\tface^{\gamma_\tface}\rho_\zface^{\gamma_\zface}\Hb^k(\sfX_\scbtop^\pm).
\end{equation}
This entails b-regularity also in $\sigma$; thus, these spaces test for $k$-fold regularity with respect to $\sigma\pa_\sigma$ and $\la x\ra\pa_x$, which span $\Vb(\sfX_\scbtop^\pm)$ over $\CI(\sfX_\scbtop^\pm)$. In order to capture partial expansions at $\zface$, we define moreover the space
\begin{equation}
\label{EqDFTHbXphg}
\begin{split}
  \Hb^{k,\bigl(\gamma_\scface,\ \gamma_\tface,\ ((\zeta,0),\gamma_\zface)\bigr)}(\sfX_\scbtop^\pm) = \Biggl\{ u(\sigma,x) &= \chi_\zface\sum_{\ell=0}^{\lceil\gamma_\zface-\zeta\rceil} \rho_\zface^{\zeta+\ell}u_\ell(x) + \tilde u \colon \\
    &\quad u_\ell\in\Hb^{k,\gamma_\tface}(\sfX),\ \tilde u\in\Hb^{k,(\gamma_\scface,\,\gamma_\tface,\,\gamma_\zface)}(\sfX_\scbtop^\pm) \Biggr\},
\end{split}
\end{equation}
where $\chi_\zface\in\CI(\sfX_\scbtop^\pm)$ equals $1$ near $\zface$ and vanishes near $\scface$. In local coordinates $\hat r:=\frac{|\sigma|}{\rho}$ (which is a local defining function of $\zface$), $\rho=\la x\ra^{-1}$, and $\omega\in\Sph^{n-1}$ near $\zface\cap\tface$, elements of this space are thus the sum of two terms: a polynomial in $\hat r$ with values in $\rho^{\gamma_\tface}\Hb^k([0,1)_\rho\times\Sph^{n-1})$ and a remainder of class $\rho^{\gamma_\tface}\hat r^{\gamma_\zface}\Hb^k([0,1)_\rho\times[0,1)_{\hat r}\times\Sph^{n-1})$.

The spaces~\eqref{EqDFTHbM} and \eqref{EqDFTHbX} have a simple relationship with product-type spaces:

\begin{lemma}[Product-type spaces]
\label{LemmaDFTProd}
  Let $k\in\N_0$ and $\beta,\gamma\in\R$. Then
  \begin{align*}
    \Hb^{k,(\gamma,\beta+\gamma,\beta)}(\sfM) &= \bigcap_{j=0}^k \Hb^{j,\beta}\bigl(\ol\R;\Hb^{k-j,\gamma}(\sfX)\bigr), \\
    \Hb^{k,(\gamma,\beta+\gamma,\beta)}(\sfX_\scbtop^\pm) &= \bigcap_{j=0}^k \Hb^{j,\beta}\bigl(\pm[0,1];\Hb^{k-j,\gamma}(\sfX)\bigr).
  \end{align*}
\end{lemma}
\begin{proof}
  Having $k$ degrees of b-regularity on $\sfM$ is equivalent to having $j$ and $k-j$ degrees of $\la t_*\ra\pa_{t_*}$- and $\la x\ra\pa_x$-regularity on $\ol\R$ and $\sfX$, respectively, for $j=0,\ldots,k$; this gives the first identity for $\beta=\gamma=0$. For the weights, we note that
  \[
    \rho_\sscri^\gamma\rho_\iota^{\beta+\gamma}\rho_\cK^\beta = \la t_*\ra^{-\beta}\la x\ra^{-\gamma}
  \]
  for the boundary defining functions~\eqref{EqDFTbdf}. The proof of the second identity is similar.
\end{proof}

At high frequencies, we can work entirely with product spaces. Analogously to~\eqref{EqSpHiNormb}, we define
\begin{equation}
\label{EqDFTSobHi}
\begin{split}
  H_{\bop^+}^{k,\beta,\gamma}(\pm[\tfrac12,\infty)\times\sfX) &:= \bigl\{ v \colon \sigma^j(\sigma\pa_\sigma)^l v \in |\sigma|^{-\gamma}L^2\bigl(\pm[\tfrac12,\infty),|\tfrac{\dd\sigma}{\sigma}|; \Hb^{k-(j+l),\beta}(\sfX) \bigr) \\
    &\quad \hspace{16em} \forall\,j,l\in\N_0,\ j+l\leq k \bigr\}.
\end{split}
\end{equation}

\begin{prop}[Fourier transform on b-Sobolev spaces on $\sfM$]
\label{PropDFTFourier}
  Let $\beta_\sscri,\beta_\iota,\beta_\cK\in\R$, and let $u\in\Hb^{k+\ell,(\beta_\sscri,\beta_\iota,\beta_\cK)}(\sfM)$ where $\ell\in\N_0$ is specified later.
  \begin{enumerate}
  \item\label{ItDFTFourierLo}{\rm (Low energies.)} Suppose that $\beta_\cK\in(0,\infty)\setminus\N$ and
    \begin{equation}
    \label{EqDFTFourierScriDecay}
      \beta_\sscri\geq\beta_\iota-\min(\beta_\cK,1-\eps)
    \end{equation}
    for some $\eps>0$; take $\ell=1$. Then
    \[
      (\cF u)|_{\sigma\in\pm[0,1]}\in\Hb^{k,(\min(\beta_\sscri,\beta_\iota)-\eps,\ \beta_\iota-1,\ ((0,0),\beta_\cK-1))}(\sfX_\scbtop^\pm).
    \]
    The restrictions of $\pa_\sigma^j\cF u$, $j=0,\ldots,\lfloor\beta_\cK-1\rfloor$, to $\zface\subset\sfX_\scbtop^\pm$ are independent of the choice of sign.
  \item\label{ItDFTFourierHi}{\rm (High energies.)} For $\ell=(\frac12-\min(\beta_\cK,\beta_\iota-\beta_\sscri))_+$, we have
    \[
      (\cF u)|_{\sigma\in\pm[\frac12,\infty)}\in H_{\bop^+}^{k,\beta_\sscri,0}(\pm[\tfrac12,\infty)_\sigma\times\sfX).
    \]
  \end{enumerate}
\end{prop}
\begin{proof}
  \pfstep{Part~\eqref{ItDFTFourierLo}.} Consider first the case $\beta_\cK\in(0,1)$. Since $\beta_\sscri\geq\beta_\iota-\beta_\cK$, we have
  \[
    u \in \bigcap_{j=0}^k \Hb^{j+1,\beta_\cK}\bigl(\ol\R;\Hb^{k-j,\beta_\iota-\beta_\cK}(\sfX)\bigr)
  \]
  by Lemma~\ref{LemmaDFTProd}. (We record the excess degree of b-regularity in the $t_*$-variable only.) Corollary~\ref{CorDFTLossy} implies that
  \[
    \cF u \in \bigcap_{j=0}^k H_{\bop^+}^{j,\beta_\cK-1}\bigl(\pm[0,\infty);\Hb^{k-j,\beta_\iota-\beta_\cK}(\sfX)\bigr),
  \]
  where we use the notation~\eqref{EqDFTHbplus}. For $|\sigma|\leq 1$, where the $|\sigma|$-decay captured by the $H_{\bop^+}$-space is vacuous, this gives
  \[
    \cF u \in \Hb^{k,(\beta_\iota-\beta_\cK,\beta_\iota-1,\beta_\cK-1)}(\sfX_\scbtop^\pm)
  \]
  by Lemma~\ref{LemmaDFTProd}. One can improve the $\scface$-order by noting that
  \[
    u\in\bigcap_{j=0}^k \Hb^{j+1,\beta_\cK-\eta}\bigl(\ol\R;\Hb^{k-j,\beta_\iota-\beta_\cK+\eta}(\sfX)\bigr)
  \]
  for all $\eta\geq 0$ such that $\beta_\sscri\geq\beta_\iota-\beta_\cK+\eta$. An almost-maximal $\eta$ such that still $\beta_\cK-\eta>0$ is thus $\eta=\min(\beta_\cK,\beta_\sscri-\beta_\iota+\beta_\cK)-\eps$, and then $\beta_\iota-\beta_\cK+\eta=\min(\beta_\iota,\beta_\sscri)-\eps$ and $\beta_\cK-\eta-1>-1$.

  Consider next the case $\beta_\cK\in(b,b+1)$, $b\in\N$. We apply what we have already shown to $u$ regarded as an element of $\Hb^{k+1,(\beta_\sscri,\beta_\iota,\beta_\cK-\eta)}$ where we give up $\eta:=\beta_\cK-1+\eps$ orders of $\cK$-decay; thus, $u$ lies a fortiori in the space $\bigcap_{j=0}^k\Hb^{j+1,\beta_\cK-\eta}(\ol\R;\Hb^{k-j,\beta_\iota-\beta_\cK+\eta}(\sfX))$, where we use that $\beta_\iota-(\beta_\cK-\eta)=\beta_\iota-(1-\eps)\leq\beta_\sscri$. Note that $\beta_\cK-\eta\in(0,1)$; the Fourier transform therefore satisfies
  \begin{equation}
  \label{EqDFTFourierEta}
    \cF u \in \Hb^{k,(\beta_\iota-\beta_\cK+\eta,\beta_\iota-1,\beta_\cK-\eta-1)}(\sfX_\scbtop^\pm)
  \end{equation}
  On the other hand, we observe that $t_*^b u\in\Hb^{k+1,(\beta_\sscri,\beta_\iota-b,\beta_\cK-b)}(\sfM)$, so
  \begin{equation}
  \label{EqDFTFourierODE}
  \begin{split}
    \sigma^b\pa_\sigma^b(\cF u)(\sigma) = (i\sigma)^b \cF(t_*^b u)(\sigma) &\in |\sigma|^b\Hb^{k,(\beta_\iota-\beta_\cK,\beta_\iota-b-1,\beta_\cK-b-1)}(\sfX_\scbtop^\pm) \\
      &= \Hb^{k,(\beta_\iota-\beta_\cK,\beta_\iota-1,\beta_\cK-1)}(\sfX_\scbtop^\pm).
  \end{split}
  \end{equation}
  The $\tface$-orders of this and~\eqref{EqDFTFourierEta} are both $\beta_\iota-1$. Note that $\sigma^b\pa_\sigma^b=\prod_{d=0}^{b-1}(\sigma\pa_\sigma-d)=\prod_{d=0}^{b-1}(\hat r\pa_{\hat r}-d)$, where we passed to the coordinates $\hat r=\frac{\sigma}{\rho}$ and $x\in\ol{\R^n}$ near $\zface$; thus, the membership~\eqref{EqDFTFourierODE} implies an ODE for $\cF u$ in $\hat r\in[0,1]$, with values in the Hilbert space $\Hb^{0,\beta_\iota-1}(\sfX)$. We integrate this towards $\hat r=0$ to extract the degree $b$ Taylor expansion of $\cF u$ at $\hat r^{-1}(0)=\zface$, with a $\rho_\tface^{\beta_\iota-1}\rho_\zface^{\beta_\cK-1}\Hb^k$ remainder. This proves the first statement of part~\eqref{ItDFTFourierLo}. The second statement is most easily checked over compact subsets of $\sfX^\circ$ and using that the Fourier transform of $u\in\Hb^{k+1,\beta_\cK}(\ol\R;L^2_\loc(\sfX^\circ))$ is $\cC^k$-regular in $\sigma\in\R$ (with values on $L^2_\loc(\sfX^\circ)$).

  \pfstep{Part~\eqref{ItDFTFourierHi}.} Setting $\beta'_\cK:=\min(\beta_\cK,\beta_\iota-\beta_\sscri)$, we have
  \[
    u \in \bigcap_{j=0}^k \Hb^{j+\ell,\beta'_\cK}\bigl(\ol\R;\Hb^{k-j,\beta_\sscri}(\sfX)\bigr).
  \]
  If $\ell+(\beta'_\cK-\frac12)\geq 0$, then by~\eqref{EqDFTL} (and using the additional $\ell$ degrees of b-regularity to balance the possibly negative Sobolev regularity order $\beta'_\cK-\frac12$), we have $\cF u\in H_{\bop^+}^{k,\beta_\sscri,0}(\pm[\frac12,\infty)\times\sfX)$.
\end{proof}

In preparation for an analogous result for the inverse Fourier transform, we need:

\begin{lemma}[Bounds for $t_*$-independent functions]
\label{LemmaDFTStat}
  Let $\beta\in\R$ and $v_0\in\Hb^{k,\beta}(\sfX)$. Then the function $v(t_*,x):=v_0(x)$ satisfies
  \[
    v \in \bigcap_{\eps>0} \Hb^{k,(\beta-\eps,\beta,-\eps)}(\sfM).
  \]
\end{lemma}
\begin{proof}
  Upon division by $\rho^\beta$, we may reduce to the case $\beta=0$. Since $\la t_*\ra\pa_{t_*}$-derivatives of $v$ vanish and $\la x\ra\pa_x$-derivatives only fall on $v_0$, it moreover suffices to consider the case $k=0$. We work in local coordinates $\rho=\la x\ra^{-1}$ and $\tau=\la t_*\ra^{-1}$ near $\pa\sfX$ and $\pa\ol\R$, respectively, and omit the spherical variables $\frac{x}{|x|}$; we use the densities $|\frac{\dd\rho}{\rho}|$ and $|\frac{\dd\tau}{\tau}|$ near $\pa\sfX$ and $\pa\ol\R$, respectively.

  Near a connected component of $\scri\cap\iota$ then, we use the local coordinates $\rho_\sscri=\frac{\rho}{\tau}$ and $\rho_\iota=\tau$ and need to estimate
  \[
    \int_0^1 \int_0^1 |\rho_\sscri^\eps v_0(\rho_\sscri\rho_\iota)|^2\,\frac{\dd\rho_\sscri}{\rho_\sscri}\,\frac{\dd\rho_\iota}{\rho_\iota} = \int_0^1 \rho^{2\eps}|v_0(\rho)|^2 \biggl(\int_\rho^1 \tau^{-2\eps} \frac{\dd\tau}{\tau}\biggr)\,\frac{\dd\rho}{\rho}.
  \]
  The inner integral is $\leq C_\eps\rho^{-2\eps}$, so overall this is bounded by $C_\eps\|v_0\|_{L^2}^2$. Near a connected component of $\iota\cap\cK$ on the other hand, we use the local coordinates $\rho_\iota=\rho$ and $\rho_\cK=\frac{\tau}{\rho}$ and estimate
  \[
    \int_0^1 \int_0^1 |\rho_\cK^\eps v_0(\rho_\iota)|^2\,\frac{\dd\rho_\iota}{\rho_\iota}\,\frac{\dd\rho_\cK}{\rho_\cK} = C_\eps\|v_0\|_{L^2}^2,
  \]
  as claimed.
\end{proof}

\begin{prop}[Inverse Fourier transforms on b-Sobolev spaces]
\label{PropDFTInv}
  Let $\beta_\cK\in(0,\infty)\setminus\N$ and $\beta_\iota,\beta_\scface,\beta_\sscri,\delta\in\R$. Let $\chi=\chi(\sigma)\in\CIc((-1,1))$ be equal to $1$ near $[-\frac12,\frac12]$. Then:
  \begin{enumerate}
  \item\label{ItDFTInvLo}{\rm (Low energies.)} Set $\ell:=\lceil(\max(\beta_\iota-\beta_\scface,\beta_\cK,0)+\eps-\frac12)_+\rceil$. If
    \[
      v_\pm\in\Hb^{k+\ell,\ \bigl(\beta_\scface,\ \beta_\iota-1,\ ((0,0),\beta_\cK-1)\bigr)}(\sfX_\scbtop^\pm)
    \]
    and $\pa_\sigma^j v_-=\pa_\sigma^j v_+$ at $\zface$ for $j=0,\ldots,\lfloor\beta_\cK-1\rfloor$, then for $v:=\chi(\sigma)\sum_\pm H(\pm\sigma)v_\pm$ we have
    \[
      \cF^{-1}v \in \bigcap_{\eps>0}\Hb^{k,(\min(\beta_\scface,\beta_\iota)-\eps,\beta_\iota,\beta_\cK)}(\sfM).
    \]
  \item\label{ItDFTInvHi}{\rm (High energies.)} Set $\ell'=\max(\lceil\delta_-+\eps,\beta_\iota-\beta_\sscri-\delta,\beta_\cK-\delta)$ and $\ell=\lceil 2\ell'+\delta\rceil$. If $v_\pm\in H_{\bop^+}^{k+\ell,\beta_\sscri,\delta-1}(\pm[\frac12,\infty)_\sigma\times\sfX)$, then
    \begin{equation}
    \label{EqDFTInvHi}
      \cF^{-1}\bigl((1-\chi(\sigma))v_\pm\bigr) \in \Hb^{k,(\beta_\sscri,\beta_\iota,\beta_\cK)}(\sfM).
    \end{equation}
  \end{enumerate}
\end{prop}

The idea behind part~\eqref{ItDFTInvHi} is that any $|\sigma|$-growth (the order $\delta-1$) can be counterbalanced by using $\bop^+$-regularity, which entails $|\sigma|$-\emph{decay}; and moreover, high b-regularity implies strong decay as $|t_*|\to\infty$ and thus at $\iota\cup\cK$. This is why we can ensure for the $\iota$- and $\cK$-orders in~\eqref{EqDFTInvHi} to be any desired fixed values.

\begin{figure}[!ht]
\centering
\includegraphics{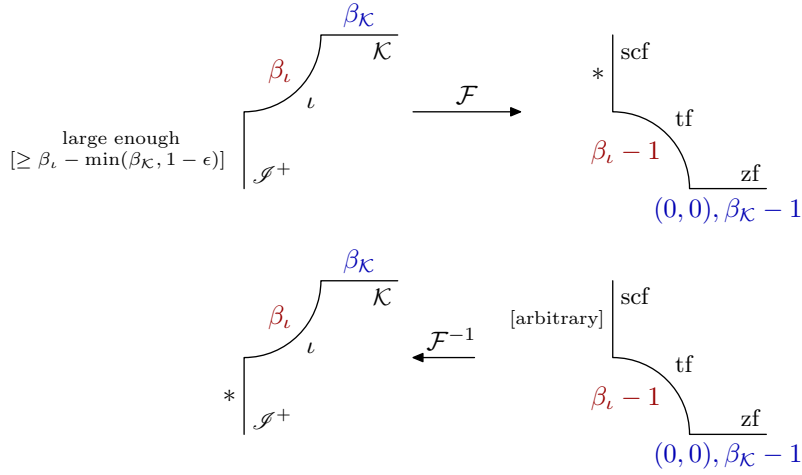}
\caption{\textit{On the top:} illustration of Proposition~\ref{PropDFTFourier}\eqref{ItDFTFourierLo}. Indicated are only the decay rates, i.e., the weights of the function spaces. The label ``$*$'' indicates that the decay rate is typically not too important (and the rate provided by Proposition~\ref{PropDFTFourier} is not sharp in any case). \textit{On the bottom:} illustration of Proposition~\ref{PropDFTInv}\eqref{ItDFTInvLo}.}
\label{FigDFT}
\end{figure}

Figure~\ref{FigDFT} summarizes the relationship between $\iota$ and $\cK$-weights on the one hand and the low-energy behavior on the Fourier transform side on the other hand.

\begin{proof}[Proof of Proposition~\usref{PropDFTInv}]
  \pfstep{Part~\eqref{ItDFTInvLo}.} Consider first the partial expansion at $\zface$ in the local coordinates $\hat r=\frac{\sigma}{\rho}$ and $x\in\ol{\R^n}$ near $\zface$. The terms in this expansion are of the form $\chi(\hat r)\hat r^j v_j(x)$ for $j\in\N_0$ where $v_j\in\Hb^{k,\beta_\iota-1}(\sfX)$; note that $\chi(\hat r)\hat r^j\in\CIc(\R)$. Writing $\phi(\hat r)=\chi(\hat r)\hat r^j$, the inverse Fourier transform is thus
  \begin{equation}
  \label{EqDFTInvExp}
    \frac{1}{2\pi}\int_\R e^{i\sigma t_*}\phi\Bigl(\frac{\sigma}{\rho}\Bigr) v_j(x)\,\dd\sigma = \rho\frac{1}{2\pi}\int_\R e^{i\hat r\,\rho t_*}\phi(\hat r)\,\dd\hat r\cdot v_j(x) = \rho\check\phi(\rho t_*)v_j(x).
  \end{equation}
  But $\check\phi(\rho t_*)$ is a Schwartz function of $\rho t_*$, and $\la\rho t_*\ra^{-1}$ is a smooth non-zero multiple of $\rho_\cK$; so since $\rho=\rho_\sscri\rho_\iota$, we have $\rho\check\phi(\rho t_*)\in\rho_\scri\rho_\iota\rho_\cK^\infty\CI(\sfM)$, and hence~\eqref{EqDFTInvExp} defines an element of $\rho_\sscri^{\beta_\iota-\eps}\rho_\iota^{\beta_\iota}\rho_\cK^\infty\Hb^k(\sfM)$ by Lemma~\ref{LemmaDFTStat}.

  It remains to consider $\cF^{-1}(\chi(\sigma)H(\sigma)v_\pm)$ where now $v_\pm\in\Hb^{k+\ell,(\gamma,\beta_\iota-1,\beta_\cK-1)}(\sfX_\scbtop^\pm)$. We only consider the ``$+$'' sign. Using a partition of unity, we may consider separately the cases that $v_+$ equals one of
  \[
    w_1 \in \Hb^{k+\ell,(\gamma,\beta_\iota-1,\infty)}(\sfX_\scbtop^+),\quad
    w_2 \in \Hb^{k+\ell,(\infty,\beta_\iota-1,\beta_\cK-1)}(\sfX_\scbtop^+).
  \]

  For $v_+=w_1$, Lemma~\ref{LemmaDFTProd} gives, for all $\eta\geq 0$,
  \[
    \chi(\sigma)w_1 \in \bigcap_{j=0}^k \Hb^{j+\ell,\beta_\iota-1-\gamma+\eta}\bigl([0,1];\Hb^{k-j,\gamma-\eta}(\sfX)\bigr).
  \]
  The $\sigma$-order is $>-1$ for all $\eta\geq\max(0,-\beta_\iota+\gamma+\eps)$ (where $\eps>0$ is arbitrary); in order to guarantee that $\beta_\iota-\gamma+\eta\notin\N_0$, we require $\eta\geq\eta_0:=\max(0,-\beta_\iota+\gamma)+\eps$. Corollary~\ref{CorDFTLossy} then implies that
  \[
    \cF^{-1}\bigl(H(\sigma)\chi(\sigma)w_1\bigr) \in \bigcap_{j=0}^k \Hb^{j,\beta_\iota-\gamma+\eta}\bigl(\ol\R;\Hb^{k-j,\gamma-\eta}(\sfX)\bigr) = \Hb^{k,(\gamma-\eta,\beta_\iota,\beta_\iota-\gamma+\eta)}(\sfM);
  \]
  the best possible $\scri$-order obtained in this fashion is thus $\gamma-\eta_0=\min(\gamma,\beta_\iota)-\eps$. According to Corollary~\ref{CorDFTLossy}, the loss $\ell$ must satisfy
  \begin{equation}
  \label{EqDFTInvLoss}
    \ell\geq\lceil(\beta_\iota-\gamma+\eta-\tfrac12)_+\rceil,
  \end{equation}
  which becomes a stronger condition when $\eta$ increases. For present purposes, however, it suffices to ensure that the $\cK$-decay order is $\geq\beta_\cK$, i.e., $\beta_\iota-\gamma+\eta\geq\beta_\cK$; the smallest $\eta\geq 0$ with this property is (up to an $\eps$-loss to guarantee avoidance of integer coincidences) given by $\eta_1:=\max(0,\beta_\cK-\beta_\iota+\gamma)+\eps$. It thus suffices to require~\eqref{EqDFTInvLoss} with $\eta=\max(\eta_0,\eta_1)$, the right-hand side of which is $\lceil(\max(\beta_\iota-\gamma,\beta_\cK,0)+\eps-\frac12)_+\rceil$. We have thus shown
  \[
    \cF^{-1}\bigl(H(\sigma)\chi(\sigma)w_1\bigr) \in \Hb^{k,(\min(\gamma,\beta_\iota)-\eps,\ \beta_\iota,\ \beta_\cK)}(\sfM).
  \]

  Turning to the case $v_+=w_2$, we note that, for all $\eta\geq 0$,
  \[
    w_2 \in \bigcap_{j=0}^k \Hb^{j+\ell,\beta_\cK-1-\eta}\bigl([0,1];\Hb^{k-j,\beta_\iota-\beta_\cK+\eta}(\sfX)\bigr);
  \]
  Provided that $\beta_\cK-\eta\in(0,\infty)\setminus\N$, Corollary~\ref{CorDFTLossy} and Lemma~\ref{LemmaDFTProd} give
  \[
    \cF^{-1}\bigl(H(\sigma)\chi(\sigma)w_2\bigr) \in \Hb^{k,(\beta_\iota-\beta_\cK+\eta,\beta_\iota,\beta_\cK-\eta)}(\sfM),
  \]
  where we need $\ell\geq\lceil(\beta_\cK-\frac12)_+\rceil$. For $\eta=0$, this gives the optimal $\cK$-order $\beta_\cK$. In order to get an almost-sharp $\scri$-order, we take $\eta=\beta_\cK-\eps$, so $\beta_\iota-\beta_\cK+\eta=\beta_\iota-\eps$ then.

  \pfstep{Part~\eqref{ItDFTInvHi}.} We first give up $\bop^+$-regularity to gain $\sigma$-decay, so
  \[
    (1-\chi(\sigma))v_\pm \in \bigcap_{j=0}^k H_{\bop^+}^{j+\ell-\ell',((0,0),\delta+\ell'-1)}\bigl( \pm[0,\infty); \Hb^{k-j,\gamma}(\sfX)\bigr)
  \]
  when $0\leq\ell'\leq\ell$. If $\delta+\ell'-1>-1$ (which holds when $\ell'\geq\lceil\delta_-+\eps\rceil$), Corollary~\ref{CorDFTLossy} gives
  \[
    \cF^{-1}\bigl((1-\chi(\sigma))v_\pm\bigr) \in \bigcap_{j=0}^k \Hb^{j,\delta+\ell'}\bigl(\ol\R;\Hb^{k-j,\gamma}(\sfX)\bigr) = \Hb^{k,(\gamma,\gamma+\delta+\ell',\delta+\ell')}(\sfM),
  \]
  provided $\ell-\ell'\geq\ell_+(\delta+\ell')=\lceil(\delta+\ell'-\frac12)_+\rceil$, which is guaranteed when $\ell-\ell'\geq\delta+\ell'$ (which is $>0$), so $\ell\geq 2\ell'+\delta$. We moreover want to choose $\ell'$ large enough so that
  \[
    \gamma+\delta+\ell' \geq \beta_\iota,\quad
    \delta+\ell' \geq \beta_\cK.
  \]
  Altogether, this yields the condition on $\ell$ in the statement of the proposition.
\end{proof}

\begin{rmk}[b-regularity losses]
\label{RmkDFTRegLoss}
  Importantly, the losses $\ell$ in Proposition~\ref{PropDFTInv}, as well as the loss of $1$ in Proposition~\ref{PropDFTFourier}, depend only on the decay orders, but \emph{not} on the b-regularity order $k$; since we will only use these two results for weights in some compact set, this means that every application gives up at most a fixed finite amount of b-regularity.
\end{rmk}

\subsection{Resolvent bounds}
\label{SsDRes}

We shall translate the resolvent bounds from Theorems~\ref{ThmSpB}, \ref{ThmSpHi}, and \ref{ThmSpLo} into bounds on the b-Sobolev spaces featuring in Propositions~\ref{PropDFTFourier} and \ref{PropDFTInv}. This merely requires establishing relationships between the variable order scattering-type spaces (with extra b-regularity) used in~\S\S\ref{SSp} and the b-Sobolev spaces of \ref{SsDFT}. We write
\[
  \mu := |\dd x|,\qquad \mu_\bop := \frac{|\dd x|}{|x|^3}
\]
for the Euclidean density (used for scattering-type spaces) and unweighted b-density on $X$, respectively. Thus,
\[
  L^2(X,\mu) = \rho^{\frac32}L^2(X,\mu_\bop),\quad \rho=|x|^{-1}.
\]
We first consider the spaces in Theorem~\ref{ThmSpB}:

\begin{lemma}[Scattering and b-Sobolev spaces]
\label{LemmaDRes}
  Let $\sfs\in\CI({}^\scop S^*X)$, $\sfr\in\CI(\ol{\Tsc^*_{\pa X}}X)$, and $k\in\N_0$, $\alpha\in\R$. Put
  \[
    \ubar\sfs:=\inf\sfs,\ \ \bar\sfs:=\sup\sfs,\qquad
    \ubar\sfr:=\inf\sfr,\ \ \bar\sfr:=\sup\sfr,\qquad
    \ubar\sfr_o:=\inf_{\cR_{\rm out}}\sfr,\ \ \bar\sfr_o:=\sup_{\cR_{\rm out}}\sfr,
  \]
  where we recall from~\eqref{EqSpOrderRinout} that $\cR_{\rm out}$ is the zero section of $\Tsc^*_{\pa X}X$. Set $\ell=\lceil\bar\sfr-\bar\sfr_o\rceil+\lceil(\bar\sfs-\lceil\bar\sfr-\bar\sfr_o\rceil-1)_+\rceil$ and $\tilde\ell:=\max(\lceil\ubar\sfr_o-\ubar\sfr\rceil,\lceil-\ubar\sfs\rceil)$. Then for all $\eps>0$, we have continuous inclusions
  \begin{subequations}
  \begin{align}
  \label{EqDResbtosc}
    \bar H_\bop^{k+\ell,\bar\sfr_o+\alpha+\frac52+\eps}(X,\mu_\bop) &\hra \bar H_{\scop;\bop}^{(\sfs-1;k),\sfr+\alpha+1}(X,\mu), \\
  \label{EqDRessctob}
    \bar H_{\scop;\bop}^{(\sfs;k),\sfr+\alpha}(X,\mu) &\hra \bar H_\bop^{k-\tilde\ell,\ubar\sfr_o+\alpha+\frac32-\eps}(X,\mu_\bop),
  \end{align}
  \end{subequations}
  where we require $k\geq\tilde\ell$ in~\eqref{EqDRessctob}.
\end{lemma}

The key features of this result are that the decay rates are (essentially) preserved and that, for fixed orders $\sfs$ and $\sfr$ the b-regularity loss is independent of $k$.

\begin{proof}[Proof of Lemma~\usref{LemmaDRes}]
  We write $H$ instead of $\bar H$ for better readability. We first claim that for every $\eps>0$, we have
  \begin{equation}
  \label{EqDResIni}
    H_{\scop;\bop}^{(\sfs;\tilde\ell),\bar r_o+\eps}(X,\mu) \hra H_\scop^{\sfs+\tilde\ell,\sfr}(X,\mu) \quad\text{if}\quad \bar\sfr_o+\tilde\ell\geq\bar\sfr.
  \end{equation}
  Let $u\in H_{\scop;\bop}^{(\sfs;\tilde\ell),\bar r_o+\eps}$. For any $\chi\in\CIc(X^\circ)$, we have $\chi u\in H^\sfs_\cp(X^\circ)\subset H_\scop^{\sfs,\bar\sfr}$, so we only need to control $(1-\chi)u$. Let $\cU,\cV\subset\Tsc^*_{\pa X}X$ be neighborhoods of $\cR_{\rm out}$ such that $\bar\sfr_o+\eps>\sfr$ on $\bar\cV$, and $\bar\cU\subset\cV$; for $A\in\Psi_\scop^{-\infty,\sfr}$ with $\WF'_\scop(A)\subset\cV$, we then have $A u\in H_\scop^{-\infty,\bar\sfr_o+\eps-\inf_{\cU}\sfr}\subset H_\scop^{-\infty,0}$. To control $u$ outside of $\bar\cU$, let $B\in\Psi_\scop^{\sfs+\tilde\ell,\sfr}$ be elliptic near $\ol{\Tsc^*_{\pa X}}X\setminus\cV$, but with $\WF'_\scop(B)\cap\bar\cU=\emptyset$. Let $D_\bop=(D_{\bop,1},\ldots,D_{\bop,N})$ be a finite spanning set of $\Vb(X)$, and write $D_\bop^\alpha=\prod_{j=1}^N D_{\bop,j}^{\alpha_j}$. We then claim that
  \begin{equation}
  \label{EqDResIniPf}
    \|B u\|_{L^2} \leq C\Biggl(\;\sum_{|\alpha|=\tilde\ell}\|D_\bop^\alpha u\|_{H_\scop^{\sfs,\bar\sfr_o+\eps}} + \|u\|_{H_\scop^{-N,-N}} \Biggr)
  \end{equation}
  for any fixed $N$. To prove this, note that for every $\varpi\in\WF'_\scop(B)$, there exists some multiindex $\alpha$ with $|\alpha|=\tilde\ell$ such that $D_\bop^\alpha\in\rho^{-\tilde\ell}\Diff_\scop^{\tilde\ell}$ is elliptic at $\varpi$. Therefore, a microlocal elliptic parametrix construction in the scattering algebra gives $\|B u\|_{L^2}\leq C(\sum_{|\alpha|=\tilde\ell} \|B Q_\alpha D_\bop^\alpha u\|_{L^2}+\|u\|_{H_\scop^{-N,-N}})$ for some $Q_\alpha\in\Psi_\scop^{-\tilde\ell,-\tilde\ell}$; but since $B Q_\alpha\in\Psi_\scop^{\sfs,\sfr-\tilde\ell}\subset\Psi_\scop^{\sfs,\bar r_o}\subset\Psi_\scop^{\sfs,\bar r_o+\eps}$ by our assumption on $\tilde\ell$, we obtain~\eqref{EqDResIniPf}. This finishes the proof of~\eqref{EqDResIni}. Additional b-regularity persists in the inclusion, so
  \begin{equation}
  \label{EqDResIni2}
    H_{\scop;\bop}^{(\sfs;k+\tilde\ell),\bar r_o+\eps} \hra H_{\scop;\bop}^{(\sfs+\tilde\ell;k),\sfr},\quad \bar\sfr_o+\tilde\ell\geq\bar\sfr.
  \end{equation}

  We conclude that for $\beta\in\R$, and large enough $k'\in\N_0$, and any variable orders $\sfs',\sfr'$, we have
  \[
    \Hb^{k',\beta}(X,\mu_\bop) = \Hb^{k',\beta-\frac32}(X,\mu) = H_{\scop;\bop}^{(0;k'),\beta-\frac32} \hra H_{\scop;\bop}^{(\sfs';k'-\lceil\bar\sfs'_+\rceil),\beta-\frac32} \overset{\text{\eqref{EqDResIni2}}}{\hra} H_{\scop;\bop}^{(\sfs'+\tilde\ell;k'-\tilde\ell-\lceil\bar\sfs'_+\rceil),\sfr'}
  \]
  provided $\beta-\frac32>\bar\sfr'_o$ (so $\beta-\frac32\geq\bar\sfr'_o+\eps$ for small $\eps>0$). Taking $\sfr'=\sfr+\alpha+1$ and $\tilde\ell=\lceil\bar\sfr'-\bar\sfr'_o\rceil=\lceil\bar\sfr-\bar\sfr_o\rceil$, and choosing $\beta>\bar\sfr_o'+\frac32$ and $k',\sfs$ such that $\sfs'+\tilde\ell=\sfs-1$ (so $\sfs'=\sfs-\lceil\bar\sfr-\bar\sfr_o\rceil-1$) and $k'-\tilde\ell-\lceil\bar\sfs'_+\rceil=k$ (so $k'=k+\lceil\bar\sfr-\bar\sfr_o\rceil+\lceil(\bar\sfs-\lceil\bar\sfr-\bar\sfr_o\rceil-1)_+\rceil$), gives~\eqref{EqDResbtosc}.

  For the proof of~\eqref{EqDRessctob}, we first note that
  \[
    H_{\scop;\bop}^{(\sfs;\tilde\ell),\sfr}(X,\mu) \hra H_\scop^{\sfs+\tilde\ell,\ubar\sfr_o-\eps}(X,\mu)\quad\forall\,\eps>0\ \ \text{if}\quad \ubar\sfr+\tilde\ell\geq\ubar\sfr_o;
  \]
  the proof is completely analogous to that of~\eqref{EqDResIni}. If moreover $\sfs+\tilde\ell\geq 0$, then
  \[
    H_{\scop;\bop}^{(\sfs;k),\sfr+\alpha}(X,\mu) \hra H_{\scop;\bop}^{(\sfs+\tilde\ell;k-\tilde\ell),\ubar\sfr_o+\alpha-\eps}(X,\mu) \subset H_\bop^{k-\tilde\ell,\ubar\sfr_o+\alpha-\eps}(X,\mu) = H_\bop^{k-\tilde\ell,\ubar\sfr_0+\alpha+\frac32-\eps}(X,\mu_\bop).\qedhere
  \]
\end{proof}

\begin{cor}[Resolvent bounds at nonzero frequencies]
\label{CorDResB}
  We use the notation of Theorem~\usref{ThmSpB}. Let $\beta_\scface<1+\ubar S$. Then there exists $\ell=\ell(P_0,\beta_\scface)\in\N_0$ such that the following holds.
  \begin{enumerate}
  \item\label{ItDResB1}{\rm (Regularity in $\sigma$.)} Suppose $\sigma_0\in\C$ with $\Im\sigma_0\geq 0$ and $\sigma_0\neq 0$ is such that $\wh{P_0}(\sigma_0)$ has trivial kernel on the space~\eqref{EqSpBNullA}. Then for all $j,k\in\N_0$ and $\eps>0$, the $j$-th derivative of $\wh{P_0}(\sigma)^{-1}$ at $\sigma=\sigma_0$ is bounded as a map
    \[
      \pa_\sigma^j\wh{P_0}(\sigma_0)^{-1} \colon \bar H_\bop^{k+j+\ell,\beta_\scface+1}(X,\mu_\bop) \to \bar H_\bop^{k,\beta_\scface-\eps}(X,\mu_\bop).
    \]
  \item\label{ItDResB2}{\rm (Joint regularity in $\sigma$ and $x$.)} Suppose $I\subset\R$ is an open interval such that $0\neq\bar I$, and $\wh{P_0}(\sigma)$ is invertible for all $\sigma\in I$. For a function $\hat f=\hat f(\sigma,x)$ on $I\times X^\circ$, set $(\wh{P_0}{}^{-1}\hat f)(\sigma,\cdot):=\wh{P_0}(\sigma)^{-1}\hat f(\sigma,\cdot)$. Write $\bar H_\bop^{k,\beta}(I\times X):=\rho^\beta\bar H_\bop^k(I\times X)$ for the Sobolev space on $I\times X$ testing for $k$ degrees of regularity with respect to $\pa_\sigma$ and b-vector fields on $X$. Then
    \begin{equation}
    \label{EqDResBjoint}
      \wh{P_0}{}^{-1} \colon \bar H_\bop^{k+\ell,\beta_\scface+1}(I\times X,|\dd\sigma|\mu_\bop) \to \bar H_\bop^{k,\beta_\scface-\eps}(I\times X,|\dd\sigma|\mu_\bop).
    \end{equation}
  \end{enumerate}
\end{cor}
\begin{proof}
  For the first part, we apply Theorem~\ref{ThmSpB} where the scattering regularity and decay orders $\sfs_\sigma$ and $\sfr_\sigma$ are induced by a stationary-$P_0$-admissible order function $\sfs\in\CI(\Stb^*M_0)$ with weights $\alpha_+,0$ such that $\beta_\scface\geq\sfr_\sigma+\alpha_++\frac32$; this is possible by~\eqref{EqSpBThrOut}. Pre- and post-composing the map~\eqref{EqSpBInvReg} with the inclusion maps~\eqref{EqDResbtosc} and \eqref{EqDRessctob} then yields the claim. (The value $\ell$ in the present result is the sum of $\ell$ and $\tilde\ell$ in Lemma~\ref{LemmaDRes}.)

  For the second part, let $\hat f\in\bar H_\bop^{k+\ell,\beta_\scface+1}(I\times X)$. We must bound, for $k'\leq k$, the norm of
  \[
    \pa_\sigma^{k'}\bigl(\wh{P_0}(\sigma)^{-1}\hat f(\sigma)\bigr) = \sum_{k''=0}^{k'} \binom{k'}{k''} \bigl(\pa_\sigma^{k''}\wh{P_0}(\sigma)^{-1}\bigr) \bigl(\pa_\sigma^{k'-k''}\hat f(\sigma)\bigr)
  \]
  in the space $L^2(I_\sigma,|\dd\sigma|;\Hb^{k-k'}(X,\mu_\bop))$. In the $k''$-th term, we note that
  \[
    \pa_\sigma^{k'-k''}\hat f\in L^2\bigl(I_\sigma;\bar H_\bop^{k-k'+k''+\ell,\beta_\scface+1}(X)\bigr),
  \]
  so part~\eqref{ItDResB1} (with $k$ and $j$ there being $k-k'$ and $k''$ in present notation) yields the desired bound.
\end{proof}

\begin{rmk}[Losses in b-regularity]
\label{RmkDResLoss}
  It is important here to translate sharp (as far as b-order is concerned) estimates for \emph{derivatives} of the resolvent on variable order spaces into (lossy) estimates solely on b-spaces. If one were to use a lossy b-estimate on $\wh{P_0}(\sigma)^{-1}$ and iterate that as around~\eqref{EqSpBInvRegPf}, one would incur $\sim k$ b-derivative losses when estimating $\pa_\sigma^j\wh{P_0}(\sigma)^{-1}$; this, however, would break tame estimates in nonlinear applications.
\end{rmk}

\bigskip

We next turn to the high-energy analogue of Corollary~\ref{CorDResB}: this involves semiclassical scattering Sobolev spaces, for which b-regularity is equivalent to scattering decay, scattering regularity, and a power of the semiclassical parameter away from the zero section over $\pa X$. Thus, we can shift these three orders by arbitrary bounded amounts upon using up a fixed finite amount of b-derivatives. Recalling the semiclassical loss function $\delta_\Gamma$ from~\eqref{EqSpHiTrLoss}, we obtain:

\begin{prop}[High-energy resolvent bounds]
\label{PropDResHi}
  We use the notation of Theorem~\usref{ThmSpHi}, so $\wh{P_0}(\sigma)$ is invertible when $\sigma\in\C$, $\Im\sigma\geq 0$, and $|\sigma|\geq C_0$; and we recall the spaces $H_{\bop^+}^{k,\beta,\gamma}$ from~\eqref{EqDFTSobHi}. Let $\beta_\scface<1+\ubar S$ and $\delta,\delta'\in\R$. Then there exists $\ell=\ell(P_0,\beta_\scface,\delta,\delta')\in\N_0$ such that for all $k\in\N_0$ and all $\eps>0$,
  \begin{equation}
  \label{EqDResHi}
    \wh{P_0}{}^{-1} \colon H_{\bop^+}^{k(1+\delta_\Gamma(-\eps))+\ell,\beta_\scface+1,\delta}(\pm[C_0,\infty)\times X) \to H_{\bop^+}^{k,\beta_\scface-\eps,\delta'}(\pm[C_0,\infty)\times X).
  \end{equation}
\end{prop}

The relative spatial decay rates of input ($\beta_\scface+1$) and output ($\beta_\scface-\eps$) are as in~\eqref{EqDResBjoint} and~\eqref{EqDResLo}. The $|\sigma|$-decay rates $\delta,\delta'$ are inconsequential, as one can shift them arbitrarily by giving up $\bop^+$-regularity. The only dangerous part is the derivative loss $k\delta_\Gamma(-\eps)$; but when $P_0$ is \emph{strongly} trapping admissible (Definition~\ref{DefSSTrapAdm})---which we show for the scalar and 1-form wave operators in Theorems~\ref{ThmA1Adm} and \ref{ThmA2Adm}---then for $\eps=\frac{1}{k}$ this loss is, in fact, independent of $k$; this is crucial for the purpose of preserving the tameness of estimates when utilizing these bounds on $\wh{P_0}{}^{-1}$ to improve decay.

\begin{cor}[Bounds for waves localized at high frequencies]
\label{CorDResHiLoc}
  Let $P_0$ be a stationary wave-type operator (Definition~\usref{DefSSAdm}) that is trapping admissible (Definitions~\usref{DefSSTrapAdm}). Let $C_0>0$ be such that $\wh{P_0}(\sigma)$ has trivial nullspace on~\eqref{EqSpBNullA} for all $\sigma\in\R$, $|\sigma|\geq C_0$. Let $\chi\in\CIc(\R)$ be equal to $1$ on $[-C_0,C_0]$, and let $\chi_T\in\CI(\R)$ be equal to $0$ on $(-\infty,1]$. Let $\beta_\sscri,\beta_\iota,\beta_\cK\in\R$ and $\beta'_\cK,\beta'_\iota\in\R$.\footnote{In practice, we will typically take $\beta_\iota'=\beta_\iota$ and $\beta'_\cK=\beta_\cK$.} Then there exist $\ell=\ell(P_0,\beta_\cK,\beta_\iota,\beta_\sscri,\beta'_\cK,\beta'_\iota)\in\N_0$ such that for all
  \[
    f\in \bigcap_{\eta>0}\Hb^{k(1+\delta_\Gamma(-\eta))+\ell,\ (\beta_\sscri+1,\ \beta_\iota+2,\ \beta_\cK)}(\Omega)^{\bullet,-},
  \]
  the function
  \[
    u = \cF^{-1} \bigl( (1-\chi(\sigma))\wh{P_0}(\sigma)^{-1}\hat f(\sigma) \bigr)
  \]
  satisfies
  \begin{equation}
  \label{EqDResHiLocu}
    \chi_T(t_*)u \in \Hb^{k,(\beta_\sscri',\beta_\iota',\beta_\cK')}(\Omega)^{\bullet,-}\quad\forall\,\beta_\sscri'<\min(\beta_\sscri,1+\ubar S).
  \end{equation}
  More generally, if we define $\sfM:=[\ol{\R_{t_*}}\times X;\pa\ol\R\times\pa X]$ as in~\eqref{EqDFTSpace}, then for all source terms $f\in\Hb^{k(1+\delta_\Gamma(-\eta))+\ell,(\beta_\sscri+1,\beta_\iota+2,\beta_\cK)}(\sfM)$ (in the notation of~\eqref{EqDFTHbM}), we have $u\in\Hb^{k,(\beta_\sscri',\beta_\iota',\beta_\cK')}(\sfM)$. 
\end{cor}

Corollary~\ref{CorDResHiLoc} asserts that the high-frequency components of $\cF^{-1}(\wh{P_0}(\sigma)^{-1}\hat f(\sigma))$ have \emph{arbitrary} decay rates at $\iota^+$ and $\cK^+$ (regardless of the decay rates of $f$), provided we give up enough b-derivatives. The reason is that the lack of decay of $f$ as $t_*\to\infty$ is only reflected in limited regularity of $\hat f$ at $\sigma=0$; but the b-regularity of $f$, which persists upon passing to the Fourier transform, is the same as standard $\sigma$-regularity away from $\sigma=0$, so upon localizing to nonzero frequencies, we can hence use any finite amount of regularity of $\hat f$ (and thus of $\wh{P_0}{}^{-1}\hat f$) to ensure any desired amount of $t_*$-decay upon inverse Fourier transforming.

\begin{proof}[Proof of Corollary~\usref{CorDResHiLoc}]
  We write $\ell'$ for the value of $\ell$ in Proposition~\ref{PropDResHi} and $\ell''$ for the value of $\ell$ in Proposition~\ref{PropDFTFourier}\eqref{ItDFTFourierHi}. Proposition~\ref{PropDFTFourier}\eqref{ItDFTFourierHi} shows that $\hat f\in H_{\bop^+}^{k(1+\delta_\Gamma(-\eta))+\ell-\ell'',\beta_\sscri+1,0}(\pm[C_0,\infty)\times X)$. By~\eqref{EqDResHi} (with $\delta=\delta'=0$), applying $\wh{P_0}{}^{-1}$ to this yields an element of $H_{\bop^+}^{k+\ell-1-\ell',(\gamma-\eps,0)}$ for $\gamma:=\min(\beta_\sscri,1+\ubar S)$. By Proposition~\ref{PropDFTInv}\eqref{ItDFTInvHi}, the inverse Fourier transform of this lies in the space $\Hb^{k+\ell-1-\ell'-\ell''',(\gamma-\eps,\beta'_\iota,\beta'_\cK)}(\sfM)$, where $\ell'''$ is number denoted $\ell$ there. Fixing $\ell=1+\ell'+\ell''+\ell'''$ gives~\eqref{EqDResHiLocu}.
\end{proof}

\bigskip

A consequence of Corollary~\ref{CorDResHiLoc} is that only the low-frequency behavior of $\wh{P_0}(\sigma)^{-1}$ has a non-trivial bearing on the late-time asymptotics of $u$. In the scattering-b-transition setting relevant for the low-energy estimates in Theorem~\ref{ThmSpLo} then, b-regularity is equivalent to scattering-b-regularity and -decay away from the zero section over $\scface\subset X_\scbtop$, completely analogously to (the proof of) Lemma~\ref{LemmaDRes}. Thus, for $\sfs\in\CI({}^\scbtop S^*X)$, $\sfr\in\CI(\ol{{}^\scbtop T^*_\scface}X)$, $k\in\N_0$, and $\alpha_+,\alpha_\tface,\alpha_\zface\in\R$, the analogue of~\eqref{EqDResbtosc} is that for any $\eps>0$, we have \emph{uniform} (in $\sigma\in\pm[0,1]$) bounds
\begin{subequations}
\begin{equation}
\label{EqDResscbt1}
  \|u\|_{\bar H_{(\scbtop,|\sigma|);\bop}^{(\sfs-1;k),(\sfr+\alpha_++1,\alpha_\tface+2,0)}(X,\mu)} \leq C\|u\|_{\bar H_{\bop,|\sigma|}^{k+\ell,(\bar\sfr_o+\alpha_++\frac52+\eps,\alpha_\tface+\frac72,0)}(X,\mu_\bop)},
\end{equation}
where we write $\bar H_{\bop,|\sigma|}^k:=\bar H_{\scbtop,|\sigma|}^{(0;k)}$, similarly for weighted spaces, and $\ell$ is as in Lemma~\ref{LemmaDRes}. (The shift of $\frac32$ in the decay orders at $\rho^{-1}(0)=\scface\cup\tface$ is due to the passage from $\mu$ to $\mu_\bop$.) Thus, the $\sigma$-dependence of the norm on the right is solely due to the $\sigma$-dependence of the weights, to wit,
\[
  \|u\|_{H_{\bop,|\sigma|}^{k,(\beta_\scface,\beta_\tface,\beta_\zface)}(X,\mu_\bop)} = \|(\rho_\scface^{-\beta_\scface}\rho_\tface^{-\beta_\tface}\rho_\zface^{-\beta_\zface})_{|\sigma|}u\|_{\Hb^k(X,\mu_\bop)},
\]
with explicit choices of the defining functions appearing here being given by~\eqref{EqDFTbdfs}. We similarly have, as the analogue of~\eqref{EqDRessctob}, the uniform bounds
\begin{equation}
\label{EqDResscbt2}
  \|u\|_{\bar H_{\bop,|\sigma|}^{k-\tilde\ell,(\ubar\sfr_o+\alpha_++\frac32-\eps,\alpha_\tface+\frac32,0)}(X,\mu_\bop)} \leq C\|u\|_{\bar H_{(\scbtop,|\sigma|);\bop}^{(\sfs;k),(\sfr+\alpha_+,\alpha_\tface,0)}(X,\mu)}
\end{equation}
\end{subequations}
for $\sigma\in\pm[0,1]$. The analogue of Corollary~\ref{CorDResB}\eqref{ItDResB2} then reads:

\begin{prop}[Uniform low-energy resolvent bounds]
\label{PropDResLo}
  We use the notation of Theorem~\usref{ThmSpLo}, including the parameters $\alpha\in\R$, with $P_0$ being $\tface$-admissible with weight $\alpha+\frac32$ lying in the indicial gap $(\beta^-,\beta^+)$ of $\wh{P_0}(0)$, and $c>0$ (with $\wh{P_0}(\sigma)$ being invertible for $|\sigma|\leq c$, $\Im\sigma\geq 0$). Set
  \[
    X_\scbtop^\pm := [ \pm[0,c]\times X; \{0\}\times\pa X ]
  \]
  (analogously to~\eqref{EqDFTXscbtpm}), and use an unweighted b-density on $X_\scbtop^\pm$ (e.g., $|\frac{\dd\sigma}{\sigma}\,\mu_\bop|$) to define b-Sobolev spaces as in~\eqref{EqDFTHbX}. Let $\beta_\scface<1+\ubar S$. Then there exists $\ell=\ell(P_0,\beta_\scface)\in\N_0$ such that for all
  \[
    \beta_\tface,\ \beta_\zface\in\R,\quad \beta_\tface-\beta_\zface\in(\beta^-,\beta^+),
  \]
  and $k\in\N_0$, $\eps>0$, we have
  \begin{equation}
  \label{EqDResLo}
    \wh{P_0}{}^{-1} \colon \Hb^{k+\ell,(\beta_\scface+1,\beta_\tface+2,\beta_\zface)}(X_\scbtop^\pm) \to \Hb^{k,(\beta_\scface-\eps,\beta_\tface,\beta_\zface)}(X_\scbtop^\pm).
  \end{equation}
\end{prop}
\begin{proof}
  For $\beta_\zface=0$, the proof is identical to that of Corollary~\ref{CorDResB}, \emph{mutatis mutandis}, except now we use $\sigma\pa_\sigma$ instead of $\pa_\sigma$. We can reduce the case of general $\beta_\zface$ to $\beta_\zface=0$ via division by $|\sigma|^{\beta_\zface}=\rho_\tface^{\beta_\zface}\rho_\zface^{\beta_\zface}$, which thus replaces $(\beta_\tface,\beta_\zface)$ by $(\beta_\tface-\beta_\zface,0)$.
\end{proof}

Note that for $\beta_\cK=-\frac12$, elements in the domain of~\eqref{EqDResLo} lie in $L^2(\pm[0,c),|\dd\sigma|;H^k_\loc(X^\circ))$, i.e., they have weight $0$ at $\sigma=0$; such spaces arise from Plancherel's theorem as in Lemma~\ref{LemmaMUetbFTb}, and for such inputs, we analyzed the output of $\wh{P_0}(\sigma)^{-1}$ in $L^2(\R_\sigma;\Hb^k(X))$-type spaces, e.g., in the proof of Theorem~\ref{ThmA1Adm}. The flexibility afforded by allowing for $\beta_\cK>-\frac12$ is thus crucial in obtaining $t_*$-decay rates stronger than $L^2$ in time (over compact subsets), which in conjunction with b-regularity and Sobolev embedding yielded only pointwise $t_*^{-\frac12}$-decay ($\sim\rho_\cK^{\frac12}$ in spatially compact sets) in~\eqref{EqFTame3} for $\aleph=0$.

\subsection{Radiation field and decay at \texorpdfstring{$\scri^+$}{null infinity}}
\label{SsDscri}

For the sake of simplicity, we study here only stationary wave-type operators $P_0$ which near $\scri^+\cap\iota^+$ are of the form
\begin{equation}
\label{EqA2RadOp}
  \frac12 \rho_\sscri^{-1}\rho_+^{-2}P_0 = -(\rho_\sscri\pa_{\rho_\sscri}-1)(\rho_\sscri\pa_{\rho_\sscri}-\rho_+\pa_{\rho_+}) + \tilde P,\quad \tilde P\in \rho_\sscri^{\frac12}\Diffb^2(M);
\end{equation}
this happens if and only if $S\in\rho\CI(X)$ in Definition~\ref{DefSSAdm}. Examples of such operators $P_0$ include the tensor wave operator on $(p,q)$-tensors on Kerr $P_0=\Box_g$ in view of~\eqref{EqSSOpMemeb}. (In fact, in this case we have $\tilde P\in\rho_\sscri\Diffb^2$ since $\Box_g$ has smooth coefficients as a b-differential operator on $M_1$.)

We shall extract a $\rho_\sscri^1$ leading order term of solutions of $P_0 u=f$ when $f=o(\rho_\sscri^2)$ at $(\scri^+)^\circ$ using ideas from \cite{HintzVasyMink4} and simple results on the integration of ODEs from \citestab{\S\S\ref*{SsTM}--\ref*{SsTInt}}. (These two works extract radiation fields in the analytically more involved, but conceptually essentially identical, setting of non-stationary wave-type equations of the form~\eqref{EqSDWAdmOp} with non-zero $\tilde p_1$ and $p_0$.) In order to capture radiation fields at $\scri^+$ efficiently, we introduce:

\begin{definition}[Partial polyhomogeneity]
\label{DefDscriPhg}
  Fix a cutoff function $\chi_\sscri\in\CI(M)$ with $\chi_\sscri=1$ near $\scri^+$ and $\supp\chi_\sscri\cap\cK^+=\emptyset$. Let $k\in\N_0$, $\beta_\sscri,\beta_+,\beta_\cK\in\R$, and $\zeta\in\R$. Then we write
  \begin{align*}
    &\Hb^{k,(((\zeta,0),\beta_\sscri),\beta_+,\beta_\cK)}(\Omega)^{\bullet,-} \\
    &\quad := \bigl\{ \chi_\sscri\rho_\sscri^\zeta u_{(1,0)} + \tilde u \colon u_{(1,0)}\in\Hb^{k,\beta_+}(\scri^+\cap\{t_*\geq 1\}),\ \tilde u\in\Hb^{k,(\beta_\sscri,\beta_+,\beta_\cK)}(\Omega)^{\bullet,-} \bigr\}.
  \end{align*}
\end{definition}

From now on, \textit{we work only in a neighborhood of $\scri^+\cap\{t_*\geq 1\}$, and thus drop the weight at $\cK^+$ from the notation.} 

\begin{prop}[Existence of and bounds for the radiation field]
\label{PropDscriRad}
  Let $\beta_\sscri,\beta_+\in\R$ with $1\neq\beta_+<\beta_\sscri$. There exists $\ell=\ell(\beta_\sscri,\beta'_\sscri)\in\N_0$ such that the following holds. Suppose $u\in\Hb^{k+\ell,(\beta_\sscri,\beta_+)}$ vanishes for $t_*\leq 1$ and satisfies $P_0 u\in\Hb^{k+\ell,(\beta'_\sscri+1,\beta_++2)}$ where $1\neq\beta'_\sscri<\frac32$ and $\beta'_\sscri\neq\beta_+$. Then
  \begin{equation}
  \label{EqA2Rad}
    u \in \Hb^{k,(((1,0),\beta'_\sscri),\beta_+)}.
  \end{equation}
\end{prop}

Thus, $u$ has a radiation field plus an error term whose weights match those of $P_0 u$ (given the usual shift of two orders because of the overall factor $\rho_\sscri\rho_+^2$ in~\eqref{EqSSOpMem}).

\begin{proof}[Proof of Proposition~\usref{PropDscriRad}]
  Suppose we have shown $u\in\Hb^{k+\ell,(((1,0),\beta_\sscri),\beta_+)}$ (which matches the assumption on $u$ if we reduce $\beta_\sscri$ so that $\beta_\sscri<1$). In view of~\eqref{EqA2RadOp}, the equation $\rho_\sscri^{-1}\rho_+^{-2}P_0 u\in\Hb^{k+\ell,(\beta'_\sscri,\beta_+)}$ can be written as
  \begin{align*}
    (\rho_\sscri\pa_{\rho_\sscri}-1)(\rho_\sscri\pa_{\rho_\sscri}-\rho_+\pa_{\rho_+})u =: h &\in \Hb^{k+\ell,(\beta_\sscri',\beta_+)} + \Hb^{k+\ell-2,(((\frac32,0),\beta_\sscri+\frac12),\beta_+)} \\
      &\subset \Hb^{k+\ell-2,(((\frac32,0),\beta''_\sscri),\beta_+)},\quad \beta''_\sscri:=\min(\beta_\sscri+\tfrac12,\beta'_\sscri)<\tfrac32.
  \end{align*}
  If $\beta''_\sscri=\beta_\sscri+\frac12$, we arrange for $\beta''_\sscri\neq\beta_+$ by decreasing $\beta''_\sscri$ by an arbitrarily small amount, if necessary. We similarly arrange that $\beta''_\sscri\neq 1$. We first integrate $\rho_\sscri\pa_{\rho_\sscri}-1$ using \citestab{Lemma~\ref*{LemmaTMIntFuchs}}, which yields
  \begin{equation}
  \label{EqDscriRadHyp}
    (\rho_\sscri\pa_{\rho_\sscri}-\rho_+\pa_{\rho_+})u \in \Hb^{k+\ell-2,(((1,0),\beta''_\sscri),\beta_+)},
  \end{equation}
  the $\rho_\sscri$ leading order term arising from the indicial root $1$ of $\rho_\sscri\pa_{\rho_\sscri}-1$.

  If $\beta''_\sscri<\beta_+$, we integrate~\eqref{EqDscriRadHyp} from $\rho_\sscri=\delta>0$ (where $u=\cO(\rho_+^{\beta_+})$ by assumption) towards $x_\sscri=0$ using \citestab{Lemma~\ref*{LemmaTIntHyp}} to deduce
  \begin{equation}
  \label{EqDscriRadImpr}
    u \in \Hb^{k+\ell-2,(((1,0),\beta''_\sscri),\beta_+)}.
  \end{equation}
  This improves the $x_\sscri$-decay of the remainder term of $v$ from $\beta_\sscri$ to $\beta''_\sscri$. If, on the other hand, $\beta''_\sscri>\beta_+$, then we integrate~\eqref{EqDscriRadHyp} from $\rho_+=t_*^{-1}=1$ (where $u$ vanishes) towards $\rho_+=0$; if $\beta_+<1$, this directly gives~\eqref{EqDscriRadImpr}, whereas if $\beta_+>1$, the radiation field of $u$ gets transported to a $\rho_+^1$ leading order term at $\iota^+$, so one initially obtains
  \begin{align*}
    &u \equiv \rho_+ u^{(1,0)}(x_\sscri,\omega) \bmod \Hb^{k+\ell-2,(((1,0),\beta''_\sscri),\beta_+)}, \\
    &\qquad u^{(1,0)}\in\Hb^{k+\ell-2,((1,0),\beta''_\sscri)}([0,1)_{\rho_\sscri}\times\Sph^2) := \rho_\sscri\CI([0,1)\times\Sph^2) + \Hb^{k+\ell-2,\beta''_\sscri}([0,1)\times\Sph^2);
  \end{align*}
  but the a priori assumption that $u=\cO(\rho_+^{\beta_+})$ at $(\iota^+)^\circ$ implies that $u^{(1,0)}$ must, in fact, vanish identically. Therefore, we have again proved~\eqref{EqDscriRadImpr}. Iterating this improvement finitely many times yields~\eqref{EqA2Rad}. The constant $\ell$ is twice the number of iterations required, and thus only depends on $\beta_\sscri$ and $\beta_\sscri'$.
\end{proof}

\subsection{Application I': power nonlinearities with smaller exponents}
\label{SsDp}

Recall from~\eqref{EqA1PowerHeur} that the only place where the restriction $p\geq 4$ on the power in the equation $\Box_g u=f+u^p$ arises is in ensuring that the $\iota^+$-decay rate of $u^p$ is (at least) two more than that of $u$ itself, thus balancing the two orders of $\iota^+$-decay lost upon inverting $\Box_g$. In the proof of Theorem~\ref{ThmA1Power}, we only recorded $1-\eps_+<1$ orders of $\iota^+$-decay on $u$; the key towards lowering $p$ to $p=3$ is thus to improve the $\iota^+$-decay order of $u$ to $>1$.

Since we need to prove $\rho_+^{1+\eps_+}$-decay for $u$ at $\iota^+$, we need to require the source term $f$ to have $\rho_+^{3+\eps_+}$-decay. Moreover, one expects anything less than $\rho_\sscri^{2+\eps_+}$-decay of $f$ at $\scri^+$ to prompt $u$ to have worse than $\rho_\sscri^{1+\eps_+}$-decay, which then propagates to $\iota^+$ and destroys the desired $\rho_+^{1+\eps_+}$-decay. Therefore, we shall work with $0<\eps_+<\eps_\sscri$ and require $\rho_\sscri^{2+\eps_\sscri}$-decay for $f$. Even for the linear scalar wave equation, however, waves with such sources have a radiation field which is a product of $r^{-1}=\rho_\sscri\rho_+$ with a function of $(t_*,\omega)$, which we must therefore capture in the space in which we seek the solution $u$ of $\Box_g u=f+u^3$. Note that in a nonlinear iteration, the presence of such a radiation field is consistent with the $\rho_\sscri^{2+\eps_\sscri}$-decay of the source term in the next iteration step since the forward operator $\Box_g u-u^p$ annihilates the radiation field to leading order (cf.\ \eqref{EqSSAdmBox} and the term $\rho\pa_\rho-1=r\pa_r+1$ in the wave operator~\eqref{EqSSAdmOp}). Motivated by these heuristics, we now assert the following extension of Theorem~\ref{ThmA1Power}:

\begin{thm}[Small data global existence: $p=3$]
\label{ThmDp}
  Let $0<\eps_\cK<\eps_+<\eps_\sscri<\frac12$ and set
  \[
    \beta_\sscri=1+\eps_\sscri,\quad
    \beta_+=1+\eps_+,\quad
    \beta_\cK=1+\eps_\cK.
  \]
  Recall the defining functions $\rho_\sscri,\rho_+,\rho_\cK$ from~\eqref{EqFWeights}. Fix a cutoff function $\chi_\sscri\in\CI(M)$ which equals $1$ near $\scri^+$ and has $\supp\chi_\sscri\cap\cK^+=\emptyset$. Then there exist $d\in\N$ and $\eps>0$ such that for all $f\in\Hb^{\infty,(\beta_\sscri+1,\beta_++2,\beta_\cK)}(\Omega,\mu_\bop)^{\bullet,-}$ (where $\mu_\bop$ is given by~\eqref{EqA1Powermub}) with
  \[
    \|f\|_{\Hb^{d,(\beta_\sscri+1,\beta_++2,\beta_\cK)}(\Omega,\mu_\bop)^{\bullet,-}}^2 = \iint_\Omega |\rho_\sscri^{-\beta_\sscri-1}\beta_+^{-\beta_+-2}\rho_\cK^{-\beta_\cK}(t_*\pa_{t_*},r\pa_x)^{\leq d}f|^2\,\dd\mu_\bop < \eps^2,
  \]
  the equation
  \begin{equation}
  \label{EqDpEq}
    \Box_g u = f + u^3
  \end{equation}
  has a unique global forward solution $u\in\Hb^{\infty,(((1,0),\beta_\sscri),\beta_+,\beta_\cK)}(\Omega,\mu_\bop)^{\bullet,-}$ on $\Omega$ (in the notation of Definition~\usref{DefDscriPhg}); that is,
  \begin{equation}
  \label{EqDpu}
  \begin{split}
    &u(t_*,x) = \chi_\sscri r^{-1}u_\sscri\Bigl(t_*,\frac{x}{|x|}\Bigr) + \tilde u(t_*,x), \\
    &\qquad u_\sscri \in \dot H_\bop^{\infty,\eps_+}([1,\infty],|\tfrac{\dd t_*}{t_*}|), \quad \tilde u\in\Hb^{\infty,(\beta_\sscri,\beta_+,\beta_\cK)}(\Omega)^{\bullet,-}.
  \end{split}
  \end{equation}
\end{thm}
\begin{proof}
  We will apply Nash--Moser iteration (Theorem~\ref{ThmA1NM}) to the map $\Phi(u):=\Box_g u-u^3-f$ using the spaces
  \begin{alignat*}{2}
    B^k_\eta &:= \bigl\{ u\in\Hb^{k,(((1,0),\beta_\sscri),\ \beta_+,\ \beta_\cK)}(\Omega,\mu_\bop)^{\bullet,-} &&\colon t_*\geq\eta\ \text{on}\ \supp u \bigr\}, \\
    \bfB^k_\eta &:= \bigl\{ f\in\Hb^{k,(\beta_\sscri+1,\beta_++2,\beta_\cK)}(\Omega,\mu_\bop)^{\bullet,-} &&\colon t_*\geq\eta\ \text{on}\ \supp f \bigr\},
  \end{alignat*}

  \pfstep{Forward mapping properties.} We first claim that $\Phi\colon B^\infty\to\bfB^\infty$. Let $u\in B^\infty$. Writing $u$ as in~\eqref{EqDpu}, we have $\Box_g\tilde u\in\bfB^\infty$ since $\Box_g\in\rho_\sscri\rho_+^2\Diffb^2$. For the contribution of $\chi_\sscri r^{-1}u_\sscri=\chi_\sscri\rho u_\sscri$, we write $\Box_g(\chi_\sscri\rho u_\sscri)=[\Box_g,\chi_\sscri]\rho u_\sscri+\chi_\sscri\Box_g(\rho u_\sscri)$. Since $\Box_g\in \rho_\sscri\rho_+^2\Diffb^2$, the first summand is of class
  \[
    \rho_\sscri^\infty\rho_+^2\rho_\cK^\infty\CI(\Omega)\cdot\rho\dot H_\bop^{\infty,\eps_+}([1,\infty]_{t_*}\times\Sph^2) \subset \Hb^{\infty,(\infty,3+\eps_+,\infty)}(\Omega)^{\bullet,-} \subset \bfB^\infty.
  \]
  For the second summand, recall from~\eqref{EqSSAdmOp} and~\eqref{EqSSAdmBox} that
  \[
    \Box_g = \pa_{t_*}\bigl(-2\rho(\rho\pa_\rho-1)+\rho^2\Diffb^1(X)\bigr) + \rho^2\Diffb^2(X) + \rho^2\CI(X)\pa_{t_*}^2.
  \]
  This maps $r^{-1}u_\sscri=\rho u_\sscri(t_*,\omega)$ into the space
  \[
    \pa_{t_*}\rho^3\Hb^{\infty,\eps_+}([1,\infty]_{t_*}\times\Sph^2) + \rho^3\Hb^{\infty,\eps_+}([1,\infty]_{t_*}\times\Sph^2) \subset \rho^3\Hb^{\infty,\eps_+}([1,\infty]_{t_*}\times\Sph^2).
  \]
  Upon multiplication by $\chi_\sscri$, this is contained in $\Hb^{\infty,(3-\eps,3+\eps_+,\infty)}(\Omega)^{\bullet,-}$ for all $\eps>0$, and thus in $\bfB^\infty$. This shows that
  \begin{equation}
  \label{EqDpBoxScri}
    \Box_g(\chi_\sscri\rho u) \in \Hb^{\infty,(\beta_\sscri+1,\beta_++2,\beta_\cK)} = \bfB^\infty.
  \end{equation}

  Next, since $u\in\Hb^{\infty,(1-\eps,\beta_+,\beta_\cK)}\subset\cC_\bop^{\infty,(1-\eps,\beta_+,\beta_\cK)}$ for all $\eps>0$, we have
  \[
    u^3 \in \Hb^{\infty,(3-3\eps,3\beta_+,3\beta_\cK)} \subset \Hb^{\infty,(\beta_\sscri+1,\beta_++2,\beta_\cK)}
  \]
  since $\beta_\sscri<2$ and, crucially, $3\beta_+\geq\beta_++2$ in view of $\beta_+\geq 1$. (This is the key gain over the analysis in Theorem~\ref{ThmA1Power}.) This shows that $\Phi(u)\in\bfB^\infty$. The estimates~\eqref{EqA1NMPhi} for $\Phi$ follow from Lemma~\ref{LemmaA1Nonlin}.

  \pfstep{Admissibility of the linearization.} The linearized operator $P_u:=\Phi'(u)=\Box_g-3 u^2$ differs from the stationary model $\Box_g$ by the zeroth order multiplication operator
  \[
    3 u^2 \in \rho_\sscri^2\cC_\bop^{\infty,(0,2\beta_+,2\beta_\cK)} = x_\sscri^2\rho_+^2\,x_\sscri^2\rho_+^{2(\beta_+-1)}\rho_\cK^{2\beta_\cK}\cC_\bop^\infty;
  \]
  this satisfies Definition~\ref{DefSDWAdm} with $\ell_\sscri=\frac12$, $\ell_+=2(\beta_+-1)>0$, $\ell_\cK=2\beta_\cK>0$. Therefore, Theorem~\ref{ThmF} is applicable, with $\alpha_\sscri$ and $\alpha_+$ there equal to $-\frac12-\frac{\eps}{2}$ and $-1-\eps$ for any $\eps\in(0,\frac12)$; using unweighted b-densities, this states that if
  \[
    P_u\dot u=h,
  \]
  then
  \begin{equation}
  \label{EqDpTame0}
    \|\dot u\|_{\Hb^{k,(1-\frac{\eps}{2},1-\eps,\frac12)}(\Omega)^{\bullet,-}} \leq C_{k,\eps}\Bigl( \|h\|_{\Hb^{k+d,(2-\frac{\eps}{2},3-\eps,\frac12)}(\Omega)^{\bullet,-}} + |u|_{k+d} \|h\|_{\Hb^{d,(2-\frac{\eps}{2},3-\eps,\frac12)}(\Omega)^{\bullet,-}}\Bigr),
  \end{equation}
  where we write $|\cdot|_k$ for the norm on $B^k$. Our subsequent arguments improve the decay of $\dot u$, with the goal of proving $\dot u\in B^\infty$ when $h\in\bfB^\infty$. They involve two types of operations:
  \begin{enumerate}
  \item multiplication operators, e.g., when controlling the ``nonlinear'' term $(P_u-\Box_g)\dot u$, for which tame estimates are provided by Lemma~\ref{LemmaA1Nonlin};
  \item the inversion of $\Box_g$, which loses only a fixed finite number of b-derivatives with each application by the results of~\S\ref{SsDRes}.
  \end{enumerate}
  With~\eqref{EqDpTame0} as the starting point, tame estimates for norms of $\dot u$ in spaces with stronger decay are therefore automatic. In the sequel, we shall thus focus exclusively on improving the weights of $\dot u$, and omit the statement of corresponding tame estimates.

  \pfstep{Improved decay at $\iota^+\cup\cK^+$ (i.e., as $t_*\to\infty$).} Given~\eqref{EqDpTame0}, we have
  \begin{equation}
  \label{EqDpNormEq}
  \begin{split}
    \Box_g\dot u = h + 3 u^2\dot u =: h' &\in \Hb^{\infty,(\beta_\sscri+1,\beta_++2,\beta_\cK)} + \cC_\bop^{\infty,(2,2\beta_+,2\beta_\cK)}\Hb^{\infty,(1-\frac{\eps}{2},1-\eps,\frac12)} \\
      &\subset \Hb^{\infty,(\beta_\sscri+1,\beta_++2,\beta_\cK)}(\Omega,\mu_\bop)^{\bullet,-}.
  \end{split}
  \end{equation}
  As shown in the proof of Theorem~\ref{ThmA1Adm}, the forward solution $\dot u$ of this equation is given by the inverse Fourier transform on the real axis, which we split as
  \begin{equation}
  \label{EqDpdotuSplit}
  \begin{split}
    &\dot u = \dot u_{\rm lo} + \dot u_{\rm hi}, \\
    &\qquad \dot u_{\rm lo} := \cF^{-1}\bigl( \chi(\sigma)\wh{\Box_g}(\sigma)^{-1}\wh{h'}(\sigma)\bigr),  \\
    &\qquad \dot u_{\rm hi} := \cF^{-1}\bigl( (1-\chi(\sigma))\wh{\Box_g}(\sigma)^{-1}\wh{h'}(\sigma)\bigr);
  \end{split}
  \end{equation}
  here $\chi\in\CIc((-1,1))$ equals $1$ on $[-\frac12,\frac12]$. Corollary~\ref{CorDResHiLoc} (with $\beta'_\iota$ and $\beta'_\cK$ equal to $\beta_+$ and $\beta_\cK$ in present notation) gives
  \begin{equation}
  \label{EqDpSolvuhi}
    \dot u_{\rm hi} \in \Hb^{\infty,(\gamma,\beta_+,\beta_\cK)}(\sfM)
  \end{equation}
  for some $\gamma\in\R$, which we allow to change from line to line in the remainder of the proof.

  At low energies on the other hand, Proposition~\ref{PropDFTFourier}\eqref{ItDFTFourierLo} (which is applicable since $(\beta_++2)-\min(\beta_\cK,1-\eps)=3+\eps_+-(1-\eps)=2+\eps_++\eps\leq 2+\eps_\sscri=\beta_\sscri+1$ for small enough $\eps>0$) shows that
  \[
    \wh{h'}|_{\sigma\in\pm[0,1]} \in \Hb^{\infty,(\gamma,\beta_++1,((0,0),\beta_\cK-1))}(X_\scbtop^\pm).
  \]
  We control $\wh{\Box_g}(\sigma)^{-1}\wh{h'}(\sigma)$ for $\sigma\in\pm[0,1]$ in two steps. \textit{The first step} is to solve away $\wh{h'}(0)$ by inverting the zero energy operator of $\Box_g$; note that
  \[
    \wh{h'}(0) \in \bar H_\bop^{\infty,\beta_++1}(X),
  \]
  and thus Theorem~\ref{ThmSp0} (expressed using an unweighted b-density, causing a shift in the weight of $\frac32$) gives
  \begin{subequations}
  \begin{equation}
  \label{EqDpSolv0a}
    \wh{\dot u}(0) = \wh{\Box_g}(0)^{-1}\wh{h'}(0) \in \bar H_\bop^{\infty,\beta_+-1}(X).
  \end{equation}
  Setting $\hat r:=\frac{|\sigma|}{\rho}$ (which is a defining function of $\zface\subset X_\scbtop$), we then estimate the remainder
  \begin{equation}
  \label{EqDpSolv0b}
    \wh{h''}(\sigma) := \wh{h'}(\sigma) - \wh{\Box_g}(\sigma)\bigl(\chi(\hat r)\wh{\dot u}(0)\bigr)
  \end{equation}
  by
  \begin{equation}
  \label{EqDpSolv0c}
    \bigl(\wh{h''}(\sigma)\bigr)\big|_{\sigma\in\pm[0,1]} \in \Hb^{\infty,(\gamma,\beta_++1,\beta_\cK-1)}(X_\scbtop^\pm)
  \end{equation}
  Indeed, since $\wh{\Box_g}(\sigma)\in\Diff_\scbtop^2(X)\subset\rho_\tface^2\Diff_\bop^2(X_\scbtop^\pm)$ (see~\eqref{EqSSSpecFamMem}), we have $\wh{\Box_g}(\sigma)(\chi(\hat r)\wh{\dot u}(0))\in\Hb^{\infty,(\infty,\beta_++1,(0,0))}(X_\scbtop^\pm)$ (by which we mean the space~\eqref{EqDFTHbXphg}, intersected over all $\gamma_\zface$), and since the restriction of this to $\zface$ cancels $\wh{h'}(0)$, we obtain~\eqref{EqDpSolv0c}.

  \textit{The second step} in controlling $\dot u_{\rm lo}$ is to apply $\wh{\Box_g}(\sigma)^{-1}$ to~\eqref{EqDpSolv0c}. Since $(\beta_++1)-(\beta_\cK-1)=2+\eps_+-\eps_\cK\in(0,1)+2$, with $(0,1)$ being the indicial gap of $\wh{\Box_g}(0)$, Proposition~\ref{PropDResLo} gives
  \begin{equation}
  \label{EqDpSolv0d}
    \wh{\Box_g}{}^{-1}\wh{h''} \in \Hb^{\infty,(\gamma,\beta_+-1,\beta_\cK-1)}(X_\scbtop^\pm).
  \end{equation}
  Recalling~\eqref{EqDpSolv0b}, we have thus shown
  \begin{equation}
  \label{EqDpSolv0e}
    \wh{\Box_g}{}^{-1}\wh{h'} \in \Hb^{\infty,(\gamma,\ \beta_+-1,\ ((0,0),\beta_\cK-1))}(X_\scbtop^\pm).
  \end{equation}
  \end{subequations}
  By Proposition~\ref{PropDFTInv}\eqref{ItDFTInvLo}, we conclude that
  \[
    \dot u_{\rm lo} \in \Hb^{\infty,(\gamma,\beta_+,\beta_\cK)}(\sfM).
  \]
  Adding this to~\eqref{EqDpSolvuhi} gives $\dot u \in \Hb^{\infty,(\gamma,\beta_+,\beta_\cK)}(\sfM)$. Since $\dot u$ is supported in $t_*\geq 1$, this implies
  \begin{equation}
  \label{EqDpuGoodt}
    \dot u \in \Hb^{\infty,(\gamma,\beta_+,\beta_\cK)}(\Omega,\mu_\bop)^{\bullet,-}.
  \end{equation}

  \pfstep{Improved decay at $\scri^+$; conclusion.} Proposition~\ref{PropDscriRad} (with $\beta_\sscri$, $\beta_\sscri'$, and $\beta_+$ there equal to $\gamma$, $\beta_\sscri$, and $\beta_+$ in present notation) improves~\eqref{EqDpuGoodt} to
  \[
    \dot u\in\Hb^{\infty,(((1,0),\beta_\sscri),\ \beta_+,\ \beta_\cK)}(\Omega,\mu_\bop)^{\bullet,-}.
  \]
  We have thus shown that the forward solution operator for $P_u$ maps $\bfB^\infty$ to $B^\infty$ (with tame estimates). Theorem~\ref{ThmA1NM} now implies the solvability of~\eqref{EqDpEq} with $u\in B^\infty$. Local uniqueness near $t_*=0$ implies that, in fact, $t_*\geq 1$ on $\supp u$. The proof is complete.
\end{proof}

\begin{rmk}[Asymptotics for linear scalar waves]
\label{RmkDpLinear}
  The arguments used in the proof of Theorem~\ref{ThmDp} imply the following result for forward solutions $u$ of the linear scalar wave equation $\Box_g u=f$: if $1<\beta_\cK<\beta_+<\beta_\sscri<2$, then
  \[
    f\in\Hb^{\infty,(\beta_\sscri+1,\beta_++2,\beta_\cK)}(\Omega,\mu_\bop)^{\bullet,-} \implies u\in\Hb^{\infty,\bigl(((1,0),\beta_\sscri),\ \beta_+,\ \beta_\cK\bigr)}(\Omega,\mu_\bop)^{\bullet,-}.
  \]
  (Indeed, the restriction $\eps_\sscri<\frac12$ in Theorem~\ref{ThmDp} was only made because of the application of Proposition~\ref{PropDscriRad}, which can be relaxed to $\eps_\sscri<1$ by exploiting the membership $\tilde P\in\rho_\sscri\Diffb^2(M)$ observed after~\eqref{EqA2RadOp}. Thus, we obtain almost-$t_*^{-2}$ pointwise decay as $t_*\to\infty$. This mapping property can be used also in the settings of Theorems~\ref{ThmA1Power}, \ref{ThmA1Q}, and \ref{ThmA1Null} to obtain almost-$t_*^{-2}$-decay (for appropriately decaying $f$). This decay rate is one order off from Price's law, which in \cite{TataruDecayAsympFlat,HintzPrice} rests on a resolvent expansion at $\sigma=0$ to one more order than what we need in~\eqref{EqDpSolv0a}--\eqref{EqDpSolv0e}.
\end{rmk}

\section{Application II: a nonlinear problem related to Maxwell's equations}
\label{SA2}

In all applications in~\S\ref{SA1}, the stationary model is the scalar wave operator for the subextremal Kerr metric $g=g_{\bhm,a}$, which does not have any zero modes. We now study wave equations on 1-forms; the 1-form Hodge d'Alembertian
\[
  \Box_g = \dd\delta_g+\delta_g\dd \in \Diff^2(M^\circ;T^*M^\circ)
\]
does have a zero mode
\begin{equation}
\label{EqA2u0}
  u_0=u_0(x) = \frac{r}{\varrho^2}(\dd t_0-a\sin^2\theta\,\dd\phi_0) + \frac{r_+-r}{\mu}\,\dd r,
\end{equation}
where we use $t_0,\phi_0$ from~\eqref{EqTs3bHt0phi0} as well as $\mu,\varrho^2$ from~\eqref{EqTsBLFunc} and $r_+$ from~\eqref{EqTsBLMfd}. This is a co-closed 1-form with spatial decay $\cO(r^{-1})$. Our goal is to prove a more precise version of Theorem~\ref{ThmI1} on the small data global existence for forward problems for the (toy) semilinear wave equation
\begin{equation}
\label{EqA2}
  \Box_g u = f+|\delta_g u|^2 u,
\end{equation}
see Theorem~\ref{ThmA2} below. The particular nonlinearity was already motivated in Remark~\ref{RmkI1Prelim}.

We will show in Theorem~\ref{ThmA2Adm} that $\Box_g$ is $1$-admissible (Definition~\ref{DefSSAlephAdm}). The bounds proved in Theorem~\ref{ThmF}, therefore, do not even imply the boundedness of the solution $u$ of $\Box_g u=f$ (rather, $t_*^{\frac12}$ growth at $(\cK^+)^\circ$, cf.\ \eqref{EqFTame3} with $\aleph=1$), which means that we must work with function spaces encoding stronger decay and more precise asymptotics in order to solve~\eqref{EqA2}. We first explain on a heuristic level what the nature of these spaces should be.

\pfstep{Heuristic I: decay of $f$ and structure of $u$ at $\cK^+$.} Given source terms $f$ with pointwise decay $t_*^{-\alpha}$, $\alpha>1$ (and spatially compact support, say), one expects the forward solution of $\Box_g u=f$ to be of the form
\begin{equation}
\label{EqA2Heur1}
  u(t_*,x) = c u_0(x) + \cO(t_*^{-\alpha+1}),
\end{equation}
at least for some range of $\alpha$. The loss of one order of decay is due to the presence of a zero mode, as arising already in the ODE example discussed around~\eqref{EqIODE}. For such $u$ then, the linearization
\begin{equation}
\label{EqA2Pu}
  P_u\dot u := \Box_g \dot u - 2 u \Re\bigl(\delta_g u\,\ol{\delta_g \dot u}\bigr) - |\delta_g u|^2 \dot u
\end{equation}
differs from $P_0=\Box_g$ by an operator that includes a term roughly of the form $u(\pa_x u)\pa_x=\cO(t_*^{-\alpha+1})\pa_x$. When $\alpha\leq 2$, this is too little decay for $P_u$ to admit a good solvability theory; we need more than $1$ order of $t_*$-decay (cf.\ \eqref{EqSDWAdmWeights} with $\aleph=1$) for Theorem~\ref{ThmF} to apply to $P_u$. But even if we take $\alpha>2$, then with $u$ having a $\cO(t_*^{-\alpha+1})$ remainder, the source term in the next step of the iteration (which is $\Box_g u-|\delta_g u|^2 u$ minus the $f$ in~\eqref{EqA2}) only has $\cO(t_*^{-\alpha+1})$ decay itself (with $\Box_g$ annihilating the first term in~\eqref{EqA2Heur1} but not improving the decay of the second). One cannot afford this loss of $1$ power of $t_*$-decay at each step.

To make progress, we need to sharpen~\eqref{EqA2Heur1} to
\begin{equation}
\label{EqA2Heur1b}
  u(t_*,x) = c u_0(x) + a(t_*)u_0(x) + \cO(t_*^{-\alpha}),\quad a=\cO(t_*^{-\alpha+1}),
\end{equation}
where the second $\cO$-term now has the same decay as $f$. (This more precise description is consistent with the introductory ODE example in~\eqref{EqIODEDecomp}.) For $\alpha=1+\eps_\cK$, $\eps_\cK>0$, then, the difference $P_u-\Box_g$ is of the schematic form (using~\eqref{EqA2Pu} and $\delta_g u_0=0$)
\[
  u \cdot a'(t_*) \iota_{\nabla t_*}u_0\,\delta_g + a'(t_*)^2 |\iota_{\nabla t_*}u_0|^2 = \cO(t_*^{-\alpha}),
\]
where we must exploit the b-regularity of $a$ in order to obtain the $\cO(t_*^{-\alpha})$-decay of first summand. Therefore, $P_u$ is now admissible; and also the source term arising from~\eqref{EqA2Heur1b} in a subsequent iteration step has acceptable $t_*$-decay since
\begin{equation}
\label{EqA2Heur1bFwd}
  \Box_g(c u_0+a(t_*)u_0)=[\Box_g,a]u_0=a'(t_*)[\Box_g,t_*]\ftrans(0)u_0+\frac12 a''(t_*)[[\Box_g,t_*],t_*]u_0
\end{equation}
decays like $a'(t_*)=\cO(t_*^{-\alpha})$.

\pfstep{Heuristic II: decay of $u$ and $f$ at $\iota^+$ and $\scri^+$.} Since we need to control $u$ at $\cK^+$ modulo a ``structureless'' $t_*^{-1-\eps_\cK}$ remainder (see~\eqref{EqA2Heur1b} with $\alpha=1+\eps_\cK$), we need to work with spaces whose structureless remainders have at least the same amount of decay at $\iota^+$ and $\scri^+$. (This is because one typically expects the decay of waves at $\cK^+$ to be inherited from that at $\iota^+$ plus contributions from zero modes, and the decay at $\iota^+$ to be inherited from that at $\scri^+$ by propagation.) So $u$ should be captured in a space with ``structureless'' remainders of size $\cO(\rho_\sscri^{1+\eps_\sscri}\rho_+^{1+\eps_+}\rho_\cK^{1+\eps_\cK})$ where $0<\eps_\cK<\eps_+<\eps_\sscri$. Correspondingly, the source terms $f$ we shall consider have the usual additional weight $\rho_\sscri\rho_+^2\rho_\cK^0$, i.e.,
\[
  f = \cO(\rho_\sscri^{2+\eps_\sscri}\rho_+^{3+\eps_+}\rho_\cK^{1+\eps_\cK}).
\]
As already remarked in~\S\ref{SsDp}, for such $f$, the solution $u$ of $\Box_g u=f$ will not decay like $\rho_\sscri^{1+\eps_\sscri}$ (which would mean $\cO(r^{-1-\eps_\sscri})$ at $(\scri^+)^\circ$), but rather has a radiation field, i.e., a $r^{-1}u_\sscri(t_*,\omega)$ leading order term at $\scri^+$ plus a faster-decaying remainder. Recall also that the $r^{-1}$ leading order term is annihilated to leading order at $\scri^+$ by the forward operator $\Box_g u-|\delta_g u|^2 u$.

Moreover, the term $c u_0(x)$ in~\eqref{EqA2Heur1b} is of class $r^{-1}\CI(X)$ (Definition~\ref{DefCSpatial}); while one might initially only want to take this term seriously at $\cK^+$, it is convenient to cut it off to a larger set, namely, to a neighborhood of $\iota^+\cup\cK^+$ (using a $t_*$-cutoff). From the perspective of the forward operator, this is, in fact, forced on us, since in the expression~\eqref{EqA2Heur1bFwd}, the term $[\Box_g,t*]\ftrans(0)\approx -2 \rho(\rho\pa_\rho-1)$ maps $u_0\in\rho\CI(X)$ to $\rho^3\CI(X)$, i.e., annihilates it to leading order at $\rho^{-1}(0)\supset\iota^+$. (In other words, $a(t_*)u_0(x)$ is not only an approximate solution of $\Box_g$ at $\cK^+$, but in fact at $\cK^+\cup\iota^+$.)

\bigskip

In the remainder of this section, we will make these heuristics rigorous; we start with the precise statement of our main result motivated by these considerations. As in~\S\ref{SsDp}, we use an unweighted b-density $\mu_\bop$ on spacetime such as~\eqref{EqA1Powermub}. We shall measure the pointwise size of a 1-form using the Euclidean metric (or any other smooth positive definite fiber metric) on the bundle $\cT^*\to M$ with frame $\dd t_*,\dd x$ (Definition~\ref{DefCTscPullback}).

\begin{thm}[Small data global existence]
\label{ThmA2}
  Let $0<\eps_\cK<\eps_+<\eps_\sscri<\frac12$, and let
  \begin{equation}
  \label{EqA2Gamma}
    \beta_\sscri=1+\eps_\sscri,\quad
    \beta_+=1+\eps_+,\quad
    \beta_\cK=1+\eps_\cK.
  \end{equation}
  Recall the defining functions $\rho_\sscri,\rho_+,\rho_\cK$ from~\eqref{EqFWeights}. Fix cutoff functions
  \begin{align*}
    &\chi_T\in\CI([1,\infty)),\quad \chi_T=0\ \text{near}\ 1,\quad \chi_T=1\ \text{near}\ \infty; \\
    &\chi_\sscri\in\CI(M),\ \text{supported in a collar neighborhood of $\scri^+$},\quad \chi_\sscri=1\ \text{near}\ \scri^+.
  \end{align*}
  Then there exist $d\in\N$ and $\eps>0$ such that for all $f\in\Hb^{\infty,(\beta_\sscri+1,\beta_++2,\beta_\cK)}(\Omega,\mu_\bop;\cT^*)^{\bullet,-}$ with
  \[
    \|f\|_{\Hb^{d,(\beta_\sscri+1,\beta_++2,\beta_\cK)}}^2 = \iint_\Omega |\rho_\sscri^{-\beta_\sscri-1}\rho_+^{-\beta_+-2}\rho_\cK^{-\beta_\cK}(t_*\pa_{t_*},r\pa_x)^{\leq d}f|^2\,\frac{\dd t_*\,\dd x}{t_* r^3} <\eps^2,
  \]
  the equation
  \[
    \Box_g u = f + |\delta_g u|^2 u
  \]
  has a unique global forward solution $u$ on $\Omega$ which is of the form
  \begin{align}
  \label{EqA2uExp}
    &u(t_*,x) = \chi_T(t_*)\bigl(c u_0(x) + a(t_*)u_0(x)\bigr) + \chi_\sscri r^{-1}u_\sscri\Bigl(t_*,\frac{x}{|x|}\Bigr) + \tilde u(t_*,x), \\
    &\qquad c\in\C,\ a\in\dot H_\bop^{\infty,\eps_\cK}([1,\infty],|\tfrac{\dd t_*}{t_*}|), \nonumber\\
    &\qquad u_\sscri\in\dot H_\bop^{\infty,\eps_+}(\scri^+\cap\{t_*\geq 1\},|\tfrac{\dd t_*}{t_*}\,\dd\slg|;\cT^*), \nonumber\\
    &\qquad \tilde u \in H_\bop^{\infty,(\beta_\sscri,\beta_+,\beta_\cK)}(\Omega,\mu_\bop;\cT^*)^{\bullet,-}. \nonumber
  \end{align}
  In particular, we have the pointwise bounds
  \begin{equation}
  \label{EqA2uBds}
    |a|\leq C t_*^{-\eps_\cK},\quad
    |r^{-1}u_\sscri|\leq C\rho_\sscri\rho_+^{1+\eps_+},\quad
    |\tilde u|\leq C\rho_\sscri^{1+\eps_\sscri}\rho_+^{1+\eps_+}\rho_\cK^{1+\eps_\cK},
  \end{equation}
  similarly for all b-derivatives.
\end{thm}

Different choices of cutoff functions lead to the same description of $u$ except for a re-definition of $\tilde u$. The cutoff function $\chi_T(t_*)$ can be omitted; we include it to emphasize that the terms involving $u_0$ only capture the late-time behavior of $u$. The cutoff function $\chi_\sscri$ is more important: without it, the term $r^{-1}u_\sscri$ would ``spill over'' to $\cK^+$ as a $\cO(t_*^{-\eps_+})$ contribution (which would not be negligible, as it decays less than $\tilde u$). In the notation of Definition~\ref{DefDscriPhg}, the final two terms in~\eqref{EqA2uExp} can be combined into a single element of $\Hb^{\infty,(((1,0),\beta_\sscri),\beta_+,\beta_\cK)}(\Omega)^{\bullet,-}$. Theorem~\ref{ThmI1} is an immediate consequence of Theorem~\ref{ThmA2}; note that the decay rate $t_*^{-\eps_\cK}r^{-1}\sim\rho_\sscri^1\rho_+^{1+\eps_\cK}\rho_\cK^{\eps_\cK}$ in~\eqref{EqI1Baby} is implied by the bounds~\eqref{EqA2uBds}.

\subsection{Admissibility of the 1-form wave operator}
\label{SsA2Adm}

Recalling Definition~\ref{DefSSAlephAdm}, we prove:

\begin{thm}[$1$-admissibility of $\Box_g$]
\label{ThmA2Adm}
  Let $\alpha_\sface\in(-\frac32,-\frac12)$ and $\delta>\alpha_\sface+\frac32$. Then the Hodge d'Alembertian (or tensor wave operator) $\Box_g$ on 1-forms is $1$-admissible with $\sface$-weight $\alpha_\sface$ and $\sface$-loss $\delta$.
\end{thm}

The structure of $\Box_g$ was described in~\eqref{EqSSAdmBox}; in particular, $\ubar S=0$. We first recall basic spectral information about $\Box_g$. Let $\sfs\in\CI(\Stb^*M_0)$ be a stationary-$\Box_g$-admissible order function with weights $\alpha_\sface,-1$ and margin $0$ (Definition~\ref{DefSSAlephAdm}); for later use, we record that by~\eqref{EqSSOrderAdmOut} this entails the threshold condition
\begin{equation}
\label{EqA2AdmThr}
  \sfs|_{\pa\cR_{\pa\cK^+,{\rm out}}}+\alpha_\sface < -\frac32.
\end{equation}

By Proposition~\ref{PropSptfAdm}, $\Box_g$ is $\tface$-admissible with all weights $\beta\in(0,1)$. Next, \cite[Theorem~3.5]{HintzGlueLocIII} (which is almost directly taken from \cite[Theorem~5.1]{AnderssonHaefnerWhitingMode}), shows that the operator $\wh{\Box_g}(\sigma)$ is invertible for $\sigma\in\C$, $\Im\sigma\geq 0$, $\sigma\neq 0$, as a map~\eqref{EqSpBMap} (with the orders there being the ones induced by $\sfs$); and the zero energy operator $\wh{\Box_g}(0)$, which is Fredholm of index $0$ as a map~\eqref{EqSp0Map} for any $\alpha\in(-\frac32,-\frac12)$ by Corollary~\ref{CorSpLoInd0}, has $1$-dimensional nullspace spanned by $u_0$ in~\eqref{EqA2u0}. (With respect to an unweighted b-density on $X$, the target space in~\eqref{EqSp0Map} has weight $\alpha+2+\frac32\in(2,3)$, which roughly speaking corresponds to pointwise $r^{-\gamma-2}$ decay where $\gamma\in(2,3)$, while elements of the domain roughly have $r^{-\gamma}$ decay.) The nullspace of its formal adjoint is (in the notation used in~\eqref{EqSp0Map})
\begin{equation}
\label{EqA2Admu0star}
  \ker_{\dot H_\bop^{-\sfs_0+1,-\alpha-2}(X;\cT^*)} \wh{\Box_g}(0)^* = \mathspan\{ u_0^* \},\quad u_0^* := \delta(r-r_+)\,\dd r.
\end{equation}
Furthermore, \cite[Lemma~3.8]{HintzGlueLocIII} and the fact that $[\Box_g,t_*]\ftrans(0)=i\pa_\sigma\wh{\Box_g}(0)$ imply that the $L^2(X^\circ;\cT^*)$-pairing (with respect to the fiber inner product induced by $g$)
\begin{equation}
\label{EqA2Adm0Pair}
  \la\pa_\sigma\wh{\Box_g}(0)u_0,u_0^*\ra = -4\pi i
\end{equation}
does not vanish.

We prepare the proof following the ideas around \cite[(3.30)]{HintzGlueLocIII} by studying the (singular) behavior of $\wh{\Box_g}(\sigma)^{-1}$ near $\sigma=0$ by passing to a rank $1$ augmentation (Grushin problem) of $\wh{\Box_g}(\sigma)$, namely
\begin{equation}
\label{EqA2AdmAug}
  \wt\Box(\sigma) := \begin{pmatrix} \wh{\Box_g}(\sigma) & \wh{\Box_g}(\sigma)(\sigma^{-1}u_0) \\ \la\cdot,f^*\ra & 0 \end{pmatrix},
\end{equation}
where $f^*\in\CIc(X^\circ;\cT^*)$ is chosen such that $\la u_0,f^*\ra\neq 0$. The $(1,2)$ entry of $\wt\Box(\sigma)$ decays at a rate consistent with the target space of $\wh{\Box_g}(0)$ (and thus of $\wh{\Box_g}(\sigma)$ at the transition face in the low-energy analysis):

\begin{lemma}[Decay of augmentation]
\label{LemmaA212}
  The $(1,2)$ entry of $\wt\Box(\sigma)$ is of class $\rho^3\CI(X;\cT^*)+\sigma\rho^3\CI(X;\cT^*)$. Its restriction to $\sigma=0$ is $\pa_\sigma\wh{\Box_g}(0)u_0$.
\end{lemma}
\begin{proof}
  Since $\wh{\Box_g}(0)u_0=0$, we have
  \[
    \wh{\Box_g}(\sigma)(\sigma^{-1}u_0) = \pa_\sigma\wh{\Box_g}(0)u_0 + \frac12\sigma\pa_\sigma^2\wh{\Box_g}(0)u_0.
  \]
  For the second term, we observe that $\pa_\sigma^2\wh{\Box_g}(0)=2 g^{0 0}\in\rho^2\CI(X)$ in the notation of~\eqref{EqSSAdmOp} and using~\eqref{EqSSSpecFam}. For the first term, we note that
  \[
    \pa_\sigma\wh{\Box_g}(0) = 2 i\rho(\rho\pa_\rho-1) - i Q \equiv 2 i\rho(\rho\pa_\rho-1) \bmod \rho^3\Diffb^1(X;\cT^*)
  \]
  is an element of $\rho\Diffb^1$ with indicial root $1$, which thus annihilates the $\rho^1$ leading order part of $u_0$. This yields the claim.
\end{proof}

\begin{lemma}[Zero energy operator of augmentation]
\label{LemmaA2Adm0Op}
  With $\sfs_0$ denoting the order induced by $\sfs$ at zero frequency as in Theorem~\usref{ThmSp0}, the operator
  \begin{equation}
  \label{EqA2Adm0Op}
    \wt\Box(0) \colon \bigl\{ u\in\bar H_\bop^{\sfs_0,\alpha_\sface}\colon \wh{\Box_g}(0)u\in\bar H_\bop^{\sfs_0-1,\alpha_\sface+2} \bigr\} \oplus C \to \bar H_\bop^{\sfs_0-1,\alpha_\sface+2} \oplus \C
  \end{equation}
  is invertible.
\end{lemma}
\begin{proof}
  Lemma~\ref{LemmaA212} shows that~\eqref{EqA2Adm0Op} is well-defined. It has Fredholm index $0$ since $\wh{\Box_g}(0)$ does. We only need to show that it is injective. Suppose that $\wt\Box(0)(u,c)=(0,0)$. Then $\wh{\Box_g}(0)u+c\pa_\sigma\wh{\Box_g}(0)u_0=0$; taking the $L^2$-inner product with $u_0^*$ from~\eqref{EqA2Admu0star}, this gives
  \[
    0 = \la u,\wh{\Box_g}(0)^*u_0^*\ra + c\la\pa_\sigma\wh{\Box_g}(0)u_0,u_0^*\ra = -4\pi i c,
  \]
  so $c=0$. Therefore $\wh{\Box_g}(0)u=0$, and hence $u=c'u_0$ for some $c'\in\C$. But then $0=\la u,f^*\ra=c'\la u_0,f^*\ra$ is a nonzero multiple of $c'$, so $c'=0$ and hence $u=0$.
\end{proof}

Therefore, Theorem~\ref{ThmSpLo} is applicable to $\wt\Box(\sigma)$ (cf.\ Remark~\ref{RmkSpLoAug}). It gives, in particular, uniform bounds on the second component of $\wt\Box(\sigma)^{-1}(\hat f(\sigma),0)$, whose inverse Fourier transform will require some care to estimate; the following result will be useful.

\begin{lemma}[Estimates for expansion terms]
\label{LemmaA2AdmExp}
  We fix the volume density $|\dd t_*|$ on $\R_{t_*}$. For $q,k\in\N_0$, write $H_\cuop^q(\ol\R):=H^q(\R)$ and
  \begin{equation}
  \label{EqA2AdmExpHcub}
    H_{\cuop;\bop}^{(q;k)}(\ol\R) := \bigl\{ u\in H_\cuop^q(\ol\R) \colon (\la t_*\ra\pa_{t_*})^j u\in H_\cuop^q(\ol\R)\ \forall\,j=0,\ldots,k \bigr\}.
  \end{equation}
  Write $\dot H_{\cuop;\bop}^{(q;k)}([1,\infty])$ for the subspace of elements supported in $\{t_*\geq 1\}$. Let $\alpha\in\R$ and $u_0\in\cA^\alpha(X)=\rho^\alpha\CI_\bop(X)$. Let $m\in\N_0$, and suppose that $a=a(t_*)$ is a distribution on $\R_{t_*}$ satisfying
  \[
    \supp a\subset\{t_*\geq 1\},\quad \pa_{t_*}^m a\in\dot H_{\cuop;\bop}^{(q;k)}([1,\infty]).
  \]
  Let $\beta_\sface\in\R$, and let $\sfs\in\CI(\Stb^*M_0)$ be a $t_*$-translation-invariant order function on $M_0$, such that
  \begin{equation}
  \label{EqA2AdmExpThr}
    \max\bigl(\beta_\sface,\ \sfs|_{\pa\cR_{\pa\cK^+,{\rm out}}}+\beta_\sface\bigr) < (-m) + \Bigl(-\frac32+\alpha\Bigr),
  \end{equation}
  where we recall Definition~\usref{DefTs3bRad}. Suppose that $q\geq\sup\sfs$. Then
  \begin{equation}
  \label{EqA2AdmExpMem}
    a(t_*)u_0(x) \in H_{\tbop;\bop}^{(\sfs;k),(\beta_\sface,-m)}(\Omega,|\dd t_*\,\dd x|)^{\bullet,-},
  \end{equation}
  with norm bounded by $C\|\pa_{t_*}^m a\|_{\dot H_{\cuop;\bop}^{(q;k)}([1,\infty])}$.
\end{lemma}
\begin{proof}
  This is a mild generalization of \cite[Lemma~2.9]{HintzGlueLocIII} (albeit with notation closer to the present paper); we give the proof for completeness and easier reference.

  We can test for the b-regularity of $a(t_*)u_0(x)$ using $t_*\pa_{t_*}$ and $r\pa_x$, which act only on $a(t_*)$ and $u_0(x)$, respectively; therefore, it suffices to consider the case $k=0$. Moreover, we may divide $u_0$ by $\rho^\alpha$ and thus reduce to the case $\alpha=0$, so $u_0\in\CI_\bop(X)\subset\Hb^{\infty,-\frac32-\eta}(X;|\dd x|)$ for all $\eta>0$; moreover, $\beta_\sface$ and $\sfs|_{\pa\cR_{\pa\cK^+,{\rm out}}}+\beta_\sface$ are $<-\frac32$, so $<-\frac32-2\eps$ for some small $\eps>0$.

  \pfstep{The case $m=0$.} We use Lemma~\ref{LemmaMUetbFT} (with $\alpha_\sface$ there equal to $\beta_\sface$) to prove~\eqref{EqA2AdmExpMem}. Concretely, writing $u(t_*,x):=a(t_*)u_0(x)$, consider first low frequencies $\sigma\in[-1,1]$: we estimate
  \[
    \|\hat u(\sigma)\|_{H_{\scbtop,|\sigma|}^{\sfs,(\sfs+\beta_\sface,\beta_\sface,0)}(X,|\dd x|)} \leq \|\hat a(\sigma)u_0\|_{H_{(\scbtop,|\sigma|);\bop}^{(\sfs;j),(\sfs+\beta_\sface,\beta_\sface,0)}} \leq C\|\hat a(\sigma)u_0\|_{H_{\scbtop,|\sigma|}^{\sfs,(s_0+\beta_\sface,\beta_\sface,0)}};
  \]
  here $s_0:=\sup_{\pa\cR_{\pa\cK^+,{\rm out}}}\sfs+\eps$, and $j\in\N_0$ (which can be inserted here since the additional b-regularity only refers to the spatial part $u_0$, which is conormal) is chosen such that $\sfs+\beta_\sface+j\geq s_0+\beta_\sface$ on $\ol{{}^\scbtop T^*_\scface}X$ (recalling~\eqref{EqMUscbtSingle} for the notation). The second estimate then follows from the fact that $j$ degrees of b-regularity imply $j$ additional degrees of decay in the complement of the zero section of ${}^\scbtop T^*X$ over $\scface$. (This was already used in the proof of Lemma~\ref{LemmaDRes}.) We can then further bound this by
  \[
    C'|\hat a(\sigma)| \|u_0\|_{\Hb^{s_+,-\frac32-\eps}} \leq C''|\hat a(\sigma)|
  \]
  where we fix $\N_0\ni s_+\geq\sup\sfs$; we use in the first estimate that b-regularity implies sc-b-transition-regularity, which implies a continuous inclusion of the corresponding Sobolev spaces. It remains to note that $\int_{-1}^1|\hat a(\sigma)|^2\,\dd\sigma\leq C_q\|a\|_{H_\cuop^q}^2$ for all $q$.

  For the high-frequency terms in Lemma~\ref{LemmaMUetbFT}, completely analogous arguments yield the estimate
  \[
    \|\hat u(\sigma)\|_{H_{\scop,|\sigma|^{-1}}^{\sfs,\sfs+\beta_\sface,\sfs}} \leq C|\sigma|^{\sup\sfs}|\hat a(\sigma)|
  \]
  for $\sigma\in\pm[1,\infty)$, with the extra power of $|\sigma|$ being due to the semiclassical order $\sfs$ on the left-hand side. This expression lies in $L^2(\pm[1,\infty)_\sigma)$ provided the Sobolev order $q$ of $a$ satisfies $q\geq\sup\sfs$.

  \pfstep{The case $m\geq 1$.} We apply Lemma~\ref{LemmaA2AdmInt} below with $u=(\pa_{t_*}^{m-\ell}a)u_0$ and $\beta=-\frac32-(\ell-1)-\eps$ for $\ell=1,2,\ldots,m-1$ to obtain~\eqref{EqA2AdmExpMem}.
\end{proof}

To complete the proof, we need:
\begin{lemma}[Integration of $\pa_{t_*}$ on 3b-Sobolev spaces]
\label{LemmaA2AdmInt}
  Let $\sfs\in\CI(\Stb^*M_0)$ be a $t_*$-translation-invariant 3b-differential order function. Let $\beta\in\R$ and $\ell>\frac12$. Then
  \begin{equation}
  \label{EqA2AdmInt}
    \pa_{t_*}u =: f \in H_\tbop^{\sfs,(\beta,-\ell+1)}(M_0),\ \ t_*\geq 1\ \text{on}\ \supp u \implies u \in H_\tbop^{\sfs,(\beta-1,-\ell)}(M_0),
  \end{equation}
  where we use the Minkowskian volume density $|\dd t\,\dd x|$ to define the underlying $L^2$-spaces.
\end{lemma}
\begin{proof}
  Consider first the case $\sfs=0$. We pass to $t=t_*+r$, so $t\geq r+1$ on $\supp u$; hence $\rho_\sface=r^{-1}$ and $\rho_\cK=\frac{r}{t}$ are local defining functions of $\sface$ and $\cK^+\subset M_0$, respectively. We drop the angular variables from the notation. We have $f(t,r)=\rho_\sface^\beta\rho_\cK^{-\ell+1}f_0=r^{-\beta-\ell+1}t^{\ell-1}f_0$ with $f_0=f_0(t,r)\in L^2$. Fix any $\alpha\in(-\frac12,\ell-1)$; then we can estimate
  \begin{align*}
    &\int_\bhm^\infty \int_{r+1}^\infty \,\biggl| r^{\beta+\ell-1}t^{-\ell} \int_{r+1}^t r^{-\beta-\ell+1}s^{\ell-1}f_0(s,r)\,\dd s\biggr|^2\,\dd t\,\dd r \\
    &\qquad = \int_\bhm^\infty \int_{r+1}^\infty \,\biggl| t^{-\ell}\int_{r+1}^t s^{\ell-1}f_0(s,r)\,\dd s\biggr|^2\,\dd t\,\dd r \\
    &\qquad \leq \int_\bhm^\infty \int_{r+1}^\infty t^{-2\ell}\biggl(\int_{r+1}^t s^{2\alpha}\,\dd s\biggr)\biggl(\int_{r+1}^t s^{2\ell-2-2\alpha}|f_0(s,r)|^2\,\dd s\biggr)\,\dd t\,\dd r \\
    &\qquad \lesssim \int_\bhm^\infty\int_{r+1}^\infty s^{2\ell-2-2\alpha}|f_0(s,r)|^2\biggl(\int_s^\infty t^{-2\ell+2\alpha+1}\,\dd t\biggr)\,\dd s\,\dd r \\
    &\qquad \lesssim \int_\bhm^\infty\int_{r+1}^\infty |f_0(s,r)|^2\,\dd s\,\dd r.
  \end{align*}

  For general $\sfs$, let $A\in\tilde\Psi_\tbop^\sfs(M_0)$ (see~\eqref{EqMUKPsdo}) be elliptic with parametrix $B\in\tilde\Psi_\tbop^{-\sfs}(M_0)$, so $I=B A+R=B A+R(B A+R)=(B+R B)A+R R=:B'A+R R$ where $R\in\tilde\Psi_\tbop^{-\infty}(M_0)$, and $B':=B+R B\in\tilde\Psi_\tbop^{-\sfs}(M_0)$; we may arrange for $A,B$ (and thus $R$ and $B'$) to be $t_*$-translation-invariant. Apply then~\eqref{EqA2AdmInt} to $A u$ and $R u$ in place of $u$: since $\pa_{t_*}(A u)=A f$ and $\pa_{t_*}(R u)=R f$ lie in $\Htb^{0,(\beta,-\ell+1)}$, we conclude that $A u,R u\in\Htb^{0,(\beta-1,-\ell)}$; this gives $u=B'(A u)+R(R u)\in\Htb^{\sfs,(\beta-1,-\ell)}$.
\end{proof}

\begin{proof}[Proof of Theorem~\usref{ThmA2Adm}]
  By density, it suffices to prove~\eqref{EqSSAlephAdmSol} for source terms $f\in\CIc((t_*,\infty)_{t_*}\times X^\circ)$ where $t_+>1$. Since the high-energy estimates of Theorem~\ref{ThmSpHi} apply to $\wh{\Box_g}(\sigma)$, the Paley--Wiener argument in the proof of Theorem~\ref{ThmA1Adm} applies \emph{verbatim} and shows that the forward solution of $\Box_g u=f$ is given by
  \[
    u(t_*,\cdot) = \frac{1}{2\pi}\int_{\Im\sigma=\gamma} e^{-i\sigma t_*} \wh{\Box_g}(\sigma)^{-1}\hat f(\sigma)\,\dd\sigma
  \]
  for large $\gamma$, and then for all $\gamma>0$ by shifting the contour and using the invertibility of $\wh{\Box_g}(\sigma)$ for $\Im\sigma>0$. (Thus, $\hat u(\sigma)=\wh{\Box_g}(\sigma)^{-1}\hat f(\sigma)$ for $\Im\sigma>0$.) Using the continuity of $\wh{\Box_g}(\sigma)^{-1}$ down to $\R\setminus\{0\}$ (which for bounded $\sigma$ follows from Theorem~\ref{ThmSpB}\eqref{ItSpBInvReg}), we can shift the contour further to
  \begin{equation}
  \label{EqA2AdmContour}
    \gamma_- \cup \gamma_0 \cup \gamma_+,\quad
    \gamma_-=(-\infty,-\tfrac12],\ \gamma_0=\tfrac12 e^{i[\pi,0]},\ \gamma_+=[\tfrac12,\infty).
  \end{equation}
  Fix a smooth cutoff function $\chi=\chi(|\sigma|)\in\CIc([0,1))$ which equals $1$ on $[0,\frac12]$.

  \pfstep{High-frequency control.} Recalling Lemma~\ref{LemmaMUetbFTb} and~\eqref{EqA1AdmLou}, we can estimate the high-frequency part
  \[
    u_{\rm hi}(t_*,\cdot) := \frac{1}{2\pi}\int_{\gamma_-\cup\gamma_+} e^{-i\sigma t_*}(1-\chi(|\sigma|))\wh{\Box_g}(\sigma)^{-1}\hat f(\sigma)\,\dd\sigma
  \]
  by
  \begin{equation}
  \label{EqA2AdmHi}
  \begin{split}
    \|u_{\rm hi}\|_{H_{\tbop;\bop}^{(\sfs;k),(\alpha_\sface,0)}}^2 &\sim \sum_{j=j_1+j_2\leq k}\int_{\R\setminus[-\frac12,\frac12]} \|\sigma^{j_1}(\sigma\pa_\sigma)^{j_2}\wh{u_{\rm hi}}(\sigma)\|_{H_{(\scop,|\sigma|^{-1});\bop}^{(\sfs_\sigma;k-j),\sfs_\sigma+\alpha_\sface,\sfs_\sigma}}^2\,\dd\sigma \\
      &\lesssim \sum_{j=j_1+j_2\leq k}\int_{\R\setminus[-\frac12,\frac12]} \bigl\| \sigma^{j_1}(\sigma\pa_\sigma)^{j_2}\bigl(\wh{\Box_g}(\sigma)^{-1}\hat f(\sigma)\bigr)\bigr\|_{H_{(\scop,|\sigma|^{-1});\bop}^{(\sfs_\sigma;k-j),\sfs_\sigma+\alpha_\sface,\sfs_\sigma}}^2\,\dd\sigma \\
      &\leq C_k\|f\|_{H_{\tbop;\bop}^{(\sfs;k),(\alpha_\sface+1,0)}}^2.
  \end{split}
  \end{equation}
  exactly as in~\eqref{EqA1AdmHiFinal}. This uses the fact that $\Box_g$ is strongly trapping admissible (Definition~\ref{DefSSTrapAdm}), and thus high-energy estimates for it lose arbitrarily small powers of $|\sigma|$ compared to the non-trapping case; this was shown in \cite[Lemma~3.3]{HintzGlueLocIII} using \cite[\S{3.2.3}]{HintzGlueLocII} and \cite[Proposition~4.1]{HintzPsdoInner} (which, in turn, relies on the fact that parallel transport of null-vectors along the trapped set is polynomially bounded in the propagation time, first shown by Marck \cite{MarckParallelNull}).

  \pfstep{Low-frequency control.} In order to control $\frac{1}{2\pi}\int_{\gamma_-\cup\gamma_0\cup\gamma_+}e^{-i\sigma t_*}\chi(|\sigma|)\wh{\Box_g}(\sigma)^{-1}\hat f(\sigma)\,\dd\sigma$, let us write
  \begin{equation}
  \label{EqA2AdmLo}
    \wt\Box(\sigma)^{-1} \bigl( \hat f(\sigma), 0 \bigr) =: \bigl( \hat u_{\rm reg}(\sigma), \hat a_1(\sigma) \bigr) \implies \hat u(\sigma) = \hat u_{\rm reg}(\sigma) + \sigma^{-1}\hat a_1(\sigma)u_0;
  \end{equation}
  note that $\hat u_{\rm reg}(\sigma)$ (in the distributional sense, i.e., when paired with any fixed $\phi\in\CIc(X^\circ;\cT^*)$) and $\hat a_1(\sigma)$ are holomorphic in $\sigma$ for $|\sigma|\leq 1$ and $\Im\sigma>0$ and continuous down to $[-1,1]\setminus\{0\}$. The bounds for the inverse of the augmented operator~\eqref{EqA2AdmAug} moreover show that they are uniformly bounded down to $[-1,1]$; in fact, writing
  \[
    u_{\rm reg}(t_*,\cdot) := \frac{1}{2\pi}\int_{\gamma_-\cup\gamma_0\cup\gamma_+} e^{-i\sigma t_*}\hat u_{\rm reg}(\sigma)\,\dd\sigma,
  \]
  we can shift the part $\gamma_0$ of the integration contour to $[-\frac12,\frac12]$ within the region where $\chi(|\sigma|)$ equals $1$ (and is thus holomorphic). The final total contour is $[-1,1]$; we thus get, analogously to~\eqref{EqA1AdmLoEst},
  \begin{equation}
  \label{EqA2AdmReg}
    \|u_{\rm reg}\|_{H_{\tbop;\bop}^{(\sfs;k),(\alpha_\sface,0)}}^2 \sim \int_{-1}^1 \| (\sigma\pa_\sigma)^j\hat u_{\rm reg}(\sigma)\|_{H_{(\scbtop,|\sigma|);\bop}^{(\sfs_\sigma;k-j),(\sfs_\sigma+\alpha_\sface,\alpha_\sface,0)}}^2\,\dd\sigma \leq C_k\|f\|_{H_{\tbop;\bop}^{(\sfs-1;k),(\alpha_\sface+2,0)}}^2.
  \end{equation}

  Finally, consider the inverse Fourier transform of the singular part of $\hat u(\sigma)$,
  \begin{equation}
  \label{EqA2AdmaDef}
    a(t_*) := \frac{1}{2\pi}\int_{\gamma_0} e^{-i\sigma t_*}\sigma^{-1}\hat a_1(\sigma)\,\dd\sigma.
  \end{equation}
  Note that by Theorem~\ref{ThmSpLo}, $\cF(\pa_{t_*}a(t_*))(\sigma)=-i\hat a_1(\sigma)$, so we can shift the contour in the inverse Fourier transform for $\pa_{t_*}a$ to $[-1,1]$ and obtain $\|\pa_{t_*}a\|_{\Hb^k(\ol\R)} \leq C_k\|f\|_{H_{\tbop;\bop}^{(\sfs-1;k),(\alpha_\sface+2,0)}}$ from the low-energy bounds for $\wt\Box(\sigma)^{-1}$, where we use the volume density $|\dd t_*|$ on $\ol{\R_{t_*}}$. Since $\cF(\pa_{t_*}a)$, as a distribution on $\R$, has compact support in $[-1,1]_\sigma$, the same bounds are valid also for all higher $t_*$-derivatives of $a$, i.e.,
  \begin{equation}
  \label{EqA2Adma}
    \|\pa_{t_*}a\|_{H_{\cuop;\bop}^{(q;k)}(\ol\R)} \leq C_{q k}\|f\|_{H_{\tbop;\bop}^{(\sfs-1;k),(\alpha_\sface+2,0)}}\quad\forall\,q\in\N_0.
  \end{equation}
  We wish to control $a$ itself. We start with the decomposition
  \begin{equation}
  \label{EqA2AdmDecomp}
    u(t_*,\cdot)=u_{\rm hi}(t_*,\cdot)+u_{\rm reg}(t_*,\cdot)+a(t_*)u_0
  \end{equation}
  and the fact that $u$ vanishes for $t_*\leq 1$. For any $\phi\in\CIc(X^\circ;\cT^*)$, this gives
  \begin{equation}
  \label{EqA2AdmDecompPhi}
    0 = \la u_{\rm hi}(t_*,\cdot)+u_{\rm reg}(t_*,\cdot),\phi\ra + a(t_*)\la u_0,\phi\ra,\quad t_*\leq 1.
  \end{equation}
  Take $\phi=f^*$, so $\la u_0,f^*\ra$ is a nonzero constant. The first summand on the right-hand side lies in\footnote{This space is defined relative to~\eqref{EqA2AdmExpHcub} in the usual fashion, cf.\ Definition~\ref{DefMUSupp}.} $\bar H_{\cuop;\bop}^{(-N;k)}([-\infty,1];|\dd t_*|)$ for some $N$ (in fact, for $-N\leq\inf\sfs$), as follows by passing to the Fourier transform and using the resolvent bounds utilized already in~\eqref{EqA2AdmHi} and \eqref{EqA2AdmReg}. But by elliptic regularity for $\pa_{t_*}$ (as a uniform differential operator on $\R$), we can upgrade this using~\eqref{EqA2Adma} to the membership
  \begin{equation}
  \label{EqA2AdmaPast}
    a|_{(-\infty,T]} \in \bar H_{\cuop;\bop}^{(\infty;k)}([-\infty,T];|\dd t_*|),\quad
    \|a|_{(-\infty,T]}\|_{\bar H_{\cuop;\bop}^{(q;k)}} \leq C_{T q k}\|f\|_{H_{\tbop;\bop}^{(\sfs-1;k),(\alpha_\sface+2,0)}}\quad\forall\,T\in\R,
  \end{equation}
  first for $T=1$ and then for all $T$ by integration (and using~\eqref{EqA2Adma}). (This should be thought of as a weak forward support property for $a$. To reiterate, while we only get $L^2$-bounds for $\pa_{t_*}a$ from its definition~\eqref{EqA2AdmaDef}, we also get $L^2$-bounds for negative times by using the forward nature of the total solution~\eqref{EqA2AdmDecomp}.) Fix now a function $\chi_T=\chi_T(t_*)\in\CI(\R_{t_*})$ which equals $1$ on $[3,\infty)$ and $0$ on $(-\infty,2]$. Then~\eqref{EqA2AdmDecomp} yields
  \begin{equation}
  \label{EqA2AdmRewrite}
    u(t_*,\cdot) = \bigl( u_{\rm hi}(t_*,\cdot) + u_{\rm reg}(t_*,\cdot) + (1-\chi_T(t_*))a(t_*)u_0 \bigr) + \chi_T(t_*)a(t_*)u_0.
  \end{equation}
  Using~\eqref{EqA2AdmHi} and \eqref{EqA2AdmReg} as well as~\eqref{EqA2AdmaPast} (for $T=3$), the first summand lies in $H_{\tbop;\bop}^{(\sfs;k),(\alpha_\sface,0)}$ with norm bounded by $\|f\|_{H_{\tbop;\bop}^{(\sfs;k),(\alpha_\sface+2,0)}}$; and, like the left-hand side and the second summand, it is supported in $\{t_*\geq 1\}$. Note that
  \[
    \|\pa_{t_*}(\chi_T a)\|_{\dot H_{\cuop;\bop}^{(q;k)}([1,\infty])} \leq C_{q k}\|f\|_{H_{\tbop;\bop}^{(\sfs-1;k),(\alpha_\sface+2,0)}}.
  \]
  We then apply Lemma~\ref{LemmaA2AdmExp} (with $\alpha=1$, $m=1$, and $\beta_\sface=-\frac32-\eps$, $\eps>0$) to obtain
  \begin{equation}
  \label{EqA2Admau0}
    \chi_T a u_0 \in H_{\tbop;\bop}^{(\sfs;k),(-\frac32-\eps,0)}(\Omega)^{\bullet,-}\quad\forall\,\eps>0,
  \end{equation}
  with norm bounded by $C_\eps\|f\|_{H_{\tbop;\bop}^{(\sfs-1;k),(\alpha_\sface+2,0)}}$. (Note that the condition~\eqref{EqA2AdmExpThr} follows from~\eqref{EqA2AdmThr} since $\beta_\sface<\alpha_\sface$.)

  \pfstep{Conclusion.} We have now established the bound
  \[
    \|u\|_{H_{\tbop;\bop}^{(\sfs;k),(-\frac32-\eps,-1)}(\Omega)^{\bullet,-}} \leq C_{\eps,k}\|f\|_{H_{\tbop;\bop}^{(\sfs;k),(\alpha_\sface+2,0)}(\Omega)^{\bullet,-}}\quad\forall\eps>0,\ k\in\N_0,
  \]
  for the forward solution of $\Box_g u=f$ when $f\in\CIc$. This implies~\eqref{EqSSAlephAdmSol} by density.
\end{proof}

\subsection{Proof of Theorem~\usref{ThmA2}}
\label{SsA2Pf}

Recall that we shall work with unweighted b-densities, unless otherwise specified. We will prove Theorem~\ref{ThmA2} using Nash--Moser iteration (Theorem~\ref{ThmA1NM}) using the spaces
\begin{align*}
  B^k_\eta &:= \C \oplus \dot H_\bop^{k,\eps_\cK}([1+\eta,\infty]) \oplus \dot H_\bop^{k,\eps_+}(\scri^+\cap\{t_*\geq 1+\eta\};\cT^*) \oplus \Hb^{k,(\beta_\sscri,\beta_+,\beta_\cK)}(\Omega_\eta;\cT^*)^{\bullet,-}, \\
  \bfB^k_\eta &:= \Hb^{k,(\beta_\sscri+1,\beta_++2,\beta_\cK)}(\Omega_\eta;\cT^*)^{\bullet,-},
\end{align*}
where $\Omega_\eta:=\Omega\cap\{t_*\geq 1+\eta\}$, $\eta\in[0,1]$, and we recall $\beta_\sscri=1+\eps_\sscri>\beta_+=1+\eps_+>\beta_\cK=1+\eps_\cK$, the usage of $\rho_\sscri=x_\sscri^2$ for $\scri^+$-weights (following~\eqref{EqFSob}), and the usage of unweighted b-densities. \textit{From now on, we drop the bundle $\cT^*$ from the notation.} We study the map
\begin{equation}
\label{EqA2PfMap}
\begin{split}
  &\Phi(c,a,u_\sscri,\tilde u):=\Box_g u-|\delta_g u|^2 u - f, \quad u=\Xi(c,a,u_\sscri,\tilde u), \\
  &\qquad \Xi(c,a,u_\sscri,\tilde u)(t_*,x):=\chi_T(t_*)\bigl(c u_0(x)+a(t_*)u_0(x)\bigr) + \chi_\sscri r^{-1} u_\sscri\Bigl(t_*,\frac{x}{|x|}\Bigr) + \tilde u(t_*,x).
\end{split}
\end{equation}

\begin{lemma}[Forward mapping properties]
\label{LemmaA2PfMap}
  If $(c,a,u_\sscri,\tilde u)\in B_\eta^\infty$, then $\Phi(u)\in\bfB_\eta^\infty$ when $u$ is given by~\eqref{EqA2PfMap}. Moreover, $\Phi\colon B^\infty\to\bfB^\infty$ is $\cC^2$ and satisfies the estimates~\eqref{EqA1NMPhi}.
\end{lemma}
\begin{proof}
  The key task is to prove the required decay properties of $\Phi(u)$; the estimates~\eqref{EqA1NMPhi} then follow easily from Lemma~\ref{LemmaA1Nonlin}. It suffices to consider $\eta=0$.

  \pfstep{Linear term: wave operator.} Since $\Box_g u_0=0$, we have
  \begin{equation}
  \label{EqA2PfBoxu}
    \Box_g u = c [\Box_g,\chi_T] u_0 + [\Box_g,\chi_T a]u_0 + \Box_g(\chi_\sscri \rho u_\sscri) + \Box_g\tilde u.
  \end{equation}
  Using~\eqref{EqSSSpecFam} (with $S\in\rho\CI(X)$ according to~\eqref{EqSSAdmBox}), we compute, for any $b=b(t_*)$,
  \begin{equation}
  \label{EqA2PfCommt}
  \begin{split}
    [\Box_g,b(t_*)] &= i b'(t_*)\pa_\sigma\wh{\Box_g}(0) - \frac12 b''(t_*)\pa_\sigma^2\wh{\Box_g}(0) \\
      &= -b'(t_*)\bigl(2\rho(\rho\pa_\rho-1-S)-Q\bigr) - b''(t_*)g^{0 0},
  \end{split}
  \end{equation}
  and recall $Q\in\rho^3\Diffb^1(X;\cT^*)$, $g^{0 0}\in\rho^2\CI(X)$. (This is equivalent to~\eqref{EqA2Heur1bFwd}.) Observe that $2\rho S+Q\in\rho^2\Diffb^1$. Since $(\rho\pa_\rho-1)u_0\in\rho^3\CI$, the first term in~\eqref{EqA2PfBoxu} thus lies in $\CIc((1,\infty)_{t_*})\rho^3\CI(X)$ and thus in $\Hb^{\infty,(\tilde\beta+1,\infty,\infty)}(\Omega)^{\bullet,-}$ provided $\tilde\beta+1<3$, which is satisfied by $\tilde\beta=\beta_\sscri$. Since
  \[
    \pa_{t_*}\colon\dot H_\bop^{\infty,\eps_\cK}([1,\infty])\to\dot H_\bop^{\infty,1+\eps_\cK}([1,\infty]),
  \]
  the second term in~\eqref{EqA2PfBoxu} lies in
  \begin{equation}
  \label{EqA2PfBoxtstar}
  \begin{split}
    &\dot H_\bop^{\infty,1+\eps_\cK}([1,\infty])\rho^3\CI(X) + \dot H_\bop^{\infty,2+\eps_\cK}([1,\infty])\rho^3\CI(X) \\
    &\quad \subset \Hb^{\infty,((2-\eps)+1,(2+\eps_\cK)+2-\eps,1+\eps_\cK)}(\Omega)^{\bullet,-}
  \end{split}
  \end{equation}
  for all $\eps>0$. Here, we use that $\rho^3\CI(X)\subset\Hb^{\infty,3-\eps}(X)$, and $\Hb^\infty(\ol{\R_{t_*}})\Hb^\infty(X)\subset\Hb^\infty(M_0)$, while the weights are $t_*^{-1-\eps_\cK}\rho^{3-\eps}=\rho_\sscri^{3-\eps}\rho_+^{1+\eps_\cK+3-\eps}\rho_\cK^{1+\eps_\cK}$ on $\Omega$. (The fact that the order at $\iota^+$ is the sum of the $t_*^{-1}$- and $r^{-1}$-decay rates comes, geometrically, from the description of $\iota^+$ as the front face of the blow-up~\eqref{EqCPXM1p}.) The weights of this space are $\geq 2\beta_\sscri+2$, $\beta_++2$, and $\beta_\cK$, respectively, so it is contained in $\bfB^\infty$.

  The membership $\Box_g(\chi_\sscri\rho u_\sscri)\in\bfB^\infty$ is proved exactly as in~\eqref{EqDpBoxScri}. Using $\Box_g\in \rho_\sscri\rho_+^2\Diffb^2(M_0)$, we deduce that the ultimate term in~\eqref{EqA2PfBoxu} lies in the space $\Hb^{\infty,(\beta_\sscri+1,\beta_++2,\beta_\cK)}(\Omega)^{\bullet,-}$ indeed.

  \pfstep{Nonlinearity.} It remains to study the nonlinear term $|\delta_g u|^2 u$. Since $\delta_g u_0=0$, we have
  \[
    \delta_g u = c [\delta_g,\chi_T(t_*)] u_0 + [\delta_g,\chi_T(t_*)a(t_*)]u_0 + \delta_g(\chi_\sscri r^{-1}u_\sscri) + \delta_g\tilde u.
  \]
  For general $b=b(t_*)$, we have $[\delta_g,b]=-b'(t_*)\iota_{\nabla t_*}$, so the first term lies in $\CIc((1,\infty)_{t_*})\rho\CI$ and the second in $\dot H_\bop^{\infty,1+\eps_\cK}([1,\infty])\rho\CI(X)$; both thus lie in $\cC_\bop^{\infty,(2\cdot 1,2+\eps_\cK,1+\eps_\cK)}(\Omega)$. For the third term, we use that
  \begin{equation}
  \label{EqA2PfLinAdmDelta}
    \delta_g\in\rho_+\Diff_\etbop^1\subset\rho_+\Diff_\bop^1
  \end{equation}
  since it is a linear combination of $\pa_{t_*}$ and $\pa_x$ with $\CI(M_0;\End(\cT^*))$-coefficients; this maps $\chi_\sscri r^{-1}u_\sscri$ into $x_\sscri^2\rho_+^2\rho_\cK^\infty\dot H_\bop^{\infty,\eps_+}([1,\infty]\times\Sph^2)\subset\cC_\bop^{\infty,(2,2+\eps_+,\infty)}(\Omega)$. Finally, $\delta_g$ maps $\tilde u\in\cC_\bop^{\infty,(2\beta_\sscri,\beta_+,\beta_\cK)}$ into the space $\cC_\bop^{\infty,(2\beta_\sscri,\beta_++1,\beta_\cK)}$. Altogether, we thus have
  \begin{equation}
  \label{EqA2PfMapDeltau}
    \delta_g u \in \cC_\bop^{\infty,(1,2+\eps_+,1+\eps_\cK)}(\Omega).
  \end{equation}
  Therefore,
  \[
    |\delta_g u|^2 u \in \cC_\bop^{\infty,(2,4+2\eps_+,2+2\eps_\cK)}\cdot\Hb^{\infty,(1-\eps,1-\eps,-\eps)} \subset \Hb^{\infty,(3-\eps,5+3\eps_+-\eps,2+2\eps_\cK-\eps)}
  \]
  lies in $\bfB^\infty$, indeed.
\end{proof}

The remainder of the proof of Theorem~\ref{ThmA2} concerns estimates for forward solutions of the linearization of $u\mapsto\Box_g u-|\delta_g u|^2 u-f$, which is given by
\begin{equation}
\label{EqA2PfLin}
  P_u \dot u = \Box_g - 2 u\Re\bigl((\delta_g u)(\delta_g \dot u)\bigr) - |\delta_g u|^2 \dot u.
\end{equation}
Our goal is to prove:

\begin{prop}[Asymptotics and tame estimates for forward solutions]
\label{PropA2PfAsy}
  There exists $d\in\N$ such that the following holds. Let $u=\Xi(c,a,u_\sscri,\tilde u)\in\Xi(B^\infty)$ (with $\Xi$ defined in~\eqref{EqA2PfMap}). Let $\eta\in[0,1]$ and $f\in\bfB_\eta^\infty$.\footnote{This $f$ can, of course, be different than the one in the statement of Theorem~\ref{ThmA2}.} Then the unique forward solution of the linear wave-type equation $P_u \dot u=f$ can be written in the form
  \[
    \dot u = \Xi(\dot c,\dot a,\dot u_\sscri,\dot{\tilde u}),
  \]
  where $(\dot c,\dot a,\dot u_\sscri,\dot{\tilde u})\in B_\eta^\infty$ obeys the tame estimates
  \[
    |(\dot c,\dot a,\dot u_\sscri,\dot{\tilde u})|_k \leq C_k\bigl( \|f\|_{k+d} + |(c,a,u_\sscri,\tilde u)|_{k+d}\|f\|_{2 d}\bigr),\quad k\geq d.
  \]
\end{prop}

Theorem~\ref{ThmA2} follows from Theorem~\ref{ThmA2Adm}, Lemma~\ref{LemmaA2PfMap}, and Proposition~\ref{PropA2PfAsy} by applying Theorem~\ref{ThmA1NM}: this produces $(c,a,u_\sscri,\tilde u)$ such that $\Phi(\Xi(c,a,u_\sscri,\tilde u))=0$, so $u=\Xi(c,a,u_\sscri,\tilde u)$ is the desired solution of the nonlinear wave equation~\eqref{EqA2}.

The proof of Proposition~\ref{PropA2PfAsy} proceeds in several steps.
\begin{itemize}
\item We first verify the assumptions of Theorem~\ref{ThmF}, which gives tame bounds for $\dot u$ on b-Sobolev spaces on $\Omega$ with polynomial weights; see Lemma~\ref{LemmaA2PfLinAdm}.
\item We prove the existence of the radiation field of $\dot u$ using Proposition~\ref{PropDscriRad}; see~\S\ref{SssA2Rad}.
\item We improve the control on $\dot u$ as $t_*\to\infty$ by re-writing $P_u\dot u=f$ as
  \begin{equation}
  \label{EqA2PfP0}
    \Box_g\dot u = f - (P_u-\Box_g)\dot u.
  \end{equation}
  (Note that $P_0=\Box_g$.) The second term on the right-hand side has more than $t_*^{-1}$ extra decay compared to $u$ itself. Upon inverting $\Box_g$ by passing to the Fourier transform and describing its resolvent in an increasingly precise fashion (for right-hand sides with increasingly high $\sigma$-regularity), in particular near zero frequency, we can then uncover the asymptotics of $\dot u$ step-by-step, beginning with almost-boundedness (\S\ref{SssA2B}), the extraction of the stationary $u_0(x)$ leading order term (\S\ref{SssA2Lot}), and finally the $\dot a(t_*)u_0(x)$ sub-leading term and the $\tilde u$ remainder (\S\ref{SssA2lot}).
\end{itemize}

\begin{lemma}[Admissibility]
\label{LemmaA2PfLinAdm}
  For $d\in\N_0$ large enough and $k\in\N_0$, the following holds. Let $(c,a,u_\sscri,\tilde u)\in B^{k+d}$ and $u:=\Xi(c,a,u_\sscri,\tilde u)$. Then $P_u$ is an admissible wave-type operator of class $((0;k),(2\ell_\sscri,\ell_+,\ell_\cK))$ for the parameters $\ell_\sscri=\frac12$, $\ell_+=1+\eps_+$, and $\ell_\cK=1+\eps_\cK$.
\end{lemma}
\begin{proof}
  The regularity accounting arises from Sobolev embedding as in the proof of Theorem~\ref{ThmA1Null}; we leave this to the reader, and work with spaces of infinite b-regularity. We need to study $P_u-\Box_g$ as an edge-3b-differential operator, and thus use powers of $x_\sscri$ as weights at $\scri^+$. The membership~\eqref{EqA2PfMapDeltau} now reads
  \[
    \delta_g u \in \cC_\bop^{\infty,(2,2+\eps_+,1+\eps_\cK)}(\Omega).
  \]
  Using~\eqref{EqA2PfLinAdmDelta} and $u\in\cC_\bop^{\infty,(2,1,0)}$, this implies
  \begin{align*}
    u\,\ol{\delta_g u}\,\delta_g(\cdot),\ u\,\delta_g u\,\ol{\delta_g(\cdot)} &\in \cC_\bop^{\infty,(4,3+\eps_+,1+\eps_\cK)}\Diff_\etbop^1, \\
    |\delta_g u|^2\times(\cdot) &\in \cC_\bop^{\infty,(4,4+2\eps_+,2+2\eps_\cK)}\Diff_\etbop^0.
  \end{align*}
  This yields~\eqref{EqSDWAdmOp} (with $P=P_u$ and $p_0=\tilde p_1=0$) provided $2\ell_\sscri+2\leq 4$, $\ell_++2\leq 3+\eps_+$, and $\ell_\cK\leq 1+\eps_\cK$.
\end{proof}

For later use, we record that the proofs of Lemmas~\ref{LemmaA2PfMap} and~\ref{LemmaA2PfLinAdm} imply the first line of
\begin{subequations}
\begin{align}
\label{EqA2PfLinDiff}
  u\in\Xi(B^\infty) \implies &P_u - \Box_g \in \cC_\bop^{\infty,(2,3+\eps_+,1+\eps_\cK)}\Diff_\bop^1 = \rho_\sscri^2\rho_+^{3+\eps_+}\rho_\cK^{1+\eps_\cK}\CI_\bop\Diff_\bop^1, \\
\label{EqA2PfLinDiffHb}
    &P_u - \Box_g \in \Hb^{\infty,(2-\eps,3+\eps_+-\eps,1+\eps_\cK)}\Diffb^1\quad\forall\,\eps>0;
\end{align}
\end{subequations}
for the membership~\eqref{EqA2PfLinDiffHb}, we only need to keep track of $L^2$-integrability of the terms with sharp $\rho_\cK^{1+\eps_\cK}$-decay, which are $c[\delta_g,\chi_T]u_0+[\delta_g,\chi_T a]u_0\in\Hb^{\infty,(1-\eps,2+\eps_\cK-\eps,1+\eps_\cK)}$ and $\delta_g\tilde u\in\Hb^{\infty,(\beta_\sscri,\beta_++1,\beta_\cK)}$; this leads to $u\,\ol{\delta_g u}\,\delta_g(\cdot)$ and $u\,\delta_gu\,\ol{\delta_g(\cdot)}$ lying in $\Hb^{\infty,(2-\eps,3+\eps_+-\eps,1+\eps_\cK)}\Diffb^1$.

Given Lemma~\ref{LemmaA2PfLinAdm}, we may now apply Theorem~\ref{ThmF}, with $\alpha_\sscri$ and $\alpha_+$ there equal to $-\frac12-\frac{\eps}{2}$ and $-1-\eps$ for any $\eps\in(0,\frac12)$, and with $\delta=\frac12-\eps+\eps'$ (where $\eps'>0$ is arbitrary) according to Theorem~\ref{ThmA2Adm}. This gives (using Lemma~\ref{LemmaA1Nonlin} for the tame bounds of $P_u-\Box_g$ in terms of $u$) the tame estimate
\begin{align*}
  &\|\dot u\|_{\Hb^{k,(1-\frac{\eps}{2},\frac12-\eps',-\frac12)}(\Omega)^{\bullet,-}} \\
  &\qquad \leq C_{k,\eps,\eps'}\Bigl( \|f\|_{\Hb^{k+d,(2-\frac{\eps}{2},3-\eps,\frac12)}(\Omega)^{\bullet,-}} + |(c,a,u_\sscri,\tilde u)|_{k+d} \|f\|_{\Hb^{d,(2-\frac{\eps}{2},3-\eps,\frac12)}(\Omega)^{\bullet,-}}\Bigr).
\end{align*}
for every $\eps>0$ (and for some fixed $d\in\N$). The shifts in the orders are due to our present usage of unweighted b-densities (as opposed to the Kerr density in Theorem~\ref{ThmF}). Thus, $\dot u$ gains one power of $\rho_\sscri$ (i.e., $2$ powers of $x_\sscri$) over $f$ at $\scri^+$, loses in the best case (for $\eps$ close to $\frac12$ and $\eps'$ close to $0$) a bit more than $2$ powers of $\rho_+$-decay, and one power of $\rho_\cK$-decay. For the purposes of proving Proposition~\ref{PropA2PfAsy}, it suffices to record
\begin{equation}
\label{EqA2PfLinTame}
  \|\dot u\|_{\Hb^{k,(1-\frac{\eps}{2},\frac12-\eps,-\frac12)}(\Omega)^{\bullet,-}} \leq C_{k,\eps}\bigl( \|f\|_{k+d} + |(c,a,u_\sscri,\tilde u)|_{k+d} \|f\|_d\bigr),\quad \eps>0.
\end{equation}

When $f\in\bfB^\infty$, our remaining task is to show $\dot u\in\Xi(B^\infty)$, together with tame estimates. The same arguments as after~\eqref{EqDpTame0} show that, with~\eqref{EqA2PfLinTame} as the starting point, tame estimates for norms of $\dot u$ in spaces with stronger decay are immediate consequences of Lemma~\ref{LemmaA1Nonlin} and the results of~\S\ref{SsDRes}. In the sequel, we shall thus focus exclusively on improving the weights of
\begin{equation}
\label{EqA2PfdotuWeak}
  \dot u\in\Hb^{\infty,(1-\frac{\eps}{2},\frac12-\eps,-\frac12)}(\Omega)^{\bullet,-},\quad P_u\dot u=f\in\bfB^\infty=\Hb^{\infty,(\beta_\sscri+1,\beta_++2,\beta_\cK)}(\Omega)^{\bullet,-},
\end{equation}
and not discuss the tameness of the corresponding quantitative estimates any further.

\subsubsection{Radiation field}
\label{SssA2Rad}

Starting with~\eqref{EqA2PfdotuWeak}, we deduce from~\eqref{EqA2PfLinDiff} (and omitting $\cK^+$-orders) that
\[
  \Box_g\dot u = f - (P_u-\Box_g)\dot u \in \Hb^{\infty,(\beta_\sscri+1,\beta_++2)} + \Hb^{\infty,(3-\frac{\eps}{2},\frac72+\eps_+-\eps)} \subset \Hb^{\infty,(\beta_\sscri+1,\beta_++2)}
\]
for $0<\eps<\frac12$, using that $\beta_\sscri=1+\eps_\sscri$, $\beta_+=1+\eps_+<\frac32$. Proposition~\ref{PropDscriRad} (with $\beta_\sscri=1-\frac{\eps}{2}$, $\beta_+=\frac12-\eps$, and $\beta'_\sscri=\beta_\sscri$) improves~\eqref{EqA2PfdotuWeak} to
\begin{equation}
\label{EqA2PfdotuWeakScri}
  \dot u \in \Hb^{\infty,(((1,0),\beta_\sscri),\frac12-\eps,-\frac12)}(\Omega)^{\bullet,-}\quad\forall\,\eps>0.
\end{equation}
Notice that for such $\dot u$, the membership~\eqref{EqA2PfLinDiff} gives
\begin{equation}
\label{EqA2PfdotuErr}
  (P_u-\Box_g)\dot u \in \bigcap_{\eps>0}\Hb^{\infty,(((3,0),\beta_\sscri+2), \frac72+\eps_+-\eps, \frac12+\eps_\cK)}(\Omega)^{\bullet,-} \subset \Hb^{\infty,(\beta_\sscri+1,\beta_++2,\frac12+\eps_\cK)}(\Omega)^{\bullet,-},
\end{equation}
where we use $\beta_\sscri+1<3$ and $\beta_++2=3+\eps_+<\frac72+\eps_+-\eps$ for $\eps<\frac12$; thus, only the weight of this expression at $\cK^+$ is weaker than that of $f$ (but stronger than that of $\dot u$ in~\eqref{EqA2PfdotuWeakScri} by $1+\eps_\cK>1$ orders).

\subsubsection{Almost-boundedness}
\label{SssA2B}

Slightly generalizing~\eqref{EqA2PfdotuWeakScri}--\eqref{EqA2PfdotuErr}, note that if $\beta'_+\geq\tfrac12-\eps$ for some $\eps>0$ and $\beta'_\cK\geq-\tfrac12$, then
\begin{align}
\label{EqA2LotRHS}
  &\dot u\in\Hb^{\infty,(((1,0),\beta_\sscri),\beta'_+,\beta'_\cK)}(\Omega)^{\bullet,-} \\
\label{EqA2LotRHS2}
  &\qquad \implies \Box_g\dot u = f-(P_u-\Box_g)\dot u =: h \in \Hb^{\infty,(\beta_\sscri+1,\ \beta_++2,\ \min(\beta'_\cK+\beta_\cK,\beta_\cK))}(\Omega)^{\bullet,-}.
\end{align}
The hypothesis on $\dot u$ was proved for $\beta'_+=\frac12-\eps$ and $\beta'_\cK=-\frac12$; we shall use equation~\eqref{EqA2LotRHS2} and the Fourier transform to improve the decay rates of $\dot u$ at $\iota^+\cup\cK^+$.

Assume~\eqref{EqA2LotRHS} to hold for some $\beta'_\cK\in[-\frac12,1-\beta_\cK)$. We split $\dot u$ into a high- and low-frequency piece similarly to~\eqref{EqDpdotuSplit}, but avoid the singularity of $\wh{\Box_g}(\sigma)^{-1}$ at $\sigma=0$ by integrating over the contour
\[
  \gamma_-\cup\gamma_0\cup\gamma_+,\quad
  \gamma_-=(-\infty,-\tfrac12],\ \gamma_0=\tfrac12 e^{i[\pi,0]},\ \gamma_+=[\tfrac12,\infty)
\]
(which is similar to~\eqref{EqA2AdmContour} but more convenient for present purposes). Let $\chi_0\in\CIc([0,1))$ be equal to $1$ on $[0,\frac12]$. We then write the solution of~\eqref{EqA2LotRHS2} as
\begin{equation}
\label{EqA2udotDecomp}
\begin{split}
  &\dot u = \dot u_{\rm lo} + \dot u_{\rm hi}, \\
  &\qquad \dot u_{\rm lo}(t_*,\cdot) := \frac{1}{2\pi}\int_{\gamma_-\cup\gamma_0\cup\gamma_+} e^{-i\sigma t_*}\chi_0(|\sigma|) \wh{\Box_g}(\sigma)^{-1}\hat h(\sigma)\,\dd\sigma,  \\
  &\qquad \dot u_{\rm hi} := \cF^{-1}\bigl( (1-\chi_0(|\sigma|))\wh{\Box_g}(\sigma)^{-1}\hat h(\sigma)\bigr).
\end{split}
\end{equation}
The high-energy piece satisfies $\dot u_{\rm hi}\in\Hb^{\infty,(\gamma',\beta_+,\beta_\cK)}(\sfM)$ in the notation of Corollary~\ref{CorDResHiLoc}.

\pfstep{Low-energy control.} For the low-energy piece, note that by a density argument it suffices to prove bounds for $\dot u_{\rm lo}$ in $\Hb$-spaces in terms of $\Hb^{\infty,(\beta_\sscri+1,\ \beta_++2,\ \beta'_\cK+\beta_\cK)}$-seminorms of $h$ for $h\in\CIc(\{t_*>1\})$; thus the holomorphicity considerations of the proof of Theorem~\ref{ThmA2Adm} are applicable. We first use Proposition~\ref{PropDFTFourier}\eqref{ItDFTFourierLo}, with $\beta_\sscri$, $\beta_\cK$, and $\beta_\iota$ there equal to $\beta_\sscri+1$, $\beta'_\cK+\beta_\cK$, and
\begin{equation}
\label{EqA1hathLoBetaI}
  \beta_\iota := \min(\beta_\sscri+1+\beta_\cK+\beta'_\cK,\ \beta_++2)
\end{equation}
in present notation. (Thus, $\beta_\iota=\beta_\sscri+\beta_\cK+\beta'_\cK$ unless $\beta'_\cK\geq 1-\beta_\cK-(\eps_+-\eps_\sscri)$.) This gives
\begin{equation}
\label{EqA2hathLo}
  \hat h|_{\sigma\in\pm[0,1]} \in \Hb^{\infty,(\gamma,\ \beta_\iota-1,\ \beta_\cK+\beta'_\cK-1)}(X_\scbtop^\pm)
\end{equation}
for some $\gamma\in\R$,\footnote{By this membership---which is trivial when $h\in\CIc$---, we mean that the $\Hb^k$-seminorm of~\eqref{EqA2hathLo} are bounded in terms of the $\Hb^{k+\ell}$-seminorm of $h$ for some $k$-independent $\ell$.} which we allow to change at each step in the remainder of the proof. We invert $\wh{\Box_g}(\sigma)$ using the inverse of the augmentation $\wt\Box(\sigma)$ defined in~\eqref{EqA2AdmAug}. Note that the relative order of~\eqref{EqA2hathLo} at $\tface$ and $\zface$ is $\beta_\iota-(\beta_\cK+\beta'_\cK)\geq\beta_\sscri+1=2+\eps_\sscri$ and is thus falls into the indicial gap $(0,1)$ of $\Box_g$ shifted by $2$. Proposition~\ref{PropDResLo} (which applies to $\wt\Box(\sigma)$, \emph{mutatis mutandis}, due to the invertibility of $\wt\Box(0)$, cf.\ Remark~\ref{RmkSpLoAug}) therefore implies that for
\begin{equation}
\label{EqA2uregDec}
  \wt\Box(\sigma)^{-1}(\hat h(\sigma),0) =: \bigl( \wh{\dot u}_{\rm reg}(\sigma), \hat a_1(\sigma) \bigr),
\end{equation}
we have $\wh{\dot u}(\sigma)=\sigma^{-1}\hat a_1(\sigma)u_0+\wh{\dot u}_{\rm reg}(\sigma)$ and
\begin{equation}
\label{EqA2ureg}
  \wh{\dot u}_{\rm reg} \in \Hb^{\infty,(\gamma,\ \beta_\iota-3,\ \beta_\cK+\beta'_\cK-1)}(X_\scbtop^\pm),\quad
  \hat a_1 \in \Hb^{\infty,\beta_\cK+\beta'_\cK-1}(\pm[0,1]_\sigma).
\end{equation}
Moreover, as discussed after~\eqref{EqA2AdmLo}, $\wh{\dot u}_{\rm reg}(\sigma)$ (in the distributional sense) and $\hat a_1(\sigma)$ are holomorphic in $\sigma$ for $|\sigma|\leq 1$ and $\Im\sigma>0$, and continuous down to $[-1,1]\setminus\{0\}$. Write
\begin{align*}
  \dot u_{\rm lo} = \dot u_{\rm reg} + a(t_*)u_0&, \\
  \qquad \dot u_{\rm reg}(t_*,\cdot)&:=\frac{1}{2\pi}\int_{\gamma_-\cup\gamma_0\cup\gamma_+} e^{-i\sigma t_*}\chi_0(|\sigma|)\wh{\dot u}_{\rm reg}(\sigma)\,\dd\sigma, \\
  \qquad a(t_*)&:=\frac{1}{2\pi}\int_{\gamma_-\cup\gamma_0\cup\gamma_+} e^{-i\sigma t_*}\chi_0(|\sigma|)\sigma^{-1}\hat a_1(\sigma)\,\dd\sigma.
\end{align*}

\pfstep{Regular part.} We can shift the contour $\gamma_0$ in the integral for $\dot u_{\rm reg}$ to the real axis; by~\eqref{EqA2ureg} and Proposition~\ref{PropDFTInv}\eqref{ItDFTInvLo}, we get
\begin{equation}
\label{EqA2uregBound}
  \dot u_{\rm reg} = \cF^{-1}\bigl(\chi_0(|\sigma|)\wh{\dot u}_{\rm reg}(\sigma)\bigr) \in \Hb^{\infty,(\gamma,\ \beta_\iota-2,\ \beta_\cK+\beta'_\cK)}(\sfM).
\end{equation}

\pfstep{Singular part.} We analyze $a$ similarly to the proof of Theorem~\ref{ThmA2Adm}: we have
\begin{equation}
\label{EqA2paa}
  \pa_{t_*}a = \cF^{-1}\hat a_1 \in \Hb^{\infty,\beta_\cK+\beta'_\cK}(\ol\R,|\tfrac{\dd t_*}{t_*}|),
\end{equation}
as follows from shifting the integration contour to the real axis and using Corollary~\ref{CorDFTLossy}. Similarly to the arguments around~\eqref{EqA2AdmDecompPhi}, we can then use the decomposition $\dot u=\dot u_{\rm hi}+\dot u_{\rm reg}+a u_0$ for $t_*<1$ to bound $a|_{(-1,1)}$ in $H^k((-1,1))$ by the $H^k((-1,1)_{t_*})$-norm of $\la\dot u_{\rm hi}(t_*,\cdot)+\dot u_{\rm reg}(t_*,\cdot),\phi\ra$ where $\phi\in\CIc(X^\circ)$ satisfies $\la u_0,\phi\ra\neq 0$. By integrating~\eqref{EqA2paa}, we obtain bounds also for $a|_{(-1,3)}\in H^\infty((-1,3))$. We then write
\begin{equation}
\label{EqA2dotuSplit}
  \dot u(t_*,\cdot) = \bigl(\dot u_{\rm hi}(t_*,\cdot) + \dot u_{\rm reg}(t_*,\cdot) + (1-\chi_T(t_*))a(t_*)u_0\bigr) + \chi_T(t_*)a(t_*)u_0,
\end{equation}
where
\begin{equation}
\label{EqA2chiT}
  \chi_T\in\CI(\R),\quad \chi_T|_{[3,\infty)}=1, \ \chi_T|_{(-\infty,2]}=0.
\end{equation}
The first summand is supported in $t_*\geq 1$ since the left-hand side and the second summand are, and it satisfies
\begin{equation}
\label{EqA2dotuSum1}
  \dot u_{\rm hi} + \dot u_{\rm reg} + (1-\chi_T)a u_0 \in \Hb^{\infty,(\gamma,\ \beta_\iota-2,\ \beta_\cK'+\beta_\cK)}(\Omega)^{\bullet,-}.
\end{equation}
We control the second summand $(\chi_T a)u_0$ using the following general result:

\begin{lemma}[Integration on $\ol\R$]
\label{LemmaA2Int}
  Let $\beta\in\R\setminus\Z$ and $k,m\in\N_0$. Let $a$ be a distribution on $\R$ with support in $[1,\infty)$ such that $\pa_{t_*}^m a\in\dot H_\bop^{k,\beta}([1,\infty],|\frac{\dd t_*}{t_*}|)$. Then $a\in\dot H_\bop^{k+m,((-m,0),\beta-m)}([1,\infty])$. Here, we write
  \begin{equation}
  \label{EqA2IntSpace}
    \dot H_\bop^{s,((\zeta,0),\gamma)}([1,\infty]) := \Biggl\{ \chi_T(t_*)\sum_{\ell=0}^{\lfloor\gamma-\zeta\rfloor} t_*^{-\zeta-\ell}a_\ell + \tilde a \colon a_\ell\in\C,\ \tilde a(t_*)\in\dot H_\bop^{s,\gamma}([1,\infty]) \Biggr\},
  \end{equation}
  where $\chi_T$ is as in~\eqref{EqA2chiT}.
\end{lemma}
\begin{proof}
  The assumption on $a$ is equivalent to
  \[
    \dot H_\bop^{k,\beta-m} \ni t_*^m\pa_{t_*}^m a = \prod_{j=0}^{m-1} (t_*\pa_{t_*}-j)a = \prod_{j=0}^{m-1} t_*^j t_*\pa_{t_*}(t_*^{-j}a).
  \]
  By induction on $m$, it thus suffices to integrate a single factor of $t_*\pa_{t_*}$ and observe that
  \[
    \beta\notin\Z,\ \ t_*\pa_{t_*}a \in \dot H_\bop^{k,\beta} \implies a \in \dot H_\bop^{k+1,((0,0),\beta)},
  \]
  which is easily proved by direct integration (or by passing to the Mellin transform in $\tau=t_*^{-1}$).
\end{proof}

In the case at hand, since $\pa_{t_*}(\chi_T a)\in\dot H_\bop^{\infty,\beta_\cK+\beta'_\cK}$, we obtain
\[
  \chi_T a \in \dot H_\bop^{\infty,\beta_\cK+\beta'_\cK-1}([1,\infty])
\]
since $\beta'_\cK+\beta_\cK-1<0$, and therefore
\begin{equation}
\label{EqA2dotuSum2}
  (\chi_T a)u_0 \in \dot H_\bop^{\infty,\beta_\cK+\beta'_\cK-1}([1,\infty])\rho\CI(X) \subset \Hb^{\infty,(1-\eps,\beta_\cK+\beta'_\cK-\eps,\beta_\cK+\beta'_\cK-1)}(\Omega)^{\bullet,-}
\end{equation}
for all $\eps>0$. Together with~\eqref{EqA2dotuSum1}, this shows that
\[
  \dot u \in \Hb^{\infty,(\gamma,\ \min(\beta_\iota-2,\beta_\cK+\beta'_\cK-\eps),\ \beta_\cK+\beta'_\cK-1)}(\Omega)^{\bullet,-}\quad\forall\,\eps>0.
\]
An application of Proposition~\ref{PropDscriRad} recovers the radiation field, so
\[
  \dot u \in \Hb^{\infty,(((1,0),\beta_\sscri),\ \min(\beta_\iota-2,\beta_\cK+\beta'_\cK-\eps),\ \beta_\cK+\beta'_\cK-1)}(\Omega)^{\bullet,-}\quad\forall\,\eps>0.
\]
This improves the $\cK^+$-decay rate of~\eqref{EqA2LotRHS} by $1-\beta_\cK=\eps_\cK>0$; the new $\iota^+$-decay rate is $\beta'_+=\min(\beta_+,\beta_\cK+\beta'_\cK-\eps)>-\frac12$ (for all small $\eps>0$). Iterating this $\lceil\frac{1/2}{\eps_\cK}\rceil$ many times yields the $\cK^+$-decay rate $\beta_\cK+\beta'_\cK-1=-\eps$ for any desired $\eps>0$, and thus
\begin{equation}
\label{EqA2AlmostB}
  \dot u \in \Hb^{\infty,(((1,0),\beta_\sscri),\ 1-2\eps,\ -\eps)}(\Omega)^{\bullet,-}\quad\forall\,\eps>0.
\end{equation}
This is an improvement over~\eqref{EqA2PfdotuWeakScri} by a factor of $t_*^{-\frac12+\eps}$.

\begin{rmk}[Stronger asymptotic description]
\label{RmkA2dotuStrong}
  A minor modification of the above arguments shows that $\dot u=(\chi_T a)u_0+\tilde u$ where $\tilde u\in\Hb^{\infty,(((1,0),\beta_\sscri,1-2\eps,1-\eps)}$ for all $\eps>0$; and this holds more generally for the forward solution of $\Box_g\dot u=h$ for \emph{any} $h\in\Hb^{\infty,(\beta_\sscri+1,3-2\eps,1-\eps)}$. (The existence of the radiation field of $\tilde u$ can be proved using the arguments leading up to~\eqref{EqA2udotuAlmost2} below.) That is, the regular part $\tilde u$ of $\dot u$ inherits all decay rates from $h$, and only the expansion term $(\chi_T a)u_0$ has weaker $t_*$-bounds. Thus, we are, in fact, already quite close to proving Proposition~\ref{PropA2PfAsy}.
\end{rmk}

\subsubsection{Extraction of the leading-order term}
\label{SssA2Lot}

We use~\eqref{EqA2LotRHS}--\eqref{EqA2LotRHS2} with $\beta'_+=1-2\eps$ and $\beta'_\cK=-\eps$ (for all $\eps>0$) from~\eqref{EqA2AlmostB}; we must thus control the forward solution of
\begin{equation}
\label{EqA2LotEq}
  \Box_g\dot u =: h \in \Hb^{\infty,(\beta_\sscri+1,\ \beta_++2,\ \beta_\cK-\eps)}(\Omega)^{\bullet,-}.
\end{equation}
We again express $\dot u$ via~\eqref{EqA2udotDecomp}, with $\dot u_{\rm hi}\in\Hb^{\infty,(\gamma',\beta_+,\beta_\cK)}$ as there. For the low-energy piece, we note that, as in~\eqref{EqA2hathLo} (where now $\beta_\iota=\beta_++2$ in~\eqref{EqA1hathLoBetaI}), we now have
\[
  \hat h|_{\sigma\in\pm[0,1]} \in \Hb^{\infty,(\gamma,\ \beta_++1,\ ((0,0),\beta_\cK-1-\eps))}(X_\scbtop^\pm)\quad\forall\,\eps>0,
\]
by Proposition~\ref{PropDFTFourier}\eqref{ItDFTFourierLo}. As in~\eqref{EqA2uregDec}, we again write
\[
  \wh{\dot u}(\sigma)=\sigma^{-1}\hat a_1(\sigma)u_0+\wh{\dot u}_{\rm reg}(\sigma),\quad \bigl(\wh{\dot u}_{\rm reg}(\sigma),\hat a_1(\sigma)\bigr)=\wt\Box(\sigma)^{-1}(\hat h(\sigma),0)
\]
for $\sigma\in\C$, $\Im\sigma\geq 0$, $|\sigma|\leq 1$. Unlike in~\S\ref{SssA2B}, we claim that the functions $\hat a_1(\sigma)$ and $\wh{\dot u}_{\rm reg}(\sigma)$ now have a well-defined limit at $\sigma=0$ plus a conormal remainder; this is analogous to~\eqref{EqDpSolv0a}--\eqref{EqDpSolv0e}. To wit, since $\hat h(0)\in\Hb^{\infty,\beta_++1}(X)=\Hb^{\infty,2+\eps_+}(X)$ (using an unweighted b-density), with $2+\eps_+$ lying in the indicial gap $(0,1)$ of $\wh{\Box_g}(0)$ shifted by $2$, $(\wh{\dot u}_{\rm reg}(0),\hat a_1(0))=\wt\Box(0)^{-1}(\hat h(0),0)$ is well-defined, with
\[
  \wh{\dot u}_{\rm reg}(0) \in \Hb^{\infty,\beta_+-1}(X),\quad \hat a_1(0)\in\C.
\]
Consider the remainder
\[
  \bigl(\wh{h'}(\sigma), \hat b_1(\sigma) \bigr) := \bigl(\hat h(\sigma), 0\bigr) - \wt\Box(\sigma) \bigl( \chi_0(\hat r)\wh{\dot u}_{\rm reg}(0), \hat a_1(0) \bigr),
\]
where $\hat r=\frac{|\sigma|}{\rho}$, and $\chi_0\in\CIc((-1,1))$ equals $1$ on $[-\frac12,\frac12]$. Similarly to~\eqref{EqDpSolv0c}, the first component satisfies
\begin{align*}
  \wh{h'}(\sigma)&=\hat h(\sigma)-\wh{\Box_g}(\sigma)\bigl(\chi_0(\hat r)\wh{\dot u}_{\rm reg}(0)\bigr) - \hat a_1(0)\wh{\Box_g}(\sigma)(\sigma^{-1}u_0) \\
    &\in \Hb^{\infty,(\gamma,\ \beta_++1,\ ((0,0),\beta_\cK-1-\eps))}(X_\scbtop^\pm)
\end{align*}
(where we use Lemma~\ref{LemmaA212} to deduce this membership for $\wh{\Box_g}(\sigma)(\sigma^{-1}u_0)$); but its restriction to $\sigma=0$ being $0$, so in fact
\[
  \bigl(\wh{h'}(\sigma)\bigr)_{\sigma\in\pm[0,1]} \in \Hb^{\infty,(\gamma,\ \beta_++1,\ \beta_\cK-1-\eps)}(X_\scbtop^\pm).
\]
Furthermore, $\hat b_1(\sigma)=-\la\chi_0(\hat r)\wh{\dot u}_{\rm reg}(0),f^*\ra$ (in the notation of~\eqref{EqA2AdmAug}, so with $f^*\in\CIc(X^\circ)$). Since $\la\wh{\dot u}_{\rm reg}(0),f^*\ra=0$, we conclude that $\hat b_1(\sigma)$ vanishes near $\sigma=0$ (and is smooth). Since $(\beta_+-1)-(\beta_\cK-1-\eps)=\eps_+-\eps_\cK-\eps\in(0,1)$ for small $\eps>0$, Proposition~\ref{PropDResLo} applies and gives
\[
  \bigl(\wt\Box(\sigma)^{-1}\bigl( \wh{h'}(\sigma), \hat b_1(\sigma) \bigr)\bigr)_{\sigma\in\pm[0,1]} \in \Hb^{\infty,(\gamma,\ \beta_+-1,\ \beta_\cK-1-\eps)}(X_\scbtop^\pm),
\]
and hence
\begin{align*}
  \bigl(\wh{\dot u}_{\rm reg}(\sigma),\hat a_1(\sigma)\bigr) &= \bigl( \chi_0(\hat r)\wh{\dot u}_{\rm reg}(0), \hat a_1(0) \bigr) + \wt\Box(\sigma)^{-1}\bigl( \wh{h'}(\sigma), \hat b_1(\sigma) \bigr) \\
    &\in \Hb^{\infty,(\gamma,\ \beta_+-1,\ ((0,0),\beta_\cK-1-\eps)}(X_\scbtop^\pm) \oplus \Hb^{\infty,((0,0),\beta_\cK-1-\eps)}(\pm[0,1]).
\end{align*}

We now pass to the inverse Fourier transform to deduce that
\[
  \dot u_{\rm reg} = \cF^{-1}\bigl(\chi_0(|\sigma|)\wh{\dot u}_{\rm reg}(\sigma)\bigr) \in \Hb^{\infty,(\gamma,\ \beta_+,\ \beta_\cK-\eps)}(\sfM)
\]
by Proposition~\ref{PropDFTInv}\eqref{ItDFTFourierLo}. On the other hand, for
\[
  a(t_*) := \frac{1}{2\pi}\int_{\gamma_-\cup\gamma_0\cup\gamma_+} e^{-i\sigma t_*}\chi_0(|\sigma|)\sigma^{-1}\hat a_1(\sigma)\,\dd\sigma,
\]
we can express $\dot u$ as in~\eqref{EqA2dotuSplit}, with the first term lying in $\Hb^{\infty,(\gamma,\beta_+,\beta_\cK-\eps)}(\Omega)^{\bullet,-}$, while for the second term we note that $\pa_{t_*}(\chi_T a)\in\dot\Hb^{\infty,\beta_\cK-\eps}([1,\infty],|\frac{\dd t_*}{t_*}|)$ and thus
\[
  \chi_T a \in \dot\Hb^{\infty,((0,0),\beta_\cK-1-\eps)}([1,\infty]) = \dot\Hb^{\infty,((0,0),\eps_\cK-\eps)}([1,\infty]),
\]
which for $t_*\geq 1$ is the sum of a constant $c\in\C$ and a decaying remainder.

Thus far, we have shown that the solution of~\eqref{EqA2LotEq} can be written in the form
\begin{subequations}
\begin{equation}
\label{EqA2udotuAlmost}
\begin{split}
  &\dot u(t_*,x) = \chi_T( c u_0(x) + a(t_*)u_0(x) \bigr) + \tilde u(t_*,x), \\
  &\qquad c\in\C,\ a\in\dot H_\bop^{\infty,\beta_\cK-\eps-1}([1,\infty]), \quad \tilde u\in\Hb^{\infty,(\gamma,\beta_+,\beta_\cK-\eps)}(\Omega)^{\bullet,-},
\end{split}
\end{equation}
for all $\eps>0$. We recover the radiation field by noting that
\[
  \Box_g\tilde u = h - \Box_g\bigl(\chi_T(c u_0(x)+a(t_*)u_0(x))\bigr) \in \Hb^{\infty,(\beta_\sscri+1,\beta_++2,\beta_\cK-\eps)}(\Omega)^{\bullet,-},
\]
where we use that $\Box_g(\chi_T a(t_*)u_0(x))\in\Hb^{\infty,(3-\eps,4+\eps_\cK-\eps,1+\eps_\cK-\eps)}(\Omega)^{\bullet,-}$, which follows from the arguments leading up to~\eqref{EqA2PfBoxtstar}. We can thus apply Proposition~\ref{PropDscriRad} to improve~\eqref{EqA2udotuAlmost} to
\begin{equation}
\label{EqA2udotuAlmost2}
  \tilde u\in\Hb^{\infty,(((1,0),\beta_\sscri),\beta_+,\beta_\cK-\eps)}(\Omega)^{\bullet,-}.
\end{equation}
\end{subequations}

\subsubsection{Precise sub-leading term and remainder; conclusion}
\label{SssA2lot}

The memberships~\eqref{EqA2udotuAlmost}--\eqref{EqA2udotuAlmost2} are only off by an $\eps$ in the $\cK^+$-decay order of $a$ and $\tilde u$ from the claimed membership of $\dot u\in\Xi(B^\infty_0)$ in Proposition~\ref{PropA2PfAsy}. In a final step, we revisit~\eqref{EqA2LotRHS2} and~\eqref{EqA2PfdotuErr}: we claim that for $\dot u$ given by~\eqref{EqA2udotuAlmost}--\eqref{EqA2udotuAlmost2}, we have
\begin{equation}
\label{EqA2lotErr}
  (P_u-\Box_g)\dot u \in \Hb^{\infty,(\beta_\sscri+1,\beta_++2,\beta_\cK)}(\Omega)^{\bullet,-}.
\end{equation}
Indeed, by~\eqref{EqA2PfLinDiff}, $P_u-\Box_g$ maps the term $\tilde u$ of $\dot u$ to $\Hb^{\infty,(3-\eps,2+2\beta_+,2\beta_\cK-\eps)}$, and $\chi_T(t_*)a(t_*)u_0\in\Hb^{\infty,(1-\eps,1+\eps_\cK-\eps,\eps_\cK-\eps)}$ to $\Hb^{\infty,(3-\eps,4+\eps_++\eps_\cK-\eps,1+2\eps_\cK-\eps)}$. For the term $\chi_T c u_0(x)$ of $\dot u$, we must now use~\eqref{EqA2PfLinDiffHb} in order to get the sharp $\cK^+$-decay order (\emph{relative to $L^2$}); thus $P_u-\Box_g$ maps $\chi_T c u_0(x)\in\cC_\bop^{\infty,(1,1,0)}$ to $\Hb^{\infty,(3-\eps,4+\eps_+-\eps,1+\eps_\cK)}$. This gives~\eqref{EqA2lotErr} (in fact, with stronger $\scri^+$- and $\iota^+$-decay), and thus
\[
  \Box_g\dot u = f-(P_u-\Box_g)\dot u \in \Hb^{\infty,(\beta_\sscri+1,\ \beta_++2,\ \beta_\cK)}.
\]
This is the same as~\eqref{EqA2LotEq} but without the $\eps$-loss at $\cK^+$. Repeating the arguments in~\S\ref{SssA2Lot} thus yields the description~\eqref{EqA2udotuAlmost}--\eqref{EqA2udotuAlmost2} for $\dot u$ with $\beta_\cK-\eps$ replaced by $\beta_\cK$. \textit{This completes the proof of Proposition~\usref{PropA2PfAsy} and thus of Theorem~\usref{ThmA2}.}

\clearpage
\bibliographystyle{alphaurl}
\newcommand{\etalchar}[1]{$^{#1}$}


\end{document}